\documentclass{article}
\usepackage{amsmath}
\usepackage{amssymb}
\usepackage{graphicx}
\usepackage{amsthm}
\usepackage{bbm}
\usepackage{cite}
\usepackage{xcolor}
\usepackage{stmaryrd}
\SetSymbolFont{stmry}{bold}{U}{stmry}{m}{n}
\usepackage{euscript}
\usepackage[arrow,curve,matrix,arc,2cell,frame]{xy}
\UseAllTwocells
\usepackage[utf8]{inputenc}
\usepackage{tikz}
\usepackage{tikz-cd}
\usepackage[new]{old-arrows}
\usetikzlibrary{calc}
\usetikzlibrary{positioning}
\usetikzlibrary{matrix,arrows}
\usetikzlibrary{decorations.markings}
\graphicspath{{pictures/},{Inkscape/}}
\usepackage{color}
\usepackage[breaklinks,unicode]{hyperref}
\usepackage{slashed}
\usepackage{mathrsfs}
\DeclareMathAlphabet{\mathscrbf}{OMS}{mdugm}{b}{n}
\def\scr{\mathscr}
\DeclareFontFamily{U}{mathb}{\hyphenchar\font45}
\DeclareFontShape{U}{mathb}{m}{n}{
      <5> <6> <7> <8> <9> <10> gen * mathb
      <10.95> mathb10 <12> <14.4> <17.28> <20.74> <24.88> mathb12
}{}
\DeclareSymbolFont{mathb}{U}{mathb}{m}{n}
\DeclareMathSymbol{\llcurly}{3}{mathb}{"CE}
\DeclareMathSymbol{\ggcurly}{3}{mathb}{"CF}
\begin{document}
%%%%%%%%%%%%%%%%%%%%%%%%%%%%%%%%%%%%%%%%%%%%%%%%%%%%%%%%%%%%%%%%%%%%%%%%
%%%%%%%%%%%%%%%%%%%%%%%%%%     Macros      %%%%%%%%%%%%%%%%%%%%%%%%%%%%%
%%%%%%%%%%%%%%%%%%%%%%%%%%%%%%%%%%%%%%%%%%%%%%%%%%%%%%%%%%%%%%%%%%%%%%%%
\def\e#1\e{\begin{equation}#1\end{equation}}
\def\ea#1\ea{\begin{align}#1\end{align}}
\def\eq#1{{\rm(\ref{#1})}}
\theoremstyle{plain}% default
\newtheorem{thm}{Theorem}[section]
\newtheorem{lem}[thm]{Lemma}
\newtheorem{prop}[thm]{Proposition}
\newtheorem{cor}[thm]{Corollary}
\newtheorem{quest}[thm]{Question}
\newtheorem{princ}[thm]{Principle}
\theoremstyle{definition}
\newtheorem{dfn}[thm]{Definition}
\newtheorem{altdfn}[thm]{Alternative Definition}
\newtheorem{ex}[thm]{Example}
\newtheorem{rem}[thm]{Remark}
\newtheorem{property}[thm]{Property}
\newtheorem{cond}[thm]{Condition}
\newtheorem{ass}[thm]{Assumption}
\newtheorem{conv}[thm]{Convention}
\numberwithin{figure}{section}
\numberwithin{table}{section}
\numberwithin{equation}{section}
\def\dim{\mathop{\rm dim}\nolimits}
\def\codim{\mathop{\rm codim}\nolimits}
\def\vdim{\mathop{\rm vdim}\nolimits}
\def\Cat{\mathop{\bf Cat}\nolimits}
\def\SCat{\mathop{\bf SCat}\nolimits}
\def\sCat{\mathop{{\mathscrbf C}\kern-.13em\bs{at}}\nolimits}
\def\sSpc{\mathop{{\mathscrbf S}\kern-.13em\bs{pc}}\nolimits}
\def\Set{{\mathop{\bf Set}\nolimits}}
\def\bKan{{\mathop{\bf Kan}}}
\def\Com{{\mathop{\bf Com}}}
\def\Pith{{\mathop{\rm Pith}}}
\def\Core{{\mathop{\rm Core}}}
\def\Gpds{{\mathop{\bf Gpds}\nolimits}}
\def\CGHaus{\mathop{\bf CGHaus}\nolimits}
\def\TR{\mathop{\bf TR}\nolimits}
\def\ES{\mathop{\bf ES}\nolimits}
\def\Sing{\mathop{\bf Sing}\nolimits}
\def\Alg{\mathop{\bf Alg}\nolimits}
\def\Der{\mathop{\rm Der}\nolimits}
\def\Tw{\mathop{\rm Tw}\nolimits}
\def\tTw{\mathop{\smash{\rm\ti Tw}\vphantom{T}}\nolimits}
\def\Str{\mathop{\rm Str}\nolimits}
\def\PSh{\mathop{\rm PSh}\nolimits}
\def\Rep{\mathop{\rm Rep}\nolimits}
\def\top{{\rm top}}
\def\lin{{\rm lin}}
\def\loc{{\rm loc}}
\def\tu{{\rm tu}}
\def\tc{{\rm tc}}
\def\rti{{\rm ti}}
\def\cu{{\rm cu}}
\def\wc{{\rm wc}}
\def\sc{{\rm sc}}
\def\un{{\rm un}}
\def\hp{{\rm hp}}
\def\hs{{\rm hs}}
\def\rfi{{\rm fi}}
\def\wf{{\rm wf}}
\def\cf{{\rm cf}}
\def\lc{{\rm lc}}
\def\Kan{{\rm Kan}}
\def\rsu{{\rm su}}
\def\lu{{\rm lu}}
\def\dg{{\rm dg}}
\def\virt{{\rm virt}}
\def\depth{\mathop{\rm depth}\nolimits}
\def\Im{\mathop{\rm Im}\nolimits}
\def\Ker{\mathop{\rm Ker}}
\def\Coker{\mathop{\rm Coker}}
\def\Ho{\mathop{\rm Ho}}
\def\ind{\mathop{\rm ind}}
\def\GL{\mathop{\rm GL}}
\def\Aut{\mathop{\rm Aut}}
\def\Stab{\mathop{\rm Stab}\nolimits}
\def\supp{\mathop{\rm supp}}
\def\rank{\mathop{\rm rank}\nolimits}
\def\Sym{\mathop{\rm Sym}\nolimits}
\def\Hom{\mathop{\rm Hom}\nolimits}
\def\Mor{\mathop{\rm Mor}\nolimits}
\def\Iso{\mathop{\rm Iso}\nolimits}
\def\bHom{\mathop{\bf Hom}\nolimits}
\def\bEqu{\mathop{\bf Equ}\nolimits}
\def\cHom{\mathop{\mathcal{H}om}\nolimits}
\def\cIso{\mathop{\mathcal{I}so}\nolimits}
\def\bcHom{\mathop{\bs{\mathcal{H}om}}\nolimits}
\def\bcEqu{\mathop{\bs{\mathcal{E}qu}}\nolimits}
\def\id{{\mathop{\rm id}\nolimits}}
\def\Id{{\mathop{\rm Id}\nolimits}}
\def\Inc{{\mathop{\rm Inc}\nolimits}}
\def\co{\mathop{\textstyle\coprod}\nolimits}
\def\Sh{\mathop{\rm Sh}}
\def\Obj{{\mathop{\rm Obj}}}
\def\AAlg{{\mathop{\bf A_\iy Alg}\nolimits}}
\def\dAAlg{{\mathop{\bf\dot A_\iy Alg}\nolimits}{}}
\def\hACat{{\mathop{\bf\hat A_\iy Cat}\nolimits}{}}
\def\hAAlg{{\mathop{\bf\hat A_\iy Alg}\nolimits}{}}
\def\cACat{{\mathop{\bf\check A_\iy Cat}\nolimits}{}}
\def\cAAlg{{\mathop{\bf\check A_\iy Alg}\nolimits}{}}
\def\tAAlg{{\mathop{\bf\ti A_\iy Alg}\nolimits}{}}
\def\ACat{{\mathop{\bf A_\iy Cat}\nolimits}}
\def\ASCat{{\mathop{\bf A_\iy SCat}\nolimits}}
\def\ADCat{{\mathop{{\bf A}^\De_\iy{\bf Cat}}\nolimits}{}}
\def\dADCat{{\mathop{{\bf\dot A}^\De_\iy{\bf Cat}}\nolimits}{}}
\def\dACat{{\mathop{\bf\dot A_\iy Cat}\nolimits}{}}
\def\tACat{{\mathop{\bf\ti A_\iy Cat}\nolimits}{}}
\def\AfAlg{{\mathop{\bf A^{fi}_\iy Alg}\nolimits}{}}
\def\hAfAlg{{\mathop{\bf\hat A^{fi}_\iy Alg}\nolimits}{}}
\def\dAfAlg{{\mathop{\bf\dot A^{fi}_\iy Alg}\nolimits}{}}
\def\dAfCat{{\mathop{\bf\dot A^{fi}_\iy Cat}\nolimits}{}}
\def\bAfCat{{\mathop{\bf\bar A^{fi}_\iy Cat}\nolimits}{}}
\def\bAfAlg{{\mathop{\bf\bar A^{fi}_\iy Alg}\nolimits}{}}
\def\tAfAlg{{\mathop{\bf\ti A^{fi}_\iy Alg}\nolimits}{}}
\def\AfCat{{\mathop{\bf A^{fi}_\iy Cat}\nolimits}{}}
\def\hAfCat{{\mathop{\bf\hat A^{fi}_\iy Cat}\nolimits}{}}
\def\tAfCat{{\mathop{\bf\ti A^{fi}_\iy Cat}\nolimits}{}}
\def\AwAlg{{\mathop{\bf A^{wf}_\iy Alg}\nolimits}{}}
\def\hAwAlg{{\mathop{\bf\hat A^{wf}_\iy Alg}\nolimits}{}}
\def\dAwAlg{{\mathop{\bf\dot A^{wf}_\iy Alg}\nolimits}{}}
\def\dAwCat{{\mathop{\bf\dot A^{wf}_\iy Cat}\nolimits}{}}
\def\bAwCat{{\mathop{\bf\bar A^{wf}_\iy Cat}\nolimits}{}}
\def\bAwAlg{{\mathop{\bf\bar A^{wf}_\iy Alg}\nolimits}{}}
\def\tAwAlg{{\mathop{\bf\ti A^{wf}_\iy Alg}\nolimits}{}}
\def\AwCat{{\mathop{\bf A^{wf}_\iy Cat}\nolimits}{}}
\def\hAwCat{{\mathop{\bf\hat A^{wf}_\iy Cat}\nolimits}{}}
\def\tAwCat{{\mathop{\bf\ti A^{wf}_\iy Cat}\nolimits}{}}
\def\NCat{{\mathop{\bf N_\iy Cat}\nolimits}}
\def\NDCat{{\mathop{\bf N^\De_\iy Cat}\nolimits}}
\def\DGA{{\mathop{\bf DGA}\nolimits}}
\def\DGCat{{\mathop{\bf DGCat}\nolimits}}
\def\hDGCat{{\mathop{\bf\hat DGCat}\nolimits}}
\def\DGMod{{\mathop{\bf DGMod}\nolimits}}
\def\Com{{\mathop{\bf Com}\nolimits}}
\def\Top{{\mathop{\bf Top}\nolimits}}
\def\Fun{{\mathop{\rm Fun}\nolimits}}
\def\tFun{{\mathop{\rm\tilde Fun}\nolimits}{}}
\def\hFun{{\mathop{\rm\hat Fun}\nolimits}{}}
\def\cFun{{\mathop{\rm\check Fun}\nolimits}{}}
\def\Man{{\mathop{\bf Man}}}
\def\Manc{{\mathop{\bf Man^c}}}
\def\Mon{{\mathop{\bf Mon}}}
\def\ul{\underline}
\def\bs{\boldsymbol}
\def\ge{\geqslant}
\def\le{\leqslant\nobreak}
\def\pr{\prec}
\def\pre{\preceq}
\def\su{\succ}
\def\ssu{\ggcurly}
\def\sue{\succeq}
\def\tl{\trianglelefteq}
\def\QE{{\mathbin{\mathcal{QE}}}}
\def\QI{{\mathbin{\mathcal{QI}}}}
\def\HE{{\mathbin{\mathcal{HE}}}}
\def\O{{\mathcal O}}
\def\R{{\mathbin{\mathbb R}}}
\def\Z{{\mathbin{\mathbb Z}}}
\def\Q{{\mathbin{\mathbb Q}}}
\def\N{{\mathbin{\mathbb N}}}
\def\C{{\mathbin{\mathbb C}}}
\def\CP{{\mathbin{\mathbb{CP}}}}
\def\K{{\mathbin{\mathbb K}}}
\def\RP{{\mathbin{\mathbb{RP}}}}
\def\fC{{\mathbin{\mathfrak C}\kern.05em}}
\def\fD{{\mathbin{\mathfrak D}}}
\def\fE{{\mathbin{\mathfrak E}}}
\def\fG{{\mathbin{\mathfrak G}}}
\def\fH{{\mathbin{\mathfrak H}}}
\def\fI{{\mathbin{\mathfrak I}}}
\def\fJ{{\mathbin{\mathfrak J}}}
\def\fK{{\mathbin{\mathfrak K}}}
\def\fL{{\mathbin{\mathfrak L}}}
\def\fS{{\mathbin{\mathfrak S}}}
\def\fT{{\mathbin{\mathfrak T}}}
\def\fU{{\mathbin{\mathfrak U}}}
\def\cA{{\mathcal A}}
\def\cB{{\mathcal B}}
\def\cC{{\mathcal C}}
\def\cD{{\mathcal D}}
\def\cE{{\mathcal E}}
\def\cF{{\mathcal F}}
\def\cG{{\mathcal G}}
\def\cH{{\mathcal H}}
\def\cI{{\mathcal I}}
\def\cJ{{\mathcal J}}
\def\cK{{\mathcal K}}
\def\cL{{\mathcal L}}
\def\cM{{\mathcal M}}
\def\cN{{\mathcal N}}
\def\cP{{\mathcal P}}
\def\cQ{{\mathcal Q}}
\def\cR{{\mathcal R}}
\def\cS{{\mathcal S}}
\def\cT{{\mathcal T}}
\def\cW{{\mathcal W}}
\def\cX{{\mathcal X}}
\def\cY{{\mathcal Y}}
\def\cZ{{\mathcal Z}}
\def\fV{{\mathbin{\mathfrak V}}}
\def\fW{{\mathbin{\mathfrak W}}}
\def\fX{{\mathbin{\mathfrak X}}}
\def\fY{{\mathbin{\mathfrak Y}}}
\def\fZ{{\mathbin{\mathfrak Z}}}
\def\fe{{\mathfrak e}}
\def\ff{{\mathfrak f}}
\def\fg{{\mathfrak g}}
\def\fh{{\mathfrak h}}
\def\oM{{\mathbin{\smash{\,\,\overline{\!\!\mathcal M\!}\,}}}}
\def\bcE{{\mathbin{\bs{\mathcal E}}}}
\def\bcM{{\mathbin{\bs{\mathcal M}}}}
\def\bcW{{\mathbin{\bs{\mathcal W}}}}
\def\bcX{{\mathbin{\bs{\mathcal X}}}}
\def\bcY{{\mathbin{\bs{\mathcal Y}}}}
\def\bcZ{{\mathbin{\bs{\mathcal Z}}}}
\def\cV{{\mathcal V}}
\def\cW{{\mathcal W}}
\def\bA{{\bs A}}
\def\bB{{\bs B}}
\def\bC{{\bs C}}
\def\bD{{\bs D}}
\def\bE{{\bs E}}
\def\bF{{\bs F}}
\def\bG{{\bs G}}
\def\bH{{\bs H}}
\def\bM{{\bs M}}
\def\bN{{\bs N}}
\def\bO{{\bs O}}
\def\bP{{\bs P}}
\def\bQ{{\bs Q}}
\def\bS{{\bs S}}
\def\bT{{\bs T}}
\def\bU{{\bs U}}
\def\bV{{\bs V}\kern -0.2em}
\def\bW{{\bs W}\kern -0.1em}
\def\bX{{\bs X}}
\def\bY{{\bs Y}\kern -0.1em}
\def\btX{{\bs{\ti X}}}
\def\bZ{{\bs Z}}
\def\al{\alpha}
\def\be{\beta}
\def\ga{\gamma}
\def\de{\delta}
\def\io{\iota}
\def\ep{\epsilon}
\def\la{\lambda}
\def\ka{\kappa}
\def\th{\theta}
\def\ze{\zeta}
\def\up{\upsilon}
\def\vp{\varphi}
\def\si{\sigma}
\def\om{\omega}
\def\De{\Delta}
\def\La{\Lambda}
\def\Om{\Omega}
\def\Io{{\rm I}}
\def\Ka{{\rm K}}
\def\Mu{{\rm M}}
\def\Nu{{\rm N}}
\def\Tau{{\rm T}}
\def\Up{\Upsilon}
\def\Al{{\rm A}}
\def\Be{{\rm B}}
\def\Ga{\Gamma}
\def\Ep{{\rm E}}
\def\Ze{{\rm Z}}
\def\Eta{{\rm H}}
\def\Si{\Sigma}
\def\Th{\Theta}
\def\Rho{{\rm P}}
\def\Eta{{\rm H}}
\def\pd{\partial}
\def\ts{\textstyle}
\def\st{\scriptstyle}
\def\sst{\scriptscriptstyle}
\def\w{\wedge}
\def\sm{\setminus}
\def\lt{\ltimes}
\def\bu{\bullet}
\def\sh{\sharp}
\def\op{\oplus}
\def\od{\odot}
\def\op{\oplus}
\def\ot{\otimes}
\def\ov{\overline}
\def\bigop{\bigoplus}
\def\bigot{\bigotimes}
\def\iy{\infty}
\def\es{\emptyset}
\def\ra{\rightarrow}
\def\rra{\rightrightarrows}
\def\Ra{\Rightarrow}
\def\Nearrow{\mathrel{\rotatebox[origin=c]{45}{$\Rightarrow$}}}
\def\RRa{\Rrightarrow}
\def\Longra{\Longrightarrow}
\def\ab{\allowbreak}
\def\longra{\longrightarrow}
\def\hookra{\hookrightarrow}
\def\dashra{\dashrightarrow}
\def\dia{\diamond}
\def\ha{{\ts\frac{1}{2}}}
\def\t{\times}
\def\ci{\circ}
\def\ti{\tilde}
\def\d{{\rm d}}
\def\md#1{\vert #1 \vert}
\def\bmd#1{\big\vert #1 \big\vert}
\def\an#1{\langle #1 \rangle}
\def\ban#1{\bigl\langle #1 \bigr\rangle}
%%%%%%%%%%%%%%%%%%%%%%%%%%%%%%%%%%%%%%%%%%%%%%%%%%%%%%%%%%%%%%%%%%%%%%%%
%%%%%%%%%%%%%%%%%%%%%%%%%%%%%%%%%%%%%%%%%%%%%%%%%%%%%%%%%%%%%%%%%%%%%%%%
\def\depth{\mathop{\rm depth}\nolimits}
\def\dim{\mathop{\rm dim}\nolimits}
\def\codim{\mathop{\rm codim}\nolimits}
\def\vdim{\mathop{\rm vdim}\nolimits}
\def\sign{\mathop{\rm sign}\nolimits}
\def\Im{\mathop{\rm Im}\nolimits}
\def\det{\mathop{\rm det}\nolimits}
\def\Res{\mathop{\rm Res}\nolimits}
\def\Area{\mathop{\rm Area}\nolimits}
\def\Type{\mathop{\rm Type}\nolimits}
\def\Turn{\mathop{\rm Turn}\nolimits}
\def\Ker{\mathop{\rm Ker}}
\def\Coker{\mathop{\rm Coker}}
\def\Shape{\mathop{\rm Shape}}
\def\cha{\mathop{\rm char}\nolimits}
\def\Spec{\mathop{\rm Spec}}
\def\Perf{\mathop{\rm Perf}}
\def\Vect{\mathop{\rm Vect}}
\def\LCon{\mathop{\rm LCon}}
\def\Flag{\mathop{\rm Flag}\nolimits}
\def\FlagSt{\mathop{\rm FlagSt}\nolimits}
\def\Ham{\mathop{\rm Ham}\nolimits}
\def\Simp{\mathop{\rm Simp}\nolimits}
\def\OFA{\mathop{\text{\rm OrFi$A_\iy$-alg}}\nolimits}
\def\SOFA{\mathop{\text{\rm StOrFi$A_\iy$-alg}}\nolimits}
\def\Lmodgr{\mathop{\text{\rm $\La_\nov$-mod}^{\rm gr}}\nolimits}
\def\Iso{\mathop{\rm Iso}\nolimits}
\def\Aut{\mathop{\rm Aut}}
\def\End{\mathop{\rm End}}
\def\Exp{\mathop{\rm Exp}\nolimits}
\def\tExp{\mathop{\smash{\,\,\,\widetilde{\!\!\!\rm Exp\!\!}}\,\,}\nolimits}
\def\hExp{\mathop{\smash{\,\,\,\widehat{\!\!\!\rm Exp\!\!}}\,\,}\nolimits}
\def\Inv{\mathop{\rm Inv}}
\def\Ho{\mathop{\rm Ho}}
\def\Hol{\mathop{\rm Hol}}
\def\PGL{\mathop{\rm PGL}\nolimits}
\def\GL{\mathop{\rm GL}\nolimits}
\def\SL{\mathop{\rm SL}}
\def\SO{\mathop{\rm SO}}
\def\SU{\mathop{\rm SU}}
\def\Sp{\mathop{\rm Sp}}
\def\Re{\mathop{\rm Re}}
\def\Ch{\mathop{\rm Ch}\nolimits}
\def\ev{{\rm ev}}
\def\ror{{\rm or}}
\def\Ev{\mathop{\rm Ev}\nolimits}
\def\rod{{\rm od}}
\def\Pin{\mathop{\rm Pin}}
\def\Spin{\mathop{\rm Spin}}
\def\Spinc{\mathop{\rm Spin^c}}
\def\Man{{\mathop{\bf Man}}}
\def\Manc{{\mathop{\bf Man^c}}}
\def\MC{{\mathop{\mathcal{MC}}\nolimits}}
\def\SF{\mathop{\rm SF}}
\def\SFa{\mathop{\rm SF_{al}}\nolimits}
\def\SFai{\mathop{\rm SF_{al}^{ind}}\nolimits}
\def\CF{\mathop{\rm CF}\nolimits}
\def\CFi{\mathop{\rm CF^{ind}}\nolimits}
\def\Tr{\mathop{\rm Tr}}
\def\U{{\mathbin{\rm U}}}
\def\vol{\mathop{\rm vol}\nolimits}
\def\inc{\mathop{\rm inc}\nolimits}
\def\iinc{\mathop{\text{\rm\i nc}}\nolimits}
\def\rin{{\rm in}}
\def\ror{{\rm or}}
\def\rgr{{\rm gr}}
\def\out{{\rm out}}
\def\go{{\rm go}}
\def\nov{{\rm nov}}
\def\CY{{\rm\sst CY}}
\def\na{{\rm na}}
\def\ad{{\rm ad}}
\def\stk{{\rm stk}}
\def\ch{\mathop{\rm ch}\nolimits}
\def\ind{\mathop{\rm ind}\nolimits}
\def\tind{{\text{\rm t-ind}}}
\def\bdim{{\mathbin{\bf dim}\kern.1em}}
\def\rk{\mathop{\rm rk}}
\def\Pic{\mathop{\rm Pic}}
\def\colim{\mathop{\rm colim}\nolimits}
\def\Stab{\mathop{\rm Stab}\nolimits}
\def\Exact{\mathop{\rm Exact}\nolimits}
\def\Crit{\mathop{\rm Crit}}
\def\Cass{\mathop{\rm Cass}}
\def\supp{\mathop{\rm supp}}
\def\Or{\mathop{\rm Or}\nolimits}
\def\rank{\mathop{\rm rank}\nolimits}
\def\Hom{\mathop{\rm Hom}\nolimits}
\def\bHom{\mathop{\bf Hom}\nolimits}
\def\Ext{\mathop{\rm Ext}\nolimits}
\def\cExt{\mathop{{\mathcal E}\mathit{xt}}\nolimits}
\def\fExact{\mathop{\mathfrak{Exact}}}
\def\id{{\mathop{\rm id}\nolimits}}
\def\Id{{\mathop{\rm Id}\nolimits}}
\def\Sch{\mathop{\bf Sch}\nolimits}
\def\Map{{\mathop{\rm Map}\nolimits}}
\def\AlgSp{\mathop{\bf AlgSp}\nolimits}
\def\Lag{\mathop{\rm Lag}\nolimits}
\def\glb{\mathop{\rm glb}\nolimits}
\def\lub{\mathop{\rm lub}\nolimits}
\def\area{\mathop{\rm area}\nolimits}
\def\tLag{\mathop{\,\,\widetilde{\smash{\!\!\rm Lag\!\!}\vphantom{-}}\,\,}\nolimits}
\def\Sta{\mathop{\bf Sta}\nolimits}
\def\AlgSp{\mathop{\bf AlgSp}\nolimits}
\def\Art{\mathop{\bf Art}\nolimits}
\def\DMSta{\mathop{\bf DMSta}\nolimits}
\def\HSta{\mathop{\bf HSta}\nolimits}
\def\dSta{\mathop{\bf dSta}\nolimits}
\def\dArt{\mathop{\bf dArt}\nolimits}
\def\VertAlg{\mathop{\bf VertAlg}\nolimits}
\def\LieAlg{\mathop{\bf LieAlg}\nolimits}
\def\TopSta{{\mathop{\bf TopSta}\nolimits}}
\def\Topho{{\mathop{\bf Top^{ho}}}}
\def\Gpds{{\mathop{\bf Gpds}\nolimits}}
\def\Aff{\mathop{\bf Aff}\nolimits}
\def\SPr{\mathop{\bf SPr}\nolimits}
\def\SSet{\mathop{\bf SSet}\nolimits}
\def\Ad{\mathop{\rm Ad}}
\def\pl{{\rm pl}}
\def\eu{{\rm eu}}
\def\mix{{\rm mix}}
\def\rint{{\rm int}}
\def\fd{{\rm fd}}
\def\wt{{\rm wt}}
\def\ran{{\rm an}}
\def\fpd{{\rm fpd}}
\def\pfd{{\rm pfd}}
\def\rbig{{\rm big}}
\def\coa{{\rm coa}}
\def\rsi{{\rm si}}
\def\ssi{{\rm ssi}}
\def\rst{{\rm st}}
\def\rop{{\rm op}}
\def\ss{{\rm ss}}
\def\vi{{\rm vi}}
\def\smq{{\rm smq}}
\def\rsm{{\rm sm}}
\def\cla{{\rm cla}}
\def\rArt{{\rm Art}}
\def\po{{\rm po}}
\def\pe{{\rm pe}}
\def\rp{{\rm rp}}
\def\spo{{\rm spo}}
\def\Kur{{\rm Kur}}
\def\dcr{{\rm dcr}}
\def\top{{\rm top}}
\def\fc{{\rm fc}}
\def\cla{{\rm cla}}
\def\num{{\rm num}}
\def\irr{{\rm irr}}
\def\red{{\rm red}}
\def\sing{{\rm sing}}
\def\virt{{\rm virt}}
\def\inv{{\rm inv}}
\def\fund{{\rm fund}}
\def\qcoh{{\rm qcoh}}
\def\coh{{\rm coh}}
\def\vect{{\rm vect}}
\def\lft{{\rm lft}}
\def\lfp{{\rm lfp}}
\def\cs{{\rm cs}}
\def\St{{\rm St}}
\def\BM{{\sst\rm BM}}
\def\dR{{\rm dR}}
\def\cdga{{\mathop{\bf cdga}\nolimits}}
\def\Rmod{\mathop{R\text{\rm -mod}}}
\def\Top{{\mathop{\bf Top}\nolimits}}
\def\modKQ{\mathop{\text{\rm mod-}\K Q}}
\def\modKQI{\text{\rm mod-$\K Q/I$}}
\def\modCQ{\mathop{\text{\rm mod-}\C Q}}
\def\modCdQ{\mathop{\text{\rm mod-}\C\dot Q}}
\def\modCQI{\text{\rm mod-$\C Q/I$}}
\def\modCtQ{\mathop{\text{\rm mod-}\C\ti Q}}
\def\ul{\underline}
\def\bs{\boldsymbol}
\def\boo{{\mathbin{\mathbbm 1}}}
\def\ds{/\!/}
\def\O{{\mathcal O}}
\def\bA{{\mathbin{\mathbb A}}}
\def\bG{{\mathbin{\mathbb G}}}
\def\bH{{\mathbin{\mathbb H}}}
\def\bL{{\mathbin{\mathbb L}}}
\def\P{{\mathbin{\mathbb P}}}
\def\K{{\mathbin{\mathbb K}}}
\def\R{{\mathbin{\mathbb R}}}
\def\bT{{\mathbin{\mathbb T}}}
\def\Z{{\mathbin{\mathbb Z}}}
\def\bP{{\mathbin{\mathbb P}}}
\def\Q{{\mathbin{\mathbb Q}}}
\def\N{{\mathbin{\mathbb N}}}
\def\C{{\mathbin{\mathbb C}}}
\def\CP{{\mathbin{\mathbb{CP}}}}
\def\KP{{\mathbin{\mathbb{KP}}}}
\def\RP{{\mathbin{\mathbb{RP}}}}
\def\fC{{\mathbin{\mathfrak C}\kern.05em}}
\def\fD{{\mathbin{\mathfrak D}}}
\def\fE{{\mathbin{\mathfrak E}}}
\def\fF{{\mathbin{\kern.1em\mathfrak F\kern.1em}}}
\def\A{{\mathbin{\mathcal A}}}
\def\cAA{{\mathbin{\hat{\mathcal A}}}}
\def\tiA{{\mathbin{\ti{\mathcal A}}}}
\def\dA{{\mathbin{\dot{\mathcal A}}}}
\def\ddA{{\mathbin{\ddot{\mathcal A}}}}
\def\dddA{{\mathbin{\dddot{\mathcal A}}}}
\def\baA{{\mathbin{\smash{\bar{\mathcal A}}}\vphantom{\mathcal A}}}
\def\grA{{\mathbin{\smash{\grave{\mathcal A}}}\vphantom{\mathcal A}}}
\def\brA{{\mathbin{\smash{\breve{\mathcal A}}}\vphantom{\mathcal A}}}
\def\acA{{\mathbin{\smash{\acute{\mathcal A}}}\vphantom{\mathcal A}}}
\def\G{{{\mathcal G}}}
\def\M{{\mathbin{\mathcal M}}}
\def\B{{\mathbin{\mathcal B}}}  
\def\haB{{\mathbin{\hat{\mathcal B}}}}
\def\tiB{{\mathbin{\ti{\mathcal B}}}}
\def\dB{{\mathbin{\dot{\mathcal B}}}}
\def\ddB{{\mathbin{\ddot{\mathcal B}}}}
\def\dddB{{\mathbin{\dddot{\mathcal B}}}}
\def\baB{{\mathbin{\smash{\bar{\mathcal B}}}\vphantom{\mathcal B}}}
\def\grB{{\mathbin{\smash{\grave{\mathcal B}}}\vphantom{\mathcal B}}}
\def\brB{{\mathbin{\smash{\breve{\mathcal B}}}\vphantom{\mathcal B}}}
\def\acB{{\mathbin{\acute{\mathcal B}}}}
\def\ovB{{\mathbin{\smash{\,\overline{\!\mathcal B}}}}}
\def\cC{{\mathbin{\mathcal C}}}
\def\cD{{\mathbin{\mathcal D}}}
\def\cE{{\mathbin{\mathcal E}}}
\def\baE{{\mathbin{\smash{\bar{\mathcal E}}}\vphantom{\mathcal E}}}
\def\dE{{\mathbin{\dot{\mathcal E}}}{}}
\def\cF{{\mathbin{\mathcal F}}}
\def\baF{{\mathbin{\smash{\bar{\mathcal F}}}\vphantom{\mathcal F}}}
\def\grF{{\mathbin{\smash{\grave{\mathcal F}}}\vphantom{\mathcal F}}}
\def\acF{{\mathbin{\smash{\ac{\mathcal F}}}\vphantom{\mathcal F}}}
\def\cG{{\mathbin{\mathcal G}}}
\def\baG{{\mathbin{\smash{\bar{\mathcal G}}}\vphantom{\mathcal G}}}
\def\grG{{\mathbin{\smash{\grave{\mathcal G}}}\vphantom{\mathcal G}}}
\def\acG{{\mathbin{\smash{\ac{\mathcal G}}}\vphantom{\mathcal G}}}
\def\baH{{\mathbin{\smash{\bar{\mathcal H}}}\vphantom{\mathcal H}}}
\def\acH{{\mathbin{\smash{\ac{\mathcal H}}}\vphantom{\mathcal H}}}
\def\cH{{\mathbin{\mathcal H}}}
\def\cI{{\mathbin{\mathcal I}}}
\def\cJ{{\mathbin{\mathcal J}}}
\def\cK{{\mathbin{\mathcal K}}}
\def\cL{{\mathbin{\mathcal L}}}
\def\cB{{\mathbin{\mathcal B}}}
\def\bcM{{\mathbin{\bs{\mathcal M}}}}
\def\cN{{\mathcal N}}
\def\cP{{\mathbin{\mathcal P}}}
\def\cQ{{\mathbin{\mathcal Q}}}
\def\cR{{\mathbin{\mathcal R}}}
\def\cS{{\mathbin{\mathcal S}}}
\def\T{{{\mathcal T}\kern .04em}}
\def\cU{{\mathbin{\mathcal U}}}
\def\cV{{\mathbin{\mathcal V}}}
\def\cW{{\mathbin{\mathcal W}}}
\def\cX{{\mathcal X}}
\def\cY{{\mathcal Y}}
\def\cZ{{\mathcal Z}}
\def\oM{{\mathbin{\smash{\,\,\overline{\!\!\mathcal M\!}\,}}\vphantom{\mathcal M}}}
\def\tiM{{\mathbin{\smash{\ti{\mathcal M}}}\vphantom{\mathcal M}}}
\def\haM{{\mathbin{\smash{\hat{\mathcal M}}}\vphantom{\mathcal M}}}
\def\acM{{\mathbin{\smash{\,\,\ac{\!\!\mathcal M\!}\,}}\vphantom{\mathcal M}}}
\def\dM{{\mathbin{\smash{\dot{\mathcal M}}}\vphantom{\mathcal M}}}
\def\ddM{{\mathbin{\smash{\ddot{\mathcal M}}}\vphantom{\mathcal M}}}
\def\dddM{{\mathbin{\smash{\dddot{\mathcal M}}}\vphantom{\mathcal M}|}}
\def\baM{{\mathbin{\smash{\bar{\mathcal M}}}\vphantom{\mathcal M}}}
\def\dbaM{{\mathbin{\smash{\dot{\bar{\mathcal M}}}}\vphantom{\mathcal M}}}
\def\grM{{\mathbin{\smash{\grave{\mathcal M}}}\vphantom{\mathcal M}}}
\def\brM{{\mathbin{\smash{\breve{\mathcal M}}}\vphantom{\mathcal M}}}
\def\cV{{\mathcal V}}
\def\baV{{\mathbin{\smash{\bar{\mathcal V}}}\vphantom{\mathcal V}}}
\def\cW{{\mathcal W}}
\def\sA{{{\mathscr A}}}
\def\sB{{{\mathscr B}}}
\def\sC{{{\mathscr C}}}
\def\sD{{{\mathscr D}}}
\def\sE{{{\mathscr E}}}
\def\sF{{{\mathscr F}}}
\def\sG{{{\mathscr G}}}
\def\sH{{{\mathscr H}}}
\def\sI{{{\mathscr I}}}
\def\sJ{{{\mathscr J}}}
\def\sK{{{\mathscr K}}}
\def\sL{{{\mathscr L}}}
\def\sM{{{\mathscr M}}}
\def\sN{{{\mathscr N}}}
\def\sO{{{\mathscr O}}}
\def\sP{{{\mathscr P}}}
\def\sQ{{{\mathscr Q}}}
\def\sR{{{\mathscr R}}}
\def\sS{{{\mathscr S}}}
\def\hsS{{\smash{\,\,\,\hat{\!\!\!\mathscr S}}\vphantom{\mathscr S}}}
\def\tsS{{\smash{\,\,\,\ti{\!\!\!\mathscr S}}\vphantom{\mathscr S}}}
\def\baS{\mathbin{\smash{\,\,\,\bar{\!\!\!\mathscr S}}\vphantom{\mathcal S}}}
\def\acS{\mathbin{\smash{\,\,\,\acute{\!\!\!\mathscr S}}\vphantom{\mathcal S}}}
\def\grS{\mathbin{\smash{\,\,\,\grave{\!\!\!\mathscr S}}\vphantom{\mathcal S}}}
\def\baB{{\mathbin{\smash{\bar{\mathcal B}}}\vphantom{\mathcal B}}}
\def\sT{{{\mathscr T}}}
\def\sU{{{\mathscr U}}}
\def\sV{{{\mathscr V}}}
\def\sW{{{\mathscr W}}}
\def\baW{\mathbin{\smash{\,\,\,\bar{\!\!\!\mathscr W}}\vphantom{\mathscr W}}}
\def\sY{{{\mathscr Y}}}
\def\b{{\mathfrak b}}
\def\fe{{\mathfrak e}}
\def\f{{\mathfrak f}}
\def\g{{\mathfrak g}}
\def\h{{\mathfrak h}}
\def\k{{\mathfrak k}}
\def\m{{\mathfrak m}}
\def\n{{\mathfrak n}}
\def\o{{\mathfrak o}}
\def\p{{\mathfrak p}}
\def\q{{\mathfrak q}}
\def\u{{\mathfrak u}}
\def\H{{\mathfrak H}}
\def\so{{\mathfrak{so}}}
\def\sp{{\mathfrak{sp}}}
\def\fW{{\mathfrak W}}
\def\fX{{\mathfrak X}}
\def\fY{{\mathfrak Y}}
\def\fZ{{\mathfrak Z}}
\def\bfb{{\bs{\mathfrak b}}}
\def\bfc{{\bs{\mathfrak c}}}
\def\bfd{{\bs{\mathfrak d}}}
\def\bfe{{\bs{\mathfrak e}}}
\def\bff{{\bs{\mathfrak f}}}
\def\bfg{{\bs{\mathfrak g}}}
\def\bfh{{\bs{\mathfrak h}}}
\def\bfU{{\bs{\mathfrak U}}}
\def\bfV{{\bs{\mathfrak V}}}
\def\bfW{{\bs{\mathfrak W}}}
\def\bfX{{\bs{\mathfrak X}}}
\def\bfY{{\bs{\mathfrak Y}}}
\def\bfZ{{\bs{\mathfrak Z}}}
\def\bM{{\bs M}}
\def\bN{{\bs N}}
\def\bO{{\bs O}}
\def\bQ{{\bs Q}}
\def\bS{{\bs S}}
\def\bU{{\bs U}}
\def\bV{{\bs V}}
\def\bW{{\bs W}\kern -0.1em}
\def\bX{{\bs X}}
\def\bY{{\bs Y}\kern -0.1em}
\def\bZ{{\bs Z}}
\def\ts{\textstyle}
\def\st{\scriptstyle}
\def\sst{\scriptscriptstyle}
\def\w{\wedge}
\def\sm{\setminus}
\def\lt{\ltimes}
\def\bu{\bullet}
\def\sh{\sharp}
\def\di{\diamond}
\def\he{\heartsuit}
\def\od{\odot}
\def\op{\oplus}
\def\ot{\otimes}
\def\hot{\mathbin{\hat\otimes}}
\def\bt{\boxtimes}
\def\bp{\boxplus}
\def\ov{\overline}
\def\bigop{\bigoplus}
\def\bigot{\bigotimes}
\def\tr{\blacktriangle}
\def\sq{{\scriptstyle\blacksquare}}
\def\iy{\infty}
\def\es{\emptyset}
\def\ra{\rightarrow}
\def\rra{\rightrightarrows}
\def\Ra{\Rightarrow}
\def\Longra{\Longrightarrow}
\def\ab{\allowbreak}
\def\longra{\longrightarrow}
\def\hookra{\hookrightarrow}
\def\dashra{\dashrightarrow}
\def\lb{\llbracket}
\def\pf{\pitchfork}
\def\rb{\rrbracket}
\def\ha{{\ts\frac{1}{2}}}
\def\t{\times}
\def\ci{\circ}
\def\ti{\tilde}
\def\ac{\acute}
\def\gr{\grave}
\def\d{{\rm d}}
\def\D{{\rm D}}
\def\md#1{\vert #1 \vert}
\def\ms#1{\vert #1 \vert^2}
\def\bmd#1{\big\vert #1 \big\vert}
\def\bms#1{\big\vert #1 \big\vert^2}
\def\an#1{\langle #1 \rangle}
\def\ban#1{\bigl\langle #1 \bigr\rangle}
\def\nm#1{\Vert #1 \Vert}
\def\bnm#1{\big\Vert #1 \big\Vert}
%%%%%%%%%%%%%%%%%%%%%%%%%%%%%%%%%%%%%%%%%%%%%%%%%%%%%%%%%%%%%%%%%%%%%%%%
%%%%%%%%%%%%%%%%%%%%%%%%    Text of book    %%%%%%%%%%%%%%%%%%%%%%%%%%%%
%%%%%%%%%%%%%%%%%%%%%%%%%%%%%%%%%%%%%%%%%%%%%%%%%%%%%%%%%%%%%%%%%%%%%%%%
\title{An $\iy$-category approach to \\ $A_\iy$-algebras and $A_\iy$-categories}
\author{Dominic Joyce}
\date{}
\maketitle

\begin{abstract}
$A_\iy$-algebras over a commutative ring $R$ are generalizations of $R$-algebras in which multiplication is only associative up to higher homotopies. Similarly, $A_\iy$-categories over $R$ are generalizations of $R$-linear categories in which composition of morphisms is only associative up to higher homotopies. They are important in the definition of Fukaya categories in Symplectic Geometry, and in other areas of mathematics.

By $\iy$-categories we mean quasicategories, a class of simplicial sets which are a model for $(\iy,1)$-categories, used in the work of Jacob Lurie.

We define and study $\iy$-categories of $A_\iy$-algebras and $A_\iy$-cat\-eg\-or\-ies. Our approach is unusual in that rather than defining them by localization, we construct our $\iy$-categories completely explicitly as simplicial sets.

Our $\iy$-categories are intended for applications to Fukaya categories in Symplectic Geometry. Given a suitable symplectic manifold $(X,\om)$, one should construct its Fukaya category $\sF(X,\om)$, a (filtered) $A_\iy$-category, and its derived Fukaya category $D^b\sF(X,\om)$, a triangulated category (or better, a stable $\iy$-category). Typically $\sF(X,\om)$ depends on many arbitrary choices. Our $\iy$-categories give a clean way to explain in what sense $\sF(X,\om)$ and $D^b\sF(X,\om)$ are independent of these choices.
\end{abstract}

\setcounter{tocdepth}{2}
\tableofcontents

\section{Introduction}
\label{ai1}

$A_\iy$-{\it algebras\/} were invented in 1963 by Stasheff \cite{Stas1,Stas2} as a tool in the study of `group-like' topological spaces in homotopy theory, and their study was mostly confined to Algebraic Topology until the 1990's. $A_\iy$-{\it categories\/} were introduced by Fukaya \cite{Fuka1} in 1993 to define Fukaya categories in Symplectic Geometry. These formed part of Kontsevich's Homological Mirror Symmetry Conjecture \cite{Kont} in 1995, and then received much attention. $A_\iy$-algebras and $A_\iy$-categories are now important in many areas of mathematics.

By $\iy$-{\it category\/} we mean quasicategories, a special kind of simplicial set which are now the most popular model for $(\iy,1)$-categories. They were introduced by Boardman and Vogt \cite{BoVo}, and developed by Joyal \cite{Joya2} and Lurie \cite{Luri2,Luri3,Luri5}. This book is heavily influenced by Lurie's philosophy and results.

$A_\iy$-algebras and $A_\iy$-categories belong in an $\iy$-categorical setting since one usually cares about them not up to isomorphism but up to `quasi-equivalence' or `homotopy equivalence', and $\iy$-categories generalize homotopy theory to a wider context. In applications it is usually necessary to make some arbitrary choices to construct an $A_\iy$-algebra or $A_\iy$-category, and different choices yield homotopy equivalent (but not necessarily isomorphic) answers.

There is a literature on $\iy$-categories of $A_\iy$-categories \cite{COS1,COS2,Faon2,OhTa,Orna1,Tana,Pasc}, reviewed in \S\ref{ai410}--\S\ref{ai411}. Most previous work has studied $\iy$-categories of the form $\ACat_R^\un[\QE_\un^{-1}]^\iy$, in which one begins with an ordinary category $\ACat_R^\un$ of (for example, strictly unital) $A_\iy$-categories over a commutative ring $R$, and then localizes at the `quasi-equivalences' $\QE_\un$ to get an $\iy$-category. Note that $\ACat_R^\un[\QE_\un^{-1}]^\iy$ is not really known explicitly: by general abstract theory, it exists as the solution to a complicated $\iy$-categorical universal property, and the challenge is to prove something useful about it. 

Our approach is very different to most of the existing literature, and indeed to the mainstream of work on $\iy$-categories. We will define our $\iy$-categories, such as the $\iy$-category $\ACat_R^{\un,\iy}$ of unital $A_\iy$-categories over $R$ in \S\ref{ai72}, {\it completely explicitly as simplicial sets}, by saying what the $N$-simplices are. This is unusual: conventional wisdom is that one uses $\iy$-categories for abstract theory, but model categories for calculations and examples.

We quote from the introduction of Land's useful book \cite{Land} on $\iy$-categories:
\begin{quotation}
\noindent `\ldots the theory of $\iy$-categories is useful in practice. However, it also has some drawbacks. In my opinion, the biggest one is the following (exaggerated on purpose):
\smallskip

\noindent\raisebox{2pt}{\begin{large}{\fontencoding{U}\fontfamily{futs}\selectfont\char 49\relax}\end{large}} It is practically impossible to construct anything by hand.'
\end{quotation}
Most of this book is a counterexample to this, as we construct by hand many $\iy$-categories and $\iy$-functors, so our theory serves as a worked example of parts of $\iy$-category theory. Our definitions are usually quite complicated.

\subsection{\texorpdfstring{Background on $A_\iy$-algebras and $A_\iy$-categories}{Background on A∞-algebras and A∞-categories}}
\label{ai1.1}

We briefly recall some basic definitions and facts for $A_\iy$-algebras and $A_\iy$-cat\-eg\-or\-ies, which are given in much more detail in Chapters \ref{ai3}--\ref{ai4}. This material is all due to other authors, attributed in Chapters \ref{ai3}--\ref{ai4}. Let $R$ be a commutative ring throughout.

\begin{dfn}
\label{ai1def1}
 An $A_\iy$-{\it category} $\sA$ over $R$ consists of:
\begin{itemize}
\setlength{\itemsep}{0pt}
\setlength{\parsep}{0pt}
\item[(a)] A set, or proper class, $\Obj\sA$ of {\it objects\/} of $\sA$.
\item[(b)] For all objects $x,y\in\Obj\sA$, a $\Z$-graded $R$-module
\begin{equation*}
\Hom(x,y)=\bigop\nolimits_{d\in\Z}\Hom^d(x,y).
\end{equation*}
We also write $\Hom_\sA(x,y)$ when we want to identify the $A_\iy$-category $\sA$.
\item[(c)] For all $k\ge 1$ and objects $x_0,x_1,\ldots,x_k$ in $\Obj\sA$, a graded $R$-linear map
\begin{equation*}
\m_k^{x_0,\ldots,x_k}:\Hom(x_0,x_1)\ot_R\cdots \ot_R\Hom(x_{k-1},x_k)\longra\Hom(x_0,x_k)
\end{equation*}
of degree $2-k$. These should satisfy for all $n\ge 1$ and objects $x_0,\ldots,x_n$
\end{itemize}
\ea
\sum_{\begin{subarray}{l} i,k\ge 0,\; j\ge 1:\\ i+j+k=n\end{subarray}}\begin{aligned}[t]&
(-1)^{i+jk}\m_{i+k+1}^{x_0,\ldots,x_i,x_{i+j},\ldots,x_n}\ci\bigl(\id_{\Hom(x_0,x_1)}\ot\cdots\ot\id_{\Hom(x_{i-1},x_i)}\ot {}\\
&\m_j^{x_i,\ldots,x_{i+j}}\ot \id_{\Hom(x_{i+j},x_{i+j+1})}\ot\cdots\ot\id_{\Hom(x_{n-1},x_n)}\bigr)=0.
\end{aligned}
\nonumber\\[-17pt]
\label{ai1eq1}
\ea

Equation \eq{ai1eq1} for $n=1$ becomes $\m_1^{x_0,x_1}\ci\m_1^{x_0,x_1}=0$, with $\m_1^{x_0,x_1}$ of degree 1. Thus $\bigl(\Hom(x_0,x_1),\m_1^{x_0,x_1}\bigr)$ is a complex, and we form its cohomology $H^*(\Hom(x_0,x_1))$. For all objects $x_0,x_1,x_2$ in $\sA$ we have a well-defined, associative, $R$-bilinear product, which we call {\it composition\/} $\ci$
\e
\ci:H^p(\Hom(x_0,x_1))\t H^q(\Hom(x_1,x_2))\longra H^{p+q}(\Hom(x_0,x_2))
\label{ai1eq2}
\e
given by
\begin{equation*}
(a_1+\Im\m^{x_0,x_1}_1)\ci(a_2+\Im\m^{x_1,x_2}_1)=\m^{x_0,x_1,x_2}_2(a_1\ot a_2)+\Im\m^{x_0,x_2}_1.
\end{equation*}

We can now define a {\it non-unital\/} $R$-linear category $H^0(\sA)$ to have objects $x$ in $\Obj\sA$, morphism spaces $\Hom_{H^0(\sA)}(x,y)=H^0(\Hom_\sA(x,y))$, and composition of morphisms $\ci:H^0(\Hom(x,y))\t H^0(\Hom(y,z))\ab\ra H^0(\Hom(x,z))$ as in \eq{ai1eq2}. We call $H^0(\sA)$ the {\it cohomology category\/} of $\sA$. Note that in general $H^0(\sA)$ is {\it not a category\/} in the usual sense, as there need not exist units $1_x\in\Hom_{H^0(\sA)}(x,x)$ for $x\in\Obj\sA$.
\end{dfn}

\begin{dfn}
\label{ai1def2}
An $A_\iy$-category $\sA$ over $R$ is called {\it unital\/} (or {\it strictly unital\/}) if for all $x\in\Obj\sA$ there exists an {\it identity\/} $\boo_x\in\Hom^0(x,x)$ satisfying
\e
\m_2^{x_0,x_0,x_1}(\boo_{x_0}\ot a)=\m_2^{x_0,x_1,x_1}(a\ot \boo_{x_1})=a
\label{ai1eq3}
\e
for all objects $x_0,x_1$ in $\sA$ and $a\in\Hom(x_0,x_1)$, and
\e
\m_{k+1}^{x_0,\ldots,x_i,x_i,\ldots,x_k}(a_1\ot\cdots\ot a_i\ot \boo_{x_i}\ot a_{i+1}\ot\! \cdots \ot a_k)=0
\label{ai1eq4}
\e
for all $k\ne 1$, $i=0,\ldots,k$, $x_0,\ldots,x_k\in\Obj\sA$ and $a_j\in\Hom(x_{j-1},x_j)$. Identities are unique if they exist, so to be unital is a property of~$\sA$.

When $\sA$ is unital, the cohomology category $H^0(\sA)$ has identities $[\boo_x]$ in $H^0(\Hom(x,x))$, so it is an honest $R$-linear category.
\end{dfn}

There are several notions of unitality for $A_\iy$-categories, of which our unital is the strongest. In \S\ref{ai41} we also define {\it cohomologically unital\/} and {\it weakly cohomologically unital\/} $A_\iy$-categories. Note that names are not consistent across the literature: some authors call our unital $A_\iy$-categories `strictly unital', and some authors use `unital' to mean either cohomologically unital, or weakly cohomologically unital, in our sense.

\begin{dfn}
\label{ai1def3}
Let $\sA,\sB$ be $A_\iy$-categories over $R$. An $A_\iy$-{\it functor\/} $\fF:\sA\ra\sB$ consists of:
\begin{itemize}
\setlength{\itemsep}{0pt}
\setlength{\parsep}{0pt}
\item[(a)] A map $\fF_\Obj:\Obj\sA\ra\Obj\sB$, sometimes just written~$\fF$.
\item[(b)] For all $k\ge 1$ and objects $x_0,x_1,\ldots,x_k$ in $\Obj\sA$, a graded $R$-linear map
\begin{align*}
\fF_k^{x_0,\ldots,x_k}:\Hom_\sA(x_0,x_1)\ot_R\cdots \ot_R\Hom_\sA(x_{k-1},x_k)&\\
\longra\Hom_\sB(\fF_\Obj(x_0),\fF_\Obj(x_k))&,
\end{align*}
of degree $1-k$. These must satisfy for all $n\ge 1$ and $x_0,\ldots,x_n\in\Obj\sA$
\end{itemize}
\ea
&\sum_{\begin{subarray}{l} i,k\ge 0,\; j\ge 1:\\ i+j+k=n \end{subarray}}
\begin{aligned}[t]&
(-1)^{i+jk}\fF_{i+k+1}^{x_0,\ldots,x_i,x_{i+j},\ldots,x_n}\ci\bigl(\id_{\Hom_\sA(x_0,x_1)}\ot\cdots\ot\id_{\Hom_\sA(x_{i-1},x_i)}\\
&{}\ot\m_j^{\sA;x_i,\ldots,x_{i+j}}\ot \id_{\Hom_\sA(x_{i+j},x_{i+j+1})}\ot\cdots\ot\id_{\Hom_\sA(x_{n-1},x_n)}\bigr)\end{aligned} 
\nonumber\\ 
&=\sum_{\begin{subarray}{l}  n=n_1+\cdots+n_k,\\ n_a\ge 1\end{subarray}}
\begin{aligned}[t]&(-1)^{\sum_{i=1}^k(n_i+1)(k-i)}\m^{\sB;\fF_\Obj(x_0),\fF_\Obj(x_{n_1}),\ldots,\fF_\Obj(x_n),}_k\bigl(\fF_{n_1}^{x_0,\ldots,x_{n_1}}\\
&{}\ot\fF_{n_2}^{x_{n_1},\ldots,x_{n_1+n_2}}\ot\cdots\ot\fF_{n_k}^{x_{n-n_k},\ldots,x_n}).\end{aligned} 
\nonumber\\[-16pt] 
\label{ai1eq5}
\ea

The {\it identity functor\/} $\id^\sA:\sA\ra\sA$ has $\id^\sA_\Obj=\id_{\Obj\sA}$, and $\id^{\sA,x}_1=\id:\Hom_\sA(x,x)\ra\Hom_\sA(x,x)$, and $\id^\sA_k=0$ for $k>1$.

When $n=1$, equation \eq{ai1eq5} becomes
\begin{align*}
&\fF^{x_0,x_1}_1\ci\m^{\sA;x_0,x_1}_1\!=\!\m_1^{\sB;\fF(x_0),\fF(x_1)}\ci\fF^{x_0,x_1}_1:\\
&\qquad \Hom_\sA(x_0,x_1)\longra \Hom_\sB(\fF(x_0),\fF(x_1)).
\end{align*}
Thus $\fF^{x_0,x_1}_1$ induces morphisms on cohomology 
\e
H^k(\fF^{x_0,x_1}_1):H^k\bigl(\Hom_\sA(x_0,x_1)\bigr)\longra H^k\bigl(\Hom_\sB(\fF(x_0),\fF(x_1))\bigr).
\label{ai1eq6}
\e

If $\fF:\sA\ra\sB$, $\fG:\sB\ra\sC$ are functors of $A_\iy$-categories over $R$, the {\it composition} $\fG\ci\fF:\sA\ra\sC$ is given by $(\fG\ci\fF)_\Obj=\fG_\Obj\ci\fF_\Obj$ and
\ea
&(\fG\ci\fF)_n^{x_0,\ldots,x_n}=
\label{ai1eq7}\\
&\sum_{\begin{subarray}{l} n=n_1+\cdots+n_k,\\ n_a\ge 1\end{subarray}}
\begin{aligned}[t]& (-1)^{\sum_{i=1}^k(n_i+1)(k-i)}\fG_k^{\fF(x_0),\fF(x_{n_1}),\fF(x_{n_1+n_2}),\ldots,\fF(x_n)} \\ 
&\bigl(\fF_{n_1}^{x_0,\ldots,x_{n_1}}\ot\fF_{n_2}^{x_{n_1},\ldots,x_{n_1+n_2}}\ot\cdots\ot\fF_{n_k}^{x_{n-n_k},\ldots,x_n}),
\end{aligned}
\nonumber
\ea
for all $n\ge 1$ and $x_0,\ldots,x_n\in\Obj\sA$. Composition is associative.

If $\sA,\sB$ are unital, we call an $A_\iy$-functor $\fF:\sA\ra\sB$ {\it unital\/} if
$\fF_1^{x,x}(\boo_x)=\boo_{\fF_\Obj(x)}$ for all objects $x$ in $\sA$, and
\begin{equation*}
\fF_{k+1}^{x_0,\ldots,x_i,x_i,\ldots,x_k}(a_1\ot\cdots\ot a_i\ot \boo_{x_i}\ot a_{i+1}\ot \cdots \ot a_k)=0
\end{equation*}
for all $k\ge 1$, $x_0,\ldots,x_k\in\Obj\sA$, $i=0,\ldots,k$ and $a_j\in\Hom_\sA(x_{j-1},x_j)$ for $j=1,\ldots,k$. 

Write $\ACat_R$ for the category with objects $A_\iy$-categories $\sA$ over $R$ and morphisms $A_\iy$-functors $\fF:\sA\ra\sB$. Write $\ACat_R^\un\subset\ACat_R$ for the non-full subcategory with objects unital $A_\iy$-categories and morphisms unital $A_\iy$-functors. In \S\ref{ai41}--\S\ref{ai42} we also define further subcategories
\e
\ACat_R^\un\subset\ACat_R^\cu\subset\ACat_R^\wc\subset\ACat_R
\label{ai1eq8}
\e
with different notions of unitality for objects and morphisms. 

Here $\ACat_R^\cu$ has objects {\it cohomologically unital\/} $A_\iy$-categories, that is, units $\boo_x\in\Hom^0(x,x)$ need not exist, but there should exist (nonunique) {\it cohomological units\/} $e_x\in\Hom^0(x,x)$ for $x\in\Obj\sA$ such that left and right multiplication by $e_x,e_y$ on the complex $(\Hom^*(x,y),\m_1^{x,y})$ using $\m_2^{x,x,y},\m_2^{x,y,y}$ are chain homotopic to the identity, and $\ACat_R^\wc$ has objects {\it weakly cohomologically unital\/} $A_\iy$-categories $\sA$, that is, there exist (unique) units $1_x\in H^0(\Hom^*(x,x),\m_1^{x,x})$ making the cohomology category $H^0(\sA)$ unital.

In both of $\ACat_R^\cu,\ACat_R^\wc$, morphisms are {\it cohomologically unital\/} $A_\iy$-functors $\fF:\sA\ra\sB$, that is, $H^0(\fF):H^0(\sA)\ra H^0(\sB)$ takes identities $[\boo_x],[e_x],1_x$ in $H^0(\sA)$ to identities in~$H^0(\sB)$.

If $\fF:\sA\ra\sB$ is a morphism in $\ACat_R^\un$, the morphisms \eq{ai1eq6} induce an $R$-linear functor $H^0(\fF_1):H^0(\sA)\ra H^0(\sB)$ of ordinary $R$-linear categories.
\end{dfn}

\begin{dfn}
\label{ai1def4}
A unital $A_\iy$-category $\sA$ over $R$ has {\it split units\/} if for all $x\in\Obj\sA$ the map $R\ra\Hom^0(x,x)$, $\la\mapsto\la\boo_x$, is injective and splits as a map of $R$-modules. An $A_\iy$-category over $R$ is called {\it h-projective\/} if for all $x,y\in\Obj\sA$, the complex $(\Hom(x,y),\m_1^{x,y})$ is h-projective (see Definition~\ref{ai4def19}).

We write $\ACat^{\un,\hs}_R\subset\ACat^{\un,\rsu}_R\subset\ACat^\un_R$ for the full subcategories of objects with split units, and of h-projective objects with split units. In the cohomologically unital case, we write $\ACat^{\cu,\hp}_R\subset\ACat^\cu_R$ for the full subcategory of h-projective objects. 

A unital $A_\iy$-category $\sA$ behaves like a `cofibrant object', in a model-the\-or\-et\-ic sense, if it is h-projective with split units. A cohomologically unital $A_\iy$-category $\sA$ behaves like a `cofibrant object' if it is h-projective.

If $R$ is a field then every $\sA$ is h-projective. If $\sA$ is unital then it has split units provided it has no {\it zero objects}, that is, no objects $z$ with~$\boo_z=0$.
\end{dfn}

\begin{dfn}
\label{ai1def5}
A {\it differential graded category}, or {\it dg-category}, over $R$ is a unital $A_\iy$-category $\sA$ over $R$ with $\m_k^{x_0,\ldots,x_k}=0$ for all $k\ge 3$ and objects $x_0,\ldots,x_k$ in $\sA$. If $\sA,\sB$ are dg-categories, a {\it dg-functor\/} $\fF:\sA\ra\sB$ is a unital $A_\iy$-functor $\fF:\sA\ra\sB$ such that $\fF_k^{x_0,\ldots,x_k}=0$ for all $k\ge 2$ and objects $x_0,\ldots,x_k$ in~$\sA$.

We write $\DGCat_R\subset\ACat^\un_R$ for the non-full subcategory with objects dg-categories and morphisms dg-functors. We also write $\DGCat^{A_\iy}_R$ for the full subcategory of $\ACat^\un_R$ whose objects are dg-categories, but whose morphisms are $A_\iy$-functors. Then~$\DGCat_R\subset\DGCat^{A_\iy}_R\subset\ACat^\un_R$.

Dg-categories are important in many areas of mathematics where $A_\iy$-cat\-eg\-or\-ies are not obviously relevant.
\end{dfn}

\begin{dfn}
\label{ai1def6}
An $A_\iy$-{\it algebra\/} is equivalent to an $A_\iy$-category $\sA$ with one object $*$, so that $\Obj\sA=\{*\}$. We use different notation for $A_\iy$-algebras: we generally write $A=\Hom_\sA(*,*)$ as a graded $R$-module, and $\m_k=\m_k^{*,\ldots,*}$, and $\m=(\m_k)_{k\ge 1}$, and then write $A_\iy$-algebras as $(A,\m)$ rather than $\sA$. We write
\begin{equation*}
\AAlg_R^\un\subset\AAlg_R^\cu\subset\AAlg_R^\wc\subset\AAlg_R
\end{equation*}
for the full subcategories of objects $\sA$ in the categories \eq{ai1eq8} with~$\Obj\sA=\{*\}$.
\end{dfn}

The theory of $A_\iy$-algebras follows from that of $A_\iy$-categories by restricting to $A_\iy$-categories with one object $*$, but it is simpler both notationally and conceptually. We discuss $A_\iy$-algebras in Chapters \ref{ai3} and \ref{ai6} before the corresponding material for $A_\iy$-categories in Chapters \ref{ai4} and~\ref{ai7}.

\begin{dfn}
\label{ai1def7}
Let $\sA,\sB$ be $A_\iy$-categories over $R$. As in \S\ref{ai43}, we can define an $A_\iy$-category $\Fun_{A_\iy}(\sA,\sB)$ whose objects are $A_\iy$-functors $\fF:\sA\ra\sB$. If $\sA,\sB$ are unital, or cohomologically unital, there are $A_\iy$-subcategories $\Fun_{A_\iy}^\un(\sA,\sB),\Fun_{A_\iy}^\cu(\sA,\sB)$ with objects (cohomologically) unital $A_\iy$-func\-tors, which are themselves (cohomologically) unital. They extend to a functor
\begin{equation*}
\Fun_{A_\iy}:(\ACat_R)^\rop\t\ACat_R\longra\ACat_R,
\end{equation*}
and similarly for $\Fun_{A_\iy}^\un,\Fun_{A_\iy}^\cu$.
\end{dfn}

\begin{dfn}
\label{ai1def8}
Let $\sA,\sB$ be (possibly unital, or cohomologically unital) $A_\iy$-categories over $R$, and $\fF,\fG:\sA\ra\sB$ be (possibly unital, or cohomologically unital) $A_\iy$-functors. As explained in \S\ref{ai44}, there are three notions of `2-morphism' $\fH:\fF\Ra\fG$, called $A_\iy$-{\it homotopies}, $A_\iy$-{\it natural transformations}, and $A_\iy$-{\it natural isomorphisms}, which can be unital, or cohomologically unital.

Here $A_\iy$-natural transformations are (closed) morphisms in $\Fun_{A_\iy}(\sA,\sB),\ab \Fun_{A_\iy}^\un(\sA,\sB),\Fun_{A_\iy}^\cu(\sA,\sB)$. $A_\iy$-natural isomorphisms are only defined in the (cohomologically) unital cases, and are $A_\iy$-natural transformations which induce isomorphisms in $H^0(\Fun_{A_\iy}^\un(\sA,\sB))$ or $H^0(\Fun_{A_\iy}^\cu(\sA,\sB))$. $A_\iy$-homo\-top\-ies have a different definition, but can be identified with special $A_\iy$-natural isomorphisms in the unital case. An $A_\iy$-natural transformation (isomorphism) $\fH:\fF\Ra\fG$ induces an ordinary natural transformation (isomorphism) of functors~$H^0(\fH):H^0(\fF)\Ra H^0(\fG)$.

A (cohomologically) unital $A_\iy$-functor $\fF:\sA\ra\sB$ is called a {\it homotopy equivalence\/} if there exist a (cohomologically) unital $A_\iy$-functor $\fG:\sB\ra\sA$ and $A_\iy$-natural isomorphisms $\fS:\fG\ci\fF\Ra\id_\sA$ and $\fT:\fF\ci\fG\Ra\id_\sB$.

We call a (cohomologically) unital $\fF:\sA\ra\sB$ a {\it quasi-equivalence\/} if the morphisms \eq{ai1eq6} are isomorphisms for all $k,x_0,x_1$, and $H^0(\fF_1):H^0(\sA)\ra H^0(\sB)$ is an equivalence of categories. All homotopy equivalences are quasi-equivalences.

If $R$ is a field, and in the unital case $\sA,\sB$ do not have zero objects (objects $z$ with $\boo_z=0$), then homotopy equivalences and quasi-equivalences coincide. 

More generally, in the unital case, if $\sA,\sB$ are h-projective with split units then $\fF:\sA\ra\sB$ is a quasi-equivalence if and only if it is a homotopy equivalence, and in the cohomologically unital case, if $\sA,\sB$ are h-projective then $\fF:\sA\ra\sB$ is a quasi-equivalence if and only if it is a homotopy equivalence.
\end{dfn}

One can argue that for $A_\iy$-algebras, $A_\iy$-homotopy is the best notion of when two $A_\iy$-morphisms are `the same', but for $A_\iy$-categories, $A_\iy$-natural isomorphism is best. Quasi-equivalence, and homotopy equivalence, are two competing notions of when two (cohomologically) unital $A_\iy$-categories are `the same'. Notice that Definition \ref{ai1def8} already makes $A_\iy$-categories into something like a 2-category.

\begin{dfn}
\label{ai1def9}
In \S\ref{ai48} we define the {\it simplicial nerve functor\/} $\sN_{A_\iy}:\ACat_R^\un\ab\ra\Cat^1_\iy$, where $\Cat^1_\iy\subset\SSet$ is the ordinary category of $\iy$-categories. It is only defined for {\it unital\/} $A_\iy$-categories $\sA$. There is a strict isomorphism of categories $H^0(\sA)\cong\Ho(\sN_{A_\iy}(\sA))$, where $\Ho(\sC)$ is the homotopy category of an $\iy$-category $\sC$. That is, $\sN_{A_\iy}(\sA)$ is an $\iy$-categorical enhancement of the cohomology category~$H^0(\sA)$.
\end{dfn}

\begin{dfn}
\label{ai1def10}
In \S\ref{ai49} we define when an $A_\iy$-category $\sA$ is {\it pretriangulated}. If $\sA$ is also (cohomologically) unital then the cohomology category $H^0(\sA)$ is a triangulated category. We write $\ADCat_R\subset\ACat_R$, $\ADCat^\un_R\subset\ACat^\un_R$, $\ADCat^\cu_R\subset\ACat^\cu_R$ for the full subcategories of pretriangulated objects.

We define the {\it twist functor\/} $\Tw:\ACat_R\ra\ACat_R$, which maps $\ACat^\un_R\ab\ra\ACat^\un_R$ and $\ACat^\cu_R\ra\ACat^\cu_R$. The essential images of these functors are $\ADCat_R,\ADCat^\un_R,\ADCat^\cu_R$. For any $\sA$ there is an embedding $\sA\hookra\Tw\sA$ as a full $A_\iy$-subcategory. We can think of $\Tw$ as being like the derived category functor $D^b:\cA\mapsto D^b\cA$ which maps an abelian category $\cA$ to its bounded derived category~$D^b\cA$.

If $\sA$ is unital and pretriangulated then the simplicial nerve $\sN_{A_\iy}(\sA)$ is a {\it stable\/ $\iy$-category}. These are the $\iy$-categorical generalization of triangulated categories, discussed in~\S\ref{ai27}.
\end{dfn}

\subsection{\texorpdfstring{Curved and filtered $A_\iy$-algebras and $A_\iy$-categories}{Curved and filtered A∞-algebras and A∞-categories}}
\label{ai1.2}

Next we very briefly describe the material of Chapter~\ref{ai5}.

\begin{dfn}
\label{ai1def11}
Let $R$ be a commutative ring. A {\it curved\/ $A_\iy$-category}, or {\it weak\/ $A_\iy$-category}, $\sA$ over $R$ consists of:
\begin{itemize}
\setlength{\itemsep}{0pt}
\setlength{\parsep}{0pt}
\item[(a)] A set, or proper class, $\Obj\sA$ of {\it objects\/} of $\sA$.
\item[(b)] For all objects $x,y\in\Obj\sA$, a $\Z$-graded $R$-module
\begin{equation*}
\Hom(x,y)=\bigop\nolimits_{d\in\Z}\Hom^d(x,y).
\end{equation*}
\item[(c)] For all $k\ge 0$ and objects $x_0,x_1,\ldots,x_k$ in $\Obj\sA$, a graded $R$-linear map
\begin{align*}
\m_k^{x_0,\ldots,x_k}&:\Hom(x_0,x_1)\ot_R\Hom(x_1,x_2)\ot_R\cdots \ot_R\Hom(x_{k-1},x_k)\\
&\quad \longra\Hom(x_0,x_k)
\end{align*}
of degree $2-k$. For all $n\ge 0$ and objects $x_0,\ldots,x_n$ these must satisfy
\end{itemize}
\ea
\sum_{\begin{subarray}{l} i,j,k\ge 0:\\ i+j+k=n\end{subarray}}\begin{aligned}[t]&
(-1)^{i+jk}\m_{i+k+1}^{x_0,\ldots,x_i,x_{i+j},\ldots,x_n}\ci\bigl(\id_{\Hom(x_0,x_1)}\ot\cdots\ot\id_{\Hom(x_{i-1},x_i)}\ot {}\\
&\m_j^{x_i,\ldots,x_{i+j}}\ot \id_{\Hom(x_{i+j},x_{i+j+1})}\ot\cdots\ot\id_{\Hom(x_{n-1},x_n)}\bigr)=0.
\end{aligned}
\nonumber\\[-17pt]
\label{ai1eq9}
\ea

This is the same as an $A_\iy$-category in Definition \ref{ai1def1}, except that we have maps $\m_k^{x_0,\ldots,x_k}$ for $k\ge 0$ rather than $k\ge 1$, and in \eq{ai1eq9} we allow $j\ge 0$ rather than $j\ge 1$ in \eq{ai1eq1}. We define {\it unital\/} curved $A_\iy$-categories as in Definition~\ref{ai1def1}.

A {\it curved\/ $A_\iy$-algebra\/} is equivalent to a curved $A_\iy$-category with one object $*$, as in Definition~\ref{ai1def6}.
\end{dfn}

Curved $A_\iy$-algebras and $A_\iy$-categories occur naturally in many places in Mathematics and Physics. The obvious way to define {\it curved\/ $A_\iy$-functors\/} of curved $A_\iy$-categories is to follow Definition \ref{ai1def3}, but defining $\fF_k$ for $k\ge 0$ rather than $k\ge 1$. However, we immediately run into a problem: equation \eq{ai1eq5} then involves an {\it infinite\/} sum, as on the right hand side we should sum over decompositions $n=n_1+\cdots+n_k$ with $k\ge 0$ and $n_i\ge 0$, rather than $n_i\ge 1$. The same applies to compositions in \eq{ai1eq7}, and the definitions of $A_\iy$-homotopies and $A_\iy$-natural transformations.

Because of this, it is natural to study curved $A_\iy$-algebras and curved $A_\iy$-categories in settings in which we have notions of {\it convergence of infinite sums}. Chapters \ref{ai5} and \ref{ai8} develop a theory of {\it filtered\/ $A_\iy$-algebras\/} and {\it filtered\/ $A_\iy$-categories}, in which the $R$-modules $\Hom^d(x,y)$ are equipped with complete filtrations over a partially ordered monoid $(\Ga,\pre)$, and convergent infinite sums in $\Hom^d(x,y)$ are well defined. Our definitions are new, but are based closely on work by Fukaya--Oh--Ohta--Ono \cite{Fuka2,Fuka3,Fuka4,Fuka5,FOOO} and De Deken and Lowen~\cite{DDLo}.

As in \cite{Fuka2,Fuka3,Fuka4,Fuka5,FOOO}, (derived) Fukaya categories of general symplectic manifolds should be constructed using filtered $A_\iy$-categories over a Novikov ring.

Cohomology categories $H^0(\sA)$ are not defined for curved or filtered $A_\iy$-categories. So we do not have a well-behaved notion of `cohomologically unital' in the filtered case. We write $\AfCat^\un_R\subset\AfCat_R$ for the categories of (unital) filtered $A_\iy$-categories over a complete filtered commutative ring~$R$.

Define subcategories $\AfCat_{R,\iy}\subset\AfCat_{R,1}\subset\AfCat_{R,0}=\AfCat_R$, such that $\AfCat_{R,1}$ has objects {\it strict\/} filtered $A_\iy$-categories $\sA$ over $R$, with $\m_0^x=0$ for all $x\in\Obj\sA$ and morphisms {\it all\/} filtered $A_\iy$-functors $\fF:\sA\ra\sB$, including those with $\fF^x_0\ne 0$, and $\AfCat_{R,\iy}$ has objects {\it strict\/} filtered $A_\iy$-categories $\sA$, with $\m^x_0=0$ for all $x\in\Obj\sA$, and morphisms {\it strict\/} filtered $A_\iy$-functors $\fF:\sA\ra\sB$, with $\fF^x_0=0$ for all $x\in\Obj\sA$. Similarly we define $\AfCat_{R,\iy}^\un\subset\AfCat_{R,1}^\un\subset\AfCat_{R,0}^\un=\AfCat^\un_R$. Now objects and morphisms of $\AfCat_{R,\iy},\AfCat_{R,\iy}^\un$ are just $A_\iy$-categories and $A_\iy$-functors, with extra filtration data. Thus forgetting filtrations gives functors
\e
\begin{split}
\Pi_{A_\iy^\rfi}^{A_\iy}&:\AfCat_{R,\iy}\longra\ACat_R,\\
\Pi_{A_\iy^\rfi}^{A_\iy}&:\AfCat^\un_{R,\iy}\longra\ACat^\un_R.
\end{split}
\label{ai1eq10}
\e

Functor categories $\Fun_{A^\rfi_\iy}(\sA,\sB),\Fun_{A^\rfi_\iy}^\un(\sA,\sB)$ can be defined in the filtered case, and surprisingly, they are $A_\iy$-categories rather than curved $A_\iy$-categories. So the definitions of $A_\iy$-natural transformation, $A_\iy$-natural isomorphism and homotopy equivalence also work in the filtered case. There are functors
\e
\begin{split}
\Str&:\AfCat_R\longra\AfCat_{R,\iy},\\
\Str&:\AfCat_R^\un\longra\AfCat_{R,\iy}^\un.
\end{split}
\label{ai1eq11}
\e
which convert a filtered $A_\iy$-category $\sA$ into a strict filtered $A_\iy$-category $\Str(\sA)$, called the {\it associated strict\/ $A_\iy$-category}. Objects of $\Str(\sA)$ are pairs $(x,b)$ of an object $x$ in $\sA$ and a `bounding cochain' $b$ for $x$.

Composing with $\Pi_{A_\iy^\rfi}^{A_\iy}$ gives functors
\e
\begin{split}
\Pi_{A_\iy^\rfi}^{A_\iy}\ci\Str&:\AfCat_R\longra\ACat_R,\\
\Pi_{A_\iy^\rfi}^{A_\iy}\ci\Str&:\AfCat^\un_R\longra\ACat^\un_R.
\end{split}
\label{ai1eq12}
\e
These are an important part of the programme of Fukaya--Oh--Ohta--Ono \cite{Fuka1,Fuka2,Fuka3,Fuka4,Fuka5,FOOO} for defining Fukaya categories of symplectic manifolds $(X,\om)$. First they define a filtered $A_\iy$-category $\sF(X,\om)$ over the Novikov ring $\La^\nov_\Q$. Then they apply the functors \eq{ai1eq12} to make it into an $A_\iy$-category.

\subsection{\texorpdfstring{Our $\iy$-categories of $A_\iy$-algebras and $A_\iy$-categories}{Our ∞-categories of A∞-algebras and A∞-categories}}
\label{ai1.3}

Chapter \ref{ai7}, the central chapter in the book, defines many different $\iy$-categories of $A_\iy$-categories over a commutative ring $R$, which fit into a diagram \eq{ai7eq1} of inclusions of $\iy$-categories:
\ea
\label{ai1eq13}\\[-19pt]
\xymatrix@!0@C=85pt@R=37pt{{\begin{subarray}{l}\ts \ACat^{\un,\iy}_R\!\!=  \\[2pt] \ts \ACat^{\un,\iy}_{R,2}\end{subarray}} \\
\ACat^{\un,\iy^{\vphantom{(}}}_{R,3}\, \ar@{^{(}->}[u] \ar@{=}[r] \ar@{^{(}->}[u] & \dACat^{\un,\iy^{\vphantom{(}}}_{R,3}\, \ar@{^{(}->}[r] & {\begin{subarray}{l}\ts \ACat^{\cu,\iy}_R\!\!=  \\[2pt] \ts \ACat^{\cu,\iy}_{R,3}\end{subarray}}\,  \ar@{^{(}->}[r] & {\begin{subarray}{l}\ts \ACat^\iy_R\!=  \\[2pt] \ts \ACat^\iy_{R,3}\end{subarray}}  \\
\ACat^{\un,\iy^{\vphantom{(}}}_{R,4}\, \ar@{=}[r] \ar@{^{(}->}[u] & \dACat^{\un,\iy^{\vphantom{(}}}_{R,4}\, \ar@{^{(}->}[r] \ar@{^{(}->}[u] & \ACat^{\cu,\iy^{\vphantom{(}}}_{R,4}\,  \ar@{^{(}->}[r] \ar@{^{(}->}[u] & \ACat^{\iy^{\vphantom{(}}}_{R,4} \ar@{^{(}->}[u] \\
\ACat^{\un,\iy^{\vphantom{(}}}_{R,5}\, \ar@{=}[r] \ar@{^{(}->}[u] & \dACat^{\un,\iy^{\vphantom{(}}}_{R,5}\, \ar@{^{(}->}[r] \ar@{^{(}->}[u] & \ACat^{\cu,\iy^{\vphantom{(}}}_{R,5}\,  \ar@{^{(}->}[r] \ar@{^{(}->}[u] & \ACat^{\iy^{\vphantom{(}}}_{R,5} \ar@{^{(}->}[u]  \\
\vdots_{\vphantom{(}} \ar@{^{(}->}[u] & \vdots_{\vphantom{(}} \ar@{^{(}->}[u] & \vdots_{\vphantom{(}} \ar@{^{(}->}[u] & \vdots_{\vphantom{(}} \ar@{^{(}->}[u] \\
\ACat^{\un,\iy^{\vphantom{(}}}_{R,\iy}\, \ar@{=}[r] \ar@{^{(}->}[u] & \dACat^{\un,\iy^{\vphantom{(}}}_{R,\iy}\, \ar@{^{(}->}[r] \ar@{^{(}->}[u] & \ACat^{\cu,\iy^{\vphantom{(}}}_{R,\iy}\,  \ar@{^{(}->}[r] \ar@{^{(}->}[u] & \ACat^{\iy^{\vphantom{(}}}_{R,\iy} \ar@{^{(}->}[u] \\
\sN(\ACat^\un_R)\, \ar@{=}[r] \ar@{^{(}->}[u] & \sN(\ACat^\un_R)\, \ar@{^{(}->}[r] \ar@{^{(}->}[u] & \sN(\ACat^\cu_R)\,  \ar@{^{(}->}[r] \ar@{^{(}->}[u] & \sN(\ACat_R).\! \ar@{^{(}->}[u] }
\nonumber
\ea

Chapter \ref{ai6} first constructs the analogous diagram \eq{ai6eq1} of $\iy$-categories of $A_\iy$-algebras, obtained from \eq{ai1eq13} by restricting to $A_\iy$-categories with one object $*$. Here in~\eq{ai1eq13}:
\begin{itemize}
\setlength{\itemsep}{0pt}
\setlength{\parsep}{0pt}
\item[(i)] The  $\iy$-categories in the first two columns, with superscripts `$\un$', have objects unital $A_\iy$-categories and 1-morphisms unital $A_\iy$-functors.
\item[(ii)] The $\iy$-categories in the third column, with superscript `$\cu$', have objects cohomologically unital $A_\iy$-categories and 1-morphisms cohomologically unital $A_\iy$-functors.
\item[(iii)] The $\iy$-categories in the fourth column have objects all $A_\iy$-categories, and 1-morphisms all $A_\iy$-functors.
\item[(iv)] The $\iy$-categories in the first and second columns, linked by `$=$', are strictly isomorphic, but with definitions written in a different way.
\item[(v)] In each row we have subscripts `$R,K$' for $K=2,3,4,\ldots,\iy$, where we call $K$ the {\it type\/} of the $\iy$-category.
\item[(vi)] $\sN(\ACat^\un_R),\ldots,\sN(\ACat_R)$ are the nerves of the ordinary categories $\ACat^\un_R,\ldots,\ACat_R$ from \S\ref{ai1.1}.
\end{itemize}

The 2-morphisms in $\ACat^{\un,\iy}_R$ of type 2 in the first row of \eq{ai1eq13} correspond to {\it unital\/ $A_\iy$-natural isomorphisms\/} in Definition \ref{ai1def8}. But in the $\iy$-categories of type $K\ge 3$ in the subsequent rows, the 2-morphisms correspond to ({\it unital\/}) $A_\iy$-{\it homotopies\/} of $A_\iy$-functors, as in Definition \ref{ai1def8}. In the author's view, $\ACat^{\un,\iy}_R$ will be the best $\iy$-category of $A_\iy$-categories for most applications, as $A_\iy$-natural isomorphisms are the best notion of when two $A_\iy$-functors are `the same'. 

The first column of \eq{ai1eq13} can only be defined for {\it unital\/} $A_\iy$-categories, as the degeneracy maps in the simplicial sets involve identities $\boo_x$. The last three columns are defined by a different method which does not need identities, but requires a more restrictive notion of 2-morphism.

To give a flavour of the definitions of the $\iy$-categories in \eq{ai1eq13}, we give a partial definition of $\ACat^{\un,\iy}_R$. For $N\ge 0$, an $N$-simplex in $\ACat^{\un,\iy}_R$ is a quadruple $\bigl((\sO^i)_{i=0}^N,(f^{ij})_{0\le i\le j\le N},(\sH^i)_{i=0}^N,(\m^I)_{\es\ne I\subseteq\{0,\ldots,N\}}\bigr)$, where:
\begin{itemize}
\setlength{\itemsep}{0pt}
\setlength{\parsep}{0pt}
\item[(a)] $\sO^i$ is a set for $i=0,\ldots,N$, thought of as the set of objects $\Obj\sA^i$ of an $A_\iy$-category $\sA^i$ at vertex $i$ of the $N$-simplex $\De^N$.
\item[(b)] $f^{ij}:\sO^i\ra\sO^j$ is a map for all $0\le i\le j\le N$, thought of as the map $\fF_\Obj^{ij}:\Obj\sA^i\ra\Obj\sA^j$ in an $A_\iy$-functor $\fF^{ij}:\sA^i\ra\sA^j$ living on the edge $[i,j]$ of $\De^N$, with $f^{ii}=\id_{\sO^i}$ for $i=0,\ldots,N$.

Suppose $\es\ne I\subseteq\{0,1,\ldots,N\}$ with $I=\{i_0,i_1,\ldots,i_M\}$ for $i_0<i_1<\cdots<i_M$. Define $f^I=f^{i_{M-1}i_M}\ci f^{i_{M-2}i_{M-1}}\ci\cdots\ci f^{i_0i_1}:\sO^{i_0}\ra\sO^{i_M}$.
\item[(c)] For all $i=0,1,\ldots,N$ and $x,y\in\sO^i$, we are given a $\Z$-graded $R$-module $\sH^i(x,y)=\bigop_{d\in\Z}\sH^i(x,y)^d$, thought of as the Hom space $\Hom_{\sA^i}(x,y)$ in the $A_\iy$-category $\sA^i$. We write $\sH^i=(\sH^i(x,y))_{x,y\in\sO^i}$.
\item[(d)] For each nonempty subset $I\subseteq\{0,\ldots,N\}$, which we write as $I=\{i_0,\ab i_1,\ab\ldots,\ab i_M\}$ for $i_0<i_1<\cdots<i_M$ with $M\ge 0$, and for all $k\ge 0$ and elements $x_0,\ldots,x_k\in \sO^{i_0}$, thought of as objects in $\sA^{i_0}$, we are given an $R$-linear map 
\begin{align*}
\m_k^{I;x_0,\ldots,x_k}:\sH^{i_0}&(x_0,x_1)\ot_R\sH^{i_0}(x_1,x_2)\ot_R\cdots \ot_R\sH^{i_0}(x_{k-1},x_k)
\nonumber\\
&\longra\sH^{i_M}(f^I(x_0),f^{i_0i_M}(x_k))
\end{align*}
graded of degree $2-M-k=3-\md{I}-k$. When $k=0$ we require that $\m^{I;x_0}_0=0$ when $\md{I}=1,2$. 

We write $\m^I_k=(\m_k^{I;x_0,\ldots,x_k})_{x_0,\ldots,x_k\in\sO^{i_0}}$ and $\m^I=(\m^I_k)_{k\ge 0}$.

We require these maps $\m_k^{I;x_0,\ldots,x_k}$ to satisfy the following identities for all $\es\ne I=\{i_0<\cdots<i_M\}\subseteq\{0,1,\ldots,N\}$, $n\ge 0$, and $x_0,\ldots,x_n\in \sO^{i_0}$, firstly if $M=0$
\e
0=\sum_{\begin{subarray}{l} k=1,\ldots,n, \\
y_0,\ldots,y_k\in\sO^{i_0} \end{subarray}}\m_k^{\{i_0\};y_0,\ldots,y_k}\ci\d_{n,k;y_0,\ldots,y_k}^{\{i_0\};x_0,\ldots,x_n},
\label{ai1eq14}
\e
and secondly if $M>0$
\ea
&0=\sum_{\begin{subarray}{l} k=1,\ldots,n \\
y_0,\ldots,y_k\in\sO^{i_0}  \end{subarray}}\m_k^{I;y_0,\ldots,y_k}\ci\d_{n,k;y_0,\ldots,y_k}^{\{i_0\};x_0,\ldots,x_n}
\label{ai1eq15}\\
&+(-1)^M\sum_{\begin{subarray}{l} k=1,\ldots,n+M-1,\\
y_0,\ldots,y_k\in\sO^{i_M} \end{subarray}}\m_k^{\{i_M\};y_0,\ldots,y_k}\ci\d_{n,k;y_0,\ldots,y_k}^{I;x_0,\ldots,x_n}.
\nonumber
\ea
Here $\d_{n,k;y_0,\ldots,y_k}^{I;x_0,\ldots,x_n}$ is a graded $R$-linear map with a complicated definition in \S\ref{ai61}--\S\ref{ai62} and \S\ref{ai71}, as a sum over labelled trees of multiple compositions of the $\m_l^{J;w_0,\ldots,w_l}$. In both \eq{ai1eq14} and \eq{ai1eq15} there are only finitely many nonzero terms, and all nonzero terms map to $\sH^{i_M}(f^I(x_0),f^{i_0i_M}(x_n))$, so \eq{ai1eq14}--\eq{ai1eq15} should be regarded as equations in this space.
\item[(e)] There should exist elements $\boo^i_x\in\sH^i(x,x)^0$ for each $i=0,\ldots,N$ and $x\in\sO^i$ satisfying analogues of \eq{ai1eq3}--\eq{ai1eq4}. These make $\sA^i$ into a {\it unital\/} $A_\iy$-category.
\item[(f)] Parts (a)--(e) imply that using the data $\bigl((\sO^i),(f^{ij}),(\sH^i),(\m^I)\bigr)$, for all $i=0,\ldots,N$ we can define a unital $A_\iy$-category $\sA^i$, for all $0\le i<j\le N$ we can define a unital $A_\iy$-functor $\fF^{ij}:\sA^i\ra\sA^j$, and for all $0\le i<j<k\le N$, using the data $\m_k^{\{i,j,k\};x_0,\ldots,x_k}$ we can define a unital $A_\iy$-natural transformation $\fH^{ijk}:\fF^{jk}\ci\fF^{ij}\Ra\fF^{ik}$. We require this to be a unital $A_\iy$-natural isomorphism.
\end{itemize}
For the $\iy$-subcategories of `type $K$' in \eq{ai1eq13}, we require that $\m_0^I=0$ if~$\md{I}\le K$.

Actually proving that the simplicial sets in \eq{ai1eq13} are $\iy$-categories is very complicated, and among the main results of this book.

We write $\ACat^{\un,\rsu,\iy}_R\subset\ACat^{\un,\iy}_R$ for the full $\iy$-subcategory whose objects have {\it split units}.

We write $\ACat^{\un,\hs,\iy}_R\subset\ACat^{\un,\iy}_R$ for the full $\iy$-subcategory whose objects are {\it h-projective with split units}.

We write $\ADCat^{\un,\iy}_R\subset\ACat^{\un,\iy}_R$ for the full $\iy$-subcategory whose objects are {\it pretriangulated}.

\begin{rem}
\label{ai1rem1}
{\bf(a)} The heart of the definition of $\ACat^{\un,\iy}_R$ is equations \eq{ai1eq14}--\eq{ai1eq15}, which generalize the $A_\iy$-equation \eq{ai1eq1} to $N$-simplices. Their apparent simplicity is deceptive, as the definition of the operators $\d_{n,k;y_0,\ldots,y_k}^{I;x_0,\ldots,x_n}$, which takes up \S\ref{ai61}, \S\ref{ai62} and \S\ref{ai71}, is extraordinarily complicated. The simplest case when $M=0$ is
\begin{align*}
&\d_{n,k;y_0,\ldots,y_k}^{\{i_0\};x_0,\ldots,x_n}=\\
&\sum_{\begin{subarray}{l} 0\le j<k, \\
y_i=x_i,\; 0\le i\le j, \\
y_i=x_{n-k+i},\; j<i\le k
\end{subarray}}\begin{aligned}[t]&
(-1)^{j+(n-k+1)(k-j-1)}\id_{\sH^{i_0}(x_0,x_1)}\ot\cdots\ot\\
&\id_{\sH^{i_0}(x_{j-1},x_j)}\ot\m_{n-k+1}^{\{i_0\};x_j,x_{j+1},\ldots,x_{j+n-k+1}}\ot\\
&\id_{\sH^{i_0}(x_{j+n-k+1},x_{j+n-k+2})}\ot\cdots\ot\id_{\sH^{i_0}(x_{n-1},x_n)}.
\end{aligned}
\end{align*}
Substituting this into \eq{ai1eq14} shows that the $\m_k^{\{i_0\};y_0,\ldots,y_k}$ satisfy \eq{ai1eq1}, consistent with the $N$-simplex above including the data of an $A_\iy$-category $\sA^i$ at each vertex $i=0,\ldots,N$.

Equations \eq{ai1eq14}--\eq{ai1eq15} have some very special properties, which are not obvious. From the full definition of $\ACat^{\un,\iy}_R$ in \S\ref{ai72}, which is only sketched above, it is fairly straightforward to show that $\ACat^{\un,\iy}_R$ is a simplicial set. However, the crucial proof in Theorem \ref{ai7thm3} that $\ACat^{\un,\iy}_R$ is an $\iy$-category is very long and difficult, and relies on deep properties of \eq{ai1eq14}--\eq{ai1eq15}.
\smallskip

\noindent{\bf(b)} Equations \eq{ai1eq14}--\eq{ai1eq15} are used to define the $\iy$-categories in the first column of \eq{ai1eq13}. For the last three columns, they are replaced by
\ea
&\sum_{\begin{subarray}{l} 0\le j\le M, \; 0\le k\le\max(n,n+j-1), \\ y_0,\ldots,y_k\in\sO^{i_j}\end{subarray}\!\!\!\!\!\!\!\!\!\!\!\!\!\!\!\!\!\!\!\!\!\!\!\!}
(-1)^j\,\m_k^{\{i_j,\ldots,i_M\};y_0,\ldots,y_k}\ci\d_{n,k;y_0,\ldots,y_k}^{\{i_0,\ldots,i_j\};x_0,\ldots,x_n}
\nonumber\\
&=\sum_{j=1,\ldots,M-1}(-1)^j\m_n^{I\sm\{i_j\};x_0,\ldots,x_n}.
\label{ai1eq16}
\ea
This only makes sense if the maps $f^{ij}:\sO^i\ra\sO^j$ in (b) above satisfy the extra condition $f^{jk}\ci f^{ij}=f^{ik}$ for all $0\le i<j<k\le N$, as otherwise different terms in \eq{ai1eq16} map to different spaces. This condition $f^{jk}\ci f^{ij}=f^{ik}$ is imposed in $\iy$-categories $\ACat_{R,K}^{*,\iy}$ in \eq{ai1eq13} of type $K$ for $K\ge 3$, but not $K=2$. This is why the last three columns in \eq{ai1eq13} have no entry with~$K=2$.

Equation \eq{ai1eq16} also has non-obvious special properties needed to show that the last three columns of \eq{ai1eq13} are $\iy$-categories.

The definitions of $\iy$-categories using \eq{ai1eq14}--\eq{ai1eq15} requires the $A_\iy$-cat\-eg\-or\-ies $\sA^i$ to be (strictly) unital, as the identities $\boo^i_x$ are used in the degeneracy maps in the simplicial set. The definitions using \eq{ai1eq16} work for non-unital $A_\iy$-categories, but require $f^{jk}\ci f^{ij}=f^{ik}$. The effect of this is that 2-morphisms in the last three columns of \eq{ai1eq13} are $A_\iy$-homotopies, not unital $A_\iy$-natural isomorphisms as in (f) above.
\smallskip

\noindent{\bf(c)} It would have been very difficult to invent equations \eq{ai1eq14}--\eq{ai1eq15} or \eq{ai1eq16} from scratch. The author derived \eq{ai1eq16} from a construction in Symplectic Geometry, as part of a project with Guillem Cazassus and Alexander Ritter \cite{CJR}. Very roughly, from a (for simplicity, say exact) symplectic manifold $(X,\om)$ with a choice of compatible almost complex structure $J$, one can construct a Fukaya $A_\iy$-category $\sF(X,\om)_J$. The author generalized this to define a categorical structure associated to a smooth family of almost complex structures $J_\si$ parametrized by points $\si$ in the topological $N$-simplex $\De^N_\top$. This was an example of an $N$-simplex in $\ACat^{\cu,\iy}_R$, and led to its definition in \S\ref{ai73}. Later, the author transformed \eq{ai1eq16} to \eq{ai1eq14}--\eq{ai1eq15}, in order to define $\iy$-categories of $A_\iy$-categories in which 2-morphisms were $A_\iy$-natural transformations.
\end{rem}

For {\it filtered\/} $A_\iy$-categories over a complete $\Ga$-filtered commutative ring $R$, as in \S\ref{ai1.2}, Chapter \ref{ai8} defines a diagram \eq{ai8eq1} of $\iy$-categories of filtered $A_\iy$-categories over $R$ similar to \eq{ai1eq13}:
\ea
\begin{gathered}
\!\!\!\!\xymatrix@C=8pt@R=12pt{
{\AfCat^{\un,\iy}_R=\AfCat^{\un,\iy}_{R,0,2}\hskip -60pt} \\
\AfCat^{\un,\iy^{\vphantom{(}}}_{R,0,3} \ar@{=}[rr] \ar@{^{(}->}[u] && \dAfCat^{\un,\iy^{\vphantom{(}}}_{R,0,3}\,\, \ar@{^{(}->}[rr]  && {\begin{subarray}{l}\ts \AfCat^\iy_R=  \\[2pt] \ts \AfCat^\iy_{R,0,3}\end{subarray}} \\
\AfCat^{\un,\iy^{\vphantom{(}}}_{R,K',L'} \ar@{=}[rr] \ar@{^{(}->}[u] && \dAfCat^{\un,\iy^{\vphantom{(}}}_{R,K',L'}\,\, \ar@{^{(}->}[rr] \ar@{^{(}->}[u] &&  \AfCat^{\iy^{\vphantom{(}}}_{R,K',L'} \ar@{^{(}->}[u] & *+[l]{\;\begin{subarray}{l}(0\le K'\le L', \\ L'\ge 3)\end{subarray}} \\
\AfCat^{\un,\iy^{\vphantom{(}}}_{R,K,L} \ar@{=}[rr] \ar@<-4ex>@{^{(}->}[u]^{K\ge K',\; L\ge L'} && \dAfCat^{\un,\iy^{\vphantom{(}}}_{R,K,L}\,\, \ar@{^{(}->}[rr] \ar@<-4ex>@{^{(}->}[u]^{K\ge K',\; L\ge L'}  && \AfCat^{\iy^{\vphantom{(}}}_{R,K,L} \ar@<-4ex>@{^{(}->}[u]^{K\ge K',\; L\ge L'} & *+[l]{\;\begin{subarray}{l}(0\le K\le L, \\ L\ge 3)\end{subarray}} \\
\AfCat^{\un,\iy^{\vphantom{(}}}_{R,K,\iy} \ar@{=}[rr]  && \dAfCat^{\un,\iy^{\vphantom{(}}}_{R,K,\iy}\,\, \ar@{^{(}->}[rr]  &&  \AfCat^{\iy^{\vphantom{(}}}_{R,K,\iy} & \ar@<3ex>@{}[d]_{(K=0,1,\iy.)} \\
\sN(\AfCat^{\un^{\vphantom{(}}}_{R,K}) \ar@{^{(}->}[u] \ar@{=}[rr] && \sN(\AfCat^{\un^{\vphantom{(}}}_{R,K})\,\, \ar@{^{(}->}[rr] \ar@{^{(}->}[u] &&  \sN(\AfCat^{\phantom\un^{\vphantom{(}}}_{R,K}) \ar@{^{(}->}[u] & }\!\!\!\!\!\!\!\end{gathered}
\label{ai1eq17}
\ea

Here the $\iy$-categories $\AfCat^{*,\iy}_{R,K,L}$ are of {\it type\/} $(K,L)$ for $0\le K\le L\le\iy$. For $N$-simplices $\bigl((\sO^i),(f^{ij}),(\sH^i),(\m^I)\bigr)$ in $\AfCat^{*,\iy}_{R,K,L}$, when $\md{I}\ne 3$, the rule is that $\m_0^{I;x_0}=0$ if $\md{I}\le K$, and $\m_0^{I;x_0}(1)$ should have strictly positive grading if~$K<\md{I}\le L$.

We also define a diagram of $\iy$-categories \eq{ai8eq2} identical to \eq{ai1eq17}, but for $\iy$-categories $\AwCat^{\un,\iy}_R,\AwCat^{\un,\iy}_{R,K,L},\ldots$ that we call {\it weakly filtered}. These have the same objects and 1-morphisms as $\AfCat^{\un,\iy}_R,\AfCat^{\un,\iy}_{R,K,L},\ldots,$ but have a weaker notion of $N$-morphism for $N\ge 2$: in $\AfCat^{\un,\iy}_R,\ldots$ the maps $\m_k^{I;x_0,\ldots,x_k}$ are required to preserve filtrations, but in $\AwCat^{\un,\iy}_R,\ldots$ the $\m_k^{I;x_0,\ldots,x_k}$ are allowed to decrease filtrations by a bounded amount.

We will compare our $\iy$-categories with $\iy$-categories of $A_\iy$-categories defined by localization. See \S\ref{ai26} for more details on the next definition.

\begin{dfn}
\label{ai1def12}
A {\it category with weak equivalences\/} $(\cC,\cW)$ is a category $\cC$ and a subset $\cW\subseteq\Mor(\cC)$ of the morphisms of $\cC$, such that $\cW$ contains all isomorphisms and satisfies the `two out of three' property. A {\it Dwyer--Kan localization of\/ $\cC$ along\/} $\cW$ is an $\iy$-category $\cC[\cW^{-1}]^\iy$ and an $\iy$-functor $\la:\sN(\cC)\ra\cC[\cW^{-1}]^\iy$, such that if $f\in\cW$ then $\la(f)$ is an equivalence in $\cC[\cW^{-1}]^\iy$, and $\la$ has the universal property that whenever $\fF:\sN(\cC)\ra\sD$ is an $\iy$-functor which maps $\cW$ to equivalences in $\sD$, then we can factor $\fF$ as 
\begin{equation*}
\xymatrix@C=50pt@R=15pt{
*+[r]{\sN(\cC)} \ar[rr]_(0.4){\Pi_\loc^\iy} \ar@/_.5pc/[drr]_\fF &  \drtwocell_{}\omit^{}\omit{^{}} & *+[l]{\cC[\cW^{-1}]^\iy} \ar@{..>}[d]^\fG \\
&& *+[l]{\sD } }
\end{equation*}
for some $\fG:\cC[\cW^{-1}]^\iy\ra\sD$ unique up to a contractible choice. Dwyer--Kan localizations $\cC[\cW^{-1}]^\iy$ always exist, and are unique up to Joyal equivalence. We can take $\cC[\cW^{-1}]^\iy$ to have the same objects as~$\cC$.
\end{dfn}

Starting from an ordinary category of $A_\iy$-categories such as $\ACat^\un_R$ or $\ACat^\cu_R$, the obvious notions of weak equivalence are {\it homotopy equivalences\/} or {\it quasi-equivalences\/} from Definition \ref{ai1def8}. If $R$ is a field, and if we restrict to $A_\iy$-categories without zero objects in the unital case, then these are the same, but in general $\{$homotopy equivalences$\}\subsetneq\{$quasi-equivalences$\}$. We will write $\HE,\QE$ to denote the classes of homotopy equivalences and quasi-equivalences in our category, with subscripts $\un,\cu,\ldots$ to indicate the category. Thus we could for example consider the $\iy$-categories
\e
\begin{aligned}
&\ACat^\un_R[\HE_\un^{-1}]^\iy, &
&\ACat^\un_R[\QE_\un^{-1}]^\iy, \\
&\ACat^\cu_R[\QE_\cu^{-1}]^\iy, &
&\DGCat_R[\QE_\dg^{-1}]^\iy.
\end{aligned}
\label{ai1eq18}
\e

We do not consider localizations of $\ACat_R$, as homotopy equivalences and quasi-equivalences are not well behaved for non-unital $A_\iy$-categories. The previous literature on $\iy$-categories of $A_\iy$-categories \cite{COS1,COS2,Faon2,OhTa,Orna1,Tana,Pasc}, discussed in \S\ref{ai410}, generally studies $\ACat^\un_R[\QE_\un^{-1}]^\iy$ and~$\ACat^\cu_R[\QE_\cu^{-1}]^\iy$.

Some important differences between our $\iy$-categories \eq{ai1eq13}, and localized $\iy$-categories such as \eq{ai1eq18}, are that our $\iy$-categories are completely explicit, and it is a lot of work to define them and prove they are $\iy$-categories. For localized $\iy$-categories, their definition and existence is trivial (given the existence of Dwyer--Kan localizations), but they are non-explicit, and the goal is to prove something useful about them.

\subsection{\texorpdfstring{$\iy$-bicategories of $A_\iy$-categories}{∞-bicategories of A∞-categories}}
\label{ai1.4}

Informally, an $(\iy,n)$-{\it category\/} is some kind of higher category with a notion of $k$-morphism for all $k=1,2,\ldots,$ such that all $k$-morphisms are invertible for $k>n$. Our $\iy$-categories are a kind of $(\iy,1)$-category. In \S\ref{ai28} we explain $\iy$-{\it bicategories}, following Lurie \cite[\S 4]{Luri4}, which are a model for $(\iy,2)$-categories. 

An $\iy$-{\it bicategory\/} $(S,T)$ is a simplicial set $S$ with a distinguished subset $T\subseteq S([2])$ of 2-simplices in $S$, called {\it thin\/} 2-simplices. These $S,T$ satisfy more complicated horn-filling conditions than $\iy$-categories do. Roughly, we think of a 2-simplex $t:\De^2\ra S$ as a diagram
\begin{equation*}
\xymatrix@C=70pt@R=15pt{
*+[l]{X\;\overset{0}{\bu}} \ar@<.8ex>[d]_f \ar@<.6ex>[dr]^h \ar@{}@<0ex>[dr]_(0.35){\eta\, \Nearrow} \\
*+[l]{Y\;\underset{1}{\bu}} \ar@<.7ex>[r]_g & *+[r]{\underset{2}{\bu}\; Z,} }
\end{equation*}
where $f:X\ra Y$, $g:Y\ra Z$ and $h:X\ra Z$ are 1-morphisms in $S$ and $\eta:g\ci f\Ra h$ is a 2-morphism. Then $t$ is thin if $\eta$ is a 2-isomorphism. We write $\Core(S,T)\subset S$ for the simplicial subset of $S$ consisting of $N$-simplices all of whose 2-simplices are thin. Then $\Core(S,T)$ is an $\iy$-category.

In \S\ref{ai72} we define our $\iy$-category $\ACat^{\un,\iy}_R$ of unital $A_\iy$-categories over $R$ by first defining a simplicial set $\hACat^{\un,\iy}_R$ of unital $A_\iy$-categories, and then defining a simplicial subset $\ACat^{\un,\iy}_R\subset\hACat^{\un,\iy}_R$ in which all 2-simplices satisfy an invertibility condition. 

Define $\cT\subset\hACat^{\un,\iy}_R([2])$ to be the subset of 2-simplices satisfying this condition. Then Theorem \ref{ai7thm4} shows that $(\hACat^{\un,\iy}_R,\cT)$ is an $\iy$-bicategory, and $\ACat^{\un,\iy}_R=\Core(\hACat^{\un,\iy}_R,\cT)$. The basic difference is that 2-morphisms in $(\hACat^{\un,\iy}_R,\cT)$ are unital $A_\iy$-natural transformations, but 2-morphisms in $\ACat^{\un,\iy}_R$ are $A_\iy$-natural isomorphisms.

Theorems \ref{ai6thm3} and \ref{ai8thm3} give analogues of this for $\iy$-bicategories of $A_\iy$-algebras, and of (weakly) filtered $A_\iy$-categories. These are all interesting as examples of $\iy$-bicategories which appear `in nature'.

As we explain in \S\ref{ai1.62} and \S\ref{ai411}, Faonte \cite[\S 3.4]{Faon2} defines another kind of $(\iy,2)$-category of unital $A_\iy$-categories over $R$. We expect, but do not prove, that our construction and Faonte's are equivalent in a suitable sense.

Apart from proving the existence of these $\iy$-bicategories, we do not develop their theory in this book. Perhaps they will find some applications in future. For example, they may be useful in problems involving $A_\iy$-natural transformations which are not $A_\iy$-natural isomorphisms. Also, the full $\iy$-bisubcategory $(\hDGCat^{\un,\iy}_R,\cT^\dg)\subset(\hACat^{\un,\iy}_R,\cT)$ with objects dg-categories is an $(\iy,2)$-category of dg-categories over $R$, which may be of interest.

\subsection{\texorpdfstring{Theorems about our $\iy$-categories of $A_\iy$-categories}{Theorems about our ∞-categories of A∞-categories}}
\label{ai1.5}

We now review some of the results of Chapter \ref{ai9}, which generalizes various results about $A_\iy$-categories to the $\iy$-categorical level.

\subsubsection{Mapping spaces}
\label{ai1.51}

If $\sC$ is an $\iy$-category and $x,y$ are objects of $\sC$, as in \S\ref{ai23} we can define the space of {\it left morphisms\/} $\Hom_\sC^L(x,y)$, a Kan complex which is a particular model for the mapping space $\Hom_\sC(x,y)$ up to homotopy equivalence. Being able to compute mapping spaces $\Hom_\sC(x,y)$ explicitly is useful in several ways, and a sign that we understand the $\iy$-category $\sC$. Here is part of Theorem~\ref{ai9thm1}.

\begin{thm}
\label{ai1thm1}
Let\/ $\sA,\sB$ be unital\/ $A_\iy$-categories over a commutative ring $R$. Then there is a canonical strict isomorphism of Kan complexes
\ea
\Hom_{\ACat^{\un,\iy}_R}^L(\sA,\sB)&\cong \sN_{A_\iy}(\Fun_{A_\iy}^\un(\sA,\sB))^\Kan.
\label{ai1eq19}
\ea
Here $\Hom_{\ACat^{\un,\iy}_R}^L(\sA,\sB)$ is the space of left morphisms in $\ACat^{\un,\iy}_R$ from\/ {\rm\S\ref{ai23},} one model for the mapping space $\Hom_{\ACat^{\un,\iy}_R}(\sA,\sB)$ in the $\iy$-category $\ACat^{\un,\iy}_R,$ and\/ $\sN_{A_\iy}(\Fun_{A_\iy}^\un(\sA,\sB))$ is the simplicial nerve from Definition\/ {\rm\ref{ai1def9}} of the $A_\iy$-category $\Fun_{A_\iy}^\un(\sA,\sB)$ from Definition\/ {\rm\ref{ai1def7},} and\/ $(\cdots)^\Kan$ denotes the maximal Kan complex in an $\iy$-category.

Analogues of\/ \eq{ai1eq19} hold for $\AfCat^{\un,\iy}_R,\AwCat^{\un,\iy}_R$ in the filtered case.
\end{thm}

The author views this as evidence $\ACat^{\un,\iy}_R,\AfCat^{\un,\iy}_R,\AwCat^{\un,\iy}_R$ are the `correct' definitions of $\iy$-categories generalizing $\ACat^\un_R,\ab\AfCat^\un_R$. In particular, \eq{ai1eq19} holds not just up to homotopy equivalence, but on the nose. 

The next theorem, explained in \S\ref{ai410}, is proved by Tanaka \cite[Th.~3.37]{Tana}, building on work by Canonaco--Ornaghi--Stellari \cite[Th.~C]{COS2} and To\"en~\cite{Toen}.

\begin{thm}
[Tanaka \cite{Tana}]
\label{ai1thm2}
Let\/ $\sA,\sB$ be unital\/ $A_\iy$-categories over a commutative ring $R,$ and\/ $\Phi:\sA'\ra\sA$ be a quasi-equivalence with\/ $\sA'$ h-projective with split units (such\/ $\Phi$ exist). Then we have a homotopy equivalence
\begin{equation*}
\Hom^L_{\ACat^\un_R[\QE_\un^{-1}]^\iy}(\sA,\sB)\simeq\sN_{A_\iy}\bigl(\Fun^\un_{A_\iy}(\sA',\sB)\bigr){}^\Kan,
\end{equation*}
where $\Hom^L_{\ACat^\un_R[\QE_\un^{-1}]^\iy}(\sA,\sB)$ is the mapping space in the localized\/ $\iy$-category $\ACat^\un_R[\QE_\un^{-1}]^\iy,$ for $\QE_\un$ the class of unitary quasi-equivalences in $\ACat^\un_R,$ and\/ $\sN_{A_\iy}\bigl(\Fun^\un_{A_\iy}(\sA',\sB)\bigr){}^\Kan$ is the maximal Kan complex in the simplicial nerve of the\/ $A_\iy$-category\/ $\Fun^\un_{A_\iy}(\sA',\sB)$.
\end{thm}

Here is Theorem \ref{ai9thm2}(d) below, which is new.

\begin{thm}
\label{ai1thm3}
Let\/ $\sA,\sB$ be unital\/ $A_\iy$-categories over a commutative ring $R$. Then we have a homotopy equivalence
\begin{align*}
\Hom_{\ACat^\un_R[\HE_\un^{-1}]^\iy}^L(\sA,\sB)&\simeq \sN_{A_\iy}(\Fun_{A_\iy}^\un(\sA,\sB))^\Kan,
\end{align*}
where $\HE_\un$ is the class of unital homotopy equivalences in $\ACat^\un_R$.
\end{thm}

\begin{rem}
\label{ai1rem4}
Although Theorems \ref{ai1thm2} and \ref{ai1thm3} look similar, there are some important differences:
\smallskip

\noindent{\bf(a)} The proof of Theorem \ref{ai1thm2} is indirect: it relies on an $\iy$-category equivalence $\ACat^\un_R[\QE_\un^{-1}]^\iy\simeq \DGCat_R[\QE_\dg^{-1}]^\iy$, and a model category structure on $\DGCat_R$ which enables computation of mapping spaces in $\DGCat_R[\QE_\dg^{-1}]^\iy$. The same method would only work for localized $\iy$-categories which are equivalent to~$\DGCat_R[\QE_\dg^{-1}]^\iy$.

The proof of Theorem \ref{ai1thm3} is more direct, and does not involve model categories, but instead uses the `Derived Mapping Space Lemma' of Arakawa--Carranza--Kapulkin \cite{ACK}, and an explicit `simplicial resolution'. The same method also works for localizations of many full subcategories of~$\ACat^\un_R$.
\smallskip

\noindent{\bf(b)} Theorem \ref{ai1thm2} involves a `cofibrant replacement' $\sA'$ of $\sA$, as in the theory of model categories. Theorem \ref{ai1thm3} does not.
\smallskip

\noindent{\bf(c)} Theorem \ref{ai1thm2} localizes quasi-equivalences, but Theorem \ref{ai1thm3} localizes homotopy equivalences. Although the previous literature focuses on localizing quasi-equivalences \cite{COS1,COS2,Faon2,OhTa,Orna1,Tana,Pasc}, we argue that homotopy equivalences are actually the more natural choice to invert. One reason is that homotopy equivalences {\it are already invertible up to\/ $A_\iy$-natural isomorphism}, so inverting them does less violence to the category; for example, we can take 1-morphisms in $\ACat^\un_R[\HE_\un^{-1}]^\iy$ to be just morphisms in $\ACat^\un_R$.
\end{rem}

\subsubsection{Equivalences of $\iy$-categories}
\label{ai1.52}

We consider two $\iy$-categories $\sC,\sD$ to be `the same' if there is a (Joyal) equivalence $\fF:\sC\ra\sD$. We summarize some of our results on equivalences between $\iy$-categories of $A_\iy$-categories. They are proved by combining the computations of mapping spaces in \S\ref{ai1.51} with the fact in Theorem \ref{ai2thm1} that $\fF:\sC\ra\sD$ is an equivalence if and only if $\fF$ induces homotopy equivalences of mapping spaces $\Hom_\sC(x,y)\ra\Hom_\sD(\fF(x),\fF(y))$, and $\Ho(\fF):\Ho(\sC)\ra\Ho(\sD)$ is an equivalence of ordinary categories. The next theorem combines parts of Theorems \ref{ai9thm3} and \ref{ai9thm4} and Corollary~\ref{ai9cor2}.

\begin{thm}
\label{ai1thm4}
Let $R$ be a commutative ring.
\smallskip

\noindent{\bf(a)} Consider the categories with weak equivalences $(\ACat^{\un,\hs}_R,\HE_{\un,\hs})$ and\/ $(\ACat^\un_R,\QE_\un)$. The inclusion $\inc:\ACat_R^{\un,\hs}\hookra\ACat^\un_R$ maps $\HE_\hs$ to $\QE_\un,$ and so induces an $\iy$-functor of localized\/ $\iy$-categories
\begin{equation*}
\inc_*:\ACat^{\un,\hs}_R[\HE_{\un,\hs}^{-1}]^\iy\longra\ACat^\un_R[\QE_\un^{-1}]^\iy.
\end{equation*}
Then $\inc_*$ is a Joyal equivalence.
\smallskip

\noindent{\bf(b)} As in\/ {\rm\eq{ai1eq13},} in\/ {\rm\S\ref{ai73}} we define an explicit\/ $\iy$-functor $\iinc^\un_R:\sN(\ACat^\un_R)\ab\ra\ACat^{\un,\iy}_R,$ which maps homotopy equivalences in $\ACat^\un_R$ to equivalences in $\ACat^{\un,\iy}_R$. Thus by the universal property of $\iy$-category localizations, there exists an $\iy$-functor $\Phi$ in the diagram
\begin{equation*}
\xymatrix@C=80pt@R=15pt{
*+[r]{\sN(\ACat^\un_R)} \ar[rr]_(0.55){\Pi_\loc^\iy} \ar@/_.5pc/[drr]_{\iinc^\un_R} &  \drtwocell_{}\omit^{}\omit{^{}} & *+[l]{\ACat^\un_R[\HE_\un^{-1}]^\iy} \ar@{..>}[d]^{\Phi} \\
&& *+[l]{\ACat^{\un,\iy}_R.} }
\end{equation*}
Then $\Phi:\ACat^\un_R[\HE_\un^{-1}]^\iy\ra\ACat^{\un,\iy}_R$ is a Joyal equivalence.

Furthermore, if\/ $\cC$ is any of the full subcategories of\/ $\ACat^\un_R$
\begin{align*}
&\DGCat_R^{A_\iy},\!\!\! & &\DGCat_R^{A_\iy,\rsu},\!\!\! & &\DGCat_R^{A_\iy,\hp},\!\!\! & &\DGCat_R^{A_\iy,\hs},\!\!\! &
&\ACat_R^{\un,\rsu}, \\ 
&\ACat_R^{\un,\hp},\!\!\! & &\ACat_R^{\un,\hs},\!\!\! & &\ADCat_R^{\un,\rsu},\!\!\! & &\ADCat_R^{\un,\hp},\!\!\! & &\ADCat_R^{\un,\hs}, 
\end{align*}
then the localization of\/ $\cC$ at its homotopy equivalences is equivalent to the corresponding full\/ $\iy$-subcategory of\/ $\ACat^{\un,\iy}_R$.
\smallskip

\noindent{\bf(c)} All the following $\iy$-categories are Joyal equivalent:
\ea
&\ACat^{\un,\hs,\iy}_R,  & &\DGCat^{A_\iy,\hs,\iy}_R, &
&\ACat^\un_R[\QE_\un^{-1}]^\iy,
\nonumber\\
&\ACat_R^\cu[\QE_\cu^{-1}]^\iy, & &\DGCat_R[\QE_\dg^{-1}]^\iy, &
&\ACat_R^{\un,\hs}[\HE_{\un,\hs}^{-1}]^\iy,  
\nonumber\\
&\DGCat_R^\hs[\HE_{\dg,\hs}^{-1}]^\iy,\!\!\! & &\DGCat_R^{A_\iy,\hs}[\HE_{A_\iy,\dg,\hs}^{-1}]^\iy.\!\!\!\!\!\!\!\!\!\!\!\!\!\!\!\!\!\!\!
\label{ai1eq20}
\ea
They are all presentable $\iy$-categories with closed symmetric monoidal structure $\ot^\iy,$ as this is known to be true for $\DGCat_R[\QE_\dg^{-1}]^\iy$.
\end{thm}

\begin{rem}
\label{ai1rem5}
The $\iy$-categories $\ACat^\un_R[\QE_\un^{-1}]^\iy$ and $\ACat_R^\cu[\QE_\cu^{-1}]^\iy$ in \eq{ai1eq20} are considered in the literature \cite{COS1,COS2,Faon2,OhTa,Orna1,Tana,Pasc} to be `the' $\iy$-category of $A_\iy$-categories over $R$. Part (c) shows these can be identified with full $\iy$-subcategory $\ACat^{\un,\hs,\iy}_R$ of our explicit $\iy$-category~$\ACat^{\un,\iy}_R$.
\end{rem}

We do not prove any analogues of Theorem \ref{ai1thm4} for {\it filtered\/} $A_\iy$-categories, as some of the tools involved do not work, for example, neither cohomology categories $H^0(\sA)$ nor quasi-equivalences are defined in the filtered case.

\subsubsection{The $\iy$-categorical simplicial nerve functor}
\label{ai1.53}

Definition \ref{ai1def9} defined the {\it simplicial nerve functor\/} $\sN_{A_\iy}:\ACat_R^\un\ab\ra\Cat^1_\iy$, where $\Cat^1_\iy\subset\SSet$ is the ordinary category of $\iy$-categories. To each unital $A_\iy$-category $\sA$ this associates an $\iy$-category $\sN_{A_\iy}(\sA)$ with homotopy category $\Ho(\sN_{A_\iy}(\sA))\cong H^0(\sA)$. That is, $\sN_{A_\iy}(\sA)$ is an $\iy$-categorical enhancement of the cohomology category $H^0(\sA)$. 

Here is part of Theorem \ref{ai9thm5} and Corollary~\ref{ai9cor3}.

\begin{thm}
\label{ai1thm5}
Let\/ $R$ be a commutative ring. By an explicit construction, we can define an $\iy$-functor $\sN^\iy_{A_\iy}:\ACat^{\un,\iy}_R\ra\sCat_\iy,$ where $\sCat_\iy$ is the $\iy$-category of\/ $\iy$-categories, such that the following diagram of\/ $\iy$-categories commutes up to homotopy:
\begin{equation*}
\xymatrix@C=180pt@R=15pt{
*+[r]{\sN(\ACat^\un_R)} \ar[r]_(0.36){\sN(\sN_{A_\iy})} \ar[d]^{\iinc^\un_R} \drtwocell_{}\omit^{}\omit{^{}} & *+[l]{\sN(\Cat^1_\iy)} \ar[d]_{I^{\sCat_\iy}} \\
*+[r]{\ACat^{\un,\iy}_R } \ar[r]^(0.66){\sN^\iy_{A_\iy}} & *+[l]{\sCat_\iy.\!}}
\end{equation*}
We call\/ $\sN^\iy_{A_\iy}$ the \begin{bfseries}$\iy$-categorical simplicial nerve functor\end{bfseries}.

This $\sN^\iy_{A_\iy}$ restricts to an $\iy$-functor
\e
\sN^\iy_{A_\iy}:\ADCat^{\un,\iy}_R\longra\sCat^\rst_\iy
\label{ai1eq21}
\e
from the $\iy$-category $\ADCat^{\un,\iy}_R$ of \begin{bfseries}pretriangulated\end{bfseries} unital\/ $A_\iy$-cat\-eg\-or\-ies over $R$  to the $\iy$-category $\sCat^\rst_\iy$ of \begin{bfseries}stable\end{bfseries} $\iy$-categories.
\end{thm}

\subsubsection{`Inner Hom' functors}
\label{ai1.54}

Using the equivalence $\ACat^\un_R[\HE_\un^{-1}]^\iy\simeq\ACat^{\un,\iy}_R$ in Theorem \ref{ai1thm4}(b), in Theorem \ref{ai9thm6} we prove:

\begin{thm}
\label{ai1thm6}
Let\/ $R$ be a commutative ring, and\/ $\Fun_{A_\iy}^\un$ be as in Definition\/ {\rm\ref{ai1def7}}. Then there exists an\/ $\iy$-functor $\Fun_{A_\iy}^{\un,\iy},$ unique up to homotopy, such that the following diagram of\/ $\iy$-categories commutes up to homotopy:
\e
\begin{gathered}
\xymatrix@!0@C=280pt@R=35pt{
*+[r]{\sN(\ACat^\un_R)^\rop\t\sN(\ACat^\un_R)} \ar[r]_(0.66){\sN(\Fun_{A_\iy}^\un)} \ar[d]^{(\iinc^\un_R)^\rop\t\iinc^\un_R} \drtwocell_{}\omit^{}\omit{^{}} & *+[l]{\sN(\ACat^\un_R)} \ar[d]_{\iinc^\un_R} \\
*+[r]{(\ACat^{\un,\iy}_R)^\rop\t\ACat^{\un,\iy}_R} \ar[r]^(0.66){\Fun_{A_\iy}^{\un,\iy}} & *+[l]{\ACat^{\un,\iy}_R.\!}}
\end{gathered}
\label{ai1eq22}
\e
\end{thm}

Since $\Fun_{A_\iy}^{\un,\iy}$ maps between explicitly-defined $\iy$-categories, it is an obvious question whether we can give explicit formulae for it. The author has not managed to do this. However, in \S\ref{ai94} we do define explicit  $\iy$-functors
\ea
\hFun_{A_\iy}^{\un,\iy}&:\sN(\ACat^\un_R)^\rop\t\ACat^{\un,\iy}_R \longra \ACat^{\un,\iy}_R,
\label{ai1eq23}\\
\cFun_{A_\iy}^{\un,\iy}&:(\ACat^{\un,\iy}_R)^\rop\t\sN(\ACat^\un_R) \longra \ACat^{\un,\iy}_R.
\label{ai1eq24}
\ea
which are homotopic to the restrictions of $\Fun_{A_\iy}^{\un,\iy}$ to the domain $\iy$-sub\-cat\-eg\-or\-ies. So for example, linking \eq{ai1eq22} and \eq{ai1eq23}, in the following diagram of\/ $\iy$-categories, the top rectangle commutes strictly, and the bottom commutes up to homotopy:
\begin{equation*}
\xymatrix@!0@C=280pt@R=35pt{
*+[r]{\sN(\ACat^\un_R)^\rop\t\sN(\ACat^\un_R)} \ar[r]_(0.66){\sN(\Fun_{A_\iy}^\un)} \ar[d]^{\id\t\iinc^\un_R} \ar@{}@<-1ex>[dr]^{=} & *+[l]{\sN(\ACat^\un_R)} \ar[d]_{\iinc^\un_R} \\
*+[r]{\sN(\ACat^\un_R)^\rop\t\ACat^{\un,\iy}_R} \ar[r]^(0.66){\hFun_{A_\iy}^{\un,\iy}} \ar[d]^{(\iinc^\un_R)^\rop\t\id} \drtwocell_{}\omit^{}\omit{^{}} & *+[l]{\ACat^{\un,\iy}_R} \ar@{=}[d]_{\id} \\
*+[r]{(\ACat^{\un,\iy}_R)^\rop\t\ACat^{\un,\iy}_R} \ar@{..>}[r]^(0.66){\Fun_{A_\iy}^{\un,\iy}} & *+[l]{\ACat^{\un,\iy}_R.\!}}
\end{equation*}

For {\it filtered\/} $A_\iy$-categories $\AfCat^{\un,\iy}_R$, $\AwCat^{\un,\iy}_R$, we do not construct analogues of \eq{ai1eq22}, as we do not have filtered analogues of the equivalence $\ACat^\un_R[\HE_\un^{-1}]^\iy\simeq\ACat^{\un,\iy}_R$. However, we do define explicit analogues of \eq{ai1eq23}--\eq{ai1eq24} in the filtered case:
\ea
\hFun_{A^\rfi_\iy}^{\un,\iy}&:\sN(\AfCat^\un_R)^\rop\t\AfCat^{\un,\iy}_R \longra \AfCat^{\un,\iy}_{R,2,2},
\label{ai1eq25}\\
\hFun_{A^\wf_\iy}^{\un,\iy}&:\sN(\AfCat^\un_R)^\rop\t\AwCat^{\un,\iy}_R \longra \AwCat^{\un,\iy}_{R,2,2},
\label{ai1eq26}\\
\cFun_{A^\rfi_\iy}^{\un,\iy}&:(\AfCat^{\un,\iy}_R)^\rop\t\sN(\AfCat^\un_R) \longra \AfCat^{\un,\iy}_{R,2,2},
\label{ai1eq27}\\
\cFun_{A^\wf_\iy}^{\un,\iy}&:(\AwCat^{\un,\iy}_R)^\rop\t\sN(\AfCat^\un_R) \longra \AwCat^{\un,\iy}_{R,2,2}.
\label{ai1eq28}
\ea

\subsubsection{Some $\iy$-functors involving filtered $A_\iy$-categories}
\label{ai1.55}

Let $R$ be a complete $\Ga$-filtered commutative ring. As in \S\ref{ai1.3}, we define $\iy$-categories $\AfCat^{\un,\iy}_R,\AfCat^{\un,\iy}_R$ of unital filtered $A_\iy$-categories over $R$, and we define $\iy$-subcategories $\AfCat^{\un,\iy}_{R,2,2}\!\subset\!\AfCat^{\un,\iy}_R$ and $\AwCat^{\un,\iy}_{R,2,2}\!\subset\!\AwCat^{\un,\iy}_R$ of type (2,2), in which the objects $\sA$ and 1-morphisms $\sF$ are {\it strict\/} (they have $\m_0^x=\fF_0^x=0$ for all objects $x$). In \S\ref{ai81}, we define $\iy$-functors
\ea
\Pi_{A_\iy^\rfi}^{A_\iy}&:\AfCat^{\un,\iy}_{R,2,2}\longra\ACat^{\un,\iy}_R,
\label{ai1eq29}\\
\Pi_{A_\iy^\wf}^{A_\iy}&:\AwCat^{\un,\iy}_{R,2,2}\longra\ACat^{\un,\iy}_R,
\label{ai1eq30}
\ea
which are analogues of \eq{ai1eq10}. They are forgetful $\iy$-functors, which forget the filtration information.

In \S\ref{ai94} we use \eq{ai1eq25}--\eq{ai1eq26} to define $\iy$-categorical analogues of \eq{ai1eq11}:
\ea
\Str^\iy=\hFun_{A^\rfi_\iy}^{\un,\iy}(\sP^\un,-)&:\AfCat^{\un,\iy}_R\longra\AfCat^{\un,\iy}_{R,2,2},
\label{ai1eq31}\\
\Str_\wf^\iy=\hFun_{A^\wf_\iy}^{\un,\iy}(\sP^\un,-)&:\AwCat^{\un,\iy}_R\longra\AwCat^{\un,\iy}_{R,2,2}.
\label{ai1eq32}
\ea
Here $\sP^\un$ is the unital filtered $A_\iy$-category with one object $*$ and $\Hom_{\sP^\un}(*,*)\ab=R$. We can compose \eq{ai1eq29}--\eq{ai1eq30} and \eq{ai1eq31}--\eq{ai1eq32} as in~\eq{ai1eq12}.

\subsubsection{The twist $\iy$-functor}
\label{ai1.56}

Let $R$ be a commutative ring. In \S\ref{ai95} we define an explicit $\iy$-functor
\e
\Tw^\iy:\ACat^{\un,\iy}_R \longra \ADCat^{\un,\iy}_R,
\label{ai1eq33}
\e
which is an $\iy$-categorical analogue of $\Tw:\ACat^\un_R\ra\ADCat^\un_R$ in Definition \ref{ai1def10}. The following diagram of $\iy$-functors strictly commutes:
\begin{equation*}
\xymatrix@C=180pt@R=15pt{
*+[r]{\sN(\ACat^\un_R)} \ar@<-1ex>@{}[dr]^{=} \ar[d]^{\iinc^\un_R} \ar[r]_(0.4){\sN(\Tw)} & *+[l]{\sN(\ADCat^\un_R)}  \ar[d]_{\iinc^\un_R} \\
*+[r]{\ACat^{\un,\iy}_R} \ar[r]^(0.6){\Tw^\iy} & *+[l]{\ADCat^{\un,\iy}_R.\!} }
\end{equation*}

Here is part of Theorem \ref{ai9thm9}:

\begin{thm}
\label{ai1thm7}
{\bf(a)} There is an $\iy$-natural isomorphism $\Tw^\iy\Ra\Tw^\iy\ci\ab\Tw^\iy$.

\smallskip

\noindent{\bf(b)} $\Tw^\iy$ is left adjoint to the fully faithful inclusion $\iy$-functor
\begin{equation*}
\inc:\ADCat^{\un,\iy}_R\longhookrightarrow\ACat^{\un,\iy}_R.
\end{equation*}
\end{thm}

Observe that if $R$ is a commutative ring we have a composition of $\iy$-functors
\begin{equation*}
\xymatrix@C=60pt{ 
\ACat^{\un,\iy}_R \ar[r]^{\Tw^\iy}_{\eq{ai1eq33}} & \ADCat^{\un,\iy}_R \ar[r]^{\sN^\iy_{A_\iy}}_{\eq{ai1eq21}} &
\sCat^\rst_\iy, }
\end{equation*}
and if $R$ is a complete $\Ga$-filtered commutative ring we have compositions
\begin{equation*}
\text{\begin{footnotesize}$
\xymatrix@C=21.5pt@R=15pt{ 
\AfCat^{\un,\iy}_R \ar[r]^{\Str^\iy}_{\eq{ai1eq31}} &
\AfCat^{\un,\iy}_{R,2,2} \ar[r]^{\Pi_{A_\iy^\rfi}^{A_\iy}}_{\eq{ai1eq29}} &
\ACat^{\un,\iy}_R \ar[r]^{\Tw^\iy}_{\eq{ai1eq33}} & \ADCat^{\un,\iy}_R \ar[r]^(0.55){\sN^\iy_{A_\iy}}_(0.55){\eq{ai1eq21}} & \sCat^\rst_\iy, \\
\AwCat^{\un,\iy}_R \ar[r]^{\Str^\iy_\wf}_{\eq{ai1eq32}} &
\AwCat^{\un,\iy}_{R,2,2} \ar[r]^{\Pi_{A_\iy^\wf}^{A_\iy}}_{\eq{ai1eq30}} &
\ACat^{\un,\iy}_R \ar[r]^{\Tw^\iy}_{\eq{ai1eq33}} & \ADCat^{\un,\iy}_R \ar[r]^(0.55){\sN^\iy_{A_\iy}}_(0.55){\eq{ai1eq21}} & \sCat^\rst_\iy. }$\end{footnotesize}}
\end{equation*}
These map (filtered) unital $A_\iy$-categories over $R$ to stable $\iy$-categories, in an $\iy$-categorical way. For a symplectic manifold $(X,\om)$, these compositions should be used to transform the Fukaya category $\sF(X,\om)$ (a unital $A_\iy$-category or unital filtered $A_\iy$-category), to the derived Fukaya category $D^b\sF(X,\om)$ (a stable $\iy$-category, so that $\Ho(D^b\sF(X,\om))$ is a triangulated category).

\subsection{Relation to other work in the literature}
\label{ai1.6}

We now compare our $\iy$-categories of $A_\iy$-algebras and $A_\iy$-categories with work by other authors.

\subsubsection{Comparison with localized $\iy$-categories of $A_\iy$-categories}
\label{ai1.61}

As we explain in \S\ref{ai410} below, Canonaco--Ornaghi--Stellari \cite{COS1,COS2}, Oh--Tanaka \cite{OhTa}, Ornaghi \cite{Orna1}, Tanaka \cite{Tana}, and Pascaleff \cite{Pasc} study $\iy$-categories of $A_\iy$-categories defined by starting with the category $\ACat^\un_R$ (or some variant), and localizing the quasi-equivalences $\QE_\un\subset\Mor(\ACat^\un_R)$ from Definition \ref{ai4def15}, as in \S\ref{ai26}, to get an $\iy$-category~$\ACat^\un_R[\QE_\un^{-1}]^\iy$.

Section \ref{ai92}, summarized in Theorem \ref{ai1thm4}, explores the relation between our $\iy$-category $\ACat^{\un,\iy}_R$ and this localized $\iy$-category $\ACat^\un_R[\QE_\un^{-1}]^\iy$, which is considered to be `the' $\iy$-category of $A_\iy$-categories by the authors of \cite{COS1,COS2,OhTa,Orna1,Tana,Pasc}. In particular, we show that $\ACat^\un_R[\QE_\un^{-1}]^\iy$ is Joyal equivalent to the full $\iy$-subcategory $\ACat^{\un,\hs,\iy}_R\subset\ACat^{\un,\iy}_R$ with objects $A_\iy$-categories which are h-projective with split units.

Theorem \ref{ai1thm1} computes mapping spaces in $\ACat^{\un,\iy}_R$. As in Theorem \ref{ai1thm2}, work in \cite{CaSt,COS2,Faon2,Tana,Toen} computes mapping spaces in $\ACat^\un_R[\QE_\un^{-1}]^\iy$, and after cofibrant replacement gives the same answer.

The authors of \cite{COS1,COS2,OhTa,Orna1,Tana,Pasc} take the view that to form an $\iy$-category of unital $A_\iy$-categories by localization, one should invert {\it quasi-equivalences}, as in Definition \ref{ai1def8}. One reason may be the Tabuada model structure \cite{Tabu} on $\DGCat_R$ from Theorem \ref{ai4thm1}, whose weak equivalences are quasi-equivalences.

Our results suggest that instead inverting {\it homotopy equivalences\/} may be a fruitful thing to do. Inverting homotopy equivalences is a more gentle thing to do, as homotopy equivalences are already invertible up to $A_\iy$-natural isomorphism; so  cofibrant replacement is unnecessary. As in Theorem \ref{ai1thm4}(b), we prove an equivalence~$\ACat^\un_R[\HE_\un^{-1}]^\iy\simeq\ACat^{\un,\iy}_R$.

\subsubsection{Comparison with Faonte \cite{Faon2}}
\label{ai1.62}

Section \ref{ai411} summarizes work of Faonte \cite[\S 3]{Faon2}. This is in one way the closest paper the author has found to our own work. For $R$ a field of characteristic zero, Faonte defines two explicit simplicial categories of unital $A_\iy$-categories over $R$, the first a category enriched in Kan complexes, and the second a category enriched in $\iy$-categories (rather than just enriched in simplicial sets). These can be interpreted as $(\iy,1)$- and $(\iy,2)$-categories of $A_\iy$-categories. 

They should be compared with our $\iy$-category $\ACat^{\un,\iy}_R$ from \S\ref{ai1.3}, and our $\iy$-bicategory $\hACat^{\un,\iy}_R$ from \S\ref{ai1.4}, respectively. In particular, Theorem \ref{ai9thm1}(a) implies that the left mapping spaces in $\ACat^{\un,\iy}_R,\hACat^{\un,\iy}_R$ agree with those in Faonte's simplicial categories not just up to homotopy equivalence, but by strict isomorphism.

Because of our Theorem \ref{ai9thm2}(d), the author expects that if Faonte's construction is done over a general ring $R$, it will yield an $(\iy,1)$-category equivalent to $\ACat^\un_R[\HE_\un^{-1}]^\iy$ rather than~$\ACat^\un_R[\QE_\un^{-1}]^\iy$.

It seems likely that under some suitable equivalences between different models of $(\iy,1)$- and $(\iy,2)$-categories, Faonte's constructions may be equivalent to our $\ACat^{\un,\iy}_R,\hACat^{\un,\iy}_R$. But we do not prove this.

\subsubsection{Comparison with Cazassus \cite{Caza}}
\label{ai1.63}

Section \ref{ai37} briefly describes work by Cazassus \cite[\S 5.2]{Caza}, who defines an explicit $\iy$-category $\mathop{\text{f-Alg}}_R$ of (nonunital) $A_\iy$-algebras over $R$. This is in a different way the closest paper the author has found to our own work since, like us in Chapter \ref{ai6}, Cazassus constructs an explicit simplicial set whose 0-simplices are $A_\iy$-algebras over $R$, and shows it is an $\iy$-category.

The author does not yet understand the exact relation between $\mathop{\text{f-Alg}}_R$ and our $\iy$-categories $\AAlg^{\un,\iy}_R,\AAlg^\iy_R,\ldots$ from Chapter \ref{ai6}. They are not equivalent. For $N\ge 2$, the extra data in an $N$-simplex in $\AAlg^{\un,\iy}_R$ which is not contained in lower-dimensional simplices is maps $\m_k^{\{0,\ldots,N\}}:(A^0)^{\ot^k}\ra A^N$ for $k\ge 0$, but the extra data in an $N$-simplex in $\mathop{\text{f-Alg}}_R$ is finite families of maps $D_{k,l}^\al:(A^0)^{\ot^k}\ra (A^N)^{\ot^l}$ for $k,l\ge 0$ and $\al\in F_{k,l}$. So it seems that $\mathop{\text{f-Alg}}_R$ is `larger' (has many more $N$-simplices) than $\AAlg^{\un,\iy}_R,\ldots.$ Also $\mathop{\text{f-Alg}}_R$ is a more linear object than $\AAlg^{\un,\iy}_R$, as the defining equations of $\mathop{\text{f-Alg}}_R$ are linear in the $D_{k,l}^\al$, but those of $\AAlg^{\un,\iy}_R$ are polynomial in the $\m_k^I$. 

The 2-morphisms in $\mathop{\text{f-Alg}}_R$ seem not correspond to either $A_\iy$-homotopies, or $A_\iy$-natural isomorphisms. But they do have the property that if $A_\iy$-functors $\f,\g:(A,\m)\ra(B,\n)$ are 2-isomorphic in $\mathop{\text{f-Alg}}_R$ then $H^0(\f)=H^0(\g)$, which is necessary to use $\mathop{\text{f-Alg}}_R$ to study functoriality of Morse cohomology.

The author expects that some of our $\iy$-categories $\AAlg^{\un,\iy}_{R,K},\ldots$ from Chapter \ref{ai6} can be embedded in $\mathop{\text{f-Alg}}_R$ as $\iy$-subcategories. It would be interesting to know whether Cazassus' `Morse complex' $\iy$-functor $CM^\bu:\sN(\Man^\rop)\ab\ra\mathop{\text{f-Alg}}_R$ from \cite[Th.~C]{Caza} factors via one of our $\iy$-categories.

\subsection{On Fukaya categories, and the origin of this book}
\label{ai1.7}

The author came to write this book as part of a long-term in-progress project \cite{CJR}, joint with Guillem Cazassus and Alexander Ritter, to develop a theory of (derived) Fukaya categories $\sF(X,\om),D^b\sF(X,\om)$ of semipositive symplectic manifolds $(X,\om)$. The theory should have certain desirable characteristics we want for applications, for example, that all compact immersed Lagrangians $L$ in $(X,\om)$ should appear as objects in~$\sF(X,\om)$.

The reader should understand that although the idea of Fukaya categories has been around for 30 years, constructing a fully general theory of (derived) Fukaya categories with all the properties one would like is a hugely long and challenging task, which is far from complete. The theories currently available either involve restrictive assumptions, for example Seidel \cite{Seid} for exact symplectic manifolds, or they are partial, for example Fukaya--Oh--Ohta--Ono \cite{FOOO} define Lagrangian Floer cohomology for embedded Lagrangians (roughly, the morphism spaces in Fukaya categories), but do not go as far as Fukaya categories.

The author hopes that this book will form part of the foundations of the fully general theory of Fukaya categories, when this emerges. We will now explain why our explicit $\iy$-categories of $A_\iy$-categories, such as $\ACat_R^{\un,\iy}$, will be helpful in this project. See \S\ref{ai96} below for more details. We will (over)simplify the discussion by not including {\it filtered\/} $A_\iy$-categories, as in Chapters \ref{ai5} and~\ref{ai8}.

Given a symplectic manifold $(X,\om)$, we want to construct its {\it Fukaya category\/} $\sF(X,\om)$, which is an $A_\iy$-category (let us say, strictly unital) over a commutative ring $R$ (this could be $\Z$ for $(X,\om)$ exact, or the Novikov ring $\La^\nov_\Q$ or $\La^\nov_\Z$ for more general $(X,\om$)). Defining $\sF(X,\om)$ involves difficult analysis of moduli spaces of $J$-holomorphic curves in $X$ with boundary in Lagrangians $L_0,\ldots,L_k$. Then by a category-theoretic construction, starting from $\sF(X,\om)$ we define the {\it derived Fukaya category\/} $D^b\sF(X,\om)$, which is a triangulated category, or (better) a stable $\iy$-category, and is what we really care about. 

To define $\sF(X,\om)$ we must choose some {\it arbitrary auxiliary data\/} $\fD$, so we write it $\sF(X,\om)_\fD$. Here $\fD$ should include a choice of almost complex structure $J$ on $X$ compatible with $\om$ (and possibly domain-dependent). It could also include families of perturbations of Lagrangians as in Seidel \cite{Seid}, choices of Morse functions, \ldots. The details depend on the set-up, and there are several different approaches in the literature.

We would like $\sF(X,\om)_\fD$ and $D^b\sF(X,\om)_\fD$ to be {\it independent of the extra data\/} $\fD$, and thus be invariants of $(X,\om)$. If $\fD_1,\fD_2$ are choices for $\fD$, in general $\sF(X,\om)_{\fD_1},D^b\sF(X,\om)_{\fD_1}$ and $\sF(X,\om)_{\fD_2},D^b\sF(X,\om)_{\fD_2}$ are not isomorphic, but they are `homotopy equivalent'.

Our explicit $\iy$-categories of $A_\iy$-categories, such as $\ACat_R^{\un,\iy}$, can be used to make precise how $\sF(X,\om)_\fD$ depends on the auxiliary data $\fD$. To do this, we should define a contractible Kan complex $\sS_{X,\om}$ whose 0-simplices are choices of data $\fD$, and an $\iy$-functor $F:\sS_{X,\om}\ra\ACat^{\un,\iy}_R$ which maps $F:\fD\mapsto\sF(X,\om)_\fD$ on 0-simplices. As a contractible Kan complex is equivalent to a point in $\iy$-category theory, this means the $\sF(X,\om)_\fD$ are independent of $\fD$ up to coherent homotopies, in an $\iy$-categorical sense.

It is important here that $F$ maps to our {\it explicit\/} $\iy$-category $\ACat_R^{\un,\iy}$, and not to a non-explicit localized $\iy$-category such as $\ACat_R^\un[\HE_\un^{-1}]^\iy$. This is because the definition of $F$ is also explicit: we generalize the construction of $\sF(X,\om)_\fD$ to a structure depending on an $N$-simplex $\De^N$ family of choices of $\fD$, and this maps to an $N$-simplex in~$\ACat_R^{\un,\iy}$.

\subsection{Remarks on Set Theory issues, and conventions}
\label{ai1.8}

\begin{rem}
\label{ai1rem2}
{\bf (On set theory foundations.)} In conventional category theory (see for example MacLane \cite{MacL}) one allows the collection of objects $\mathop{\rm Obj}\scr C$ of a category $\scr C$ to be a {\it proper class}, which may be larger than a set. A category for which $\mathop{\rm Obj}\scr C$ is a set is called a {\it small category}. The distinction between sets and classes is necessary in Set Theory to avoid paradoxes such as `the set of all sets which are not members of themselves' (since the family of all sets is not a set, it is a proper class).

Now in the theory of $\iy$-categories in Chapter \ref{ai2}, an $\iy$-category ${\scr C}_\iy$ is a simplicial set, so its family of objects $\mathop{\rm Obj}{\scr C}_\iy$ must be a set. Thus, to promote categories to $\iy$-categories, as we wish to do, we must restrict our attention to small categories. However, if for example we just defined $\AAlg_R$ to be the category of all $A_\iy$-algebras over a ring $R$, it would not be a small category.

There is a standard way of dealing with this problem: as in Grothendieck--Verdier \cite[\S 0--\S 1]{GrVe} and Bourbaki \cite{Bour}, we fix a `Grothendieck universe' $\sU$, a very large set in which all of mathematics can be performed, and we assume all sets in our theory are subsets of $\sU$ (that is, they are `$\sU$-small'). See Lurie \cite[\S 1.2.15, \S 5.4]{Luri2} for a discussion of these issues in the context of $\iy$-categories.

This book raises no new Set Theory problems not already discussed in \cite{Bour,GrVe,Luri2} and many other places. So we will just fix a Grothendieck universe $\sU$, assume implicitly from now on that all structures in this book are $\sU$-small, and not worry about the Set Theory issues. (Actually, for some purposes this is an oversimplification: there are constructions which require one to pass from a universe $\sU$ to a larger one $\sU'$. We will ignore this.)
\end{rem}

\begin{rem}
\label{ai1rem3}
{\bf(On sign, order, and degree conventions.)} There are multiple different conventions in the literature on $A_\iy$-algebras and $A_\iy$-categories on {\it signs\/}: for example $(-1)^{i+jk}$ in the defining equation \eq{ai3eq3} for $A_\iy$-algebras
\begin{equation*}
\sum_{i,k\ge 0,\; j\ge 1:\; i+j+k=n}
(-1)^{i+jk}\m_{i+k+1}\ci\bigl(\id_A^{\ot^i}\ot\m_j\ot \id_A^{\ot^k}\bigr)=0,
\end{equation*}
and on {\it orders\/}: for example, in \S\ref{ai41} we define $A_\iy$-categories using morphisms
\begin{equation*}
\m_k:\Hom(x_0,x_1)\ot\Hom(x_1,x_2)\ot\cdots \ot\Hom(x_{k-1},x_k)\longra\Hom(x_0,x_k),
\end{equation*}
but Seidel \cite[\S 1]{Seid} defines them using morphisms
\begin{equation*}
\m_k:\Hom(x_{k-1},x_k)\ot\Hom(x_{k-2},x_{k-1})\ot\cdots \ot\Hom(x_0,x_1)\longra\Hom(x_0,x_k),
\end{equation*}
and {\it degrees}, for example whether we regard $\m_k$ as having degree $2-k$ or 1.

Changing order convention is straightforward: we just redefine $\smash{\,\,\,\widetilde{\!\!\!\Hom\!\!\!}\,\,\,}(x,y)=\Hom(y,x)$, reversing the order in which we regard morphisms as acting, and make no other changes. This is like taking an `opposite category'.

Signs in $A_\iy$-theory (and degrees, which affect signs) are always painful. I have tried to choose sign conventions which give the simplest signs in important formulae, and which are consistent throughout the book. For $A_\iy$-algebras, our signs in \S\ref{ai32}--\S\ref{ai33} follow Keller \cite{Kell1}. Beyond this, most of the signs were computed by the author, and usually differ from those in our primary sources.
\end{rem}

\noindent{\it Use of generative AI.} This project began in 2022-23, and was 90\% finished before AI began to challenge human mathematicians. In late 2026 I had completed proofs of all the major theorems except Theorems \ref{ai9thm3} and \ref{ai9thm4}. I had several conversations with ChatGPU Edu about $\iy$-categories and issues in the proofs, in general terms, without uploading any part of this document. 

ChatGPT was really useful at helping me understand colimits in $\iy$-cat\-eg\-or\-ies and Arakawa--Carranza--Kapulkin \cite{ACK} better than I had before. It also suggested references which helped me finish the proofs: for Theorem \ref{ai9thm3} it pointed out Holstein \cite{Hols}, and for Theorem \ref{ai9thm4} it located results I wanted in Lurie \cite{Luri2,Luri5}, including the `last vertex map'. Apart from these, all results, proofs, and mistakes in this book were written by a human.
\smallskip

\noindent{\it Acknowledgements.} I would like to thank Mohammed Abouzaid, Denis Auroux, Lukas Brantner, Alberto Canonaco, Guillem Cazassus, Jason Lotay, Mattia Ornaghi, Alexander Ritter, Paolo Stellari, Hiro Lee Tanaka, and ChatGPT Edu for helpful conversations. This research was supported by EPSRC grant EP/T012749/1. For the purpose of open access, the author has applied a CC BY public copyright licence to any Author Accepted Manuscript (AAM) version arising from this submission.

\section{\texorpdfstring{Simplicial sets and $\iy$-categories}{Simplicial sets and ∞-categories}}
\label{ai2}

\subsection{Simplicial sets}
\label{ai21}

For an elementary introduction to simplicial sets we recommend Friedman \cite{Frie}, and for foundational material see Goerss--Jardine \cite[\S I]{GoJa}. 

\begin{dfn}
\label{ai2def1}
The {\it simplex category\/} $\De$ has objects sets $[N]=\{0,1,\ldots,N\}$ for $N\ge 0$ and morphisms $f:[N]\ra[M]$ the order-preserving functions, that is, if $0\le i\le j\le N$ then $f(i)\le f(j)$.

A {\it simplicial set\/} is a functor $S:\De^\rop\ra\Set$, where $\Set$ is the category of sets. A {\it morphism of simplicial sets\/} $f:S\ra T$, or just a {\it map of simplicial sets}, is a natural transformation of functors $f:S\Ra T$. This makes simplicial sets into a category $\SSet$.	It is complete and cocomplete, and (co)limits may be defined objectwise in $\De^\rop$ by taking (co)limits in~$\Set$.

Define morphisms in $\De$, the {\it face maps\/} $\de^{N,i}:[N-1]\ra[N]$ and {\it degeneracy maps\/} $\si^{N,i}:[N+1]\ra[N]$ for $i=0,\ldots,N$ by
\begin{equation*}
\de^{N,i}(j)=\begin{cases} j, & j<i, \\ j+1, & j\ge i, \end{cases}\quad 
\si^{N,i}(j)=\begin{cases} j, & j\le i, \\ j-1, & j>i.
\end{cases}
\end{equation*}
That is, $\de^{N,i}$ misses $i$, and $\si^{N,i}$ repeats $i$. These satisfy the identities
\e
\begin{aligned}
\de^{N,j}\ci\de^{N-1,i}&=\de^{N,i}\ci\de^{N-1,j-1},& 0&\le i<j\le N,	\\
\si^{N,j}\ci\si^{N+1,i}&=\si^{N,i}\ci\si^{N+1,j+1},& 0&\le i\le j\le N,	\\
\si^{N,j}\ci\de^{N+1,i}&=\de^{N,i}\ci\si^{N-1,j-1},& 0&\le i<j\le N,	\\
\si^{N-1,j}\ci\de^{N,i}&=\de^{N-1,i-1}\ci\si^{N-2,j},& 0&\le j<j+1<i\le N,	\\
\si^{N,j}\ci\de^{N+1,i}&=\id, & i&=j\text{ or }j+1.
\end{aligned}
\label{ai2eq1}
\e
The category $\De$ is generated by the $\de^{N,i},\si^{N,i}$ subject only to the relations \eq{ai2eq1}.

Now let $S:\De^\rop\ra\Set$ be a simplicial set. Write $S_N=S([N])$, and call $S_N$ the set of $N$-{\it simplices\/} of $S$. Define {\it face maps\/} $d_{N,i}=S(\de^{N,i}):S_N\ra S_{N-1}$ and {\it degeneracy maps\/} $s_{N,i}=S(\si^{N,i}):S_N\ra S_{N+1}$. Reversing directions of morphisms, these satisfy
\e
\begin{aligned}
d_{N-1,i}\ci d_{N,j}&=d_{N-1,j-1}\ci d_{N,i},& 0&\le i<j\le N, \\
s_{N+1,i}\ci s_{N,j}&=s_{N+1,j+1}\ci s_{N,i},& 0&\le i\le j\le N,	\\
d_{N+1,i}\ci s_{N,j}&=s_{N-1,j-1}\ci d_{N,i},& 0&\le i<j\le N, \\
d_{N,i}\ci s_{N-1,j}&=s_{N-2,j}\ci d_{N-1,i-1},& 0&\le j<j+1<i\le N,	\\
d_{N+1,i}\ci s_{N,j}&=\id, & i&=j\text{ or }j+1.
\end{aligned}
\label{ai2eq2}
\e
As $\De$ is generated by the $\de^{N,i},\si^{N,i}$ subject to \eq{ai2eq1}, to define a simplicial set $S$ it is enough to give sets $S_N=S(\De_N)$ for $N\ge 0$ and maps $d_{N,i}:S_N\ra S_{N-1}$ and $s_{N,i}:S_N\ra S_{N+1}$ for $0\le i\le N$ satisfying \eq{ai2eq2}, and all the other morphisms $S(f)$ in $S$ for $f:[k]\ra[l]$ can be written as compositions of the $d_{N,i},s_{N,i}$. We also write $d^S_{N,i},s^S_{N,i}$ for $d_{N,i},s_{N,i}$ if we need to specify the simplicial set $S$.

This gives us a way to draw a picture of a simplicial set:
\e
\begin{gathered}
\xymatrix@C=60pt{
S_0\,\,
\ar@<0ex>[r]^{s_{0,0}}
&
\,\,S_1\,\,
\ar@<-3ex>[l]_{d_{1,0}}
\ar@<3ex>[l]_{d_{1,1}}
\ar@<-3ex>[r]^{s_{1,0}}
\ar@<3ex>[r]^{s_{1,1}}
& 
\,\,S_2\,\,
\ar@<-6ex>[l]_{d_{2,0}}
\ar@<0ex>[l]_{d_{2,1}}
\ar@<6ex>[l]_{d_{2,2}}
\ar@<-6ex>[r]^{s_{2,0}}
\ar@<0ex>[r]^{s_{2,1}}
\ar@<6ex>[r]^{s_{2,2}}
&
\,\,S_3\qquad \cdots. 
\ar@<-9ex>[l]_{d_{3,0}}
\ar@<-3ex>[l]_{d_{3,1}}
\ar@<3ex>[l]_{d_{3,2}}
\ar@<9ex>[l]_{d_{3,3}} }\!\!\!\!\!\!\!\!\!\!
\end{gathered}
\label{ai2eq3}
\e

An $N$-simplex $\Si\in S_N$ is called {\it degenerate\/} if it lies in the image of $s_{N-1,i}:S_{N-1}\ra S_N$ for some $i=0,\ldots,N-1$, and {\it nondegenerate\/} otherwise.

For each $\Si\in S_N$ there exist $0\le M\le N$, a surjective morphism $\al:[N]\ra[M]$ in $\De$, and a nondegenerate $M$-simplex $\ti\Si\in S_M$ with $\Si=S(\al)(\ti\Si)$, and $M,\al,\ti\Si$ are unique. Here $\Si$ is nondegenerate if $M=N$, $\al=\id_{[N]}$ and~$\Si=\ti\Si$.

A map of simplicial sets $f:S\ra T$ is equivalent to maps $f_N:S_N\ra T_N$ for $N=0,1,2,\ldots$ commuting with the face and degeneracy maps~$d_{N,i},s_{N,i}$.

If $S,T$ are simplicial sets, we call $S$ a {\it simplicial subset\/} of $T$, written $S\subseteq T$, if $S_N\subseteq T_N$ for $N=0,1,\ldots$ and the face and degeneracy maps satisfy $d_{N,i}^S=d_{N,i}^T\vert_{S_N}$ and $s_{N,i}^S=s_{N,i}^T\vert_{S_N}$ for all $N,i$. Inclusion $\inc:S\hookra T$ is a map of simplicial sets.

Now let $\sC$ be a category. A {\it simplicial object in\/} $\sC$ is a functor $F:\De^\rop\ra\sC$. Thus a simplicial set is a simplicial object in $\Set$. One often abbreviates `simplicial object in things' to `simplicial thing'. To pass from classical Algebraic Geometry to Derived Algebraic Geometry, one often replaces objects by simplicial objects. For example, a {\it simplicial commutative\/ $\K$-algebra\/} (simplicial object in $\Alg_\K$) is a functor $A:\De^\rop\ra\Alg_\K$, and these are the appropriate notion of `derived commutative algebra' used to define derived $\K$-schemes.

We can draw a picture of a simplicial object in $\sC$ as a diagram \eq{ai2eq3}, where now the $S_i$ are objects in $\sC$ and the $d_{N,i},s_{N,i}$ are morphisms in $\sC$ satisfying~\eq{ai2eq2}.
\end{dfn}

\begin{dfn}
\label{ai2def2}
{\bf(i)} The {\it topological $N$-simplex\/} $\De^N_\top$ is
\e
\De^N_\top=\bigl\{(x_0,\ldots,x_N)\in \R^{N+1}:x_i\ge 0,\;\> x_0+\cdots+x_N=1\bigr\}.
\label{ai2eq4}
\e
If $f:[N]\ra[M]$ is order-preserving we define $f_\top:\De^M_\top\ra \De^N_\top$ by
\begin{align*}
f_\top(x_0,\ldots,x_M)=(y_0,\ldots,y_N),\quad 
y_i=\ts\sum_{j\in\{0,\ldots,M\}:f(j)=i}x_j.
\end{align*}
This defines a functor $G:\De^\rop\ra\CGHaus$ mapping $[N]\mapsto\De^N_\top$, $f\mapsto f_\top$, where $\CGHaus\subset\Top$ is the category of compactly generated Hausdorff topological spaces.
\smallskip

\noindent{\bf(ii)} Let $X$ be a compactly-generated, Hausdorff topological space. A {\it triangulation\/} of $X$ is data $(f_i^N:\De^N_\top\ra X)_{i\in I^N}$ for $N\ge 0$, where $I^N$ is an indexing set, such that $f_i^N$ is a homeomorphism with a closed subset $\Im f_i^N\subseteq X$, and for each $j=0,\ldots,N$ the restriction of $f_i^N$ to the face $x_j=0$ of $\De^N_\top$ is $f_{i'}^{N-1}$ for some unique $i'\in I^{N-1}$, that is,
\begin{equation*}
f_i^N(x_0,\ldots,x_{j-1},0,x_{j+1},\ldots,x_N)=f_{i'}^{N-1}(x_0,\ldots,x_{j-1},x_{j+1},\ldots,x_N),	
\end{equation*}
and $X=\coprod_{N\ge 0}\coprod_{i\in I^N}f_i^N((\De^n_\top)^\ci)$ is a stratification of $X$ into disjoint locally closed sets. For example, think of a topological surface divided into triangles.
\smallskip

\noindent{\bf(iii)} Let $S$ be a simplicial set. We define a compactly-generated, Hausdorff topological space $X_S$ with a triangulation, the {\it topological realization\/} of $S$, by
\begin{align*}
X_S=\bigl(\ts\coprod_{N\ge 0}S([N])\t\De^N_\top\bigr)\big/\sim,
\end{align*}
where the equivalence relation $\sim$ is generated by
\begin{equation*}
(\Si,f_\top(x_0,\ldots,x_M))\simeq (S(f)\Si,(x_0,\ldots,x_M)),
\end{equation*}
for all $f:[N]\ra[M]$ in $\De$, $\Si\in S([N])$ and $(x_0,\ldots,x_M)\in \De^M_\top$. The maps $f_i^N:\De^N_\top\ra X$ in the triangulation of $X$ come from the inclusions $\De^N_\top\ra \{\Si\}\t\De^N_\top$ for {\it nondegenerate\/} $N$-simplices $\Si\in S_N$. In this way we can define a {\it topological realization functor\/} $\TR:\SSet\ra\CGHaus\subset\Top$. It is a right Kan extension
\begin{equation*}
\xymatrix@C=130pt@R=15pt{ 
\De^\rop \ar[r]_(0.6){\text{Yoneda embedding}} \ar@/_.3pc/[dr] _G & \SSet \ar@{..>}[d]^\TR \\ & \CGHaus. }	
\end{equation*}

We use $\CGHaus$ rather than $\Top$ here as it has better categorical properties: both $\SSet$ and $\CGHaus$ are Cartesian closed, and $\TR$ is product-preserving as a functor to $\CGHaus$, but not as a functor to~$\Top$. 
\end{dfn}

\begin{ex}
\label{ai2ex1}
Let $\De^N$ be the simplicial set represented by the object $[N]$, so that by definition $\De^N:\De^\rop\ra\Set$ has $\De^N([M])=\Hom_\De([M],[N])$. By the Yoneda Lemma, if $S$ is a simplicial set the set of maps $\De^N\ra S$ (that is, morphisms of simplicial sets) is in natural bijection with $S_N$. The topological realization $\TR(\De^N)$ is naturally homeomorphic to $\De^N_\top$ in~\eq{ai2eq4}.
\end{ex}

\begin{ex}
\label{ai2ex2}
The next example will be important in \S\ref{ai22}. Let the simplicial set $\De^N$ be as in Example \ref{ai2ex1}, and let $0\le i\le N$. The $i^{\rm th}$ $N$-{\it horn\/} is the simplicial subset $\La^N_i\subset\De^N$ such that $\La^N_i([M])$ is the subset of order-preserving functions $f:[M]\ra[N]$ with $[N]\sm\{i\}\not\subseteq f([M])$. The topological realization $\TR(\La^N_i)$ is naturally homeomorphic to the subset of $\De^N_\top$ in \eq{ai2eq4} obtained by deleting the interior of $\De^N_\top$, and the interior of the $(N-1)$-dimensional boundary face $\De^{N-1}_{i,\top}$ opposite the $i^{\rm th}$ vertex. So for example $\TR(\La^2_i)$ consists of two sides of the triangle~$\De^2_\top$.
\end{ex}

\begin{ex}
\label{ai2ex3}
There is a good notion of {\it product\/} of simplicial sets. Considering simplicial sets as functors $S,T:\De^\rop\ra\Set$, the product $S\t T$ just takes products in $\Set$, so that $(S\t T)([N])=S([N])\t T([N])$. This satisfies the universal property of products in a category.

The picture of topological realizations of products in $\SSet$ is more interesting. The topological realization $\TR(S\t T)$ is canonically homeomorphic to the topological product $\TR(S)\t\TR(T)$, but the triangulation of $\TR(S\t T)$ is constructed from those of $\TR(S)$ and $\TR(T)$ by dividing products $\De^M_\top\t\De^N_\top$ canonically into $\binom{M+N}{M}$ copies of $\De^{M+N}_\top$.
\end{ex}

\begin{ex}
\label{ai2ex4}
There is a functor $\Sing:\Top\ra\SSet$ which maps a space $X$ to the functor $F:\De^\rop\ra\Set$ with $F([N])=\Map_{C^0}(\De^N_\top,X)$ for objects $[N]$ and $F(f)=-\ci f_\top:\Map_{C^0}(\De^N_\top,X)\ra\Map_{C^0}(\De^M_\top,X)$ for morphisms $f:[M]\ra[N]$, so that $f_\top:\De^M_\top\ra\De^N_\top$. It is right adjoint to~$\TR:\SSet\ra\Top$.

If $X$ is a topological space, there is a natural map $\eta_X:\TR\ci\Sing(X)\ra X$ coming from the adjunction. It turns out that $\eta_X$ is a {\it weak homotopy equivalence\/} (i.e.\ it induces isomorphisms $\pi_n(\TR\ci\Sing(X),x')\ra \pi_n(X,\eta_X(x'))$ on homotopy groups for all $n\ge 0$ and $x'\in \TR\ci\Sing(X)$). Thus by the Whitehead Theorem, if $X$ is a CW complex then $\eta_X$ is a homotopy equivalence.
\end{ex}

\begin{rem}
\label{ai2rem1}
One moral of Example \ref{ai2ex4} is that simplicial sets up to (a suitable notion of) homotopy equivalence are more-or-less the same as topological spaces (or nice topological spaces, such as CW complexes) up to homotopy equivalence.

Theories of $\iy$-categories have their roots in homotopy theory of topological spaces.  In an $\iy$-category (however one defines them, there are several approaches), one wants notions of homotopies and homotopy equivalences between morphisms, higher homotopies, and so on. Simplicial sets provide a model for homotopy theory of topological spaces which is very convenient for developing a theory of $\iy$-categories, as in Joyal \cite{Joya2} and Lurie~\cite{Luri2,Luri3,Luri5}.
\end{rem}

The next definition follows Lurie~\cite[\S 1.2.1]{Luri2}.

\begin{dfn}
\label{ai2def3}
Define a functor $\Si:\De\ra\De$ by $\Si([N])=[N]$ on objects, and if $f:[N]\ra[M]$ is a morphism in $\De$ then $\Si(f):[N]\ra[M]$ is $\Si(f):i\mapsto M-f(N-i)$. Then $\Si\ci\Si=\id_\De$, a strict equality of functors.

If $S$ is a simplicial set, define the {\it opposite simplicial set\/} $S^\rop$ by $S^\rop=S\ci\Si^\rop$ in functors $\De^\rop\ra\Set$, that is, $S^\rop$ is the composition $\De^\rop\,{\buildrel\Si^\rop\over\longra}\,\De^\rop\,{\buildrel Sp\over\longra}\,\Set$. Clearly $(S^\rop)^\rop=S$, a strict equality of simplicial sets. This extends to a functor $(-)^\rop:\SSet\ra\SSet$ by composing morphisms in $\SSet$, which are natural transformations of functors $\De^\rop\ra\Set$, with $\Si^\rop$.

If we think of a simplicial set $S$ as being built of $N$-simplices $\De^N$ with vertices numbered $0,1,\ldots,N$, then $S^\rop$ is built of the same $N$-simplices, but with the numbering of the vertices changed, so that vertex $i$ in $S$ gets numbered $N-i$ in $S^\rop$. In Definition \ref{ai2def1}, $S^\rop$ has the same sets $S_N^\rop=S_N$, and face and degeneracy maps $d^\rop_{N,i}=d_{N,N-i}$, $s^\rop_{N,i}=s_{N,N-i}$.
\end{dfn}

\subsection{\texorpdfstring{Kan complexes and $\iy$-categories}{Kan complexes and ∞-categories}}
\label{ai22}

Groth \cite{Grot} and Joyal \cite{Joya1} give surveys of $\iy$-categories, and Land \cite{Land} is an introductory book. They were introduced by Boardman and Vogt \cite{BoVo} under the name {\it weak Kan complexes}. Joyal \cite{Joya2} and Lurie's monumental works \cite{Luri2,Luri3,Luri5} develop the theory in detail. 

\begin{dfn}
\label{ai2def4}
Let the simplicial sets $\De^N,\La^N_i$ with inclusion morphism $\inc:\La^N_i\hookra\De^N$ be as in Examples \ref{ai2ex1}--\ref{ai2ex2}. A simplicial set $S$ is called a {\it Kan complex\/} if every map of simplicial sets $f:\La^N_i\ra S$ can be extended to an $N$-simplex $\ti f:\De^N\ra S$ for all $0\le i\le N$. That is, there should exist (not necessarily unique) $\ti f$ in the commutative diagram in~$\SSet$:
\e
\begin{gathered}
\xymatrix@C=100pt@R=15pt{ *+[r]{\La^N_{i_{\phantom{(}}}} \ar@{^{(}->}[d]^\inc \ar@/^.5pc/[dr]^f \\
*+[r]{\De^N} \ar@{..>}[r]^{\ti f} & S. }
\end{gathered}
\label{ai2eq5}
\e

A simplicial set $S$ is called an $\iy$-{\it category\/} (or {\it quasicategory}, or {\it inner Kan complex}, or {\it weak Kan complex\/}) if every map of simplicial sets $f:\La^N_i\ra S$ can be extended to an $N$-simplex $\ti f:\De^N\ra S$ for all $0<i<N$.

We will write $\Cat^1_\iy\subset\SSet$ for the full subcategory of $\iy$-categories in $\SSet$. This is an ordinary category. Note that it is different to the $\iy$-category of $\iy$-categories $\sCat_\iy$ discussed in \S\ref{ai24}. If $\sC,\sD$ are $\iy$-categories, a map of simplicial sets $F:\sC\ra\sD$ is called an $\iy$-{\it functor}, and these are the morphisms in~$\Cat^1_\iy$.
\end{dfn}

\begin{rem}
\label{ai2rem2}
{\bf(a)} There are actually many notions of infinity category, beginning with model categories in the 1960s. See Riehl \cite{Rieh} for a survey. Our use of `$\iy$-category' to mean weak Kan complex is becoming standard, but the term is also used with different meanings in the literature.
\smallskip

\noindent{\bf(b)} Definition \ref{ai2def4} begs the question: Why should we think of $\iy$-categories as some kind of category? Where are the objects, morphisms, $N$-morphisms, \ldots? The rest of this section will try to make this clearer.
\end{rem}

\begin{ex}
\label{ai2ex5}
Let $\cC$ be a small category. We will define a simplicial set $\sN(\cC)$ called the {\it nerve\/} of $\cC$. We define the $N$-simplices $\sN(\cC)_N$ for $N\ge 0$ to be the set of $N$-tuples of composable morphisms $C_0\,{\buildrel f_1\over\longra}\, C_1\,{\buildrel f_2\over\longra}\,\cdots \,{\buildrel f_N\over\longra}\,C_N$ in $\cC$, that is, the set of $(2N+1)$-tuples $(C_0,\ldots,C_N,f_1,\ldots,f_N)$ where $f_i:C_{i-1}\ra C_i$ is a morphism in $\cC$ for~$1\le i\le N$. 

The face map $d_{N,i}:\sN(\cC)_N\ra \sN(\cC)_{N-1}$ for $0<i<N$ composes $f_{i-1}$ and $f_i$, so that $d_{N,i}(C_0,\ldots,f_N)=(C_0,\ldots,C_{i-1},C_{i+1},\ldots,C_N,f_1,\ldots,f_i\ci f_{i-1},\ldots,f_N)$, and $d_{N,0}$ omits $C_0,f_1$, and $d_{N,N}$ omits $C_N,f_N$.

The degeneracy map $s_{N,i}:\sN(\cC)_N\ra \sN(\cC)_{N+1}$ for $0\le i\le N$ acts on $(C_0,\ldots,f_N)$ by inserting $\id_{C_i}:C_i\ra C_i$ between $f_i$ and~$f_{i+1}$.

Then $\sN(\cC)$ is a simplicial set, and an $\iy$-category. Also $\sN(\cC)$ is a Kan complex if and only if $\cC$ is a {\it groupoid\/} (i.e.\ all morphisms in $\cC$ are isomorphisms). We can recover the category $\cC$ from $\sN(\cC)$ (objects are 0-simplices, morphisms are 1-simplices, 2-simplices give composition, identities are given by~$s_{0,0}$). 

In fact $\sN(\cC)$ has the extra property that the lifts $\ti f$ of $f$ in \eq{ai2eq5} exist and are {\it unique\/} for all $0<i<N$ and $f$. A simplicial set $S$ is isomorphic to a nerve $\sN(\cC)$ if and only if it is an $\iy$-category for which lifts $\ti f$ have this uniqueness property. Thus, we can see $\iy$-categories (or Kan complexes) as generalizations of small categories (or groupoids) in which the lifts $\ti f$ in \eq{ai2eq5} need not be unique.

The nerve defines a functor $\sN:\Cat^1\ra\Cat^1_\iy$, where $\Cat^1$ is the ordinary category of (small) ordinary categories, and $\Cat^1_\iy$ the ordinary category of $\iy$-categories from Definition~\ref{ai2def4}.
\end{ex}

\begin{dfn}
\label{ai2def5}
Let $\sC$ be an $\iy$-category. An {\it object\/} $X$ of $\sC$ is a 0-simplex in $\sC$, or equivalently, a map of simplicial sets $X:\De^0\ra\sC$. A {\it morphism}, or 1-{\it morphism}, $f:X\ra Y$ of $\sC$ is a 1-simplex $f$ in $\sC$ with $d_{1,0}(f)=X$ and $d_{1,1}(f)=Y$, where $d_{N,i}$ are the degeneracy maps in $\sC$ as in \eq{ai2eq3}, or equivalently, a map of simplicial sets $f:\De^1\ra\sC$ with $f\vert_{\De^{\{0\}}}=X$ and~$f\vert_{\De^{\{1\}}}=Y$.

For an object $X$, the {\it identity morphism\/} is $\id_X=s_{0,0}(X)$, where $s_{N,i}$ are the face maps in $\sC$ as in~\eq{ai2eq3}.

Let $f:X\ra Y$ and $g:Y\ra Z$ be morphisms in $\sC$. A morphism $h:X\ra Z$ is a {\it composition\/} $g\ci f$ if there exists a 2-simplex $\eta$ in $\sC$ of the form
\e
\begin{gathered}
\xymatrix@C=70pt@R=15pt{
*+[l]{X\;\overset{0}{\bu}} \ar@<.8ex>[d]_f \ar@<.6ex>[dr]^{h=\,\text{`}g\ci f\text{'}} \ar@{}@<-1ex>[dr]_(0.3)\eta \\
*+[l]{Y\;\underset{1}{\bu}} \ar@<.7ex>[r]_g & *+[r]{\underset{2}{\bu}\; Z.} }
\end{gathered}
\label{ai2eq6}
\e
The horn-filling condition for $\La^2_1$ guarantees such $h$ exist, but {\it they are not unique}. That is, in an $\iy$-category composition of morphisms is not uniquely defined, but there can be many candidates $h$ for $g\ci f$. Because of this non-uniqueness, composition is not associative (though it is `up to homotopy').

For fixed objects $X,Y$, define a binary relation $\approx$ on morphisms $f,g:X\ra Y$ by $f\approx g$ 
if there exists a 2-simplex $\eta$ in $\sC$ of the form
\begin{equation*}
\xymatrix@C=70pt@R=15pt{
*+[l]{X\;\overset{0}{\bu}} \ar@<.8ex>[d]_{\id_X} \ar@<.6ex>[dr]^{f} \ar@{}@<-1ex>[dr]_(0.3)\eta \\
*+[l]{X\;\underset{1}{\bu}} \ar@<.7ex>[r]_g & *+[r]{\underset{2}{\bu}\; Y.} }
\end{equation*}
The $\iy$-category axioms imply $\approx$ is an equivalence relation (exercise), and also that $f\approx g$ if and only if there exists a 2-simplex $\eta$ in $\sC$ of the form
\begin{equation*}
\xymatrix@C=70pt@R=15pt{
*+[l]{X\;\overset{0}{\bu}} \ar@<.8ex>[d]_{g} \ar@<.6ex>[dr]^{f} \ar@{}@<-1ex>[dr]_(0.3)\eta \\
*+[l]{Y\;\underset{1}{\bu}} \ar@<.7ex>[r]_{\id_Y} & *+[r]{\underset{2}{\bu}\; Y.} }
\end{equation*}

As in \cite[\S 1.1]{Grot}, \cite[\S 1.2]{Land}, the {\it homotopy category\/} $\Ho(\sC)$ is the ordinary category whose objects are objects of $\sC$, and whose morphisms $[f]:X\ra Y$ are $\approx$-equivalence classes $[f]$ of morphisms $f:X\ra Y$ in $\sC$. Identity morphisms are $\id_X=[\id_X]:X\ra X$, and for morphisms $[f]:X\ra Y$, $[g]:Y\ra Z$, the composition is $[g]\ci[f]=[h]$ for $h$ as in \eq{ai2eq6}. This is well defined.

A morphism $f:X\ra Y$ in $\sC$ is called an {\it equivalence\/} if $[f]:X\ra Y$ is an isomorphism in~$\Ho(\sC)$.

If $\sC,\sD$ are $\iy$-categories then an $\iy$-functor $F:\sC\ra\sD$ induces an ordinary functor $\Ho(F):\Ho(\sC)\ra\Ho(\sD)$. As in Cisinski \cite[Def.~3.9.3]{Cisi} and Land \cite[Def.~2.3.3]{Land}, an $\iy$-functor $F:\sC\ra\sD$ is called {\it essentially surjective\/} if $\Ho(F):\Ho(\sC)\ra\Ho(\sD)$ is essentially surjective, that is, every object in $\Ho(\sD)$ is isomorphic to an object in the image of~$\Ho(F)$.

Let $\sC,\sD$ be $\iy$-categories, and $F,G:\sC\ra\sD$ be $\iy$-functors. An $\iy$-{\it natural transformation\/} $\eta:F\Ra G$ is a morphism of simplicial sets $\eta:\sC\t\De^1\ra\sC$ with $\eta\vert_{\sC\t\{0\}}=F$ and $\eta\vert_{\sC\t\{1\}}=G$. Then $\eta$ induces a natural transformation $\Ho(\eta):\Ho(F)\Ra\Ho(G)$ of ordinary functors $\Ho(\sC)\ra\Ho(\sD)$. We call $\eta$ an $\iy$-{\it natural isomorphism\/} if for each object $x$ in $\sC$, $\eta(\{x\}\t\De^1)$ is an equivalence in $\sD$. Then $\Ho(\eta):\Ho(F)\Ra\Ho(G)$ is a natural isomorphism.

If $\cC$ is an ordinary category, so that its nerve $\sN(\cC)$ is an $\iy$-category, there is a canonical strict isomorphism of categories $\Ho(\sN(\cC))\cong\cC$.

For any $\iy$-category $\sC$, there is a unique $\iy$-functor
\e
\Pi_{\Ho}:\sC\longra\sN(\Ho(\sC))
\label{ai2eq7}
\e
which maps $X\mapsto X$ on objects and $f\mapsto[f]$ on morphisms (1-simplices).
\end{dfn}

The next definition follows Lurie~\cite[\S 1.2.1]{Luri2}.

\begin{dfn}
\label{ai2def6}
Let $\sC$ be an $\iy$-category. The {\it opposite $\iy$-category\/} $\sC^\rop$ is the opposite simplicial set of $\sC$ in the sense of Definition \ref{ai2def3}. Then $\sC^\rop$ is an $\iy$-category, since taking opposites turns the horn-filling condition for $\La^N_i$ into the horn-filling condition for $\La^N_{N-i}$. We have~$(\sC^\rop)^\rop=\sC$.

If $f:X\ra Y$ is a morphism in $\sC$ then $f:Y\ra X$ is a morphism in $\sC^\rop$, as $f:\De^1\ra\sC$ corresponds to $f^\rop:(\De^1)^\rop\ra\sC^\rop$, and changing $\De^1\mapsto(\De^1)^\rop$ exchanges vertices 0,1 of $\De^1$. Thus, opposite $\iy$-categories are the analogues of opposite categories in ordinary category theory.

If $\fF:\sC\ra\sD$ is an $\iy$-functor then $\fF^\rop:\sC^\rop\ra\sD^\rop$ is an $\iy$-functor.
\end{dfn}

\begin{dfn}
\label{ai2def7}
An $\iy$-category $\sC$ is called an $\iy$-{\it groupoid\/} if all morphisms in $\sC$ are equivalences, \cite[Def.~1.2.24]{Land}. As in Land \cite[Lem.~1.2.29 \& Cor.~2.1.12]{Land}, it is a nontrivial theorem that $\iy$-groupoids are the same as Kan complexes.

If $\sC$ is an $\iy$-category, define the {\it maximal Kan complex\/} $\sC^\Kan$ in $\sC$ to be the simplicial subset of $\sC$ such that an $N$-simplex $S:\De^N\ra\sC$ in $\sC$ lies in $\sC^\Kan$ if and only if every 1-simplex in $S$ is an equivalence. Then $\sC^\Kan$ is an $\iy$-groupoid, and thus a Kan complex, and it contains every simplicial subset of $\sC$ which is a Kan complex, so it is the maximal Kan complex in $\sC$. 

If $F:\sC\ra\sD$ is an $\iy$-functor of $\iy$-categories then $F(\sC^\Kan)\subseteq\sD^\Kan\subseteq\sD$, and we write $F^\Kan=F\vert_{\sC^\Kan}:\sC^\Kan\ra\sD^\Kan$.
\end{dfn}

\subsection{\texorpdfstring{Mapping spaces in $\iy$-categories}{Mapping spaces in ∞-categories}}
\label{ai23}

As one model for $(\iy,1)$-categories $\sC$ is categories enhanced in topological spaces, one should expect mapping spaces $\Hom(X,Y)$ in $\sC$ to have the structure not just of sets, but of topological spaces up to homotopy equivalence. In the world of simplicial sets, topological spaces correspond to Kan complexes. So in $\iy$-categories (quasicategories), we expect mapping spaces to be Kan complexes, natural up to homotopy equivalence. There are several different, but homotopy-equivalent, definitions of the Kan complex $\Hom(X,Y)$. The next definition and proposition follow Lurie \cite[\S 1.2.2, Prop.~1.2.2.3 \& Cor.~4.2.1.8]{Luri2} and~\cite[\S 4.6]{Luri5}.

\begin{dfn}
\label{ai2def8}
Let $\sC$ be an $\iy$-category, and let $X,Y$ be objects (0-simplices) in $\sC$. Define a simplicial set $\Hom_\sC^L(X,Y)$, the space of {\it left morphisms\/} from $X$ to $Y$ in $\sC$, by taking $\Hom_{\SSet}(\De^N,\Hom_\sC^L(X,Y))$ to be the set of all morphisms $Z:\De^{N+1}\ra\sC$ such that $Z\vert_{\De^{\{0\}}}$ is the 0-simplex $X$, and $Z\vert_{\De^{\{1,2,\ldots,N+1\}}}$ is the constant $N$-simplex at $Y$. The face and degeneracy maps in $\Hom_\sC^L(X,Y)$ are induced from those in $\sC$ by extending an order-preserving function $f:[N]\ra[M]$ to an order-preserving function $\ti f:[N+1]\ra[M+1]$ with $\ti f(0)=0$ and $\ti f(i)=f(i-1)+1$ for~$i=1,\ldots,N+1$.

Similarly, define a simplicial set $\Hom_\sC^R(X,Y)$, the space of {\it right morphisms\/} from $X$ to $Y$ in $\sC$, by taking $\Hom_{\SSet}(\De^N,\Hom_\sC^R(X,Y))$ to be the set of all morphisms $Z:\De^{N+1}\ra\sC$ such that $Z\vert_{\De^{\{0,1,\ldots,N\}}}$ is a constant $N$-simplex at $X$, and $Z\vert_{\De^{\{N+1\}}}$ is the 0-simplex $Y$. The face and degeneracy maps in $\Hom_\sC^R(X,Y)$ are induced from those in $\sC$ by extending an order-preserving function $f:[N]\ra[M]$ to an order-preserving function $\ti f:[N+1]\ra[M+1]$ with $\ti f\vert_{[N]}=f$ and~$f(N+1)=M+1$.
\end{dfn}

\begin{prop}
\label{ai2prop1}
In Definition\/ {\rm\ref{ai2def8},} $\Hom_\sC^L(X,Y)$ and\/ $\Hom_\sC^R(X,Y)$ are Kan complexes, which are homotopy equivalent.
\end{prop}

When we only care about mapping spaces up to homotopy equivalence, rather than the particular simplicial set, we may write $\Hom_\sC(X,Y)$ rather than $\Hom_\sC^L(X,Y)$ or $\Hom_\sC^R(X,Y)$. The next definition follows Cisinski \cite[Def.~3.9.3]{Cisi} and Land \cite[Def.~2.3.2]{Land}.

\begin{dfn}
\label{ai2def9}
Let $F:\sC\ra\sD$ be an $\iy$-functor of $\iy$-categories, and $X,Y$ be objects in $\sC$. Then composition of morphisms $\De^{N+1}\ra\sC\,{\buildrel F\over\longra}\,\sD$ induces a morphism of Kan complexes $F_{X,Y}^L:\Hom_\sC^L(X,Y)\ra \Hom_\sD^L(F(X),F(Y))$. We call $F$ {\it fully faithful\/} if $F_{X,Y}^L$ is a homotopy equivalence for all $X,Y\in\sC$.
\end{dfn}

\subsection{\texorpdfstring{The $\iy$-category of $\iy$-categories}{The ∞-category of ∞-categories}}
\label{ai24}

The $\iy$-category $\sCat_\iy$ is discussed in Lurie \cite[\S 3]{Luri2} and Land~\cite[Def.~2.1.13]{Land}.

\begin{dfn}
\label{ai2def10}
In \cite[Def.~3.0.0.1]{Luri2}, Lurie defines the $\iy$-{\it category of\/ $\iy$-categories\/} $\sCat_\iy$. It is constructed as the simplicial nerve $\sN(\sCat^\De_\iy)$ of a simplicial category $\sCat^\De_\iy$ of $\iy$-categories. It has the properties:
\begin{itemize}
\setlength{\itemsep}{0pt}
\setlength{\parsep}{0pt}
\item[(a)] 0-simplices of $\sCat_\iy$ are (small) $\iy$-categories $\sC,\sD$, and 1-simplices of $\sCat_\iy$ are $\iy$-functors $F:\sC\ra\sD$.
\item[(b)] The ordinary category of $\iy$-categories $\Cat^1_\iy$ from Definition \ref{ai2def4} is a subcategory of $\sCat^\De_\iy$ (taking only 0-simplices in mapping spaces in $\sCat^\De_\iy$). The inclusion $\Cat^1_\iy\hookra\sCat^\De_\iy$ induces an inclusion of $\iy$-subcategories $I^{\sCat_\iy}:\sN(\Cat^1_\iy)\hookra\sCat_\iy$.
\item[(c)] For $\iy$-categories $\sC,\sD$, the right morphisms $\Hom_{\sCat_\iy}^R(\sC,\sD)$ may be identified with the maximal Kan complex $\Fun(\sC,\sD)^\Kan$ in the $\iy$-category $\Fun(\sC,\sD)$ of $\iy$-functors $F:\sC\ra\sD$, where $N$-simplices in $\Fun(\sC,\sD)$ correspond to maps of simplicial sets $\De^N\t\sC\ra\sD$.
\item[(d)] As in \cite[\S 3.3.3]{Luri2}, $\sCat^\De_\iy$ may be identified with the subcategory of fibrant-cofibrant objects in a simplicial model category ${\cal S}{\rm et}_\De^+$ of `marked simplicial sets' \cite[\S 3.1]{Luri2}. This implies nice properties of $\sCat_\iy$, for example, all small limits and colimits exist in~$\sCat_\iy$.
\end{itemize}
\end{dfn}

\begin{dfn}
\label{ai2def11}
Let $\sC,\sD$ be $\iy$-categories. An $\iy$-functor $F:\sC\ra\sD$ is called a {\it Joyal equivalence\/} \cite[Def.~2.1.14]{Land} if the corresponding morphism $[F]:\sC\ra\sD$ in $\Ho(\sCat_\iy)$ is an isomorphism. This is the $\iy$-category analogue of equivalence of categories.
\end{dfn}

Here is Cisinski \cite[Th.~3.9.7]{Cisi} and Land \cite[Th.~2.3.20]{Land}.

\begin{thm}
\label{ai2thm1}
Let\/ $F:\sC\ra\sD$ be an $\iy$-functor. Then $F$ is a Joyal equivalence if and only if it is fully faithful and essentially surjective, in the sense of Definitions\/ {\rm\ref{ai2def5}} and\/~{\rm\ref{ai2def9}}.
\end{thm}

In \S\ref{ai93} we will construct an $\iy$-functor $\ACat^{\un,\iy}_R\ra\sCat_\iy$ from our $\iy$-category of $A_\iy$-categories $\ACat^{\un,\iy}_R$. We do this using Lurie's `straightening-unstraightening equivalence' \cite[\S 3.2]{Luri2}, \cite[\S 3]{Land}, which for an $\iy$-category $\sD$ identifies $\iy$-functors $\sD\ra\sCat_\iy$ with `cocartesian fibrations'~$p:\sC\ra\sD$.

The next definition follows Lurie \cite[\S 2.3--\S 2.4]{Luri2} and Land~\cite[\S 1.3 \& \S 3.1]{Land}, see also~\cite[\S 4.1, \S 4.2 \& \S 5.1]{Luri5}.

\begin{dfn}
\label{ai2def12}
A map of simplicial sets $p:S\ra T$ is called an {\it inner fibration\/} if for all $0<i<N$, any commutative square of maps of simplicial sets `$\ra$' as follows can be completed with a map `$\dashra$' in a commutative diagram
\e
\begin{gathered}
\xymatrix@C=80pt@R=15pt{
*+[r]{\La^N_{i_{\phantom{(}}}} \ar[r] \ar@{^{(}->}[d] & *+[l]{S} \ar[d]_p \\
*+[r]{\De^N} \ar@{..>}[ur] \ar[r] & *+[l]{T.\!} }
\end{gathered}
\label{ai2eq8}
\e
If $T=\De^0$ is a point, this is equivalent to $S$ being an $\iy$-category. More generally it implies that the fibres of $p:S\ra T$ are $\iy$-categories. If $T$ is an $\iy$-category, then $S$ must be an $\iy$-category.

Similarly, we call $p:S\ra T$ a {\it left fibration\/} if we can complete \eq{ai2eq8} for all {\it left horns\/} $\La^N_i$ with $0\le i<N$, for $N\ge 1$, and a {\it right fibration\/} if we can complete \eq{ai2eq8} for all {\it right horns\/} $\La^N_i$ with $0<i\le N$, for~$N\ge 1$.

Let $p:S\ra T$ be an inner fibration. A morphism $f:\De^1\ra S$ is called $p$-{\it cocartesian\/} if for all $N\ge 2$ and all commutative diagrams of morphisms `$\ra$' as follows can be completed with a map `$\dashra$' in a commutative diagram
\e
\begin{gathered}
\xymatrix@C=80pt@R=15pt{
*+[r]{\De^1\,} \ar@<1ex>@/^.8pc/[rr]^f  \ar@{_{(}->}[r]_{\De^{\{0,1\}}} & \La^N_{0_{\phantom{(}}} \ar[r] \ar@{^{(}->}[d] & *+[l]{S} \ar[d]_p \\
& *+[r]{\De^N} \ar@{..>}[ur] \ar[r] & *+[l]{T.\!} }
\end{gathered}
\label{ai2eq9}
\e

Let $p:S\ra T$ be an inner fibration. We say that $p$ is a {\it cocartesian fibration\/} if any commutative square of maps of simplicial sets `$\ra$' as follows can be completed with a map `$\dashra$' in a commutative diagram
\e
\begin{gathered}
\xymatrix@C=80pt@R=15pt{
*+[r]{\De^0_{\phantom{(}}} \ar[r] \ar@{^{(}->}[d]_{\De^{\{0\}}} & *+[l]{S} \ar[d]_p \\
*+[r]{\De^1} \ar@{..>}[ur]^{\ti f} \ar[r]^f & *+[l]{T,\!} }
\end{gathered}
\label{ai2eq10}
\e
such that $\ti f$ is $p$-cocartesian.

If $p$ is a cocartesian fibration and we are given an outer commutative square \eq{ai2eq10}, then the space of $p$-cocartesian lifts $\ti f$ making \eq{ai2eq10} commute is a contractible Kan complex. That is, $\ti f$ is unique up to homotopy in a strong sense.
\end{dfn}

By Lurie \cite[Prop.~5.1.4.15]{Luri5} and Land \cite[Lem.~3.1.22]{Land}, we have

\begin{prop}
\label{ai2prop2}
A map of simplicial sets\/ $p:S\ra T$ is a left fibration if and only if it is cocartesian and every edge of\/ $S$ is\/ $p$-cocartesian.
\end{prop}

The next theorem comes from Lurie \cite[\S 3.2--\S 3.3]{Luri2} and Land~\cite[\S 3.3]{Land}.

\begin{thm}
\label{ai2thm2}
There is a \begin{bfseries}universal cocartesian fibration\end{bfseries} $u:\sF\ra\sCat_\iy$ with the following properties:
\begin{itemize}
\setlength{\itemsep}{0pt}
\setlength{\parsep}{0pt}
\item[{\bf(a)}] Every $\iy$-category $\sC$ is an object in $\sCat_\iy,$ and so corresponds to a morphism $\pi_\sC:\De^0\ra\sCat_\iy$. There is a homotopy Cartesian square of\/ $\iy$-categories
\begin{equation*}
\xymatrix@C=90pt@R=15pt{
*+[r]{\sC}  \drtwocell_{}\omit^{}\omit{^{}} \ar[r]_(0.7){\Phi_\sC} \ar[d] & *+[l]{\sF} \ar[d]_u \\
*+[r]{\De^0} \ar[r]^(0.3){\pi_\sC} & *+[l]{\sCat_\iy.\!} }
\end{equation*}
That is, the fibre $\sF_\sC$ of\/ $u:\sF\ra\sCat_\iy$ over each object\/ $\sC$ in $\sCat_\iy$ is (Joyal equivalent to) $\sC$. 
\item[{\bf(b)}] Let\/ $p:\sC\ra\sD$ be a cocartesian fibration of\/ $\iy$-categories. Then there is an $\iy$-functor $f:\sD\ra\sCat_\iy,$ unique up to a contractible choice, which fits into a homotopy Cartesian square of\/ $\iy$-categories
\e
\begin{gathered}
\xymatrix@C=90pt@R=15pt{
*+[r]{\sC}  \drtwocell_{}\omit^{}\omit{^{}} \ar[r] \ar[d]^p & *+[l]{\sF} \ar[d]_u \\
*+[r]{\sD} \ar[r]^(0.35)f & *+[l]{\sCat_\iy.\!} }
\end{gathered}
\label{ai2eq11}
\e
We call\/ $f$ the \begin{bfseries}straightening\end{bfseries} of\/~$p$.
\item[{\bf(c)}] Let\/ $f:\sD\ra\sCat_\iy$ be an $\iy$-functor of\/ $\iy$-categories. Then there is a cocartesian fibration $p:\sC\ra\sD,$ unique up to a contractible choice, which fits into a homotopy Cartesian square \eq{ai2eq11}. We call\/ $p$ the \begin{bfseries}unstraightening\end{bfseries} of\/~$f$.
\end{itemize}
In\/ {\bf(b)\rm,\bf(c)\rm,} part\/ {\bf(a)} implies that for all objects $X\in\sD,$ the image $f(X)\in\sCat_\iy,$ considered as an $\iy$-category, is Joyal equivalent to the fibre $p^{-1}(X)$.
\end{thm}

The next example summarizes parts of Lurie \cite[\S 3.2.5]{Luri2}.

\begin{ex}
\label{ai2ex6}
Let $\cD$ be an ordinary category, and $G:\cD\ra\Cat^1_\iy$ a functor, where $\Cat^1_\iy\subset\SSet$ is the ordinary category of $\iy$-categories from Definition \ref{ai2def4}. Definition \ref{ai2def10}(b) gives an $\iy$-functor $I^{\sCat_\iy}:\sN(\Cat^1_\iy)\hookra\sCat_\iy$. So we can consider the composition $\iy$-functor $I^{\sCat_\iy}\ci\sN(G):\sN(\cD)\ra\sCat_\iy$. Theorem \ref{ai2thm2} then gives a cocartesian fibration $p:\sC\ra\sN(\cD)$, the {\it unstraightening\/} of $I^{\sCat_\iy}\ci\sN(G)$, uniquely up to a contractible choice.

Lurie \cite[Def.~3.2.5.2]{Luri2} writes down an explicit simplicial set $N_G(\cD)$ and an $\iy$-functor $\pi:N_G(\cD)\ra\sN(\cD)$, the {\it relative simplicial nerve\/} of $G$, such that $\pi:N_G(\cD)\ra\sN(\cD)$ is a possible choice for the unstraightening~$p:\sC\ra\sN(\cD)$.

If $N\ge 0$, an $N$-simplex $(S,(\tau(I))_{\es\ne I\subseteq[N]})$ of $N_G(\cD)$ consists of the data:
\begin{itemize}
\setlength{\itemsep}{0pt}
\setlength{\parsep}{0pt}
\item[(i)] an $N$-simplex $S$ of $\sN(\cD)$. Here $S$ is equivalent to objects $X_0,\ldots,X_N$ in $\cD$ and morphisms $f^{ij}:X_i\ra X_j$ in $\cD$ for $0\le i\le j\le N$, with $f^{ii}=\id_{X_i}$ and $f^{jk}\ci f^{ij}=f^{ik}$ for $0\le i\le j\le k\le N$.
\item[(ii)] For all $\es\ne I\subseteq\{0,\ldots,N\}$, a map of simplicial sets $\tau(I):\De^I\ra G(X_{\max I})$, where $\De^I$ is the simplex with ordered vertices $(I,\le)$, considered as a simplicial set, so that~$\De^I\cong\De^{\md{I}-1}$.
\item[(iii)] For all $\es\ne J\subseteq I\subseteq\{0,\ldots,N\}$, the following diagram of simplicial sets must commute:
\e
\begin{gathered}
\xymatrix@C=150pt@R=15pt{
*+[r]{\De^J} \ar[r]_{\tau(J)} \ar[d]^{\De^{\inc(J\hookra I)}} & *+[l]{G(X_{\max J})} \ar[d]_{G(f^{(\max J)(\max I)})} \\
*+[r]{\De^I} \ar[r]^{\tau(I)} & *+[l]{G(X_{\max I}),\!}
}
\end{gathered}
\label{ai2eq12}
\e
with $\De^{\inc(J\hookra I)}:\De^J\ra\De^I$ the inclusion of simplices induced by~$J\hookra I$.
\end{itemize}
These $N$-simplices have the obvious functoriality under morphisms $\al:[\bar N]\ra[N]$ in $\De^\rop$, and so define a simplicial set $N_G(\cD)$. The morphism $\pi:N_G(\cD)\ra\sN(\cD)$ acts on $N$-simplices by $\pi:(S,(\tau(I))_{\es\ne I\subseteq[N]})\mapsto S$.

If $X$ is an object in $\cD$, and thus a 0-simplex in $\sN(\cD)$, consider the fibre $\pi^{-1}(X)\subseteq N_G(\cD)$. An $N$-simplex $(S,(\tau(I))_{\es\ne I\subseteq[N]})$ in $\pi^{-1}(X)$ has $S$ the constant $N$-simplex at $X$, and the maps $\tau(I):\De^I\ra G(X)$ for $I\subseteq[N]$ are uniquely determined by the $N$-simplex $\tau([N])$ in $G(X)$. Hence there is a canonical strict isomorphism of simplicial sets~$\pi^{-1}(X)\cong G(X)$.
\end{ex}

\subsection{\texorpdfstring{The $\iy$-category of spaces, mapping space functors, and Yoneda embedding}{The ∞-category of spaces, mapping space functors, and Yoneda embedding}}
\label{ai25}

\begin{dfn}
\label{ai2def13}
We will write $\sSpc$ for the $\iy$-{\it category of spaces}. Objects of $\sSpc$ are `spaces' (roughly, topological spaces up to homotopy), or $\iy$-groupoids. It is a complete and cocomplete $\iy$-category.

There are several different ways to define $\sSpc$ in the literature, which yield equivalent $\iy$-categories. We follow Lurie \cite[\S 1.2.16]{Luri2}, \cite[\S 3.1.6]{Luri5} and Land \cite[Def.~1.3.41]{Land}, and define $\sSpc$ to be the simplicial nerve $\sN(\bKan^\De)$ of the simplicial category $\bKan^\De$ of Kan complexes, which is the full simplicial subcategory with objects Kan complexes in the simplicial category $\mathop{\bf SSet^\De}$ of simplicial sets. Then $\sSpc$ is a full $\iy$-subcategory of $\sCat_\iy$ in Definition~\ref{ai2def10}.

Thus, objects of $\sSpc$ are Kan complexes. A morphism $f:X\ra Y$ in $\sSpc$ is an equivalence in $\sSpc$ if and only if it is a homotopy equivalence of Kan complexes.

The ordinary category $\bKan$ of Kan complexes, a full subcategory $\bKan\subset\Cat^1_\iy\subset\SSet$, is a subcategory of $\bKan^\De$ (taking only 0-simplices in mapping spaces in $\bKan^\De$). Thus as for $I^{\sCat_\iy}$ in Definition \ref{ai2def10}(b), the inclusion $\bKan\hookra\bKan^\De$ induces an inclusion of $\iy$-subcategories $I^\bKan:\sN(\bKan)\hookra\sSpc$.

Land \cite[Def.~1.2.72]{Land} defines an $\iy$-category with objects CW-complexes (nice topological spaces), and shows \cite[Obs.~1.3.43]{Land} it is equivalent to $\sSpc$. Cisinski \cite[Def.~5.2.3]{Cisi} gives an alternative definition, yielding an equivalent $\iy$-category~\cite[Rem.s 5.4.11 \& 7.8.10]{Cisi}.
\end{dfn}

\begin{dfn}
\label{ai2def14}
Let $\sC$ be an $\iy$-category. Then there is an $\iy$-functor $\bHom_\sC:\sC^\rop\t\sC\ra\sSpc$, with a difficult definition in Lurie \cite[\S 5.2.1, Prop.~5.2.11]{Luri3} (see also \cite[\S 5.1.3]{Luri2}), Cisinski \cite[\S 5.8.1]{Cisi}, and Land \cite[Def.~4.2.5]{Land}, such that for objects $X,Y\in\sC$ we have a homotopy equivalence $\bHom_\sC(X,Y)\simeq \Hom_\sC^L(X,Y)$ (i.e.\ an equivalence in $\sSpc$), where $\Hom_\sC^L(X,Y)$ is the Kan complex in Definition \ref{ai2def8}. That is, $\bHom_\sC$ makes mapping spaces $\Hom_\sC^L(X,Y)$ (or just $\Hom_\sC(X,Y)$) in \S\ref{ai23} functorial in~$X,Y$.
\end{dfn}

\begin{dfn}
\label{ai2def15}
Let $\sC$ be an $\iy$-category. As in Lurie \cite[\S 5.1]{Luri2} define the $\iy$-category $\PSh(\sC)$ of {\it presheaves\/} on $\sC$ to be the $\iy$-category $\Fun(\sC^\rop,\sSpc)$ of $\iy$-functors $G:\sC^\rop\ra\sSpc$. Then $\PSh(\sC)$ is complete and cocomplete \cite[Cor.~5.1.2.4]{Luri2}, and (co)limits in $\PSh(\sC)$ may be computed objectwise in~$\sC^\rop$.

If $F:\sC\ra\sD$ is an $\iy$-functor there is a pullback $\iy$-functor $F^*:\PSh(\sD)\ra\PSh(\sC)$ mapping $G\mapsto G\ci F$. This $F^*$ has left and right adjoints, given by left and right Kan extension. We denote the left adjoint by $F_!:\PSh(\sC)\ra\PSh(\sD)$.
\end{dfn}

\begin{dfn}
\label{ai2def16}
The {\it Yoneda functor\/} $\sY_\sC:\sC\ra\PSh(\sC)$ maps $\sY_\sC:X\mapsto\bHom_\sC(-,X)$. As in Lurie \cite[\S 5.1.3]{Luri2}, Cisinski \cite[Th.~5.8.13]{Cisi}, and Land \cite[Prop.~4.2.11]{Land}, this is a full and faithful $\iy$-functor, so we call it the {\it Yoneda embedding}. Objects of $\PSh(\sC)$ in the essential image of $\sY_\sC$ are called {\it representable}. Write $\Rep(\sC)$ for the full $\iy$-subcategory of representable objects in $\PSh(\sC)$. Then $\sY_\sC:\sC\ra\Rep(\sC)$ is an equivalence of $\iy$-categories. 

If $F:\sC\ra\sD$ is an $\iy$-functor then by Lurie \cite[Prop.~5.2.6.3]{Luri2} the following commutes up to natural equivalence:
\e
\begin{gathered}
\xymatrix@C=100pt@R=15pt{ *+[r]{\sC} \drtwocell_{}\omit^{}\omit{^{}} \ar[r]_(0.33){\sY_\sC} \ar[d]^F & *+[l]{\PSh(\sC)} \ar[d]_{F_!} \\
*+[r]{\sD} \ar[r]^(0.6){\sY_\sD} & *+[l]{\PSh(\sD).\!} }
\end{gathered}
\label{ai2eq13}
\e
\end{dfn}

In Theorem \ref{ai2thm2} we explained Lurie's `straightening/unstraightening correspondence' between cocartesian fibrations of $\iy$-categories $p:\sC\ra\sD$ and $\iy$-functors $f:\sD\ra\sCat_\iy$. As above $\sSpc$ is a full $\iy$-subcategory of $\sCat_\iy$. Definition \ref{ai2def12} defined left fibrations $p:\sC\ra\sD$, which are special examples of cocartesian fibrations by Proposition \ref{ai2prop2}. If $p$ is a cocartesian fibration then the fibres $p^{-1}(X)$ for $X\in\sD$  are $\iy$-categories, but if $p$ is a left fibration then the $p^{-1}(X)$ are Kan complexes. Here is Land \cite[Th.~3.3.16]{Land}.

\begin{thm}
\label{ai2thm3}
The straightening/unstraightening correspondence in Theorem\/ {\rm\ref{ai2thm2}(b),(c)} restricts to a correspondence between \begin{bfseries}left fibrations\end{bfseries} $p:\sC\ra\sD$ and\/ $\iy$-functors $f:\sD\ra\sSpc\subset\sCat_\iy$.
\end{thm}

\begin{ex}
\label{ai2ex7}
Let $\sC$ be an $\iy$-category, and let $X$ be an object (0-simplex) in $\sC$. Define a simplicial set $\sC_{X/}$ by taking $\Hom_{\SSet}(\De^N,\sC_{X/})$ to be the set of all morphisms $Z:\De^{N+1}\ra\sC$ such that $Z\vert_{\De^{\{0\}}}$ is the the 0-simplex $X$. The face and degeneracy maps in $\Hom_\sC^L(X,Y)$ are induced from those in $\sC$ by extending an order-preserving function $f:[N]\ra[M]$ to an order-preserving function $\ti f:[N+1]\ra[M+1]$ with $\ti f(0)=0$ and $\ti f(i)=f(i-1)+1$ for $i=1,\ldots,N+1$. Define a map $p_\sC:\sC_{X/}\ra\sC$ to map an $N$-simplex $Z$ in $\sC_{X/}$, thought of as an $(N+1)$-simplex $Z:\De^{N+1}\ra\sC$ with $Z\ci\De^{\{0\}}=X$, to the $N$-simplex $Z\ci\De^{\{1,2,\ldots,N+1\}}$ in $\sC$.

By Lurie \cite[Cor.~2.1.2.2]{Luri2}, \cite[Prop.~4.3.6.1]{Luri5}, $\sC_{X/}$ is an $\iy$-category, called the {\it undercategory of\/} $X$, and $p_\sC:\sC_{X/}\ra\sC$ is a left fibration. From Definition \ref{ai2def8}, we see that the fibre of $p_\sC$ over $Y\in\sC$ is the left morphisms $\Hom^L_\sC(X,Y)$, a Kan complex. Theorem \ref{ai2thm3} says that $p_\sC$ corresponds to an $\iy$-functor $F_X:\sC\ra\sSpc$, the straightening of $p_\sC$, such that $F_X(Y)$ is homotopy equivalent to $\Hom^L_\sC(X,Y)$. As observed by Lurie \cite[\S 2.1.2]{Luri2}, this is a way to make mapping spaces $\Hom_\sC(X,-)$ into an $\iy$-functor $\sC\ra\sSpc$.
\end{ex}

\subsection{\texorpdfstring{Defining $\iy$-categories via localization}{Defining ∞-categories via localization}}
\label{ai26}

There are many situations in mathematics where one is given a category $\cC$ and a class of morphisms $\cW$ in $\cC$ (usually called {\it weak equivalences\/}), and one wants to define a {\it localized category\/} $\cC[\cW^{-1}]$ by inverting the morphisms in $\cW$. That is, we should have a functor $\Pi:\cC\ra\cC[\cW^{-1}]$ which maps morphisms in $\cW$ to isomorphisms, and this should be universal amongst functors $\Pi:\cC\ra\cD$ mapping $\cW$ to isomorphisms.

For example, the derived category $D^b\cA$ of an abelian category $\cA$ is obtained by inverting quasi-isomorphisms in a category of complexes $\Com^b(\sA)$. Homotopy theory is the study of the `homotopy category' obtained by localizing topological spaces $\Top$ along weak homotopy equivalences.

It turns out that there is a natural $\iy$-category $\cC[\cW^{-1}]^\iy$, whose homotopy category $\Ho(\cC[\cW^{-1}]^\iy)$ is $\cC[\cW^{-1}]$. Often it is helpful to work in $\cC[\cW^{-1}]^\iy$ rather than in $\cC[\cW^{-1}]$. For example, derived categories $D^b\cA$ suffer from `nonfunctoriality of the cone', as in Gelfand and Manin \cite[\S IV.1.7]{GeMa}, which is repaired at the $\iy$-category level. This is an important reason for studying $\iy$-categories.

We will discuss localizations of ordinary categories with weak equivalences. In fact localization works more generally (one can localize any $\iy$-category at any set of morphisms \cite[Rem.~1.3.4.2]{Luri3}, \cite[\S 6.3.2]{Luri5}, \cite[Prop.~2.4.8]{Land}), but this is all we need. See Lurie \cite[\S 1.3.4]{Luri3}, \cite[\S 5.2.7]{Luri2} and Land~\cite[\S 2.4]{Land}.

\begin{dfn}
\label{ai2def17}
A {\it category with weak equivalences\/} $(\cC,\cW)$ is an ordinary category $\cC$ and a subset $\cW\subseteq\Mor(\cC)$ of the morphisms of $\cC$, satisfying
\begin{itemize}
\setlength{\itemsep}{0pt}
\setlength{\parsep}{0pt}
\item[(a)] All isomorphisms in $\cC$ lie in $\cW$. 
\item[(b)] (The `two out of three' property.) If $f:X\ra Y$, $g:Y\ra Z$ are morphisms in $\cC$ and two out of three of $f,g,g\ci f$ lie in $\cW$, then the third lies in~$\cW$.
\end{itemize}

Let $(\cC,\cW)$ be a category with weak equivalences. A {\it Dwyer--Kan localization of\/ $\cC$ along\/} $\cW$ is an $\iy$-category $\sL$ and an $\iy$-functor $\la:\sN(\cC)\ra\sL$, such that if $f\in\cW$ then $\la(f)$ is an equivalence in $\sL$ (that is, on homotopy categories $\Ho(\la):\cC\cong\Ho(\sN(\cC))\ra\Ho(\sL)$ maps morphisms in $\cW$ to isomorphisms), and $\la$ has the universal property that whenever $\fF:\sN(\cC)\ra\sD$ is an $\iy$-functor which maps $\cW$ to equivalences in $\sD$, then we can factor $\fF$ as 
\e
\begin{gathered}
\xymatrix@C=50pt@R=15pt{
*+[r]{\sN(\cC)} \ar[rr]_(0.6)\la \ar@/_.5pc/[drr]_\fF &  \drtwocell_{}\omit^{}\omit{^{}} & *+[l]{\sL} \ar@{..>}[d]^\fG \\
&& *+[l]{\sD } }
\end{gathered}
\label{ai2eq14}
\e
for some $\fG:\sL\ra\sD$ unique up to a contractible choice. (See \cite[\S 2.4]{Land} for a precise statement.)

Dwyer--Kan localizations $\sL$ always exist, and are unique up to Joyal equivalence. There is an explicit construction, by first taking the `Dwyer--Kan hammock localization' to give a simplicial category, and then taking the simplicial nerve of this. We will use the notation $\cC[\cW^{-1}]^\iy$ for this explicit localization, and $\Pi_\loc^\iy:\sN(\cC)\ra\cC[\cW^{-1}]^\iy$ for the localization $\iy$-functor.

As in Definition \ref{ai2def15} we have a pullback $\iy$-functor
\begin{equation*}
(\Pi_\loc^\iy)^*:\PSh(\cC[\cW^{-1}]^\iy)\longra\PSh(\sN(\cC)).
\end{equation*}
The universal property of $\cC[\cW^{-1}]^\iy$ implies that $(\Pi_\loc^\iy)^*$ is full and faithful, with essential image those $\iy$-functors $\sN(\cC)^\rop\ra\sSpc$ which map morphisms in $\cW^\rop$ to equivalences in $\sSpc$. Therefore $\cC[\cW^{-1}]^\iy$ is equivalent to a full $\iy$-subcategory of $\PSh(\sN(\cC))$, the essential image of $\Rep(\cC[\cW^{-1}]^\iy)$ under~$(\Pi_\loc^\iy)^*$.
\end{dfn}

\begin{rem}
\label{ai2rem3}
In Definition \ref{ai2def17}, we can also consider the localization $\cC[\cW^{-1}]$ of $\cC$ along $\cW$ just in ordinary categories, not in $\iy$-categories. It turns out that there is an equivalence of categories $\cC[\cW^{-1}]\simeq\Ho(\cC[\cW^{-1}]^\iy)$, canonical up to natural isomorphism. To see why, suppose that $f:\cC\ra\cD$ is a functor of ordinary categories which maps $\cW$ to isomorphisms in $\cD$. Then $f$ should factor via the localization $\cC[\cW^{-1}]$ as in
\begin{equation*}
\xymatrix@C=50pt@R=15pt{
*+[r]{\cC} \ar[rr]_(0.4){\Pi_\loc} \ar@/_.5pc/[drr]_f &  \drtwocell_{}\omit^{}\omit{^{}} & *+[l]{\cC[\cW^{-1}]} \ar@{..>}[d]^g \\
&& *+[l]{\cD, } }
\end{equation*}
where $g$ is natural up to natural isomorphism.

Now consider the homotopy commutative $\iy$-category diagram
\begin{equation*}
\xymatrix@C=37pt@R=15pt{
*+[r]{\sN(\cC)} \ar@<1ex>@/^.8pc/[rrrr]^(0.1){\sN(\la)} \ar[rr]_(0.45){\Pi_\loc^\iy} \ar@/_.5pc/[drr]_{\sN(f)} &  \drtwocell_{}\omit^{}\omit{^{}} & {\cC[\cW^{-1}]^\iy} \ar@{..>}[d]^\fG \ar[rr]_(0.3){\Pi_{\Ho}} && *+[l]{\sN(\Ho(\cC[\cW^{-1}]^\iy))} \ar@{..>}[dll]^{\sN(g)} \\
&& {\sN(\cD),} }
\end{equation*}
where $\Pi_{\Ho}$ is as in \eq{ai2eq7}. Here $\sN(f)$ factors via $\fG$, uniquely up to a contractible choice, by the universal property of $\cC[\cW^{-1}]^\iy$. But as the target of $\fG$ is $\sN(\cD)$, $\fG$ factors via $\Pi_{\Ho}$, and any $\iy$-functor $\sN(\Ho(\cC[\cW^{-1}]^\iy))\ra\sN(\cD)$ is $\sN(g)$ for some functor $g:\Ho(\cC[\cW^{-1}]^\iy)\ra\cD$. It follows that $\Ho(\cC[\cW^{-1}]^\iy)$ satisfies the universal property for $\cC[\cW^{-1}]$, and~$\cC[\cW^{-1}]\simeq\Ho(\cC[\cW^{-1}]^\iy)$.

Thus, we can consider $\cC[\cW^{-1}]^\iy$ as an $\iy$-categorical enhancement of the ordinary category localization~$\cC[\cW^{-1}]$.
\end{rem}

The next proposition, which is easy to prove, is used in \cite{Pasc,COS2} to prove Theorem \ref{ai4thm8} below.

\begin{prop}
\label{ai2prop3}
Suppose $(\cC,\cW)$ and\/ $(\cC',\cW')$ are categories with weak equivalences, and we are given adjoint functors
\begin{equation*}
\xymatrix@C=80pt{ \cC \ar@<.5ex>[r]^F & \cC' \ar@<.5ex>[l]^G }
\end{equation*}
satisfying 
\begin{itemize}
\setlength{\itemsep}{0pt}
\setlength{\parsep}{0pt}
\item[{\bf(i)}] $F(\cW)\subseteq\cW'$ and\/ $G(\cW')\subseteq\cW$.
\item[{\bf(ii)}] Write $\ze:\Id_\cC\Ra G\ci F$ and\/ $\eta:F\ci G\Ra\Id_\cC$ for the unit and counit of the adjunction. Then $\ze(X)\in\cW$ and\/ $\eta(X')\in\cW'$ for all\/ $X\in\cC$ and\/~$X'\in\cC'$.
\end{itemize}
Then $F,G$ induce quasi-inverse Joyal equivalences of localized\/ $\iy$-categories:
\begin{equation*}
\xymatrix@C=80pt{ \cC[\cW^{-1}]^\iy \ar@<.5ex>[r]^{F_*} & \cC'[\cW^{\prime -1}]^\iy \ar@<.5ex>[l]^{G_*} }
\end{equation*}
\end{prop}

The next theorem follows from Land \cite[Lem.~2.4.9 \& Th.~2.4.10]{Land}. It describes $\sCat_\iy,\sSpc$ in Definitions \ref{ai2def10} and \ref{ai2def13} as localizations.

\begin{thm}
\label{ai2thm4}
{\bf(a)} Let\/ $I^{\sCat_\iy}:\sN(\Cat^1_\iy)\hookra\sCat_\iy$ be as in Definition\/ {\rm\ref{ai2def10}(b)}. Then $I^{\sCat_\iy}$ maps Joyal equivalences to equivalences in $\sCat_\iy$ by Definition\/ {\rm\ref{ai2def11}}. Hence as in \eq{ai2eq14} we have a canonical factorization
\e
\begin{gathered}
\xymatrix@C=50pt@R=15pt{
*+[r]{\sN(\Cat^1_\iy)} \ar[rr]_(0.45){\Pi_\loc^\iy} \ar@/_.5pc/[drr]_{I^{\sCat_\iy}} &  \drtwocell_{}\omit^{}\omit{^{}} & *+[l]{\Cat^1_\iy[{\rm Joy}^{-1}]^\iy} \ar@{..>}[d]^\Pi \\
&& *+[l]{\sCat_\iy,\! } }
\end{gathered}
\label{ai2eq15}
\e
where ${\rm Joy}\subset\Mor(\Cat^1_\iy)$ is the set of Joyal equivalences. Then $\Pi$ is a Joyal equivalence of\/~$\iy$-categories.
\smallskip

\noindent{\bf(b)} The category of simplicial sets $\SSet$ has the \begin{bfseries}Quillen model category structure\end{bfseries}, in which the weak equivalences are weak homotopy equivalences (morphisms of simplicial sets which induce weak homotopy equivalences of the topological realizations). All objects are cofibrant, and fibrant objects are Kan complexes. Write $\bKan$ for the full subcategory of Kan complexes in $\SSet$. Then weak homotopy equivalences in $\bKan$ are homotopy equivalences. 

We have Joyal equivalences of\/ $\iy$-categories
\e
\SSet[{\cal{WE}}^{-1}]^\iy\simeq \bKan[{\cal{HE}}^{-1}]^\iy\simeq\sSpc,
\label{ai2eq16}
\e
where $\cal{WE}$ are the weak homotopy equivalences in $\SSet,$ and\/ ${\cal{HE}}$ are the homotopy equivalences in $\bKan$. Thus as for \eq{ai2eq15} we have a canonical factorization
\e
\begin{gathered}
\xymatrix@C=50pt@R=15pt{
*+[r]{\sN(\bKan)} \ar[rr]_(0.45){\Pi_\loc^\iy} \ar@/_.5pc/[drr]_{I^\bKan} &  \drtwocell_{}\omit^{}\omit{^{}} & *+[l]{\bKan[{\cal{HE}}^{-1}]^\iy} \ar@{..>}[d]^\Pi \\
&& *+[l]{\sSpc,\! } }
\end{gathered}
\label{ai2eq17}
\e
where $I^\bKan$ is as in Definition\/ {\rm\ref{ai2def13},} and\/ $\Pi$ is a Joyal equivalence.
\end{thm}

In several important situations there are tools for describing the morphisms in $\cC[\cW^{-1}]$, and the homotopy types of the mapping spaces in $\cC[\cW^{-1}]^\iy$. One such situation is the theory of model categories (or simplicial model categories).

\begin{rem}
\label{ai2rem4}
In this remark we use terms from the theory of model categories with little explanation. See Hirschhorn \cite{Hirs} and Hovey \cite{Hove} for books on model categories, and Dwyer--Spalinski \cite{DwSp} and Hess \cite{Hess} for surveys.

A {\it model category\/} $\cC$ is a category equipped with distinguished classes of morphisms: {\it weak equivalences\/} $\cW$, {\it fibrations\/} ${\mathcal F}\text{\it ib}$, and {\it cofibrations\/} ${\mathcal C}\text{\it ofib}$, satisfying some axioms. The main purpose of the theory of model categories is to study the `homotopy category' $\cC[\cW^{-1}]$ and its $\iy$-categorical enhancement $\cC[\cW^{-1}]^\iy$. 

A central example is when $\cC$ is the category of topological spaces $\Top$, and $\cW$ is weak homotopy equivalences, and ${\mathcal F}\text{\it ib}$ is Serre fibrations. Then, following Quillen in the 1960s, much of classical homotopy theory can be rewritten in this language, and so generalized to any model category.

An object $X$ in a model category $\cC$ is called {\it fibrant\/} if $\pi:X\ra *$ is a fibration, where $*$ is the terminal object in $\cC$, and {\it cofibrant\/} if $i:\es\ra X$ is a cofibration, where $\es$ is the initial object in $\cC$. A {\it fibrant replacement\/} for $X$ is a weak equivalence $i:X\ra X'$ with $X'$ fibrant. A {\it cofibrant replacement\/} for $X$ is a weak equivalence $i:X'\ra X$ with $X'$ cofibrant.
\smallskip

\noindent{\bf(a)} Suppose $X,Y$ are objects of $\cC$, and $i:X'\ra X$ is a cofibrant replacement, and $j:Y\ra Y'$ is a fibrant replacement. Then as in \cite[\S 5]{DwSp}, \cite[Th.~2.3.2]{Hess}, \cite[\S 7.5]{Hirs}, and \cite[Th.~1.2.10]{Hove} the set of morphisms $X\ra Y$ in $\cC[\cW^{-1}]$ is given by
\begin{equation*}
\Hom_{\cC[\cW^{-1}]}(X,Y)\cong \Hom_\cC(X',Y')/\approx,
\end{equation*}
where $\approx$ is the equivalence relation of `homotopy' on morphisms $f':X'\ra Y'$ in $\cC$, as in \cite[\S 4]{DwSp}, \cite[Def.~2.2.4]{Hess}, \cite[Def.~7.3.2]{Hirs}, and~\cite[Def.~1.2.4]{Hove}.\smallskip

\noindent{\bf(b)} Suppose $X,Y$ are objects of $\cC$. Define a category $\cM(X,Y)$ whose objects $(U,V,f,i,j)$ are diagrams in $\cC$
\begin{equation*}
\xymatrix@C=40pt{ X & U \ar[r]^f \ar[l]_i & V & Y \ar[l]_j }
\end{equation*}
with $i,j$ weak equivalences, and whose morphisms $(U,V,f,i,j)\ra(U',\ab V',\ab f',\ab i',\ab j')$ are commutative diagrams in $\cC$ with $k,l$ weak equivalences
\begin{equation*}
\xymatrix@C=40pt@R=15pt{ X \ar[d]^{\id_X} & U \ar[d]^k \ar[r]_f \ar[l]^i & V \ar[d]^l & Y \ar[d]_{\id_Y} \ar[l]^j \\ 
X & U' \ar[r]^{f'} \ar[l]_{i'} & V' & Y.\! \ar[l]_{j'} }
\end{equation*}
Then (in modern language) Dwyer--Kan \cite[\S 6.2, \S 8.4]{DwKa} show that the mapping space $\Hom_{\cC[\cW^{-1}]^\iy}^L(X,Y)$ in the $\iy$-category $\cC[\cW^{-1}]^\iy$ is homotopy equivalent to the nerve $\sN(\cM(X,Y))$, considered as a simplicial set.

Here we can think of $(U,V,f,i,j)$ as a morphism $j^{-1}\ci f\ci i^{-1}:X\ra Y$ in a category $\cC[\cW^{-1}]$ in which the weak equivalences $i,j$ have been inverted.

Dugger \cite{Dugg} shows that if $Y$ is fibrant then $\Hom_{\cC[\cW^{-1}]^\iy}^L(X,Y)$ is also homotopy equivalent to $\sN(\cM(X,Y)_{j=\id})$, where $\cM(X,Y)_{j=\id}\subset \cM(X,Y)$ is the full subcategory of objects with $V=Y$ and $j=\id_Y$, and dually, if $X$ is cofibrant then $\Hom_{\cC[\cW^{-1}]^\iy}^L(X,Y)\simeq\sN(\cM(X,Y)_{i=\id})$, where $\cM(X,Y)_{i=\id}\subset \cM(X,Y)$ is the full subcategory of objects with $U=X$ and~$i=\id_X$.

\smallskip

\noindent{\bf(c)} A second way to compute mapping spaces $\Hom_{\cC[\cW^{-1}]^\iy}^L(X,Y)$ is the theory of {\it framings\/} in Hirschhorn \cite[\S 16]{Hirs} and Hovey \cite[\S 5]{Hove}. As in \cite[Def.~16.1.2]{Hirs}, a {\it simplicial resolution\/} of an object $Y$ in $\cC$ is a functor $S:\De^\rop\ra\cC$ (that is, a simplicial object in $\cC$) with an morphism $i:Y\,{\buildrel\cong\over\longra}\,S([0])$ in $\cW$ such that $S(f)\in\cW$ for all morphisms $f$ in $\De^\rop$, and $S$ is a fibrant replacement for the constant functor $\De^\rop\ra\{Y\}\subset\cC$ in the Reedy model structure on $\Fun(\De^\rop,\cC)$.

Let $X$ be another object in $\cC$, and $j:X'\ra X$ a cofibrant replacement. Then as in \cite[\S 17.2]{Hirs}, \cite[Th.~5.4.9]{Hove}, the mapping space $\Hom_{\cC[\cW^{-1}]^\iy}^L(X,Y)$ is homotopy equivalent to the simplicial set $\Hom_\cC(X',S(-)):\De^\rop\ra\Set$, which is a Kan complex by \cite[Cor.~16.5.3(2)]{Hirs}. This Kan complex is called the {\it right homotopy function complex\/} in~\cite[\S 17.2]{Hirs}.
\smallskip

\noindent{\bf(d)} Write $\cC_\cf\subset\cC$ for the full subcategory of cofibrant objects, and $\cW_\cf$ for the class of weak equivalences between cofibrant objects. Then the inclusion $\cC_\cf\hookra\cC$ induces an $\iy$-functor $\cC^{\vphantom{-1}}_\cf[\cW_\cf^{-1}]^\iy\ra\cC[\cW^{-1}]^\iy$. Cisinski \cite[Th.~7.5.18]{Cisi} implies that this is an equivalence of $\iy$-categories.

Thus for some purposes, it is helpful to restrict to the subcategory $\cC_\cf$, which makes cofibrant replacements $X'$ above unnecessary. In many interesting cases, all objects are fibrant, and then fibrant replacements are also unnecessary.
\end{rem}

The next definition and theorem generalize Remark \ref{ai2rem4}(c) to a case when $\cC$ need not be a model category. We will prove the theorem using the `Derived Mapping Space Lemma' of Arakawa--Carranza--Kapulkin \cite{ACK}, which is based on a theorem on `left calculi of fractions' in Cisinski~\cite[Def.~7.2.6 \& Th.~7.2.8]{Cisi}. 

To compare Definition \ref{ai2def18} and Theorem \ref{ai2thm5} to the model category case, we should imagine $\cC,\cW$ below to be $\cC_\cf,\cW_\cf$ in Remark~\ref{ai2rem4}(d).

\begin{dfn}
\label{ai2def18}
Let $(\cC,\cW)$ be a category with weak equivalences, as in Definition \ref{ai2def17}. Suppose $S:\De^\rop\ra\cC$ is a functor, that is, $S$ is a simplicial object in $\cC$. Then there is a natural functor $F_S:\cC^\rop\ra\SSet$ mapping $X\mapsto F_S(X)=\Hom_\cC(X,S(-)):\De^\rop\ra\Set$ for objects $X$ in $\cC$, and mapping $\ga\mapsto -\ci \ga:F_S(X_2)\ra F_S(X_1)$ for morphisms $\ga:X_1\ra X_2$ in~$\cC$.

Let $Y$ be an object of $\cC$. A {\it simplicial resolution\/} $(S,i)$ of $Y$ in $(\cC,\cW)$ is a functor $S:\De^\rop\ra\cC$ and a morphism $i:Y\ra S([0])$ in $\cW$, satisfying
\begin{itemize}
\setlength{\itemsep}{0pt}
\setlength{\parsep}{0pt}
\item[(i)] If $f:[N]\ra[M]$ is a morphism in $\De$ then $S(f)\in\cW$.
\item[(ii)] $F_S(\ga):F_S(X_2)\ra F_S(X_1)$ is a weak homotopy equivalence of simplicial sets for all morphisms $\ga:X_1\ra X_2$ in~$\cW$.
\end{itemize}

Using the notation of Example \ref{ai2ex5} and Definition \ref{ai2def16} we have an $\iy$-functor
\begin{equation*}
\sY_{\sN(\cC)}\ci \sN(S):\sN(\De^\rop)\longra\PSh(\sN(\cC)).
\end{equation*}
Since $\PSh(\sN(\cC))$ is cocomplete, we may form the $\iy$-categorical colimit
\e
\Ga_{\cC,S}:=\colim_{\sN(\De^\rop)}\sY_{\sN(\cC)}\ci \sN(S)\in\PSh(\sN(\cC)).
\label{ai2eq18}
\e
For colimits in $\iy$-categories see Lurie \cite[\S 4.2]{Luri2}, \cite[\S 7]{Luri5} and Land \cite[\S 4.3]{Land}. Here $\Ga_{\cC,S}$ is an $\iy$-functor $\sN(\cC)^\rop\ra\sSpc$. Since colimits in $\PSh(\sN(\cC))$ can be computed objectwise, for each $X\in\cC$ we have an equivalence in $\sSpc$
\e
\Ga_{\cC,S}(X)\simeq \colim_{\sN(\De^\rop)}\bHom_\cC(X,\sN(S)).
\label{ai2eq19}
\e

Now the functor $\bHom_\cC(X,\sN(S))$ in \eq{ai2eq19} is equivalent to the composition
\begin{equation*}
\xymatrix@C=26pt{  \sN(\De^\rop)\,\, \ar[r]^(0.52){\sN(F_S(X))} & \sN(\Set)\, \ar@{^{(}->}[r]^(0.52){\sN(\inc)} & \sN(\SSet) \ar[r]^(0.42){\Pi^\iy_\loc} & \SSet[{\cal{WE}}^{-1}]^\iy \ar[r]^(0.66){\eq{ai2eq16}}_(0.66)\simeq & \sSpc,  }
\end{equation*}
where $\inc:\Set\hookra\SSet$ maps a set $S$ to the constant simplicial set $\De^\rop\ra S\in\Set$. It follows from Lurie \cite[Ex.~A.2.9.31]{Luri2} that the $\iy$-categorical colimit of this composition is the simplicial set $F_S(X)$ (up to weak homotopy equivalence), considered as an object of $\SSet[{\cal{WE}}^{-1}]^\iy$. That is, we have
\e
\Ga_{\cC,S}(X)\simeq F_S(X).
\label{ai2eq20}
\e
Similarly, if $\ga:X_1\ra X_2$ is a morphism in $\cC$ then $\Ga_{\cC,S}(f):\Ga_{\cC,S}(X_2)\ra \Ga_{\cC,S}(X_1)$ is $F_S(\ga):F_S(X_2)\ra F_S(X_1)$, modulo weak homotopy equivalences of $F_S(X_1)$ and $F_S(X_2)$. In particular, if $\ga:X_1\ra X_2$ lies in $\cW$ then $\Ga_{\cC,S}(f):\Ga_{\cC,S}(X_2)\ra \Ga_{\cC,S}(X_1)$ is a weak homotopy equivalence by (ii) above, that is, $\Ga_{\cC,S}(f)$ is an equivalence in~$\sSpc$.

The rest of the definition follows Arakawa--Carranza--Kapulkin \cite[Def.~2.1]{ACK} in our special case. Let $[n]$ be an object in $\De$, and $\pi_n:[n]\ra[0]$ the unique morphism. Write $\Phi_{[n]}$ for the composition of morphisms in $\PSh(\sN(\cC))$
\begin{equation*}
\xymatrix@C=80pt{ 
\sY_{\sN(\cC)}(Y) \ar[r]^{\sY_{\sN(\cC)}(i\ci S(\pi_n))}  &  \sY_{\sN(\cC)}(S([n])) \ar[r]^{\colim_{[n]}}  & \Ga_{\cC,S}, }
\end{equation*}
where $\colim_{[n]}$ is the morphism to the colimit in \eq{ai2eq18}, which is natural up to homotopy. Then $\Phi_{[n]}$ is independent of $[n]$ up to homotopy, by properties of colimits.

Using the $\iy$-functor $\Pi_\loc^\iy:\sN(\cC)\ra\cC[\cW^{-1}]^\iy$, from Definition \ref{ai2def15} we have $\iy$-functors $(\Pi_\loc^\iy)^*:\PSh(\cC[\cW^{-1}]^\iy)\ra\PSh(\sN(\cC))$ and its left adjoint $(\Pi_\loc^\iy)_!:\PSh(\sN(\cC))\ra\PSh(\cC[\cW^{-1}]^\iy)$. Thus we have a 2-commutative diagram
\e
\begin{gathered}
\xymatrix@C=170pt@R=17pt{ 
*+[r]{\sY_{\sN(\cC)}(Y)} \drtwocell_{}\omit^{}\omit{} \ar[r]_(0.4){\Phi_{[n]}} \ar[d]^\eta & *+[l]{\Ga_{\cC,S}} \ar[d]_\simeq^\eta \\
*+[r]{(\Pi_\loc^\iy)^*\ci(\Pi_\loc^\iy)_!(\sY_{\sN(\cC)}(Y))} \ar[r]^(0.54){(\Pi_\loc^\iy)^*\ci(\Pi_\loc^\iy)_!(\Phi_{[n]})}  & *+[l]{(\Pi_\loc^\iy)^*\ci(\Pi_\loc^\iy)_!(\Ga_{\cC,S}),} }
\end{gathered}
\label{ai2eq21}
\e
where $\eta:\Id\Ra(\Pi_\loc^\iy)^*\ci(\Pi_\loc^\iy)_!$ is the unit of the adjunction. 

We showed above that if $\ga:X_1\ra X_2$ lies in $\cW$ then $\Ga_{\cC,S}(f)$ is an equivalence in $\sSpc$. From Definition \ref{ai2def17}, this implies that $\Ga_{\cC,S}(f)$ lies in the essential image of $(\Pi_\loc^\iy)^*:\PSh(\cC[\cW^{-1}]^\iy)\ra\PSh(\sN(\cC))$, which is full and faithful. So by properties of adjoint functors we see that the right hand column of \eq{ai2eq21} is an equivalence in~$\PSh(\sN(\cC))$.

Applying equation \eq{ai2eq13} to $F=\Pi_\loc^\iy$, and noting that $\Pi_\loc^\iy(Y)=Y$ as $\sN(\cC)$ and $\cC[\cW^{-1}]^\iy$ have the same objects, gives an equivalence in $\PSh(\cC[\cW^{-1}]^\iy)$
\e
\sY_{\cC[\cW^{-1}]^\iy}(Y)\simeq (\Pi_\loc^\iy)_!(\sY_{\sN(\cC)}(Y)).
\label{ai2eq22}
\e
Combining \eq{ai2eq21}--\eq{ai2eq22} and the fact that the right hand column of \eq{ai2eq21} is an equivalence, we see that there exists a morphism $\Psi$ unique up to equivalence which makes the following commute up to equivalence:
\e
\begin{gathered}
\xymatrix@C=175pt@R=17pt{ 
*+[r]{(\Pi_\loc^\iy)^*\ci \sY_{\cC[\cW^{-1}]^\iy}(Y)} \drtwocell_{}\omit^{}\omit{} \ar[r]_(0.4){\Psi} \ar[d]^{(\Pi_\loc^\iy)^*\eq{ai2eq22}}_\simeq & *+[l]{\Ga_{\cC,S}} \ar[d]_\simeq^\eta \\
*+[r]{(\Pi_\loc^\iy)^*\!\ci\!(\Pi_\loc^\iy)_!(\sY_{\sN(\cC)}(Y))} \ar[r]^(0.53){(\Pi_\loc^\iy)^*\ci(\Pi_\loc^\iy)_!(\Phi_{[n]})}  & *+[l]{(\Pi_\loc^\iy)^*\!\ci\!(\Pi_\loc^\iy)_!(\Ga_{\cC,S}).} }
\end{gathered}
\label{ai2eq23}
\e

Evaluating $\Psi$ at an object $X$ in $\cC$ and noting \eq{ai2eq21} and
\begin{align*}
\bigl((\Pi_\loc^\iy)^*\ci \sY_{\cC[\cW^{-1}]^\iy}(Y)\bigr)(X)&\simeq \bigl(\sY_{\cC[\cW^{-1}]^\iy}(Y)\bigr)(X)=\bHom_{\cC[\cW^{-1}]^\iy}(X,Y)\\
&\simeq\Hom_{\cC[\cW^{-1}]^\iy}^L(X,Y),
\end{align*}
gives a morphism in $\sSpc$
\e
\Psi(X):\Hom_{\cC[\cW^{-1}]^\iy}^L(X,Y)\longra F_S(X).
\label{ai2eq24}
\e
\end{dfn}

\begin{thm}
\label{ai2thm5}
In Definition\/ {\rm\ref{ai2def18},} the morphism $\Psi$ in \eq{ai2eq23} is an equivalence in $\PSh(\sN(\cC)),$ so that\/ $\Psi(X)$ in \eq{ai2eq24} is an equivalence in $\sSpc$ for all\/ $X\in\cC$.
\end{thm}

\begin{proof}
We apply Arakawa--Carranza--Kapulkin \cite[Def.~2.1 \& Th.~2.2]{ACK} in a special case, which is similar to their proof of Remark \ref{ai2rem4}(c) above in~\cite[Cor.~3.1]{ACK}. 

They begin with an $\iy$-category $\sC$ with a suitable $\iy$-subcategory $\sW\subset\sC$ to be localized, an object $Y\in\sC$, a small $\iy$-category $\sI$, and an $\iy$-functor $Y_\bu:\sI\ra\sC\downarrow Y$ to the undercategory of $Y$ in $\sC$. We take $\sC=\sN(\cC)$, $\sW=\sN(\cW)$, $I=\sN(\De^\rop)$, and $Y_\bu$ to be the lift of $\sN(S):\sN(\De^\rop)\ra\sN(\cC)$ to $\sN(\cC)\downarrow Y$ using the morphisms $i\ci S(\pi):Y\ra S([n])$ for all $[n]$ in $\De^\rop$. 

They impose conditions \cite[Def.~2.1(R1)--(R3)]{ACK} for $S$ to be an `$\cI$-resolution' of $Y$. Here (R1) is automatic as $\De^\rop$ has an initial object, (R2) follows from Definition \ref{ai2def18}(i), and (R3) follows from $\Ga_{\cC,S}(f):\Ga_{\cC,S}(X_2)\ra \Ga_{\cC,S}(X_1)$ an equivalence in $\sSpc$ if $\ga:X_1\ra X_2$ is a morphism in $\cW$, as proved in Definition \ref{ai2def18}. Thus \cite[Th.~2.2]{ACK} says that $\Psi$ in \eq{ai2eq23} is an equivalence.
\end{proof}

\subsection{\texorpdfstring{Stable $\iy$-categories}{Stable ∞-categories}}
\label{ai27}

Stable $\iy$-categories were introduced by Lurie \cite{Luri1}, \cite[\S 1]{Luri3}, see Groth~\cite[\S 5.1]{Grot}.

\begin{dfn}
\label{ai2def19}
Let $\sC$ be an $\iy$-category. A {\it zero object\/} in $\sC$ is an object $0\in\sC$ which is both initial and terminal (in the $\iy$-categorical sense). 

Suppose $\sC$ has a zero object. A {\it triangle\/} in $\sC$ is a 2-commutative diagram
\e
\begin{gathered}
\xymatrix@C=100pt@R=15pt{
*+[r]{X} \drtwocell_{}\omit^{}\omit{^{}} \ar[r]_(0.7)f \ar[d] & *+[l]{Y} \ar[d]_g \\ *+[r]{0} \ar[r] & *+[l]{Z,\!} }
\end{gathered}
\label{ai2eq25}
\e
where $X,Y,Z$ are objects in $\sC$ and $f,g$ are morphisms. We call \eq{ai2eq25} an {\it exact triangle\/} if it is Cartesian (in the $\iy$-categorical sense), and then we say that $X$ is the {\it kernel\/} of $g$. We call \eq{ai2eq25} a {\it coexact triangle\/} if it is co-Cartesian, and then we say that $Z$ is the {\it cokernel\/} of~$f$. 

We call $\sC$ {\it stable\/} if every morphism in $\sC$ has a kernel and a cokernel, and a triangle \eq{ai2eq25} is exact if and only if it is coexact. Write $\sCat_\iy^\rst\subset\sCat_\iy$ for the full $\iy$-subcategory with objects stable $\iy$-categories.
\end{dfn}

If $\sC$ is a stable $\iy$-category then its homotopy category $\Ho(\sC)$ has the structure of a triangulated category, as in Gelfand--Manin \cite{GeMa}. Thus stable $\iy$-categories are an $\iy$-categorical version of triangulated categories. Note however that to be triangulated is an extra structure on a category $\cT$, but to be stable is a property of an $\iy$-category $\sC$.

It has been apparent since the introduction of triangulated categories by Grothendieck and Verdier that there were problems with the theory -- for example, the `nonfunctoriality of the cone' in triangulated categories, as in Gelfand and Manin \cite[\S IV.1.7]{GeMa} -- and several attempts have been made to fix these problems by somehow enhancing the triangulated category structure.

There is a strong case for believing that stable $\iy$-categories are the `true' theory to which triangulated categories are an imperfect approximation. Stable $\iy$-categories have many wonderful properties. 

\subsection{\texorpdfstring{$\iy$-bicategories}{∞-bicategories}}
\label{ai28}

Informally, an $(\iy,n)$-{\it category\/} is some kind of higher category with a notion of $k$-morphism for all $k=1,2,\ldots,$ such that all $k$-morphisms are invertible for $k>n$. This is not a precise definition, and there are many different models for $(\iy,n)$-categories. Our $\iy$-categories in \S\ref{ai22}--\S\ref{ai27} are a kind of $(\iy,1)$-category, and Kan complexes are a kind of $(\iy,0)$-category.

We will explain $\iy$-{\it bicategories}, which are a particularly nice model for $(\iy,2)$-categories, following Lurie \cite[\S 4]{Luri4} and Gagna--Harpaz--Lanari \cite{GHL}. A slightly different notion appears in Lurie \cite[\S 5.4]{Luri5}. There it is called $\iy$-bicategory, conflicting with the notation of \cite{Luri4}, but we will call them {\it strong $\iy$-bicategories}. The next definition follows Lurie \cite[\S 3.1]{Luri4}. 

\begin{dfn}
\label{ai2def20}
A {\it scaled simplicial set\/} is a pair $(S,T)$ where $S$ is a simplicial set, and $T\subseteq S([2])$ is a subset of the 2-simplices of $S$ which contains every degenerate 2-simplex in $S$. A 2-simplex $t:\De^2\ra S$ is called {\it thin\/} if it lies in~$T$.

If $(S,T),(S',T')$ are scaled simplicial sets, a {\it morphism\/} $f:(S,T)\ra(S',T')$ is a morphism of simplicial sets $f:S\ra S'$ such that $f$ maps 2-simplices in $T$ to 2-simplices in $T'$. Then scaled simplicial sets form a category $\SSet^{\rm sc}$.
\end{dfn}

Lurie \cite[Th.~4.2.7]{Luri4} defines a combinatorial model category structure on $\SSet^{\rm sc}$, in which every object is cofibrant. In \cite[Th.~0.0.3]{Luri4} he shows that the model category $\SSet^{\rm sc}$ is part of a diagram of Quillen equivalences with other model categories that are agreed to be models for $(\iy,2)$-categories (see Gagna--Harpaz--Lanari \cite[Fig.~1]{GHL} for a larger diagram).

Lurie \cite[Def.~4.1.1]{Luri4} defines a class of scaled simplicial sets satisfying horn-filling conditions that he calls {\it weak\/ $\iy$-bicategories}. In \cite[Def.~4.2.8]{Luri4} he defines a scaled simplicial set to be an $\iy$-{\it bicategory\/} if it is a fibrant object in his model structure. As in \cite[Prop.~1.30]{GHL}, one can show using \cite[Prop.~3.1.13]{Luri4} that every $\iy$-bicategory is a weak $\iy$-bicategory. Gagna--Harpaz--Lanari \cite[Th.~5.1]{GHL} prove the converse, that every weak $\iy$-bicategory is an $\iy$-bicategory. The next definition follows \cite[Def.~4.1.1]{Luri4}, but using the name $\iy$-bicategory.

\begin{dfn}
\label{ai2def21}
A scaled simplicial set $(S,T)$ is called an $\iy$-{\it bicategory\/} if:
\begin{itemize}
\setlength{\itemsep}{0pt}
\setlength{\parsep}{0pt}
\item[(a)] For all $N\ge 3$ and $0<i<N$, a commutative diagram of morphisms `$\ra$' in which $t$ is a thin 2-simplex, as follows:
\e
\begin{gathered}
\xymatrix@C=120pt@R=15pt{ *+[r]{\La^N_{i_{\phantom{(}}}} \ar@{^{(}->}[d]^\inc \ar[dr]_(0.3)f 
& \De^2 \ar[l]^(0.25){\De^{\{i-1,i,i+1\}}} \ar[d]_t
\\
*+[r]{\De^N} \ar@{..>}[r]^{\ti f} & S, }
\end{gathered}
\label{ai2eq26}
\e
can be completed with a morphism $\ti f$ as shown. Here $\inc:\La^N_i\hookra\De^N$ is as in Examples \ref{ai2ex1}--\ref{ai2ex2}, and $\De^{\{i-1,i,i+1\}}:\De^2\ra\De^N$ is the 2-simplex in $\De^N$ with vertices $i-1,i,i+1$. This lies in $\La^N_i\subset\De^N$ as $N\ge 3$.
\item[(b)] Every commutative diagram of morphisms `$\ra$'
\e
\begin{gathered}
\xymatrix@C=100pt@R=15pt{ *+[r]{\La^2_{i_{\phantom{(}}}} \ar@{^{(}->}[d]^\inc \ar@/^.5pc/[dr]^f \\
*+[r]{\De^2} \ar@{..>}[r]^{\ti f} & S }
\end{gathered}
\label{ai2eq27}
\e
can be completed with a morphism $\ti f$ as shown, with $\ti f$ a thin 2-simplex.
\item[(c)] Suppose $f:\De^4\ra S$ is a morphism and $f\ci\De^{\{i,j,k\}}:\De^2\ra S$ is thin for $\{i,j,k\}=\{0,2,4\},\{1,2,3\},\{0,1,3\},\{1,3,4\},\{0,1,2\}$. Then $f\ci\De^{\{i,j,k\}}:\De^2\ra S$ is also thin for~$\{i,j,k\}=\{0,3,4\},\{0,1,4\}$.
\item[(d)] For $N\ge 3$, every commutative diagram of morphisms `$\ra$'
\e
\begin{gathered}
\xymatrix@C=70pt@R=15pt{ *+[r]{\La^N_{0_{\phantom{(}}}} \ar@{^{(}->}[d]^\inc \ar[drr]_(0.3)f 
& \De^2 \ar[l]_{\De^{\{0,1,N\}}} \ar[r]^{s_{1,0}} & \De^1 \ar[d]_g
\\
*+[r]{\De^N} \ar@{..>}[rr]^(0.45){\ti f} && S, }
\end{gathered}
\label{ai2eq28}
\e
can be completed with a morphism $\ti f$ as shown. Here $s_{1,0}:\De^2\ra\De^1$ is the {\it degeneracy map\/} acting on vertices by $0\mapsto 0$, $1\mapsto 0$, $2\mapsto 1$.

This is an {\it outer horn-filling condition\/}: it says that an outer horn $f:\La^N_0\ra S$ can be filled provided the 2-simplex $f\ci\De^{\{0,1,N\}}:\De^2\ra S$ is {\it left-degenerate\/} (factors via $s_{1,0}$). 
\end{itemize}

We can think of $\iy$-bicategory $(S,T)$ as a kind of category in a similar way to Definition \ref{ai2def5}. An {\it object\/} $X$ of $(S,T)$ is a 0-simplex in $S$, or equivalently, a map of simplicial sets $X:\De^0\ra S$. A {\it morphism}, or 1-{\it morphism}, $f:X\ra Y$ of $(S,T)$ is a 1-simplex $f$ in $S$ with $f\vert_{\De^{\{0\}}}=X$ and $f\vert_{\De^{\{1\}}}=Y$. If $f:X\ra Y$ and $g:Y\ra Z$ are 1-morphisms, to define the `composition' $g\ci f:X\ra Z$ we choose a thin 2-simplex $t:\De^2\ra S$ with $t\ci\De^{\{0,1\}}=f$ and $t\ci\De^{\{1,2\}}=g$, which exists by (b). Then we set $g\ci f=t\ci\De^{\{0,2\}}$. As in $\iy$-categories, composition is not uniquely defined, as $t$ need not be unique.

Let $\eta:\De^2\ra S$ be a morphism, and set $X=\eta\ci\De^{\{0\}}$, $Y=\eta\ci\De^{\{1\}}$, $Z=\eta\ci\De^{\{2\}}$, $f=\eta\ci\De^{\{0,1\}}$, $g=\eta\ci\De^{\{1,2\}}$, $h=\eta\ci\De^{\{0,2\}}$, so that $f:X\ra Y$, $g:Y\ra Z$, $h:X\ra $ are 1-morphisms. We think of $\eta$ as giving a 2-{\it morphism\/} $g\ci f\Ra h$ in $(S,T)$, which is a 2-{\it isomorphism\/} if $\eta$ is thin. 

Although part (d) above is not symmetric under taking opposite simplicial sets $S\mapsto S^\rop$, Gagna--Harpaz--Lanari \cite[Cor.~5.5]{GHL} show that $(S,T)$ is an $\iy$-bicategory if and only if $(S^\rop,T)$ is an $\iy$-bicategory.
\end{dfn}

\begin{dfn}
\label{ai2def22}
Let $(S,T)$ be an $\iy$-bicategory. Define the {\it core\/} $\Core(S,T)$ of $(S,T)$ to be the simplicial subset $\Core(S,T)\subseteq S$ such that $f:\De^N\ra S$ is an $N$-simplex of $\Core(S,T)$ if and only if all compositions $\De^2\,{\buildrel g\over\longra}\,\De^N\,{\buildrel f\over\longra}\,S$ are thin 2-simplices in $S$. It follows from Definition \ref{ai2def21}(a),(b) that $\Core(S,T)$ is an $\iy$-category. Here our notation `core' comes from Gagna--Harpaz--Lanari \cite[Def.~1.22]{GHL}. It is called the {\it underlying\/ $\iy$-category\/} of $(S,T)$ in Lurie \cite[Rem.~4.1.3]{Luri4}, and the {\it pith\/} of $(S,T)$ in Lurie~\cite[\S 5.4.5]{Luri5}. 

In effect, $\Core(S,T)$ is the $(\iy,2)$-subcategory of $(S,T)$ that only includes 2-morphisms which are 2-isomorphisms, and is thus an $(\iy,1)$-subcategory.
\end{dfn}

\begin{ex}
Let $(S,T)$ be a scaled simplicial set in which every 2-simplex in $S$ is thin. Then by Lurie \cite[Rem.~4.1.2]{Luri4}, $(S,T)$ is an $\iy$-bicategory if and only if $S$ is an $\iy$-category, and then $\Core(S,T)=S$.
\end{ex}

\begin{dfn}
\label{ai2def23}
An $\iy$-bicategory $(S,T)$ is {\it strong\/} if $T$ is the set of all 2-simplices $t:\De^2\ra S$ for which Definition \ref{ai2def21}(a) holds. Then $T$ is determined by $S$, so we can just refer to $S$ as a strong $\iy$-bicategory. Lurie \cite[\S 5.4]{Luri5} develops a theory of `$\iy$-bicategories', which are equivalent to our strong $\iy$-bicategories.
\end{dfn}

If $\cC$ is an ordinary category, Example \ref{ai2ex5} defined an $\iy$-category $\sN(\cC)$, the nerve of $\sC$. In a similar way, if $\cC$ is a weak 2-category (bicategory), Lurie \cite[\S 2.3]{Luri5} defines the {\it Duskin nerve\/} $\sN^D(\cC)$, a simplicial set, and shows it is a strong $\iy$-bicategory~\cite[Prop.~5.4.1.5]{Luri5}.

In Theorems \ref{ai6thm4}, \ref{ai7thm4} and \ref{ai8thm3} below we will construct $\iy$-bicategories of $A_\iy$-algebras and $A_\iy$-categories. As in Remark \ref{ai7rem4}(f), the author expects that these are {\it not\/} strong $\iy$-bicategories. This is our reason for choosing the $\iy$-bicategory theory from \cite[\S 4]{Luri4} rather than from~\cite[\S 5.4]{Luri5}.

\section{\texorpdfstring{The classical theory of $A_\iy$-algebras}{The classical theory of A∞-algebras}}
\label{ai3}

$A_\iy$-algebras were introduced by Stasheff \cite{Stas1,Stas2}. Some good references are Fukaya et al.\ \cite[\S 3.2--\S 3.3]{FOOO}, Seidel \cite{Seid}, and two survey papers by Keller \cite{Kell1,Kell2}. Kontsevich and Soibelman \cite{KoSo} give a geometric point of view on $A_\iy$-algebras. We follow the grading and sign conventions of Keller~\cite{Kell1}.

\subsection{Sign conventions}
\label{ai31}

We will use the following sign convention for tensor products of graded $R$-linear and $R$-multilinear operators throughout the book.

\begin{conv}
\label{ai3conv1}
Suppose $R$ is a commutative ring, $n\ge 1$, and $p_1,\ldots,p_n\ge 0$, and that $B_{i,j}=\bigop_{d\in\Z}B_{i,j}^d$ and $C_i=\bigop_{d\in\Z}C_i^d$ are $\Z$-graded $R$-modules for $i=1,\ldots,n$ and $j=1,\ldots,p_i$, and that 
\begin{equation*}
g_i:B_{i,1}\t B_{i,2}\t\cdots\t B_{i,p_i}\longra C_i
\end{equation*}
is an $R$-multilinear graded map of degree $m_i$ for $i=1,\ldots,n$, that is, $g_i$ maps $B_{i,1}^{d_1}\t \cdots\t B_{i,p_i}^{d_{p_i}}\ra C_i^{m_i+\sum_{j=1}^{p_i}d_j}$ for $d_1,\ldots,d_{p_i}\in\Z$. Equivalently, we can regard $g_i$ as an $R$-linear map graded of degree $m_i$
\begin{equation*}
g_i:B_{i,1}\ot_R B_{i,2}\ot_R\cdots\ot_R B_{i,p_i}\longra C_i,
\end{equation*}
by identifying $g_i(b_{i,1},\ldots,b_{i,p_i})=g_i(b_{i,1}\ot\cdots\ot b_{i,p_i})$. We often use the same symbol for $R$-multilinear maps, and $R$-linear maps of tensor products, in this~way.

When $p_i=0$, we consider $g_i$ to be an element of $C_i^{m_i}$, or equivalently, an $R$-linear map $R\ra C_i^{m_i}$ which recovers the element of $C_i^{m_i}$ as $g_i(1)$.

Define
\begin{equation*}
g_1\ot\cdots\ot g_n:B_{1,1}\ot_R\cdots\ot_R B_{1,p_1}\ot_R B_{2,1}\ot_R\cdots \ot_R B_{n,p_n}
\longra C_1\ot_R\cdots\ot_RC_n
\end{equation*}
to be the $R$-linear map graded of degree $m_1+\cdots+m_n$ satisfying
\e
\begin{split}
&(g_1\ot\cdots\ot g_n)(b_{1,1}\ot\cdots\ot b_{1,p_1}\ot b_{2,1}\ot\cdots\ot b_{n,1}\cdots\ot b_{n,p_n})\\
&=(-1)^{\sum_{\begin{subarray}{l} 1\le i<k\le n,\\ j=1,\ldots,p_i\end{subarray}}d_{i,j}m_k}g_1(b_{1,1},\ldots,b_{1,p_1})\ot\cdots\ot g_n(b_{n,1},\ldots,b_{n,p_n})
\end{split}
\label{ai3eq1}
\e
whenever $b_{i,j}\in B_{i,j}$ is pure of degree $d_{i,j}$, that is, $b_{i,j}\in B_{i,j}^{d_{i,j}}$, for~$d_{i,j}\in\Z$. 

The important point in \eq{ai3eq1} is the sign $(-1)^{\sum_{i<k}d_{i,j}m_k}$, which we consider to be {\it part of the definition\/} of $g_1\ot\cdots\ot g_n$. Here is how to think about it: we split $B_{i,j}$ as a direct sum $B_{i,j}=B_{i,j}^{\rm even}\op B_{i,j}^{\rm odd}$ of `even' and `odd' parts, where $B_{i,j}^{\rm even}=\bigop_{d\in 2\Z}B_{i,j}^d$ and $B_{i,j}^{\rm odd}=\bigop_{d\in 2\Z+1}B_{i,j}^d$. We also call $g_i$ `even' or `odd' according to whether $m_i=\deg g_i$ is even or odd. 

The rule, known as the {\it Koszul sign convention}, is that whenever we exchange the order of two `odd' elements or operators, we should change the sign by $-1$. In equation \eq{ai3eq1}, the $g_i,b_{i,j}$ are ordered $g_1,\ldots,g_n,b_{1,1},\ldots,b_{n,p_n}$ on the left hand side, but ordered $g_1,b_{1,1},\ldots,b_{1,p_1},g_2,b_{2,1},\ldots,g_n,b_{n,1},\ldots,b_{n,p_n}$ on the right hand side. So the order of $b_{i,j}$ and $g_k$ is different on the two sides if $i<k$, which gives the sign change $(-1)^{d_{i,j}m_k}$, and the overall sign~$(-1)^{\sum_{i<k}d_{i,j}m_k}$.

Now suppose $A_{i,j,k}$ is a $\Z$-graded $R$-module for $i=1,\ldots,n$, $j=1,\ldots,p_i$ and $k=1,\ldots q_{i,j}$, and we are given $R$-multilinear maps
\begin{equation*}
f_{i,j}:A_{i,j,1}\t\cdots\t A_{i,j,q_{i,j}}\longra B_{i,j},
\end{equation*}
graded of degree $l_{i,j}$ for $i=1,\ldots,n$, $j=1,\ldots,p_i$. Equivalently, $f_{i,j}$ is an $R$-linear map graded of degree $l_{i,j}$
\begin{equation*}
f_{i,j}:A_{i,j,1}\ot\cdots\ot A_{i,j,p_{i,j}}\longra B_{i,j}.
\end{equation*}
Then we can show from \eq{ai3eq1} that
\ea
&\bigl(g_1\ci(f_{1,1}\ot\cdots\ot f_{1,p_1})\bigr)\ot\cdots\ot \bigl(g_n\ci (f_{n,1}\ot\cdots\ot f_{n,p_n})\bigr)
\label{ai3eq2}\\
&=(-1)^{\sum_{\begin{subarray}{l} 1\le i<k\le n, \\ j=1,\ldots,p_i\end{subarray}}l_{i,j}m_k}(g_1\ot\cdots\ot g_n)
\nonumber\\
&\qquad \ci (f_{1,1}\ot\cdots\ot f_{1,p_1}\ot f_{2,1}\ot\cdots\ot f_{n,1}\ot\cdots\ot f_{n,p_n}).
\nonumber
\ea
This is consistent with the principle that exchanging the order of $f_{i,j}$ and $g_k$ changes signs by $(-1)^{l_{i,j}m_k}$.
\end{conv}

Of course, we could have defined $g_1\ot\cdots\ot g_n$ without the sign $(-1)^{\sum_{i<k}d_{i,j}m_k}$ in \eq{ai3eq1}. But then the signs in the rest of the book would have to be be different, and more complicated. Convention \ref{ai3conv1} simplifies the signs later.

\subsection{\texorpdfstring{$A_\iy$-algebras, and units in $A_\iy$-algebras}{A∞-algebras, and units in A∞-algebras}}
\label{ai32}

Our sign conventions for $A_\iy$-algebras follow Keller~\cite{Kell1}.

\begin{dfn}
\label{ai3def1}
Let $R$ be a commutative ring. An $A_\iy$-{\it algebra} $(A,\m)$ over $R$ consists of:
\begin{itemize}
\setlength{\itemsep}{0pt}
\setlength{\parsep}{0pt}
\item[(a)] A $\Z$-graded $R$-module $A=\bigop_{d\in\Z}A^d$; and
\item[(b)] Graded $R$-linear maps $\m_k:A^{\ot^k}\ra A$ of degree $2-k$ for $k\ge 1$.
Write $\m=(\m_k)_{k\ge 1}$. These must satisfy the equation in maps $A^{\ot^n}\ra A$ for all $n\ge 1$
\end{itemize}
\e
\sum_{i,k\ge 0,\; j\ge 1:\; i+j+k=n}
(-1)^{i+jk}\m_{i+k+1}\ci\bigl(\id_A^{\ot^i}\ot\m_j\ot \id_A^{\ot^k}\bigr)=0.
\label{ai3eq3}
\e

By abuse of notation, we also think of $\m_k$ as a graded $R$-multilinear map $\m_k:{\buildrel{\ulcorner\,\,\,\text{$k$ copies } \,\,\,\urcorner} \over {\vphantom{m}\smash{A\t\cdots\t A}}}\ra A$ of degree $2-k$, with $\m_k(a_1,\ldots,a_k)=\m_k(a_1\ot\cdots\ot a_k)$. Then by Convention \ref{ai3conv1}, equation \eq{ai3eq3} is equivalent to the equation for all $n\ge 1$ and all pure $a_1,\ldots,a_n$ in $A$
\e
\begin{gathered}
\sum_{i,k\ge 0,\; j\ge 1:\; i+j+k=n}
\begin{aligned}[t]
&(-1)^{i+jk+j\sum_{l=1}^i\deg a_l}\,\m_{i+k+1}(a_1,\ldots,a_i,\\
&\qquad\m_j(a_{i+1},\ldots,a_{i+j}),a_{i+j+1},\ldots,a_n)=0.\end{aligned}
\end{gathered}
\label{ai3eq4}
\e

An $A_\iy$-algebra $(A,\m)$ is called {\it minimal\/} if $\m_1=0$.
\end{dfn}

\begin{rem}
\label{ai3rem1}
{\bf(a)} Our $A_\iy$-algebras are also called {\it non-unital\/} $A_\iy$-algebras, as in Seidel \cite[\S 1]{Seid} for instance. See Definitions \ref{ai3def2}--\ref{ai3def4} for units in $A_\iy$-algebras.
\smallskip

\noindent{\bf(b)} There is an alternative way to write the definition of $A_\iy$-algebras using {\it bar complexes}, which is simpler in some ways. See~\S\ref{ai36}.
\end{rem}

If $(A,\m)$ is an $A_\iy$-algebra then \eq{ai3eq3} for $n=1$ becomes $\m_1\ci\m_1=0$, with $\m_1$ of degree 1. Thus $(A,\m_1)$ is a complex, and we form its cohomology $H^*(A)$~by
\begin{equation*}
H^p(A)=\frac{\Ker\bigl(\m_1:A^p\ra A^{p+1}\bigr)}{\Im\bigl(\m_1:A^{p-1}\ra A^p\bigr)}\,.
\end{equation*}
Then $\m_k$ for $k>1$ induce various operations on $H^*(A)$. For example, \eq{ai3eq4} when $k=2$ yields
\begin{equation*}
-\m_2(\m_1(a_1),a_2)-(-1)^{\deg a_1}\m_2(a_1,\m_1(a_2))+\m_1(\m_2(a_1,a_2))=0.
\end{equation*}
This implies that the bilinear product $\bu:H^p(A)\t H^q(A)\ra H^{p+q}(A)$ given by
\e
(a_1+\Im\m_1)\bu(a_2+\Im\m_1)=\m_2(a_1,a_2)+\Im\m_1
\label{ai3eq5}
\e
is well-defined. Then \eq{ai3eq4} when $k=3$ implies that $\bu$ is {\it associative}. Thus $(H^*(A),\bu)$ is an associative $R$-algebra without unit.

Here are three ways to include units in $A_\iy$-algebras.

\begin{dfn}
\label{ai3def2}
Let $(A,\m)$ be an $A_\iy$-algebra over $R$. A {\it unit\/} for $(A,\m)$ is an element $\boo\in A^0$ satisfying
\e
\begin{split}
\m_2(\boo\ot a)=\m_2(a\ot \boo)&=a \;\>\text{for all }a\in A, \\
\m_k(a_1\!\ot\!\cdots\!\ot\! a_{i-1}\!\ot\! \boo\!\ot\! a_{i+1}\!\ot\! \cdots\!\ot\! a_k)&=0\; \text{for all $a_j\!\in\! A$, $1\!\le\! i\!\le\! k$, $k\!\ne\! 2$.}
\end{split}
\label{ai3eq6}
\e
We call $(A,\m)$ {\it unital\/} if it has a unit. If $\boo,\boo'$ are units for $(A,\m)$ then $\boo=\m_2(\boo,\boo')=\boo'$, so a unit $\boo$ is unique if it exists. Thus, to be unital is a property of an $A_\iy$-algebra $(A,\m)$, not extra data.
\end{dfn}

\begin{dfn}
\label{ai3def3}
We call an $A_\iy$-algebra $(A,\m)$ over $R$ {\it weakly cohomologically unital\/} if the associative algebra $H^*(A)$ has a unit, that is, if there exists $1_A$ in $H^0(A)$ such that $1_A\bu\al=\al\bu 1_A=\al$ for all $\al\in H^*(A)$. This $1_A$ is unique if it exists, so to be cohomologically unital is a property of $(A,\m)$, not extra data.
\end{dfn}

\begin{dfn}
\label{ai3def4}
A {\it cohomological unit\/} for an $A_\iy$-algebra $(A,\m)$ over $R$ is $e_A$ in $A^0$ such that $\m_1(e_A)=0$, and the following maps of chain complexes are chain homotopic to the identity:
\begin{align*}
\Phi&:(A,\m_1)\longra(A,\m_1),& \Phi&:a\longmapsto \m_2(e_A\ot a), \\
\Psi&:(A,\m_1)\longra(A,\m_1),& \Psi&:a\longmapsto \m_2(a\ot e_A).
\end{align*}
These imply that the cohomology class $1_A=[e_A]$ in $H^0(A)$ is a unit in $H^*(A)$.

We call $(A,\m)$ {\it cohomologically unital\/} if it has a cohomological unit~$e_A$.

If $R$ is a field then weakly cohomologically unital and cohomologically unital coincide, since $\Phi,\Psi$ are chain homotopic to the identity if and only if they induce the identity maps on cohomology. 

A unit $\boo_A$ in $(A,\m)$ is a cohomological unit. Thus, $(A,\m)$ unital implies $(A,\m)$ cohomologically unital, implies $(A,\m)$ weakly cohomological unital.

If $e_A'\in A^0$ is any other representative for $1_A$ in $H^0(A)$, then $e_A'=e_A+\m_1(f)$ for $f\in A^{-1}$, and defining $\Th,\Xi:A^*\ra A^{*-1}$ by $\Th(a)=\m_2(f\ot a)$, $\Xi(a)=-\m_2(a\ot f)$, we see from \eq{ai3eq4} that
\begin{equation*}
\Phi'(a)\!=\!\Phi(a)\!+\!\m_1\ci\Th(a)\!-\!\Th\ci\m_1(a),\;\> \Psi'(a)\!=\!\Psi(a)\!+\!\m_1\ci\Xi(a)\!-\!\Xi\ci\m_1(a).
\end{equation*}
Thus, $\Phi',\Psi'$ are chain homotopic to the identity if and only if $\Phi,\Psi$ are, and $e_A'$ is a cohomological unit if and only if $e_A$ is. Hence, to be cohomologically unital is a property of $(A,\m)$, not of the particular choice of element~$e_A$.
\end{dfn}

\begin{rem}
\label{ai3rem2}
Readers are warned that the names used for different kinds of unitality vary widely. Our notion of unital is often called `strictly unital'. Our notion of cohomologically unital is sometimes just called `unital'. Some references on different forms of units (for both $A_\iy$-algebras and $A_\iy$-categories) are Keller \cite[\S 3.5]{Kell1}, Seidel \cite[\S 2]{Seid}, Canonaco--Ornaghi--Stellari \cite[\S 1.3]{COS2}, Lyubashenko--Manzyuk \cite{LyMa2}, and Tanaka \cite[\S 1]{Tana}. Seidel \cite[\S 2]{Seid} and Fukaya et al.\ \cite[\S 3.3]{FOOO} also define {\it homotopy units\/} in $A_\iy$-algebras, but we do not discuss them.

Although weakly cohomologically unital $A_\iy$-algebras seem a very natural definition, the additional chain homotopy condition in cohomologically unital tends to be needed in proving theorems (see e.g.\ \cite{LyMa2}, \cite[Th.~1.1]{Tana}).
\end{rem}

\subsection{\texorpdfstring{Morphisms between $A_\iy$-algebras}{Morphisms between A∞-algebras}}
\label{ai33}

Here is the notion of morphism of $A_\iy$-algebras.

\begin{dfn}
\label{ai3def5}
Let $(A,\m),(B,\n)$ be $A_\iy$-algebras over $R$. An $A_\iy$-{\it morphism}, or just {\it morphism}, $\f:(A,\m)\ra(B,\n)$ is $\f=(\f_k)_{k\ge 1}$, where $\f_k:A^{\ot^k}\ra B$ for $k=1,2,\ldots$ are graded $R$-linear maps of degree $1-k$, satisfying in maps $A^{\ot^n}\ra B$ for all $n\ge 1$
\e
\begin{split}
&\sum_{i,k\ge 0,\; j\ge 1:\; i+j+k=n}(-1)^{i+jk}\f_{i+k+1}\bigl(
\id_A^{\ot^i}\ot\m_j\ot\id_A^{\ot^k}\bigr) \\ 
&=\sum_{n=n_1+\cdots+n_k,\; n_a\ge 1\!\!\!\!\!\!\!\!}
(-1)^{\sum_{i=1}^k(n_i+1)(k-i)}\n_k\bigl(\f_{n_1}\ot\f_{n_2}\ot\cdots\ot\f_{n_k}\bigr).
\end{split}
\label{ai3eq7}
\e

The {\it identity morphism\/} $\id_A:(A,\m)\ra(A,\m)$ is $\id_A=(\id_{A,k})_{k\ge 1}$ with $\id_{A,1}=\id_A$ and $\id_{A,k}=0$ for $k>1$.

We call a morphism $\f:(A,\m)\ra(B,\n)$ {\it strict\/} if $\f_k=0$ for $k\ne 1$, an $A_\iy$-{\it isomorphism} if $\f_1:A\ra B$ is an isomorphism of $R$-modules, and a {\it strict\/ $A_\iy$-isomorphism} if it is both strict and an $A_\iy$-isomorphism. When $n=1$, \eq{ai3eq7} becomes $\f_1\ci\m_1=\n_1\ci\f_1:A\ra B$. Thus we have a morphism of complexes
\e
\f_1:(A,\m_1)\longra(B,\n_1),
\label{ai3eq8}
\e
which induces a morphism of cohomology groups 
\e
H^*(\f_1):H^*(A,\m_1)\longra H^*(B,\n_1).
\label{ai3eq9}
\e

Following \cite[Def.~3.2.10]{FOOO}, we call a morphism $\f:(A,\m)\ra(B,\n)$ a {\it weak homotopy equivalence\/} if \eq{ai3eq8} is a chain homotopy equivalence. This is automatic if $\f$ is an $A_\iy$-isomorphism. As in \cite[Rem.~3.2.11]{FOOO}, if $R$ is a field then $\f$ is a weak homotopy equivalence if and only if \eq{ai3eq9} is an isomorphism.

If $(A,\m),(B,\n),(C,{\mathfrak o})$ are $A_\iy$-algebras and $\f:(A,\m)\ra(B,\n)$, $\g:(B,\n)\ra(C,{\mathfrak o})$ are morphisms, the {\it composition} $\g\ci\f:(A,\m) \ra(C,{\mathfrak o})$ is given by
\e
(\g\ci\f)_n=\sum_{n=n_1+\cdots+n_k,\; n_a\ge 1\!\!\!\!\!\!\!\!\!\!\!\!\!\!}(-1)^{\sum_{i=1}^k(n_i+1)(k-i)}\g_k\bigl(\f_{n_1}\ot\f_{n_2}\ot\cdots\ot\f_{n_k}).
\label{ai3eq10}
\e
Composition is associative. Write $\AAlg_R$ for the category with objects $A_\iy$-algebras $(A,\m)$ over $R$ and morphisms $\f:(A,\m)\ra(B,\n)$.

If $(A,\m),(B,\n)$ have units $\boo_A,\boo_B$, as in Definition \ref{ai3def2}, we call a morphism $\f:(A,\m)\ra(B,\n)$ {\it unital\/} if
\e
\begin{gathered}
\f_1(\boo_A)=\boo_B \quad\text{and}\quad \f_k(a_1\ot\cdots\ot a_{i-1}\ot\boo_A\ot a_{i+1}\ot\cdots\ot a_k)=0\\
\text{for all $a_j\in A$, $i=1,\ldots,k$, and $k\ne 1$.}
\end{gathered}
\label{ai3eq11}
\e
Unital $A_\iy$-algebras over $R$ and unital morphisms form a non-full subcategory $\AAlg_R^\un\subset\AAlg_R$.

Let $(A,\m),(B,\n)$ be (weakly) cohomologically unital $A_\iy$-algebras, as in Definitions \ref{ai3def3} and \ref{ai3def4}, and $\f:(A,\m)\ra(B,\n)$ be a morphism. We call $\f$ ({\it weakly\/}) {\it cohomologically unital\/} if $H^0(\f_1)(1_A)=1_B$. Then (weakly) cohomologically unital $A_\iy$-algebras over $R$ and (weakly) cohomologically unital morphisms form non-full subcategories $\AAlg_R^\cu\subset\AAlg_R^\wc\subset\AAlg_R$. When $(A,\m),(B,\n)$ are unital, unital morphisms $\f:(A,\m)\ra(B,\n)$ are cohomologically unital, so~$\AAlg_R^\un\subset\AAlg_R^\cu$.
\end{dfn}

\subsection{\texorpdfstring{Homotopies between $A_\iy$-morphisms, and $A_\iy$-algebras}{Homotopies between A∞-morphisms, and A∞-algebras}}
\label{ai34}

Now let $(A,\m),(B,\n)$ be $A_\iy$-algebras, and $\f,\g:(A,\m)\ra(B,\n)$ be morphisms. We will define the notion of {\it homotopy} $\H$ from $\f$ to $\g$. Our definition is based on Keller \cite[\S 3.7]{Kell1}. Fukaya et al.\ \cite[\S 4]{FOOO} use a different, more complicated definition, involving `models of $[0,1]\t B$', but they show in \cite[Prop.~4.2.40]{FOOO} that the two definitions yield the same notion of whether $\f,\g$ are homotopic.

\begin{dfn}
\label{ai3def6}
Let $(A,\m),(B,\n)$ be $A_\iy$-algebras over $R$, and $\f,\g:(A,\m)\ra(B,\n)$ be morphisms. A {\it homotopy} $\H:\f\Ra\g$ from $\f$ to $\g$ is $\H=(\H_k)_{k\ge 1}$, where $\H_k:A^{\ot^k}\ra B$ for $k=1,2,\ldots$ are graded $R$-linear maps of degree $-k$, satisfying
\ea
&\f_n-\g_n=
\label{ai3eq12}\\
&\sum_{\begin{subarray}{l} n=n_1+\cdots+n_k, \\ 
n_i\ge 1,\; 1\le l\le k \end{subarray}}
\begin{aligned}[t]&
(-1)^{k+1}\cdot(-1)^{\sum_{i=1}^{l-1}(n_i+1)}\cdot(-1)^{\sum_{i=1}^k(n_i+1)(k-i)}\\
&\n_k\bigl(\f_{n_1}\ot\f_{n_2}\ot\cdots\ot\f_{n_{l-1}}\ot
\H_{n_l}\ot\g_{n_{l+1}}\ot\cdots\ot\g_{n_k}\bigr)
\end{aligned}
\nonumber\\
&+\sum_{i,k\ge 0,\; j\ge 1:\; i+j+k=n}
(-1)^{i+jk}\H_{i+k+1}\ci\bigl(\id_A^{\ot^i}\ot\m_j\ot \id_A^{\ot^k}\bigr). 
\nonumber
\ea

When $n=1,$ equation \eq{ai3eq12} becomes
\begin{equation*}
\f_1-\g_1=\n_1\ci\H_1+\H_1\ci\m_1.
\end{equation*}
This implies that $\H_1:A^*\ra B^{*-1}$ is a chain homotopy between the morphisms of chain complexes $\f_1,\g_1:(A,\m_1)\ra(B,\n_1)$. Thus for the morphisms \eq{ai3eq9} on cohomology we have
\begin{equation*}
H^*(\f)=H^*(\g):H^*(A)\longra H^*(B).
\end{equation*}

If $\f,\g:(A,\m)\ra(B,\n)$ are {\it unital\/} morphisms of {\it unital\/} $A_\iy$-algebras, as in Definitions \ref{ai3def2} and \ref{ai3def5}, we call a homotopy $\H:\f\Ra\g$ {\it unital\/} if for all $1\le i\le k$ and $a_j\in A$ we have \e
\H_k(a_1\ot\cdots\ot a_{i-1}\ot\boo_A\ot a_{i+1}\ot\cdots\ot a_k)=0.
\label{ai3eq13}
\e

For ({\it weakly\/}) {\it cohomologically unital\/} $A_\iy$-algebras and morphisms, as in Definitions \ref{ai3def3}, \ref{ai3def4} and \ref{ai3def5}, there seems to be no useful notion of (weakly) cohomologically unital homotopy between (weakly) cohomologically unital morphisms, one should allow all homotopies.
\end{dfn}

\begin{dfn}
\label{ai3def7}
Let $(A,\m),(B,\n)$ be $A_\iy$-algebras, and $\f:(A,\m)\ra(B,\n)$ a morphism. A {\it homotopy inverse} for $\f$ is a morphism $\g:(B,\n)\ra(A,\m)$ such that $\g\ci\f:(A,\m)\ra(A,\m)$ is homotopic to $\id_A:(A,\m)\ra(A,\m)$, and $\f\ci\g:(B,\n)\ra(B,\n)$ is homotopic to $\id_B:(B,\n)\ra(B,\n)$. If $\f$ has a homotopy inverse, we call $\f$ a {\it homotopy equivalence}, and we call $(A,\m),(B,\n)$ {\it homotopy equivalent}. If $\f:(A,\m)\ra(B,\n)$ is a homotopy equivalence then $H^*(\f):H^*(A)\ra H^*(B)$ is an isomorphism. 
\end{dfn}

The next theorem is proved by Fukaya et al.\ \cite[Cor.~4.2.44, Th.~4.2.45]{FOOO}, see also Keller \cite[\S 3.7]{Kell1}, Lef\`evre-Hasegawa \cite[Cor.~1.3.1.3]{LeHa} and Seidel \cite[\S 1(h)]{Seid}.

\begin{thm}
\label{ai3thm1}
Let\/ $(A,\m),(B,\n)$ be $A_\iy$-algebras over $R$. Then
\begin{itemize}
\setlength{\itemsep}{0pt}
\setlength{\parsep}{0pt}
\item[{\rm(a)}] Homotopy is an equivalence relation on morphisms $\f,\g:(A,\m)\ra(B,\n)$.
\item[{\rm(b)}] Homotopy equivalence is an equivalence relation on $A_\iy$-algebras.
\item[{\rm(c)}] A morphism $\f:(A,\m)\ra(B,\n)$ is a homotopy equivalence if and only if it is a weak homotopy equivalence in the sense of Definition\/ {\rm\ref{ai3def5}}. Hence, any $A_\iy$-isomorphism is a homotopy equivalence.
\end{itemize}
\end{thm}

When $R$ is a field, Kadeishvili \cite{Kade} proves:

\begin{thm}
\label{ai3thm2}
Let\/ $(A,\m)$ be an $A_\iy$-algebra over a field\/ $R$. Then there exists a \begin{bfseries}minimal\end{bfseries}\/ $A_\iy$-algebra $(B,\n)$ homotopy equivalent to $(A,\m),$ in the sense of Definition\/ {\rm\ref{ai3def1}}. As graded\/ $R$-modules we have $B\cong H^*(A)$.
\end{thm}

In practice, homotopy equivalence is a more useful notion of when two $A_\iy$-algebras are `the same' than (strict) $A_\iy$-isomorphism. We are interested in properties of $A_\iy$-algebras which are invariant under homotopy. Constructions of $A_\iy$-algebras generally depend on some arbitrary choices, and different choices yield homotopy equivalent but not (strictly) isomorphic $A_\iy$-algebras.

\begin{dfn}
\label{ai3def8}
The {\it homotopy category\/} $\Ho(\AAlg_R)$ of $\AAlg_R$ is the category with objects $A_\iy$-algebras $(A,\m)$ over $R$ and morphisms homotopy classes $[\f]$ of morphisms $\f:(A,\m)\ra(B,\n)$, with composition $[\g]\ci[\f]=[\g\ci\f]$.

For cohomologically unital $A_\iy$-categories, we define the {\it homotopy category\/} $\Ho(\AAlg_R^\cu)$ of $\AAlg_R^\cu$ in the same way. 

For unital $A_\iy$-categories, we define the {\it homotopy category\/} $\Ho(\AAlg_R^\un)$ of $\AAlg_R^\un$ by taking equivalence classes of morphisms under {\it unital\/} homotopies.
\end{dfn}

\subsection{\texorpdfstring{$A_\iy$-algebras and differential graded algebras}{A∞-algebras and differential graded algebras}}
\label{ai35}

\begin{dfn}
\label{ai3def9}
A {\it differential graded algebra\/} ({\it dga\/}) over $R$ is a quadruple $(A,\cdot,\boo,\d)$, where $A=\bigop_{d\in\Z}A^d$ is a $\Z$-graded $R$-module, and $\cdot:A\t A\ra A$ is a graded associative $R$-bilinear product of degree 0, and $\boo$ is an identity for $\cdot$, and $\d:A\ra A$ is graded $R$-bilinear of degree 1 with $\d^2=0$, $\d(a\cdot b)=(\d a)\cdot b+(-1)^{\deg a}a\cdot(\d b)$ for pure $a,b\in A$, and $\d\boo=0$. We usually write dgas as $A$, leaving the data $\cdot,\boo,\d$ implicit, and sometimes write this as~$\cdot_A,\boo_A,\d_A$.

Define graded linear $\m_k:A^{\ot^k}\ra A$ for $k\ge 0$ by
\begin{equation*}
\m_1(a)=\d a,\; \m_2(a\ot b)=a\cdot b,\; \m_k=0,\; k\ne 1,2. 
\end{equation*}
Then $(A,\m)$ is an $A_\iy$-algebra. It has a unit $\boo$ in the sense of Definition~\ref{ai3def2}.

If $A,B$ are dgas over $R$, a {\it dga morphism}, or {\it morphism} $f:A\ra B$ is an $R$-linear map $f:A\ra B$ which preserves all the structure, so that $f$ is graded of degree 0 with $f(a\cdot a')=f(a)\cdot f(a')$, $\d_B f(a)=f(\d_A a)$, and $f(\boo_A)=\boo_B$. 

Then regarding $A,B$ as $A_\iy$-algebras, defining $\f_1=f$ and $\f_k=0$ for $k>1$, we find that $\f=(\f_k)_{k\ge 1}$ is an $A_\iy$-algebra morphism $\f:A\ra B$, which is unital in the sense of Definition \ref{ai3def5}. Note however that if $\f=(\f_k)_{k\ge 1}:A\ra B$ is a unital $A_\iy$-morphism, then $\f$ comes from a dga morphism $f=\f_1:A\ra B$ if and only if $\f_k=0$ for~$k>1$.

Write $\DGA_R$ for the category of dgas over $R$, with morphisms dga morphisms $f:A\ra B$. Also write $\DGA_R^{A_\iy}$ for the category of dgas over $R$, with morphisms unital $A_\iy$-algebra morphisms $f:A\ra B$. Then $\DGA_R\subset \DGA_R^{A_\iy}$ is a non-full subcategory. Also $\DGA_R^{A_\iy}\subset\AAlg_R^\un$ is a full subcategory.

As in Definition \ref{ai3def8}, we write $\Ho(\DGA_R^{A_\iy})\subset\Ho(\AAlg_R^\un)$ for the {\it homotopy category\/} of $\DGA_R^{A_\iy}$, whose morphisms are unital homotopy equivalence classes of morphisms in~$\DGA_R^{A_\iy}$.

A dga morphism $f:A\ra B$ is a {\it quasi-isomorphism\/} if $f:(A,\d_A)\ra(B,\d_B)$ is a quasi-isomorphism of complexes, that is, $H^k(f):H^k(A,\d_A)\ra H^k(B,\d_B)$ is an isomorphism for all $k\in\Z$. Also $f$ is a {\it chain homotopy equivalence\/} if $f:(A,\d_A)\ra(B,\d_B)$ is a chain homotopy equivalence of complexes. If $R$ is a field then every quasi-isomorphism is a chain homotopy equivalence.
\end{dfn}

\subsection{\texorpdfstring{Defining $A_\iy$-algebras using bar complexes}{Defining A∞-algebras using bar complexes}}
\label{ai36}

See Getzler--Jones \cite[\S 1]{GeJo} and Keller \cite[\S 3.6]{Kell1} for the material of this section.

\begin{dfn}
\label{ai3def10}
Let $(A,\m)$ be an $A_\iy$-algebra over $R$. Write $A[1]$ for $A$ with grading shifted by 1, that is, $A[1]^k=A^{k+1}$. Write $S[1]:A\ra A[1]$ for the identity map $\id_A$, thought of as  an $R$-linear map graded of degree $1$. Write $S[-1]:A[1]\ra A$ for its inverse, which is also $\id_A$. We call $S[\pm 1]$ {\it shift maps}.

Define $R$-linear maps $\bar\m_k:A[1]^{\ot^k}\ra A[1]$ for $k\ge 1$ by
\e
\bar\m_k=(-1)^{k(k-1)/2}S[1]\ci \m_k\ci S[-1]^{\ot^k}.
\label{ai3eq14}
\e
Since $\m_k$ is graded of degree $2-k$, $\bar\m_k$ is graded of degree 1. As $(S[-1]^{\ot^k})\ci(S[1]^{\ot^k})=(-1)^{k(k-1)/2}\id$ by \eq{ai3eq2}, equation \eq{ai3eq14} is equivalent to
\begin{equation*}
\m_k=S[-1]\ci\bar\m_k\ci S[1]^{\ot^k},
\end{equation*}

Here is a subtle point: $\m_k$ and $\bar\m_k$ are essentially the same map, as $S[\pm 1]$ are $\id_A$. However, when we translate equations in the $\m_k$ such as \eq{ai3eq3} into equations in the $\bar\m_k$, {\it there may be sign changes}, because of the Koszul sign rule in Convention \ref{ai3conv1}. For example, if $i,k\ge 0$ and $j\ge 1$ with $i+j+k=n$ then
\ea
&\bar\m_{i+k+1}\ci\bigl[\id_{A[1]}^{\ot^i}\ot\bar\m_j\ot \id_{A[1]}^{\ot^k}\bigr]
\label{ai3eq15}\\
&=(-1)^{(i+k+1)(i+k)/2}(-1)^{j(j-1)/2}S[1]\ci\m_{i+k+1}\ci S[-1]^{\ot^{i+k+1}}\ci
\nonumber\\
&\quad\bigl[(S[1]\ci\id_A\ci S[-1])^{\ot^i}\ot(S[1]\ci\m_j\ci S[-1]^{\ot^j})\ot (S[1]\ci\id_A\ci S[-1])^{\ot^k}\bigr]
\nonumber\\
&=(-1)^{n(n-1)/2+(j+1)(i+k)}(-1)^kS[1]\ci\m_{i+k+1}\ci \bigl[(S[-1]\ci S[1]\ci\id_A\ci S[-1])^{\ot^i} 
\nonumber\\
&\quad\ot(S[-1]\ci S[1]\ci\m_j\ci S[-1]^{\ot^j})\ot (S[-1]\ci S[1]\ci\id_A\ci S[-1])^{\ot^k}\bigr]
\nonumber\\
&=(-1)^{n(n-1)/2+i+ij+jk}S[1]\!\ci\!\m_{i+k+1}\!\ci\! \bigl[S[-1]^{\ot^i}\!\ot\!(\m_j\!\ci\! S[-1]^{\ot^j})\!\ot\!(S[-1])^{\ot^k}\bigr]
\nonumber\\
&=(-1)^{n(n-1)/2+i+ij+jk}(-1)^{i(2-j)}S[1]\!\ci\!\m_{i+k+1}\!\ci\!\bigl(\id_A^{\ot^i}\!\ot\!\m_j\!\ot\! \id_A^{\ot^k}\bigr)\!\ci\! S[-1]^{\ot^{i+j+k}}
\nonumber\\
&=(-1)^{n(n-1)/2}(-1)^{i+jk}S[1]\ci\m_{i+k+1}\ci\bigl(\id_A^{\ot^i}\ot\m_j\ot \id_A^{\ot^k}\bigr)\ci S[-1]^{\ot^{i+j+k}},
\nonumber
\ea
using \eq{ai3eq2} in the second and fourth steps.

Combining \eq{ai3eq3} and \eq{ai3eq15} shows that the $\bar\m_k$ satisfy the identities for $n\ge 1$
\e
\begin{gathered}
\sum_{i,k\ge 0,\; j\ge 1:\; i+j+k=n}
\bar\m_{i+k+1}\ci\bigl(\id_{A[1]}^{\ot^i}\ot\bar\m_j\ot \id_{A[1]}^{\ot^k}\bigr)=0,
\end{gathered}
\label{ai3eq16}
\e
as in Keller \cite[\S 3.6]{Kell1}. Note the sign change compared to \eq{ai3eq3}.

Now let $\f:(A,\m)\ra(B,\n)$ be a morphism of $A_\iy$-algebras. Define $R$-linear maps $\bar\f_k:A[1]^{\ot^k}\ra B[1]$ for $k\ge 1$ by
\begin{equation*}
\bar\f_k=(-1)^{k(k-1)/2}S_B[1]\ci\f_k\ci S_A[-1]^{\ot^k}.
\end{equation*}
Note that although $\f_k$ is graded of degree $1-k$, $\bar\f_k$ is graded of degree 0. Then as for \eq{ai3eq3} and \eq{ai3eq16}, equation \eq{ai3eq7} is equivalent to the equation for $n\ge 1$
\e
\begin{split}
&\sum_{i,k\ge 0,\; j\ge 1:\; i+j+k=n}\bar\f_{i+k+1}\bigl(
\id_{A[1]}^{\ot^i}\ot\bar\m_j\ot\id_{A[1]}^{\ot^k}\bigr) \\ 
&=\sum_{n=k_1+\cdots+k_l,\; k_a\ge 1}\bar\n_l\bigl(\bar\f_{k_1}\ot\bar\f_{k_2}\ot\cdots\ot\bar\f_{k_l}\bigr).
\end{split}
\label{ai3eq17}
\e

Next let $\f,\g:(A,\m)\ra(B,\n)$ be morphisms of $A_\iy$-algebras, and $\H:\f\Ra\g$ be a homotopy. Define $R$-linear maps $\bar\H_k:A[1]^{\ot^k}\ra B[1]$ for $k\ge 1$ by
\begin{equation*}
\bar\H_k=(-1)^{k(k-1)/2}S_B[1]\ci\H_k\ci S_A[-1]^{\ot^k}.
\end{equation*}
Then $\bar\H_k$ is graded of degree $-1$, and equation \eq{ai3eq12} is equivalent to the equation for $n\ge 1$
\begin{align*}
&\bar\f_n-\bar\g_n=\\
&\sum_{\begin{subarray}{l} l,m\ge 0,\; i_1,\ldots,i_l,j,k_1,\ldots,k_m\ge 1: \\
n=i_1+\cdots+i_l+j+k_1+\cdots+k_m\end{subarray}\!\!\!\!\!\!}
\bar\n_{l+m+1}\ci\bigl(\bar\f_{i_1}\ot\cdots\ot\bar\f_{i_l}\ot
\bar\H_j\ot\bar\g_{k_1}\ot\cdots\ot
\bar\g_{k_m}\bigr)\\
&+\sum_{i,k\ge 0,\; j\ge 1:\; i+j+k=n}\bar\H_{i+k+1}\bigl(\id_{A[1]}^{\ot^i}\ot\bar\m_j\ot\id_{A[1]}^{\ot^k}\bigr). 
\end{align*}

We can also use labelled trees $\bs\Ga$ as in \S\ref{ai61} to represent compositions of the shifted operators $\bar\m_j,\bar\f_j,\bar\H_j,\ldots$ above. Then we will write the compositions as $\bar C(\bs\Ga)$ rather than~$C(\bs\Ga)$.
\end{dfn}

As in Keller \cite[\S 3.6]{Kell1}, \S\ref{ai32}--\S\ref{ai34} can be rewritten in terms of {\it bar complexes}.

\begin{dfn}
\label{ai3def11}
Let $(A,\m)$ be an $A_\iy$-algebra over $R$. Write $A[1]$ for $A$ with grading shifted by 1, that is, $A[1]^k=A^{k+1}$. The {\it tensor coalgebra\/} $T(A[1])$ of $A[1]$ is $T(A[1])=\bigop_{n\ge 0}A^{\ot^n}$, where we write $A^{\ot^0}=R$. It is graded in the obvious way, so that
\begin{equation*}
T(A[1])^d=\bigop_{d_1+\cdots+d_n=d}A[1]^{d_1}\ot\cdots\ot A[1]^{d_n}=\bigop_{d_1+\cdots+d_n=d}A^{d_1+1}\ot\cdots\ot A^{d_n+1}.
\end{equation*}
 It has a {\it coproduct\/} $\co:T(A[1])\ra T(A[1])\ot T(A[1])$ given by
\begin{equation*}
\co(a_1\ot\cdots\ot a_n)=\ts\sum_{k=0}^n (a_1\ot\cdots\ot a_k)\ot(a_{k+1}\ot\cdots\ot a_n),
\end{equation*}
taking the $k=0$ and $k=n$ terms to be $1\ot(a_1\ot\cdots\ot a_n)$ and $(a_1\ot\cdots\ot a_n)\ot 1$. 

It has an increasing filtration $T(A[1])_{\le 0}\subset T(A[1])_{\le 1}\subset T(A[1])_{\le 2}\subset\cdots,$ with $T(A[1])_{\le k}=\bigop_{n=0}^kA^{\ot^n}$. The coproduct maps
\begin{equation*}
\co:T(A[1])_{\le k}\longra\sum_{i+j=k}T(A[1])_{\le i}\ot T(A[1])_{\le j}.
\end{equation*}

Define $R$-linear maps $\bar\m_k:A[1]^{\ot^k}\ra A[1]$ from $\m_k$ as in Definition \ref{ai3def10}. Define an $R$-linear map $\bar\d:T(A[1])\ra T(A[1])$ by, for all $n\ge 0$
\e
\bar\d\vert_{A[1]^{\ot^n}}=\sum_{i,k\ge 0,\; j\ge 1:\; i+j+k=n}\id_{A[1]}^{\ot^i}\ot\bar\m_j\ot \id_{A[1]}^{\ot^k}.
\label{ai3eq18}
\e
Then $\bar\d$ is a graded $R$-linear map of degree 1, and equation \eq{ai3eq16} is equivalent to $\bar\d\ci\bar\d=0$, so that $\bigl(T(A[1]),\bar\d\bigr)$ is a complex, the {\it bar complex} of $(A,\m)$. Also $\bar\d$ is a derivation for the coproduct $\co$, that is,
\e
\co\ci\bar\d=\bigl(\bar\d\ot\id_{T(A[1])}+\id_{T(A[1])}\ot\bar\d\bigr)\ci\co
\label{ai3eq19}
\e
so that $\bigl(T(A[1]),\co,\bar\d\bigr)$ is a {\it differential graded coalgebra}. And $\bar\d$ preserves the filtration, that is,
\e
\bar\d\bigl(T(A[1])_{\le k}\bigr)\subseteq T(A[1])_{\le k}.
\label{ai3eq20}
\e
\end{dfn}

\begin{dfn}
\label{ai3def12}
Let $\f:(A,\m)\ra(B,\n)$ be a morphism of $A_\iy$-algebras over $R$, as in Definition \ref{ai3def5}. Use the notation of Definition \ref{ai3def10}. We can rewrite \eq{ai3eq7} in terms of the bar complexes of $(A,\m)$ and $(B,\n)$: define an $R$-linear map $\bar\f:T(A[1])\ra T(B[1])$ by, for $n>0$
\e
\bar\f\vert_{A[1]^{\ot^n}}=\sum_{n=k_1+\cdots+k_l,\; k_a\ge 1}
\bar\f_{k_1}\ot\bar\f_{k_2}\ot\cdots\ot\bar\f_{k_l}.
\label{ai3eq21}
\e
When $n=0$, we set $\bar f(\la)=\la$ for $\la\in R=A^{\ot^0}=B^{\ot^0}$. Then $\bar\f$ is graded of degree 0, and \eq{ai3eq17} is equivalent to
\e
\bar\d_B\ci\bar\f=\bar\f\ci\bar\d_A:T(A[1])\longra T(B[1]),
\label{ai3eq22}
\e
that is, $\bar\f$ is a morphism of bar complexes $\bigl(T(A[1]),\bar\d_A\bigr)\ra\bigl(T(B[1]),\bar\d_B\bigr)$. It also intertwines the coproducts $\co_A,\co_B$ on $T(A[1]),T(B[1])$, that is
\begin{equation*}
\co_B\ci\bar f=(\bar f\ot\bar f)\ci\co_A,
\end{equation*}
and it preserves filtrations $T(A[1])_{\le k},T(B[1])_{\le k}$, namely
\e
\bar f\bigl(T(A[1])_{\le k}\bigr)\subseteq T(B[1])_{\le k}.
\label{ai3eq23}
\e

For identity morphisms $\id_A:(A,\m)\ra(A,\m)$ we have $\bar\id_A=\id_{T(A[1])}$ on bar complexes. For morphisms $\f:(A,\m)\ra(B,\n)$, $\g:(B,\n)\ra(C,{\mathfrak o})$, equation \eq{ai3eq10} implies that $\overline{(\g\ci\f)}=\bar\g\ci\bar\f$ on bar complexes.
\end{dfn}

\begin{dfn}
\label{ai3def13}
Let $\f,\g:(A,\m)\ra(B,\n)$ be morphisms of $A_\iy$-algebras over $R$ and $\H:\f\Ra\g$ be a homotopy, as in Definition \ref{ai3def6}.  Use the notation of Definition \ref{ai3def10}. Define an $R$-linear map $\bar\H:T(A[1])\ra T(B[1])$ by, for $n\ge 0$
\e
\bar\H\vert_{A[1]^{\ot^n}}=\sum_{\begin{subarray}{l} l,m\ge 0,\; i_1,\ldots,i_l,j,k_1,\ldots,k_m\ge 1: \\
n=i_1+\cdots+i_l+j+k_1+\cdots+k_m\end{subarray}\!\!\!\!\!\!}
\bar\f_{i_1}\ot\cdots\ot\bar\f_{i_l}\ot
\bar\H_j\ot\bar\g_{k_1}\ot\cdots\ot
\bar\g_{k_m}.
\label{ai3eq24}
\e
Then $\bar\H$ is a graded $R$-linear map of degree $-1$, and \eq{ai3eq12} is equivalent to
\e
\bar\f-\bar\g=\bar\d_B\ci\bar\H+\bar\H\ci\bar\d_A,
\label{ai3eq25}
\e
so that $\bar\f$ and $\bar\g$ are homotopic as morphisms of chain complexes in the usual sense.
Also $\bar\H$ satisfies
\e
\co_B\ci\bar\H=(\bar\f\ot\bar\H+\bar\H\ot\bar\g)\ci\co_A,
\label{ai3eq26}
\e
and $\bar\H$ preserves filtrations, that is
\e
\bar\H\bigl(T(A[1])_{\le k}\bigr)\subseteq T(B[1])_{\le k}.
\label{ai3eq27}
\e
\end{dfn}

These are equivalent to the previous definitions, in the following sense:

\begin{prop}
\label{ai3prop1}
{\bf(a)} Let\/ $R$ be a commutative ring and\/ $A$ be a $\Z$-graded\/ $R$-module. Then $A_\iy$-algebra structures $(A,\m)$ on $A$ are equivalent under \eq{ai3eq18} to $R$-linear differentials $\bar\d:T(A[1])\ra T(A[1])$ which are graded of degree $+1$ and satisfy $\bar\d\ci\bar\d=0$ and equations\/ {\rm\eq{ai3eq19}--\eq{ai3eq20}}. 
\smallskip

\noindent{\bf(b)} Let\/ $(A,\m),(B,\n)$ be $A_\iy$-algebras over $R$. Then morphisms $\f:(A,\m)\ra(B,\n)$ are equivalent under \eq{ai3eq21} to $R$-linear maps $\bar\f:T(A[1])\ra T(B[1])$ which are graded of degree $0$ and satisfy\/ {\rm\eq{ai3eq22}--\eq{ai3eq23}}.
\smallskip

\noindent{\bf(c)} Let\/ $\f,\g:(A,\m)\ra(B,\n)$ be morphisms of\/ $A_\iy$-algebras over $R$. Then homotopies $\H:\f\Ra\g$ are equivalent under \eq{ai3eq24} to $R$-linear maps $\bar\H:T(A[1])\ra T(B[1])$ which are graded of degree $-1$ and satisfy\/ {\rm\eq{ai3eq25}--\eq{ai3eq27}}.
\end{prop}

\begin{rem}
\label{ai3rem3}
{\bf(i)} Note that the definitions of $A_\iy$-algebras, morphisms and homotopies using bar complexes are shorter and simpler than the usual definitions. One could argue that bar complexes should be the primary definition.

\smallskip

\noindent{\bf(ii)} The conventional notation for coproducts in bar complexes is $\De$, not $\co$. We change this to avoid confusion with $\De$ in simplicial sets, as in Chapter~\ref{ai2}.
\end{rem}

\subsection{\texorpdfstring{Mazuir and Cazassus' work on $\iy$-categories of $A_\iy$-algebras}{Mazuir and Cazassus' work on ∞-categories of A∞-algebras}}
\label{ai37}
 
As part of a project to study the higher structure of Morse theory of a compact manifold \cite{Mazu1,Mazu2}, Mazuir \cite[\S 2.4]{Mazu2} discusses several attempts to construct an $\iy$-category of $A_\iy$-algebras, but he does not succeed in this goal. 

If $(A,\m)$ and $(B,\n)$ are $A_\iy$-algebras over $R$, Mazuir \cite[Th.~1]{Mazu2} defines a Kan complex $\cHom((A,\m),(B,\n))$, which may be related to mapping spaces $\Hom^L((A,\m),(B,\n))$ in our $\iy$-categories of $A_\iy$-algebras in Chapter~\ref{ai6}.

Cazassus \cite[\S 5.2]{Caza} defines an $\iy$-category $\mathop{\text{f-Alg}}_R$ of $A_\iy$-algebras over $R$. Part of his aim is to understand the sense in which Morse (co)homology of smooth manifolds $X,Y$ is functorial under smooth maps $\phi:X\ra Y$. Singular (co)homology, for example, is functorial in the obvious sense on (co)chains. But to define the Morse cochain complex $CM^\bu(X,R)$ of $X$ (which can be given the structure of an $A_\iy$-algebra) one must choose a Morse function $f:X\ra\R$ and a Palais--Smale pseudo-gradient $v$ for $f$ (and more choices for the $A_\iy$-structure), and these choices cannot be made functorial at the 1-categorical level.

Cazassus defines an $\iy$-functor $CM^\bu:\sN(\Man^\rop)\ra\mathop{\text{f-Alg}}_R$ which restores the functoriality at the $\iy$-categorical level. Here $CM^\bu$ makes an arbitrary choice of $f,v$ as above on each 0-simplex $X$ of $\sN(\Man^\rop)$, and on $N$-simplices for $N\ge 1$ makes more complicated choices of Morse-type data. It is necessary to define an $\iy$-category $\mathop{\text{f-Alg}}_R$ of $A_\iy$-algebras as the target of this functor.

Cazassus also defines $\iy$-categories $\mathop{\text{f-Coalg}}_R$ of `$A_\iy$-coalgebras', which are the target of an $\iy$-functor for Morse homology, and $\mathop{\text{f-Bialg}}_R$ of `$A_\iy$-bialgebras', which are the target of an $\iy$-functor for Morse cohomology of Lie groups. All these are intended for eventual applications in (equivariant) Lagrangian Floer cohomology and (equivariant) Fukaya categories, since Floer theories are (morally) infinite-dimensional Morse theories.

To define $\mathop{\text{f-Alg}}_R$, Cazassus defines an $R$-linear category $(R,u)\text{-mod}$ of {\it ascending foresty $R$-modules}, or $(R,u)$-{\it modules}, with objects graded $R$-modules $A,B$ and morphisms certain collections of graded $R$-linear maps $A^{\ot k}\ra B[l]$. He forms the dg-category $\Ch_{(R,u)\text{-mod}}$ of complexes of $(R,u)$-modules. As in Lurie \cite[\S 1.3.1]{Luri3} this has a {\it dg-nerve\/} $\sN_\dg(\Ch_{(R,u)\text{-mod}})$, which is an $\iy$-category.

Cazassus defines a faithful (but not full) inclusion functor $I:\AAlg_R\hookra\Ch_{(R,u)\text{-mod}}$. Then he sets $\mathop{\text{f-Alg}}_R$ to be the $\iy$-subcategory of $\sN_\dg(\Ch_{(R,u)\text{-mod}})$ whose 0- and 1-simplices lie in the image of $I$, but whose $N$-simplices for $N\ge 2$ are all $N$-simplices in $\sN_\dg(\Ch_{(R,u)\text{-mod}})$ with these 0- and 1-simplices.

\section{\texorpdfstring{The classical theory of $A_\iy$-categories}{The classical theory of A∞-categories}}
\label{ai4}

\subsection{\texorpdfstring{The definition of $A_\iy$-categories, and units}{The definition of A∞-categories, and units}}
\label{ai41}

Some references on $A_\iy$-categories are Fukaya \cite[\S 1]{Fuka3} and Seidel \cite{Seid}.

\begin{dfn}
\label{ai4def1}
Let $R$ be a commutative ring. An $A_\iy$-{\it category} $\sA$ over $R$ consists of:
\begin{itemize}
\setlength{\itemsep}{0pt}
\setlength{\parsep}{0pt}
\item[(a)] A set, or proper class, $\Obj\sA$ of {\it objects\/} of $\sA$.
\item[(b)] For all objects $x,y\in\Obj\sA$, a $\Z$-graded $R$-module
\begin{equation*}
\Hom(x,y)=\bigop\nolimits_{d\in\Z}\Hom^d(x,y).
\end{equation*}
We also write $\Hom_\sA(x,y)$ when we want to identify the $A_\iy$-category $\sA$.
\item[(c)] For all $k\ge 1$ and objects $x_0,x_1,\ldots,x_k$ in $\Obj\sA$, a graded $R$-linear map
\begin{equation*}
\!\!\!\!\!\m_k:\Hom(x_0,x_1)\ot_R\Hom(x_1,x_2)\ot_R\cdots \ot_R\Hom(x_{k-1},x_k)\longra\Hom(x_0,x_k)
\end{equation*}
of degree $2-k$. We also write $\m_k^\sA$, or $\m_k^{x_0,\ldots,x_k}$, or $\m_k^{\sA;x_0,\ldots,x_k}$ for $\m_k$ when we want to identify $x_0,\ldots,x_k$ or $\sA$. We write $\m=(\m_k)_{k\ge 1}$. 

Generalizing \eq{ai3eq3}, we require that for all $n\ge 1$ and objects $x_0,\ldots,x_n$\end{itemize}
\ea
\sum_{\begin{subarray}{l} i,k\ge 0,\; j\ge 1:\\ i+j+k=n\end{subarray}}\begin{aligned}[t]&
(-1)^{i+jk}\m_{i+k+1}^{x_0,\ldots,x_i,x_{i+j},\ldots,x_n}\ci\bigl(\id_{\Hom(x_0,x_1)}\ot\cdots\ot\id_{\Hom(x_{i-1},x_i)}\ot {}\\
&\m_j^{x_i,\ldots,x_{i+j}}\ot \id_{\Hom(x_{i+j},x_{i+j+1})}\ot\cdots\ot\id_{\Hom(x_{n-1},x_n)}\bigr)=0.
\end{aligned}
\nonumber\\[-17pt]
\label{ai4eq1}
\ea
By abuse of notation, we also regard $\m_k$ in (c) as a graded $R$-multilinear map 
\begin{equation*}
\m_k:\Hom(x_0,x_1)\t\Hom(x_1,x_2)\t\cdots \t\Hom(x_{k-1},x_k)\longra\Hom(x_0,x_k).
\end{equation*}
Then generalizing \eq{ai3eq4}, by Convention \ref{ai3conv1}, equation \eq{ai4eq1} is equivalent to the equation for all $n\ge 1$ and all pure $a_i\in\Hom(x_{i-1},x_i)$ for $i=1,\ldots,n$
\e
\begin{gathered}
\sum_{i,k\ge 0,\; j\ge 1:\; i+j+k=n}
\begin{aligned}[t]
&(-1)^{i+jk+j\sum_{l=1}^i\deg a_l}\,\m_{i+k+1}(a_1,\ldots,a_i,\\
&\qquad\m_j(a_{i+1},\ldots,a_{i+j}),a_{i+j+1},\ldots,a_n)=0.\end{aligned}
\end{gathered}
\label{ai4eq2}
\e

An $A_\iy$-category $\sA$ for which $\Obj\sA$ is a set (rather than a proper class) is called a {\it small\/} $A_\iy$-category.
\end{dfn}

\begin{rem}
\label{ai4rem1}
{\bf(On Set Theory foundations.)} As in Remark \ref{ai1rem2}, to study $A_\iy$-categories $\sA$ using $\iy$-categories, we will need $\Obj\sA$ to be a set, so $\sA$ is a small $A_\iy$-category. In this book we fix a Grothendieck universe $\sU$, and assume implicitly from now on that $\Obj\sA$ is a subset of~$\sU$.
\end{rem}

\begin{rem}
\label{ai4rem2}
{\bf(Relations between $A_\iy$-categories and $A_\iy$-algebras.)} If $\Obj\sA=\{*\}$ is a single point then Definition \ref{ai4def1} reduces to Definition \ref{ai3def1}. That is, $A_\iy$-algebras are essentially the same as $A_\iy$-categories with one object.
\end{rem}

\begin{rem}
\label{ai4rem3}
{\bf(Other remarks.) (a)} Note that $A_\iy$-categories are not actually categories, as by \eq{ai4eq1} the composition maps
\begin{equation*}
\m_2^{x,y,z}:\Hom(x,y)\t\Hom(y,z)\longra\Hom(x,z)
\end{equation*}
need not be associative, and also we do not require the existence of identities (units). In Definition \ref{ai4def4} we will show that a (cohomologically) unital $A_\iy$-category $\sA$ has a {\it cohomology category\/} $H^0(\sA)$ which is an honest category.

\smallskip

\noindent{\bf(b)} $A_\iy$-categories are important in Symplectic Geometry: Fukaya categories of symplectic manifolds are initially defined as $A_\iy$-categories \cite{Fuka3,Seid}. They also have applications elsewhere. For example, dg-categories are special cases of $A_\iy$-categories, and to study the homotopy theory of dg-categories it is helpful to regard them as $A_\iy$-categories, with $A_\iy$-functors between them.
\end{rem}

If $\sA$ is an $A_\iy$-category then \eq{ai4eq1} for $n=1$ becomes $\m_1^{x_0,x_1}\ci\m_1^{x_0,x_1}=0$, with $\m_1^{x_0,x_1}$ of degree 1. Thus $\bigl(\Hom(x_0,x_1),\m_1^{x_0,x_1}\bigr)$ is a complex, and we form its cohomology $H^*(\Hom(x_0,x_1))$ by
\e
H^p(\Hom(x_0,x_1))=\frac{\Ker\bigl(\m_1^{x_0,x_1}:\Hom^p(x_0,x_1)\ra \Hom^{p+1}(x_0,x_1)\bigr)}{\Im\bigl(\m_1^{x_0,x_1}:\Hom^{p-1}(x_0,x_1)\ra \Hom^p(x_0,x_1)\bigr)}\,.
\label{ai4eq3}
\e
As in \eq{ai3eq5}, for all objects $x_0,x_1,x_2$ in $\sA$ we have a well-defined, associative $R$-bilinear product, which we call {\it composition\/} $\ci$
\e
\ci:H^p(\Hom(x_0,x_1))\t H^q(\Hom(x_1,x_2))\longra H^{p+q}(\Hom(x_0,x_2))
\label{ai4eq4}
\e
given by
\begin{equation*}
(a_1+\Im\m^{x_0,x_1}_1)\ci(a_2+\Im\m^{x_1,x_2}_1)=\m^{x_0,x_1,x_2}_2(a_1,a_2)+\Im\m^{x_0,x_2}_1.
\end{equation*}

Here are the analogues of Definitions \ref{ai3def2}--\ref{ai3def4}.

\begin{dfn}
\label{ai4def2}
Let $\sA$ be an $A_\iy$-category. We say that $\sA$ is {\it unital\/} if for all $x\in\Obj\sA$ there exists $\boo_x\in\Hom^0(x,x)$, called the {\it identity\/} for $x$, satisfying
\e
\m_2^{x_0,x_0,x_1}(\boo_{x_0}\ot a)=\m_2^{x_0,x_1,x_1}(a\ot \boo_{x_1})=a
\label{ai4eq5}
\e
for all objects $x_0,x_1$ in $\sA$ and $a\in\Hom(x_0,x_1)$, and
\e
\m_{k+1}^{x_0,\ldots,x_i,x_i,\ldots,x_k}(a_1\ot\cdots\ot a_i\ot \boo_{x_i}\ot a_{i+1}\ot\! \cdots \ot a_k)=0
\label{ai4eq6}
\e
for all $k\ne 1$, $i=0,\ldots,k$, $x_0,\ldots,x_k\in\Obj\sA$ and $a_j\in\Hom(x_{j-1},x_j)$.

If $\boo_x,\boo_x'$ are alternative choices for $\boo_x$ then \eq{ai4eq5} gives $\boo_x=\m_2^{x,x}(\boo_x,\boo_x')=\boo_x'$, so identities are unique if they exist, and to be unital is a property of $\sA$, not extra data.
\end{dfn}

\begin{dfn}
\label{ai4def3}
Let $\sA$ be an $A_\iy$-category over a commutative ring $R$. We call $\sA$ {\it weakly cohomologically unital\/} if for all $x\in\Obj\sA$ there exists $1_x$ in $H^0(\Hom(x,x))$, called the {\it cohomological identity\/} for $x$, satisfying 
\e
1_{x_0}\ci \al=\al\ci 1_{x_1}=\al
\label{ai4eq7}
\e
for all objects $x_0,x_1$ in $\sA$ and $\al\in H^*(\Hom(x_0,x_1))$, where $\ci$ is as in \eq{ai4eq4}.

Again, if $1_x,1_x'$ are alternative choices for $1_x$ then \eq{ai4eq7} gives $1_x=1_x\ci 1_x'=1_x'$, so cohomological identities are unique if they exist, and to be cohomologically unital is a property of $\sA$, not extra data.

Let $\sA$ be weakly cohomologically unital. The {\it cohomology category\/} $H^0(\sA)$ has objects $x\in\Obj\sA$, morphisms $\Hom_{H^0(\sA)}(x,y)=H^0(\Hom_\sA(x,y))$ for $x,y\in\Obj\sA$, composition of morphisms $\ci:H^0(\Hom(x,y))\!\t\! H^0(\Hom(y,z))\ab\ra H^0(\Hom(x,z))$ as in \eq{ai4eq4}, and identities $1_x\in H^0(\Hom(x,x))$. Then $H^0(\sA)$ is an ordinary $R$-linear category.
\end{dfn}

\begin{dfn}
\label{ai4def4}
Let $\sA$ be an $A_\iy$-category over a commutative ring $R$, and $x$ an object in $\sA$. A {\it cohomological unit for\/} $x$ is an element $e_x$ in $\Hom_\sA^0(x,x)$ such that $\m_1^{x,x}(e_x)=0$, and for all objects $y$ in $\sA$ the following maps of chain complexes are chain homotopic to the identity:
\begin{align*}
\Phi&:\bigl(\Hom_\sA(x,y),\m_1^{x,y}\bigr)\longra\bigl(\Hom_\sA(x,y),\m_1^{x,y}\bigr),& \Phi&:a\longmapsto \m_2^{x,y}(e_x\ot a), \\
\Psi&:\bigl(\Hom_\sA(y,x),\m_1^{y,x}\bigr)\longra\bigl(\Hom_\sA(y,x),\m_1^{y,x}\bigr),& \Psi&:a\longmapsto \m_2^{y,x}(a\ot e_x).
\end{align*}
Then the cohomology class $1_x=[e_x]$ in $H^0(\Hom(x,x))$ satisfies~\eq{ai4eq7}.

If $e_x'\in\Hom_\sA^0(x,x)$ is any other representative for $1_x$ in $H^0(\Hom(x,x))$, then $e_x'=e_x+\m_1^{x,x}(f)$ for $f\in\Hom_\sA^{-1}(x,x)$, and defining $\Th,\Xi$ by $\Th(a)=\m_2^{x,y}(f\ot a)$, $\Xi(a)=-\m_2^{y,x}(a\ot f)$, we see from \eq{ai3eq4} that
\begin{equation*}
\Phi'(a)\!=\!\Phi(a)\!+\!\m_1^{x,y}\ci\Th(a)\!-\!\Th\!\ci\!\m_1^{x,y}(a),\; \Psi'(a)\!=\!\Psi(a)\!+\!\m_1^{y,x}\ci\Xi(a)\!-\!\Xi\!\ci\!\m_1^{y,x}(a).
\end{equation*}
Thus, $\Phi',\Psi'$ are chain homotopic to the identity if and only if $\Phi,\Psi$ are, and $e_x'$ is a cohomological unit if and only if $e_x$ is. Thus, for $x$ to have a cohomologically unit is a property of $\sA$, not of the particular choice of element~$e_x$.

We call $\sA$ {\it cohomologically unital\/} if every $x\in\sA$ has a cohomological unit.

Clearly, $\sA$ unital implies $\sA$ cohomologically unital, implies $\sA$ weakly cohomologically unital. If $R$ is a field then weakly cohomologically unital and cohomologically unital coincide, since $\Phi,\Psi$ are chain homotopic to the identity if and only if they induce the identity maps on cohomology. 
\end{dfn}

Readers are warned that the names used for different kinds of unitality in $A_\iy$-categories vary widely. Our notion of unital is often called `strictly unital'. Our notion of cohomologically unital is sometimes just called `unital'. 

The next definition follows Fukaya~\cite[Def.~6.22]{Fuka3}.

\begin{dfn}
\label{ai4def5}
Suppose $\sA$ is a unital or cohomologically unital $A_\iy$-category over $R$, and $x,y$ are objects in $\sA$, and $f\in\Hom^0(x,y)$ with $\m_1^{x,y}(f)=0$, so that $f$ defines a morphism $[f]:x\ra y$ in the cohomology category $H^0(\sA)$. We call $f$ a {\it homotopy equivalence\/} if $[f]$ is an isomorphism in $H^0(\sA)$. We say that $x,y$ are {\it homotopy equivalent\/} if there exists a homotopy equivalence between them.

Homotopy equivalence is an equivalence relation on objects in~$\sA$.
\end{dfn}

\subsection{\texorpdfstring{Functors between $A_\iy$-categories}{Functors between A∞-categories}}
\label{ai42}

Here is the analogue of Definition \ref{ai3def5}.

\begin{dfn}
\label{ai4def6}
Let $\sA,\sB$ be $A_\iy$-categories over a commutative ring $R$. An $A_\iy$-{\it functor}, or just {\it functor}, $\fF:\sA\ra\sB$ consists of:
\begin{itemize}
\setlength{\itemsep}{0pt}
\setlength{\parsep}{0pt}
\item[(a)] A map $\fF_\Obj:\Obj\sA\ra\Obj\sB$. When this will not cause confusion, for $x\in\Obj\sA$ we write $\fF(x)\in\Obj\sB$ rather than $\fF_\Obj(x)$.
\item[(b)] For all $k\ge 1$ and objects $x_0,x_1,\ldots,x_k$ in $\Obj\sA$, a graded $R$-linear map
\begin{align*}
\fF_k:\Hom_\sA(x_0,x_1)\ot_R\Hom_\sA(x_1,x_2)\ot_R\cdots \ot_R\Hom_\sA(x_{k-1},x_k)&\\
\longra\Hom_\sB(\fF_\Obj(x_0),\fF_\Obj(x_k))&,
\end{align*}
of degree $1-k$. We also write $\fF_k^{x_0,\ldots,x_k}$ for $\fF_k$ when we want to identify $x_0,\ldots,x_k$. Generalizing \eq{ai3eq7}, these must satisfy for all $n\ge 1$ and objects $x_0,\ldots,x_n$ in~$\Obj\sA$
\end{itemize}
\ea
&\sum_{\begin{subarray}{l} i,k\ge 0,\; j\ge 1:\\ i+j+k=n \end{subarray}}
\begin{aligned}[t]&
(-1)^{i+jk}\fF_{i+k+1}^{x_0,\ldots,x_i,x_{i+j},\ldots,x_n}\ci\bigl(\id_{\Hom_\sA(x_0,x_1)}\ot\cdots\ot\id_{\Hom_\sA(x_{i-1},x_i)}\\
&{}\ot\m_j^{\sA;x_i,\ldots,x_{i+j}}\ot \id_{\Hom_\sA(x_{i+j},x_{i+j+1})}\ot\cdots\ot\id_{\Hom_\sA(x_{n-1},x_n)}\bigr)\end{aligned} 
\nonumber\\ 
&=\sum_{\begin{subarray}{l}  n=n_1+\cdots+n_k,\\ n_a\ge 1\end{subarray}}
\begin{aligned}[t]&(-1)^{\sum_{i=1}^k(n_i+1)(k-i)}\m^{\sB;\fF_\Obj(x_0),\fF_\Obj(x_{n_1}),\ldots,\fF_\Obj(x_n),}_k\bigl(\fF_{n_1}^{x_0,\ldots,x_{n_1}}\\
&{}\ot\fF_{n_2}^{x_{n_1},\ldots,x_{n_1+n_2}}\ot\cdots\ot\fF_{n_k}^{x_{n-n_k},\ldots,x_n}).\end{aligned} 
\nonumber\\[-16pt] 
\label{ai4eq8}
\ea

The {\it identity functor\/} $\id^\sA:\sA\ra\sA$ has $\id^\sA_\Obj=\id_{\Obj\sA}$, and $\id^{\sA,x}_1=\id:\Hom_\sA(x,x)\ra\Hom_\sA(x,x)$, and $\id^\sA_k=0$ for $k>1$.

When $n=1$, equation \eq{ai4eq8} becomes
\begin{align*}
&\fF^{x_0,x_1}_1\ci\m^{\sA;x_0,x_1}_1\!=\!\m_1^{\sB;\fF(x_0),\fF(x_1)}\ci\fF^{x_0,x_1}_1:\\
&\qquad \Hom_\sA(x_0,x_1)\longra \Hom_\sB(\fF(x_0),\fF(x_1)).
\end{align*}
Thus we have a morphism of complexes
\begin{equation*}
\fF^{x_0,x_1}_1:(\Hom_\sA(x_0,x_1),\m^{\sA;x_0,x_1}_1)\longra(\Hom_\sB(\fF(x_0),\fF(x_1)),\m_1^{\sB;\fF(x_0),\fF(x_1)}),
\end{equation*}
which induces a morphism of cohomology groups 
\e
H^*(\fF^{x_0,x_1}_1):H^*\bigl(\Hom_\sA(x_0,x_1)\bigr)\longra H^*\bigl(\Hom_\sB(\fF(x_0),\fF(x_1))\bigr).
\label{ai4eq9}
\e

If $\fF:\sA\ra\sB$, $\fG:\sB\ra\sC$ are functors of $A_\iy$-categories over $R$, the {\it composition} $\fG\ci\fF:\sA\ra\sC$ is given by $(\fG\ci\fF)_\Obj=\fG_\Obj\ci\fF_\Obj$ and
\ea
&(\fG\ci\fF)_n^{x_0,\ldots,x_n}=
\label{ai4eq10}\\
&\sum_{\begin{subarray}{l} n=n_1+\cdots+n_k,\\ n_a\ge 1\end{subarray}}
\begin{aligned}[t]& (-1)^{\sum_{i=1}^k(n_i+1)(k-i)}\fG_k^{\fF(x_0),\fF(x_{n_1}),\fF(x_{n_1+n_2}),\ldots,\fF(x_n)} \\ 
&\bigl(\fF_{n_1}^{x_0,\ldots,x_{n_1}}\ot\fF_{n_2}^{x_{n_1},\ldots,x_{n_1+n_2}}\ot\cdots\ot\fF_{n_k}^{x_{n-n_k},\ldots,x_n}),
\end{aligned}
\nonumber
\ea
for all $n\ge 1$ and $x_0,\ldots,x_n\in\Obj\sA$. Composition is associative.

Write $\ACat_R$ for the category with objects $A_\iy$-categories $\sA$ over $R$ and morphisms functors $\fF:\sA\ra\sB$.

If $\sA,\sB$ are unital $A_\iy$-categories over $R$, as in Definition \ref{ai4def2}, we call a functor $\fF:\sA\ra\sB$ {\it unital\/} if
\e
\fF_1^{x,x}(\boo_x)=\boo_{\fF_\Obj(x)}
\label{ai4eq11}
\e
for all objects $x$ in $\sA$, and
\e
\fF_{k+1}^{x_0,\ldots,x_i,x_i,\ldots,x_k}(a_1\ot\cdots\ot a_i\ot \boo_{x_i}\ot a_{i+1}\ot \cdots \ot a_k)=0
\label{ai4eq12}
\e
for all $k\ge 1$, $x_0,\ldots,x_k\in\Obj\sA$, $i=0,\ldots,k$ and $a_j\in\Hom_\sA(x_{j-1},x_j)$ for $j=1,\ldots,k$. Unital $A_\iy$-categories over $R$ and unital functors form a non-full subcategory $\ACat_R^\un\subset\ACat_R$.

Let $\sA,\sB$ be (weakly) cohomologically unital $A_\iy$-categories, as in Definitions \ref{ai4def3} and \ref{ai4def4}, and $\fF:\sA\ra\sB$ be a functor. We call $\fF$ ({\it weakly\/}) {\it cohomologically unital\/} if $H^0(\fF_1^x)(1_x)=1_{\fF_\Obj(x)}$ for all $x\in\Obj\sA$. Then (weakly) cohomologically unital $A_\iy$-categories over $R$ and (weakly) cohomologically unital functors form non-full subcategories $\ACat_R^\cu\subset\ACat_R^\wc \subset\ACat_R$. When $\sA,\sB$ are unital, unital functors $\fF:\sA\ra\sB$ are cohomologically unital, so~$\ACat_R^\un\subset\ACat_R^\cu$.

If $\fF:\sA\ra\sB$ is a (cohomologically) unital functor of (cohomologically) unital $A_\iy$-categories over $R$, the morphisms \eq{ai4eq9} induce an $R$-linear functor
\begin{equation*}
H^0(\fF_1):H^0(\sA)\longra H^0(\sB)
\end{equation*}
of ordinary $R$-linear categories.
\end{dfn}

\subsection{\texorpdfstring{The $A_\iy$-category of $A_\iy$-functors}{The A∞-category of A∞-functors}}
\label{ai43}

Let $\sA,\sB$ be $A_\iy$-categories over $R$. Then one can define an $A_\iy$-category $\Fun_{A_\iy}(\sA,\sB)$ whose objects are $A_\iy$-functors $\fF:\sA\ra\sB$. This idea seems to be due to Kontsevich circa 1998, and accounts can be found in Fukaya \cite[Th.-Def.~7.55]{Fuka3}, Lef\`evre-Hasegawa \cite[\S 8]{LeHa}, Lyubashenko \cite[\S 5]{Lyub}, and Seidel \cite[\S 1d, \S 2e]{Seid}. See also Keller \cite[\S 5]{Kell2}, and Bespalov--Lyubashenko--Manzyuk \cite[Prop.~8.15]{BLM} who give a version for $A_\iy$-multifunctors $\fF:\sA_1\t\cdots\t\sA_n\ra\sB$. Our presentation is based on \cite{Fuka3,LeHa}, but these are written in terms of bar complexes as in \S\ref{ai36}, so we have used the method of \S\ref{ai36} to compute the signs in \eq{ai4eq15} and \eq{ai4eq17} from the (much simpler) signs in~\cite{Fuka3,LeHa}.

\begin{dfn}
\label{ai4def7}
Let $\sA,\sB$ be $A_\iy$-categories over a commutative ring $R$. We will define an $A_\iy$-category $\Fun_{A_\iy}(\sA,\sB)$. Objects of $\Fun_{A_\iy}(\sA,\sB)$ are $A_\iy$-functors $\fF:\sA\ra\sB$, as in \S\ref{ai42}. For objects $\fF,\fG$, we define a graded $R$-module
\begin{equation*}
\Hom_{\Fun_{A_\iy}(\sA,\sB)}(\fF,\fG)=\Hom_{A_\iy}(\fF,\fG)=\bigop_{l\in\Z}\Hom^l_{A_\iy}(\fF,\fG)
\end{equation*}
as follows. An element $\fT$ of $\Hom^l_{A_\iy}(\fF,\fG)$ consists of, for all $k\ge 0$ and objects $x_0,\ldots,x_k$ in $\sA$, a graded $R$-linear map of degree $l-k$
\ea
\fT_k=\fT^{x_0,\ldots,x_k}_k&:\Hom_\sA(x_0,x_1)\ot_R\Hom_\sA(x_1,x_2)\ot_R\cdots \ot_R\Hom_\sA(x_{k-1},x_k)
\nonumber\\
&\qquad\longra\Hom_\sB(\fF_\Obj(x_0),\fG_\Obj(x_k)).
\label{ai4eq13}
\ea
Note that here we allow $k=0$, in contrast to Definitions \ref{ai4def1} and \ref{ai4def6} when we consider $\m_k,\fF_k$ for $k\ge 1$, and $\fT_0^{x_0}$ maps $R\ra\Hom_\sB^l(\fF(x_0),\fG(x_0))$. Thus $\fT_0^{x_0}$ is equivalent to an element $\fT_0^{x_0}(1)\in \Hom_\sB^l(\fF(x_0),\fG(x_0))$.

Such $\fT$ are called {\it pre natural transformations from $\fF$ to\/} $\fG$ \cite[Def.~7.49]{Fuka3}. We write the $A_\iy$-operations on the $A_\iy$-category $\Fun_{A_\iy}(\sA,\sB)$ as $\n_p$ for $p\ge 1$. Firstly, when $p=1$ we define
\e
\n_1=\n_1^{\fF,\fG}:\Hom_{A_\iy}(\fF,\fG)\longra \Hom_{A_\iy}(\fF,\fG)
\label{ai4eq14}
\e
by $\n_1(\fT)=\fU$, where if $n\ge 0$ and $x_0,\ldots,x_n$ are objects in $\sA$ then we define
\ea
&\fU_n^{x_0,\ldots,x_n}=
\label{ai4eq15}\\
&\sum_{\begin{subarray}{l} 
0=l_0\le l_1\le \cdots\le l_k=n,\\
1\le j\le k: \\
l_i>l_{i-1},\; i\ne 0,j \end{subarray}}
\begin{aligned}[t]&
(-1)^{(k+1)(\deg\fT+1)}\cdot(-1)^{\sum_{i=1}^{j-1}(l_i-l_{i-1}+1)(\deg\fT+1)}\cdot{}\\
&(-1)^{\sum_{i=1}^k(l_i-l_{i-1}+1)(k-i)} \cdot
\m^{\begin{subarray}{l} \sB;\fF(x_{l_0}),\fF(x_{l_1}),\fF(x_{l_2}),\ldots,\fF(x_{l_{j-1})},\\
\fG(x_{l_j}),\fG(x_{l_{j+1}}),\ldots,\fG(x_{l_k})\end{subarray}}_k\!\ci{} \\ 
&\bigl(\fF_{l_1-l_0}^{x_{l_0},\ldots,x_{l_1}}\ot\fF_{l_2-l_1}^{x_{l_1},\ldots,x_{l_2}}\ot\cdots\ot\fF_{l_{j-1}-l_{j-2}}^{x_{l_{j-2}},\ldots,x_{l_{j-1}}} \\
&{}\ot\fT_{l_j-l_{j-1}}^{x_{l_{j-1}},\ldots,x_{l_j}}\ot\fG_{l_{j+1}-l_j}^{x_{l_j},\ldots,x_{l_{j+1}}}\ot\cdots\ot\fG_{l_k-l_{k-1}}^{x_{l_{k-1}},\ldots,x_{l_k}}\bigr)
\end{aligned}
\nonumber\\
&+\sum_{\begin{subarray}{l} i,k\ge 0,\; j\ge 1:\\ i+j+k=n \end{subarray}}\begin{aligned}[t]&(-1)^{\deg\fT}\cdot(-1)^{i+jk}\cdot{}\\
&\fT^{x_0,\ldots,x_i,x_{i+j},\ldots,x_n}_{i+k+1}\ci\bigl(\id_{\Hom_\sA(x_0,x_1)}\ot\cdots\ot\id_{\Hom_\sA(x_{i-1},x_i)}\\
&{}\ot\m_j^{\sA;x_i,\ldots,x_{i+j}}\ot \id_{\Hom_\sA(x_{i+j},x_{i+j+1})}\ot\cdots\ot\id_{\Hom_\sA(x_{n-1},x_n)}\bigr). \end{aligned}
\nonumber
\ea

Secondly, let $p\ge 2$ and $\fF^0,\fF^1,\ldots,\fF^p:\sA\ra\sB$ be $A_\iy$-functors. We define the $A_\iy$-operation $\n_p$ in $\Fun_{A_\iy}(\sA,\sB)$,
\e
\begin{split}
\n_p=\n_p^{\fF^0,\ldots,\fF^p}:\Hom_{A_\iy}(\fF^0,\fF^1)\ot_R\cdots\ot_R\Hom_{A_\iy}(\fF^{p-1},\fF^p)&\\
\longra \Hom_{A_\iy}(\fF^0,\fF^p)&
\end{split}
\label{ai4eq16}
\e
by, if $\fT^i\in\Hom_{A_\iy}(\fF^{i-1},\fF^i)$ for $i=1,\ldots,p$, then $\n_p^{\fF^0,\ldots,\fF^p}(\fT^1\ot\cdots\ot\fT^p)\!=\!\fU$ in $\Hom_{A_\iy}(\fF^0,\fF^p)$, where if $n\ge 0$ and $x_0,\ldots,x_n$ are objects in $\sA$ then we define
\ea
&\fU_n^{x_0,\ldots,x_n}=
\label{ai4eq17}\\
&\sum_{\begin{subarray}{l} 
0=l_0\le l_1\le \\
\cdots\le l_k=n,\\
1\le a_1<a_2< \\
\cdots<a_p\le k: \\
l_j>l_{j-1},\\ j\ne 0,a_i 
\end{subarray}}
\begin{aligned}[t]&
(-1)^{\sum_{i=1}^p(k+1)(\deg\fT^i+1)}\cdot(-1)^{\sum_{i=1}^p\sum_{j=1}^{a_i-1}(l_j-l_{j-1}+1)(\deg\fT^i+1)}\cdot{}\\
&(-1)^{\sum_{i=1}^k(l_i-l_{i-1}+1)(k-i)} \cdot{} \\
&\m_k^{\begin{subarray}{l} \sB;\fF^0(x_{l_0}),\fF^0(x_{l_1}),\ldots,\fF^0(x_{l_{a_1-1}}),
\fF^1(x_{l_{a_1}}),\fF^1(x_{l_{a_1+1}}),\\
\ldots,\fF^1(x_{l_{a_2-1}}),\ldots,\fF^p(x_{l_{a_p}}),\fF^p(x_{l_{a_p+1}}),\ldots,\fF^p(x_{l_k})
\end{subarray}}\ci\\
&\bigl(\fF_{l_1-l_0}^{0;x_{l_0},\ldots,x_{l_1}}\ot\fF_{l_2-l_1}^{0;x_{l_1},\ldots,x_{l_2}}\ot\cdots\ot\fF_{l_{a_1-1}-l_{a_1-2}}^{0;x_{l_{a_1-2}},\ldots,x_{l_{a_1-1}}}\\
&{}\ot\fT_{l_{a_1}-l_{a_1-1}}^{1;x_{l_{a_1-1}},\ldots,x_{l_{a_1}}}\ot\fF_{l_{a_1+1}-l_{a_1}}^{1;x_{l_{a_1}},\ldots,x_{l_{a_1+1}}}\!\ot\!\cdots\!\ot\!\fF_{l_{a_2-1}-l_{a_2-2}}^{1;x_{l_{a_2-2}},\ldots,x_{l_{a_2-1}}} \\
&{}\ot
\fT_{l_{a_2}-l_{a_2-1}}^{2;x_{l_{a_2-1}},\ldots,x_{l_{a_2}}}\ot \cdots \ot\fT_{l_{a_p}-l_{a_p-1}}^{p;x_{l_{a_p-1}},\ldots,x_{l_{a_p}}} {}\\
&\ot\fF_{l_{a_p+1}-l_{a_p}}^{p;x_{l_{a_p}},\ldots,x_{l_{a_p+1}}}\ot\cdots\ot\fF_{l_k-l_{k-1}}^{p;x_{l_{k-1}},\ldots,x_{l_k}}\bigr).
\end{aligned}
\nonumber
\ea
When $p=1$, equation \eq{ai4eq17} for $(\fF^0,\fF^1)=(\fF,\fG)$ reduces to the first sum in \eq{ai4eq15}, but omits the second sum, which is why we treat $p=1$ separately. As in \cite{Fuka2,LeHa,Lyub,Seid}, the above make $\Fun_{A_\iy}(\sA,\sB)$ into an $A_\iy$-category over~$R$.

Next suppose $\sA,\sB$ are {\it cohomologically unital\/} $A_\iy$-categories. Following \cite[\S 2e]{Seid}, define $\Fun_{A_\iy}^\cu(\sA,\sB)$ to be the full $A_\iy$-subcategory of $\Fun_{A_\iy}(\sA,\sB)$ with objects cohomologically unital $A_\iy$-functors $\fF:\sA\ra\sB$. Then $\Fun_{A_\iy}^\cu(\sA,\sB)$ is a cohomologically unital $A_\iy$-category by Lyubashenko~\cite[Prop.~7.7]{Lyub}.
\end{dfn}

\begin{dfn}
\label{ai4def8}
Suppose $\sA,\sB$ are {\it unital\/} $A_\iy$-categories over $R$. Following \cite[\S 8.1.3]{LeHa}, define a non-full $A_\iy$-subcategory $\Fun_{A_\iy}^\un(\sA,\sB)$ of $\Fun_{A_\iy}(\sA,\sB)$ to have objects unital $A_\iy$-functors $\fF:\sA\ra\sB$, and if $\fF,\fG:\sA\ra\sB$ are unital functors, we define
\begin{equation*}
\Hom_{\Fun_{A_\iy}^\un(\sA,\sB)}(\fF,\fG)=\Hom_{A_\iy}^\un(\fF,\fG)\subset \Hom_{\Fun_{A_\iy}(\sA,\sB)}(\fF,\fG)
\end{equation*}
to be the subcomplex of $\fT$ in $\Hom_{A_\iy}(\fF,\fG)$ such that
\e
\fT_{k+1}^{x_0,\ldots,x_i,x_i,\ldots,x_k}(a_1\ot\cdots\ot a_i\ot \boo_{x_i}\ot a_{i+1}\ot \cdots \ot a_k)=0
\label{ai4eq18}
\e
for all $0\le i\le k$, $x_0,\ldots,x_k\in\Obj\sA$ and $a_j\in\Hom(x_{j-1},x_j)$ for $j=1,\ldots,k$. As in \cite[\S 8.1.3]{LeHa}, these are closed under the operations $\n_p^{\fF^0,\ldots,\fF^p}$, so $\Fun_{A_\iy}^\un(\sA,\sB)$ is an $A_\iy$-category.

In fact $\Fun_{A_\iy}^\un(\sA,\sB)$ is a {\it unital\/} $A_\iy$-category. For each unital $A_\iy$-functor $\fF:\sA\ra\sB$, the identity is $\boo_\fF=\fS$, where
\e
\begin{split}
&\fS_0^{x_0}:R\longra \Hom_\sB(\fF(x_0),\fF(x_0)), \quad\fS^{x_0}_0:\la\longmapsto \la\boo_{\fF(x_0)},\\
&\fS_k^{x_0,\ldots,x_k}=0\quad\text{for all $k\ge 1$ and $x_0,\ldots,x_k\in\Obj\sA$.}
\end{split}
\label{ai4eq19}
\e
We can check this satisfies \eq{ai4eq5}--\eq{ai4eq6}. For \eq{ai4eq5}, if $x_0,x_1\in\Obj\sA$ and $\fT\in \Hom_{A_\iy}(\fF^0,\fF^1)$ then \eq{ai4eq17} and \eq{ai4eq5} for $\sB$ yield
\begin{align*}
\n_2^{\fF^0,\fF^0,\fF^1}(\boo_{\fF^0}\ot\fT)^{x_0,\ldots,x_n}&=\m_2^{\fF^0,\fF^0,\fF^1}\bigl(\boo_{\fF^0(x_0)}\ot\fT_n^{x_0,\ldots,x_n}\bigr)=\fT_n^{x_0,\ldots,x_n},\\
\n_2^{\fF^0,\fF^1,\fF^1}(\fT\ot\boo_{\fF^1})^{x_0,\ldots,x_n}&=\m_2^{\fF^0,\fF^1,\fF^1}\bigl(\fT_n^{x_0,\ldots,x_n}\ot\boo_{\fF^1(x_n)}\bigr)=\fT_n^{x_0,\ldots,x_n}.
\end{align*}
\end{dfn}

\begin{rem}
\label{ai4rem4}
Definition \ref{ai4def7} shows that in the category $\ACat_R$, the mapping spaces $\Hom_{\ACat_R}(\sA,\sB)$ can be given the structure of $A_\iy$-categories $\Fun_{A_\iy}(\sA,\sB)$. That is, $\ACat_R$ is morally `a category enriched over itself'.

This is relevant to a central question for this book: {\it what kind of higher category should\/ $A_\iy$-categories form?\/} By Lyubashenko \cite[\S 7]{Lyub}, $A_\iy$-categories should be at least a 2-category. It is clear that $\ACat_R$ should really be an $(\iy,\iy)$-category, since if $\ACat_R$ has (possibly non-invertible) $n$-morphisms, then $\Fun_{A_\iy}(\sA,\sB)$ should have (possibly non-invertible) $n$-morphisms, and thus $\ACat_R$ must have (possibly non-invertible) $(n+1)$-morphisms.

We do not really understand $(\iy,\iy)$-categories, but we do have a good theory of $(\iy,1)$-categories, namely $\iy$-categories (quasicategories). To make $\ACat_R$ into an $\iy$-category, we allow only invertible $n$-morphisms for $n\ge 2$, which roughly corresponds to restricting (in Definition \ref{ai7def4} below) to the subcategory of $\Fun_{A_\iy}(\sA,\sB)$ with morphisms homotopy equivalences.

In Theorem \ref{ai7thm4} we will also define an $\iy$-{\it bicategory\/} of unital $A_\iy$-categories over $R$, which is a kind of $(\iy,2)$-category, as in~\S\ref{ai28}.
\end{rem}

Next we discuss functoriality of $\Fun_{A_\iy}(\sA,\sB)$ under $A_\iy$-functors $\fE:\sA_1\ra\sA_2$ or $\fG:\sB_1\ra\sB_2$, following Lef\`evre-Hasegawa \cite[\S 8.1.2]{LeHa}, Lyubashenko \cite[Prop.~6.2]{Lyub}, and Fukaya~\cite[Prop.-Def.~8.41]{Fuka3}.

\begin{dfn}
\label{ai4def9}
{\bf(a)} Let $\sA_1,\sA_2,\sB$ be $A_\iy$-categories over a commutative ring $R$, and $\fE:\sA_1\ra\sA_2$ be an $A_\iy$-functor. We define an $A_\iy$-functor $\fE^*:\Fun_{A_\iy}(\sA_2,\sB)\ra\Fun_{A_\iy}(\sA_1,\sB)$ to act on objects by $\fE^*_\Obj:\fF\mapsto\fF\ci\fE$ for $\fF:\sA_2\ra\sB$ an $A_\iy$-functor, and if $p\ge 1$ and $\fF_0,\ldots,\fF_p$ are objects in $\Fun_{A_\iy}(\sA_2,\sB)$, we define
\begin{align*}
(\fE^*)_p^{\fF_0,\ldots,\fF_p}&:\Hom_{A_\iy}(\fF_0,\fF_1)\ot\cdots\ot\Hom_{A_\iy}(\fF_{p-1},\fF_p)\\
&\qquad \longra \Hom_{A_\iy}(\fF_0\ci\fE,\fF_p\ci\fE)
\end{align*}
by $(\fE^*)_p^{\fF_0,\ldots,\fF_p}=0$ if $p\ne 1$, and
\ea
&\bigl((\fE^*)_1^{\fF_0,\fF_1}(\fT)\bigr){}_n^{x_0,\ldots,x_n}=
\label{ai4eq20}\\
&\sum_{0=l_0<l_1<\cdots<l_k=n}
\begin{aligned}[t]& (-1)^{\sum_{i=1}^k(l_i-l_{i-1}+1)(k-i)}\fT_k^{\fE(x_{l_0}),\fE(x_{l_1}),\fE(x_{l_2}),\ldots,\fE(x_{l_k})} \\ 
&\bigl(\fE_{l_1-l_0}^{x_{l_0},\ldots,x_{l_1}}\ot\fE_{l_2-l_1}^{x_{l_1},\ldots,x_{l_2}}\ot\cdots\ot\fE_{l_k-l_{k-1}}^{x_{l_{k-1}},\ldots,x_{l_k}}\bigr).
\end{aligned}
\nonumber
\ea
This defines an $A_\iy$-functor as in \cite{Fuka3,LeHa,Lyub}, contravariantly functorial in~$\fE$.
\smallskip

\noindent{\bf(b)} Let $\sA,\sB_1,\sB_2$ be $A_\iy$-categories over $R$, and $\fG:\sB_1\ra\sB_2$ be an $A_\iy$-functor. We define an $A_\iy$-functor $\fG_*:\Fun_{A_\iy}(\sA,\sB_1)\ra\Fun_{A_\iy}(\sA,\sB_2)$ to act on objects by $(\fG_*)_\Obj:\fF\mapsto\fG\ci\fF$ for $\fF:\sA\ra\sB_1$ an $A_\iy$-functor, and if $p\ge 1$ and $\fF_0,\ldots,\fF_p$ are objects in $\Fun_{A_\iy}(\sA,\sB_1)$, we define
\ea
&(\fG_*)_p^{\fF_0,\ldots,\fF_p}:\Hom_{A_\iy}(\fF_0,\fF_1)\ot\cdots\ot\Hom_{A_\iy}(\fF_{p-1},\fF_p)
\label{ai4eq21}\\
&\qquad \longra \Hom_{A_\iy}(\fG\ci\fF_0,\fG\ci\fF_p)\qquad\text{by}
\nonumber\\
&\bigl((\fG_*)_p^{\fF_0,\ldots,\fF_p}(\fT_1\ot\cdots\ot\fT_p)\bigr){}_n^{x_0,\ldots,x_n}=
\nonumber\\
&\sum_{\begin{subarray}{l} 
0=l_0\le l_1\le \\
\cdots\le l_k=n,\\
1\le a_1<a_2< \\
\cdots<a_p\le k: \\
l_j>l_{j-1},\\ j\ne 0,a_i 
\end{subarray}}
\begin{aligned}[t]&
(-1)^{\sum_{i=1}^p(k+1)(\deg\fT^i+1)}\cdot(-1)^{\sum_{i=1}^p\sum_{j=1}^{a_i-1}(l_j-l_{j-1}+1)(\deg\fT^i+1)}\cdot{}\\
&(-1)^{\sum_{i=1}^k(l_i-l_{i-1}+1)(k-i)} \cdot{}\\
&\fG_k^{\begin{subarray}{l} \fF^0(x_{l_0}),\fF^0(x_{l_1}),\ldots,\fF^0(x_{l_{a_1-1}}),\fF^1(x_{l_{a_1}}),\fF^1(x_{l_{a_1+1}}),\\ 
\ldots,\fF^1(x_{l_{a_2-1}}),\ldots,\fF^p(x_{l_{a_p}}),\fF^p(x_{l_{a_p+1}}),\ldots,\fF^p(x_{l_k})
\end{subarray}}\ci\\
&\bigl(\fF_{l_1-l_0}^{0;x_{l_0},\ldots,x_{l_1}}\ot\fF_{l_2-l_1}^{0;x_{l_1},\ldots,x_{l_2}}\ot\cdots\ot\fF_{l_{a_1-1}-l_{a_1-2}}^{0;x_{l_{a_1-2}},\ldots,x_{l_{a_1-1}}} \\
&\ot\fT_{l_{a_1}-l_{a_1-1}}^{1;x_{l_{a_1-1}},\ldots,x_{l_{a_1}}}\ot\fF_{l_{a_1+1}-l_{a_1}}^{1;x_{l_{a_1}},\ldots,x_{l_{a_1+1}}}\!\ot\!\cdots\!\ot\!\fF_{l_{a_2-1}-l_{a_2-2}}^{1;x_{l_{a_2-2}},\ldots,x_{l_{a_2-1}}} \\
&{}\ot
\fT_{l_{a_2}-l_{a_2-1}}^{2;x_{l_{a_2-1}},\ldots,x_{l_{a_2}}}\ot \cdots \ot\fT_{l_{a_p}-l_{a_p-1}}^{p;x_{l_{a_p-1}},\ldots,x_{l_{a_p}}} \\
&\ot\fF_{l_{a_p+1}-l_{a_p}}^{p;x_{l_{a_p}},\ldots,x_{l_{a_p+1}}}\ot\cdots\ot\fF_{l_k-l_{k-1}}^{p;x_{l_{k-1}},\ldots,x_{l_k}}\bigr).
\end{aligned}
\nonumber
\ea
This defines an $A_\iy$-functor as in \cite{Fuka3,LeHa,Lyub}, covariantly functorial in~$\fG$.
\smallskip

\noindent{\bf(c)} Combining Definition \ref{ai4def7} and {\bf(a)\rm,\bf(b)} gives a functor
\e
\Fun_{A_\iy}:(\ACat_R)^\rop\t\ACat_R\longra \ACat_R.
\label{ai4eq22}
\e
acting on objects $(\sA,\sB)\mapsto\Fun_{A_\iy}(\sA,\sB)$ and on morphisms by $(\fE^\rop,\fG)\mapsto \fE^*\ci\fG_*$. Here $\fE^*$ and $\fG_*$ commute, so their order does not matter.
\smallskip

\noindent{\bf(d)} The analogue of {\bf(a)\rm--\bf(c)} works for cohomologically unital $A_\iy$-categories, cohomologically unital $A_\iy$-fun\-ct\-ors, and $\Fun_{A_\iy}^\cu(-,-)$ in Definition \ref{ai4def7}, giving a functor
\e
\Fun_{A_\iy}^\cu:(\ACat^\cu_R)^\rop\t\ACat^\cu_R\longra \ACat^\cu_R.
\label{ai4eq23}
\e

\noindent{\bf(e)} The analogue of {\bf(a)\rm--\bf(c)} works for unital $A_\iy$-categories, unital $A_\iy$-fun\-ct\-ors, and $\Fun_{A_\iy}^\un(-,-)$ in Definition \ref{ai4def8}, giving a functor
\e
\Fun_{A_\iy}^\un:(\ACat^\un_R)^\rop\t\ACat^\un_R\longra \ACat^\un_R.
\label{ai4eq24}
\e
\end{dfn}

\begin{rem}
\label{ai4rem5}
Bespalov--Lyubashenko--Manzyuk \cite[\S 8.13, Prop.~9.5]{BLM} show that $\ACat_R$ and $\ACat^\un_R$ are {\it closed categories}, with internal Hom functors $\Fun_{A_\iy},\Fun_{A_\iy}^\un$ in \eq{ai4eq22}--\eq{ai4eq24}.
\end{rem}

We explain two simple examples of functor $A_\iy$-categories, which will be important in \S\ref{ai57}--\S\ref{ai58} and \S\ref{ai93} below.

\begin{dfn}
\label{ai4def10}
Let $R$ be a commutative ring. Define an $A_\iy$-category $\sP$ over $R$ to have $\Obj\sP=\{*\}$, that is, $\sP$ has one object $*$, and set $\Hom_{\sP}(*,*)=0$ with maps $\m_k^{*,\ldots,*}=0$ for $k\ge 1$. We call $\sP$ the ({\it non-unital\/}) {\it point\/ $A_\iy$-category}.
\end{dfn}

\begin{prop}
\label{ai4prop1}
Let\/ $\sA$ be any\/ $A_\iy$-category over $R$. Then there is a canonical isomorphism of\/ $A_\iy$-categories $\sA\cong\Fun_{A_\iy}(\sP,\sA),$ where $\Fun_{A_\iy}(\sP,\sA)$ is as in Definition\/~{\rm\ref{ai4def7}}.
\end{prop}

\begin{proof}
Objects of $\Fun_{A_\iy}(\sP,\sA)$ are $A_\iy$-functors $\fF:\sP\ra\sA$. Such $\fF$ has a map $\fF_\Obj:\Obj\sP=\{*\}\ra\Obj\sA$, and so determines an object $x=\fF(*)$ in $\sA$. Since $\Hom_\sP(*,*)=0$ we have $\fF_k^{*,\ldots,*}=0$ for all $k\ge 1$. Hence objects $\fF$ in $\Fun_{A_\iy}(\sP,\sA)$ are in 1-1 correspondence with objects $x$ in $\sA$, and we write $\fF^x$ for the functor with~$\fF^x(*)=x$.

For objects $x,y$ in $\sA$, consider the graded $R$-module $\Hom_{A_\iy}(\fF^x,\fF^y)\ab=\ab\Hom_{\Fun_{A_\iy}(\sP,\sA)}(\fF^x,\fF^y)$. Elements of this are $(\fT_k^{*,\ldots,*})_{k\ge 0}$. As $\Hom_{\sP}(*,*)=0$ we have $\fT_k^{*,\ldots,*}=0$ for $k\ge 1$, so by \eq{ai4eq13} the only data is 
\begin{equation*}
\fT_0^{*,*}\in\Hom_R\bigl(R,\Hom_\sA(x,y)\bigr)\cong \Hom_\sA(x,y).
\end{equation*}
This gives a canonical isomorphism of graded $R$-modules
\begin{equation*}
\Hom_{\Fun_{A_\iy}(\sP,\sA)}(\fF^x,\fF^y)\cong\Hom_\sA(x,y).
\end{equation*}
It is now easy to check that these isomorphisms also identify the maps $\m_k^{x_0,\ldots,x_k}$ in $\sA$ and $\Fun_{A_\iy}(\sP,\sA)$.
\end{proof}

Here are the unital analogues of Definition \ref{ai4def10} and Proposition~\ref{ai4prop1}.

\begin{dfn}
\label{ai4def11}
Let $R$ be a commutative ring. Define a unital $A_\iy$-category $\sP^\un$ over $R$ to have $\Obj\sP^\un=\{*\}$, that is, $\sP^\un$ has one object $*$, and set $\Hom_{\sP^\un}(*,*)=R\cdot\boo^{\sP^\un}_*$ with $\boo^{\sP^\un}_*$ in degree 0, and define the maps $\m_k^{*\ldots,*}$ by
\begin{equation*}
\m_2^{*,*}(\boo^{\sP^\un}_*\ot\boo^{\sP^\un}_*)=\boo^{\sP^\un}_*, \quad \m_k^{*\ldots,*}=0, \quad k\ne 2.
\end{equation*}
We call $\sP^\un$ the {\it unital point\/ $A_\iy$-category}. It is a dg-category in \S\ref{ai45} below.
\end{dfn}

\begin{prop}
\label{ai4prop2}
Let\/ $\sA$ be any unital\/ $A_\iy$-category over $R$. Then there is a canonical isomorphism of unital\/ $A_\iy$-categories $\sA\cong\Fun^\un_{A_\iy}(\sP^\un,\sA),$ where $\Fun^\un_{A_\iy}(\sP^\un,\sA)$ is as in Definitions\/~{\rm\ref{ai4def7}--\ref{ai4def8}}.
\end{prop}

\begin{proof}
Objects of $\Fun^\un_{A_\iy}(\sP^\un,\sA)$ are unital $A_\iy$-functors $\fF:\sP^\un\ra\sA$. Such $\fF$ has a map $\fF_\Obj:\Obj\sP^\un=\{*\}\ra\Obj\sA$, and so determines an object $x=\fF(*)$ in $\sA$. Since $\Hom_\sP^\un(*,*)=R\cdot\boo^{\sP^\un}_*$, by \eq{ai4eq11}--\eq{ai4eq12}, the maps $\fF_k^{*,\ldots,*}$ are determined uniquely by $\fF_1^{*,*}(\boo^{\sP^\un}_*)=\boo^\sA_x$ and $\fF_k^{*,\ldots,*}=0$ for $k\ne 1$. Hence objects $\fF$ in $\Fun^\un_{A_\iy}(\sP^\un,\sA)$ are in 1-1 correspondence with objects $x$ in $\sA$, and we write $\fF^x$ for the functor with~$\fF^x(*)=x$.

For objects $x,y$ in $\sA$, consider the graded $R$-module $\Hom_{A_\iy}^\un(\fF^x,\fF^y)$. Elements of this are $(\fT_k^{*,\ldots,*})_{k\ge 0}$. As $\Hom_{\sP^\un}(*,*)\ab=R\cdot\boo^{\sP^\un}_*$, we have $\fT_k^{*,\ldots,*}=0$ for $k\ge 1$ by \eq{ai4eq18}, so by \eq{ai4eq13} the only data is 
\begin{equation*}
\fT_0^{*,*}\in\Hom_R\bigl(R,\Hom_\sA(x,y)\bigr)\cong \Hom_\sA(x,y).
\end{equation*}
This gives a canonical isomorphism of graded $R$-modules
\begin{equation*}
\Hom_{\Fun_{A_\iy}^\un(\sP^\un,\sA)}(\fF^x,\fF^y)\cong\Hom_\sA(x,y).
\end{equation*}
When $x=y$, this isomorphism identifies $\boo^{\Fun^\un_{A_\iy}(\sP^\un,\sA)}_{\fF^x}$ in \eq{ai4eq19} with $\boo^\sA_x$ in $\Hom_\sA^0(x,x)$. It is now easy to check that these isomorphisms also identify the maps $\m_k^{x_0,\ldots,x_k}$ in $\sA$ and $\Fun^\un_{A_\iy}(\sP^\un,\sA)$.
\end{proof}

\begin{rem}
\label{ai4rem6}
Let $\sA$ be a non-unital $A_\iy$-category over $R$. For some applications it may be useful to consider the $A_\iy$-category $\Fun_{A_\iy}(\sP^\un,\sA)$. It should be regarded as the {\it idempotent completion\/} of $\sA$, as on the level of cohomology categories, an object $\fF$ in $\Fun_{A_\iy}(\sP^\un,\sA)$ gives an object $x=\fF(*)$ in $\sA$ and an element $i_\fF=[\fF_1^{*,*}(\boo_*^{\sP^\un})]$ in $H^0(\Hom_\sA(x,x))$ with $i_\fF^2=i_\fF$, that is, $i_\fF$ is an idempotent in the non-unital $R$-algebra $H^0(\Hom_\sA(x,x))$. This $i_\fF$ serves as a cohomological unit for $\fF$ in~$\Fun_{A_\iy}(\sP^\un,\sA)$.
\end{rem}

\begin{rem}
\label{ai4rem7}
Let $\sA,\sB,\sC$ be $A_\iy$-categories over $R$. It seems natural to ask whether there should exist some kind of `$A_\iy$-bifunctor'
\e
\la_{\sA\sB\sC}:\Fun_{A_\iy}(\sA,\sB)\t\Fun_{A_\iy}(\sB,\sC)\longra\Fun_{A_\iy}(\sA,\sC),
\label{ai4eq25}
\e
which acts on objects by $\la_{\sA\sB\sC}:(\fF,\fG)\mapsto\fG\ci\fF$, and similarly for the  (cohomologically) unital versions $\Fun_{A_\iy}^\un(\cdots),\Fun_{A_\iy}^\cu(\cdots)$. Lyubashenko \cite{Lyub} provides a construction of this kind, which we outline without giving details. 

Recall from \S\ref{ai36} that if $(A,\m)$ is an $A_\iy$-algebra over $R$, then the bar complex construction gives a dg-coalgebra $\bigl(T(A[1]),\co,\bar\d\bigr)$. This gives an equivalence between $A_\iy$-algebras and a certain class of dg-coalgebras, and one can define an $A_\iy$-structure $\m$ on a graded $R$-module $A$ as a suitable graded coderivation $\bar\d$ on $(T(A[1]),\co)$. Similarly, Lyubashenko \cite[\S 2]{Lyub} gives an equivalence between $A_\iy$-categories $\sA$ over $R$ and a certain class of {\it counital dg-cocategories}, and in fact he defines $A_\iy$-categories in terms of dg-cocategories. 

We will write $T\sA[1]$ for the counital dg-cocategory corresponding to $\sA$. $A_\iy$-functors $\fF:\sA\ra\sB$ are identified with suitable functors $T\sA[1]\ra T\sB[1]$, and similarly for $A_\iy$-natural transformations.

There is a good notion of {\it tensor product\/} $\ot$ of counital dg-cocategories. Note that this does not give a definition of tensor product of $A_\iy$-categories, since $(T\sA[1])\ot(T\sB[1])$ need not be of the form~$T\sC[1]$.

Let $\sA_1,\ldots,\sA_k,\sB$ be $A_\iy$-categories over $R$. Following Lyubashenko \cite[5.3]{Lyub}, we can define an $A_\iy$-{\it multifunctor\/} $\sA_1\t\cdots\t\sA_k\ra\sB$ to be a suitable kind of functor of counital dg-cocategories $(T\sA_1[1])\ot\cdots\ot(T\sA_k[1])\ra T\sB[1]$. When $k=1$ this recovers ordinary $A_\iy$-functors.

Given $A_\iy$-categories $\sA,\sB,\sC$ over $R$, Lyubashenko \cite[\S 4, p.~14 \& \S 6, p.~21]{Lyub}
defines a functor of counital dg-cocategories 
\e
\La_{\sA\sB\sC}:T\Fun_{A_\iy}(\sA,\sB)[1]\ot T\Fun_{A_\iy}(\sB,\sC)[1]\longra T\Fun_{A_\iy}(\sA,\sC)[1],
\label{ai4eq26}
\e
which is by definition an $A_\iy$-bifunctor \eq{ai4eq25}. These $\la_{\sA\sB\sC},\La_{\sA\sB\sC}$ have strict associativity and unit properties.

If $\sA,\sB,\sC$ are (cohomologically) unital, then $\la_{\sA\sB\sC},\La_{\sA\sB\sC}$ above restrict to bifunctors on $\Fun_{A_\iy}^\un(\cdots),\Fun_{A_\iy}^\cu(\cdots)$.
\end{rem}

\subsection{\texorpdfstring{$A_\iy$-homotopies, $A_\iy$-natural transformations, and $A_\iy$-natural isomorphisms}{A∞-homotopies, A∞-natural transformations, and A∞-natural isomorphisms}}
\label{ai44}

We now consider different notions of when two $A_\iy$-functors $\fF,\fG:\sA\ra\sB$ are `the same'. The next definition follows Seidel \cite[\S 1h]{Seid}, and is the analogue of homotopies of morphisms of $A_\iy$-algebras in Definition \ref{ai3def6}.

\begin{dfn}
\label{ai4def12}
Let $\sA,\sB$ be $A_\iy$-categories over $R$, and $\fF,\fG:\sA\ra\sB$ be $A_\iy$-functors which act the same way on objects, $\fF_\Obj=\fG_\Obj$. An $A_\iy$-{\it homotopy}, or just {\it homotopy}, $\H:\fF\Ra\fG$ consists of for all $k\ge 1$ and objects $x_0,x_1,\ldots,x_k$ in $\sA$, a graded $R$-linear map
\begin{align*}
\fH_k^{x_0,\ldots,x_k}:\Hom_\sA(x_0,x_1)\ot_R\Hom_\sA(x_1,x_2)\ot_R\cdots \ot_R\Hom_\sA(x_{k-1},x_k)&\\
\longra\Hom_\sB(\fF_\Obj(x_0),\fF_\Obj(x_k))&,
\end{align*}
of degree $-k$, satisfying the analogue of \eq{ai3eq12}: for all $n\ge 1$ and objects $x_0,x_1,\ldots,x_n$ in $\sA$, we require that
\ea
&\fF_n^{x_0,\ldots,x_n}-\fG_n^{x_0,\ldots,x_n}=
\label{ai4eq27}\\
&\sum_{\!\!\!\!\!\!\begin{subarray}{l} 
0=l_0<l_1\le \cdots< l_k=n,\\
1\le j\le k \end{subarray}}
\begin{aligned}[t]&(-1)^{k+1}\cdot(-1)^{\sum_{i=1}^{j-1}(l_i-l_{i-1}+1)}\cdot(-1)^{\sum_{i=1}^k(l_i-l_{i-1}+1)(k-i)}\\
&\m^{\sB;\fF(x_{l_0}),\fF(x_{l_1}),\ldots,\fF(x_{l_k})}_k\ci\bigl[\fF^{x_{l_0},\ldots,x_{l_1}}_{l_1-l_0}\!\ot\!\cdots\!\ot\!\fF^{x_{l_{j-2}},\ldots,x_{l_{j-1}}}_{l_{j-1}-l_{j-2}}\\
&{}\ot \fH^{x_{l_{j-1}},\ldots,x_{l_j}}_{l_j-l_{j-1}} \ot\fG^{x_{l_j},\ldots,x_{l_{j+1}}}_{l_{j+1}-l_j}\ot\cdots\ot\fG^{x_{l_{k-1}},\ldots,x_{l_k}}_{l_k-l_{k-1}}\bigr]
\end{aligned}
\nonumber\\
&+\sum_{\begin{subarray}{l} i,k\ge 0,\; j\ge 1:\\ i+j+k=n \end{subarray}}\begin{aligned}[t]&
(-1)^{i+jk}\fH^{x_0,x_1,\ldots,x_i,x_{i+j},x_{i+j+1},\ldots,x_n}_{i+k+1}\ci\\
&\bigl(\id_{\Hom_\sA(x_0,x_1)}\ot\cdots\ot \id_{\Hom_\sA(x_{i-1},x_i)}\ot\m^{\sA;x_i,x_{i+1},\ldots,x_{i+j}}_j \\
&\ot \id_{\Hom_\sA(x_{i+j},x_{i+j+1})}\ot\cdots\ot \id_{\Hom_\sA(x_{n-1},x_n)}\bigr).\end{aligned}
\nonumber
\ea

When $n=1,$ equation \eq{ai4eq27} becomes
\begin{equation*}
\fF_1^{x_0,x_1}-\fG_1^{x_0,x_1}=\m^{\sB;\fF(x_0),\fF(x_1)}_1\ci\H^{x_0,x_1}_1+\H^{x_0,x_1}_1\ci\m^{\sA;x_0,x_1}_1.
\end{equation*}
This implies that $\H^{x_0,x_1}_1$ is a chain homotopy between the morphisms of chain complexes $\fF_1^{x_0,x_1},\fG_1^{x_0,x_1}:\bigl(\Hom_\sA(x_0,x_1),\m^{\sA;x_0,x_1}_1)\!\ra\!(\Hom_\sB(\fF(x_0),\fF(x_1)),\ab\m^{\sB;\fF(x_0),\fF(x_1)}_1)$. Thus for the morphisms \eq{ai4eq9} on cohomology we have
\begin{equation*}
H^*(\fF^{x_0,x_1}_1)\!=\!H^*(\fG^{x_0,x_1}_1):H^*\bigl(\Hom_\sA(x_0,x_1)\bigr)\!\longra\! H^*\bigl(\Hom_\sB(\fF(x_0),\fF(x_1))\bigr).
\end{equation*}
Hence in functors of cohomology categories we have
\begin{equation*}
H^0(\fF_1)=H^0(\fG_1):H^0(\sA)\longra H^0(\sB).
\end{equation*}

If $\fF,\fG$ are unital $A_\iy$-functors of unital $A_\iy$-categories $\sA,\sB$, we call a homotopy $\H:\fF\Ra\fG$ {\it unital\/} if for all $k\ge 0$, $i=0,\ldots,k$, $x_0,\ldots,x_k\in\Obj\sA$ and $a_j\in\Hom_\sA(x_{j-1},x_j)$ we have
\e
\fH_{k+1}^{x_0,\ldots,x_i,x_i,\ldots,x_k}(a_1\ot\cdots\ot a_i\ot \boo_{x_i}\ot a_{i+1}\ot \cdots \ot a_k)=0.
\label{ai4eq28}
\e
\end{dfn}

\begin{dfn}
\label{ai4def13}
Let $\sA,\sB$ be $A_\iy$-categories over $R$, and $\fF,\fG:\sA\ra\sB$ be $A_\iy$-functors. Using the notation of Definition \ref{ai4def7}, an $A_\iy$-{\it natural transformation\/} $\fT:\fF\Ra\fG$ is an element $\fT\in\Hom_{A_\iy}^0(\fF,\fG)$ with $\n_1^{\fF,\fG}(\fT)=0$. This is the most natural notion of morphism $\fF\ra\fG$ in the $A_\iy$-category $\Fun_{A_\iy}(\sA,\sB)$.

Note that $\fT$ is data $(\fT_k)_{k\ge 0}$, where $\fT_0$ gives $\fT_0^x(1)\in\Hom_\sB^0(\fF(x),\fG(x))$ for all $x\in\Obj\sA$, and the $n=0$ equation in $\n_1^{\fF,\fG}(\fT)=0$ from \eq{ai4eq15} gives $\m_1^{\sB;\fF(x),\fG(x)}(\fT_0^x)=0$. In cohomology this gives a morphism $[\fT_0^x(1)]:H^0(\fF_1)(x)\ra H^0(\fG_1)(x)$ in the (nonunital) category $H^0(\sB)$. 

The $n=0$ equation in $\n_1^{\fF,\fG}(\fT)=0$ from \eq{ai4eq15} becomes
\begin{align*}
&\fT_1^{x_0,x_1}\ci\m_1^{\sA;x_0,x_1}+\m_1^{\sB;\fF(x_0),\fG(x_1)}\ci\fT_1^{x_0,x_1}\\
&-\m_2^{\sB;\fF(x_0),\fF(x_1),\fG(x_1)}\bigl(\fF_1^{x_0,x_1}\ot\fT^{x_1}_0(1)\bigr)\\
&+\m_2^{\sB;\fF(x_0),\fG(x_0),\fG(x_1)}\bigl(\fT^{x_0}_0(1)\ot\fG_1^{x_0,x_1}\bigr)=0.
\end{align*}
By applying this to a cochain representative of a morphism $\al:x_0\ra x_1$ in $H^0(\sA)$, we find the following commutes in the (nonunital) category~$H^0(\sB)$:
\begin{equation*}
\xymatrix@C=150pt@R=15pt{
*+[r]{H^0(\fF_1)(x_0)} \ar[r]_(0.5){[\fT^{x_0}_0(1)]} \ar[d]^{H^0(\fF_1)(\al)} & *+[l]{H^0(\fG_1)(x_0)} \ar[d]_{H^0(\fG_1)(\al)} \\
*+[r]{H^0(\fF_1)(x_1)} \ar[r]^(0.5){[\fT^{x_1}_0(1)]}  & *+[l]{H^0(\fG_1)(x_1).\!} }
\end{equation*}
Thus, we can define a natural transformation $H^0(\fT_0):H^0(\fF_1)\Ra H^0(\fG_1)$ of (nonunital) functors $H^0(\sA)\ra H^0(\sB)$ by $H^0(\fT_0)(x)=[\fT_0^x(1)]$.

Now suppose $\sA,\sB,\fF,\fG$ are unital, so that $\fF,\fG$ are objects in the $A_\iy$-subcategory $\Fun_{A_\iy}^\un(\sA,\sB)\subset\Fun_{A_\iy}(\sA,\sB)$. We call $\fT:\fF\Ra\fG$ a {\it unital\/} $A_\iy$-natural transformation if $\fT\in\Hom_{A_\iy}^{\un,0}(\fF,\fG)\subset\Hom_{A_\iy}^0(\fF,\fG)$. Then $H^0(\fT_0)$ is a natural transformation of honest categories and functors, not just nonunital categories and nonunital functors.

A unital $A_\iy$-natural transformation $\fT:\fF\Ra\fG$ is called a ({\it unital\/}) $A_\iy$-{\it natural isomorphism\/} if it is a homotopy equivalence in the unital $A_\iy$-category $\Fun_{A_\iy}^\un(\sA,\sB)$, in the sense of Definition \ref{ai4def5}, that is, $[\fT]$ is an isomorphism in $H^0(\Fun_{A_\iy}^\un(\sA,\sB))$. Unital $A_\iy$-natural isomorphism is an equivalence relation on unital functors~$\sA\ra\sB$.

Similarly, suppose $\sA,\sB,\fF,\fG$ are cohomologically unital, so that $\fF,\fG$ are objects in the  cohomologically unital $A_\iy$-category $\Fun_{A_\iy}^\cu(\sA,\sB)$. An $A_\iy$-natural transformation $\fT:\fF\Ra\fG$ is called a ({\it cohomologically unital\/}) $A_\iy$-{\it natural isomorphism\/} if it is a homotopy equivalence in~$\Fun_{A_\iy}^\cu(\sA,\sB)$.
\end{dfn}

\begin{rem}
\label{ai4rem8}
Fukaya \cite[Def.~8.1]{Fuka3} calls an $A_\iy$-natural isomorphism $\fT:\fF\Ra\fG$ a {\it homotopy equivalence\/} between $\fF$ and $\fG$. Lef\`evre-Hasegawa \cite[Def.~8.2.2.1]{LeHa} calls $A_\iy$-natural isomorphisms {\it weak equivalences}.
\end{rem}

Lef\`evre-Hasegawa \cite[Prop.~8.2.2.3]{LeHa} proves:

\begin{prop}
\label{ai4prop3}
Let\/ $\sA,\sB$ be unital\/ $A_\iy$-categories over $R,$ and\/ $\fF,\fG:\sA\ra\sB$ be unital\/ $A_\iy$-functors, and\/ $\fT:\fF\Ra\fG$ be a unital\/ $A_\iy$-natural transformation. Then $\fT$ is an $A_\iy$-natural  isomorphism if and only if\/ $H^0(\fT_0):H^0(\fF_1)\Ra H^0(\fG_1)$ is a natural isomorphism of ordinary functors.
\end{prop}

The next proposition, outlined in Seidel \cite[\S 2c]{Seid}, explains that $A_\iy$-homo\-top\-ies are special $A_\iy$-natural isomorphisms. We give a detailed proof.

\begin{prop}
\label{ai4prop4}
Let\/ $\sA,\sB$ be unital\/ $A_\iy$-categories over $R,$ and\/ $\fF,\fG:\sA\ra\sB$ be unital\/ $A_\iy$-functors which act the same way on objects, $\fF_\Obj=\fG_\Obj,$ and\/ $\fH:\fF\Ra\fG$ be a unital\/ $A_\iy$-homotopy. Define $\fT\in\Hom_{A_\iy}^{\un,0}(\fF,\fG)$ by
\e
\begin{split}
&\fT_0^{x_0}:R\longra \Hom^0_\sB(\fF(x_0),\fG(x_0))\quad\text{is given by}\quad \fT_0^{x_0}:\la\longmapsto\la\boo^\sB_{\fF(x_0)},\\
&\text{and}\quad \fT_k^{x_0,\ldots,x_k}=\H_k^{x_0,\ldots,x_k},\quad \text{all\/ $k\ge 1$ and\/ $x_0,\ldots,x_k\in\Obj\sA$.}
\end{split}
\label{ai4eq29}
\e
Then $\fT:\fF\Ra\fG$ is an $A_\iy$-natural isomorphism. This gives a 1-1 correspondence between unital\/ $A_\iy$-homotopies and\/ $A_\iy$-natural isomorphisms\/ $\fT:\fF\Ra\fG$ with\/ $\fT_0^{x_0}:\la\mapsto\la\boo^\sB_{\fF(x_0)}$ for all\/~$x_0\in\Obj\sA$.
\end{prop}

\begin{proof}
We claim that for $\fT$ as in \eq{ai4eq29}, the equation $\n_1^{\fF,\fG}(\fT)=0$ from \eq{ai4eq15} is equivalent to \eq{ai4eq27}. This holds as $\n_1^{\fF,\fG}(\fT)$ in \eq{ai4eq15} is the sum of terms from $\fT_k^{x_0,\ldots,x_k}=\H_k^{x_0,\ldots,x_k}$ for $k\ge 1$, which is just the right hand side of \eq{ai4eq27}, and terms in $\fT_0^{x_0}$. The terms in $\fT_0^{x_0}$ come only from the first sum in \eq{ai4eq15}, and as $\sB$ is unital, when $\fT_0^{x_0}=\boo^\sB_{\fF(x_0)}$ is substituted into $\m^{\sB;\fF(x_0),\ldots,\fF(x_n)}_k$ we get zero unless $k=2$. Thus there are only two possibly nonzero terms in $\fT_0^{x_0}$, which come from $k=2$, $j=1$, $l_0=l_1=0$, $l_2=n$ and from $k=2$, $j=2$, $l_0=0$, $l_1=l_2=n$ in the first sum in \eq{ai4eq15}, which are
\begin{align*}
&\m_2^{\sB;\fF(x_0),\fF(x_1)}\ci(\fT_0^{x_0}(1)\ot\fG_n^{x_0,\ldots,x_n})
-\m_2^{\sB;\fF(x_0),\fF(x_1)}\ci(\fF_n^{x_0,\ldots,x_n}\ot\fT_0^{x_n}(1))\\
&=\fG_n^{x_0,\ldots,x_n}-\fF_n^{x_0,\ldots,x_n},
\end{align*}
which is minus the right hand side of \eq{ai4eq27}. Hence \eq{ai4eq27} for $\fH$ is equivalent to $\n_1^{\fF,\fG}(\fT)=0$, and $\fT$ is an $A_\iy$-natural transformation.

We will prove directly that $\fT$ is an $A_\iy$-natural isomorphism. First we claim that there is a unique $\fU\in\Hom_{A_\iy}^{\un,0}(\fG,\fF)$ such that $\fU_0^{x_0}:\la\mapsto\la\boo^\sB_{\fF(x_0)}$ for all $x_0\in\Obj\sA$ and
\e
\n_2^{\fF,\fG,\fF}(\fT\ot\fU)=\boo_\fF\quad\text{in $\Hom_{A_\iy}^{\un,0}(\fF,\fF)$.}
\label{ai4eq30}
\e
This holds as by \eq{ai4eq15}, the $(\cdots)_n^{x_0,\ldots,x_n}$ component in \eq{ai4eq30} is of the form
\ea
&\m_2^{\fF(x_0),\fF(x_0),\fF(x_n)}(\fT_0^{x_0}(1)\ot\fU_n^{x_0,\ldots,x_n})+(\text{terms in $\fU_k^{x_0',\ldots,x_k'}$ for $k<n$})
\nonumber\\
&\qquad =\begin{cases} \boo^\sB_{\fF(x_0)}, & n=0, \\ 0, & n>0.
\end{cases}
\label{ai4eq31}
\ea
But $\fT_0^{x_0}(1)=\boo^\sB_{\fF(x_0)}$, so $\m_2^{\fF(x_0),\fF(x_0),\fF(x_n)}(\fT_0^{x_0}(1)\ot\fU_n^{x_0,\ldots,x_n})=\fU_n^{x_0,\ldots,x_n}$. Thus we can use \eq{ai4eq31} to determine $\fU_n^{x_0,\ldots,x_n}$ uniquely by induction on $n=0,1,2,\ldots.$

The $A_\iy$-equation \eq{ai4eq1} for $n=2$ and objects $\fF,\fG,\fF$ applied to $\fT\ot\fU$ yields
\begin{align*}
0&=\n_1^{\fF,\fF}\ci\n_2^{\fF,\fG,\fF}(\fT\ot\fU)-\n_2^{\fF,\fG,\fF}\bigl(\n_1^{\fF,\fG}(\fT)\ot\fU\bigr)
-\n_2^{\fF,\fG,\fF}\bigl(\fT\ot\n_1^{\fG,\fF}(\fU)\bigr)\\
&=\n_1^{\fF,\fF}(\boo_\fF)-0-\n_2^{\fF,\fG,\fF}\bigl(\fT\ot\n_1^{\fG,\fF}(\fU)\bigr)=-\n_2^{\fF,\fG,\fF}\bigl(\fT\ot\n_1^{\fG,\fF}(\fU)\bigr).
\end{align*}
Hence the uniqueness argument for $\fU$ above shows that $\n_1^{\fG,\fF}(\fU)=0$.

Similarly, there is a unique $\fS\in\Hom_{A_\iy}^{\un,0}(\fG,\fF)$ such that $\fS_0^{x_0}:\la\mapsto\la\boo^\sB_{\fF(x_0)}$ for all $x_0\in\Obj\sA$ and $\n_1^{\fG,\fF}(\fS)=0$ and
\begin{equation*}
\n_2^{\fG,\fF,\fG}(\fS\ot\fT)=\boo_\fG\quad\text{in $\Hom_{A_\iy}^{\un,0}(\fG,\fG)$.}
\end{equation*}
Now equation \eq{ai4eq1} for $n=3$ and objects $\fG,\fF,\fG,\fF$ applied to $\fS\ot\fT\ot\fU$ yields  
\ea
0&=\n_1^{\fG,\fF}\ci\n_3^{\fG,\fF,\fG,\fF}(\fS\ot\fT\ot\fU)
+\n_2^{\fG,\fG,\fF}\bigl(\n_2^{\fG,\fF,\fG}(\fS\ot\fT)\ot\fU\bigr)
\nonumber\\
&-\n_2^{\fG,\fF,\fF}\bigl(\fS\ot\n_2^{\fF,\fG,\fF}(\fT\ot\fU)\bigr)
+\n_3^{\fG,\fF,\fG,\fF}\bigl(\n_1^{\fG,\fF}(\fS)\ot\fT\ot\fU\bigr)
\nonumber\\
&+\n_3^{\fG,\fF,\fG,\fF}\bigl(\fS\ot\n_1^{\fG,\fF}(\fT)\ot\fU\bigr)
+\n_3^{\fG,\fF,\fG,\fF}\bigl(\fS\ot\fT\ot\n_1^{\fG,\fF}(\fU)\bigr)
\nonumber\\
&=\n_1^{\fG,\fF}\ci\n_3^{\fG,\fF,\fG,\fF}(\fS\ot\fT\ot\fU)+\fU-\fS+0+0+0.
\label{ai4eq32}
\ea
In the cohomology category $H^0(\Fun_{A_\iy}^\un(\sA,\sB))$ we now have morphisms $[\fT]:\fF\ra\fG$ and $[\fS],[\fU]:\fG\ra\fF$ with $[\fT]\ci[\fU]=1_\fF$ and $[\fS]\ci[\fT]=1_\fG$. But \eq{ai4eq32} shows that $[\fS]=[\fU]$. Thus $[\fS]$ is the inverse morphism to $[\fT]$, so $[\fT]$ is an isomorphism, and $\fT$ is an $A_\iy$-natural isomorphism.
\end{proof}

\begin{rem}
\label{ai4rem9}
Both $A_\iy$-homotopies, and $A_\iy$-natural transformations, are defined for non-unital $A_\iy$-categories. However, the relation in Proposition \ref{ai4prop4} between $A_\iy$-homotopies and $A_\iy$-natural transformations only makes sense for unital $A_\iy$-categories.
\end{rem}

The next definition follows \cite[Def.~8.5]{Fuka3} and~\cite[Def.~9.2.0.1]{LeHa}.

\begin{dfn}
\label{ai4def14}
Let $\sA,\sB$ be unital $A_\iy$-categories over $R,$ and $\fF:\sA\ra\sB$ be a unital $A_\iy$-functor. We call $\fF$ a ({\it unital\/}) {\it homotopy equivalence\/} if there exists a unital $A_\iy$-functor $\fG:\sB\ra\sA$ and unital $A_\iy$-natural isomorphisms $\fS:\fG\ci\fF\Ra\id_\sA$ and $\fT:\fF\ci\fG\Ra\id_\sB$. Then we also call $\sA,\sB$ ({\it unital\/}) {\it homotopy equivalent}. 

Similarly, if $\sA,\sB$ are cohomologically unital $A_\iy$-categories over $R,$ and $\fF:\sA\ra\sB$ is a cohomologically unital $A_\iy$-functor, we call $\fF$ a ({\it cohomologically unital\/}) {\it homotopy equivalence\/} if there exists a cohomologically unital $A_\iy$-functor $\fG:\sB\ra\sA$ and $A_\iy$-natural isomorphisms $\fS:\fG\ci\fF\Ra\id_\sA$ and $\fT:\fF\ci\fG\Ra\id_\sB$. Then we also call $\sA,\sB$ ({\it cohomologically unital\/}) {\it homotopy equivalent}. 

Compositions of homotopy equivalences are homotopy equivalences. Homotopy equivalence is an equivalence relation on unital, or cohomologically unital, $A_\iy$-categories.
\end{dfn}

The next definition follows \cite[Def.~1.13]{COS2} and~\cite[\S 1.4]{Orna1}.

\begin{dfn}
\label{ai4def15}
Let $\sA,\sB$ be unital, or cohomologically unital, $A_\iy$-categories over a commutative ring $R$, and $\fF:\sA\ra\sB$ be a unital, or cohomologically unital, $A_\iy$-functor. We call $\fF$ a {\it quasi-equivalence\/} if:
\begin{itemize}
\setlength{\itemsep}{0pt}
\setlength{\parsep}{0pt}
\item[(a)] For all objects $x,y$ in $\sA$ and $k\in\Z$, the following is an isomorphism
\begin{equation*}
H^k(\fF_1^{x,y}):H^k\bigl(\Hom_\sA(x,y)\bigr)\longra H^k\bigl(\Hom_\sB(\fF(x),\fF(y))\bigr).
\end{equation*}
\item[(b)] $H^0(\fF_1):H^0(\sA)\ra H^0(\sB)$ is an equivalence of categories.
\end{itemize}
All homotopy equivalences are quasi-equivalences. Theorem \ref{ai4thm3} below gives sufficient conditions for a quasi-equivalence to be a homotopy equivalence. 

If $R$ is a field, and in the unital case $\sA,\sB$ do not have zero objects (objects $z$ with $\boo_z=0$), then homotopy equivalences and quasi-equivalences coincide.
\end{dfn}

\begin{rem}
\label{ai4rem10}
Homotopy equivalence and quasi-equivalence are two competing notions of when two (cohomologically) unital $A_\iy$-categories are `the same'. One could argue that for applications  --- e.g.\ in comparing two constructions of Fukaya category for a symplectic manifold --- homotopy equivalence is better. 

However, there is a consensus in the mathematical community (see for example \S\ref{ai410} and \cite{COS1,COS2,Faon2,OhTa,Orna1,Tana,Pasc}) that to form localized ($\iy$-)categories of $A_\iy$-categories, one should invert quasi-equivalences, not homotopy equivalences. One reason for this is that in the Tabuada model structure on dg-categories (see Theorem \ref{ai4thm1}) the weak equivalences are quasi-equivalences.

One lesson from our results in \S\ref{ai92} is that localizing homotopy equivalences of $A_\iy$-categories is also an interesting and fruitful thing to do.
\end{rem}

Definitions \ref{ai4def7}--\ref{ai4def8} defined $A_\iy$-categories $\Fun_{A_\iy}(\sA,\sB),\ab\Fun_{A_\iy}^\cu(\sA,\sB),\ab\Fun_{A_\iy}^\un(\sA,\sB)$, and Definition \ref{ai4def9} explained their functoriality under $A_\iy$-functors in $\sA,\sB$. We now consider their functoriality under $A_\iy$-natural transformations and $A_\iy$-natural isomorphisms. The only reference the author can find to the content of Definition \ref{ai4def16} in the literature is Lyubashenko \cite[Prop.~8.4]{Lyub}, who appears to claim that equivalence classes $[\fL_*]$ of $\fL_*:(\fG_1)_*\Ra(\fG_2)_*$ in part (b) below are part of a strict 2-functor $\Fun_{A_\iy}(\sA,-)$ defined on 2-categorical enhancements of~$\ACat_R,\ACat^\cu_R$.

\begin{dfn}
\label{ai4def16}
{\bf(a)} Let $R$ be a commutative ring, and $\sA_1,\sA_2,\sB$ be $A_\iy$-categories over $R$, and $\fE_1,\fE_2:\sA_1\ra\sA_2$ be $A_\iy$-functors, and $\fK:\fE_1\Ra\fE_2$ be an $A_\iy$-natural transformation. Definition \ref{ai4def9} constructed $A_\iy$-functors $\fE_1^*,\fE_2^*:\Fun_{A_\iy}(\sA_2,\sB)\ab\ra\Fun_{A_\iy}(\sA_1,\sB)$. We will define an $A_\iy$-natural transformation $\fK^*:\fE_1^*\Ra\fE_2^*$. To do this, for all $p\ge 1$ and $A_\iy$-functors $\fF_0,\fF_1,\ldots,\fF_p:\sA_2\ra\sB$, define a graded $R$-linear map of degree $-p$
\begin{align*}
(\fK^*)_p^{\fF_0,\ldots,\fF_p}&:\Hom_{A_\iy}(\fF_0,\fF_1)\ot\cdots\ot \Hom_{A_\iy}(\fF_{p-1},\fF_p)\\
&\quad \longra\Hom_{A_\iy}(\fF_0\ci\fE_1,\fF_p\ci\fE_2)
\end{align*}
by $(\fK^*)_p^{\fF_0,\ldots,\fF_p}=0$ if $p\ne 1$, and in a similar way to \eq{ai4eq20}
\ea
&\bigl((\fK^*)_1^{\fF_0,\fF_1}(\fT)\bigr){}_n^{x_0,\ldots,x_n}=
\label{ai4eq33}\\
&\sum_{\begin{subarray}{l}
0=l_0\le l_1\le\cdots\le l_k=n,\\
1\le j\le k: \\
\text{$l_{i-1}<l_i$ if $i\ne j$} \end{subarray}}
\begin{aligned}[t]& (-1)^{k+1}\cdot(-1)^{\sum_{i=1}^{j-1}(l_i-l_{i-1}+1)}\cdot(-1)^{\sum_{i=1}^k(l_i-l_{i-1}+1)(k-i)}\\
&\fT_k^{\fE_1(x_{l_0}),\fE_1(x_{l_1}),\ldots,\fE_1(x_{l_{j-1}}),\fE_2(x_{l_j}),\fE_2(x_{l_{j+1}}),\ldots,\fE_2(x_{l_k})} \\ 
&\bigl((\fE_1)_{l_1-l_0}^{x_{l_0},\ldots,x_{l_1}}\ot(\fE_1)_{l_2-l_1}^{x_{l_1},\ldots,x_{l_2}}\ot\cdots\ot(\fE_1)_{l_{j-1}-l_{j-2}}^{x_{l_{j-2}},\ldots,x_{l_{j-1}}}\\
&\ot\fK_{l_j-l_{j-1}}^{x_{l_{j-1}},\ldots,x_{l_j}}\ot
(\fE_2)_{l_{j+1}-l_j}^{x_{l_j},\ldots,x_{l_{j+1}}}\ot\cdots\ot(\fE_2)_{l_k-l_{k-1}}^{x_{l_{k-1}},\ldots,x_{l_k}}\bigr).
\end{aligned}
\nonumber
\ea

We can prove using the identity $\n_1^{\fE_1,\fE_2}(\fK)=0$ that $\n_1^{\fE_1^*,\fE_2^*}(\fK^*)=0$, so $\fK^*:\fE_1^*\Ra\fE_2^*$ is an $A_\iy$-natural transformation. 

If $\sA_1,\sA_2,\sB,\fE_1,\fE_2$ are cohomologically unital then $\fE_1^*,\fE_2^*$ restrict to $\fE_1^*,\fE_2^*:\Fun_{A_\iy}^\cu(\sA_2,\sB)\ab\ra\Fun_{A_\iy}^\cu(\sA_1,\sB)$, where $\Fun_{A_\iy}^\cu(\sA_i,\sB)\subset\Fun_{A_\iy}(\sA_i,\sB)$ is the full $A_\iy$-subcategory with objects cohomologically unital $A_\iy$-functors, and $\fK^*$ above restricts to an $A_\iy$-natural transformation of $A_\iy$-functors $\Fun_{A_\iy}^\cu(\sA_2,\sB)\ab\ra\Fun_{A_\iy}^\cu(\sA_1,\sB)$ by taking $\fF_0,\ldots,\fF_p$ above to be cohomologically unital.

If $\sA_1,\sA_2,\sB,\fE_1,\fE_2,\fK$ are unital then $\fE_1^*,\fE_2^*$ restrict to $\fE_1^*,\fE_2^*:\Fun_{A_\iy}^\un(\sA_2,\ab\sB)\ab\ra\Fun_{A_\iy}^\un(\sA_1,\ab\sB)$. We replace $\Hom_{A_\iy}(\cdots)$ by $\Hom^\un_{A_\iy}(\cdots)$, and take $\fF_0,\ab\ldots,\ab\fF_p$ to be unital, and then $\fK^*$ above defines an $A_\iy$-natural transformation of unital $A_\iy$-functors $\Fun_{A_\iy}^\un(\sA_2,\sB)\ab\ra\Fun_{A_\iy}^\un(\sA_1,\sB)$. 

To see \eq{ai4eq33} defines an element of $\Hom_{A_\iy}^\un(\fF_0\ci\fE_1,\fF_p\ci\fE_2)\subset\Hom_{A_\iy}(\fF_0\ci\fE_1,\fF_1\ci\fE_2)$ in this case, note that if $1\le j\le n$ and $x_{j-1}=x_j$, and we evaluate \eq{ai4eq33} on an element of the form $a_1\ot\cdots\ot a_{j-1}\ot\boo_{x_j}\ot a_{j+1}\ot\cdots\ot a_n$ for $a_i\in\Hom_{\sA}(x_{i-1},x_i)$ then in each term in the sum in \eq{ai4eq33}, either $\boo_{x_j}$ is an input in $\fK_{n_l}^{x_{m_{l-1}},\ldots,x_{m_l}}$, so it gives zero as $\fK$ is unital, or $\boo_{x_j}$ is an input in
$(\fE_a)_{n_l}^{x_{m_{l-1}},\ldots,x_{m_l}}$ for $a=1,2$ and $n_l\ge 2$, so it gives zero as $\fE_a$ is unital, or $\boo_{x_j}$ is an input in $(\fE_a)_1^{x_j,x_j}$ for $a=1,2$, and then $(\fE_a)_1^{x_j,x_j}(\boo_{x_j})=\boo_{\fE_a(x_j)}$ is an input in $\fT_k^{\fE_1(x_{m_0}),\ldots,\fE_2(x_{m_k})}$, so it gives zero as $\fT$ lies in $\Hom_{A_\iy}^\un(\fF_0,\fF_1)$.

In the cohomologically unital, or unital, cases, we can also show that if $\fK$ is a (cohomologically) unital $A_\iy$-natural isomorphism then $\fK^*$ is also a (cohomologically) unital $A_\iy$-natural isomorphism. This is because then there is another $A_\iy$-natural transformation $\fJ:\fE_2\Ra\fE_1$ such that $[\fJ]$ is inverse to $[\fK]$ in $H^0(\Fun^*(\sA_1,\sA_2))$, and we show that $\fJ^*:\fE_2^*\Ra\fE_1^*$ has $[\fJ^*]$ inverse to $[\fK^*]$ in $H^0\bigl(\Fun^*\bigl(\Fun_{A_\iy}^*(\sA_2,\sB),\Fun_{A_\iy}^*(\sA_1,\sB)\bigr)\bigr)$.
\smallskip

\noindent{\bf(b)} Let $\sA,\sB_1,\sB_2$ be $A_\iy$-categories over $R$, and $\fG_1,\fG_2:\sB_1\ra\sB_2$ be $A_\iy$-functors, and $\fL:\fG_1\Ra\fG_2$ be an $A_\iy$-natural transformation. Definition \ref{ai4def9} constructed $A_\iy$-functors $(\fG_1)_*,(\fG_2)_*:\Fun_{A_\iy}(\sA,\sB_1)\ab\ra\Fun_{A_\iy}(\sA,\sB_2)$. We will define an $A_\iy$-natural transformation $\fL_*:(\fG_1)_*\Ra(\fG_2)_*$. To do this, for all $p\ge 1$ and $A_\iy$-functors $\fF_0,\fF_1,\ldots,\fF_p:\sA\ra\sB_1$, in a very similar way to \eq{ai4eq21} define a graded $R$-linear map of degree $-p$
\begin{align*}
&(\fL_*)_p^{\fF_0,\ldots,\fF_p}:\Hom_{A_\iy}(\fF_0,\fF_1)\ot\cdots\ot\Hom_{A_\iy}(\fF_{p-1},\fF_p)\\
&\qquad \longra \Hom_{A_\iy}(\fG_1\ci\fF_0,\fG_2\ci\fF_p)\qquad\text{by}
\nonumber\\
&\bigl((\fL_*)_p^{\fF_0,\ldots,\fF_p}(\fT_1\ot\cdots\ot\fT_p)\bigr){}_n^{x_0,\ldots,x_n}=\\
&\sum_{\begin{subarray}{l} 
0=l_0\le l_1\le \\
\cdots\le l_k=n,\\
1\le a_1<a_2< \\
\cdots<a_p\le k: \\
l_j>l_{j-1},\\ j\ne 0,a_i 
\end{subarray}}
\begin{aligned}[t]&
(-1)^{\sum_{i=1}^p(k+1)(\deg\fT^i+1)}\cdot(-1)^{\sum_{i=1}^p\sum_{j=1}^{a_i-1}(l_j-l_{j-1}+1)(\deg\fT^i+1)}\cdot{}\\
&(-1)^{\sum_{i=1}^k(l_i-l_{i-1}+1)(k-i)} \cdot{} \\
&\fL_k^{\begin{subarray}{l} \fF^0(x_{l_0}),\fF^0(x_{l_1}),\ldots,\fF^0(x_{l_{a_1-1}}),\fF^1(x_{l_{a_1}}),\fF^1(x_{l_{a_1+1}}),\\
\ldots,\fF^1(x_{l_{a_2-1}}),\ldots,\fF^p(x_{l_{a_p}}),\fF^p(x_{l_{a_p+1}}),\ldots,\fF^p(x_{l_k})
\end{subarray}}\ci\\
&\bigl(\fF_{l_1-l_0}^{0;x_{l_0},\ldots,x_{l_1}}\ot\fF_{l_2-l_1}^{0;x_{l_1},\ldots,x_{l_2}}\ot\cdots\ot\fF_{l_{a_1-1}-l_{a_1-2}}^{0;x_{l_{a_1-2}},\ldots,x_{l_{a_1-1}}} \\
&\ot\fT_{l_{a_1}-l_{a_1-1}}^{1;x_{l_{a_1-1}},\ldots,x_{l_{a_1}}}\ot\fF_{l_{a_1+1}-l_{a_1}}^{1;x_{l_{a_1}},\ldots,x_{l_{a_1+1}}}\!\ot\!\cdots\!\ot\!\fF_{l_{a_2-1}-l_{a_2-2}}^{1;x_{l_{a_2-2}},\ldots,x_{l_{a_2-1}}} \\
&{}\ot
\fT_{l_{a_2}-l_{a_2-1}}^{2;x_{l_{a_2-1}},\ldots,x_{l_{a_2}}}\ot \cdots \ot\fT_{l_{a_p}-l_{a_p-1}}^{p;x_{l_{a_p-1}},\ldots,x_{l_{a_p}}} \\
&\ot\fF_{l_{a_p+1}-l_{a_p}}^{p;x_{l_{a_p}},\ldots,x_{l_{a_p+1}}}\ot\cdots\ot\fF_{l_k-l_{k-1}}^{p;x_{l_{k-1}},\ldots,x_{l_k}}\bigr).
\end{aligned}
\end{align*}

Again, we can prove using the identity $\n_1^{\fG_1,\fG_2}(\fL)=0$ that $\n_1^{\fG_1^*,\fG_2^*}(\fL^*)=0$, so $\fL^*:\fG_1^*\Ra\fG_2^*$ is an $A_\iy$-natural transformation. 

By the same reasoning as in {\bf(a)}, we can restrict $\fL^*$ to the cohomologically unital, and unital, cases, and if $\fL$ is a (cohomologically) unital $A_\iy$-natural isomorphism then $\fL^*$ is also a (cohomologically) unital $A_\iy$-natural isomorphism.
\end{dfn}

\begin{cor}
\label{ai4cor1}
{\bf(a)} Let\/ $\sA,\sA_1,\sA_2,\sB,\sB_1,\sB_2$ be unital\/ $A_\iy$-categories over\/ $R,$ and\/ $\fE:\sA_1\ra\sA_2,$ $\fG:\sB_1\ra\sB_2$ be unital homotopy equivalences. Then $\fE^*:\Fun^\un_{A_\iy}(\sA_2,\sB)\ab\ra\Fun^\un_{A_\iy}(\sA_1,\sB)$ and\/ $\fG_*:\Fun^\un_{A_\iy}(\sA,\sB_1)\ra\Fun^\un_{A_\iy}(\sA,\sB_2)$ from Definition\/ {\rm\ref{ai4def9}} are unital homotopy equivalences.
\smallskip

\noindent{\bf(b)} The analogue of\/ {\bf(a)} holds for the cohomologically unital case.
\end{cor}

\begin{proof}
For the first part of (a), as $\fE:\sA_1\ra\sA_2$ is a unital homotopy equivalence there exists a unital $A_\iy$-functor $\fF:\sA_2\ra\sA_1$ and unital $A_\iy$-natural isomorphisms $\fJ:\fF\ci\fE\Ra\id_{\sA_1}$ and $\fK:\fE\ci\fF\Ra\id_{\sA_2}$. Then we have a unital $A_\iy$-functor $\fF^*:\Fun_{A_\iy}(\sA_1,\sB)\ab\ra\Fun_{A_\iy}(\sA_2,\sB)$ and unital $A_\iy$-natural isomorphisms $\fJ^*:\fE^*\ci\fF^*=(\fF\ci\fE)^*\Ra\id_{\Fun_{A_\iy}(\sA_1,\sB)}$ and $\fK:\fF^*\ci\fE^*=(\fE\ci\fF)^*\Ra\id_{\Fun_{A_\iy}(\sA_2,\sB)}$. Hence $\fE^*$ is a unital homotopy equivalence. The rest of the proof is essentially the same.
\end{proof}

\subsection{\texorpdfstring{$A_\iy$-categories and dg-categories}{A∞-categories and dg-categories}}
\label{ai45}

We generalize the relation between $A_\iy$-algebras and dg-algebras in \S\ref{ai35}.

\begin{dfn}
\label{ai4def17}
Let $R$ be a commutative ring. A {\it differential graded category}, or {\it dg-category}, $\sD$ over $R$ consists of:
\begin{itemize}
\setlength{\itemsep}{0pt}
\setlength{\parsep}{0pt}
\item[(a)] A set $\Obj\sD$ of {\it objects\/} of $\sD$. (More generally $\Obj\sD$ could be a proper class, but as in Remark \ref{ai1rem2} we restrict to $\sU$-small sets in this book.)

\item[(b)] For all objects $x,y\in\Obj\sD$, a $\Z$-graded $R$-module of {\it morphisms}
\begin{equation*}
\Hom_\sD(x,y)=\Hom(x,y)=\bigop\nolimits_{d\in\Z}\Hom^d(x,y),
\end{equation*}
and an $R$-linear differential $\d:\Hom^d(x,y)\ra\Hom^{d+1}(x,y)$ with $\d^2=0$.
\item[(c)] For all $x\in\Obj\sD$ a unit $\boo_x$ (or $\boo_x^\sD$) in $\Hom^0(x,x)$ with $\d\boo_x=0$.
\item[(d)] For all $x,y,z\in\Obj\sD$, an $R$-linear {\it composition map\/} graded of degree~zero
\begin{equation*}
\mu^{x,y,z}:\Hom(x,y)\ot_R\Hom(y,z)\longra\Hom(x,z).
\end{equation*}
This should satisfy for all $w,x,y,z\in\Obj\sD$ and $a\in\Hom^k(w,\ab x)$, $b\in\Hom^l(x,y)$, and $c\in\Hom^m(y,z)$
\end{itemize}
\e
\begin{gathered}
\mu^{w,w,x}(\boo_w\ot a)=\mu^{w,x,x}(a\ot\boo_x)=a,\\
\d\ci\mu^{w,x,y}(a\ot b)=\mu^{w,x,y}((\d a)\ot b)+(-1)^k\mu^{w,x,y}(a\ot(\d b)),\\
\mu^{w,y,z}\bigl(\mu_2^{w,x,y}(a\ot b)\ot c\bigr)=\mu^{w,x,z}\bigl(a\ot \mu^{x,y,z}(b\ot c)\bigr).
\end{gathered}
\label{ai4eq34}
\e

Note that (unlike $A_\iy$-categories), $\sD$ is an honest category, with morphisms $\Hom(x,y)$ (or alternatively $\Hom^0(x,y)$), units $\boo_x$, and composition $\mu^{x,y,z}$, which is associative by~\eq{ai4eq34}.

Let $\sD,\sE$ be dg-categories over $R$. A {\it dg-functor\/} $\fF:\sD\ra\sE$ consists of a map $\fF_\Obj:\Obj\sD\ra\Obj\sE$, often written $x\mapsto\fF(x)$ rather than $x\mapsto\fF_\Obj(x)$, and for all $x,y\in\Obj\sD$ an $R$-linear map graded of degree zero
\begin{equation*}
\fF^{x,y}:\Hom_\sD(x,y)\longra\Hom_\sD(\fF(x),\fF(y)),
\end{equation*}
satisfying for all $x,y,z\in\Obj\sD$
\begin{gather*}
\fF^{x,y}(\boo^\sD_x)=\boo^\sE_{\fF(x)},\quad \d_\sE\ci\fF^{x,y}=\fF^{x,y}\ci\d_\sD,\\
\fF^{x,z}\ci\mu_\sD^{x,y,z}=\mu_\sE^{\fF(x),\fF(y),\fF(z)}\ci(\fF^{x,y}\ot\fF^{y,z}).
\end{gather*}
Again, dg-functors are honest functors of categories. Identity functors and compositions of functors, are defined in the obvious way. Then dg-categories over $R$ and dg-functors form a category~$\DGCat_R$.

Any dg-category can be made into a unital $A_\iy$-category by defining $\m_1^{x,y}=\d:\Hom(x,y)\ra\Hom(x,y)$, and $\m_2^{x,y,z}=\mu^{x,y,z}$, and $\m_k^{x_0,\ldots,x_k}=0$ for $k\ge 3$. This identifies dg-categories with unital $A_\iy$-categories with $\m_k=0$ for~$k\ge 3$.

Similarly, any dg-functor $\fF:\sD\ra\sE$ can be made into a unital $A_\iy$-functor of $A_\iy$-categories by setting $\fF_1^{x,y}=\fF^{x,y}$ and $\fF_k^{x_0,\ldots,x_k}=0$ for $k\ge 2$. This identifies dg-functors with unital $A_\iy$-functors $\fF$ (between $A_\iy$-categories coming from dg-categories) with $\fF_k=0$ for $k\ge 2$. Thus we have an inclusion $\DGCat_R\subset\ACat^\un_R$ as a non-full subcategory.

We also write $\DGCat^{A_\iy}_R$ for the full subcategory of $\ACat^\un_R$ whose objects are dg-categories, but whose morphisms are $A_\iy$-functors rather than dg-functors. Then~$\DGCat_R\subset\DGCat^{A_\iy}_R\subset\ACat^\un_R$.

We call a dg-functor $\fF:\sD\ra\sE$ a {\it quasi-equivalence\/} if it is a quasi-equivalence of $A_\iy$-categories, as in Definition~\ref{ai4def15}.
\end{dfn}

\begin{rem}
\label{ai4rem11}
Dg-categories have many applications (for example, to derived categories), and have received more attention than $A_\iy$-categories. Note that one can consider $A_\iy$-functors between dg-categories, and for some purposes these are better than dg-functors. For example, a claim of Kontsevich proved in \cite[Th.~C]{COS2} says that if $\sA,\sB$ are h-projective dg-categories over $R$ with split units, then the internal Hom object $\cHom(\sA,\sB)$ in $\DGCat_R$ is~$\Fun^\un_{A_\iy}(\sA,\sB)$.

See \cite[\S 1.8]{BLM} for discussion of why $A_\iy$-categories and $A_\iy$-functors are important in the study of dg-categories.
\end{rem}

\begin{rem}
\label{ai4rem12}
If $\sA$ is a unital $A_\iy$-category over $R$, and $\sB$ is a dg-category over $R$, then $\Fun_{A_\iy}^\un(\sA,\sB)$ in Definition \ref{ai4def8} is a dg-category over $R$. This follows from \eq{ai4eq16}--\eq{ai4eq17}, as $\m_k^\sB=0$ for $k\ge 3$ implies that the $A_\iy$-operations $\n_p$ in $\Fun_{A_\iy}^\un(\sA,\sB)$ have $\n_p=0$ for $p\ge 3$.
\end{rem}

Tabuada \cite{Tabu} proves:

\begin{thm}
\label{ai4thm1}
Let\/ $R$ be a commutative ring. There is a combinatorial model category structure on the category $\DGCat_R$ of dg-categories over $R,$ whose weak equivalences are quasi-equivalences. All objects in $\DGCat_R$ are fibrant.
\end{thm}

By properties of $\iy$-category localizations of model categories, and using Lurie \cite[Prop.~1.3.4.22]{Luri3} to deduce presentability, this implies:

\begin{cor}
\label{ai4cor2}
The $\iy$-category $\DGCat_R[\QE_\dg^{-1}]^\iy$ localizing quasi-eq\-uiv\-al\-ences in $\DGCat_R,$ as in {\rm\S\ref{ai26},} has all small homotopy limits and homotopy colimits. Also $\DGCat_R[\QE_\dg^{-1}]^\iy$ is \begin{bfseries}presentable\end{bfseries} in the sense of Lurie\/ {\rm\cite[\S 5]{Luri2}, \cite[\S 9.5]{Luri5}}.
\end{cor}

\begin{rem}
\label{ai4rem13}
For more information on model categories, see Hovey \cite{Hove}. In \S\ref{ai26} we discussed $\iy$-category localizations $\cC[\cW^{-1}]^\iy$ of a category with weak equivalences $(\cC,\cW)$. If $\cC$ is a model category and the weak equivalences $\cW$ come from the model category structure, then the theory of model categories provides tools for describing $\cC[\cW^{-1}]^\iy$ and $\Ho(\cC[\cW^{-1}]^\iy)$. So Theorem \ref{ai4thm1} helps us understand $\DGCat_R[\QE_\dg^{-1}]^\iy$. We return to this in \S\ref{ai410} and~\S\ref{ai92}.

In contrast to Theorem \ref{ai4thm1}, Canonaco--Ornaghi--Stellari \cite[\S 1.5]{COS1} show that $\ACat^\un_R$ {\it has no model structure}, as not all equalizers exist in $\ACat^\un_R$.
\end{rem}

The next theorem is primarily due to To\"en \cite{Toen}, although \cite{Toen} predates the theory of $\iy$-categories, and we have restated it in $\iy$-categorical language. See also Canonaco--Stellari \cite{CaSt}, who give a new proof of the results of \cite{Toen}, and Tanaka \cite[\S 3]{Tana}, who gives an alternative $\iy$-categorical proof of the existence of $\ot^\iy,\cHom$ using results of Haugseng and Gepner--Haugseng.

\begin{thm}
\label{ai4thm2}
{\bf(a)} There is a \begin{bfseries}symmetric monoidal structure\end{bfseries} $\ot^\iy$ on the $\iy$-category $\DGCat_R[\QE_\dg^{-1}]^\iy$. On objects $\sC,\sD$ it is given by
\begin{equation*}
\sC\ot^\iy\sD\cong \sC\ot^\bL\sD\cong\sC'\ot\sD',
\end{equation*}
where $\ot:\DGCat_R\t\DGCat_R\ra\DGCat_R$ is the classical tensor product of dg-categories, and\/ $\ot^\bL$ is its left derived tensor product, and $\sC',\sD'$ are cofibrant replacements for $\sC,\sD$ in the Tabuada model structure.
\smallskip

\noindent{\bf(b)} The monoidal structure $\ot^\iy$ is \begin{bfseries}closed\end{bfseries}. That is, there exist internal mapping objects $\cHom(\sC,\sD)$ in $\DGCat_R[\QE_\dg^{-1}]^\iy$ such that for any other object\/ $\sB,$ we have homotopy equivalences of mapping spaces
\begin{align*}
&\Hom_{\DGCat_R[\QE_\dg^{-1}]^\iy}\bigl(\sB,\cHom(\sC,\sD)\bigr)\\
&\qquad \simeq\Hom_{\DGCat_R[\QE_\dg^{-1}]^\iy}\bigl(\sB\ot^\iy\sC,\sD).
\end{align*}

\noindent{\bf(c)} The mapping object\/ $\cHom(\sC,\sD)$ may be described explicitly as a dg-category as the full subcategory of `right quasi-representable objects' in the dg-category of\/ $\sC\ot^\iy\sD^\rop$-modules (dg-functors from $\sC\ot^\iy\sD^\rop$ to complexes of\/ $R$-modules).

\smallskip

\noindent{\bf(d)} For all objects $\sC,\sD$ in $\DGCat_R[\QE_\dg^{-1}]^\iy,$ we have a homotopy equivalence
\begin{equation*}
\Hom_{\DGCat_R[\QE_\dg^{-1}]^\iy}(\sC,\sD)\simeq\sN_\dg\bigl(\cHom(\sC,\sD)\bigr){}^\Kan,
\end{equation*}
where $\sN_\dg(-)$ is the \begin{bfseries}dg-nerve\end{bfseries} $\sN_\dg:\DGCat_R\ra\Cat^1_\iy$ defined by Lurie\/ {\rm\cite[\S 1.3.1]{Luri3},} an $\iy$-category, and\/ $(\cdots)^\Kan$ is the maximal Kan complex inside it.
\end{thm}

\subsection{\texorpdfstring{Unital $A_\iy$-categories with split units}{Unital A∞-categories with split units}}
\label{ai46}

To have {\it split units\/} is a condition we can impose on a unital $A_\iy$-category. It becomes important when we want to compare unital $A_\iy$-categories with dg-categories, or with cohomologically unital $A_\iy$-categories. It is also part of a notion of `cofibrant' unital $A_\iy$-category, as explained in \S\ref{ai47}. The next definition follows Tanaka \cite[Def.~3.33]{Tana}, Ornaghi \cite[Def.~1.13]{Orna1} and Canonaco--Ornaghi--Stellari~\cite[Def.~3.4]{COS2}.

\begin{dfn}
\label{ai4def18}
We say that a unital $A_\iy$-category $\sA$ over a commutative ring $R$ has {\it split units\/} if for all objects $x$ in $\sA$ the map $R\ra\Hom^0(x,x)$, $\la\mapsto\la\boo_x$, is injective and splits as a map of $R$-modules (not necessarily as a map of complexes).
\end{dfn}

To have split units is a homotopy-invariant property:

\begin{prop}
\label{ai4prop5}
Let\/ $\fF:\sA\ra\sB$ be a unital homotopy equivalence of unital\/ $A_\iy$-categories over $R$. Then $\sA$ has split units if and only if\/ $\sB$ has split units.
\end{prop}

\begin{proof}
Suppose $\sB$ has split units, and let $x\in\Obj_\sA$. Consider the $R$-linear morphism
\begin{equation*}
\fF_1^{x,x}:\Hom_\sA^0(x,x)\longra\Hom_\sB^0(\fF(x),\fF(x))=R\cdot\boo^\sB_{\fF(x)}\op \Hom_\sB^0(\fF(x),\fF(x))',
\end{equation*}
where the right hand side is split as $\sB$ has split units.

As $\fF_1^{x,x}(\boo^\sA_x)=\boo^\sB_{\fF(x)}$, using $\fF_1^{x,x}$ we see that the map $R\ra\Hom^0_\sA(x,x)$, $\la\mapsto\la\boo^\sA_x$ is injective, and we have a splitting
\begin{equation*}
\Hom_\sA^0(x,x)=R\cdot\boo^\sA_x\op (\fF_1^{x,x})^{-1}\bigl(\Hom_\sB^0(\fF(x),\fF(x))'\bigr).
\end{equation*}
Thus $\sA$ has split units. Applying the same argument to a quasi-inverse $\fG:\sB\ra\sA$ for $\fF$, we see that if $\sA$ has split units then so does~$\sB$.
\end{proof}

\begin{rem}
\label{ai4rem14}
Some natural unital $A_\iy$-categories $\sA$ contain a `zero object' $z\in\Obj\sA$ with $\Hom(z,z)=0$, and $\Hom(x,z)=\Hom(z,x)=0$ for all $x$ in $\Obj\sA$. Then $\boo_z=0$, so the map $R\ra\Hom^0(z,z)$, $\la\mapsto\la\boo_z$ is not injective, and $\sA$ does not have split units, and nor does any unital $A_\iy$-category $\sB$ homotopy equivalent to $\sA$. This holds even if $R$ is a field.
\end{rem}

Here is Canonaco--Ornaghi--Stellari \cite[Lem.~4.1]{COS2}.

\begin{lem}
\label{ai4lem1}
Let\/ $\fF:\sA\ra\sB$ be a cohomologically unital functor between unital\/ $A_\iy$-categories $\sA,\sB$ over a commutative ring $R$. If\/ $\sA$ has split units then there exists a unital $A_\iy$-functor $\fF':\sA\ra\sB$ and a cohomologically unital\/ $A_\iy$-natural isomorphism~$\fS:\fF\Ra\fF'$.
\end{lem}

\begin{ex}
\label{ai4ex1}
Let $\sA,\sB$ be the unital dg-categories over the commutative ring $R=\Z_4$ with one object $*$ and, regarding $\Z_2$ as a module over $\Z_4$,
\begin{align*}
\Hom_\sA(*,*)=\Z_2\an{e_\sA}\op\Z_2\an{\boo^\sA_*}, \quad \Hom_\sB(*,*)=\Z_4\an{e_\sB}\op\Z_4\an{\boo^\sB_*}, 
\end{align*}
where $e_\sA,e_\sB$ have degree $-1$ and $\m_1^{\sA;*,*}(e_\sA)=\boo^\sA_*$, $\m_1^{\sB;*,*}(e_\sB)=\boo^\sB_*$. Then $H^*(\Hom_\sA(*,*))=H^*(\Hom_\sB(*,*))=0$. There is a cohomologically unital $A_\iy$-functor $\fF:\sA\ra\sB$ with $\fF(*)=*$ and $\fF_k^{*,\ldots,*}=0$ for all $k\ge 1$. 

However, there are no unital $A_\iy$-functors $\fF':\sA\ra\sB$, as we cannot have $\fF_1^{\prime *,*}(\boo^\sA_*)=\boo^\sB_*$, since any $\Z_4$-linear morphism $\Z_2\ra\Z_4$ is zero. Here $\sA$ does not have split units. This shows that the condition that $\sA$ has split units in Lemma \ref{ai4lem1} is necessary.
\end{ex}

As in Canonaco--Ornaghi--Stellari \cite[Rem.~1.24 \& Cor.~2.6]{COS1}, \cite[Rem.~5.4]{COS2}, the next result follows from Lef\`evre-Hasegawa \cite[Lem.~8.2.1.3]{LeHa} and Lemma \ref{ai4lem1}. Lef\`evre-Hasegawa takes $R$ to be a field, but his argument works for general $R$ provided $\sA$ has split units. Here $\Fun_{A_\iy}^\cu(\sA,\sB)$ is unital, not merely cohomologically unital, as $\sA,\sB$ are unital.

\begin{prop}
\label{ai4prop6}
Let\/ $R$ be a commutative ring, and $\sA,\sB$ be unital\/ $A_\iy$-categories over $R,$ where\/ $\sA$ has split units. Then the inclusion
\begin{equation*}
\inc:\Fun_{A_\iy}^\un(\sA,\sB)\longhookrightarrow\Fun_{A_\iy}^\cu(\sA,\sB)
\end{equation*}
is a quasi-equivalence of cohomologically unital\/ $A_\iy$-categories.
\end{prop}

Lyubashenko \cite[Th.~8.8]{Lyub} proves part (a) of the next theorem (see also Lef\`evre-Hasegawa \cite[Th.~9.2.0.4]{LeHa}, Seidel \cite[Th.~2.9]{Seid}, or Fukaya \cite[Th.~8.6]{Fuka3} when $R$ is a field, or for free $R$-modules). We can deduce part (b) from (a) and Proposition \ref{ai4prop6} as follows. Given unital $\fF:\sA\ra\sB$ satisfying (i),(ii), part (a) gives a cohomologically unital quasi-inverse $\fG':\sB\ra\sA$. Proposition \ref{ai4prop6} in $\Fun_{A_\iy}^*(\sB,\sA)$ shows $\fG'$ is homotopic to unital $\fG:\sB\ra\sA$. Then we have $A_\iy$-natural isomorphisms $\fS:\fG\ci\fF\Ra\id_\sA$ and $\fT:\fF\ci\fG\Ra\id_\sB$, and using Proposition \ref{ai4prop6} in $\Fun_{A_\iy}^*(\sA,\sA)$ and $\Fun_{A_\iy}^*(\sB,\sB)$, we can replace these by unital $A_\iy$-natural isomorphisms, so $\fF$ is a unital homotopy equivalence.

\begin{thm}
\label{ai4thm3}
{\bf(a)} Let\/ $\fF:\sA\ra\sB$ be a cohomologically unital\/ $A_\iy$-functor of cohomologically unital\/ $A_\iy$-categories $\sA,\sB$ over $R$. Then $\fF$ is a cohomologically unital homotopy equivalence if and only if:
\begin{itemize}
\setlength{\itemsep}{0pt}
\setlength{\parsep}{0pt}
\item[{\bf(i)}] For all objects\/ $x,y$ in $\sA,$ the following is a chain homotopy equivalence of chain complexes
\end{itemize}
\e
\fF_1^{x,y}:\bigl(\Hom_\sA(x,y),\m_1^{\sA;x,y}\bigr)\longra \bigl(\Hom_\sB(\fF(x),\fF(y)),\m_1^{\sB;\fF(x),\fF(y)}\bigr).
\label{ai4eq35}
\e
\begin{itemize}
\setlength{\itemsep}{0pt}
\setlength{\parsep}{0pt}
\item[] If\/ $R$ is a field, this holds if\/ $H^*(\fF_1^{x,y})$ is an isomorphism on cohomology.
\item[{\bf(ii)}] $H^0(\fF_1):H^0(\sA)\ra H^0(\sB)$ is an equivalence of categories.
\end{itemize}

\noindent{\bf(b)} Suppose instead that\/ $\sA,\sB,\fF$ are unital, and\/ $\sA,\sB$ have split units. Then $\fF$ is a unital homotopy equivalence if and only if\/ {\bf(i)\rm,\bf(ii)} hold.
\end{thm}

The next theorem follows from Canonaco--Ornaghi--Stellari \cite[Prop.~3.5]{COS2}. (See \cite[Prop.~2.1]{COS1} and \cite[Cor.~1.4]{BLM} for related results.)

\begin{thm}
\label{ai4thm4}
Let\/ $R$ be a commutative ring. For every $\sA\in\ACat^\un_R$ there exists a unital quasi-equivalence $\fF:\sA\ra\sB$ in the sense of Definition\/ {\rm\ref{ai4def15}} with $\sB\in\DGCat_R\subset\ACat^\un_R,$ where we may take $\fF_\Obj:\Obj\sA\ra\Obj\sB$ to be a bijection. If\/ $\sA$ has split units, we may take $\fF$ to be a unital homotopy equivalence in the sense of Definition\/ {\rm\ref{ai4def14},} and then $\sB$ also has split units by Proposition\/~{\rm\ref{ai4prop5}}.
\end{thm}

\subsection{\texorpdfstring{H-projective $A_\iy$-categories}{H-projective A∞-categories}}
\label{ai47}

As in Remark \ref{ai4rem13}, $\ACat^\un_R$ and $\ACat^\cu_R$ are not model categories. However, the authors of \cite{COS1,COS2,Faon2,OhTa,Orna1,Tana,Pasc} treat $\ACat^\un_R,\ACat^\cu_R$ {\it as though they were model categories}, in which weak equivalences are quasi-equivalences, and all objects are fibrant. This is true of the closely related category $\DGCat_R$ by Theorem~\ref{ai4thm1}.

When working in model categories, to compute things like $\iy$-category mapping spaces, one needs `cofibrant replacements'. In \cite{COS2,OhTa,Orna1,Tana} the substitute for cofibrant objects in $\ACat^\cu_R$ are {\it h-projective\/} cohomologically unital $A_\iy$-categories. Although this is not explicit in  \cite{COS2,OhTa,Orna1,Tana}, the appropriate substitute for cofibrant objects in $\ACat^\un_R$ are {\it h-projective\/} unital $A_\iy$-categories {\it with split units}, as in \S\ref{ai46}. The next definition follows Ornaghi \cite[\S 6.1]{Orna1} and Tanaka~\cite[Rem.s 1.6, 1.7, Def.~3.6 \& \S 3.11]{Tana}.

\begin{dfn}
\label{ai4def19}
Let $R$ be a commutative ring. A complex of $R$-modules $(\cA,\d)$ is called {\it h-projective} (short for {\it homotopically projective\/}) if whenever $(\cB,\d)$ is a complex of $R$-modules with $H^*(\cB,\d)=0$ then the morphism complex $(\Hom(\cA,\cB),\d)$ has~$H^*(\Hom(\cA,\cB),\d)=0$.

If $(\cA,\d)$ is a complex of {\it free\/} $R$-modules it is automatically h-projective. If $R$ is a field then all complexes $(\cA,\d)$ are h-projective. An $A_\iy$-category or dg-category $\sA$ over $R$ is called {\it h-projective\/} if for all $x,y\in\Obj\sA$, the complex $(\Hom(x,y),\m_1^{x,y})$ is h-projective. Oh--Tanaka \cite[Def.~1.2]{OhTa} call h-projective $A_\iy$-categories {\it cofibrant}. Note that $\sA$ is automatically h-projective if $R$ is a field, or if all Hom spaces $\Hom_\sA(x,y)$ are free $R$-modules.

Define full subcategories
\e
\begin{split}
\DGCat_R^\rsu,\DGCat_R^\hp,\DGCat_R^\hs&\subset\DGCat_R,\\
\DGCat_R^{A_\iy,\rsu},\DGCat_R^{A_\iy,\hp},\DGCat_R^{A_\iy,\hs}&\subset\DGCat^{A_\iy}_R,\\
\ACat_R^{\un,\rsu},\ACat_R^{\un,\hp},\ACat_R^{\un,\hs}&\subset\ACat_R^\un,\\
\ACat_R^{\cu,\hp}&\subset\ACat^\cu_R,
\end{split}
\label{ai4eq36}
\e
where `$\rsu$' means objects have split units, `$\hp$' means objects are h-projective, and `$\hs$' means objects are h-projective with split units.
\end{dfn}

The following proposition comes from Tanaka~\cite[Prop~3.8]{Tana}.

\begin{prop}
\label{ai4prop7}
Let\/ $f:(\cA,\d)\ra(\cB,\d)$ be a morphism of h-projective complexes over a commutative ring $R$ such that\/ $H^k(f):H^k(\cA,\d)\ra H^k(\cB,\d)$ is an isomorphism for all\/ $k\in\Z$. Then $f$ is a chain homotopy equivalence.
\end{prop}

As in Ornaghi \cite[Cor.~2.5]{Orna1}, Theorem \ref{ai4thm3} and Proposition \ref{ai4prop7} yield:

\begin{cor}
\label{ai4cor3}
{\bf(a)} Let\/ $R$ be a commutative ring and\/ $\fF:\sA\ra\sB$ be a quasi-equivalence of h-projective cohomologically unital\/ $A_\iy$-categories over $R,$ as in Definition\/ {\rm\ref{ai4def15}}. Then $\fF$ is a cohomologically unital homotopy equivalence, as in Definition\/~{\rm\ref{ai4def14}}.
\smallskip

\noindent{\bf(b)} Let\/ $\fF:\sA\ra\sB$ be a unital quasi-equivalence of h-projective unital\/ $A_\iy$-categories over $R$ with split units, as in Definitions\/ {\rm\ref{ai4def14}--\ref{ai4def15}}. Then $\fF$ is a unital homotopy equivalence, as in Definition\/~{\rm\ref{ai4def14}}.
\end{cor}

\begin{rem}
\label{ai4rem15}
We said above that a cohomologically unital $A_\iy$-category is regarded as `cofibrant' if it is h-projective, and a unital $A_\iy$-category is regarded as `cofibrant' if it is h-projective and has split units. Corollary \ref{ai4cor3} justifies this: between cofibrant objects, quasi-equivalences agree with homotopy equivalences.
\end{rem}

To be h-projective is a homotopy-invariant property.

\begin{prop}
\label{ai4prop8}
Let\/ $R$ be a commutative ring and\/ $\fF:\sA\ra\sB$ be a (cohomologically) unital homotopy equivalence of (cohomologically) unital\/ $A_\iy$-categories over $R$. Then $\sA$ is h-projective if and only if $\sB$ is h-projective.
\end{prop}

\begin{proof}
Note that if $\phi:(\cC,\d)\ra(\cD,\d)$ is a chain homotopy equivalence of complexes of $R$-modules then $(\cC,\d)$ is h-projective if and only if $(\cD,\d)$ is h-projective. Suppose $\sB$ is h-projective. By Theorem \ref{ai4thm3}(a), equation \eq{ai4eq35} is a chain homotopy equivalence. As $\sB$ is h-projective, the right hand side of \eq{ai4eq35} is h-projective, so the left hand side is h-projective for all $x,y\in\Obj\sA$, and thus $\sA$ is h-projective. Applying the same argument to a quasi-inverse $\fG:\sB\ra\sA$ for $\fF$, we see that if $\sA$ is h-projective then so is~$\sB$.
\end{proof}

Ornaghi \cite[Th.s B$'$, 5.1 \& 5.2]{Orna1} (see also Tanaka \cite[Rem.~1.8]{Tana}) proves:

\begin{thm}
\label{ai4thm5}
Let\/ $\sA$ be a unital\/ $A_\iy$-category over a commutative ring $R$. Then there exists an h-projective unital\/ $A_\iy$-category $\sA'$ over $R$ and a unital quasi-equivalence $\Phi:\sA'\ra\sA$. We call\/ $\Phi$ an \begin{bfseries}h-projective resolution of\/\end{bfseries}~$\sA$.

We may take $\sA'$ to have the same objects as $\sA,$ and to have split units. If\/ $\sA$ is a dg-category, the above holds with\/ $\sA'$ a dg-category and\/ $\Phi$ a dg-functor.
\end{thm}

\subsection{The simplicial nerve construction}
\label{ai48}

Example \ref{ai2ex5} defined the nerve functor $\sN:\Cat^1\ra\Cat^1_\iy$ from categories to $\iy$-categories. Faonte \cite[Def.~2.8]{Faon1} defines a functor $\sN_{A_\iy}:\ACat_R^\un\ra\Cat^1_\iy$, called the {\it simplicial nerve}, from the category of unital $A_\iy$-categories over $R$ to the ordinary category of $\iy$-categories $\Cat^1_\iy$ from Definition \ref{ai2def4}. (Note that this is different to the $\iy$-category of $\iy$-categories $\sCat_\iy$ in \S\ref{ai24}.) 

Faonte works with $R$ a field, but his construction works just as well with $R$ a commutative ring. As in \cite[Cor.~2.10]{Faon1}, the restriction $\sN_{A_\iy}\vert_{\DGCat_R}$ of $\sN_{A_\iy}$ to the subcategory $\DGCat_R\subset\ACat_R^\un$ of dg-categories over $R$ from Definition \ref{ai4def17} is the {\it dg-nerve\/} $\sN_\dg:\DGCat_R\ra\Cat^1_\iy$ defined by Lurie \cite[\S 1.3.1]{Luri3}. Le Grignou \cite{LeGr} studies the simplicial nerve using operads.

\begin{dfn}
\label{ai4def20}
Let $R$ be a commutative ring. For each $N\ge 0$, define a unital $A_\iy$-category $\sC^\De_{[N]}$ over $R$ to have objects $\Obj\sC^\De_{[N]}=\{0,1,\ldots,N\}$ and Hom spaces
\begin{equation*}
\Hom(i,j)=\begin{cases} R\cdot\om_{ij}, & i\le j, \\
0, & i>j,
\end{cases}
\end{equation*}
with $\deg\om_{ij}=0$, and units $\boo_i=\om_{ii}$, and maps $\m_k^{i_0,\ldots,i_k}$ given by
\e
\m_2^{i_0,i_1,i_2}(\om_{i_0i_1}\!\ot\!\om_{i_1i_2})\!=\!\om_{i_0i_2},\; 0\!\le\! i_0\!\le\! i_1\!\le\! i_2\!\le\! N,\; \m_k^{i_0,\ldots,i_k}\!=\!0,\; k\!\ne\! 2.
\label{ai4eq37}
\e

If $\al:[N]\ra[M]$ is a morphism in $\De$ in Definition \ref{ai2def1}, so that $\al:\{0,\ldots,N\}\ab\ra\{0,\ldots,M\}$ is an order-preserving function, we define a unital $A_\iy$-functor $\sC^\De_\al:\sC^\De_{[N]}\ra \sC^\De_{[M]}$ by $\Obj\sC^\De_\al=\al$, and to have maps.
\begin{equation*}
(\sC^\De_\al)_1^{i_0,i_1}:\om_{i_0i_1}\longmapsto \om_{\al(i_0)\al(i_1)},\;\> 0\le i_0\le i_1\le N,\quad (\sC^\De_\al)_k^{i_0,\ldots,i_k}=0,\;\> k\ne 1.
\end{equation*}
This defines a functor $\sC^\De:\De\ra\ACat^\un_R$ mapping $[N]\mapsto\sC^\De_{[N]}$ on objects and $\al\mapsto\sC^\De_\al$ on morphisms. In fact $\sC^\De_{[N]}$ is a dg-category, and $\sC^\De_\al$ is a dg-functor, so $\sC^\De$ is a functor~$\De\ra\DGCat_R\subset\ACat^\un_R$.

Noting that $\SSet$ is the category of functors $S:\De^\rop\ra\Set$, define a functor $\sN_{A_\iy}:\ACat_R^\un\ra\SSet$ by $\sN_{A_\iy}(-)=\Hom_{\ACat_R^\un}(\sC^\De,-)$. That is, $\sN_{A_\iy}$ maps an object $\sA$ to the functor $\sN_{A_\iy}(\sA):\De^\rop\ra\Set$, $\sN_{A_\iy}(\sA)(-)=\Hom_{\ACat_R^\un}(\sC^\De(-),\sA)$, and $\sN_{A_\iy}$ maps a morphism $\fF:\sA\ra\sB$ to the natural transformation $\sN_{A_\iy}(\fF):\sN_{A_\iy}(\sA)\Ra\sN_{A_\iy}(\sB)$, $\sN_{A_\iy}(\fF)=-\ci\fF:\Hom_{\ACat_R^\un}(\sC^\De(-),\sA)\ra \Hom_{\ACat_R^\un}(\sC^\De(-),\sB)$.

Faonte shows \cite[Prop.~2.15]{Faon1} that $\sN_{A_\iy}$ maps to~$\Cat^1_\iy\subset\SSet$.
\end{dfn}

The next proposition, which is easy to prove, shows that we should regard $\sN_{A_\iy}(\sA)$ as an $\iy$-categorical enhancement of~$H^0(\sA)$.

\begin{prop}
\label{ai4prop9}
Let $\sA$ be a unital\/ $A_\iy$-category over $R$. Then there is a strict isomorphism of ordinary categories $H^0(\sA)\cong\Ho(\sN_{A_\iy}(\sA)),$ where $H^0(\sA)$ is the cohomology category of\/ $\sA$ from Definition\/ {\rm\ref{ai4def4},} and\/ $\Ho(\sN_{A_\iy}(\sA))$ is the homotopy category from Definition\/ {\rm\ref{ai2def5}} of\/ $\sN_{A_\iy}(\sA)$ above. This extends to a natural isomorphism $H^0\cong\Ho\ci\sN_{A_\iy}$ of functors\/~$\ACat_R^\un\ra\Cat^1$.

This implies that homotopy equivalences $f:x\ra y$ in $\sA$ in the sense of Definition\/ {\rm\ref{ai4def5}} are the same thing as equivalences in the $\iy$-category $\sN(\sA)$ in the sense of Definition\/~{\rm\ref{ai2def5}}.
\end{prop}

Faonte \cite[Prop.s 2.9 \& 2.11]{Faon1} describes the $\iy$-category $\sN_{A_\iy}(\sA)$ explicitly. We have rewritten it for compatibility with our sign and order conventions.

\begin{prop}
\label{ai4prop10}
Let\/ $\sA$ be a unital\/ $A_\iy$-category over $R,$ so that\/ $\sN_{A_\iy}(\sA)$ is an $\iy$-category. Then
\begin{itemize}
\setlength{\itemsep}{0pt}
\setlength{\parsep}{0pt}
\item[{\bf(a)}] An $N$-simplex in $\sN_{A_\iy}(\sA)$ is equivalent to data $\bigl((x_i)_{i=0}^N,(f_k^{i_0,\ldots,i_k})^{\bs i}_{k\ge 1}\bigr),$ where $x_0,\ldots,x_N$ are objects in $\sA,$ and for all\/ $k\ge 1$ and\/ $0\ab\le\ab i_0\ab\le\ab i_1\ab\le\ab\cdots\ab\le\ab i_k\ab\le N,$ we have $f_k^{i_0,\ldots,i_k}\in\Hom_\sA^{1-k}(x_{i_0},x_{i_k}),$ satisfying
\e
f_1^{i_0i_0}=\boo_{x_{i_0}}, \quad f_{k+1}^{i_0,i_1,\ldots,i_{j-1},i_j,i_j,i_{j+1},\ldots,i_k}=0, \;\> 0\le j\le k,\;\> k\ne 1,
\label{ai4eq38}
\e
and for all\/ $n\ge 1$ and\/ $0\le i_0\le i_1\le\cdots\le i_n\le N$ 
\ea
&\sum_{0<j<n}(-1)^{j-1}f_{n-1}^{i_0,i_1,\ldots,i_{j-1},i_{j+1},\ldots,i_n}
\label{ai4eq39}\\
&=\sum_{\begin{subarray}{l} k\ge 1,\; 0=l_0<l_1<\cdots <l_k=n\end{subarray}}
\begin{aligned}[t]&(-1)^{\sum_{i=1}^k(l_i-l_{i-1}+1)(k-i)}\m^{\sA;x_{l_0},x_{l_1},\ldots,x_{l_k}}_k\ci{}\\
&\bigl(f_{l_1-l_0}^{i_{l_0},\ldots,i_{l_1}}\ot f_{l_2-l_1}^{i_{l_1},\ldots,i_{l_2}}\ot\cdots\ot f_{l_k-l_{k-1}}^{i_{l_{k-1}},\ldots,i_{l_k}}).\end{aligned} 
\nonumber
\ea
\item[{\bf(b)}] If\/ $\al:[\ti N]\ra[N]$ is a morphism in $\De,$ so $\al:\{0,\ldots,\ti N\}\ra\{0,\ldots,N\}$ is a nondecreasing function, then $\sN_{A_\iy}(\sA)(\al)$ maps $\bigl((x_i)_{i=0}^N,(f_k^{i_0,\ldots,i_k})^{\bs i}_{k\ge 1}\bigr)$ to $\bigl((\ti x_i)_{i=0}^{\ti N},(\ti f_k^{i_0,\ldots,i_k})^{\bs i}_{k\ge 1}\bigr),$ where $\ti x_i=x_{\al(i)}$ and\/ $\ti f_k^{i_0,\ldots,i_k}=f_k^{\al(i_0),\ldots,\al(i_k)}$ for all\/ $k\ge 1$ and\/ $0\le i_0\le i_1\le\cdots\le i_k\le\ti N$.
\end{itemize}
\end{prop}

Observe in (a) that if $i_{j-1}=i_j$ for any $j=1,\ldots,k$ then $f_k^{i_0,\ldots,i_k}$ is uniquely determined by \eq{ai4eq38}. Thus we can restrict to data $f_k^{i_0,\ldots,i_k}$ for $0\ab\le\ab i_0\ab<\ab i_1\ab<\ab\cdots\ab<\ab i_k\ab\le N$. This gives the following way to rewrite Proposition \ref{ai4prop10}, which will be useful in Chapter \ref{ai9}. In \eq{ai4eq40}, the last two cases are deduced from~\eq{ai4eq38}. 

\begin{cor}
\label{ai4cor4}
Let\/ $\sA$ be a unital\/ $A_\iy$-category over $R,$ so that\/ $\sN_{A_\iy}(\sA)$ is an $\iy$-category. Then
\begin{itemize}
\setlength{\itemsep}{0pt}
\setlength{\parsep}{0pt}
\item[{\bf(a)}] An $N$-simplex in $\sN_{A_\iy}(\sA)$ is equivalent to data $\bigl((x_i)_{i=0}^N,(f_k^{i_0,\ldots,i_k})^{\bs i}_{k\ge 1}\bigr),$ where $x_0,\ldots,x_N$ are objects in $\sA,$ and for all\/ $k\ge 1$ and\/ $0\ab\le\ab i_0\ab<\ab i_1\ab<\ab\cdots\ab<\ab i_k\ab\le N,$ we have $f_k^{i_0,\ldots,i_k}\in\Hom_\sA^{1-k}(x_{i_0},x_{i_k}),$ satisfying \eq{ai4eq39} for all\/ $1\le n\le N$ and\/ $0\le i_0< i_1<\cdots< i_n\le N$.
\item[{\bf(b)}] If\/ $\al:[\ti N]\ra[N]$ is a morphism in $\De,$ so $\al:\{0,\ldots,\ti N\}\ra\{0,\ldots,N\}$ is a nondecreasing function, then $\sN_{A_\iy}(\sA)(\al)$ maps $\bigl((x_i)_{i=0}^N,(f_k^{i_0,\ldots,i_k})^{\bs i}_{k\ge 1}\bigr)$ to $\bigl((\ti x_i)_{i=0}^{\ti N},(\ti f_k^{i_0,\ldots,i_k})^{\bs i}_{k\ge 1}\bigr),$ where $\ti x_i=x_{\al(i)}$ and for all\/ $1\le k\le N$ and\/ $0\le\ab i_0\ab<\ab i_1<\ab\cdots\ab<i_k\le\ti N$ we have
\e
\ti f_k^{i_0,\ldots,i_k}=\begin{cases}
f_k^{\al(i_0),\ldots,\al(i_k)} & \text{if\/ $\al\vert_{\{i_0,\ldots,i_k\}}$ is injective,} \\
\boo_{x^i} & \text{if\/ $k=1$ and\/ $\al(i_0)=\al(i_1)=i,$} \\
0 & \text{otherwise.}
\end{cases}
\label{ai4eq40}
\e
\end{itemize}
\end{cor}

\begin{rem}
\label{ai4rem16}
The simplicial nerve construction does not work nicely for {\it nonunital\/} $A_\iy$-categories, because we need units in \eq{ai4eq38} and \eq{ai4eq40}.
\end{rem}

\begin{thm}
\label{ai4thm6}
Let\/ $\sA$ be a unital\/ $A_\iy$-category over a commutative ring $R,$ and\/ $x$ an object in $\sA$. Write\/ $\sN_{A_\iy}(\sA)^\Kan$ for the maximal Kan complex in $\sN_{A_\iy}(\sA),$ as in Definition\/ {\rm\ref{ai2def7},} and\/ $\sN_{A_\iy}(\sA)^\Kan_x$ for the connected component of\/ $\sN_{A_\iy}(\sA)^\Kan$ containing $x$. Then the homotopy groups of\/ $\sN_{A_\iy}(\sA)^\Kan_x$ are given by
\e
\pi_1\bigl(\sN_{A_\iy}(\sA)^\Kan_x\bigr)\cong\bigl\{[f]\,\vert\,\text{$[f]:x\ra x$ is an isomorphism in $H^0(\sA)$}\bigr\},
\label{ai4eq41}
\e
with group operation composition in $H^0(\sA),$ and for $n\ge 2$
\e
\pi_n\bigl(\sN_{A_\iy}(\sA)^\Kan_x\bigr)\cong H^{1-n}\bigl(\Hom_\sA(x,x),\m_1^{x,x}\bigr),
\label{ai4eq42}
\e
with group operation addition in cohomology.

\end{thm}

\begin{proof}
As $\sN_{A_\iy}(\sA)^\Kan_x$ is a Kan complex, elements of $\pi_1\bigl(\sN_{A_\iy}(\sA)^\Kan_x\bigr)$ are homotopy classes of morphisms $f:x\ra x$ in $\sN_{A_\iy}(\sA)^\Kan_x$. Such $f$ are morphisms $f:x\ra x$ in $\sN_{A_\iy}(\sA)$ which are equivalences, and hence morphisms $f:x\ra x$ in $\sA$ for which $[f]:x\ra x$ is an isomorphism in~$H^0(\sA)$. 

Two such morphisms $f,f':x\ra x$ are homotopic if there exists a 2-simplex $\eta$ in $\sN_{A_\iy}(\sA)^\Kan_x$ of the form
\begin{equation*}
\xymatrix@C=70pt@R=15pt{
*+[l]{x\;\overset{0}{\bu}} \ar@<.8ex>[d]_{\id_x} \ar@<.6ex>[dr]^{f'} \ar@{}@<-1ex>[dr]_(0.3)\eta \\
*+[l]{x\;\underset{1}{\bu}} \ar@<.7ex>[r]_{f} & *+[r]{\underset{2}{\bu}\; x.} }
\end{equation*}
By Corollary \ref{ai4cor4}(a), we can show such $\eta$ is equivalent to $\eta_2^{012}\in\Hom_\sA^{-1}(x,x)$ satisfying $f'=f+\m_1^{x,x}\ci\eta_2^{012}$. Thus $f,f'$ are homotopic if and only if $[f]=[f']$ in $\Ho(\sA)$. This proves \eq{ai4eq41} as sets. To see that the group operation is composition in $H^0(\sA)$, note that for all isomorphisms $f,f':x\ra x$ in $\sA$ there is a 2-simplex $\ze$ in $\sN_{A_\iy}(\sA)^\Kan_x$ as follows, with $\ze_2^{012}=0$
\begin{equation*}
\xymatrix@C=70pt@R=15pt{
*+[l]{x\;\overset{0}{\bu}} \ar@<.8ex>[d]_{f} \ar@<.6ex>[dr]^{f'\ci f} \ar@{}@<-1ex>[dr]_(0.3)\ze \\
*+[l]{x\;\underset{1}{\bu}} \ar@<.7ex>[r]_{f'} & *+[r]{\underset{2}{\bu}\; x.} }
\end{equation*}

Now let $n\ge 2$. Then elements of $\pi_n\bigl(\sN_{A_\iy}(\sA)^\Kan_x\bigr)$ are homotopy classes of $n$-simplices $S:\De^n\ra\sN_{A_\iy}(\sA)^\Kan_x$ such that $S\vert_{\pd\De^n}$ is constant at $x$. By Corollary \ref{ai4cor4}(a), $S$ is equivalent to data $\bigl((x_i)_{i=0}^n,(S_k^{i_0,\ldots,i_k})^{\bs i}_{k\ge 1}\bigr)$ in which $x_i=x$ for all $i$, and as $S\vert_{\pd\De^n}\equiv x$, the only possibly nonzero $S_k^{i_0,\ldots,i_k}$ is $S_n^{01\cdots n}$ in $\Hom_\sA^{1-n}(x,x)$. Then \eq{ai4eq39} becomes $\m_1^{x,x}(S_n^{01\cdots n})=0$, so we have a cohomology class $[S_n^{01\cdots n}]\in H^{1-n}\bigl(\Hom_\sA(x,x),\m_1^{x,x}\bigr)$.

Two such $S,S'$ are homotopic if there is an $(n+1)$-simplex $T$ in $\sN_{A_\iy}(\sA)^\Kan_x$ such that two prescribed $\De^n$-faces of $T$ are $S,S'$, and the rest of $T\vert_{\pd\De^{n+1}}$ is constant at $x$. By Corollary \ref{ai4cor4}(a), we find that $T$ is equivalent to $T_{n+1}^{01\cdots n+1}$ in $\Hom_\sA^{-n}(x,x)$ satisfying $\m_1^{x,x}(T_{n+1}^{01\cdots n+1})=S_n^{\prime 01\cdots n}-S_n^{01\cdots n}$. Thus $S,S'$ are homotopic if and only if $[S_n^{01\cdots n}]=[S_n^{\prime 01\cdots n}]$ in $H^{1-n}\bigl(\Hom_\sA(x,x),\m_1^{x,x}\bigr)$. 

This proves \eq{ai4eq42} as sets. To see the group operation is addition, $[S_n^{01\cdots n}]\cdot[S_n^{\prime 01\cdots n}]=[S_n^{01\cdots n}+S_n^{\prime 01\cdots n}]$, consider an $(n+1)$-simplex $T$ in $\sN_{A_\iy}(\sA)^\Kan_x$ which has three prescribed $\De^n$-faces $S,S',S+S'$, and the rest of $T\vert_{\pd\De^{n+1}}$ is constant at $x$, and $T_{n+1}^{01\cdots n+1}=0$, and check that \eq{ai4eq39} holds.
\end{proof}

Note that the next corollary only needs $\fF$ to be a quasi-equivalence, not a homotopy equivalence.

\begin{cor}
\label{ai4cor5}
Suppose $\fF:\sA\ra\sB$ is a unital quasi-equivalence of unital\/ $A_\iy$-categories over a commutative ring $R$. Then the following is a homotopy equivalence:
\begin{equation*}
\sN_{A_\iy}(\fF)^\Kan:\sN_{A_\iy}(\sA)^\Kan\longra\sN_{A_\iy}(\sB)^\Kan.
\end{equation*}
\end{cor}

\begin{proof}
Connected components of $\sN_{A_\iy}(\sA)^\Kan$ and $\sN_{A_\iy}(\sB)^\Kan$ correspond to isomorphism classes of objects in $H^0(\sA)$ and $H^0(\sB)$, so $\sN_{A_\iy}(\fF)^\Kan$ induces a 1-1 correspondence of connected components as $H^0(\fF):H^0(\sA)\ra H^0(\sB)$ is an equivalence of categories.

Let $x\in\Obj\sA$, so that $\sN_{A_\iy}(\fF)^\Kan$ maps $\sN_{A_\iy}(\sA)^\Kan_x\ra \sN_{A_\iy}(\sB)^\Kan_{\fF(x)}$. Then $\sN_{A_\iy}(\fF)^\Kan$ gives an isomorphism $\pi_1\bigl(\sN_{A_\iy}(\sA)^\Kan_x\bigr)\ab\ra \pi_1\bigl(\sN_{A_\iy}(\sB)^\Kan_{\fF(x)}\bigr)$ by \eq{ai4eq41}, as $H^0(\fF)$ is an equivalence of categories. Also $\sN_{A_\iy}(\fF)^\Kan$ induces an isomorphism $\pi_n\bigl(\sN_{A_\iy}(\sA)^\Kan_x\bigr)\ab\ra \pi_n\bigl(\sN_{A_\iy}(\sB)^\Kan_{\fF(x)}\bigr)$ for $n\ge 2$ by \eq{ai4eq42}, as $\fF_1^{x,x}:\bigl(\Hom_\sA(x,x),\m_1^{x,x}\bigr)\ra \bigl(\Hom_\sB(\fF(x),\fF(x)),\ab\m_1^{\fF(x),\fF(x)}\bigr)$ is a quasi-isomorphism. Hence $\sN_{A_\iy}(\fF)^\Kan$ is a weak homotopy equivalence, so it is a homotopy equivalence by the Whitehead Theorem. 
\end{proof}

\subsection{\texorpdfstring{Pretriangulated $A_\iy$-categories and the twist functor}{Pretriangulated A∞-categories and the twist functor}}
\label{ai49}

\begin{dfn}
\label{ai4def21}
As in Gelfand and Manin \cite{GeMa}, a {\it triangulated category\/} $\cT$ is an additive category equipped with a `translation functor' and a class of `distinguished triangles' satisfying some axioms. They appear in many places in mathematics. {\it Derived categories\/} $D^b\cC$ of abelian categories $\cC$, such as derived categories of coherent sheaves $D^b\coh X$, are important examples. 

Let $R$ be a commutative ring, and $\sA$ an $A_\iy$-category over $R$, which should be weakly cohomologically unital, or cohomologically unital, or unital. There is a notion of when $\sA$ is {\it pretriangulated}, defined in Bespalov--Lyubashenko--Manzyuk \cite{BLM} and Seidel \cite[\S 3h]{Seid} (who calls them {\it triangulated\/}). They are defined directly in \cite{BLM,Seid} in terms of the existence of translations and distinguished triangles. We will not write this out, but instead give an indirect definition in Definition \ref{ai4def22} using the twist functor.

If $\sA$ is a pretriangulated $A_\iy$-category, its cohomology category $H^0(\sA)$ is triangulated. Note that although to be triangulated is additional structure on a category, to be pretriangulated is a property of an $A_\iy$-category. We write
\begin{equation*}
\ADCat_R^\un\subset\ACat_R^\un, \;\> \ADCat_R^\wc\subset\ACat_R^\wc,\;\> \ADCat_R^\cu\subset\ACat_R^\cu,
\end{equation*}
for the full subcategories whose objects are pretriangulated. For the subcategories $\ACat_R^{\un,\rsu},\ldots,\ACat_R^{\cu,\hp}$ in \eq{ai4eq36}, we also write $\ADCat_R^{\un,\rsu},\ab\ldots,\ab\ADCat_R^{\cu,\hp}$ for the full subcategories of these categories  whose objects are pretriangulated.
\end{dfn}

The next theorem is proved by Faonte \cite[Th.~3.18]{Faon1} for pretriangulated dg-categories over a field, and by Ornaghi \cite[Th.~5.1]{Orna3} in the general case.

\begin{thm}
\label{ai4thm7}
Let\/ $\sA$ be a pretriangulated unital\/ $A_\iy$-category over a commutative ring $R$. Then the simplicial nerve $\sN_{A_\iy}(\sA)$ from Definition\/ {\rm\ref{ai4def20}} is a \begin{bfseries}stable $\iy$-category\end{bfseries} in the sense of\/ {\rm\S\ref{ai27}}. The strict isomorphism $H^0(\sA)\cong\Ho(\sN_{A_\iy}(\sA))$ in Proposition\/ {\rm\ref{ai4prop9}}
is an isomorphism of triangulated categories.
\end{thm}

\begin{dfn}
\label{ai4def22}
Let $R$ be a commutative ring. The {\it twist functor\/} is a functor
\begin{equation*}
\Tw:\ACat_R\longra\ACat_R,
\end{equation*}
defined in Seidel \cite[\S 3]{Seid}, Keller \cite[\S 7.6]{Kell1}, Lef\`evre-Hasegawa \cite[\S 7.2]{LeHa}, and Bespalov--Lyubashenko--Manzyuk \cite[\S 12]{BLM}. We will construct it explicitly in \S\ref{ai58}, using filtered $A_\iy$-categories; doing it this way will make it easy to generalize $\Tw$ to $\iy$-categories in Chapter \ref{ai9}. Here we explain it informally, and review its properties.

The twist functor preserved the properties of being weakly cohomologically unital, cohomologically unital, and unital, and if $\sA$ is at least weakly cohomologically unital then $\Tw\sA$ is pretriangulated. Thus $\Tw$ restricts to functors
\e
\begin{gathered}
\Tw:\ACat^\wc_R\longra\ADCat^\wc_R, \qquad \Tw:\ACat^\cu_R\longra\ADCat^\cu_R,\\
\Tw:\ACat^\un_R\longra\ADCat^\un_R.
\end{gathered}
\label{ai4eq43}
\e

For each $A_\iy$-category $\sA$ there is an inclusion $\inc_\sA:\sA\hookra\Tw\sA$, defined explicitly in \S\ref{ai58}, which is a strict isomorphism with a full $A_\iy$-subcategory, so that $H^0(\inc_\sA):H^0(\sA)\hookra H^0(\Tw\sA)$ is full and faithful. Taken over all $\sA$ these form a natural transformation $\inc:\Id\Ra\Tw$ of functors $\ACat_R\ra\ACat_R$. Composing on the left or right with $\Tw$ gives two natural transformations $\id_{\Tw}*\inc,\inc*\id_{\Tw}:\Tw\Ra\Tw\ci\Tw$. For each $\sA$ these give $A_\iy$-functors $\Tw(\inc_\sA),\inc_{\Tw(\sA)}:\Tw\sA\hookra\Tw\Tw\sA$. These are full and faithful embeddings, strictly $A_\iy$-naturally isomorphic to each other, and they are homotopy equivalences for unital $\sA$ \cite[Prop.~12.11]{BLM}. A quasi-inverse $\Tw\Tw\sA\ra\Tw\sA$ is constructed in~\cite[\S 12.2.3]{BLM}. 

As in \cite[Lem.~3.33]{Seid} and \cite[Prop.~12.10]{BLM}, a weakly cohomological $A_\iy$-category $\sA$ is pretriangulated if and only if $\inc_\sA:\sA\hookra\Tw\sA$ is a quasi-equivalence (this happens provided $\inc_\sA$ is essentially surjective, and then $\inc_\sA$ is also a homotopy equivalence). Thus, we can define $\sA$ to be {\it pretriangulated\/} if $\inc_\sA:\sA\hookra\Tw\sA$ is a homotopy equivalence. If $\sA,\sB$ are homotopy equivalent then $\sA$ is pretriangulated if and only if $\sB$ is.

The definition of $\Tw\sA$ is similar to the definition of the bounded derived category $D^b\cC$ of an abelian category $\cC$. Objects of $\Tw\sA$ are {\it twisted complexes\/} in $\sA$. Basically, we make $\Tw\sA$ by starting from $\sA$ and adding many more objects: shifts $x[k]$ for $x\in\Obj\sA$ and $k\in\Z$, direct sums $x\op y$ of objects $x,y$, and iterated extensions of objects $x_1,\ldots,x_k$, to make $\Tw\sA$ closed under shifts $[k]$ and extensions in exact triangles.
\end{dfn}

\begin{rem}
\label{ai4rem17}
Here is an important application of the twist functor. If $(X,\om)$ is a symplectic manifold, one can try to define the {\it Fukaya category\/} $\sF(X,\om)$, a cohomologically unital $A_\iy$-category with objects Lagrangian submanifolds $L$ in $(X,\om)$ (plus extra data). Then the {\it derived Fukaya category\/} is $D^b\sF(X,\om)=H^0(\Tw(\sF(X,\om)))$, a triangulated category.

Kontsevich's Homological Mirror Symmetry Conjecture \cite{Kont} explains Mirror Symmetry of Calabi--Yau manifolds $X,\check X$ in terms of equivalences of triangulated categories $D^b\sF(X,\om)\simeq D^b\coh\check X$, $D^b\sF(\check X,\check\om)\simeq D^b\coh X$.
\end{rem}

\subsection{\texorpdfstring{$\iy$-categories of $A_\iy$-categories as localizations}{∞-categories of A∞-categories as localizations}}
\label{ai410}

We now discuss work by Canonaco--Ornaghi--Stellari \cite{COS1,COS2}, Faonte \cite{Faon2}, Oh--Tanaka \cite{OhTa}, Ornaghi \cite{Orna1}, Tanaka \cite{Tana}, and Pascaleff \cite{Pasc}, which studies $\iy$-categories of $A_\iy$-categories (or dg-categories) defined by starting with the category $\ACat^\un_R$ (or some variant), and localizing the quasi-equivalences $\QE_\un\subset\Mor(\ACat^\un_R)$ from Definition \ref{ai4def15}, using $\iy$-category localization from \S\ref{ai26}, to get an $\iy$-category $\ACat^\un_R[\QE_\un^{-1}]^\iy$. So for them (unlike us) the definition of their $\iy$-category is trivial, and they are concerned with its properties. 

The next theorem is proved by Canonaco--Ornaghi--Stellari \cite[Th.~A \& Rem.~3.13]{COS2}. Some of the equivalences in \eq{ai4eq44} were proved earlier by Pascaleff \cite[Cor.~5.2]{Pasc} when $R$ is a field, and by Tanaka~\cite[Th.~1.1]{Tana}.

\begin{thm}
\label{ai4thm8}
For any commutative ring $R,$ the inclusions of subcategories 
\begin{align*}
&\DGCat_R\longhookrightarrow\DGCat^{A_\iy}_R\longhookrightarrow\ACat_R^\un\longhookrightarrow\ACat_R^\cu\longhookrightarrow\ACat_R^\wc
\end{align*}
induce equivalences of\/ $\iy$-categories
\ea
&\DGCat_R[\QE_\dg^{-1}]^\iy\,{\buildrel\simeq\over\longra}\,\DGCat^{A_\iy}_R[\QE_{\dg,A_\iy}^{-1}]^\iy\,{\buildrel\simeq\over\longra}\,
\ACat_R^\un[\QE_\un^{-1}]^\iy\,{\buildrel\simeq\over\longra} 
\nonumber\\
&\ACat_R^\cu[\QE_\cu^{-1}]^\iy\,{\buildrel\simeq\over\longra}\,\ACat_R^\wc[\QE_\wc^{-1}]^\iy.
\label{ai4eq44}
\ea
\end{thm}

Combining Theorem \ref{ai4thm8} with Corollary \ref{ai4cor2} and Theorem \ref{ai4thm2} yields:

\begin{cor}
\label{ai4cor6}
{\bf(a)} The $\iy$-category $\ACat^\un_R[\QE_\un^{-1}]^\iy$ has all small homotopy limits and homotopy colimits.
\smallskip

\noindent{\bf(b)} The $\iy$-category $\ACat^\un_R[\QE_\un^{-1}]^\iy$ is presentable.
\smallskip

\noindent{\bf(c)} There is a symmetric monoidal structure $
\ot^\iy$ on $\ACat^\un_R[\QE_\un^{-1}]^\iy$.

\smallskip

\noindent{\bf(d)} The monoidal structure $\ot^\iy$ is closed, so there exist internal mapping objects $\cHom(\sC,\sD)$ in $\ACat^\un_R[\QE_\un^{-1}]^\iy$ with homotopy equivalences
\begin{align*}
&\Hom_{\ACat^\un_R[\QE_\un^{-1}]^\iy}\bigl(\sB,\cHom(\sC,\sD)\bigr)\\
&\qquad \simeq\Hom_{\ACat^\un_R[\QE_\un^{-1}]^\iy}\bigl(\sB\ot^\iy\sC,\sD).
\end{align*}

The analogues hold for $\ACat_R^\cu[\QE_\cu^{-1}]^\iy$.
\end{cor}

We noted in \S\ref{ai47} that the authors of \cite{COS1,COS2,Faon2,OhTa,Orna1,Tana,Pasc} treat $\ACat^\un_R$ {\it as though it was a model category}, in which weak equivalences are quasi-equiv\-al\-en\-ces, and cofibrant objects are h-projective unital $A_\iy$-categories with split units. The next three theorems use this idea, since Theorems \ref{ai4thm9}--\ref{ai4thm10} are like using cofibrant replacements in a model category, and Theorem \ref{ai4thm11} is an analogue of Remark \ref{ai2rem4}(d) for model categories.

Here is Canonaco--Ornaghi--Stellari \cite[Th.~C]{COS2}, which proves a claim by Kontsevich. A version was claimed earlier by Faonte \cite[Th.~1.2]{Faon2} for $R$ a field of characteristic zero. However, Canonaco et al.\ \cite[p.~4]{COS1} point out errors in \cite{Faon2}. In particular, \cite[Th.~1.2]{Faon2} is false as stated, as it has no cofibrant replacement of $\sC$ by $\sC'$, but this is necessary even for $R$ a field if $\sC$ has zero objects.

\begin{thm}
\label{ai4thm9}
Let\/ $\sC,\sD$ be dg-categories over a commutative ring $R$. Then as in Theorem\/ {\rm\ref{ai4thm2}(b)} there is an internal mapping object $\cHom(\sC,\sD)$ in the $\iy$-category $\DGCat_R[\QE_\dg^{-1}]^\iy,$ which is explicitly described in Theorem\/~{\rm\ref{ai4thm2}(b)}.

We can give an alternative description of\/ $\cHom(\sC,\sD)$ as
\begin{equation*}
\cHom(\sC,\sD)\simeq \Fun^\un_{A_\iy}(\sC',\sD).
\end{equation*}
Here $\Phi:\sC'\ra\sC$ is a quasi-equivalence in $\DGCat_R$ such that\/ $\sC'$ is h-projective with split units, which exists by Theorem\/ {\rm\ref{ai4thm5},} and\/ $\Fun^\un_{A_\iy}(\sC',\sD)$ is the $A_\iy$-category from\/ {\rm\S\ref{ai43},} which is a dg-category in this case as $\sC',\sD$ are dg-categories. Hence by Theorem\/ {\rm\ref{ai4thm2}(d)} we have a homotopy equivalence
\begin{equation*}
\Hom_{\DGCat_R[\QE_\dg^{-1}]^\iy}(\sC,\sD)\simeq\sN_{A_\iy}\bigl(\Fun^\un_{A_\iy}(\sC',\sD)\bigr){}^\Kan,
\end{equation*}
where $\cN_{A_\iy}(-)$ is as in\/ {\rm\S\ref{ai48},} and coincides with\/ $\cN_\dg(-)$ for dg-categories.
\end{thm}

The next theorem is taken from Tanaka \cite[Th.s 3.36 \& 3.37]{Tana}. It generalizes Theorem \ref{ai4thm9} to $A_\iy$-categories. 

\begin{thm}
\label{ai4thm10}
Let\/ $\sA,\sB$ be unital\/ $A_\iy$-categories over a commutative ring $R,$ and let\/ $\Phi:\sA'\ra\sA$ be an h-projective resolution of\/ $\sA,$ as in Theorem\/ {\rm\ref{ai4thm5},} such that\/ $\sA'$ has split units, as in Definition\/ {\rm\ref{ai4def18}}. Then the internal mapping object $\cHom(\sA,\sB)$ in $\ACat^\un_R[\QE_\un^{-1}]^\iy$ given by Corollary\/ {\rm\ref{ai4cor6}(d)} satisfies
\begin{equation*}
\cHom(\sA,\sB)\simeq \Fun^\un_{A_\iy}(\sA',\sB),
\end{equation*}
where $\Fun^\un_{A_\iy}(\sA',\sB)$ is as in\/ {\rm\S\ref{ai43}}. Also, we have a homotopy equivalence
\e
\Hom_{\ACat^\un_R[\QE_\un^{-1}]^\iy}(\sA,\sB)\simeq\sN_{A_\iy}\bigl(\Fun^\un_{A_\iy}(\sA',\sB)\bigr){}^\Kan,
\label{ai4eq45}
\e
where $\Hom_{\ACat^\un_R[\QE_\un^{-1}]^\iy}(\sA,\sB)$ is the mapping space as in\/ {\rm\S\ref{ai23}} in the $\iy$-category $\ACat^\un_R[\QE_\un^{-1}]^\iy,$ and\/ $\sN_{A_\iy}\bigl(\Fun^\un_{A_\iy}(\sA',\sB)\bigr){}^\Kan$ is the maximal Kan complex in the simplicial nerve of the\/ $A_\iy$-category\/ $\Fun^\un_{A_\iy}(\sA',\sB),$ as in\/~{\rm\S\ref{ai48}}.
\end{thm}

If a dg-category $\sA$ is cofibrant in the Tabuada model structure from Theorem \ref{ai4thm1} then it is h-projective with split units. Because of this, by the argument in Tanaka \cite[Cons.~3.20 \& Rem.s 3.21 \& 3.22]{Tana} we can show:

\begin{thm}
\label{ai4thm11}
Let $R$ be a commutative ring. Then the inclusion of subcategories $\DGCat_R^\hs\hookra\DGCat_R$ from Definition\/ {\rm\ref{ai4def19}} induces an equivalence of\/ $\iy$-categories
\e
\DGCat_R^\hs[\HE_{\dg,\hs}^{-1}]^\iy\,{\buildrel\simeq\over\longra}\, \DGCat_R[\QE_\dg^{-1}]^\iy.
\label{ai4eq46}
\e
\end{thm}

\begin{rem}
\label{ai4rem18}
If it were true (it is not) that $\ACat^\un_R,\ACat^\cu_R$ were model categories in which weak equivalences are quasi-equivalences, and cofibrant objects in $\ACat^\un_R$ are h-projective with split units, and cofibrant objects in $\ACat^\cu_R$ are h-projective, then by Remark \ref{ai2rem4}(d) and Corollary \ref{ai4cor3} we could deduce that the inclusions of subcategories $\ACat_R^{\un,\hs}\ab\hookra\ACat_R^\un$ and $\ACat_R^{\cu,\hp}\hookra\ACat_R^\cu$ from Definition \ref{ai4def19} induce equivalences
\ea
\ACat_R^{\un,\hs}[\HE_{\un,\hs}^{-1}]^\iy&\,{\buildrel\simeq\over\longra}\, \ACat_R^\un[\QE_\un^{-1}]^\iy,
\label{ai4eq47}\\
\ACat_R^{\cu,\hp}[\HE_{\cu,\hp}^{-1}]^\iy&\,{\buildrel\simeq\over\longra}\, \ACat_R^\cu[\QE_\cu^{-1}]^\iy.
\label{ai4eq48}
\ea

Following the philosophy of \cite{COS1,COS2,Faon2,OhTa,Orna1,Tana,Pasc}, one might expect \eq{ai4eq47}--\eq{ai4eq48} to be true even though $\ACat^\un_R,\ACat^\cu_R$ are not model categories. The author cannot find a proof of these in the literature. But we will prove \eq{ai4eq47} ourselves in Theorem \ref{ai9thm3} below.
\end{rem}

\subsection[\texorpdfstring{Faonte's $(\iy,1)$- and $(\iy,2)$-categories of $A_\iy$-categories}{Faonte's (∞,1)- and (∞,2)-categories of A∞-categories}]{Faonte's \kern -.1em $(\iy,1)$- \kern -.1em and \kern -.1em $(\iy,2)$-categories \kern -.1em of \kern -.1em $A_\iy$-categories}
\label{ai411}

We now review some ideas due to Faonte \cite{Faon2}. We note two things: firstly, Faonte works over a field $\K$ of characteristic zero. However, the parts of \cite{Faon2} we discuss work over a general commutative ring $R$, and we will state them like this. Secondly, Canonaco--Ornaghi--Stellari \cite[Intro]{COS1} point out errors in \cite{Faon2}, which affect \cite[\S A.4]{Faon2} and \cite[Th.~1.2]{Faon2}, which is false as stated.

\begin{dfn}
\label{ai4def23}
Let $\sA,\sB,\sC$ be unital $A_\iy$-categories over $R$. In Remark \ref{ai4rem7} we explained how work of Lyubashenko \cite{Lyub} defines an `$A_\iy$-bifunctor'
\begin{equation*}
\ka_{\sA\sB\sC}^\un:\Fun_{A_\iy}^\un(\sA,\sB)\t\Fun_{A_\iy}^\un(\sB,\sC)\longra\Fun_{A_\iy}^\un(\sA,\sC),
\end{equation*}
which acts on objects as $(\fF,\fG)\mapsto\fG\ci\fF$. 

Faonte \cite[\S 3.4, Def.~3.18]{Faon2} uses $\ka_{\sA\sB\sC}^\un$ (actually, the corresponding functor of counital dg-categories, as in \eq{ai4eq26}) to construct an explicit $\iy$-functor
\e
\begin{split}
\la_{\sA\sB\sC}&:\sN_{A_\iy}\bigl(\Fun^\un_{A_\iy}(\sA,\sB)\bigr)\t\sN_{A_\iy}\bigl(\Fun^\un_{A_\iy}(\sB,\sC)\bigr)\\
&\qquad \longra \sN_{A_\iy}\bigl(\Fun^\un_{A_\iy}(\sA,\sC)\bigr).
\end{split}
\label{ai4eq49}
\e
There are no arbitrary choices involved, $\la_{\sA\sB\sC}$ is defined as an explicit map on $N$-simplices. Restricting to maximal Kan complexes $(-)^\Kan$ gives a morphism 
\e
\begin{split}
\la_{\sA\sB\sC}^\Kan&:\sN_{A_\iy}\bigl(\Fun^\un_{A_\iy}(\sA,\sB)\bigr){}^\Kan\t\sN_{A_\iy}\bigl(\Fun^\un_{A_\iy}(\sB,\sC)\bigr){}^\Kan\\
&\qquad \longra \sN_{A_\iy}\bigl(\Fun^\un_{A_\iy}(\sA,\sC)\bigr){}^\Kan.
\end{split}
\label{ai4eq50}
\e
Note that spaces $\sN_{A_\iy}\bigl(\Fun^\un_{A_\iy}(\sA,\sB)\bigr){}^\Kan$ are models for mapping spaces in $\iy$-categories in Theorems \ref{ai4thm9} and \ref{ai4thm10}.
\end{dfn}

Faonte \cite[Prop.~3.19]{Faon2} proves:

\begin{prop}
\label{ai4prop11}
Let\/ $\sA,\sB,\sC,\sD$ be unital\/ $A_\iy$-categories over $R$. Then the $\iy$-functors $\la_{\sA\sB\sC}$ in \eq{ai4eq49} satisfy the identities
\ea
&\la_{\sA\sB\sD}\ci\bigl(\id_{\cdots(\sA,\sB)}\t\la_{\sB\sC\sD}\bigr)=
\la_{\sA\sC\sD}\ci\bigl(\la_{\sA\sB\sC}\t\id_{\cdots(\sC,\sD)}\bigr):
\nonumber\\
&\sN_{A_\iy}\bigl(\Fun^\un_{A_\iy}(\sA,\sB)\bigr)\t\sN_{A_\iy}\bigl(\Fun^\un_{A_\iy}(\sB,\sC)\bigr)
\t\sN_{A_\iy}\bigl(\Fun^\un_{A_\iy}(\sC,\sD)\bigr)
\nonumber\\
&\qquad \longra \sN_{A_\iy}\bigl(\Fun^\un_{A_\iy}(\sA,\sD)\bigr),
\label{ai4eq51}\\
&\la_{\sA\sA\sB}\vert_{\{\id_\sA\}\t\sN_{A_\iy}(\Fun^\un_{A_\iy}(\sA,\sB))}=
\id_{\sN_{A_\iy}(\Fun^\un_{A_\iy}(\sA,\sB))}:
\label{ai4eq52}\\
&*\t \sN_{A_\iy}\bigl(\Fun^\un_{A_\iy}(\sA,\sB)\bigr)\cong \sN_{A_\iy}\bigl(\Fun^\un_{A_\iy}(\sA,\sB)\bigr)\longra\sN_{A_\iy}\bigl(\Fun^\un_{A_\iy}(\sA,\sB)\bigr),
\nonumber\\
&\la_{\sA\sB\sB}\vert_{\sN_{A_\iy}(\Fun^\un_{A_\iy}(\sA,\sB))\t\{\id_\sB\}}=
\id_{\sN_{A_\iy}(\Fun^\un_{A_\iy}(\sA,\sB))}:
\label{ai4eq53}\\
&\sN_{A_\iy}\bigl(\Fun^\un_{A_\iy}(\sA,\sB)\bigr)\t *\cong \sN_{A_\iy}\bigl(\Fun^\un_{A_\iy}(\sA,\sB)\bigr)\longra\sN_{A_\iy}\bigl(\Fun^\un_{A_\iy}(\sA,\sB)\bigr).
\nonumber
\ea
Here for\/ \eq{ai4eq52} we restrict to the $0$-simplex $\id_\sA$ in $\sN_{A_\iy}\bigl(\Fun^\un_{A_\iy}(\sA,\sA)\bigr),$ regarded as a point\/ $*,$ and similarly for \eq{ai4eq53}. Note that\/ {\rm\eq{ai4eq51}--\eq{ai4eq53}} hold not just up to homotopy, but strictly as maps of simplicial sets.
\end{prop}

Recall as in \cite[Def.~1.2.39]{Land}, \cite[\S 1.1.4]{Luri2} that a {\it simplicial category\/} $\sS$ is a category enriched in simplicial sets. That is, for objects $X,Y$ in $\sS$, $\Hom_\sS(X,Y)$ is a simplicial set. Composition maps $\la_{XYZ}:\Hom_\sS(X,Y)\t\Hom_\sS(Y,Z)\ra\Hom_\sS(X,Z)$ are maps of simplicial sets. The identity $\id_X$ is a 0-simplex in $\Hom_\sS(X,X)$. Associativity and identity properties hold strictly in maps of simplicial sets. Write $\SCat$ for the category of simplicial categories.

Proposition \ref{ai4prop11} implies that we can make unital $A_\iy$-categories over $R$ into a simplicial category $\sS$ in two different ways:
\begin{itemize}
\setlength{\itemsep}{0pt}
\setlength{\parsep}{0pt}
\item[(i)] We can take $\Hom_\sS(\sA,\sB)=\sN_{A_\iy}\bigl(\Fun^\un_{A_\iy}(\sA,\sB)\bigr){}^\Kan$, with composition maps \eq{ai4eq50} and identities $\id_\sA\in\sN_{A_\iy}\bigl(\Fun^\un_{A_\iy}(\sA,\sA)\bigr){}^\Kan$. We will write $\ASCat_R^{\un,\Kan}$ for this simplicial category.
\item[(ii)] We can take $\Hom_\sS(\sA,\sB)=\sN_{A_\iy}\bigl(\Fun^\un_{A_\iy}(\sA,\sB)\bigr)$, with composition maps \eq{ai4eq49} and identities $\id_\sA\in\sN_{A_\iy}\bigl(\Fun^\un_{A_\iy}(\sA,\sA)\bigr)$. We will write $\ASCat_R^\un$ for this simplicial category, so $\ASCat_R^{\un,\Kan}\subset\ASCat_R^\un$.
\end{itemize}
Considered as higher categories, 2-morphisms in $\ASCat_R^\un$ are $A_\iy$-natural transformations as in \S\ref{ai44}, but 2-morphisms in $\ASCat_R^{\un,\Kan}$ are $A_\iy$-natural isomorphisms. We discuss (i),(ii) in the next two remarks.

\begin{rem}
\label{ai4rem19}
{\bf(a)} As in Faonte \cite[\S 3.1]{Faon2}, we should regard a simplicial category $\sS$ with all $\Hom_\sS(X,Y)$ Kan complexes, such as $\ASCat_R^{\un,\Kan}$, as a kind of $(\iy,1)$-category. Here $n$-morphisms in $\sS$ for $n\ge 1$ are $(n-1)$-morphisms in $\Hom_\sS(X,Y)$, so taking $\Hom_\sS(X,Y)$ a Kan complex means that $n$-morphisms in $\sS$ are invertible for $n\ge 2$, as we want in an $(\iy,1)$-category.
\smallskip

\noindent{\bf(b)} Bergner \cite{Berg} defined a model structure on $\SCat$. There is a Quillen equivalence between $\SCat$ with the Bergner model structure and $\SSet$ with the Joyal model structure. This induces the first in a chain of equivalences of $\iy$-categories, which give different models for the $\iy$-category of $\iy$-categories:
\begin{equation*}
\xymatrix@C=30pt{ \SCat[\cW_{\rm Berg}^{-1}]^\iy \ar[r]^\simeq & \SSet[\cW_{\rm Joy}^{-1}]^\iy \ar[r]^(0.45)\simeq & 
\Cat^1_\iy[{\rm Joy}^{-1}]^\iy \ar[r]^(0.55)\simeq & \sCat_\iy. }
\end{equation*}

The fibrant objects in the Bergner model structure are simplicial categories whose Hom spaces are Kan complexes, such as $\ASCat_R^{\un,\Kan}$.
\smallskip

\noindent{\bf(c)} Faonte works with $R$ a field $\K$ of characteristic zero. When $\sC,\sD$ are dg-categories, he claims \cite[\S 4]{Faon2} that the Hom spaces $\sN_{A_\iy}\bigl(\Fun^\un_{A_\iy}(\sC,\sD)\bigr){}^\Kan$ in $\ASCat_R^{\un,\Kan}$ are homotopic to those in $\DGCat_R[\QE_\dg^{-1}]^\iy$. In fact this may be false if $\sC$ has zero objects, because of the error in \cite[Th.~1.2]{Faon2} we noted before Theorem \ref{ai4thm9}. He then claims without proof \cite[\S 4]{Faon2} that the full simplicial subcategory $\mathop{\bf DGSCat}_R^\Kan\subset\ASCat_R^{\un,\Kan}$ with objects dg-categories should be equivalent as an $(\iy,1)$-category to~$\DGCat_R[\QE_\dg^{-1}]^\iy$.

\smallskip

\noindent{\bf(d)} The following observation may be new. As in \S\ref{ai1.62}, because of our Theorems \ref{ai9thm1}(a), \ref{ai9thm2}(d), \ref{ai9thm4} below, the author expects that Faonte's construction of $\ASCat_R^{\un,\Kan}$ works for general rings $R$, and gives an $(\iy,1)$-category which is a model for $\ACat_R^\un[\HE^{-1}]^\iy$, {\it not\/} $\ACat_R^\un[\QE^{-1}]^\iy$ or $\DGCat_R[\QE_\dg^{-1}]^\iy$, and should be equivalent as an $(\iy,1)$-category to our~$\ACat^{\un,\iy}_R$.
\end{rem}

\begin{rem}
\label{ai4rem20}
{\bf(a)} As in Faonte \cite[\S 3.4]{Faon2}, we should regard a simplicial category $\sS$ with all $\Hom_\sS(X,Y)$ $\iy$-categories, such as $\ASCat_R^\un$, as a kind of $(\iy,2)$-category. Here $n$-morphisms in $\sS$ for $n\ge 1$ are $(n-1)$-morphisms in $\Hom_\sS(X,Y)$, so taking $\Hom_\sS(X,Y)$ an $\iy$-category means that $n$-morphisms in $\sS$ are invertible for $n\ge 3$, as we want in an $(\iy,2)$-category. We should think of 2-morphisms in $\ASCat_R^\un$ as $A_\iy$-natural transformations, which need not be invertible unless they are $A_\iy$-natural isomorphisms.
\smallskip

\noindent{\bf(b)} In \cite[Def.~3.20]{Faon2}, Faonte applies a construction to convert $\ASCat_R^\un$ into a `fibrant Segal category in the model category $({\rm SSet},{\rm Joy})$', which is a widely accepted notion of $(\iy,2)$-category.
\smallskip

\noindent{\bf(c)} Any $(\iy,2)$-category has an associated $(\iy,1)$-category, in which one restricts to only invertible 2-morphisms. The $(\iy,1)$-category associated to $\ASCat_R^\un$ is~$\ASCat_R^{\un,\Kan}$.
\smallskip

\noindent{\bf(d)} In Theorem \ref{ai7thm4} we will construct an $\iy$-{\it bicategory\/} $(\hACat_R^{\un,\iy},\cT)$ of unital $A_\iy$-categories, in the sense of \S\ref{ai28}. As in Remark \ref{ai7rem4}(d), we expect this to be closely connected to Faonte's $(\iy,2)$-category.
\end{rem}

\section{\texorpdfstring{Filtered $A_\iy$-algebras and $A_\iy$-categories}{Filtered A∞-algebras and A∞-categories}}
\label{ai5}

As we explain in \S\ref{ai51}, {\it curved\/ $A_\iy$-algebras\/} and {\it curved\/ $A_\iy$-categories\/} are defined as for $A_\iy$-algebras and $A_\iy$-categories in \S\ref{ai32} and \S\ref{ai41}, but allowing operations $\m_k,\m_k^{x_0,\ldots,x_k}$ for $k\ge 0$ rather than $k\ge 1$. To define morphisms, composition, homotopies, \ldots, equations such as \eq{ai3eq7}, \eq{ai3eq10}, \eq{ai3eq12}, \eq{ai4eq8}, \ldots, become {\it infinite\/} sums. So the theory must be set in a context in which convergence of infinite sums makes sense. This leads to {\it filtered\/ $A_\iy$-algebras\/} and {\it filtered\/ $A_\iy$-categories}, written using complete filtered modules with a notion of convergence.

There are three main notions of filtered $A_\iy$-algebra and filtered $A_\iy$-category in the literature, due to Fukaya et al.\ \cite{Fuka2,Fuka3,Fuka4,Fuka5,FOOO}, Positselski \cite{Posi2} and De Deken and Lowen \cite{DDLo}, reviewed in \S\ref{ai51}. We will give definitions of filtered $A_\iy$-algebras and filtered $A_\iy$-categories which are new, but generalize all three of these.

\subsection{\texorpdfstring{Literature review on curved and filtered $A_\iy$-algebras and $A_\iy$-categories}{Literature review on curved and filtered A∞-algebras and A∞-categories}}
\label{ai51}

\begin{dfn}
\label{ai5def1}
Let $R$ be a commutative ring. A {\it curved\/ $A_\iy$-algebra}, or {\it weak\/ $A_\iy$-algebra}, $(A,\m)$ over $R$ consists of:
\begin{itemize}
\setlength{\itemsep}{0pt}
\setlength{\parsep}{0pt}
\item[(a)] A $\Z$-graded $R$-module $A=\bigop_{d\in\Z}A^d$; and
\item[(b)] Graded $R$-linear maps $\m_k:A^{\ot^k}\ra A$ of degree $2-k$ for $k\ge 0$. When $k=0$ we take $\m_0$ to be an $R$-linear map $\m_0:R\ra A^2$. Write $\m=(\m_k)_{k\ge 0}$. These must satisfy the equation in maps $A^{\ot^n}\ra A$ for all $n\ge 0$
\end{itemize}
\e
\sum_{i,j,k\ge 0,\; i+j+k=n}
(-1)^{i+jk}\m_{i+k+1}\ci\bigl(\id_A^{\ot^i}\ot\m_j\ot \id_A^{\ot^k}\bigr)=0.
\label{ai5eq1}
\e

This is the same as an $A_\iy$-algebra in Definition \ref{ai3def1}, except that we have maps $\m_k$ for $k\ge 0$ rather than $k\ge 1$, and in \eq{ai5eq1} we allow $j\ge 0$ rather than $j\ge 1$ in \eq{ai3eq3}. We say that $(A,\m)$ is {\it strict\/} if $\m_0=0$. Then $(A,\m)$ is an ordinary $A_\iy$-algebra. We define {\it units\/} $\boo\in A^0$ for curved $A_\iy$-algebras $(A,\m)$ as for $A_\iy$-algebras in \eq{ai3eq6}. We call a unital curved $A_\iy$-algebra $(A,\m)$ a {\it curved dg-algebra\/} or {\it weak dg-algebra\/} if $\m_k=0$ for $k\ge 3$. 
\end{dfn}

\begin{dfn}
\label{ai5def2}
Let $R$ be a commutative ring. A {\it curved\/ $A_\iy$-category}, or {\it weak\/ $A_\iy$-category}, $\sA$ over $R$ consists of:
\begin{itemize}
\setlength{\itemsep}{0pt}
\setlength{\parsep}{0pt}
\item[(a)] A set, or proper class, $\Obj\sA$ of {\it objects\/} of $\sA$.
\item[(b)] For all objects $x,y\in\Obj\sA$, a $\Z$-graded $R$-module
\begin{equation*}
\Hom(x,y)=\bigop\nolimits_{d\in\Z}\Hom^d(x,y).
\end{equation*}
We also write $\Hom_\sA(x,y)$ when we want to identify the $A_\iy$-category $\sA$.
\item[(c)] For all $k\ge 0$ and objects $x_0,x_1,\ldots,x_k$ in $\Obj\sA$, a graded $R$-linear map
\begin{align*}
\m_k^{x_0,\ldots,x_k}&:\Hom(x_0,x_1)\ot_R\Hom(x_1,x_2)\ot_R\cdots \ot_R\Hom(x_{k-1},x_k)\\
&\quad \longra\Hom(x_0,x_k)
\end{align*}
of degree $2-k$. For all $n\ge 0$ and objects $x_0,\ldots,x_n$ these must satisfy
\end{itemize}
\ea
\sum_{\begin{subarray}{l} i,j,k\ge 0:\\ i+j+k=n\end{subarray}}\begin{aligned}[t]&
(-1)^{i+jk}\m_{i+k+1}^{x_0,\ldots,x_i,x_{i+j},\ldots,x_n}\ci\bigl(\id_{\Hom(x_0,x_1)}\ot\cdots\ot\id_{\Hom(x_{i-1},x_i)}\ot {}\\
&\m_j^{x_i,\ldots,x_{i+j}}\ot \id_{\Hom(x_{i+j},x_{i+j+1})}\ot\cdots\ot\id_{\Hom(x_{n-1},x_n)}\bigr)=0.
\end{aligned}
\nonumber\\[-17pt]
\label{ai5eq2}
\ea

This is the same as an $A_\iy$-category in Definition \ref{ai4def1}, except that we have maps $\m_k^{x_0,\ldots,x_k}$ for $k\ge 0$ rather than $k\ge 1$, and in \eq{ai5eq2} we allow $j\ge 0$ rather than $j\ge 1$ in \eq{ai4eq2}. We define {\it unital\/} curved $A_\iy$-categories and {\it identities\/} $\boo_x$ as in Definition \ref{ai4def2}. We call a unital curved $A_\iy$-category $\sA$ a {\it curved dg-category\/} or {\it weak dg-category\/} if $\m_k^{x_0,\ldots,x_k}=0$ for all $k\ge 3$ and~$x_0,\ldots,x_k$.
\end{dfn}

\begin{rem}
\label{ai5rem1}
{\bf(a)} Curved $A_\iy$-algebras/dg-algebras/$A_\iy$-categories/dg-cat\-eg\-or\-ies occur naturally in many places in Mathematics and Physics. For example, curved $A_\iy$-algebras appear in the study of partial differential operators in Calderbank--Diemer \cite[\S 6]{CaDi}, of Poisson geometry in Cattaneo--Felder \cite[\S 3.1]{CaFe}, of differential forms on loop spaces in Getzler--Jones--Petrack \cite{GeJo,GJP}, of algebraic Chern--Weil theory in Lazarev--Schedler \cite{LaSc} and Nikolov--Zahariev \cite{NiZa}, of Koszul algebras in Floystad \cite{Floy} and Positselski \cite{Posi1}, and of derived categories in Nicol\'as \cite{Nico} and Positselski \cite{Posi2}, and in deformation theory in de Kleijn--Wierstra \cite{dKWi}. Curved $A_\iy$-categories appear in the Lagrangian Floer theory of Fukaya--Oh--Ohta--Ono~\cite{FOOO}. 

Curved dg-algebras and dg-categories appear in the study of noncommutative geometry in Block--Daenzer \cite{BlDa}, Kajiura \cite{Kaji}, Schwarz \cite{Schw1,Schw2}, and Tang \cite{Tang}, and of Landau--Ginzburg models in Caldararu--Tu \cite{CaTu} and Kapustin--Li~\cite{KaLi1,KaLi2}.

\smallskip

\noindent{\bf(b)} If $(A,\m)$ is an $A_\iy$-algebra then \eq{ai3eq4} for $k=1$ gives $\m_1\ci\m_1=0$, so $(A,\m_1)$ is a complex and we can form its cohomology $H^*(A)$ as in \eq{ai3eq5}.

However, if $(A,\m)$ is only a curved $A_\iy$-algebra then \eq{ai5eq1} for $k=1$ yields
\begin{equation*}
\m_1\ci\m_1=-\m_2\ci(\m_0(1)\ot\id_A)+\m_2\ci(\id_A\ot\m_0(1)).
\end{equation*}
So we may no longer have $\m_1\ci\m_1=0$, and we cannot form the cohomology $H^*(A)$. We regard $\m_0$ as the {\it obstruction} to $(A,\m_1)$ being a complex.
\smallskip

\noindent{\bf(c)} Part {\bf(b)} means that definitions involving $H^*(A,\m_1)$, such as cohomological units and quasi-isomorphisms, do not make sense for curved $A_\iy$-algebras and curved $A_\iy$-categories.
\smallskip

\noindent{\bf(d)} The obvious way to define {\it morphisms\/} of curved $A_\iy$-algebras is to follow Definition \ref{ai3def5}, but defining $\f_k$ for $k\ge 0$ rather than $k\ge 1$. However, we immediately run into a problem: equation \eq{ai3eq7} then involves an {\it infinite\/} sum, as on the right hand side we should sum over decompositions $n=n_1+\cdots+n_k$ with $k\ge 0$ and $n_i\ge 0$, rather than $n_i\ge 1$. The same applies to compositions in \eq{ai3eq10} and homotopies in \eq{ai3eq12}. Because of this, it is natural to study curved $A_\iy$-algebras and curved $A_\iy$-categories in settings in which we have notions of {\it convergence of infinite sums}. Parts {\bf(e)\rm--\bf(g)} describe three cases in the literature.

\smallskip

\noindent{\bf(e)} There is a massive project by Fukaya--Oh--Ohta--Ono over 2000--2025, including \cite{Fuka2,Fuka3,Fuka4,Fuka5,FOOO}, which constructs Lagrangian Floer cohomology and Fukaya categories for symplectic manifolds $(X,\om)$. They define notions of {\it filtered\/ $A_\iy$-algebra\/} \cite[\S 3.2]{FOOO} and {\it filtered\/ $A_\iy$-category\/} \cite[\S 1]{Fuka3}, \cite[\S 2]{Fuka5}, where \cite{FOOO} associates a filtered $A_\iy$-algebra to each Lagrangian $L$ in $(X,\om)$. These involve complete filtered modules over filtered {\it Novikov rings\/} such as $\La^\nov_\K$ and $\La^{0,\nov}_\K$ in Example \ref{ai5ex3} below, so convergence of infinite sums in such modules makes sense.
\smallskip

\noindent{\bf(f)} Positselski \cite{Posi2} studies curved $A_\iy$-algebras $(A,\m)$ which are complete modules over a pro-Artinian local ring $R$ with (maximal) ideal ${\mathfrak i}\subset R$, such that $\m_0(1)\in{\mathfrak i}\,A^2\subset A^2$. An infinite sum $\sum_{i\in I}a_i$ in $A^d$ is convergent if its projection to $A^d/{\mathfrak i}^NA^d$ has finitely many nonzero terms for all~$N\ge 1$.

Filtered $A_\iy$-algebras over $\La^{0,\nov}_\K$ in Fukaya et al.\ \cite{FOOO} fit into Positselski's framework, but filtered $A_\iy$-algebras over $\La^\nov_\K$ do not.
\smallskip

\noindent{\bf(g)} De Deken and Lowen \cite{DDLo} define a notion of filtered $A_\iy$-category (which they call {\it filtered $cA_\iy$-categories}, where the `$c$' stands for `curved'), intended for applications in noncommutative geometry. It is different to that in Fukaya \cite[\S 1]{Fuka3}, \cite[\S 2]{Fuka5}. They fix a totally ordered monoid $(\Ga,\le)$ and a complete $\Ga$-filtered commutative ring $R$, and their filtered $A_\iy$-categories involve complete $\Ga$-filtered $R$-modules. Convergence of infinite sums in such modules makes sense. 

They require that $0\in\Ga$ is the minimal element in the order $\Ga$. Because of this, filtered $A_\iy$-algebras over $\La^{0,\nov}_\K$ in Fukaya et al.\ \cite{FOOO} fit into De Deken--Lowen's framework with $(\Ga,\le)=([0,\iy),\le)$, but filtered $A_\iy$-algebras over $\La^\nov_\K$ with $(\Ga,\le)=(\R,\le)$ do not.
\smallskip

\noindent{\bf(h)} My notion of filtered $A_\iy$-algebra and $A_\iy$-category in \S\ref{ai53}--\S\ref{ai54} are new, but include {\bf(e)\rm--\bf(g)} as special cases. Although my primary motivation comes from Symplectic Geometry, I hope this work will be useful in other areas too, as in {\bf(a)}. As a sample application, in \S\ref{ai58} I will write the `twist construction' for $A_\iy$-categories from \S\ref{ai49} in terms of filtered $A_\iy$-categories.
\smallskip

\noindent{\bf(i)} Here are some references on filtered $A_\iy$-theory which are not very relevant to our project. Biran--Cornea--Shelukhin \cite[\S 2]{BCS} discuss `weakly filtered' $A_\iy$-categories and $A_\iy$-functors. Fukaya \cite[\S 4]{Fuka6} defines `filtered $A_\iy$-functors with energy loss'. In both \cite{BCS,Fuka6}, the $A_\iy$-categories and $A_\iy$-functors are ordinary (not curved) $A_\iy$-categories $\sA$ and $A_\iy$-functors $\fF$, with $\m_0^x=0$ and $\fF_0^x=0$. They are equipped with $\R$-filtrations, but the $A_\iy$-maps in $\sA$ and $\fF$ need not preserve filtrations, they may decrease filtrations by the `energy cost' or `energy loss', a measure of the `distance' between two $A_\iy$-categories.
\end{rem}

\subsection{Complete filtered rings and modules}
\label{ai52}

\begin{dfn}
\label{ai5def3}
Let $\Ga$ be a commutative monoid, with identity $0\in\Ga$ and monoid operation $+$, and $\pre$ be a partial order on $\Ga$. We will use the following notation: we write $\ga\pr\de$ to mean $\ga\pre\de$ and $\ga\ne\de$, and we write $\ga\sue\de$ and $\ga\su\de$ to mean $\de\pre\ga$ and $\de\pr\ga$ respectively. For $\ga\in\Ga$ and $n\in\N$ we write $n\ga={\buildrel {\ulcorner\,\,\,\text{$n$ copies }\,\,\,\urcorner}\over{\ga+\cdots+\ga}}$, with $0\ga=0$. If $\ga\in\Ga$ with $\ga\su 0$, we write $\ga\ssu 0$ if for all $\de\in\Ga$ there exists $n\gg 0$ in $\N$ with~$n\ga\sue\de$.

We call $(\Ga,\pre)$ a {\it partially-ordered monoid}, or {\it po-monoid}, if:
\begin{itemize}
\setlength{\itemsep}{0pt}
\setlength{\parsep}{0pt}
\item[(a)] If $\al\pre\be$ and $\ga\pre\de$ in $\Ga$ then $\al+\ga\pre\be+\de$.
\item[(b)] If $\ga,\de\in\Ga$ there is a unique {\it least upper bound\/} $\lub(\ga,\de)$ and a unique {\it greatest lower bound\/} $\glb(\ga,\de)$, such that $\ga\pre\lub(\ga,\de)$, $\de\pre\lub(\ga,\de)$, $\ga\sue\glb(\ga,\de)$, $\de\sue\glb(\ga,\de)$, and if $\ep\in\Ga$ with $\ga\pre\ep$, $\de\pre\ep$ then $\lub(\ga,\de)\pre\ep$, and if $\ze\in\Ga$ with $\ga\sue\ze$, $\de\sue\ze$ then~$\glb(\ga,\de)\sue\ze$.

By induction on the size of the set, least upper bounds and greatest lower bounds also exist for all nonempty finite subsets of~$\Ga$.
\item[(c)] There exists $\ga\in\Ga$ such that $\ga\ssu 0$.
\item[(d)] If $\ga\ssu 0$ and $\de\ssu 0$ then $\glb(\ga,\de)\ssu 0$.
\end{itemize}

If $(\ga_n)_{n=0}^\iy$ is a sequence in $\Ga$, we say that $\ga_n\ra\iy$ in $\Ga$ as $n\ra\iy$ if for all $\de\in\Ga$ there exists $N\ge 0$ such that $\ga_n\sue\de$ for all $n\ge N$. Then $\ga\ssu 0$ if and only if $n\ga\ra\iy$ as $n\ra\iy$.
\end{dfn}

\begin{ex}
\label{ai5ex1}
{\bf(a)} The trivial example is $\Ga=\{0,\iy\}$ with $0+0=0$, $0+\iy=\iy$, $\iy+\iy=\iy$ and $0\pr\iy$.
\smallskip

\noindent{\bf(b)} We can take $\Ga=\R,\Z,\Q,[0,\iy)$ or $\N$ with the usual total order $\le$. Then $\lub(\ga,\de)=\max(\ga,\de)$, and $\glb(\ga,\de)=\min(\ga,\de)$, and $\ga\ssu 0$ if and only if~$\ga>0$.

\smallskip

\noindent{\bf(c)} Let $(\Ga_1,\pre_1)$ and $(\Ga_2,\pre_2)$ be po-monoids. Define a partial order $\pre_{12}$ on the product monoid $\Ga_1\t\Ga_2$ by $(\ga_1,\ga_2)\pre_{12}(\de_1,\de_2)$ if $\ga_1\pre_1\de_1$ and $\ga_2\pre_2\de_2$. It is easy to show that $(\Ga_1\t\Ga_2,\pre_{12})$ is a po-monoid.
\end{ex}

\begin{dfn}
\label{ai5def4}
Let $(\Ga,\pre)$ be a po-monoid. A $\Ga$-{\it filtered ring\/} $R$ is a commutative ring $R$ equipped with a decreasing filtration $(F^\ga R)_{\ga\in\Ga}$ satisfying:
\begin{itemize}
\setlength{\itemsep}{0pt}
\setlength{\parsep}{0pt}
\item[(i)] $F^\ga R$ is an additive subgroup of $R$, with $F^\de R\subseteq F^\ga R$ if $\ga\pre \de$ in~$\Ga$.
\item[(ii)] $1\in F^0R$.
\item[(iii)] $F^\ga R\cdot F^\de R\subseteq F^{\ga+\de}R$ for all $\ga,\de\in\Ga$.
\item[(iv)] $\bigcup_{\ga\in\Ga}F^\ga R=R$ and $\bigcap_{\ga\in\Ga}F^\ga R=\{0\}$.
\end{itemize}

If $\ga\in\Ga$, we write $F^{\ssu\ga}R=\bigcup_{\de\in\Ga:\de\ssu 0}F^{\ga+\de}R$. Observe that $F^0R$ is a subring of $R$, and $F^{\ssu 0}R$ is an ideal in $F^0R$, so the quotient $F^0R/F^{\ssu 0}R$ is a commutative ring. For any $\ga\in\Ga$,  $F^\ga R/F^{\ssu \ga}R$ is a module over~$F^0R/F^{\ssu 0}R$.

A $\Ga$-{\it filtered module}, or just {\it filtered module}, $A$ over $R$ is an $R$-module $A$ equipped with a decreasing filtration $(F^\ga A)_{\ga\in\Ga}$ satisfying:
\begin{itemize}
\setlength{\itemsep}{0pt}
\setlength{\parsep}{0pt}
\item[(v)] $F^\ga A$ is an additive subgroup of $A$, with $F^\de A\subseteq F^\ga A$ if $\ga\pre \de$ in~$\Ga$.
\item[(vi)] $F^\ga R\cdot F^\de A\subseteq F^{\ga+\de}A$ for all $\ga,\de\in\Ga$.
\item[(vii)] $\bigcup_{\ga\in\Ga}F^\ga A=A$ and $\bigcap_{\ga\in\Ga}F^\ga A=\{0\}$.
\end{itemize}
If $\ga\in\Ga$, we write $F^{\ssu\ga}A=\bigcup_{\de\in\Ga:\de\ssu 0}F^{\ga+\de}A$. Then  $F^{\ssu\ga}A\subseteq F^\ga A$, so we can take the quotient $F^\ga A/F^{\ssu\ga}A$. It is a module over~$F^0R/F^{\ssu 0}R$.

Let $(a_n)_{n=0}^\iy$ be a sequence in $A$. We say that $(a_n)_{n=0}^\iy$ is {\it Cauchy\/} if for all $\ga\in\Ga$ there exists $N\ge 0$ such that $a_m-a_n\in F^\ga A$ for all $m,n\ge N$. We say that $(a_n)_{n=0}^\iy$ is {\it convergent}, with {\it limit\/} $a\in A$ if for all $\ga\in\Ga$ there exists $N\ge 0$ such that $a_n-a\in F^\ga A$ for all $n\ge N$. We write this as $a_n\ra a$ as $n\ra\iy$. The condition $\bigcap_{\ga\in\Ga}F^\ga A=\{0\}$ implies that limits are unique if they exist.

We can interpret convergence in topological terms, by making $A$ into a topological space, where the topology on $A$ is generated by open sets $\{a+F^\ga A:a\in A,$ $\ga\in\Ga\}$. The condition $\bigcap_{\ga\in\Ga}F^\ga A=\{0\}$ implies that $A$ is Hausdorff.

We say that the filtered $R$-module $A$ is {\it complete\/} if every Cauchy sequence $(a_n)_{n=0}^\iy$ is convergent. We call a filtered ring $R$ {\it complete\/} if it is complete as a module over itself. Every filtered $R$-module $A$ (or filtered ring $R$) has a natural completion $\bar A$ (or $\bar R$), where points of $\bar A$ are $\sim$-equivalence classes $[(a_n)_{n=0}^\iy]$ of Cauchy sequences $(a_n)_{n=0}^\iy$ in $A$, with equivalence relation $(a_n)_{n=0}^\iy\sim(a'_n)_{n=0}^\iy$ if $a_n-a_n'\ra 0$ as~$n\ra\iy$.

We can also define convergence and limits of infinite sums in $A$. A sum $\sum_{i\in I}\al_i$ with $\al_i\in A$ is called {\it Cauchy\/} if for all $\ga\in\Ga$ there are only finitely many $i\in I$ with $\al_i\notin F^\ga A$. A sum $\sum_{i\in I}\al_i$ in $A$ is called {\it convergent}, with {\it limit\/} $\al\in A$, if for all $\ga\in\Ga$ we have
\begin{equation*}
\al-\sum_{i\in I:\; \al_i\notin F^\ga A}\al_i\in F^\ga A.
\end{equation*}
Limits are unique if they exist. If $A$ is complete then every Cauchy sum is convergent.

Suppose $(\ga_n)_{n=0}^\iy$ is a sequence in $\Ga$ with $\ga_n\ra\iy$ as $n\ra\iy$, as in Definition \ref{ai5def3}, and $\sum_{n=0}^\iy\al_n$ is an infinite sum in $A$ with $\al_n\in F^{\ga_n}A$ for all $n\ge 0$. Then $\sum_{n=0}^\iy\al_n$ is Cauchy, and convergent if $A$ is complete.

If $A,B$ are filtered modules over $R$ then $A\ot_RB$ is a filtered module, such that $F^\ga(A\ot B)$ is the subgroup of elements of $A\ot_RB$ spanned by elements $r\,a\ot b$ for $r\in F^{\ga_1}R$, $a\in F^{\ga_2}A$ and $b\in F^{\ga_3}B$ with $\ga_1+\ga_2+\ga_3\sue\ga$. If $A,B$ are complete then $A\ot_RB$ need not be complete, but it has a unique completion, written $A\hat\ot_RB$.
\end{dfn}

\begin{rem}
\label{ai5rem2}
Definitions \ref{ai5def3} and \ref{ai5def4} are quite similar to parts of De Deken--Lowen \cite[\S 2.1--\S 2.5]{DDLo}. Our definition is more general, in that they consider filtrations over $(\Ga,\pre)$ in which $\pre$ is a total order, not a partial order, and require that $0\in\Ga$ is minimal in $\pre$. So \cite{DDLo} is a special case of our theory.
\end{rem}

\begin{ex}
\label{ai5ex2}
Take the po-monoid $\Ga=\{0,\iy\}$ as in Example \ref{ai5ex1}(a). Then any commutative ring $R$ has a unique $\Ga$-filtration satisfying Definition \ref{ai5def4}(i)--(iv), with $F^0R=R$ and $F^\iy R=\{0\}$. Similarly, every $R$-module $A$ has a unique $\Ga$-filtration satisfying Definition \ref{ai5def4}(v)--(vii), with $F^0A=A$ and $F^\iy A=\{0\}$. We have $F^{\ssu 0}R\!=\!F^{\ssu 0}A\!=\!\{0\}$. Then $R,A$ are complete. An infinite sum $\sum_{i\in I}a_i$ in $A$ is convergent if and only if there are only finitely many $i\in I$ with~$a_i\ne 0$.
\end{ex}

The next example is based on Fukaya et al.\ \cite[\S 1.7]{FOOO}.

\begin{ex}
\label{ai5ex3}
Let $\K$ be a commutative ring. The {\it universal Novikov ring\/ $\La^\nov_\K$ over\/} $\K$ is the formal power series ring
\begin{align*}
\La^\nov_\K=\bigl\{\ts\sum_{i=1}^N a_iq^{\la_i}:\,&\text{$N\in\{0,1,2,\ldots,\iy\}$, $a_i\in\K\sm\{0\}$, $\la_i\in\R$,}\\
&\text{$\la_1<\la_2<\cdots,$ $\la_i\ra\iy$ as $i\ra\iy$ if $N=\iy$}\bigr\}.
\end{align*}
It is a $\K$-algebra under addition and multiplication of formal power series, with unit $1q^0$. It has a decreasing $\R$-filtration $(F^\ga\La_\K)_{\ga\in\R}$ given by
\begin{align*}
F^\ga\La^\nov_\K=\bigl\{\ts\sum_{i=1}^N a_iq^{\la_i}:\,&\text{$N\in\{0,1,2,\ldots,\iy\}$, $a_i\in\K\sm\{0\}$, $\la_i\in\R$,}\\
&\text{$\ga\le\la_1<\la_2<\cdots,$ $\la_i\ra\iy$ as $i\ra\iy$ if $N=\iy$}\bigr\}.
\end{align*}
Then $\La^\nov_\K$ is a complete filtered ring over $(\R,\le)$. If $\K$ is a field then $\La^\nov_\K$ is a field. Also $\La^{0,\nov}_\K:=F^0\La^\nov_\K$ is a complete filtered ring over~$([0,\iy),\le)$.

Filtered $A_\iy$-algebras over both $\La^\nov_\K$ and $\La^{0,\nov}_\K$ are used in the Lagrangian Floer cohomology theory of Fukaya--Oh--Ohta--Ono \cite{Fuka2,Fuka3,Fuka4,Fuka5,FOOO}. The operations $\m_k$ in the filtered $A_\iy$-algebras are defined using sums of counts of $J$-holomorphic curves $\Si$ weighted by $q^{\area(\Si)}$ in $\La_\K$. Although there may be infinitely many possible homology classes for $\Si$, for any $\ga\in\R$ there are only finitely many homology classes for $\Si$ with $\area(\Si)<\ga$. So the sum over $\Si$ is convergent in a complete $\La^\nov_\K$-module. 
\end{ex}

\subsection{\texorpdfstring{Filtered $A_\iy$-algebras}{Filtered A∞-algebras}}
\label{ai53}

Next we define {\it filtered\/ $A_\iy$-algebras}. 

\begin{dfn}
\label{ai5def5}
Let $(\Ga,\pre)$ be a po-monoid and $R$ a complete $\Ga$-filtered ring, as in Definitions \ref{ai5def3} and \ref{ai5def4}. A $\Ga$-{\it filtered\/ $A_\iy$-algebra}, or just {\it filtered\/ $A_\iy$-algebra}, $(A,\m)$ over $R$ consists of:
\begin{itemize}
\setlength{\itemsep}{0pt}
\setlength{\parsep}{0pt}
\item[(a)] A $\Z$-graded $\Ga$-filtered $R$-module $A=\bigop_{d\in\Z}A^d$, where $A^d$ is a complete filtered $R$-module for each $d\in\Z$. (We do not require $A$ to be complete.) We say that $A$ is {\it degree-wise complete}, meaning that $A^d$ is complete for each degree~$d\in\Z$.
\item[(b)] Graded $R$-linear maps $\m_k:A^{\ot^k}\ra A$ of degree $2-k$ for $k\ge 0$. When $k=0$ we take $\m_0$ to be an $R$-linear map $\m_0:R\ra A^2$. Sometimes we identify $\m_0$ with $\m_0(1)\in A^2$. Write $\m=(\m_k)_{k\ge 0}$. 

We require $\m_k:A^{\ot^k}\ra A$ to be {\it filtration-preserving} for all $k\ge 1$, that is, $\m_k(F^\ga(A^{\ot^k}))\subseteq F^\ga A$ for all $\ga\in\Ga$. Equivalently, if $a_i\in F^{\ga_i}A$ for $i=1,\ldots,k$, $\ga_i\in\Ga$, then $\m_k(a_1\ot\cdots\ot a_k)\in F^{\ga_1+\cdots+\ga_k}A$.

We require $\m_0:R\ra A^2$ to satisfy $\m_0(1)\in F^{\ssu 0}A^2$. That is, $\m_0(1)\in F^\ga A^2$ for some $\ga\in\Ga$ with $\ga\ssu 0$, so $\ga\su 0$ and $n\ga\ra 0$ in $\Ga$ as~$n\ra\iy$.

As in \eq{ai5eq1}, these must satisfy the equation in maps $A^{\ot^n}\ra A$ for all $n\ge 0$
\end{itemize}
\e
\sum_{i,j,k\ge 0,\; i+j+k=n}
(-1)^{i+jk}\m_{i+k+1}\ci\bigl(\id_A^{\ot^i}\ot\m_j\ot \id_A^{\ot^k}\bigr)=0.
\label{ai5eq3}
\e

Note that a filtered $A_\iy$-algebra $(A,\m)$ is a curved $A_\iy$-algebra in the sense of Definition \ref{ai5def1}, with extra structure.

A filtered $A_\iy$-algebra $(A,\m)$ is called {\it strict\/} if $\m_0=0$. Then $(A,\m)$ is an ordinary $A_\iy$-algebra, forgetting the filtration structure.

A {\it unit\/} for $(A,\m)$ is an element $\boo\in F^0A^0$ satisfying \eq{ai3eq6}. Then we call $(A,\m)$ {\it unital}. Units are unique if they exist.

We will not define cohomologically unital filtered $A_\iy$-algebras.
\end{dfn}

\begin{ex}
\label{ai5ex4}
If $\Ga=\{0,\iy\}$ as in Example \ref{ai5ex1}(a), then as in Example \ref{ai5ex2}, $\Ga$-filtered rings $R$ and modules $A$ are equivalent to ordinary rings $R$ and modules $A$, all $R$-linear maps are filtration-preserving, and $F^{\ssu 0}A=0$. Thus in Definition \ref{ai5def5} we have $\m_0=0$, and $\Ga$-filtered $A_\iy$-algebras are equivalent to $A_\iy$-algebras.
\end{ex}

\begin{rem}
\label{ai5rem3}
{\bf(a)} When $\Ga=(\R,\le)$ and $R=\La^\nov_\K$ in Example \ref{ai5ex3}, or $\Ga=([0,\iy),\le)$ and $R=\La^{0,\nov}_\K$ in Example \ref{ai5ex3}, and $A$ is the completion of a free $\La^\nov_\K$- or $\La^{0,\nov}_\K$-module, then Definition \ref{ai5def5} reduces to filtered $A_\iy$-algebras as defined in Fukaya--Oh--Ohta--Ono \cite[\S 3.2]{FOOO}. Fukaya et al.\ associate filtered $A_\iy$-algebras of this kind to suitable Lagrangians $L$ in a symplectic manifold $(X,\om)$, as part of their theory of Lagrangian Floer cohomology.
\smallskip

\noindent{\bf(b)} Let $R$ be a pro-Artinian local ring with maximal ideal ${\mathfrak i}\subset R$. Make $R$ into a filtered ring over $(\N,\le)$ by setting $F^nR={\mathfrak i}^n$ for $n\in\N$. Then $\N$-filtered $A_\iy$-algebras $(A,\m)$ over $R$ in Definition \ref{ai5def5}, with filtration $F^nA={\mathfrak i}^n\cdot A$, are the curved $A_\iy$-algebras over $R$ studied by  Positselski~\cite{Posi2}.
\end{rem}

We can generalize much of \S\ref{ai33}--\S\ref{ai34} to the filtered case.

\begin{dfn}
\label{ai5def6}
Let $(\Ga,\pre)$ be a po-monoid, $R$ a complete $\Ga$-filtered ring, and $(A,\m),(B,\n)$ be filtered $A_\iy$-algebras over $R$. A {\it filtered\/ $A_\iy$-morphism\/} $\f:(A,\m)\ra(B,\n)$ is $\f=(\f_k)_{k\ge 0}$, where $\f_k:A^{\ot^k}\ra B$ for $k=0,1,2,\ldots$ are graded $R$-linear maps of degree $1-k$, satisfying:
\begin{itemize}
\setlength{\itemsep}{0pt}
\setlength{\parsep}{0pt}
\item[(i)] $\f_k$ is filtration-preserving for $k\ge 1$, that is, $\f_k(F^\ga(A^{\ot^k}))\subseteq F^\ga B$ for all~$\ga\in\Ga$.
\item[(ii)] $\f_0:R\ra B^1$ has $\f_0(1)\in F^{\ssu 0}B^1$. That is, $\f_0(1)\in F^\ga B^1$ for some $\ga\in\Ga$ with $\ga\ssu 0$, so $\ga\su 0$ and $n\ga\ra 0$ in $\Ga$ as~$n\ra\iy$.
\item[(iii)] Generalizing \eq{ai3eq7}, in maps $A^{\ot^n}\ra B$ for all $n\ge 0$ we have
\e
\begin{split}
&\sum_{i,j,k\ge 0:\; i+j+k=n}(-1)^{i+jk}\f_{i+k+1}\bigl(
\id_A^{\ot^i}\ot\m_j\ot\id_A^{\ot^k}\bigr) \\ 
&=\sum_{k\ge 0,\; n=n_1+\cdots+n_k,\; n_a\ge 0\!\!\!\!\!\!\!\!}
(-1)^{\sum_{i=1}^k(n_i+1)(k-i)}\n_k\bigl(\f_{n_1}\ot\f_{n_2}\ot\cdots\ot\f_{n_k}\bigr).
\end{split}
\label{ai5eq4}
\e

Note that the second line of \eq{ai5eq4} is an {\it infinite\/} sum, as for example for any $k\ge 1$ we can take $n_1=n$ and $n_a=0$ for $a=2,3,\ldots,k$. For a term in this sum with $k>n$, we see that $n_a=0$ for at least $k-n$ values of $a=1,\ldots,k$, so $\f_{n_1}\ot\cdots\ot\f_{n_k}$ contains at least $n-k$ factors $\f_0(1)$.
\end{itemize}

Here is how to interpret \eq{ai5eq4} in terms of convergent sums. Let $a_i\in F^{\ga_i}A^{d_i}$, $\ga_i\in\Ga$, $d_i\in\Z$ for $i=1,\ldots,n$. Evaluating \eq{ai5eq4} on $a_1\ot\cdots\ot a_n$ using Convention \ref{ai3conv1} gives an equation in $B^d$, for $d=2-n+\sum_{i=1}^nd_i$
\ea
&\sum_{\begin{subarray}{l} i,j,k\ge 0:\\ 
i+j+k=n \end{subarray}}
\begin{aligned}[t]&
(-1)^{i+jk+j\sum_{l=1}^id_i}\f_{i+k+1}\bigl(a_1\ot\cdots\ot a_i
\ot {}\\
&\m_j(a_{i+1}\ot\cdots\ot a_{i+j})\ot a_{i+j+1}\ot\cdots a_n\bigr) 
\end{aligned}
\label{ai5eq5}\\ 
&=\sum_{\begin{subarray}{l} k\ge 0,\; n= \\ n_1+\cdots+n_k,\\ n_a\ge 0 \end{subarray}}\begin{aligned}[t]&
(-1)^{\sum_{i=1}^k(n_i+1)(k-i+\sum_{a=1}^{n_1+\cdots+n_{i-1}}d_a)}\n_k\bigl(\f_{n_1}(a_1\ot\cdots\ot a_{n_1})\\
&{}\ot\f_{n_2}(a_{n_1+1}\ot\cdots\ot a_{n_1+n_2})\ot\cdots\\
&{}\ot\f_{n_k}(a_{n_1+\cdots+n_{k-1}+1}\ot\cdots\ot a_n)\bigr).
\end{aligned}
\nonumber
\ea

We will show that the second line of \eq{ai5eq5} is a convergent infinite sum in the complete $R$-module $B^d$ in the sense of Definition \ref{ai5def4}, and so is well defined. Consider the term in the second line of \eq{ai5eq5} from $n=n_1+\cdots+n_k$ with $k>n$. This includes at least $k-n$ terms $\f_0(1)\in F^\ga B^1$. Thus by $a_i\in F^{\ga_i}A^{d_i}$ and the filtration-preserving properties of $\n_k$ and $\f_{n_i}$ for $n_i\ge 1$, we see see that
\e
\begin{split}
&\pm\n_k\bigl(\f_{n_1}(a_1\ot\cdots\ot a_{n_1})\ot\cdots\ot\f_{n_k}(a_{n_1+\cdots+n_{k-1}+1}\ot\cdots\ot a_n)\bigr)\\
&\qquad \in F^{(k-n)\ga+\ga_1+\cdots+\ga_n}B^d.
\end{split}
\label{ai5eq6}
\e
Since $k\ga\ra\iy$ in $\Ga$ as $k\ra\iy$, we see that $(k-n)\ga+\ga_1+\cdots+\ga_n\ra\iy$ in $\Ga$ as $k\ra\iy$. Hence for any $\de\in\Ga$ there exists $K\ge 0$ such that $(k-n)\ga+\ga_1+\cdots+\ga_n\sue\de$ for all $k\ge K$. Thus \eq{ai5eq6} lies in $F^\de B^d$ if $k\ge K$. As there are only finitely many terms in the sum with $k<K$, only finitely many terms in the sum in the second line of \eq{ai5eq5} do not lie in $F^\de B^d$ for each $\de\in\Ga$, so the sum is convergent.

We have shown that \eq{ai5eq4} makes sense as a convergent sum when evaluated on $a_1\ot\cdots\ot a_n$ in $A^n$ for $a_i\in A^{d_i}$, $i=1,\ldots,n$. Extending by linearity, equation \eq{ai5eq4} makes sense as an equation in $R$-linear maps $A^{\ot^n}\ra B$. 

If $(A,\m),(B,\n)$ are strict, that is, $\m_0=\n_0=0$, we call a filtered $A_\iy$-morphism $\f:(A,\m)\ra(B,\n)$ {\it strict\/} if $\f_0=0$. Then $\f$ is an ordinary $A_\iy$-morphism of $A_\iy$-algebras.

The {\it identity morphism\/} $\id_A:(A,\m)\ra(A,\m)$ is $\id_A=(\id_{A,k})_{k\ge 0}$ with $\id_{A,1}=\id_A$ and $\id_{A,k}=0$ for $k\ge 1$.

If $(A,\m),(B,\n),(C,{\mathfrak o})$ are $\Ga$-filtered $A_\iy$-algebras over $R$ and $\f:(A,\m)\ra(B,\n)$, $\g:(B,\n)\ra(C,{\mathfrak o})$ are filtered $A_\iy$-morphisms, generalizing \eq{ai3eq10}, the {\it composition} $\g\ci\f:(A,\m) \ra(C,{\mathfrak o})$ is given by
\begin{equation*}
(\g\ci\f)_n=\sum_{n=n_1+\cdots+n_k,\; n_a\ge 0\!\!\!\!\!\!\!\!\!\!\!\!\!\!}(-1)^{\sum_{i=1}^k(n_i+1)(k-i)}\g_k\bigl(\f_{n_1}\ot\f_{n_2}\ot\cdots\ot\f_{n_k}).
\end{equation*}
This is an infinite sum, but it is convergent in the same sense as \eq{ai5eq4}, as if $k>n$ then $\g_k\bigl(\f_{n_1}\ot\cdots\ot\f_{n_k})$ contains at least $k-n$ factors of $\f_0(1)$. Composition is associative. Write $\AfAlg_R$ for the category with objects filtered $A_\iy$-algebras $(A,\m)$ over $R$ and morphisms filtered $A_\iy$-morphisms~$\f:(A,\m)\ra(B,\n)$.

If $(A,\m),(B,\n)$ are unital filtered $A_\iy$-algebras, as in Definition \ref{ai5def5}, we call a filtered $A_\iy$-morphism $\f:(A,\m)\ra(B,\n)$ {\it unital\/} if \eq{ai3eq11} holds. Unital filtered $A_\iy$-algebras over $R$ and unital filtered $A_\iy$-morphisms form a non-full subcategory~$\AfAlg_R^\un\subset\AfAlg_R$.

Define subcategories $\AfAlg_{R,\iy}\subset\AfAlg_{R,1}\subset\AfAlg_{R,0}=\AfAlg_R$, such that $\AfAlg_{R,1}$ has objects {\it strict\/} filtered $A_\iy$-algebras $(A,\m)$ with $\m_0=0$ and morphisms {\it all\/} filtered $A_\iy$-morphisms $\f:(A,\m)\ra(B,\n)$, including those with $\f_0\ne 0$, and $\AfAlg_{R,\iy}$ has objects {\it strict\/} filtered $A_\iy$-algebras $(A,\m)$, with $\m_0=0$, and morphisms {\it strict\/} filtered $A_\iy$-morphisms $\f$, with~$\f_0=0$.

Objects of $\AfAlg_{R,1},\AfAlg_{R,\iy}$ are ordinary $A_\iy$-algebras (with extra filtration data), and morphisms of $\AfAlg_{R,\iy}$ are $A_\iy$-morphisms (which preserve filtrations), but morphisms of $\AfAlg_{R,1}$ are not $A_\iy$-morphisms. Thus there is a forgetful functor $\AfAlg_{R,\iy}\ra\AAlg_R$ which forgets filtrations. This gives a diagram of categories and functors
\e
\begin{gathered}
\xymatrix@C=35pt@R=15pt{
\AfAlg_{R,\iy}\,\,  \ar@{^{(}->}[r] \ar[d] &  \AfAlg_{R,1}\,\, \ar@{^{(}->}[r]  &  {\AfAlg_{R,0}=\AfAlg_R}  \\
{\AAlg_R.\,\,}  }
\end{gathered}
\label{ai5eq7}
\e
We define $\AfAlg^\un_{R,\iy}\subset\AfAlg^\un_{R,1}\subset\AfAlg^\un_{R,0}=\AfAlg^\un_R$ in the same way. For the categories $\AfAlg_{R,K},\AfAlg^\un_{R,K}$ for $K=0,1,\iy$ we call $K$ the {\it type}. It is the analogue of the type $K$ of $\iy$-categories $\AAlg^{\un,\iy}_{R,K},\ldots$ in Chapters~\ref{ai6}--\ref{ai8}.
\end{dfn}

The next definition follows Fukaya et al.~\cite[Def.~4.2.35 \& Prop.~4.2.40]{FOOO}.

\begin{dfn}
\label{ai5def7}
Let $(\Ga,\pre)$ be a po-monoid, $R$ a complete $\Ga$-filtered ring, $(A,\m),(B,\n)$ be filtered $A_\iy$-algebras over $R$, and $\f,\g:(A,\m)\ra(B,\n)$ be filtered $A_\iy$-morphisms. A {\it filtered homotopy\/} $\H:\f\Ra\g$ from $\f$ to $\g$ is $\H=(\H_k)_{k\ge 0}$, where $\H_k:A^{\ot^k}\ra B$ for $k=0,1,2,\ldots$ are graded $R$-linear maps of degree $-k$, satisfying:
\begin{itemize}
\setlength{\itemsep}{0pt}
\setlength{\parsep}{0pt}
\item[(i)] $\H_k$ is filtration-preserving for $k\ge 1$, that is, $\H_k(F^\ga(A^{\ot^k}))\subseteq F^\ga B$ for all~$\ga\in\Ga$.
\item[(ii)] $\H_0:R\ra B^0$ has $\H_0(1)\in F^{\ssu 0}B^0$. That is, $\H_0(1)\in F^\ga B^0$ for some $\ga\in\Ga$ with $\ga\ssu 0$, so $\ga\su 0$ and $n\ga\ra 0$ in $\Ga$ as~$n\ra\iy$.
\item[(iii)] Generalizing \eq{ai3eq12}, in maps $A^{\ot^n}\ra B$ for all $n\ge 0$ we have
\begin{align*}
&\f_n-\g_n=
\\
&\sum_{\begin{subarray}{l} n=n_1+\cdots+n_k, \\ 
n_i\ge 0,\; 1\le l\le k \end{subarray}}
\begin{aligned}[t]&
(-1)^{k+1}\cdot(-1)^{\sum_{i=1}^{l-1}(n_i+1)}\cdot(-1)^{\sum_{i=1}^k(n_i+1)(k-i)}\\
&\n_k\bigl(\f_{n_1}\ot\f_{n_2}\ot\cdots\ot\f_{n_{l-1}}\ot
\H_{n_l}\ot\g_{n_{l+1}}\ot\cdots\ot\g_{n_k}\bigr)
\end{aligned}
\\
&+\sum_{i,j,k\ge 0,\; i+j+k=n}
(-1)^{i+jk}\H_{i+k+1}\ci\bigl(\id_A^{\ot^i}\ot\m_j\ot \id_A^{\ot^k}\bigr). 
\end{align*}
The second line is an infinite sum, which converges as for~\eq{ai5eq4}.
\end{itemize}

If $(A,\m),(B,\n),\f,\g$ are unital, we call $\H$ {\it unital\/} if \eq{ai3eq13} holds.

Let $(A,\m),(B,\n)$ be filtered $A_\iy$-algebras over $R$, and $\f:(A,\m)\ra(B,\n)$ a filtered $A_\iy$-morphism. A {\it filtered homotopy inverse} for $\f$ is a filtered $A_\iy$-morphism $\g:(B,\n)\ra(A,\m)$ such that $\g\ci\f:(A,\m)\ra(A,\m)$ is filtered homotopic to $\id_A:(A,\m)\ra(A,\m)$, and $\f\ci\g:(B,\n)\ra(B,\n)$ is filtered homotopic to $\id_B:(B,\n)\ra(B,\n)$. If $\f$ has a filtered homotopy inverse, we call $\f$ a {\it filtered homotopy equivalence}, and we call $(A,\m),(B,\n)$ {\it filtered homotopy equivalent}.
\end{dfn}

Here is the analogue of Theorem \ref{ai3thm1}(a),(b), whose proof generalizes to the filtered case. (See Seidel \cite[\S 1(h)]{Seid} for a proof of (a) which extends easily.)

\begin{thm}
\label{ai5thm1}
Let\/ $(A,\m),(B,\n)$ be filtered\/ $A_\iy$-algebras over $R$. Then
\begin{itemize}
\setlength{\itemsep}{0pt}
\setlength{\parsep}{0pt}
\item[{\rm(a)}] Filtered homotopy is an equivalence relation on filtered\/ $A_\iy$-morphisms $\f,\g:(A,\m)\ra(B,\n)$.
\item[{\rm(b)}] Filtered homotopy equivalence is an equivalence relation on filtered\/ $A_\iy$-algebras.
\end{itemize}
\end{thm}

\begin{rem}
\label{ai5rem4}
Let $(A,\m)$ be a $\Ga$-filtered $A_\iy$-algebra. Then $F^0A/F^{\ssu 0}A$ is a module over the commutative ring $F^0R/F^{\ssu 0}R$, and $\m_k$ for $k\ge 0$ maps $(F^0A)^{\ot^k}\ra F^0A$, and this descends to an $F^0R/F^{\ssu 0}R$-linear map
\begin{equation*}
\dot\m_k:(F^0A/F^{\ssu 0}A)^{\ot^k}\longra F^0A/F^{\ssu 0}A,
\end{equation*}
with $\dot\m_0\!=0$. Equation \eq{ai5eq3} holds for the $\dot\m_k$, so as $\dot\m_0=0$ we see that $(F^0A/F^{\ssu 0}A,\dot\m=(\dot\m_k)_{k\ge 1})$ is an ordinary $A_\iy$-algebra over~$F^0R/F^{\ssu 0}R$.

Similarly, if we choose $\ga\su 0$ with $\m_0(1)\in F^\ga A^2$ then $F^0A/F^\ga A$ becomes an ordinary $A_\iy$-algebra $(F^0A/F^\ga A,\ddot\m=(\ddot\m_k)_{k\ge 1})$ over~$F^0R/F^\ga R$.

Filtered $A_\iy$-morphisms $\f$ induce $A_\iy$-morphisms $\dot\f$ for the $(F^0A/F^{\ssu 0}A,\dot\m)$, and (if we can choose $\ga$ `small enough' in $\Ga$) $\ddot\f$ for the $(F^0A/F^\ga A,\ddot\m)$.

We could use this to make definitions for $\Ga$-filtered $A_\iy$-algebras $(A,\m)$, and then try and prove theorems. For example, we could define $(A,\m)$ to be cohomologically unital if $(F^0A/F^{\ssu 0}A,\dot\m)$ or $(F^0A/F^\ga A,\ddot\m)$ is, in the sense of Definition \ref{ai3def4}, or we could define $\f$ to be a weak homotopy equivalence if $\dot f$ or $\ddot f$ is, in the sense of Definition \ref{ai3def5}. This approach is taken by Fukaya et al.\ \cite{Fuka2,Fuka3,Fuka4,Fuka5,Fuka6,FOOO} with their notion of filtered $A_\iy$-algebra. For example, in \cite[Th.~4.2.45]{FOOO} they prove an analogue of of Theorem \ref{ai3thm1}(c) for `gapped filtered $A_\iy$-algebras'.

We will not do this, as it is not helpful for the main results that we want to prove, but it may be useful in further developments of the theory.
\end{rem}

\subsection{\texorpdfstring{Filtered $A_\iy$-categories}{Filtered A∞-categories}}
\label{ai54}

Next we extend much of \S\ref{ai41}--\S\ref{ai42} to the filtered case. The next definition combines Definitions \ref{ai4def1} and \ref{ai5def5}. It generalizes definitions in Fukaya \cite[\S 1]{Fuka3}, \cite[\S 2]{Fuka5} and De Deken--Lowen~\cite[\S 2]{DDLo}.

\begin{dfn}
\label{ai5def8}
Let $(\Ga,\pre)$ be a po-monoid and $R$ a complete $\Ga$-filtered ring, as in Definitions \ref{ai5def3} and \ref{ai5def4}. A $\Ga$-{\it filtered\/ $A_\iy$-category}, or just {\it filtered\/ $A_\iy$-category}, $\sA$ over $R$ consists of:
\begin{itemize}
\setlength{\itemsep}{0pt}
\setlength{\parsep}{0pt}
\item[(a)] A set, or proper class, $\Obj\sA$ of {\it objects\/} of $\sA$.
\item[(b)] For all objects $x,y\in\Obj\sA$, a $\Z$-graded $\Ga$-filtered $R$-module
\begin{equation*}
\Hom(x,y)=\bigop\nolimits_{d\in\Z}\Hom^d(x,y),
\end{equation*}
where $\Hom^d(x,y)$ is a complete filtered $R$-module for each $d\in\Z$. (We do not require $\Hom(x,y)$ to be complete.) We say that $\Hom(x,y)$ is {\it degree-wise complete}. We also write $\Hom_\sA(x,y)$ when we want to identify the filtered $A_\iy$-category $\sA$.
\item[(c)] For all $k\ge 0$ and objects $x_0,x_1,\ldots,x_k$ in $\Obj\sA$, a graded $R$-linear map
\begin{align*}
&\m_k^{x_0,\ldots,x_k}:\Hom(x_0,x_1)\ot_R\Hom(x_1,x_2)\ot_R\cdots \ot_R\Hom(x_{k-1},x_k)\\
&\qquad\qquad \longra\Hom(x_0,x_k)
\end{align*}
of degree $2-k$. When $k=0$, this is an $R$-linear map $\m_0^{x_0}:R\ra\Hom^2(x_0,x_0)$. We write $\m_k=(\m_k^{x_0,\ldots,x_k})_{x_0,\ldots,x_k\in\Obj\sA}$ and $\m=(\m_k)_{k\ge 0}$.

We require $\m_k^{x_0,\ldots,x_k}$ to be {\it filtration-preserving} for all $k\ge 1$ and $x_0,\ldots,x_k$ in $\Obj\sA$, that is, $\m_k^{x_0,\ldots,x_k}(F^\ga(\Hom(x_0,x_1)\ot\cdots\ot\Hom(x_{k-1},x_k)))\subseteq F^\ga \Hom(x_0,x_k)$ for all $\ga\in\Ga$.

We require that $\m_0^{x_0}(1)\in F^{\ssu 0}\Hom^2(x_0,x_0)$ for all $x_0$ in $\Obj\sA$. That is, $\m_0^{x_0}(1)\in F^\ga\Hom^2(x_0,x_0)$ for some $\ga\in\Ga$ with $\ga\ssu 0$, so $\ga\su 0$ and $n\ga\ra 0$ in $\Ga$ as $n\ra\iy$. This $\ga$ may depend on~$x_0$.

Generalizing \eq{ai4eq1}, we require that for all $n\ge 0$ and objects $x_0,\ldots,x_n$
\end{itemize}
\ea
\sum_{\begin{subarray}{l} i,j,k\ge 0:\\ i+j+k=n\end{subarray}}\begin{aligned}[t]&
(-1)^{i+jk}\m_{i+k+1}^{x_0,\ldots,x_i,x_{i+j},\ldots,x_n}\ci\bigl(\id_{\Hom(x_0,x_1)}\ot\cdots\ot\id_{\Hom(x_{i-1},x_i)}\ot {}\\
&\m_j^{x_i,\ldots,x_{i+j}}\ot \id_{\Hom(x_{i+j},x_{i+j+1})}\ot\cdots\ot\id_{\Hom(x_{n-1},x_n)}\bigr)=0.
\end{aligned}
\nonumber\\[-17pt]
\label{ai5eq8}
\ea

We say that $\sA$ is {\it strict\/} if $\m_0^{x_0}=0$ for all $x_0\in\Obj\sA$. Then $\sA$ is an ordinary $A_\iy$-category, forgetting the filtration structure.

We call $\sA$ {\it unital\/} if for all $x\in\Obj\sA$ there exists $\boo_x\in F^0\Hom^0(x,x)$, called the {\it identity\/} for $x$, satisfying \eq{ai4eq5}--\eq{ai4eq6}. Identities $\boo_x$ are unique if they exist.
\end{dfn}

The next definition combines Definitions \ref{ai4def6} and~\ref{ai5def6}.

\begin{dfn}
\label{ai5def9}
Let $(\Ga,\pre)$ be a po-monoid, $R$ a complete $\Ga$-filtered ring, and $\sA,\sB$ be filtered $A_\iy$-categories over $R$. A {\it filtered\/ $A_\iy$-functor\/} $\fF:\sA\ra\sB$ consists of:
\begin{itemize}
\setlength{\itemsep}{0pt}
\setlength{\parsep}{0pt}
\item[(a)] A map $\fF_\Obj:\Obj\sA\ra\Obj\sB$. When this will not cause confusion, for $x\in\Obj\sA$ we write $\fF(x)\in\Obj\sB$ rather than $\fF_\Obj(x)$.
\item[(b)] For all $k\ge 0$ and objects $x_0,x_1,\ldots,x_k$ in $\Obj\sA$, a graded $R$-linear map
\begin{align*}
\fF_k^{x_0,\ldots,x_k}:\Hom_\sA(x_0,x_1)\ot_R\Hom_\sA(x_1,x_2)\ot_R\cdots \ot_R\Hom_\sA(x_{k-1},x_k)&\\
\longra\Hom_\sB(\fF_\Obj(x_0),\fF_\Obj(x_k))&,
\end{align*}
of degree $1-k$. When $k=0$, this is an $R$-linear map $\fF_0^{x_0}:R\ra\Hom_\sB^1(\fF(x_0),\fF(x_0))$. 

We require $\fF_k^{x_0,\ldots,x_k}$ to be {\it filtration-preserving} for all $k\ge 1$ and $x_0,\ldots,x_k$ in $\Obj\sA$, that is, $\fF_k^{x_0,\ldots,x_k}(F^\ga(\Hom_\sA(x_0,x_1)\ot\cdots\ot\Hom_\sA(x_{k-1},x_k)))\subseteq F^\ga \Hom_\sB(\fF(x_0),\fF(x_k))$ for all $\ga\in\Ga$.

We require that $\fF_0^{x_0}(1)\in F^{\ssu 0}\Hom^1_\sB(\fF(x_0),\fF(x_0))$ for all $x_0$ in $\Obj\sA$.

Generalizing \eq{ai4eq8} and \eq{ai5eq4}, the $\fF_k^{x_0,\ldots,x_k}$ must satisfy for all $n\ge 0$ and $x_0,\ldots,x_n$ in $\Obj\sA$
\end{itemize}
\ea
&\sum_{\begin{subarray}{l} i,j,k\ge 0:\\ i+j+k=n \end{subarray}}
\begin{aligned}[t]&
(-1)^{i+jk}\fF_{i+k+1}^{x_0,\ldots,x_i,x_{i+j},\ldots,x_n}\ci\bigl(\id_{\Hom_\sA(x_0,x_1)}\ot\cdots\ot\id_{\Hom_\sA(x_{i-1},x_i)}\\
&{}\ot\m_j^{\sA;x_i,\ldots,x_{i+j}}\ot \id_{\Hom_\sA(x_{i+j},x_{i+j+1})}\ot\cdots\ot\id_{\Hom_\sA(x_{n-1},x_n)}\bigr)\end{aligned} 
\nonumber\\ 
&=\sum_{\begin{subarray}{l}  n=n_1+\cdots+n_k,\\ n_a\ge 0\end{subarray}}
\begin{aligned}[t]&(-1)^{\sum_{i=1}^k(n_i+1)(k-i)}\m^{\sB;\fF_\Obj(x_0),\fF_\Obj(x_{n_1}),\ldots,\fF_\Obj(x_n),}_k\bigl(\fF_{n_1}^{x_0,\ldots,x_{n_1}}\\
&{}\ot\fF_{n_2}^{x_{n_1},\ldots,x_{n_1+n_2}}\ot\cdots\ot\fF_{n_k}^{x_{n-n_k},\ldots,x_n}).\end{aligned} 
\nonumber\\[-16pt] 
\label{ai5eq9}
\ea
\begin{itemize}
\setlength{\itemsep}{0pt}
\setlength{\parsep}{0pt}
\item[] The second line of \eq{ai5eq9} is an {\it infinite\/} sum, which converges in a similar way to \eq{ai5eq4}, since $\fF_0^{x_j}(1)\in F^{\ssu 0}\Hom^1_\sB(\fF(x_j),\fF(x_j))$.
\end{itemize}

If $\sA,\sB$ are strict, we call a filtered $A_\iy$-functor $\fF:\sA\ra\sB$ {\it strict\/} if $\fF^{x_0}_0=0$ for all $x_0\in\Obj\sA$. Then $\fF$ is an ordinary $A_\iy$-functor of $A_\iy$-categories.

The {\it identity functor\/} $\id^\sA:\sA\ra\sA$ has $\id^\sA_\Obj=\id_{\Obj\sA}$, and $\id^{\sA,x}_1=\id:\Hom_\sA(x,x)\ra\Hom_\sA(x,x)$, and $\id^\sA_k=0$ for $k\ne 1$.

If $\fF:\sA\ra\sB$, $\fG:\sB\ra\sC$ are filtered $A_\iy$-functors, generalizing \eq{ai4eq10} the {\it composition} $\fG\ci\fF:\sA\ra\sC$ is given by $(\fG\ci\fF)_\Obj=\fG_\Obj\ci\fF_\Obj$ and
\ea
&(\fG\ci\fF)_n^{x_0,\ldots,x_n}=
\label{ai5eq10}\\
&\sum_{\begin{subarray}{l} n=n_1+\cdots+n_k,\\ n_a\ge 0\end{subarray}}
\begin{aligned}[t]& (-1)^{\sum_{i=1}^k(n_i+1)(k-i)}\fG_k^{\fF(x_0),\fF(x_{n_1}),\fF(x_{n_1+n_2}),\ldots,\fF(x_n)} \\ 
&\bigl(\fF_{n_1}^{x_0,\ldots,x_{n_1}}\ot\fF_{n_2}^{x_{n_1},\ldots,x_{n_1+n_2}}\ot\cdots\ot\fF_{n_k}^{x_{n-n_k},\ldots,x_n}),
\end{aligned}
\nonumber
\ea
for all $n\ge 0$ and $x_0,\ldots,x_n\in\Obj\sA$. Again, the second line of \eq{ai5eq10} is an {\it infinite\/} sum, which converges as for \eq{ai5eq9}. Composition is associative.

Write $\AfCat_R$ for the category with objects filtered $A_\iy$-categories $\sA$ over $R$ and morphisms filtered $A_\iy$-functors $\fF:\sA\ra\sB$.

If $\sA,\sB$ are unital filtered $A_\iy$-categories over $R$, as in Definition \ref{ai4def2}, we call a filtered $A_\iy$-functor $\fF:\sA\ra\sB$ {\it unital\/} if \eq{ai4eq11}--\eq{ai4eq12} hold. Unital filtered $A_\iy$-categories over $R$ and unital filtered $A_\iy$-functors form a non-full subcategory $\AfCat_R^\un\subset\AfCat_R$.

Define subcategories $\AfCat_{R,\iy}\subset\AfCat_{R,1}\subset\AfCat_{R,0}=\AfCat_R$, such that $\AfCat_{R,1}$ has objects {\it strict\/} filtered $A_\iy$-categories $\sA$ over $R$, with $\m_0^x=0$ for all $x\in\Obj\sA$ and morphisms {\it all\/} filtered $A_\iy$-functors $\fF:\sA\ra\sB$, including those with $\fF^x_0\ne 0$, and $\AfCat_{R,\iy}$ has objects {\it strict\/} filtered $A_\iy$-categories $\sA$, with $\m^x_0=0$ for all $x\in\Obj\sA$, and morphisms {\it strict\/} filtered $A_\iy$-functors $\fF:\sA\ra\sB$, with $\fF^x_0=0$ for all~$x\in\Obj\sA$.

Objects of $\AfCat_{R,1},\AfCat_{R,\iy}$ are ordinary $A_\iy$-categories (with extra filtration data), and morphisms of $\AfCat_{R,\iy}$ are $A_\iy$-functors (which preserve filtrations), but morphisms of $\AfCat_{R,1}$ are not $A_\iy$-functors. Thus there is a forgetful functor $\Pi_{A_\iy^\rfi}^{A_\iy}:\AfCat_{R,\iy}\ra\ACat_R$ which forgets filtrations. As in \eq{ai5eq7}, this gives a diagram of categories and functors
\e
\begin{gathered}
\xymatrix@C=35pt@R=15pt{
\AfCat_{R,\iy}\,\,  \ar@{^{(}->}[r] \ar[d]_{\Pi_{A_\iy^\rfi}^{A_\iy}} &  \AfCat_{R,1}\,\, \ar@{^{(}->}[r]  &  {\AfCat_{R,0}=\AfCat_R}  \\
{\ACat_R.\,\,}  }
\end{gathered}
\label{ai5eq11}
\e
We define $\AfCat^\un_{R,\iy}\subset\AfCat^\un_{R,1}\subset\AfCat^\un_{R,0}=\AfCat^\un_R$ in the same way. For the categories $\AfCat_{R,K},\AfCat^\un_{R,K}$ for $K=0,1,\iy$ we call $K$ the {\it type}. (The choice of labels $K=0,1,\iy$ correspond to the subscript $K$ in the $\iy$-categories $\AfCat^\iy_{R,K,L},\AfCat^{\un,\iy}_{R,K,L}$ in Chapter \ref{ai8}, see equation~\eq{ai8eq1}.)
\end{dfn}

\subsection{\texorpdfstring{The $A_\iy$-category of filtered $A_\iy$-functors}{The A∞-category of filtered A∞-functors}}
\label{ai55}

When $\sA,\sB$ are $A_\iy$-categories over $R$, section \ref{ai43} defined the $A_\iy$-category $\Fun_{A_\iy}(\sA,\sB)$ of $A_\iy$-functors $\fF:\sA\ra\sB$. There is an analogue of this for filtered $A_\iy$-categories, which is discussed by De Deken--Lowen \cite[\S 3.4--\S 3.5]{DDLo} and Fukaya \cite[Def.~2.18, Th.~2.19 \& Prop.~2.20]{Fuka5} for their notions of filtered $A_\iy$-categories. 

\begin{dfn}
\label{ai5def10}
{\bf(a)} Let $(\Ga,\pre)$ be a po-monoid, $R$ a complete $\Ga$-filtered ring, and $\sA,\sB$ be filtered $A_\iy$-categories over $R$. We will define a strict filtered $A_\iy$-category $\Fun_{A^\rfi_\iy}(\sA,\sB)$ over $R$. Objects of $\Fun_{A^\rfi_\iy}(\sA,\sB)$ are filtered $A_\iy$-functors $\fF:\sA\ra\sB$, as in \S\ref{ai54}. For objects $\fF,\fG$, we define a $\Ga$-filtered graded $R$-module
\begin{equation*}
\Hom_{\Fun_{A^\rfi_\iy}(\sA,\sB)}(\fF,\fG)=\Hom_{A^\rfi_\iy}(\fF,\fG)=\bigop_{l\in\Z}\Hom^l_{A^\rfi_\iy}(\fF,\fG)
\end{equation*}
as follows. As in \eq{ai4eq13}, an element $\fT$ of $\Hom^l_{A^\rfi_\iy}(\fF,\fG)$ consists of, for all $k\ge 0$ and objects $x_0,\ldots,x_k$ in $\sA$, a graded $R$-linear map of degree $l-k$
\begin{align*}
\fT_k=\fT^{x_0,\ldots,x_k}_k&:\Hom_\sA(x_0,x_1)\ot_R\Hom_\sA(x_1,x_2)\ot_R\cdots \ot_R\Hom_\sA(x_{k-1},x_k)\\
&\qquad\longra\Hom_\sB(\fF_\Obj(x_0),\fG_\Obj(x_k)),
\end{align*}
such that there exists $\de\in\Ga$ such that for all $k\ge 0,$ $x_0,\ldots,x_k\in\Obj\sA$ and $\ga\in\Ga$ we have
\e
\begin{split}
&\fT^{x_0,\ldots,x_k}_k\bigl(F^\ga(\Hom_\sA(x_0,x_1)\ot\cdots \ot\Hom_\sA(x_{k-1},x_k))\bigr) \\
&\qquad \subseteq F^{\ga+\de}\Hom_\sB(\fF(x_0),\fG(x_k)).
\end{split}
\label{ai5eq12}
\e

We define the $\Ga$-filtration on $\Hom^l_{A^\rfi_\iy}(\fF,\fG)$ by $\fT\in F^\de\Hom^l_{A^\rfi_\iy}(\fF,\fG)$ for $\de\in\Ga$ if \eq{ai5eq12} holds for all $k\ge 0,$ $x_0,\ldots,x_k\in\Obj\sA$ and $\ga\in\Ga$.

The $A_\iy$-operations $\n_p^{\fF^0,\ldots,\fF^p}$ for $p\ge 0$ are defined by $\n_0^{\fF_0}=0$, and $\n_p^{\fF^0,\ldots,\fF^p}$ for $p\ge 1$ are as in \eq{ai4eq14}--\eq{ai4eq17}, except that we include terms $\m_0^{\sA;x_0},\m_0^{\sB;y_0},\ab\fF_0^{i;x_0}$, so in \eq{ai4eq15} and \eq{ai4eq17} we omit the conditions $l_i>l_{i-1}$, $l_j>l_{j-1}$. Then the second line of \eq{ai4eq15}, and \eq{ai4eq17}, become infinite sums, which converge in the usual way. As in De Deken--Lowen \cite[\S 3.4.2]{DDLo}, this makes $\Fun_{A^\rfi_\iy}(\sA,\sB)$ into a strict filtered $A_\iy$-category over $R$, and thus an $A_\iy$-category over $R$.
\smallskip

\noindent{\bf(b)} Let $\sA_1,\sA_2,\sB$ be filtered $A_\iy$-categories over $R$, and $\fE:\sA_1\ra\sA_2$ be a filtered $A_\iy$-functor. Generalizing Definition \ref{ai4def9}(a), we define a strict filtered $A_\iy$-functor $\fE^*:\Fun_{A^\rfi_\iy}(\sA_2,\sB)\ra\Fun_{A^\rfi_\iy}(\sA_1,\sB)$ to act on objects by $\fE^*_\Obj:\fF\mapsto\fF\ci\fE$ for $\fF:\sA_2\ra\sB$ a filtered $A_\iy$-functor, and if $p\ge 0$ and $\fF_0,\ldots,\fF_p$ are objects in $\Fun_{A_\iy}(\sA_2,\sB)$, we define
\begin{align*}
(\fE^*)_p^{\fF_0,\ldots,\fF_p}&:\Hom_{A^\rfi_\iy}(\fF_0,\fF_1)\ot\cdots\ot\Hom_{A^\rfi_\iy}(\fF_{p-1},\fF_p)\\
&\qquad \longra \Hom_{A^\rfi_\iy}(\fF_0\ci\fE,\fF_p\ci\fE)
\end{align*}
by $(\fE^*)_p^{\fF_0,\ldots,\fF_p}=0$ if $p\ne 1$, and $(\fE^*)_1^{\fF_0,\fF_1}$ is as in \eq{ai4eq20}, except that we include terms $\fE_0^{x_j}$, so we sum over $n_a\ge 0$, and do not require $k\le n$. Then \eq{ai4eq20} becomes an infinite sum, which converges in the usual way. As in De Deken--Lowen \cite[Prop.~3.49]{DDLo}, this defines a strict filtered $A_\iy$-functor.
\smallskip

\noindent{\bf(c)} Let $\sA,\sB_1,\sB_2$ be filtered $A_\iy$-categories over $R$, and $\fG:\sB_1\ra\sB_2$ be a filtered $A_\iy$-functor. Generalizing Definition \ref{ai4def9}(b), we define a strict filtered $A_\iy$-functor $\fG_*:\Fun_{A_\iy}(\sA,\sB_1)\ra\Fun_{A_\iy}(\sA,\sB_2)$ to act on objects by $(\fG_*)_\Obj:\fF\mapsto\fG\ci\fF$ for $\fF:\sA\ra\sB_1$ a filtered $A_\iy$-functor, and if $p\ge 0$ and $\fF_0,\ldots,\fF_p$ are objects in $\Fun_{A_\iy}(\sA,\sB_1)$, we define
\begin{align*}
&(\fG_*)_p^{\fF_0,\ldots,\fF_p}:\Hom_{A^\rfi_\iy}(\fF_0,\fF_1)\ot\cdots\ot\Hom_{A^\rfi_\iy}(\fF_{p-1},\fF_p)
\\
&\qquad \longra \Hom_{A^\rfi_\iy}(\fG\ci\fF_0,\fG\ci\fF_p)
\end{align*}
by $(\fG_*)_0^{\fF_0}=0$, and $(\fG_*)_p^{\fF_0,\ldots,\fF_p}$ for $p\ge 1$ are as in \eq{ai4eq21}, except that we include terms $\fF_0^{i;x_j}$ and $\fG_0^{y}$, so we sum over $n_j\ge 0$, and do not require $k\le n$. Then \eq{ai4eq21} becomes an infinite sum, which converges in the usual way. As in De Deken--Lowen \cite[Prop.~3.50]{DDLo}, this defines a strict filtered $A_\iy$-functor.
\smallskip

\noindent{\bf(d)} Generalizing Definition \ref{ai4def9}(c), combining {\bf(a)\rm--\bf(c)} gives a functor
\begin{equation*}
\Fun_{A^\rfi_\iy}:(\AfCat_R)^\rop\t\AfCat_R\longra \AfCat_{R,\iy}.
\end{equation*}
acting on objects $(\sA,\sB)\mapsto\Fun_{A^\rfi_\iy}(\sA,\sB)$ and on morphisms by $(\fE^\rop,\fG)\mapsto \fE^*\ci\fG_*$. Here $\fE^*$ and $\fG_*$ commute, so their order does not matter. Note that the target $\AfCat_{R,\iy}$ is as in Definition \ref{ai5def9}, with {\it strict\/} objects and morphisms, and has a forgetful functor $\AfCat_{R,\iy}\ra\ACat_R$ forgetting filtrations.
\smallskip

\noindent{\bf(e)} The analogue of {\bf(a)\rm--\bf(d)} works for unital filtered $A_\iy$-categories and unital filtered $A_\iy$-functors, giving a functor
\begin{equation*}
\Fun_{A^\rfi_\iy}^\un:(\AfCat^\un_R)^\rop\t\AfCat^\un_R\longra \AfCat^\un_{R,\iy}.
\end{equation*}

\noindent{\bf(f)} Define $\Ga_0=\{\ga\in\Ga:\ga\sue 0\}$, so that $(\Ga_0,\pre)$ is a po-submonoid of $(\Ga,\pre)$, and $F^0R$ is a complete $(\Ga_0,\pre)$-filtered ring. Define $F^0\Fun_{A^\rfi_\iy}(\sA,\sB)$ to be the strict filtered $A_\iy$-category over $F^0R$ which has the same objects $\fF,\fG$ as $\Fun_{A^\rfi_\iy}(\sA,\sB)$, and has
\begin{equation*}
\Hom_{F^0\Fun_{A^\rfi_\iy}(\sA,\sB)}(\fF,\fG)=F^0\Hom_{\Fun_{A^\rfi_\iy}(\sA,\sB)}(\fF,\fG),
\end{equation*}
as a $\Ga_0$-filtered $F^0R$-module. The $A_\iy$-operations $\n_p^{\fF^0,\ldots,\fF^p}$ for $p\ge 0$ are the restrictions of those for $\Fun_{A^\rfi_\iy}(\sA,\sB)$ in {\bf(a)} to $F^0\Hom_{\Fun_{A^\rfi_\iy}(\sA,\sB)}(\fF^i,\fF^j)\subset \Hom_{\Fun_{A^\rfi_\iy}(\sA,\sB)}(\fF^i,\fF^j)$, and the $A_\iy$-equations also hold by restriction.

For $\sA,\sB$ unital we define a strict filtered $A_\iy$-category $F^0\Fun^\un_{A^\rfi_\iy}(\sA,\sB)$ over $F^0R$ by applying $F^0(\cdots)$ to $\Fun^\un_{A^\rfi_\iy}(\sA,\sB)$ in the same way.

As we explain in \S\ref{ai56}, the point of this is that morphisms in $\Fun_{A^\rfi_\iy}(\sA,\sB)$ are a kind of $A_\iy$-natural transformation; if we want these to preserve filtrations, we should restrict to morphisms in~$F^0\Fun_{A^\rfi_\iy}(\sA,\sB)$.
\end{dfn}

\begin{rem}
\label{ai5rem5}
We emphasize the surprising fact in Definition \ref{ai5def10} that although $\sA,\sB$ are general filtered $A_\iy$-categories (with $\m_0^{x_0}\ne 0$), $\Fun_{A^\rfi_\iy}(\sA,\sB)$ is a {\it strict\/} filtered $A_\iy$-category (with $\n_0^{\fF^0}=0$), so in particular $\Fun_{A^\rfi_\iy}(\sA,\sB)$ is an $A_\iy$-category. Similarly, although $\fE,\fG$ are general filtered $A_\iy$-functors, $\fE^*$ and $\fG_*$ are {\it strict\/} filtered $A_\iy$-functors, so in particular they are $A_\iy$-functors.

This is known to De Deken--Lowen \cite[\S 3.4--\S 3.5]{DDLo}. But Fukaya \cite[Th.~2.19]{Fuka5} says that $\Fun_{A_\iy^\rfi}(\sA,\sB)$ is strict if $\sB$ is strict, which is an unnecessary condition.

The reason why $\Fun_{A^\rfi_\iy}(\sA,\sB)$ is strict is that, in a similar equation to \eq{ai4eq14}--\eq{ai4eq17}, the appropriate definition of $\n_0^{\fF^0}$ in $\Fun_{A^\rfi_\iy}(\sA,\sB)$ is the right hand side of \eq{ai5eq9} for $\fF^0$ minus the left hand side, and so is zero.

This means we can use cohomology categories such as $H^0(\Fun_{A^\rfi_\iy}^\un(\sA,\sB))$, even though cohomology categories $H^0(\sA)$ for general filtered $A_\iy$-categories $\sA$ are not defined.
\end{rem}

\subsection[\texorpdfstring{Filtered $A_\iy$-homotopies, filtered $A_\iy$-natural transformations, \\ and filtered $A_\iy$-natural isomorphisms}{Filtered A∞-homotopies, filtered A∞-natural transformations, and filtered A∞-natural isomorphisms}]{Filtered $A_\iy$-homotopies, filtered $A_\iy$-natural transformations, and filtered $A_\iy$-natural isomorphisms}\label{ai56}

Here is the analogue of Definitions \ref{ai4def12} and \ref{ai5def7}, which is equivalent to Fukaya \cite[Def.~13.5]{Fuka5} as in~\cite[Prop.~4.2.20]{FOOO}.

\begin{dfn}
\label{ai5def11}
Let $(\Ga,\pre)$ be a po-monoid, $R$ a complete $\Ga$-filtered ring, $\sA,\sB$ be filtered $A_\iy$-categories over $R$, and $\fF,\fG:\sA\ra\sB$ be filtered $A_\iy$-functors which act the same way on objects, $\fF_\Obj=\fG_\Obj$. A {\it filtered\/ $A_\iy$-homotopy\/} $\H:\fF\Ra\fG$ consists of for all $k\ge 0$ and objects $x_0,x_1,\ldots,x_k$ in $\sA$, a graded $R$-linear map of degree $-k$
\begin{align*}
\fH_k^{x_0,\ldots,x_k}:\Hom_\sA(x_0,x_1)\ot_R\Hom_\sA(x_1,x_2)\ot_R\cdots \ot_R\Hom_\sA(x_{k-1},x_k)&\\
\longra\Hom_\sB(\fF_\Obj(x_0),\fF_\Obj(x_k))&.
\end{align*}

We require $\fH_k^{x_0,\ldots,x_k}$ to be {\it filtration-preserving} for all $k\ge 1$ and $x_0,\ldots,x_k$ in $\Obj\sA$, that is, $\fH_k^{x_0,\ldots,x_k}(F^\ga(\Hom_\sA(x_0,x_1)\ot\cdots\ot\Hom_\sA(x_{k-1},x_k)))\subseteq F^\ga \Hom_\sB(\fF(x_0),\fF(x_k))$ for all $\ga\in\Ga$.

We require that $\fH_0^{x_0}(1)\in F^{\ssu 0}\Hom^0_\sB(\fF(x_0),\fF(x_0))$ for all $x_0$ in $\Obj\sA$.

Generalizing \eq{ai4eq27}, the $\fH_k^{x_0,\ldots,x_k}$ must satisfy for all $n\ge 0$ and objects $x_0,x_1,\ldots,x_n$ in $\sA$
\ea
&\fF_n^{x_0,\ldots,x_n}-\fG_n^{x_0,\ldots,x_n}=
\label{ai5eq13}\\
&\sum_{\!\!\!\!\!\!\begin{subarray}{l} 
n=n_1+\cdots+n_k, \\
1\le l\le k,\; n_i\ge 0. \\
\text{Set $p_j=\sum_{i=1}^jn_i$,} \\ 
0\le j\le k \end{subarray}}
\begin{aligned}[t]&(-1)^{k+1}\cdot(-1)^{\sum_{i=1}^{l-1}(n_i+1)}\cdot(-1)^{\sum_{i=1}^k(n_i+1)(k-i)}\\
&\m^{\sB;\fF(x_{p_0}),\fF(x_{p_1}),\ldots,\fF(x_{p_k})}_k\ci\bigl[\fF^{x_{p_0},\ldots,x_{p_1}}_{n_1}\!\ot\!\cdots\!\ot\!\fF^{x_{p_{l-2}},\ldots,x_{p_{l-1}}}_{n_{l-1}}\\
&{}\ot \fH^{x_{p_{l-1}},\ldots,x_{p_l}}_{n_l} \ot\fG^{x_{p_l},\ldots,x_{p_{l+1}}}_{n_{l+1}}\ot\cdots\ot\fG^{x_{p_{k-1}},\ldots,x_{p_k}}_{n_k}\bigr]
\end{aligned}
\nonumber\\
&+\sum_{\begin{subarray}{l} i,j,k\ge 0:\\ i+j+k=n \end{subarray}}\begin{aligned}[t]&
(-1)^{i+jk}\fH^{x_0,x_1,\ldots,x_i,x_{i+j},x_{i+j+1},\ldots,x_n}_{i+k+1}\ci\\
&\bigl(\id_{\Hom_\sA(x_0,x_1)}\ot\cdots\ot \id_{\Hom_\sA(x_{i-1},x_i)}\ot\m^{\sA;x_i,x_{i+1},\ldots,x_{i+j}}_j \\
&\ot \id_{\Hom_\sA(x_{i+j},x_{i+j+1})}\ot\cdots\ot \id_{\Hom_\sA(x_{n-1},x_n)}\bigr).\end{aligned}
\nonumber
\ea

If $\sA,\sB,\fF,\fG$ are strict, we call $\fH$ {\it strict\/} if $\fH^{x_0}_0=0$ for all $x_0\in\Obj\sA$. Then $\fF$ is an ordinary $A_\iy$-homotopy of $A_\iy$-functors.

If $\fF,\fG$ are unital filtered $A_\iy$-functors of unital filtered $A_\iy$-categories $\sA,\sB$, we call a filtered $A_\iy$-homotopy $\H:\fF\Ra\fG$ {\it unital\/} if \eq{ai4eq28} holds for all $i,k,a_j,x_l$.
\end{dfn}

Here is the analogue of Definition \ref{ai4def13}. Part (a) follows Fukaya~\cite[Lem.-Def.~2.17]{Fuka5}.

\begin{dfn}
\label{ai5def12}
Let $(\Ga,\pre)$ be a po-monoid, $R$ a complete $\Ga$-filtered ring, $\sA,\sB$ be filtered $A_\iy$-categories over $R$, and $\fF,\fG:\sA\ra\sB$ be filtered $A_\iy$-functors. Use the notation of Definition \ref{ai5def10}.
\begin{itemize}
\setlength{\itemsep}{0pt}
\setlength{\parsep}{0pt}
\item[{\bf(a)}] A {\it filtered\/ $A_\iy$-natural transformation\/} $\fT:\fF\Ra\fG$ is an element $\fT\in F^0\Hom_{A^\rfi_\iy}^0(\fF,\fG)$ with $\n_1^{\fF,\fG}(\fT)=0$. 

Then $[\fT]:\fF\ra\fG$ is a morphism in the (non-unital) cohomology category $H^0(F^0\Fun_{A^\rfi_\iy}(\sA,\sB))$ of the (non-unital) $A_\iy$-category $F^0\Fun_{A^\rfi_\iy}(\sA,\sB)$ over $F^0R$ from Definition \ref{ai5def10}(f).
\item[{\bf(b)}] A {\it weakly filtered\/ $A_\iy$-natural transformation\/} $\fT:\fF\Ra\fG$ is an element $\fT\in \Hom_{A^\rfi_\iy}^0(\fF,\fG)$ with $\n_1^{\fF,\fG}(\fT)=0$. 

Note that this implies that $\fT\in F^\de\Hom_{A^\rfi_\iy}^0(\fF,\fG)$ for some $\de\in\Ga$, and by replacing $\de$ by $\glb(0,\de)$ we can suppose that $\de\pre 0$. Then each $\fT_k^{x_0,\ldots,x_k}$ in $\fT$ maps $F^\ga(\cdots)\ra F^{\ga+\de}(\cdots)$ for all $\ga\in\Ga$, as in~\eq{ai5eq12}.

Then $[\fT]:\fF\ra\fG$ is a morphism in the (non-unital) cohomology category $H^0(\Fun_{A^\rfi_\iy}(\sA,\sB))$ of the (non-unital) $A_\iy$-category $\Fun_{A^\rfi_\iy}(\sA,\sB)$.
\end{itemize}

Now suppose $\sA,\sB,\fF,\fG$ are unital, so that $\fF,\fG$ are objects in the unital $A_\iy$-subcategory $\Fun_{A^\rfi_\iy}^\un(\sA,\sB)\subset\Fun_{A^\rfi_\iy}(\sA,\sB)$. 
\begin{itemize}
\setlength{\itemsep}{0pt}
\setlength{\parsep}{0pt}
\item[{\bf(c)}] A {\it unital filtered\/ $A_\iy$-natural transformation\/} $\fT:\fF\Ra\fG$ is an element $\fT\in F^0\Hom_{A^\rfi_\iy}^{\un,0}(\fF,\fG)$ with $\n_1^{\fF,\fG}(\fT)=0$. 

Then $[\fT]:\fF\ra\fG$ is a morphism in the unital cohomology category $H^0(F^0\Fun_{A^\rfi_\iy}^\un(\sA,\sB))$ of the unital $A_\iy$-category $F^0\Fun_{A^\rfi_\iy}^\un(\sA,\sB)$.

We call $\fT$ a {\it filtered\/ $A_\iy$-natural isomorphism\/} if it is a homotopy equivalence in the unital $A_\iy$-category $F^0\Fun_{A^\rfi_\iy}^\un(\sA,\sB)$ over $F^0R$, in the sense of Definition \ref{ai4def5}, that is, $[\fT]$ is an isomorphism in $H^0(F^0\Fun_{A^\rfi_\iy}^\un(\sA,\sB))$.

\item[{\bf(d)}] A {\it unital weakly filtered\/ $A_\iy$-natural transformation\/} $\fT:\fF\Ra\fG$ is an element $\fT\in \Hom_{A^\rfi_\iy}^{\un,0}(\fF,\fG)$ with $\n_1^{\fF,\fG}(\fT)=0$. 

Then $[\fT]:\fF\ra\fG$ is a morphism in the unital cohomology category $H^0(\Fun_{A^\rfi_\iy}^\un(\sA,\sB))$ of the unital $A_\iy$-category $\Fun_{A^\rfi_\iy}^\un(\sA,\sB)$.

We call $\fT$ a {\it weakly filtered\/ $A_\iy$-natural isomorphism\/} if it is a homotopy equivalence in the unital $A_\iy$-category $\Fun_{A^\rfi_\iy}^\un(\sA,\sB)$ over $R$, that is, $[\fT]$ is an isomorphism in~$H^0(\Fun_{A^\rfi_\iy}^\un(\sA,\sB))$.
\end{itemize}

(Weakly) filtered $A_\iy$-natural isomorphism is an equivalence relation on unital filtered $A_\iy$-functors~$\fF:\sA\ra\sB$.
\end{dfn}

\begin{rem}
\label{ai5rem6}
{\bf(a)} In our definitions of filtered $A_\iy$-categories $\sA$ and $A_\iy$-functors $\fF$, we required the $A_\iy$-maps $\m_k^{x_0,\ldots,x_k},\fF_k^{x_0,\ldots,x_k}$ to preserve filtrations, and to increase them when $k=0$. This was necessary to avoid divergent sums.

For the filtered versions of $A_\iy$-natural transformations and $A_\iy$-natural isomorphisms, we can relax these conditions without causing divergent sums. In Definition \ref{ai5def12}(a), requiring that $\fT\in F^0\Hom_{A^\rfi_\iy}^0(\fF,\fG)\subset\Hom_{A^\rfi_\iy}^0(\fF,\fG)$ means that the maps $\fT_k^{x_0,\ldots,x_k}$ preserve filtrations. But in Definition \ref{ai5def12}(b), we must have $\fT\in F^\ep\Hom_{A^\rfi_\iy}^0(\fF,\fG)$ for some $\ep\in\Ga$, and then $\fT_k^{x_0,\ldots,x_k}$ maps $F^\ga(\cdots)\ra F^{\ga+\ep}(\cdots)$ for all $k,x_0,\ldots,x_k$. So for {\it weakly\/} filtered $A_\iy$-natural transformations, the maps $\fT_k^{x_0,\ldots,x_k}$ need not preserve filtrations, but can decrease them by a bounded amount~$\ep$.

We outline an application of weakly filtered $A_\iy$-natural isomorphisms in Symplectic Geometry, to justify the definition. In the ongoing project of Fukaya--Oh--Ohta--Ono, if $(X,\om)$ is a suitable symplectic manifold, they aim (though this may not yet be complete?) to define a {\it filtered Fukaya category\/} $\sF^\rfi(X,\om)$, a filtered $A_\iy$-category over $\La^\nov_\Q$ in Example \ref{ai5ex3}, with objects suitable Lagrangians $L$ in $(X,\om)$. Although \cite{FOOO} is not written in terms of $A_\iy$-categories, the author expects that some of their results \cite[Th.~G(4), Th.~4.1.5, Def.~5.2.1, \& Th.~5.3.65]{FOOO} may be interpreted as follows. 

Let $\Phi:(X,\om)\ra(X,\om)$ be a symplectomorphism Hamiltonian isotopic to the identity. Then one should construct a filtered $A_\iy$-functor $\Phi_*:\sF^\rfi(X,\om)\ra\sF^\rfi(X,\om)$ acting on objects by $L\mapsto\Phi(L)$. We might expect there to exist a filtered $A_\iy$-natural isomorphism $\fH:\Id\Ra\Phi_*$. However, it seems clear from \cite{FOOO} that the maps $\fH_k^{x_0,\ldots,x_k}$ should not be filtration-preserving; instead, following `weakly filtered $A_\iy$-bimodule homomorphisms' in \cite[Def.~5.2.1]{FOOO}, we should allow $\fH_k^{x_0,\ldots,x_k}$ to map $F^\la(\cdots)\mapsto F^{\la-\ep}(\cdots)$ for some $\ep\ge 0$ independent of $k,x_0,\ldots,x_k$. So $\fH$ should be a {\it weakly\/} filtered $A_\iy$-natural isomorphism.
\smallskip

\noindent{\bf(b)} We do not define a `weak' version of filtered $A_\iy$-homotopy in Definition \ref{ai5def11}. This is because, as in \cite[Def.~13.5]{Fuka5}, filtered $A_\iy$-homotopies can be defined as examples of filtered $A_\iy$-functors, so they should satisfy the same filtration conditions as filtered $A_\iy$-functors.
\smallskip

\noindent{\bf(c)} The analogues of Proposition \ref{ai4prop3}, Theorem \ref{ai4thm3}, and Definition \ref{ai4def15} do not make sense in the filtered case, as filtered $A_\iy$-categories $\sA$ do not have cohomology categories~$H^0(\sA)$.
\end{rem}

Here is the analogue of Proposition \ref{ai4prop4}.

\begin{prop}
\label{ai5prop1}
Let\/ $(\Ga,\pre)$ be a po-monoid, $R$ a complete $\Ga$-filtered ring, $\sA,\sB$ be  unital filtered\/ $A_\iy$-categories over $R,$ and\/ $\fF,\fG:\sA\ra\sB$ be unital filtered\/ $A_\iy$-functors which act the same way on objects, $\fF_\Obj=\fG_\Obj,$ and\/ $\fH:\fF\Ra\fG$ be a unital filtered\/ $A_\iy$-homotopy. Define $\fT\in\Hom_{A^\rfi_\iy}^0(\fF,\fG)$ by
\begin{equation*}
\fT_k^{x_0,\ldots,x_k}=\begin{cases} \boo^\sB_{\fF(x_0)}+\fH_0^{x_0}, & k=0, \\
\H_k^{x_0,\ldots,x_k}, & k\ge 1,
\end{cases}
\end{equation*}
for all\/ $k\ge 0$ and\/ $x_0,\ldots,x_k\in\Obj\sA$. Then $\fT:\fF\Ra\fG$ is a filtered $A_\iy$-natural isomorphism. This gives a 1-1 correspondence between unital filtered\/ $A_\iy$-homotopies $\fH:\fF\Ra\fG$ and filtered\/ $A_\iy$-natural isomorphisms\/ $\fT:\fF\Ra\fG$ with\/ $\fT_0^{x_0}(1)-\boo^\sB_{\fF(x_0)}\in F^{\ssu 0}\Hom_\sB^0(\fF(x_0),\fF(x_0))$ for all\/~$x_0\in\Obj\sA$.
\end{prop}

\begin{proof}
Essentially the same proof as that of Proposition \ref{ai4prop4} works in the filtered case, but with the following modification. Equation \eq{ai4eq31} no longer holds, and should be replaced by
\begin{align*}
&\m_2^{\fF(x_0),\fF(x_0),\fF(x_n)}\bigl(\fT_0^{x_0}(1)\ot\fU_n^{x_0,\ldots,x_n}\bigr)+\bigl(\text{terms in $\fU_k^{x_0',\ldots,x_k'}$ for $k<n$}\bigr) \\
&+\bigl(\text{terms in $\fU_{n'}^{x'_0,\ldots,x'_{n'}}$, $n'\ge n$, and at least one $\fF_0^x(1),\fG_0^x(1),\fT_0^x(1)-\boo^\sB_{\fF(x)}$}\bigr) \\
&\qquad =\begin{cases} \boo^\sB_{\fF(x_0)}, & n=0, \\ 0, & n>0.
\end{cases}
\end{align*}
Since $\fF_0^x(1),\fG_0^x(1),\fT_0^x(1)-\boo^\sB_{\fF(x)}$ lie in $F^{\ssu 0}\Hom_\sB(\fF(x),\fF(x))$, the terms on the second line increase the filtration of $\fU_{n'}^{x'_0,\ldots,x'_{n'}}$. We then show there is a unique $\fU$ satisfying \eq{ai4eq30} by a double induction, in which we determine $\fU_n^{x_0,\ldots,x_n}$ modulo $F^\ga(\cdots)$ by induction firstly on $n$, and secondly on increasing~$\ga$.
\end{proof}

Here is the analogue of Definition \ref{ai4def14}.

\begin{dfn}
\label{ai5def13}
Let $(\Ga,\pre)$ be a po-monoid, $R$ a complete $\Ga$-filtered ring, $\sA,\sB$ be unital filtered $A_\iy$-categories over $R$, and $\fF:\sA\ra\sB$ be a unital filtered $A_\iy$-functor. We call $\fF$ a {\it homotopy equivalence\/} if there exists a unital filtered $A_\iy$-functor $\fG:\sB\ra\sA$ and filtered $A_\iy$-natural isomorphisms $\fS:\fG\ci\fF\Ra\id_\sA$ and $\fT:\fF\ci\fG\Ra\id_\sB$. Then we also call $\sA,\sB$ {\it homotopy equivalent}. 

We call $\fF$ a {\it weak homotopy equivalence}, and $\sA,\sB$ {\it weakly homotopy equivalent}, if the same happens with $\fS,\fT$ weakly filtered $A_\iy$-natural isomorphisms. 

Compositions of (weak) homotopy equivalences are (weak) homotopy equivalences. (Weak) homotopy equivalence is an equivalence relation on unital filtered $A_\iy$-categories.
\end{dfn}

\subsection{Bounding cochains, and their applications}
\label{ai57}

{\it Bounding cochains\/} are a method of modifying filtered $A_\iy$-algebras or filtered $A_\iy$-categories to make them strict. They were introduced by Fukaya et al.\ \cite[\S 3.6]{FOOO} for $A_\iy$-algebras, and by Fukaya \cite[Def.~1.18]{Fuka3}, \cite[Def.~2.5]{Fuka5} for $A_\iy$-categories, and are a vital part of their programme for constructing Lagrangian Floer cohomology and Fukaya categories. We first explain them following \cite{Fuka3,Fuka5,FOOO}, and then in Proposition \ref{ai5prop2} we give the author's perspective on them.

\begin{dfn}
\label{ai5def14}
Let $(\Ga,\pre)$ be a po-monoid, $R$ a complete $\Ga$-filtered ring, and $(A,\m)$ be a filtered $A_\iy$-algebra over $R$. Let $b\in F^{\ssu 0}A^1$. Define graded $R$-linear maps $\m_k^b:A^{\ot^k}\ra A$ for $k=0,1,\ldots,$ of degree $2-k$, by
\ea
&\m_k^b(a_1\ot\cdots\ot a_k)=
\label{ai5eq14}\\
&\sum_{n_0,\ldots,n_k\ge
0}\m_{k+n_0+\cdots+n_k}\bigl(b^{\ot^{n_0}}\ot a_1\ot b^{\ot^{n_1}}\ot a_2\ot\cdots\ot b^{\ot^{n_{k-1}}}\ot a_k\ot b^{\ot^{n_k}}\bigr).
\nonumber
\ea
This is an infinite sum, but converges as $b\in F^{\ssu 0}A^1$. It is filtration-preserving. Write $\m^b=(\m_k^b)_{k\ge 0}$. Then as in \cite[\S 3.6.1]{FOOO}, $(A,\m^b=(\m^b_k)_{k\ge 0})$ is another filtered $A_\iy$-algebra over $R$, with $A_\iy$-maps $\m_k$ `twisted by 
~$b$'.

If $\boo$ is a unit for $(A,\m)$ then it is also a unit for $(A,\m^b)$.

There is a filtered $A_\iy$-isomorphism $\fI^b:(A,\m)\ra(A,\m^b)$ given by $\fI^b_0=-b$, $\fI^b_1=\id_A$ and $\fI^b_k=0$ if $k\ge 2$. It has an inverse filtered $A_\iy$-isomorphism $\fJ^b:(A,\m^b)\ra(A,\m)$ given by $\fJ^b_0=b$, $\fJ^b_1=\id_A$ and $\fJ^b_k=0$ if~$k\ge 2$.
 
We call $b$ a {\it bounding cochain\/} for $(A,\m)$ if $\m_0^b=0$, i.e.~if
\e
\ts\sum_{k\ge 0}\m_k(b,\ldots,b)=0\quad\text{in $A^2$.}
\label{ai5eq15}
\e
This is the {\it Maurer--Cartan equation}, or {\it Batalin--Vilkovisky master equation}. If $b$ is a bounding cochain then $(A,\m^b)$ is a {\it strict\/} filtered $A_\iy$-algebra, and thus an $A_\iy$-algebra.
\end{dfn}

\begin{dfn}
\label{ai5def15}
Let $(\Ga,\pre)$ be a po-monoid, $R$ a complete $\Ga$-filtered ring, and $\sA$ a filtered $A_\iy$-category over $R$. Suppose $x$ is an object of $\sA$. A {\it bounding cochain\/ $b$ for\/} $x$ is an element $b$ in $F^{\ssu 0}\Hom_\sA^1(x,x)$ satisfying, as in \eq{ai5eq15}
\e
\ts\sum_{k\ge 0}\m_k^{x,\ldots,x}(b,\ldots,b)=0\quad\text{in $\Hom_\sA^2(x,x)$.}
\label{ai5eq16}
\e

\noindent{\bf(b)} We will define a strict filtered $A_\iy$-category $\Str(\sA)$ over $R$ called the {\it associated strict category\/} to $\sA$. Objects of $\Str(\sA)$ are pairs $(x,b)$ of an object $x$ in $\sA$ and a bounding cochain $b$ for $x$ in $\sA$. For objects $(x,b),(x',b')$ we define $\Hom_{\Str(\sA)}\bigl((x,b),(x',b')\bigr)= \Hom_\sA(x,x')$ as a $\Ga$-filtered $\Z$-graded $R$-module. 

If $k\ge 0$ and $(x_0,b_0),\ldots,(x_k,b_k)$ are objects in $\Str(\sA)$, then generalizing \eq{ai5eq14} we define the $A_\iy$-operations $\n_k^{(x_0,b_0),\ldots,(x_k,b_k)}$ in $\Str(\sA)$ by
\ea
&\n_k^{(x_0,b_0),\ldots,(x_k,b_k)}(a_1\ot\cdots\ot a_k)=
\label{ai5eq17}\\
&\sum_{n_0,\ldots,n_k\ge
0}\begin{aligned}[t]&\m_{k+n_0+\cdots+n_k}^{(n_0+1)x_0,(n_1+1)x_1,\ldots,(n_k+1)x_k}\bigl(b_0^{\ot^{n_0}}\ot a_1\ot b_1^{\ot^{n_1}}\ot a_2\ot\cdots\\
&\qquad \qquad\ot b_{k-1}^{\ot^{n_{k-1}}}\ot a_k\ot b_k^{\ot^{n_k}}\bigr).
\end{aligned}
\nonumber
\ea
Here by $(n_0+1)x_0,(n_1+1)x_1,\ldots$ we mean the list with $x_0$ repeated $n_0+1$ times, followed by $x_1$ repeated $n_1+1$ times, and so on. Equation \eq{ai5eq17} is an infinite sum in $\Hom_{\Str(\sA)}((x_0,b_0),(x_k,b_k))=\Hom_\sA(x_0,x_k)$, which converges as $b_i\in F^{\ssu 0}\Hom_\sA^1(x_i,x_i)$. As in Fukaya \cite[Def.~1.18]{Fuka3}, \cite[Def.~2.5]{Fuka5}, this defines a filtered $A_\iy$-category $\Str(\sA)$ over~$R$. 
\smallskip

\noindent{\bf(c)} If $\sA$ is a {\it unital\/} filtered $A_\iy$-category over $R$, with units $\boo_x^\sA\in \Hom_\sA^0(x,x)$, then $\Str(\sA)$ is also unital, with units $\boo_{(x,b)}^{\Str(\sA)}=\boo_x^\sA$.
\end{dfn}

\begin{rem}
\label{ai5rem7}
Here is how bounding cochains and associated strict categories are used in Fukaya--Oh--Ohta--Ono's programme. Suppose $(X,\om)$ is a symplectic manifold, satisfying some assumptions. They perform the following steps:
\begin{itemize}
\setlength{\itemsep}{0pt}
\setlength{\parsep}{0pt}
\item[(i)] For the po-monoid $(\Ga,\pre)=(\R,\le)$ and $\R$-filtered ring $R=\La_\Q^\nov$ from Example \ref{ai5ex3}, after making some additional choices such as an almost complex structure $J$ on $X$, they construct a $\La_\Q^\nov$-filtered $A_\iy$-category $\sF(X,\om)$, whose objects are suitable Lagrangians $L$ in $(X,\om)$. 

This construction is very long and difficult, with many challenging analytic issues. Most of the work is done in \cite{FOOO}, although this is written in terms of filtered $A_\iy$-algebras rather than filtered $A_\iy$-categories.
\item[(ii)] They pass to the associated strict category $\Str(\sF(X,\om))$, an $A_\iy$-category over $\La_\Q^\nov$, whose objects are pairs $(L,b)$ of a Lagrangian $L$ in $(X,\om)$ and a bounding cochain $b$ for the filtered $A_\iy$-algebra of~$L$.
\item[(iii)] They apply the `twist construction' $\Tw\ci\Str(\sF(X,\om))$ from \S\ref{ai49} to get a pretriangulated $A_\iy$-category.
\item[(iv)] The {\it derived Fukaya category\/} is $D^b\sF(X,\om)=H^0(\Tw\ci\Str(\sF(X,\om)))$. It is a triangulated category, of great interest in Symplectic Geometry. For example, Kontsevich's Homological Mirror Symmetry Conjecture \cite{Kont} explains Mirror Symmetry of Calabi--Yau manifolds $X,\check X$ in terms of equivalences of triangulated categories $D^b\sF(X,\om)\simeq D^b\coh\check X$, $D^b\sF(\check X,\check\om)\simeq D^b\coh X$.
\end{itemize}
\end{rem}

We will now interpret associated strict categories in terms of filtered $A_\iy$-functor categories in a way which seems very natural, but the author has not found it in the literature, even in Fukaya \cite{Fuka5} which defines both ideas. Our point of view will make it easy to generalize associated strict categories to the $\iy$-categorical setting. Recall that Definition \ref{ai4def10} defined the point $A_\iy$-category $\sP$ and Proposition \ref{ai4prop1} showed that $\Fun_{A_\iy}(\sP,\sA)\cong\sA$ for any $A_\iy$-category $\sA$, and Definition \ref{ai4def11} and Proposition \ref{ai4prop2} did the same for unital $A_\iy$-categories. In the filtered case, instead we find that~$\Fun_{A^\rfi_\iy}(\sP,\sA)\cong\Str(\sA)$.

\begin{prop}
\label{ai5prop2}
Let\/ $(\Ga,\pre)$ be a po-monoid, $R$ a complete $\Ga$-filtered ring, and\/ $\sA$ be a filtered\/ $A_\iy$-category over $R$. Let\/ $\sP$ and\/ $\sP^\un$ be the non-unital and unital point\/ $A_\iy$-categories over $R$ as in Definitions\/ {\rm\ref{ai4def10}} and\/ {\rm\ref{ai4def11},} made into strict filtered\/ $A_\iy$-categories in the obvious way.
\smallskip

\noindent{\bf(a)} A filtered\/ $A_\iy$-functor $\fF:\sP\ra\sA$ is equivalent to a pair $(x,b)$ where $x\in\Obj\sA$ and\/ $b$ is a bounding cochain for $x$ in $\sA$.
\smallskip

\noindent{\bf(b)} The strict filtered\/ $A_\iy$-category $\Fun_{A^\rfi_\iy}(\sP,\sA)$ from Definition\/ {\rm\ref{ai5def10}} may be identified with the associated strict category\/ $\Str(\sA)$ from Definition\/~{\rm\ref{ai5def15}(b)}.

Define a functor $\Str=\Fun_{A^\rfi_\iy}(\sP,-):\AfCat_R\ra\AfCat_{R,\iy}$. This makes associated strict $A_\iy$-categories functorial.
\smallskip

\noindent{\bf(c)} If\/ $\sA$ is unital then a unital filtered\/ $A_\iy$-functor $\fF:\sP^\un\ra\sA$ is equivalent to a pair $(x,b)$ where $x\in\Obj\sA$ and\/ $b$ is a bounding cochain for $x$ in~$\sA$.
\smallskip

\noindent{\bf(d)} If\/ $\sA$ is unital then the unital strict filtered\/ $A_\iy$-category $\Fun_{A^\rfi_\iy}^\un(\sP^\un,\sA)$ from Definition\/ {\rm\ref{ai5def10}} may be identified with the associated strict category\/ $\Str(\sA)$ from Definition\/ {\rm\ref{ai5def15}(b)}. The restriction of the functor $\Str$ in {\bf(b)} to $\AfCat^\un_R$ is naturally isomorphic to $\Fun^\un_{A^\rfi_\iy}(\sP^\un,-):\AfCat^\un_R\ra\AfCat^\un_{R,\iy}$.
\end{prop}

\begin{proof}
A filtered $A_\iy$-functor $\fF:\sP\ra\sA$ has a map $\fF_\Obj:\Obj\sP=\{*\}\ra\Obj\sA$, and so determines an object $x=\fF(*)$ in $\sA$. For the maps $\fF_k^{*,\ldots,*}$ for $k\ge 0$, we see that $\fF_k^{*,\ldots,*}=0$ for $k\ge 1$ as $\Hom_\sP(*,*)=0$, and $\fF_0^*:R\ra\Hom_\sA^1(x,x)$ is equivalent to a choice of $b=\fF_0^*(1)$ in $\Hom_\sA^1(x,x)$. By definition $b\in F^{\ssu 0}\Hom_\sA^1(x,x)$, and equation \eq{ai5eq9} for $n=0$ is equivalent to \eq{ai5eq16}. Thus $b$ is a bounding cochain, and (a) follows. Part (b) is then easy to check.

For the unital case, let $\fF:\sP^\un\ra\sA$ be a unital filtered $A_\iy$-functor, and set $x=\fF(*)$ in $\Obj\sA$. Since $\Hom_{\sP^\un}(*,*)=R\cdot\boo^{\sP^\un}_*$, equations \eq{ai4eq11}--\eq{ai4eq12} imply that $\fF_1^{*,*}$ maps $\la\boo^{\sP^\un}_*\mapsto \la\boo^\sA_x$ for $\la\in R$ and $\fF_k^{*,\ldots,*}=0$ for $k\ge 2$. Again we set $b=\fF_0^*(1)$ in $F^{\ssu 0}\Hom_\sA^1(x,x)$. Then \eq{ai5eq9} for $n=0$ is equivalent to \eq{ai5eq16}, so $b$ is a bounding cochain. We can show that \eq{ai5eq9} for $n\ge 1$ holds automatically, and impose no extra conditions on $(x,b)$. Part (c) follows, and (d) is then easy to check.
\end{proof}

\subsection{Explicit definition of the twist functor}
\label{ai58}

Let $R$ be a commutative ring. As in Definition \ref{ai4def22}, the {\it twist functor\/}
\e
\Tw:\ACat_R\longra\ADCat_R,
\label{ai5eq18}
\e
is defined in Seidel \cite[\S 3]{Seid}, Keller \cite[\S 7.6]{Kell1}, Lef\`evre-Hasegawa \cite[\S 7.2]{LeHa}, and Bespalov--Lyubashenko--Manzyuk \cite[\S 12]{BLM}. We will describe a construction of $\Tw$ which passes through $\Z$-filtered $A_\iy$-categories, and uses the functor $\sA\mapsto\Fun_{A^\rfi_\iy}(\sP,\sA)$ defining associated strict categories in Proposition \ref{ai5prop2}(a),(b) as an intermediate step. This will be easy to extend to $\iy$-categories in Chapter~\ref{ai9}.

\begin{dfn}
\label{ai5def16}
Let $R$ be a commutative ring. Make $R$ into a complete filtered ring over the po-monoid $(\Z,\le)$ by setting $F^\la R=R$ if $\la\le 0$ in $\Z$ and $F^\la R=0$ if $\la>0$ in $\Z$. We will define $\Tw$ in \eq{ai5eq18} as the composition of three functors
\e
\begin{gathered}
\xymatrix@C=140pt@R=17pt{
*+[r]{\ACat_R} \ar[r]_\Si \ar[d]^\Tw & *+[l]{\AfCat_{R,\iy}} \ar[d]_{\Fun_{A^\rfi_\iy}(\sP,-)} \\
*+[r]{\ACat_R}  & *+[l]{\AfCat_{R,\iy}.\!} \ar[l]_{\text{forget filtrations}} }
\end{gathered}
\label{ai5eq19}
\e

The restriction of $\Tw$ to $\ACat_R^\un$ in \eq{ai4eq43} may equivalently be defined as the composition
\e
\begin{gathered}
\xymatrix@C=140pt@R=17pt{
*+[r]{\ACat^\un_R} \ar[r]_\Si \ar[d]^\Tw & *+[l]{\AfCat^\un_{R,\iy}} \ar[d]_{\Fun^\un_{A^\rfi_\iy}(\sP^\un,-)} \\
*+[r]{\ACat^\un_R}  & *+[l]{\AfCat^\un_{R,\iy},\!} \ar[l]_{\text{forget filtrations}} }
\end{gathered}
\label{ai5eq20}
\e
where $\sP^\un$ is as in Definition~\ref{ai4def11}.

Here for each $\sA$ in $\ACat_R$ we define $\Si(\sA)$ in $\AfCat_R$ as follows:
\begin{itemize}
\setlength{\itemsep}{0pt}
\setlength{\parsep}{0pt}
\item[(a)] Objects of $\Si(\sA)$ are triples $\bigl((x_1,\ldots,x_l),(s_1,\ldots,s_l),(f_1,\ldots,f_l)\bigr)$ where $l\ge 0$, and $x_1,\ldots,x_l\in\Obj\sA$, and $s_i,f_i\in\Z$ for $i=1,\ldots,l$. We call $s_i$ the {\it shift\/} of $x_i$ and $f_i$ the {\it filtration weight\/} of $x_i$. As a shorthand we write $(\bs x,\bs s,\bs f)=\bigl((x_1,\ldots,x_l),(s_1,\ldots,s_l),(f_1,\ldots,f_l)\bigr)$.

\item[(b)] Let $\bigl((x_1,\ldots,x_l),(s_1,\ldots,s_l),(f_1,\ldots,f_l)\bigr)$ and $\bigl((x'_1,\ldots,x'_{l'})$, $(s'_1,\ldots,\ab s'_{l'})$, $(f'_1,\ab\ldots,f'_{l'})\bigr)$ be objects of $\Si(\sA)$. Define a graded $R$-module
\begin{align*}
&\Hom_{\Si(\sA)}\bigl((\bs x,\bs s,\bs f),(\bs x',\bs s',\bs f')\bigr)=\bigop_{d\in\Z}\Hom_{\Si(\sA)}^d\bigl((\bs x,\bs s,\bs f),(\bs x',\bs s',\bs f')\bigr)\;\>\text{by} \\
&\Hom_{\Si(\sA)}^d\bigl((\bs x,\bs s,\bs f),(\bs x',\bs s',\bs f')\bigr)=\bigop_{i=1,\ldots,l,\;
j=1,\ldots,l'\!\!\!\!\!\!\!\!\!\!\!\!\!\!\!\!\!\!\!} \Hom_\sA^{d+s_i-s'_j}(x_i,x'_j).
\end{align*}
We write elements of $\Hom_{\Si(\sA)}\bigl((\bs x,\bs s,\bs f),(\bs x',\bs s',\bs f')\bigr)$ as $(\fT_{ij})_{i=1,\ldots,l}^{j=1,\ldots,l'}$ in the obvious way. Define a $\Z$-filtration on $\Hom_{\Si(\sA)}^d\bigl((\bs x,\bs s,\bs f),(\bs x',\bs s',\bs f')\bigr)$ by, for $\la\in\Z$
\begin{equation*}
F^\la\Hom_{\Si(\sA)}^d\bigl((\bs x,\bs s,\bs f),(\bs x',\bs s',\bs f')\bigr)=\bigop_{\begin{subarray}{l} i=1,\ldots,l,\; j=1,\ldots,l': \\ f_i-f_j\ge\la \end{subarray}\!\!\!\!\!\!\!\!\!\!\!\!\!\!\!\!\!\!}
\Hom_\sA^{d+s_i-s'_j}(x_i,x'_j).
\end{equation*}
\item[(c)] Let $(\bs x^0,\bs s^0,\bs f^0),\ldots,(\bs x^k,\bs s^k,\bs f^k)$ be objects of $\Si(\sA)$. For $k\ge 0$, define a graded $R$-linear map of degree $2-k$
\begin{align*}
\n_k^{(\bs x^0,\bs s^0,\bs f^0),\ldots,(\bs x^k,\bs s^k,\bs f^k)}:\,&\Hom_{\Si(\sA)}\bigl((\bs x^0,\bs s^0,\bs f^0),(\bs x^1,\bs s^1,\bs f^1)\bigr)\ot_R\cdots \ot_R\\
&\Hom_{\Si(\sA)}\bigl((\bs x^{k-1},\bs s^{k-1},\bs f^{k-1}),(\bs x^k,\bs s^k,\bs f^k)\bigr)\\
&\quad \longra\Hom_{\Si(\sA)}\bigl((\bs x^0,\bs s^0,\bs f^0),(\bs x^k,\bs s^k,\bs f^k)\bigr)
\end{align*}
by $\n_0^{(\bs x^0,\bs s^0,\bs f^0)}=0$ for $k=0$, and for $k\ge 1$ define, for $j^0=1,\ldots,l^0$ and $j^k=1,\ldots,l^k$
\end{itemize}
\begin{align*}
&\Bigl(\n_k^{(\bs x^0,\bs s^0,\bs f^0),\ldots,(\bs x^k,\bs s^k,\bs f^k)}\bigl((\fT^1_{i^0i^1})_{i^0=1,\ldots,l^0}^{i^1=1,\ldots,l^1}\ot\cdots\ot (\fT^k_{i^{k-1}i^k})_{i^{k-1}=1,\ldots,l^{k-1}}^{i^k=1,\ldots,l^k}\bigr)\Bigr)_{j^0j^k} \\
&=\sum_{\begin{subarray}{l} j^1=1,\ldots,l^1,\ldots, \\
j^{k-1}=1,\ldots,l^{k-1}
\end{subarray}}
\m_k^{\sA;x^0_{j^0},x^1_{j^1},\ldots,x^{k-1}_{j^{k-1}},x^k_{j^k}}\bigl(\fT^1_{j^0j^1}\ot\cdots\ot \fT^k_{j^{k-1}j^k}\bigr).
\end{align*}

It is now easy to show that this data makes $\Si(\sA)$ into a filtered $A_\iy$-category over $R$, which is strict as $\m_0^{(\bs x^0,\bs s^0,\bs f^0)}=0$ for all $(\bs x^0,\bs s^0,\bs f^0)$. In particular, equation \eq{ai4eq1} for the $\m_k^{x_0,\ldots,x_k}$ implies \eq{ai5eq8} for the~$\n_k^{(\bs x^0,\bs s^0,\bs f^0),\ldots,(\bs x^k,\bs s^k,\bs f^k)}$.

Next let $\fF:\sA\ra\sB$ be a morphism in $\ACat_R$. Define a filtered $A_\iy$-functor $\Si(\fF):\Si(\sA)\ra\Si(\sB)$ by:
\begin{itemize}
\setlength{\itemsep}{0pt}
\setlength{\parsep}{0pt}
\item[(d)] $\Si(\fF)$ acts on objects by $\bigl((x_1,\ldots,x_l),(s_1,\ldots,s_l),(f_1,\ldots,f_l)\bigr)\mapsto \ab\bigl((\fF(x_1),\ab\ldots,\ab\fF(x_l)),(s_1,\ldots,s_l),(f_1,\ldots,f_l)\bigr)$.
\item[(e)] Let $(\bs x^0,\bs s^0,\bs f^0),\ldots,(\bs x^k,\bs s^k,\bs f^k)$ be objects of $\Si(\sA)$. For $k\ge 0$, define a graded $R$-linear map of degree $1-k$
\begin{align*}
&\Si(\fF_k)^{(\bs x^0,\bs s^0,\bs f^0),\ldots,(\bs x^k,\bs s^k,\bs f^k)}:\Hom_{\Si(\sA)}\bigl((\bs x^0,\bs s^0,\bs f^0),(\bs x^1,\bs s^1,\bs f^1)\bigr)\\
&\qquad\qquad \ot_R\cdots \ot_R\Hom_{\Si(\sA)}\bigl((\bs x^{k-1},\bs s^{k-1},\bs f^{k-1}),(\bs x^k,\bs s^k,\bs f^k)\bigr)\\
&\qquad\qquad\quad \longra\Hom_{\Si(\sA)}\bigl((\fF(\bs x^0),\bs s^0,\bs f^0),(\fF(\bs x^k),\bs s^k,\bs f^k)\bigr)
\end{align*}
by $\Si(\fF_0)^{(\bs x^0,\bs s^0,\bs f^0)}=0$ for $k=0$, and for $k\ge 1$ define
\end{itemize}
\begin{align*}
&\Bigl(\Si(\fF_k)_k^{(\bs x^0,\bs s^0,\bs f^0),\ldots,(\bs x^k,\bs s^k,\bs f^k)}\bigl((\fT^1_{i^0i^1})_{i^0=1,\ldots,l^0}^{i^1=1,\ldots,l^1}\!\ot\!\cdots\!\ot\!(\fT^k_{i^{k-1}i^k})_{i^{k-1}=1,\ldots,l^{k-1}}^{i^k=1,\ldots,l^k}\bigr)\Bigr)_{j^0j^k} \\
&=\sum_{\begin{subarray}{l} j^1=1,\ldots,l^1,\ldots, \\
j^{k-1}=1,\ldots,l^{k-1}
\end{subarray}}
\fF_k^{\sA;x^0_{j^0},x^1_{j^1},\ldots,x^{k-1}_{j^{k-1}},x^k_{j^k}}\bigl(\fT^1_{j^0j^1}\ot\cdots\ot \fT^k_{j^{k-1}j^k}\bigr),
\end{align*}
\begin{itemize}
\setlength{\itemsep}{0pt}
\setlength{\parsep}{0pt}
\item[] for $j^0=1,\ldots,l^0$ and $j^k=1,\ldots,l^k$. It is now easy to show that this data makes $\Si(\fF)$ into a filtered $A_\iy$-functor over $R$, which is strict as  $\Si(\fF_0)^{(\bs x^0,\bs s^0,\bs f^0)}=0$ for all $(\bs x^0,\bs s^0,\bs f^0)$. In particular, equation \eq{ai4eq8} for the $\fF_k^{x_0,\ldots,x_k}$ implies equation \eq{ai5eq9} for the~$\Si(\fF)_k^{(\bs x^0,\bs s^0,\bs f^0),\ldots,(\bs x^k,\bs s^k,\bs f^k)}$.
\end{itemize}
As $\fF\mapsto\Si(\fF)$ is compatible with compositions and identities, this defines a functor $\Si:\ACat_R\ra\AfCat_{R,\iy}$, where the target is $\AfCat_{R,\iy}\subset\AfCat_R$ as $\Si(\sA)$ and $\Si(\fF)$ are both strict.

The right hand functor $\Fun_{A^\rfi_\iy}(\sP,-)$ in \eq{ai5eq19} is as in Proposition \ref{ai5prop2}, where $\sP$ is the non-unital point $A_\iy$-category over $R$ from Definition \ref{ai4def10}, made into a $\Z$-filtered $A_\iy$-category over $R$ in the unique way. Note that this functor maps $\sA\mapsto\Str(\sA)$, where $\Str(\sA)$ is the associated strict category from Definition \ref{ai5def15}(b). The bottom functor in \eq{ai5eq17} forgets filtration data, as in the vertical functor in \eq{ai5eq11}. This completes \eq{ai5eq17}, and the definition of~$\Tw$.

Objects of $\Tw\sA$, called {\it twisted complexes}, are pairs $\bigl((\bs x,\bs s,\bs f),\bs b\bigr)$ where $(\bs x,\bs s,\bs f)$ is an object in $\Si(\sA)$ and $\bs b$ is a bounding cochain for $(\bs x,\bs s,\bs f)$. Explicitly, $\bs b=(b_{ij})_{i=1,\ldots,l}^{j=1,\ldots,l}$ with $b_{ij}\in \Hom_\sA^{1+s_i-s'_j}(x_i,x_j)$, where the condition that $\bs b\in F^{\ssu 0}\bigl((\bs x,\bs s,\bs f),(\bs x,\bs s,\bs f)\bigr)$ means that $b_{ij}=0$ if $f_i-f_j\le 0$. This $\bs b$ must satisfy \eq{ai5eq16}, which has only finitely many nonzero terms.

The natural transformation $\inc:\Id\Ra\Tw$ in Definition \ref{ai4def22} is defined as follows: for each $\sA$ in $\ACat_R$, the $A_\iy$-functor $\inc_\sA:\sA\hookra\Tw\sA$ acts on objects by $x\mapsto\bigl(((x),(0),(0)),(0)\bigr)$, and $(\inc_\sA)_k^{x_0,\ldots,x_k}$ is $\id_{\Hom_\sA(x_0,x_1)}$ if $k=1$ and $0$ if $k\ne 1$. This is a strict isomorphism with a full $A_\iy$-subcategory of~$\Tw\sA$.
\end{dfn}

\section{\texorpdfstring{$\iy$-categories of $A_\iy$-algebras}{∞-categories of A∞-algebras}}
\label{ai6}

We now define many different $\iy$-categories of $A_\iy$-algebras over a commutative ring $R$, which fit into a diagram of inclusions of $\iy$-categories:
\ea
\label{ai6eq1}\\[-19pt]
\xymatrix@!0@C=85pt@R=37pt{{\begin{subarray}{l}\ts \AAlg^{\un,\iy}_R\!\!=  \\[2pt] \ts \AAlg^{\un,\iy}_{R,2}\end{subarray}} \\
\AAlg^{\un,\iy^{\vphantom{(}}}_{R,3}\, \ar@{^{(}->}[u] \ar@{=}[r] \ar@{^{(}->}[u] & \dAAlg^{\un,\iy^{\vphantom{(}}}_{R,3}\, \ar@{^{(}->}[r] & {\begin{subarray}{l}\ts \AAlg^{\cu,\iy}_R\!\!=  \\[2pt] \ts \AAlg^{\cu,\iy}_{R,3}\end{subarray}}\,  \ar@{^{(}->}[r] & {\begin{subarray}{l}\ts \AAlg^\iy_R\!=  \\[2pt] \ts \AAlg^\iy_{R,3}\end{subarray}}  \\
\AAlg^{\un,\iy^{\vphantom{(}}}_{R,4}\, \ar@{=}[r] \ar@{^{(}->}[u] & \dAAlg^{\un,\iy^{\vphantom{(}}}_{R,4}\, \ar@{^{(}->}[r] \ar@{^{(}->}[u] & \AAlg^{\cu,\iy^{\vphantom{(}}}_{R,4}\,  \ar@{^{(}->}[r] \ar@{^{(}->}[u] & \AAlg^{\iy^{\vphantom{(}}}_{R,4} \ar@{^{(}->}[u] \\
\AAlg^{\un,\iy^{\vphantom{(}}}_{R,5}\, \ar@{=}[r] \ar@{^{(}->}[u] & \dAAlg^{\un,\iy^{\vphantom{(}}}_{R,5}\, \ar@{^{(}->}[r] \ar@{^{(}->}[u] & \AAlg^{\cu,\iy^{\vphantom{(}}}_{R,5}\,  \ar@{^{(}->}[r] \ar@{^{(}->}[u] & \AAlg^{\iy^{\vphantom{(}}}_{R,5} \ar@{^{(}->}[u]  \\
\vdots_{\vphantom{(}} \ar@{^{(}->}[u] & \vdots_{\vphantom{(}} \ar@{^{(}->}[u] & \vdots_{\vphantom{(}} \ar@{^{(}->}[u] & \vdots_{\vphantom{(}} \ar@{^{(}->}[u] \\
\AAlg^{\un,\iy^{\vphantom{(}}}_{R,\iy}\, \ar@{=}[r] \ar@{^{(}->}[u] & \dAAlg^{\un,\iy^{\vphantom{(}}}_{R,\iy}\, \ar@{^{(}->}[r] \ar@{^{(}->}[u] & \AAlg^{\cu,\iy^{\vphantom{(}}}_{R,\iy}\,  \ar@{^{(}->}[r] \ar@{^{(}->}[u] & \AAlg^{\iy^{\vphantom{(}}}_{R,\iy} \ar@{^{(}->}[u] \\
\sN(\AAlg^\un_R)\, \ar@{=}[r] \ar@{^{(}->}[u] & \sN(\AAlg^\un_R)\, \ar@{^{(}->}[r] \ar@{^{(}->}[u] & \sN(\AAlg^\cu_R)\,  \ar@{^{(}->}[r] \ar@{^{(}->}[u] & \sN(\AAlg_R).\! \ar@{^{(}->}[u] }
\nonumber
\ea
We comment on these:
\begin{itemize}
\setlength{\itemsep}{0pt}
\setlength{\parsep}{0pt}
\item[(i)] The $\iy$-categories in the first and second columns, with superscripts `$\un$', have objects unital $A_\iy$-algebras and 1-morphisms unital $A_\iy$-morphisms.
\item[(ii)] The $\iy$-categories in the third column, with superscript `$\cu$', have objects cohomologically unital $A_\iy$-algebras and 1-morphisms cohomologically unital $A_\iy$-morphisms.
\item[(iii)] The $\iy$-categories in the fourth column have objects all $A_\iy$-algebras, and 1-morphisms all $A_\iy$-morphisms.
\item[(iv)] The $\iy$-categories in the first and second columns, linked by `$=$', are strictly isomorphic, but with definitions written in a different way. See Theorem \ref{ai6thm6} for the isomorphism $I^\un_{R,K}:\AAlg^{\un,\iy}_{R,K}\ra\dAAlg^{\un,\iy}_{R,K}$.
\item[(v)] In each row we have subscripts `$R,K$' for $K=2,3,4,\ldots,\iy$, where we call $K$ the {\it type\/} of the $\iy$-category.
\item[(vi)] $\sN(\AAlg^\un_R),\ldots,\sN(\AAlg_R)$ are the nerves of the ordinary categories $\AAlg^\un_R,\ldots,\AAlg_R$ from Definition~\ref{ai3def5}.
\end{itemize}

Here is what we mean by `type'. Recall that $A_\iy$-algebras, $A_\iy$-morphisms, and $A_\iy$-homotopies  in \S\ref{ai3}, and $A_\iy$-categories and $A_\iy$-functors in \S\ref{ai41}--\S\ref{ai42}, involved families of maps $\m_k,\f_k,\H_k$ for $k\ge 1$. However, the $A_\iy$-category $\Fun_{A_\iy}(\sA,\sB)$ in \S\ref{ai43}, $A_\iy$-natural transformations in \S\ref{ai44}, and filtered $A_\iy$-algebras and $A_\iy$-categories in Chapter \ref{ai5}, involved families of maps $\fT_k,\m_k,\f_k$ for~$k\ge 0$. 

The $N$-simplices in the $\iy$-categories in \eq{ai6eq1} all involve families of maps $\m_k^I$ for $\es\ne I\subseteq\{0,\ldots,N\}$. For $\iy$-categories of type $K$, the rule is that the maps $\m_k^I$ are defined for $k\ge 1$ if $\md{I}\le K$, and for $k\ge 0$ if $\md{I}>K$, that is, we set $\m_0^I=0$ if $\md{I}\le K$. In effect this means that $N$-morphisms in $\iy$-categories of type $K$ include maps $\m_k$ for $k\ge 1$ if $n<K$, and maps $\m_k$ for $k\ge 0$ if~$N\ge K$.

When $N=2$, in the $\iy$-categories $\AAlg^{\un,\iy}_R$ of type 2 in the first row of \eq{ai6eq1}, the 2-morphisms correspond to ({\it cohomologically\/}) {\it unital\/ $A_\iy$-natural isomorphisms\/} in Definition \ref{ai4def13}, regarding $A_\iy$-algebras as $A_\iy$-categories with one object. But in the $\iy$-categories of type $K\ge 3$ in the subsequent rows, the 2-morphisms correspond to ({\it unital\/}) {\it homotopies\/} of $A_\iy$-algebras, as in \S\ref{ai34}. This suggests that $\iy$-categories of type $K\ge 3$ may be better for studying $A_\iy$-algebras, but when we pass to $A_\iy$-categories in Chapter \ref{ai7}, type 2 should be preferred.

We give two different, but related, ways to define $\iy$-categories of $A_\iy$-algebras. The first, in \S\ref{ai63}, gives the first column of \eq{ai6eq1}, and the second, in \S\ref{ai64}, gives the last three columns. The first method works for types $K\ge 2$, but only for unital $A_\iy$-algebras, as the degeneracy maps involve units. The second method only works for types $K\ge 3$, but for general $A_\iy$-algebras.

We postpone most discussion of the motivation for the definitions, and their relation to the literature, to Chapter \ref{ai7}. Here we just note that our definition of $\AAlg^\iy_R$ was a byproduct of an ongoing project with Guillem Cazassus and Alexander Ritter, on a new construction of Fukaya categories of symplectic manifolds $(X,\om)$, which we hope will eventually appear in~\cite{CJR}.

Our definitions of $\AAlg^{\un,\iy}_R,\ldots,\AAlg^\iy_{R,\iy}$ will be written using sums over a class of graphs $\Ga$ called {\it rooted planar trees}. The vertices and edges of $\Ga$ will be labelled by data including $R$-linear graded operators $\m_k^I$. Section \ref{ai61} defines graphical notation suitable for describing (filtered) $A_\iy$-algebras, $A_\iy$-functors, and homotopies, as in \S\ref{ai3} and \S\ref{ai52}--\S\ref{ai54}. Then \S\ref{ai62} extends the notation to {\it simplices\/} of morphisms. The $\iy$-categories in \eq{ai6eq1} are defined in~\S\ref{ai63}--\S\ref{ai64}.

\subsection{Graphical notation for compositions of operators}
\label{ai61}

In working with $A_\iy$-algebras, we meet many formulae involving multiple compositions of the maps $\m_k,\f_k,\g_k,\H_k$ of \S\ref{ai32} and \S\ref{ai34}. We will explain a way to draw pictures of these.

\begin{figure}[htb]
\centering
\includegraphics{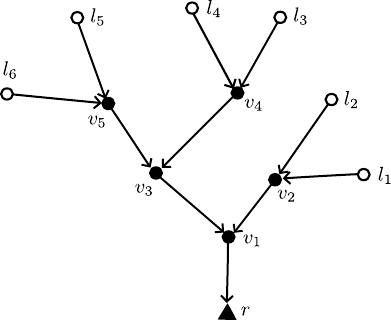}    
\caption{A rooted planar tree}
\label{ai6fig1}
\end{figure}

\begin{dfn}
\label{ai6def1}
We will define a class of finite directed graphs $\Ga$ called {\it rooted planar trees}. (For brevity we often shorten this to {\it tree}.)

\begin{itemize}
\setlength{\itemsep}{0pt}
\setlength{\parsep}{0pt}
\item[(i)] $\Ga$ has finitely many (and at least two) {\it vertices\/} $v$ and finitely many (and at least one) directed {\it edges\/} $e$, which go between two vertices. An edge $e$ starting at vertex $v$ and ending at vertex $w$ is drawn~$\overset{v}{\bu}\,{\buildrel e\over\longra}\,\overset{w}{\bu}$.
\item[(ii)] We require that $\Ga$ is a {\it tree}, that is, $\Ga$ is connected and simply-connected.
\item[(iii)] There is one distinguished vertex $r$ called the {\it root vertex}, drawn `$\tr$'. (This is why we call $\Ga$ {\it rooted}.) There should be exactly one edge $e$ ending at $r$, and no edges starting at~$r$.
\item[(iv)] As $\Ga$ is a tree, there is a unique shortest path in $\Ga$ from any vertex $v$ to the root vertex $r$. We require that all edges in $\Ga$ should be oriented in the direction of the shortest path to~$r$.

This implies that every non-root vertex $v$ has exactly one edge $e$ starting at $v$, and finitely many edges ending at~$v$.
\item[(v)] The vertices $v$ in $\Ga$ which have no edges ending at $v$ are divided into two classes, {\it leaves}, drawn as `$\ci$' and often labelled $l_i$, and {\it buds}, drawn as `$\sq$' and often labelled $b_i$. The distinction between leaves and buds is part of the data of~$\Ga$.

Vertices which are not leaves, buds, or the root vertex, are called {\it internal vertices}, drawn as `$\bu$', and often labelled~$v_i$.
\item[(vi)] When we say that $\Ga$ is {\it planar}, we mean that $\Ga$ is drawn in the plane $\R^2$, and we care about $\Ga$ up to continuous deformations in $\R^2$.

Since $\Ga$ is a tree, being planar is equivalent to the data of a {\it cyclic order\/} of the edges starting and ending at each vertex $v$. That is, if $\Ga$ is drawn in the plane then set of the edges starting and ending at $v$ has a natural cyclic order going anticlockwise around $v$.
\item[(vii)] Let $v$ be an internal vertex, one outgoing edge starting at $v$ and $k\ge 1$ incoming edges ending at $v$. We sometimes write these edges as $e_{v,0},\ab e_{v,1},\ab\ldots,\ab e_{v,k}$, where $e_{v,0}$ is the outgoing edge, and $e_{v,1},\ldots,e_{v,k}$ are the incoming edges, and the edges are in the anticlockwise cyclic order $e_{v,0},e_{v,1},\ldots,e_{v,k}$ around $v$. This gives a numbering $0,1,\ldots,k$ of the edges around $v$.
\item[(viii)] Let there be $n\ge 0$ leaves in $\Ga$. We sometimes write these as $l_1,\ldots,l_n$, such that the root $r$ and the leaves appear in the anticlockwise cyclic order $r,l_1,\ldots,l_n$, going round the outside of $\Ga$ in $\R^2$.

We call $\Ga$ an $n$-{\it leaved\/} rooted planar tree if it has $n$ leaves.
\item[(ix)] Let there be $m\ge 0$ buds in $\Ga$. We sometimes write these as $b_1,\ldots,b_m$, such that the root $r$ and the buds appear in the anticlockwise cyclic order $r,b_1,\ldots,b_m$, going round the outside of $\Ga$ in $\R^2$.
\item[(x)] We sometimes write the $m+n$ leaves and buds together as $u_1,\ldots,u_{m+n}$, such that the root $r$ and the leaves and buds appear in the anticlockwise cyclic order $r,u_1,\ldots,u_{m+n}$, going round the outside of $\Ga$ in $\R^2$. Then $u_1,\ldots,u_{m+n}$ consists of the leaves $l_1,\ldots,l_n$ and the buds $b_1,\ldots,b_m$. For example, if $m=2$ and $n=3$ we might have~$(u_1,\ldots,u_5)=(b_1,l_1,l_2,b_2,l_3)$. 
\end{itemize}
A {\it rooted planar tree\/} (or just {\it tree\/}) $\Ga$ means the data above: a finite set of vertices $V$, a finite set of edges $E$, head and tail maps $h,t:V\ra E$, a division of vertices with no incoming edges into leaves and buds, a distinguished root $r\in V$, and cyclic orders of edges at each internal vertex, satisfying the conditions above. We illustrate a 6-leaved rooted planar tree in Figure~\ref{ai6fig1}.
\end{dfn}

\begin{figure}[htb]
\centering
\includegraphics{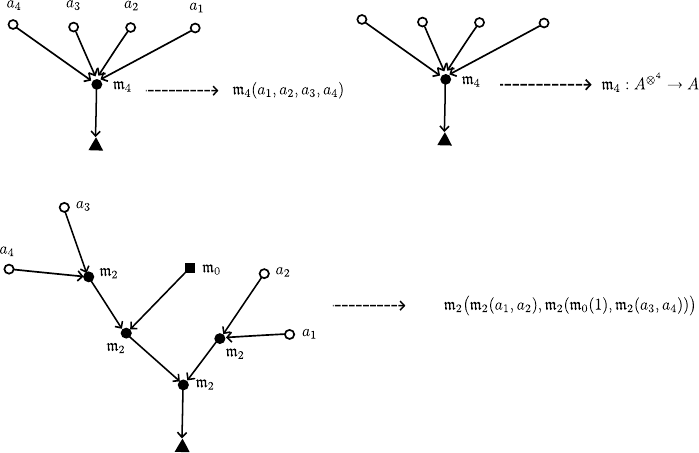}    
\caption{Examples of multiple compositions from labelled trees}
\label{ai6fig2}
\end{figure}

\begin{dfn}
\label{ai6def2}
Suppose that $R$ is a commutative ring, and $A$ is a graded $R$-module, and $\m_k:A^{\ot^k}\ra A$ are graded $R$-linear maps of some degree for $k\ge 0$ (or possibly $k\ge 1$, but then we just define $\m_0=0$), where when $k=0$ by convention we write $A^{\ot^0}=R$, so that $\m_0:R\ra A$ and often we identify $\m_0$ with $\m_0(1)\in A$. For example, such maps arise in $A_\iy$-algebras in Definition \ref{ai3def1}, and curved $A_\iy$-algebras in Definition \ref{ai5def1}. We can use rooted planar trees $\Ga$ to represent multiple compositions of the $\m_k$ operators as follows:
\begin{itemize}
\setlength{\itemsep}{0pt}
\setlength{\parsep}{0pt}
\item[(i)] Label each internal vertex $v$ with $k\ge 1$ incoming edges by the operator~$\m_k$.
\item[(ii)] Label each bud $b$ with $k$ incoming edges by the element~$\m_0(1)\in A$.
\item[(iii)] To get an expression like $\m_3(a_1,\m_2(a_2,a_3),a_4)$ as in \eq{ai3eq4} involving evaluation on particular elements $a_1,a_2,\ldots$ in $A$, label the leaves $l_i$ of $\Ga$ by elements of $A$. When the leaves are ordered $l_1,\ldots,l_n$ as in Definition \ref{ai6def1}(viii), we usually label $l_i$ by $a_i$ for $a_1,\ldots,a_n$ in $A$.
\item[(iv)] To get an expression like $\m_3(\id_A\ot\m_2\ot\id_A):A^{\ot^4}\ra A$ as in \eq{ai3eq3}, with no evaluation on particular elements of $A$, do not label the leaves $l_i$.
\item[(v)] The root $r$ represents the output, the composition of the whole formula.
\end{itemize}

We then build up the multiple composition by starting at the leaves, which represent the domain $(a_1,\ldots,a_n)\in A^n$ or $A^{\ot^n}$ of the formula, and following the directed edges, applying the operator $\m_k$ at each internal vertex to its $k$ inputs coming along the $k$ incoming edges, using the numbering $e_{v,1},\ldots,e_{v,k}$ of the incoming edges in Definition \ref{ai6def1}(vii) to determine the order of the $k$ inputs into $\m_k$, and substituting $\m_0(1)$ as the output at each bud. Some examples of this are shown in Figure~\ref{ai6fig2}.

We sometimes write bold $\bs\Ga$ to mean a rooted planar tree $\Ga$ together with choices of suitable labels at vertices (and later, also edges) which determine a multiple composition in an $A_\iy$-algebra. Then we call $\bs\Ga$ a {\it labelled tree}. Note that we also use $\bs\Ga$ in \S\ref{ai62} and \S\ref{ai71} to mean $\Ga$ together with other kinds of~labels.

We will write $C(\bs\Ga)$, or $C(\bs\Ga,a_1,\ldots,a_n)$, for these multiple compositions of operators $\m_k$, and elements $a_1,\ldots,a_n$.

In cases when we only have operators $\m_k$ for $k\ge 1$, not $k\ge 0$, such as the definitions of $A_\iy$-algebras and $A_\iy$-categories, we restrict to trees with no buds.

\end{dfn}

\begin{dfn}
\label{ai6def3}
Let $\Ga$ be a rooted planar tree. Define a total order $\pr$ on the internal vertices $v$ of $\Ga$ such that if $v,w$ are internal vertices then $v\pr w$ if~either:
\begin{itemize}
\setlength{\itemsep}{0pt}
\setlength{\parsep}{0pt}
\item[(i)] There is a directed path $\overset{w}{\bu}\,{\buildrel e_1\over\longra}\,\bu\,{\buildrel e_2\over\longra}\,\bu\,\cdots\,\bu\,{\buildrel e_k\over\longra}\,\overset{v}{\bu}$ in $\Ga$ for $k>0$; or
\item[(ii)] The shortest paths from $v,w$ to the root $r$ first intersect at an internal vertex $\hat v$ with $v\ne\hat v\ne w$, and writing the directed paths in $\Ga$ from $v,w$ to $v'$ as $\overset{v}{\bu}\,{\buildrel e_1\over\longra}\,\bu\,{\buildrel e_2\over\longra}\,\bu\,\cdots\,\bu\,{\buildrel e_k\over\longra}\,\overset{\hat v}{\bu}$ and $\overset{w}{\bu}\,{\buildrel e'_1\over\longra}\,\bu\,{\buildrel e'_2\over\longra}\,\bu\,\cdots\,\bu\,{\buildrel e'_{k'}\over\longra}\,\overset{\hat v}{\bu}$ for $k,k'>0$, then the incoming edge $e_k$ at $\hat v$ is smaller than $e'_{k'}$ in the numbering of Definition \ref{ai6def1}(vii). That is, writing $e_k=e_{\hat v,j}$ and $e'_{k'}=e_{\hat v,j'}$, then~$j<j'$.
\end{itemize}
By definition $v\not\pr v$, that is, $\pr$ is a total order like $<$, not like~$\le$.

The importance of $\pr$ is that when we write out a multiple composition $C(\bs\Ga)$ in a form such as $\m_2\ci\bigl((\m_3\ci(\m_1\ot\m_2\ot\m_5))\ot\m_4\bigr)$, the operators $\m_k=\m_k^v$ associated to the internal vertices $v$ appear in the order $\pr$ on vertices. This is important in determining signs in formulae, as in Convention~\ref{ai3conv1}.

We can also extend $\pr$ to both internal vertices $v$ and leaves $l$ of $\Ga$ by the same definition, and then $\pr$ determines the order of both operators $\m_k=\m_k^v$ associated to internal vertices $v$ and also elements $a_i\in A$ associated to leaves $l_i$ of $\Ga$ in expressions such as~$\m_2\bigl(\m_2(\m_1(a_1),\m_2(a_2,a_3)),\m_3(a_4,a_5,a_6)\bigr)$.
\end{dfn}

\begin{dfn}
\label{ai6def4}
Let $\f:(A,\m)\ra(B,\n)$ be a morphism of $A_\iy$-algebras over $R$. Then we can also use labelled trees $\bs\Ga$ to describe multiple compositions of maps $\f_i,\m_j,\n_k$, as in \eq{ai3eq7}. Internal vertices of $\Ga$ can be labelled with operators $\f_i,\m_j$ or $\n_k$. In this case, it can also be helpful to label the {\it edges\/} of $\Ga$ by the $A_\iy$-algebra $A$ or $B$, with the rules:
\begin{itemize}
\setlength{\itemsep}{0pt}
\setlength{\parsep}{0pt}
\item[(i)] If an internal vertex $v$ is labelled by $\m_k$, all edges at $v$ are labelled~$A$. 
\item[(ii)] If an internal vertex $v$ is labelled by $\n_k$, all edges at $v$ are labelled~$B$. 
\item[(iii)] If an internal vertex $v$ is labelled by $\f_k$, all incoming edges at $v$ are labelled $A$, and the outgoing edge is labelled~$B$.
\item[(iv)] Leaves $l_i$ may be labelled by an element of $A$ or $B$ depending on whether the outgoing edge is labelled by $A$ or $B$.
\end{itemize}
These ensure that compositions make sense, that is, the inputs to $\f_i,\m_j,\n_k$ lie in $A$ or $B$ as appropriate.

Given another morphism $\g:(B,\n)\ra(C,{\mathfrak o})$, we can also use labelled trees $\bs\Ga$ to describe multiple compositions of maps $\f_i,\g_j,\m_k,\n_l,{\mathfrak o}_m$, in the obvious way. Also, if $\f,\g:(A,\m)\ra(B,\n)$ be morphisms of $A_\iy$-algebras, and $\H:\f\Ra\g$ is a homotopy, we can use labelled trees $\bs\Ga$ to describe multiple compositions of maps $\f_i,\g_j,\H_k,\m_k,\n_l$.
\end{dfn}

\begin{ex}
\label{ai6ex1}
{\bf(a)} We represent the defining equation \eq{ai3eq3} of $A_\iy$-algebras graphically in Figure \ref{ai6fig3}.
\begin{figure}[htb]
\centering
\includegraphics{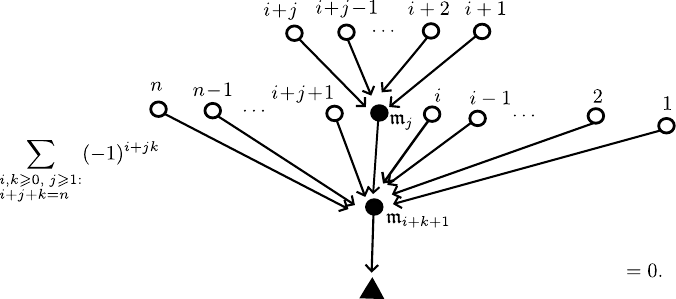}    
\caption{Equation \eq{ai3eq3}: defining $A_\iy$-algebras $(A,\m)$}
\label{ai6fig3}
\end{figure}

\noindent{\bf(b)} We represent the defining equation \eq{ai3eq7} of morphisms of $A_\iy$-algebras graphically in Figure~\ref{ai6fig4}.
\begin{figure}[htb]
\centering
\includegraphics{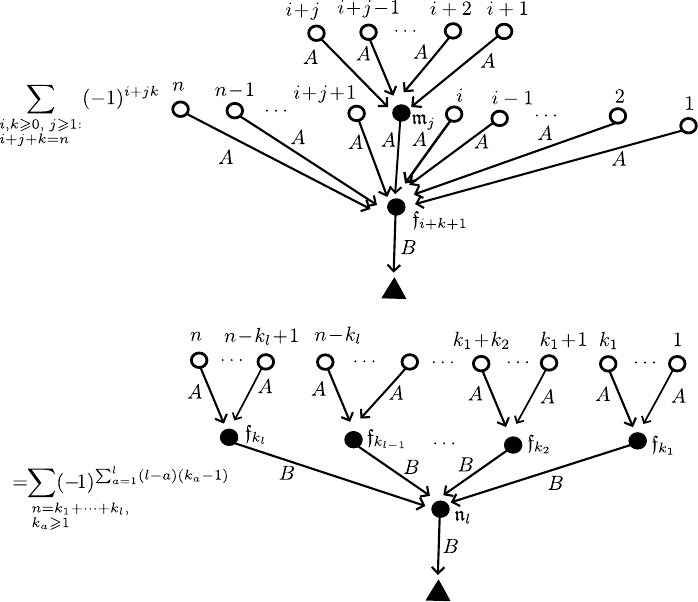}    
\caption{Equation \eq{ai3eq7}: morphisms of $A_\iy$-algebras $\f:(A,\m)\ra(B,\n)$}
\label{ai6fig4}
\end{figure}
\end{ex}

\subsection{Labelling graphs for simplices of morphisms}
\label{ai62}

In our definition of $\iy$-categories of $A_\iy$-algebras in \S\ref{ai63}--\S\ref{ai64}, the $N$-simplices involve graded $R$-modules $A^0,\ldots,A^N$ and $R$-multilinear maps $\m_k^I:(A^{\min I})^{\ot^k}\ab\ra A^{\max I}$ for all subsets $\es\ne I\subseteq\{0,1,\ldots,N\}$ and $k\ge 0$. We extend the graphical notation of \S\ref{ai61} to include multiple compositions of such maps~$\m_k^I$.

\begin{dfn}
\label{ai6def5}
Let $N\ge 0$. A $\{0,\ldots,N\}$-{\it labelled tree\/} $\bs\Ga$ is a rooted planar tree $\Ga$, as in Definition \ref{ai6def1}, together with a nonempty subset $\es\ne I_v\subseteq\{0,\ldots,N\}$ for each vertex $v$ in $\Ga$, satisfying: 
\begin{itemize}
\setlength{\itemsep}{0pt}
\setlength{\parsep}{0pt}
\item[(a)] If $l$ is a leaf of $\Ga$ then $I_l=\{0\}$.
\item[(b)] The root $r$ of $\Ga$ has $I_r=\{N\}$.
\item[(c)] Let $\overset{v}{\bu}\,{\buildrel e\over\longra}\,\overset{w}{\bu}$ be a directed edge in $\Ga$. Then $\max(I_v)=\min(I_w)$.
\end{itemize}

When we draw a $\{0,\ldots,N\}$-labelled tree, we write the subset $I_v$ next to the vertex $v$. We may omit labels of root and leaves as they are automatic by (a)--(b). For each edge $\overset{v}{\bu}\,{\buildrel e\over\longra}\,\overset{w}{\bu}$, we may also label $e$ by $I_e:=\max(I_v)=\min(I_w)$ in $\{0,\ldots,n\}$. Figure \ref{ai6fig5} gives an example of a $\{0,\ldots,5\}$-labelled tree.

\begin{figure}[htb]
\centering
\includegraphics{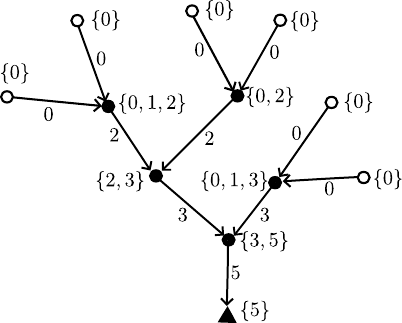}    
\caption{A $\{0,\ldots,5\}$-labelled tree}
\label{ai6fig5}
\end{figure}

We do not need to take the labelling set to be $\{0,1,\ldots,N\}$. For example, we could replace it by $\{a_0,a_1,\ldots,a_N\}$ for $a_0<a_1<\cdots<a_N$ in $\N$, and get a notion of $\{a_0,a_1,\ldots,a_N\}$-{\it labelled tree}, where above we require that $\es\ne I_v\subseteq\{a_0,a_1,\ldots,a_N\}$, and in (a),(b) we take $I_l=\{a_0\}$ and $I_r=\{a_N\}$. Or we could replace $\{0,1,\ldots,N\}$ by any finite totally-ordered set.

Now let $R$ be a commutative ring, and suppose we are given the data:
\begin{itemize}
\setlength{\itemsep}{0pt}
\setlength{\parsep}{0pt}
\item[(i)] $A^i=\bigop_{d\in\Z}A^{i,d}$ is a $\Z$-graded $R$-module for each $i=0,\ldots,N$.
\item[(ii)] For each nonempty subset $I\subseteq\{0,\ldots,N\}$, which we write as $I=\{i_0,\ab i_1,\ab\ldots,\ab i_M\}$ for $i_0<i_1<\cdots<i_M$ with $M\ge 0$, and each $k\ge 0$, we are given an $R$-linear map 
\begin{equation*}
\m^I_k:(A^{i_0})^{\ot^k}\longra A^{i_M},
\end{equation*}
graded of degree $2-M-k$. Here when $k=0$, by convention $(A^{i_0})^{\ot^0}=R$, and we sometimes identify $\m_0^I:R\ra A^{i_M}$ with $\m_0^I(1)\in A^{i_M}$. Equivalently, we may regard $\m^I_k$ as an $R$-multilinear map $\m^I_k:{\buildrel {\ulcorner\,\,\,\text{$k$ copies } \,\,\,\urcorner} \over {A^{i_0}\t\cdots\t A^{i_0}}}\longra A^{i_M}$.
\end{itemize}

If $\bs\Ga$ is an $n$-leaved $\{0,\ldots,N\}$-labelled tree, we define an $R$-linear map
\e
C(\bs\Ga):(A^0)^{\ot^n}\longra A^N
\label{ai6eq2}
\e
as in Definition \ref{ai6def2}, by associating the operator $\m^{I_v}_k$ to each internal vertex $v$ with $k\ge 1$ incoming edges, and $\m^{I_b}_0(1)$ to each bud $b$, and taking the composition over the diagram. (Note that we do {\it not\/} associate operators $\m^{\{0\}}_0$ with the leaves, or $\m^{\{N\}}_1$ with the root.) Conditions (a)--(c) above ensure that this composition is well-defined, that is, that the inputs to each $\m_k^I$ in the formula do lie in $A^{i_0}$ for the correct $i_0$, and the output lies in $A^N$. So, for example, Figure \ref{ai6fig5} yields the $R$-linear map
\begin{equation*}
\m_2^{\{3,5\}}\bigl(\m_2^{\{0,1,3\}}\ot\m_2^{\{2,3\}}(\m_2^{\{0,2\}}\ot\m_2^{\{0,1,2\}})\bigr):(A^0)^{\ot^6}\longra A^5.
\end{equation*}
If we are given elements $a_1,\ldots,a_n\in A^0$ we can insert these as labels at the leaves $l_1,\ldots,l_n$, and form $C(\bs\Ga)(a_1\ot\cdots\ot a_n)\in A^N$.

Since $\m_k^I$ has degree $3-\md{I}-k$, an easy calculation shows $C(\bs\Ga)$ has degree
\e
\deg C(\bs\Ga)=1-n+\sum_{\text{$v$ internal vertex or bud of $\Ga$}}(2-\md{I_v}).
\label{ai6eq3}
\e
\end{dfn}

{\it Balanced\/} $\{0,\ldots,N\}$-labelled trees will be important in our theory, as our definition of $\AAlg^\iy_R$ involves sums over balanced trees.

\begin{dfn}
\label{ai6def6}
Let $\bs\Ga$ be a $\{0,\ldots,N\}$-labelled tree with $n\ge 0$ leaves and $m\ge 0$ buds. As in Definition \ref{ai6def1}(x), let us write the $m+n$ leaves and buds together as $u_1,\ldots,u_{m+n}$, such that the root $r$ and the leaves and buds appear in the anticlockwise cyclic order $r,u_1,\ldots,u_{m+n}$, going round the outside of $\Ga$ in $\R^2$. We call $\bs\Ga$ {\it balanced\/} if the following hold:
\begin{itemize}
\setlength{\itemsep}{0pt}
\setlength{\parsep}{0pt}
\item[(a)] For every internal vertex $v$ in $\Ga$ we have $\md{I_v}\ge 2$.
\item[(b)] For every bud $b$ in $\Ga$ we have $\md{I_b}\ge 3$.
\item[(c)] As extra data, each bud $b$ of $\Ga$ should also be labelled by a nonempty subset $P_b$ of $\{0,\ldots,N\}$ called the {\it prefix\/} of $b$, which should satisfy $\min P_b=0$ and $\max P_b=\min I_b$.
\item[(d)] For each $i=1,\ldots,m+n-1$, consider the shortest paths in $\Ga$ from the leaves or buds $u_i$ and $u_{i+1}$ to the root $r$. Let $w$ be the unique internal vertex that lies on both these paths and is closest to $l_i,l_{i+1}$. Write the shortest paths from $u_i,u_{i+1}$ to $w$ in $\Ga$ as $\overset{u_i=v^0}{\ci}\longra\overset{v^1}{\bu}\longra\overset{v^2}{\bu}\cdots\overset{v^k}{\bu}\longra\overset{w}{\bu}$ and $\overset{u_{i+1}=\ti v^0}{\ci}\longra\overset{\ti v^1}{\bu}\longra\overset{\ti v^2}{\bu}\cdots\overset{\ti v^{\ti k}}{\bu}\longra\overset{w}{\bu}$, for~$k,\ti k\ge 0$. 

Write $I_{v^j}=\{i^j_0,i^j_1,\ldots,i^j_{a^j}\}$ for $i^j_0<i^j_1<\cdots<i^j_{a^j}$, $j=0,\ldots,k$, and $I_{\ti v^j}=\{\ti\imath^j_0,\ti\imath^j_1,\ldots,\ti\imath^j_{\ti a^j}\}$ for $\ti\imath^j_0<\ti\imath^j_1<\cdots<\ti\imath^j_{\ti a^j}$, $j=0,\ldots,\ti k$. Then $a^j,\ti a^j\ge 1$ by (i), and Definition \ref{ai6def5}(c) gives
\begin{equation*}
i^j_{a^j}=i^{j+1}_0,\quad \ti\imath^j_{\ti a^j}=\ti\imath^{j+1}_0,\quad i^k_{a^k}=\ti\imath^{\ti k}_{\ti a^{\ti k}}=\min(I_w).
\end{equation*}

If $u_i$ is a bud, write the prefix $P_{u_i}$ as $P_{u_i}=\{p_0,\ldots,p_l\}$ with $0=p_0<p_1<\cdots<p_l=\min I_{u_i}=i^0_0$. Otherwise write $l=-1$.

Similarly, if $u_{i+1}$ is a bud, write $P_{u_{i+1}}=\{\ti p_0,\ldots,\ti p_{\ti l}\}$ with $0=\ti p_0<\ti p_1<\cdots<\ti p_{\ti l}=\min I_{u_{i+1}}=\ti\imath^0_0$. Otherwise write $\ti l=0$ and $\ti p_0=0$.

We require that in strictly increasing subsequences of $\{0,\ldots,N\}$ we have
\end{itemize}
\ea
&(p_0,p_1,\ldots,p_l,i^0_{a^0}=i^1_1,i^1_{a^1}=i^2_0,\ldots, i^k_{a^k})=
\label{ai6eq4}\\
&(\ti p_0,\ti p_1,\ldots,\ti p_{\ti l}\!=\!\ti\imath_0^0,\ti\imath^0_1,\ti\imath^0_2,\ldots,\ti\imath^0_{\ti a^0}\!=\!\ti\imath^1_0,\ti\imath^1_1,\ldots,\ti\imath^1_{\ti a^1}=\ti\imath^2_0,\ldots,\ti\imath^{\ti k-1}_{\ti a^{\ti k-1}}\!=\!\ti\imath^{\ti k}_0,\ldots, \ti\imath^{\ti k}_{\ti a^{\ti k}}).
\nonumber
\ea

Suppose $\bs\Ga$ is balanced. Write the shortest paths from $u_{m+n},u_1$ to the root $r$ in $\Ga$ as $\overset{u_{m+n}}{\ci}\longra\overset{v^1}{\bu}\longra\overset{v^2}{\bu}\cdots\overset{v^k}{\bu}\longra\overset{r}{\tr}$ and $\overset{u_1}{\ci}\longra\overset{\ti v^1}{\bu}\longra\overset{\ti v^2}{\bu}\cdots\overset{\ti v^{\ti k}}{\bu}\longra\overset{r}{\tr}$, for $k,\ti k\ge 0$. If $u_{m+n}$ is a bud, write the prefix $P_{u_{m+n}}$ as $P_{u_{m+n}}=\{p_0,\ldots,p_l\}$ with $0=p_0<p_1<\cdots<p_l=\min I_{u_i}=i^0_0$. Otherwise write $l=-1$. Similarly, if $u_1$ is a bud, write $P_{u_1}=\{\ti p_0,\ldots,\ti p_{\ti l}\}$ with $0=\ti p_0<\ti p_1<\cdots<\ti p_{\ti l}=\min I_{u_1}=\ti\imath^0_0$. Otherwise write $\ti l=0$ and $\ti p_0=0$. Define the {\it right and left sequences\/} $R(\bs\Ga)$ and $L(\bs\Ga)$ of $\bs\Ga$ to be 
\ea
\begin{split}
R(\bs\Ga)&=\{p_0,p_1,\ldots,p_l,i^0_{a^0}=i^1_1,i^1_{a^1}=i^2_0,\ldots, i^k_{a^k}\},\\
L(\bs\Ga)&=\{\ti p_0,\ti p_1,\ldots,\ti p_{\ti l}\!=\!\ti\imath_0^0,\ti\imath^0_1,\ti\imath^0_2,\ldots,\ti\imath^0_{\ti a^0}\!=\!\ti\imath^1_0,\ti\imath^1_1,\ldots,\ti\imath^1_{\ti a^1}=\ti\imath^2_0,
\end{split}
\label{ai6eq5}\\
&\qquad\qquad\qquad \ldots,\ti\imath^{\ti k-1}_{\ti a^{\ti k-1}}\!=\!\ti\imath^{\ti k}_0,\ldots, \ti\imath^{\ti k}_{\ti a^{\ti k}}\}.
\nonumber
\ea
Using \eq{ai6eq3}, we can show that if $\bs\Ga$ is balanced then
\e
\deg C(\bs\Ga)=1-n+\md{R(\bs\Ga)}-\md{L(\bs\Ga)}.
\label{ai6eq6}
\e

For example, Figure \ref{ai6fig5} is not balanced, as \eq{ai6eq4} for leaves $l_2,l_3$ requires that $\{0,3\}=\{0,2,3\}$, and \eq{ai6eq4} for leaves $l_4,l_5$ requires that $\{0,2\}=\{0,1,2\}$, both of which are false. Figure \ref{ai6fig6} gives an example of a balanced $\{0,\ldots,8\}$-labelled tree, with the values of $R(\bs\Ga)$ and $L(\bs\Ga)$.
\begin{figure}[htb]
\centering
\includegraphics{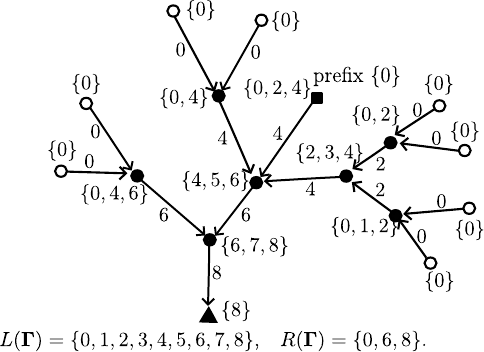}    
\caption{A balanced $\{0,\ldots,8\}$-labelled tree}
\label{ai6fig6}
\end{figure}

Let $\bs\Ga$ be a balanced $n$-leaved $\{0,\ldots,N\}$-labelled tree, with leaves numbered $l_1,\ldots,l_n$ as in Definition \ref{ai6def1}(viii). For $i=1,\ldots,n$, let the shortest path from $l_i$ to the root $r$ in $\Ga$ be $\overset{l_n}{\ci}\,{\buildrel e_1\over\longra}\,\bu\,{\buildrel e_2\over\longra}\,\bu\,\cdots\,\bu\,{\buildrel e_k\over\longra}\,\overset{r}{\tr}$. Then $0=I_{e_1}<I_{e_2}<\cdots<I_{e_k}=N$ by (a), where $I_{e_i}$ is as in Definition \ref{ai6def5}. Define the {\it edge sequence\/} of $l_i$ to be $\ES(l_i)=\{I_{e_1},I_{e_2},\ldots,I_{e_k}\}$. We can show from (d) and \eq{ai6eq5} that
\begin{equation*}
L(\bs\Ga)\supseteq\ES(l_1)\supseteq \ES(l_2)\supseteq \ES(l_3)\supseteq\cdots\supseteq \ES(l_n)\subseteq R(\bs\Ga).
\end{equation*}
In particular, this implies that $R(\bs\Ga)\subseteq \ES(l_i)$ for all $i=1,\ldots,n$.
\end{dfn}

\begin{rem}
\label{ai6rem1}
Condition Definition \ref{ai6def6}(d) probably seems mysterious at this stage. However, in \S\ref{ai71} we extend our notation to maps $\m_k^{I,x_0,\ldots,x_k}$ depending on a sequence of objects $x_0,\ldots,x_k$, as for $A_\iy$-categories in Definition \ref{ai4def1}. Then part (d) will turn out to be necessary to match the domains of inputs and outputs of edges, so that the multiple compositions make sense.
\end{rem}

The next lemma is elementary.

\begin{lem}
\label{ai6lem1}
Let\/ $\bs\Ga$ be a balanced\/ $\{0,\ldots,N\}$-labelled tree. For $v$ an internal vertex or bud in $\Ga,$ let the \begin{bfseries}interior\end{bfseries} $\rint(I_v)$ of\/ $I_v$ be
\begin{equation*}
\rint(I_v)=I_v\sm\{\min(I_v),\max(I_v)\}.
\end{equation*}
Then $\rint(I_v)\subseteq L(\bs\Ga)\sm R(\bs\Ga)$. Furthermore, for each\/ $i$ in $L(\bs\Ga)\sm R(\bs\Ga),$ there is a \begin{bfseries}unique\end{bfseries} internal vertex or bud\/ $v$ of\/ $\Ga$ with\/~$i\in\rint(I_v)$.
\end{lem}

We define an operation of {\it sawing\/} a balanced $\{0,\ldots,N\}$-labelled tree.

\begin{dfn}
\label{ai6def7}
Let $\bs\Ga$ be a balanced $\{0,\ldots,N\}$-labelled tree, with right sequence $R(\bs\Ga)=\{i_0,i_1,\ldots,i_d\}$ with $0=i_0<i_1<\cdots<i_d=N$. Suppose $d>1$, and let $0<c<d$. We will explain how to {\it saw $\bs\Ga$ at level\/}~$i_c$. 

As in Definition \ref{ai6def5}, each edge $e$ of $\Ga$ is labelled by an element $I_e$ in $\{0,\ldots,N\}$. Form a new, disconnected labelled graph $\bs{\ti\Ga}$ by modifying $\bs\Ga$ as follows: for each edge $\overset{v}{\bu}\,{\buildrel e\over\longra}\,\overset{w}{\bu}$ of $\Ga$ with $I_e=i_c$, we `break $e$ in the middle' (or `saw through the branch $e$'), replacing $e$ by two new edges $\overset{v}{\bu}\,{\buildrel e'\over\longra}\,\overset{r_e}{\bu}$ and $\overset{l_e}{\bu}\,{\buildrel e''\over\longra}\,\overset{w}{\bu}$, with two new vertices $r_e$ and $l_e$. Here $r_e$ is the root of a component of $\ti\Ga$, and $l_e$ is a leaf of another component of $\ti\Ga$. We label $r_e$ and $l_e$ by~$I_{r_e}=I_{l_e}=\{i_c\}$.

As $R(\bs\Ga)\subseteq \ES(l_i)$, the shortest path in $\Ga$ from any leaf $l_i$ to the root $r$ contains a unique edge $e$ with $I_e=i_c$. So cutting $\Ga$ into $\ti\Ga$ disconnects $r$ from every leaf $l_i$ of~$\Ga$.

Write $\ti\Ga_0$ for the connected component of $\ti\Ga$ containing $r$. Then $\bs{\ti\Ga_0}$ is an $\{i_c,\ldots,N\}$-labelled tree, with leaves $l_e$ for each edge $e$ in $\Ga$ with $I_e=i_c$. Number this set of edges $e$ in $\Ga$ as $e_1,\ldots,e_k$, such that the leaves of $\ti\Ga_0$ in the order in Definition \ref{ai6def1}(viii) are $l_{e_1},l_{e_2},\ldots,l_{e_k}$. Define $\ti\Ga_i$ to be the connected component of $\ti\Ga$ containing $r_{e_i}$ for $i=1,\ldots,k$. Write $\bs{\ti\Ga}_i$ for the labelled tree corresponding to $\ti\Ga_i$ for $i=0,\ldots,k$. 

To make $\bs{\ti\Ga}_0,\ldots,\bs{\ti\Ga}_k$ into {\it balanced\/} labelled trees, we must explain how to define the extra data of the {\it prefixes\/} $\ti P_b$ of a bud $b$ in some $\ti\Ga_i$ from Definition \ref{ai6def6}(c), from the prefix $P_b$ of the corresponding bud $b$ in $\Ga$. If $b$ lies in $\ti\Ga_i$ for $i=1,\ldots,k$, we just set $\ti P_b=P_b$. If $b$ lies in $\ti\Ga_0$, write $P_b=\{j_0,\ldots,j_l\}$ for $0=j_0<j_1<\cdots<j_l=\min P_v$. It is easy to see that $I_e=i_c=j_m$ for some $m=0,\ldots,l$. Then we define the prefix of $b$ in $\ti\Ga_0$ to be $\ti P_b=\{j_m,\ldots,j_l\}$.

We can now show that:
\begin{itemize}
\setlength{\itemsep}{0pt}
\setlength{\parsep}{0pt}
\item[(i)] $\bs{\ti\Ga}$ is the disjoint union of $\bs{\ti\Ga}_0,\bs{\ti\Ga}_1,\ldots,\bs{\ti\Ga}_k$.
\item[(ii)] $\bs{\ti\Ga}_0$ is a balanced $\{i_c,\ldots,N\}$-labelled tree, with prefixes as above.
\item[(iii)] $\bs{\ti\Ga}_1,\ldots,\bs{\ti\Ga}_k$ are balanced $\{0,\ldots,i_c\}$-labelled trees, with prefixes as above.
\item[(iv)] Write $L(\bs\Ga)=\{\ti\imath_0,\ti\imath_1,\ldots,\ti\imath_{\ti d}\}$ with $0=\ti\imath_0<\ti\imath_1<\cdots<\ti\imath_d=N$ and $\ti\imath_{\ti c}=i_c$. Then
\end{itemize}
\begin{align*}
L(\bs{\ti\Ga}_0)&=\{\ti\imath_{\ti c},\ti\imath_1,\ldots,\ti\imath_{\ti d}\}, & R(\bs{\ti\Ga}_0)&=\{i_c,i_{c+1},\ldots,i_d\},\\
L(\bs{\ti\Ga}_1)&=\{\ti\imath_0,\ti\imath_1,\ldots,\ti\imath_{\ti c}\}, & R(\bs{\bs\Ga}_k)&=\{i_0,i_1,\ldots,i_c\},\\
R(\bs{\ti\Ga}_i)&=L(\bs{\ti\Ga}_{i+1}),\quad i=1,2,\ldots,k-1.\!\!\!\!\!\!\!\!\!\!
\end{align*}

Note that as the added vertices are leaves and roots, the set of internal vertices of $\Ga$ is the disjoint union of the sets of internal vertices of $\ti\Ga_0,\ldots,\ti\Ga_k$, and the set of buds of $\Ga$ is the disjoint union of the sets of buds of~$\ti\Ga_0,\ldots,\ti\Ga_k$.

We claim that the compositions satisfy
\e
\begin{split}
C(\bs\Ga)=\,&\biggl(\prod_{\begin{subarray}{l} \text{$v$ an internal vertex or bud of $\ti\Ga_i$, $i>0$,} \\ \text{$w$ an internal vertex or bud of $\ti\Ga_0$,} \\ \text{$v\pr w$ in order $\pr$ of $\Ga$,} \\
\text{$v$ has $l$ incoming edges,} \\ \text{$w$ has $m$ incoming edges} \end{subarray}} 
(-1)^{(3-\md{I_v}-l)(3-\md{I_w}-m)}\biggr) \\
&{}\cdot C(\bs{\ti\Ga}_0)\ci\bigl(C(\bs{\ti\Ga}_1)\ot C(\bs{\ti\Ga}_2)\ot\cdots\ot C(\bs{\ti\Ga}_k)\bigr),
\end{split}
\label{ai6eq7}
\e
where the total order $\pr$ on the internal vertices of $\Ga$ is as in Definition \ref{ai6def3}. To see this, note that in $C(\bs\Ga)$, the operators $\m_{k_v}^{I_v}$ for the internal vertices $v$ of $\Ga$ appear in order $\pr$. In $C(\bs{\ti\Ga}_0)\ci\bigl(C(\bs{\ti\Ga}_1)\ot C(\bs{\ti\Ga}_2)\ot\cdots\ot C(\bs{\ti\Ga}_k)\bigr)$, the operators $\m_{k_v}^{I_v}$ appear firstly in the order of the tree $\ti\Ga_0,\ti\Ga_1,\ldots,\ti\Ga_k$ containing $v$, and then for a fixed $\ti\Ga_i$ in the order $\pr$ restricted to $\ti\Ga_i$. Thus, two vertices $v,w$ appear in a different order on the two sides of \eq{ai6eq7} if and only if $v$ lies in $\ti\Ga_i$ for $i>0$, and $w$ lies in $\ti\Ga_0$, and $v\pr w$ (or vice versa). As the degrees of $\m_{k_v}^{I_v}$ and $\m_{k_w}^{I_w}$ are $3-\md{I_v}-l$, $3-\md{I_w}-m$, the factor $(-1)^{(3-\md{I_v}-l)(3-\md{I_w}-m)}$ is that prescribed in \eq{ai3eq2} for exchanging the order of $\m_{k_v}^{I_v}$ and~$\m_{k_w}^{I_w}$.

Figure \ref{ai6fig7} shows the result of sawing Figure \ref{ai6fig6} at level~6.
\begin{figure}[htb]
\centering
\includegraphics{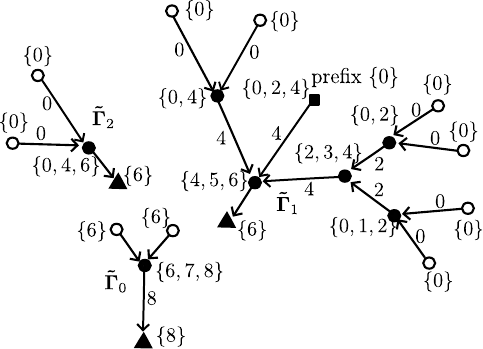}    
\caption{Sawing Figure \ref{ai6fig6} at level 6}
\label{ai6fig7}
\end{figure}
\end{dfn}

\begin{dfn}
\label{ai6def8}
Let $R$ be a commutative ring, and suppose as in Definition \ref{ai6def5} that we are given the data:
\begin{itemize}
\setlength{\itemsep}{0pt}
\setlength{\parsep}{0pt}
\item[(i)] $A^i=\bigop_{d\in\Z}A^{i,d}$ is a $\Z$-graded $R$-module for each $i=0,\ldots,N$.
\item[(ii)] For each nonempty subset $I\subseteq\{0,\ldots,N\}$, which we write as $I=\{i_0,\ab i_1,\ab\ldots,\ab i_M\}$ for $i_0<i_1<\cdots<i_M$ with $M\ge 0$, and each $k\ge 0$, we are given an $R$-linear map $\m^I_k:(A^{i_0})^{\ot^k}\ra A^{i_M}$, graded of degree $2-M-k=3-\md{I}-k$.

Actually, because of Definition \ref{ai6def6}(b), when $k=0$ we only use maps $\m^I_0$ for $\md{I}\ge 3$, so we need not define $\m^I_0$ when $\md{I}=1,2$. It is convenient to set $\m^I_0=0$ when~$\md{I}=1,2$.
\end{itemize}

For each nonempty subset $I=\{i_0,i_1,\ldots,\ab i_M\}\subseteq\{0,\ldots,N\}$ as in (ii) and all $k,n\ge 0$, define an $R$-linear map
\e
\d_{n,k}^I:(A^{i_0})^{\ot^n}\longra(A^{i_M})^{\ot^k}
\label{ai6eq8}
\e
as follows:
\begin{itemize}
\setlength{\itemsep}{0pt}
\setlength{\parsep}{0pt}
\item[(a)] Noting that by convention $(A^i)^{\ot^0}=R$, when $k=0$ set
\e
\d_{n,0}^I=\begin{cases} 
\id_R:(A^{i_0})^{\ot^0}\longra (A^{i_M})^{\ot^0}, & \text{$n=0$ and $M=1$,} \\
0, & \text{otherwise.}
\end{cases}
\label{ai6eq9}
\e
\item[(b)] If $M=0$, $I=\{i_0\}$, and $k\ge 1$, set 
\ea
&\d_{n,k}^{\{i_0\}}=
\label{ai6eq10}\\
&\begin{cases} \sum_{j=0}^{k-1}(-1)^{j+(n-k+1)(k-j-1)}
\id_{A^{i_0}}^{\ot^j}\ot\m_{n-k+1}^{\{i_0\}}\ot \id_{A^{i_0}}^{\ot^{k-j-1}}, & n\ge k, \\
0, & n<k. \end{cases}
\nonumber
\ea
\item[(c)] If $M>0$ and $k\ge 1$, set 
\end{itemize}
\e
\d_{n,k}^I=
\sum_{\begin{subarray}{l} \text{isomorphism classes} \\ \text{of $\bs\Ga_1,\ldots,\bs\Ga_k$:}\\ \text{$\bs\Ga_i$ is a balanced} \\ 
\text{$I$-labelled tree} \\ \text{with $n_i\ge 0$ leaves,} \\ \text{$1\le i\le k$,} \\
\text{$n=n_1+\cdots+n_k,$} \\ \text{$L(\bs\Ga_1)=I$,} \\ \text{$R(\bs\Ga_k)=\{i_0,i_M\}$,} \\
\text{$R(\bs\Ga_i)=L(\bs\Ga_{i+1})$,} \\ \text{$i=1,\ldots,k-1$}\end{subarray}}
\begin{aligned}[t]
&(-1)^{(k+1)(M+n+1)}\\
&{}\cdot\biggl(\prod_{\begin{subarray}{l} \text{$v$ internal vertex or bud of\/ $\bs\Ga_i,$ $1\le i\le k,$} \\ \text{$w$ internal vertex or bud of\/ $\bs\Ga_{i'},$ $1\le i'\le k,$} \\ 
\text{either $i<i',$ or $i=i'$ and} \\ \text{$v\pr w$ in order $\pr$ of $\Ga_i$,} \\ \text{$v$ has $p$ incoming edges including} \\ \text{$q$ leaves, $w$ has $r$ incoming edges}\end{subarray}\!\!\!\!\!\!\!\!\!\!\!\!\!\!\!\!\!\!\!\!\!\!\!\!\!\!\!\!\!\!\!\!\!\!\!\!\!\!\!\!\!\!\!\!\!\!\!\!\!\!} (-1)^{\begin{subarray}{l} N(v,w)+q(r+1) \\ +(p+1)(\md{I_w}+r+1)\end{subarray}}\biggr)\\
&{}\cdot C(\bs\Ga_1)\ot C(\bs\Ga_2)\ot\cdots\ot C(\bs\Ga_k),
\end{aligned}
\label{ai6eq11}
\e
\begin{itemize}
\setlength{\itemsep}{0pt}
\setlength{\parsep}{0pt}
\item[] where we define
\e
N(v,w)=\bmd{\bigl\{(r,s)\in\rint(I_v)\t\rint(I_w):r<s\bigr\}}.
\label{ai6eq12}
\e
with $\rint(I_v)=I_v\sm\{\min(I_v),\max(I_v)\}$ as in Lemma \ref{ai6lem1}.
\end{itemize}
It follows from \eq{ai6eq6} and \eq{ai6eq9}--\eq{ai6eq11} that $\d_{n,k}^I$ is graded of degree
\begin{equation*}
\deg\d_{n,k}^I=k-n+2-\md{I}.
\end{equation*}
\end{dfn}

\begin{rem}
\label{ai6rem2}
The conditions Definition \ref{ai6def6}(a),(b) in the definition of {\it balanced\/} labelled trees $\bs\Ga$, which say that $\md{I_v}\ge 2$ if $v$ is an internal vertex and $\md{I_b}\ge 3$ if $b$ is a bud, are both necessary to ensure that \eq{ai6eq11} is a finite sum.
\end{rem}

\begin{ex}
\label{ai6ex2}
{\bf(i)} If $k=1$ then for any $I,n$ we have
\e
\d^I_{n,1}=\begin{cases} \m_n^I, & \text{$n\ge 1$ or $\md{I}\ge 3$,} \\
0, & \text{$n=0$ and $\md{I}=1,2$.}
\end{cases}
\label{ai6eq13}
\e
This is because for $n\ge 1$ the only balanced $I$-labelled graph $\bs\Ga$ with $L(\bs\Ga)=I$ and $R(\bs\Ga)=\{i_0,i_m\}$ is the graph with only one internal vertex $v$ with $I_v=I$ and $n$ leaves, and all the signs in \eq{ai6eq11} are trivial. If $n=0$ then $v$ is a bud, and we require $\md{I_v}\ge 3$ by Definition~\ref{ai6def6}(b).
\smallskip

\noindent{\bf(ii)} If $I=\{0,1\}$ and $k\ge 1$ then after simplifying signs
\e
\d_{n,k}^{\{0,1\}}=\sum_{n=n_1+\cdots+n_k,\; n_a\ge 1}\begin{aligned}[t]
&(-1)^{\sum_{i=1}^k(n_i+1)(k-i)} \\ 
&\quad \m^{\{0,1\}}_{n_1}\ot\m^{\{0,1\}}_{n_2}\ot\cdots\ot\m^{\{0,1\}}_{n_k}.\end{aligned}
\label{ai6eq14}
\e

\noindent{\bf(iii)} If $I=\{0,1,2\}$ and $k\ge 1$ then after simplifying signs
\begin{equation*}
\d_{n,k}^{\{0,1,2\}}=\sum_{\begin{subarray}{l} 
n=n_1+\cdots+n_k, \\ 1\le l\le k \\
n_i=p^1_i+\cdots+p^{m_i}_i,\\ i=1,\ldots,l-1, \\
p_i^j\ge 1,\; m_i\ge 1, \\
n_i\ge 1,\; i\ne l, \\
n_l\ge 0 \end{subarray}}
\begin{aligned}[t]
&(-1)^{k+1}\cdot(-1)^{\sum_{i=1}^{l-1}(n_i+1)}\cdot{} \\
&(-1)^{\sum_{i=1}^k(n_i+1)(k-i)}\cdot\prod_{i=1}^{l-1}
(-1)^{\sum_{a=1}^{m_i}(p_i^a+1)(m_i-a)}\\
&\m^{\{1,2\}}_{m_1}\bigl(\m^{\{0,1\}}_{p^1_1}\!\ot\!\m^{\{0,1\}}_{p^2_1}\!\ot\!\cdots\!\ot\!\m^{\{0,1\}}_{p^{m_1}_1}\bigr)\\
&\ot\cdots\ot\m^{\{1,2\}}_{m_{l-1}}\bigl(\m^{\{0,1\}}_{p^1_{l-1}}\!\ot\!\m^{\{0,1\}}_{p^2_{l-1}}\!\ot\!\cdots\!\ot\!\m^{\{0,1\}}_{p^{m_{l-1}}_{l-1}}\bigr)\\
&\ot\m^{\{0,1,2\}}_{n_l}\ot\m^{\{0,2\}}_{n_{l+1}}\ot\cdots\ot\m^{\{0,2\}}_{n_k}.
\end{aligned}
\end{equation*}
\end{ex}

\begin{lem}
\label{ai6lem2}
In Definition\/ {\rm\ref{ai6def8},} $\d_{n,k}^I=0$ unless $n\ge \min(k,k+2-\md{I})$.
\end{lem}

\begin{proof}
In \eq{ai6eq9}--\eq{ai6eq10} this is obvious. Suppose $k>n$. For the sum in \eq{ai6eq11} to be nonempty, there must exist $\bs\Ga_1,\ldots,\bs\Ga_k$ as in \eq{ai6eq11}, where at least $k-n$ of the $\bs\Ga_i$ have no leaves, and so must have at least one bud $b$. Then $\md{I_b}\ge 3$ by Definition \ref{ai6def6}(b), so $\md{\rint(I_b)}\ge 1$. But the elements of $\rint(I_b)$ for the $\ge k-n$ buds $b$ are disjoint elements of $I\sm\{i_0,i_M\}$. Thus~$k-n\ge\md{I}-2$.
\end{proof}

The next theorem is proved in \S\ref{ai101}. 

\begin{thm}
\label{ai6thm1}
In Definition\/ {\rm\ref{ai6def8},} let\/ $k=k'+k''$ with\/ $k',k''\ge 1$. Then
\begin{itemize}
\setlength{\itemsep}{0pt}
\setlength{\parsep}{0pt}
\item[{\bf(a)}] If\/ $M=0$ and\/ $n\ge k$ then
\e
\d_{n,k}^{\{i_0\}}=(-1)^{k''(k'+n)}\d_{n-k'',k'}^{\{i_0\}}\ot\id_{A^{i_0}}^{\ot^{k''}}+
(-1)^{k'}\id_{A^{i_0}}^{\ot^{k'}}\ot\d_{n-k',k''}^{\{i_0\}}.
\label{ai6eq15}
\e
If\/ $M=0$ and\/ $n<k$ then $\d_{n,k}^{\{i_0\}}=0$ by \eq{ai6eq10}.
\item[{\bf(b)}] If\/ $M>0$ then
\end{itemize}
\ea
&\d_{n,k}^I
\label{ai6eq16}\\
&=\!\!\!\sum_{\begin{subarray}{l} 
n=n'+n'', \; n',n''\ge 0, \\ 
0<l\le M,\; 0=j_0<j_1<\cdots<j_l=M, \\
k_0,\ldots,k_l\ge 0, \; k_0=n', \; k_l=k', \\
k_{a-1}\ge k_a+1-j_a+j_{a-1}, \; a=1,\ldots,l
\end{subarray}}\begin{aligned}[t]&
(-1)^{k''(k'+l+M)+n'(k''+l+1)+\sum_{a=1}^{l-1}(j_a-a)}\\
&\bigl(\d_{k_{l-1},k_l}^{\{i_{j_{l-1}},\ldots,i_{j_l}\}}\ci\d_{k_{l-2},k_{l-1}}^{\{i_{j_{l-2}},\ldots,i_{j_{l-1}}\}}\ci\cdots\ci\d_{k_0,k_1}^{\{i_{j_0},\ldots,i_{j_1}\}}\bigr)\\
&{}\ot\d_{n'',k''}^{\{i_{j_0},i_{j_1},i_{j_2},\ldots,i_{j_l}\}},
\end{aligned}
\nonumber
\ea
where the sum is finite. Here the condition $k_{a-1}\ge k_a+1-j_a+j_{a-1}$ is unnecessary, as if it does not hold then $\d_{k_{a-1},k_a}^{\{i_{j_{a-1}},\ldots,i_{j_a}\}}=0$ by Lemma\/ {\rm\ref{ai6lem2},} but without it the sum\/ \eq{ai6eq16} would contain infinitely many zero terms.
\end{thm}

\subsection{\texorpdfstring{The $\iy$-categories $\AAlg^{\un,\iy}_R,\AAlg^{\un,\iy}_{R,K}$}{The ∞-categories A∞-Algᵘⁿ∞ᵣ,A∞-Algᵘⁿ∞ᵣₖ}}
\label{ai63}

We now define an $\iy$-category $\AAlg^{\un,\iy}_R$ of unital $A_\iy$-algebras. We first define a simplicial set $\hAAlg^{\un,\iy}_R$ which is {\it not\/} an $\iy$-category, and then define our $\iy$-category as a simplicial subset~$\AAlg^{\un,\iy}_R\subset\hAAlg^{\un,\iy}_R$. Theorem \ref{ai6thm4} extends $\hAAlg^{\un,\iy}_R$ to an $\iy$-{\it bicategory\/} $(\hAAlg^{\un,\iy}_R,\cT)$ in the sense of \S\ref{ai28}, with~$\AAlg^{\un,\iy}_R=\Core(\hAAlg^{\un,\iy}_R,\cT)$.

\begin{dfn}
\label{ai6def9}
Use the notation of Chapter \ref{ai2} and \S\ref{ai61}--\S\ref{ai62}. Let $R$ be a commutative ring. We will define a simplicial set $\hAAlg^{\un,\iy}_R$ of unital $A_\iy$-algebras over~$R$.

\smallskip

\noindent{\bf(a)} For $N\ge 0$, define the $N$-simplices $(\hAAlg^{\un,\iy}_R)_N$ in $\hAAlg^{\un,\iy}_R$ to be the set of pairs $\bigl((A^i)_{i=0}^N,(\m^I)_{\es\ne I\subseteq\{0,\ldots,N\}}\bigr)$, where:
\begin{itemize}
\setlength{\itemsep}{0pt}
\setlength{\parsep}{0pt}
\item[(i)] $A^i=\bigop_{d\in\Z}A^{i,d}$ is a $\Z$-graded $R$-module for each $i=0,\ldots,N$.
\item[(ii)] For each nonempty subset $I\subseteq\{0,\ldots,N\}$, which we write as $I=\{i_0,\ab i_1,\ab\ldots,\ab i_M\}$ for $i_0<i_1<\cdots<i_M$ with $M\ge 0$, we have $\m^I=(\m_k^I)_{k\ge 0}$, where $\m_k^I$ is an $R$-linear map 
\e
\m^I_k:(A^{i_0})^{\ot^k}\longra A^{i_M},
\label{ai6eq17}
\e
graded of degree $2-M-k=3-\md{I}-k$. Here when $k=0$, by convention $(A^{i_0})^{\ot^0}=R$, and we sometimes identify $\m_0^I:R\ra A^{i_M}$ with~$\m_0^I(1)\in A^{i_M}$.

In fact when $k=0$ we only want maps $\m^I_0$ for $\md{I}\ge 3$, so we require that $\m^I_0=0$ when $\md{I}=1,2$. Alternatively, we can say that $\m_k^I$ in \eq{ai6eq17} should only be defined for $k\ge 1$ when $\md{I}=1,2$ and for $k\ge 0$ when $\md{I}\ge 3$, and then we should exclude terms in \eq{ai6eq18} and other equations below involving $\m^I_0$ for~$\md{I}=1,2$.

We require that for all $I=\{i_0<i_1<\cdots<i_M\}\subseteq\{0,\ldots,N\}$ as above and all $n\ge 0$, the maps $\m^J_l$ should satisfy
\e
\begin{aligned}
\sum_{k=1}^n\m_k^{\{i_0\}}\ci\d_{n,k}^{\{i_0\}}&=0 && \text{if $M=0$,} \\
\sum_{k=1}^n\m_k^I\ci\d_{n,k}^{\{i_0\}}+(-1)^M\sum_{k=1}^{n+M-1}\m_k^{\{i_M\}}\ci\d_{n,k}^I&=0 && \text{if $M>0$,}
\end{aligned}
\label{ai6eq18}
\e
where $\d_{n,k}^I$ is defined from the $\m_l^J$ as in \eq{ai6eq8}--\eq{ai6eq11}. Here the upper limits $k=n$ and $k=n+M-1$ in the sums in \eq{ai6eq18} are unnecessary as $\d_{n,k}^{\{i_0\}}=0$ if $k>n$ and $\d_{n,k}^I=0$ when $k>n+M-1$ for $M>0$ by Lemma \ref{ai6lem2}, so we could replace the sums by $\sum_{k=0}^\iy$, but we write it as in \eq{ai6eq18} to get finite sums. The lower limits $k=1$ are as $\m_0^{\{i\}}=0$.
\item[(iii)] As in Example \ref{ai6ex3}(a), parts (i)--(ii) imply that $(A^i,\m^{\{i\}})$ is an $A_\iy$-algebra for each $i=0,\ldots,N$. We require this $A_\iy$-algebra $(A^i,\m^{\{i\}})$ to be {\it unital}, with unique unit $\boo^i\in A^{i,0}$. That is, we have
\e
\begin{split}
&\m_2^{\{i\}}\ci(\boo^i\ot\id_{A^i})=\m_2^{\{i\}}\ci(\id_{A^i}\ot\boo^i)=\id_{A^i},\\
&\m_k^{\{i\}}\ci(\id_{A^i}^{\ot^{l-1}}\ot\boo^i\ot\id_{A^i}^{\ot^{k-l}})=0,\quad \text{all $1\le l\le k$, $k\ne 2$.}
\end{split}
\label{ai6eq19}
\e
\item[(iv)] As in Example \ref{ai6ex3}(b), parts (i)--(iii) imply that $\m^{\{i,j\}}:(A^i,\m^{\{i\}})\ra(A^j,\m^{\{j\}})$ is a morphism of unital $A_\iy$-algebras for all $0\le i<j\le N$. We require this to be a {\it unital\/} morphism. That is, we have
\e
\begin{split}
&\m_1^{\{i,j\}}(\boo^i)=\boo^j,\\
&\m_k^{\{i,j\}}\ci(\id_{A^i}^{\ot^{l-1}}\ot\boo^i\ot\id_{A^i}^{\ot^{k-l}})=0,\quad \text{all $1\le l\le k$, $k\ne 1$.}
\end{split}
\label{ai6eq20}
\e
\item[(v)] For all $I\subseteq\{0,\ldots,N\}$ with $\md{I}\ge 3$, which we write as $I=\{i_0,\ab i_1,\ab\ldots,\ab i_M\}$ for $i_0<i_1<\cdots<i_M$ with $M\ge 2$, we have
\e
\m_k^I\ci(\id_{A^{i_0}}^{\ot^{l-1}}\ot\boo^{i_0}\ot\id_{A^{i_0}}^{\ot^{k-l}})=0, \quad \text{all $1\le l\le k$.}
\label{ai6eq21}
\e
\end{itemize}

\noindent{\bf(b)} If $\al:[\bar N]\ra[N]$ is a morphism in $\De$, that is, $\al:\{0,\ldots,\bar N\}\ra\{0,\ldots,N\}$ is a nondecreasing function, define
\begin{equation*}
\hAAlg^{\un,\iy}_R(\al):(\hAAlg^{\un,\iy}_R)_N\longra(\hAAlg^{\un,\iy}_R)_{\bar N}
\end{equation*}
to map $\bigl((A^i)_{i=0}^N,(\m^I)_{\es\ne I\subseteq\{0,\ldots,N\}}\bigr)$ to $\bigl((A^{\al(i)})_{i=0}^{\bar N},(\bar\m{}^{\bar I})_{\es\ne \bar I\subseteq\{0,\ldots,\bar N\}}\bigr)$, where $\bar\m{}^{\bar I}\ab =(\bar\m{}^{\bar I}_k)_{k\ge 0}$ with
\e
\bar\m_k^{\bar I}=\begin{cases}
\m_k^{\al(\bar I)}, & \text{if $\al\vert_{\bar I}$ is injective,} \\
\id_{A^i}, & \text{if $k\!=\!1$ and $\bar I\!=\!\{\bar\imath_0,\bar\imath_1\}$ with $\bar\imath_0\!\ne\! \bar\imath_1$, $\al(\bar\imath_0)\!=\!\al(\bar\imath_1)\!=\!i$,} \\
\boo^{\bar\imath_2}, & \text{if $k=0$, $\bar I=\{\bar\imath_0\!<\!\bar\imath_1\!<\!\bar\imath_2\}$ and $\al\vert_{\bar I}$ is not injective,} \\
0, & \text{otherwise.}
\end{cases}
\label{ai6eq22}
\e
Here by $\bar\m_0^{\{\bar\imath_0,\bar\imath_1,\bar\imath_2\}}=\boo^{\bar\imath_2}$ we mean that $\bar\m_0^{\{\bar\imath_0,\bar\imath_1,\bar\imath_2\}}:R\ra A^{\bar\imath_2}$ maps $\la\mapsto\la \boo^{\bar\imath_2}$.

Proposition \ref{ai6prop1} shows that $\hAAlg^{\un,\iy}_R(\al)$ does map to $(\hAAlg^{\un,\iy}_R)_{\bar N}$. It is easy to check that if $\be:[\check N]\ra[\bar N]$ is another morphism in $\De$ then
\begin{align*}
\hAAlg^{\un,\iy}_R(\al\ci\be)&=\hAAlg^{\un,\iy}_R(\be)\ci\hAAlg^{\un,\iy}_R(\al):\\
&(\hAAlg^{\un,\iy}_R)_N\longra(\hAAlg^{\un,\iy}_R)_{\check N}.
\end{align*} 
Also $\hAAlg^{\un,\iy}_R(\id_{[N]})=\id_{(\hAAlg^{\un,\iy}_R)_N}$. Hence $\hAAlg^{\un,\iy}_R:\De^\rop\ra\Set$ is a functor, that is, $\hAAlg^{\un,\iy}_R$ is a simplicial set.
\end{dfn}

\begin{rem}
\label{ai6rem3}
{\bf(On Set Theory foundations.)}
Recall from Remark \ref{ai1rem2} that we fix a Grothendieck universe $\sU$, and implicitly assume all sets in our theory are $\sU$-small. Without this convention, the $N$-simplices $(\hAAlg^\iy_R)_N$ above would be proper classes, rather than sets, as they are built out of the proper class of all graded $R$-modules $A$. This would cause Set Theory problems.
\end{rem}

The next proposition will be proved in~\S\ref{ai102}.

\begin{prop}
\label{ai6prop1}
In Definition\/ {\rm\ref{ai6def9}(b),} $\bigl((A^{\al(i)})_{i=0}^{\bar N},(\bar\m{}^{\bar I})_{\es\ne \bar I\subseteq\{0,\ldots,\bar N\}}\bigr)$ lies in $(\hAAlg^{\un,\iy}_R)_{\bar N}$.
\end{prop}

The next theorem will be proved in~\S\ref{ai103}.

\begin{thm}
\label{ai6thm2}
Let\/ $\bigl((A^i)_{i=0}^N,(\m^I)_{I\subseteq\{0,\ldots,N\}}\bigr)$ be an $N$-simplex in $\hAAlg^{\un,\iy}_R$ as in Definition\/ {\rm\ref{ai6def9},} and define maps\/ $\d_{n,k}^I$ from the $\m_k^I$ as in\/ {\rm\eq{ai6eq8}--\eq{ai6eq11}}. Then for each\/ $\es\ne I=\{i_0<i_1<\cdots<i_M\}\subseteq\{0,\ldots,N\}$ and all\/ $k,n\ge 0,$ in $R$-linear maps $(A^{i_0})^{\ot^n}\ra(A^{i_M})^{\ot^k}$ we have
\e
\begin{aligned}
\sum_{l=k}^n\d_{l,k}^{\{i_0\}}\ci\d_{n,l}^{\{i_0\}}&=0 && \text{if\/ $M=0,$} \\
\sum_{l=1}^n\d_{l,k}^I\ci\d_{n,l}^{\{i_0\}}+(-1)^M\sum_{l=k}^{n+M-1}\d_{l,k}^{\{i_M\}}\ci\d_{n,l}^I&=0 && \text{if\/ $M>0$.}
\end{aligned}
\label{ai6eq23}
\e
Here the limits in the sums in \eq{ai6eq23} are unnecessary and could be replaced by $\sum_{l=0}^\iy\cdots,$ as for $l$ outside these limits the summands are zero by {\rm\eq{ai6eq8}--\eq{ai6eq11}} and Lemma\/ {\rm\ref{ai6lem2}}. 
\end{thm}

\begin{ex}
\label{ai6ex3}
{\bf(a)} Let $\bigl((A^0),(\m^{\{0\}})\bigr)$ be a 0-simplex of $\hAAlg^{\un,\iy}_R$. Then \eq{ai6eq10} and \eq{ai6eq18} for $I=\{0\}$ say that
\begin{equation*}
\sum_{i,k\ge 0,\; j\ge 1:\; i+j+k=n}
(-1)^{i+jk}\m^{\{0\}}_{i+k+1}\ci\bigl(\id_{A^0}^{\ot^i}\ot\m^{\{0\}}_j\ot \id_{A^0}^{\ot^k}\bigr)=0.
\end{equation*}
Thus by \eq{ai3eq3}, $(A^0,\m^{\{0\}})$ is an $A_\iy$-algebra over $R$. Definition \ref{ai6def9}(a)(iii) says that $(A^0,\m^{\{0\}})$ is unital.
\smallskip

\noindent{\bf(b)} A 1-simplex $\bigl((A^0,A^1),(\m^{\{0\}},\m^{\{1\}},\m^{\{0,1\}})\bigr)$ of $\hAAlg^{\un,\iy}_R$ consists of unital $A_\iy$-algebras $(A^0,\m^{\{0\}}),(A^1,\m^{\{1\}})$ over $R$ by {\bf(a)}, and $R$-linear maps 
\begin{equation*}
\m^{\{0,1\}}_k:(A^0)^{\ot^k}\longra A^1
\end{equation*}
for all $k\ge 1$, graded of degree $1-k$, satisfying by \eq{ai6eq10}, \eq{ai6eq11}, \eq{ai6eq18} for~$n\ge 1$
\begin{align*}
&\sum_{i,k\ge 0,\; j\ge 1:\; i+j+k=n}(-1)^{i+jk}\m^{\{0,1\}}_{i+k+1}\bigl(
\id_{A^0}^{\ot^i}\ot\m^{\{0\}}_j\ot\id_{A^0}^{\ot^k}\bigr) \\ 
&-\sum_{n=n_1+\cdots+n_k,\; n_a\ge 1\!\!\!\!\!\!\!\!\!\!\!\!\!\!\!\!}(-1)^{\sum_{i=1}^k(n_i+1)(k-i)}\m^{\{1\}}_k\bigl(\m^{\{0,1\}}_{n_1}\ot\m^{\{0,1\}}_{n_2}\ot\cdots\ot\m^{\{0,1\}}_{n_k}\bigr)=0,
\end{align*}
after simplifying signs. By \eq{ai3eq7}, this says $\m^{\{0,1\}}$ is a morphism $(A^0,\m^{\{0\}})\ra(A^1,\m^{\{1\}})$ of $A_\iy$-algebras over $R$. Definition \ref{ai6def9}(a)(iv) says that $\m^{\{0,1\}}$ is unital. Thus 1-simplices of $\hAAlg^{\un,\iy}_R$ are equivalent to unital morphisms of unital $A_\iy$-algebras over~$R$.
\smallskip

\noindent{\bf(c)} A 2-simplex of $\hAAlg^{\un,\iy}_R$ consists of unital $A_\iy$-algebras $(A^0,\m^{\{0\}}),\ab(A^1,\ab\m^{\{1\}}),\ab(A^2,\m^{\{2\}})$ over $R$ as in {\bf(a)}, and unital morphisms $\m^{\{i,j\}}:(A^i,\m^{\{i\}})\ra(A^j,\m^j)$ for $(i,j)=(0,1),(0,2),(1,2)$ as in {\bf(b)}, and $R$-linear maps 
\begin{equation*}
\m^{\{0,1,2\}}_k:(A^0)^{\ot^k}\longra A^2
\end{equation*}
for $k=0,1,2,\ldots,$ graded of degree $-k$, satisfying \eq{ai6eq21} and for all $n\ge 0$
\ea
&\sum_{i,k\ge 0,\; j\ge 1:\; i+j+k=n}
(-1)^{i+jk}\m^{\{0,1,2\}}_{i+k+1}\ci\bigl(\id_{A^0}^{\ot^i}\ot\m^{\{0\}}_j\ot \id_{A^0}^{\ot^k}\bigr)
\label{ai6eq24}\\
&+\sum_{\begin{subarray}{l} 
1\le l\le k\le n+1, \; n=n_1+\cdots+n_k, \\
n_i=p^1_i+\cdots+p^{m_i}_i,\; i=1,\ldots,l-1, \\
p_i^j\ge 1,\; m_i\ge 1,\; n_l\ge 0,\;
n_i\ge 1,\; i\ne l\end{subarray}}
\begin{aligned}[t]&(-1)^{k+1}\cdot(-1)^{\sum_{i=1}^{l-1}(n_i+1)}\cdot{} \\
&(-1)^{\sum_{i=1}^k(n_i+1)(k-i)}\cdot\prod_{i=1}^{l-1}
(-1)^{\sum_{a=1}^{m_i}(p_i^a+1)(m_i-a)}\\
&\m^{\{2\}}_k\bigl[\m^{\{1,2\}}_{m_1}\bigl(\m^{\{0,1\}}_{p^1_1}\!\ot\!\m^{\{0,1\}}_{p^2_1}\!\ot\!\cdots\!\ot\!\m^{\{0,1\}}_{p^{m_1}_1}\bigr)\\
&\ot\cdots\ot\m^{\{1,2\}}_{m_{l-1}}\bigl(\m^{\{0,1\}}_{p^1_{l-1}}\!\ot\!\m^{\{0,1\}}_{p^2_{l-1}}\!\ot\!\cdots\!\ot\!\m^{\{0,1\}}_{p^{m_{l-1}}_{l-1}}\bigr)\\
&\ot\m^{\{0,1,2\}}_{n_l}\ot\m^{\{0,2\}}_{n_{l+1}}\ot\cdots\ot\m^{\{0,2\}}_{n_k}\bigr]=0,
\end{aligned}
\nonumber
\ea
after simplifying signs. If we write $\f=\m^{\{1,2\}}\ci\m^{\{0,1\}}$ and $\g=\m^{\{0,2\}}$ as morphisms $\f,\g:(A^0,\m^{\{0\}})\ra(A^2,\m^{\{2\}})$, then we can rewrite \eq{ai6eq24} as
\ea
&\sum_{i,k\ge 0,\; j\ge 1:\; i+j+k=n}
(-1)^{i+jk}\m^{\{0,1,2\}}_{i+k+1}\ci\bigl(\id_{A^0}^{\ot^i}\ot\m^{\{0\}}_j\ot \id_{A^0}^{\ot^k}\bigr)
\label{ai6eq25}\\
&+\sum_{\begin{subarray}{l} 
1\le l\le k\le n+1, \; n=n_1+\cdots+n_k, \\
n_i=p^1_i+\cdots+p^{m_i}_i,\; i=1,\ldots,l-1, \\
p_i^j\ge 1,\; m_i\ge 1,\; n_l\ge 0,\;
n_i\ge 1,\; i\ne l\end{subarray}}
\begin{aligned}[t]&(-1)^{k+1}\cdot(-1)^{\sum_{i=1}^{l-1}(n_i+1)}\cdot(-1)^{\sum_{i=1}^k(n_i+1)(k-i)}\\
&\m^{\{2\}}_k\bigl[\f_{n_1}\ot\cdots\ot\f_{n_{l-1}}\ot \\
&\m^{\{0,1,2\}}_{n_l}\ot\g_{n_{l+1}}\ot\cdots\ot\g_{n_k}\bigr]=0.
\end{aligned}
\nonumber
\ea
Regarding $(A^0,\m^{\{0\}}),(A^2,\m^{\{2\}})$ as unital $A_\iy$-categories with one object, this is a {\it unital\/ $A_\iy$-natural transformation\/} $\eta:\f\Ra\g$ as in Definition~\ref{ai4def13}.

In particular, if $\bigl((A^i)_{i=0}^2,(\m^I)_{\es\ne I\subseteq\{0,1,2\}}\bigr)$ is a 2-simplex in $\hAAlg^{\un,\iy}_R$ with $\m^{\{0,1,2\}}_0=\boo^2$ and $\m^{\{0,1,2\}}_k=0$ for $k>0$, so that $\eta=\id_{\m^{\{0,2\}}}$, then $\m^{\{0,2\}}=\m^{\{1,2\}}\ci\m^{\{0,1\}}$ as morphisms of $A_\iy$-algebras $(A^0,\m^0)\ab\ra(A^2,\m^2)$. This is how the notion of {\it composition of morphisms\/} of $A_\iy$-algebras in Definition \ref{ai3def5} is encoded in our theory.

When $n=0$, equation \eq{ai6eq24} becomes
\begin{equation*}
\m^{\{2\}}_1\ci\m^{\{0,1,2\}}_0=0:R\longra A^{2,1}.
\end{equation*}
Thus $\m^{\{0,1,2\}}_0(1)\in A^{2,0}$  defines a cohomology class
\e
[\m^{\{0,1,2\}}_0(1)]\in H^0(A^2,\m^{\{2\}}_1).
\label{ai6eq26}
\e
Note that $H^*(A^2,\m^{\{2\}}_1)$ is a graded unital $R$-algebra, as in~\S\ref{ai32}.

When $n=1$, equation \eq{ai6eq24} becomes
\begin{align*}
&\m_1^{\{0,1,2\}}\ci\m_1^{\{0\}}+\m_1^{\{2\}}\ci\m_1^{\{0,1,2\}} \\
&-\m_2^{\{2\}}\bigl((\m_1^{\{1,2\}}\ci\m_1^{\{0,1\}})\ot\m^{\{0,1,2\}}_0(1)\bigr)
+\m_2^{\{2\}}\bigl(\m^{\{0,1,2\}}_0(1)\ot\m_1^{\{0,2\}}\bigr)=0.
\end{align*}
By taking cohomology, we can deduce that
\begin{align*}
&\bigl(H^*(\m_1^{\{1,2\}})\ci H^*(\m_1^{\{0,1\}})\bigr)\cdot [\m^{\{0,1,2\}}_0(1)]
=[\m^{\{0,1,2\}}_0(1)]\cdot H^*(\m_1^{\{0,2\}}) \\
&\text{in maps}\qquad H^*(A^0,\m_1^{\{0\}})\longra H^*(A^2,\m_1^{\{2\}}),
\end{align*}
where we first map $H^*(A^0)\ra H^*(A^2)$, and then multiply on the right or left with $[\m^{\{0,1,2\}}_0(1)]$ in the $R$-algebra $H^*(A^2)$.

\noindent{\bf(d)} For a 3-simplex, equation \eq{ai6eq18} for $I=\{0,1,2,3\}$ and $n=0$ becomes
\begin{align*}
&\m_1^{\{3\}}\ci\m_1^{\{0,1,2,3\}}+\m_2^{\{3\}}\bigl((\m_1^{\{2,3\}}\ci\m_0^{\{0,1,2\}}(1))\ot\m_0^{\{0,2,3\}}(1)\bigr) \\
&-\m_2^{\{3\}}\bigl(\m_0^{\{1,2,3\}}(1)\ot\m_0^{\{0,1,3\}}(1)\bigr)=0.
\end{align*}
By taking cohomology, we see that in the $R$-algebra $H^0(A^3,\m_1^{\{3\}})$ we have
\begin{equation*}
H^*(\m_1^{\{2,3\}})([\m_0^{\{0,1,2\}}(1)])\cdot [\m_0^{\{0,2,3\}}(1)]=
[\m_0^{\{1,2,3\}}(1)]\cdot [\m_0^{\{0,1,3\}}(1)].
\end{equation*}
\end{ex}

\begin{dfn}
\label{ai6def10}
Define a simplicial subset $\AAlg^{\un,\iy}_R\subset\hAAlg^{\un,\iy}_R$, by an $N$-simplex $\bigl((A^i)_{i=0}^N,(\m^I)_{\es\ne I\subseteq\{0,\ldots,N\}}\bigr)$ in $\hAAlg^{\un,\iy}_R$ lies in $\AAlg^{\un,\iy}_R$ if and only if whenever $0\le i_0<i_1<i_2\le N$ then the cohomology class $[\m^{\{i_0,i_1,i_2\}}_0(1)]$ from \eq{ai6eq26} is invertible in the unital $R$-algebra~$H^0(A^{i_2},\m^{\{i_2\}}_1)$.

If $\al:[\bar N]\ra[N]$ is a morphism in $\De$ then $\hAAlg^{\un,\iy}_R(\al)$ in Definition \ref{ai6def9}(b) maps $(\AAlg^{\un,\iy}_R)_N\ra(\AAlg^{\un,\iy}_R)_{\bar N}$, so $\AAlg^{\un,\iy}_R$ is a simplicial set. We also write $\AAlg^{\un,\iy}_R=\AAlg^{\un,\iy}_{R,2}$ to fit with the $\iy$-categories $\AAlg^{\un,\iy}_{R,K}$ in Definition~\ref{ai6def11}.
\end{dfn}

The next theorem is proved in~\S\ref{ai104}.

\begin{thm}
\label{ai6thm3}
In Definition\/ {\rm\ref{ai6def10},} $\AAlg^{\un,\iy}_R$ is an $\iy$-category.
\end{thm}

In \S\ref{ai28}, following Lurie \cite[\S 4]{Luri4}, we defined $\iy$-{\it bicategories}, a kind of $(\iy,2)$-category. The next theorem will be proved in~\S\ref{ai105}.

\begin{thm}
\label{ai6thm4}
In the simplicial set\/ $\hAAlg^{\un,\iy}_R$ in Definition\/ {\rm\ref{ai6def9},} define a $2$-simplex to be \begin{bfseries}thin\end{bfseries} if the cohomology class $[\m^{\{0,1,2\}}_0(1)]$ from \eq{ai6eq26} is invertible in the unital $R$-algebra $H^0(A^2,\m^{\{2\}}_1),$ that is, if it is a $2$-simplex in $\AAlg^{\un,\iy}_R\subset\hACat^{\un,\iy}_R$. Write $\cT\subset\hAAlg^{\un,\iy}_R([2])$ for the set of thin $2$-simplices. Then $(\hAAlg^{\un,\iy}_R,\cT)$ is an \begin{bfseries}$\iy$-bicategory\end{bfseries} in the sense of\/ {\rm\S\ref{ai28}}. We have $\AAlg^{\un,\iy}_R=\Core(\hAAlg^{\un,\iy}_R,\cT),$ as in Definition\/~{\rm\ref{ai2def22}}.
\end{thm}

\begin{rem}
\label{ai6rem4}
Theorem \ref{ai6thm3} follows trivially from Theorem \ref{ai6thm4}, as the core of an $\iy$-bicategory is an $\iy$-category. But we use the proof of Theorem \ref{ai6thm3} in that of Theorem \ref{ai6thm4}, so we prove Theorem \ref{ai6thm3} first.
\end{rem}

\begin{dfn}
\label{ai6def11}
Let $K\in\{3,4,5,\ldots,\iy\}$. Define simplicial subsets
\begin{gather*}
\AAlg^{\un,\iy}_{R,K}\subset\AAlg^{\un,\iy}_R\subset\hAAlg^{\un,\iy}_R
\end{gather*}
such that an $N$-simplex $\bigl((A^i)_{i=0}^N,(\m^I)_{\es\ne I\subseteq\{0,\ldots,N\}}\bigr)$ in $\hAAlg^{\un,\iy}_R$ in Definition \ref{ai6def9}(a) lies in $\AAlg^{\un,\iy}_{R,K}$ if for all $0\le i_0<i_1<i_2\le N$ we have $\m_0^{\{i_0,i_1,i_2\}}:\la\mapsto\la\boo^{i_2}$, and $\m_0^I=0$ for all $\es\ne I\subseteq\{0,\ldots,N\}$ with $\md{I}\le K$ and~$\md{I}\ne 3$. 

This condition is preserved by the maps in Definition \ref{ai6def9}(b), and so defines a simplicial subset $\AAlg^{\un,\iy}_{R,K}\subset\hAAlg^{\un,\iy}_R$. As $[\m_0^{\{i_0,i_1,i_2\}}(1)]=[\boo^{i_2}]$ is invertible in $H^0(A^{i_2},\m_1^{\{i_2\}})$ we see that $\AAlg^{\un,\iy}_{R,K}\subset\AAlg^{\un,\iy}_R$. Clearly $\AAlg^{\un,\iy}_{R,K}\subset\AAlg^{\un,\iy}_{R,L}$ if $L\le K$ in $\{3,4,5,\ldots,\iy\}$, giving the inclusion relations in the first column of \eq{ai6eq1}. Theorem \ref{ai6thm7} will show that $\AAlg^{\un,\iy}_{R,K}$ is an $\iy$-category. We say that $\AAlg^{\un,\iy}_{R,K}$ has {\it type\/}~$K$. 
\end{dfn}

\begin{rem}
\label{ai6rem5} 
It is not difficult to rewrite the definitions of $\hAAlg^{\un,\iy}_R,\ab\AAlg^{\un,\iy}_R,\ab\AAlg^{\un,\iy}_{R,K}$ above in terms of {\it bar complexes}, as in \S\ref{ai36}. Implicitly we do this already, as the operators in the bar complex picture are essentially $\d^I=\coprod_{k,n}\d^I_{n,k}$ as in Definition \ref{ai6def8}, but with degrees shifted by 1, and Theorem \ref{ai6thm2} gives identities satisfied by these $\d^I$. In the bar complex version some of the signs become simpler, for example in \eq{ai6eq11}. But other than this, the author did not find the exercise very illuminating, so we will not write it out here.
\end{rem}

\subsection{\texorpdfstring{Other $\iy$-categories of $A_\iy$-algebras}{Other ∞-categories of A∞-algebras}}
\label{ai64}

We wish to generalize the definition of $\AAlg^{\un,\iy}_{R,K}$ in \S\ref{ai63} for $K\ge 3$ to $A_\iy$-algebras that are not unital, but only cohomologically unital, or non-unital. The main difficulty in doing this is that the degeneracy maps in the simplicial set $\hAAlg^{\un,\iy}_R$ in \eq{ai6eq22} use units $\boo^{\bar\imath_2}$ in~$\bar\m_0^{\{\bar\imath_0,\bar\imath_1,\bar\imath_2\}}$.

Our idea is simple. Suppose $\bigl((A^i)_{i=0}^N,(\m^I)_{\es\ne I\subseteq\{0,\ldots,N\}}\bigr)$ is an $N$-simplex in $\AAlg^{\un,\iy}_{R,K}$ in Definitions \ref{ai6def9} and \ref{ai6def10}. We define new maps $\dot\m_k^I$ by
\e
\dot\m_k^{\{i_0,\ldots,i_M\}}=\begin{cases} \m_0^{\{i_0,i_1,i_2\}}-\boo^{i_2}, & \text{$k=0$ and $M=2$,} \\
\m_k^{\{i_0,\ldots,i_M\}}, & \text{otherwise,} \end{cases}
\label{ai6eq27}
\e
where $\boo^{i_2}$ means the map $R\ra A^{i_2,0}$, $\la\mapsto\la\boo^{i_2}$. Then we rewrite the definition of $\AAlg^{\un,\iy}_{R,K}$ in terms of the maps $\dot\m_k^I$ in Definitions \ref{ai6def12}, to give $\dAAlg^{\un,\iy}_{R,K}$. The result no longer uses units in degeneracy maps, so is suitable for generalizing to the non-unital case.

We will do this only for $\AAlg^{\un,\iy}_{R,K}$ for $K\ge 3$ in Definition \ref{ai6def11}, and not for $\AAlg^{\un,\iy}_R=\AAlg^{\un,\iy}_{R,2}$ in Definition \ref{ai6def10}, so in \eq{ai6eq27} we automatically have $\m_0^{\{i_0,i_1,i_2\}}=\boo^{i_2}$ and $\dot\m_0^{\{i_0,i_1,i_2\}}=0$. In fact the approach of this section would also suffice to define an $\iy$-category $\dAAlg^{\un,\iy}_R=\dAAlg^{\un,\iy}_{R,2}$. However, the generalization of this to $A_\iy$-categories does not work, and as $A_\iy$-categories are our primary interest, we will not bother.

Here is the analogue of Definitions \ref{ai6def9} and~\ref{ai6def11}.

\begin{dfn}
\label{ai6def12}
Let $R$ be a commutative ring and $K\in\{3,4,5,\ldots,\iy\}$. We will define a simplicial set $\dAAlg^{\un,\iy}_{R,K}$ of unital $A_\iy$-algebras over~$R$.
\smallskip

\noindent{\bf(a)} For $N\ge 0$, define the $N$-simplices $(\dAAlg^{\un,\iy}_{R,K})_N$ in $\dAAlg^{\un,\iy}_{R,K}$ to be the set of pairs $\bigl((A^i)_{i=0}^N,(\m^I)_{\es\ne I\subseteq\{0,\ldots,N\}}\bigr)$, where:
\begin{itemize}
\setlength{\itemsep}{0pt}
\setlength{\parsep}{0pt}
\item[(i)] $A^i=\bigop_{d\in\Z}A^{i,d}$ is a $\Z$-graded $R$-module for each $i=0,\ldots,N$.
\item[(ii)] For each nonempty subset $I\subseteq\{0,\ldots,N\}$, which we write as $I=\{i_0,\ab i_1,\ab\ldots,\ab i_M\}$ for $i_0<i_1<\cdots<i_M$ with $M\ge 0$, we have $\m^I=(\m_k^I)_{k\ge 0}$, where $\m_k^I$ is an $R$-linear map $\m^I_k:(A^{i_0})^{\ot^k}\ra A^{i_M}$, graded of degree $2-M-k=3-\md{I}-k$. Here when $k=0$, by convention $(A^{i_0})^{\ot^0}=R$, and we sometimes identify $\m_0^I:R\ra A^{i_M}$ with $\m_0^I(1)\in A^{i_M}$. When $k=0$ we require that $\m_0^I=0$ for all $\es\ne I\subseteq\{0,\ldots,N\}$ with~$\md{I}\le K$.

We require that for all $I=\{i_0<i_1<\cdots<i_M\}\subseteq\{0,\ldots,N\}$ as above and all $n\ge 0$, the maps $\m^J_l$ should satisfy
\e
\sum_{j=0}^M\sum_{0\le k\le\max(n,n+j-1)\!\!\!\!\!\!\!\!\!\!\!\!\!\!\!\!\!\!}(-1)^j\,\m_k^{\{i_j,\ldots,i_M\}}\ci\d_{n,k}^{\{i_0,\ldots,i_j\}}=\sum_{j=1}^{M-1}(-1)^j\m_n^{I\sm\{i_j\}},
\label{ai6eq28}
\e
where $\d_{n,k}^{\{i_0,\ldots,i_j\}}$ is defined from the $\m_k^I$ as in \eq{ai6eq8}--\eq{ai6eq11}, using sums over $\{0,\ldots,N\}$-labelled rooted trees $\bs\Ga_i$. Here the upper limit on $k$ is not really necessary, as if it does not hold then $\d_{n,k}^{\{i_0,\ldots,i_j\}}=0$ by Lemma~\ref{ai6lem2}.
\item[(iii)] Definition \ref{ai6def9}(iii)--(v) hold.
\end{itemize}

\noindent{\bf(b)} If $\al:[\bar N]\ra[N]$ is a morphism in $\De$, that is, $\al:\{0,\ldots,\bar N\}\ra\{0,\ldots,N\}$ is a nondecreasing function, define
\begin{equation*}
\dAAlg^{\un,\iy}_{R,K}(\al):(\dAAlg^{\un,\iy}_{R,K})_N\longra(\dAAlg^{\un,\iy}_{R,K})_{\bar N}
\end{equation*}
to map $\bigl((A^i)_{i=0}^N,(\m^I)_{\es\ne I\subseteq\{0,\ldots,N\}}\bigr)$ to $\bigl((A^{\al(i)})_{i=0}^{\bar N},(\bar\m{}^{\bar I})_{\es\ne \bar I\subseteq\{0,\ldots,\bar N\}}\bigr)$, where $\bar\m{}^{\bar I}\ab =(\bar\m{}^{\bar I}_k)_{k\ge 0}$ with
\e
\bar\m_k^{\bar I}=\begin{cases}
\m_k^{\al(\bar I)} & \text{if $\al\vert_{\bar I}$ is injective,} \\
\id_{A^i} & \text{if $k\!=\!1$ and $\bar I\!=\!\{\bar\imath_0,\bar\imath_1\}$ with $\bar\imath_0\!\ne\! \bar\imath_1$, $\al(\bar\imath_0)\!=\!\al(\bar\imath_1)\!=\!i$,} \\
0 & \text{otherwise.}
\end{cases}
\label{ai6eq29}
\e

Proposition \ref{ai6prop2} shows $\dAAlg^{\un,\iy}_{R,K}(\al)$ does map to $(\dAAlg^{\un,\iy}_{R,K})_{\bar N}$, and thus $\dAAlg^{\un,\iy}_{R,K}$ is a simplicial set as in Definition \ref{ai6def9}. Theorem \ref{ai6thm7} will show $\dAAlg^{\un,\iy}_{R,K}$ is an $\iy$-category. We say that $\dAAlg^{\un,\iy}_{R,K}$ has {\it type\/} $K$. Clearly $\dAAlg^{\un,\iy}_{R,K}\subset\dAAlg^{\un,\iy}_{R,L}$ if $L\le K$ in $\{3,4,5,\ldots,\iy\}$, giving the inclusion relations in the second column of~\eq{ai6eq1}.
\end{dfn}

Note that \eq{ai6eq18} and \eq{ai6eq28} are different, and \eq{ai6eq22} and \eq{ai6eq29} are different. The following analogue of Proposition \ref{ai6prop1} is proved in~\S\ref{ai106}.

\begin{prop}
\label{ai6prop2}
In Definition\/ {\rm\ref{ai6def12}(b),} $\bigl((A^{\al(i)})_{i=0}^{\bar N},(\bar\m{}^{\bar I})_{\es\ne \bar I\subseteq\{0,\ldots,\bar N\}}\bigr)$ lies in $(\dAAlg^{\un,\iy}_{R,K})_{\bar N}$.
\end{prop}

The next theorem, an analogue of Theorem \ref{ai6thm2}, will be proved in \S\ref{ai107} using equation~\eq{ai6eq28}.

\begin{thm}
\label{ai6thm5}
Let\/ $\bigl((A^i),(\m^I)\bigr)$ be an $N$-simplex in\/ $\dAAlg^{\un,\iy}_{R,K},$ as in Definition\/ {\rm\ref{ai6def12},} and define maps\/ $\d_{n,k}^I$ from the $\m_k^I$ as in\/ {\rm\eq{ai6eq8}--\eq{ai6eq11}}. Then for each\/ $\es\ne I=\{i_0<\cdots<i_M\}\subseteq\{0,\ldots,N\}$ and all\/ $k,n\ge 0,$ in $R$-linear maps $(A^{i_0})^{\ot^n}\ra(A^{i_M})^{\ot^k}$ we have
\e
\sum_{j=0}^M\sum_{0\le l\le \max(n,n+j-1)\!\!\!\!\!\!}
(-1)^j\,\d_{l,k}^{\{i_j,\ldots,i_M\}}\ci\d_{n,l}^{\{i_0,\ldots,i_j\}}=\sum_{j=1}^{M-1}(-1)^j\,\d_{n,k}^{I\sm\{i_j\}}.
\label{ai6eq30}
\e
Here the upper limit on $l$ is not really necessary, as if it does not hold then $\d_{n,l}^{\{i_0,\ldots,i_j\}}=0$ by Lemma\/ {\rm\ref{ai6lem2}}.
\end{thm}

Here is the analogue of Example~\ref{ai6ex3}.

\begin{ex}
\label{ai6ex4}
{\bf(a)} A 0-simplex $\bigl((A^0),(\m^{\{0\}})\bigr)$ of $\dAAlg^{\un,\iy}_{R,K}$ is a unital $A_\iy$-algebra over $R$.
\smallskip

\noindent{\bf(b)} A 1-simplex $\bigl((A^0,A^1),(\m^{\{0\}},\m^{\{1\}},\m^{\{0,1\}})\bigr)$ of $\dAAlg^{\un,\iy}_{R,K}$ is equivalent to a unital morphism $\m^{\{0,1\}}:(A^0,\m^{\{0\}})\ra(A^1,\m^{\{1\}})$ of unital $A_\iy$-algebras.
\smallskip

\noindent{\bf(c)} A 2-simplex of $\dAAlg^{\un,\iy}_{R,K}$ is equivalent to unital $A_\iy$-algebras $(A^0,\m^{\{0\}}),\ab\ab(A^1,\ab\m^{\{1\}}),\ab(A^2,\m^{\{2\}})$ as in {\bf(a)}, and unital morphisms $\m^{\{i,j\}}:(A^i,\m^{\{i\}})\ab\ra(A^j,\m^j)$ for $(i,j)=(0,1),(0,2),(1,2)$ as in {\bf(b)}, and $R$-linear maps 
\begin{equation*}
\m^{\{0,1,2\}}_k:(A^0)^{\ot^k}\longra A^2
\end{equation*}
for $k=1,2,\ldots,$ graded of degree $-k$, with $\m^{\{0,1,2\}}_0=0$, satisfying for all $n\ge 0$
\ea
&\sum_{i,k\ge 0,\; j\ge 1:\; i+j+k=n}
(-1)^{i+jk}\m^{\{0,1,2\}}_{i+k+1}\ci\bigl(\id_{A^0}^{\ot^i}\ot\m^{\{0\}}_j\ot \id_{A^0}^{\ot^k}\bigr)
\nonumber\\
&-\sum_{n=n_1+\cdots+n_k,\; n_a\ge 1\!\!\!\!\!\!\!\!\!\!\!\!\!\!\!\!}(-1)^{\sum_{i=1}^k(n_i+1)(k-i)}\m^{\{1,2\}}_k\bigl(\m^{\{0,1\}}_{n_1}\ot\m^{\{0,1\}}_{n_2}\ot\cdots\ot\m^{\{0,1\}}_{n_k}\bigr)
\nonumber\\
&+\sum_{\begin{subarray}{l} 
n=n_1+\cdots+n_k, \; 1\le l\le k \\
n_i=p^1_i+\cdots+p^{m_i}_i,\; i=1,\ldots,l-1, \\
\text{$p_i^j\ge 1$, $m_i\ge 1$, all $i,j$,}\\
n_i\ge 1,\; i=1,\ldots,k \end{subarray}}
\begin{aligned}[t]&(-1)^{k+1}\cdot(-1)^{\sum_{i=1}^{l-1}(n_i+1)}\cdot{} \\
&(-1)^{\sum_{i=1}^k(n_i+1)(k-i)}\cdot\prod_{i=1}^{l-1}
(-1)^{\sum_{a=1}^{m_i}(p_i^a+1)(m_i-a)}\\
&\m^{\{2\}}_k\bigl[\m^{\{1,2\}}_{m_1}\bigl(\m^{\{0,1\}}_{p^1_1}\!\ot\!\m^{\{0,1\}}_{p^2_1}\!\ot\!\cdots\!\ot\!\m^{\{0,1\}}_{p^{m_1}_1}\bigr)\\
&\ot\cdots\ot\m^{\{1,2\}}_{m_{l-1}}\bigl(\m^{\{0,1\}}_{p^1_{l-1}}\!\ot\!\m^{\{0,1\}}_{p^2_{l-1}}\!\ot\!\cdots\!\ot\!\m^{\{0,1\}}_{p^{m_{l-1}}_{l-1}}\bigr)\\
&\ot\m^{\{0,1,2\}}_{n_l}\ot\m^{\{0,2\}}_{n_{l+1}}\ot\cdots\ot\m^{\{0,2\}}_{n_k}\bigr]
\end{aligned}
\nonumber\\
&=-\m^{\{0,2\}}_n.
\label{ai6eq31}
\ea
Note that we require $n_l\ge 1$ in the third sum rather than $n_l\ge 0$, as $\m^{\{0,1,2\}}_0=0$. Equation \eq{ai6eq31} is trivial when $n=0$. By equations \eq{ai3eq10}, \eq{ai3eq12}, \eq{ai3eq13}, \eq{ai6eq21}, and \eq{ai6eq31}, $\m^{\{0,1,2\}}=(\m^{\{0,1,2\}}_k)_{k\ge 1}$ is a {\it unital homotopy\/} $\m^{\{1,2\}}\ci \m^{\{0,1\}}\Ra\m^{\{0,2\}}$. Hence, a 2-simplex in $\dAAlg^{\un,\iy}_{R,K}$ is a 2-commutative diagram of $A_\iy$-algebras
\begin{align*}
\xymatrix@C=70pt@R=13pt{
\drrtwocell_{}\omit^{}\omit{_{\,\,\,\,\,\,\,\,\,\,\,\,\,\,\,\,\,\,\m^{\{0,1,2\}}}} & (A^1,\m^{\{1\}}) \ar[dr]^{\m^{\{1,2\}}}  \\
(A^0,\m^{\{0\}}) \ar[ur]^{\m^{\{0,1\}}} \ar[rr]^(0.35){\m^{\{0,2\}}} &&  (A^2,\m^{\{2\}}).}
\end{align*}

In particular, if $\bigl((A^i)_{i=0}^2,(\m^I)_{\es\ne I\subseteq\{0,1,2\}}\bigr)$ is a 2-simplex in $\dAAlg^\iy_{R,K}$ and $\m^{\{0,1,2\}}=0$ then $\m^{\{0,2\}}=\m^{\{1,2\}}\ci\m^{\{0,1\}}$ as morphisms of $A_\iy$-algebras $(A^0,\m^0)\ab\ra(A^2,\m^2)$. This is how the notion of {\it composition of morphisms\/} of $A_\iy$-algebras in Definition \ref{ai3def5} is encoded in our theory.
\end{ex}

The next theorem is proved in~\S\ref{ai108}.

\begin{thm}
\label{ai6thm6}
There is a strict isomorphism of\/ $\iy$-categories
\e
I^\un_{R,K}:\AAlg^{\un,\iy}_{R,K}\longra\dAAlg^{\un,\iy}_{R,K},\quad K\in\{3,4,5,\ldots,\iy\},
\label{ai6eq32}
\e
from Definitions\/ {\rm\ref{ai6def11}} and\/ {\rm\ref{ai6def12},} which acts on $N$-simplices by
\begin{equation*}
I^\un_{R,K}:\bigl((A^i)_{i=0}^N,(\m^I)_{\es\ne I\subseteq\{0,\ldots,N\}}\bigr)\longmapsto
\bigl((A^i)_{i=0}^N,(\dot\m^I)_{\es\ne I\subseteq\{0,\ldots,N\}}\bigr),
\end{equation*}
where $\dot\m^I_k$ is defined from $\m^I_k$ by
\e
\dot\m_k^{\{i_0,\ldots,i_M\}}=\begin{cases} 0, & \text{$k=0$ and\/ $M=2$,} \\
\m_k^{\{i_0,\ldots,i_M\}}, & \text{otherwise.} \end{cases}
\label{ai6eq33}
\e
This is equivalent to \eq{ai6eq27} in this case, as $\m_0^{\{i_0,i_1,i_2\}}=\boo^{i_2}$.
\end{thm}

\begin{rem}
\label{ai6rem6}
One could argue that for $A_\iy$-algebras, homotopies in \S\ref{ai34} are the best notion of when two $A_\iy$-morphisms are `the same'. In contrast, for $A_\iy$-categories, $A_\iy$-natural isomorphisms in \S\ref{ai44} (which also make sense for $A_\iy$-algebras, considered as $A_\iy$-categories with one object) are probably the better notion. So Example \ref{ai6ex4} and Theorem \ref{ai6thm6} give a reason to prefer the $\iy$-categories $\AAlg^{\un,\iy}_{R,K},\dAAlg^{\un,\iy}_{R,K}$ for $K\ge 3$ to~$\AAlg^{\un,\iy}_R$.
\end{rem}

We now define the remaining $\iy$-categories in~\eq{ai6eq1}.

\begin{dfn}
\label{ai6def13}
Let $R$ be a commutative ring and $K\in\{3,4,5,\ldots,\iy\}$. Generalizing Definition \ref{ai6def12}, define a simplicial set $\AAlg^{\cu,\iy}_{R,K}$ of cohomologically unital $A_\iy$-algebras over $R$ such that the $N$-simplices $\bigl((A^i)_{i=0}^N,(\m^I)_{\es\ne I\subseteq\{0,\ldots,N\}}\bigr)$ for $N\ge 0$ satisfy Definition \ref{ai6def12}(a)(i)--(ii), but instead of Definition \ref{ai6def12}(a)(iii), we require that
\begin{itemize}
\setlength{\itemsep}{0pt}
\setlength{\parsep}{0pt}
\item[(iii$)'$] As in Example \ref{ai6ex4}(a), Definition \ref{ai6def12}(a)(i)--(ii) imply that $(A^i,\m^{\{i\}})$ is an $A_\iy$-algebra for each $i=0,\ldots,N$. We require this $A_\iy$-algebra $(A^i,\m^{\{i\}})$ to be {\it cohomologically unital}, as in Definition \ref{ai3def4}.
\item[(iv$)'$] As in Example \ref{ai6ex4}(b), Definition \ref{ai6def12}(a)(i)--(ii) mean $\m^{\{i,j\}}:(A^i,\m^{\{i\}})\ab\ra(A^j,\m^{\{j\}})$ is a morphism of $A_\iy$-algebras for all $0\le i<j\le N$. We require this to be a {\it cohomologically unital\/} morphism, as in Definition~\ref{ai3def5}.
\end{itemize}
If $\al:[\bar N]\ra[N]$ is a morphism in $\De$, define $\AAlg^{\cu,\iy}_{R,K}(\al):(\AAlg^{\cu,\iy}_{R,K})_N\ab\ra(\AAlg^{\cu,\iy}_{R,K})_{\bar N}$ as in Definition \ref{ai6def12}(b). The proof of Proposition \ref{ai6prop2} also works in the cohomologically unital case, and shows $\AAlg^{\cu,\iy}_{R,K}$ is a simplicial set. It contains $\AAlg^{\un,\iy}_{R,K}$ in Definition \ref{ai6def12} as a simplicial subset.

Theorem \ref{ai6thm5} also holds in $\AAlg^{\cu,\iy}_{R,K}$, as the proof does not use unitality.

Theorem \ref{ai6thm7} will show that $\AAlg^{\cu,\iy}_{R,K}$ is an $\iy$-category. We say that $\AAlg^{\cu,\iy}_{R,K}$ has {\it type\/} $K$. We also write~$\AAlg^{\cu,\iy}_R=\AAlg^{\cu,\iy}_{R,3}$.
\end{dfn}

\begin{dfn}
\label{ai6def14}
Let $R$ be a commutative ring and $K\in\{3,4,5,\ldots,\iy\}$. Generalizing Definition \ref{ai6def12}, define a simplicial set $\AAlg^\iy_{R,K}$ of (possibly non-unital) $A_\iy$-algebras over $R$ whose $N$-simplices $\bigl((A^i)_{i=0}^N,(\m^I)_{\es\ne I\subseteq\{0,\ldots,N\}}\bigr)$ for $N\ge 0$ satisfy Definition \ref{ai6def12}(a)(i)--(ii), but need not satisfy (a)(iii). If $\al:[\bar N]\ra[N]$ is a morphism in $\De$, define $\AAlg^\iy_{R,K}(\al):(\AAlg^\iy_{R,K})_N\ab\ra(\AAlg^\iy_{R,K})_{\bar N}$ as in Definition \ref{ai6def12}(b). The proof of Proposition \ref{ai6prop2} also works in the non-unital case, and shows $\AAlg^\iy_{R,K}$ is a simplicial set. It contains $\AAlg^{\cu,\iy}_{R,K}$ in Definition \ref{ai6def13} as a simplicial subset.

Theorem \ref{ai6thm5} also holds in $\AAlg^\iy_{R,K}$, as the proof does not use unitality.

Theorem \ref{ai6thm7} will show that $\AAlg^\iy_{R,K}$ is an $\iy$-category. We say that $\AAlg^\iy_{R,K}$ has {\it type\/} $K$. We also write~$\AAlg^\iy_R=\AAlg^\iy_{R,3}$.
\end{dfn}

The next theorem, one of our main results, is a special case of Theorem~\ref{ai7thm7}.

\begin{thm}
\label{ai6thm7}
The following simplicial sets defined in Definitions\/ {\rm\ref{ai6def11}, \ref{ai6def12}, \ref{ai6def13},} and\/ {\rm\ref{ai6def14},} for all\/ $K\in\{3,4,5,\ldots,\iy\},$ are all\/ $\iy$-categories:
\begin{equation*}
\AAlg^{\un,\iy}_{R,K},\dAAlg^{\un,\iy}_{R,K},\AAlg^{\cu,\iy}_{R,K},\AAlg^\iy_{R,K}.
\end{equation*}
\end{thm}

The next definition is an analogue of Lurie \cite[Rem.~1.3.1.9]{Luri3} for dg-cat\-e\-gor\-ies. It completes our diagram \eq{ai6eq1} of inclusions of $\iy$-categories of $A_\iy$-algebras.

\begin{dfn}
\label{ai6def15}
Let $\AAlg^\un_R$ be the ordinary category of unital $A_\iy$-algebras over $R$ from Definition \ref{ai3def5}, and $\sN(\AAlg^\un_R)$ be its nerve from Example \ref{ai2ex5}, which is an $\iy$-category. Then $N$-simplices in $\sN(\AAlg^\un_R)$ are of the form $\bigl((A^i,\m^i)_{i=0}^N,(\f^i)_{i=1}^N\bigr)$, where $(A^i,\m^i)$ is a unital $A_\iy$-algebra over $R$ for $i=0,\ldots,N$, and $\f^i:(A^{i-1},\m^{i-1})\ra(A^i,\m^i)$ is a unital $A_\iy$-morphism for~$i=1,\ldots,N$.

Define an $\iy$-functor $\inc^\un_R:\sN(\AAlg^\un_R)\ra\dAAlg^{\un,\iy}_{R,\iy}$ to act on $N$-simplices taking $\bigl((A^i,\m^i)_{i=0}^N,(\f^i)_{i=1}^N\bigr)$ to $\bigl((A^i)_{i=0}^N,(\m^I)_{\es\ne I\subseteq\{0,\ldots,N\}}\bigr)$, where we set $\m^{\{i\}}=\m^i$ for $i=0,\ldots,N$, and $\m^{\{i,j\}}=\f^j\ci\f^{j-1}\ci\cdots\ci\f^{i+1}:(A^i,\m^{\{i\}})\ra(A^j,\m^{\{j\}})$ if $0\le i<j\le N$, and $\m^I_k=0$ if $\md{I}\ge 3$ and $k\ge 0$.

We claim $\bigl((A^i)_{i=0}^N,(\m_k^I)_{k\ge 0,\; \es\ne I\subseteq\{0,\ldots,N\}}\bigr)$ is an $N$-simplex of $\dAAlg^{\un,\iy}_{R,\iy}$. The condition that $\m_0^I=0$ for all $\es\ne I\subseteq\{0,\ldots,N\}$ with~$\md{I}\le \iy$ is clear. We must verify \eq{ai6eq28} for all $\es\ne I\subseteq\{0,\ldots,N\}$ and $n\ge 0$. Divide into cases:
\begin{itemize}
\setlength{\itemsep}{0pt}
\setlength{\parsep}{0pt}
\item[(i)] If $I=\{i\}$, this follows from $(A^i,\m^{\{i\}})$ being an $A_\iy$-algebra, as in Example~\ref{ai6ex4}(a).
\item[(ii)] If $I=\{i,j\}$ with $i<j$, this follows from $\m^{\{i,j\}}:(A^i,\m^{\{i\}})\ra (A^j,\m^{\{j\}})$ being a morphism of $A_\iy$-algebras, as in Example~\ref{ai6ex4}(b).
\item[(iii)] If $I=\{i,j,k\}$ with $i<j<k$, this follows from 
\begin{equation*}
\m^{\{i,k\}}=\f^k\ci\f^{k-1}\ci\cdots\ci \f^{j+1}\ci\f^j\ci\cdots\ci\f^{i+1}=\m^{\{j,k\}}\ci\m^{\{i,j\}}
\end{equation*}
and $\m^{\{i,j,k\}}=0$, as in Example~\ref{ai6ex4}(c).
\item[(iv)] If $\md{I}\ge 4$ then on the left hand side of \eq{ai6eq28} we have 
$\m_k^{\{i_j,\ldots,i_M\}}=0$ if $j<M-1$, as then $\md{\{i_j,\ldots,i_M\}}\ge 3$, and 
$\d_{n,k}^{\{i_0,\ldots,i_j\}}=0$ if $j\ge M-1\ge 2$, as every term of $\d_{n,k}^{\{i_0,\ldots,i_j\}}$ must contain a factor $\m_k^J$ with $\md{J}\ge 3$, and so must be zero. Also the right hand side of \eq{ai6eq28} is zero as $\md{I\sm\{i_j\}}\ge 3$.
\end{itemize}
Definition \ref{ai6def12}(a)(iii) follows from Definition~\ref{ai3def5}.

It is now easy to see that $\inc^\un_R$ is a morphism of simplicial sets, that is, an $\iy$-functor. It is a strict isomorphism with the $\iy$-subcategory of $\dAAlg^{\un,\iy}_{R,\iy}$ consisting of $N$-simplices $\bigl((A^i)_{i=0}^N,(\m_k^I)_{k\ge 0,\; \es\ne I\subseteq\{0,\ldots,N\}}\bigr)$ with $\m_k^I=0$ if $\md{I}\ge 3$ and $k\ge 0$. In the same way we define $\iy$-functors $\inc^\cu_R:\sN(\AAlg^\cu_R)\hookra\AAlg^{\cu,\iy}_{R,\iy}$ for cohomologically unital $A_\iy$-algebras, and $\inc_R:\sN(\AAlg_R)\ab\hookra\AAlg^\iy_{R,\iy}$ for general $A_\iy$-algebras.

We can also define an $\iy$-functor $\iinc^\un_R:\sN(\AAlg^\un_R)\ra\AAlg^{\un,\iy}_{R,\iy}$ which is identified with $\inc^\un_R$ under the isomorphism $I^\un_{R,\iy}$ in \eq{ai6eq32}. The difference is that as in \eq{ai6eq27}, when we define $\bigl((A^i)_{i=0}^N,(\m^I)_{\es\ne I\subseteq\{0,\ldots,N\}}\bigr)$ above we define $\m_0^{\{i_0,i_1,i_2\}}=\boo^{i_2}$ rather than~$\m_0^{\{i_0,i_1,i_2\}}=0$. 
\end{dfn}

\section{\texorpdfstring{$\iy$-categories of $A_\iy$-categories}{∞-categories of A∞-categories}}
\label{ai7}

We now define many different $\iy$-categories of $A_\iy$-categories over a commutative ring $R$, which fit into a diagram of inclusions of $\iy$-categories:
\ea
\label{ai7eq1}\\[-19pt]
\xymatrix@!0@C=85pt@R=37pt{{\begin{subarray}{l}\ts \ACat^{\un,\iy}_R\!\!=  \\[2pt] \ts \ACat^{\un,\iy}_{R,2}\end{subarray}} \\
\ACat^{\un,\iy^{\vphantom{(}}}_{R,3}\, \ar@{^{(}->}[u] \ar@{=}[r] \ar@{^{(}->}[u] & \dACat^{\un,\iy^{\vphantom{(}}}_{R,3}\, \ar@{^{(}->}[r] & {\begin{subarray}{l}\ts \ACat^{\cu,\iy}_R\!\!=  \\[2pt] \ts \ACat^{\cu,\iy}_{R,3}\end{subarray}}\,  \ar@{^{(}->}[r] & {\begin{subarray}{l}\ts \ACat^\iy_R\!=  \\[2pt] \ts \ACat^\iy_{R,3}\end{subarray}}  \\
\ACat^{\un,\iy^{\vphantom{(}}}_{R,4}\, \ar@{=}[r] \ar@{^{(}->}[u] & \dACat^{\un,\iy^{\vphantom{(}}}_{R,4}\, \ar@{^{(}->}[r] \ar@{^{(}->}[u] & \ACat^{\cu,\iy^{\vphantom{(}}}_{R,4}\,  \ar@{^{(}->}[r] \ar@{^{(}->}[u] & \ACat^{\iy^{\vphantom{(}}}_{R,4} \ar@{^{(}->}[u] \\
\ACat^{\un,\iy^{\vphantom{(}}}_{R,5}\, \ar@{=}[r] \ar@{^{(}->}[u] & \dACat^{\un,\iy^{\vphantom{(}}}_{R,5}\, \ar@{^{(}->}[r] \ar@{^{(}->}[u] & \ACat^{\cu,\iy^{\vphantom{(}}}_{R,5}\,  \ar@{^{(}->}[r] \ar@{^{(}->}[u] & \ACat^{\iy^{\vphantom{(}}}_{R,5} \ar@{^{(}->}[u]  \\
\vdots_{\vphantom{(}} \ar@{^{(}->}[u] & \vdots_{\vphantom{(}} \ar@{^{(}->}[u] & \vdots_{\vphantom{(}} \ar@{^{(}->}[u] & \vdots_{\vphantom{(}} \ar@{^{(}->}[u] \\
\ACat^{\un,\iy^{\vphantom{(}}}_{R,\iy}\, \ar@{=}[r] \ar@{^{(}->}[u] & \dACat^{\un,\iy^{\vphantom{(}}}_{R,\iy}\, \ar@{^{(}->}[r] \ar@{^{(}->}[u] & \ACat^{\cu,\iy^{\vphantom{(}}}_{R,\iy}\,  \ar@{^{(}->}[r] \ar@{^{(}->}[u] & \ACat^{\iy^{\vphantom{(}}}_{R,\iy} \ar@{^{(}->}[u] \\
\sN(\ACat^\un_R)\, \ar@{=}[r] \ar@{^{(}->}[u] & \sN(\ACat^\un_R)\, \ar@{^{(}->}[r] \ar@{^{(}->}[u] & \sN(\ACat^\cu_R)\,  \ar@{^{(}->}[r] \ar@{^{(}->}[u] & \sN(\ACat_R).\! \ar@{^{(}->}[u] }
\nonumber
\ea
This is the exact analogue of our diagram \eq{ai6eq1} of $\iy$-categories of $A_\iy$-algebras, and reduces to \eq{ai6eq1} if we restrict to $A_\iy$-categories with one object $*$, which are equivalent to $A_\iy$-algebras. 

As in Chapter \ref{ai6}, we note that:
\begin{itemize}
\setlength{\itemsep}{0pt}
\setlength{\parsep}{0pt}
\item[(i)] The $\iy$-categories in the first two columns, with superscripts `$\un$', have objects unital $A_\iy$-categories and 1-morphisms unital $A_\iy$-functors.
\item[(ii)] The $\iy$-categories in the third column, with superscript `$\cu$', have objects cohomologically unital $A_\iy$-categories and 1-morphisms cohomologically unital $A_\iy$-functors.
\item[(iii)] The $\iy$-categories in the fourth column have objects all $A_\iy$-categories, and 1-morphisms all $A_\iy$-functors.
\item[(iv)] The $\iy$-categories in the first and second columns, linked by `$=$', are strictly isomorphic, but with definitions written in a different way.
\item[(v)] In each row we have subscripts `$R,K$' for $K=2,3,4,\ldots,\iy$, where we call $K$ the {\it type\/} of the $\iy$-category.
\item[(vi)] $\sN(\ACat^\un_R),\ldots,\sN(\ACat_R)$ are the nerves of the ordinary categories $\ACat^\un_R,\ldots,\ACat_R$ from Definition~\ref{ai4def6}.
\end{itemize}

The 2-morphisms in $\ACat^{\un,\iy}_R$ of type 2 in the first row of \eq{ai7eq1} correspond to {\it unital\/ $A_\iy$-natural isomorphisms\/} in Definition \ref{ai4def13}. But in the $\iy$-categories of type $K\ge 3$ in the subsequent rows, the 2-morphisms correspond to ({\it unital\/}) {\it homotopies\/} of $A_\iy$-functors, as in Definition \ref{ai4def12}. In the author's view, $\ACat^{\un,\iy}_R$ will be the best $\iy$-category of $A_\iy$-categories for most applications, as $A_\iy$-natural isomorphisms are the best notion of when two $A_\iy$-functors are `the same'. It is a pity that this only works for unital $A_\iy$-categories.

\subsection{Labelling graphs by objects in categories}
\label{ai71}

Sections \ref{ai61}--\ref{ai62} set up notation for defining $\iy$-categories of $A_\iy$-algebras, involving $R$-linear maps $\m^I_k:(A^{i_0})^{\ot^k}\ra A^{i_M}$. To extend this to $A_\iy$-categories, as in Chapter \ref{ai4} we must extend these to maps $\m^{I;x_0,\ldots,x_k}_k$ depending on objects $x_0,\ldots,x_k$ in an $A_\iy$-category. When we deal with {\it simplices\/} of morphisms, as in \S\ref{ai62}, it will turn out that for $I=\{i_0<i_1<\cdots<i_M\}$ we must work with graded $R$-linear maps of the form
\begin{align*}
\m_k^{I;x_0,\ldots,x_k}:\Hom_{\sA^{i_0}}(x_0,x_1)\ot_R\Hom_{\sA^{i_0}}(x_1,x_2)\ot_R\cdots \ot_R\Hom_{\sA^{i_0}}(x_{k-1},x_k)&\\
\longra\Hom_{\sA^{i_M}}\bigl(f^{i_{M-1}i_M}\ci f^{i_{M-2}i_{M-1}}\ci\cdots\ci f^{i_0i_1}(x_0),f^{i_0i_M}(x_k)\bigr),&
\end{align*}
where $\sA^i$ are $A_\iy$-categories and $f^{ij}=\fF_\Obj^{ij}:\Obj\sA^i\ra\Obj\sA^j$ for $A_\iy$-functors $\fF^{ij}:\sA^i\ra\sA^j$.

We now extend the graphical notation of \S\ref{ai61}--\S\ref{ai62} to cover multiple compositions of such maps $\m_k^{I;x_0,\ldots,x_k}$. To do this we add labels of ordered pairs $(y_e,z_e)$ of objects in $\sA^i$ to each edge $e$ in $\Ga$, such that $e$ represents an input/output in $\Hom_{\sA^i}(y_e,z_e)$. These pairs $(y_e,z_e)$ must satisfy consistency conditions to match the domains and ranges of the maps~$\m_k^{I;x_0,\ldots,x_k}$.

We first define maps $C(\bs\Ga)^{x_0,\ldots,x_n}$, generalizing $C(\bs\Ga)$ in Definition \ref{ai6def5}. We restrict to {\it balanced\/} $\bs\Ga$, as in Definition \ref{ai6def6}, as being balanced automatically implies the required consistency conditions on pairs~$(y_e,z_e)$.

\begin{dfn}
\label{ai7def1}
Let $R$ be a commutative ring. Let $N\ge 0$, and suppose we are given data $(\sO^i)_{i=0}^N,(f^{ij})_{0\le i\le j\le N},(\sH^i)_{i=0}^N,(\m^I)_{\es\ne I\subseteq\{0,\ldots,N\}}$, where:
\begin{itemize}
\setlength{\itemsep}{0pt}
\setlength{\parsep}{0pt}
\item[(i)] $\sO^i$ is a set for $i=0,\ldots,N$. Elements of $\sO^i$ are called {\it objects}. Later we will take $\sO^i=\Obj\sA^i$ for an $A_\iy$-category $\sA^i$.
\item[(ii)] $f^{ij}:\sO^i\ra\sO^j$ is a map for all $0\le i\le j\le N$. We require these maps to satisfy $f^{ii}=\id_{\sO^i}$ for $i=0,\ldots,N$. Later we will take $f^{ij}=\fF_\Obj^{ij}$ for an $A_\iy$-functor $\fF^{ij}:\sA^i\ra\sA^j$. For the $\iy$-category $\ACat^{\un,\iy}_R$ we do {\it not\/} require that $f^{jk}\ci f^{ij}=f^{ik}$, so for the present we do not assume this, though for the other $\iy$-categories in \eq{ai7eq1} we do require it.

Let $\es\ne I\subseteq\{0,1,\ldots,N\}$ with $I=\{i_0,i_1,\ldots,i_M\}$ for $i_0<i_1<\cdots<i_M$. Define
\end{itemize}
\e
f^I:\sO^{i_0}\ra\sO^{i_M}\;\>\text{by}\;\>
f^I=\begin{cases}
f^{i_0i_0}=\id_{\sO^{i_0}}, & M=0, \\
f^{i_{M-1}i_M}\!\ci\! f^{i_{M-2}i_{M-1}}\!\ci\!\cdots\!\ci\! f^{i_0i_1}, & M>0.
\end{cases}
\label{ai7eq2}
\e
\begin{itemize}
\setlength{\itemsep}{0pt}
\setlength{\parsep}{0pt}
\item[(iii)] For all $i=0,1,\ldots,N$ and $x,y\in\sO^i$, we are given a $\Z$-graded $R$-module $\sH^i(x,y)=\bigop_{d\in\Z}\sH^i(x,y)^d$. We write $\sH^i=(\sH^i(x,y))_{x,y\in\sO^i}$. Later we will take $\sH^i(x,y)=\Hom_{\sA^i}(x,y)$ for the $A_\iy$-category $\sA^i$.
\item[(iv)] For each nonempty subset $I\subseteq\{0,\ldots,N\}$, which we write as $I=\{i_0,\ab i_1,\ab\ldots,\ab i_M\}$ for $i_0<i_1<\cdots<i_M$ with $M\ge 0$, and for all $k\ge 0$ and elements $x_0,\ldots,x_k\in \sO^{i_0}$, we are given an $R$-linear map 
\ea
&\m_k^{I;x_0,\ldots,x_k}:\sH^{i_0}(x_0,x_1)\ot_R\sH^{i_0}(x_1,x_2)\ot_R\cdots \ot_R\sH^{i_0}(x_{k-1},x_k)
\nonumber\\
&\quad\longra\sH^{i_M}(f^I(x_0),f^{i_0i_M}(x_k)),
\label{ai7eq3}
\ea
graded of degree $2-M-k=3-\md{I}-k$.

Actually, for the case $k=0$, our definition of $\ACat^{\un,\iy}_R$ only uses maps $\m_0^{I;x_0}$ when $\md{I}\ge 3$. So we could say we want maps $\m_k^{I;x_0,\ldots,x_k}$ when either $\md{I}\ge 3$ and $k\ge 0$, or when $\md{I}=1,2$ and $k\ge 1$. It is convenient to set $\m_0^{I;x_0}=0$ when $\md{I}=1,2$. For the $\iy$-category of {\it filtered\/} $A_\iy$-categories in Chapter \ref{ai8} we will use maps $\m_0^{I;x_0}$ for all $\md{I}\ge 1$.
\end{itemize}

Let $\es \ne I\subseteq\{0,\ldots,N\}$, written as $I=\{i_0,\ab i_1,\ab\ldots,\ab i_M\}$ for $i_0<i_1<\cdots<i_M$. Let $\bs\Ga$ be a balanced $n$-leaved $I$-labelled tree, as in Definition \ref{ai6def6}. Let $x_0,\ldots,x_n\in\sO^{i_0}$. Generalizing $C(\bs\Ga)$ in \eq{ai6eq2}, we will define an $R$-linear map
\ea
C(\bs\Ga)^{x_0,\ldots,x_n}&:\sH^{i_0}(x_0,x_1)\ot_R\sH^{i_0}(x_1,x_2)\ot_R\cdots \ot_R\sH^{i_0}(x_{n-1},x_n)
\label{ai7eq4}\\
&\qquad\qquad \longra\sH^{i_M}(f^{L(\bs\Ga)}(x_0),f^{R(\bs\Ga)}(x_n)),
\nonumber
\ea
where $L(\bs\Ga)\supseteq R(\bs\Ga)$ are the left and right sequences of $\bs\Ga$ in Definition~\ref{ai6def6}.

We do this by associating operators $\m^{I_v;y_0,\ldots,y_k}_k$ from \eq{ai7eq3} to each internal vertex or bud $v$ in $\bs\Ga$ with $k$ incoming edges, and taking the composition over the diagram, as in Definition \ref{ai6def2}. The new feature (compared to Definition \ref{ai6def5}, which involved maps $\m^{I_v}_k$ rather than $\m^{I_v;y_0,\ldots,y_k}_k$) is that we have to explain how to define the objects $y_0,\ldots,y_k$ in $\sO^{\min I_v}$, in such a way that the maps $\m^{I_v;y_0,\ldots,y_k}_k$ are composable, that is, the output $\sH^{i_M}(f^I(x_0),f^{i_0i_M}(x_k))$ from one map \eq{ai7eq3} matches an input $\sH^i(y_a,y_{a+1})$ of another map \eq{ai7eq3} into which it is substituted.

We associate additional {\it object labels\/} to the tree $\Ga$: to each directed edge $\overset{v}{\bu}\,{\buildrel e\over\longra}\,\overset{w}{\bu}$, we associate an ordered pair $(y_e,z_e)$ with $y_e,z_e\in \sO^{\max(I_v)}=\sO^{\min(I_w)}$, and to each bud $b$ we associate an object $x_{j_b}$, according to the following rules:
\begin{itemize}
\setlength{\itemsep}{0pt}
\setlength{\parsep}{0pt}
\item[(a)] Let $l_j$ be the $j^{\rm th}$ leaf of $\Ga$ for $j=1,\ldots,n$, as in Definition \ref{ai6def1}(viii), and let $\overset{l_j}{\ci}\,{\buildrel e\over\longra}\,\overset{w}{\bu}$ be the unique edge in $\Ga$ starting at $l_j$. We set $(y_e,z_e)=(x_{j-1},x_j)$, for $x_k\in\sO^{i_0}$ as given.
\item[(b)] Let $r$ be the root of $\Ga$, and let $\overset{v}{\bu}\,{\buildrel e\over\longra}\,\overset{r}{\tr}$ be the unique edge in $\Ga$ ending at $r$. We set $(y_e,z_e)=(f^{L(\bs\Ga)}(x_0),f^{R(\bs\Ga)}(x_n))$.
\item[(c)] Let $v$ be an internal vertex of $\Ga$. As in Definition \ref{ai6def1}(vii), write the edges of $\Ga$ in anticlockwise cyclic order as $e_{v,0},\ab e_{v,1},\ab\ldots,\ab e_{v,k}$, where $e_{v,0}$ is the outgoing edge. We require that
\begin{itemize}
\setlength{\itemsep}{0pt}
\setlength{\parsep}{0pt}
\item[(i)] $z_{e_{v,i}}=y_{e_{v,i+1}}$ in $\sO^{\min I_v}$ for all $1\le i<k$.
\item[(ii)] $y_{e_{v,0}}=f^{I_v}(y_{e_{v,1}})$ and $z_{e_{v,0}}=f^{\min I_v\max I_v}(z_{e_{v,k}})$ in $\sO^{\max I_v}$.
\end{itemize}
\item[(d)] Let $b$ be a bud of $\Ga$. As in Definition \ref{ai6def1}(x), write the leaves and buds of $\Ga$ as $u_1,\ldots,u_{m+n}$, ordered anticlockwise around the outside of $\Ga$ starting at the root $r$. Then $b$ lies in the list $u_1,\ldots,u_{m+n}$, and the leaves $l_1,\ldots,l_n$ appear in this order in the list $u_1,\ldots,u_{m+n}$. Let $j_b=0,\ldots,n$ be unique such that $j_b=0$ when $b$ lies before $l_1$ in $u_1,\ldots,u_{m+n}$, and $b$ lies between $l_{j_b}$ and $l_{j_b+1}$ in $u_1,\ldots,u_{m+n}$ when $1\le j_b<n$, and $j_b=n$ when $b$ lies after $l_n$ in $u_1,\ldots,u_{m+n}$. We label $b$ by the object $x_{j_b}$. Let $\overset{b}{\sq}\,{\buildrel e\over\longra}\,\overset{v}{\bu}$ be the unique edge in $\Ga$ starting at $b$, and $P_b$ be the prefix of $b$. We require that
\begin{equation*}
y_e= f^{I_b}\ci f^{P_b}(x_{j_b}),\quad
z_e=f^{\min I_b\max I_b}\ci f^{P_b}(x_{j_b}).
\end{equation*}
\end{itemize}

We can now show that there are unique labels $(y_e,z_e)$ for each edge $e$ and $x_{j_b}$ for each bud $b$ satisfying (a)--(d). This works because the condition \eq{ai6eq4} in the definition of $\bs\Ga$ balanced is exactly what is required to ensure that $z_{e_{v,i}}=y_{e_{v,i+1}}$ as in (c)(i), and the definitions of $L(\bs\Ga),R(\bs\Ga)$ are what is required to describe the target of \eq{ai7eq4} at the root $r$. We call these labels $(y_e,z_e),x_{j_b}$ {\it object labels}, and $\bs\Ga$ with its object labels an {\it object-labelled $n$-leaved $I$-labelled tree}. 

Figure \ref{ai7fig1} shows the tree $\bs\Ga$ of Figure \ref{ai6fig6} labelled with objects.

\begin{figure}[htb]
\centering
\includegraphics{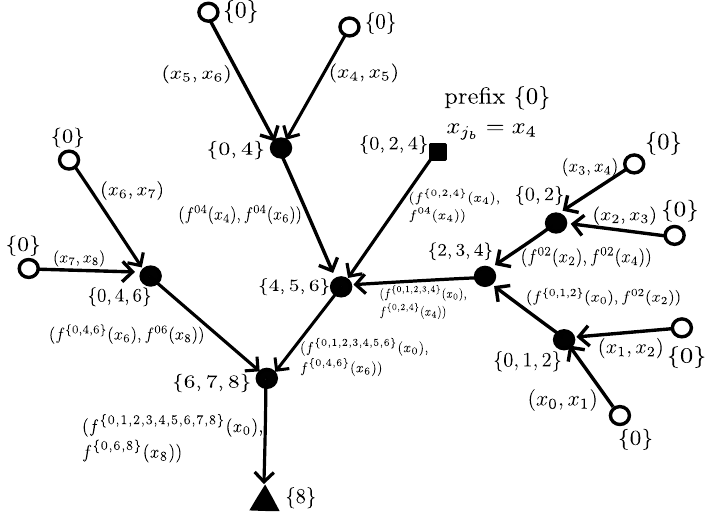}    
\caption{An object-labelled, balanced $\{0,\ldots,8\}$-labelled tree}
\label{ai7fig1}
\end{figure}

Now we have defined the labels $(y_e,z_e)$ for each edge $e$ of $\Ga$, to each internal vertex $v$ of $\Ga$ with ordered edges $e_{v,0},\ab e_{v,1},\ab\ldots,\ab e_{v,k}$ as in (c), we associate the list of objects $w_0,\ldots,w_k$ in $\sO^{\min I_v}$ where $w_0=y_{e_{v,1}}$, $w_1=z_{e_{v,1}}=y_{e_{v,2}}$, \ldots, $w_{k-1}=z_{e_{k-1}}=y_{e_{v,k}}$, $w_k=z_{e_{v,k}}$, and then we associate the $R$-linear map $\m_k^{I_v;w_0,\ldots,w_k}$ to $v$. To each bud $b$ with label $x_{j_b}$ and prefix $P_b$ we associate the element
\e
\m_0^{I_b,f^{P_b}(x_{j_b})}(1)\in \sH^{\max I_b}\bigl(f^{I_b}\ci f^{P_b}(x_{j_b}),f^{\min I_b\max I_b}\ci f^{P_b}(x_{j_b})\bigr).
\label{ai7eq5}
\e

We define $C(\bs\Ga)^{x_0,\ldots,x_n}$ in \eq{ai7eq4} to be the composition over all internal vertices $v$ of $\Ga$ of the maps $\m_k^{I_v;w_0,\ldots,w_k}$, and all buds $b$ of the elements \eq{ai7eq5}, as in Definition \ref{ai6def2}. As in Definition \ref{ai6def5}, $C(\bs\Ga)^{x_0,\ldots,x_n}$ is graded of degree~\eq{ai6eq3}.
\end{dfn}

\begin{dfn}
\label{ai7def2}
Let $R$, $N$, $(\sO^i)_{i=0}^N,(f^{ij})_{0\le i\le j\le N}$, $(\sH^i)_{i=0}^N$, $(\m^I)_{\es\ne I\subseteq\{0,\ldots,N\}}$ be as in Definition \ref{ai7def1}. Choose $\es\ne I\subseteq\{0,\ldots,N\}$, which we write as $I=\{i_0,\ab i_1,\ab\ldots,\ab i_M\}$ for $i_0<i_1<\cdots<i_M$ with $M\ge 0$. Suppose $k,n\ge 0$, and $x_0,\ldots,x_n\in \sA^{i_0}$, and $y_0,\ldots,y_k\in \sA^{i_M}$. Generalizing Definition \ref{ai6def8}, we will define an $R$-linear map graded of degree $k-n+2-\md{I}$
\ea
&\d_{n,k;y_0,\ldots,y_k}^{I;x_0,\ldots,x_n}:
\sH^{i_0}(x_0,x_1)\ot_R\sH^{i_0}(x_1,x_2)\ot_R\cdots \ot_R\sH^{i_0}(x_{n-1},x_n)
\nonumber\\
&\longra \sH^{i_M}(y_0,y_1)\ot_R\sH^{i_M}(y_1,y_2)\ot_R\cdots \ot_R\sH^{i_M}(y_{k-1},y_k).
\label{ai7eq6}
\ea
Here if $n=0$ the domain $\sH^{i_0}(x_0,x_1)\ot\cdots\ot\sH^{i_0}(x_{n-1},x_n)$ is $R$, as it is a tensor product of zero factors which is $R$ by convention, and similarly if $k=0$ the target $\sH^{i_M}(y_0,y_1)\ot\cdots\ot\sH^{i_M}(y_{k-1},y_k)$ is~$R$.
\smallskip

\noindent{\bf(a)} Suppose $k=n=0$, and $M=1$ so $I=\{i_0,i_1\}$, and $y_0=f^{i_0i_1}(x_0)$. As in \eq{ai6eq9}, we define
\e
\d_{0,0;f^{i_0i_1}(x_0)}^{\{i_0,i_1\};x_0}=\id_R:R\longra R.
\label{ai7eq7}
\e
In all other cases when $k=0$ --- if $n\ne 0$, or $M\ne 1$, or $n=0$, $M=1$ and $y_0\ne f^{i_0i_1}(x_0)$ --- we define $\d_{n,0;y_0}^{I;x_0,\ldots,x_n}=0$, as in \eq{ai6eq9}.
\smallskip

\noindent{\bf(b)} If $M=0$ and $1\le k\le n$, generalizing \eq{ai6eq10} we define
\ea
&\d_{n,k;y_0,\ldots,y_k}^{\{i_0\};x_0,\ldots,x_n}=
\label{ai7eq8}\\
&\sum_{\begin{subarray}{l} 0\le j<k, \\
y_i=x_i,\; 0\le i\le j, \\
y_i=x_{n-k+i},\; j<i\le k
\end{subarray}}\begin{aligned}[t]&
(-1)^{j+(n-k+1)(k-j-1)}\id_{\sH^{i_0}(x_0,x_1)}\ot\cdots\ot\\
&\id_{\sH^{i_0}(x_{j-1},x_j)}\ot\m_{n-k+1}^{\{i_0\};x_j,x_{j+1},\ldots,x_{j+n-k+1}}\ot\\
&\id_{\sH^{i_0}(x_{j+n-k+1},x_{j+n-k+2})}\ot\cdots\ot\id_{\sH^{i_0}(x_{n-1},x_n)}.
\end{aligned}
\nonumber
\ea
If $M=0$ and $k>n$ we define $\d_{n,k;y_0,\ldots,y_k}^{\{i_0\};x_0,\ldots,x_n}=0$.
\smallskip

\noindent{\bf(c)} If $M>0$ and $k\ge 1$ then generalizing \eq{ai6eq11}, define 
\ea
&\d_{n,k;y_0,\ldots,y_k}^{I;x_0,\ldots,x_n}=
\label{ai7eq9}\\
&\sum_{\begin{subarray}{l} 0=l_0\le l_1\le \cdots\le l_k=n, \\
\text{isomorphism classes} \\ \text{of $\bs\Ga_1,\ldots,\bs\Ga_k$:}\\ \text{$\bs\Ga_i$ is a balanced} \\ 
\text{$I$-labelled tree} \\ \text{with $l_i-l_{i-1}$} \\ \text{leaves, $1\le i\le k$,} \\
\text{$L(\bs\Ga_1)=I$,} \\ \text{$R(\bs\Ga_k)=\{i_0,i_M\}$,} \\
\text{$R(\bs\Ga_i)=L(\bs\Ga_{i+1})$,} \\ \text{$i=1,\ldots,k-1$,} \\
y_i=f^{R(\bs\Ga_i)}(x_{l_i}), \\
i=1,\ldots,k,\; y_0=f^I(x_0) \end{subarray}}
\begin{aligned}[t]
&(-1)^{(k+1)(M+n+1)}\cdot{}\\
&\biggl(\prod_{\begin{subarray}{l} \text{$v$ internal vertex or bud of\/ $\bs\Ga_i,$ $1\le i\le k,$} \\ \text{$w$ internal vertex or bud of\/ $\bs\Ga_{i'},$ $1\le i'\le k,$} \\ 
\text{either $i<i',$ or $i=i'$ and} \\ \text{$v\pr w$ in order $\pr$ of $\Ga_i$,} \\ \text{$v$ has $p$ incoming edges including} \\ \text{$q$ leaves, $w$ has $r$ incoming edges}\end{subarray}\!\!\!\!\!\!\!\!\!\!\!\!\!\!\!\!\!\!\!\!\!\!\!\!\!\!\!\!\!\!\!\!\!\!\!\!\!\!\!} (-1)^{\begin{subarray}{l} N(v,w)+q(r+1) \\ +(p+1)(\md{I_w}+r+1)\end{subarray}}\biggr)\cdot{}\\
& C(\bs\Ga_1)^{x_{l_0},\ldots,x_{l_1}}\!\ot\! C(\bs\Ga_2)^{x_{l_1},\ldots,x_{l_2}}\!\ot\!\cdots\!\ot\! C(\bs\Ga_k)^{x_{l_{k-1}},\ldots,x_{l_k}},
\end{aligned}
\nonumber
\ea
where we use the notation of \S\ref{ai61}--\S\ref{ai62}, and $N(v,w)$ is as in \eq{ai6eq12}.
\end{dfn}

Here is the analogue of Example \ref{ai6ex2}(i),(ii).

\begin{ex}
\label{ai7ex1}
{\bf(i)} If $k=1$ then for any $I,n$, as in \eq{ai6eq13} we have
\ea
\d^{I;x_0,\ldots,x_n}_{n,1;y_0,y_1}=\begin{cases} \m_n^{I;x_0,\ldots,x_n}, & \text{$n\ge 1$ or $\md{I}\ge 3$, and $y_0=f^I(x_0)$, $y_1=f^{i_0i_M}(x_n)$,} \\
0, & \text{otherwise.}
\end{cases}
\nonumber\\[-15pt]
\label{ai7eq10}
\ea

\noindent{\bf(ii)} If $I=\{0,1\}$ and $k\ge 1$ then as in \eq{ai6eq14} we have
\ea
&\d_{n,k;y_0,\ldots,y_k}^{\{0,1\};x_0,\ldots,x_n}=
\label{ai7eq11}\\
&\sum_{\begin{subarray}{l} n=n_1+\cdots+n_k,\; n_a\ge 1. \\
\text{Set $p_j=\sum_{i=0}^jn_i$, $j=0,\ldots,k$.} \\
\text{Then $y_j=f^{01}(x_{p_j})$, $j=0,\ldots,k$}\end{subarray}}\begin{aligned}[t]
&(-1)^{\sum_{i=1}^k(n_i+1)(k-i)}\m^{\{0,1\};x_{p_0},\ldots,x_{p_1}}_{n_1}\ot{} \\ 
&\m^{\{0,1\};x_{p_1},\ldots,x_{p_2}}_{n_2}\ot\cdots\ot\m^{\{0,1\};x_{p_{k-1}},\ldots,x_{p_k}}_{n_k}.
\end{aligned}
\nonumber
\ea
\end{ex}

\begin{rem}
\label{ai7rem1}
In the defining equations \eq{ai7eq7}--\eq{ai7eq9} for $\d_{n,k;y_0,\ldots,y_k}^{I;x_0,\ldots,x_n}$, to get a nonzero term the $y_j$ are determined by the $x_i$ by equations of the form $y_j=f^J(x_i)$ for some $J\subseteq I$ with $i_0,i_M\in J$, and if the $y_j$ are not of this form then $\d_{n,k;y_0,\ldots,y_k}^{I;x_0,\ldots,x_n}=0$. Thus, for fixed $I,n,k,x_0,\ldots,x_n$, there are only finitely many lists $y_0,\ldots,y_k$ such that~$\d_{n,k;y_0,\ldots,y_k}^{I;x_0,\ldots,x_n}\ne 0$. 

Below, we will give equations involving $\sum_{y_0,\ldots,y_k}(\cdots)\ci\d_{n,k;y_0,\ldots,y_k}^{I;x_0,\ldots,x_n}$, such as \eq{ai7eq13}, \eq{ai7eq15} and \eq{ai7eq16}. Although these involve apparently infinite sums over objects $y_0,\ldots,y_k$, there are only finitely many choices for which $\d_{n,k;y_0,\ldots,y_k}^{I;x_0,\ldots,x_n}\ne 0$, so they are in effect finite sums, with no convergence issues.
\end{rem}

Here is the analogue of Lemma \ref{ai6lem2}, proved in the same way.

\begin{lem}
\label{ai7lem1}
$\d_{n,k;y_0,\ldots,y_k}^{I;x_0,\ldots,x_n}=0$ above unless $n\ge \min(k,k+2-\md{I})$.
\end{lem}

Here is the analogue of Theorem \ref{ai6thm1}, proved in \S\ref{ai101}. 

\begin{thm}
\label{ai7thm1}
In Definition\/ {\rm\ref{ai7def2},} let\/ $k=k'+k''$ with\/ $k',k''\ge 1$. Then
\smallskip

\noindent{\bf(a)} If\/ $M=0$ then
\ea
&\d_{n,k;y_0,\ldots,y_k}^{\{i_0\};x_0,\ldots,x_n}
\label{ai7eq12}\\
&=\begin{cases}
{\begin{subarray}{l}\ts (-1)^{k''(k'+n)}\d_{n-k'',k';y_0,\ldots,y_{k'}}^{\{i_0\};x_0,\ldots,x_{n-k''}}\ot{} \\ \ts
\id_{\sH^{i_0}(y_{k'},y_{k'+1})}\ot\cdots\ot\id_{\sH^{i_0}(y_{k-1},y_k)}
\end{subarray}}
& \text{if\/ $y_i=x_{n-k+i},$ $k'\le i\le k,$} \\
0, & \text{otherwise}
\end{cases}
\nonumber\\
&+\begin{cases}{\begin{subarray}{l}\ts  (-1)^{k'}\id_{\sH^{i_0}(y_0,y_1)}\ot\cdots\ot{} \\
\ts \id_{\sH^{i_0}(y_{k'-1},y_{k'})}\ot\d_{n-k',k'';y_{k'},\ldots,y_k}^{\{i_0\};x_{k'},\ldots,x_n} \end{subarray}}
& \text{if\/ $y_i=x_i,$ $0\le i\le k',$} \\
0, & \text{otherwise.}
\end{cases}
\nonumber
\ea

\noindent{\bf(b)} If\/ $M>0$ then
\ea
&\d_{n,k;y_0,\ldots,y_k}^{I;x_0,\ldots,x_n}
\label{ai7eq13}\\
&=\!\!\!\sum_{\begin{subarray}{l} 
n=n'+n'', \; n',n''\ge 0, \\ 
0<l\le M,\\
0=j_0<j_1<\cdots<j_l=M, \\
k_0,\ldots,k_l\ge 0, \; k_0=n', \; k_l=k', \\
w_0^i,\ldots,w_{k_i}^i\in\sO^i,\; 0<i<l
\end{subarray}}\begin{aligned}[t]&
(-1)^{k''(k'+l+M)+n'(k''+l+1)+\sum_{a=1}^{l-1}(j_a-a)}\\
&\bigl(\d_{k_{l-1},k_l;y_0,\ldots,y_{k'}}^{\{i_{j_{l-1}},\ldots,i_{j_l}\};w_0^{l-1},\ldots,w_{k_{l-1}}^{l-1}}\ci\d_{k_{l-2},k_{l-1};w_0^{l-1},\ldots,w_{k_{l-1}}^{l-1}}^{\{i_{j_{l-2}},\ldots,i_{j_{l-1}}\};w_0^{l-2},\ldots,w_{k_{l-2}}^{l-2}}\\
&{}\ci\cdots\ci\d_{k_0,k_1;w_0^1,\ldots,w_{k_1}^1}^{\{i_{j_0},\ldots,i_{j_1}\};x_0,\ldots,x_{n'}}\bigr)\\
&{}\ot\d_{n'',k'';y_{k'},y_{k'+1},\ldots,y_k}^{\{i_{j_0},i_{j_1},i_{j_2},\ldots,i_{j_l}\};x_{n'},x_{n'+1},\ldots,x_n},
\end{aligned}
\nonumber
\ea
where there are only finitely many nonzero terms in the sum.
\end{thm}

\subsection{\texorpdfstring{The $\iy$-categories $\ACat^{\un,\iy}_R,\ACat^{\un,\iy}_{R,K}$}{The ∞-categories A∞-Catᵘⁿ∞ᵣ,A∞-Catᵘⁿ∞ᵣₖ}}
\label{ai72}

This section is the analogue of \S\ref{ai63}. We first define a simplicial set $\hACat^{\un,\iy}_R$ of unital $A_\iy$-categories over $R$, and then define our $\iy$-category as a subset $\ACat^{\un,\iy}_R\!\subset\!\hACat^{\un,\iy}_R$. Theorem \ref{ai7thm4} extends $\hACat^{\un,\iy}_R$ to an $\iy$-{\it bi\-cat\-eg\-ory\/} $(\hACat^{\un,\iy}_R,\cT)$, as in \S\ref{ai28}, with~$\ACat^{\un,\iy}_R\!=\!\Core(\hACat^{\un,\iy}_R,\cT)$.

\begin{dfn}
\label{ai7def3}
Let $R$ be a commutative ring. We will define a simplicial set $\hACat^{\un,\iy}_R$ of unital $A_\iy$-categories over~$R$.
\smallskip

\noindent{\bf(a)} For $N\ge 0$, define the $N$-simplices $(\hACat^{\un,\iy}_R)_N$ in $\hACat^{\un,\iy}_R$ to be the set of quadruples $\bigl((\sO^i)_{i=0}^N,(f^{ij})_{0\le i\le j\le N},(\sH^i)_{i=0}^N,(\m^I)_{\es\ne I\subseteq\{0,\ldots,N\}}\bigr)$, where:
\begin{itemize}
\setlength{\itemsep}{0pt}
\setlength{\parsep}{0pt}
\item[(i)] $\sO^i$ is a set for $i=0,\ldots,N$, thought of as the set of objects $\Obj\sA^i$ of a (small) $A_\iy$-category $\sA^i$ at vertex $i$ of the $N$-simplex $\De^N$.
\item[(ii)] $f^{ij}:\sO^i\ra\sO^j$ is a map for all $0\le i\le j\le N$, thought of as the map $\fF_\Obj^{ij}:\Obj\sA^i\ra\Obj\sA^j$ in an $A_\iy$-functor $\fF^{ij}:\sA^i\ra\sA^j$ living on the edge $[i,j]$ of $\De^N$. We require these maps to satisfy $f^{ii}=\id_{\sO^i}$ for $i=0,\ldots,N$. We do {\it not\/} require that $f^{jk}\ci f^{ij}=f^{ik}$.

Suppose $\es\ne I\subseteq\{0,1,\ldots,N\}$ with $I=\{i_0,i_1,\ldots,i_M\}$ for $i_0<i_1<\cdots<i_M$. Define $f^I:\sO^{i_0}\ra\sO^{i_M}$ as in \eq{ai7eq2}.
\item[(iii)] For all $i=0,1,\ldots,N$ and $x,y\in\sO^i$, we are given a $\Z$-graded $R$-module $\sH^i(x,y)=\bigop_{d\in\Z}\sH^i(x,y)^d$, thought of as the Hom space $\Hom_{\sA^i}(x,y)$ in the $A_\iy$-category $\sA^i$. We write $\sH^i=(\sH^i(x,y))_{x,y\in\sO^i}$.
\item[(iv)] For each nonempty subset $I\subseteq\{0,\ldots,N\}$, which we write as $I=\{i_0,\ab i_1,\ab\ldots,\ab i_M\}$ for $i_0<i_1<\cdots<i_M$ with $M\ge 0$, and for all $k\ge 0$ and elements $x_0,\ldots,x_k\in \sO^{i_0}$, thought of as objects in $\sA^{i_0}$, we are given an $R$-linear map 
\ea
\m_k^{I;x_0,\ldots,x_k}:\sH^{i_0}&(x_0,x_1)\ot_R\sH^{i_0}(x_1,x_2)\ot_R\cdots \ot_R\sH^{i_0}(x_{k-1},x_k)
\nonumber\\
&\longra\sH^{i_M}(f^I(x_0),f^{i_0i_M}(x_k))
\label{ai7eq14}
\ea
graded of degree $2-M-k=3-\md{I}-k$. Here when $k=0$, by convention the empty tensor product on the top line of \eq{ai7eq14} is $R$, so $\m_0^{I;x_0}$ maps
\begin{equation*}
\m_0^{I;x_0}:R\longra\sH^{i_M}(f^I(x_0),f^{i_0i_M}(x_0))^{2-M}.
\end{equation*}
We sometimes identify $\m_0^{I;x_0}$ with $\m_0^{I;x_0}(1)\!\in\!\sH^{i_M}(f^I(x_0),f^{i_0i_M}(x_0))^{2-M}$.

In fact when $k=0$ we only want maps $\m^{I;x_0}_0$ for $\md{I}\ge 3$, so we require that $\m^{I;x_0}_0=0$ when $\md{I}=1,2$. Alternatively, we can say that $\m_k^{I;x_0,\ldots,x_k}$ in \eq{ai7eq14} should only be defined for $k\ge 1$ when $\md{I}=1,2$ and for $k\ge 0$ when $\md{I}\ge 3$, and then we should exclude terms in \eq{ai7eq16} and other equations below involving $\m^{I;x_0}_0$ for~$\md{I}=1,2$.

We write $\m^I_k=(\m_k^{I;x_0,\ldots,x_k})_{x_0,\ldots,x_k\in\sO^{i_0}}$ and $\m^I=(\m^I_k)_{k\ge 0}$.

We require these maps $\m_k^{I;x_0,\ldots,x_k}$ to satisfy the following identities generalizing \eq{ai6eq18} for all $\es\ne I=\{i_0,\ab i_1,\ab\ldots,\ab i_M\}\subseteq\{0,1,\ldots,N\}$, $n\ge 0$, and $x_0,\ldots,x_n\in \sO^{i_0}$, firstly if $M=0$
\e
0=\sum_{\begin{subarray}{l} k=1,\ldots,n, \\
y_0,\ldots,y_k\in\sO^{i_0} \end{subarray}}\m_k^{\{i_0\};y_0,\ldots,y_k}\ci\d_{n,k;y_0,\ldots,y_k}^{\{i_0\};x_0,\ldots,x_n},
\label{ai7eq15}
\e
and secondly if $M>0$
\ea
&0=\sum_{\begin{subarray}{l} k=1,\ldots,n \\
y_0,\ldots,y_k\in\sO^{i_0}  \end{subarray}}\m_k^{I;y_0,\ldots,y_k}\ci\d_{n,k;y_0,\ldots,y_k}^{\{i_0\};x_0,\ldots,x_n}
\label{ai7eq16}\\
&+(-1)^M\sum_{\begin{subarray}{l} k=1,\ldots,n+M-1,\\
y_0,\ldots,y_k\in\sO^{i_M} \end{subarray}}\m_k^{\{i_M\};y_0,\ldots,y_k}\ci\d_{n,k;y_0,\ldots,y_k}^{I;x_0,\ldots,x_n}.
\nonumber
\ea
Here $\d_{n,k;y_0,\ldots,y_k}^{I;x_0,\ldots,x_n}$ is defined using the $\m_l^{J;w_0,\ldots,w_l}$ as in Definition \ref{ai7def2}. In both \eq{ai7eq15} and \eq{ai7eq16} there are only finitely many nonzero terms, and all nonzero terms map to $\sH^{i_M}(f^I(x_0),f^{i_0i_M}(x_n))$, so \eq{ai7eq15}--\eq{ai7eq16} should be regarded as equations in this space.
\item[(v)] There should exist elements $\boo^i_x\in\sH^i(x,x)^0$ for each $i=0,\ldots,N$ and $x\in\sO^i$, which are uniquely determined by \eq{ai7eq17}, satisfying
\begin{gather}
\begin{aligned}
&\m_2^{\{i\};x_0,x_0,x_1}(\boo^i_{x_0}\ot a)=\m_2^{\{i\};x_0,x_1,x_1}(a\ot \boo^i_{x_1})=a \\
&\text{for all $i=0,\ldots,N$, $x_0,x_1\in\sO^i$ and $a\in\sH^i(x_0,x_1)$.}
\end{aligned}
\label{ai7eq17}\\
\m_1^{\{i,j\};x,x}(\boo_x^i)=\boo_{f^{ij}(x)}^j\quad\text{for all $0\le i<j\le N$ and $x\in\sO^i$,}
\label{ai7eq18}\\
\begin{aligned}
&\m_{k+1}^{I;x_0,\ldots,x_j,x_j,\ldots,x_k}(a_1\ot\cdots\ot a_j\ot \boo^{i_0}_{x_j}\ot a_{j+1}\ot \cdots \ot a_k)=0\\
&\text{for all $I=\{i_0<i_1<\cdots<i_M\}\subseteq\{0,\ldots,N\}$, $0\le j\le k$,}\\
&\text{$x_0,\ldots,x_k\in\sO^{i_0}$, and $a_l\in\sH^{i_0}(x_{l-1},x_l)$, $l=1,\ldots,k$,} \\
&\text{provided $(M,k)\ne (0,1),(1,0)$. }
\end{aligned}
\label{ai7eq19}
\end{gather}
The exclusion $(M,k)\ne (0,1),(1,0)$ in \eq{ai7eq19} is because \eq{ai7eq17} prescribes the case $(M,k)=(0,1)$, and \eq{ai7eq18} the case~$(M,k)=(1,0)$.
\end{itemize}

\noindent{\bf(b)} If $\al:[\bar  N]\ra[N]$ is a morphism in $\De$, that is, $\al:\{0,\ldots,\bar  N\}\ra\{0,\ldots,N\}$ is nondecreasing, let $\hACat^{\un,\iy}_R(\al):(\hACat^{\un,\iy}_R)_N\!\ra\!\hACat^{\un,\iy}_{R,\bar  N}$ map $\bigl((\sO^i)_{i=0}^N,\ab(f^{ij})_{0\le i\le j\le N},(\sH^i)_{i=0}^N,(\m^I)_{\es\ne I\subseteq\{0,\ldots,N\}}\bigr)$ to $\bigl((\bar \sO^i)_{i=0}^{\bar  N}$, $(\bar  f^{ij})_{0\le i\le j\le\ti N},$ $(\bar\sH{}^i)_{i=0}^{\bar N},\ab(\bar\m{}^{\bar I})_{\es\ne\bar I\subseteq\{0,\ldots,\bar N\}}\bigr)$, where $\bar\sO^i=\sO^{\al(i)}$ for $i=0,\ldots,\bar N$, and $\bar f^{ij}=f^{\al(i)\al(j)}$ for $0\le i\le j\le\bar N$, and $\bar\sH^i=\sH^{\al(i)}$ for $i=0,\ldots,\bar N$. To define $\bar\m{}^{\bar I}$, if $\bar I=\{\bar\imath_0<\bar\imath_1<\cdots<\bar\imath_{\bar M}\}\subseteq\{0,1,\ldots,\bar N\}$, $k\ge 0$ and $x_0,\ldots,x_k\in\bar\sO^{\bar\imath_0}=\sO^{\al(\bar\imath_0)}$, generalizing \eq{ai6eq22} we define
\e
\bar\m_k^{\bar I;x_0,\ldots,x_k}=\begin{cases}
\m_k^{\al(\bar I);x_0,\ldots,x_k}, & \text{if $\al\vert_{\bar I}$ is injective,} \\
\id_{\sH^i(x_0,x_1)}, & {\begin{subarray}{l} \ts \text{if $k\!=\!1$ and $\bar I\!=\!\{\bar\imath_0,\bar\imath_1\}$ with} \\ \ts \text{$\bar\imath_0\!\ne\! \bar\imath_1$, $\al(\bar\imath_0)\!=\!\al(\bar\imath_1)\!=\!i$,} \end{subarray}} \\[8pt]
\boo^{\bar\imath_2}_{f^{\bar\imath_0\bar\imath_2}(x_0)}, & {\begin{subarray}{l}\ts \text{if $k=0$, $\bar I=\{\bar\imath_0\!<\!\bar\imath_1\!<\!\bar\imath_2\}$} \\ \ts\text{and $\al\vert_{\bar I}$ is not injective,}\end{subarray}} \\[3pt]
0, & \text{otherwise.}
\end{cases}
\label{ai7eq20}
\e

We will show in Proposition \ref{ai7prop1} that $\hACat^{\un,\iy}_R(\al)$ maps to $\hACat^{\un,\iy}_{R,\bar N}$. Clearly, if $\be:[\hat N]\ra[\bar N]$ is another morphism in $\De$ then
\begin{align*}
\hACat^{\un,\iy}_R(\al\ci\be)&=\hACat^{\un,\iy}_R(\be)\ci\hACat^{\un,\iy}_R(\al):\\
&(\hACat^{\un,\iy}_R)_N\longra\hACat^{\un,\iy}_{R,\hat N}.
\end{align*} 
Also $\hACat^{\un,\iy}_R(\id_{[N]})=\id_{(\hACat^{\un,\iy}_R)_N}$. Hence $\hACat^{\un,\iy}_R:\De^\rop\ra\Set$ is a functor, that is, $\hACat^{\un,\iy}_R$ is a simplicial set.
\end{dfn}

\begin{rem}
\label{ai7rem2}
{\bf(On Set Theory foundations.)}
Recall from Remark \ref{ai1rem2} that we fix a Grothendieck universe $\sU$, and implicitly assume all sets in our theory are $\sU$-small. As in Remark \ref{ai6rem3}, without this convention, the $N$-simplices $(\hACat^{\un,\iy}_R)_N$ above would be proper classes, rather than sets. This would cause Set Theory problems.
\end{rem}

The following analogue of Proposition \ref{ai6prop1} will be proved in~\S\ref{ai102}.

\begin{prop}
\label{ai7prop1}
The data $\bigl((\bar\sO^i)_{i=0}^{\bar N},\!(\bar f^{ij})_{0\le i\le j\le\bar N},\!(\bar\sH{}^i)_{i=0}^{\bar N},\!(\bar\m{}^{\bar I})_{\es\ne\bar I\subseteq\{0,\ldots,\bar N\}}\bigr)$ in Definition\/ {\rm\ref{ai7def3}(b)} lies in $(\hACat^{\un,\iy}_R)_{\bar N}$.
\end{prop}

The following analogue of Theorem \ref{ai6thm2} will be proved in~\S\ref{ai103}.

\begin{thm}
\label{ai7thm2}
Let\/ $\bigl((\sO^i)_{i=0}^N,(f^{ij})_{0\le i\le j\le N},(\sH^i)_{i=0}^N,(\m^I)_{\es\ne I\subseteq\{0,\ldots,N\}}\bigr)$ be an $N$-simplex in $\hACat^{\un,\iy}_R$ as in Definition\/ {\rm\ref{ai7def3},} and define maps\/ $\d_{n,k;y_0,\ldots,y_k}^{I;x_0,\ldots,x_n}$ from the $\m_k^{I;x_0,\ldots,x_k}$ as in\/ {\rm\eq{ai7eq6}--\eq{ai7eq9}}. Then for each nonempty subset\/ $I\subseteq\{0,\ldots,N\},$ written $I=\{i_0,\ab i_1,\ab\ldots,\ab i_M\}$ for $i_0<i_1<\cdots<i_M,$ and all\/ $k,n\ge 0,$ $x_0,\ldots,x_n\in\sO^{i_0}$ and\/ $y_0,\ldots,y_k\in\sO^{i_M},$ in $R$-linear maps as in \eq{ai7eq6} we~have
\e
\begin{aligned}
&\sum_{\begin{subarray}{l} k\le l\le n, \\ w_0,\ldots,w_l\in\sO^{i_0} \end{subarray}}^n\d_{l,k;y_0,\ldots,y_k}^{\{i_0\};w_0,\ldots,w_l}\ci\d_{n,l;w_0,\ldots,w_l}^{\{i_0\};x_0,\ldots,x_n}=0 && \text{if\/ $M=0,$} \\
&\begin{aligned}[t]& \sum_{\begin{subarray}{l} 1\le l\le n, \\ w_0,\ldots,w_l\in\sO^{i_0} \end{subarray}} \d_{l,k;y_0,\ldots,y_k}^{I;w_0,\ldots,w_l}\ci\d_{n,l;w_0,\ldots,w_l}^{\{i_0\};x_0,\ldots,x_n}\\
&+(-1)^M\sum_{\begin{subarray}{l} k\le l\le n+M-1, \\  w_0,\ldots,w_l\in\sO^{i_M} \end{subarray}}\d_{l,k}^{\{i_M\}}\ci\d_{n,l;w_0,\ldots,w_l}^{I;x_0,\ldots,x_n}&=0\end{aligned} && \text{if\/ $M>0,$}
\end{aligned}
\label{ai7eq21}
\e
where there are only finitely many nonzero terms in the sums.
\end{thm}

Here is the analogue of Example \ref{ai6ex3}.

\begin{ex}
\label{ai7ex2}
{\bf(a)} A 0-simplex $(\sO^0,f^{00}=\id_{\sO^0},\sH^0,\m^{\{0\}})$ of $\hACat^{\un,\iy}_R$ is equivalent to a (small) unital $A_\iy$-category $\sA^0$ over $R$, as in Definition \ref{ai4def2}, with objects $\Obj\sA^0=\sO^0$, Hom spaces $\Hom_{\sA^0}(x,y)=\sH^0(x,y)$, morphisms $(\m_k^{\{0\};x_0,\ldots,x_k})_{k\ge 1}$, and units~$\boo^0_x\in \Hom_{\sA^0}^0(x,x)=\sH^0(x,x)^0$.
\smallskip

\noindent{\bf(b)} A 1-simplex of $\hACat^{\un,\iy}_R$ is equivalent to a unital $A_\iy$-functor $\fF^{01}:\sA^0\ra\sA^1$ of unital $A_\iy$-categories, as in Definition \ref{ai4def6}, with object map $\fF_\Obj^{01}=f^{01}$ and morphisms $\fF^{01;x_0,\ldots,x_k}_k=\m_k^{\{0,1\};x_0,\ldots,x_k}$. 
\smallskip

\noindent{\bf(c)} A 2-simplex of $\hACat^{\un,\iy}_R$ is equivalent to unital $A_\iy$-categories $\sA^0,\ab\sA^1,\ab\sA^2$ over $R$ as in {\bf(a)}, and unital $A_\iy$-functors $\fF^{ij}:\sA^i\ra\sA^j$ for $(i,j)=(0,1),\ab(0,2),\ab(1,2)$ as in {\bf(b)}, and $R$-linear maps for all $k\ge 0$ and $x_0,\ldots,x_k\in\sO^0$
\begin{align*}
\m_k^{\{0,1,2\};x_0,\ldots,x_k}:\sH^0&(x_0,x_1)\ot_R\sH^0(x_1,x_2)\ot_R\cdots \ot_R\sH^0(x_{k-1},x_k)
\nonumber\\
&\longra\sH^2(f^{12}\ci f^{01}(x_0),f^{02}(x_k))
\end{align*}
graded of degree $-k$, satisfying \eq{ai7eq19} and for all $n\ge 0$ and $x_0,\ldots,x_n\in\sO^0$
\ea
&\sum_{\begin{subarray}{l} i,k\ge 0,\; j\ge 1:\\ i+j+k=n \end{subarray}}\begin{aligned}[t]&
(-1)^{i+jk}\m^{\{0,1,2\};x_0,x_1,\ldots,x_i,x_{i+j},x_{i+j+1},\ldots,x_n}_{i+k+1}\ci\\
&\bigl(\id_{\sH^0(x_0,x_1)}\ot\cdots\ot \id_{\sH^0(x_{i-1},x_i)}\ot\m^{\{0\};x_i,x_{i+1},\ldots,x_{i+j}}_j\\
&\ot \id_{\sH^0(x_{i+j},x_{i+j+1})}\ot\cdots\ot \id_{\sH^0(x_{n-1},x_n)}\bigr)\end{aligned}
\label{ai7eq22}\\
&+\sum_{\!\!\!\!\!\!\begin{subarray}{l} 
1\le l\le k\le n+1, \\ 
n=n_1+\cdots+n_k, \\
n_l\ge 0,\; n_i\ge 1,\; i\ne l. \\
\text{Set $p_j=\sum_{i=1}^jn_i$,} \\ 
0\le j\le k \end{subarray}}
\begin{aligned}[t]&(-1)^{k+1}\cdot(-1)^{\sum_{i=1}^{l-1}(n_i+1)}\cdot(-1)^{\sum_{i=1}^k(n_i+1)(k-i)}\\
&\m^{\{2\};f^{\{0,1,2\}}(x_{p_0}),f^{\{0,1,2\}}(x_{p_1}),\ldots,f^{\{0,1,2\}}(x_{p_{l-1}}),f^{02}(x_{p_l}),\ldots,f^{02}(x_{p_k})}_k\\
&\ci\bigl[(\m^{\{1,2\}}\!\ci\!\m^{\{0,1\}})^{x_{p_0},\ldots,x_{p_1}}_{n_1}\!\ot\!\cdots\!\ot\!(\m^{\{1,2\}}\!\ci\!\m^{\{0,1\}})^{x_{p_{l-2}},\ldots,x_{p_{l-1}}}_{n_{l-1}}\\
&\ot \m^{\{0,1,2\};x_{p_{l-1}},\ldots,x_{p_l}}_{n_l}\ot\m^{\{0,2\};x_{p_l},\ldots,x_{p_{l+1}}}_{n_{l+1}}\ot\\
&\cdots\ot\m^{\{0,2\};x_{p_{k-1}},\ldots,x_{p_k}}_{n_k}\bigr]=0.
\end{aligned}
\nonumber
\ea
This is an analogue of \eq{ai6eq24} and \eq{ai6eq25} which we have simplified by writing in terms of the composition $\m^{\{1,2\}}\ci\m^{\{0,1\}}$ of the data $\m^{\{0,1\}},\m^{\{1,2\}}$, corresponding to working with the composition of $A_\iy$-functors $\fF^{12}\ci\fF^{01}$. Comparing \eq{ai4eq15} and \eq{ai7eq22} and using \eq{ai7eq19}, we see that $\m^{\{0,1,2\}}$ is equivalent to a {\it unital\/ $A_\iy$-natural transformation\/} $\fH^{012}:\fF^{12}\ci\fF^{01}\Ra\fF^{02}$, in the sense of Definition~\ref{ai4def13}.
\end{ex}

Here is the analogue of Definition~\ref{ai6def10}.

\begin{dfn}
\label{ai7def4}
Define a simplicial subset $\ACat^{\un,\iy}_R\subset\hACat^{\un,\iy}_R$, by an $N$-simplex $\bigl((\sO^i)_{i=0}^N,(f^{ij})_{0\le i\le j\le N},(\sH^i)_{i=0}^N,(\m^I)_{\es\ne I\subseteq\{0,\ldots,N\}}\bigr)$ in $\hACat^{\un,\iy}_R$ lies in $\ACat^{\un,\iy}_R$ if and only if whenever $0\le i_0<i_1<i_2\le N$, any (hence all) of the following four equivalent conditions hold:
\begin{itemize}
\setlength{\itemsep}{0pt}
\setlength{\parsep}{0pt}
\item[(i)] For all objects $x_0\in\sO^{i_0}$, the following morphism is an isomorphism in the cohomology category $H^0(\sA^{i_2})$:
\begin{equation*}
[\m^{\{i_0,i_1,i_2\};x_0}_0(1)]:f^{i_1i_2}\ci f^{i_0i_1}(x_0)\longra f^{i_2i_0}(x_0).
\end{equation*}
\item[(ii)] $[\m^{\{i_0,i_1,i_2\}}_0(1)]:H^0(\fF_1^{\{i_1,i_2\}}\ci\fF_1^{\{i_0,i_1\}})\Ra H^0(\fF_1^{\{i_0,i_2\}})$ is a natural isomorphism of functors $H^0(\sA^{i_0})\ra H^0(\sA^{i_2})$.
\item[(iii)] The unital $A_\iy$-natural transformation $\fH^{i_0i_1i_2}:\fF^{i_1i_2}\ci\fF^{i_0i_1}\Ra\fF^{i_0i_2}$ in Example \ref{ai7ex2}(c) is a unital $A_\iy$-natural isomorphism, in the sense of Definition~\ref{ai4def13}.
\item[(iv)] $[\m^{\{i_0,i_1,i_2\}}]:\fF^{i_1i_2}\ci\fF^{i_0i_1}\ra\fF^{i_0i_2}$ is an isomorphism in the cohomology category~$H^0(\Fun_{A_\iy}^\un(\sA^{i_0},\sA^{i_2}))$.
\end{itemize}

Here (i) and (ii) are equivalent by definition of natural isomorphism of ordinary functors, and (iii) and (iv) are equivalent by definition of unital $A_\iy$-natural isomorphisms in Definition \ref{ai4def13}, and (ii) and (iii) are equivalent by Proposition \ref{ai4prop3}, so (i)--(iv) are equivalent. Part (i) is the analogue of the condition that $[\m^{\{i_0,i_1,i_2\}}_0(1)]$ is invertible for $A_\iy$-algebras in Definition \ref{ai6def10}. If $\al:[\bar N]\ra[N]$ is a morphism in $\De$ then $\hACat^{\un,\iy}_R(\al)$ in Definition \ref{ai7def3}(b) maps $(\ACat^{\un,\iy}_R)_N\ra(\ACat^{\un,\iy}_R)_{\bar N}$, so $\ACat^{\un,\iy}_R$ is a simplicial set.
\end{dfn}

The following analogue of Theorem \ref{ai6thm3}, one of our main results, is proved in \S\ref{ai104}. The proof is long and difficult.

\begin{thm}
\label{ai7thm3}
In Definition\/ {\rm\ref{ai7def4},} $\ACat^{\un,\iy}_R$ is an $\iy$-category.
\end{thm}

\begin{rem}
\label{ai7rem3}
In the generalization to {\it filtered\/} $A_\iy$-categories in Chapter \ref{ai8}, the analogues of Definition \ref{ai7def4}(i),(ii) do not make sense as we can have $(\m_1^{\{i_2\}})^2\ne 0$, so the cohomology category $H^0(\sA^{i_2})$ is not defined. However, the analogues of Definition \ref{ai7def4}(iii),(iv) do make sense. We have written the proof of Theorem \ref{ai7thm3} using Definition \ref{ai7def4}(iv), so that it extends to the filtered case.
\end{rem}

\begin{dfn}
\label{ai7def5}
Let $\bigl((\sO^i),(f^{ij}),(\sH^i),(\m^I)\bigr)$ be an $N$-simplex in $\ACat^{\un,\iy}_R$ as in Definitions \ref{ai7def3} and \ref{ai7def4}, for $N\ge 2$. If $N>2$, we call $\bigl((\sO^i),\ldots,(\m^I)\bigr)$ {\it special\/} if $\m^{\{0,\ldots,N\}}=0$. If $N=2$ we call $\bigl((\sO^i),\ldots,(\m^I)\bigr)$ {\it special\/} if in the notation of Example \ref{ai7ex2}(c), we have $\fF^{02}=\fF^{21}\ci\fF^{01}$ and~$\fH^{012}=\id_{\fF^{02}}$.
\end{dfn}

In the proof of Theorem \ref{ai7thm3} in \S\ref{ai104}, to prove $\ACat^{\un,\iy}_R$ is an $\iy$-category, as in Definition \ref{ai2def4}, for every $0<i<N$ and simplicial set map $f:\La^N_i\ra\ACat^{\un,\iy}_R$, we show that $f$ can be extended to an $N$-simplex $\ti f:\De^N\ra\ACat^{\un,\iy}_R$ as in \eq{ai2eq5}. The proof also yields:

\begin{prop}
\label{ai7prop2}
Let\/ $f:\La^N_i\ra\ACat^{\un,\iy}_R$ be a map of simplicial sets $0<i<N$. Then we can extend $f$ to a \begin{bfseries}unique\end{bfseries} map $\ti f:\De^N\ra\ACat^{\un,\iy}_R$ as in {\rm\eq{ai2eq5},} such that\/ $\ti f$ is a \begin{bfseries}special\end{bfseries}\/ $N$-simplex in the sense of Definition\/~{\rm\ref{ai7def5}}.
\end{prop}

In \S\ref{ai28}, following Lurie \cite[\S 4]{Luri4}, we defined $\iy$-{\it bicategories}, a kind of $(\iy,2)$-category. The next theorem will be proved in~\S\ref{ai105}.

\begin{thm}
\label{ai7thm4}
In the simplicial set\/ $\hACat^{\un,\iy}_R$ in Definition\/ {\rm\ref{ai7def3},} define a $2$-simplex to be \begin{bfseries}thin\end{bfseries} if the unital\/ $A_\iy$-natural transformation $\fH^{012}:\fF^{12}\ci\fF^{01}\Ra\fF^{02}$ in Example\/ {\rm\ref{ai7ex2}(c)} is an\/ $A_\iy$-natural isomorphism, that is, if it is a $2$-simplex in $\ACat^{\un,\iy}_R\subset\hACat^{\un,\iy}_R$. Write $\cT\subset\hACat^{\un,\iy}_R([2])$ for the set of thin $2$-simplices. Then $(\hACat^{\un,\iy}_R,\cT)$ is an \begin{bfseries}$\iy$-bicategory\end{bfseries} in the sense of\/ {\rm\S\ref{ai28}}. We have $\ACat^{\un,\iy}_R=\Core(\hACat^{\un,\iy}_R,\cT),$ as in Definition\/~{\rm\ref{ai2def22}}.
\end{thm}

\begin{rem}
\label{ai7rem4}
{\bf(a)} Theorem \ref{ai7thm3} follows trivially from Theorem \ref{ai7thm4}, as the core of an $\iy$-bicategory is an $\iy$-category. But we use the proof of Theorem \ref{ai7thm3} in that of Theorem \ref{ai7thm4}, so we prove Theorem \ref{ai7thm3} first.
\smallskip

\noindent{\bf(b)} Other than proving that $(\hACat^{\un,\iy}_R,\cT)$ is an $\iy$-bicategory, we will not take these ideas much further in this book; our focus will be on $\iy$-categories, not $(\iy,2)$-categories. But maybe someone will find applications for this $\iy$-bicategory in future. At the least, the author thinks $(\hACat^{\un,\iy}_R,\cT)$ is interesting as an example of an $\iy$-bicategory which appears `in nature'.
\smallskip

\noindent{\bf(c)} The most basic difference between $\ACat^{\un,\iy}_R$ and $(\hACat^{\un,\iy}_R,\cT)$ is that $\ACat^{\un,\iy}_R$ is an $\iy$-category ($(\iy,1)$-category) in which 2-morphisms are $A_\iy$-natural isomorphisms, and $(\hACat^{\un,\iy}_R,\cT)$ is an $\iy$-bicategory ($(\iy,2)$-category) in which 2-morphisms are unital $A_\iy$-natural transformations.
\smallskip

\noindent{\bf(d)} In \S\ref{ai411} we discussed Faonte's proposal \cite{Faon2} for defining an $(\iy,2)$-category of unital $A_\iy$-categories over $R$, as a fibrant Segal category in the model category $({\rm SSet},{\rm Joy})$, at least when $R$ is a field of characteristic zero. Although we will not prove this, we expect that under the equivalences between different models of $(\iy,2)$-category in Gagna--Harpaz--Lanari \cite{GHL}, our $\iy$-bicategory $(\hACat^{\un,\iy}_R,\cT)$ should be equivalent to Faonte's $(\iy,2)$-category. As evidence for this, Theorem \ref{ai9thm1}(a) below shows left mapping spaces in $(\hACat^{\un,\iy}_R,\cT)$ are given by
\begin{equation*}
\Hom_{\hACat^{\un,\iy}_R}^L(\sA,\sB)\cong \sN_{A_\iy}(\Fun_{A_\iy}^\un(\sA,\sB)),
\end{equation*}
where the right hand side is the mapping space in Faonte's $(\iy,2)$-category.
\smallskip

\noindent{\bf(e)} We will use Theorem \ref{ai7thm4} in the proof of Theorem \ref{ai9thm5} in~\S\ref{ai115}.
\smallskip

\noindent{\bf(f)} In Definition \ref{ai2def23} we defined {\it strong\/ $\iy$-bicategories}, which are just called $\iy$-bicategories in Lurie \cite[\S 5.4]{Luri5}. The author expects that $(\hACat^{\un,\iy}_R,\cT)$ is {\it not\/} a strong $\iy$-bicategory. That is, there may be 2-simplices $t:\De^2\ra\hACat^{\un,\iy}_R$ which satisfy Definition \ref{ai2def21}(a), but for which the unital $A_\iy$-natural transformation $\fH^{012}:\fF^{12}\ci\fF^{01}\Ra\fF^{02}$ in Example \ref{ai7ex2}(c) is not an $A_\iy$-natural isomorphism. Using Definition \ref{ai2def21}(a), it is not difficult to show that $\fH^{012}$ has a left inverse in the cohomology category $H^0(\Fun_{A_\iy}^\un(\sA^0,\ab\sA^2))$, but it does not seem possible to show that $\fH^{012}$ has a right inverse unless $\fF^{01}$ or $\fF^{12}$ is a homotopy equivalence. This is our reason for choosing the $\iy$-bicategory theory from \cite[\S 4]{Luri4} rather than from~\cite[\S 5.4]{Luri5}.
\end{rem}

Here is the analogue of Definition \ref{ai6def11}.

\begin{dfn}
\label{ai7def6}
Let $K\in\{3,4,5,\ldots,\iy\}$. Define simplicial subsets
\begin{gather*}
\ACat^{\un,\iy}_{R,K}\subset\ACat^{\un,\iy}_R\subset\hACat^{\un,\iy}_R
\end{gather*}
such that an $N$-simplex $\bigl((\sO^i)_{i=0}^N,(f^{ij})_{0\le i\le j\le N},(\sH^i)_{i=0}^N,(\m^I)_{\es\ne I\subseteq\{0,\ldots,N\}}\bigr)$ in $\hACat^{\un,\iy}_R$ in Definition \ref{ai7def3}(a) lies in $\ACat^{\un,\iy}_{R,K}$ if $f^{ik}=f^{jk}\ci f^{ij}$ for all $0\le i\le j\le k\le N$, and for all $0\le i_0<i_1<i_2\le N$ and $x_0\in\sO^{x_0}$ we have $\m_0^{\{i_0,i_1,i_2\};x_0}:\la\mapsto\la\boo^{i_2}_{f^{i_0i_2}(x_0)}$, and $\m_0^{I;x_0}=0$ for all $\es\ne I\subseteq\{0,\ldots,N\}$ with $\md{I}\le K$ and $\md{I}\ne 3$ and~$x_0\in\sO^{\min I}$. 

These conditions are preserved by the maps in Definition \ref{ai7def3}(b), and so define a simplicial subset $\ACat^{\un,\iy}_{R,K}\subset\hACat^{\un,\iy}_R$. As $[\m_0^{\{i_0,i_1,i_2\};x_0}(1)]=[\boo^{i_2}_{f^{i_0i_2}(x_0)}]$ is invertible in $H^0(\sA^{i_2})$ for all $x_0\in\sO^{i_0}$ we see that $\ACat^{\un,\iy}_{R,K}\subset\ACat^{\un,\iy}_R$ by Definition \ref{ai7def4}(i). Clearly $\ACat^{\un,\iy}_{R,K}\subset\ACat^{\un,\iy}_{R,L}$ if $L\le K$ in $\{3,4,5,\ldots,\iy\}$, giving the inclusion relations in the first column of \eq{ai7eq1}. Theorem \ref{ai7thm7} will show that $\ACat^{\un,\iy}_{R,K}$ is an $\iy$-category. We say that $\ACat^{\un,\iy}_{R,K}$ has {\it type\/}~$K$. 
\end{dfn}

\begin{rem}
\label{ai7rem5}
{\bf(a)} Note the condition that $f^{ik}=f^{jk}\ci f^{ij}$ in Definition \ref{ai7def6}, that is, the $A_\iy$-functors $\fF^{ik}$ and $\fF^{jk}\ci\fF^{ij}$ must act the same way on objects. We need this as $\m_0^{\{i_0,i_1,i_2\};x_0}$ maps $R\ra\sH^{i_2}(f^{i_1i_2}\ci f^{i_0i_1}(x_0),f^{i_0i_2}(x_0))$, so defining $\m_0^{\{i_0,i_1,i_2\};x_0}:\la\mapsto\la\boo^{i_2}_{f^{i_0i_2}(x_0)}$ only makes sense if~$f^{i_0i_2}(x_0)=f^{i_1i_2}\ci f^{i_0i_1}(x_0)$.

This condition is trivial for $A_\iy$-algebras, when we take $f^{ij}=\id_{\{*\}}$ as in Definition \ref{ai7def7}, but for $A_\iy$-categories, it is {\it very restrictive}. It means that we do not have an interesting notion of natural isomorphism of functors in $\ACat^{\un,\iy}_{R,K}$ for $K\ge 3$. In Example \ref{ai7ex2}(c), a 2-simplex in $\ACat^{\un,\iy}_{R,K}$ defines unital $A_\iy$-categories $\sA^0,\ab\sA^1,\ab\sA^2$ over $R$ as in {\bf(a)}, and unital $A_\iy$-functors $\fF^{ij}:\sA^i\ra\sA^j$ for $(i,j)=(0,1),(0,2),(1,2)$, such that the cohomology functors $H^0(\fF^{12}_1\ci\fF^{01}_1),H^0(\fF^{02}_1):H^0(\sA^0)\ra H^0(\sA^2)$ are not just naturally isomorphic, as in Definition \ref{ai7def4}(ii), but actually equal.
\smallskip

\noindent{\bf(b)} In \cite{CJR}, with Guillem Cazassus and Alexander Ritter, the author is planning an application  in Symplectic Geometry of $\ACat^\iy_{R,\iy}$ in \S\ref{ai73} below, which is a non-unital version of $\ACat^{\un,\iy}_{R,\iy}$. The aim is to give a new construction of the Fukaya $A_\iy$-category of (say exact) Lagrangians $L$ in a symplectic manifold $(X,\om)$. To do this, we must make some arbitrary choices, such as an almost complex structure $J$ on $X$ compatible with $\om$. These arbitrary choices are organized into a contractible Kan complex $\sK$, and our definition of Fukaya category yields a map of simplicial sets~$\sK\ra\ACat^\iy_{R,\iy}$. 

The point we want to make here is that though different choices $\ka$ in $\sK$ yield different $A_\iy$-categories $\sA_\ka$, the set of objects in $\sA_\ka$, which is a set $\sO$ of suitable Lagrangians $L$ in $(X,\om)$, is independent of $\ka$. So it is natural in this case to require $A_\iy$-functors $\sA_{\ka_1}\ra\sA_{\ka_2}$ to induce the same map $\id_\sO$ on objects. 

\smallskip

\noindent{\bf(c)} The condition $f^{ik}=f^{jk}\ci f^{ij}$ simplifies parts of Definition \ref{ai7def3}: the maps $f^I:\sO^{i_0}\ra\sO^{i_M}$ in Definition \ref{ai7def3}(a)(ii) satisfy~$f^I=f^{i_0i_M}$.
\end{rem}

\begin{dfn}
\label{ai7def7}
For $\AAlg^{\un,\iy}_R$ as in \S\ref{ai63}, define an $\iy$-functor $I_{\AAlg^\un}^{\ACat^\un}:\AAlg^{\un,\iy}_R\ra\ACat^{\un,\iy}_R$ to act on $N$-simplices by
\begin{align*}
&I_{\AAlg^\un}^{\ACat^\un}([N]):\bigl((A^i)_{i=0}^N,(\m^I)_{\es\ne I\subseteq\{0,\ldots,N\}}\bigr)\\
&\quad\longmapsto\bigl((\{*\})_{i=0}^N,(\id_{\{*\}})_{0\le i\le j\le N},(A^i)_{i=0}^N,(\m^I)_{\es\ne I\subseteq\{0,\ldots,N\}}\bigr).
\end{align*}
That is, we take $\sO^i$ to be the one point set $\{*\}$, and $f^{ij}=\id_{\{*\}}$, and $\sH^i(*,*)=A^i$, and $\m^{I;*,\ldots,*}_k=\m^I_k$. This just turns $A_\iy$-algebras into $A_\iy$-categories with one object $*$. As the definition of $\ACat^{\un,\iy}_R$ extends that of $\AAlg^{\un,\iy}_R$ from one object to many, this is a strict isomorphism with the full $\iy$-subcategory of $\ACat^{\un,\iy}_R$ with objects unital $A_\iy$-categories $\sA$ with one object~$*$.

This restricts to an $\iy$-functor $I_{\AAlg^\un}^{\ACat^\un}:\AAlg^{\un,\iy}_{R,K}\!\ra\!\ACat^{\un,\iy}_{R,K}$ for all $K\in\{3,4,5,\ldots,\iy\}$.

Results for $\AAlg^{\un,\iy}_R$ such as Proposition \ref{ai6prop1} and Theorems \ref{ai6thm2} and \ref{ai6thm3} can be deduced from their analogues for $\ACat^{\un,\iy}_R$ by restricting to the image of $I_{\AAlg^\un}^{\ACat^\un}$, and changing notation appropriately.
\end{dfn}

\subsection{\texorpdfstring{Other $\iy$-categories of $A_\iy$-categories}{Other ∞-categories of A∞-categories}}
\label{ai73}

This section is the analogue of \S\ref{ai64}. We first define $\iy$-categories $\dACat^{\un,\iy}_{R,K}$ isomorphic to $\ACat^{\un,\iy}_{R,K}$ for $K=3,4,\ldots,\iy$, extending Definition~\ref{ai6def12}.

Let $\bigl((\sO^i)_{i=0}^N,(f^{ij})_{0\le i\le j\le N},(\sH^i)_{i=0}^N,(\m^I)_{\es\ne I\subseteq\{0,\ldots,N\}}\bigr)$ be an $N$-simplex in $\ACat^{\un,\iy}_{R,K}\!$ in Definition \ref{ai7def6}. This will correspond to an $N$-simplex $\bigl((\sO^i)_{i=0}^N,\ab(f^{ij})_{0\le i\le j\le N},(\sH^i)_{i=0}^N,(\dot\m^I)_{\es\ne I\subseteq\{0,\ldots,N\}}\bigr)$ in $\dACat^{\un,\iy}_{R,K}$, where we define
\e
\dot\m_k^{\{i_0,\ldots,i_M\};x_0,\ldots,x_k}=\begin{cases} \m_0^{\{i_0,i_1,i_2\};x_0}-\boo^{i_2}_{f^{i_0i_2}(x_0)}, & \text{$k=0$ and $M=2$,}  \\
\m_k^{\{i_0,\ldots,i_M\};x_0,\ldots,x_k}, & \text{otherwise,} \end{cases}
\label{ai7eq23}
\e
extending \eq{ai6eq27}. This only makes sense if $f^{i_0i_2}(x_0)=f^{i_1i_2}\ci f^{i_0i_1}(x_0)$. So we only do it for the $\ACat^{\un,\iy}_{R,K}$ for $K=3,4,\ldots,\iy$, which include the condition $f^{ik}=f^{jk}\ci f^{ij}$ for all $0\le i\le j\le k\le N$ as in Definitions \ref{ai7def6} and Remark \ref{ai7rem5}. The construction of this section does not work for $\ACat^{\un,\iy}_R=\ACat^{\un,\iy}_{R,2}$ in Definitions \ref{ai7def3} and \ref{ai7def4}, which allow~$f^{ik}\ne f^{jk}\ci f^{ij}$.

\begin{dfn}
\label{ai7def8}
Let $R$ be a commutative ring and $K\in\{3,4,5,\ldots,\iy\}$. We will define a simplicial set $\dACat^{\un,\iy}_{R,K}$ of unital $A_\iy$-categories over~$R$.
\smallskip

\noindent{\bf(a)} For $N\ge 0$, define the $N$-simplices $(\dACat^{\un,\iy}_{R,K})_N$ in $\dACat^{\un,\iy}_{R,K}$ to be the set of quadruples $\bigl((\sO^i)_{i=0}^N,(f^{ij})_{0\le i\le j\le N},(\sH^i)_{i=0}^N,(\m^I)_{\es\ne I\subseteq\{0,\ldots,N\}}\bigr)$, where:
\begin{itemize}
\setlength{\itemsep}{0pt}
\setlength{\parsep}{0pt}
\item[(i)] $\sO^i$ is a set for $i=0,\ldots,N$, thought of as the set of objects $\Obj\sA^i$ of a (unital)  $A_\iy$-category $\sA^i$ at vertex $i$ of the $N$-simplex $\De^N$.
\item[(ii)] $f^{ij}:\sO^i\ra\sO^j$ is a map for all $0\le i\le j\le N$, thought of as the map $\fF_\Obj^{ij}:\Obj\sA^i\ra\Obj\sA^j$ in a (unital) $A_\iy$-functor $\fF^{ij}:\sA^i\ra\sA^j$ living on the edge $[i,j]$ of $\De^N$. We require these maps to satisfy $f^{ii}=\id_{\sO^i}$ for $i=0,\ldots,N$ and $f^{jk}\ci f^{ij}=f^{ik}$ for $0\le i\le j\le k\le N$.
\item[(iii)] For all $i=0,1,\ldots,N$ and $x,y\in\sO^i$, we are given a $\Z$-graded $R$-module $\sH^i(x,y)=\bigop_{d\in\Z}\sH^i(x,y)^d$, thought of as the Hom space $\Hom_{\sA^i}(x,y)$ in the $A_\iy$-category $\sA^i$. We write $\sH^i=(\sH^i(x,y))_{x,y\in\sO^i}$.
\item[(iv)] For each nonempty subset $I\subseteq\{0,\ldots,N\}$, which we write as $I=\{i_0,\ab i_1,\ab\ldots,\ab i_M\}$ for $i_0<i_1<\cdots<i_M$ with $M\ge 0$, and for all $k\ge 0$ and elements $x_0,\ldots,x_k\in \sO^{i_0}$, thought of as objects in $\sA^{i_0}$, we are given an $R$-linear map
\ea
\m_k^{I;x_0,\ldots,x_k}:\sH^{i_0}&(x_0,x_1)\ot_R\sH^{i_0}(x_1,x_2)\ot_R\cdots \ot_R\sH^{i_0}(x_{k-1},x_k)
\nonumber\\
&\longra\sH^{i_M}(f^{i_0i_M}(x_0),f^{i_0i_M}(x_k))
\label{ai7eq24}
\ea
graded of degree $2-M-k=3-\md{I}-k$. Here when $k=0$, by convention the empty tensor product on the top line of \eq{ai7eq24} is $R$, so $\m_0^{I;x_0}$ maps
\begin{equation*}
\m_0^{I;x_0}:R\longra\sH^{i_M}(f^{i_0i_M}(x_0),f^{i_0i_M}(x_0))^{2-M}.
\end{equation*}
When $k=0$ we require that $\m_0^{I;x_0}=0$ for all $\es\ne I\subseteq\{0,\ldots,N\}$ with $\md{I}\le K$ and $x_0\in\sO^{\min I}$.

We write $\m^I_k=(\m_k^{I;x_0,\ldots,x_k})_{x_0,\ldots,x_k\in\sO^{i_0}}$ and $\m^I=(\m^I_k)_{k\ge 0}$.

We require these maps $\m_k^{I;x_0,\ldots,x_k}$ to satisfy the following identity generalizing \eq{ai6eq28} for all $I=\{i_0<\cdots<i_M\}\subseteq\{0,\ldots,N\}$, $n\ge 0$ and objects~$x_0,\ldots,x_n\in\sO^{i_0}$
\ea
&\sum_{\begin{subarray}{l} 0\le j\le M, \; 0\le k\le\max(n,n+j-1), \\ y_0,\ldots,y_k\in\sO^{i_j}\end{subarray}\!\!\!\!\!\!\!\!\!\!\!\!\!\!\!\!\!\!\!\!\!\!\!\!}
(-1)^j\,\m_k^{\{i_j,\ldots,i_M\};y_0,\ldots,y_k}\ci\d_{n,k;y_0,\ldots,y_k}^{\{i_0,\ldots,i_j\};x_0,\ldots,x_n}
\nonumber\\
&=\sum_{j=1,\ldots,M-1}(-1)^j\m_n^{I\sm\{i_j\};x_0,\ldots,x_n},
\label{ai7eq25}
\ea
where $\d_{n,k;y_0,\ldots,y_k}^{\{i_0,\ldots,i_j\};x_0,\ldots,x_n}$ is defined from the $\m_k^{J;w_0,\ldots,w_k}$ as in \eq{ai7eq7}--\eq{ai7eq9}. In \eq{ai7eq25} there are only finitely many nonzero terms, which all map to $\sH^{i_M}(f^{i_0i_M}(x_0),f^{i_0i_M}(x_n))$, so \eq{ai7eq25} is an equation in this space.
\item[(v)] Definition \ref{ai7def3}(a)(v) holds.
\end{itemize}

\noindent{\bf(b)} If $\al:[\bar  N]\ra[N]$ is a morphism in $\De$, that is, $\al:\{0,\ldots,\bar  N\}\ra\{0,\ldots,N\}$ is nondecreasing, let $\dACat^{\un,\iy}_{R,K}(\al):(\dACat^{\un,\iy}_{R,K})_N\!\ra\!(\dACat^{\un,\iy}_{R,K})_{\bar N}$ map $\bigl((\sO^i),\ab(f^{ij}),(\sH^i),(\m^I)\bigr)$ to $\bigl((\bar \sO^i),(\bar  f^{ij}),(\bar\sH{}^i),(\bar\m{}^{\bar I})\bigr)$, where $\bar\sO^i=\sO^{\al(i)}$ for $i=0,\ldots,\bar N$, and $\bar f^{ij}=f^{\al(i)\al(j)}$ for $0\le i\le j\le\bar N$, and $\bar\sH^i=\sH^{\al(i)}$ for $i=0,\ldots,\bar N$. To define $\bar\m{}^{\bar I}$, if $\bar I=\{\bar\imath_0<\bar\imath_1<\cdots<\bar\imath_{\bar M}\}\subseteq\{0,1,\ldots,\bar N\}$, $k\ge 0$ and $x_0,\ldots,x_k\in\bar\sO^{\bar\imath_0}=\sO^{\al(\bar\imath_0)}$, generalizing \eq{ai6eq22} we define
\e
\bar\m_k^{\bar I;x_0,\ldots,x_k}=\begin{cases}
\m_k^{\al(\bar I);x_0,\ldots,x_k} & \text{if $\al\vert_{\bar I}$ is injective,} \\
\id_{\sH^i(x_0,x_1)} & {\begin{subarray}{l} \ts \text{if $k\!=\!1$ and $\bar I\!=\!\{\bar\imath_0,\bar\imath_1\}$ with} \\ \ts \text{$\bar\imath_0\!\ne\! \bar\imath_1$, $\al(\bar\imath_0)\!=\!\al(\bar\imath_1)\!=\!i$,} \end{subarray}} \\[8pt]
0 & \text{otherwise.}
\end{cases}
\label{ai7eq26}
\e
Note that in contrast to \eq{ai7eq20}, there is no special case when $k=0$, $\md{\bar I}=3$ and $\al\vert_{\bar I}$ is not injective, and $\bar\m_0^{\{\bar\imath_0,\bar\imath_1,\bar\imath_2\};x_0}=0$ in this case.

Proposition \ref{ai7prop3} shows $\dACat^{\un,\iy}_{R,K}(\al)$ does map to $(\dACat^{\un,\iy}_{R,K})_{\bar N}$, and thus $\dACat^{\un,\iy}_{R,K}$ is a simplicial set as in Definitions \ref{ai6def12} and \ref{ai7def3}. Theorem \ref{ai7thm7} will show that $\dACat^{\un,\iy}_{R,K}$ is an $\iy$-category. We say that $\dACat^{\un,\iy}_{R,K}$ has {\it type\/} $K$. Clearly $\dACat^{\un,\iy}_{R,K}\subset\dACat^{\un,\iy}_{R,L}$ if $L\le K$ in $\{3,4,5,\ldots,\iy\}$, giving the inclusion relations in the second column of~\eq{ai7eq1}.
\end{dfn}

Here is the analogue of Propositions \ref{ai6prop2} and \ref{ai7prop1}, proved in~\S\ref{ai106}.

\begin{prop}
\label{ai7prop3}
The data $\bigl((\bar\sO^i)_{i=0}^{\bar N},\!(\bar f^{ij})_{0\le i\le j\le\bar N},\!(\bar\sH{}^i)_{i=0}^{\bar N},\!(\n{}^{\bar I})_{\es\ne\bar I\subseteq\{0,\ldots,\bar N\}}\bigr)$ in Definition\/ {\rm\ref{ai7def8}(b)} lies in $(\dACat^{\un,\iy}_{R,K})_{\bar N}$.
\end{prop}

The next theorem, an analogue of Theorems \ref{ai6thm5} and \ref{ai7thm2}, will be proved in \S\ref{ai107} using equation~\eq{ai7eq25}.

\begin{thm}
\label{ai7thm5}
Let\/ $\bigl((\sO^i)_{i=0}^N,(f^{ij})_{0\le i\le j\le N},(\sH^i)_{i=0}^N,(\m^I)_{\es\ne I\subseteq\{0,\ldots,N\}}\bigr)$ be an $N$-simplex in $\dACat^{\un,\iy}_{R,K}$ as in Definition\/ {\rm\ref{ai7def8},} and define maps\/ $\d_{n,k;y_0,\ldots,y_k}^{I;x_0,\ldots,x_n}$ from the $\m_k^{I;x_0,\ldots,x_k}$ as in\/ {\rm\eq{ai7eq6}--\eq{ai7eq9}}. Then for each nonempty subset\/ $I\subseteq\{0,\ldots,N\},$ written $I=\{i_0,\ab i_1,\ab\ldots,\ab i_M\}$ for $i_0<i_1<\cdots<i_M,$ and all\/ $k,n\ge 0,$ $x_0,\ldots,x_n\in\sO^{i_0}$ and\/ $y_0,\ldots,y_k\in\sO^{i_M},$ in $R$-linear maps as in \eq{ai7eq6} we~have
\ea
&\sum_{j=0}^M\sum_{\begin{subarray}{l} 0\le l\le \max(n,n+j-1), \\ w_0,\ldots,w_l\in\sO^{i_j} \end{subarray}}(-1)^j\d_{l,k;y_0,\ldots,y_k}^{\{i_j,\ldots,i_M\};w_0,\ldots,w_l}\ci\d_{n,l;w_0,\ldots,w_l}^{\{i_0,\ldots,i_j\};x_0,\ldots,x_n}
\nonumber\\
&=\sum_{j=1}^{M-1}(-1)^j\,\d_{n,k;y_0,\ldots,y_k}^{I\sm\{i_j\};x_0,\ldots,x_n},
\label{ai7eq27}
\ea
where there are only finitely many nonzero terms in the sums.
\end{thm}

Here is the analogue of Examples \ref{ai6ex4} and~\ref{ai7ex2}.

\begin{ex}
\label{ai7ex3}
{\bf(a)} A 0-simplex $(\sO^0,f^{00}=\id_{\sO^0},\sH^0,\m^{\{0\}})$ of $\dACat^{\un,\iy}_{R,K}$ is equivalent to a unital $A_\iy$-category $\sA^0$ over $R$, with objects $\Obj\sA^0=\sO^0$, Hom spaces $\Hom_{\sA^0}(x,y)=\sH^0(x,y)$, morphisms $(\m_k^{\{0\};x_0,\ldots,x_k})_{k\ge 1}$, and units~$\boo^0_x\in \Hom_{\sA^0}^0(x,x)=\sH^0(x,x)^0$.
\smallskip

\noindent{\bf(b)} A 1-simplex of $\dACat^{\un,\iy}_{R,K}$ is equivalent to unital $A_\iy$-categories $\sA^0,\sA^1$ over $R$, and a unital $A_\iy$-functor $\fF^{01}:\sA^0\ra\sA^1$, with $\fF_\Obj^{01}=f^{01}$ and morphisms $\fF^{01;x_0,\ldots,x_k}_k=\m_k^{\{0,1\};x_0,\ldots,x_k}$. 
\smallskip

\noindent{\bf(c)} A 2-simplex of $\dACat^{\un,\iy}_{R,K}$ is equivalent to unital $A_\iy$-categories $\sA^0,\ab\sA^1,\ab\sA^2$ over $R$ as in {\bf(a)}, and unital $A_\iy$-functors $\fF^{ij}:\sA^i\ra\sA^j$ for $(i,j)=(0,1),\ab(0,2),\ab(1,2)$ as in {\bf(b)}, and $R$-linear maps for all $k\ge 1$ and $x_0,\ldots,x_k\in\sO^0$
\begin{align*}
\m_k^{\{0,1,2\};x_0,\ldots,x_k}:\sH^0&(x_0,x_1)\ot_R\sH^0(x_1,x_2)\ot_R\cdots \ot_R\sH^0(x_{k-1},x_k)
\nonumber\\
&\longra\sH^2(f^{02}(x_0),f^{02}(x_k))
\end{align*}
graded of degree $-k$, with $\m^{\{0,1,2\};x_0}_0=0$, satisfying \eq{ai7eq19} and
\ea
&\sum_{\begin{subarray}{l} i,k\ge 0,\; j\ge 1:\\ i+j+k=n \end{subarray}}\begin{aligned}[t]&
(-1)^{i+jk}\m^{\{0,1,2\};x_0,x_1,\ldots,x_i,x_{i+j},x_{i+j+1},\ldots,x_n}_{i+k+1}\ci\\
&\bigl(\id_{\sH^0(x_0,x_1)}\ot\cdots\ot \id_{\sH^0(x_{i-1},x_i)}\ot\m^{\{0\};x_i,x_{i+1},\ldots,x_{i+j}}_j\\
&\ot \id_{\sH^0(x_{i+j},x_{i+j+1})}\ot\cdots\ot \id_{\sH^0(x_{n-1},x_n)}\bigr)\end{aligned}
\label{ai7eq28}\\
&-(\m^{\{1,2\}}\!\ci\!\m^{\{0,1\}})^{x_0,\ldots,x_n}_n
\nonumber\\
&+\sum_{\!\!\!\!\!\!\begin{subarray}{l} 
1\le l\le k\le n+1, \\ 
n=n_1+\cdots+n_k, \\
n_i\ge 1,\; i=1,\ldots,k. \\
\text{Set $p_j=\sum_{i=1}^jn_i$,} \\ 
0\le j\le k \end{subarray}}
\begin{aligned}[t]&(-1)^{k+1}\cdot(-1)^{\sum_{i=1}^{l-1}(n_i+1)}\cdot(-1)^{\sum_{i=1}^k(n_i+1)(k-i)}\\
&\m^{\{2\};f^{\{0,1,2\}}(x_{p_0}),f^{\{0,1,2\}}(x_{p_1}),\ldots,f^{\{0,1,2\}}(x_{p_{l-1}}),f^{02}(x_{p_l}),\ldots,f^{02}(x_{p_k})}_k\\
&\ci\bigl[(\m^{\{1,2\}}\!\ci\!\m^{\{0,1\}})^{x_{p_0},\ldots,x_{p_1}}_{n_1}\!\ot\!\cdots\!\ot\!(\m^{\{1,2\}}\!\ci\!\m^{\{0,1\}})^{x_{p_{l-2}},\ldots,x_{p_{l-1}}}_{n_{l-1}}\\
&\ot \m^{\{0,1,2\};x_{p_{l-1}},\ldots,x_{p_l}}_{n_l}\ot\m^{\{0,2\};x_{p_l},\ldots,x_{p_{l+1}}}_{n_{l+1}}\ot\\
&\cdots\ot\m^{\{0,2\};x_{p_{k-1}},\ldots,x_{p_k}}_{n_k}\bigr]
\end{aligned}
\nonumber\\
&=-\m^{\{0,2\};x_0,\ldots,x_n}_n.
\nonumber
\ea

By equations \eq{ai4eq27}, \eq{ai4eq28}, \eq{ai7eq19}, and \eq{ai7eq28}, $\m^{\{0,1,2\}}=(\m^{\{0,1,2\}}_k)_{k\ge 1}$ is a {\it unital\/ $A_\iy$-homotopy\/} $\fF^{\{1,2\}}\ci \fF^{\{0,1\}}\Ra\fF^{\{0,2\}}$. Thus, a 2-simplex in $\dACat^{\un,\iy}_{R,K}$ is a 2-commutative diagram of unital $A_\iy$-categories, unital $A_\iy$-functors, and a unital $A_\iy$-homotopy
\begin{align*}
\xymatrix@C=70pt@R=13pt{
\drrtwocell_{}\omit^{}\omit{_{\,\,\,\,\,\,\,\,\,\,\,\,\,\,\,\,\,\,\m^{\{0,1,2\}}}} & \sA^1 \ar[dr]^{\fF^{\{1,2\}}}  \\
\sA^0 \ar[ur]^{\fF^{\{0,1\}}} \ar[rr]^(0.35){\fF^{\{0,2\}}} &&  \sA^2.}
\end{align*}

Note as in \S\ref{ai44} that $A_\iy$-homotopies are less general than $A_\iy$-natural transformations and $A_\iy$-natural isomorphisms, which we met in Example \ref{ai7ex2}(c) and Definition~\ref{ai7def4}(iii).
\end{ex}

The following analogue of Theorem \ref{ai6thm6} is proved in~\S\ref{ai108}.

\begin{thm}
\label{ai7thm6}
There is a strict isomorphism of\/ $\iy$-categories
\e
I^\un_{R,K}:\ACat^{\un,\iy}_{R,K}\longra\dACat^{\un,\iy}_{R,K},\quad K\in\{3,4,5,\ldots,\iy\},
\label{ai7eq29}
\e
from Definitions\/ {\rm\ref{ai7def6}} and\/ {\rm\ref{ai7def8},} which acts on $N$-simplices by
\begin{align*}
I^\un_{R,K}:&\bigl((\sO^i)_{i=0}^N,(f^{ij})_{0\le i\le j\le N},(\sH^i)_{i=0}^N,(\m^I)_{\es\ne I\subseteq\{0,\ldots,N\}}\bigr)\longmapsto\\
&\bigl((\sO^i)_{i=0}^N,(f^{ij})_{0\le i\le j\le N},(\sH^i)_{i=0}^N,(\dot\m^I)_{\es\ne I\subseteq\{0,\ldots,N\}}\bigr),
\end{align*}
where $\dot\m^{I;x_0,\ldots,x_k}_k$ is defined from $\m^{I;x_0,\ldots,x_k}_k$ by
\begin{equation*}
\dot\m_k^{\{i_0,\ldots,i_M\};x_0,\ldots,x_k}=\begin{cases} 0, & \text{$k=0$ and\/ $M=2$,} \\
\m_k^{\{i_0,\ldots,i_M\};x_0,\ldots,x_k}, & \text{otherwise.} \end{cases}
\end{equation*}
\end{thm}

Here are the analogues of Definitions \ref{ai6def13} and~\ref{ai6def14}.

\begin{dfn}
\label{ai7def9}
Let $R$ be a commutative ring and $K\in\{3,4,5,\ldots,\iy\}$. Generalizing Definition \ref{ai7def8}, define a simplicial set $\ACat^{\cu,\iy}_{R,K}$ of cohomologically unital $A_\iy$-categories over $R$ such that the $N$-simplices $\bigl((\sO^i)_{i=0}^N,(f^{ij})_{0\le i\le j\le N},\ab(\sH^i)_{i=0}^N,(\m^I)_{\es\ne I\subseteq\{0,\ldots,N\}}\bigr)$ for $N\ge 0$ satisfy Definition \ref{ai7def8}(a)(i)--(iv), but instead of Definition \ref{ai7def8}(a)(v), we require that
\begin{itemize}
\setlength{\itemsep}{0pt}
\setlength{\parsep}{0pt}
\item[(v$)'$] As in Example \ref{ai7ex3}(a), Definition \ref{ai7def8}(a)(i)--(iv) imply that $\sA^i\!=\!(\sO^i,\ab\sH^i,\ab\m^{\{i\}})$ is an $A_\iy$-category for each $i=0,\ldots,N$. We require this $A_\iy$-category to be {\it cohomologically unital}, as in Definition \ref{ai4def4}.
\item[(vi$)'$] As in Example \ref{ai7ex3}(b), Definition \ref{ai7def8}(a)(i)--(iv) mean $\fF^{ij}\!=\!(f^{ij},\m^{\{i,j\}}):\sA^i\ra\sA^j$ is an $A_\iy$-functor for all $0\le i<j\le N$. We require this to be a {\it cohomologically unital\/} $A_\iy$-functor, as in Definition~\ref{ai4def6}.
\end{itemize}
If $\al:[\bar N]\ra[N]$ is a morphism in $\De$, define $\ACat^{\cu,\iy}_{R,K}(\al):(\ACat^{\cu,\iy}_{R,K})_N\ab\ra(\ACat^{\cu,\iy}_{R,K})_{\bar N}$ as in Definition \ref{ai7def8}(b). The proof of Proposition \ref{ai7prop3} also works in the cohomologically unital case, and shows $\ACat^{\cu,\iy}_{R,K}$ is a simplicial set. It contains $\ACat^{\un,\iy}_{R,K}$ in Definition \ref{ai7def8} as a simplicial subset.

Theorem \ref{ai7thm5} also holds in $\ACat^{\cu,\iy}_{R,K}$, as the proof does not use unitality.

Theorem \ref{ai7thm7} will show that $\ACat^{\cu,\iy}_{R,K}$ is an $\iy$-category. We say that $\ACat^{\cu,\iy}_{R,K}$ has {\it type\/} $K$. We also write~$\ACat^{\cu,\iy}_R=\ACat^{\cu,\iy}_{R,3}$.

Following Definition \ref{ai7def7} we can define $\iy$-functors $I_{\AAlg^\cu}^{\ACat^\cu}:\AAlg^{\cu,\iy}_{R,K}\ab\ra\ACat^{\cu,\iy}_{R,K}$ for all $K\in\{3,4,5,\ldots,\iy\}$, for $\AAlg^{\cu,\iy}_{R,K}$ as in Definition \ref{ai6def13}, which are strict isomorphisms with the full $\iy$-subcategory of $\ACat^{\cu,\iy}_{R,K}$ with objects cohomologically unital $A_\iy$-categories $\sA$ with one object~$*$.
\end{dfn}

\begin{dfn}
\label{ai7def10}
Let $R$ be a commutative ring and $K\in\{3,4,5,\ldots,\iy\}$. Generalizing Definition \ref{ai7def8}, define a simplicial set $\ACat^\iy_{R,K}$ of (possibly non-unital) $A_\iy$-categories over $R$ whose $N$-simplices $\bigl((\sO^i)_{i=0}^N,(f^{ij})_{0\le i\le j\le N},\ab(\sH^i)_{i=0}^N,\ab(\m^I)_{\es\ne I\subseteq\{0,\ldots,N\}}\bigr)$ for $N\ge 0$ satisfy Definition \ref{ai7def8}(a)(i)--(iv), but need not satisfy (a)(v). If $\al:[\bar N]\ra[N]$ is a morphism in $\De$, define $\ACat^\iy_{R,K}(\al):(\ACat^\iy_{R,K})_N\ab\ra(\ACat^\iy_{R,K})_{\bar N}$ as in Definition \ref{ai7def8}(b). The proof of Proposition \ref{ai7prop3} also works in the non-unital case, and shows $\ACat^\iy_{R,K}$ is a simplicial set. It contains $\ACat^{\cu,\iy}_{R,K}$ in Definition \ref{ai7def9} as a simplicial subset.

Theorem \ref{ai7thm5} also holds in $\ACat^\iy_{R,K}$, as the proof does not use unitality.

Theorem \ref{ai7thm7} will show that $\ACat^\iy_{R,K}$ is an $\iy$-category. We say that $\ACat^\iy_{R,K}$ has {\it type\/} $K$. We also write~$\ACat^\iy_R=\ACat^\iy_{R,3}$.

Following Definition \ref{ai7def7} we can define $\iy$-functors $I_{\AAlg}^{\ACat}:\AAlg^\iy_{R,K}\ab\ra\ACat^\iy_{R,K}$ for all $K\in\{3,4,5,\ldots,\iy\}$, for $\AAlg^\iy_{R,K}$ as in Definition \ref{ai6def14}, which are strict isomorphisms with the full $\iy$-subcategory of $\ACat^\iy_{R,K}$ with objects $A_\iy$-categories $\sA$ with one object~$*$.
\end{dfn}

The following generalization of Theorem \ref{ai6thm7}, one of our main results, is proved in \S\ref{ai109}. The proof is easier than that of Theorem~\ref{ai7thm3}.

\begin{thm}
\label{ai7thm7}
The following simplicial sets defined in Definitions\/ {\rm \ref{ai7def6}, \ref{ai7def8}, \ref{ai7def9},} and\/ {\rm\ref{ai7def10},} for all\/ $K\in\{3,4,5,\ldots,\iy\},$ are all\/ $\iy$-categories:
\e
\begin{gathered}
\ACat^{\un,\iy}_{R,K},\dACat^{\un,\iy}_{R,K},\ACat^{\cu,\iy}_{R,K},\ACat^\iy_{R,K}.
\end{gathered}
\label{ai7eq30}
\e
\end{thm}

Here is the analogue of Definition \ref{ai6def15}. It completes our diagram \eq{ai7eq1} of inclusions of $\iy$-categories of $A_\iy$-categories.

\begin{dfn}
\label{ai7def11}
Let $\ACat^\un_R$ be the category from Definition \ref{ai4def6}, and $\sN(\ACat^\un_R)$ be its nerve from Example \ref{ai2ex5}, which is an $\iy$-category. Then $N$-simplices in $\sN(\ACat^\un_R)$ are of the form $\bigl((\sA^i)_{i=0}^N,(\fF^i)_{i=1}^N\bigr)$, where $\sA^i$ is a unital $A_\iy$-category over $R$ for $i=0,\ldots,N$, and $\fF^i:\sA^{i-1}\ra\sA^i$ is a unital $A_\iy$-functor for~$i=1,\ldots,N$.

Define an $\iy$-functor $\inc^\un_R:\sN(\ACat^\un_R)\ra\dACat^{\un,\iy}_{R,\iy}$ to act on $N$-simplices for $N\ge 0$ mapping
\begin{equation*}
\bigl((\sA^i)_{i=0}^N,(\fF^i)_{i=1}^N\bigr)\longmapsto\bigl((\sO^i)_{i=0}^N,(f^{ij})_{0\le i\le j\le N},(\sH^i)_{i=0}^N,(\m^I)_{\es\ne I\subseteq\{0,\ldots,N\}}\bigr),
\end{equation*}
where we set $\sO^i=\Obj\sA^i$, and writing $\fF^{ii}=\id_{\sA^i}:\sA^i\ra\sA^i$ and $\fF^{ij}=\fF^j\ci\fF^{j-1}\ci\cdots\ci\fF^{i+1}:\sA^i\ra\sA^j$ for $0\le i<j\le N$, we define $f^{ij}=\fF_\Obj^{ij}:\sO^i\ra\sO^j$ for $0\le i\le j\le N$. We set $\sH^i(x,y)=\Hom_{\sA^i}(x,y)$ for $x,y\in\sO^i=\Obj\sA^i$. We take $\m^{\{i\};x_0,\ldots,x_k}_k$ to be the map $\m^{x_0,\ldots,x_k}_k$ in $\sA^i$ for $i=0,\ldots,N$, and $\m^{\{i,j\};x_0,\ldots,x_k}_k$ to be the map $\fF^{ij;x_0,\ldots,x_k}_k$ in $\fF^{ij}$ for $0\le i<j\le N$, and $\m^{I;x_0,\ldots,x_k}_k=0$ if $\md{I}\ge 3$ and~$k\ge 0$.

We claim $\bigl((\sO^i)_{i=0}^N,(f^{ij})_{0\le i\le j\le N},(\sH^i)_{i=0}^N,(\m^I)_{\es\ne I\subseteq\{0,\ldots,N\}}\bigr)$ is an $N$-sim\-pl\-ex of $\dACat^{\un,\iy}_{R,\iy}$. We can verify \eq{ai7eq25} for all $\es\ne I\subseteq\{0,\ldots,N\}$, $n\ge 0$ and $x_0,\ldots,x_n\in\sO^{\min I}$ in a very similar way to Definition \ref{ai6def15}(i)--(iv). Definition \ref{ai7def8}(a)(v) follows as the $\sA^i,\fF^{ij}$ are unital.

It is now easy to see that $\inc^\un_R$ is a morphism of simplicial sets, that is, an $\iy$-functor. It is a strict isomorphism with the $\iy$-subcategory of $\dACat^{\un,\iy}_{R,\iy}$ consisting of $N$-simplices $\bigl((\sO^i)_{i=0}^N$, $(f^{ij})_{0\le i\le j\le N}$, $(\sH^i)_{i=0}^N$, $(\m^I)_{\es\ne I\subseteq\{0,\ldots,N\}}\bigr)$ with $\m_k^{I;x_0,\ldots,x_k}=0$ if $\md{I}\ge 3$ and $k\ge 0$. Similarly we define $\iy$-functors $\inc^\cu_R:\sN(\ACat^\cu_R)\hookra\ACat^{\cu,\iy}_{R,\iy}$ for cohomologically unital $A_\iy$-categories, and $\inc_R:\sN(\ACat_R)\hookra\ACat^\iy_{R,\iy}$ for general $A_\iy$-categories.

We can also define an $\iy$-functor $\iinc^\un_R:\sN(\ACat^\un_R)\ra\ACat^{\un,\iy}_{R,\iy}$ which is identified with $\inc^\un_R$ under the isomorphism $I^\un_{R,\iy}$ in \eq{ai7eq29}. The difference is that as in \eq{ai7eq23}, when we define $\bigl((\sO^i)_{i=0}^N,(f^{ij})_{0\le i\le j\le N},(\sH^i)_{i=0}^N,\ab(\m^I)_{\es\ne I\subseteq\{0,\ldots,N\}}\bigr)$ above we set $\m_0^{\{i_0,i_1,i_2\};x_0}=\boo^{i_2}_{f^{i_0i_2}(x_0)}$, not~$\m_0^{\{i_0,i_1,i_2\};x_0}=0$.

Similarly we define an $\iy$-functor~$\iinc^\un_R:\sN(\ACat^\un_R)\ra\ACat^{\un,\iy}_{R,\iy}$.
\end{dfn}

\begin{rem}
\label{ai7rem6}
Composing $\iinc^\un_R$ with $\ACat^{\un,\iy}_{R,\iy}\hookra\ACat^{\un,\iy}_R$ gives an $\iy$-functor $\iinc^\un_R:\sN(\ACat^\un_R)\ra\ACat^{\un,\iy}_R$, which is a strict isomorphism with an $\iy$-subcategory of $\ACat^{\un,\iy}_R$.

It follows from Example \ref{ai7ex2}(c) and Definition \ref{ai7def4} that if morphisms $\fF,\fG:\sA\ra\sB$ in $\ACat^\un_R$ are unitally $A_\iy$-naturally isomorphic, then their images under $\iinc^\un_R$ are homotopic in $\ACat^{\un,\iy}_R$, that is $[\iinc^\un_R(\fF)]=[\iinc^\un_R(\fG)]$ in $\Ho(\ACat^{\un,\iy}_R)$. Therefore $\iinc^\un_R$ maps unital homotopy equivalences in $\ACat^\un_R$, as in Definition \ref{ai4def14}, to equivalences in~$\ACat^{\un,\iy}_R$.

Thus by the universal property of $\iy$-category localizations in \S\ref{ai26}, there exists an $\iy$-functor $\Phi$ in the 2-commutative diagram
\begin{equation*}
\xymatrix@C=80pt@R=15pt{
*+[r]{\sN(\ACat^\un_R)} \ar[rr]_(0.55){\Pi_\loc^\iy} \ar@/_.5pc/[drr]_{\iinc^\un_R} &  \drtwocell_{}\omit^{}\omit{^{}} & *+[l]{\ACat^\un_R[\HE_\un^{-1}]^\iy} \ar@{..>}[d]^{\Phi} \\
&& *+[l]{\ACat^{\un,\iy}_R,} }
\end{equation*}
where $\HE_\un$ is the class of unital homotopy equivalences in $\ACat^\un_R$. Theorem \ref{ai9thm4} will show that $\Phi$ is an equivalence of $\iy$-categories.
\end{rem}

\begin{dfn}
\label{ai7def12}
Define $\iy$-subcategories 
\begin{equation*}
\ACat^{\un,\rsu,\iy}_R\subset\ACat^{\un,\iy}_R, \quad 
\ACat^{\un,\hs,\iy}_R\subset\ACat^{\un,\iy}_R, 
\end{equation*}
to be the full $\iy$-subcategories with objects that are unital $A_\iy$-categories with {\it split units\/} for superscript `$\rsu$', as in \S\ref{ai46}, and objects that are {\it h-projective\/} unital $A_\iy$-categories with {\it split units\/} for superscript `$\hs$', as in~\S\ref{ai46}--\S\ref{ai47}.
\end{dfn}

\begin{dfn}
\label{ai7def13}
Define $\iy$-subcategories 
\begin{align*}
\ADCat^{\un,\iy}_R&\subset\ACat^{\un,\iy}_R,& \ADCat^{\un,\iy}_{R,K}&\subset\ACat^{\un,\iy}_{R,K}, \\
\ADCat^{\cu,\iy}_R&\subset\ACat^{\cu,\iy}_R,& \ADCat^{\cu,\iy}_{R,K}&\subset\ACat^{\cu,\iy}_{R,K}, \\
\ADCat^\iy_R&\subset\ACat^\iy_R,& \ADCat^\iy_{R,K}&\subset\ACat^\iy_{R,K}, 
\end{align*}
for $K\in\{3,4,\ldots,\iy\}$ to be the full $\iy$-subcategories whose objects are {\it pretriangulated\/} $A_\iy$-categories over $R$, as in Definition~\ref{ai4def21}.
\end{dfn}

\begin{dfn}
\label{ai7def14}
Define $\iy$-categories $\DGCat_R^{*,\iy}$ in a commutative diagram of inclusions of $\iy$-subcategories:
\e
\begin{gathered}
\xymatrix@C=40pt@R=15pt{ 
\DGCat_{R_{\vphantom{(}}}^{\hs,\iy}\,\, \ar@{^{(}->}[r]^(0.46){\text{$\not\simeq$ all $R$}} \ar@{^{(}->}[d]^{\not\simeq} & \DGCat_{R_{\vphantom{(}}}^{A_\iy,\hs,\iy}\,\, \ar@{^{(}->}[d]^{\not\simeq} \ar@{^{(}->}[r]^\simeq &  \ACat_{R_{\vphantom{(}}}^{\un,\hs,\iy} \ar@{^{(}->}[d]^{\not\simeq} \\
\DGCat_{R_{\vphantom{(}}}^{\rsu,\iy}\,\, \ar@{^{(}->}[r]^(0.46){\text{$\not\simeq$ all $R$}} \ar@{^{(}->}[d]^{\not\simeq} & \DGCat_{R_{\vphantom{(}}}^{A_\iy,\rsu,\iy}\,\, \ar@{^{(}->}[d]^{\not\simeq} \ar@{^{(}->}[r]^\simeq &  \ACat_{R_{\vphantom{(}}}^{\un,\rsu,\iy} \ar@{^{(}->}[d]^{\not\simeq} \\
\DGCat_R^\iy\,\, \ar@{^{(}->}[r]^(0.46){\text{$\not\simeq$ all $R$}} & \DGCat_R^{A_\iy,\iy}\,\, \ar@{^{(}->}[r] &  \ACat_R^{\un,\iy}, }
\end{gathered}
\label{ai7eq31}
\e
such that:
\begin{itemize}
\setlength{\itemsep}{0pt}
\setlength{\parsep}{0pt}
\item[(i)] Objects of $\DGCat_R^{*,\iy}$ are dg-categories over $R$, which are h-projective with split units in the top row of \eq{ai7eq31} with superscript `$\hs$', and have split units in the middle row of \eq{ai7eq31} with superscript `$\rsu$'.
\item[(ii)] In the left hand column of \eq{ai7eq31}, 1-morphisms in $\DGCat_R^{*,\iy}$ are dg-functors $\fF:\sA\ra\sB$ between dg-categories. In the middle column of \eq{ai7eq31}, with superscript $A_\iy$, 1-morphisms in $\DGCat_R^{A_\iy,*,\iy}$ are unital $A_\iy$-functors $\fF:\sA\ra\sB$ between dg-categories.
\item[(iii)] For $N\ge 2$, the $N$-simplices in each $\DGCat_R^{*,\iy}$ are all $N$-simplices in $\ACat_R^{\un,\iy}$ which have 0- and 1-simplices as prescribed in (i),(ii). 
\end{itemize}
Since the classes of 1-morphisms in (ii) are closed under composition, these do define $\iy$-subcategories. See Theorem \ref{ai7thm8} for the meaning of the labels `$\simeq$', `$\not\simeq$' on the inclusions in \eq{ai7eq31}.
\end{dfn}

% todo work out if final morphism in \eq{ai7eq31} is equivalence or not. Counterexample?

\begin{thm}
\label{ai7thm8}
In equation {\rm\eq{ai7eq31},} the inclusions marked `$\,\simeq$' are Joyal equivalences for all commutative rings $R$. For the inclusions marked `$\,\not\simeq$', there exist commutative rings $R$ (e.g.\ $R=\Z$) for which they are not Joyal equivalences. For those marked `$\,\not\simeq$ all\/ $R$', they are not Joyal equivalences for all\/ $R$.

If\/ $R$ is a field, the upper vertical morphisms in \eq{ai7eq31} are identity morphisms, so they are Joyal equivalences in this case, even though marked `$\,\not\simeq$'.
\end{thm}

\begin{proof}
By Theorem \ref{ai2thm1}, an inclusion in \eq{ai7eq31} is a Joyal equivalence if it is {\it essentially surjective}, as in Definition \ref{ai2def5}, and {\it fully faithful}, as in Definition \ref{ai2def9}. The two inclusions $\DGCat_R^{A_\iy,\hs,\iy}\hookra\ACat_R^{\un,\hs,\iy}$ and $\DGCat_R^{A_\iy,\rsu,\iy}\hookra\ab\ACat_R^{\un,\rsu,\iy}$ are essentially surjective by Theorem \ref{ai4thm4} and Proposition \ref{ai4prop8}. They are trivially fully faithful as they are inclusions of full $\iy$-subcategories. 

For the inclusions marked `$\not\simeq$' in \eq{ai7eq31}, Example \ref{ai7ex4} shows they are not Joyal equivalences when $R=\Z$. The last part holds as in~\S\ref{ai46}--\S\ref{ai47}.
\end{proof}

\begin{ex}
\label{ai7ex4}
{\bf(a)} Take $R=\Z$ and let $\sA$ be the dg-category over $R$ with one object $*$ and $\Hom_\sA(*,*)=\Z_2\cdot\id_*$. Then $\sA$ does not have split units, so does not lie in the essential image of the vertical inclusions  between the second and third rows of \eq{ai7eq31} by Proposition \ref{ai4prop5}. Thus these are not Joyal equivalences.

\smallskip

\noindent{\bf(b)} Take $R=\Z$ and let $\sA$ be the dg-category over $R$ with one object $*$ and $\Hom_\sA(*,*)=\Z\cdot\id_*\op\Z_2\cdot\ep$ with $\ep^2=0$. Then $\sA$ is not h-projective, so does not lie in the essential image of the vertical inclusions  between the first and second rows of \eq{ai7eq31} by Proposition \ref{ai4prop8}. Hence these are not Joyal equivalences.
\smallskip

\noindent{\bf(c)} Let $R$ be any ring. Take $\sA$ to be the dg-category over $R$ with one object $*$ and $\Hom_\sA(*,*)=R\cdot\id_*\op R\cdot x$ with $\deg x=0$ and $\d x=0$. Take $\sB$ to be the dg-category over $R$ with one object $*$ and $\Hom_\sB(*,*)=R\cdot\id_*\op R\cdot y$ with $\deg y=-1$ and $\d y=0$. Define $\fF:\sA\ra\sB$ to be the unital $A_\iy$-functor with $\fF_\Obj(*)=*$, $\fF_1(\id_*)=\id_*$, $\fF_1(x)=0$ and $\fF_2(x,x)=y$. Then $\fF$ is not $A_\iy$-natural isomorphic to any dg-functor $\fG:\sA\ra\sB$. This means that all three horizontal inclusions between the first and second columns of \eq{ai7eq31} are not fully faithful, and so are not Joyal equivalences.
\end{ex}

\section{\texorpdfstring{$\iy$-categories of filtered $A_\iy$-algebras and \\ filtered $A_\iy$-categories}{∞-categories of filtered A∞-algebras and filtered A∞-categories}}
\label{ai8}

Filtered $A_\iy$-categories were discussed in Chapter \ref{ai5}. Let $(\Ga,\pre)$ be a po-monoid, and $R$ a complete $\Ga$-filtered ring. We will define many $\iy$-categories of filtered $A_\iy$-categories over $R$, which fit into a diagram of inclusions of $\iy$-categories:
\ea
\begin{gathered}
\!\!\!\!\xymatrix@C=8pt@R=12pt{
{\begin{subarray}{l}\ts \AfCat^{\un,\iy}_R=  \\[2pt] \ts \AfCat^{\un,\iy}_{R,0,2}\end{subarray}} \\
\AfCat^{\un,\iy^{\vphantom{(}}}_{R,0,3} \ar@{=}[rr] \ar@{^{(}->}[u] && \dAfCat^{\un,\iy^{\vphantom{(}}}_{R,0,3}\,\, \ar@{^{(}->}[rr]  && {\begin{subarray}{l}\ts \AfCat^\iy_R=  \\[2pt] \ts \AfCat^\iy_{R,0,3}\end{subarray}} \\
\AfCat^{\un,\iy^{\vphantom{(}}}_{R,K',L'} \ar@{=}[rr] \ar@{^{(}->}[u] && \dAfCat^{\un,\iy^{\vphantom{(}}}_{R,K',L'}\,\, \ar@{^{(}->}[rr] \ar@{^{(}->}[u] &&  \AfCat^{\iy^{\vphantom{(}}}_{R,K',L'} \ar@{^{(}->}[u] & *+[l]{\;\begin{subarray}{l}(0\le K'\le L', \\ L'\ge 3)\end{subarray}} \\
\AfCat^{\un,\iy^{\vphantom{(}}}_{R,K,L} \ar@{=}[rr] \ar@<-4ex>@{^{(}->}[u]^{K\ge K',\; L\ge L'} && \dAfCat^{\un,\iy^{\vphantom{(}}}_{R,K,L}\,\, \ar@{^{(}->}[rr] \ar@<-4ex>@{^{(}->}[u]^{K\ge K',\; L\ge L'}  && \AfCat^{\iy^{\vphantom{(}}}_{R,K,L} \ar@<-4ex>@{^{(}->}[u]^{K\ge K',\; L\ge L'} & *+[l]{\;\begin{subarray}{l}(0\le K\le L, \\ L\ge 3)\end{subarray}} \\
\AfCat^{\un,\iy^{\vphantom{(}}}_{R,K,\iy} \ar@{=}[rr]  && \dAfCat^{\un,\iy^{\vphantom{(}}}_{R,K,\iy}\,\, \ar@{^{(}->}[rr]  &&  \AfCat^{\iy^{\vphantom{(}}}_{R,K,\iy} & \ar@<3ex>@{}[d]_{(K=0,1,\iy).} \\
\sN(\AfCat^{\un^{\vphantom{(}}}_{R,K}) \ar@{^{(}->}[u] \ar@{=}[rr] && \sN(\AfCat^{\un^{\vphantom{(}}}_{R,K})\,\, \ar@{^{(}->}[rr] \ar@{^{(}->}[u] &&  \sN(\AfCat^{\phantom\un^{\vphantom{(}}}_{R,K}) \ar@{^{(}->}[u] & }\!\!\!\!\!\!\!\end{gathered}
\label{ai8eq1}
\ea
We then define $\iy$-categories of filtered $A_\iy$-algebras by restricting to the full $\iy$-subcategories of filtered $A_\iy$-categories with one object $*$.

We comment on these:
\begin{itemize}
\setlength{\itemsep}{0pt}
\setlength{\parsep}{0pt}
\item[(i)] The objects of the $\iy$-categories in the first two columns, with superscript `$\un$', are {\it unital\/} filtered $A_\iy$-categories.
\item[(ii)] The objects in the third column are all filtered $A_\iy$-categories. 

(In contrast to Chapter \ref{ai7} we do not define $\iy$-categories of {\it cohomologically unital\/} filtered $A_\iy$-categories, as we do not have a good definition of these.)
\item[(iii)] The $\iy$-categories in the first two columns, linked by `$=$', are strictly isomorphic, but with definitions written in a different way.
\item[(iv)] We have subscripts `$R,K,L$' for $K,L\in\{0,1,2,\ldots,\iy\}$ with $K\le L$, and $L\ge 2$ for the first column and $L\ge 3$ for the second and third. We call $(K,L)$ the {\it type\/} of the $\iy$-category, generalizing type $K$ in Chapters \ref{ai6}--\ref{ai7}. We have inclusions $\AfCat^{*,\iy}_{R,K,L}\hookra\AfCat^{*,\iy}_{R,K',L'}$ if $K\ge K'$ and~$L\ge L'$.

As in Chapters \ref{ai6}--\ref{ai7}, the $N$-simplices in the $\iy$-categories in \eq{ai8eq1} all involve families of maps $\m_k^{I;x_0,\ldots,x_k}$ for $\es\ne I\subseteq\{0,\ldots,N\}$, $k\ge 0$ and $x_0,\ldots,x_k$ in $\Obj\sA^{\min I}$, where $\m_k^{I;x_0,\ldots,x_k}$ is filtration-preserving. 

For type $(K,L)$ and $\md{I}\ne 3$, the rule is that $\m_0^{I;x_0}=0$ if $\md{I}\le K$, and $\m_0^{I;x_0}(1)\in F^{\ssu 0}\Hom_{\sA^{i_M}}(f^I(x_0),f^{i_0i_M}(x_0))$ if $K<\md{I}\le L$, and $\m_0^{I;x_0}(1)\in F^0\Hom_{\sA^{i_M}}(f^I(x_0),f^{i_0i_M}(x_0))$ if~$L<\md{I}$.
\item[(v)] $\sN(\AfCat^\un_{R,K}),\sN(\AfCat_{R,K})$ are the nerves of the ordinary categories $\AfCat^\un_{R,K},\AfCat_{R,K}$ from Definition \ref{ai5def9}, for~$K=0,1,\iy$.
\end{itemize}

Now 2-morphisms in $\AfCat^{\un,\iy}_R$ are filtered $A_\iy$-natural isomorphisms in the sense of Definition \ref{ai5def12}(c). We define a second diagram of $\iy$-categories:
\ea
\begin{gathered}
\!\!\!\!\xymatrix@C=8pt@R=12pt{
{\AwCat^{\un,\iy}_R=\AwCat^{\un,\iy}_{R,0,2}\hskip -60pt} \\
\AwCat^{\un,\iy^{\vphantom{(}}}_{R,0,3} \ar@{=}[rr] \ar@{^{(}->}[u] && \dAwCat^{\un,\iy^{\vphantom{(}}}_{R,0,3}\,\, \ar@{^{(}->}[rr]  && {\begin{subarray}{l}\ts \AwCat^\iy_R=  \\[2pt] \ts \AwCat^\iy_{R,0,3}\end{subarray}} \\
\AwCat^{\un,\iy^{\vphantom{(}}}_{R,K',L'} \ar@{=}[rr] \ar@{^{(}->}[u] && \dAwCat^{\un,\iy^{\vphantom{(}}}_{R,K',L'}\,\, \ar@{^{(}->}[rr] \ar@{^{(}->}[u] &&  \AwCat^{\iy^{\vphantom{(}}}_{R,K',L'} \ar@{^{(}->}[u] & *+[l]{\;\begin{subarray}{l}(0\le K'\le L', \\ L'\ge 3)\end{subarray}} \\
\AwCat^{\un,\iy^{\vphantom{(}}}_{R,K,L} \ar@{=}[rr] \ar@<-4ex>@{^{(}->}[u]^{K\ge K',\; L\ge L'} && \dAwCat^{\un,\iy^{\vphantom{(}}}_{R,K,L}\,\, \ar@{^{(}->}[rr] \ar@<-4ex>@{^{(}->}[u]^{K\ge K',\; L\ge L'}  && \AwCat^{\iy^{\vphantom{(}}}_{R,K,L} \ar@<-4ex>@{^{(}->}[u]^{K\ge K',\; L\ge L'} & *+[l]{\;\begin{subarray}{l}(0\le K\le L, \\ L\ge 3)\end{subarray}} \\
\AwCat^{\un,\iy^{\vphantom{(}}}_{R,K,\iy} \ar@{=}[rr]  && \dAwCat^{\un,\iy^{\vphantom{(}}}_{R,K,\iy}\,\, \ar@{^{(}->}[rr]  &&  \AwCat^{\iy^{\vphantom{(}}}_{R,K,\iy} & \ar@<3ex>@{}[d]_{(K=0,1,\iy.)} \\
\sN(\AfCat^{\un^{\vphantom{(}}}_{R,K}) \ar@{^{(}->}[u] \ar@{=}[rr] && \sN(\AfCat^{\un^{\vphantom{(}}}_{R,K})\,\, \ar@{^{(}->}[rr] \ar@{^{(}->}[u] &&  \sN(\AfCat^{\phantom\un^{\vphantom{(}}}_{R,K}) \ar@{^{(}->}[u] &}\!\!\!\!\!\!\!\end{gathered}
\label{ai8eq2}
\ea
such that 2-morphisms in $\AwCat^{\un,\iy}_R$ are {\it weakly\/} filtered $A_\iy$-natural isomorphisms in the sense of Definition \ref{ai5def12}(d). There are inclusions of $\iy$-categories $\AfCat^{*,\iy}_{R,K,L}\hookra\AwCat^{*,\iy}_{R,K,L}$. The difference with the $\iy$-categories in \eq{ai8eq1} is that (iv) above should be replaced by:
\begin{itemize}
\setlength{\itemsep}{0pt}
\setlength{\parsep}{0pt}
\item[(iv$)'$]  We have subscripts `$R,K,L$' for $K,L\in\{0,1,2,\ldots,\iy\}$ with $K\le L$, and $L\ge 2$ for the first column and $L\ge 3$ for the second and third. We call $(K,L)$ the {\it type\/} of the $\iy$-category, generalizing type $K$ in Chapters \ref{ai6}--\ref{ai7}. We have inclusions $\AwCat^{*,\iy}_{R,K,L}\hookra\AwCat^{*,\iy}_{R,K',L'}$ if $K\ge K'$ and~$L\ge L'$.

The $N$-simplices in the $\iy$-categories in \eq{ai8eq2} all involve families of maps $\m_k^{I;x_0,\ldots,x_k}$ for $\es\ne I\subseteq\{0,\ldots,N\}$, $k\ge 0$ and $x_0,\ldots,x_k$ in $\Obj\sA^{\min I}$, where $\m_k^{I;x_0,\ldots,x_k}$ is filtration-preserving if $\md{I}=1,2$, and $\m_k^{I;x_0,\ldots,x_k}$ decreases filtrations by a bounded amount $\ep\in\Ga$ if $\md{I}\ge 3$. 

For type $(K,L)$ and $\md{I}\ne 3$, the rule is that $\m_0^{I;x_0}=0$ if $\md{I}\le K$, and $\m_0^{I;x_0}(1)\in F^{\ssu 0}\Hom_{\sA^{i_M}}(f^I(x_0),f^{i_0i_M}(x_0))$ if $K<\md{I}\le L$, and $\m_0^{I;x_0}(1)\in \Hom_{\sA^{i_M}}(f^I(x_0),f^{i_0i_M}(x_0))$ if~$L<\md{I}$.
\end{itemize}
The objects (0-simplices) and morphisms (1-simplices) in \eq{ai8eq1} and \eq{ai8eq2} are the same, the two differ in $N$-simplices for~$N\ge 2$.

In brief, in \eq{ai8eq1} all maps $\m_k^{I;x_0,\ldots,x_k}$ must be filtration-preserving, but in \eq{ai8eq2} the $\m_k^{I;x_0,\ldots,x_k}$ can decrease filtrations for $\md{I}>L\ge 2$. See Remark \ref{ai5rem6}(a) for an application in Symplectic Geometry which needs weakly filtered $A_\iy$-natural isomorphisms, and so should be carried out in the $\iy$-category~$\AwCat^{\un,\iy}_R$.

\subsection[\texorpdfstring{The $\iy$-categories $\AfCat^{\un,\iy}_R,\AfCat^{\un,\iy}_{R,K,L},\AwCat^{\un,\iy}_R$, \\ $\AwCat^{\un,\iy}_{R,K,L}$}{The ∞-categories Aᶠⁱ∞-Catᵘⁿ∞ᵣ,Aᶠⁱ∞-Catᵘⁿ∞ᵣₖₗ,Aʷᶠ∞-Catᵘⁿ∞ᵣ,Aʷᶠ∞-Catᵘⁿ∞ᵣₖₗ}]{The $\iy$-categories $\AfCat^{\un,\iy}_R,\AfCat^{\un,\iy}_{R,K,L},\AwCat^{\un,\iy}_R,\ab\AwCat^{\un,\iy}_{R,K,L}$}
\label{ai81}

This section is the analogue of \S\ref{ai72}. Here is a change in the filtered case.

\begin{rem}
\label{ai8rem1}
In \S\ref{ai62} and \S\ref{ai71} we defined {\it balanced\/} $\{0,\ldots,N\}$-labelled trees $\bs\Ga$ in Definition \ref{ai6def6}, and then in Definitions \ref{ai6def8} and \ref{ai7def2}, given collections of $R$-linear maps $\m_k^I$ or $\m_k^{I;x_0,\ldots,x_k}$, we defined operators $\d_{n,k}^I$ or $\d_{n,k;y_0,\ldots,y_k}^{I;x_0,\ldots,x_n}$ as sums over balanced $\bs\Ga_1,\ldots,\bs\Ga_k$ of multiple compositions of $\m_l^J$ or~$\m_l^{J;x_{j_0},\ldots,x_{j_l}}$.

Definition \ref{ai6def6}(a),(b) require that for $\bs\Ga$ to be balanced, every internal vertex $v$ in $\Ga$ should have $\md{I_v}\ge 2$, but every bud $b$ in $\Ga$ should have $\md{I_b}\ge 3$. Now a bud is just an internal vertex with no incoming edges, corresponding to an operator $\m_0^{I_v}$ or $\m_0^{I_v;x}$. In Chapters \ref{ai6} and \ref{ai7} we required that $\m_0^{I_v}=\m_0^{I_v;x}=0$ if $\md{I_v}=1$ or 2. This is the reason for the condition that $\md{I_b}\ge 3$ in Definition \ref{ai6def6}(b), as if $\md{I_b}=2$ then~$\m_0^{I_b}=\m_0^{I_b;x}=0$.

In the filtered case we no longer require $\m_0^{I_v}=\m_0^{I_v;x}=0$ if $\md{I_v}=1$ or 2, but only that $\m_0^{I_v}\in F^{\ssu 0}A^{\max I_v}$ or $\m_0^{I_v;x}\in F^{\ssu 0}\sH^{\max I_v}(f^{I_v}(x),f^{\min I_v\max I_v}(x))$. Therefore we must relax the definition of balanced, to allow $\md{I_b}\ge 2$ for buds $b$ instead of $\md{I_b}\ge 3$. We call such $\{0,\ldots,N\}$-labelled trees {\it filtered balanced}.

In the filtered case, we should define operators $\d_{n,k}^I$ or $\d_{n,k;y_0,\ldots,y_k}^{I;x_0,\ldots,x_n}$ as in Definitions \ref{ai6def8} and \ref{ai7def2}, but using filtered balanced labelled trees $\bs\Ga_1,\ldots,\bs\Ga_k$.

The proofs of Lemmas \ref{ai6lem2} and \ref{ai7lem1} use the condition that $\md{I_b}\ge 3$ for buds $b$ in $\Ga$, so these lemmas no longer hold in the filtered case. Instead, we should replace them with Lemma \ref{ai8lem1} below.

Because of the failure of Lemmas \ref{ai6lem2} and \ref{ai7lem1}, some upper limits in sums such as $k=n+M-1$ in \eq{ai6eq18}, \eq{ai6eq23}, \eq{ai7eq16}, \eq{ai7eq21} and $k=\max(n,n+j-1)$ in \eq{ai6eq28}, \eq{ai6eq30}, \eq{ai7eq25}, \eq{ai7eq27} are no longer correct in the filtered case. Instead, we should sum these up to $k=\iy$, to give an infinite sum which is convergent in the filtered sense, as explained in Definition~\ref{ai8def1}. Also limits $k=n$ in \eq{ai6eq18}, \eq{ai7eq15}, \eq{ai7eq16} should be replaced by $k=n+1$, to allow for $\m_0^{i_0}\ne 0$ or~$\m_0^{i_0;x_0}\ne 0$.
\end{rem}

\begin{lem}
\label{ai8lem1}
{\bf(a)} Suppose we define operators $\d_{n,k}^I$ using maps $\m_l^J$ as in Definition\/ {\rm\ref{ai6def8},} but no longer requiring $\m_0^I=0$ for $\md{I}=1,2,$ and using filtered balanced\/ $\bs\Ga_i$ rather than balanced $\bs\Ga_i$ in \eq{ai6eq11}. If $\md{I}=1$ then $\d_{n,k}^I=0$ unless $k\le n+1$. If\/ $\md{I}\ge 2$ and\/ $k>n+\md{I}-2$ then every term in \eq{ai6eq11} involves at least\/ $k-n-\md{I}+2$ factors $\m_0^{I_b}$ with\/~$\md{I_b}=2$.
\smallskip

\noindent{\bf(b)} The analogue of\/ {\bf(a)} holds for operators $\d_{n,k;y_0,\ldots,y_k}^{I;x_0,\ldots,x_n}$ as in Definition\/ {\rm\ref{ai7def2},} with at least\/ $k-n-\md{I}+2$ factors $\m_0^{I_b;x}$ with\/~$\md{I_b}=2$.
\end{lem}

Here is the filtered analogue of Definition \ref{ai7def3}, used to define~\eq{ai8eq1}.

\begin{dfn}
\label{ai8def1}
Let $(\Ga,\pre)$ be a po-monoid and $R$ a complete $\Ga$-filtered ring, as in Definitions \ref{ai5def3} and \ref{ai5def4}. We will define a simplicial set $\hAfCat^{\un,\iy}_R$ of unital $\Ga$-filtered $A_\iy$-categories over~$R$.
\smallskip

\noindent{\bf(a)} For $N\ge 0$, define the $N$-simplices $(\hAfCat^{\un,\iy}_R)_N$ in $\hAfCat^{\un,\iy}_R$ to be the set of quadruples $\bigl((\sO^i)_{i=0}^N,(f^{ij})_{0\le i\le j\le N},(\sH^i)_{i=0}^N,(\m^I)_{\es\ne I\subseteq\{0,\ldots,N\}}\bigr)$, where:
\begin{itemize}
\setlength{\itemsep}{0pt}
\setlength{\parsep}{0pt}
\item[(i)] $\sO^i$ is a set for $i=0,\ldots,N$, thought of as the set of objects $\Obj\sA^i$ of a filtered $A_\iy$-category $\sA^i$ at vertex $i$ of the $N$-simplex $\De^N$.
\item[(ii)] $f^{ij}:\sO^i\ra\sO^j$ is a map for all $0\le i\le j\le N$, thought of as the map $\fF_\Obj^{ij}:\Obj\sA^i\ra\Obj\sA^j$ in a filtered $A_\iy$-functor $\fF^{ij}:\sA^i\ra\sA^j$ living on the edge $[i,j]$ of $\De^N$. We require these maps to satisfy $f^{ii}=\id_{\sO^i}$ for $i=0,\ldots,N$. We do {\it not\/} require that $f^{jk}\ci f^{ij}=f^{ik}$.

Suppose $\es\ne I\subseteq\{0,1,\ldots,N\}$ with $I=\{i_0,i_1,\ldots,i_M\}$ for $i_0<i_1<\cdots<i_M$. Define $f^I:\sO^{i_0}\ra\sO^{i_M}$ as in \eq{ai7eq2}.
\item[(iii)] For all $i=0,1,\ldots,N$ and $x,y\in\sO^i$, we are given a $\Z$-graded $\Ga$-filtered $R$-module $\sH^i(x,y)=\bigop_{d\in\Z}\sH^i(x,y)^d$, thought of as the Hom space $\Hom_{\sA^i}(x,y)$ in the $A_\iy$-category $\sA^i$. 

We require that $\sH^i(x,y)^d$ is a complete filtered $R$-module for each $d\in\Z$ (this is called {\it degree-wise complete\/}), but we do not require $\sH^i(x,y)$ to be complete. We write $\sH^i=(\sH^i(x,y))_{x,y\in\sO^i}$.
\item[(iv)] For each nonempty subset $I\subseteq\{0,\ldots,N\}$, which we write as $I=\{i_0,\ab i_1,\ab\ldots,\ab i_M\}$ for $i_0<i_1<\cdots<i_M$ with $M\ge 0$, and for all $k\ge 0$ and elements $x_0,\ldots,x_k\in \sO^{i_0}$, thought of as objects in $\sA^{i_0}$, we are given an $R$-linear map 
\ea
\m_k^{I;x_0,\ldots,x_k}:\sH^{i_0}&(x_0,x_1)\ot_R\sH^{i_0}(x_1,x_2)\ot_R\cdots \ot_R\sH^{i_0}(x_{k-1},x_k)
\nonumber\\
&\longra\sH^{i_M}(f^I(x_0),f^{i_0i_M}(x_k))
\label{ai8eq3}
\ea
graded of degree $2-M-k=3-\md{I}-k$. Here when $k=0$, by convention the empty tensor product on the top line of \eq{ai8eq3} is $R$, so $\m_0^{I;x_0}$ maps
\begin{equation*}
\m_0^{I;x_0}:R\longra\sH^{i_M}(f^I(x_0),f^{i_0i_M}(x_0))^{2-M}.
\end{equation*}
We sometimes identify $\m_0^{I;x_0}$ with $\m_0^{I;x_0}(1)\!\in\!\sH^{i_M}(f^I(x_0),f^{i_0i_M}(x_0))^{2-M}$.

We require $\m_k^{I;x_0,\ldots,x_k}$ to be {\it filtration-preserving} for all $k\!\ge\! 1$ and $x_0,\ldots,x_k$ in $\sO^{i_0}$, that is, for all $\ga\in\Ga$ we have
\begin{align*}
&\m_k^{x_0,\ldots,x_k}\bigl(F^\ga(
\sH^{i_0}(x_0,x_1)\ot\cdots \ot\sH^{i_0}(x_{k-1},x_k))\bigr)\\
&\qquad \subseteq F^\ga\bigl(\sH^{i_M}(f^I(x_0),f^{i_0i_M}(x_k))\bigr).
\end{align*}

We require that $\m_0^{I;x_0}(1)\in F^{\ssu 0}\sH^{i_M}(f^I(x_0),f^{i_0i_M}(x_0))$ for all $x_0\in\sO^{i_0}$ if $\md{I}=1,2$, and $\m_0^{I;x_0}(1)\in F^0\sH^{i_M}(f^I(x_0),f^{i_0i_M}(x_0))$ for all $x_0\in\sO^{i_0}$ if~$\md{I}\ge 3$.

We write $\m^I_k=(\m_k^{I;x_0,\ldots,x_k})_{x_0,\ldots,x_k\in\sO^{i_0}}$ and $\m^I=(\m^I_k)_{k\ge 0}$.

We require these maps $\m_k^{I;x_0,\ldots,x_k}$ to satisfy the following identities generalizing \eq{ai7eq15}--\eq{ai7eq16} for all $I=\{i_0<\ab i_1<\ab\ldots<\ab i_M\}\subseteq\{0,1,\ldots,N\}$, $n\ge 0$, and $x_0,\ldots,x_n\in \sO^{i_0}$, firstly if $M=0$
\e
0=\sum_{\begin{subarray}{l} k=1,\ldots,n+1, \\
y_0,\ldots,y_k\in\sO^{i_0} \end{subarray}}\m_k^{\{i_0\};y_0,\ldots,y_k}\ci\d_{n,k;y_0,\ldots,y_k}^{\{i_0\};x_0,\ldots,x_n},
\label{ai8eq4}
\e
and secondly if $M>0$
\ea
&0=\sum_{\begin{subarray}{l} k=1,\ldots,n+1 \\
y_0,\ldots,y_k\in\sO^{i_0}  \end{subarray}}\m_k^{I;y_0,\ldots,y_k}\ci\d_{n,k;y_0,\ldots,y_k}^{\{i_0\};x_0,\ldots,x_n}
\label{ai8eq5}\\
&+(-1)^M\sum_{k\ge 0,\; y_0,\ldots,y_k\in\sO^{i_M}}\m_k^{\{i_M\};y_0,\ldots,y_k}\ci\d_{n,k;y_0,\ldots,y_k}^{I;x_0,\ldots,x_n}.
\nonumber
\ea
Here $\d_{n,k;y_0,\ldots,y_k}^{I;x_0,\ldots,x_n}$ is defined using the $\m_l^{J;w_0,\ldots,w_l}$ as in Definition \ref{ai7def2}, but modified as in Remark \ref{ai8rem1} to sum over {\it filtered balanced\/} trees $\bs\Ga_i$ in~\eq{ai7eq9}. 

The second line of \eq{ai8eq5} is an infinite sum: for each value of $k$, there are only finitely many possible $y_0,\ldots,y_k$ with $\d_{n,k;y_0,\ldots,y_k}^{I;x_0,\ldots,x_n}\ne 0$, but we can have $\d_{n,k;y_0,\ldots,y_k}^{I;x_0,\ldots,x_n}\ne 0$ for infinitely many $k\gg 0$. We claim that it converges in the filtered sense. 

To justify this, note that by Lemma \ref{ai8lem1}(b), if $k>n+M-1$ then $\d_{n,k;y_0,\ldots,y_k}^{I;x_0,\ldots,x_n}$ contains at least $k-n-M+1$ factors $\m_0^{I_b;x}$ with $\md{I_b}=2$. Here $I_b=\{i_p,i_q\}$ for $0\le p<q\le M$ and $x=f^P(x_j)$ for $j=0,\ldots,n$ and $P\subseteq \{i_0,\ldots,i_p\}$ is the prefix, with $i_0,i_p\in P$. Thus there are only finitely many possibilities for $\m_0^{I_b;x}$, and each lies in $F^{\ssu 0}\sH^{i_b}(f^{I_b}(x),f^{I_b}(x))$, so there exists $\de\ssu 0$ in $\Ga$ with~$\m_0^{I_b;x}\in F^\de\sH^{i_b}(f^{I_b}(x),f^{I_b}(x))$.

By Definition \ref{ai5def3}(b) we can take the greatest lower bound $\ep\in\Ga$ of these $\de\in\Ga$ over all possibilities for $\m_0^{I_b;x}$, and by Definition \ref{ai5def3}(d) we have $\ep\ssu 0$. That is, there exists a single $\ep\ssu 0$ in $\Ga$ such that for all choices of $\m_0^{I_b;x}$ we have $\m_0^{I_b;x}\in F^\ep\sH^{i_b}(f^{I_b}(x),f^{I_b}(x))$. Therefore, if $k>n+M-1$ then $\d_{n,k;y_0,\ldots,y_k}^{I;x_0,\ldots,x_n}$ in \eq{ai8eq5} increases filtrations by at least $(k-n-M+1)\ep$, that is, it maps $F^\ga(\cdots)\ra F^{\ga+(k-n-M+1)\ep}(\cdots)$. Since $\ga+(k-n-M+1)\ep\ra\iy$ in $\Ga$ as $k\ra\iy$, it follows that when applied to any element of $\sH^{i_0}(x_0,x_1)\ot\cdots\ot\sH^{i_0}(x_{n-1},x_n)$, equation \eq{ai8eq5} is a convergent sum.
\item[(v)] There should exist elements $\boo^i_x\in F^0\sH^i(x,x)^0$ for each $i=0,\ldots,N$ and $x\in\sO^i$ satisfying \eq{ai7eq17}--\eq{ai7eq19}. They are uniquely determined by~\eq{ai7eq17}.
\end{itemize}

\noindent{\bf(b)} If $\al:[\bar  N]\ra[N]$ is a morphism in $\De$, we define
\begin{equation*}
\hAfCat^{\un,\iy}_R(\al):(\hAfCat^{\un,\iy}_R)_N\longra\hAfCat^{\un,\iy}_{R,\bar N}
\end{equation*}
as in Definition \ref{ai7def3}(b). Essentially the same proof as Proposition \ref{ai7prop1} shows that $\hAfCat^{\un,\iy}_R(\al)$ maps to $\hAfCat^{\un,\iy}_{R,\bar N}$. As in Definition \ref{ai7def3}(b), $\hAfCat^{\un,\iy}_R:\De^\rop\ra\Set$ is a functor, that is, $\hAfCat^{\un,\iy}_R$ is a simplicial set.
\end{dfn}

We also give an analogue with {\it weakly\/} filtered $N$-morphisms for $N\ge 2$, used to define~\eq{ai8eq2}.

\begin{dfn}
\label{ai8def2}
Let $(\Ga,\pre)$ be a po-monoid and $R$ a complete $\Ga$-filtered ring, as in Definitions \ref{ai5def3} and \ref{ai5def4}. We will define a simplicial set $\hAwCat^{\un,\iy}_R$ of unital $\Ga$-filtered $A_\iy$-categories over~$R$.

For $N\ge 0$, the $N$-simplices $(\hAwCat^{\un,\iy}_R)_N$ in $\hAwCat^{\un,\iy}_R$ are defined as for $(\hAfCat^{\un,\iy}_R)_N$ in Definition \ref{ai8def1}(a), but with the following changes to (a)(iv): there should exist $\de\pre 0$ in $\Ga$ independent of $I,k,x_0,\ldots,x_k$ such that the following holds. 

The $R$-linear maps $\m_k^{I;x_0,\ldots,x_k}$ in \eq{ai8eq3} should be filtration-preserving for $\md{I}=1,2$ and $k\ge 1$, but for $\md{I}\ge 3$ and $k\ge 1$ they need not preserve filtrations, but should instead map $F^\ga(\cdots)\ra F^{\ga+\de}(\cdots)$ for all~$\ga\in\Ga$.

We require that $\m_0^{I;x_0}(1)\in F^{\ssu 0}\sH^{i_M}(f^I(x_0),f^{i_0i_M}(x_0))$ for all $x_0\in\sO^{i_0}$ if $\md{I}=1,2$, and $\m_0^{I;x_0}(1)\in F^\de\sH^{i_M}(f^I(x_0),f^{i_0i_M}(x_0))$ for all $x_0\in\sO^{i_0}$ if~$\md{I}\ge 3$ (this replaces $\m_0^{I;x_0}(1)\in F^0\sH^{i_M}(f^I(x_0),f^{i_0i_M}(x_0))$).

Each term in $\d_{n,k;y_0,\ldots,y_k}^{I;x_0,\ldots,x_n}$ contains at most $M-1$ factors $\m_l^{J;x'_0,\ldots,x'_l}$ with $\md{J}\ge 3$. Thus the argument in Definition \ref{ai8def1}(a)(iv) that \eq{ai8eq5} converges in the filtered sense works here, but we should include an extra shift in filtrations $F^\ga(\cdots)\mapsto F^{\ga+(M-1)\de}(\cdots)$ when $M\ge 2$, and conclude that the infinite sums converge because $\ga+(M-1)\de+(k-n-M+1)\ep\ra\iy$ in $\Ga$ as~$k\ra\iy$.

The rest of the definition of $\hAwCat^{\un,\iy}_R$ is as in Definition \ref{ai8def1}. Note that $\hAfCat^{\un,\iy}_R\subset\hAwCat^{\un,\iy}_R$.
\end{dfn}

The following analogue of Theorem \ref{ai7thm2} will be proved in~\S\ref{ai103}.

\begin{thm}
\label{ai8thm1}
Let\/ $\bigl((\sO^i)_{i=0}^N,(f^{ij})_{0\le i\le j\le N},(\sH^i)_{i=0}^N,(\m^I)_{\es\ne I\subseteq\{0,\ldots,N\}}\bigr)$ be an $N$-simplex in $\hAfCat^{\un,\iy}_R$ or $\hAwCat^{\un,\iy}_R$ as in Definitions\/ {\rm\ref{ai8def1}} and\/ {\rm\ref{ai8def2},} and define maps\/ $\d_{n,k;y_0,\ldots,y_k}^{I;x_0,\ldots,x_n}$ from the $\m_k^{I;x_0,\ldots,x_k}$ as in Remark\/ {\rm\ref{ai8rem1}}. Then for each nonempty subset\/ $I\subseteq\{0,\ldots,N\},$ written $I=\{i_0,\ab i_1,\ab\ldots,\ab i_M\}$ for $i_0<i_1<\cdots<i_M,$ and all\/ $k,n\ge 0,$ $x_0,\ldots,x_n\in\sO^{i_0}$ and\/ $y_0,\ldots,y_k\in\sO^{i_M},$ in $R$-linear maps as in {\rm\eq{ai7eq6},} generalizing \eq{ai7eq21} we have
\e
\begin{aligned}
&\sum_{\begin{subarray}{l} k-1\le l\le n+1, \\ w_0,\ldots,w_l\in\sO^{i_0} \end{subarray}}^n\d_{l,k;y_0,\ldots,y_k}^{\{i_0\};w_0,\ldots,w_l}\ci\d_{n,l;w_0,\ldots,w_l}^{\{i_0\};x_0,\ldots,x_n}=0 && \text{if\/ $M=0,$} \\
&\begin{aligned}[t]& \sum_{\begin{subarray}{l} 1\le l\le n+1, \\ w_0,\ldots,w_l\in\sO^{i_0} \end{subarray}} \d_{l,k;y_0,\ldots,y_k}^{I;w_0,\ldots,w_l}\ci\d_{n,l;w_0,\ldots,w_l}^{\{i_0\};x_0,\ldots,x_n}\\
&+(-1)^M\sum_{l\ge 0, \;  w_0,\ldots,w_l\in\sO^{i_M}}\d_{l,k;y_0,\ldots,y_k}^{\{i_M\};w_0,\ldots,w_l}\ci\d_{n,l;w_0,\ldots,w_l}^{I;x_0,\ldots,x_n}=0\end{aligned} && \text{if\/ $M>0,$}
\end{aligned}
\label{ai8eq6}
\e
where the final sum may have infinitely many nonzero terms, but converges in the filtered sense. 
\end{thm}

Here is the analogue of Example \ref{ai7ex2}.

\begin{ex}
\label{ai8ex1}
{\bf(a)} A 0-simplex $(\sO^0,f^{00}=\id_{\sO^0},\sH^0,\m^{\{0\}})$ of $\hAfCat^{\un,\iy}_R$ is equivalent to a unital filtered $A_\iy$-category $\sA^0$ over $R$, as in Definition \ref{ai5def8}, with objects $\Obj\sA^0=\sO^0$, Hom spaces $\Hom_{\sA^0}(x,y)=\sH^0(x,y)$, morphisms $(\m_k^{\{0\};x_0,\ldots,x_k})_{k\ge 0}$, and units~$\boo^0_x\in \Hom_{\sA^0}^0(x,x)=\sH^0(x,x)^0$.
\smallskip

\noindent{\bf(b)} A 1-simplex of $\hAfCat^{\un,\iy}_R$ is equivalent to a unital filtered $A_\iy$-functor $\fF^{01}:\sA^0\ra\sA^1$ of unital filtered $A_\iy$-categories, as in Definition \ref{ai5def9}, with object map $\fF_\Obj^{01}=f^{01}$ and morphisms $\fF^{01;x_0,\ldots,x_k}_k=\m_k^{\{0,1\};x_0,\ldots,x_k}$. 
\smallskip

\noindent{\bf(c)} A 2-simplex of $\hAfCat^{\un,\iy}_R$ is equivalent to unital filtered $A_\iy$-categories $\sA^0,\ab\sA^1,\ab\sA^2$ over $R$ as in {\bf(a)}, and unital filtered $A_\iy$-functors $\fF^{ij}:\sA^i\ra\sA^j$ for $(i,j)=(0,1),\ab(0,2),\ab(1,2)$ as in {\bf(b)}, and $R$-linear maps for all $k\ge 0$ and $x_0,\ldots,x_k\in\sO^0$
\begin{align*}
\m_k^{\{0,1,2\};x_0,\ldots,x_k}:\sH^0&(x_0,x_1)\ot_R\sH^0(x_1,x_2)\ot_R\cdots \ot_R\sH^0(x_{k-1},x_k)
\nonumber\\
&\longra\sH^2(f^{12}\ci f^{01}(x_0),f^{02}(x_k)),
\end{align*}
graded of degree $-k$, which are filtration-preserving, satisfying \eq{ai7eq19} and for all $n\ge 0$ and $x_0,\ldots,x_n\in\sO^0$, a generalization of \eq{ai7eq22}
\ea
&\sum_{\begin{subarray}{l} i,j,k\ge 0:\\ i+j+k=n \end{subarray}}\begin{aligned}[t]&
(-1)^{i+jk}\m^{\{0,1,2\};x_0,x_1,\ldots,x_i,x_{i+j},x_{i+j+1},\ldots,x_n}_{i+k+1}\ci\\
&\bigl(\id_{\sH^0(x_0,x_1)}\ot\cdots\ot \id_{\sH^0(x_{i-1},x_i)}\ot\m^{\{0\};x_i,x_{i+1},\ldots,x_{i+j}}_j\\
&\ot \id_{\sH^0(x_{i+j},x_{i+j+1})}\ot\cdots\ot \id_{\sH^0(x_{n-1},x_n)}\bigr)\end{aligned}
\label{ai8eq7}\\
&+\sum_{\!\!\!\!\!\!\begin{subarray}{l} 
n=n_1+\cdots+n_k, \\
1\le l\le k, \; n_i\ge 0. \\
\text{Set $p_j=\sum_{i=1}^jn_i$,} \\ 
0\le j\le k \end{subarray}}
\begin{aligned}[t]&(-1)^{k+1}\cdot(-1)^{\sum_{i=1}^{l-1}(n_i+1)}\cdot(-1)^{\sum_{i=1}^k(n_i+1)(k-i)}\\
&\m^{\{2\};f^{\{0,1,2\}}(x_{p_0}),f^{\{0,1,2\}}(x_{p_1}),\ldots,f^{\{0,1,2\}}(x_{p_{l-1}}),f^{02}(x_{p_l}),\ldots,f^{02}(x_{p_k})}_k\\
&\ci\bigl[(\m^{\{1,2\}}\!\ci\!\m^{\{0,1\}})^{x_{p_0},\ldots,x_{p_1}}_{n_1}\!\ot\!\cdots\!\ot\!(\m^{\{1,2\}}\!\ci\!\m^{\{0,1\}})^{x_{p_{l-2}},\ldots,x_{p_{l-1}}}_{n_{l-1}}\\
&\ot \m^{\{0,1,2\};x_{p_{l-1}},\ldots,x_{p_l}}_{n_l}\ot\m^{\{0,2\};x_{p_l},\ldots,x_{p_{l+1}}}_{n_{l+1}}\ot\\
&\cdots\ot\m^{\{0,2\};x_{p_{k-1}},\ldots,x_{p_k}}_{n_k}\bigr]=0.
\end{aligned}
\nonumber
\ea

Comparing \eq{ai4eq15} and \eq{ai8eq7}, we see $\m^{\{0,1,2\}}$ is equivalent to a {\it unital filtered\/ $A_\iy$-natural transformation\/} $\fH^{012}:\fF^{12}\ci\fF^{01}\Ra\fF^{02}$, as in Definition~\ref{ai5def12}(a).
\smallskip

\noindent{\bf(d)} In the analogues of {\bf(a)\rm--\bf(c)} for $\hAwCat^{\un,\iy}_R$, {\bf(a)\rm,\bf(b)} are unchanged, but in {\bf(c)}, instead of taking $\m_k^{\{0,1,2\};x_0,\ldots,x_k}$ to be filtration-preserving, as in Definition \ref{ai8def2} there should exist $\de\pre 0$ in $\Ga$ independent of $k,x_0,\ldots,x_k$ such that $\m_k^{\{0,1,2\};x_0,\ldots,x_k}$ maps $F^\ga(\cdots)\ra F^{\ga+\de}(\cdots)$ for all $\ga\in\Ga$. Hence $\m^{\{0,1,2\}}$ is equivalent to a {\it weakly filtered\/ $A_\iy$-natural transformation\/} $\fF^{12}\ci\fF^{01}\Ra\fF^{02}$ in the sense of Definition~\ref{ai5def12}(b).
\end{ex}

Here is the analogue of Definition~\ref{ai7def4}.

\begin{dfn}
\label{ai8def3}
Define a simplicial subset $\AfCat^{\un,\iy}_R\subset\hAfCat^{\un,\iy}_R$, by an $N$-simplex $\bigl((\sO^i)_{i=0}^N,(f^{ij})_{0\le i\le j\le N},(\sH^i)_{i=0}^N,(\m^I)_{\es\ne I\subseteq\{0,\ldots,N\}}\bigr)$ in $\hAfCat^{\un,\iy}_R$ lies in $\AfCat^{\un,\iy}_R$ if and only if whenever $0\le i_0<i_1<i_2\le N$, either (hence both) of the following two equivalent conditions hold:
\begin{itemize}
\setlength{\itemsep}{0pt}
\setlength{\parsep}{0pt}
\item[(a)] The unital filtered $A_\iy$-natural transformation $\m^{\{i_0,i_1,i_2\}}:\fF^{i_1i_2}\ci\fF^{i_0i_1}\Ra\fF^{i_0i_2}$ in Example \ref{ai8ex1}(c) is a filtered $A_\iy$-natural isomorphism, in the sense of Definition~\ref{ai5def12}(c).
\item[(b)] $[\m^{\{i_0,i_1,i_2\}}]:\fF^{i_1i_2}\ci\fF^{i_0i_1}\ra\fF^{i_0i_2}$ is an isomorphism in the cohomology category~$H^0(F^0\Fun_{A^\rfi_\iy}^\un(\sA^{i_0},\sA^{i_2}))$.
\end{itemize}

Here (a),(b) are equivalent by Definition \ref{ai5def12}(c). If $\al:[\bar N]\ra[N]$ is a morphism in $\De$ then $\hAfCat^{\un,\iy}_R(\al)$ maps $(\AfCat^{\un,\iy}_R)_N\ra(\AfCat^{\un,\iy}_R)_{\bar N}$, so $\AfCat^{\un,\iy}_R$ is a simplicial set.

Similarly, define a simplicial subset $\AwCat^{\un,\iy}_R\subset\hAwCat^{\un,\iy}_R$, by an $N$-simplex $\bigl((\sO^i)_{i=0}^N,(f^{ij})_{0\le i\le j\le N},(\sH^i)_{i=0}^N,(\m^I)_{\es\ne I\subseteq\{0,\ldots,N\}}\bigr)$ in $\hAwCat^{\un,\iy}_R$ lies in $\AwCat^{\un,\iy}_R$ if and only if whenever $0\le i_0<i_1<i_2\le N$, either (hence both) of the following two equivalent conditions hold:
\begin{itemize}
\setlength{\itemsep}{0pt}
\setlength{\parsep}{0pt}
\item[(c)] The unital weakly filtered $A_\iy$-natural transformation $\m^{\{i_0,i_1,i_2\}}:\fF^{i_1i_2}\ci\fF^{i_0i_1}\Ra\fF^{i_0i_2}$ in Example \ref{ai8ex1}(d) is a weakly filtered $A_\iy$-natural isomorphism, in the sense of Definition~\ref{ai5def12}(d).
\item[(d)] $[\m^{\{i_0,i_1,i_2\}}]:\fF^{i_1i_2}\ci\fF^{i_0i_1}\ra\fF^{i_0i_2}$ is an isomorphism in the cohomology category~$H^0(\Fun_{A^\rfi_\iy}^\un(\sA^{i_0},\sA^{i_2}))$.
\end{itemize}
\end{dfn}

\begin{rem}
\label{ai8rem2}
Definition \ref{ai7def4} gave four equivalent conditions (i)--(iv) for an $N$-simplex in $\hACat^{\un,\iy}_R$ to 
lie in $\ACat^{\un,\iy}_R$, where (i)--(ii) are equivalent to (iii)--(iv) by Proposition \ref{ai4prop3}. In the (weakly) filtered case, Definition \ref{ai7def4}(i)--(ii) do not make sense, and we have no analogue of Proposition \ref{ai4prop3}, but (a)--(b) and (c)--(d) above are analogues of Definition \ref{ai7def4}(iii)--(iv).
\end{rem}

The following analogue of Theorem \ref{ai7thm3}, one of our main results, is proved in \S\ref{ai104}. The proof is long and difficult.

\begin{thm}
\label{ai8thm2}
$\AfCat^{\un,\iy}_R$ and\/ $\AwCat^{\un,\iy}_R$ above are $\iy$-categories.
\end{thm}

The analogues of Definition \ref{ai7def5} and Proposition \ref{ai7prop2} also work in the (weakly) filtered cases.

Here is the analogue of Theorem \ref{ai7thm4}. It will be proved in~\S\ref{ai105}.

\begin{thm}
\label{ai8thm3}
In the simplicial sets\/ $\hAfCat^{\un,\iy}_R$ and\/ $\hAwCat^{\un,\iy}_R$ in Definitions\/ {\rm\ref{ai8def1}} and\/ {\rm\ref{ai8def2},} define a $2$-simplex to be \begin{bfseries}thin\end{bfseries} if the unital (weakly) filtered\/ $A_\iy$-natural transformation $\fH^{012}:\fF^{12}\ci\fF^{01}\Ra\fF^{02}$ in Example\/ {\rm\ref{ai8ex1}(c)} is a (weakly) filtered\/ $A_\iy$-natural isomorphism, that is, if it is a $2$-simplex in $\AfCat^{\un,\iy}_R\subset\hAfCat^{\un,\iy}_R$ or\/ $\AwCat^{\un,\iy}_R\subset\hAwCat^{\un,\iy}_R$.

Write $\cT^\rfi\subset\hAfCat^{\un,\iy}_R([2])$ and\/ $\cT^\wf\subset\hAwCat^{\un,\iy}_R([2])$ for the sets of thin $2$-simplices. Then $(\hAfCat^{\un,\iy}_R,\cT^\rfi)$ and\/ $(\hAwCat^{\un,\iy}_R,\cT^\wf)$ are \begin{bfseries}$\iy$-bicategories\end{bfseries}, as in\/ {\rm\S\ref{ai28}}. We have $\AfCat^{\un,\iy}_R=\Core(\hAfCat^{\un,\iy}_R,\cT^\rfi)$ and\/ $\AwCat^{\un,\iy}_R=\Core(\hAwCat^{\un,\iy}_R,\cT^\wf),$ as in Definition\/~{\rm\ref{ai2def22}}.
\end{thm}

Here is the analogue of Definition~\ref{ai7def6}.

\begin{dfn}
\label{ai8def4}
Let $K,L\in\{0,1,2,\ldots,\iy\}$ with $K\le L$ and $L\ge 2$. Define a simplicial subset $\AfCat^{\un,\iy}_{R,K,L}\subset\AfCat^{\un,\iy}_R$ such that an $N$-simplex $\bigl((\sO^i)_{i=0}^N,(f^{ij})_{0\le i\le j\le N},(\sH^i)_{i=0}^N,(\m^I)_{\es\ne I\subseteq\{0,\ldots,N\}}\bigr)$ in $\AfCat^{\un,\iy}_R$ in Definition \ref{ai8def1}(a) lies in $\AfCat^{\un,\iy}_{R,K,L}$ if:
\begin{itemize}
\setlength{\itemsep}{0pt}
\setlength{\parsep}{0pt}
\item[(i)] If $I=\{i_0<i_1<\cdots<i_M\}\subseteq\{0,\ldots,N\}$ with $\md{I}\ne 3$ and $x_0\in\sO^{i_0}$, we have $\m_0^{I;x_0}\!=\!0$ if $\md{I}\!\le\! K$, and $\m_0^{I;x_0}(1)\!\in\! F^{\ssu 0}\Hom_{\sA^{i_M}}(f^I(x_0),f^{i_0i_M}(x_0))$ if $K<\md{I}\le L$.
\item[(ii)] If $L\ge 3$ then $f^{ik}=f^{jk}\ci f^{ij}$ for all $0\le i\le j\le k\le N$, 
\item[(iii)] If $L\ge 3$, $I=\{i_0<i_1<i_2\}\subseteq\{0,\ldots,N\}$ and $x_0\in\sO^{i_0}$ then $\m_0^{I;x_0}(1)=\boo^{i_2}_{f^{i_0i_2}(x_0)}$ if $K\ge 3$ (this makes sense by (ii) as $f^I(x_0)=f^{i_0i_2}(x_0)$), and $\m_0^{I;x_0}(1)-\boo^{i_2}_{f^{i_0i_2}(x_0)}\in F^{\ssu 0}\Hom_{\sA^{i_2}}(f^{i_0i_2}(x_0),f^{i_0i_2}(x_0))$ if~$K\le 2$.
\end{itemize}

These conditions are preserved by the maps in Definition \ref{ai8def1}(b), and so define a simplicial subset $\AfCat^{\un,\iy}_{R,K,L}\subset\AfCat^{\un,\iy}_R$. When $K=0$ and $L=2$ part (i) is included in Definition \ref{ai8def1}(a)(iv) and (ii),(iii) are trivial, so $\AfCat^{\un,\iy}_{R,0,2}=\AfCat^{\un,\iy}_R$. Clearly $\AfCat^{\un,\iy}_{R,K,L}\subseteq\AfCat^{\un,\iy}_{R,K',L'}$ if $K\ge K'$ and $L\ge L'$, as in \eq{ai8eq1}. Theorem \ref{ai8thm6} will show that $\AfCat^{\un,\iy}_{R,K,L}$ is an $\iy$-category. We say that $\AfCat^{\un,\iy}_{R,K,L}$ has {\it type\/}~$(K,L)$. 

We also define $\iy$-subcategories $\AwCat^{\un,\iy}_{R,K,L}\subset\AwCat^{\un,\iy}_R$ as in \eq{ai8eq2} exactly as above, by the same (i)--(iii), but applied to $N$-simplices in $\AwCat^{\un,\iy}_R$, plus the extra condition which will be needed in the proof of Theorem~\ref{ai8thm6}:
\begin{itemize}
\setlength{\itemsep}{0pt}
\setlength{\parsep}{0pt}
\item[(iv)] If $L\ge 3$, $I=\{i_0<i_1<i_2\}\subseteq\{0,\ldots,N\},$ $n\ge 1$, and $x_0,\ldots,x_n\in\sO^{i_0}$, then $\m_n^{I;x_0,\ldots,x_n}$ is filtration-preserving.
\end{itemize}
\end{dfn}

\begin{ex}
\label{ai8ex2}
If $K\le 2$ and $L\ge 3$ then modifying Example \ref{ai8ex1}(c), a 2-simplex in $\AfCat^{\un,\iy}_{R,K,L}$ consists of unital filtered $A_\iy$-categories $\sA^0,\sA^1,\sA^2$ over $R$ (which are strict if $K>0$) and unital filtered $A_\iy$-functors $\fF^{ij}:\sA^i\ra\sA^j$ for $(i,j)=(0,1),\ab(0,2),\ab(1,2)$ (which are strict if $K>1$), and a filtered $A_\iy$-natural isomorphism $\m^{\{0,1,2\}}:\fF^{12}\ci\fF^{01}\Ra\fF^{02}$, such that $\m_0^{\{0,1,2\};x_0}(1)-\boo^{i_2}_{f^{02}(x_0)}$ lies in $F^{\ssu 0}\Hom_{\sA^2}(f^{02}(x_0),f^{02}(x_0))$ for all $x_0\in\sO^0$. 

By Proposition \ref{ai5prop1}, such $\m^{\{0,1,2\}}$ are in 1-1 correspondence with {\it unital filtered\/ $A_\iy$-homotopies\/} $\fH:\fF^{12}\ci\fF^{01}\Ra\fF^{02}$, as in Definition \ref{ai5def11}. Thus, 2-morphisms in $\AfCat^{\un,\iy}_{R,0,3}$ correspond to unital filtered $A_\iy$-homotopies.
\end{ex}

\begin{dfn}
\label{ai8def5}
Apart from the extra data of filtrations, the important difference between the definitions of $\ACat^{\un,\iy}_R$ in Definitions \ref{ai7def3} and \ref{ai7def4} and $\AfCat^{\un,\iy}_R$ in Definitions \ref{ai8def1} and \ref{ai8def3} is that in $N$-simplices in $\ACat^{\un,\iy}_R$ we have $\m_0^{I;x_0}=0$ if $\md{I}=1,2$, but in $N$-simplices in $\AfCat^{\un,\iy}_R$ we allow $0\ne \m_0^{I;x_0}\in F^{\ssu 0}(\cdots)$ if $\md{I}=1,2$. However, in the $\iy$-subcategory $\AfCat^{\un,\iy}_{R,2,2}\subset\AfCat^{\un,\iy}_R$ we require that $\m_0^{I;x_0}=0$ if $\md{I}=1,2$. 

Thus, $N$-simplices of $\AfCat^{\un,\iy}_{R,2,2}$ are $N$-simplices of $\ACat^{\un,\iy}_R$ with extra filtration data, and similarly for $\AwCat^{\un,\iy}_{R,2,2}$. Hence forgetting filtrations gives $\iy$-functors
\ea
\Pi_{A_\iy^\rfi}^{A_\iy}&:\AfCat^{\un,\iy}_{R,2,2}\longra\ACat^{\un,\iy}_R,
\label{ai8eq8}\\
\Pi_{A_\iy^\wf}^{A_\iy}&:\AwCat^{\un,\iy}_{R,2,2}\longra\ACat^{\un,\iy}_R.
\label{ai8eq9}
\ea
\end{dfn}

\subsection{\texorpdfstring{Other $\iy$-categories of filtered $A_\iy$-categories}{Other ∞-categories of filtered A∞-categories}}
\label{ai82}

Here is the filtered analogue of Definition~\ref{ai7def8}.

\begin{dfn}
\label{ai8def6}
Let $(\Ga,\pre)$ be a po-monoid, $R$ a complete $\Ga$-filtered ring, and $K,L\in\{0,1,2,\ldots,\iy\}$ with $K\le L$ and $L\ge 3$. We will define a simplicial set $\dAfCat^{\un,\iy}_{R,K,L}$ of unital  filtered $A_\iy$-categories over~$R$.
\smallskip

\noindent{\bf(a)} For $N\ge 0$, define the $N$-simplices $(\dAfCat^{\un,\iy}_{R,K,L})_N$ in $\dAfCat^{\un,\iy}_{R,K,L}$ to be the set of quadruples $\bigl((\sO^i)_{i=0}^N,(f^{ij})_{0\le i\le j\le N},(\sH^i)_{i=0}^N,(\m^I)_{\es\ne I\subseteq\{0,\ldots,N\}}\bigr)$, where
\begin{itemize}
\setlength{\itemsep}{0pt}
\setlength{\parsep}{0pt}
\item[(i)] $\sO^i$ is a set for $i=0,\ldots,N$, thought of as the set of objects $\Obj\sA^i$ of a filtered $A_\iy$-category $\sA^i$ at vertex $i$ of the $N$-simplex $\De^N$.
\item[(ii)] $f^{ij}:\sO^i\ra\sO^j$ is a map for all $0\le i\le j\le N$, thought of as the map $\fF_\Obj^{ij}:\Obj\sA^i\ra\Obj\sA^j$ in a filtered $A_\iy$-functor $\fF^{ij}:\sA^i\ra\sA^j$ living on the edge $[i,j]$ of $\De^N$. We require these maps to satisfy $f^{ii}=\id_{\sO^i}$ for $i=0,\ldots,N$ and $f^{jk}\ci f^{ij}=f^{ik}$ for $0\le i\le j\le k\le N$.
\item[(iii)] For all $i=0,1,\ldots,N$ and $x,y\in\sO^i$, we are given a $\Z$-graded $\Ga$-filtered $R$-module $\sH^i(x,y)=\bigop_{d\in\Z}\sH^i(x,y)^d$, thought of as the Hom space $\Hom_{\sA^i}(x,y)$ in the $A_\iy$-category $\sA^i$. 

We require that $\sH^i(x,y)^d$ is a complete filtered $R$-module for each $d\in\Z$ (this is called {\it degree-wise complete\/}), but we do not require $\sH^i(x,y)$ to be complete. We write $\sH^i=(\sH^i(x,y))_{x,y\in\sO^i}$.
\item[(iv)] For each nonempty subset $I\subseteq\{0,\ldots,N\}$, which we write as $I=\{i_0,\ab i_1,\ab\ldots,\ab i_M\}$ for $i_0<i_1<\cdots<i_M$ with $M\ge 0$, and for all $k\ge 0$ and elements $x_0,\ldots,x_k\in \sO^{i_0}$, thought of as objects in $\sA^{i_0}$, we are given an $R$-linear map
\ea
\m_k^{I;x_0,\ldots,x_k}:\sH^{i_0}&(x_0,x_1)\ot_R\sH^{i_0}(x_1,x_2)\ot_R\cdots \ot_R\sH^{i_0}(x_{k-1},x_k)
\nonumber\\
&\longra\sH^{i_M}(f^{i_0i_M}(x_0),f^{i_0i_M}(x_k))
\label{ai8eq10}
\ea
graded of degree $2-M-k=3-\md{I}-k$.

We require $\m_k^{I;x_0,\ldots,x_k}$ to be {\it filtration-preserving} for all $k\!\ge\! 1$ and $x_0,\ldots,x_k$ in $\sO^{i_0}$, that is, for all $\ga\in\Ga$ we have
\begin{align*}
&\m_k^{x_0,\ldots,x_k}\bigl(F^\ga(
\sH^{i_0}(x_0,x_1)\ot\cdots \ot\sH^{i_0}(x_{k-1},x_k))\bigr)\\
&\qquad \subseteq F^\ga\bigl(\sH^{i_M}(f^{i_0i_M}(x_0),f^{i_0i_M}(x_k))\bigr).
\end{align*}

For all $x_0\in\sO^{i_0}$ we require that $\m_0^{I;x_0}=0$ if $\md{I}\le K$, and $\m_0^{I;x_0}(1)\!\in\! F^{\ssu 0}\Hom_{\sA^{i_M}}(f^{i_0i_M}(x_0),f^{i_0i_M}(x_0))$ if $K<\md{I}\le L$, and $\m_0^{I;x_0}(1)$ lies in $F^0\sH^{i_M}(f^{i_0i_M}(x_0),f^{i_0i_M}(x_0))$ if~$\md{I}>L$.

We require these maps $\m_k^{I;x_0,\ldots,x_k}$ to satisfy the following identity generalizing \eq{ai7eq25} for all $I=\{i_0<\cdots<i_M\}\subseteq\{0,\ldots,N\}$, $n\ge 0$ and objects~$x_0,\ldots,x_n\in\sO^{i_0}$
\ea
&\sum_{0\le j\le M, \; k\ge 0, \; y_0,\ldots,y_k\in\sO^{i_j}\!\!\!\!\!\!\!\!\!\!\!\!\!\!\!\!\!\!\!\!\!\!\!\!}
(-1)^j\,\m_k^{\{i_j,\ldots,i_M\};y_0,\ldots,y_k}\ci\d_{n,k;y_0,\ldots,y_k}^{\{i_0,\ldots,i_j\};x_0,\ldots,x_n}
\nonumber\\
&=\sum_{j=1,\ldots,M-1}(-1)^j\m_n^{I\sm\{i_j\};x_0,\ldots,x_n}.
\label{ai8eq11}
\ea
Here $\d_{n,k;y_0,\ldots,y_k}^{I;x_0,\ldots,x_n}$ is defined using the $\m_l^{J;w_0,\ldots,w_l}$ as in Definition \ref{ai7def2}, but modified as in Remark \ref{ai8rem1} to sum over {\it filtered balanced\/} trees $\bs\Ga_i$ in \eq{ai7eq9}. Equation \eq{ai8eq11}  converges in the filtered sense as for~\eq{ai8eq5}.
\item[(v)] Definition \ref{ai8def1}(a)(v) holds.
\end{itemize}

\noindent{\bf(b)} If $\al:[\bar  N]\ra[N]$ is a morphism in $\De$, we define
\begin{equation*}
\dAfCat^{\un,\iy}_{R,K,L}(\al):(\dAfCat^{\un,\iy}_{R,K,L})_N\longra(\dAfCat^{\un,\iy}_{R,K,L})_{\bar N}
\end{equation*}
as in Definition \ref{ai7def8}(b). Essentially the same proof as Proposition \ref{ai7prop3} shows that $\dAfCat^{\un,\iy}_{R,K,L}(\al)$ maps to $(\dAfCat^{\un,\iy}_{R,K,L})_{\bar N}$. As in Definition \ref{ai7def8}(b), $\dAfCat^{\un,\iy}_{R,K,L}:\De^\rop\ra\Set$ is a functor, i.e.\ $\dAfCat^{\un,\iy}_{R,K,L}$ is a simplicial set.

Theorem \ref{ai8thm6} will show that $\dAfCat^{\un,\iy}_{R,K,L}$ is an $\iy$-category. We say that $\dAfCat^{\un,\iy}_{R,K,L}$ has {\it type\/} $(K,L)$. Clearly $\dAfCat^{\un,\iy}_{R,K,L}\subseteq\dAfCat^{\un,\iy}_{R,K',L'}$ if $K\ge K'$ and $L\ge L'$, as in~\eq{ai8eq1}.
\end{dfn}

As for Definition \ref{ai8def2}, we also give an analogue with {\it weakly\/} filtered $N$-morphisms for $N\ge 2$, used to define~\eq{ai8eq2}.

\begin{dfn}
\label{ai8def7}
Let $(\Ga,\pre)$ be a po-monoid, $R$ a complete $\Ga$-filtered ring, and $K,L\in\{0,1,2,\ldots,\iy\}$ with $K\le L$ and $L\ge 3$. We will define a simplicial set $\dAwCat^{\un,\iy}_{R,K,L}$ of unital  filtered $A_\iy$-categories over~$R$.

For $N\ge 0$, the $N$-simplices $(\dAwCat^{\un,\iy}_{R,K,L})_N$ in $\dAwCat^{\un,\iy}_{R,K,L}$ are defined as for $(\dAfCat^{\un,\iy}_{R,K,L})_N$ in Definition \ref{ai8def6}(a), but with the following changes to (a)(iv): there should exist $\de\pre 0$ in $\Ga$ independent of $I,k,x_0,\ldots,x_k$ such that the following holds. 

The $R$-linear maps $\m_k^{I;x_0,\ldots,x_k}$ in \eq{ai8eq10} should be filtration-preserving for $\md{I}=1,2$ and $k\ge 1$, but for $\md{I}\ge 3$ and $k\ge 1$ they need not preserve filtrations, but should instead map $F^\ga(\cdots)\ra F^{\ga+\de}(\cdots)$ for all~$\ga\in\Ga$.

We require that $\m_0^{I;x_0}\!=\!0$ if $\md{I}\!\le\! K$, and $\m_0^{I;x_0}(1)\!\in\! F^{\ssu 0}\Hom_{\sA^{i_M}}(f^{i_0i_M}(x_0),\ab f^{i_0i_M}(x_0))$ if $K<\md{I}\le L$, and $\m_0^{I;x_0}(1)\in F^\de\sH^{i_M}(f^{i_0i_M}(x_0),f^{i_0i_M}(x_0))$ if $\md{I}>L$ (this replaces $\m_0^{I;x_0}(1)\in F^0\sH^{i_M}(f^{i_0i_M}(x_0),f^{i_0i_M}(x_0))$), for all~$x_0\in\sO^{i_0}$. 

Then \eq{ai8eq11} converges in the filtered sense, by the argument in Definition \ref{ai8def2}. We also require that Definition \ref{ai8def4}(iv) holds. The rest of the definition of $\dAwCat^{\un,\iy}_{R,K,L}$ follows Definition \ref{ai8def6}. Theorem \ref{ai8thm6} will show $\dAwCat^{\un,\iy}_{R,K,L}$ is an $\iy$-category. Note that~$\dAfCat^{\un,\iy}_{R,K,L}\subset\dAwCat^{\un,\iy}_{R,K,L}$.
\end{dfn}

The next theorem, the analogue of Theorem \ref{ai7thm5}, will be proved in \S\ref{ai107} using equation~\eq{ai8eq11}.

\begin{thm}
\label{ai8thm4}
Let\/ $\bigl((\sO^i)_{i=0}^N,(f^{ij})_{0\le i\le j\le N},(\sH^i)_{i=0}^N,(\m^I)_{\es\ne I\subseteq\{0,\ldots,N\}}\bigr)$ be an $N$-simplex in $\dAfCat^{\un,\iy}_{R,K,L}$ or $\dAwCat^{\un,\iy}_{R,K,L}$ as in Definitions\/ {\rm\ref{ai8def6}} and\/ {\rm\ref{ai8def7},} and define maps\/ $\d_{n,k;y_0,\ldots,y_k}^{I;x_0,\ldots,x_n}$ from the $\m_k^{I;x_0,\ldots,x_k}$ as in\/ {\rm\eq{ai7eq6}--\eq{ai7eq9}}. Then for each nonempty subset\/ $I\subseteq\{0,\ldots,N\},$ written $I=\{i_0,\ab i_1,\ab\ldots,\ab i_M\}$ for $i_0<i_1<\cdots<i_M,$ and all\/ $k,n\ge 0,$ $x_0,\ldots,x_n\in\sO^{i_0}$ and\/ $y_0,\ldots,y_k\in\sO^{i_M},$ in $R$-linear maps as in \eq{ai7eq6} we~have
\ea
&\sum_{j=0}^M\sum_{l\ge 0, \; w_0,\ldots,w_l\in\sO^{i_j}}(-1)^j\d_{l,k;y_0,\ldots,y_k}^{\{i_j,\ldots,i_M\};w_0,\ldots,w_l}\ci\d_{n,l;w_0,\ldots,w_l}^{\{i_0,\ldots,i_j\};x_0,\ldots,x_n}
\nonumber\\
&=\sum_{j=1}^{M-1}(-1)^j\,\d_{n,k;y_0,\ldots,y_k}^{I\sm\{i_j\};x_0,\ldots,x_n},
\label{ai8eq12}
\ea
where the first line is an infinite sum which converges in the filtered sense.
\end{thm}

Here is the analogue of Examples \ref{ai7ex3} and~\ref{ai8ex1}.

\begin{ex}
\label{ai8ex3}
{\bf(a)} A 0-simplex $(\sO^0,f^{00}=\id_{\sO^0},\sH^0,\m^{\{0\}})$ of $\dAfCat^{\un,\iy}_{R,K,L}$ is equivalent to a unital filtered $A_\iy$-category $\sA^0$ over $R$, which is strict if $K>0$, with objects $\Obj\sA^0=\sO^0$, Hom spaces $\Hom_{\sA^0}(x,y)=\sH^0(x,y)$, morphisms $(\m_k^{\{0\};x_0,\ldots,x_k})_{k\ge 1}$, and units~$\boo^0_x\in \Hom_{\sA^0}^0(x,x)=\sH^0(x,x)^0$.
\smallskip

\noindent{\bf(b)} A 1-simplex of $\dAfCat^{\un,\iy}_{R,K,L}$ is equivalent to a unital filtered $A_\iy$-functor $\fF^{01}:\sA^0\ra\sA^1$, which is strict if $K>1$, with $\fF_\Obj^{01}=f^{01}$ and morphisms $\fF^{01;x_0,\ldots,x_k}_k=\m_k^{\{0,1\};x_0,\ldots,x_k}$. 
\smallskip

\noindent{\bf(c)} A 2-simplex of $\dAfCat^{\un,\iy}_{R,K,L}$ is equivalent to unital filtered $A_\iy$-categories $\sA^0,\ab\sA^1,\ab\sA^2$ over $R$ as in {\bf(a)}, which are strict if $K>0$, and unital filtered $A_\iy$-functors $\fF^{ij}:\sA^i\ra\sA^j$ for $(i,j)=(0,1),\ab(0,2),\ab(1,2)$ as in {\bf(b)}, which are strict if $K>1$, and $R$-linear maps for all $k\ge 0$ and $x_0,\ldots,x_k\in\sO^0$
\begin{align*}
\m_k^{\{0,1,2\};x_0,\ldots,x_k}:\sH^0&(x_0,x_1)\ot_R\sH^0(x_1,x_2)\ot_R\cdots \ot_R\sH^0(x_{k-1},x_k)
\nonumber\\
&\longra\sH^2(f^{02}(x_0),f^{02}(x_k))
\end{align*}
graded of degree $-k$, satisfying \eq{ai7eq19} and
\ea
&\sum_{\begin{subarray}{l} i,j,k\ge 0:\\ i+j+k=n \end{subarray}}\begin{aligned}[t]&
(-1)^{i+jk}\m^{\{0,1,2\};x_0,x_1,\ldots,x_i,x_{i+j},x_{i+j+1},\ldots,x_n}_{i+k+1}\ci\\
&\bigl(\id_{\sH^0(x_0,x_1)}\ot\cdots\ot \id_{\sH^0(x_{i-1},x_i)}\ot\m^{\{0\};x_i,x_{i+1},\ldots,x_{i+j}}_j\\
&\ot \id_{\sH^0(x_{i+j},x_{i+j+1})}\ot\cdots\ot \id_{\sH^0(x_{n-1},x_n)}\bigr)\end{aligned}
\label{ai8eq13}\\
&-(\m^{\{1,2\}}\!\ci\!\m^{\{0,1\}})^{x_0,\ldots,x_n}_n
\nonumber\\
&+\sum_{\!\!\!\!\!\!\begin{subarray}{l} 
1\le l\le k, \\ 
n=n_1+\cdots+n_k, \\
n_i\ge 0,\; i=1,\ldots,k. \\
\text{Set $p_j=\sum_{i=1}^jn_i$,} \\ 
0\le j\le k \end{subarray}}
\begin{aligned}[t]&(-1)^{k+1}\cdot(-1)^{\sum_{i=1}^{l-1}(n_i+1)}\cdot(-1)^{\sum_{i=1}^k(n_i+1)(k-i)}\\
&\m^{\{2\};f^{\{0,1,2\}}(x_{p_0}),f^{\{0,1,2\}}(x_{p_1}),\ldots,f^{\{0,1,2\}}(x_{p_{l-1}}),f^{02}(x_{p_l}),\ldots,f^{02}(x_{p_k})}_k\\
&\ci\bigl[(\m^{\{1,2\}}\!\ci\!\m^{\{0,1\}})^{x_{p_0},\ldots,x_{p_1}}_{n_1}\!\ot\!\cdots\!\ot\!(\m^{\{1,2\}}\!\ci\!\m^{\{0,1\}})^{x_{p_{l-2}},\ldots,x_{p_{l-1}}}_{n_{l-1}}\\
&\ot \m^{\{0,1,2\};x_{p_{l-1}},\ldots,x_{p_l}}_{n_l}\ot\m^{\{0,2\};x_{p_l},\ldots,x_{p_{l+1}}}_{n_{l+1}}\ot\\
&\cdots\ot\m^{\{0,2\};x_{p_{k-1}},\ldots,x_{p_k}}_{n_k}\bigr]
\end{aligned}
\nonumber\\
&=-\m^{\{0,2\};x_0,\ldots,x_n}_n,
\nonumber
\ea
with $\m_0^{\{0,1,2\};x_0}=0$ if $K>2$.

By Definition \ref{ai5def11} and equations \eq{ai5eq13} and \eq{ai8eq13}, $\m^{\{0,1,2\}}=(\m^{\{0,1,2\}}_k)_{k\ge 0}$ is a {\it unital filtered\/ $A_\iy$-homotopy\/} $\fF^{\{1,2\}}\ci \fF^{\{0,1\}}\Ra\fF^{\{0,2\}}$, which is strict if $K>2$. Hence, a 2-simplex in $\dAfCat^{\un,\iy}_{R,K,L}$ is a 2-commutative diagram of (strict if $K>0$) unital filtered $A_\iy$-categories, (strict if $K>1$) unital filtered $A_\iy$-functors, and a (strict if $K>2$) unital filtered $A_\iy$-homotopy
\begin{align*}
\xymatrix@C=70pt@R=13pt{
\drrtwocell_{}\omit^{}\omit{_{\,\,\,\,\,\,\,\,\,\,\,\,\,\,\,\,\,\,\m^{\{0,1,2\}}}} & \sA^1 \ar[dr]^{\fF^{\{1,2\}}}  \\
\sA^0 \ar[ur]^{\fF^{\{0,1\}}} \ar[rr]^(0.35){\fF^{\{0,2\}}} &&  \sA^2.}
\end{align*}
\end{ex}

The following analogue of Theorem \ref{ai7thm6} is proved in~\S\ref{ai108}.

\begin{thm}
\label{ai8thm5}
Let\/ $K,L\in\{0,1,2,\ldots,\iy\}$ with\/ $K\le L$ and\/ $L\ge 3$. Then there is a strict isomorphism of\/ $\iy$-categories
\e
I^{\un,\rfi}_{R,K,L}:\AfCat^{\un,\iy}_{R,K,L}\longra\dAfCat^{\un,\iy}_{R,K,L}
\label{ai8eq14}
\e
from Definitions\/ {\rm\ref{ai8def4}} and\/ {\rm\ref{ai8def6},} which acts on $N$-simplices by
\e
\begin{split}
I^{\un,\rfi}_{R,K,L}:&\bigl((\sO^i)_{i=0}^N,(f^{ij})_{0\le i\le j\le N},(\sH^i)_{i=0}^N,(\m^I)_{\es\ne I\subseteq\{0,\ldots,N\}}\bigr)\longmapsto\\
&\bigl((\sO^i)_{i=0}^N,(f^{ij})_{0\le i\le j\le N},(\sH^i)_{i=0}^N,(\dot\m^I)_{\es\ne I\subseteq\{0,\ldots,N\}}\bigr),
\end{split}
\label{ai8eq15}
\e
where $\dot\m^{I;x_0,\ldots,x_k}_k$ is defined from $\m^{I;x_0,\ldots,x_k}_k$ by
\e
\dot\m_k^{\{i_0,\ldots,i_M\};x_0,\ldots,x_k}=\begin{cases} \m_0^{\{i_0,i_1,i_2\};x_0}-\boo^{i_2}_{f^{i_0i_2}(x_0)}, & \text{$k=0$ and\/ $M=2,$}  \\
\m_k^{\{i_0,\ldots,i_M\};x_0,\ldots,x_k}, & \text{otherwise.} \end{cases}
\label{ai8eq16}
\e

Similarly, there is a strict isomorphism of\/ $\iy$-categories
\e
I^{\un,\wf}_{R,K,L}:\AwCat^{\un,\iy}_{R,K,L}\longra\dAwCat^{\un,\iy}_{R,K,L}
\label{ai8eq17}
\e
from Definitions\/ {\rm\ref{ai8def4}} and\/~{\rm\ref{ai8def7}}.
\end{thm}

Here is the analogue of Definition~\ref{ai7def10}.

\begin{dfn}
\label{ai8def8}
Let $(\Ga,\pre)$ be a po-monoid, $R$ a complete $\Ga$-filtered ring, and $K,L\in\{0,1,2,\ldots,\iy\}$ with $K\le L$ and $L\ge 3$. Generalizing Definition \ref{ai8def6}, define a simplicial set $\AfCat^\iy_{R,K,L}$ of (possibly non-unital) filtered $A_\iy$-categories over $R$ whose $N$-simplices $\bigl((\sO^i)_{i=0}^N,(f^{ij})_{0\le i\le j\le N},\ab(\sH^i)_{i=0}^N,\ab(\m^I)_{\es\ne I\subseteq\{0,\ldots,N\}}\bigr)$ for $N\ge 0$ satisfy Definition \ref{ai8def6}(a)(i)--(iv), but need not satisfy (a)(v). If $\al:[\bar N]\ra[N]$ is a morphism in $\De$, define $\AfCat^\iy_{R,K,L}(\al):(\AfCat^\iy_{R,K,L})_N\ab\ra(\AfCat^\iy_{R,K,L})_{\bar N}$ as in Definition \ref{ai8def6}(b). Essentially the same proof as Proposition \ref{ai7prop3} shows $\AfCat^\iy_{R,K,L}$ is a simplicial set.

Theorem \ref{ai8thm4} also holds in $\AfCat^\iy_{R,K,L}$, as the proof does not use unitality.

Theorem \ref{ai8thm6} will show that $\AfCat^\iy_{R,K,L}$ is an $\iy$-category. We say that $\AfCat^\iy_{R,K,L}$ has {\it type\/} $(K,L)$. We also write~$\AfCat^\iy_R=\AfCat^\iy_{R,0,3}$.

We also define $\iy$-categories $\AwCat^\iy_{R,K,L}$ as in \eq{ai8eq2} exactly as above,  but 
modifying the definition of $\dAwCat^{\un,\iy}_{R,K,L}$ in Definition \ref{ai8def7} rather than $\dAfCat^{\un,\iy}_{R,K,L}$ in Definition \ref{ai8def6}. As in Definition \ref{ai8def7}, we require that Definition \ref{ai8def4}(iv) holds.
\end{dfn}

The following generalization of Theorem \ref{ai7thm7} is proved in~\S\ref{ai1010}.

\begin{thm}
\label{ai8thm6}
The following simplicial sets defined in Definitions\/ {\rm \ref{ai8def4}, \ref{ai8def6}, \ref{ai8def7},} and\/ {\rm\ref{ai8def8},} for all\/ $K,L$ with\/ $L\ge 3,$ are all\/ $\iy$-categories:
\begin{gather*}
\AfCat^{\un,\iy}_{R,K,L},\dAfCat^{\un,\iy}_{R,K,L},\AfCat^\iy_{R,K,L},\\
\AwCat^{\un,\iy}_{R,K,L},\dAwCat^{\un,\iy}_{R,K,L},\AwCat^\iy_{R,K,L}.
\end{gather*}
\end{thm}

Here is the analogue of Definition \ref{ai7def11}. It completes our diagrams \eq{ai8eq1}--\eq{ai8eq2} of inclusions of $\iy$-categories of filtered $A_\iy$-categories.

\begin{dfn}
\label{ai8def9}
Suppose $(\Ga,\pre)$ is a po-monoid, and $R$ is a complete $\Ga$-filtered ring, and $K=0,1$ or $\iy$. Let $\AfCat^\un_{R,K}$ be the ordinary category of (possibly strict) unital filtered $A_\iy$-categories over $R$ from Definition \ref{ai5def9}, and $\sN(\AfCat^\un_{R,K})$ be its nerve from Example \ref{ai2ex5}, which is an $\iy$-category. 

Then $N$-simplices in $\sN(\AfCat^\un_{R,K})$ are of the form $\bigl((\sA^i)_{i=0}^N,(\fF^i)_{i=1}^N\bigr)$, where $\sA^i$ is a (strict if $K=1,\iy$) unital filtered $A_\iy$-category over $R$ for $i=0,\ldots,N$, and $\fF^i:\sA^{i-1}\ra\sA^i$ is a (strict if $K=\iy$) unital filtered $A_\iy$-functor for~$i=1,\ldots,N$.

Define an $\iy$-functor $\inc^{\un,\rfi}_{R,K}:\sN(\AfCat^\un_{R,K})\hookra\dAfCat^{\un,\iy}_{R,K,\iy}$ to act on $N$-simplices for $N\ge 0$ mapping
\begin{equation*}
\bigl((\sA^i)_{i=0}^N,(\fF^i)_{i=1}^N\bigr)\longmapsto\bigl((\sO^i)_{i=0}^N,(f^{ij})_{0\le i\le j\le N},(\sH^i)_{i=0}^N,(\m^I)_{\es\ne I\subseteq\{0,\ldots,N\}}\bigr),
\end{equation*}
where we set $\sO^i=\Obj\sA^i$, and writing $\fF^{ii}=\id_{\sA^i}:\sA^i\ra\sA^i$ and $\fF^{ij}=\fF^j\ci\fF^{j-1}\ci\cdots\ci\fF^{i+1}:\sA^i\ra\sA^j$ for $0\le i<j\le N$, we define $f^{ij}=\fF_\Obj^{ij}:\sO^i\ra\sO^j$ for $0\le i\le j\le N$. We set $\sH^i(x,y)=\Hom_{\sA^i}(x,y)$ for $x,y\in\sO^i=\Obj\sA^i$. We take $\m^{\{i\};x_0,\ldots,x_k}_k$ to be the map $\m^{x_0,\ldots,x_k}_k$ in $\sA^i$ for $i=0,\ldots,N$, and $\m^{\{i,j\};x_0,\ldots,x_k}_k$ to be the map $\fF^{ij;x_0,\ldots,x_k}_k$ in $\fF^{ij}$ for $0\le i<j\le N$, and $\m^{I;x_0,\ldots,x_k}_k=0$ if $\md{I}\ge 3$ and~$k\ge 0$.

We can prove as in Definition \ref{ai7def11} that $\inc^{\un,\rfi}_{R,K}$ is an $\iy$-functor. It is a strict isomorphism with the $\iy$-subcategory of $\dAfCat^{\un,\iy}_{R,K,\iy}$ consisting of $N$-simplices $\bigl((\sO^i)_{i=0}^N$, $(f^{ij})_{0\le i\le j\le N}$, $(\sH^i)_{i=0}^N$, $(\m^I)_{\es\ne I\subseteq\{0,\ldots,N\}}\bigr)$ with $\m_k^{I;x_0,\ldots,x_k}\!=\!0$ if $\md{I}\!\ge\!3$ and $k\!\ge\!0$. Similarly we define $\iy$-functors $\inc^\rfi_{R,K}:\sN(\AfCat_{R,K})\ab\hookra\AfCat^\iy_{R,K,\iy}$ for general filtered $A_\iy$-categories. 

We also define an $\iy$-functor $\iinc^{\un,\rfi}_{R,K,\iy}:\sN(\AfCat^\un_{R,K})\ra\AfCat^{\un,\iy}_{R,K,\iy}$ which is identified with $\inc^{\un,\rfi}_{R,K}$ under the isomorphism $I^{\un,\rfi}_{R,K,\iy}$ in \eq{ai8eq14}. The difference is that as in \eq{ai7eq23}, when we define $\bigl((\sO^i)_{i=0}^N,(f^{ij})_{0\le i\le j\le N},(\sH^i)_{i=0}^N,\ab(\m^I)_{\es\ne I\subseteq\{0,\ldots,N\}}\bigr)$ above we set $\m_0^{\{i_0,i_1,i_2\};x_0}=\boo^{i_2}_{f^{i_0i_2}(x_0)}$, not~$\m_0^{\{i_0,i_1,i_2\};x_0}=0$.
\end{dfn}

\begin{dfn}
\label{ai8def10}
Let $(\Ga,\pre)$ be a po-monoid and $R$ a complete $\Ga$-filtered ring. For all $K,L$, define $\iy$-subcategories
\e
\begin{aligned}
\AfAlg^{\un,\iy}_R&\subset\AfCat^{\un,\iy}_R, & \AfAlg^{\un,\iy}_{R,K,L}&\subset\AfCat^{\un,\iy}_{R,K,L},\\
\dAfAlg^{\un,\iy}_{R,K,L}&\subset\dAfCat^{\un,\iy}_{R,K,L}, & \AfAlg^\iy_R&\subset\AfCat^\iy_R, \\
\AfAlg^\iy_{R,K,L}&\subset\AfCat^\iy_{R,K,L}, & \AwAlg^{\un,\iy}_R&\subset\AwCat^{\un,\iy}_R, \\
\AwAlg^{\un,\iy}_{R,K,L}&\subset\AwCat^{\un,\iy}_{R,K,L}, & \dAwAlg^{\un,\iy}_{R,K,L}&\subset\dAwCat^{\un,\iy}_{R,K,L}, \\
\AwAlg^\iy_R&\subset\AwCat^\iy_R, & \AwAlg^\iy_{R,K,L}&\subset\AwCat^\iy_{R,K,L},
\end{aligned}
\label{ai8eq18}
\e
to be the full $\iy$-subcategories of filtered $A_\iy$-categories $\sA$ with $\Obj_\sA=\{*\}$. These are equivalent to {\it filtered\/ $A_\iy$-algebras\/} over $R$, so \eq{ai8eq18} are filtered analogues of the $\iy$-categories $\AAlg^{\un,\iy}_R,\ldots,\AAlg^\iy_{R,K}$ from~\S\ref{ai63}.

For $N$-simplices $\bigl((\sO^i)_{i=0}^N,(f^{ij})_{0\le i\le j\le N},(\sH^i)_{i=0}^N,(\m^I)_{\es\ne I\subseteq\{0,\ldots,N\}}\bigr)$ in these $\iy$-categories we have $\sO^i=\{*\}$ and $f^{ij}=\id_{\{*\}}$ for all $i,j$, and we use the notation $\sH^i(*,*)=A^i$ and $\m_k^{I;*,\ldots,*}=\m_k^I$, so that $N$-simplices may alternatively be written $\bigl((A^i)_{i=0}^N,(\m^I)_{\es\ne I\subseteq\{0,\ldots,N\}}\bigr)$, as in Chapter~\ref{ai6}.
\end{dfn}

\section{\texorpdfstring{$\iy$-category constructions for $A_\iy$-categories}{∞-category constructions for A∞-categories}}
\label{ai9}

\subsection{\texorpdfstring{Describing mapping spaces in $\ACat^{\un,\iy}_R$}{Describing mapping spaces in A∞-Algᵘⁿ∞ᵣ}}
\label{ai91}

The next theorem, proved in \S\ref{ai111}, generalizes equation \eq{ai4eq45} in Theorem \ref{ai4thm10} to our $\iy$-categories $\ACat^{\un,\iy}_R,\AfCat^{\un,\iy}_R,\AwCat^{\un,\iy}_R$. 

\begin{thm}
\label{ai9thm1}
{\bf(a)} Let\/ $\sA,\sB$ be unital\/ $A_\iy$-categories over a commutative ring $R$. Then there is a canonical strict isomorphism of Kan complexes
\e
\Hom_{\ACat^{\un,\iy}_R}^L(\sA,\sB)\cong \sN_{A_\iy}(\Fun_{A_\iy}^\un(\sA,\sB))^\Kan.
\label{ai9eq1}
\e
Here $\Hom_{\ACat^{\un,\iy}_R}^L(\sA,\sB)$ is the space of left morphisms in $\ACat^{\un,\iy}_R$ from Definition\/ {\rm\ref{ai2def8},} one model for the mapping space $\Hom_{\ACat^{\un,\iy}_R}(\sA,\sB)$ in the $\iy$-category $\ACat^{\un,\iy}_R,$ and\/ $\sN_{A_\iy}(\Fun_{A_\iy}^\un(\sA,\sB))$ is the simplicial nerve from Definition\/ {\rm\ref{ai4def20}} of the unital\/ $A_\iy$-category $\Fun_{A_\iy}^\un(\sA,\sB)$ from Definitions\/ {\rm\ref{ai4def7}--\ref{ai4def8},} and\/ $(\cdots)^\Kan$ denotes the maximal Kan complex in an $\iy$-category.

We prove\/ {\rm\eq{ai9eq1}} by first showing that
\e
\Hom_{\hACat^{\un,\iy}_R}^L(\sA,\sB)\cong \sN_{A_\iy}(\Fun_{A_\iy}^\un(\sA,\sB)).
\label{ai9eq2}
\e
Here we use the left morphisms $\Hom_\sC^L(X,Y)$ defined as in Definition\/ {\rm\ref{ai2def8},} but in the simplicial set\/ $\hACat^{\un,\iy}_R$ from Definition\/ {\rm\ref{ai7def3},} which is not an\/ $\iy$-category, so that its left morphisms need not be a Kan complex. As $(\hACat^{\un,\iy}_R,\cT)$ is an $\iy$-bicategory by Theorem\/ {\rm\ref{ai7thm4},} its left morphisms are an $\iy$-category.
\smallskip

\noindent{\bf(b)} Now let\/ $(\Ga,\pre)$ be a po-monoid and\/ $R$ a complete $\Ga$-filtered ring, as in\/ {\rm\S\ref{ai52},} and suppose\/ $\sA,\sB$ are unital filtered\/ $A_\iy$-categories over $R$. Then generalizing {\rm\eq{ai9eq1},} there are canonical strict isomorphisms of Kan complexes
\ea
\Hom_{\AfCat^{\un,\iy}_R}^L(\sA,\sB)&\cong \sN_{A_\iy}(F^0\Fun_{A^\rfi_\iy}^\un(\sA,\sB))^\Kan,
\label{ai9eq3}\\
\Hom_{\AwCat^{\un,\iy}_R}^L(\sA,\sB)&\cong \sN_{A_\iy}(\Fun_{A^\rfi_\iy}^\un(\sA,\sB))^\Kan,
\label{ai9eq4}
\ea
where $\Fun_{A^\rfi_\iy}^\un(\sA,\sB),F^0\Fun_{A^\rfi_\iy}(\sA,\sB)$ are as in Definition\/ {\rm\ref{ai5def10}(e),(f)}. Note that as in Remark\/ {\rm\ref{ai5rem5},} $\Fun_{A^\rfi_\iy}^\un(\sA,\sB)$ is an honest unital\/ $A_\iy$-category over $R$ (it is not curved, it has $\m_0^{\fF}=0$ for all objects $\fF$), and\/ $F^0\Fun_{A^\rfi_\iy}^\un(\sA,\sB)$ is a unital\/ $A_\iy$-cat\-eg\-ory over $F^0R,$ so $\sN_{A_\iy}(\Fun_{A^\rfi_\iy}^\un(\sA,\sB)),\sN_{A_\iy}(F^0\Fun_{A^\rfi_\iy}^\un(\sA,\sB))$ make sense. The analogues of\/ \eq{ai9eq2} also hold.
\end{thm}

The author views Theorem \ref{ai9thm1} as evidence that $\ACat^{\un,\iy}_R,\ab\AfCat^{\un,\iy}_R,\ab\AwCat^{\un,\iy}_R$ are the `right' definitions of $\iy$-categories generalizing $\ACat^\un_R,\ab\AfCat^\un_R$. In particular, \eq{ai9eq1}--\eq{ai9eq4} hold not just up to homotopy equivalence, but on the nose. Combining Definition \ref{ai7def14} and Theorem \ref{ai9thm1}(a), we can also describe mapping spaces in the $\iy$-category $\DGCat_R^\iy$ from Definition~\ref{ai7def14}:

\begin{cor}
\label{ai9cor1}
Let\/ $\sA,\sB$ be dg-categories over a commutative ring $R$. Then there is a canonical strict isomorphism of Kan complexes
\begin{equation*}
\Hom_{\DGCat_R^\iy}^L(\sA,\sB)\cong \sN_{A_\iy}(\Fun_{A_\iy}^\un(\sA,\sB)^\dg)^\Kan.
\end{equation*}
Here $\Fun_{A_\iy}^\un(\sA,\sB)^\dg$ is the full\/ $A_\iy$-subcategory of\/ $\Fun_{A_\iy}^\un(\sA,\sB)$ with objects dg-functors $\fF:\sA\ra\sB$.
\end{cor}

\subsection{\texorpdfstring{Comparing our $\iy$-categories with localizations}{Comparing our ∞-categories with localizations}}
\label{ai92}

Definition \ref{ai7def11} gives an $\iy$-functor $\iinc^\un_R:\sN(\ACat^\un_R)\ra\ACat^{\un,\iy}_R$, which maps unital homotopy equivalences in $\ACat^\un_R$ in Definition \ref{ai4def6} to equivalences in $\ACat^{\un,\iy}_R$ by Remark \ref{ai7rem6}. Thus by the universal property of $\iy$-category localizations in \S\ref{ai26}, there exists an $\iy$-functor $\Phi$ in the diagram
\e
\begin{gathered}
\xymatrix@C=80pt@R=15pt{
*+[r]{\sN(\ACat^\un_R)} \ar[rr]_(0.55){\Pi_\loc^\iy} \ar@/_.5pc/[drr]_{\iinc^\un_R} &  \drtwocell_{}\omit^{}\omit{^{}} & *+[l]{\ACat^\un_R[\HE_\un^{-1}]^\iy} \ar@{..>}[d]^{\Phi} \\
&& *+[l]{\ACat^{\un,\iy}_R,} }
\end{gathered}
\label{ai9eq5}
\e
unique up to homotopy, where $\HE_\un$ is the class of unital homotopy equivalences in $\ACat^\un_R$, in the sense of Definition \ref{ai4def14}. 

Theorem \ref{ai9thm4} will show that $\Phi$ is an equivalence of $\iy$-categories. To do this, we will show that $\Phi$ induces homotopy equivalences of mapping spaces. Mapping spaces in $\ACat^{\un,\iy}_R$ are given in Theorem \ref{ai9thm1}. To describe mapping spaces in $\ACat^\un_R[\HE_\un^{-1}]^\iy$, in the next definition and theorem we will construct {\it simplicial resolutions\/} of objects in $(\ACat^\un_R,\HE_\un)$ as in Definition \ref{ai2def18}, and then apply Theorem~\ref{ai2thm5}.

\begin{dfn}
\label{ai9def1}
Let $R$ be a commutative ring, and $\sB$ an object in $\ACat^\un_R$. We will define a functor $S_\sB:\De^\rop\ra\ACat^\un_R$ and an isomorphism $i_\sB:\sB\ra S_\sB([0])$. For each $N\ge 0$, Definition \ref{ai4def20} gives a unital $A_\iy$-category $\sC^\De_{[N]}$ over $R$. This has objects $0,\ldots,N$, and for all $0\le i\le j\le N$ there is a distinguished closed degree 0 morphism $\om_{ij}:i\ra j$ in $\sC^\De_{[N]}$, with~$\om_{ii}=\boo_i$.

Definitions \ref{ai4def7}--\ref{ai4def8} define a unital $A_\iy$-category $\Fun_{A_\iy}^\un(\sC^\De_{[N]},\sB)$ whose objects are unital $A_\iy$-functors $\fF:\sC^\De_{[N]}\ra\sB$. For the object $[N]$ in $\De^\rop$, define $S_\sB([N])\in\ACat^\un_R$ to be the full $A_\iy$-subcategory of $\Fun_{A_\iy}^\un(\sC^\De_{[N]},\sB)$ with objects those $\fF:\sC^\De_{[N]}\ra\sB$ such that if $0\le i<j\le N$ then $\fF(\om_{ij})$ is a homotopy equivalence in $\sB$, in the sense of Definition~\ref{ai4def5}.

If $\al:[N]\ra[M]$ is a morphism in $\De$ then Definition \ref{ai4def20} gives a unital $A_\iy$-functor $\sC^\De_\al:\sC^\De_{[N]}\ra \sC^\De_{[M]}$. By Definition \ref{ai4def9}(a),(e) this induces a unital $A_\iy$-functor $(\sC^\De_\al)^*:\Fun_{A_\iy}^\un(\sC^\De_{[M]},\sB)\ra \Fun_{A_\iy}^\un(\sC^\De_{[N]},\sB)$. Since $\sC^\De_\al(\om_{ij})=\om_{\al(i),\al(j)}$, this $(\sC^\De_\al)^*$ restricts to a unital $A_\iy$-functor $S_\sB([M])\ra S_\sB([N])$. Define $S_\sB(\al)=(\sC^\De_\al)^*\vert_{S_\sB([M])}:S_\sB([M])\ra S_\sB([N])$. The functoriality of the $\sC^\De_\al$ in Definition \ref{ai4def20} and of pullbacks imply that $S_\sB:\De^\rop\ra\ACat^\un_R$ is a functor.

Now $\sC^\De_{[0]}$ is isomorphic to the unital point $A_\iy$-category $\sP^\un$ in Definition \ref{ai4def11}, so Proposition \ref{ai4prop2} gives a canonical isomorphism $\sB\ra\Fun_{A_\iy}^\un(\sC^\De_{[0]},\sB)$. As $S_\sB([0])=\Fun_{A_\iy}^\un(\sC^\De_{[0]},\sB)$, we define $i_\sB:\sB\ra S_\sB([0])$ to be this.

Now suppose that $\sB$ lies in $\DGCat_R\subset\ACat^\un_R$. As in Definition \ref{ai4def20}, the $A_\iy$-categories $\sC^\De_{[N]}$ are actually dg-categories. Hence $\Fun_{A_\iy}^\un(\sC^\De_{[N]},\sB)$ is a dg-category, as in Remark \ref{ai4rem12}. It is immediate from Definition \ref{ai4def9}(a) that the $A_\iy$-functors $(\sC^\De_\al)^*$ above are dg-functors. Thus $S_\sB$ is a functor $S_\sB:\De^\rop\ra\DGCat_R\subset\ACat^\un_R$.\end{dfn}

The next theorem will be proved in \S\ref{ai112}. Note that (a),(b) refer to two different, though closely related, notions of simplicial resolution.

\begin{thm}
\label{ai9thm2}
Work in the situation of Definition\/~{\rm\ref{ai9def1}}.
\smallskip

\noindent{\bf(a)} $(S_\sB,i_\sB)$ is a \begin{bfseries}simplicial resolution\end{bfseries} of\/ $\sB$ in\/ $(\ACat^\un_R,\HE_\un),$ in the sense of Definition\/~{\rm\ref{ai2def18}}.
\smallskip

\noindent{\bf(b)} Suppose $\sB$ lies in $\DGCat_R\subset\ACat^\un_R,$ and is fibrant in the Tabuada model structure on $\DGCat_R$ from Theorem\/ {\rm\ref{ai4thm1}}. Then $(S_\sB,i_\sB)$ is a \begin{bfseries}simplicial resolution\end{bfseries} of\/ $\sB$ in the model category\/ $\DGCat_R,$ in the sense of Remark\/ {\rm\ref{ai2rem4}(c)}. This simplicial resolution was constructed by Holstein\/~{\rm\cite[\S 3]{Hols}}.

\smallskip

\noindent{\bf(c)} Define $F^\un_{S_\sB}:(\ACat^\un_R)^\rop\ra\SSet$ as in Definition\/ {\rm\ref{ai2def18}}. Then for each unital\/ $A_\iy$-category $\sA$ over $R$ there is a canonical strict isomorphism of simplicial sets
\e
F^\un_{S_\sB}(\sA)\cong\sN_{A_\iy}(\Fun_{A_\iy}^\un(\sA,\sB))^\Kan.
\label{ai9eq6}
\e

If\/ $\fF:\sA_1\ra\sA_2$ is a morphism in $\ACat^\un_R$ then the isomorphism \eq{ai9eq6} identifies $F^\un_{S_\sB}(\fF):F^\un_{S_\sB}(\sA_2)\ra F^\un_{S_\sB}(\sA_1)$ with
\e
\!\!\sN_{A_\iy}(\fF^*)^\Kan:\sN_{A_\iy}(\Fun_{A_\iy}^\un(\sA_2,\sB))^\Kan\!\longra\!\sN_{A_\iy}(\Fun_{A_\iy}^\un(\sA_1,\sB))^\Kan,\!\!
\label{ai9eq7}
\e
where $\fF^*:\Fun_{A_\iy}^\un(\sA_2,\sB)\ra\Fun_{A_\iy}^\un(\sA_1,\sB)$ is as in Definition\/~{\rm\ref{ai4def9}(a),(e)}. 
\smallskip

\noindent{\bf(d)} Combining {\bf(a)\rm,\bf(c)} and Theorem\/ {\rm\ref{ai2thm5}} gives a homotopy equivalence for all unital\/ $A_\iy$-categories $\sA,\sB$ over $R$
\e
\Hom_{\ACat^\un_R[\HE_\un^{-1}]^\iy}^L(\sA,\sB)\simeq \sN_{A_\iy}(\Fun_{A_\iy}^\un(\sA,\sB))^\Kan.
\label{ai9eq8}
\e

\noindent{\bf(e)} Suppose\/ $\cC$ is any of the full subcategories of\/ $\ACat^\un_R$
\e
\begin{aligned}
&\DGCat_R^{A_\iy}, & &\DGCat_R^{A_\iy,\rsu}, & &\DGCat_R^{A_\iy,\hp}, & &\DGCat_R^{A_\iy,\hs}, \\
&\ACat_R^{\un,\rsu}, & &\ACat_R^{\un,\hp}, & &\ACat_R^{\un,\hs}, & &\ADCat_R^\un,\\
&\ADCat_R^{\un,\rsu}, & &\ADCat_R^{\un,\hp}, & &\ADCat_R^{\un,\hs}, \end{aligned}
\label{ai9eq9}
\e
from Definitions\/ {\rm \ref{ai4def17}, \ref{ai4def19},} and\/ {\rm\ref{ai4def21},} and let\/ $\sB$ be an object of\/ $\cC$. Then $S_\sB$ maps to $\cC\subset\ACat^\un_R,$ so that\/ $(S_\sB,i_\sB)$ is also a simplicial resolution of\/ $\sB$ in\/ $(\cC,\HE_\un\cap\cC),$ in the sense of Definition\/ {\rm\ref{ai2def18}}. Hence the computation of mapping spaces \eq{ai9eq8} also holds in the localized\/ $\iy$-category $\cC[(\HE_\un\cap\cC)^{-1}]^\iy$. This implies that the natural morphism $\cC[(\HE_\un\cap\cC)^{-1}]^\iy\ra \ACat^\un_R[\HE_\un^{-1}]^\iy$ is fully faithful.
\end{thm}

\begin{rem}
\label{ai9rem2}
{\bf(a)} The simplicial resolution $(S_\sB,i_\sB)$ in Definition \ref{ai9def1} is {\it functorial in\/} $\sB$ over $\ACat^\un_R$. We can define a functor $S:\De^\rop\t\ACat^\un_R\ra \ACat^\un_R$ and a natural isomorphism $i:\Id_{\ACat^\un_R}\Ra S\vert_{[0]\t\ACat^\un_R}$ with $S\vert_{\De^\rop\t\sB}=S_\sB$ and $i\vert_\sB=i_\sB$ for all $\sB$. But we will not need this.
\smallskip

\noindent{\bf(b)} Although it is a `classical' result, our proof of Theorem \ref{ai9thm2}(a) uses Theorems \ref{ai7thm3} and \ref{ai9thm1}(a) on our new $\iy$-category $\ACat^{\un,\iy}_R$.
\end{rem}

The next theorem will be proved in \S\ref{ai113}. It proves equation \eq{ai4eq47} in Remark~\ref{ai4rem18}.

\begin{thm}
\label{ai9thm3}
Let\/ $R$ be a commutative ring. Then inclusions of subcategories between the categories $\DGCat_R^\hs,\ab\DGCat_R,\ab\DGCat^{A_\iy,\hs}_R,\ab\DGCat^{A_\iy}_R,\ab\ACat^{\un,\hs}_R,\ab\ACat_R^\un$ from {\rm\S\ref{ai4},} and between the corresponding classes of homotopy equivalences $\HE_{\cdots}$ or quasi-equivalences $\QE_{\cdots},$ induce a homotopy commutative diagram of\/ $\iy$-functors of localizations:
\e
\begin{gathered}
\xymatrix@C=150pt@R=15pt{
*+[r]{\DGCat_R^\hs[\HE_{\dg,\hs}^{-1}]^\iy}
\ar[r]^\simeq_{\eq{ai4eq46}} \ar@{.>}[d]
& *+[l]{\DGCat_R[\QE_\dg^{-1}]^\iy} \ar[d]^\simeq_{\eq{ai4eq44}}
\\
*+[r]{\DGCat^{A_\iy,\hs}_R[\HE_{\dg,A_\iy,\hs}^{-1}]^\iy} \ar@{.>}[r] \ar@{.>}[d]
& *+[l]{\DGCat^{A_\iy}_R[\QE_{\dg,A_\iy}^{-1}]^\iy} \ar[d]^\simeq_{\eq{ai4eq44}}
\\
*+[r]{\ACat^{\un,\hs}_R[\HE_{\un,\hs}^{-1}]^\iy} \ar@{.>}[r] &
*+[l]{\ACat_R^\un[\QE_\un^{-1}]^\iy.}
}
\end{gathered}
\label{ai9eq10}
\e
The top and right morphisms in the diagram, written `$\,\ra$', are Joyal equivalences by Theorems\/ {\rm\ref{ai4thm8}} and\/ {\rm\ref{ai4thm11},} as in equations \eq{ai4eq44} and\/~\eq{ai4eq46}. 

The other morphisms in {\rm\eq{ai9eq10},} written `$\,\dashra$', are also Joyal equivalences.
\end{thm}

The next theorem will be proved in \S\ref{ai114}.

\begin{thm}
\label{ai9thm4}
Let\/ $R$ be a commutative ring. Then the $\iy$-functor from equation \eq{ai9eq5} is a Joyal equivalence:
\e
\Phi:\ACat^\un_R[\HE_\un^{-1}]^\iy\longra\ACat^{\un,\iy}_R.
\label{ai9eq11}
\e

The analogue holds for all the full subcategories of\/ $\ACat^\un_R$ in {\rm\eq{ai9eq9},} giving  Joyal equivalences of their localizations with the full\/ $\iy$-subcategories of\/ $\ACat^{\un,\iy}_R$ with the same objects. So for example we have equivalences
\e
\begin{aligned}
\ACat^{\un,\hs}_R[\HE_{\un,\hs}^{-1}]^\iy&\simeq \ACat^{\un,\hs,\iy}_R,\\
\ACat^{\un,\rsu}_R[\HE_{\un,\rsu}^{-1}]^\iy&\simeq \ACat^{\un,\rsu,\iy}_R,\\
\ADCat^\un_R[\HE_{\un,\De}^{-1}]^\iy&\simeq\ADCat^{\un,\iy}_R.
\end{aligned}
\label{ai9eq12}
\e 
\end{thm}

Combining Theorems \ref{ai4thm8}, \ref{ai4thm11}, \ref{ai7thm8}, \ref{ai9thm3}, \ref{ai9thm4} and Corollary \ref{ai4cor6}, we deduce:

\begin{cor}
\label{ai9cor2}
{\bf(a)} Let\/ $R$ be a commutative ring. Then all the following $\iy$-categories are Joyal equivalent:
\ea
&\ACat^{\un,\hs,\iy}_R, & &\DGCat^{A_\iy,\hs,\iy}_R, & &\ACat^\un_R[\QE_\un^{-1}]^\iy, 
\nonumber \\
&\ACat_R^\cu[\QE_\cu^{-1}]^\iy, &
&\ACat_R^\wc[\QE_\wc^{-1}]^\iy, & 
&\DGCat_R[\QE_\dg^{-1}]^\iy, 
\nonumber \\
&\DGCat^{A_\iy}_R[\QE_{\dg,A_\iy}^{-1}]^\iy, &
&\ACat_R^{\un,\hs}[\HE_{\un,\hs}^{-1}]^\iy,\!\!\!\!
\nonumber \\
&\DGCat_R^\hs[\HE_{\dg,\hs}^{-1}]^\iy,\!\!\! & &\DGCat_R^{A_\iy,\hs}[\HE_{A_\iy,\dg,\hs}^{-1}]^\iy.\!\!\!\!\!\!\!\!\!\!\!\!\!\!\!\!\!
\label{ai9eq13}
\ea
They satisfy the conclusions of Corollary\/ {\rm\ref{ai4cor6},} so they are presentable and closed symmetric monoidal. 
\smallskip

\noindent{\bf(b)} The following $\iy$-categories are Joyal equivalent:
\e
\begin{aligned}
&\DGCat^{A_\iy,\rsu,\iy}_R, & &\ACat^{\un,\rsu,\iy}_R, \\
&\DGCat^{A_\iy,\rsu}_R[\HE_{\dg,A_\iy,\rsu}^{-1}]^\iy, & &\ACat^{\un,\rsu}_R[\HE_{\un,\rsu}^{-1}]^\iy.
\end{aligned}
\label{ai9eq14}
\e
\end{cor}

\begin{rem}
\label{ai9rem3}
The $\iy$-categories $\ACat^\un_R[\QE_\un^{-1}]^\iy$ and $\ACat_R^\cu[\QE_\cu^{-1}]^\iy$ in \eq{ai9eq13} are considered in the literature \cite{COS1,COS2,Faon2,OhTa,Orna1,Tana,Pasc} to be `the' $\iy$-category of $A_\iy$-categories over $R$, as in \S\ref{ai49}. We have shown that these can be identified with the full $\iy$-subcategory $\ACat^{\un,\hs,\iy}_R$ of our explicit $\iy$-category $\ACat^{\un,\iy}_R$ from Chapter~\ref{ai7}.
\end{rem}

\begin{rem}
\label{ai9rem4}
We did not treat {\it filtered\/} $A_\iy$-categories from Chapter \ref{ai8} in the material above. This is because the method of Definition \ref{ai9def1} does not work to define simplicial resolutions in the filtered case; to define the isomorphism $i_\sB:\sB\ra S_\sB([0])$ we used that $\sC^\De_{[0]}\cong\sP^\un$ and $\Fun_{A_\iy}^\un(\sP^\un,\sB)\cong\sB$, but as in Proposition \ref{ai5prop2}, in the filtered case we have $\Fun_{A^\rfi_\iy}^\un(\sP^\un,\sB)\cong\Str(\sB)$. The author experimented with other methods to define simplicial resolutions in $(\AfCat^\un_R,\HE_{\un,\rfi})$, but did not succeed. It seems plausible that
\begin{equation*}
\AfCat^\un_R[\HE_{\un,\rfi}^{-1}]^\iy\simeq\AfCat^{\un,\iy}_R,
\end{equation*}
but the author does not have a proof of this.
\end{rem}

\subsection{\texorpdfstring{The $\iy$-categorical simplicial nerve functor}{The ∞-categorical simplicial nerve functor}}
\label{ai93}

This section will construct an $\iy$-functor $\sN^\iy_{A_\iy}:\ACat^{\un,\iy}_R\ra\sCat_\iy$ which generalizes the simplicial nerve functor $\sN_{A_\iy}:\ACat^\un_R\ra\sCat_\iy$ from \S\ref{ai48}. Recall that Definition \ref{ai4def11} defined the unital point $A_\iy$-category $\sP^\un$, and Proposition \ref{ai4prop2} showed that $\sA\cong\Fun^\un_{A_\iy}(\sP^\un,\sA)$ for any unital $A_\iy$-category $\sA$. Combining \eq{ai9eq2} and Proposition \ref{ai4prop2}, we see that if we could define an $\iy$-functor $\sN^\iy_{A_\iy}:\ACat^{\un,\iy}_R\ra\sCat_\iy$ mapping $\sA\mapsto \Hom_{\hACat^{\un,\iy}_R}^L(\sP^\un,\sA)$, then on objects it would map $\sA\mapsto \sN_{A_\iy}(\Fun_{A_\iy}^\un(\sP^\un,\sA))\cong\sN_{A_\iy}(\sA)$, as we want. We do this using the `straightening construction' of Theorem~\ref{ai2thm2}.

\begin{dfn}
\label{ai9def2}
Let the simplicial sets $\ACat^{\un,\iy}_R\subset\hACat^{\un,\iy}_R$ be as in Definitions \ref{ai7def3} and \ref{ai7def4}, for $R$ a commutative ring. Then $\sP^\un$ in Definition \ref{ai4def11} is a 0-simplex in $\ACat^{\un,\iy}_R$. Define a simplicial set $\NCat^{\un,\iy}_R$ and a map of simplicial sets $p:\NCat^{\un,\iy}_R\ra \ACat^{\un,\iy}_R$, such that $N$-simplices of $\NCat^{\un,\iy}_R$ are $(N+1)$-simplices $f:\De_{N+1}\ra\hACat^{\un,\iy}_R$ such that $f(\De^{\{0\}})=\sP^\un$ and $f(\De^{\{1,2,\ldots,N+1\}})$ is an $N$-simplex in $\ACat^{\un,\iy}_R\subset\hACat^{\un,\iy}_R$, and $p$ maps $f\mapsto f(\De^{\{1,2,\ldots,N+1\}})$. Maps between $N$-simplices corresponding to morphisms $\al:[\bar N]\ra[N]$ in $\De$, to make $\NCat^{\un,\iy}_R$ into a simplicial set, are defined as for $\Hom_\sC^L(x,y)$ in Definition~\ref{ai2def8}.

This is an example of a well known construction: as in \cite[\S 1.2.9]{Luri2} we have a simplicial set $(\hACat^{\un,\iy}_R)_{\sP^\un/}$ with a morphism $\Pi:(\hACat^{\un,\iy}_R)_{\sP^\un/}\ra \hACat^{\un,\iy}_R$, the {\it undercategory\/} of $\sP^\un$ in $\hACat^{\un,\iy}_R$, and $p$ is the restriction of $\Pi$ to the preimage of $\ACat^{\un,\iy}_R\subset\hACat^{\un,\iy}_R$.
\end{dfn}

The next theorem will be proved in \S\ref{ai115}.

\begin{thm}
\label{ai9thm5}
{\bf(a)} In Definition\/ {\rm\ref{ai9def2},} $p:\NCat^{\un,\iy}_R\ra \ACat^{\un,\iy}_R$ is a \begin{bfseries}cocartesian fibration\end{bfseries} in the sense of Definition\/ {\rm\ref{ai2def12}}. Therefore Theorem\/ {\rm\ref{ai2thm2}(b)} defines an $\iy$-functor $\sN^\iy_{A_\iy}:\ACat^{\un,\iy}_R\ra\sCat_\iy,$ the \begin{bfseries}straightening\end{bfseries} of\/ $p,$ uniquely up to a contractible choice. We call $\sN^\iy_{A_\iy}$ the \begin{bfseries}$\iy$-categorical simplicial nerve functor\end{bfseries}, as it extends $\sN_{A_\iy}$ in {\rm\S\ref{ai48}} to $\iy$-categories.
\smallskip

\noindent{\bf(b)} Let\/ $\sA$ be an object of\/ $\ACat^{\un,\iy}_R$. Then there is a strict isomorphism of $\iy$-categories $p^{-1}(\sA)\cong\sN_{A_\iy}(\sA),$ where $\sN_{A_\iy}(\sA)$ is as in Definition\/ {\rm\ref{ai4def20}}. Hence by Theorem\/ {\rm\ref{ai2thm2}(a),} $\sN^\iy_{A_\iy}(\sA)$ is Joyal equivalent to~$\sN_{A_\iy}(\sA)$.
\smallskip

\noindent{\bf(c)} The following diagram of\/ $\iy$-categories commutes up to homotopy:
\e
\begin{gathered}
\xymatrix@C=180pt@R=15pt{
*+[r]{\sN(\ACat^\un_R)} \ar[r]_(0.36){\sN(\sN_{A_\iy})} \ar[d]^{\iinc^\un_R} \drtwocell_{}\omit^{}\omit{^{}} & *+[l]{\sN(\Cat^1_\iy)} \ar[d]_{I^{\sCat_\iy}} \\
*+[r]{\ACat^{\un,\iy}_R } \ar[r]^(0.66){\sN^\iy_{A_\iy}} & *+[l]{\sCat_\iy,\!}}
\end{gathered}
\label{ai9eq15}
\e
using the notation of Definitions\/ {\rm\ref{ai2def10}, \ref{ai4def20}} and\/ {\rm\ref{ai7def11}}. We prove this by showing that the pullback of\/ $p:\NCat^{\un,\iy}_R\ra\ACat^{\un,\iy}_R$ to $\sN(\ACat^\un_R)$ is strictly isomorphic to the relative simplicial nerve\/ $\pi:N_{\sN_{A_\iy}}(\ACat^\un_R)\ra\sN(\ACat^\un_R)$ from Example\/ {\rm\ref{ai2ex6},} which is an explicit model for the \begin{bfseries}unstraightening\end{bfseries} of\/ $I^{\sCat_\iy}\ci\sN(\sN_{A_\iy})$ in Theorem\/~{\rm\ref{ai2thm2}(c)}.
\end{thm}

Combining Theorems \ref{ai4thm7} and \ref{ai9thm5} yields:

\begin{cor}
\label{ai9cor3}
In Theorem\/ {\rm\ref{ai9thm5},} $\sN^\iy_{A_\iy}$ restricts to
\e
\sN^\iy_{A_\iy}:\ADCat^{\un,\iy}_R\longra\sCat^\rst_\iy
\label{ai9eq16}
\e
from the $\iy$-category $\ADCat^{\un,\iy}_R$ of \begin{bfseries}pretriangulated\end{bfseries} unital\/ $A_\iy$-cat\-eg\-or\-ies over $R$ in Definition\/ {\rm\ref{ai7def13}} to the $\iy$-category $\sCat^\rst_\iy$ of \begin{bfseries}stable\end{bfseries} $\iy$-categories in Definition\/~{\rm\ref{ai2def19}}.
\end{cor}

\begin{rem}
\label{ai9rem5}
{\bf(a)} As in Remark \ref{ai4rem17}, if $(X,\om)$ is a suitable symplectic manifold, the {\it derived Fukaya category\/} is $D^b\sF(X,\om)=H^0(\Tw(\sF(X,\om)))$, a triangulated category. Kontsevich's Homological Mirror Symmetry Conjecture \cite{Kont} explains Mirror Symmetry of Calabi--Yau manifolds $X,\check X$ in terms of equivalences of triangulated categories $D^b\sF(X,\om)\simeq D^b\coh\check X$.

As in \S\ref{ai27}, stable $\iy$-categories are the $\iy$-categorical version of triangulated categories, and are better-behaved than triangulated categories in some ways.

Corollary \ref{ai9cor3} is important as it will allow construction of derived Fukaya categories as stable $\iy$-categories, in an $\iy$-functorial way.
\smallskip

\noindent{\bf(b)} We can interpret $\sN^\iy_{A_\iy}$ in Theorem \ref{ai9thm5} in terms of the $\iy$-bicategory $(\hACat^{\un,\iy}_R,\cT)$ in Theorem \ref{ai7thm4}, as mapping $\sA\!\mapsto\!\Hom_{\hACat^{\un,\iy}_R}(\sP^\un,\sA)$. Probably Theorem \ref{ai9thm5}(a) could be proved using the theory of $\iy$-bicategories, but we give a direct proof.
\smallskip

\noindent{\bf(c)} An alternative, less explicit way to construct the $\iy$-functor $\sN^\iy_{A_\iy}$ in Theorem \ref{ai9thm5} is via the homotopy commutative diagram:
\begin{equation*}
\xymatrix@C=145pt@R=15pt{
*+[r]{\sN(\ACat^\un_R)} \ar[r]_(0.36){\sN(\sN_{A_\iy})} \ar[d]^{\Pi_\loc^\iy} \ar@<-1.4ex>@/_.6pc/[dd]_{\iinc^\un_R} \drtwocell_{}\omit^{}\omit{^{}} & *+[l]{\sN(\Cat^1_\iy)} \ar[d]_{\Pi_\loc^\iy} \ar@<1.4ex>@/^.6pc/[dd]^{I^{\sCat_\iy}} \\
*+[r]{\ACat^\un_R[\HE_\un^{-1}]^\iy} \ar@{..>}[r]^(0.6){\sN(\sN_{A_\iy})^\iy_\loc} \ar[d]^{\eq{ai9eq11}}_(0.4)\simeq \drtwocell_{}\omit^{}\omit{^{}} & *+[l]{\Cat^1_\iy[{\rm Joy}^{-1}]^\iy} \ar[d]_{\eq{ai2eq15}}^(0.3)\simeq \\
*+[r]{\ACat^{\un,\iy}_R } \ar@{..>}[r]^(0.66){\sN^\iy_{A_\iy}} & *+[l]{\sCat_\iy.\!}}
\end{equation*}
Here we first show that $\sN_{A_\iy}$ takes homotopy equivalences to Joyal equivalences, then the middle horizontal morphism exists by properties of localizations, and the bottom horizontal morphism exists as the bottom vertical morphisms are equivalences. Properties of localizations then imply that $\sN^\iy_{A_\iy}$ is characterized uniquely up to homotopy by the fact that \eq{ai9eq15} commutes.
\smallskip

\noindent{\bf(d)} It would be interesting to strengthen Theorem \ref{ai9thm5} to make $\sN^\iy_{A_\iy}$ map to $R$-{\it linear\/} $\iy$-categories, in a suitable sense. For the restriction \eq{ai9eq16} to $\ADCat^{\un,\iy}_R$, these could be tensored over the stable $\iy$-category~$D^b\text{mod-}R$.
\end{rem}

\subsection{\texorpdfstring{An `Inner Hom functor' for $\ACat^{\un,\iy}_R$}{An 'Inner Hom functor' for A∞-Catᵘⁿ∞ᵣ}}
\label{ai94}

\begin{thm}
\label{ai9thm6}
Let\/ $R$ be a commutative ring. Then there exist\/ $\iy$-functors $\tFun_{A_\iy}^{\un,\iy},\Fun_{A_\iy}^{\un,\iy},$ unique up to homotopy, such that the following diagrams of\/ $\iy$-categories commute up to homotopy:
\ea
\begin{gathered}
\xymatrix@!0@C=280pt@R=40pt{
*+[r]{\sN(\ACat^\un_R)^\rop\t\sN(\ACat^\un_R)} \ar[r]_(0.66){\sN(\Fun_{A_\iy}^\un)} \ar[d]^(0.4){(\Pi_\loc^\iy)^\rop\t\Pi_\loc^\iy} \drtwocell_{}\omit^{}\omit{^{}} & *+[l]{\sN(\ACat^\un_R)} \ar[d]_{\Pi_\loc^\iy} \\
*+[r]{\;\>\begin{subarray}{l} \ts (\ACat^\un_R[\HE_\un^{-1}]^\iy)^\rop \\ \ts {}\t\ACat^\un_R[\HE_\un^{-1}]^\iy\end{subarray}} \ar@{..>}[r]^(0.61){\tFun_{A_\iy}^{\un,\iy}} & *+[l]{\ACat^\un_R[\HE_\un^{-1}]^\iy,\!}}
\end{gathered}&
\label{ai9eq17}
\allowdisplaybreaks\\
\begin{gathered}
\xymatrix@!0@C=280pt@R=35pt{
*+[r]{\sN(\ACat^\un_R)^\rop\t\sN(\ACat^\un_R)} \ar[r]_(0.66){\sN(\Fun_{A_\iy}^\un)} \ar[d]^{(\iinc^\un_R)^\rop\t\iinc^\un_R} \drtwocell_{}\omit^{}\omit{^{}} & *+[l]{\sN(\ACat^\un_R)} \ar[d]_{\iinc^\un_R} \\
*+[r]{(\ACat^{\un,\iy}_R)^\rop\t\ACat^{\un,\iy}_R} \ar@{..>}[r]^(0.66){\Fun_{A_\iy}^{\un,\iy}} & *+[l]{\ACat^{\un,\iy}_R,\!}}
\end{gathered}&
\label{ai9eq18}
\ea
Here the ordinary functor $\Fun_{A_\iy}^\un$ is as in {\rm\eq{ai4eq24},} and\/ $\iinc^\un_R$ is as in Definition\/~{\rm\ref{ai7def11}}.
\end{thm}

\begin{proof}
Corollary \ref{ai4cor1} implies that $\Fun_{A_\iy}^\un$ maps $(\HE_\un)^\rop\t\HE_\un\ra\HE_\un$. Thus by the universal property of localizations of $\iy$-categories, there exists an $\iy$-functor $\tFun_{A_\iy}^\un$, unique up to homotopy, making \eq{ai9eq17} commute up to homotopy. Equation \eq{ai9eq18} then follows from \eq{ai9eq17} by the equivalence \eq{ai9eq11}, which identifies $\Pi_\loc^\iy$ with $\iinc^\un_R$ up to homotopy.
\end{proof}

\begin{rem}
\label{ai9rem6}
{\bf(a)} The $\iy$-functor $\Fun_{A_\iy}^{\un,\iy}$ is an analogue of the internal mapping objects $\cHom(\sC,\sD)$ in $\DGCat_R[\QE_\dg^{-1}]^\iy,\ACat^\un_R[\QE_\un^{-1}]^\iy$ discussed in Theorems \ref{ai4thm2}, \ref{ai4thm9} and \ref{ai4thm10}, and should coincide with them on restriction to the $\iy$-subcategory $\ACat^{\un,\hs,\iy}_R$ which is equivalent to $\DGCat_R[\QE_\dg^{-1}]^\iy,\ab\ACat^\un_R[\QE_\un^{-1}]^\iy$ by Corollary~\ref{ai9cor2}(a). 

However, our $\iy$-functors are both more general -- they are defined on larger $\iy$-categories -- and are easier to construct. One reason is that we are inverting homotopy equivalences, not quasi-equivalences, so we can apply Corollary~\ref{ai4cor1}.

\smallskip

\noindent{\bf(b)} It is an obvious question whether we can define a symmetric monoidal structure $\ot^\iy$ on $\ACat^{\un,\iy}_R$ for which $\Fun_{A_\iy}^{\un,\iy}$ is the associated internal mapping object, as for $\DGCat_R[\QE_\dg^{-1}]^\iy$ in Theorem \ref{ai4thm2}. This is possible for the $\iy$-subcategory $\ACat^{\un,\hs,\iy}_R$ by Theorem \ref{ai4thm2} and Corollary \ref{ai9cor2}(a). But the author does not know whether $\ot^\iy$ can be extended to $\ACat^{\un,\iy}_R$. See Saneblidze--Umble \cite{SaUm1,SaUm2}, Markl--Shnider \cite{MaSh} and Ornaghi \cite{Orna2} for results related to this question, on defining tensor products of $A_\iy$-algebras and $A_\iy$-categories.
\end{rem}

Another interesting question is whether we can find an explicit formula for the $\iy$-functor $\Fun_{A_\iy}^{\un,\iy}$ in \eq{ai9eq18}, written in terms of the presentation of $\ACat^{\un,\iy}_R$ in \S\ref{ai72} as an explicit simplicial set. The author has not managed to find such a general formula. 

However, we can write down explicit formulae for $\Fun_{A_\iy}^{\un,\iy}$ on the $\iy$-sub\-cat\-eg\-ories $(\sN(\ACat^\un_R))^\rop\t\ACat^{\un,\iy}_R$ and~$(\ACat^{\un,\iy}_R)^\rop\t\sN(\ACat^\un_R)$.

\begin{dfn}
\label{ai9def3}
{\bf(a)} Let $R$ be a commutative ring. We will define an $\iy$-functor
\e
\hFun_{A_\iy}^{\un,\iy}:\sN(\ACat^\un_R)^\rop\t\ACat^{\un,\iy}_R \longra \ACat^{\un,\iy}_R.
\label{ai9eq19}
\e
Form the opposite $\iy$-category $\sN(\ACat^\un_R)^\rop$ as in Definitions \ref{ai2def3} and \ref{ai2def6}. Then $N$-simplices of the product $\iy$-category $\sN(\ACat^\un_R)^\rop\t\ACat^{\un,\iy}_R$ are pairs 
\begin{align*}
(S,T):=\bigl(&((\sO^i)_{i=0}^N,(f^{ij})_{0\le i\le j\le N},(\sH^i)_{i=0}^N,(\m^I)_{\es\ne I\subseteq\{0,\ldots,N\}}),\\
&((\sP^i)_{i=0}^N,(g^{ij})_{0\le i\le j\le N},(\sI^i)_{i=0}^N,(\n^I)_{\es\ne I\subseteq\{0,\ldots,N\}})\bigr)
\end{align*}
where $S$ is an $N$-simplex in $\sN(\ACat^\un_R)\subset\ACat^\un_R$ and $T$ is an $N$-simplex in $\ACat^\un_R$, as in \S\ref{ai72}. (Here we regard $\iinc^\un_R$ in Definition \ref{ai7def11} as embedding  $\sN(\ACat^\un_R)$ as a simplicial subset of $\ACat^\un_R$.) A morphism $\al:[\bar N]\ra[N]$ in $\De$ acts on $(S,T)$ as 
\begin{equation*}
\al:(S,T)\longmapsto\bigl(\ACat^{\un,\iy}_R(\al^\rop)(S),\ACat^{\un,\iy}_R(\al)(T)\bigr),
\end{equation*}
where $\ACat^{\un,\iy}_R(\al)$ is as in Definition \ref{ai7def3}(b), and $\al^\rop(i)=N-\al(\bar N-i)$.

For $i=0,\ldots,N$, write $\sA^i,\sB^i$ for the $A_\iy$-categories over $R$ at vertex $i$ of $S,T$, so that $\Obj\sA^i=\sO^i,$ $\Obj\sB^i=\sP^i$, $\Hom_{\sA^i}(x_0,x_1)=\sH^i(x_0,x_1)$, $\Hom_{\sB^i}(y_0,y_1)=\sI^i(y_0,y_1)$, $\m_k^{\sA^i;x_0,\ldots,x_k}=\m_k^{\{i\};x_0,\ldots,x_k}$, and $\m_k^{\sB^i;y_0,\ldots,y_k}=\n_k^{\{i\};y_0,\ldots,y_k}$. For $0\le i\le j\le N$, write $\fF^{ij}:\sA^i\ra\sA^j$ and $\fG^{ij}:\sB^i\ra\sB^j$ for the $A_\iy$-functors at edges $[i,j]$ of $S,T$, so that $\fF_\Obj^{ij}=f^{ij}$, $\fG_\Obj^{ij}=g^{ij}$, $\fF_k^{ij;x_0,\ldots,x_k}=\m_k^{\{i,j\};x_0,\ldots,x_k}$, and~$\fG_k^{ij;y_0,\ldots,y_k}=\n_k^{\{i,j\};y_0,\ldots,y_k}$. 

As $S$ is an $N$-simplex in $\sN(\ACat^\un_R)\subset\ACat^\un_R$, we have $\fF^{ik}=\fF^{jk}\ci\fF^{ij}$ for $0\le i<j<k\le N$, and $\m_0^{\{i,j,k\};x}=\boo^k_{\fF^{ik}(x)}$, and $\m^{I;x_0,\ldots,x_k}_k=0$ if $\md{I}\ge 4$ or $\md{I}=3$ and~$k\ge 1$.

To deal with $N$-simplices $S$ in the opposite $\iy$-categories  $(\sN(\ACat^\un_R))^\rop$ or $(\ACat^{\un,\iy}_R)^\rop$ we will regard $S$ as an $N$-simplex in $\sN(\ACat^\un_R)$ or $\ACat^{\un,\iy}_R$, but we will consistently reverse the order of the vertices $0,1,\ldots,N$ in $\De^N$, replacing vertex $i$ by vertex $N-i$. So our formulae will involve $\sO^{N-i},\ab f^{(N-i)(N-j)},\ab\sA^{N-i},\ab\fF^{(N-i)(N-j)}$, and so on. We write $\d_{k,l;x_0,\ldots,x_l}^{\sA,I;w_0,\ldots,w_k}$ for the operators $\d_{k,l;x_0,\ldots,x_l}^{I;w_0,\ldots,w_k}$ in Definition \ref{ai7def2} defined using the data in the $N$-simplex~$S$.

Given such an $N$-simplex $(S,T)$ for $N\ge 0$ we will define another $N$-simplex 
\e
\hFun_{A_\iy}^{\un,\iy}(S,T):=\bigl((\sQ^i)_{i=0}^N,(h^{ij})_{0\le i\le j\le N},(\sJ^i)_{i=0}^N,(\o^I)_{\es\ne I\subseteq\{0,\ldots,N\}}\bigr)
\label{ai9eq20}
\e
in $\ACat^{\un,\iy}_R$. Theorem \ref{ai9thm7} shows that these define an $\iy$-functor~\eq{ai9eq19}.

Define the data in \eq{ai9eq20} as follows:
\begin{itemize}
\setlength{\itemsep}{0pt}
\setlength{\parsep}{0pt}
\item[{\bf(i)}] For $i=0,\ldots,N$, $\sQ^i$ is the set of unital $A_\iy$-functors $\fH^i:\sA^{N-i}\ra\sB^i$. 
\item[{\bf(ii)}] For $0\le i\le j\le N$, $h^{ij}:\sQ^i\ra\sQ^j$ maps $h^{ij}:\fH^i\mapsto\fG^{ij}\ci\fH^i\ci\fF^{(N-j)(N-i)}$.

Suppose $\es\ne I\subseteq\{0,1,\ldots,N\}$ with $I=\{i_0,i_1,\ldots,i_M\}$ for $i_0<i_1<\cdots<i_M$. Define $h^I:\sQ^{i_0}\ra\sQ^{i_M}$ using the $h^{ij}$ as in~\eq{ai7eq2}.
\item[{\bf(iii)}] For $i=0,\ldots,N$ and $\fH,\fH'$ in $\sQ^i$, define $\sJ^i(\fH,\fH')=\Hom^\un_{A_\iy}(\fH,\fH')$ as a graded $R$-module, where $\Hom^\un_{A_\iy}(-,-)$ is as in Definition~\ref{ai4def8}.
\item[{\bf(iv)}] Suppose $I=\{i_0,i_1,\ldots,i_M\}$ with $0\le i_0<i_1<\cdots<i_M\le N$, and $p\ge 0$, and $\fH^0,\ldots,\fH^p\in\sQ^{i_0}$. We will define an $R$-linear map
\ea
\o_p^{I;\fH^0,\ldots,\fH^p}:\sJ^{i_0}&(\fH^0,\fH^1)\ot_R\sJ^{i_0}(\fH^1,\fH^2)\ot_R\cdots \ot_R\sJ^{i_0}(\fH^{p-1},\fH^p)
\nonumber\\
&\longra\sJ^{i_M}(h^I(\fH^0),h^{i_0i_M}(\fH^p))
\label{ai9eq21}
\ea
graded of degree $2-M-p=3-\md{I}-p$. If $\fT^j\in\sJ^{i_0,d_j}(\fH^{j-1},\fH^j)$ for $j=1,\ldots,p$ then we will construct
\e
\begin{split}
&\o_p^{I;\fH^0,\ldots,\fH^p}(\fT^1\ot\cdots\ot\fT^p)=\fU\\
&\text{in $\sJ^{i_M,d_1+\cdots+d_p+2-M-p}(h^I(\fH^0),h^{i_0i_M}(\fH^p))$.}
\end{split}
\label{ai9eq22}
\e
To define $\fU$ in \eq{ai9eq22}, for all $k\ge 0$ and $w_0,\ldots,w_k\in\sO^{N-i_M}$, we define
\end{itemize} 
\begin{align*}
&\fU_k^{w_0,\ldots,w_k}:\Hom_{\sA^{N-i_M}}(w_0,w_1)\ot_R\Hom_{\sA^{N-i_M}}(w_1,w_2)\ot_R\cdots \\
&\ot_R\Hom_{\sA^{N-i_M}}(w_{k-1},w_k)\longra\Hom_{\sB^{i_M}}(h^I(\fH^0)(w_0),h^{i_0i_M}(\fH^p)(w_k))
\end{align*}
\begin{itemize}
\setlength{\itemsep}{0pt}
\setlength{\parsep}{0pt}
\item[] by, if $M=0$,
\end{itemize}
\ea
&\fU_k^{w_0,\ldots,w_k}=
\sum_{\begin{subarray}{l} 
m:\; p\le m\le k+p,\\
y_0,\ldots,y_m\in\sP^{i_0}
\end{subarray}}
\begin{aligned}[t]&
\n_m^{\{i_0\};y_0,\ldots,y_m}\ci{} \\
&\bD_{k,m;y_0,\ldots,y_m}^{p;w_0,\ldots,w_k}(\fH^0,\ldots,\fH^p,\fT^1\ot\cdots\ot\fT^p)
\end{aligned}
\label{ai9eq23}\\
&+(-1)^{d_1+\cdots+d_p}\sum_{\begin{subarray}{l} 
l:\; 0\le l\le k, \\
x_0,\ldots,x_l\in\sO^{N-i_0}, \\
\end{subarray}}
\begin{aligned}[t]&
\bD_{l,1;\fH^0(x_0),\fH^p(x_l)}^{p;x_0,\ldots,x_l}(\fH^0,\ldots,\fH^p,\fT^1\ot\cdots\ot\fT^p)\\
&{}\ci\d^{\sA,\{N-i_0\};w_0,\ldots,w_k}_{k,l;x_0,\ldots,x_l},
\end{aligned}
\nonumber
\ea
\begin{itemize}
\setlength{\itemsep}{0pt}
\setlength{\parsep}{0pt}
\item[] and if $M>0$
\end{itemize}
\ea
&\fU_k^{w_0,\ldots,w_k}=
\sum_{\begin{subarray}{l} 
l,m:\\
0\le l\le k,\; p\le m\le l+p,\\
x_0,\ldots,x_l\in\sO^{N-i_0}, \\
y_0,\ldots,y_m\in\sP^{i_0}
\end{subarray}}
\begin{aligned}[t]&
\n_m^{I;y_0,\ldots,y_m}\ci{} \\
&\bD_{l,m;y_0,\ldots,y_m}^{p;x_0,\ldots,x_l}(\fH^0,\ldots,\fH^p,\fT^1\ot\cdots\ot\fT^p)\\
&{}\ci\d^{\sA,\{N-i_M,N-i_0\};w_0,\ldots,w_k}_{k,l;x_0,\ldots,x_l}.
\end{aligned}
\label{ai9eq24}
\ea

Here $\d^{\sA,\{N-i_M,N-i_0\};w_0,\ldots,w_k}_{k,l;x_0,\ldots,x_l}$ is $\d^{I;w_0,\ldots,w_k}_{k,l;x_0,\ldots,x_l}$ from \S\ref{ai71} for the data $S=((\sO^i),\ab(f^{ij}),\ab(\sH^i),\ab(\m^I))$ and $I=\{N-i_M,N-i_0\}$, and if $l,m,p\ge 0$, and $x_0,\ldots,x_l\in\sO^{N-i_0}$, and $y_0,\ldots,y_m\in\sP^{i_0}$, and $\fH^0,\ldots,\fH^p\in\sQ^{i_0}$, we define
\ea
&\bD_{l,m;y_0,\ldots,y_m}^{p;x_0,\ldots,x_l}(\fH^0,\ldots,\fH^p,\fT^1\ot\cdots\ot\fT^p)
\label{ai9eq25}\\
&=\sum_{\begin{subarray}{l} 
l=b_1+\cdots+b_m,\; \\
0=c_0<c_1< \\
\cdots<c_{p+1}=m+1: \\
b_{c_i}\ge 0, \; 1\le i\le p, \\
b_j\ge 1, \; j\ne c_i, \\
\text{Set $s_j=\sum_{i=1}^jb_i$.}\\
\text{Require $y_j=\fH^i(x_{s_j})$,} \\ 
0\le j\le m,\; c_i\le j<c_{i+1}.
\end{subarray}}
\begin{aligned}[t]&
(-1)^{\sum_{i=1}^p(m+1)(d_i+1)}\cdot
(-1)^{\sum_{i=1}^p\sum_{j=1}^{c_i-1}(b_j+1)(d_i+1)}\cdot{} \\
&(-1)^{\sum_{i=1}^m(b_i+1)(m-i)}\cdot{}\\
&\fH_{b_1}^{0;x_{s_0},\ldots,x_{s_1}}\ot\fH_{b_2}^{0;x_{s_1},\ldots,x_{s_2}}\ot\cdots\ot\fH_{b_{c_1-1}}^{0;x_{s_{c_1-2}},\ldots,x_{s_{c_1-1}}}\ot{} \\
&\fT_{b_{c_1}}^{1;x_{s_{c_1-1}},\ldots,x_{s_{c_1}}}\!\ot\!\fH_{b_{c_1+1}}^{1;x_{s_{c_1}},\ldots,x_{s_{c_1+1}}}\!\ot\!\cdots\!\ot\!\fH_{b_{c_2-1}}^{1;x_{s_{c_2-2}},\ldots,x_{s_{c_2-1}}} \\
&{}\ot\!
\fT_{b_{c_2}}^{2;x_{s_{c_2-1}},\ldots,x_{s_{c_2}}}\ot \cdots \ot\fT_{b_{c_p}}^{p;x_{s_{c_p-1}},\ldots,x_{s_{c_p}}} \\
&\ot\!\fH_{b_{c_p+1}}^{p;x_{s_{c_p}},\ldots,x_{s_{c_p+1}}}\ot\cdots\ot\fH_{b_m}^{p;x_{s_{m-1}},\ldots,x_{s_m}}.
\end{aligned}
\nonumber
\ea

This defines $\hFun_{A_\iy}^{\un,\iy}$ in \eq{ai9eq19}. Note that when $N=0$, the 0-simplex $\hFun_{A_\iy}^{\un,\iy}(S,T)$ is the unital $A_\iy$-category $\Fun_{A_\iy}^\un(\sA^0,\sB^0)$ from Definition~\ref{ai4def8}.

\smallskip

\noindent{\bf(b)} Similarly we define an $\iy$-functor
\e
\cFun_{A_\iy}^{\un,\iy}:(\ACat^{\un,\iy}_R)^\rop\t\sN(\ACat^\un_R) \longra \ACat^{\un,\iy}_R.
\label{ai9eq26}
\e
Most of the definition is the same as in {\bf(a)}. The differences are that as now $T$ is an $N$-simplex in $\sN(\ACat^\un_R)\subset\ACat^\un_R$, we have $\fG^{ik}=\fG^{jk}\ci\fG^{ij}$ for $0\le i\le j\le k\le N$, and $\n_0^{ijk;x}=\boo^k_{\fG^{ik}(x)}$, and $\n^{I;x_0,\ldots,x_k}_k=0$ if $\md{I}\ge 4$ or $\md{I}=3$ and $k\ge 1$, and we leave \eq{ai9eq23} unchanged, but replace \eq{ai9eq24} by
\ea
&\fU_k^{w_0,\ldots,w_k}=(-1)^{(M+1)(d_1+\cdots+d_p)}\cdot{}
\label{ai9eq27}\\
&\sum_{\begin{subarray}{l} 
l,m:\\
0\le l\le k,\; p\le m\le l+p,\\
x_0,\ldots,x_l\in\sO^{N-i_0}, \\
y_0,\ldots,y_m\in\sP^{i_0}
\end{subarray}}
\begin{aligned}[t]&
\n_m^{\{i_0,i_M\};y_0,\ldots,y_m}\ci{} \\
&\bD_{l,m;y_0,\ldots,y_m}^{p;x_0,\ldots,x_l}(\fH^0,\ldots,\fH^p,\fT^1\ot\cdots\ot\fT^p)\\
&{}\ci\d_{k,l;x_0,\ldots,x_l}^{\sA,\{N-i:i\in I\};w_0,\ldots,w_k}.
\end{aligned}
\nonumber
\ea
\end{dfn}

The next theorem will be proved in \S\ref{ai116}.

\begin{thm}
\label{ai9thm7}
The $\iy$-functors in \eq{ai9eq19} and\/ \eq{ai9eq26} are well defined. Relating \eq{ai9eq18} with\/ {\rm\eq{ai9eq19},} in the following diagram of\/ $\iy$-categories, the top rectangle commutes strictly, and the bottom commutes up to homotopy:
\ea
\begin{gathered}
\xymatrix@!0@C=280pt@R=35pt{
*+[r]{\sN(\ACat^\un_R)^\rop\t\sN(\ACat^\un_R)} \ar[r]_(0.66){\sN(\Fun_{A_\iy}^\un)} \ar[d]^{\id\t\iinc^\un_R} \ar@{}@<-1ex>[dr]^{=} & *+[l]{\sN(\ACat^\un_R)} \ar[d]_{\iinc^\un_R} \\
*+[r]{\sN(\ACat^\un_R)^\rop\t\ACat^{\un,\iy}_R} \ar[r]^(0.66){\hFun_{A_\iy}^{\un,\iy}} \ar[d]^{(\iinc^\un_R)^\rop\t\id} \drtwocell_{}\omit^{}\omit{^{}} & *+[l]{\ACat^{\un,\iy}_R} \ar@{=}[d]_{\id} \\
*+[r]{(\ACat^{\un,\iy}_R)^\rop\t\ACat^{\un,\iy}_R} \ar@{..>}[r]^(0.66){\Fun_{A_\iy}^{\un,\iy}} & *+[l]{\ACat^{\un,\iy}_R.\!}}
\end{gathered}
\label{ai9eq28}
\ea
The analogue holds relating \eq{ai9eq18} with\/ \eq{ai9eq26}.
\end{thm}

We did not treat {\it filtered\/} $A_\iy$-categories from Chapter \ref{ai8} in the material above, because as in Remark \ref{ai9rem4} the author does not know how to write $\AfCat^{\un,\iy}_R,\ab\AwCat^{\un,\iy}_R$ as localizations, so we have no filtered analogues of \eq{ai9eq17}. However, we can generalize Definition \ref{ai9def3} to the filtered case.

\begin{dfn}
\label{ai9def4}
Let $(\Ga,\pre)$ be a po-monoid and $R$ a complete $\Ga$-filtered ring. We will define $\iy$-functors
\ea
\hFun_{A^\rfi_\iy}^{\un,\iy}&:\sN(\AfCat^\un_R)^\rop\t\AfCat^{\un,\iy}_R \longra \AfCat^{\un,\iy}_{R,2,2},
\label{ai9eq29}\\
\hFun_{A^\wf_\iy}^{\un,\iy}&:\sN(\AfCat^\un_R)^\rop\t\AwCat^{\un,\iy}_R \longra \AwCat^{\un,\iy}_{R,2,2}.
\label{ai9eq30}\\
\cFun_{A^\rfi_\iy}^{\un,\iy}&:(\AfCat^{\un,\iy}_R)^\rop\t\sN(\AfCat^\un_R) \longra \AfCat^{\un,\iy}_{R,2,2},
\label{ai9eq31}\\
\cFun_{A^\wf_\iy}^{\un,\iy}&:(\AwCat^{\un,\iy}_R)^\rop\t\sN(\AfCat^\un_R) \longra \AwCat^{\un,\iy}_{R,2,2}.
\label{ai9eq32}
\ea
To define \eq{ai9eq29} and \eq{ai9eq31} we follow Definition \ref{ai9def3} with the following differences: $\sA^i,\ab\sB^i,\ab\fH^i,\ab\fF^{ij},\ab\fG^{ij}$ are {\it filtered\/} unital $A_\iy$-categories and unital $A_\iy$-functors, and we replace $\ACat^\un_R,\ACat^{\un,\iy}_R$ by $\AfCat^\un_R,\AfCat^{\un,\iy}_R$, and $\Hom^\un_{A_\iy}(\fH,\fH')$ by $\Hom^\un_{A^\rfi_\iy}(\fH,\fH')$, and we change the limits of sums in \eq{ai9eq23}--\eq{ai9eq25} and \eq{ai9eq27}:
\begin{itemize}
\setlength{\itemsep}{0pt}
\setlength{\parsep}{0pt}
\item[(i)] In \eq{ai9eq23} the first sum is over $m\ge p$ and the second over $l\ge 0$.
\item[(ii)] In \eq{ai9eq24} and \eq{ai9eq27} we sum over $l\ge 0$ and $m\ge p$. 
\item[(iii)] In \eq{ai9eq25} we do not require $b_j\ge 1$, $j\ne c_i$.
\end{itemize}
Thus these equations become infinite sums, which converge in the filtered sense. 

Then for \eq{ai9eq29} the 0-simplex at vertex $i$ of the $N$-simplex we define is the unital filtered $A_\iy$-category $\Fun_{A_\iy^\rfi}^\un(\sA^{N-i},\sB^i)$, and the 1-simplex on edge $i\ra j$ is the unital filtered $A_\iy$-functor $(\fF^{(N-j)(N-i)})^*\ci\fG^{ij}_*:\Fun_{A_\iy^\rfi}^\un(\sA^{N-i},\sB^i)\ra\Fun_{A_\iy^\rfi}^\un(\sA^{N-j},\sB^j)$, as in \S\ref{ai55}. Note as in Definition \ref{ai5def10} and Remark \ref{ai5rem5} that these are {\it strict\/} filtered $A_\iy$-categories and {\it strict\/} filtered $A_\iy$-functors, so that $\m_0^{I;x_0}=0$ for $\md{I}=1,2$. Hence the target $\iy$-category in \eq{ai9eq29} is $\AfCat^{\un,\iy}_{R,2,2}$ from Definition \ref{ai8def4}. Equations \eq{ai9eq30}--\eq{ai9eq32} are similar.
\end{dfn}

The proof of the next theorem is the same as that of the first part of Theorem \ref{ai9thm7} in \S\ref{ai116}, except that the limits in sums are slightly different, which leads to infinite sums convergent in the filtered sense.

\begin{thm}
\label{ai9thm8}
The $\iy$-functors in {\rm\eq{ai9eq29}--\eq{ai9eq32}} are well defined. 
\end{thm}

\begin{dfn}
\label{ai9def5}
Let $(\Ga,\pre)$ be a po-monoid and $R$ a complete $\Ga$-filtered ring. Let $\sP^\un$ be the unital point $A_\iy$-category over $R$ from Definition \ref{ai4def11}, made into a strict filtered $A_\iy$-category in the obvious way. Following Definition \ref{ai5def15} and Proposition \ref{ai5prop2} and using \eq{ai9eq29}--\eq{ai9eq30}, define $\iy$-functors
\ea
\Str^\iy=\hFun_{A^\rfi_\iy}^{\un,\iy}(\sP^\un,-)&:\AfCat^{\un,\iy}_R\longra\AfCat^{\un,\iy}_{R,2,2},
\label{ai9eq33}\\
\Str_\wf^\iy=\hFun_{A^\wf_\iy}^{\un,\iy}(\sP^\un,-)&:\AwCat^{\un,\iy}_R\longra\AwCat^{\un,\iy}_{R,2,2}.
\label{ai9eq34}
\ea
By Proposition \ref{ai5prop2}, these map an object $\sA$ to the {\it associated strict filtered\/ $A_\iy$-category\/} $\Str(\sA)$ from Definition \ref{ai5def15}, and are $\iy$-categorical analogues of the functor $\Str$ in Proposition \ref{ai5prop2}(b),(d).
\end{dfn}

\begin{rem}
\label{ai9rem7}
{\bf(a)} Note that we can compose $\hFun_{A^\rfi_\iy}^{\un,\iy},\ldots,\cFun_{A^\wf_\iy}^{\un,\iy}$ in \eq{ai9eq29}--\eq{ai9eq32} with $\Pi_{A_\iy^\rfi}^{A_\iy},\Pi_{A_\iy^\wf}^{A_\iy}$ in \eq{ai8eq8}--\eq{ai8eq9} to map to $\ACat^{\un,\iy}_R$. The same applies to $\Str^\iy$ and $\Str_\wf^\iy$ in~\eq{ai9eq33}--\eq{ai9eq34}.
\smallskip

\noindent{\bf(b)} It seems likely that \eq{ai9eq29}--\eq{ai9eq30} and \eq{ai9eq31}--\eq{ai9eq32} extend (at least up to homotopy) to $\iy$-functors
\begin{align*}
\Fun_{A^\rfi_\iy}^{\un,\iy}&:(\AfCat^{\un,\iy}_R)^\rop\t\AfCat^{\un,\iy}_R \longra \AfCat^{\un,\iy}_{R,2,2}, \\
\Fun_{A^\wf_\iy}^{\un,\iy}&:(\AwCat^{\un,\iy}_R)^\rop\t\AwCat^{\un,\iy}_R \longra \AwCat^{\un,\iy}_{R,2,2},
\end{align*}
but we will not prove this.
\end{rem}

\subsection{\texorpdfstring{The twist $\iy$-functor}{The twist ∞-functor}}
\label{ai95}

We now give an $\iy$-categorical version of the twist functor $\Tw$ in \S\ref{ai49} and \S\ref{ai58}, in the unital case, by generalizing Definition~\ref{ai5def16}.

\begin{dfn}
\label{ai9def6}
Let $R$ be a commutative ring. We will define an explicit $\iy$-functor
\e
\Tw^\iy:\ACat^{\un,\iy}_R \longra \ADCat^{\un,\iy}_R,
\label{ai9eq35}
\e
where $\ACat^{\un,\iy}_R$ is as in \S\ref{ai72} and $\ADCat^{\un,\iy}_R$ as in Definition \ref{ai7def13}.

Make $R$ into a complete filtered ring over the po-monoid $(\Z,\le)$ by setting $F^\la R=R$ if $\la\le 0$ in $\Z$ and $F^\la R=0$ if $\la>0$ in $\Z$. Generalizing \eq{ai5eq20}, we will define $\Tw^\iy$ in \eq{ai9eq35} as the composition of $\iy$-functors
\e
\!\!\xymatrix@!0@C=84pt{
\ACat^{\un,\iy}_R \ar[r]^(0.5){\Si^\iy} & \AfCat^{\un,\iy}_R \ar[r]^(0.52){\raisebox{8pt}{$\st \hFun^\un_{A^\rfi_\iy}(\sP^\un,-)$}} & \AfCat^{\un,\iy}_{R,2,2} \ar[r]^{\Pi_{A_\iy^\rfi}^{A_\iy}} & \ACat^{\un,\iy}_R.  }\!\!
\label{ai9eq36}
\e
Here $\Si^\iy$ is an $\iy$-categorical generalization of $\Si$ in Definition \ref{ai5def16} which we will define below. Also $\sP^\un$ is the unital point $A_\iy$-category from Definition \ref{ai4def11}, with $\Z$-filtration on $\Hom_{\sP^\un}(*,*)=R$ as above, and $\hFun^\un_{A^\rfi_\iy}(\sP^\un,-)$ is as in \eq{ai9eq29}, restricted to the 0-simplex $\sP^\un$ in $\sN(\AfCat^\un_R)^\rop$, and $\Pi_{A_\iy^\rfi}^{A_\iy}$ is as in~\eq{ai8eq8}.

It will follow by comparing the definition with Definition \ref{ai5def16} that the following diagram of $\iy$-functors strictly commutes:
\e
\begin{gathered}
\xymatrix@C=180pt@R=15pt{
*+[r]{\sN(\ACat^\un_R)} \ar@<-1ex>@{}[dr]^{=} \ar[d]^{\iinc^\un_R} \ar[r]_(0.4){\sN(\Tw)} & *+[l]{\sN(\ACat^\un_R)}  \ar[d]_{\iinc^\un_R} \\
*+[r]{\ACat^{\un,\iy}_R} \ar[r]^(0.6){\Tw^\iy} & *+[l]{\ACat^{\un,\iy}_R.\!} }
\end{gathered}
\label{ai9eq37}
\e
Since $\Tw$ maps $\ACat^\un_R\ra\ADCat^\un_R\subset\ACat^\un_R$, it follows that $\Tw^\iy$ maps to $\ADCat^{\un,\iy}_R\subset\ACat^{\un,\iy}_R$, as in \eq{ai9eq35}, so we can rewrite \eq{ai9eq37} as
\e
\begin{gathered}
\xymatrix@C=180pt@R=15pt{
*+[r]{\sN(\ACat^\un_R)} \ar@<-1ex>@{}[dr]^{=} \ar[d]^{\iinc^\un_R} \ar[r]_(0.4){\sN(\Tw)} & *+[l]{\sN(\ADCat^\un_R)}  \ar[d]_{\iinc^\un_R} \\
*+[r]{\ACat^{\un,\iy}_R} \ar[r]^(0.6){\Tw^\iy} & *+[l]{\ADCat^{\un,\iy}_R.\!} }
\end{gathered}
\label{ai9eq38}
\e

To define $\Si^\iy$ in \eq{ai9eq36}, suppose $N\ge 0$ and $\bigl((\sO^i)_{i=0}^N,\ab(f^{ij})_{0\le i\le j\le N},\ab(\sH^i)_{i=0}^N,\ab(\m^I)_{\es\ne I\subseteq\{0,\ldots,N\}}\bigr)$ is an $N$-simplex in $\ACat^{\un,\iy}_R\subset\hACat^{\un,\iy}_R$, as in Definitions \ref{ai7def3} and \ref{ai7def4}. Define 
\e
\begin{split}
&\Si^\iy\bigl((\sO^i)_{i=0}^N,(f^{ij})_{0\le i\le j\le N},(\sH^i)_{i=0}^N,(\m^I)_{\es\ne I\subseteq\{0,\ldots,N\}}\bigr)\\
&\quad =\bigl((\sP^i)_{i=0}^N,(g^{ij})_{0\le i\le j\le N},(\sI^i)_{i=0}^N,(\n^I)_{\es\ne I\subseteq\{0,\ldots,N\}}\bigr),
\end{split}
\label{ai9eq39}
\e
where, generalizing Definition \ref{ai5def16}, we set:
\begin{itemize}
\setlength{\itemsep}{0pt}
\setlength{\parsep}{0pt}
\item[(a)] Elements of $\sP^i$ are triples $\bigl((x_1,\ldots,x_l),(s_1,\ldots,s_l),(f_1,\ldots,f_l)\bigr)$ where $l\ge 0$, and $x_1,\ldots,x_l\in\sO^i$, and $s_i,f_i\in\Z$ for $i=1,\ldots,l$. As a shorthand we write $(\bs x,\bs s,\bs f)=\bigl((x_1,\ldots,x_l),(s_1,\ldots,s_l),(f_1,\ldots,f_l)\bigr)$.
\item[(b)] For $0\le i\le j\le N$, define $g^{ij}:\sP^i\ra\sP^j$ by $g^{ij}\bigl((x_1,\ldots,x_l),\ab(s_1,\ldots,s_l),\ab(f_1,\ldots,f_l)\bigr)=\bigl((f^{ij}(x_1),\ldots,f^{ij}(x_l)),(s_1,\ldots,s_l),(f_1,\ldots,f_l)\bigr)$.

When $I=\{i_0<i_1<\cdots<i_M\}\subseteq\{0,\ldots,N\}$, we define $g^I:\sP^{i_0}\ra\sP^{i_M}$ from the $g^{ij}$ as in~\eq{ai7eq2}.
\item[(c)] Let $\bigl((x_1,\ldots,x_l),(s_1,\ldots,s_l),(f_1,\ldots,f_l)\bigr)$ and $\bigl((x'_1,\ldots,x'_{l'})$, $(s'_1,\ldots,\ab s'_{l'})$, $(f'_1,\ab\ldots,f'_{l'})\bigr)$ be elements of $\sP^i$. Define a graded $R$-module
\begin{align*}
&\sI^i\bigl((\bs x,\bs s,\bs f),(\bs x',\bs s',\bs f')\bigr)=\bigop_{d\in\Z}\sI^i\bigl((\bs x,\bs s,\bs f),(\bs x',\bs s',\bs f')\bigr){}^d\quad\text{by} \\
&\sI^i\bigl((\bs x,\bs s,\bs f),(\bs x',\bs s',\bs f')\bigr){}^d=\bigop_{i=1,\ldots,l,\;
j=1,\ldots,l'\!\!\!\!\!\!\!\!\!\!\!\!\!\!\!\!\!\!\!} \sH^i(x_i,x'_j)^{d+s_i-s'_j}.
\end{align*}
We write elements of $\sI^i\bigl((\bs x,\bs s,\bs f),(\bs x',\bs s',\bs f')\bigr)$ as $(\fT_{ij})_{i=1,\ldots,l}^{j=1,\ldots,l'}$ in the obvious way. Define a $\Z$-filtration on $\sI^i\bigl((\bs x,\bs s,\bs f),(\bs x',\bs s',\bs f')\bigr)^d$ by, for $\la\in\Z$
\begin{equation*}
F^\la\sI^i\bigl((\bs x,\bs s,\bs f),(\bs x',\bs s',\bs f')\bigr){}^d=\bigop_{\begin{subarray}{l} i=1,\ldots,l,\; j=1,\ldots,l': \\ f_i-f_j\ge\la \end{subarray}\!\!\!\!\!\!\!\!\!\!\!\!\!\!\!\!\!\!}
\sH(x_i,x'_j)^{d+s_i-s'_j}.
\end{equation*}
\item[(d)] For all $I=\{i_0<i_1<\cdots<i_M\}\subseteq\{0,\ldots,N\}$, $k\ge 0$ and elements $x_0,\ldots,x_k\in \sO^{i_0}$ and elements $(\bs x^0,\bs s^0,\bs f^0),\ldots,(\bs x^k,\bs s^k,\bs f^k)$ of $\sI^{i_0}$, define a graded $R$-linear map of degree $2-k$
\begin{align*}
\n_k^{I;(\bs x^0,\bs s^0,\bs f^0),\ldots,(\bs x^k,\bs s^k,\bs f^k)}:\,&\sI^{i_0}\bigl((\bs x^0,\bs s^0,\bs f^0),(\bs x^1,\bs s^1,\bs f^1)\bigr)\ot_R\cdots \ot_R\\
&\sI^{i_0}\bigl((\bs x^{k-1},\bs s^{k-1},\bs f^{k-1}),(\bs x^k,\bs s^k,\bs f^k)\bigr)\\
&\quad \longra\sI^{i_M}\bigl(g^I(\bs x^0,\bs s^0,\bs f^0),g^{i_0i_M}(\bs x^k,\bs s^k,\bs f^k)\bigr)
\end{align*}
by, if $k=0$ and $j^0,\ti\jmath^0=1,\ldots,l^0$ then
\e
\bigl(\n_0^{I;(\bs x^0,\bs s^0,\bs f^0)}(1)\bigr)_{j^0\ti\jmath^0}=\de_{j^0\ti\jmath^0}\,\m_0^{I;x^0_{j^0}}(1),
\label{ai9eq40}
\e
and if $k\ge 1$ and $j^0=1,\ldots,l^0$ and $j^k=1,\ldots,l^k$ then
\end{itemize}
\begin{align*}
&\Bigl(\n_k^{I;(\bs x^0,\bs s^0,\bs f^0),\ldots,(\bs x^k,\bs s^k,\bs f^k)}\bigl((\fT^1_{i^0i^1})_{i^0=1,\ldots,l^0}^{i^1=1,\ldots,l^1}\ot\cdots\ot (\fT^k_{i^{k-1}i^k})_{i^{k-1}=1,\ldots,l^{k-1}}^{i^k=1,\ldots,l^k}\bigr)\Bigr)_{j^0j^k} \\
&=\sum_{\begin{subarray}{l} j^1=1,\ldots,l^1,\ldots, \\
j^{k-1}=1,\ldots,l^{k-1}
\end{subarray}}
\m_k^{I;x^0_{j^0},x^1_{j^1},\ldots,x^{k-1}_{j^{k-1}},x^k_{j^k}}\bigl(\fT^1_{j^0j^1}\ot\cdots\ot \fT^k_{j^{k-1}j^k}\bigr).
\end{align*}

It is now straightforward to show that \eq{ai9eq39} is an $N$-simplex in $\hAfCat^{\un,\iy}_R$, and then in $\AfCat^{\un,\iy}_R$. In particular, the $(\cdots)^{j^0j^n}$ component of \eq{ai8eq4}--\eq{ai8eq5} for \eq{ai9eq39} is the sum over $j^1=1,\ldots,l^1,\ldots,j^{n-1}=1,\ldots,l^{n-1}$ of \eq{ai7eq15}--\eq{ai7eq16} for $\bigl((\sO^i),(f^{ij}),(\sH^i),(\m^I)\bigr)$ at $x^0_{j^0},\ldots,x^n_{j^n}$. By \eq{ai9eq40}, the isomorphism conditions on $\n_0^{\{i_0,i_1,i_2\}}$ in Definition \ref{ai8def3} follows from a product of the isomorphism conditions on $\m_0^{\{i_0,i_1,i_2\}}$ in Definition \ref{ai7def4}. The above definitions are clearly compatible with the actions of morphisms $\al:[\bar N]\ra[N]$ in $\De$ on the simplicial sets, so $\Si^\iy$ is an $\iy$-functor.
\end{dfn}

\begin{rem}
\label{ai9rem8}
Combining equations \eq{ai9eq11}, \eq{ai9eq12} and \eq{ai9eq38}, we see that there exists an $\iy$-functor $\tTw^\iy$ making the following commute up to homotopy:
\e
\begin{gathered}
\xymatrix@C=140pt@R=15pt{
*+[r]{\sN(\ACat^\un_R)} \ar[r]_(0.36){\sN(\Tw)} \ar[d]^{\iinc^\un_R} \ar@<-1.4ex>@/_.6pc/[dd]_{\Pi_\loc^\iy} \ar@<-1ex>@{}[dr]^{=} & *+[l]{\sN(\ADCat^\un_R)} \ar[d]_{\iinc^\un_R} \ar@<1.4ex>@/^.6pc/[dd]^{\Pi_\loc^\iy} \\
*+[r]{\ACat^{\un,\iy}_R} \ar[r]^(0.66){\Tw^\iy} \ar[d]^{\eq{ai9eq11}}_(0.4)\simeq \drtwocell_{}\omit^{}\omit{^{}} & *+[l]{\ADCat^{\un,\iy}_R\!} \ar[d]_{\eq{ai9eq12}}^(0.3)\simeq \\
*+[r]{\ACat^\un_R[\HE_\un^{-1}]^\iy} \ar@{..>}[r]^(0.57){\tTw^\iy}   & *+[l]{\ADCat^\un_R[\HE_{\un,\De}^{-1}]^\iy.\!} 
}
\end{gathered}
\label{ai9eq41}
\e
Therefore $\Tw$ in \S\ref{ai49} maps $\HE_\un^{-1}\ra\HE_\un^{-1}$, so the outer rectangle
\e
\begin{gathered}
\xymatrix@C=180pt@R=15pt{
*+[r]{\sN(\ACat^\un_R)} \drtwocell_{}\omit^{}\omit{^{}} \ar[d]^{\Pi_\loc^\iy} \ar[r]_(0.4){\sN(\Tw)} & *+[l]{\sN(\ADCat^\un_R)}  \ar[d]_{\Pi_\loc^\iy} \\
*+[r]{\ACat^\un_R[\HE_\un^{-1}]^\iy} \ar[r]^(0.6){\tTw^\iy} & *+[l]{\ADCat^\un_R[\HE_{\un,\De}^{-1}]^\iy} }
\end{gathered}
\label{ai9eq42}
\e
is obtained by localization. Presumably we could have proved directly that $\Tw$ maps $\HE_\un^{-1}\ra\HE_\un^{-1}$, and so obtained \eq{ai9eq42} without first constructing $\Tw^\iy$. But $\Tw^\iy$ has the advantage of being explicit.
\end{rem}

\begin{thm}
\label{ai9thm9}
Let\/ $R$ be a commutative ring, and\/ $\Tw^\iy,\tTw^\iy$ be as in \eq{ai9eq35} and\/ \eq{ai9eq41}. Then:
\smallskip

\noindent{\bf(a)} There are $\iy$-natural isomorphisms of\/ $\iy$-functors $\eta:\Tw^\iy\Ra\Tw^\iy\ci\ab\Tw^\iy$ and\/~$\ti\eta:\tTw^\iy\Ra \tTw^\iy\ci\tTw^\iy,$ in the sense of Definition\/~{\rm\ref{ai2def5}}.
\smallskip

\noindent{\bf(b)} $\Tw^\iy$ and\/ $\tTw^\iy$ are left adjoint to the fully faithful inclusion $\iy$-functors
\begin{align*}
\inc&:\ADCat^{\un,\iy}_R\longhookrightarrow\ACat^{\un,\iy}_R,\\
\inc&:\ADCat^\un_R[\HE_{\un,\De}^{-1}]^\iy\longhookrightarrow\ACat^\un_R[\HE_\un^{-1}]^\iy.
\end{align*}
\end{thm}

\begin{proof}
Recall from \S\ref{ai49} that natural transformation $\inc:\Id\Ra\Tw$ of functors $\ACat^\un_R\ra\ACat^\un_R$. Composing on the left and right with $\Tw$ give two natural transformations $\id_{\Tw}*\inc,\inc*\id_{\Tw}:\Tw\Ra\Tw\ci\Tw$ of functors $\ACat^\un_R\ra\ACat^\un_R$. For each object $\sA$ this gives two $A_\iy$-functors $\Tw(\inc_\sA),\inc_{\Tw\sA}:\Tw\sA\hookra\Tw\Tw\sA$, both homotopy equivalences.

By properties of localizations, there is a homotopy commutative diagram
\begin{equation*}
\xymatrix@C=180pt@R=15pt{
*+[r]{\sN(\ACat^\un_R)\t\De^1} \drtwocell_{}\omit^{}\omit{^{}} \ar[d]^{\Pi_\loc^\iy\t\id_{\De^1}} \ar[r]_(0.4){\sN(\inc)} & *+[l]{\sN(\ACat^\un_R)}  \ar[d]_{\Pi_\loc^\iy} \\
*+[r]{\ACat^\un_R[\HE_\un^{-1}]^\iy\t\De^1} \ar@{..>}[r]^(0.6){{\rm\ti{\text{\i}}nc}^\iy} & *+[l]{\ACat^\un_R[\HE_\un^{-1}]^\iy,\!} }
\end{equation*}
where $\sN(\inc):\sN(\Id_{\ACat^\un_R})\Ra\sN(\Tw)$ and ${\rm\ti{\text{\i}}nc}^\iy:\Id_{\ACat^\un_R[\HE_\un^{-1}]^\iy}\Ra\tTw^\iy$ are $\iy$-natural transformations, in the sense of Definition \ref{ai2def5}. Composing on the left and right with $\tTw^\iy$ gives two $\iy$-natural transformations $\id_{\tTw^\iy}*{\rm\ti{\text{\i}}nc}^\iy,{\rm\ti{\text{\i}}nc}^\iy *\id_{\tTw^\iy}:\tTw^\iy\Ra\tTw^\iy\ci\tTw^\iy$. Up to equivalences, these act on objects $\sA$ by $\Pi_\loc^\iy(\Tw(\inc_\sA)),\Pi_\loc^\iy(\inc_{\Tw\sA})$ above. Since $\Tw(\inc_\sA),\inc_{\Tw\sA}$ are homotopy equivalences, their images in $\ACat^\un_R[\HE_\un^{-1}]^\iy$ are equivalences. Hence $\id_{\tTw^\iy}*{\rm\ti{\text{\i}}nc}^\iy,{\rm\ti{\text{\i}}nc}^\iy *\id_{\tTw^\iy}$ are $\iy$-natural isomorphisms.

For (a), we define $\ti\eta=\id_{\tTw^\iy}*{\rm\ti{\text{\i}}nc}^\iy:\tTw^\iy\Ra\tTw^\iy\ci\tTw^\iy$. Then by \eq{ai9eq41}, there is an $\iy$-natural isomorphism~$\eta:\Tw^\iy\Ra\Tw^\iy\ci\Tw^\iy$.

For (b), Lurie \cite[Prop.~5.2.74]{Luri2} proves the following result. Suppose $\sC$ is an $\iy$-category, $F:\sC\ra\sC$ is an $\iy$-functor, and $\sD\subset\sC$ is the essential image of $\tau$. Then the following are equivalent: 
\begin{itemize}
\setlength{\itemsep}{0pt}
\setlength{\parsep}{0pt}
\item[(i)] There exists an $\iy$-natural transformation $\al:\Id_\sC\Ra F$ such that for each object $x$ in $\sC$, the morphisms $F(\al(x)),\al(F(x)):F(x)\ra F\ci F(x)$ are equivalences in $\sC$, that is, $\id_F*\al,\al*\id_F:F\Ra F\ci F$ are $\iy$-natural isomorphisms.
\item[(ii)] Regarded as an $\iy$-functor $F:\sC\ra\sD$, $F$ is left adjoint to the inclusion $\inc:\sD\hookra\sC$.
\end{itemize}
Apply this with $\sC\!=\!\ACat^\un_R[\HE_\un^{-1}]^\iy$, $F\!=\!\tTw^\iy$, $\sD\!=\!\ADCat^\un_R[\HE_{\un,\De}^{-1}]^\iy$, and $\al={\rm\ti{\text{\i}}nc}^\iy$. Then (i) holds from above, so by (ii) $\tTw^\iy$ is left adjoint to $\inc:\ADCat^\un_R[\HE_{\un,\De}^{-1}]^\iy\hookra\ACat^\un_R[\HE_\un^{-1}]^\iy$. The analogue for $\Tw^\iy$ follows by equation~\eq{ai9eq41}.
\end{proof}

\subsection{On future applications to Symplectic Geometry}
\label{ai96}

Finally we outline how our theory could be used to define and study Fukaya categories in Symplectic Geometry. Let $(X,\om)$ be a suitable symplectic manifold. Symplectic Geometers would like to construct the {\it derived Fukaya category\/} $D^b\sF(X,\om)$, which should be a triangulated category. These have important applications, e.g.\ Kontsevich's Homological Mirror Symmetry Conjecture \cite{Kont} explains Mirror Symmetry of Calabi--Yau manifolds $X,\check X$ in terms of equivalences of triangulated categories $D^b\sF(X,\om)\simeq D^b\coh\check X$, $D^b\sF(\check X,\check\om)\simeq D^b\coh X$.

See Auroux \cite{Auro} for an introduction to (derived) Fukaya categories, Seidel \cite{Seid} for a definition of $D^b\sF(X,\om)$ under strong restrictions on $(X,\om)$, and Fukaya--Oh--Ohta--Ono \cite{Fuka3,FOOO} for the foundations of a more general construction.

In the world of $\iy$-categories one should replace triangulated categories by stable $\iy$-categories, as in \S\ref{ai27}. Here is how we envisage using our theory to produced derived Fukaya categories $D^b\sF(X,\om)$ as stable $\iy$-categories. We fix a suitable symplectic manifold $(X,\om)$.
\smallskip

\noindent{\bf(a)} Take the po-monoid $(\Ga,\pre)$ to be $(\R,\le)$, and the complete $\R$-filtered ring $R$ to be the Novikov ring $\La^\nov_\Q$ or $\La^\nov_\Z$ from Example \ref{ai5ex3}, where we use $\La^\nov_\Q$ if the relevant counts of $J$-holomorphic curves in $(X,\om)$ lie in $\Q$, and $\La^\nov_\Z$ if they lie in $\Z$ (typically this requires $(X,\om)$ `semipositive'). After choosing some auxiliary extra data $\fD$ on $(X,\om)$, as in Fukaya--Oh--Ohta--Ono \cite{Fuka3,FOOO} we should define the {\it Fukaya category\/} $\sF(X,\om)_\fD$, a unital filtered $A_\iy$-category $\sF(X,\om)_\fD$ over $R$, and so an object of $\AfCat^{\un,\iy}_R$. Objects of $\sF(X,\om)_\fD$ are compact Lagrangians in $(X,\om)$ (maybe equipped with orientations, gradings, \ldots).

The auxiliary data $\fD$ should include a choice of almost complex structure $J$ on $X$ compatible with $\om$ (and possibly domain-dependent). It could also include families of perturbations of Lagrangians as in Seidel \cite{Seid}, choices of Morse functions, \ldots. The details depend on the set-up, and there are several different approaches in the literature.
\smallskip

\noindent{\bf(b)} Having constructed the unital filtered $A_\iy$-category $\sF(X,\om)_\fD$, we define the {\it derived Fukaya category\/} $D^b\sF(X,\om)_\fD$, a stable $\iy$-category to be the composition of $\iy$-functors from \eq{ai8eq8}, \eq{ai9eq16}, \eq{ai9eq33}, and \eq{ai9eq35} applied to~$\sF(X,\om)_\fD$:
\ea
\xymatrix@C=15pt{ \AfCat^{\un,\iy}_R \ar@<-1ex>@/_.6pc/[rrrr]_{\Der_R^\rfi} \ar[r]^{\Str^\iy} & \AfCat^{\un,\iy}_{R,2,2} \ar[r]^(0.45){\Pi_{A_\iy^\rfi}^{A_\iy}} & \ACat^{\un,\iy}_R \ar[r]^{\Tw^\iy} & \ADCat^{\un,\iy}_R \ar[r]^(0.59){\sN^\iy_{A_\iy}} & \sCat^\rst_\iy. }
\nonumber\\[-12pt]
\label{ai9eq43}
\ea
We write $\Der_R^\rfi$ for this composition. So far we have not used $\iy$-categories in an interesting way, we could have defined $D^b\sF(X,\om)_\fD$ using ordinary categories.
\smallskip

\noindent{\bf(c)} We now want to consider how $\sF(X,\om)_\fD$ and $D^b\sF(X,\om)_\fD$ depend on the auxiliary data $\fD$ in {\bf(a)}. To do this, we should define a contractible Kan complex $\sS_{X,\om}$ whose 0-simplices are choices of data $\fD$, and an $\iy$-functor $F:\sS_{X,\om}\ra \AfCat^{\un,\iy}_R$ which maps $F:\fD\mapsto\sF(X,\om)_\fD$ on 0-simplices.

To see why this might be reasonable, consider the (oversimplified) case in which $\fD$ is just a choice of an arbitrary almost complex structure $J$ on $X$ compatible with $\om$. Then we take $\sS_{X,\om}$ to be the simplicial set whose $N$-simplices $\{J_{\bs x}:\bs x\in\De^N_\top\}$ are smooth families of almost complex structures $J_{\bs x}$ on $X$ compatible with $\om$ parametrized by points $\bs x$ in the topological $N$-simplex $\De^N_\top$ of \eq{ai2eq4}, thought of as a smooth manifold with corners. Equivalently, $N$-simplices are smooth maps $\De^N_\top\ra\cJ_{X,\om}$, where $\cJ_{X,\om}$ is the infinite-dimensional manifold of almost complex structures on $X$ compatible with $\om$. Since $\cJ_{X,\om}$ is contractible, $\sS_{X,\om}$ is a contractible Kan complex.

It is important here that $F$ maps to our {\it explicit\/} $\iy$-category $\AfCat^{\un,\iy}_R$, whose $N$-simplices involve collections of maps $\m_k^{I;x_0,\ldots,x_k}$ satisfying \eq{ai8eq4}--\eq{ai8eq5}, and not to a non-explicit localized $\iy$-category such as $\ACat_R^\un[\HE_\un^{-1}]^\iy$. The definition of $F$ will be a generalization of the definition of $\sF(X,\om)_\fD$ in {\bf(a)}, which reduces to this for 0-simplices. For each $N$-simplex in $\sS_{X,\om}$, we arrange to define collections of maps $\m_k^{I;x_0,\ldots,x_k}$ by defining virtual chains for moduli spaces of $J_{\bs x}$-holomorphic curves for $\bs x\in\De^N_\top$, and prove they satisfy equations \eq{ai8eq4}--\eq{ai8eq5}. Guillem Cazassus, Alexander Ritter and the author are working on a construction of this kind~\cite{CJR}.
\smallskip

\noindent{\bf(d)} We have an $\iy$-functor $\Der_R^\rfi\ci F:\sS_{X,\om}\ra\sCat^\rst_\iy$ which maps $\fD\mapsto D^b\sF(X,\om)_\fD$ on 0-simplices. As $\sS_{X,\om}$ is a {\it connected\/} Kan complex, if $\fD_1,\fD_2$ are choices of auxiliary data we have a Joyal equivalence $D^b\sF(X,\om)_{\fD_1}\simeq D^b\sF(X,\om)_{\fD_2}$. As $\sS_{X,\om}$ is a {\it contractible\/} Kan complex, these equivalences are coherent up to a contractible choice in an $\iy$-categorical sense.

Note that the $\iy$-category formalism gives a very clean way to say how $D^b\sF(X,\om)_\fD$ depends on the auxiliary data~$\fD$.
\smallskip

\noindent{\bf(e)} An alternative way to define $D^b\sF(X,\om)_\fD$ for auxiliary data $\fD$ is by the Cartesian diagram of simplicial sets
\begin{equation*}
\xymatrix@C=150pt@R=15pt{ *+[r]{D^b\sF(X,\om)_\fD} \ar[r] \ar[d] &  *+[l]{\NDCat^{\un,\iy}_R} \ar[d]_{p^\De} \\ 
*+[r]{*} \ar[r]^{\Tw^\iy\ci\Pi_{A_\iy^\rfi}^{A_\iy}\ci\Str^\iy\ci F\ci\fD}  &  *+[l]{\ADCat^{\un,\iy}_R.} }
\end{equation*}
Here $p^\De:\NDCat^{\un,\iy}_R\ra\ADCat^{\un,\iy}_R$ is the restriction to $\ADCat^{\un,\iy}_R\subset\ACat^{\un,\iy}_R$ of the cocartesian fibration $p:\NCat^{\un,\iy}_R\ra \ACat^{\un,\iy}_R$ from Definition \ref{ai9def2} and Theorem \ref{ai9thm5}, and $*=\De^0$ is the point simplicial set, and we regard $\fD$ as a morphism $\fD:*\ra\sS_{X,\om}$. The idea here is that by Theorem \ref{ai9thm5}(b), the final morphism $\sN^\iy_{A_\iy}$ in \eq{ai9eq43} maps an object $\sA$ to the fibre of $p^\De$ over~$\sA$.

Define the {\it big derived Fukaya category} $D^b\sF(X,\om)^\rbig$ by the Cartesian diagram of simplicial sets
\e
\begin{gathered}
\xymatrix@C=150pt@R=15pt{ *+[r]{D^b\sF(X,\om)^\rbig} \ar[r] \ar[d]^{\Pi^\rbig} &  *+[l]{\NDCat^{\un,\iy}_R} \ar[d]_{p^\De} \\ *+[r]{\sS_{X,\om}} \ar[r]^{\Tw^\iy\ci\Pi_{A_\iy^\rfi}^{A_\iy}\ci\Str^\iy\ci F}  &  *+[l]{\ADCat^{\un,\iy}_R.} }
\end{gathered}
\label{ai9eq44}
\e
As $p^\De$ is a cocartesian fibration and $\sS_{X,\om}$ is a Kan complex, $D^b\sF(X,\om)^\rbig$ is an $\iy$-category. Since $\sS_{X,\om}$ is contractible, $D^b\sF(X,\om)^\rbig$ is unchanged up to equivalence by replacing $\sS_{X,\om}$ by a point $\fD$ of $\sS_{X,\om}$. That is, for any $\fD$ in $\sS_{X,\om}$, there is an embedding of $\iy$-categories $D^b\sF(X,\om)_\fD\hookra D^b\sF(X,\om)^\rbig$ which is a Joyal equivalence.

Thus $D^b\sF(X,\om)^\rbig$ is an alternative model for the derived Fukaya category of $(X,\om)$ which has the advantage that it is {\it independent of choices\/} $\fD$. It is constructed from not one, but a huge number, of equivalent $A_\iy$-categories.
\smallskip

\noindent{\bf(f)} Here is a (speculative) application of these ideas. Suppose $G$ is a Lie group with an action on $X$ by symplectomorphisms which admits a moment map $\mu:X\ra\g^*$. Define a simplicial group $G^\De$ -- a group in simplicial sets -- such that $N$-simplices of $G^\De$ are smooth maps $S:\De^N_\top\ra G$, regarding $\De^N_\top$ as a smooth manifold with corners. Then we should write down a simplicial action of $G^\De$ on $\sS_{X,\om}$. The $\iy$-functor $\Tw^\iy\ci\Pi_{A_\iy^\rfi}^{A_\iy}\ci\Str^\iy\ci F:\sS_{X,\om}\ra\ADCat^{\un,\iy}_R$ should be invariant under $G^\De$ up to coherent isomorphisms. Thus there is a simplicial action of $G^\De$ on $D^b\sF(X,\om)^\rbig$ such that $\Pi^\rbig$ in \eq{ai9eq44} is $G^\De$-equivariant.

Define the {\it big $G$-equivariant derived Fukaya category\/} $D^b\sF(X,\om)^\rbig_G$ to be the simplicial set of $G^\De$-equivariant sections $s:\sS_{X,\om}\ra D^b\sF(X,\om)^\rbig$ of $\Pi^\rbig$, that is, $\Pi^\rbig\ci s=\id_{\sS_{X,\om}}$. We claim that this is a stable $\iy$-category. We expect that under good conditions, $D^b\sF(X,\om)^\rbig_G$ will roughly be a direct sum of stable $\iy$-subcategories $D^b\sF(X,\om)^\rbig_{G,\eta}$ for each $\eta\in Z(\g^*)\subset\g^*$, which should roughly coincide with the derived Fukaya category of the symplectic quotient $X/\!/_\eta G=\mu^{-1}(\eta)/G$ when this is nonsingular.

Note that this construction would be impossible without $\iy$-categories.

%
% todo Could write homotopy (co)limits section
%
%\subsection{\texorpdfstring{Homotopy limits and colimits of $A_\iy$-categories}{Homotopy limits and colimits of A∞-categories}}
%\label{ai97}
%
%
%{\bf Notes:}
%
%Corollary \ref{ai4cor6}(a), \ref{ai9cor2}.
%
%
%
%Ganatra--Pardon--Shende \cite{GPS1,GPS2}.
%
%Pushouts of $A_\iy$-categories, other colimits, in \cite[\S A.4]{GPS2}. 
%
%\cite[Rem.~A.18]{GPS2}
%
%
%\cite[\S 3.1.3]{GPS1} quotients of $A_\iy$-categories and localization.
%
%Lyubashenko--Manzyuk \cite{LyMa1}
%
%Lyubashenko--Ovsienko \cite{LyOv}.
%
%Canonaco--Ornaghi--Stellari \cite[Lem.~1.28]{COS1} show that not all equalizers exist in $\ACat^\un_R$, for $R$ a field, as an ordinary category.
%
%Oh--Tanaka \cite{OhTa}: prove that quotients and localizations of $A_\iy$-categories satisfy $\iy$-categorical %universal property.
%
%Question: also true in my $\iy$-category without cofibrant assumptions?
%
%
%

\section{Proofs of results in Chapters \ref{ai6}, \ref{ai7} and \ref{ai8}}
\label{ai10}

\subsection{Proofs of Theorems \ref{ai6thm1} and \ref{ai7thm1}}
\label{ai101}

We prove Theorem \ref{ai7thm1}(a), which implies Theorem \ref{ai6thm1}(a) by restricting to $A_\iy$-categories $\sA$ with $\sO=\{*\}$. If $n<k$ then by Definition \ref{ai7def2}(b) both sides of \eq{ai7eq12} are zero. If $n\ge k$ we have
\ea
&\d_{n,k;y_0,\ldots,y_k}^{\{i_0\};x_0,\ldots,x_n}=\sum_{\begin{subarray}{l} 0\le j<k, \\
y_i=x_i,\; 0\le i\le j, \\
y_i=x_{n-k+i},\; j<i\le k
\end{subarray}}\begin{aligned}[t]&
(-1)^{j+(n-k+1)(k-j-1)}\id_{\sH^{i_0}(x_0,x_1)}\ot\cdots\ot\\
&\id_{\sH^{i_0}(x_{j-1},x_j)}\ot\m_{n-k+1}^{\{i_0\};x_j,x_{j+1},\ldots,x_{j+n-k+1}}\ot\\
&\id_{\sH^{i_0}(x_{j+n-k+1},x_{j+n-k+2})}\ot\cdots\ot\id_{\sH^{i_0}(x_{n-1},x_n)}
\end{aligned}
\nonumber
\allowdisplaybreaks\\
&=(-1)^{k''(k'+n)}\cdot{}
\nonumber\\
&\biggl(\sum_{\begin{subarray}{l} 0\le j<k', \\
y_i=x_i,\; 0\le i\le j, \\
y_i=x_{n-k+i},\; j<i\le k
\end{subarray}}
\begin{aligned}[t]&(-1)^{\begin{subarray}{l}j+(k'-j-1)\\((n-k'')-k'+1)\end{subarray}}
\id_{\sH^{i_0}(x_0,x_1)}\ot\cdots\ot\id_{\sH^{i_0}(x_{j-1},x_j)}\\
&\ot\!\m_{(n-k'')-k'+1}^{\{i_0\};x_j,x_{j+1},\ldots,x_{j+n-k+1}}\!\ot\!
\id_{\sH^{i_0}(x_{j+n-k+1},x_{j+n-k+2})}\\
&\ot\cdots\ot\id_{\sH^{i_0}(x_{n-k''-1},x_{n-k''})}
\end{aligned}\biggr)
\nonumber\\
&\ot \id_{\sH^{i_0}(y_{k'},y_{k'+1})}\ot\cdots\ot\id_{\sH^{i_0}(y_{k-1},y_k)}
\nonumber\\
&+(-1)^{k'}\id_{\sH^{i_0}(y_0,y_1)}\ot\cdots\ot\id_{\sH^{i_0}(y_{k'-1},y_{k'})}
\nonumber\\
&\ot\biggl(
\sum_{\begin{subarray}{l} 0\le\ti\jmath<k'', \\
y_i=x_i,\; 0\le i\le k'+\ti\jmath, \\
y_i=x_{n-k+i},\\ 
k'+\ti\jmath<i\le k
\end{subarray}}
\begin{aligned}[t]&(-1)^{\begin{subarray}{l}\ti\jmath+(k''-\ti\jmath-1)\\((n-k')-k''+1)\end{subarray}}
\id_{\sH^{i_0}(x_{k'},x_{k'+1})}\ot\cdots\ot\\
&\id_{\sH^{i_0}(x_{k'+\ti\jmath-1},x_{k'+\ti\jmath})}\ot\m_{(n-k')-k''+1}^{\{i_0\};x_{k'+\ti\jmath},x_{k'+\ti\jmath+1},\ldots,x_{\ti\jmath+n-k''+1}}\\
&\ot\id_{\sH^{i_0}(x_{\ti\jmath+n-k''+1},x_{\ti\jmath+n-k''+2})}\ot\cdots\ot\id_{\sH^{i_0}(x_{n-1},x_n)}
\end{aligned}\biggr)
\nonumber
\allowdisplaybreaks\\
&=\begin{cases}
{\begin{subarray}{l}\ts (-1)^{k''(k'+n)}\d_{n-k'',k';y_0,\ldots,y_{k'}}^{\{i_0\};x_0,\ldots,x_{n-k''}}\ot{} \\ \ts
\id_{\sH^{i_0}(y_{k'},y_{k'+1})}\ot\cdots\ot\id_{\sH^{i_0}(y_{k-1},y_k)}
\end{subarray}}
& \text{if\/ $y_i=x_{n-k+i},$ $k'\le i\le k,$} \\
0, & \text{otherwise}
\end{cases}
\nonumber\\
&+\begin{cases}{\begin{subarray}{l}\ts  (-1)^{k'}\id_{\sH^{i_0}(y_0,y_1)}\ot\cdots\ot{} \\
\ts \id_{\sH^{i_0}(y_{k'-1},y_{k'})}\ot\d_{n-k',k'';y_{k'},\ldots,y_k}^{\{i_0\};x_{k'},\ldots,x_n} \end{subarray}}
& \text{if\/ $y_i=x_i,$ $0\le i\le k',$} \\
0, & \text{otherwise,}
\end{cases}
\label{ai10eq1}
\ea
where in the first and third steps we use \eq{ai7eq8}, and in the second step we split the sum $\sum_{j=0}^{k-1}$ into two sums $\sum_{j=0}^{k'-1}$ and $\sum_{j=k'}^{k'+k''-1}$, replace $j$ by $\ti\jmath=j-k'$ in the second sum, and rewrite the signs.

For Theorems \ref{ai6thm1}(b) and \ref{ai7thm1}(b), for simplicity we will prove equation \eq{ai6eq16} in Theorem \ref{ai6thm1} rather than equation \eq{ai7eq13} in Theorem \ref{ai7thm1}, that is, we will omit objects $x_0,\ldots,x_n,\ab y_0,\ab\ldots,\ab y_k$ from operators $\d_{n,k;y_0,\ldots,y_k}^{I;x_0,\ldots,x_n}$. We could include objects as in equation \eq{ai10eq1}, but this adds notational complexity without changing the important ideas. With $I=\{i_0<\cdots<i_M\}$ for $M>0$, consider the long equation
\ea
&\d_{n,k}^I=
\sum_{\begin{subarray}{l} \text{isomorphism classes} \\ \text{of $\bs\Ga_1,\ldots,\bs\Ga_k$:}\\ \text{$\bs\Ga_i$ is a balanced} \\ 
\text{$I$-labelled tree} \\ \text{with $n_i\ge 0$ leaves,} \\ \text{$1\le i\le k$,} \\
\text{$n=n_1+\cdots+n_k,$} \\ \text{$L(\bs\Ga_1)=I$,} \\ \text{$R(\bs\Ga_k)=\{i_0,i_M\}$,} \\
\text{$R(\bs\Ga_i)=L(\bs\Ga_{i+1})$,} \\ \text{$i=1,\ldots,k-1$}\end{subarray}}
\begin{aligned}[t]
&(-1)^{(k+1)(M+n+1)}\\
&{}\cdot\biggl(\prod_{\begin{subarray}{l} \text{$v$ internal vertex or bud of\/ $\bs\Ga_i,$ $1\le i\le k,$} \\ \text{$w$ internal vertex or bud of\/ $\bs\Ga_{i'},$ $1\le i'\le k,$} \\ 
\text{either $i<i',$ or $i=i'$ and} \\ \text{$v\pr w$ in order $\pr$ of $\Ga_i$,} \\ \text{$v$ has $p$ incoming edges including} \\ \text{$q$ leaves, $w$ has $r$ incoming edges}\end{subarray}\!\!\!\!\!\!\!\!\!\!\!\!\!\!\!\!\!\!\!\!\!\!\!\!\!\!\!\!\!\!\!\!\!\!\!\!\!\!\!} (-1)^{\begin{subarray}{l} N(v,w)+q(r+1) \\ +(p+1)(\md{I_w}+r+1)\end{subarray}}\biggr)\\
&{}\cdot C(\bs\Ga_1)\ot C(\bs\Ga_2)\ot\cdots\ot C(\bs\Ga_k)
\end{aligned}
\nonumber
\allowdisplaybreaks\\
&=\sum_{\begin{subarray}{l}  0<l\le M, \; n=n'+n'',\; n',n''\ge 0, \\
0=j_0<j_1<\cdots<j_l=M \end{subarray}}
\begin{aligned}[t]&(-1)^{(k+1)(M+n+1)}\cdot (-1)^{\sum_{a=1}^{l-1}(j_a-a)} \\
&\;\cdot(-1)^{n'(n''-k'')+(n'-k')(k''-n''+1-l)}\\
&\;\cdot(-1)^{(k''+1)(l+n''+1)}\end{aligned}
\allowdisplaybreaks
\nonumber\\
&{}\cdot\biggl\{\sum_{\begin{subarray}{l} \text{isomorphism classes} \\ \text{of $\bs\Ga_1,\ldots,\bs\Ga_{k'}$:}\\ \text{$\bs\Ga_i$ is a balanced} \\ 
\text{$I$-labelled tree} \\ \text{with $n_i\ge 0$ leaves,} \\ \text{$1\le i\le k'$,} \\
\text{$n'=n_1+\cdots+n_{k'},$} \\ \text{$L(\bs\Ga_1)=I$,} \\ \text{$R(\bs\Ga_{k'})=\{j_0,j_1,\ldots,j_l\}$,} \\
\text{$R(\bs\Ga_i)=L(\bs\Ga_{i+1})$,} \\ \text{$i=1,\ldots,k'-1$}\end{subarray}}
\begin{aligned}[t]
&\biggl(\prod_{\begin{subarray}{l} \text{$v$ internal vertex or bud of\/ $\bs\Ga_i,$ $1\le i\le k',$} \\ \text{$w$ internal vertex or bud of\/ $\bs\Ga_{i'},$ $1\le i'\le k',$} \\ 
\text{either $i<i',$ or $i=i'$ and} \\ \text{$v\pr w$ in order $\pr$ of $\Ga_i$,} \\ \text{$v$ has $p$ incoming edges including} \\ \text{$q$ leaves, $w$ has $r$ incoming edges}\end{subarray}\!\!\!\!\!\!\!\!\!\!\!\!\!\!\!\!\!\!\!\!\!\!\!\!\!\!\!\!\!\!\!\!\!\!\!\!\!\!\!} (-1)^{\begin{subarray}{l} N(v,w)+q(r+1) \\ +(p+1)(\md{I_w}+r+1)\end{subarray}}\biggr)\\
&{}\cdot C(\bs\Ga_1)\ot C(\bs\Ga_2)\ot\cdots\ot C(\bs\Ga_{k'})\biggr\}
\end{aligned}
\allowdisplaybreaks
\nonumber\\
&{}\ot\biggl\{\sum_{\begin{subarray}{l} \text{isomorphism classes} \\ \text{of $\bs\Ga_{k'+1},\ldots,\bs\Ga_k$:}\\ \text{$\bs\Ga_i$ is a balanced} \\ 
\text{$I$-labelled tree} \\ \text{with $n_i\ge 0$ leaves,} \\ \text{$k'<i\le k'$,} \\
\text{$n''=n_{k'+1}+\cdots+n_k,$} \\ \text{$L(\bs\Ga_{k'+1})=\{j_0,j_1,\ldots,j_l\}$,} \\ \text{$R(\bs\Ga_k)=\{i_0,i_M\}$,} \\
\text{$R(\bs\Ga_i)=L(\bs\Ga_{i+1})$,} \\ \text{$i=k'+1,\ldots,k-1$}\end{subarray}}
\begin{aligned}[t]
&(-1)^{(k''+1)(l+n''+1)}\\
&\biggl(\prod_{\begin{subarray}{l} \text{$v$ internal vertex or bud of\/ $\bs\Ga_i,$ $k'<i\le k,$} \\ \text{$w$ internal vertex or bud of\/ $\bs\Ga_{i'},$ $k'<i'\le k,$} \\ 
\text{either $i<i',$ or $i=i'$ and} \\ \text{$v\pr w$ in order $\pr$ of $\Ga_i$,} \\ \text{$v$ has $p$ incoming edges including} \\ \text{$q$ leaves, $w$ has $r$ incoming edges}\end{subarray}\!\!\!\!\!\!\!\!\!\!\!\!\!\!\!\!\!\!\!\!\!\!\!\!\!\!\!\!\!\!\!\!\!\!\!\!\!\!\!} (-1)^{\begin{subarray}{l} N(v,w)+q(r+1) \\ +(p+1)(\md{I_w}+r+1)\end{subarray}}\biggr)\\
&{}\cdot C(\bs\Ga_{k'+1})\ot C(\bs\Ga_{k'+2})\ot\cdots\ot C(\bs\Ga_k)\biggr\}
\end{aligned}
\allowdisplaybreaks
\nonumber\\
&=\sum_{\begin{subarray}{l} 
0<l\le M, \; n=n'+n'',\; n',n''\ge 0,\\
0=j_0<j_1<\cdots<j_l=M \end{subarray}\!\!\!\!\!\!\!\!\!\!\!\!\!\!\!\!\!\!\!\!\!\!\!}\bigl((-1)^{(k'+k''+1)(l+M+n')+k'k''+ln'+n'+\sum_{a=1}^{l-1}(j_a-a)}\bigr)_{(a)}
\nonumber\\
&\biggl\{\sum_{\begin{subarray}{l} \text{isomorphism classes} \\ 
\text{of $\bs\Ga_1,\ldots,\bs\Ga_{k'}$: $\bs\Ga_i$ is a} \\ 
\text{balanced $I$-labelled tree} \\ \text{with $n_i$ leaves,} \\ \text{$i=1,\ldots,k'$,} \\
\text{$n'=n_1+\cdots+n_{k'},$} \\ \text{$L(\bs\Ga_1)=I$,} \\ \text{$R(\bs\Ga_{k'})=\{j_0,j_1,\ldots,j_l\}$,} \\
\text{$R(\bs\Ga_i)=L(\bs\Ga_{i+1})$,} \\ \text{$i=1,\ldots,k'-1$,} \\
\text{trees $\bs\Ga_{i,a}^c$} \\ \text{obtained by} \\ \text{sawing $\bs\Ga_i$ at levels} \\ \text{$j_1,j_2,\ldots,j_{l-1}$.} \\
\text{Write $k_c$ for the} \\ \text{total number of} \\ \text{leaves of trees $\bs\Ga_{i,a}^{c+1}$,} \\ \text{$1\le i\le k'$, $1\le a\le b_i^c$,} \\
\text{$0\le c<l$, and $k_l=k'$.} \\
\text{Then $k_0,\ldots,k_l\ge 0$,} \\
\text{$k_0=n'$, $k_l=k'$,} \\
\text{$k_{a-1}\ge k_a+1-j_a+j_{a-1}$,} \\
\text{$a=1,\ldots,l$}
\end{subarray}}\begin{aligned}[t]
&\biggl(\prod_{\begin{subarray}{l} \text{$v$ internal vertex or bud of\/ $\bs\Ga_i,$ $1\le i\le k',$} \\ \text{$w$ internal vertex or bud of\/ $\bs\Ga_{i'},$ $1\le i'\le k',$} \\ 
\text{either $i<i',$ or $i=i'$ and} \\ \text{$v\pr w$ in order $\pr$ of $\Ga_i$,} \\ \text{$v$ has $p$ incoming edges including} \\ \text{$q$ leaves, $w$ has $r$ incoming edges}\end{subarray}\!\!\!\!\!\!\!\!\!\!\!\!\!\!\!\!\!\!\!\!\!\!\!\!\!\!\!\!\!\!\!\!\!\!\!\!\!\!\!} (-1)^{\begin{subarray}{l} N(v,w)+q(r+1) \\ +(p+1)(\md{I_w}+r+1)\end{subarray}}\biggr)_{(b)}\\
&{}\cdot\biggl(\prod_{\begin{subarray}{l} \text{$v$ internal vertex or bud of\/ $\bs\Ga_i,$ $1\le i\le k',$} \\ \text{$w$ internal vertex or bud of\/ $\bs\Ga_{i'},$ $1\le i'\le k',$} \\ 
\text{either $i<i',$ or $i=i'$ and} \\ \text{$v\pr w$ in order $\pr$ of $\Ga_i$,} \\ \text{$v$ has $p$ incoming edges,} \\ \text{$w$ has $r$ incoming edges,} \\ \text{$\max I_v\le j_c\le\min I_w$, some $1\le c<l$}
\end{subarray}\!\!\!\!\!\!\!\!\!\!\!\!\!\!\!\!\!\!\!\!\!\!\!\!\!\!\!\!\!\!\!\!\!\!\!\!\!\!\!\!\!\!\!\!\!\!\!\!\!\!\!} (-1)^{(3-\md{I_v}-p)(3-\md{I_w}-r)}\biggr)_{(c)}\\
&{}\cdot\Bigl(\bigot_{i=1}^{k'}C(\bs\Ga^l_{i,1})\Bigr)\ci\Bigl(\bigot_{i=1}^{k'}\bigot_{a=1}^{b^{l-1}}C(\bs\Ga^{l-1}_{i,a})\Bigr) \\
&{}\ci\cdots\ci\Bigl(\bigot_{i=1}^{k'}\bigot_{a=1}^{b^2}C(\bs\Ga^2_{i,a})\Bigr)\ci\Bigl(\bigot_{i=1}^{k'}\bigot_{a=1}^{b^1}C(\bs\Ga^1_{i,a})\Bigr)\biggr\}\\
&{}\ot\d_{n'',k''}^{\{i_{j_0},i_{j_1},i_{j_2},\ldots,i_{j_l}\}}\end{aligned}
\allowdisplaybreaks
\nonumber\\
&=\sum_{\begin{subarray}{l} n=n'+n'', \; n',n''\ge 0, \\ 
0<l\le M,\; 0=j_0<j_1<\cdots<j_l=M, \\
k_0,\ldots,k_l\ge 0, \; k_0=n', \; k_l=k', \\
k_{a-1}\ge k_a+1-j_a+j_{a-1}, \; a=1,\ldots,l
\end{subarray}\!\!\!\!\!\!\!\!\!\!\!\!\!\!\!\!\!\!\!\!\!\!\!\!\!\!\!\!\!\!\!\!\!\!\!\!\!\!}\bigl((-1)^{(k'+k''+1)(l+M+n')+k'k''+ln'+n'}
\cdot(-1)^{\sum_{a=1}^{l-1}(j_a-a)}\bigr)_{(d)}
\nonumber
\allowdisplaybreaks\\
&\biggl\{\sum_{\begin{subarray}{l} \text{$\bs\Ga_1,\ldots,\bs\Ga_{k'}$} \\ \text{and $\bs\Ga_{i,a}^c$ as} \\ \text{above, with} \\ \text{the given} \\
\text{values of} \\ \text{$k_0,\ldots,k_l$}\end{subarray}}
\biggl(\prod_{\begin{subarray}{l} \text{$v$ internal vertex or bud of $\bs\Ga_{i,a}^c,$ $1\le i\le k',$ $1\le c\le l$, $1\le a\le b^c$} \\ \text{$w$ internal vertex or bud of $\bs\Ga_{i',a'}^c,$ $1\le i'\le k',$ $1\le a'\le b^c$} \\ 
\text{either $i<i',$ or $i=i'$ and $v\pr w$ in order $\pr$ of $\Ga_i$,} \\ \text{$v$ has $p$ incoming edges including $q$ leaves in $\bs\Ga_{i,a}^c,$} \\ \text{$w$ has $r$ incoming edges}\end{subarray}\!\!\!\!\!\!\!\!\!\!\!\!\!\!\!\!\!\!\!\!\!\!\!\!\!\!\!\!\!\!\!\!\!\!\!\!\!\!\!\!\!\!\!\!\!\!\!\!\!\!\!\!} (-1)^{\begin{subarray}{l} N(v,w)+q(r+1) \\ +(p+1)(\md{I_w}+r+1)\end{subarray}}\biggr)_{(e)}
\allowdisplaybreaks
\nonumber\\
&\qquad\qquad\;\>\>
\biggl(\prod_{\begin{subarray}{l} \text{$v$ internal vertex or bud of $\bs\Ga_{i,a}^c,$ $1\le i\le k',$ $1\le c\le l$, $1\le a\le b^c$} \\ \text{$w$ internal vertex or bud of $\bs\Ga_{i',a'}^{c'},$ $1\le i'\le k',$ $1\le c'<c,$ $1\le a'\le b^{c'}$} \\ 
\text{either $i<i',$ or $i=i'$ and $v\pr w$ in order $\pr$ of $\Ga_i$,} \\ \text{$v$ has $p$ incoming edges, $w$ has $r$ incoming edges}\end{subarray}\!\!\!\!\!\!\!\!\!\!\!\!\!\!\!\!\!\!\!\!\!\!\!\!\!\!\!\!\!\!\!\!\!\!\!\!\!\!\!\!\!\!\!\!\!\!\!\!\!\!\!\!\!\!\!\!\!\!} (-1)^{(p+1)(\md{I_w}+r+1)}\biggr)_{(f)}
\allowdisplaybreaks
\nonumber\\
&\qquad\qquad\;\>\>\biggl(\prod_{\begin{subarray}{l} \text{$v$ internal vertex or bud of $\bs\Ga_{i,a}^c,$ $1\le i\le k',$ $1\le c<l$, $1\le a\le b^c$} \\ \text{$w$ internal vertex or bud of $\bs\Ga_{i',a'}^{c'},$ $1\le i'\le k',$ $c<c'\le l,$ $1\le a'\le b^{c'}$} \\ 
\text{either $i<i',$ or $i=i'$ and $v\pr w$ in order $\pr$ of $\Ga_i$,} \\ \text{$v$ has $p$ incoming edges, $w$ has $r$ incoming edges}\end{subarray}\!\!\!\!\!\!\!\!\!\!\!\!\!\!\!\!\!\!\!\!\!\!\!\!\!\!\!\!\!\!\!\!\!\!\!\!\!\!\!\!\!\!\!\!\!\!\!\!\!\!\!\!\!\!\!\!\!\!\!\!\!\!\!\!\!\!\!\!\!\!\!\!\!\!\!\!\!\!\!\!\!\!\!\!\!\!} (-1)^{\begin{subarray}{l} (\md{I_v}-2)(\md{I_w}-2)+(p-1)(r+1) \\ +(p+1)(\md{I_w}+r+1)\end{subarray}}\biggr)_{(g)}
\allowdisplaybreaks
\nonumber\\
&\qquad\qquad\;\>\>{}\cdot\biggl(\prod_{\begin{subarray}{l} \text{$v$ internal vertex or bud of $\bs\Ga_{i,a}^c,$ $1\le i\le k',$ $1\le c<l$, $1\le a\le b^c$} \\ \text{$w$ internal vertex or bud of $\bs\Ga_{i',a'}^{c'},$ $1\le i'\le k',$ $c<c'\le l,$ $1\le a'\le b^{c'}$} \\ 
\text{either $i<i',$ or $i=i'$ and $v\pr w$ in order $\pr$ of $\Ga_i$,} \\ \text{$v$ has $p$ incoming edges, $w$ has $r$ incoming edges} \end{subarray}\!\!\!\!\!\!\!\!\!\!\!\!\!\!\!\!\!\!\!\!\!\!\!\!\!\!\!\!\!\!\!\!\!\!\!\!\!\!\!\!\!\!\!\!\!\!\!\!\!\!\!\!\!\!\!\!\!\!\!\!\!\!\!\!\!\!\!\!\!} (-1)^{(3-\md{I_v}-p)(3-\md{I_w}-r)}\biggr)_{(h)}
\allowdisplaybreaks
\nonumber\\
&\qquad\qquad\;\>\>{}\cdot\Bigl(\bigot_{i=1}^{k'}C(\bs\Ga^l_{i,1})\Bigr)\ci\Bigl(\bigot_{i=1}^{k'}\bigot_{a=1}^{b^{l-1}}C(\bs\Ga^{l-1}_{i,a})\Bigr) 
\allowdisplaybreaks
\nonumber\\
&\qquad\qquad\;\>\>{}\ci\cdots\ci\Bigl(\bigot_{i=1}^{k'}\bigot_{a=1}^{b^2}C(\bs\Ga^2_{i,a})\Bigr)\ci\Bigl(\bigot_{i=1}^{k'}\bigot_{a=1}^{b^1}C(\bs\Ga^1_{i,a})\Bigr)\biggr\}\allowdisplaybreaks
\nonumber\\
&\qquad\qquad\;\>\>{}\ot\d_{n'',k''}^{\{i_{j_0},i_{j_1},i_{j_2},\ldots,i_{j_l}\}}
\allowdisplaybreaks
\label{ai10eq2}\\
&=\!\!\!\sum_{\begin{subarray}{l} 
n=n'+n'', \; n',n''\ge 0, \\ 
0<l\le M,\; 0=j_0<j_1<\cdots<j_l=M, \\
k_0,\ldots,k_l\ge 0, \; k_0=n', \; k_l=k', \\
k_{a-1}\ge k_a+1-j_a+j_{a-1}, \; a=1,\ldots,l
\end{subarray}}\begin{aligned}[t]&
\bigl((-1)^{k''(k'+l+M)+n'(k''+l+1)}\cdot(-1)^{\sum_{a=1}^{l-1}(j_a-a)}\bigr)_{(i)}\\
&\bigl(\d_{k_{l-1},k_l}^{\{i_{j_{l-1}},\ldots,i_{j_l}\}}\ci\d_{k_{l-2},k_{l-1}}^{\{i_{j_{l-2}},\ldots,i_{j_{l-1}}\}}\ci\cdots\ci\d_{k_0,k_1}^{\{i_{j_0},\ldots,i_{j_1}\}}\bigr)\\
&{}\ot\d_{n'',k''}^{\{i_{j_0},i_{j_1},i_{j_2},\ldots,i_{j_l}\}}.
\end{aligned}
\nonumber
\ea

Here the first step is \eq{ai6eq11}. In the second step we split the sum over $\bs\Ga_1,\ldots,\bs\Ga_k$ into separate sums over $\bs\Ga_1,\ldots,\bs\Ga_{k'}$ and $\bs\Ga_{k'+1},\ldots,\bs\Ga_k$, and rearrange the signs. 

For the initial signs, the second sign in the second step comes from
\ea
&\prod_{\begin{subarray}{l} \text{$v$ internal vertex or bud of $\bs\Ga_i,$ $1\le i\le k',$} \\ \text{$w$ internal vertex or bud of\/ $\bs\Ga_{i'},$ $k'<i'\le k,$}\end{subarray}\!\!\!\!\!\!\!\!\!\!\!\!\!\!\!\!\!\!\!\!\!\!\!\!\!\!\!\!\!\!\!\!\!\!\!\!\!\!\!\!\!}(-1)^{N(v,w)}
=\sum_{\begin{subarray}{l} r\in \{0,\ldots,M\}\sm\{j_0,\ldots,j_l\}, \\ s\in\{j_0,\ldots,j_l\}\sm\{0,M\}:\; r<s\end{subarray}\!\!\!\!\!\!\!\!\!\!\!\!\!\!\!\!\!\!\!\!\!\!\!\!\!\!\!\!\!\!\!\!\!\!\!\!\!\!\!\!\!\!\!\!}(-1)=(-1)^{\sum_{a=1}^{l-1}(j_a-a)}.
\label{ai10eq3}
\ea
This holds as by Lemma \ref{ai6lem1}, the possible values of $r\in\rint(I_v)$ for $v$ an internal vertex or bud of $\bs\Ga_i,$ $1\le i\le k'$ are $\{0,\ldots,M\}\sm\{j_0,\ldots,j_l\}$, and each value of $r$ in $\{0,\ldots,M\}\sm\{j_0,\ldots,j_l\}$ occurs for exactly one such $v$. Similarly, the possible values of $s\in\rint(I_w)$ for $w$ an internal vertex or bud of $\bs\Ga_{i'},$ $k'<i'\le k,$ are $\{j_0,\ldots,j_l\}\sm\{0,M\}$, and each value of $s$ occurs for exactly one such $w$. The first step of \eq{ai10eq3} follows using \eq{ai6eq12}, and the second step is a combinatorial calculation.

The third sign in the second step comes from
\begin{align*}
\prod_{\begin{subarray}{l} \text{$v$ internal vertex or bud of\/ $\bs\Ga_i,$ $1\le i\le k',$} \\ \text{$w$ internal vertex or bud of\/ $\bs\Ga_{i'},$ $k'<i'\le k,$} \\ 
\text{$v$ has $p$ incoming edges including} \\ \text{$q$ leaves, $w$ has $r$ incoming edges}\end{subarray}\!\!\!\!\!\!\!\!\!\!\!\!\!\!\!\!\!\!\!\!\!\!\!\!\!\!\!\!\!\!\!\!\!\!\!\!\!\!\!\!\!\!\!\!\!\!\!\!\!\!\!\!\!\!\!\!\!\!\!\!\!\!\!\!} (-1)^{q(r+1)+(p+1)(\md{I_w}+r+1)}
=(-1)^{n'(n''-k'')+(n'-k')(k''-n''+1-l)},
\end{align*}
which holds as
\begin{align*}
\sum_{\begin{subarray}{l} \text{$v$ internal vertex or bud of $\bs\Ga_i,$} \\ \text{$1\le i\le k',$ $v$ has $q$ leaves}\end{subarray}\!\!\!\!\!\!\!\!}q&=n', & \sum_{\begin{subarray}{l} \text{$v$ internal vertex or bud of $\bs\Ga_i,$} \\ \text{$1\le i\le k',$ $v$ has $p$ incoming edges}\end{subarray}\!\!\!\!\!\!\!\!\!\!\!\!\!\!\!\!\!\!\!\!\!\!\!\!\!\!\!\!\!\!\!\!\!\!\!\!\!\!\!\!\!\!\!\!\!\!\!}(p-1)&=n'-k',\\
\sum_{\begin{subarray}{l} \text{$w$ internal vertex or bud of $\bs\Ga_{i'},$} \\ \text{$k'<i'\le k,$} \\
\text{$w$ has $r$ incoming edges}\end{subarray}\!\!\!\!\!\!\!\!\!\!\!\!\!\!\!\!\!\!\!\!\!\!\!\!\!\!\!\!\!\!\!\!\!\!\!\!\!\!\!\!\!\!\!\!\!\!\!}(r-1)&=n''\!-\!k'', & \sum_{\begin{subarray}{l} \text{$w$ internal vertex or bud of $\bs\Ga_{i'},$} \\ \text{$k'<i'\le k,$} \\
\text{$w$ has $r$ incoming edges}\end{subarray}\!\!\!\!\!\!\!\!\!\!\!\!\!\!\!\!\!\!\!\!\!\!\!\!\!\!\!\!\!\!\!\!\!\!\!\!\!\!\!\!\!\!\!}(3\!-\!\md{I_w}\!-\!r)&=k''\!-\!n''\!+\!1\!-\!l.
\end{align*}
The fourth sign in the second step cancels the same sign three lines later.

In the third step of \eq{ai10eq2}, for the sum over $\bs\Ga_1,\ldots,\bs\Ga_{k'}$, we saw the $\bs\Ga_i$ at levels $j_1,j_2,\ldots,j_{l-1}$, in this order, using \eq{ai6eq7} inductively. The notation is that sawing $\bs\Ga_i$ at level $j_1$ gives trees $\bs\Ga^1_{i,a}$ for $a=0,1,\ldots,b^1$, where $\bs\Ga^1_{i,0}$ is a $\{j_1,j_1+1,\ldots,j_l=M\}$-labelled tree and $\bs\Ga^1_{i,a}$ is a $\{j_0,\ldots,j_1\}$-labelled tree for $a=1,2,\ldots,b_i^1$. Then sawing $\bs\Ga_{i,0}^1$ at level $j_2$ gives trees $\bs\Ga^2_{i,a}$ for $a=0,1,\ldots,b_i^2$, where $\bs\Ga^2_{i,0}$ is a $\{j_2,j_2+1,\ldots,j_l\}$-labelled tree and $\bs\Ga^2_{i,a}$ is a $\{j_1,\ldots,j_2\}$-labelled tree for $a=1,2,\ldots,b_i^2$. Inductively, sawing $\bs\Ga_{i,0}^{c-1}$ at level $j_c$ gives trees $\bs\Ga^c_{i,a}$ for $a=0,1,\ldots,b_i^c$, where $\bs\Ga^c_{i,0}$ is a $\{j_c,j_c+1,\ldots,j_l\}$-labelled tree and $\bs\Ga^c_{i,a}$ is a $\{j_{c-1},\ldots,j_c\}$-labelled tree for $a=1,2,\ldots,b_i^c$, for $c=2,3,\ldots,l-1$. We also write $b_i^l=1$ and $\bs\Ga^l_{i,1}=\bs\Ga_{i,0}^{l-1}$.

The extra sign factor with a product $(-1)^{(3-\md{I_v}-p)(\md{I_w}-r)}$ comes from two sources: the factors from $v,w$ with $i=i'$ and $v\pr w$ in the order $\pr$ of $\Ga_i$ come from the sign in \eq{ai6eq7} from cutting $\bs\Ga_i$ at level $j_c$ for some $1\le c<l$. The factors from $v,w$ with $i<i'$ come from an additional rearrangement of the order of the factors $C(\bs\Ga_{i,a}^c)$, in which we order them first by $c=l,l-1,\ldots,1$ and then by $i=1,\ldots,k'$, rather than ordering them first by $i=1,\ldots,k'$ and then by $c=l,l-1,\ldots,1$. Also, we use \eq{ai6eq11} for the sum over $\bs\Ga_{k'+1},\ldots,\bs\Ga_k'$, and we simplify the initial signs.

In the third step we also define numbers $k_0,\ldots,k_l\ge 0$ with $k_0=n'$ and $k_l=k'$, associated to each choice of $\bs\Ga_1,\ldots,\bs\Ga_{k'}$, such that $k_c$ is the total number of leaves of trees $\bs\Ga_{i,a}^{c+1}$ over all $i=1,\ldots,k'$ and $a=1,\ldots,b_i^{c+1}$, or equivalently, the total number of roots of trees $\bs\Ga_{i,a}^c$ over all $i=1,\ldots,k'$ and $a=1,\ldots,b_i^c$, although these are not used in the third step. The inequality $k_{a-1}\ge k_a+1-j_a+j_{a-1}$ holds as in the proof of Lemma \ref{ai6lem2}: if $k_a>k_{a-1}$ then there are at least $k_a-k_{a-1}$ trees $\bs\Ga_{i,a}^c$ with no leaves, and each must contain a bud $b$ with $\md{I_b}\ge 3$, so $\md{\rint(I_b)}\ge 1$, and the disjoint union of $\rint(I_v)$ over internal vertices and buds $v$ of trees $\bs\Ga_{i,a}^c$ is $\{i_{j_{a-1}},\ldots,i_{j_a}\}\sm\{i_{j_{a-1}},i_{j_a}\}$, which has size $j_a-j_{a-1}-1$, so that $k_a-k_{a-1}\le j_a-j_{a-1}-1$.

In the fourth step of \eq{ai10eq2} we rewrite by first fixing values of $k_0,\ldots,k_l$ in the first sum in the third step, and then summing over $\bs\Ga_1,\ldots,\bs\Ga_{k'}$ and $\bs\Ga^c_{i,a}$ with these values of $k_0,\ldots,k_l$. We also rewrite the signs. To help explain this we have labelled the signs in the third--fifth steps (a)--(i). Then in the fourth step:
\begin{itemize}
\setlength{\itemsep}{0pt}
\setlength{\parsep}{0pt}
\item[(i)] Signs $(e),(f),(g)$ are $(b)$ split into three factors, with $(b)=(e)\cdot(f)\cdot(g)$. Roughly, what we do is to regard $v$ as a vertex of some $\bs\Ga_{i,a}^c$ rather than $\bs\Ga_i$, and $w$ as a vertex of some $\bs\Ga_{i',a'}^{c'}$, and then $(e)$ is the terms in $(b)$ with $c'=c$, and $(f)$ the terms in $(b)$ with $c'<c$, and $(g)$ the terms in $(b)$ with $c'>c$. But there are some subtleties, see (ii),(iii).
\item[(ii)] Recall the definition of $N(v,w)$ in \eq{ai6eq12}. In $(f)$ we have $I_v\!\subseteq\!\{j_{c-1},\ldots,j_c\}$ and $I_w\subseteq\{j_{c'-1},\ldots,j_{c'}\}$ with $c'<c$, so $N(v,w)=0$ and we omit the $N(v,w)$ term from $(f)$. In $(g)$ we have $I_v\subseteq\{j_{c-1},\ldots,j_c\}$ and $I_w\subseteq\{j_{c'-1},\ldots,j_{c'}\}$ with $c<c'$, so $N(v,w)=(\md{I_v}-2)(\md{I_w}-2)$. 
\item[(iii)] We have rewritten the terms in the number of leaves $q$ in $(b)$ and~$(e)$--$(g)$. 

In $(b)$, $q$ is the number of leaves of $\bs\Ga_i$ at $v$. But in $(e)$, $q$ is the number of leaves of $\bs\Ga_{i,a}^c$ at $v$; if $c=1$ these are leaves of $\bs\Ga_i$, but if $c>1$ they are just internal edges of $\bs\Ga_i$ at $v$ labelled by $j_{c-1}$.

Fix a vertex $v$ in some $\Ga_{i,a}^{c'}$ with $q$ leaves in $\Ga_{i,a}^{c'}$. We claim that
\ea
\sum_{\begin{subarray}{l}
\text{$\dot v$ internal vertex or bud in $\Ga_{i,\dot a}^1$:} \\
\text{$\dot v\ne v$, $\dot v$ has $\dot q$ leaves in $\Ga$,} \\ \text{there is a directed path} \\ \text{$\overset{\dot v}{\bu}\longra\cdots\longra\,\overset{v}{\bu}$ in $\Ga_i$} \\ \text{ending in a leaf of $v$ in $\Ga_{i,a}^{c'}$} \end{subarray}\!\!\!\!\!\!\!\!}\dot q
=q+\sum_{\begin{subarray}{l}
\text{$\ddot v$ internal vertex or bud in $\Ga_{i,\ddot a}^{\ddot c}$:} \\
\text{$\ddot v\ne v$, $\ddot v$ has $\ddot p$ incoming edges,} \\ \text{there is a directed path} \\ \text{$\overset{\ddot v}{\bu}\longra\cdots\longra\,\overset{v}{\bu}$ in $\Ga_i$} \\ \text{ending in a leaf of $v$ in $\Ga_{i,a}^{c'}$} \end{subarray}\!\!\!\!\!\!\!\!\!\!\!\!\!\!\!\!\!\!\!}(\ddot p-1).
\label{ai10eq4}
\ea
To see this, note that the subgraph of $\Ga_i$ with vertices $\ddot v$ as in the right hand side of \eq{ai10eq4} is the disjoint union of $q$ trees, each with a root corresponding to a leaf of $v$ in $\Ga_{i,a}^c$. For any tree $\Ga$
\begin{equation*}
(\text{number of leaves of $\Ga$})=1+\sum_{\text{$v$ internal vertex or bud of $\Ga$}\!\!\!\!\!\!\!\!\!\!\!\!\!\!\!\!\!\!\!\!\!\!\!\!\!\!\!\!\!\!\!\!\!\!\!\!\!\!\!\!\!\!\!\!\!\!\!\!\!\!\!\!\!\!}(\text{number of incoming edges at $v$}\,-1), 
\end{equation*}
and summing this over the $q$ trees above gives \eq{ai10eq4}.

Next fix a vertex $w$ of some $\bs\Ga_{i',a'}^{c'}$ with $r$ incoming edges. We claim that
\ea
&\biggl(\prod_{\begin{subarray}{l} \text{$\dot v$ internal vertex or bud of $\bs\Ga_{i,\dot a}^1,$ $1\le i\le k',$ $1\le\dot a\le b_i^1$} \\  
\text{either $i<i',$ or $i=i'$ and $v\pr w$ in order $\pr$ of $\Ga_i$,} \\ \text{$\dot v$ has $\dot q$ leaves in $\Ga_i$}\end{subarray}\!\!\!\!\!\!\!\!\!\!\!\!\!\!\!\!\!\!\!\!\!\!\!\!\!\!\!\!\!\!\!} (-1)^{\dot q(r+1)}\biggr)
\label{ai10eq5}\\
&=\biggl(\prod_{\begin{subarray}{l} \text{$v$ internal vertex or bud of $\bs\Ga_{i,a}^{c'},$ $1\le i\le k',$ $1\le a\le b^{c'}$} \\  
\text{either $i<i',$ or $i=i'$ and $\dot v\pr w$ in order $\pr$ of $\Ga_i$,} \\ \text{$v$ has $q$ leaves in $\Ga_{i,a}^{c'}$}\end{subarray}\!\!\!\!\!\!\!\!\!\!\!\!\!\!\!\!\!\!\!\!\!\!\!\!\!\!\!\!\!\!\!\!\!\!} (-1)^{q(r+1)}\biggr)
\nonumber\\
&\cdot\biggl(\prod_{\begin{subarray}{l} \text{$\ddot v$ internal vertex or bud of $\bs\Ga_{i,\ddot a}^{\ddot c},$ $1\le\ddot c<c'$, $1\le i\le k',$ $1\le\ddot a\le b_i^{\ddot c}$} \\  
\text{either $i<i',$ or $i=i'$ and $\ddot v\pr w$ in order $\pr$ of $\Ga_i$,} \\ \text{$\ddot v$ has $\ddot p$ incoming edges}\end{subarray}\!\!\!\!\!\!\!\!\!\!\!\!\!\!\!\!\!\!\!\!\!\!\!\!\!\!\!\!\!\!\!\!\!\!} (-1)^{(\ddot p-1)(r+1)}\biggr).
\nonumber
\ea
This follows by raising \eq{ai10eq4} to the power $(-1)^{r+1}$, and taking the product over all $v$ satisfying the conditions in the second line.

Now take the product of \eq{ai10eq5} over all $w$ and $\bs\Ga_{i',a'}^{c'}$. The left hand side gives the terms in $(-1)^{q(r+1)}$ in $(b)$. The first term on the right hand side gives the terms in $(-1)^{q(r+1)}$ in $(e)$. The second term on the right hand side gives the terms in $(-1)^{(p-1)(r+1)}$ in $(g)$, explaining these terms.

For $(f)$, as $1\le c'<c$ and $v$ can only have leaves of $\Ga$ if $c=1$ there is no contribution from the $(-1)^{q(r+1)}$ term in $(b)$.
\item[(iv)] Sign $(h)$ is just $(c)$  written in terms of $\bs\Ga_{i,a}^c,\bs\Ga_{i',a'}^{c'}$ rather than $\bs\Ga_i,\bs\Ga_{i'}$.
\end{itemize}

In the fifth step of \eq{ai10eq2} we observe that for each $c=1,2,\ldots,l$ the sets of trees $\bs\Ga_{i,a}^c$, in the orders in which they appear in the tensor products in the fourth step, actually depend only on the data $j_{c-1},j_c,k_{c-1},k_c$ in the first sum in the fourth step: we need sets of trees $\bs\Ga_{i,a}^c$ such that there are $k_c$ trees, and $k_{c-1}$ total leaves, and the right and left sequences $R(\bs\Ga_{i,a}^c),L(\bs\Ga_{i,a}^c)$ satisfy the appropriate matching conditions. Furthermore, these are exactly the sets of trees, orders, and matching conditions appearing in the definition of $\d_{k_{c-1},k_c}^{\{i_{j_{c-1}},\ldots,i_{j_c}\}}$ in~\eq{ai6eq11}.

We will also show that after some complicated cancellations, the product of signs in the fourth step equals the product of the sign in the fifth step, with the signs in the definitions of $\d_{k_{l-1},k_l}^{\{i_{j_{l-1}},\ldots,i_{j_l}\}},\ab\ldots,\ab\d_{k_0,k_1}^{\{i_{j_0},\ldots,i_{j_1}\}}$ in \eq{ai6eq11}. That is:
\begin{itemize}
\setlength{\itemsep}{0pt}
\setlength{\parsep}{0pt}
\item[(v)] As $(g),(h)$ in the fourth step are sums over the same set of $v,w,p,r$, and
\begin{equation*}
\!\!\!\!\!\!\!\!\!\!\!(-1)^{\begin{subarray}{l} (\md{I_v}-2)(\md{I_w}-2)+(p-1)(r+1) \\ +(p+1)(\md{I_w}+r+1)\end{subarray}}\cdot(-1)^{(3-\md{I_v}-p)(3-\md{I_w}-r)}=(-1)^{(r+1)(\md{I_v}+p+1)},
\end{equation*}
we see that
\e
(g)\cdot(h)=\prod_{\begin{subarray}{l} \text{$v$ internal vertex or bud of $\bs\Ga_{i,a}^c,$ $1\le i\le k',$ $1\le c<l$, $1\le a\le b^c$} \\ \text{$w$ internal vertex or bud of $\bs\Ga_{i',a'}^{c'},$ $1\le i'\le k',$ $c<c'\le l,$ $1\le a'\le b^{c'}$} \\ 
\text{either $i<i',$ or $i=i'$ and $v\pr w$ in order $\pr$ of $\Ga_i$,} \\ \text{$v$ has $p$ incoming edges, $w$ has $r$ incoming edges} \end{subarray}\!\!\!\!\!\!\!\!\!\!\!\!\!\!\!\!\!\!\!\!\!\!\!\!\!\!\!\!\!\!\!\!\!\!\!\!\!\!\!\!\!\!\!\!\!\!\!\!\!\!\!\!\!\!} (-1)^{(r+1)(\md{I_v}+p+1)}.
\label{ai10eq6}
\e
Let us exchange $v,a,i,c,p$ with $w,a',i',c',r$ in \eq{ai10eq6}, respectively, and compare the result with $(f)$ in the fourth step of \eq{ai10eq2}. The differences are that $(f)$ has the condition that either $i<i'$, or $i=i'$ and $v\pr w$, but the exchanged version of \eq{ai10eq6} has the condition that either $i'<i$, or $i'=i$ and $w\pr v$. Also both $(f)$ and the exchanged version of \eq{ai10eq6} have $c'<c$, which excludes $v=w$. Thus we deduce the first step in
\ea
&(f)\cdot(g)\cdot(h)=\prod_{\!\!\!\!\!\!\!\begin{subarray}{l} \text{$v$ internal vertex or bud of $\bs\Ga_{i,a}^c,$ $1\le i\le k',$ $1\le c\le l$, $1\le a\le b^c$} \\ \text{$w$ internal vertex or bud of $\bs\Ga_{i',a'}^{c'},$ $1\le i'\le k',$ $1\le c'<c,$ $1\le a'\le b^{c'}$} \\ 
\text{$v$ has $p$ incoming edges, $w$ has $r$ incoming edges}\end{subarray}\!\!\!\!\!\!\!\!\!\!\!\!\!\!\!\!\!\!\!\!\!\!\!\!\!\!\!\!\!\!\!\!\!\!\!\!\!\!\!\!\!\!\!\!\!\!\!\!\!\!\!\!\!\!\!\!\!\!\!} (-1)^{(p+1)(\md{I_w}+r+1)}
\nonumber\\
&=\prod_{1\le c'<c\le l}(-1)^{(k_{c-1}-k_c)(j_{c'}-j_{c'-1}+k_{c'-1}-k_{c'}+1)}
\nonumber\\
&=\prod_{c'=1}^l(-1)^{(k_{c'}-k')(j_{c'}-j_{c'-1}+k_{c'-1}-k_{c'}+1)}
\nonumber\\
&=\Bigl(\prod_{c=1}^l(-1)^{(k_c+1)(j_c-j_{c-1}+k_{c-1}+1)}\Bigr)
\cdot\Bigl(\prod_{c=1}^l(-1)^{(k'+1)(j_c-j_{c-1}+k_{c-1}-k_c+1)}\Bigr)
\nonumber\\
&=\Bigl(\prod_{c=1}^l(-1)^{(k_c+1)(j_c-j_{c-1}+k_{c-1}+1)}\Bigr)\cdot(-1)^{(k'+1)(M+n'-k'+l)}.\nonumber\\
&=\Bigl(\prod_{c=1}^l(-1)^{(k_c+1)(j_c-j_{c-1}+k_{c-1}+1)}\Bigr)\cdot
(-1)^{(k'+k''+1)(l+M+n')+k'k''+ln'+n'}
\nonumber\\
&\qquad \cdot(-1)^{k''(k'+l+M)+n'(k''+l+1)}.
\label{ai10eq7}
\ea
Here in the second step of \eq{ai10eq7} we use that
\begin{align*}
\prod_{\begin{subarray}{l} \text{$v$ internal vertex or bud of $\bs\Ga_{i,a}^c,$} \\ \text{$1\le i\le k',$ $1\le a\le b^c$} \\ \text{$v$ has $p$ incoming edges}\end{subarray}}(-1)^{p+1}&=(-1)^{k_{c-1}-k_c},\\
\prod_{\begin{subarray}{l} \text{$w$ internal vertex or bud of $\bs\Ga_{i',a'}^{c'},$ $1\le i'\le k',$  $1\le a'\le b^{c'}$} \\ 
\text{$w$ has $r$ incoming edges}\end{subarray}\!\!\!\!\!\!\!\!\!\!\!\!\!\!\!\!\!\!\!\!\!\!\!\!\!\!\!\!\!\!\!\!\!\!\!\!\!\!\!\!\!\!\!\!\!\!\!\!\!\!\!\!} (-1)^{\md{I_w}+r+1}&=(-1)^{j_{c'}-j_{c'-1}+k_{c'-1}-k_{c'}+1},
\end{align*}
and in the third that $\sum_{c=c'+1}^l(k_{c-1}-k_c)=k_{c'}-k'$ as $k_l=k'$.
\item[(vi)] The first sign in the last step of \eq{ai10eq7} is the product over $c=1,\ldots,l$ of the first signs in \eq{ai6eq11} of $\d_{k_{c-1},k_c}^{\{i_{j_{c-1}},\ldots,i_{j_c}\}}$. Sign $(e)$ is the product over $c=1,\ldots,l$ of the second signs in \eq{ai6eq11} of $\d_{k_{c-1},k_c}^{\{i_{j_{c-1}},\ldots,i_{j_c}\}}$. The second sign in the last step of \eq{ai10eq7} is the first sign in $(d)$. The third sign in the last step of \eq{ai10eq7} is the first sign in $(i)$. The second signs in $(d)$ and $(i)$ agree. 
\end{itemize}
This proves equation \eq{ai10eq2}, and hence \eq{ai6eq16}. The claims on finiteness after \eq{ai6eq16} are straightforward. This completes Theorems \ref{ai6thm1} and~\ref{ai7thm1}.

\subsection{Proofs of Propositions \ref{ai6prop1} and \ref{ai7prop1}}
\label{ai102}

We prove Proposition \ref{ai7prop1}, which implies Proposition \ref{ai6prop1}. It is easy to see that $\bigl((\bar \sO^i)_{i=0}^{\bar N},(\bar  f^{ij})_{0\le i\le j\le\ti N},(\bar\sH{}^i)_{i=0}^{\bar N},\ab(\bar\m{}^{\bar I})_{\es\ne\bar I\subseteq\{0,\ldots,\bar N\}}\bigr)$ satisfies all of Definition \ref{ai7def3}(a)(i)--(v), except perhaps equations \eq{ai7eq15}--\eq{ai7eq16}. So we must show that the $\bar\m{}^{\bar I;x_0,\ldots,x_k}_k$ satisfy \eq{ai7eq15}--\eq{ai7eq16} for all $\es\ne \bar I\subseteq\{0,\ldots,\bar N\}$, which we write as $\bar I=\{\bar\imath_0,\ab \bar\imath_1,\ab\ldots,\ab \bar\imath_{\bar M}\}$ for $\bar\imath_0<\bar\imath_1<\cdots<\bar\imath_{\bar M}$ with $\bar M\ge 0$, and $n\ge 0$, $x_0,\ldots,x_n\in\bar\sO^{\bar\imath_0}$. Write $\bar\d_{n,k;y_0,\ldots,y_k}^{\bar I;x_0,\ldots,x_n}$ for the operators $\d_{n,k;y_0,\ldots,y_k}^{I;x_0,\ldots,x_n}$ in Definition \ref{ai7def2} defined using $\bigl((\bar \sO^i)_{\cdots},(\bar f^{ij})_{\cdots},(\bar\sH{}^i)_{\cdots},\ab(\bar\m{}^{\bar I})_{\cdots}\bigr)$. 

Divide into the following cases:
\begin{itemize}
\setlength{\itemsep}{0pt}
\setlength{\parsep}{0pt}
\item[(a)] $\al\vert_{\bar I}$ is injective, so that $\md{\al(\bar I)}=\bar M+1$.
\item[(b)] $\bar M=1$ and $\al(\bar\imath_0)=\al(\bar\imath_1)$.
\item[(c)] $\bar M\ge 2$ and $\md{\al(\bar I)}=\bar M$, and $\al(\bar\imath_{j-1})=\al(\bar\imath_j)$ for $1<j<M$.
\item[(d)] $\bar M\ge 2$ and $\md{\al(\bar I)}=\bar M$, and $\al(\bar\imath_0)=\al(\bar\imath_1)$.
\item[(e)] $\bar M\ge 2$ and $\md{\al(\bar I)}=\bar M$, and $\al(\bar\imath_{\bar M-1})=\al(\bar\imath_{\bar M})$.
\item[(f)] $\bar M=3$ and $\md{\al(\bar I)}=1$ or 2.
\item[(g)] $\bar M\ge 4$ and $\md{\al(\bar I)}<\bar M$.
\end{itemize}

In case (a), by \eq{ai7eq20} we have $\bar\m_k^{\bar J;x_0,\ldots,x_k}=\m_k^{\al(\bar J);x_0,\ldots,x_k}$ for all $\es\ne \bar J\subseteq \bar I$, so \eq{ai7eq15}--\eq{ai7eq16} for $\bigl((\bar \sO^i)_{\cdots},(\bar f^{ij})_{\cdots},(\bar\sH{}^i)_{\cdots},\ab(\bar\m{}^{\bar I})_{\cdots}\bigr)$ and $\bar I$ are equivalent to \eq{ai7eq15}--\eq{ai7eq16} for $\bigl((\sO^i)_{\cdots},\ab(f^{ij})_{\cdots},(\sH^i)_{\cdots},(\m^I)_{\cdots}\bigr)$ and~$I=\al(\bar I)$.

In case (b), from \eq{ai7eq9} and \eq{ai7eq20} we can show that 
\begin{equation*}
\bar\d_{n,n;y_0,\ldots,y_n}^{\{\bar\imath_0,\bar\imath_1\};x_0,\ldots,x_n}=\begin{cases}
\id, & x_i=y_i,\; i=0,\ldots,n, \\ 0, & \text{otherwise,}
\end{cases}
\end{equation*}
and $\bar\d_{n,k;y_0,\ldots,y_k}^{\{\bar\imath_0,\bar\imath_1\};x_0,\ldots,x_n}=0$ for $k<n$. Hence by \eq{ai7eq10} and \eq{ai7eq20},
equation \eq{ai7eq16} for $\bigl((\bar \sO^i)_{\cdots}$, $(\bar f^{ij})_{\cdots}$, $(\bar\sH{}^i)_{\cdots}$, $(\bar\m{}^{\bar I})_{\cdots}\bigr)$ holds as
\begin{align*}
&\sum_{y_0,y_1\in\bar\sO^0}\bar\m_1^{\{\bar\imath_0,\bar\imath_1\};y_0,y_1}\ci\bar\d_{n,1;y_0,y_1}^{\{\bar\imath_0\};x_0,\ldots,x_n}-\sum_{y_0,\ldots,y_n\in\bar\sO^1}\bar\m_n^{\{\bar\imath_1\};y_0,\ldots,y_n}\ci\bar\d_{n,n;y_0,\ldots,y_n}^{\{\bar\imath_0,\bar\imath_1\};x_0,\ldots,x_n}\\
&=\id\ci\m_n^{\{\al(\bar\imath_0)\};x_0,\ldots,x_n}-\m_n^{\{\al(\bar\imath_0)\};x_0,\ldots,x_n}\ci\id=0.
\end{align*}

For cases (c)--(g), we begin with some discussion. As $\al\vert_{\bar I}$ is not injective, there exists at least one $j=1,\ldots,M$ with $\al(\bar\imath_{j-1})=\al(\bar\imath_j)$. The first term $\sum_{k,y_i}\bar\m_k^{\bar I;y_0,\ldots,y_k}\ci\bar\d_{n,k;y_0,\ldots,y_k}^{\{\bar\imath_0\};x_0,\ldots,x_n}$ in \eq{ai7eq16} for $\bar I$ is zero as $\bar\m_k^{\bar I;y_0,\ldots,y_k}=0$ for all $k\ge 1$ by \eq{ai7eq20}. In the second term $(-1)^{\bar M}\sum_{k,y_i}\bar\m_k^{\{\bar\imath_{\bar M}\};y_0,\ldots,y_k}\ci\bar\d_{n,k;y_0,\ldots,y_k}^{\bar I;x_0,\ldots,x_n}$, equation \eq{ai7eq9} defines $\bar\d_{n,k;y_0,\ldots,y_k}^{\bar I;x_0,\ldots,x_n}$ as a sum over sets of labelled trees $\bs\Ga_1,\ldots,\bs\Ga_k$. By considering the least $i=1,\dots,k$ such that one of $\bar\imath_{j-1}$ or $\bar\imath_j$ does not lie in $R(\bs\Ga_i)$ we can show that $\bs\Ga_i$ has a vertex $v$ with $\bar\imath_{j-1},\bar\imath_j\in I_v$ and $\md{I_v}\ge 3$. Then $\al\vert_{I_v}$ is not injective, so \eq{ai7eq20} shows that $\bar\m_k^{I_v;y_0,\ldots,y_k}=0$ unless $k=0$, so that $v$ is a bud, and $\md{I_v}=3$, and then $\bar\m_0^{I_v;y_0}=\boo^{\bar\jmath_2}_{\bar f{}^{\bar\jmath_0\bar\jmath_2}(y_0)}$ for~$I_v=\{\ti\jmath_0<\ti\jmath_1<\ti\jmath_2\}$.

Suppose $v$ lies on an edge $\overset{v}{\sq}\,{\buildrel e\over\longra}\,\overset{w}{\bu}$ in $\bs\Ga_i$, with $w$ an internal vertex. Then in $C(\bs\Ga_i)^{x_{l_{i-1}},\ldots,x_{l_i}}$ in \eq{ai7eq4}, $\bar\m_0^{I_v;y_0}$ is substituted as an input in $\bar\m_l^{I_w;z_0,\ldots,z_l}$, where $l\ge 1$ is the number of incoming edges at $w$. But $\bar\m_l^{I_w;z_0,\ldots,z_l}$ satisfies Definition \ref{ai7def3}(a)(v), and $\md{I_w}\ge 2$ by definition of $\bs\Ga_i$. Thus substituting $\bar\m_0^{I_v;y_0}=\boo^{\bar\jmath_2}_{\bar f{}^{\bar\jmath_0\bar\jmath_2}(y_0)}$ into $\bar\m_l^{I_w;z_0,\ldots,z_l}$ gives zero unless $l=1$, $\md{I_w}=2$ and $z_0=z_1=\bar f{}^{\bar\jmath_0\bar\jmath_2}(y_0)$, and then $\bar\m_1^{I_w;z_0,z_0}\ci \bar\m_0^{I_v;y_0}=\bar\m_1^{I_w;z_0,z_0}(\boo^{\max I_v}_{z_0})=\boo^{\max I_w}_{\bar f{}^{I_w}(z_0)}$. Continuing the argument, we see that $C(\bs\Ga_i)^{x_{l_{i-1}},\ldots,x_{l_i}}$ can only be nonzero if $\bs\Ga_i$ is of the form
\e
\begin{gathered}
\xymatrix@C=17pt@R=-4pt{
\overset{\!\!\!\!\!\!\{\bar\imath_{a_0},\bar\imath_{a_1},\bar\imath_{a_2}\}\!\!\!\!\!\!}{\sq} 
\ar[r] & 
\overset{\!\!\!\!\!\!\{\bar\imath_{a_2},\bar\imath_{a_3}\}\!\!\!\!\!\!}{\bu}  \ar[r] & \overset{\!\!\!\!\!\!\{\bar\imath_{a_3},\bar\imath_{a_4}\}\!\!\!\!\!\!}{\bu}  \ar[r] &  \cdots \ar[r] & 
\overset{\!\!\!\!\!\!\{\bar\imath_{a_{m-1}},\bar\imath_{a_m}=\bar\imath_{\bar M}\}\!\!\!\!\!\!}{\bu} \ar[r] & \tr & \!\!\!\!\!\!\!\!\!\!\!\!\bs\Ga_i, }
\end{gathered}
\label{ai10eq8}
\e
where $0\le a_0<a_1<\cdots<a_m=\bar M$ and where $\{j-1,j\}\subset\{a_0,a_1,a_2\}$, and then~$C(\bs\Ga_i)^{x_{l_{i-1}},\ldots,x_{l_i}}=\boo^{\bar\imath_{\bar M}}_z$.

Thus, any nonzero term $\pm C(\bs\Ga_1)^{x_{l_0},\ldots,x_{l_1}}\ot\cdots\ot C(\bs\Ga_k)^{x_{l_{k-1}},\ldots,x_{l_k}}
$ in \eq{ai7eq9} for $\bar\d_{n,k;y_0,\ldots,y_k}^{\bar I;x_0,\ldots,x_n}$ must contain a factor $C(\bs\Ga_i)^{x_{l_{i-1}},\ldots,x_{l_i}}=\boo^{\bar\imath_{\bar M}}_z$. But then when we apply $\bar\m_k^{\{\bar\imath_{\bar M}\};y_0,\ldots,y_k}$, by \eq{ai7eq19} we get zero unless $k=2$. This proves that the only way to get a nonzero term in $\bar\m_k^{\{\bar\imath_{\bar M}\}}\ci\bar\d_{n,k}^{\bar I}$ in (c)--(g) is when $k=2$, and for one of the graphs $\bs\Ga_1,\bs\Ga_2$ to be of the form \eq{ai10eq8}. We can then enumerate the possible $\bs\Ga_1,\bs\Ga_2$ by hand in each case.

In case (c), it is now easy to show that there are only two possibly nonzero contributions to \eq{ai7eq16}, coming from two sets of trees $\bs\Ga_1,\bs\Ga_2$ in equation \eq{ai7eq9} for $\bar\d_{n,2;y_0,y_1,y_2}^{\bar I;x_0,\ldots,x_n}$ illustrated as follows:
\begin{align*}
\xymatrix@C=17pt@R=-4pt{
\overset{\!\!\!\!\!\!\{\bar\imath_{j-2},\bar\imath_{j-1},\bar\imath_j\}\!\!\!\!\!\!}{\sq} 
\ar[r] & 
\overset{\!\!\!\!\!\!\{\bar\imath_j,\bar\imath_{j+1}\}\!\!\!\!\!\!}{\bu}  \ar[r] & \overset{\!\!\!\!\!\!\{\bar\imath_{j+1},\bar\imath_{j+2}\}\!\!\!\!\!\!}{\bu}  \ar[r] &  \cdots \ar[r] & 
\overset{\!\!\!\!\!\!\{\bar\imath_{\bar M-1},\bar\imath_{\bar M}\}\!\!\!\!\!\!}{\bu} \ar[r] & \tr & \!\!\!\!\!\!\!\!\!\!\!\!\bs\Ga_1 \\
&&& \!\!\!\!{\begin{subarray}{l} \ts 1\; \ci \\[-3pt] \ts \kern .91em\vdots \\[2pt] \ts n\; \ci\end{subarray}}
\ar[r] \ar@<1ex>[r]\ar@<-1ex>[r]
& \overset{\!\!\{\bar\imath_0,\ldots,\bar\imath_{\bar M}\}\sm\{\bar\imath_{j-1}\}\!\!}{\bu} \ar[r] & \tr & \!\!\!\!\!\!\!\!\!\!\!\!\bs\Ga_2 }\\
\xymatrix@C=17pt@R=-4pt{
\overset{\!\!\!\!\!\!\{\bar\imath_{j-1},\bar\imath_j,\bar\imath_{j+1}\}\!\!\!\!\!\!}{\sq} 
\ar[r] & 
\overset{\!\!\!\!\!\!\{\bar\imath_{j+1},\bar\imath_{j+2}\}\!\!\!\!\!\!}{\bu}  \ar[r] & \overset{\!\!\!\!\!\!\{\bar\imath_{j+2},\bar\imath_{j+3}\}\!\!\!\!\!\!}{\bu}  \ar[r] &  \cdots \ar[r] & 
\overset{\!\!\!\!\!\!\{\bar\imath_{\bar M-1},\bar\imath_{\bar M}\}\!\!\!\!\!\!}{\bu} \ar[r] & \tr & \!\!\!\!\!\!\!\!\!\!\!\!\bs\Ga_1 \\
&&& \!\!\!\!{\begin{subarray}{l} \ts 1\; \ci \\[-3pt] \ts \kern .91em\vdots \\[2pt] \ts n\; \ci\end{subarray}}
\ar[r] \ar@<1ex>[r]\ar@<-1ex>[r]
& \overset{\!\!\{\bar\imath_0,\ldots,\bar\imath_{\bar M}\}\sm\{\bar\imath_j\}\!\!}{\bu} \ar[r] & \tr & \!\!\!\!\!\!\!\!\!\!\!\!\bs\Ga_2 }
\end{align*}
Thus equation \eq{ai7eq16} for $\bar I$ reduces to
\begin{align*}
&(-1)^{\bar M}\cdot(-1)^{j-1+\bar M+1}\bar\m_2^{\{\bar\imath_{\bar M}\};y_0,y_1,y_2}(\boo_{\bar f{}^{\bar I}(x_0)}^{\bar\imath_{\bar M}}\ot \bar\m_n^{\{\bar\imath_0,\ldots,\bar\imath_{\bar M}\}\sm\{\bar\imath_{j-1}\};x_0,\ldots,x_n})\\
&+(-1)^{\bar M}\cdot(-1)^{j+\bar M+1}\bar\m_2^{\{\bar\imath_{\bar M}\};y_0,y_1,y_2}(\boo_{\bar f{}^{\bar I}(x_0)}^{\bar\imath_{\bar M}}\ot \bar\m_n^{\{\bar\imath_0,\ldots,\bar\imath_{\bar M}\}\sm\{\bar\imath_j\};x_0,\ldots,x_n})\\
&=(-1)^j\m_n^{\al(\bar I);x_0,\ldots,x_n}+(-1)^{j+1}\m_n^{\al(\bar I);x_0,\ldots,x_n}=0.
\end{align*}

In case (d), there are again only two possibly nonzero contributions to \eq{ai7eq16}, coming from two sets of trees $\bs\Ga_1,\bs\Ga_2$ in \eq{ai7eq9} for $\bar\d_{n,2;y_0,y_1,y_2}^{\bar I;x_0,\ldots,x_n}$ as follows:
\begin{align*}
\xymatrix@C=17pt@R=-4pt{
\overset{\!\!\!\!\!\!\{\bar\imath_0,\bar\imath_1,\bar\imath_2\}\!\!\!\!\!\!}{\sq} 
\ar[r] & 
\overset{\!\!\!\!\!\!\{\bar\imath_2,\bar\imath_3\}\!\!\!\!\!\!}{\bu}  \ar[r] & \overset{\!\!\!\!\!\!\{\bar\imath_3,\bar\imath_4\}\!\!\!\!\!\!}{\bu}  \ar[r] &  \cdots \ar[r] & 
\overset{\!\!\!\!\!\!\{\bar\imath_{\bar M-1},\bar\imath_{\bar M}\}\!\!\!\!\!\!}{\bu} \ar[r] & \tr & \!\!\!\!\!\!\!\!\!\!\!\!\bs\Ga_1 \\
&&& \!\!\!\!{\begin{subarray}{l} \ts 1\; \ci \\[-3pt] \ts \kern .91em\vdots \\[2pt] \ts n\; \ci\end{subarray}}
\ar[r] \ar@<1ex>[r]\ar@<-1ex>[r]
& \overset{\!\!\{\bar\imath_0,\ldots,\bar\imath_{\bar M}\}\sm\{\bar\imath_1\}\!\!}{\bu} \ar[r] & \tr & \!\!\!\!\!\!\!\!\!\!\!\!\bs\Ga_2 }\\
\xymatrix@C=17pt@R=14pt{
&& {\begin{subarray}{l} \ts 1\vphantom{\{}\; \ci \\[-3pt] \ts \kern .91em\vdots \\[2pt] \ts n\vphantom{\{}\; \ci\end{subarray}} \ar@<3ex>[r]  \ar@<1ex>[r] \ar@<-1ex>[r] \ar@<-3ex>[r] & {\begin{subarray}{l} \ts \{\bar\imath_0,\bar\imath_1\} \; \bu \\[-3pt] \ts \kern 3.4em\vdots \\[2pt] \ts \{\bar\imath_0,\bar\imath_1\}\; \bu\end{subarray}}
\ar@<3ex>[r]  \ar@<1ex>[r] \ar@<-1ex>[r] \ar@<-3ex>[r]
& \overset{\!\!\{\bar\imath_1,\bar\imath_2,\ldots,\bar\imath_{\bar M}\}\!\!}{\bu} \ar[r] & \tr & \!\!\!\!\!\!\!\!\!\!\!\!\bs\Ga_1 \\
&&&& \overset{\!\!\!\!\!\!\{\bar\imath_0,\bar\imath_1,\bar\imath_{\bar M}\}\!\!\!\!\!\!}{\sq} \ar[r] & \tr & \!\!\!\!\!\!\!\!\!\!\!\!\bs\Ga_2
}
\end{align*}
Thus equation \eq{ai7eq16} for $\bar I$ reduces to
\begin{align*}
&(-1)^{\bar M}\cdot(-1)^{\bar M}\bar\m_2^{\{\bar\imath_{\bar M}\};y_0,y_1,y_2}(\boo_{\bar f{}^{\bar I}(x_0)}^{\bar\imath_{\bar M}}\ot \bar\m_n^{\{\bar\imath_0,\ldots,\bar\imath_{\bar M}\}\sm\{\bar\imath_1\};x_0,\ldots,x_n})\\
&+(-1)^{\bar M}\cdot(-1)^{\bar M+1}\bar\m_2^{\{\bar\imath_{\bar M}\};y_0,y_1,y_2}\bigl((\bar\m_n^{\{\bar\imath_1,\ldots,\bar\imath_{\bar M}\};x_0,\ldots,x_n}\\
&\qquad {}\ci(\ts\bigot_{i=1}^n\bar\m_1^{\{\bar\imath_0,\bar\imath_1\};x_{i-1},x_i}))\ot \boo_{\bar f{}^{\bar\imath_0\bar\imath_{\bar M}}(x_n)}^{\bar\imath_{\bar M}}\bigr)\\
&=\m_n^{\al(\bar I);x_0,\ldots,x_n}-\m_n^{\al(\bar I);x_0,\ldots,x_n}\ci(\id^{\ot n})=0.
\end{align*}

Case (e) is similar, with two sets of trees $\bs\Ga_1,\bs\Ga_2$ in \eq{ai7eq9} for $\bar\d_{n,2;y_0,y_1,y_2}^{\bar I;x_0,\ldots,x_n}$:
\begin{align*}
\xymatrix@C=17pt@R=-4pt{
&& \overset{\!\!\!\!\!\!\{\bar\imath_{\bar M-2},\bar\imath_{\bar M-1},\bar\imath_{\bar M}\}\!\!\!\!\!\!}{\sq} \ar[r] & \tr & \!\!\!\!\!\!\!\!\!\!\!\!\bs\Ga_1 \\
& \!\!\!\!{\begin{subarray}{l} \ts 1\; \ci \\[-3pt] \ts \kern .91em\vdots \\[2pt] \ts n\; \ci\end{subarray}}
\ar[r] \ar@<1ex>[r]\ar@<-1ex>[r]
& \overset{\!\!\{\bar\imath_0,\ldots,\bar\imath_{\bar M}\}\sm\{\bar\imath_{\bar M-1}\}\!\!}{\bu} \ar[r] & \tr & \!\!\!\!\!\!\!\!\!\!\!\!\bs\Ga_2 }\\
\xymatrix@C=17pt@R=-4pt{
\!\!\!\!{\begin{subarray}{l} \ts 1\; \ci \\[-3pt] \ts \kern .91em\vdots \\[2pt] \ts n\; \ci\end{subarray}}
\ar[r] \ar@<1ex>[r]\ar@<-1ex>[r] &
\overset{\!\!\{\bar\imath_0,\ldots,\bar\imath_{\bar M-1}\}\!\!}{\bu}
\ar[r] & \overset{\!\!\{\bar\imath_{\bar M-1},\bar\imath_{\bar M}\}\!\!}{\bu} \ar[r] & \tr & \!\!\!\!\!\!\!\!\!\!\!\!\bs\Ga_1
\\
&& \overset{\!\!\!\!\!\!\{\bar\imath_0,\bar\imath_{\bar M-1},\bar\imath_{\bar M}\}\!\!\!\!\!\!}{\sq} \ar[r] & \tr & \!\!\!\!\!\!\!\!\!\!\!\!\bs\Ga_2
}
\end{align*}
Thus equation \eq{ai7eq16} for $\bar I$ reduces to
\begin{align*}
&(-1)^{\bar M}\bar\m_2^{\{\bar\imath_{\bar M}\};y_0,y_1,y_2}(\boo_{\bar f{}^{\bar I}(x_0)}^{\bar\imath_{\bar M}}\ot \bar\m_n^{\{\bar\imath_0,\ldots,\bar\imath_{\bar M}\}\sm\{\bar\imath_{\bar M-1}\};x_0,\ldots,x_n})\\
&+(-1)^{\bar M}\cdot(-1)\bar\m_2^{\{\bar\imath_{\bar M}\};y_0,y_1,y_2}\bigl((\bar\m_1^{\{\bar\imath_{\bar M-1},\bar\imath_{\bar M}\};\bar f{}^{\bar I}(x_0),\bar f{}^{\bar\imath_0\bar\imath_{\bar M}}(x_n)}\\
&\qquad \ci \bar\m_n^{\{\bar\imath_0,\ldots,\bar\imath_{\bar M-1}\};x_0,\ldots,x_n})\ot \boo_{\bar f{}^{\bar\imath_0\bar\imath_{\bar M}}(x_n)}^{\bar\imath_{\bar M}}\bigr)\\
&=(-1)^{\bar M}\m_n^{\al(\bar I);x_0,\ldots,x_n}-(-1)^{\bar M}\id\ci\m_n^{\al(\bar I);x_0,\ldots,x_n}=0.
\end{align*}

In case (f), the only way to get a nonzero term is for both $\bs\Ga_1,\bs\Ga_2$ to be of the form \eq{ai10eq8}, which forces $n=0$. There are two such sets of trees:
\begin{align*}
\xymatrix@C=17pt@R=-4pt{
\overset{\!\!\!\!\!\!\{\bar\imath_1,\bar\imath_2,\bar\imath_3\}\!\!\!\!\!\!}{\sq} \ar[r] & \tr & \!\!\!\!\!\!\!\!\!\!\!\!\bs\Ga_1 \\
\overset{\!\!\{\bar\imath_0,\bar\imath_1,\bar\imath_3\}\!\!}{\sq} \ar[r] & \tr & \!\!\!\!\!\!\!\!\!\!\!\!\bs\Ga_2 }\qquad\qquad\qquad
\xymatrix@C=17pt@R=-4pt{
\overset{\!\!\{\bar\imath_0,\bar\imath_1,\bar\imath_2\}\!\!}{\bu}
\ar[r] & \overset{\!\!\{\bar\imath_2,\bar\imath_3\}\!\!}{\sq} \ar[r] & \tr & \!\!\!\!\!\!\!\!\!\!\!\!\bs\Ga_1
\\
& \overset{\!\!\!\!\!\!\{\bar\imath_0,\bar\imath_2,\bar\imath_3\}\!\!\!\!\!\!}{\sq} \ar[r] & \tr & \!\!\!\!\!\!\!\!\!\!\!\!\bs\Ga_2
}
\end{align*}
Thus equation \eq{ai7eq16} for $\bar I$ reduces to
\begin{align*}
&-\bar\m_2^{\{\bar\imath_3\};y_0,y_1,y_2}(\boo_{\bar f{}^{\bar I}(x_0)}^{\bar\imath_3}\ot \boo_{\bar f{}^{\bar I}(x_0)}^{\bar\imath_3})\\
&+\bar\m_2^{\{\bar\imath_3\};y_0,y_1,y_2}\bigl((\bar\m_1^{\{\bar\imath_2,\bar\imath_3\};\bar f{}^{\bar\imath_0\bar\imath_2}(x_0),\bar f{}^{\bar\imath_0\bar\imath_2}(x_0)}(\boo_{\bar f{}^{\bar\imath_0\bar\imath_2}(x_0)}^{\bar\imath_2})\ot \boo_{\bar f{}^{\bar I}(x_0)}^{\bar\imath_3}\bigr)\\
&=-\boo_{\bar f{}^{\bar I}(x_0)}^{\bar\imath_3}+\boo_{\bar f{}^{\bar I}(x_0)}^{\bar\imath_3}=0.
\end{align*}

In case (g) there are no nonzero terms, so equation \eq{ai7eq16} for $\bar I$ holds trivially. This completes the proof of Proposition~\ref{ai7prop1}.

\subsection{Proofs of Theorems \ref{ai6thm2}, \ref{ai7thm2}, and \ref{ai8thm1}}
\label{ai103}

We will prove Theorem \ref{ai8thm1}. Theorem \ref{ai7thm2} is a special case of this with $\Ga=\{0,\iy\}$, as in Example \ref{ai5ex2}. Theorem \ref{ai6thm2} follows from Theorem \ref{ai7thm2} by restricting to $A_\iy$-categories $\sA$ with $\sO=\{*\}$. The only differences in the proofs of Theorems \ref{ai7thm2} and \ref{ai8thm1} are in the ranges of variables in the sums in \eq{ai10eq9}, and that in the filtered case some sums in \eq{ai10eq9} are infinite and convergent, but in the non-filtered case there are only finitely many nonzero terms.

Observe that when $k=0$ equation \eq{ai8eq6} follows from Definition \ref{ai7def2}(a), since $\d_{l,0;y_0}^{\{i_M\};w_0,\ldots,w_l}=0$, and $\d_{l,0;y_0}^{I;w_0,\ldots,w_l}=0$ unless $M=1$ and $l=0$, in which case $\d_{n,0;w_0}^{\{i_0\};x_0,\ldots,x_n}=0$. When $k=1$ equation \eq{ai8eq6} follows from \eq{ai7eq10} and \eq{ai8eq4}--\eq{ai8eq5}. When $M=0$ equation \eq{ai8eq6} follows from the filtered $A_\iy$-category analogue of $\bar\d\ci\bar\d=0$ for the bar complex of an $A_\iy$-algebra in \S\ref{ai36}. So it remains to prove \eq{ai8eq6} when $M>0$ and $k\ge 2$. We do this by a double induction on $M$ and $k$. For given $\ti M>0$ and $\ti k\ge 2$, our inductive hypothesis is
\begin{itemize}
\setlength{\itemsep}{0pt}
\setlength{\parsep}{0pt}
\item[$(*)^{\ti M}_{\ti k}$] Equation \eq{ai8eq6} holds for $M>0$ and $k,n\ge 0$ provided either $M<\ti M$, or $M=\ti M$ and $k\le\ti k$.
\end{itemize}
Let $\ti M>0$ and $\ti k\ge 2$ be given. If $\ti k>2$ we suppose that $(*)^{\ti M}_{\ti k-1}$ holds. If $\ti k=2$ we suppose that $(*)^{\ti M-1}_k$ holds for all $k\ge 2$, unless $\ti M=1$, when we make no assumption. For any $\ti n\ge 0$, we will now prove \eq{ai8eq6} with $\ti M,\ti k,\ti n$ in place of $M,k,n$, in the long equation
\ea
& \sum_{\begin{subarray}{l} 1\le \ti l\le \ti n+1, \\ w_0,\ldots,w_{\ti l}\in\sO^{i_0} \end{subarray}} \d_{\ti l,\ti k;y_0,\ldots,y_{\ti k}}^{I;w_0,\ldots,w_{\ti l}}\ci\d_{\ti n,\ti l;w_0,\ldots,w_{\ti l}}^{\{i_0\};x_0,\ldots,x_{\ti n}}
\label{ai10eq9}\\
&+(-1)^{\ti M}\sum_{\ti l\ge 0, \;  w_0,\ldots,w_{\ti l}\in\sO^{i_{\ti M}}}\d_{\ti l,\ti k;y_0,\ldots,y_{\ti k}}^{\{i_{\ti M}\};w_0,\ldots,w_{\ti l}}\ci\d_{\ti n,\ti l;w_0,\ldots,w_{\ti l}}^{I;x_0,\ldots,x_{\ti n}}
\nonumber\\
&=\sum_{\ti l=1}^{\ti n+1}
\begin{aligned}[t]&
\sum_{\begin{subarray}{l} 
\ti l=l'+l'', \; l',l''\ge 0, \\ 
0<m\le \ti M,\\ 
0=j_0<j_1<\cdots<j_m=\ti M, \\
k_0,\ldots,k_m\ge 0, \\
k_0=l', \; k_m=\ti k-1, \\
z_0^a,\ldots,z_{k_a}^a\in\sO^{i_{j_a}},\; 0<a<m
\end{subarray}}
\Bigl(\begin{aligned}[t]&
(-1)^{\ti k-1+m+\ti M+l'm+\sum_{a=1}^{m-1}(j_a-a)}\\
&\bigl(\d_{k_{m-1},k_m;y_0,\ldots,y_{\ti k-1}}^{\{i_{j_{m-1}},\ldots,i_{j_m}\};z_0^{m-1},\ldots,z_{k_{m-1}}^{m-1}}\\
&\ci\d_{k_{m-2},k_{m-1};z_0^{m-1},\ldots,z_{k_{m-1}}^{m-1}}^{\{i_{j_{m-2}},\ldots,i_{j_{m-1}}\};z_0^{m-2},\ldots,z_{k_{m-2}}^{m-2}}\\
&{}\ci\cdots\ci\d_{k_0,k_1;z_0^1,\ldots,z_{k_1}^1}^{\{i_{j_0},\ldots,i_{j_1}\};w_0,\ldots,w_{l'}}\bigr)\\
&{}\ot\d_{l'',1;y_{\ti k-1},y_{\ti k}}^{\{i_{j_0},i_{j_1},i_{j_2},\ldots,i_{j_m}\};w_{l'},w_{l'+1},\ldots,w_l}
\Bigr)\end{aligned}\\
&\ci\Bigl((-1)^{l''(l'+\ti n)}
\de_{\begin{subarray}{l}
w_i=x_{\ti n-\ti l+i}, \\ l'\le i\le \ti l\end{subarray}}\cdot
\d_{\ti n-l'',l';w_0,\ldots,w_{l'}}^{\{i_0\};x_0,\ldots,x_{\ti n-l''}}\ot\id^{\ot^{l''}}\\
&+(-1)^{l'}\de_{\begin{subarray}{l} w_i=x_i, \\ 0\le i\le l'\end{subarray}}\cdot
\id^{\ot^{l'}}\ot\d_{\ti n-l',l'';w_{l'},\ldots,w_{\ti l}}^{\{i_0\};x_{l'},\ldots,x_{\ti n}}
\Bigr)\end{aligned}
\allowdisplaybreaks
\nonumber\\
&+\sum_{\ti l\ge 0}\begin{aligned}[t]&
(-1)^{\ti M}\Bigl((-1)^{\ti k-1+\ti l}
\de_{\begin{subarray}{l}y_i=w_{\ti l-\ti k+i}, \\ \ti k-1\le i\le \ti k\end{subarray}}\cdot
\d_{\ti l-1,\ti k-1;y_0,\ldots,y_{\ti k-1}}^{\{i_{\ti M}\};w_0,\ldots,w_{\ti l-1}}\ot\id\Bigr)\ci{}\\
&\biggl(
\sum_{\begin{subarray}{l} 
\ti n=n'+n'', \; n',n''\ge 0, \\ 
0<m\le \ti M,\\ 
0=j_0<j_1<\cdots<j_m=\ti M, \\
k_0,\ldots,k_m\ge 0, \\ 
k_0=n', \; k_m=\ti l-1,\\
z_0^a,\ldots,z_{k_a}^a\in\sO^{i_{j_a}},\; 0<a<m
\end{subarray}}
\begin{aligned}[t]&
(-1)^{\ti l-1+m+\ti M+n'm+\sum_{a=1}^{m-1}(j_a-a)}\\
&\bigl(\d_{k_{m-1},k_m;w_0,\ldots,w_{\ti k-1}}^{\{i_{j_{m-1}},\ldots,i_{j_m}\};z_0^{m-1},\ldots,z_{k_{m-1}}^{m-1}}\\
&\ci\d_{k_{m-2},k_{m-1};z_0^{m-1},\ldots,z_{k_{m-1}}^{m-1}}^{\{i_{j_{m-2}},\ldots,i_{j_{m-1}}\};z_0^{m-2},\ldots,z_{k_{m-2}}^{m-2}}\\
&{}\ci\cdots\ci\d_{k_0,k_1;z_0^1,\ldots,z_{k_1}^1}^{\{i_{j_0},\ldots,i_{j_1}\};x_0,\ldots,x_{n'}}\bigr)\\
&{}\ot\d_{n'',1;w_{\ti l-1},w_{\ti l}}^{\{i_{j_0},i_{j_1},i_{j_2},\ldots,i_{j_m}\};x_{n'},x_{n'+1},\ldots,x_n}\biggr)\end{aligned}
\end{aligned}
\allowdisplaybreaks
\nonumber\\
&+\sum_{\ti l\ge 0}\begin{aligned}[t]&
(-1)^{\ti M}\Bigl((-1)^{\ti k-1}\de_{\begin{subarray}{l} y_i=w_i, \\ 0\le i\le \ti k-1\end{subarray}}\cdot
\id^{\ot^{\ti k-1}}\ot\d_{\ti l-\ti k+1,1;y_{\ti k-1},y_{\ti k}}^{\{i_{\ti M}\};w_{\ti k-1},\ldots,w_{\ti l}}\Bigr)\ci{}\\
&\biggl(\sum_{\begin{subarray}{l} 
\ti n=n'+n'', \; n',n''\ge 0, \\ 
0<m\le \ti M,\\ 
0=j_0<j_1<\cdots<j_m=\ti M, \\
k_0,\ldots,k_m\ge 0, \\ 
k_0=n', \; k_m=\ti k-1,\\
z_0^a,\ldots,z_{k_a}^a\in\sO^{i_{j_a}},\; 0<a<m
\end{subarray}}
\begin{aligned}[t]&
(-1)^{\begin{subarray}{l}(\ti l-\ti k+1)(\ti k-1+m+\ti M) \\ +n'(\ti l-\ti k+m)+\sum_{a=1}^{m-1}(j_a-a)\end{subarray}}\\
&\bigl(\d_{k_{m-1},k_m;w_0,\ldots,w_{\ti k-1}}^{\{i_{j_{m-1}},\ldots,i_{j_m}\};z_0^{m-1},\ldots,z_{k_{m-1}}^{m-1}}\\
&\ci\d_{k_{m-2},k_{m-1};z_0^{m-1},\ldots,z_{k_{m-1}}^{m-1}}^{\{i_{j_{m-2}},\ldots,i_{j_{m-1}}\};z_0^{m-2},\ldots,z_{k_{m-2}}^{m-2}}\\
&{}\ci\cdots\ci\d_{k_0,k_1;z_0^1,\ldots,z_{k_1}^1}^{\{i_{j_0},\ldots,i_{j_1}\};x_0,\ldots,x_{n'}}\bigr)\\
&{}\ot\d_{n'',\ti l-\ti k+1;w_{\ti k-1},\ldots,w_{\ti l}}^{\{i_{j_0},i_{j_1},i_{j_2},\ldots,i_{j_m}\};x_{n'},x_{n'+1},\ldots,x_n}\biggr)\end{aligned}
\end{aligned}
\allowdisplaybreaks
\nonumber\\
&=\sum_{\begin{subarray}{l} 
\ti n=n'+n'', \; n',n''\ge 0, \\ 
0<m\le \ti M,\\ 
0=j_0<j_1<\cdots<j_m=\ti M, \\
k_0,\ldots,k_m\ge 0, \\ 
k_0=n', \; k_m=\ti k-1,\\
z_0^a,\ldots,z_{k_a}^a\in\sO^{i_{j_a}},\; 0<a<m
\end{subarray}}
\begin{aligned}[t]&
(-1)^{1+\ti k+\ti M+m+n'+mn'+\sum_{a=1}^{m-1}(j_a-a)}\\
&\bigl(\d_{k_{m-1},k_m;y_0,\ldots,y_{\ti k-1}}^{\{i_{j_{m-1}},\ldots,i_{j_m}\};z_0^{m-1},\ldots,z_{k_{m-1}}^{m-1}}\\
&\ci\d_{k_{m-2},k_{m-1};z_0^{m-1},\ldots,z_{k_{m-1}}^{m-1}}^{\{i_{j_{m-2}},\ldots,i_{j_{m-1}}\};z_0^{m-2},\ldots,z_{k_{m-2}}^{m-2}}\\
&{}\ci\cdots\ci\d_{k_0,k_1;z_0^1,\ldots,z_{k_1}^1}^{\{i_{j_0},\ldots,i_{j_1}\};x_0,\ldots,x_{n'}}\bigr)\ot{}\\
&\biggl[(-1)^m\sum_{\begin{subarray}{l} l\ge 0, \\ v_0,\ldots,v_l\in\sO^{i_0} \end{subarray}}
\begin{aligned}[t]&\d_{l,1;y_{\ti k-1},y_{\ti k}}^{\{i_{\ti M}\};v_0,\ldots,v_l}\ci{}\\
&\d_{n'',l;v_0,\ldots,v_l}^{\{i_{j_0},i_{j_1},i_{j_2},\ldots,i_{j_m}\};x_{n'},x_{n'+1},\ldots,x_n}
\end{aligned}\\
&+\sum_{\begin{subarray}{l}
l=1,\ldots,n''+1,\\ v_0,\ldots,v_l\in\sO^{i_0} 
\end{subarray}}\begin{aligned}[t]& \d_{l,1;y_{\ti k-1},y_{\ti k}}^{\{i_{j_0},i_{j_1},i_{j_2},\ldots,i_{j_m}\};v_0,\ldots,v_l}\ci{}\\
&\d_{n'',l;v_0,\ldots,v_l}^{\{i_0\};x_{n'},x_{n'+1},\ldots,x_n}\end{aligned}\biggr]\end{aligned}
\allowdisplaybreaks
\nonumber\\
&+\sum_{\begin{subarray}{l} 
\ti n=n'+n'', \; n',n''\ge 0, \\ 
0<m\le \ti M,\\ 
0=j_0<j_1<\cdots<j_m=\ti M, \\
k_0,\ldots,k_m\ge 0, \; k_m=\ti k-1, \\
z_0^a,\ldots,z_{k_a}^a\in\sO^{i_{j_a}},\; 0\le a<m
\end{subarray}}
\begin{aligned}[t]&(-1)^{1+\ti k+\ti M+n'm+\sum_{a=1}^{m-1}(j_a-a)}\\
&\bigl(\d_{k_{m-1},k_m;y_0,\ldots,y_{\ti k-1}}^{\{i_{j_{m-1}},\ldots,i_{j_m}\};z_0^{m-1},\ldots,z_{k_{m-1}}^{m-1}}\\
&\ci\d_{k_{m-2},k_{m-1};z_0^{m-1},\ldots,z_{k_{m-1}}^{m-1}}^{\{i_{j_{m-2}},\ldots,i_{j_{m-1}}\};z_0^{m-2},\ldots,z_{k_{m-2}}^{m-2}}\\
&{}\ci\cdots\ci\d_{k_0,k_1;z_0^1,\ldots,z_{k_1}^1}^{\{i_{j_0},\ldots,i_{j_1}\};z_0^0,\ldots,z_{k_0}^0}
\ci \d_{n',k_0;z_0^0,\ldots,z_{k_0}^0}^{\{i_0\};x_0,\ldots,x_{n'}}\bigr)\\
&\ot\d_{n'',1;y_{\ti k-1},y_{\ti k}}^{\{i_{j_0},i_{j_1},i_{j_2},\ldots,i_{j_m}\};x_{n'},x_{n'+1},\ldots,x_n}\end{aligned}
\allowdisplaybreaks
\nonumber\\
&+\sum_{\begin{subarray}{l} 
\ti n=n'+n'', \; n',n''\ge 0, \\ 
0<m\le \ti M,\\ 
0=j_0<j_1<\cdots<j_m=\ti M, \\
k_0,\ldots,k_m\ge 0, \; k_0=n',\\
z_0^a,\ldots,z_{k_a}^a\in\sO^{i_{j_a}},\; 0<a\le m
\end{subarray}}
\begin{aligned}[t]&(-1)^{1+\ti k+\ti M+n'm+\sum_{a=1}^{m-1}(j_a-a)}\cdot\bigl((-1)^{\ti M-m+1}\\
&\d_{k_m,\ti k-1;y_0,\ldots,y_{\ti k-1}}^{\{i_{\ti M}\};z_0^m,\ldots,z_{k_m}^m}\ci\d_{k_{m-1},k_m;z_0^m,\ldots,z_{k_m}^m}^{\{i_{j_{m-1}},\ldots,i_{j_m}\};z_0^{m-1},\ldots,z_{k_{m-1}}^{m-1}}\\
&\ci\d_{k_{m-2},k_{m-1};z_0^{m-1},\ldots,z_{k_{m-1}}^{m-1}}^{\{i_{j_{m-2}},\ldots,i_{j_{m-1}}\};z_0^{m-2},\ldots,z_{k_{m-2}}^{m-2}}\\
&{}\ci\cdots\ci\d_{k_0,k_1;z_0^1,\ldots,z_{k_1}^1}^{\{i_{j_0},\ldots,i_{j_1}\};x_0,\ldots,x_{n'}}\bigr)\\
&{}\ot\d_{n'',1;y_{\ti k-1},y_{\ti k}}^{\{i_{j_0},i_{j_1},i_{j_2},\ldots,i_{j_m}\};x_{n'},x_{n'+1},\ldots,x_n}\end{aligned}
\allowdisplaybreaks
\nonumber\\
&+\sum_{\begin{subarray}{l} 
\ti n=n'+n'', \; n',n''\ge 0, \\ 
0<m\le \ti M,\\ 
0=j_0<j_1<\cdots<j_m=\ti M, \\
k_0,\ldots,k_m\ge 0, \\ 
k_0=n', \; k_m=\ti k-1, \\
z_0^a,\ldots,z_{k_a}^a\in\sO^{i_{j_a}},\\ 
0<a<m, \\
0<l<m,\; \hat k\ge 0, \\
v_0,\ldots,v_{\hat k}\in\sO^{i_{j_l}}
\end{subarray}}
\begin{aligned}[t]&(-1)^{1+\ti k+\ti M+n'm+\sum_{a=1}^{m-1}(j_a-a)}\\
&\bigl(\d_{k_{m-1},k_{m};y_0,\ldots,y_{\ti k-1}}^{\{i_{j_{m-1}},\ldots,i_{j_m}\};z_0^{m-1},\ldots,z_{k_{m-1}}^{m-1}}\ci\cdots\ci\d_{k_l,k_{l+1};z_0^{l+1},\ldots,z_{k_{l+1}}^{l+1}}^{\{i_{j_l},\ldots,i_{j_{l+1}}\};z_0^l,\ldots,z_{k_l}^l}\\
&\ci\bigl[(-1)^{l+j_l}\d_{\hat k,k_l;z_0^l,\ldots,z_{k_l}^l}^{\{i_{j_l}\};v_0,\ldots,v_{\hat k}}-(-1)^{l+j_l}\d_{\hat k,k_l;z_0^l,\ldots,z_{k_l}^l}^{\{i_{j_l}\};v_0,\ldots,v_{\hat k}}\bigr]\\
&\ci\d_{k_{l-1},\hat k;v_0,\ldots,v_{\hat k}}^{\{i_{j_{l-1}},\ldots,i_{j_l}\};z_0^{l-1},\ldots,z_{k_{l-1}}^{l-1}}\!\ci\!\d_{k_{l-2},k_{l-1};z_0^{l-1},\ldots,z_{k_{l-1}}^{l-1}}^{\{i_{j_{l-2}},\ldots,i_{j_{l-1}}\};z_0^{l-2},\ldots,z_{k_{l-2}}^{l-2}}\!\ci\!\cdots\\
&{}\ci\d_{k_0,k_1;z_0^1,\ldots,z_{k_1}^1}^{\{i_{j_0},\ldots,i_{j_1}\};x_0,\ldots,x_{n'}}\bigr)
\ot\d_{n'',1;y_{\ti k-1},y_{\ti k}}^{\{i_{j_0},i_{j_1},i_{j_2},\ldots,i_{j_m}\};x_{n'},x_{n'+1},\ldots,x_n}\end{aligned}
\nonumber\\
&=0.
\nonumber
\ea

Here in the first step, the first term rewrites the first term in the first line of \eq{ai10eq9}, first by applying \eq{ai7eq13} to $\d_{\ti l,\ti k;y_0,\ldots,y_{\ti k}}^{I;w_0,\ldots,w_{\ti l}}$ with $k'=\ti k-1$ and $k''=1$, which yields a sum over splittings $\ti l=l'+l''$, and then applying \eq{ai7eq12} to $\d_{\ti n,\ti l;w_0,\ldots,w_{\ti l}}^{\{i_0\};x_0,\ldots,x_{\ti n}}$ splitting $\ti l=l'+l''$. The second and third terms rewrite the second term in the first line of \eq{ai10eq9}, first by applying \eq{ai7eq12} to $\d_{\ti l,\ti k;y_0,\ldots,y_{\ti k}}^{\{i_M\};w_0,\ldots,w_{\ti l}}$ with $k'=\ti k-1$ and $k''=1$. We split this into two terms, and then apply \eq{ai7eq13} to $\d_{\ti n,\ti l;w_0,\ldots,w_{\ti l}}^{I;x_0,\ldots,x_{\ti n}}$, firstly with $\ti l=(\ti l-1)+1$, and secondly with~$\ti l=(\ti k-1)+(\ti l-\ti k+1)$.

In the first three terms in the second step, we write out the compositions explicitly and simplify. In the first and second terms we include extra signs from changing the order of operators $\d_{n,k;y_0,\ldots,y_k}^{I;x_0,\ldots,x_n}$, as in \S\ref{ai31}. We also change variables in the sums. The fourth term in the second step is clearly zero.

For the third step, observe that the first term in the second step is zero by \eq{ai8eq6} with $k=1$ applied to $[\cdots]$. We claim that the final three terms in the second step also add up to zero. This holds because:
\begin{itemize}
\setlength{\itemsep}{0pt}
\setlength{\parsep}{0pt}
\item[(i)] The third term, involving $\d_{k_0,k_1;z_0^1,\ldots,z_{k_1}^1}^{\{i_{j_0},\ldots,i_{j_1}\};z_0^0,\ldots,z_{k_0}^0}\ci \d_{n',k_0;z_0^0,\ldots,z_{k_0}^0}^{\{i_0\};x_0,\ldots,x_{n'}}$, and the second half of the case $l=1$ in the fifth term in the second step, involving $-(-1)^{1+j_1}\d_{\hat k,k_1;z_0^1,\ldots,z_{k_1}^1}^{\{i_{j_1}\};v_0,\ldots,v_{\hat k}}\ci\d_{k_0,\hat k;v_0,\ldots,v_{\hat k}}^{\{i_{j_0},\ldots,i_{j_1}\};x_0,\ldots,x_{n'}}$, sum to zero.
\item[(ii)] For $1<\ti l<m$, the first half of the case $l=\ti l-1$ in the fifth term, involving $\d_{k_{\ti l-1},k_{\ti l};z_0^{\ti l},\ldots,z_{k_{\ti l}}^{\ti l}}^{\{i_{j_{\ti l-1}},\ldots,i_{j_{\ti l}}\};z_0^{\ti l-1},\ldots,z_{k_{\ti l-1}}^{\ti l-1}}\ci (-1)^{\ti l-1+j_{\ti l-1}}\d_{\hat k,k_{\ti l-1};z_0^{\ti l-1},\ldots,z_{k_{\ti l-1}}^{\ti l-1}}^{\{i_{j_{\ti l-1}}\};v_0,\ldots,v_{\hat k}}$, and the second half of the case $l=\ti l$ in the fifth term in the second step, involving $-(-1)^{\ti l+j_{\ti l}}\d_{\hat k,k_{\ti l};z_0^{\ti l},\ldots,z_{k_{\ti l}}^{\ti l}}^{\{i_{j_{\ti l}}\};v_0,\ldots,v_{\hat k}}
\ci\d_{k_{\ti l-1},\hat k;v_0,\ldots,v_{\hat k}}^{\{i_{j_{\ti l-1}},\ldots,i_{j_{\ti l}}\};z_0^{\ti l-1},\ldots,z_{k_{\ti l-1}}^{\ti l-1}}$, sum to zero.
\item[(iii)] The fourth term in the second step, involving $(-\!1)^{\ti M-m+1}\d_{k_m,\ti k-1;y_0,\ldots,y_{\ti k-1}}^{\{i_{\ti M}\};z_0^m,\ldots,z_{k_m}^m}\ab\ci\d_{k_{m-1},k_m;z_0^m,\ldots,z_{k_m}^m}^{\{i_{j_{m-1}},\ldots,i_{j_m}\};z_0^{m-1},\ldots,z_{k_{m-1}}^{m-1}}$, and the first half of the case $l=m-1$ in the fifth term in the second step, involving $\d_{k_{m-1},k_m;z_0^m,\ldots,z_{k_m}^m}^{\{i_{j_{m-1}},\ldots,i_{j_m}\};z_0^{m-1},\ldots,z_{k_{m-1}}^{m-1}}\ci(-1)^{m-1+j_{m-1}}\d_{\hat k,k_{m-1};z_0^{m-1},\ldots,z_{k_{m-1}}^{m-1}}^{\{i_{j_{m-1}}\};v_0,\ldots,v_{\hat k}}$, sum to zero.
\end{itemize}
In each case, we can show this holds using \eq{ai8eq6} and the inductive hypothesis, where the `involved' terms indicated above are those appearing in \eq{ai8eq6}. Some relabelling of the variables $k_a,\hat k,n'$ is required, for example, in (i), in the first term we fix $k_1,n'$ and sum over $k_0$, and in the second term we fix $k_0,k_1$ and sum over $\hat k$, where $n'$ in the first term is relabelled to $k_0$ in the second term.

This proves \eq{ai10eq9}, and Theorem \ref{ai8thm1} follows by induction.

\subsection{Proofs of Theorems \ref{ai6thm3}, \ref{ai7thm3}, and \ref{ai8thm2}}
\label{ai104}

We will prove Theorem \ref{ai8thm2} for $\AfCat^{\un,\iy}_R$. Working with $\AfCat^{\un,\iy}_R$ involves using the $A_\iy$-category $F^0\Fun_{A^\rfi_\iy}^\un(\sA^{\ti\imath_0},\ab\sA^{\ti\imath_{\ti M}})$ over $F^0R$ from Definition \ref{ai5def10}(f), with morphisms $F^0\Hom_{A^\rfi_\iy}(\fF,\fG)$. The proof for $\AwCat^{\un,\iy}_R$ is essentially the same, except that we omit $F^0(\cdots)$ throughout.

Theorem \ref{ai7thm3} is a special case of Theorem \ref{ai8thm2} for $\AfCat^{\un,\iy}_R$ with $\Ga=\{0,\iy\}$, as in Example \ref{ai5ex2}. The only differences in the proofs of Theorems \ref{ai7thm3} and \ref{ai8thm2} are in the ranges of variables in the sums (for Theorem \ref{ai8thm2} we allow larger ranges, but the extra terms are automatically zero in the case of Theorem \ref{ai7thm3}), and that in the filtered case some sums are infinite and convergent, but in the non-filtered case there are only finitely many nonzero terms. Theorem \ref{ai6thm3} follows from Theorem \ref{ai7thm3} by restricting to $A_\iy$-categories $\sA$ with~$\sO=\{*\}$. 

First we set up some notation and prove three preliminary propositions.

\begin{dfn}
\label{ai10def1}
In the proof we will consider the following situation:

Suppose data $\bigl((\sO^i)_{i=0}^N,(f^{ij})_{0\le i\le j\le N},(\sH^i)_{i=0}^N,(\m^I)_{\es\ne I\subseteq\{0,\ldots,N\}}\bigr)$ is given as in Definitions \ref{ai8def1}(a)(i)-(v) and \ref{ai8def3} for some $N\ge 2$, po-monoid $(\Ga,\pre)$  and complete $\Ga$-filtered commutative ring $R$, {\it except\/} that we do not assume \eq{ai8eq4}--\eq{ai8eq5} hold for all $\es\ne I=\{i_0<i_1<\cdots<i_M\}\subseteq\{0,\ldots,N\}$ and $n\ge 0$, but only for certain $I$ and all $n\ge 0$, where we will specify explicitly which $I$ for each result. If we assume \eq{ai8eq4}--\eq{ai8eq5} hold for $I$, we will also assume \eq{ai8eq4}--\eq{ai8eq5} hold for all~$\es\ne I'\subsetneq I$.

We always require the filtration conditions on $\m_0^{I;x_0}$ in Definition \ref{ai8def1}(a)(iv) hold, as these are needed to make infinite sums such as \eq{ai10eq14} below converge.

We always assume that \eq{ai8eq4}--\eq{ai8eq5} hold when $I=\{i\}$ for all $i=0,\ldots,N$ and $n\ge 0$. This implies that the data $\sO^i,\sH^i,\m^{\{i\}}$ defines a $\Ga$-filtered $A_\iy$-category $\sA^i$ over $R$ as in Definition \ref{ai8def1} and Example \ref{ai8ex1}(a), and we require $\sA^i$ to be {\it unital\/} with identities $\boo^i_x$ for $x\in\sO^i$ as in Definition \ref{ai8def1}(a)(v), so that we always require \eq{ai7eq17} holds.

In cases when we require \eq{ai8eq4}--\eq{ai8eq5} hold when $I=\{i,j\}$ for all $0\le i<j\le N$ and $n\ge 0$, the data $f^{ij},\m^{\{i,j\}}$ defines a filtered $A_\iy$-functor $\fF^{ij}:\sA^i\ra\sA^j$ as in Definition \ref{ai8def1} and Example \ref{ai8ex1}(b), and then we always assume that $\fF^{ij}$ is unital, so that \eq{ai7eq18} holds.

We always assume the unitality condition \eq{ai7eq19} holds for all $I,k$ with $(M,k)\ne(0,1),(1,0)$, whether or not we assume \eq{ai8eq4}--\eq{ai8eq5} hold for these~$I$.

In cases when we require \eq{ai8eq4}--\eq{ai8eq5} to hold for $I=\{i_0,i_1,i_2\}$, so that 
$\m^{\{i_0,i_1,i_2\}}$ is equivalent to a filtered $A_\iy$-natural transformation $\fH^{i_0i_1i_2}:\fF^{i_1i_2}\ci\fF^{i_0i_1}\Ra\fF^{i_0i_2}$ as in Example \ref{ai8ex1}(c), which is unital as \eq{ai7eq19} holds, we require that $\fH^{i_0i_1i_2}$ should be a filtered $A_\iy$-natural isomorphism, as in Definition~\ref{ai8def3}(a). \end{dfn}

\begin{dfn}
\label{ai10def2}
Continue in the situation of Definition~\ref{ai10def1}.

Let $\es\ne\ti I=\{\ti\imath_0<\ti\imath_1<\cdots< \ti\imath_{\ti M}\}\subseteq\{0,\ldots,N\}$ for $\ti M\ge 2$. We assume \eq{ai8eq4}--\eq{ai8eq5} hold for all $\es\ne I\subsetneq\ti I$ and $x_0,
\ldots,x_n$ for $n\ge 0$, but we do {\it not\/} yet assume \eq{ai8eq5} when $I=\ti I$. Then we have unital filtered $A_\iy$-functors $\fF^{\ti\imath_0\ti\imath_{\ti M}}:\sA^{\ti\imath_0}\ra\sA^{\ti\imath_{\ti M}}$ and $\fF^{\ti\imath_{j-1}\ti\imath_j}:\sA^{\ti\imath_{j-1}}\ra\sA^{\ti\imath_j}$. Generalizing equation \eq{ai7eq2}, define $\fF^{\ti I}=\fF^{\ti\imath_{\ti M-1}\ti\imath_{\ti M}}\ci\fF^{\ti\imath_{\ti M-2}\ti\imath_{\ti M-1}}\ci\cdots\ci\fF^{\ti\imath_0\ti\imath_1}:\sA^{\ti\imath_0}\ra\sA^{\ti\imath_{\ti M}}$. 

Note that \eq{ai8eq4}--\eq{ai8eq5} hold for $\es\ne I\subsetneq\ti I$ imply \eq{ai8eq6} for $\es\ne I\subsetneq\ti I$, by Theorem \ref{ai8thm1}. As in Definition \ref{ai5def10} we have a (filtered, strict) unital $A_\iy$-category $F^0\Fun_{A^\rfi_\iy}^\un(\sA^{\ti\imath_0},\ab\sA^{\ti\imath_{\ti M}})$ over $F^0R$, and $\fF^{\ti I},\fF^{\ti\imath_0\ti\imath_{\ti M}}$ are objects in this, with morphism complex 
\e
\bigl(F^0\Hom_{A^\rfi_\iy}^\un(\fF^{\ti I},\fF^{\ti\imath_0\ti\imath_{\ti M}}),\n_1=\n_1^{\fF^{\ti I},\fF^{\ti\imath_0\ti\imath_{\ti M}}}\bigr),
\label{ai10eq10}
\e
for $\n_1$ as in \eq{ai4eq15} omitting the condition $l_i>l_{i-1}$, and we can take its cohomology $H^*\bigl(F^0\Hom_{A^\rfi_\iy}^\un(\fF^{\ti I},\fF^{\ti\imath_0\ti\imath_{\ti M}}),\n_1\bigr)$.

Observe that (without assuming \eq{ai8eq5} holds for $I=\ti I$), we have
\e
\m^{\ti I}=\bigl(\m_k^{\ti I;x_0,\ldots,x_k})^{x_0,\ldots,x_k\in\sO^{\ti\imath_0}}_{k\ge 0}\in F^0\Hom_{A^\rfi_\iy}^\un(\fF^{\ti I},\fF^{\ti\imath_0\ti\imath_{\ti M}})^{2-\ti M}
\label{ai10eq11}
\e
by Definition \ref{ai5def10}, where $\m^{\ti I}$ lies in $F^0\Hom_{A^\rfi_\iy}^\un(\cdots)\subset F^0\Hom_{A^\rfi_\iy}(\cdots)$ as we assume \eq{ai7eq19} holds and~$\ti M\ge 2$.

Let $x_0,\ldots,x_n\in\sO^{\ti\imath_0}$ and $y_0,\ldots,y_k\in\sO^{\ti\imath_{\ti M}}$ for $n,k\ge 0$. Then Definition \ref{ai7def2} defines an operator $\d_{n,k;y_0,\ldots,y_k}^{\ti I;x_0,\ldots,x_n}$ in \eq{ai7eq6} and \eq{ai7eq9} in terms of the data $\m_l^{I;w_0,\ldots,w_l}$ for $I\subseteq\ti I$. Write $\D_{n,k;y_0,\ldots,y_k}^{\ti I;x_0,\ldots,x_n}$ for the sum of all terms in the definition \eq{ai7eq9} of $\d_{n,k;y_0,\ldots,y_k}^{\ti I;x_0,\ldots,x_n}$ which do {\it not\/} involve any operator $\m_l^{\ti I;w_0,\ldots,w_l}$.

Note that $\D_{n,1;y_0,y_1}^{\ti I;x_0,\ldots,x_n}=0$ when $k=1$, as then $\d_{n,1;y_0,y_1}^{\ti I;x_0,\ldots,x_n}$ is either $\m_n^{\ti I;x_0,\ldots,x_n}$ or 0, depending on $y_0,y_1$. Also $\D_{n,k;y_0,\ldots,y_k}^{\ti I;x_0,\ldots,x_n}=0$ when $\ti M=2$, as every term in $\d_{n,k;y_0,\ldots,y_k}^{\ti I;x_0,\ldots,x_n}$ involves some~$\m_l^{\{\ti\imath_0,\ti\imath_1,\ti\imath_2\};w_0,\ldots,w_l}$.

For $n\ge 0$ and $x_0,\ldots,x_n\in\sO^{\ti\imath_0}$, define
\e
{\rm E}_n^{\ti I;x_0,\ldots,x_n}
=\sum_{\begin{subarray}{l} k=1,\ldots,n+\ti M-1,\\
y_0,\ldots,y_k\in\sO^{\ti\imath_{\ti M}} \end{subarray}\!\!\!\!\!\!\!\!\!\!\!}\m_k^{\{\ti\imath_{\ti M}\};y_0,\ldots,y_k}\ci\D_{n,k;y_0,\ldots,y_k}^{\ti I;x_0,\ldots,x_n}
\label{ai10eq12}
\e
in $\sH^{\ti\imath_{\ti M}}(f^{\ti I}(x_0),f^{\ti\imath_0\ti\imath_{\ti M}}(x_n))^{3-\ti M-n}$. By definition of $\D_{n,k;y_0,\ldots,y_k}^{\ti I;x_0,\ldots,x_n}$, there are only finitely many nonzero terms in this sum. We claim that
\e
{\rm E}^{\ti I}=\bigl({\rm E}_n^{\ti I;x_0,\ldots,x_n})^{x_0,\ldots,x_n\in\sO^{\ti\imath_0}}_{n\ge 0}\in F^0\Hom_{A^\rfi_\iy}^\un(\fF^{\ti I},\fF^{\ti\imath_0\ti\imath_{\ti M}})^{3-\ti M}.
\label{ai10eq13}
\e
It is obvious that ${\rm E}^{\ti I}$ lies in $F^0\Hom_{A^\rfi_\iy}(\fF^{\ti I},\fF^{\ti\imath_0\ti\imath_{\ti M}})^{3-\ti M}$. We can deduce the additional equation \eq{ai4eq18} defining $F^0\Hom_{A^\rfi_\iy}^\un(\cdots)\subset F^0\Hom_{A^\rfi_\iy}(\cdots)$ from the fact that the $\m_k^{I;w_0,\ldots,w_k}$ satisfy \eq{ai7eq18}--\eq{ai7eq20} and~$\ti M\ge 2$.

Note that ${\rm E}^{\ti I}=0$ when $\ti M=2$, as $\D_{n,k;y_0,\ldots,y_k}^{\ti I;x_0,\ldots,x_n}=0$.\end{dfn}

\begin{prop}
\label{ai10prop1}
Work in the situation of Definition\/ {\rm\ref{ai10def2}}. Then
\smallskip

\noindent{\bf(a)} We have
\ea
&\d_{n,k;y_0,\ldots,y_k}^{\ti I;x_0,\ldots,x_n}=\D_{n,k;y_0,\ldots,y_k}^{\ti I;x_0,\ldots,x_n}
\label{ai10eq14}\\
&+\sum_{\begin{subarray}{l}
0=l_0\le l_1\le \cdots\le l_k=n,\\
1\le j\le k: \\
y_i=f^{\ti I}(x_{l_i}),\; 0\le i<j,\\
y_i=f^{\ti\imath_0\ti\imath_{\ti M}}(x_{l_i}),\; j\le i\le k
\end{subarray}}
\begin{aligned}[t]&
(-1)^{(k+1)(\ti M+1)}\cdot(-1)^{\sum_{i=1}^{j-1}(l_i-l_{i-1}+1)(\ti M+1)}\\
&\cdot(-1)^{\sum_{i=1}^k(l_i-l_{i-1}+1)(k-i)} \\
&(\fF^{\ti I})_{l_1-l_0}^{x_{l_0},\ldots,x_{l_1}}
\ot 
(\fF^{\ti I})_{l_2-l_1}^{x_{l_1},\ldots,x_{l_2}}\ot\cdots\\
&\ot (\fF^{\ti I})_{l_{j-1}-l_{j-2}}^{x_{l_{j-2}},\ldots,x_{l_{j-1}}}
\ot\m_{l_j-l_{j-1}}^{\ti I;x_{l_{j-1}},\ldots,x_{l_j}}
\ot\\
&(\fF^{\ti\imath_0\ti\imath_{\ti M}})_{l_{j+1}-l_j}^{x_{l_j},\ldots,x_{l_{j+1}}}\ot\cdots\ot(\fF^{\ti\imath_0\ti\imath_{\ti M}})_{l_k-l_{k-1}}^{x_{l_{k-1}},\ldots,x_{l_k}}.
\end{aligned}
\nonumber
\ea
That is, by definition $\D_{n,k;y_0,\ldots,y_k}^{\ti I;x_0,\ldots,x_n}$ is the sum of terms in $\d_{n,k;y_0,\ldots,y_k}^{\ti I;x_0,\ldots,x_n}$ in \eq{ai7eq9} which do not involve operators $\m_l^{\ti I;w_0,\ldots,w_l},$ and the second line of\/ \eq{ai10eq14} is the sum of terms in $\d_{n,k;y_0,\ldots,y_k}^{\ti I;x_0,\ldots,x_n}$ which do involve operators~$\m_l^{\ti I;w_0,\ldots,w_l}$.
\smallskip

\noindent{\bf(b)} The $\D_{n,k;y_0,\ldots,y_k}^{\ti I;x_0,\ldots,x_n}$ satisfy an analogue of\/ {\rm\eq{ai7eq13},} for $k=k'+k''$ with\/~{\rm$k,k''\ge 1$:}
\ea
&\D_{n,k;y_0,\ldots,y_k}^{\ti I;x_0,\ldots,x_n}
\label{ai10eq15}\\
&=\sum_{n=n'+n'', \; n',n''\ge 0\!\!\!\!\!\!\!\!\!\!\!\!\!\!\!\!\!\!\!\!\!\!\!\!\!\!\!\!\!\!\!\!\!\!\!\!\!\!\!\!}
(-1)^{k''(k'+\ti M+n'+1)}\D_{n',k';y_0,\ldots,y_{k'}}^{\ti I;x_0,\ldots,x_{n'}}\!\!\ot\!\d_{n'',k'';y_{k'},y_{k'+1},\ldots,y_k}^{\{\ti\imath_0,\ti\imath_M\};x_{n'},x_{n'+1},\ldots,x_n}
\nonumber\\
&+\!\!\!\sum_{\begin{subarray}{l} 
n=n'+n'', \; n',n''\ge 0, \\ 
1<l<\ti M,\\
0=j_0<j_1<\cdots<j_l=\ti M, \\
k_0,\ldots,k_l\ge 0, \; k_0=n', \; k_l=k', \\
w_0^i,\ldots,w_{k_i}^i\in\sO^{\ti\imath_{j_i}},\; 0<i<l
\end{subarray}}\begin{aligned}[t]&
(-1)^{k''(k'+l+\ti M)+n'(k''+l+1)+\sum_{a=1}^{l-1}(j_a-a)}\\
&\bigl(\d_{k_{l-1},k_l;y_0,\ldots,y_{k'1}}^{\{\ti\imath_{j_{l-1}},\ldots,\ti\imath_{j_l}\};w_0^{l-1},\ldots,w_{k_{l-1}}^{l-1}}\ci\d_{k_{l-2},k_{l-1};w_0^{l-1},\ldots,w_{k_{l-1}}^{l-1}}^{\{\ti\imath_{j_{l-2}},\ldots,\ti\imath_{j_{l-1}}\};w_0^{l-2},\ldots,w_{k_{l-2}}^{l-2}}\\
&{}\ci\cdots\ci\d_{k_0,k_1;w_0^1,\ldots,w_{k_1}^1}^{\{\ti\imath_{j_0},\ldots,\ti\imath_{j_1}\};x_0,\ldots,x_{n'}}\bigr)\\
&{}\ot\d_{n'',k'';y_{k'},y_{k'+1},\ldots,y_k}^{\{\ti\imath_{j_0},\ti\imath_{j_1},\ti\imath_{j_2},\ldots,\ti\imath_{j_l}\};x_{n'},x_{n'+1},\ldots,x_n}
\end{aligned}
\nonumber\\
&+\!\!\!\sum_{\begin{subarray}{l} 
n=n'+n'', \; n',n''\ge 0, \\ 
k_0,\ldots,k_{\ti M}\ge 0, \\ k_0=n', \; k_{\ti M}=k', \\
w_0^i,\ldots,w_{k_i}^i\in\sO^{\ti\imath_i},\; 0<i<\ti M
\end{subarray}}\begin{aligned}[t]&
(-1)^{k''k'+n'(k''+\ti M+1)}\\
&\bigl(\d_{k_{\ti M-1},k_{\ti M};y_0,\ldots,y_{k'}}^{\{\ti\imath_{\ti M-1},\ti\imath_{\ti M}\};w_0^{\ti M-1},\ldots,w_{k_{\ti M-1}}^{\ti M-1}}\ci\d_{k_{\ti M-2},k_{\ti M-1};w_0^{\ti M-1},\ldots,w_{k_{\ti M-1}}^{\ti M-1}}^{\{\ti\imath_{\ti M-2},\ti\imath_{\ti M-1}\};w_0^{\ti M-2},\ldots,w_{k_{\ti M-2}}^{\ti M-2}}\\
&{}\ci\cdots\ci\d_{k_0,k_1;w_0^1,\ldots,w_{k_1}^1}^{\{\ti\imath_0,\ti\imath_1\};x_0,\ldots,x_{n'}}\bigr)\ot\D_{n'',k'';y_{k'},y_{k'+1},\ldots,y_k}^{\ti I;x_{n'},x_{n'+1},\ldots,x_n}.
\end{aligned}
\nonumber
\ea
\noindent{\bf(c)} The $\D_{n,k;y_0,\ldots,y_k}^{\ti I;x_0,\ldots,x_n}$ satisfy an analogue of the second equation of\/~{\rm\eq{ai8eq6}:}
\ea
&0=\sum_{\begin{subarray}{l} 1\le l\le n+1, \\ w_0,\ldots,w_l\in\sO^{\ti\imath_0} \end{subarray}}
\D_{l,k;y_0,\ldots,y_k}^{\ti I;w_0,\ldots,w_l}\ci\d_{n,l;w_0,\ldots,w_l}^{\{\ti\imath_0\};x_0,\ldots,x_n}
\label{ai10eq16}\\
&+(-1)^{\ti M}\sum_{l\ge 0, \; w_0,\ldots,w_l\in\sO^{\ti\imath_{\ti M}}}
\d_{l,k;y_0,\ldots,y_k}^{\{\ti\imath_{\ti M}\};w_0,\ldots,w_l}\ci\D_{n,l;w_0,\ldots,w_l}^{\ti I;x_0,\ldots,x_n}
\nonumber\\
&+\sum_{\begin{subarray}{l}
0\le p\le q\le n, \\
1\le j\le k, \\
k_0,\ldots,k_{\ti M}\ge 0, \\ k_0=p, \; k_{\ti M}=j-1, \\
w_0^i,\ldots,w_{k_i}^i\in\sO^{\ti\imath_i},\\ 
0<i<\ti M, \; m\ge 1\\
z_0,\ldots,z_m\in\sO^{\ti\imath_{\ti M}} 
\end{subarray}}
\begin{aligned}[t]&
(-1)^{(j+k+p+1)\ti M}\cdot(-1)^{p+jk+jq+kq}\cdot{}\\
&\bigl(\d_{k_{\ti M-1},k_{\ti M};y_0,\ldots,y_{j-1}}^{\{\ti\imath_{\ti M-1},\ti\imath_{\ti M}\};w_0^{\ti M-1},\ldots,w_{k_{\ti M-1}}^{\ti M-1}}\ci\d_{k_{\ti M-2},k_{\ti M-1};w_0^{\ti M-1},\ldots,w_{k_{\ti M-1}}^{\ti M-1}}^{\{\ti\imath_{\ti M-2},\ti\imath_{\ti M-1}\};w_0^{\ti M-2},\ldots,w_{k_{\ti M-2}}^{\ti M-2}}\\
&{}\ci\cdots\ci\d_{k_0,k_1;w_0^1,\ldots,w_{k_1}^1}^{\{\ti\imath_0,\ti\imath_1\};x_0,\ldots,x_p}\bigr)\ot
\bigl(\d_{m,1;y_{j-1},y_j}^{\{\ti\imath_{\ti M}\};z_0,\ldots,z_m}\ci\D_{q-p,m;z_0,\ldots,z_m}^{\ti I;x_p,\ldots,x_q}\bigr)\\
&{}\ot\d_{n-q,k-j;y_j,\ldots,y_k}^{\{\ti\imath_0,\ti\imath_M\};x_q,\ldots,x_n}\bigr).
\end{aligned}
\nonumber
\ea

An equivalent way to write \eq{ai10eq16} is
\ea
&0=\sum_{\begin{subarray}{l} 1\le l\le n+1, \\ w_0,\ldots,w_l\in\sO^{\ti\imath_0} \end{subarray}}
\D_{l,k;y_0,\ldots,y_k}^{\ti I;w_0,\ldots,w_l}\ci\d_{n,l;w_0,\ldots,w_l}^{\{\ti\imath_0\};x_0,\ldots,x_n}
\label{ai10eq17}\\
&+(-1)^{\ti M}\sum_{l\ge 0, \; w_0,\ldots,w_l\in\sO^{\ti\imath_{\ti M}}}
\d_{l,k;y_0,\ldots,y_k}^{\{\ti\imath_{\ti M}\};w_0,\ldots,w_l}\ci\D_{n,l;w_0,\ldots,w_l}^{\ti I;x_0,\ldots,x_n}
\nonumber\\
&-\sum_{\begin{subarray}{l}
0=l_0\le l_1\le \cdots\le l_k=n,\\
1\le j\le k: \\
y_i=f^{\ti I}(x_{l_i}),\; 0\le i<j,\\
y_i=f^{\ti\imath_0\ti\imath_{\ti M}}(x_{l_i}),\; j\le i\le k, \\
m=1,\ldots,l_j-l_{j-1}+\ti M-1,\\
z_0,\ldots,z_m\in\sO^{\ti\imath_{\ti M}} 
\end{subarray}}
\begin{aligned}[t]&
(-1)^{k\ti M}\cdot(-1)^{\sum_{i=1}^{j-1}(l_i-l_{i-1}+1)\ti M}\\
&\cdot(-1)^{\sum_{i=1}^k(l_i-l_{i-1}+1)(k-i)} 
\\
&\bigl(\fF_{l_1-l_0}^{\ti I;x_{l_0},\ldots,x_{l_1}}\ot\fF_{l_2-l_1}^{\ti I;x_{l_1},\ldots,x_{l_2}}\ot\cdots\ot\fF_{l_{j-1}-l_{j-2}}^{\ti I;x_{l_{j-2}},\ldots,x_{l_{j-1}}} \\
&{}\ot
\bigl(\m_m^{\{\ti\imath_{\ti M}\};z_0,\ldots,z_m}\ci\D_{l_j-l_{j-1},m;z_0,\ldots,z_m}^{\ti I;x_{l_{j-1}},\ldots,x_{l_j}}\bigr)\\
&{}\ot\fF_{l_{j+1}-l_j}^{\ti\imath_0\ti\imath_{\ti M};x_{l_j},\ldots,x_{l_{j+1}}}\ot\cdots\ot\fF_{l_k-l_{k-1}}^{\ti\imath_0\ti\imath_{\ti M};x_{l_{k-1}},\ldots,x_{l_k}}\bigr).
\end{aligned}
\nonumber
\ea

Here the proofs of\/ {\rm\eq{ai10eq14}--\eq{ai10eq15}} do not use {\rm\eq{ai8eq4}--\eq{ai8eq5}}. Those of\/ {\rm\eq{ai10eq16}--\eq{ai10eq17}} do, but not for $I=\ti I$.
\end{prop}

\begin{proof}
For (a), it is not difficult to show that the only terms in \eq{ai7eq9} involving some $\m_l^{\ti I;w_0,\ldots,w_l}$ have $k\ge 1$, and for $0\le l_0\le l_1\le\cdots\le l_k=n$ and some $1\le j\le k$ have labelled graphs $\bs\Ga_1,\ldots,\bs\Ga_k$ as in the following diagrams:
\begin{gather}
\xymatrix@C=17pt@R=6pt{
\!\!\!\!{\begin{subarray}{l} \ts l_{i-1}\!+\!1\; \ci \\[-3pt] \ts \kern .91em\vdots \\[2pt] \ts \,l_{i\phantom{-1}}\phantom{+1}\; \ci\end{subarray}}  
\ar[r] \ar@<2ex>[r] \ar@<-2ex>[r]
&{\begin{subarray}{l} 
\ts \overset{\{\ti\imath_0,\ti\imath_1\}}{\bu}  \\
\ts \overset{\{\ti\imath_0,\ti\imath_1\}}{\bu}  \\
\ts \overset{\{\ti\imath_0,\ti\imath_1\}}{\bu}  
\end{subarray}}
\ar[r] \ar@<2ex>[r]\ar@<-2ex>[r] & 
{\begin{subarray}{l} 
\ts \overset{\{\ti\imath_1,\ti\imath_2\}}{\bu}  \\
\ts \overset{\{\ti\imath_1,\ti\imath_2\}}{\bu}  \\
\ts \overset{\{\ti\imath_1,\ti\imath_2\}}{\bu}  
\end{subarray}}
\ar[r] \ar@<2ex>[r]\ar@<-2ex>[r] &  \cdots \ar[r] \ar@<1.5ex>[r]\ar@<-1.5ex>[r] & 
\overset{\{\ti\imath_{\ti M-1},\ti\imath_{\ti M}\}}{\bu}  
\ar[r] & \tr & \!\!\!\!\!\!\bs\Ga_i,\; 1\le i<j, }
\nonumber\\
\xymatrix@C=17pt@R=6pt{ {\begin{subarray}{l} \ts l_{j-1}\!+\!1\; \ci \\[-3pt] \ts \kern .91em\vdots \\[2pt] \ts \,l_{j\phantom{-1}}\phantom{+1}\; \ci\end{subarray}}
\ar[r] \ar@<1.5ex>[r]\ar@<-1.5ex>[r]
& \overset{\ti I}{\bu} \ar[r] & \tr & \!\!\!\!\!\!\!\!\!\!\!\!\bs\Ga_j \\
{\begin{subarray}{l} \ts l_{i-1}\!+\!1\; \ci \\[-3pt] \ts \kern .91em\vdots \\[2pt] \ts \,l_{i\phantom{-1}}\phantom{+1}\; \ci\end{subarray}}
\ar[r] \ar@<1.5ex>[r]\ar@<-1.5ex>[r]
& \overset{\{\ti\imath_0,\ti\imath_{\ti M}\}}{\bu} \ar[r] & \tr &
\bs\Ga_i,\; j<i\le k.  }
\label{ai10eq18}
\end{gather}
Here in the first line, we allow trees $\bs\Ga_i$ in which we can have different numbers of vertices in each column with labels $\{\ti\imath_{k-1},\ti\imath_k\}$, and these vertices can have arbitrary numbers of incoming edges, so there are many possible $\bs\Ga_i$.

The second line of \eq{ai10eq14} is the result of summing the terms in \eq{ai7eq9} from labelled trees $\bs\Ga_1,\ldots,\bs\Ga_k$ as in \eq{ai10eq18}. For fixed $i=1,\ldots,j-1$ and $l_{i-1},l_i$, summing over all trees $\bs\Ga_i$ as in the first line of \eq{ai10eq18} has the effect of computing $(\fF^{\ti I})_{l_i-l_{i-1}}^{x_{l_{i-1}},\ldots,x_{l_i}}$ in the multiple composition $\fF^{\ti I}=\fF^{\ti\imath_{\ti M-1}\ti\imath_{\ti M}}\ci\fF^{\ti\imath_{\ti M-2}\ti\imath_{\ti M-1}}\ci\cdots\ci\fF^{\ti\imath_0\ti\imath_1}$, where compositions of filtered $A_\iy$-functors are defined as sums as in \eq{ai5eq10}. Some of the signs $(-1)^{\cdots}$ in \eq{ai7eq9}, those coming from both $v,w$ vertices in $\bs\Ga_i$, are the signs in \eq{ai5eq10}. The remaining signs have been simplified and written in~\eq{ai10eq14}.

For (b),  we modify the proof in \S\ref{ai101} of equation \eq{ai7eq13} in Theorem \ref{ai7thm1}. Consider \eq{ai7eq13} with $M=\ti M$, and delete all terms in both sides which involve an operator $\m_l^{\ti I;w_0,\ldots,w_l}$. The left hand side yields $\D_{n,k;y_0,\ldots,y_k}^{\ti I;x_0,\ldots,x_n}$ by Definition \ref{ai10def2}. The right hand side of \eq{ai7eq13} involves a sum over $0<l\le\ti M$ and $0=j_0<\cdots<j_l=\ti M$. Observe that terms with $1<l<\ti M$ involve no $\m_l^{\ti I;w_0,\ldots,w_l}$, since $\ti I$ is split between the two halves of the formula. We have written the sum over $1<l<\ti M$ as the second line of \eq{ai10eq15}. This leaves the terms coming from $l=1$ and $l=\ti M$ in \eq{ai7eq13}. These appear in the first and third lines of \eq{ai10eq15}, and are modified by replacing $\d^{\ti I;\cdots}_{\cdots}$ by~$\D^{\ti I;\cdots}_{\cdots}$.

For (c), to prove \eq{ai10eq16} we modify the proof in \S\ref{ai103} of equation \eq{ai8eq6} in Theorem \ref{ai8thm1}. The proof is by induction on $k=0,1,\ldots.$ The first case $k=0$ is trivial. Let $\ti k\ge 1$ be given, and suppose by induction that \eq{ai10eq16} holds for all $0\le k<\ti k$ and all $n\ge 0$. For any $n\ge 0$ and $x_0,\ldots,x_n\in\sO^{\ti\imath_0}$, $y_0,\ldots,y_{\ti k}\in\sO^{\ti\imath_{\ti M}}$, we will now prove \eq{ai10eq16} with $k=\ti k$, in the long analogue of~\eq{ai10eq9}:
\ea
&\sum_{\begin{subarray}{l} 1\le l\le n+1, \\ w_0,\ldots,w_l\in\sO^{\ti\imath_0} \end{subarray}}
\D_{l,\ti k;y_0,\ldots,y_{\ti k}}^{\ti I;w_0,\ldots,w_l}\ci\d_{n,l;w_0,\ldots,w_l}^{\{\ti\imath_0\};x_0,\ldots,x_n}{}_{\color{red}(a1)\color{black}}
\label{ai10eq19}\\
&+(-1)^{\ti M}\sum_{l\ge 0, \; w_0,\ldots,w_l\in\sO^{\ti\imath_{\ti M}}}
\d_{l,\ti k;y_0,\ldots,y_{\ti k}}^{\{\ti\imath_{\ti M}\};w_0,\ldots,w_l}\ci\D_{n,l;w_0,\ldots,w_l}^{\ti I;x_0,\ldots,x_n}{}_{\color{red}(a2)\color{black}}
\nonumber\\
&+\sum_{\begin{subarray}{l}
0\le p\le q\le n, \\
1\le j\le \ti k, \\
k_0,\ldots,k_{\ti M}\ge 0, \\ k_0=p, \; k_{\ti M}=j-1, \\
w_0^i,\ldots,w_{k_i}^i\in\sO^{\ti\imath_i},\\ 
0<i<\ti M, \; m\ge 1\\
z_0,\ldots,z_m\in\sO^{\ti\imath_{\ti M}} \\
{\qquad \color{red}(a3)\color{black}}
\end{subarray}}
\begin{aligned}[t]&
(-1)^{(j+\ti k+p+1)\ti M}\cdot(-1)^{p+j\ti k+jq+\ti kq}\cdot{}\\
&\bigl(\d_{k_{\ti M-1},k_{\ti M};y_0,\ldots,y_{j-1}}^{\{\ti\imath_{\ti M-1},\ti\imath_{\ti M}\};w_0^{\ti M-1},\ldots,w_{k_{\ti M-1}}^{\ti M-1}}\ci\d_{k_{\ti M-2},k_{\ti M-1};w_0^{\ti M-1},\ldots,w_{k_{\ti M-1}}^{\ti M-1}}^{\{\ti\imath_{\ti M-2},\ti\imath_{\ti M-1}\};w_0^{\ti M-2},\ldots,w_{k_{\ti M-2}}^{\ti M-2}}\\
&{}\ci\cdots\ci\d_{k_0,k_1;w_0^1,\ldots,w_{k_1}^1}^{\{\ti\imath_0,\ti\imath_1\};x_0,\ldots,x_p}\bigr)\ot
\bigl(\d_{m,1;y_{j-1},y_j}^{\{\ti\imath_{\ti M}\};z_0,\ldots,z_m}\ci\D_{q-p,m;z_0,\ldots,z_m}^{\ti I;x_p,\ldots,x_q}\bigr)\\
&{}\ot\d_{n-q,\ti k-j;y_j,\ldots,y_{\ti k}}^{\{\ti\imath_0,\ti\imath_M\};x_q,\ldots,x_n}
\end{aligned}
\allowdisplaybreaks
\nonumber\\
&=\sum_{\begin{subarray}{l} 1\le l\le n+1, \\ w_0,\ldots,w_l\in\sO^{\ti\imath_0} \end{subarray}}
\biggl[\sum_{l=l'+l'', \; l',l''\ge 0
\!\!\!\!\!\!\!\!\!\!\!\!\!\!\!\!\!\!\!\!\!\!\!\!\!\!\!\!\!\!\!\!\!\!\!\!\!}
(-1)^{\ti k+\ti M+l'}\D_{l',\ti k-1;y_0,\ldots,y_{\ti k-1}}^{\ti I;w_0,\ldots,w_{l'}}\!\!\ot\!\d_{l'',1;y_{\ti k-1},y_{\ti k}}^{\{\ti\imath_0,\ti\imath_M\};w_{l'},\ldots,w_l}{}_{\color{red}(b1)\color{black}}
\allowdisplaybreaks
\nonumber\\
&+\begin{aligned}[t]&
\sum_{\begin{subarray}{l} 
l=l'+l'', \; l',l''\ge 0, \\ 
1<m<\ti M,\\ 
0=j_0<j_1<\cdots<j_m=\ti M, \\
k_0,\ldots,k_m\ge 0, \\
k_0=l', \; k_m=\ti k-1, \\
z_0^a,\ldots,z_{k_a}^a\in\sO^{i_{j_a}},\; 0<a<m
\end{subarray}}
\Bigl(\begin{aligned}[t]&
(-1)^{\ti k-1+m+\ti M+l'm+\sum_{a=1}^{m-1}(j_a-a)}\\
&\bigl(\d_{k_{m-1},k_m;y_0,\ldots,y_{\ti k-1}}^{\{\ti\imath_{j_{m-1}},\ldots,\ti\imath_{j_m}\};z_0^{m-1},\ldots,z_{k_{m-1}}^{m-1}}\\
&\ci\d_{k_{m-2},k_{m-1};z_0^{m-1},\ldots,z_{k_{m-1}}^{m-1}}^{\{\ti\imath_{j_{m-2}},\ldots,\ti\imath_{j_{m-1}}\};z_0^{m-2},\ldots,z_{k_{m-2}}^{m-2}}\\
&{}\ci\cdots\ci\d_{k_0,k_1;z_0^1,\ldots,z_{k_1}^1}^{\{\ti\imath_{j_0},\ldots,\ti\imath_{j_1}\};w_0,\ldots,w_{l'}}\bigr)\\
&{}\ot\d_{l'',1;y_{\ti k-1},y_{\ti k}}^{\{\ti\imath_{j_0},\ti\imath_{j_1},\ti\imath_{j_2},\ldots,\ti\imath_{j_m}\};w_{l'},w_{l'+1},\ldots,w_l}
\Bigr){}_{\color{red}(b2)\color{black}}\biggr]\end{aligned}\\
&\ci\Bigl((-1)^{l''(l'+n)}
\de_{\begin{subarray}{l}
w_i=x_{n-l+i}, \\ l'\le i\le l\end{subarray}}\cdot
\d_{n-l'',l';w_0,\ldots,w_{l'}}^{\{\ti\imath_0\};x_0,\ldots,x_{n-l''}}\ot\id^{\ot^{l''}}{}_{\color{red}(i)\color{black}}\\
&+(-1)^{l'}\de_{\begin{subarray}{l} w_i=x_i, \\ 0\le i\le l'\end{subarray}}\cdot
\id^{\ot^{l'}}\ot\d_{n-l',l'';w_{l'},\ldots,w_l}^{\{\ti\imath_0\};x_{l'},\ldots,x_n}{}_{\color{red}(ii)\color{black}}\Bigr)\end{aligned}
\allowdisplaybreaks
\nonumber\\
&+\sum_{\begin{subarray}{l} l\ge 0, \\ w_0,\ldots,w_l\in\sO^{\ti\imath_{\ti M}} \end{subarray}}
\begin{aligned}[t]&
(-1)^{\ti M}\Bigl((-1)^{\ti k-1+l}
\de_{\begin{subarray}{l}y_i=w_{l-\ti k+i}, \\ \ti k-1\le i\le \ti k\end{subarray}}\cdot
\d_{l-1,\ti k-1;y_0,\ldots,y_{\ti k-1}}^{\{\ti\imath_{\ti M}\};w_0,\ldots,w_{l-1}}\ot\id\Bigr)\ci{}\\
&\biggl[
\sum_{n=n'+n'', \; n',n''\ge 0 \!\!\!\!\!\!\!\!\!\!\!\!\!\!\!\!\!\!\!\!\!\!\!\!\!\!\!\!\!\!\!\!\!\!\!\!\!}
(-1)^{l+\ti M+n'}\D_{n',l-1;w_0,\ldots,w_{l-1}}^{\ti I;x_0,\ldots,x_{n'}}\!\!\ot\!\d_{n'',1;w_{l-1},w_l}^{\{\ti\imath_0,\ti\imath_M\};x_{n'},x_{n'+1},\ldots,x_n}{}_{\color{red}(b3)\color{black}}
\end{aligned}
\allowdisplaybreaks
\nonumber\\
&+\sum_{\begin{subarray}{l} 
n=n'+n'', \; n',n''\ge 0, \\ 
1<m<\ti M,\\ 
0=j_0<j_1<\cdots<j_m=\ti M, \\
k_0,\ldots,k_m\ge 0, \\ 
k_0=n', \; k_m=l-1,\\
z_0^a,\ldots,z_{k_a}^a\in\sO^{i_{j_a}},\; 0<a<m
\end{subarray}}
\begin{aligned}[t]&
(-1)^{l-1+m+\ti M+n'm+\sum_{a=1}^{m-1}(j_a-a)}\\
&\bigl(\d_{k_{m-1},k_m;w_0,\ldots,w_{\ti k-1}}^{\{\ti\imath_{j_{m-1}},\ldots,\ti\imath_{j_m}\};z_0^{m-1},\ldots,z_{k_{m-1}}^{m-1}}\\
&\ci\d_{k_{m-2},k_{m-1};z_0^{m-1},\ldots,z_{k_{m-1}}^{m-1}}^{\{\ti\imath_{j_{m-2}},\ldots,\ti\imath_{j_{m-1}}\};z_0^{m-2},\ldots,z_{k_{m-2}}^{m-2}}\\
&{}\ci\cdots\ci\d_{k_0,k_1;z_0^1,\ldots,z_{k_1}^1}^{\{\ti\imath_{j_0},\ldots,\ti\imath_{j_1}\};x_0,\ldots,x_{n'}}\bigr)\\
&{}\ot\d_{n'',1;w_{l-1},w_l}^{\{\ti\imath_{j_0},\ti\imath_{j_1},\ti\imath_{j_2},\ldots,\ti\imath_{j_m}\};x_{n'},x_{n'+1},\ldots,x_n}{}_{\color{red}(b4)\color{black}}
\biggr]\end{aligned}
\allowdisplaybreaks
\nonumber\\
&+\sum_{\begin{subarray}{l} l\ge 0, \\ w_0,\ldots,w_l\in\sO^{\ti\imath_{\ti M}} \end{subarray}}
\begin{aligned}[t]&
(-1)^{\ti M}\Bigl((-1)^{\ti k-1}\de_{\begin{subarray}{l} y_i=w_i, \\ 0\le i\le \ti k-1\end{subarray}}\cdot
\id^{\ot^{\ti k-1}}\ot\d_{l-\ti k+1,1;y_{\ti k-1},y_{\ti k}}^{\{\ti\imath_{\ti M}\};w_{\ti k-1},\ldots,w_l}\Bigr)\ci{}\\
&\biggl[\sum_{\begin{subarray}{l} n=n'+n'', \\ n',n''\ge 0 \end{subarray}}\begin{aligned}[t] &
(-1)^{(l-\ti k+1)(\ti k+\ti M+n')} \\
&\D_{n',\ti k-1;y_0,\ldots,y_{\ti k-1}}^{\ti I;x_0,\ldots,x_{n'}}\ot\d_{n'',l-\ti k+1;w_{\ti k-1},\ldots,w_l}^{\{\ti\imath_0,\ti\imath_{\ti M}\};x_{n'},\ldots,x_n}{}_{\color{red}(b5)\color{black}}
\end{aligned}
\end{aligned}
\allowdisplaybreaks
\nonumber\\
&+\sum_{\begin{subarray}{l} 
n=n'+n'', \; n',n''\ge 0, \\ 
1<m<\ti M,\\ 
0=j_0<j_1<\cdots<j_m=\ti M, \\
k_0,\ldots,k_m\ge 0, \\ 
k_0=n', \; k_m=\ti k-1,\\
z_0^a,\ldots,z_{k_a}^a\in\sO^{i_{j_a}},\; 0<a<m
\end{subarray}}
\begin{aligned}[t]&
(-1)^{\begin{subarray}{l}(l-\ti k+1)(\ti k-1+m+\ti M) \\ +n'(l-\ti k+m)+\sum_{a=1}^{m-1}(j_a-a)\end{subarray}}\\
&\bigl(\d_{k_{m-1},k_m;w_0,\ldots,w_{\ti k-1}}^{\{\ti\imath_{j_{m-1}},\ldots,\ti\imath_{j_m}\};z_0^{m-1},\ldots,z_{k_{m-1}}^{m-1}}\\
&\ci\d_{k_{m-2},k_{m-1};z_0^{m-1},\ldots,z_{k_{m-1}}^{m-1}}^{\{\ti\imath_{j_{m-2}},\ldots,\ti\imath_{j_{m-1}}\};z_0^{m-2},\ldots,z_{k_{m-2}}^{m-2}}\\
&{}\ci\cdots\ci\d_{k_0,k_1;z_0^1,\ldots,z_{k_1}^1}^{\{\ti\imath_{j_0},\ldots,\ti\imath_{j_1}\};x_0,\ldots,x_{n'}}\bigr)\\
&{}\ot\d_{n'',l-\ti k+1;w_{\ti k-1},\ldots,w_l}^{\{\ti\imath_{j_0},\ti\imath_{j_1},\ti\imath_{j_2},\ldots,\ti\imath_{j_m}\};x_{n'},x_{n'+1},\ldots,x_n}{}_{\color{red}(b6)\color{black}}
\end{aligned}
\allowdisplaybreaks
\nonumber\\
&+\sum_{\begin{subarray}{l} 
n=n'+n'', \; n',n''\ge 0, \\ 
k_0,\ldots,k_{\ti M}\ge 0, \\ k_0=n', \; k_{\ti M}=\ti k-1, \\
w_0^i,\ldots,w_{k_i}^i\in\sO^{\ti\imath_i},\; 0<i<\ti M
\end{subarray}}\begin{aligned}[t]&
(-1)^{(l-\ti k+1)(\ti k-1)+n'(l-\ti k+\ti M)}\\
&\bigl(\d_{k_{\ti M-1},k_{\ti M};y_0,\ldots,y_{\ti k-1}}^{\{\ti\imath_{\ti M-1},\ti\imath_{\ti M}\};w_0^{\ti M-1},\ldots,w_{k_{\ti M-1}}^{\ti M-1}}\ci\d_{k_{\ti M-2},k_{\ti M-1};w_0^{\ti M-1},\ldots,w_{k_{\ti M-1}}^{\ti M-1}}^{\{\ti\imath_{\ti M-2},\ti\imath_{\ti M-1}\};w_0^{\ti M-2},\ldots,w_{k_{\ti M-2}}^{\ti M-2}}\\
&{}\ci\cdots\ci\d_{k_0,k_1;w_0^1,\ldots,w_{k_1}^1}^{\{\ti\imath_0,\ti\imath_1\};x_0,\ldots,x_{n'}}\bigr)\ot\D_{n'',l-\ti k+1;w_{\ti k-1},\ldots,w_l}^{\ti I;x_{n'},x_{n'+1},\ldots,x_n}{}_{\color{red}(b7)\color{black}}\biggr]\end{aligned}
\allowdisplaybreaks
\nonumber\\
&+\sum_{\begin{subarray}{l}
n=n'+n'', \; n',n''\ge 0, \\ 
0\le p\le q\le n', \\
1\le j\le \ti k-1, \\
k_0,\ldots,k_{\ti M}\ge 0, \\ k_0=p, \; k_{\ti M}=j-1, \\
w_0^i,\ldots,w_{k_i}^i\in\sO^{\ti\imath_i},\\ 
0<i<\ti M, \; m\ge 1\\
z_0,\ldots,z_m\in\sO^{\ti\imath_{\ti M}} 
\end{subarray}}
\begin{aligned}[t]&
(-1)^{(j+\ti k+p+1)\ti M}\cdot(-1)^{p+j\ti k+jq+\ti kq}\cdot(-1)^{\ti k-j-1+n'-q}\cdot{}\\
&\bigl(\d_{k_{\ti M-1},k_{\ti M};y_0,\ldots,y_{j-1}}^{\{\ti\imath_{\ti M-1},\ti\imath_{\ti M}\};w_0^{\ti M-1},\ldots,w_{k_{\ti M-1}}^{\ti M-1}}\ci\d_{k_{\ti M-2},k_{\ti M-1};w_0^{\ti M-1},\ldots,w_{k_{\ti M-1}}^{\ti M-1}}^{\{\ti\imath_{\ti M-2},\ti\imath_{\ti M-1}\};w_0^{\ti M-2},\ldots,w_{k_{\ti M-2}}^{\ti M-2}}\\
&{}\ci\cdots\ci\d_{k_0,k_1;w_0^1,\ldots,w_{k_1}^1}^{\{\ti\imath_0,\ti\imath_1\};x_0,\ldots,x_p}\bigr)\ot
\bigl(\d_{m,1;y_{j-1},y_j}^{\{\ti\imath_{\ti M}\};z_0,\ldots,z_m}\ci\D_{q-p,m;z_0,\ldots,z_m}^{\ti I;x_p,\ldots,x_q}\bigr)\\
&{}\ot\d_{n'-q,\ti k-j-1;y_j,\ldots,y_{\ti k-1}}^{\{\ti\imath_0,\ti\imath_M\};x_q,\ldots,x_{n'}}
\ot\d_{n'',1;y_{\ti k-1},y_{\ti k}}^{\{\ti\imath_0,\ti\imath_{\ti M}\};x_{n'},\ldots,x_n}{}_{\color{red}(b8)\color{black}}
\end{aligned}
\allowdisplaybreaks
\nonumber\\
&+\sum_{\begin{subarray}{l}
0\le p\le n, \\
k_0,\ldots,k_{\ti M}\ge 0, \\ k_0=p, \; k_{\ti M}=\ti k-1, \\
w_0^i,\ldots,w_{k_i}^i\in\sO^{\ti\imath_i},\\ 
0<i<\ti M, \; m\ge 1\\
z_0,\ldots,z_m\in\sO^{\ti\imath_{\ti M}} 
\end{subarray}}
\begin{aligned}[t]&
(-1)^{(p+1)\ti M}\cdot(-1)^{p+\ti k}\cdot{}\\
&\bigl(\d_{k_{\ti M-1},k_{\ti M};y_0,\ldots,y_{j-1}}^{\{\ti\imath_{\ti M-1},\ti\imath_{\ti M}\};w_0^{\ti M-1},\ldots,w_{k_{\ti M-1}}^{\ti M-1}}\ci\d_{k_{\ti M-2},k_{\ti M-1};w_0^{\ti M-1},\ldots,w_{k_{\ti M-1}}^{\ti M-1}}^{\{\ti\imath_{\ti M-2},\ti\imath_{\ti M-1}\};w_0^{\ti M-2},\ldots,w_{k_{\ti M-2}}^{\ti M-2}}\\
&{}\ci\cdots\ci\d_{k_0,k_1;w_0^1,\ldots,w_{k_1}^1}^{\{\ti\imath_0,\ti\imath_1\};x_0,\ldots,x_p}\bigr)\ot
\bigl(\d_{m,1;y_{\ti k-1},y_{\ti k}}^{\{\ti\imath_{\ti M}\};z_0,\ldots,z_m}\ci\D_{n-p,m;z_0,\ldots,z_m}^{\ti I;x_p,\ldots,x_n}\bigr){}_{\color{red}(b9)\color{black}}
\end{aligned}
\allowdisplaybreaks
\nonumber\\
&=\sum_{\begin{subarray}{l} 
n=n'+n'', \; n',n''\ge 0, \\ 
1<m<\ti M,\\ 
0=j_0<j_1<\cdots<j_m=\ti M, \\
k_0,\ldots,k_m\ge 0, \\ 
k_0=n', \; k_m=\ti k-1,\\
z_0^a,\ldots,z_{k_a}^a\in\sO^{i_{j_a}},\; 0<a<m
\end{subarray}}
\begin{aligned}[t]&
(-1)^{1+\ti k+\ti M+m+n'+mn'+\sum_{a=1}^{m-1}(j_a-a)}\\
&\bigl(\d_{k_{m-1},k_m;y_0,\ldots,y_{\ti k-1}}^{\{\ti\imath_{j_{m-1}},\ldots,\ti\imath_{j_m}\};z_0^{m-1},\ldots,z_{k_{m-1}}^{m-1}}\\
&\ci\d_{k_{m-2},k_{m-1};z_0^{m-1},\ldots,z_{k_{m-1}}^{m-1}}^{\{\ti\imath_{j_{m-2}},\ldots,\ti\imath_{j_{m-1}}\};z_0^{m-2},\ldots,z_{k_{m-2}}^{m-2}}\\
&{}\ci\cdots\ci\d_{k_0,k_1;z_0^1,\ldots,z_{k_1}^1}^{\{\ti\imath_{j_0},\ldots,\ti\imath_{j_1}\};x_0,\ldots,x_{n'}}\bigr)\ot{}\\
&\biggl[(-1)^m\sum_{\begin{subarray}{l} l\ge 0, \\ v_0,\ldots,v_l\in\sO^{\ti\imath_0} \end{subarray}}
\begin{aligned}[t]&\d_{l,1;y_{\ti k-1},y_{\ti k}}^{\{\ti\imath_{\ti M}\};v_0,\ldots,v_l}\ci{}\\
&\d_{n'',l;v_0,\ldots,v_l}^{\{\ti\imath_{j_0},\ti\imath_{j_1},\ti\imath_{j_2},\ldots,\ti\imath_{j_m}\};x_{n'},x_{n'+1},\ldots,x_n}{}_{\color{red}(c1)\color{black}}
\end{aligned}\\
&+\sum_{\begin{subarray}{l}
l=1,\ldots,n''+1,\\ v_0,\ldots,v_l\in\sO^{\ti\imath_0} 
\end{subarray}}\begin{aligned}[t]& \d_{l,1;y_{\ti k-1},y_{\ti k}}^{\{\ti\imath_{j_0},\ti\imath_{j_1},\ti\imath_{j_2},\ldots,\ti\imath_{j_m}\};v_0,\ldots,v_l}\ci{}\\
&\d_{n'',l;v_0,\ldots,v_l}^{\{\ti\imath_0\};x_{n'},x_{n'+1},\ldots,x_n}{}_{\color{red}(c2)\color{black}}\end{aligned}
\biggr]\end{aligned}
\allowdisplaybreaks
\nonumber\\
&+\sum_{\begin{subarray}{l} 
n=n'+n'', \; n',n''\ge 0, \\ 
1<m<\ti M,\\ 
0=j_0<j_1<\cdots<j_m=\ti M, \\
k_0,\ldots,k_m\ge 0, \; k_m=\ti k-1, \\
z_0^a,\ldots,z_{k_a}^a\in\sO^{i_{j_a}},\; 0\le a<m
\end{subarray}}
\begin{aligned}[t]&(-1)^{1+\ti k+\ti M+n'm+\sum_{a=1}^{m-1}(j_a-a)}\\
&\bigl(\d_{k_{m-1},k_m;y_0,\ldots,y_{\ti k-1}}^{\{\ti\imath_{j_{m-1}},\ldots,\ti\imath_{j_m}\};z_0^{m-1},\ldots,z_{k_{m-1}}^{m-1}}\\
&\ci\d_{k_{m-2},k_{m-1};z_0^{m-1},\ldots,z_{k_{m-1}}^{m-1}}^{\{\ti\imath_{j_{m-2}},\ldots,\ti\imath_{j_{m-1}}\};z_0^{m-2},\ldots,z_{k_{m-2}}^{m-2}}\\
&{}\ci\cdots\ci\d_{k_0,k_1;z_0^1,\ldots,z_{k_1}^1}^{\{\ti\imath_{j_0},\ldots,\ti\imath_{j_1}\};z_0^0,\ldots,z_{k_0}^0}
\ci \d_{n',k_0;z_0^0,\ldots,z_{k_0}^0}^{\{\ti\imath_0\};x_0,\ldots,x_{n'}}\bigr)\\
&\ot\d_{n'',1;y_{\ti k-1},y_{\ti k}}^{\{\ti\imath_{j_0},\ti\imath_{j_1},\ti\imath_{j_2},\ldots,\ti\imath_{j_m}\};x_{n'},x_{n'+1},\ldots,x_n}{}_{\color{red}(c3)\color{black}}
\end{aligned}
\allowdisplaybreaks
\nonumber\\
&+\sum_{\begin{subarray}{l} 
n=n'+n'', \; n',n''\ge 0, \\ 
1<m<\ti M,\\ 
0=j_0<j_1<\cdots<j_m=\ti M, \\
k_0,\ldots,k_m\ge 0, \; k_0=n',\\
z_0^a,\ldots,z_{k_a}^a\in\sO^{i_{j_a}},\; 0<a\le m
\end{subarray}}
\begin{aligned}[t]&(-1)^{1+\ti k+\ti M+n'm+\sum_{a=1}^{m-1}(j_a-a)}\cdot\bigl((-1)^{\ti M-m+1}\\
&\d_{k_m,\ti k-1;y_0,\ldots,y_{\ti k-1}}^{\{\ti\imath_{\ti M}\};z_0^m,\ldots,z_{k_m}^m}\ci\d_{k_{m-1},k_m;z_0^m,\ldots,z_{k_m}^m}^{\{\ti\imath_{j_{m-1}},\ldots,\ti\imath_{j_m}\};z_0^{m-1},\ldots,z_{k_{m-1}}^{m-1}}\\
&\ci\d_{k_{m-2},k_{m-1};z_0^{m-1},\ldots,z_{k_{m-1}}^{m-1}}^{\{\ti\imath_{j_{m-2}},\ldots,\ti\imath_{j_{m-1}}\};z_0^{m-2},\ldots,z_{k_{m-2}}^{m-2}}\\
&{}\ci\cdots\ci\d_{k_0,k_1;z_0^1,\ldots,z_{k_1}^1}^{\{\ti\imath_{j_0},\ldots,\ti\imath_{j_1}\};x_0,\ldots,x_{n'}}\bigr)\\
&{}\ot\d_{n'',1;y_{\ti k-1},y_{\ti k}}^{\{\ti\imath_{j_0},\ti\imath_{j_1},\ti\imath_{j_2},\ldots,\ti\imath_{j_m}\};x_{n'},x_{n'+1},\ldots,x_n}{}_{\color{red}(c4)\color{black}}
\end{aligned}
\allowdisplaybreaks
\nonumber\\
&+\!\!\!\sum_{\begin{subarray}{l} 
n=n'+n'', \; n',n''\ge 0, \\ 
1<m<\ti M,\\ 
0=j_0<j_1<\cdots<j_m=\ti M, \\
k_0,\ldots,k_m\ge 0, \\ 
k_0=n', \; k_m=\ti k-1, \\
z_0^a,\ldots,z_{k_a}^a\in\sO^{i_{j_a}},\\ 
0<a<m \\
0<l<m,\; \hat k\ge 0, \\
v_0,\ldots,v_{\hat k}\in\sO^{\ti\imath_{j_l}}
\end{subarray}}
\begin{aligned}[t]&(-1)^{1+\ti k+\ti M+n'm+\sum_{a=1}^{m-1}(j_a-a)}\\
&\bigl(\d_{k_{m-1},k_{m};y_0,\ldots,y_{\ti k-1}}^{\{\ti\imath_{j_{m-1}},\ldots,\ti\imath_{j_m}\};z_0^{m-1},\ldots,z_{k_{m-1}}^{m-1}}\ci\cdots\ci\d_{k_l,k_{l+1};z_0^{l+1},\ldots,z_{k_{l+1}}^{l+1}}^{\{\ti\imath_{j_l},\ldots,\ti\imath_{j_{l+1}}\};z_0^l,\ldots,z_{k_l}^l}\\
&\ci\bigl[(-1)^{l+j_l}\d_{\hat k,k_l;z_0^l,\ldots,z_{k_l}^l}^{\{\ti\imath_{j_l}\};v_0,\ldots,v_{\hat k}}-(-1)^{l+j_l}\d_{\hat k,k_l;z_0^l,\ldots,z_{k_l}^l}^{\{\ti\imath_{j_l}\};v_0,\ldots,v_{\hat k}}\bigr]\\
&\ci\d_{k_{l-1},\hat k;v_0,\ldots,v_{\hat k}}^{\{\ti\imath_{j_{l-1}},\ldots,\ti\imath_{j_l}\};z_0^{l-1},\ldots,z_{k_{l-1}}^{l-1}}\!\ci\!\d_{k_{l-2},k_{l-1};z_0^{l-1},\ldots,z_{k_{l-1}}^{l-1}}^{\{\ti\imath_{j_{l-2}},\ldots,\ti\imath_{j_{l-1}}\};z_0^{l-2},\ldots,z_{k_{l-2}}^{l-2}}\!\ci\!\cdots\\
&{}\ci\d_{k_0,k_1;z_0^1,\ldots,z_{k_1}^1}^{\{\ti\imath_{j_0},\ldots,\ti\imath_{j_1}\};x_0,\ldots,x_{n'}}\bigr)
\ot\d_{n'',1;y_{\ti k-1},y_{\ti k}}^{\{\ti\imath_{j_0},\ti\imath_{j_1},\ti\imath_{j_2},\ldots,\ti\imath_{j_m}\};x_{n'},x_{n'+1},\ldots,x_n}{}_{\color{red}(c5)\color{black}}
\end{aligned}
\allowdisplaybreaks
\nonumber\\
&+\sum_{1\le p\le n}\begin{aligned}[t]&
(-1)^{\ti k+\ti M}
\D_{p,\ti k-1;y_0,\ldots,y_{\ti k-1}}^{\ti I;x_0,\ldots,x_p}\!\!\ot\!
\biggl[
\sum_{1\le l\le n-p,\; w_0,\ldots,w_l\in\sO^{\ti\imath_0}\!\!\!\!\!\!\!\!\!\!\!\!\!\!\!\!\!\!\!\!\!\!\!\!\!\!\!\!\!\!\!\!\!\!\!\!\!\!\!\!\!\!\!\!\!\!\!\!\!\!}\d_{l,1;y_{\ti k-1},y_{\ti k}}^{\{\ti\imath_0,\ti\imath_M\};w_0,\ldots,w_l}
\ci\d_{n-p,l;w_0,\ldots,w_l}^{\{\ti\imath_0\};x_p,\ldots,x_n}{}_{\color{red}(c6)\color{black}}\\
&-\sum_{1\le l\le n-p,\; w_0,\ldots,w_l\in\sO^{\ti\imath_{\ti M}}\!\!\!\!\!\!\!\!\!\!\!\!\!\!\!\!\!\!\!\!\!\!\!\!\!\!\!\!\!\!\!\!\!\!\!\!\!\!}\d_{l,1;y_{\ti k-1},y_{\ti k}}^{\{\ti\imath_{\ti M}\};w_0,\ldots,w_l}\ci\d_{n-p,l;w_0,\ldots,w_l}^{\{\ti\imath_0,\ti\imath_M\};x_p,\ldots,x_n}{}_{\color{red}(c7)\color{black}}
\biggr]
\end{aligned}
\allowdisplaybreaks
\nonumber\\
&+\sum_{\begin{subarray}{l} 
n=n'+n'', \\ n',n''\ge 0 \end{subarray}}
\begin{aligned}[t]&
(-1)^{\ti k+\ti M+n'-1}
\biggl[\sum_{\begin{subarray}{l} 1\le l\le n', \\ w_0,\ldots,w_l\in\sO^{\ti\imath_0} \end{subarray}}
\D_{l,\ti k-1;y_0,\ldots,y_{\ti k-1}}^{\ti I;w_0,\ldots,w_l}\ci\d_{n',l;w_0,\ldots,w_l}^{ 
\{\ti\imath_0\};x_0,\ldots,x_{n'}}{}_{\color{red}(c8)\color{black}}\end{aligned}
\allowdisplaybreaks
\nonumber\\
&+(-1)^{\ti M}\sum_{\begin{subarray}{l} \ti k-1\le m\le n'+\ti M-1, \\ z_0,\ldots,z_m\in\sO^{\ti\imath_{\ti M}} \end{subarray}}
\d_{m,\ti k-1;y_0,\ldots,y_{\ti k-1}}^{\{\ti\imath_{\ti M}\};z_0,\ldots,z_m}\ci\D_{n',m;z_0,\ldots,z_m}^{\ti I;x_0,\ldots,x_{n'}}{}_{\color{red}(c9)\color{black}}
\allowdisplaybreaks
\nonumber\\
&+\sum_{\begin{subarray}{l}
n=n'+n'', \; n',n''\ge 0, \\ 
0\le p\le q\le n', \\
1\le j\le k-1, \\
k_0,\ldots,k_{\ti M}\ge 0, \\ k_0=p, \; k_{\ti M}=j-1, \\
w_0^i,\ldots,w_{k_i}^i\in\sO^{\ti\imath_i},\\ 
0<i<\ti M, \; m\ge 1\\
z_0,\ldots,z_m\in\sO^{\ti\imath_{\ti M}} 
\end{subarray}}
\begin{aligned}[t]&
(-1)^{(j+\ti k+p)\ti M}\cdot(-1)^{p+j(\ti k-1)+jq+(\ti k-1)q}\cdot{}\\
&\bigl(\d_{k_{\ti M-1},k_{\ti M};y_0,\ldots,y_{j-1}}^{\{\ti\imath_{\ti M-1},\ti\imath_{\ti M}\};w_0^{\ti M-1},\ldots,w_{k_{\ti M-1}}^{\ti M-1}}\ci\d_{k_{\ti M-2},k_{\ti M-1};w_0^{\ti M-1},\ldots,w_{k_{\ti M-1}}^{\ti M-1}}^{\{\ti\imath_{\ti M-2},\ti\imath_{\ti M-1}\};w_0^{\ti M-2},\ldots,w_{k_{\ti M-2}}^{\ti M-2}}\\
&{}\ci\cdots\ci\d_{k_0,k_1;w_0^1,\ldots,w_{k_1}^1}^{\{\ti\imath_0,\ti\imath_1\};x_0,\ldots,x_p}\bigr)\ot
\bigl(\d_{m,1;y_{j-1},y_j}^{\{\ti\imath_{\ti M}\};z_0,\ldots,z_m}\ci\D_{q-p,m;z_0,\ldots,z_m}^{\ti I;x_p,\ldots,x_q}\bigr)\\
&{}\ot\d_{n'-q,\ti k-j-1;y_j,\ldots,y_{\ti k-1}}^{\{\ti\imath_0,\ti\imath_M\};x_q,\ldots,x_{n'}}
{}_{\color{red}(c10)\color{black}}\biggr]\ot\d_{n'',1;y_{\ti k-1},y_{\ti k}}^{\{\ti\imath_0,\ti\imath_{\ti M}\};x_{n'},\ldots,x_n}\end{aligned}
\nonumber\\
&=0.
\nonumber
\ea

Here to make the equation easier to navigate, we have labelled terms in the first step $(a1),(a2)$, and so on, where $(b1),(b2)$ are split into $(i),(ii)$. To prove the first step of \eq{ai10eq19}, note that:
\begin{itemize}
\setlength{\itemsep}{0pt}
\setlength{\parsep}{0pt}
\item[(a)] Term $(a1)$ is rewritten as $((b1)+(b2))((i)+(ii))$ by first applying \eq{ai10eq15} to $\D_{l,\ti k;y_0,\ldots,y_{\ti k}}^{\ti I;w_0,\ldots,w_l}$ with $k'=\ti k-1$ and $k''=1$ and noting that $\D_{n'',1;y_{\ti k-1},y_{\ti k}}^{\ti I;x_{n'},x_{n'+1},\ldots,x_n}=0$, which yields a sum over splittings $l=l'+l''$, and then applying \eq{ai7eq12} to $\d_{n,l;w_0,\ldots,w_l}^{\{\ti\imath_0\};x_0,\ldots,x_n}$ splitting $l=l'+l''$. 
\item[(b)] Term $(a2)$ is rewritten as $(b3)+\cdots+(b7)$ by applying \eq{ai7eq12} to $\d_{l,\ti k;y_0,\ldots,y_{\ti k}}^{\{\ti\imath_M\};w_0,\ldots,w_l}$ with $k'=\ti k-1$ and $k''=1$, splitting this into two terms, and then applying \eq{ai10eq15} to $\D_{n,l;w_0,\ldots,w_l}^{\ti I;x_0,\ldots,x_n}$, firstly splitting $l=(l-1)+1$ giving $(b3)+(b4)$, and secondly splitting $l=(\ti k-1)+(l-\ti k+1)$, giving $(b5)+(b6)+(b7)$.
\item[(c)] Term $(a3)$ is split into the sum over $j<\ti k$, which is rewritten by decomposing $\d_{n-q,\ti k-j;y_j,\ldots,y_{\ti k}}^{\{\ti\imath_0,\ti\imath_M\};x_q,\ldots,x_n}$ using \eq{ai7eq13} with $k'=\ti k-j-1$ and $k''=1$ to give $(b8)$, and the sum over $j=\ti k$ (forcing $q=n$), which is~$(b9)$.
\end{itemize}

To prove the second step of \eq{ai10eq19}, observe that:
\begin{itemize}
\setlength{\itemsep}{0pt}
\setlength{\parsep}{0pt}
\item[(d)] By writing out compositions, with appropriate sign changes upon exchanging order of factors, and relabelling some variables, we see that $(c1)=(b6)$, $(c2)=(b2)(ii)$, $(c3)=(b2)(i)$, $(c4)=(b4)$, $(c6)=(b1)(ii)$, $(c7)=(b5)$, $(c8)=(b1)(i)$, $(c9)=(b3)$, $(c10)=(b8)$.
\item[(e)] Clearly $(c5)=0$.
\item[(f)] $(b7)$ and $(b9)$ cancel out. 
\end{itemize}

To prove the third step of \eq{ai10eq19}, note that:
\begin{itemize}
\setlength{\itemsep}{0pt}
\setlength{\parsep}{0pt}
\item[(g)] $(c1)+(c2)=0$ by \eq{ai8eq6} for $\es\ne I\subsetneq\ti I$.
\item[(h)] $(c3)+(c4)+(c5)=0$ as in the final step of \eq{ai10eq9}.
\item[(i)] $(c6)+(c7)=0$ by \eq{ai8eq6} for $\es\ne I\subsetneq\ti I$.
\item[(j)] $(c8)+(c9)+(c10)=0$ by \eq{ai10eq16} with $k=\ti k-1$, which holds by induction.
\end{itemize}
Thus \eq{ai10eq19} holds. This completes the inductive step, so \eq{ai10eq16} holds by induction.

We then deduce \eq{ai10eq17} from \eq{ai10eq16} using the identities
\ea
&\d_{n,k;y_0,y_1,\ldots,y_k}^{\{\ti\imath_0,\ti\imath_M\};x_0,x_1,\ldots,x_n}
\label{ai10eq20}\\
&=\sum_{\begin{subarray}{l}
0=l_0\le l_1\le \cdots\le l_k=n:\\
y_i=f^{\ti\imath_0\ti\imath_{\ti M}}(x_{l_i}),\; 0\le i\le k
\end{subarray}}
\begin{aligned}[t]&
(-1)^{\sum_{i=1}^k(l_i-l_{i-1}+1)(k-i)} \\
&\fF_{l_1-l_0}^{\ti\imath_0\ti\imath_{\ti M};x_{l_0},\ldots,x_{l_1}}\ot\cdots\ot\fF_{l_k-l_{k-1}}^{\ti\imath_0\ti\imath_{\ti M};x_{l_{k-1}},\ldots,x_{l_k}},
\end{aligned}
\nonumber\\
&\sum_{\begin{subarray}{l} 
k_0,\ldots,k_{\ti M}\ge 0, \\ k_0=n, \; k_{\ti M}=k, \\
w_0^i,\ldots,w_{k_i}^i\in\sO^{\ti\imath_i},\; 0<i<\ti M
\end{subarray}}\begin{aligned}[t]
&\d_{k_{\ti M-1},k_{\ti M};y_0,\ldots,y_k}^{\{\ti\imath_{\ti M-1},\ti\imath_{\ti M}\};w_0^{\ti M-1},\ldots,w_{k_{\ti M-1}}^{\ti M-1}}\ci\d_{k_{\ti M-2},k_{\ti M-1};w_0^{\ti M-1},\ldots,w_{k_{\ti M-1}}^{\ti M-1}}^{\{\ti\imath_{\ti M-2},\ti\imath_{\ti M-1}\};w_0^{\ti M-2},\ldots,w_{k_{\ti M-2}}^{\ti M-2}}\\
&{}\ci\cdots\ci\d_{k_0,k_1;w_0^1,\ldots,w_{k_1}^1}^{\{\ti\imath_0,\ti\imath_1\};x_0,\ldots,x_n}
\end{aligned}
\nonumber\\
&=\sum_{\begin{subarray}{l} 
0=l_0\le l_1\le \cdots\le l_k=n:\\
y_i=f^{\ti I}(x_{l_i}),\; 0\le i\le k
\end{subarray}}\begin{aligned}[t]&
(-1)^{\sum_{i=1}^k(l_i-l_{i-1}+1)(k-i)} \\
&\fF_{l_1-l_0}^{\ti I;x_{l_0},\ldots,x_{l_1}}\ot\fF_{l_2-l_1}^{\ti I;x_{l_1},\ldots,x_{l_2}}\ot\cdots\ot\fF_{l_k-l_{k-1}}^{\ti I;x_{l_{k-1}},\ldots,x_{l_k}}.
\end{aligned}
\label{ai10eq21}
\ea
Here \eq{ai10eq20} follows from \eq{ai7eq11} and Example \ref{ai7ex2}(b), and \eq{ai10eq21} can be proved by composing $k$ copies of \eq{ai10eq20}, using $\fF^{\ti I}=\fF^{\ti\imath_{\ti M-1}\ti\imath_{\ti M}}\ci\cdots\ci\fF^{\ti\imath_0\ti\imath_1}$ and \eq{ai4eq10}, and computing signs. This proves Proposition~\ref{ai10prop1}.
\end{proof}

\begin{prop}
\label{ai10prop2}
Work in the situation of Definition\/ {\rm\ref{ai10def2}}. Then:
\smallskip

\noindent{\bf(a)} Equation \eq{ai8eq5} for $I=\ti I$ and\/ $n,x_0,\ldots,x_n$ is equivalent to
\e
\bigl(\n_1(\m^{\ti I})+{\rm E}^{\ti I}\bigr)_n^{x_0,\ldots,x_n}=0.
\label{ai10eq22}
\e
Here $\m^{\ti I},{\rm E}^{\ti I}$ are as in {\rm\eq{ai10eq11}, \eq{ai10eq13},} and\/ $\n_1$ as in\/ \eq{ai4eq15} and\/~\eq{ai10eq10}.
\smallskip

\noindent{\bf(b)} Using that\/ \eq{ai8eq4}--\eq{ai8eq5} hold for all\/ $\es\ne I\subsetneq\ti I,$ we show that
\e
\n_1({\rm E}^{\ti I})=0\quad\text{in $F^0\Hom_{A^\rfi_\iy}^\un(\fF^{\ti I},\fF^{\ti\imath_0\ti\imath_{\ti M}})^{2-\ti M}$.}
\label{ai10eq23}
\e
Hence we can form the cohomology class
\e
[{\rm E}^{\ti I}]\quad\text{in}\quad H^{3-\ti M}\bigl(F^0\Hom_{A^\rfi_\iy}^\un(\fF^{\ti I},\fF^{\ti\imath_0\ti\imath_{\ti M}}),\n_1\bigr).
\label{ai10eq24}
\e
Thus, leaving $\m^I$ for $\es\ne I\subsetneq\ti I$ fixed, we can change the definition of\/ $\m^{\ti I}$ if necessary to make \eq{ai8eq5} hold with\/ $I=\ti I$ and for all $x_0,\ldots,x_n,$ $n\ge 0,$ if and only if the cohomology class \eq{ai10eq24} is zero.
\end{prop}

\begin{proof}
For (a), for $n\ge 0$ and $x_0,\ldots,x_n\in\sO^{\ti\imath_0}$ we have
\begin{align*}
&\bigl(\n_1(\m^{\ti I})\bigr)_n^{x_0,\ldots,x_n}\\
&=\sum_{\begin{subarray}{l}
0=l_0\le l_1\le\cdots\le l_k=n,\\
1\le j\le k: \\
y_i=f^{\ti I}(x_{l_i}),\; 0\le i<j,\\
y_i=f^{\ti\imath_0\ti\imath_{\ti M}}(x_{l_i}),\; j\le i\le k \end{subarray}}
\begin{aligned}[t]&
(-1)^{(k+1)(\ti M+1)}\cdot(-1)^{\sum_{i=1}^{j-1}(l_i-l_{i-1}+1)(\ti M+1)}\\
&\cdot(-1)^{\sum_{i=1}^k(l_i-l_{i-1}+1)(k-i)} \\
&\m^{\{\ti\imath_{\ti M}\};y_0,\ldots,y_k}_k\ci{} \\ 
&\bigl(\fF_{l_1-l_0}^{\ti I;x_{l_0},\ldots,x_{l_1}}\ot\fF_{l_2-l_1}^{\ti I;x_{l_1},\ldots,x_{l_2}}\ot\cdots\ot\fF_{l_{j-1}-l_{j-2}}^{\ti I;x_{l_{j-2}},\ldots,x_{l_{j-1}}} \\
&{}\ot
\m_{l_j-l_{j-1}}^{\ti I;x_{l_{j-1}},\ldots,x_{l_j}}\\
&{}\ot\fF_{l_{j+1}-l_j}^{\ti\imath_0\ti\imath_{\ti M};x_{l_j},\ldots,x_{l_{j+1}}}\ot\cdots\ot\fF_{l_k-l_{k-1}}^{\ti\imath_0\ti\imath_{\ti M};x_{l_{k-1}},\ldots,x_{l_k}}\bigr)
\end{aligned}\\
&+\sum_{\begin{subarray}{l} i,j,k\ge 0:\\ i+j+k=n \end{subarray}}(-1)^{\ti M}(-1)^{i+jk}\m^{\ti I;x_0,\ldots,x_i,x_{i+j},\ldots,x_n}_{i+k+1}\ci\bigl(\id^{\ot^i}\ot\m_j^{\{\ti\imath_0\};x_i,\ldots,x_{i+j}}\ot \id^{\ot^k}\bigr)
\allowdisplaybreaks
\nonumber\\
&=\sum_{k\ge 0, \; y_0,\ldots,y_k\in\sO^{\ti\imath_{\ti M}}}
\m^{\{\ti\imath_{\ti M}\};y_0,\ldots,y_k}_k\ci\bigl(\d_{n,k;y_0,\ldots,y_k}^{\ti I;x_0,\ldots,x_n}-\D_{n,k;y_0,\ldots,y_k}^{\ti I;x_0,\ldots,x_n}\bigr)\\
&+\sum_{k\ge 0, \; y_0,\ldots,y_k\in\sO^{i_0}}
(-1)^{\ti M}\m_k^{\ti I;y_0,\ldots,y_k}\ci\d_{n,k;y_0,\ldots,y_k}^{\{\ti\imath_0\};x_0,\ldots,x_n}\\
&=-\bigl({\rm E}^{\ti I}\bigr)_n^{x_0,\ldots,x_n},
\end{align*}
using equation \eq{ai4eq15} in the first step, \eq{ai10eq14} and \eq{ai7eq8} in the second, and \eq{ai8eq5} with $I=\ti I$ and \eq{ai10eq12} in the third. This proves~\eq{ai10eq22}.

For (b), for $n\ge 0$ and $x_0,\ldots,x_n\in\sO^{\ti\imath_0}$ consider the equation
\begin{align*}
&\bigl(\n_1({\rm E}^{\ti I})\bigr)_n^{x_0,\ldots,x_n}\\
&=\sum_{\begin{subarray}{l}
0=l_0\le l_1\le\cdots\le l_k=n,\\
1\le j\le k: \\
y_i=f^{\ti I}(x_{l_i}),\; 0\le i<j,\\
y_i=f^{\ti\imath_0\ti\imath_{\ti M}}(x_{l_i}),\; j\le i\le k, \\
m\ge 0,\; z_0,\ldots,z_m\in\sO^{\ti\imath_{\ti M}} 
\end{subarray}}
\begin{aligned}[t]&
(-1)^{(k+1)\ti M}\cdot(-1)^{\sum_{i=1}^{j-1}(l_i-l_{i-1}+1)\ti M}\cdot{}\\
&(-1)^{\sum_{i=1}^k(l_i-l_{i-1}+1)(k-i)}\cdot \m^{\{\ti\imath_{\ti M}\};y_0,\ldots,y_k}_k\ci{} \\ 
&\bigl(\fF_{l_1-l_0}^{\ti I;x_{l_0},\ldots,x_{l_1}}\ot\fF_{l_2-l_1}^{\ti I;x_{l_1},\ldots,x_{l_2}}\ot\cdots\ot\fF_{l_{j-1}-l_{j-2}}^{\ti I;x_{l_{j-2}},\ldots,x_{l_{j-1}}} \\
&{}\ot
\bigl(\m_m^{\{\ti\imath_{\ti M}\};z_0,\ldots,z_m}\ci\D_{l_j-l_{j-1},m;z_0,\ldots,z_m}^{\ti I;x_{l_{j-1}},\ldots,x_{l_j}}\bigr)\\
&{}\ot\fF_{l_{j+1}-l_j}^{\ti\imath_0\ti\imath_{\ti M};x_{l_j},\ldots,x_{l_{j+1}}}\ot\cdots\ot\fF_{l_k-l_{k-1}}^{\ti\imath_0\ti\imath_{\ti M};x_{l_{k-1}},\ldots,x_{l_k}}\bigr)
\end{aligned}
\nonumber\\
&+\sum_{\begin{subarray}{l} 
k\ge 0,\; y_0,\ldots,y_k\in\sO^{\ti\imath_{\ti M}}, \\
l\ge 0,\; z_0,\ldots,z_l\in\sO^{\ti\imath_0}
\end{subarray}}
(-1)^{\ti M+1}\m_k^{\{\ti\imath_{\ti M}\};y_0,\ldots,y_k}\ci\D_{m,k;y_0,\ldots,y_k}^{\ti I;z_0,\ldots,z_m}\ci \d^{\{i_0\};x_0,\ldots,x_n}_{n,m;z_0,\ldots,z_m}
\allowdisplaybreaks
\nonumber\\
&=\sum_{\begin{subarray}{l} 
k\ge 0,\; y_0,\ldots,y_k\in\sO^{\ti\imath_{\ti M}}, \\
l\ge 0,\; z_0,\ldots,z_l\in\sO^{\ti\imath_0}
\end{subarray}}
 \m^{\{\ti\imath_{\ti M}\};y_0,\ldots,y_k}_k\ci\d_{l,k;y_0,\ldots,y_k}^{\{\ti\imath_{\ti M}\};w_0,\ldots,w_l}\ci\D_{n,l;w_0,\ldots,w_l}^{\ti I;x_0,\ldots,x_n}=0,
\end{align*}
using equations \eq{ai4eq15}, \eq{ai7eq8} and \eq{ai10eq12} in the first step, \eq{ai10eq17} in the second, and \eq{ai7eq15} in the third. This proves \eq{ai10eq23}. The last part is immediate.
\end{proof}

For our next proposition we need a different set-up to Definition \ref{ai10def2}. In the definitions of ${\rm G}_{n,k;y_0,\ldots,y_k}^{\ti I;x_0,\ldots,x_n},{\rm H}_{n,k;y_0,\ldots,y_k}^{\ti I;x_0,\ldots,x_n}$ below, the exclusion of terms in $\ti I\sm\{\ti\imath_0\}$ is not useful for the proof of Theorem \ref{ai8thm2}, but it will be used in \S\ref{ai1054} to complete the proof of Theorem~\ref{ai8thm3}.

\begin{dfn}
\label{ai10def3}
Continue in the situation of Definition~\ref{ai10def1}.

Let $\es\ne\ti I=\{\ti\imath_0<\ti\imath_1<\cdots< \ti\imath_{\ti M}\}\subseteq\{0,\ldots,N\}$ for $\ti M\ge 3$. We assume \eq{ai8eq4}--\eq{ai8eq5} hold for all $\es\ne I\subseteq\ti I$ with $I\ne\ti I$ and $I\ne\ti I\sm\{\ti\imath_a\}$ for $a=0,\ldots,\ti M-1$ and all $x_0,\ldots,x_n$ for $n\ge 0$, but we do {\it not\/} yet assume \eq{ai8eq5} when $I=\ti I$ or when $I=\ti I\sm\{\ti\imath_a\}$ for $0\le a<\ti M$.

Use the notation $\fF^{\ti\imath_a\ti\imath_b}$ and $\fF^{\ti J}$ for $\es\ne\ti J\subseteq\ti I$ as in Definition \ref{ai10def2}.

Observe that for all $0<a<\ti M$, Definition \ref{ai10def2} holds with $\ti I\sm\{\ti\imath_a\}$ in place of $\ti I$. Thus Proposition \ref{ai10prop2}(b) defines a cohomology class 
\begin{equation*}
\bigl[{\rm E}^{\ti I\sm\{\ti\imath_a\}}\bigr]\quad\text{in}\quad H^{4-\ti M}\bigl(F^0\Hom_{A^\rfi_\iy}^\un(\fF^{\ti I\sm\{\ti\imath_a\}},\fF^{\ti\imath_0\ti\imath_{\ti M}}),\n_1\bigr),
\end{equation*}
where $\fF^{\ti I\sm\{\ti\imath_a\}}:\sA^{\ti\imath_0}\ra\sA^{\ti\imath_{\ti M}}$ is given by
\begin{equation*}
\fF^{\ti I\sm\{\ti\imath_a\}}=\fF^{\ti\imath_{\ti M-1}\ti\imath_{\ti M}}\ci\cdots\ci\fF^{\ti\imath_{a+1}\ti\imath_{a+2}}\ci\fF^{\ti\imath_{a-1}\ti\imath_{a+1}}\ci\fF^{\ti\imath_{a-2}\ti\imath_{a-1}}\ci\cdots\ci\fF^{\ti\imath_0\ti\imath_1}. 
\end{equation*}
We will regard $\bigl[{\rm E}^{\ti I\sm\{\ti\imath_a\}}\bigr]$ as a degree $4-\ti M$ morphism
\e
\bigl[{\rm E}^{\ti I\sm\{\ti\imath_a\}}\bigr]:\fF^{\ti I\sm\{\ti\imath_a\}}\longra\fF^{\ti\imath_0\ti\imath_{\ti M}}
\label{ai10eq25}
\e
in the graded cohomology category $H^*\bigl(F^0\Fun_{A^\rfi_\iy}^\un(\sA^{\ti\imath_0},\ab\sA^{\ti\imath_{\ti M}})\bigr)$.

Similarly, when $a=0$ we have a degree $4-\ti M$ morphism
\begin{equation*}
\bigl[{\rm E}^{\ti I\sm\{\ti\imath_0\}}\bigr]:\fF^{\ti I\sm\{\ti\imath_0\}}\longra\fF^{\ti\imath_1\ti\imath_{\ti M}}.
\end{equation*}

As in \eq{ai10eq11}, the data $\m_k^{\{\ti\imath_{a-1},\ti\imath_a,\ti\imath_{a+1}\};x_0,\ldots,x_k}$ for all $k\ge 0$ and $x_0,\ldots,x_k$ in $\sO^{\ti\imath_{a-1}}$ defines an element
\begin{equation*}
\m^{\{\ti\imath_{a-1},\ti\imath_a,\ti\imath_{a+1}\}}\in F^0\Hom_{A^\rfi_\iy}^\un(\fF^{\{\ti\imath_{a-1},\ti\imath_a,\ti\imath_{a+1}\}},\fF^{\ti\imath_{a-1}\ti\imath_{a+1}})^0.
\end{equation*}
Equation \eq{ai8eq5} for $I=\{\ti\imath_{a-1},\ti\imath_a,\ti\imath_{a+1}\}$ is equivalent to $\n_1(\m^{\{\ti\imath_{a-1},\ti\imath_a,\ti\imath_{a+1}\}})=0$, so we have a cohomology class
\begin{equation*}
\bigl[\m^{\{\ti\imath_{a-1},\ti\imath_a,\ti\imath_{a+1}\}}\bigr]\in H^0\bigl(F^0\Hom_{A^\rfi_\iy}^\un(\fF^{\{\ti\imath_{a-1},\ti\imath_a,\ti\imath_{a+1}\}},\fF^{\ti\imath_{a-1}\ti\imath_{a+1}}),\n_1\bigr).
\end{equation*}
This is degree 0 morphism 
\e
\bigl[\m^{\{\ti\imath_{a-1},\ti\imath_a,\ti\imath_{a+1}\}}\bigr]:\fF^{\{\ti\imath_{a-1},\ti\imath_a,\ti\imath_{a+1}\}}\longra\fF^{\ti\imath_{a-1}\ti\imath_{a+1}}
\label{ai10eq26}
\e
in the graded cohomology category $H^*\bigl(F^0\Fun_{A^\rfi_\iy}^\un(\sA^{\ti\imath_{a-1}},\ab\sA^{\ti\imath_{a+1}})\bigr)$.

As in Definitions \ref{ai8def3}(a) and \ref{ai10def1}, we require $\m^{\{\ti\imath_{a-1},\ti\imath_a,\ti\imath_{a+1}\}}$ to be a filtered $A_\iy$-natural isomorphism, in the sense of Definition \ref{ai5def12}(c). This means that \eq{ai10eq26} is an isomorphism in~$H^*\bigl(F^0\Fun_{A^\rfi_\iy}^\un(\sA^{\ti\imath_{a-1}},\ab\sA^{\ti\imath_{a+1}})\bigr)$.

Now Definition \ref{ai4def9}(a),(b) and Definition \ref{ai5def10}(b),(c) explain how to define $A_\iy$-functors between $A_\iy$-categories $F^0\Fun_{A^\rfi_\iy}^\un(\cdots)$ by precomposing with $A_\iy$-functors $\fE:\sA_1\ra\sA_2$ and postcomposing with $A_\iy$-functors $\fG:\sB_1\ra\sB_2$. We will apply this with
\begin{align*}
\fE&=\fF^{\{\ti\imath_0,\ldots,\ti\imath_{a-1}\}}=\fF^{\ti\imath_{a-2}\ti\imath_{a-1}}\ci\cdots\ci\fF^{\ti\imath_0\ti\imath_1}:
\sA^{\ti\imath_0}\longra\sA^{\ti\imath_{a-1}},\\
\fG&=\fF^{\{\ti\imath_{a+1},\ldots,\ti\imath_{\ti M}\}}=\fF^{\ti\imath_{\ti M-1}\ti\imath_{\ti M}}\ci\cdots\ci\fF^{\ti\imath_{a+1}\ti\imath_{a+2}}:
\sA^{\ti\imath_{a+1}}\longra\sA^{\ti\imath_{\ti M}}.
\end{align*}
This yields 
\begin{align*}
&(\fF^{\{\ti\imath_0,\ldots,\ti\imath_{a-1}\}})^*\ci
(\fF^{\{\ti\imath_{a+1},\ldots,\ti\imath_{\ti M}\}})_*
(\m^{\{\ti\imath_{a-1},\ti\imath_a,\ti\imath_{a+1}\}})\in F^0\Hom_{A^\rfi_\iy}^\un(\fF^{\ti I},\fF^{\ti I\sm\{\ti\imath_a\}})^0.
\end{align*}
Since $\fE^*,\fG_*$ are $A_\iy$-functors, this is closed under $\n_1$, and defines a degree 0 morphism in $H^*\bigl(F^0\Fun_{A^\rfi_\iy}^\un(\sA^{\ti\imath_0},\ab\sA^{\ti\imath_{\ti M}})\bigr)$
\e
\bigl[(\fF^{\{\ti\imath_0,\ldots,\ti\imath_{a-1}\}})^*\ci
(\fF^{\{\ti\imath_{a+1},\ldots,\ti\imath_{\ti M}\}})_*
(\m^{\{\ti\imath_{a-1},\ti\imath_a,\ti\imath_{a+1}\}})\bigr]:\fF^{\ti I}\longra\fF^{\ti I\sm\{\ti\imath_a\}}.
\label{ai10eq27}
\e
As \eq{ai10eq26} is an isomorphism, so is \eq{ai10eq27}.

For each $a=1,\ldots,\ti M-1$, we can now take the composition of \eq{ai10eq27} and \eq{ai10eq25}, to get a degree $4-\ti M$ morphism in $H^*\bigl(F^0\Fun_{A^\rfi_\iy}^\un(\sA^{\ti\imath_0},\ab\sA^{\ti\imath_{\ti M}})\bigr)$:
\begin{align*}
&\bigl[{\rm E}^{\ti I\sm\{\ti\imath_a\}}\bigr]\ci\bigl[(\fF^{\{\ti\imath_0,\ldots,\ti\imath_{a-1}\}})^*\ci
(\fF^{\{\ti\imath_{a+1},\ldots,\ti\imath_{\ti M}\}})_*
(\m^{\{\ti\imath_{a-1},\ti\imath_a,\ti\imath_{a+1}\}})\bigr]:\\
&\qquad\qquad \fF^{\ti I}\longra\fF^{\ti\imath_0\ti\imath_{\ti M}}.
\end{align*}

Let $x_0,\ldots,x_n\in\sO^{\ti\imath_0}$ and $y_0,\ldots,y_k\in\sO^{\ti\imath_{\ti M}}$ for $n,k\ge 0$. Then Definition \ref{ai7def2} defines an operator $\d_{n,k;y_0,\ldots,y_k}^{\ti I;x_0,\ldots,x_n}$ in \eq{ai7eq6} and \eq{ai7eq9} in terms of the data $\m_l^{I;w_0,\ldots,w_l}$ for $I\subseteq\ti I$. In a similar way to $\D_{n,k;y_0,\ldots,y_k}^{\ti I;x_0,\ldots,x_n}$ in Definition \ref{ai10def2}, write ${\rm G}_{n,k;y_0,\ldots,y_k}^{\ti I;x_0,\ldots,x_n}$ for the sum of all terms in the definition \eq{ai7eq9} of $\d_{n,k;y_0,\ldots,y_k}^{\ti I;x_0,\ldots,x_n}$ which do {\it not\/} involve any operator $\m_l^{I;w_0,\ldots,w_l}$ for $I=\ti I$ or $I=\ti I\sm\{\ti\imath_a\}$ for~$a=0,\ldots,\ti M-1$. 

As for ${\rm E}_n^{\ti I;x_0,\ldots,x_n}$ in \eq{ai10eq12}, for $n\ge 0$ and $x_0,\ldots,x_n\in\sO^{\ti\imath_0}$, define
\e
{\rm H}_n^{\ti I;x_0,\ldots,x_n}
=\sum_{\begin{subarray}{l} k=1,\ldots,n+\ti M-1,\\
y_0,\ldots,y_k\in\sO^{\ti\imath_{\ti M}} \end{subarray}\!\!\!\!\!\!\!\!\!\!\!}\m_k^{\{\ti\imath_{\ti M}\};y_0,\ldots,y_k}\ci{\rm G}_{n,k;y_0,\ldots,y_k}^{\ti I;x_0,\ldots,x_n}
\label{ai10eq28}
\e
in $\sH^{\ti\imath_{\ti M}}(f^{\ti I}(x_0),f^{\ti\imath_0\ti\imath_{\ti M}}(x_n))^{3-\ti M-n}$. As for \eq{ai10eq13}, we have
\e
{\rm H}^{\ti I}=\bigl({\rm H}_n^{\ti I;x_0,\ldots,x_n})^{x_0,\ldots,x_n\in\sO^{\ti\imath_0}}_{n\ge 0}\in F^0\Hom_{A^\rfi_\iy}^\un(\fF^{\ti I},\fF^{\ti\imath_0\ti\imath_{\ti M}})^{3-\ti M}.
\label{ai10eq29}
\e
\end{dfn}

\begin{prop}
\label{ai10prop3}
Work in the situation of Definition\/ {\rm\ref{ai10def3}}. Then:
\smallskip

\noindent{\bf(a)} For each $a=1,\ldots,\ti M-1,$ we have
\ea
{\rm E}^{\ti I}&=(-1)^{\ti M-a-1}\n_2\bigl( (\fF^{\{\ti\imath_0,\ldots,\ti\imath_{a-1}\}})^*\ci
(\fF^{\{\ti\imath_{a+1},\ldots,\ti\imath_{\ti M}\}})_*
(\m^{\{\ti\imath_{a-1},\ti\imath_a,\ti\imath_{a+1}\}}),
\m^{\ti I\sm\{\ti\imath_a\}}\bigr)
\nonumber\\
&+\text{terms not involving\/ $\m^{\ti I\sm\{\ti\imath_a\}}$}.
\label{ai10eq30}
\ea

\noindent{\bf(b)} We have
\e
\begin{split}
{\rm E}^{\ti I}&=\n_2\bigl((\fF^{\ti\imath_0\ti\imath_1})^*(\m^{\ti I\sm\{\ti\imath_0\}}),\m^{\{\ti\imath_0,\ti\imath_1,\ti\imath_{\ti M}\}}\bigr)\\
&+\text{terms not involving\/ $\m^{\ti I\sm\{\ti\imath_0\}}$}.
\end{split}
\label{ai10eq31}
\e

\noindent{\bf(c)} In $F^0\Hom_{A^\rfi_\iy}^\un(\fF^{\ti I},\fF^{\ti\imath_0\ti\imath_{\ti M}})^{4-\ti M},$ we have
\ea
&\n_1({\rm H}^{\ti I})=-\n_2\bigl((\fF^{\ti\imath_0\ti\imath_1})^*({\rm E}^{\ti I\sm\{\ti\imath_0\}}),\m^{\{\ti\imath_0,\ti\imath_1,\ti\imath_{\ti M}\}}\bigr)
\label{ai10eq32}\\
&+\sum_{a=1}^{\ti M-1}(-1)^{\ti M-a}\n_2\bigl((\fF^{\{\ti\imath_0,\ldots,\ti\imath_{a-1}\}})^*\!\ci\!(\fF^{\{\ti\imath_{a+1},\ldots,\ti\imath_{\ti M}\}})_*
(\m^{\{\ti\imath_{a-1},\ti\imath_a,\ti\imath_{a+1}\}}),{\rm E}^{\ti I\sm\{\ti\imath_a\}}\bigr).
\nonumber
\ea
Passing to cohomology,  in\/ $H^{4-\ti M}\bigl(F^0\Hom_{A^\rfi_\iy}^\un(\fF^{\ti I},\fF^{\ti\imath_0\ti\imath_{\ti M}}),\n_1\bigr)$ we have
\ea
&-(\fF^{\ti\imath_0\ti\imath_1})^*([{\rm E}^{\ti I\sm\{\ti\imath_0\}}])\ci[\m^{\{\ti\imath_0,\ti\imath_1,\ti\imath_{\ti M}\}}]+{}
\label{ai10eq33}\\
&\sum_{a=1}^{\ti M-1}(-1)^{\ti M-a}\bigl[{\rm E}^{\ti I\sm\{\ti\imath_a\}}\bigr]\!\ci\!\bigl[(\fF^{\{\ti\imath_0,\ldots,\ti\imath_{a-1}\}})^*\ci
(\fF^{\{\ti\imath_{a+1},\ldots,\ti\imath_{\ti M}\}})_*
(\m^{\{\ti\imath_{a-1},\ti\imath_a,\ti\imath_{a+1}\}})\bigr]\!=\!0.
\nonumber
\ea
\end{prop}

\begin{proof}
For (a), by Definition \ref{ai10def1}, the terms in $\D_{n,k;y_0,\ldots,y_k}^{\ti I;x_0,\ldots,x_n}$ involving $\m^{\ti I\sm\{\ti\imath_a\}}$ are exactly the terms in $\d_{n,k;y_0,\ldots,y_k}^{\ti I;x_0,\ldots,x_n}$ involving $\m^{\ti I\sm\{\ti\imath_a\};z_0,\ldots,z_l}_l$. From \eq{ai7eq9} we can show that these must arise from collections of graphs $\bs\Ga_1,\ldots,\bs\Ga_k$ of the form
\begin{gather*}
\xymatrix@C=17pt@R=6pt{
\!\!\!\!{\begin{subarray}{l} \ts \ci \\[-3pt] \ts \kern 1.3pt \vdots \\[2pt] \ts \ci\end{subarray}}  
\ar[r] \ar@<2ex>[r] \ar@<-2ex>[r]
&{\begin{subarray}{l} 
\ts \overset{\{\ti\imath_0,\ti\imath_1\}}{\bu}  \\
\ts \overset{\{\ti\imath_0,\ti\imath_1\}}{\bu}  \\
\ts \overset{\{\ti\imath_0,\ti\imath_1\}}{\bu}  
\end{subarray}}
\ar[r] \ar@<2ex>[r]\ar@<-2ex>[r] & 
{\begin{subarray}{l} 
\ts \overset{\{\ti\imath_1,\ti\imath_2\}}{\bu}  \\
\ts \overset{\{\ti\imath_1,\ti\imath_2\}}{\bu}  \\
\ts \overset{\{\ti\imath_1,\ti\imath_2\}}{\bu}  
\end{subarray}}
\ar[r] \ar@<2ex>[r]\ar@<-2ex>[r] &  \cdots \ar[r] \ar@<2ex>[r]\ar@<-2ex>[r] & 
\overset{\{\ti\imath_{\ti M-1},\ti\imath_{\ti M}\}}{\bu}  
\ar[r] & \tr & \!\!\!\!\!\!{\begin{subarray}{l} \ts\bs\Ga_i,\\ \ts 1\!\le\! i\!<\!l,\end{subarray}} }
\\
\xymatrix@C=7pt@R=6pt{
\!\!\!\!{\begin{subarray}{l} \ts \ci \\[-3pt] \ts \kern 1.3pt \vdots \\[2pt] \ts \ci\end{subarray}}  
\ar[r] \ar@<2ex>[r] \ar@<-2ex>[r]
&{\begin{subarray}{l} 
\ts \overset{\{\ti\imath_0,\ti\imath_1\}}{\bu}  \\
\ts \overset{\{\ti\imath_0,\ti\imath_1\}}{\bu}  \\
\ts \overset{\{\ti\imath_0,\ti\imath_1\}}{\bu}  
\end{subarray}}
\ar[r] \ar@<2ex>[r]\ar@<-2ex>[r] &  \cdots \ar[r] \ar@<2ex>[r]\ar@<-2ex>[r] & 
{\begin{subarray}{l} 
\ts \overset{\{\ti\imath_{a-2},\ti\imath_{a-1}\}}{\bu}  \\
\ts \overset{\{\ti\imath_{a-2},\ti\imath_{a-1}\}}{\bu}  \\
\ts \overset{\{\ti\imath_{a-2},\ti\imath_{a-1}\}}{\bu}  
\end{subarray}}
\ar[r] \ar@<2ex>[r]\ar@<-2ex>[r] & {\begin{subarray}{l} 
\ts \overset{\{\ti\imath_{a-1},\ti\imath_a,\ti\imath_{a+1}\}}{\bu}  \\
\ts \overset{\{\ti\imath_{a-1},\ti\imath_a,\ti\imath_{a+1}\}}{\bu}  \\
\ts \overset{\{\ti\imath_{a-1},\ti\imath_a,\ti\imath_{a+1}\}}{\bu}  
\end{subarray}}
\ar[r] \ar@<2ex>[r]\ar@<-2ex>[r] & {\begin{subarray}{l} 
\ts \overset{\{\ti\imath_{a+1},\ti\imath_{a+2}\}}{\bu}  \\
\ts \overset{\{\ti\imath_{a+1},\ti\imath_{a+2}\}}{\bu}  \\
\ts \overset{\{\ti\imath_{a+1},\ti\imath_{a+2}\}}{\bu}  
\end{subarray}}
\ar[r] \ar@<2ex>[r]\ar@<-2ex>[r] &
 \cdots \ar[r] \ar@<2ex>[r]\ar@<-2ex>[r] & 
\overset{\{\ti\imath_{\ti M-1},\ti\imath_{\ti M}\}}{\bu}  
\ar[r] & \tr & \!\!\!\!\!\!\!\!\!\!\!\!\!\!\!\!\!\!\!\!\!{\begin{subarray}{l} \phantom{(} \\ \phantom{(} \\[20pt] \ts\bs\Ga_l,\end{subarray}} }\!\!\!\!\!\!\!
\\
\xymatrix@C=6.5pt@R=6pt{
\!\!\!\!{\begin{subarray}{l} \ts \ci \\[-3pt] \ts \kern 1.3pt \vdots \\[2pt] \ts \ci\end{subarray}}  
\ar[r] \ar@<2ex>[r] \ar@<-2ex>[r]
&{\begin{subarray}{l} 
\ts \overset{\{\ti\imath_0,\ti\imath_1\}}{\bu}  \\
\ts \overset{\{\ti\imath_0,\ti\imath_1\}}{\bu}  \\
\ts \overset{\{\ti\imath_0,\ti\imath_1\}}{\bu}  
\end{subarray}}
\ar[r] \ar@<2ex>[r]\ar@<-2ex>[r] &  \!\cdots\! \ar[r] \ar@<2ex>[r]\ar@<-2ex>[r] & 
{\begin{subarray}{l} 
\ts \overset{\{\ti\imath_{a-2},\ti\imath_{a-1}\}}{\bu}  \\
\ts \overset{\{\ti\imath_{a-2},\ti\imath_{a-1}\}}{\bu}  \\
\ts \overset{\{\ti\imath_{a-2},\ti\imath_{a-1}\}}{\bu}  
\end{subarray}}
\ar[r] \ar@<2ex>[r]\ar@<-2ex>[r] & {\begin{subarray}{l} 
\ts \overset{\{\ti\imath_{a-1},\ti\imath_{a+1}\}}{\bu}  \\
\ts \overset{\{\ti\imath_{a-1},\ti\imath_{a+1}\}}{\bu}  \\
\ts \overset{\{\ti\imath_{a-1},\ti\imath_{a+1}\}}{\bu}  
\end{subarray}}
\ar[r] \ar@<2ex>[r]\ar@<-2ex>[r] & {\begin{subarray}{l} 
\ts \overset{\{\ti\imath_{a+1},\ti\imath_{a+2}\}}{\bu}  \\
\ts \overset{\{\ti\imath_{a+1},\ti\imath_{a+2}\}}{\bu}  \\
\ts \overset{\{\ti\imath_{a+1},\ti\imath_{a+2}\}}{\bu}  
\end{subarray}}
\ar[r] \ar@<2ex>[r]\ar@<-2ex>[r] &
 \!\cdots\! \ar[r] \ar@<2ex>[r]\ar@<-2ex>[r] & 
\overset{\{\ti\imath_{\ti M-1},\ti\imath_{\ti M}\}}{\bu}  
\ar[r] & \tr & \!\!\!\!\!\!\!\!\!\!\!\!\!\!\!\!\!\!{\begin{subarray}{l} \phantom{(} \\ \phantom{(} \\[20pt] \ts\quad \bs\Ga_i,\\ \ts {}\!\! l\!<\! i\!<\!m,\end{subarray}} }\!\!\!\!\!\!\!\!\!\!
\\
\xymatrix@C=17pt@R=6pt{ {\begin{subarray}{l} \ts \ci \\[-3pt] \ts \kern 1.3pt \vdots \\[2pt] \ts \ci\end{subarray}}
\ar[r] \ar@<2ex>[r]\ar@<-2ex>[r]
& \overset{\ti I\sm\{\ti\imath_a\}}{\bu} \ar[r] & \tr & \!\!\!\!\!\!\!\!\!\!\!\!\bs\Ga_m, \\
{\begin{subarray}{l} \ts \ci \\[-3pt] \ts \kern 1.3pt \vdots \\[2pt] \ts \ci\end{subarray}}
\ar[r] \ar@<2ex>[r]\ar@<-2ex>[r]
& \overset{\{\ti\imath_0,\ti\imath_{\ti M}\}}{\bu} \ar[r] & \tr &
\bs\Ga_i,\; m<i\le k,  }
\end{gather*}
for $1\le l<m\le k$. Here the graph $\bs\Ga_l$ becomes the term $(\fF^{\{\ti\imath_0,\ldots,\ti\imath_{a-1}\}})^*\ci(\fF^{\{\ti\imath_{a+1},\ldots,\ti\imath_{\ti M}\}})_*(\m^{\{\ti\imath_{a-1},\ti\imath_a,\ti\imath_{a+1}\}})$ in \eq{ai10eq30}, and the graph $\bs\Ga_m$ becomes the term $\m^{\ti I\sm\{\ti\imath_a\}}$ in \eq{ai10eq30}. In the definition \eq{ai4eq17} of $\n_2$ in \eq{ai10eq30}, the graphs $\bs\Ga_i$ for $1\le i<l$ become the $\fF^0$ terms for $\fF^0=\fF^{\ti I}$, and the graphs $\bs\Ga_i$ for $l<i<m$ become the $\fF^1$ terms for $\fF^1=\fF^{\ti I\sm\{\ti\imath_a\}}$, and the graphs $\bs\Ga_i$ for $m<i\le k$ become the $\fF^2$ terms for $\fF^2=\fF^{\ti\imath_0\ti\imath_{\ti M}}$. Also $\m_k^{\{\ti\imath_{\ti M}\};y_0,\ldots,y_k}$ in \eq{ai10eq12} is part of~$\n_2$.

In this way we see that the terms in ${\rm E}^{\ti I}$ involving $\m^{\ti I\sm\{\ti\imath_a\}}$ agree with the expression $\n_2(\cdots)$ in \eq{ai10eq30}, up to sign. Comparing signs in the definitions, we find that the only difference is a factor $(-1)^{N(v,w)}=(-1)^{\ti M-a-1}$ in \eq{ai7eq9}, where $v$ is the vertex labelled $\{\ti\imath_{a-1},\ti\imath_a,\ti\imath_{a+1}\}$ in $\bs\Ga_l$, and $w$ is the vertex labelled $\ti I\sm\{\ti\imath_a\}$ in $\bs\Ga_m$, so that $N(v,w)=\ti M-a-1$. This proves~\eq{ai10eq30}.

For (b), the terms in $\D_{n,k;y_0,\ldots,y_k}^{\ti I;x_0,\ldots,x_n}$ involving $\m^{\ti I\sm\{\ti\imath_0\}}$ are exactly the terms in $\d_{n,k;y_0,\ldots,y_k}^{\ti I;x_0,\ldots,x_n}$ involving $\m^{\ti I\sm\{\ti\imath_0\};z_0,\ldots,z_l}_l$. From \eq{ai7eq9} we can show that these must arise from collections of graphs $\bs\Ga_1,\ldots,\bs\Ga_k$ of the form
\begin{gather*}
\xymatrix@C=17pt@R=6pt{
\!\!\!\!{\begin{subarray}{l} \ts \ci \\[-3pt] \ts \kern 1.3pt \vdots \\[2pt] \ts \ci\end{subarray}}  
\ar[r] \ar@<2ex>[r] \ar@<-2ex>[r]
&{\begin{subarray}{l} 
\ts \overset{\{\ti\imath_0,\ti\imath_1\}}{\bu}  \\
\ts \overset{\{\ti\imath_0,\ti\imath_1\}}{\bu}  \\
\ts \overset{\{\ti\imath_0,\ti\imath_1\}}{\bu}  
\end{subarray}}
\ar[r] \ar@<2ex>[r]\ar@<-2ex>[r] & 
{\begin{subarray}{l} 
\ts \overset{\{\ti\imath_1,\ti\imath_2\}}{\bu}  \\
\ts \overset{\{\ti\imath_1,\ti\imath_2\}}{\bu}  \\
\ts \overset{\{\ti\imath_1,\ti\imath_2\}}{\bu}  
\end{subarray}}
\ar[r] \ar@<2ex>[r]\ar@<-2ex>[r] &  \cdots \ar[r] \ar@<2ex>[r]\ar@<-2ex>[r] & 
\overset{\{\ti\imath_{\ti M-1},\ti\imath_{\ti M}\}}{\bu}  
\ar[r] & \tr & \!\!\!\!\!\!{\bs\Ga_i,\; 1\!\le\! i\!<\!l,} }
\\
\xymatrix@C=20pt@R=6pt{
\!\!\!\!{\begin{subarray}{l} \ts \ci \\[-3pt] \ts \kern 1.3pt \vdots \\[2pt] \ts \ci\end{subarray}}  
\ar[r] \ar@<2ex>[r] \ar@<-2ex>[r]
&{\begin{subarray}{l} 
\ts \overset{\{\ti\imath_0,\ti\imath_1\}}{\bu}  \\
\ts \overset{\{\ti\imath_0,\ti\imath_1\}}{\bu}  \\
\ts \overset{\{\ti\imath_0,\ti\imath_1\}}{\bu}  
\end{subarray}}
\ar[r] \ar@<2ex>[r]\ar@<-2ex>[r]
& \overset{\ti I\sm\{\ti\imath_0\}}{\bu} \ar[r] & \tr & \!\!\!\!\!\!\!\bs\Ga_l, }\!\!\!\!\!\!\!
\\
\xymatrix@C=20pt@R=6pt{
\!\!\!\!{\begin{subarray}{l} \ts \ci \\[-3pt] \ts \kern 1.3pt \vdots \\[2pt] \ts \ci\end{subarray}}  
\ar[r] \ar@<2ex>[r] \ar@<-2ex>[r]
&{\begin{subarray}{l} 
\ts \overset{\{\ti\imath_0,\ti\imath_1\}}{\bu}  \\
\ts \overset{\{\ti\imath_0,\ti\imath_1\}}{\bu}  \\
\ts \overset{\{\ti\imath_0,\ti\imath_1\}}{\bu}  
\end{subarray}}
\ar[r] \ar@<2ex>[r]\ar@<-2ex>[r]
& \overset{\{\ti\imath_1,\ti\imath_{\ti M}\}}{\bu} \ar[r] & \tr & \!\!\!\!\!\!\!\bs\Ga_i,\; l\!<\! i\!<\!m, }\!\!\!\!\!\!\!
\\
\xymatrix@C=17pt@R=6pt{ {\begin{subarray}{l} \ts \ci \\[-3pt] \ts \kern 1.3pt \vdots \\[2pt] \ts \ci\end{subarray}}
\ar[r] \ar@<2ex>[r]\ar@<-2ex>[r]
& \overset{\{\ti\imath_0,\ti\imath_1,\ti\imath_{\ti M}\}}{\bu} \ar[r] & \tr & \!\!\!\!\!\!\!\!\!\!\!\!\bs\Ga_m, \\
{\begin{subarray}{l} \ts \ci \\[-3pt] \ts \kern 1.3pt \vdots \\[2pt] \ts \ci\end{subarray}}
\ar[r] \ar@<2ex>[r]\ar@<-2ex>[r]
& \overset{\{\ti\imath_0,\ti\imath_{\ti M}\}}{\bu} \ar[r] & \tr &
\bs\Ga_i,\; m<i\le k,  }
\end{gather*}
for $1\le l<m\le k$. Here the graph $\bs\Ga_l$ becomes the term $(\fF^{\ti\imath_0\ti\imath_1})^*(\m^{\ti I\sm\{\ti\imath_0\}})$ in \eq{ai10eq31}, and the graph $\bs\Ga_m$ becomes the term $\m^{\{\ti\imath_0,\ti\imath_1,\ti\imath_{\ti M}\}}$ in \eq{ai10eq31}. In the definition \eq{ai4eq17} of $\n_2$ in \eq{ai10eq31}, the graphs $\bs\Ga_i$ for $1\le i<l$ become the $\fF^0$ terms for $\fF^0=\fF^{\ti I}$, and the graphs $\bs\Ga_i$ for $l<i<m$ become the $\fF^1$ terms for $\fF^1=\fF^{\{\ti\imath_0,\ti\imath_1,\ti\imath_{\ti M}\}}$, and the graphs $\bs\Ga_i$ for $m<i\le k$ become the $\fF^2$ terms for $\fF^2=\fF^{\ti\imath_0\ti\imath_{\ti M}}$. Also $\m_k^{\{\ti\imath_{\ti M}\};y_0,\ldots,y_k}$ in \eq{ai10eq12} is part of~$\n_2$.

In this way we see that the terms in ${\rm E}^{\ti I}$ involving $\m^{\ti I\sm\{\ti\imath_0\}}$ agree with the expression $\n_2(\cdots)$ in \eq{ai10eq31}, up to sign. Comparing signs in the definitions, we find they also agree. This proves~\eq{ai10eq31}.

For (c), we first claim that
\ea
&{\rm G}_{n,k;y_0,\ldots,y_k}^{\ti I;x_0,\ldots,x_n}=
\label{ai10eq34}\\
&\sum_{\begin{subarray}{l} 0=l_0\le l_1\le \cdots\le l_k=n, \\
\text{isomorphism classes} \\ \text{of $\bs\Ga_1,\ldots,\bs\Ga_k$:}\\ \text{$\bs\Ga_i$ is a balanced} \\ 
\text{$I$-labelled tree} \\ \text{with $l_i-l_{i-1}$} \\ \text{leaves, $1\le i\le k$,} \\
\text{$L(\bs\Ga_1)=\ti I$,} \\ \text{$R(\bs\Ga_k)=\{\ti\imath_0,\ti\imath_{\ti M}\}$,} \\
\text{$R(\bs\Ga_i)=L(\bs\Ga_{i+1})$,} \\ \text{$i=1,\ldots,k-1$,} \\
y_i=f^{R(\bs\Ga_i)}(x_{l_i}), \\
i=1,\ldots,k,\; y_0=f^{\ti I}(x_0). \\
\color{blue}
\text{For all $i=1,\ldots,k$,} \\
\color{blue}
\text{$(L(\bs\Ga_i),R(\bs\Ga_i))$ is not} \\
\color{blue}
\text{$(\ti I,\{\ti\imath_0,\ti\imath_{\ti M}\})$ or} \\
\color{blue}
\text{$(\ti I,\{\ti\imath_0,\ti\imath_1,\ti\imath_{\ti M}\})$ or} \\
\color{blue}
\text{$(\ti I\sm\{\ti\imath_a\},\{\ti\imath_0,\ti\imath_{\ti M}\})$,}\\
\color{blue}
\text{$1<a<\ti M$.}
\color{black}
 \end{subarray}}
\begin{aligned}[t]
&(-1)^{(k+1)(\ti M+n+1)}\cdot{}\\
&\biggl(\prod_{\begin{subarray}{l} \text{$v$ internal vertex or bud of\/ $\bs\Ga_i,$ $1\le i\le k,$} \\ \text{$w$ internal vertex or bud of\/ $\bs\Ga_{i'},$ $1\le i'\le k,$} \\ 
\text{either $i<i',$ or $i=i'$ and} \\ \text{$v\pr w$ in order $\pr$ of $\Ga_i$,} \\ \text{$v$ has $p$ incoming edges including} \\ \text{$q$ leaves, $w$ has $r$ incoming edges}\end{subarray}\!\!\!\!\!\!\!\!\!\!\!\!\!\!\!\!\!\!\!\!\!\!\!\!\!\!\!\!\!\!\!\!\!\!\!\!\!\!\!} (-1)^{\begin{subarray}{l} N(v,w)+q(r+1) \\ +(p+1)(\md{I_w}+r+1)\end{subarray}}\biggr)\cdot{}\\
& C(\bs\Ga_1)^{x_{l_0},\ldots,x_{l_1}}\!\ot\! C(\bs\Ga_2)^{x_{l_1},\ldots,x_{l_2}}\!\ot\!\cdots\!\ot\! C(\bs\Ga_k)^{x_{l_{k-1}},\ldots,x_{l_k}}.
\end{aligned}
\nonumber
\ea
This is the same as \eq{ai7eq9} for $\d_{n,k;y_0,\ldots,y_k}^{\ti I;x_0,\ldots,x_n}$, except that we have added the final conditions in the sum that $(L(\bs\Ga_i),R(\bs\Ga_i))\ne(\ti I,\{\ti\imath_0,\ti\imath_{\ti M}\})$, $(\ti I,\{\ti\imath_0,\ti\imath_1,\ti\imath_{\ti M}\})$, $(\ti I\sm\{\ti\imath_a\},\{\ti\imath_0,\ti\imath_{\ti M}\})$. Equation \eq{ai10eq34} holds since a term in \eq{ai7eq9} for $\d_{n,k;y_0,\ldots,y_k}^{\ti I;x_0,\ldots,x_n}$ contains an $\m_l^{\ti I;w_0,\ldots,w_l}$ if and only if $(L(\bs\Ga_i),R(\bs\Ga_i))=(\ti I,\{\ti\imath_0,\ti\imath_{\ti M}\})$ for some $i$, and contains an $\m_l^{\ti I\sm\{\ti\imath_0\};w_0,\ldots,w_l}$ if and only if $(L(\bs\Ga_i),R(\bs\Ga_i))=(\ti I,\{\ti\imath_0,\ti\imath_1,\ti\imath_{\ti M}\})$ for some $i$, and contains an $\m_l^{\ti I\sm\{\ti\imath_a\};w_0,\ldots,w_l}$ for $0<a<M$ if and only if $(L(\bs\Ga_i),R(\bs\Ga_i))=(\ti I\sm\{\ti\imath_a\},\{\ti\imath_0,\ti\imath_{\ti M}\})$ for some~$i$.

To prove \eq{ai10eq32}, for $n\ge 0$ and $x_0,\ldots,x_n\in\sO^{\ti\imath_0}$ consider the equation
\ea
&\bigl(\n_1({\rm H}^{\ti I})\bigr){}_n^{x_0,\ldots,x_n}
\label{ai10eq35}\\
&=\sum_{\begin{subarray}{l} 
0=a_0\le a_1\le \\
\cdots\le a_y=n,\\
1\le x\le y \\
\color{red}(b1)\color{black} \end{subarray}}
\begin{aligned}[t]&
(-1)^{(y+1)\ti M}\cdot(-1)^{\sum_{i=1}^{x-1}(a_i-a_{i-1}+1)\ti M}\cdot{}\\
&(-1)^{\sum_{i=1}^y(a_i-a_{i-1}+1)(y-i)} \cdot
\m^{\{\ti\imath_{\ti M}\};\ldots}_y\ci{} \\ 
&\bigl(\fF_{a_1-a_0}^{\ti I;x_{a_0},\ldots,x_{a_1}}\ot\fF_{a_2-a_1}^{\ti I;x_{a_1},\ldots,x_{a_2}}\ot\cdots\ot\fF_{\ti I;a_{x-1}-a_{x-2}}^{x_{a_{x-2}},\ldots,x_{a_{x-1}}} \\
&{}\ot\Biggl\{
\sum_{\begin{subarray}{l} a_{x-1}=l_0\le l_1\le \cdots\le l_t=a_x, \\
 \text{$\bs\Ga_1,\ldots,\bs\Ga_t$ with} \\
\text{usual conditions:} \\
\text{$(L(\bs\Ga_i),R(\bs\Ga_i))$ is} \\
\text{not $(\ti I,\{\ti\imath_0,\ti\imath_{\ti M}\})$ or} \\
\text{$(\ti I,\{\ti\imath_0,\ti\imath_1,\ti\imath_{\ti M}\})$ or} \\
\text{$(\ti I\sm\{\ti\imath_a\},\{\ti\imath_0,\ti\imath_{\ti M}\})$,}\\
\text{all $i$, $1<a<\ti M$.}
 \end{subarray}}
\begin{aligned}[t]
&(-1)^{(t+1)(\ti M+a_x-a_{x-1}+1)}\cdot{}\\
&\biggl(\prod_{\begin{subarray}{l} \text{$v,w,p,q,r$ usual conditions}
\end{subarray}\!\!\!\!\!\!\!\!\!\!\!\!\!\!\!\!\!\!\!\!\!\!\!\!\!\!\!\!\!\!\!\!\!\!\!\!\!\!\!\!\!\!\!\!\!\!\!\!\!\!\!\!\!\!\!} (-1)^{\begin{subarray}{l} N(v,w)+q(r+1) \\ +(p+1)(\md{I_w}+r+1)\end{subarray}}\biggr)\cdot{}\\
&\m_t^{\{\ti\imath_{\ti M}\},\ldots}\ci \bigl(C(\bs\Ga_1)^{x_{l_0},\ldots,x_{l_1}} \\
& {}\ot\!\cdots\!\ot\! C(\bs\Ga_t)^{x_{l_{t-1}},\ldots,x_{l_t}}\bigr)
\end{aligned}\Biggr\}\\
&{}\ot\fF_{a_{x+1}-a_x}^{\ti\imath_0\ti\imath_{\ti M};x_{a_x},\ldots,x_{a_{x+1}}}\ot\cdots\ot\fF_{a_y-a_{y-1}}^{\ti\imath_0\ti\imath_{\ti M};x_{a_{y-1}},\ldots,x_{a_y}}\bigr)
\end{aligned}
\allowdisplaybreaks
\nonumber\\
&+\sum_{\begin{subarray}{l} x,z\ge 0, \\ y\ge 1:\\ x+y+z=n \\
\color{red}(b2)\color{black} \end{subarray}}
\begin{aligned}[t]&(-1)^{\ti M+1}\cdot(-1)^{x+yz}\cdot{}\\
&\Biggl\{
\sum_{\begin{subarray}{l} 1\le j\le k, \\
0=l_0\le l_1\le \cdots\le l_{j-1}\le x, \\
x+y\le l_j\le\cdots\le l_k=n \\
\text{$\bs\Ga_1,\ldots,\bs\Ga_k$ with} \\
\text{usual conditions:} \\
\text{$(L(\bs\Ga_i),R(\bs\Ga_i))$ is} \\
\text{not $(\ti I,\{\ti\imath_0,\ti\imath_{\ti M}\})$ or} \\
\text{$(\ti I,\{\ti\imath_0,\ti\imath_1,\ti\imath_{\ti M}\})$ or} \\
\text{$(\ti I\sm\{\ti\imath_a\},\{\ti\imath_0,\ti\imath_{\ti M}\})$,}\\
\text{all $i$, $1<a<\ti M$.}
 \end{subarray}}
\begin{aligned}[t]
&(-1)^{(k+1)(\ti M+x+z+2)}\cdot{}\\
&\biggl(\prod_{\begin{subarray}{l} \text{$v,w,p,q,r$ usual conditions}
\end{subarray}\!\!\!\!\!\!\!\!\!\!\!\!\!\!\!\!\!\!\!\!\!\!\!\!\!\!\!\!\!\!\!\!\!\!\!\!\!\!\!\!\!\!\!\!\!\!\!\!\!\!\!\!\!\!\!} (-1)^{\begin{subarray}{l} N(v,w)+q(r+1) \\ +(p+1)(\md{I_w}+r+1)\end{subarray}}\biggr)\cdot{}\\
&\m_k^{\{\ti\imath_{\ti M}\},\ldots}\ci \bigl(C(\bs\Ga_1)^{x_{l_0},\ldots,x_{l_1}}\ot{} \\
& \cdots\!\ot\! C(\bs\Ga_j)^{x_{l_{j-1}},\ldots,x_x,x_{x+y},\ldots,x_{l_j}} \\
& {}\ot\!\cdots\!\ot\! C(\bs\Ga_k)^{x_{l_{k-1}},\ldots,x_{l_k}}\bigr)
\end{aligned}\Biggr\}
\\
&\ci\bigl(\id^{\ot^x}\ot\m_y^{\{\ti\imath_0\};x_x,\ldots,x_{x+y}}\ot \id^{\ot^z}\bigr) 
\end{aligned}
\allowdisplaybreaks
\nonumber\\
&=\sum_{\begin{subarray}{l} k\ge 1,\\ s,u\ge 0, \\ t\ge 1: \\ s+t+u=k \\
\color{red}(c1)\color{black} \end{subarray}}
\begin{aligned}[t]&
(-1)^{s+tu}\cdot\m_{s+u+1}^{\{\ti\imath_{\ti M}\};\ldots}\ci
\bigl(\id^{\ot^s}\ot\m_t^{\{\ti\imath_{\ti M}\};\ldots}\ot \id^{\ot^u}\bigr)\\
&{}\ci\Biggl\{
\sum_{\begin{subarray}{l} 0=l_0\le l_1\le \cdots\le l_k=n, \\
 \text{$\bs\Ga_1,\ldots,\bs\Ga_k$ with} \\
\text{usual conditions:} \\
\text{$(L(\bs\Ga_i),R(\bs\Ga_i))$ is} \\
\text{not $(\ti I,\{\ti\imath_0,\ti\imath_{\ti M}\})$ or} \\
\text{$(\ti I,\{\ti\imath_0,\ti\imath_1,\ti\imath_{\ti M}\})$ or} \\
\text{$(\ti I\sm\{\ti\imath_a\},\{\ti\imath_0,\ti\imath_{\ti M}\})$,}\\
\text{all $i$, $1<a<\ti M$.} \\
\color{blue} L(\bs\Ga_{s+1})=\ti I, \\
\color{blue} R(\bs\Ga_{s+t})=\{\ti\imath_0,\ti\imath_{\ti M}\}.
\color{black} \end{subarray}}
\begin{aligned}[t]
&(-1)^{(k+1)(\ti M+n+1)}\cdot{}\\
&\biggl(\prod_{\begin{subarray}{l} \text{$v,w,p,q,r$ usual conditions}
\end{subarray}\!\!\!\!\!\!\!\!\!\!\!\!\!\!\!\!\!\!\!\!\!\!\!\!\!\!\!\!\!\!\!\!\!\!\!\!\!\!\!\!\!\!\!\!\!\!\!\!\!\!\!\!\!\!\!} (-1)^{\begin{subarray}{l} N(v,w)+q(r+1) \\ +(p+1)(\md{I_w}+r+1)\end{subarray}}\biggr)\cdot{}\\
&C(\bs\Ga_1)^{x_{l_0},\ldots,x_{l_1}}\!\ot\!\cdots\!\ot\! C(\bs\Ga_k)^{x_{l_{k-1}},\ldots,x_{l_k}}
\end{aligned}\Biggr\}
\end{aligned}
\allowdisplaybreaks
\nonumber\\
&+\sum_{\begin{subarray}{l} m\ge 1,\\ s,u\ge 0, \\ t\ge 1: \\ s+t+u=m \\
\color{red}(c2)\color{black} \end{subarray}}
\begin{aligned}[t]&
(-1)^{s+tu}\cdot\m_{s+u+1}^{\{\ti\imath_{\ti M}\};\ldots}\ci
\bigl(\id^{\ot^s}\ot\m_t^{\{\ti\imath_{\ti M}\};\ldots}\ot \id^{\ot^u}\bigr)\\
&{}\ci\Biggl\{
\sum_{\begin{subarray}{l} 0=l_0\le l_1\le \cdots\le l_m=n, \\
 \text{$\bs\Ga_1,\ldots,\bs\Ga_m$ with} \\
\text{usual conditions:} \\
\text{$(L(\bs\Ga_i),R(\bs\Ga_i))$ is} \\
\text{not $(\ti I,\{\ti\imath_0,\ti\imath_{\ti M}\})$ or} \\
\text{$(\ti I,\{\ti\imath_0,\ti\imath_1,\ti\imath_{\ti M}\})$ or} \\
\text{$(\ti I\sm\{\ti\imath_a\},\{\ti\imath_0,\ti\imath_{\ti M}\})$,}\\
\text{all $i$, $1<a<\ti M$.} \\
\color{blue} \text{$(L(\bs\Ga_{s+1}), R(\bs\Ga_{s+t}))$} \\
\color{blue} \text{is not $(\ti I,\{\ti\imath_0,\ti\imath_{\ti M}\})$ or} \\
\color{blue} \text{$(\ti I,\{\ti\imath_0,\ti\imath_1,\ti\imath_{\ti M}\})$ or} \\
\color{blue} \text{$(\ti I\sm\{\ti\imath_a\},\{\ti\imath_0,\ti\imath_{\ti M}\})$,}\\
\color{blue} \text{$1<a<\ti M$.} 
\color{black} \end{subarray}}
\begin{aligned}[t]
&(-1)^{(m+1)(\ti M+n+1)}\cdot{}\\
&\biggl(\prod_{\begin{subarray}{l} \text{$v,w,p,q,r$ usual conditions}
\end{subarray}\!\!\!\!\!\!\!\!\!\!\!\!\!\!\!\!\!\!\!\!\!\!\!\!\!\!\!\!\!\!\!\!\!\!\!\!\!\!\!\!\!\!\!\!\!\!\!\!\!\!\!\!\!\!\!} (-1)^{\begin{subarray}{l} N(v,w)+q(r+1) \\ +(p+1)(\md{I_w}+r+1)\end{subarray}}\biggr)\cdot{}\\
&C(\bs\Ga_1)^{x_{l_0},\ldots,x_{l_1}}\!\ot\!\cdots\!\ot\! C(\bs\Ga_m)^{x_{l_{m-1}},\ldots,x_{l_m}}
\end{aligned}\Biggr\}
\end{aligned}
\allowdisplaybreaks
\nonumber\\
&=-\sum_{\begin{subarray}{l} k\ge 1,\\ s,u\ge 0, \\ t\ge 1: \\ s+t+u=k \\
\color{red}(d1)\color{black} \end{subarray}}
\begin{aligned}[t]&
(-1)^{s+tu}\cdot\m_{s+u+1}^{\{\ti\imath_{\ti M}\};\ldots}\ci
\bigl(\id^{\ot^s}\ot\m_t^{\{\ti\imath_{\ti M}\};\ldots}\ot \id^{\ot^u}\bigr)\\
&{}\ci\Biggl\{
\sum_{\begin{subarray}{l} 0=l_0\le l_1\le \cdots\le l_k=n, \\
 \text{$\bs\Ga_1,\ldots,\bs\Ga_k$ with} \\
\text{usual conditions:} \\
\text{$(L(\bs\Ga_i),R(\bs\Ga_i))$ is} \\
\text{not $(\ti I,\{\ti\imath_0,\ti\imath_{\ti M}\})$ or} \\
\text{$(\ti I,\{\ti\imath_0,\ti\imath_1,\ti\imath_{\ti M}\})$ or} \\
\text{$(\ti I\sm\{\ti\imath_a\},\{\ti\imath_0,\ti\imath_{\ti M}\})$,}\\
\text{all $i$, $1<a<\ti M$.} \\
\color{blue} (L(\bs\Ga_{s+1}),R(\bs\Ga_{s+t})) \\
\color{blue} \ne (\ti I,\{\ti\imath_0,\ti\imath_{\ti M}\}).
\color{black} \end{subarray}}
\begin{aligned}[t]
&(-1)^{(k+1)(\ti M+n+1)}\cdot{}\\
&\biggl(\prod_{\begin{subarray}{l} \text{$v,w,p,q,r$ usual conditions}
\end{subarray}\!\!\!\!\!\!\!\!\!\!\!\!\!\!\!\!\!\!\!\!\!\!\!\!\!\!\!\!\!\!\!\!\!\!\!\!\!\!\!\!\!\!\!\!\!\!\!\!\!\!\!\!\!\!\!} (-1)^{\begin{subarray}{l} N(v,w)+q(r+1) \\ +(p+1)(\md{I_w}+r+1)\end{subarray}}\biggr)\cdot{}\\
&C(\bs\Ga_1)^{x_{l_0},\ldots,x_{l_1}}\!\ot\!\cdots\!\ot\! C(\bs\Ga_k)^{x_{l_{k-1}},\ldots,x_{l_k}}
\end{aligned}\Biggr\}
\end{aligned}
\allowdisplaybreaks
\nonumber\\
&+\sum_{\begin{subarray}{l} k\ge 1,\\ s,u\ge 0, \\ t\ge 1: \\ s+t+u=k \\
\color{red}(d2)\color{black} \end{subarray}}
\begin{aligned}[t]&
(-1)^{s+tu}\cdot\m_{s+u+1}^{\{\ti\imath_{\ti M}\};\ldots}\ci
\bigl(\id^{\ot^s}\ot\m_t^{\{\ti\imath_{\ti M}\};\ldots}\ot \id^{\ot^u}\bigr)\\
&{}\ci\Biggl\{
\sum_{\begin{subarray}{l} 0=l_0\le l_1\le \cdots\le l_k=n, \\
 \text{$\bs\Ga_1,\ldots,\bs\Ga_k$ with} \\
\text{usual conditions:} \\
\text{$(L(\bs\Ga_i),R(\bs\Ga_i))$ is} \\
\text{not $(\ti I,\{\ti\imath_0,\ti\imath_{\ti M}\})$ or} \\
\text{$(\ti I,\{\ti\imath_0,\ti\imath_1,\ti\imath_{\ti M}\})$ or} \\
\text{$(\ti I\sm\{\ti\imath_a\},\{\ti\imath_0,\ti\imath_{\ti M}\})$,}\\
\text{all $i$, $1<a<\ti M$.} \\
\color{blue} \text{$(L(\bs\Ga_{s+1}), R(\bs\Ga_{s+t}))$} \\
\color{blue} \text{$\ne (\ti I,\{\ti\imath_0,\ti\imath_{\ti M}\})$ or} \\
\color{blue} \text{$(\ti I,\{\ti\imath_0,\ti\imath_1,\ti\imath_{\ti M}\})$ or} \\
\color{blue} \text{$(\ti I\sm\{\ti\imath_a\},\{\ti\imath_0,\ti\imath_{\ti M}\})$,}\\
\color{blue} \text{$1<a<\ti M$.} 
\color{black} \end{subarray}}
\begin{aligned}[t]
&(-1)^{(k+1)(\ti M+n+1)}\cdot{}\\
&\biggl(\prod_{\begin{subarray}{l} \text{$v,w,p,q,r$ usual conditions}
\end{subarray}\!\!\!\!\!\!\!\!\!\!\!\!\!\!\!\!\!\!\!\!\!\!\!\!\!\!\!\!\!\!\!\!\!\!\!\!\!\!\!\!\!\!\!\!\!\!\!\!\!\!\!\!\!\!\!} (-1)^{\begin{subarray}{l} N(v,w)+q(r+1) \\ +(p+1)(\md{I_w}+r+1)\end{subarray}}\biggr)\cdot{}\\
&C(\bs\Ga_1)^{x_{l_0},\ldots,x_{l_1}}\!\ot\!\cdots\!\ot\! C(\bs\Ga_k)^{x_{l_{k-1}},\ldots,x_{l_k}}
\end{aligned}\Biggr\}
\end{aligned}
\allowdisplaybreaks
\nonumber\\
&=-\sum_{\begin{subarray}{l} k\ge 1,\\ s,u\ge 0, \\ t\ge 1: \\ s+t+u=k \end{subarray}}
\begin{aligned}[t]&
(-1)^{s+tu}\cdot\m_{s+u+1}^{\{\ti\imath_{\ti M}\};\ldots}\ci
\bigl(\id^{\ot^s}\ot\m_t^{\{\ti\imath_{\ti M}\};\ldots}\ot \id^{\ot^u}\bigr)\\
&{}\ci\Biggl\{
\sum_{\begin{subarray}{l} 0=l_0\le l_1\le \cdots\le l_k=n, \\
 \text{$\bs\Ga_1,\ldots,\bs\Ga_k$ with} \\
\text{usual conditions:} \\
\text{$(L(\bs\Ga_i),R(\bs\Ga_i))$ is} \\
\text{not $(\ti I,\{\ti\imath_0,\ti\imath_{\ti M}\})$ or} \\
\text{$(\ti I,\{\ti\imath_0,\ti\imath_1,\ti\imath_{\ti M}\})$ or} \\
\text{$(\ti I\sm\{\ti\imath_a\},\{\ti\imath_0,\ti\imath_{\ti M}\})$,}\\
\text{all $i$, $1<a<\ti M$.} \\
\color{blue} \text{$(L(\bs\Ga_{s+1}), R(\bs\Ga_{s+t}))$} \\
\color{blue} \text{$=(\ti I,\{\ti\imath_0,\ti\imath_1,\ti\imath_{\ti M}\})$} \; \color{red}(e1) \\
\color{blue} \text{or $(\ti I\sm\{\ti\imath_a\},\{\ti\imath_0,\ti\imath_{\ti M}\})$,}\\
\color{blue} \text{some $1<a<\ti M$.\; \color{red}(e2)} 
\color{black} \end{subarray}}
\begin{aligned}[t]
&(-1)^{(k+1)(\ti M+n+1)}\cdot{}\\
&\biggl(\prod_{\begin{subarray}{l} \text{$v,w,p,q,r$ usual conditions}
\end{subarray}\!\!\!\!\!\!\!\!\!\!\!\!\!\!\!\!\!\!\!\!\!\!\!\!\!\!\!\!\!\!\!\!\!\!\!\!\!\!\!\!\!\!\!\!\!\!\!\!\!\!\!\!\!\!\!} (-1)^{\begin{subarray}{l} N(v,w)+q(r+1) \\ +(p+1)(\md{I_w}+r+1)\end{subarray}}\biggr)\cdot{}\\
&C(\bs\Ga_1)^{x_{l_0},\ldots,x_{l_1}}\!\ot\!\cdots\!\ot\! C(\bs\Ga_k)^{x_{l_{k-1}},\ldots,x_{l_k}}
\end{aligned}\Biggr\}
\end{aligned}
\allowdisplaybreaks
\nonumber\\
&=-\n_2\bigl((\fF^{\ti\imath_0\ti\imath_1})^*({\rm E}^{\ti I\sm\{\ti\imath_0\}}),\m^{\{\ti\imath_0,\ti\imath_1,\ti\imath_{\ti M}\}}\bigr)_{\color{red}(f1)}
\nonumber\\
&+\sum_{\begin{subarray}{l} a=1 \\ \color{red}(f2)\end{subarray}}^{\ti M-1}\begin{aligned}&(-1)^{\ti M-a}\\[-3pt]
&\n_2\bigl((\fF^{\{\ti\imath_0,\ldots,\ti\imath_{a-1}\}})^*\ci(\fF^{\{\ti\imath_{a+1},\ldots,\ti\imath_{\ti M}\}})_*
(\m^{\{\ti\imath_{a-1},\ti\imath_a,\ti\imath_{a+1}\}}),{\rm E}^{\ti I\sm\{\ti\imath_a\}}\bigr)_n^{x_0,\ldots,x_n}.
\end{aligned}
\nonumber
\ea

Here throughout we have omitted some conditions in the sums to save space, and focus on what is important. When we say `$\bs\Ga_1,\ldots,\bs\Ga_k$ with usual conditions' or `$v,w,p,q,r$ usual conditions', we mean the conditions on these in \eq{ai7eq9}, adapted as appropriate. Also we have omitted some objects, and conditions on objects: when we write terms like $\m^{\{\ti\imath_{\ti M}\};\ldots}_t$, we mean $\m^{\{\ti\imath_{\ti M}\};y_0,\ldots,y_t}_t$ for some $y_0,\ldots,y_t\in\sO^{\ti\imath_{\ti M}}$, and these $y_i$ should be determined by conditions like those in \eq{ai10eq34}, so that the compositions make sense. To make the equation easier to navigate, we have labelled terms in the second step $(b1),(b2)$, and so on. Extra new conditions in the sums are highlighted in blue.

In the first step of \eq{ai10eq35} we use \eq{ai4eq15}, \eq{ai10eq28} and \eq{ai10eq34}. The terms $\{\cdots\}$ in the first step are~$({\rm H}^{\ti I})_{\cdots}^{\cdots}$.

For the second step, we first show that $(b1)=(c1)$. Now $(b1)$ contains $\m_y^{\{\ti\imath_{\ti M}\},\ldots},\m_t^{\{\ti\imath_{\ti M}\},\ldots}$. We bring these to the front as $\m_y^{\{\ti\imath_{\ti M}\},\ldots}\ci\bigl(\id^{\ot^{x-1}}\ot\m_t^{\{\ti\imath_{\ti M}\},\ldots}\ot\id^{\ot^{y-x}}\bigr)$, and change variables from $x,y$ to $s,u,k$ with $s=x-1$, $u=y-x$ and $k=t+y-1$, making this $\m_{s+u+1}^{\{\ti\imath_{\ti M}\};\ldots}\ci\bigl(\id^{\ot^s}\ot\m_t^{\{\ti\imath_{\ti M}\};\ldots}\ot \id^{\ot^u}\bigr)$ as in~$(c1)$.

Up to sign, the term $\fF_{a_1-a_0}^{\ti I;x_{a_0},\ldots,x_{a_1}}\ot\fF_{a_2-a_1}^{\ti I;x_{a_1},\ldots,x_{a_2}}\ot\cdots\ot\fF_{\ti I;a_{x-1}-a_{x-2}}^{x_{a_{x-2}},\ldots,x_{a_{x-1}}}\ot C(\bs\Ga_1)^{x_{l_0},\ldots,x_{l_1}}\ot\cdots\ot C(\bs\Ga_t)^{x_{l_{t-1}},\ldots,x_{l_t}}\ot \fF_{a_{x+1}-a_x}^{\ti\imath_0\ti\imath_{\ti M};x_{a_x},\ldots,x_{a_{x+1}}}\ot\cdots\ot\fF_{a_y-a_{y-1}}^{\ti\imath_0\ti\imath_{\ti M};x_{a_{y-1}},\ldots,x_{a_y}}$ in $(b1)$ becomes the term $C(\bs\Ga_1)^{x_{l_0},\ldots,x_{l_1}}\ot\cdots\ot C(\bs\Ga_k)^{x_{l_{k-1}},\ldots,x_{l_k}}$ in $(c1)$. This works because of the extra conditions $L(\bs\Ga_{s+1})=\ti I$, $R(\bs\Ga_{s+t})=\{\ti\imath_0,\ti\imath_{\ti M}\}$ in $(c1)$, which ensure that the sum over trees $\bs\Ga_1,\ldots,\bs\Ga_s$ is of the form $\fF_{a_1-a_0}^{\ti I;x_{a_0},\ldots,x_{a_1}}\ot\fF_{a_2-a_1}^{\ti I;x_{a_1},\ldots,x_{a_2}}\ot\cdots\ot\fF_{\ti I;a_{x-1}-a_{x-2}}^{x_{a_{x-2}},\ldots,x_{a_{x-1}}}$ with $x-1=s$, and the sum over trees $\bs\Ga_{s+t+1},\ldots,\bs\Ga_k$ is of the form $\fF_{a_{x+1}-a_x}^{\ti\imath_0\ti\imath_{\ti M};x_{a_x},\ldots,x_{a_{x+1}}}\ot\cdots\ot\fF_{a_y-a_{y-1}}^{\ti\imath_0\ti\imath_{\ti M};x_{a_{y-1}},\ldots,x_{a_y}}$ with $y-x-1=k-s-t$. Computing signs, we find that~$(b1)=(c1)$.

Next we show that $(b2)=(c2)$. To do this, we start with $(b2)$, and fix some of the data in the sum, while summing over the choices of the remaining data. The data we fix is $j,k,l_0,\ldots,l_k,\bs\Ga_1,\ldots,\bs\Ga_{j-1},\bs\Ga_{j+1},\ldots,\bs\Ga_k$. The data we sum over is $x,y,z$ and $\bs\Ga_j$. Here $x,y,z$ must satisfy $x+y+z=n$ and $l_{j-1}\le x\le x+y\le l_j$. Also as $\bs\Ga_{j-1},\bs\Ga_{j+1}$ are fixed and $R(\bs\Ga_{j-1})=L(\bs\Ga_j)$, $R(\bs\Ga_j)=L(\bs\Ga_{j+1})$, we see that $L(\bs\Ga_j),R(\bs\Ga_j)$ are fixed. (This happens even if $j=1$ or $j=k$, as $L(\bs\Ga_1)=\ti I$, $R(\bs\Ga_k)=\{\ti\imath_0,\ti\imath_{\ti M}\}$.) Let us write
\e
L(\bs\Ga_j)=\{\ti\imath_{b_0},\ldots,\ti\imath_{b_K}\},\quad R(\bs\Ga_j)=\{\ti\imath_{b_{c_0}},\ldots,\ti\imath_{b_{c_L}}\},
\label{ai10eq36}
\e
where $0=b_0<b_1<\cdots<b_K=\ti M$ and $0=c_0<c_1<\cdots<c_L=K$.

Now consider the equation
\ea
&\sum_{\begin{subarray}{l}
x,z\ge 0, \; y\ge 1:\\ x+y+z=n \\
l_{j-1}\le x\le x+y\le l_j,\\
\text{$\bs\Ga_j$ balanced with} \\
\text{$l_j-l_{j-1}-y$ leaves,} \\
\text{\eq{ai10eq36} holds}
\end{subarray}}
\begin{aligned}[t]
&\biggl(\prod_{\begin{subarray}{l} \text{$v,w$ vertices in $\bs\Ga_j$} \\
\text{$v,w,p,q,r$ usual conditions}
\end{subarray}\!\!\!\!\!\!\!\!\!\!\!\!\!\!\!\!\!\!\!\!\!\!\!\!\!\!\!\!\!\!\!\!\!\!\!\!\!\!\!\!\!\!\!\!\!\!\!\!\!\!\!\!\!\!\!} (-1)^{\begin{subarray}{l} N(v,w)+q(r+1) \\ +(p+1)(\md{I_w}+r+1)\end{subarray}}\biggr)\cdot{}\\
&C(\bs\Ga_j)^{x_{l_{j-1}},\ldots,x_x,x_{x+y},\ldots,x_{l_j}}\\
&{}\ci\bigl((-1)^{(x-l_{j-1})+y(z-n+l_j)}\cdot{}\\
&\id^{\ot^{x-l_{j-1}}}\ot\m_y^{\{\ti\imath_0\};x_x,\ldots,x_{x+y}}\ot \id^{\ot^{l_j-x-y}}\bigr)
\end{aligned}
\allowdisplaybreaks
\label{ai10eq37}\\
&=\sum_{\begin{subarray}{l} 
k_0,\ldots,k_{L}\ge 0, \\ k_L=1, \\
w_0^i,\ldots,w_{k_i}^i\in\sO^{\ti\imath_{b_{c_i}}},\\ 
0\le i\le L \end{subarray}}
\begin{aligned}[t]&\bigl\{(-1)^{\sum_{i=1}^{L-1}(c_i-i)}
\d_{k_{L-1},k_{L};w_0^L,\ldots,w_{k_L}^L}^{\{\ti\imath_{b_{c_{L-1}}},\ti\imath_{b_{c_{L-1}+1}},\ldots,\ti\imath_{b_{c_L}}\};w_0^{L-1},\ldots,w_{k_{L-1}}^{L-1}}\\
&\ci\d_{k_{L-2},k_{L-1};w_0^{L-1},\ldots,w_{k_{L-1}}^{L-1}}^{\{\ti\imath_{b_{c_{L-2}}},\ti\imath_{b_{c_{L-2}+1}},\ldots,\ti\imath_{b_{c_{L-1}}}\};w_0^{L-2},\ldots,w_{k_{L-2}}^{L-2}}\\
&{}\ci\cdots\ci\d_{k_0,k_1;w_0^1,\ldots,w_{k_1}^1}^{\{\ti\imath_{b_{c_0}},\ti\imath_{b_{c_0+1}},\ldots,\ti\imath_{b_{c_1}}\};w_0^0,\ldots,w_{k_0}^0}\bigr\}\\
&\ci\d_{l_j-l_{j-1},k_0;w_0^0,\ldots,w_{k_0}^0}^{\{\ti\imath_0\};x_{l_{j-1}},\ldots,x_{l_j}}
\end{aligned}
\allowdisplaybreaks
\nonumber\\
&=\sum_{\begin{subarray}{l} 
k_1,\ldots,k_L\ge 0, \\
w_0^i,\ldots,w_{k_i}^i\in\sO^{\ti\imath_{b_{c_i}}},\\ 
1\le i\le L \end{subarray}}
\begin{aligned}[t]
&(-1)^{K-L}\cdot\m_{k_L}^{\{\ti\imath_{\ti M}\};w_0^L,\ldots,w_{k_L}^L}\\
&{}\ci\bigl\{(-1)^{\sum_{i=1}^{L-1}(c_i-i)}\d_{k_{L-1},k_{L};w_0^L,\ldots,w_{k_L}^L}^{\{\ti\imath_{b_{c_{L-1}}},\ti\imath_{b_{c_{L-1}+1}},\ldots,\ti\imath_{b_{c_L}}\};w_0^{L-1},\ldots,w_{k_{L-1}}^{L-1}}\\
&\ci\d_{k_{L-2},k_{L-1};w_0^{L-1},\ldots,w_{k_{L-1}}^{L-1}}^{\{\ti\imath_{b_{c_{L-2}}},\ti\imath_{b_{c_{L-2}+1}},\ldots,\ti\imath_{b_{c_{L-1}}}\};w_0^{L-2},\ldots,w_{k_{L-2}}^{L-2}}\\
&{}\ci\cdots\ci\d_{l_j-l_{j-1},k_1;w_0^1,\ldots,w_{k_1}^1}^{\{\ti\imath_{b_{c_0}},\ti\imath_{b_{c_0+1}},\ldots,\ti\imath_{b_{c_1}}\};x_{l_{j-1}},\ldots,x_{l_j}}\bigr\}\\
\end{aligned}
\allowdisplaybreaks
\nonumber\\
&=\sum_{\begin{subarray}{l} l_{j-1}=m_0\le m_1\le \cdots\le m_t=l_j, \\
\text{$\bs\Ga'_1,\ldots,\bs\Ga'_t$ with} \\
\text{usual conditions:} \\
\text{$L(\bs\Ga'_1)=\{\ti\imath_{b_0},\ldots,\ti\imath_{b_K}\}$,} \\
\text{$R(\bs\Ga'_t)=\{\ti\imath_{b_{c_0}},\ldots,\ti\imath_{b_{c_L}}\}$}
\end{subarray}}
\begin{aligned}[t]
&(-1)^{K-L}(-1)^{(t+1)(K-L+l_j-l_{j-1})}\cdot{}\\
&\biggl(\prod_{\begin{subarray}{l} \text{$v,w,p,q,r$ usual conditions}
\end{subarray}\!\!\!\!\!\!\!\!\!\!\!\!\!\!\!\!\!\!\!\!\!\!\!\!\!\!\!\!\!\!\!\!\!\!\!\!\!\!\!\!\!\!\!\!\!\!\!\!\!\!\!\!\!\!\!} (-1)^{\begin{subarray}{l} N(v,w)+q(r+1) \\ +(p+1)(\md{I_w}+r+1)\end{subarray}}\biggr)\cdot{}\\
&\m_t^{\{\ti\imath_{\ti M}\},\ldots}\ci \bigl(C(\bs\Ga'_1)^{x_{m_0},\ldots,x_{m_1}} \\
& {}\ot\cdots\ot C(\bs\Ga'_t)^{x_{m_{t-1}},\ldots,x_{m_t}}\bigr).
\end{aligned}
\nonumber
\ea
Here the first line is a part of $(b2)$, summed over the data that we have not fixed. To prove the first step, first note from \eq{ai7eq8} that
\begin{align*}
&\sum_{x,y,z}(-1)^{(x-l_{j-1})+y(z-n+l_j)}\cdot\bigl(\id^{\ot^{x-l_{j-1}}}\ot\m_y^{\{\ti\imath_0\};x_x,\ldots,x_{x+y}}\ot \id^{\ot^{l_j-x-y}}\bigr)\\
&=\sum_{k_0,\; w_0^0,\ldots,w_{k_0}^0\in\sO^{\ti\imath_{b_{c_0}}}}\d_{l_j-l_{j-1},k_0;w_0^0,\ldots,w_{k_0}^0}^{\{\ti\imath_0\};x_{l_{j-1}},\ldots,x_{l_j}},
\end{align*}
setting $k_0=l_j-l_{j-1}-y+1$ and $(w_0^0,\ldots,w_{k_0}^0)=(x_{l_{j-1}},\ldots,x_x,x_{x+y},\ldots,x_{l_j})$. Then for these fixed $k_0,(w_0^0,\ldots,w_{k_0}^0)$, we can show following the proof of Theorem \ref{ai7thm1}(b) with $k'=1$ in \S\ref{ai101} that $\sum_{\bs\Ga_j}\pm C(\bs\Ga_j)^{w_0^0,\ldots,w_{k_0}^0}$ is the sum of terms $\bigl\{\d_{k_{L-1},k_{L};\cdots}^{\cdots}\ci\cdots\ci\d_{k_0,k_1;\cdots}^{\cdots}\bigr\}$ in the second step of~\eq{ai10eq37}.

For the second step of \eq{ai10eq37}, we apply equation \eq{ai8eq6} $L$ times. The first application gives
\begin{align*}
&\sum_{k_0\ge 0,\; w_0^0,\ldots,w_{k_0}^0\in\sO^{\ti\imath_{b_{c_0}}}\!\!\!\!\!\!\!\!\!\!\!\!\!}\d_{k_0,k_1;w_0^1,\ldots,w_{k_1}^1}^{\{\ti\imath_{b_{c_0}},\ti\imath_{b_{c_0+1}},\ldots,\ti\imath_{b_{c_1}}\};w_0^0,\ldots,w_{k_0}^0}\ci\d_{l_j-l_{j-1},k_0;w_0^0,\ldots,w_{k_0}^0}^{\{\ti\imath_0\};x_{l_{j-1}},\ldots,x_{l_j}}\\
&=(-1)^{c_1-c_0-1}\sum_{k'_1\ge 0,\; w_0^{\prime 1},\ldots,w_{k_1'}^{\prime 1}\in\sO^{\ti\imath_{b_{c_1}}}\!\!\!\!\!\!\!\!\!\!\!\!\!\!\!\!\!\!\!\!\!\!\!\!\!\!}
\d_{k'_1,k_1;w_0^1,\ldots,w_{k_1}^1}^{\{\ti\imath_{b_{c_1}}\};w_0^{\prime 1},\ldots,w_{k_1'}^{\prime 1}}\ci\d_{l_j-l_{j-1},k'_1;w_0^{\prime 1},\ldots,w_{k_1'}^{\prime 1}}^{\{\ti\imath_{b_{c_0}},\ti\imath_{b_{c_0+1}},\ldots,\ti\imath_{b_{c_1}}\};x_{l_{j-1}},\ldots,x_{l_j}}.
\end{align*}
In this way we move the initial $\d_{\cdots}^{\{\ti\imath_0\};\cdots}$ through each factor $\{\d_{\cdots}^{\{\ti\imath_{b_{c_{i-1}}},\ldots,\ti\imath_{b_{c_i}}\};\cdots}$ for $i=1,\ldots,L$ until it becomes $\d_{k_L,1;\cdots}^{\{\ti\imath_{\ti M}\};\cdots}$, which we write as $\m_{k_L}^{\{\ti\imath_{\ti M}\};\cdots}$ in the third line of \eq{ai10eq37}. The overall sign change is 
\begin{equation*}
(-1)^{c_1-c_0-1}\cdot(-1)^{c_2-c_1-1}\cdots (-1)^{c_L-c_{L-1}-1}=(-1)^{c_L-c_0-L}=(-1)^{K-L}.
\end{equation*}
Note that the conditions on $L(\bs\Ga_j),R(\bs\Ga_j)$ imply that we only use \eq{ai8eq6} for $\es\ne I\subseteq\ti I$ which we have assumed \eq{ai8eq6} holds in Definition \ref{ai10def3}. For the third step of \eq{ai10eq37}, we follow the proof of Theorem \ref{ai7thm1}(b) in \S\ref{ai101}, as in the first step. 

To prove that $(b2)=(c2)$ in \eq{ai10eq35}, for each choice of fixed data $j,k,l_0,\ab\ldots,\ab l_k,\ab\bs\Ga_1,\ab\ldots,\bs\Ga_{j-1},\bs\Ga_{j+1},\ldots,\bs\Ga_k$ in $(b2)$, we use \eq{ai10eq37} to replace the sum over the remaining data $x,y,z,\bs\Ga_j$ of the terms $\pm C(\bs\Ga_j)\ci(\id\ot\m_t^{\{\ti\imath_0\};\cdots}\ot \id)$ in $(b2)$ by the last line of \eq{ai10eq37}. This replaces a sum over sequences of graphs $\bs\Ga_1,\ldots,\bs\Ga_k$ in $(b2)$ by a sum over sequences of graphs $\bs\Ga_1,\ldots,\bs\Ga_{j-1},\bs\Ga'_1,\ldots,\bs\Ga'_t,\bs\Ga_{j+1},\ldots,\bs\Ga_k$, since $\bs\Ga_j$ is replaced by $\bs\Ga'_1,\ldots,\bs\Ga'_t$.

In $(c2)$ we change the notation $\bs\Ga_1,\ldots,\bs\Ga_{j-1},\bs\Ga'_1,\ldots,\bs\Ga'_t,\bs\Ga_{j+1},\ldots,\bs\Ga_k$ to $\bs\Ga_1,\ab\ldots,\ab\bs\Ga_s,\ab\bs\Ga_{s+1},\ab\ldots,\bs\Ga_{s+t},\bs\Ga_{t+1},\ldots,\bs\Ga_m$, so that $s=j-1$ and $m=k+t-1$. Thus $(L(\bs\Ga_j),R(\bs\Ga_j))$ in $(b2)$ becomes $(L(\bs\Ga_{s+1}), R(\bs\Ga_{s+t}))$ in $(c2)$. Therefore the condition in $(b2)$ that $(L(\bs\Ga_j),R(\bs\Ga_j))\ne(\ti I,\{\ti\imath_0,\ti\imath_{\ti M}\}),(\ti I,\{\ti\imath_0,\ti\imath_1,\ti\imath_{\ti M}\}),(\ti I\sm\{\ti\imath_a\},\{\ti\imath_0,\ti\imath_{\ti M}\})$ becomes the new condition in $(c2)$ that $(L(\bs\Ga_{s+1}), R(\bs\Ga_{s+t}))\ne(\ti I,\{\ti\imath_0,\ti\imath_{\ti M}\}),(\ti I,\{\ti\imath_0,\ti\imath_1,\ti\imath_{\ti M}\}),(\ti I\sm\{\ti\imath_a\},\{\ti\imath_0,\ti\imath_{\ti M}\})$.

There are several signs to keep track of in passing from $(b2)$ to $(c2)$:
\begin{itemize}
\setlength{\itemsep}{0pt}
\setlength{\parsep}{0pt}
\item[(i)] In $(b2)$ we have a sign $(-1)^{\ti M+1}\cdot(-1)^{x+yz}\cdot(-1)^{(k+1)(\ti M+x+z+2)}$.
\item[(ii)] The first line of \eq{ai10eq37} has sign $(-1)^{(x-l_{j-1})+y(z-n+l_j)}$.
\item[(iii)] The last line of \eq{ai10eq37} has sign $(-1)^{K-L}(-1)^{(t+1)(K-L+l_j-l_{j-1})}$.
\item[(iv)] To put $C(\bs\Ga_j)$ and $\id\ot\m_y^{\{\ti\imath_0\},\ldots}\ot\id$ next to each other in $(b2)$, in order to apply \eq{ai10eq37}, we should commute $\m_y^{\{\ti\imath_0\},\ldots}$ with $C(\bs\Ga_{j+1})\ot\cdots\ot C(\bs\Ga_k)$. As in \S\ref{ai31}, this gives a sign $(-1)^{y(j+k+n+l_j+L+1)}$.
\item[(v)] After substituting \eq{ai10eq37} into $(b2)$, to get $(c2)$ we must commute $\m_t^{\{\ti\imath_{\ti M}\},\ldots}$ with $C(\bs\Ga_1)\ot\cdots\ot C(\bs\Ga_{j-1})$. This gives a sign $(-1)^{t(j-1+l_{j-1}+\ti M-K)}$.
\item[(vi)] In $(c2)$ we have a sign $(-1)^{s+tu}\cdot(-1)^{(m+1)(\ti M+n+1)}$.
\item[(vii)] Consider the difference in signs $\prod_{v,w}(-1)^{N(v,w)+\cdots}$ with $v$ in $\bs\Ga_1,\ldots,\bs\Ga_{j-1}$ and $w$ in $\bs\Ga_j$ in $(b2)$, and signs $\prod_{v,w}(-1)^{N(v,w)+\cdots}$ with $v$ in $\bs\Ga_1,\ldots,\bs\Ga_{j-1}$ and $w$ in $\bs\Ga_1',\ldots,\bs\Ga_t'$ in $(c2)$. Since $\sum_w(r+1)$ changes by $y+t$ mod 2 between $\bs\Ga_j$ and $\bs\Ga_1',\ldots,\bs\Ga_t'$,  and $\sum_vq=l_{j-1}$ and $\sum_v(p+1)=l_{j-1}+j+1$ mod 2, this gives a sign $(-1)^{(j+1)(y+t)}$.
\item[(viii)] Consider the difference in signs $\prod_{v,w}(-1)^{N(v,w)+\cdots}$ with $v$ in $\bs\Ga_j$ and $w$ in $\bs\Ga_{j+1},\ldots,\bs\Ga_k$ in $(b2)$, and signs $\prod_{v,w}(-1)^{N(v,w)+\cdots}$ with $v$ in $\bs\Ga_1',\ldots,\bs\Ga_t'$ and $w$ in $\bs\Ga_{j+1},\ldots,\bs\Ga_k$ in $(c2)$. Since $\sum_wq$ changes by $y-1$ mod 2 and $\sum_v(p+1)$ changes by $y+t$ mod 2 between $\bs\Ga_j$ and $\bs\Ga_1',\ldots,\bs\Ga_t'$, and $\sum_w\md{I_w}=L+1\mod 2$, and $\sum_w(r+1)=n-l_j+k-j\mod 2$, this gives a sign $(-1)^{(y-1)(n-l_j+k-j)}(-1)^{(y+t)(L+1+n-l_j+k-j)}$.
\end{itemize}
Taking the product of the signs in (i)--(viii), using $x+y+z=n$, $s=j-1$, $m=k+t+1$, $s+t+u=m$, and simplifying, yields 1. Hence~$(b2)=(c2)$.

For the third step of \eq{ai10eq35}, to show $(c1)=(d1)$ we note that for $k\ge 0$ and $y_0,\ldots,y_k\in\sO^{\ti\imath_{\ti M}}$ fixed, by \eq{ai5eq8} or \eq{ai7eq8} and \eq{ai8eq4} we have
\e
\sum_{s,u\ge 0, \; t\ge 1: \; s+t+u=k\!\!\!\!\!\!\!\!\!\!\!\!\!\!\!\!\!\!\!\!\!}
(-1)^{s+tu}\cdot\m_{s+u+1}^{\{\ti\imath_{\ti M}\};\ldots}\ci
\bigl(\id^{\ot^s}\ot\m_t^{\{\ti\imath_{\ti M}\};\ldots}\ot \id^{\ot^u}\bigr)^{y_0,\ldots,y_k}=0.
\label{ai10eq38}
\e
But in $(c1)$ we do not sum over all $s,t,u$, but only those with $L(\bs\Ga_{s+1})=\ti I$, $R(\bs\Ga_{s+t})=\{\ti\imath_0,\ti\imath_{\ti M}\}$. Subtracting the sum over all $s,t,u$ and using \eq{ai10eq38} shows that $(c1)=(d1)$. Also $(d2)$ is $(c2)$, with no change.

The fourth step of \eq{ai10eq35} is immediate, comparing conditions on $(L(\bs\Ga_{s+1}),\ab R(\bs\Ga_{s+t}))$ in the three sums.

For the fifth step of \eq{ai10eq35}, we split the fourth step into the sum $(e1)$ of those terms with $(L(\bs\Ga_{s+1}), R(\bs\Ga_{s+t}))=(\ti I,\{\ti\imath_0,\ti\imath_1,\ti\imath_{\ti M}\})$, which yield the first line $(f1)$ of the fifth step, and the sum $(e2)$ over $a=1,\ldots,\ti M-1$ of those terms with $(L(\bs\Ga_{s+1}),R(\bs\Ga_{s+t}))=(\ti I\sm\{\ti\imath_a\},\{\ti\imath_0,\ti\imath_{\ti M}\})$, which yield the sum over $a=1,\ldots,\ti M-1$ in the second line $(f2)$ of the fifth step.

To prove that $(e1)=(f1)$, consider the trees $\bs\Ga_{s+t+1},\ldots,\bs\Ga_k$ in a term in the sum in the fourth step with $(L(\bs\Ga_{s+1}),R(\bs\Ga_{s+t}))=(\ti I,\{\ti\imath_0,\ti\imath_1,\ti\imath_{\ti M}\})$. We have $L(\bs\Ga_{s+t+1})=R(\bs\Ga_{s+t})=\{\ti\imath_0,\ti\imath_1,\ti\imath_{\ti M}\}$, and $R(\bs\Ga_k)=\{\ti\imath_0,\ti\imath_{\ti M}\}$. Since $L(\bs\Ga_{s+t+1})\supseteq L(\bs\Ga_{s+t+2})\supseteq L(\bs\Ga_k)$, we see that there is a unique $i=s+t+1,\ldots,k$ such that $L(\bs\Ga_i)=\{\ti\imath_0,\ti\imath_1,\ti\imath_{\ti M}\}$ and $R(\bs\Ga_i)=\{\ti\imath_0,\ti\imath_{\ti M}\}$. The correspondence between the sum $(e1)$, and the definition \eq{ai4eq17} of $\n_2(\cdots)_n^{x_0,\ldots,x_n}$ in $(f1)$, works like this:
\begin{itemize}
\setlength{\itemsep}{0pt}
\setlength{\parsep}{0pt}
\item[(A)] As $L(\bs\Ga_i)=\{\ti\imath_0,\ti\imath_1,\ti\imath_{\ti M}\}$, $R(\bs\Ga_i)=\{\ti\imath_0,\ti\imath_{\ti M}\}$, $\bs\Ga_i$ has only one internal vertex $v$ with $I_v=\{\ti\imath_0,\ti\imath_1,\ti\imath_{\ti M}\}$, and $C(\bs\Ga_i)^{x_{l_{i-1}},\ldots,x_{l_i}}=\m_{l_i-l_{i-1}}^{\{\ti\imath_0,\ti\imath_1,\ti\imath_{\ti M}\};x_{l_{i-1}},\ldots,x_{l_i}}$. These terms $C(\bs\Ga_i)$ give the second input $\m^{\{\ti\imath_0,\ti\imath_1,\ti\imath_{\ti M}\}}$ in $\n_2(\cdots)$ in $(f1)$.
\item[(B)] The sum of $\pm \m_t^{\{\ti\imath_{\ti M}\};\ldots}\ci\bigl( C(\bs\Ga_{s+1})\ot\cdots\ot C(\bs\Ga_{s+t})\bigr)$ in $(e1)$ over choices of $t,\bs\Ga_{s+1},\ldots,\bs\Ga_{s+t}$, where the sign $\pm$ is $\prod_{v,w}(-1)^{N(v,w)+\cdots}$ taken over $v,w$ both lying in $\bs\Ga_{s+1},\ldots,\bs\Ga_{s+t}$, is $(\fF^{\ti\imath_0\ti\imath_1})^*({\rm E}^{\ti I\sm\{\ti\imath_0\}})$ in $(f1)$.

Here as $(L(\bs\Ga_{s+1}),R(\bs\Ga_{s+t}))=(\ti I,\{\ti\imath_0,\ti\imath_1,\ti\imath_{\ti M}\})$, for each $j=s+1,\ldots,s+t$, the sequences $L(\bs\Ga_j),R(\bs\Ga_j)$ both start $\ti\imath_0,\ti\imath_1,\ldots,$ so all vertices $v$ in $\bs\Ga_j$ with leaves must have $I_v=\{\ti\imath_0,\ti\imath_1\}$. By `sawing' the trees $\bs\Ga_{s+1},\ldots,\bs\Ga_{s+t}$ at level $\ti\imath_1$, as in Definition \ref{ai6def7}, we break the sum above into a factor coming from the interval $[\ti\imath_0,\ti\imath_1]$, which is the operation $(\fF^{\ti\imath_0\ti\imath_1})^*$, and a factor coming from $[\ti\imath_1,\ti\imath_{\ti M}]$, which is~${\rm E}^{\ti I\sm\{\ti\imath_0\}}$.
\item[(C)] The definition \eq{ai4eq17} of $\n_2$ involves products of $\fF^0,\fF^1,\fF^2$ terms. Here we have $\fF^0=\fF^{\ti I}$, $\fF^1=\fF^{\{\ti\imath_0,\ti\imath_1,\ti\imath_{\ti M}\}}$ and $\fF^2=\fF^{\ti\imath_0\ti\imath_{\ti M}}$.

The $\fF^0$ terms correspond to the sum of $\pm C(\bs\Ga_1)\ot\cdots\ot C(\bs\Ga_s)$ over choices of $s,\bs\Ga_1,\ldots,\bs\Ga_s$ in $(e1)$.

The $\fF^1$ terms correspond to the sum of $\pm C(\bs\Ga_{s+t+1})\ot\cdots\ot C(\bs\Ga_{i-1})$ over choices of $i-s-t,\bs\Ga_{s+t+1},\ldots,\bs\Ga_{i-1}$ in $(e1)$.

The $\fF^2$ terms correspond to the sum of $\pm C(\bs\Ga_{i+1})\ot\cdots\ot C(\bs\Ga_k)$ over choices of $k-i,\bs\Ga_{i+1},\ldots,\bs\Ga_k$ in $(e1)$.
\item[(D)] Including the overall minus signs, the signs in $(e1)$, and in the definition of $\n_2(\cdots)_n^{x_0,\ldots,x_n}$, turn out to agree.
\end{itemize}

To prove that $(e2)=(f2)$, consider the trees $\bs\Ga_1,\ldots,\bs\Ga_s$ in a term in the sum in the fourth step with $(L(\bs\Ga_{s+1}),R(\bs\Ga_{s+t}))=(\ti I\sm\{\ti\imath_a\},\{\ti\imath_0,\ti\imath_{\ti M}\})$. We have $L(\bs\Ga_1)=\ti I$, and $R(\bs\Ga_s)=L(\bs\Ga_{s+1})=\ti I\sm\{\ti\imath_a\}$. Since $L(\bs\Ga_1)\supseteq L(\bs\Ga_2)\supseteq L(\bs\Ga_{s+1})$, we see that there is a unique $i=1,\ldots,s$ such that $L(\bs\Ga_i)=\ti I$ and $R(\bs\Ga_i)=\ti I\sm\{\ti\imath_a\}$. The correspondence between the sum $(e2)$, and the definition \eq{ai4eq17} of $\n_2(\cdots)_n^{x_0,\ldots,x_n}$ in $(f2)$, now works like this:
\begin{itemize}
\setlength{\itemsep}{0pt}
\setlength{\parsep}{0pt}
\item[(A$)'$] The sum of $\pm C(\bs\Ga_i)$ in $(e2)$ over possible choices of $\bs\Ga_i$ with $L(\bs\Ga_i)=\ti I$, $R(\bs\Ga_i)=\ti I\sm\{\ti\imath_a\}$, where the sign $\pm$ is $\prod_{v,w}(-1)^{N(v,w)+\cdots}$ taken over vertices $v,w$ both lying in $\bs\Ga_i$, is $(\fF^{\{\ti\imath_0,\ldots,\ti\imath_{a-1}\}})^*\ab\ci\ab(\fF^{\{\ti\imath_{a+1},\ldots,\ti\imath_{\ti M}\}})_*\ab (\m^{\{\ti\imath_{a-1},\ti\imath_a,\ti\imath_{a+1}\}})$ in $(f2)$.
\item[(B$)'$] The sum of $\pm \m_t^{\{\ti\imath_{\ti M}\};\ldots}\ci\bigl( C(\bs\Ga_{s+1})\ot\cdots\ot C(\bs\Ga_{s+t})\bigr)$ in $(e2)$ over choices of $t,\bs\Ga_{s+1},\ldots,\bs\Ga_{s+t}$, where the sign $\pm$ is $\prod_{v,w}(-1)^{N(v,w)+\cdots}$ taken over $v,w$ both lying in $\bs\Ga_{s+1},\ldots,\bs\Ga_{s+t}$, is ${\rm E}^{\ti I\sm\{\ti\imath_a\}}$ in $(f2)$.
\item[(C$)'$] The definition \eq{ai4eq17} of $\n_2$ involves products of $\fF^0,\fF^1,\fF^2$ terms. Here we have $\fF^0=\fF^{\ti I}$, $\fF^1=\fF^{\ti I\sm\{\ti\imath_a\}}$ and $\fF^2=\fF^{\ti\imath_0\ti\imath_{\ti M}}$.

The $\fF^0$ terms correspond to the sum of $\pm C(\bs\Ga_1)\ot\cdots\ot C(\bs\Ga_{i-1})$ over choices of $i,\bs\Ga_1,\ldots,\bs\Ga_{i-1}$ in $(e2)$.

The $\fF^1$ terms correspond to the sum of $\pm C(\bs\Ga_{i+1})\ot\cdots\ot C(\bs\Ga_s)$ over choices of $s-i,\bs\Ga_{i+1},\ldots,\bs\Ga_s$ in $(e2)$.

The $\fF^2$ terms correspond to the sum of $\pm C(\bs\Ga_{s+t+1})\ot\cdots\ot C(\bs\Ga_k)$ over choices of $u=k-s-t,\bs\Ga_{s+t+1},\ldots,\bs\Ga_k$ in $(e2)$.
\item[(D$)'$] The signs in $(e2)$, and in the definition of $\n_2(\cdots)_n^{x_0,\ldots,x_n}$, turn out to agree, with the following exception. Consider signs $(-1)^{N(v,w)}$ in $(e2)$, in which $v$ is the vertex in $\bs\Ga_i$ labelled by $\{\ti\imath_{a-1},\ti\imath_a,\ti\imath_{a+1}\}$, and $w$ is any vertex in $\bs\Ga_{s+1},\ldots,\bs\Ga_{s+t}$. From the definition \eq{ai6eq12} of $N(v,w)$, we see that these signs $(-1)^{N(v,w)}$ count each pair $(\ti\imath_a,\ti\imath_b)$ for $b=a+1,a+2,\ldots,\ti M-1$ exactly once, giving a total sign~$(-1)^{\ti M-a-1}$.
\end{itemize}

This finishes the proof of \eq{ai10eq35}, which yields \eq{ai10eq32}. Equation \eq{ai10eq33} is immediate. The proof of Proposition \ref{ai10prop3} is complete.
\end{proof}

We can now prove Theorem \ref{ai8thm2} for $\AfCat^{\un,\iy}_R$. Definition \ref{ai8def3} shows that $\AfCat^{\un,\iy}_R$ is a simplicial set. By Definition \ref{ai2def4}, to prove $\AfCat^{\un,\iy}_R$ is an $\iy$-category we must show that every map of simplicial sets $f:\La^N_i\ra \AfCat^{\un,\iy}_R$ can be extended to an $N$-simplex $\ti f:\De^N\ra\AfCat^{\un,\iy}_R$ for all $N\ge 2$ and $0<i<N$. Unwinding the definitions, this means:
\begin{itemize}
\setlength{\itemsep}{0pt}
\setlength{\parsep}{0pt}
\item[$(*)_i^N$] Suppose we are given data $(\sO^i)_{i=0}^N,(f^{ij})_{0\le i\le j\le N},(\sH^i)_{i=0}^N$ and $\m^I$ for all $\es\ne I\subseteq\{0,\ldots,N\}$ with $I\ne\{0,\ldots,N\},\{0,\ldots,N\}\sm\{i\}$ (note this exclusion) as in Definitions \ref{ai8def1}(a)(i)--(v) and \ref{ai8def3}, such that equations \eq{ai7eq17}--\eq{ai7eq19} and \eq{ai8eq4}--\eq{ai8eq5} hold for all $\es\ne I\subseteq\{0,\ldots,N\}$ with $I\ne\{0,\ldots,N\},\{0,\ldots,N\}\sm\{i\}$ and~$n\ge 0$. 

Also, as in Definition \ref{ai8def3}(b), for every subset $\{i_0<i_1<i_2\}$ in $\{0,\ldots,N\}$ with $\{i_0,i_1,i_2\}\ne \{0,\ab\ldots,\ab N\},\{0,\ldots,N\}\sm\{i\}$ we assume $[\m^{\{i_0,i_1,i_2\}}]:\fF^{i_1i_2}\ci\fF^{i_0i_1}\ra\fF^{i_0i_2}$ is an isomorphism in $H^0(F^0\Fun_{A^\rfi_\iy}^\un(\sA^{i_0},\sA^{i_2}))$.

We must show that there exist (not necessarily unique) choices of $\m^I$ for $I=\{0,\ab\ldots,\ab N\},\{0,\ldots,N\}\sm\{i\}$, such that \eq{ai7eq17}--\eq{ai7eq19} and \eq{ai8eq4}--\eq{ai8eq5} also hold for $I=\{0,\ldots,N\},\{0,\ldots,N\}\sm\{i\}$ and~$n\ge 0$.

Also, if $N=2$ we must choose $\m^{\{0,1,2\}}$ with $[\m^{\{0,1,2\}}]$ an isomorphism in $H^0(F^0\Fun_{A^\rfi_\iy}^\un(\sA^0,\sA^2))$, and if $N=3$ we must choose $\m^{\{0,1,2,3\}\sm\{i\}}$ with $[\m^{\{0,1,2,3\}\sm\{i\}}]$ an isomorphism in $H^0(F^0\Fun_{A^\rfi_\iy}^\un(\sA^0,\sA^3))$.
\end{itemize}

Since $N=2,3$ are special cases, we first prove $(*)_i^N$ for $N\ge 4$, and then explain the modifications for $N=2,3$. Let $N\ge 4$ and $0<i<N$ be given, and suppose the hypotheses of $(*)_i^N$ hold. For consistency with Definition \ref{ai10def1} (which assumed $\m^I$ is defined for all $I$), we first make arbitrary choices of $\m^I$ for $I=\{0,\ldots,N\},\{0,\ldots,N\}\sm\{i\}$ (e.g.\ $\m^I=0$), which satisfy \eq{ai7eq19}, but which initially need not satisfy \eq{ai8eq5}. We will then explain how to modify these choices to make \eq{ai8eq5} hold for~$I=\{0,\ldots,N\},\{0,\ldots,N\}\sm\{i\}$.

Observe that the hypotheses of Definition \ref{ai10def3} hold with $\ti M=N$ and $\ti I=\{0,1,
\ldots,N\}$. Thus for $a=1,\ldots,N-1$ we have cohomology classes
\begin{equation*}
\bigl[{\rm E}^{\{0,\ldots,N\}\sm\{a\}}\bigr]\quad\text{in}\quad H^{4-N}\bigl(F^0\Hom_{A^\rfi_\iy}^\un(\fF^{\{0,\ldots,N\}\sm\{a\}},\fF^{0N}),\n_1\bigr),
\end{equation*}
and we also have a cohomology class
\begin{equation*}
\bigl[{\rm E}^{\{1,\ldots,N\}}\bigr]\quad\text{in}\quad H^{4-N}\bigl(F^0\Hom_{A^\rfi_\iy}^\un(\fF^{\{1,\ldots,N\}},\fF^{1N}),\n_1\bigr),
\end{equation*}
and Proposition \ref{ai10prop3}(c) shows that \eq{ai10eq33} holds.

Now if $a=1,\ldots,N-1$ with $a\ne i$ then the hypotheses of Definition \ref{ai10def2} hold with $\ti M=N-1$ and $\ti I=\{0,1,\ldots,N\}\sm\{a\}$, and as \eq{ai8eq5} holds for $I=\{0,\ldots,N\}\sm\{a\}$, Proposition \ref{ai10prop2} implies that $\bigl[{\rm E}^{\{0,\ldots,N\}\sm\{a\}}\bigr]=0$. Also \eq{ai8eq5} holds for $I=\{1,\ldots,N\}$, so $\bigl[{\rm E}^{\{1,\ldots,N\}}\bigr]=0$. Therefore \eq{ai10eq33} reduces to
\e
\!\!(-1)^{N-i}\bigl[{\rm E}^{\{0,\ldots,N\}\sm\{i\}}\bigr]\!\ci\!\bigl[(\fF^{\{0,\ldots,i-1\}})^*\!\ci\!(\fF^{\{i+1,\ldots,N\}})_*(\m^{\{i-1,i,i+1\}})\bigr]\!=\!0.\!
\label{ai10eq39}
\e
Since $\bigl[(\fF^{\{0,\ldots,i-1\}})^*\ci(\fF^{\{i+1,\ldots,N\}})_*(\m^{\{i-1,i,i+1\}})\bigr]$ is an isomorphism as in Definition \ref{ai10def3}, this implies that~$\bigl[{\rm E}^{\{0,\ldots,N\}\sm\{i\}}\bigr]=0$.

Proposition \ref{ai10prop2} with $\ti M=N-1$ and $\ti I=\{0,1,\ldots,N\}\sm\{i\}$ now implies that we can change the definition of $\m^{\{0,\ldots,N\}\sm\{i\}}$ (which was an arbitrary choice above, so can be changed) to make \eq{ai7eq19} and \eq{ai8eq5} hold with $\ti I=\{0,1,\ldots,N\}\sm\{i\}$, for all $n\ge 0$ and~$x_0,\ldots,x_n\in\sO^0$.

We have now arranged that \eq{ai7eq17}--\eq{ai7eq19} and \eq{ai8eq4}--\eq{ai8eq5} hold for all $\es\ne I\subseteq\{0,\ldots,N\}$ with $I\ne\{0,\ldots,N\}$ and $n\ge 0$. Hence the hypotheses of Definition \ref{ai10def2} hold with $\ti M=N$ and $\ti I=\{0,1,\ldots,N\}$, so Proposition \ref{ai10prop2} shows that $\n_1\bigl({\rm E}^{\{0,\ldots,N\}}\bigr)=0$, and we have a cohomology class
\begin{equation*}
[{\rm E}^{\{0,\ldots,N\}}]\quad\text{in}\quad H^{3-N}\bigl(F^0\Hom_{A^\rfi_\iy}^\un(\fF^{\{0,\ldots,N\}},\fF^{0N}),\n_1\bigr).
\end{equation*}

We now make a further modification of $\m^{\{0,\ldots,N\}\sm\{i\}}$. As in Proposition \ref{ai10prop2}(a), to make \eq{ai8eq5} hold with $\ti I=\{0,1,\ldots,N\}\sm\{i\}$, we had to choose $\m^{\{0,\ldots,N\}\sm\{i\}}$ such that $\n_1(\m^{\{0,\ldots,N\}\sm\{i\}})+{\rm E}^{\{0,\ldots,N\}\sm\{i\}}=0$. Thus we may change $\m^{\{0,\ldots,N\}\sm\{i\}}$ by adding on an element of $\Ker(\n_1)$, that is, a representative of a cohomology class in $H^{3-N}\bigl(F^0\Hom_{A^\rfi_\iy}^\un(\fF^{\{0,\ldots,N\}\sm\{i\}},\fF^{0N}),\n_1\bigr)$, and \eq{ai8eq5} with $\ti I=\{0,1,\ldots,N\}\sm\{i\}$ will still hold.

We change $\m^{\{0,\ldots,N\}\sm\{i\}}$ by adding on a representative of
\begin{equation*}
(-1)^{N-i}[{\rm E}^{\{0,\ldots,N\}}]\ci[(\fF^{\{0,\ldots,i-1\}})^*\ci
(\fF^{\{i+1,\ldots,N\}})_*(\m^{\{i-1,i,i+1\}})]^{-1},
\end{equation*}
where the inverse exists as \eq{ai10eq27} is an isomorphism. Since $\n_2$ in \eq{ai10eq30} is used to define composition in $H^*\bigl(F^0\Fun_{A^\rfi_\iy}^\un(\sA^{\ti\imath_0},\ab\sA^{\ti\imath_{\ti M}})\bigr)$, by Proposition \ref{ai10prop3}(a), the effect of this on the cohomology class $[{\rm E}^{\{0,\ldots,N\}}]$ is to add on
\begin{align*}
&(-1)^{N-i-1}\bigl((-1)^{N-i}[{\rm E}^{\{0,\ldots,N\}}]\!\ci\![(\fF^{\{0,\ldots,i-1\}})^*\!\ci\!
(\fF^{\{i+1,\ldots,N\}})_*(\m^{\{i-1,i,i+1\}})]^{-1}\bigr)\\
&\ci[(\fF^{\{0,\ldots,i-1\}})^*\ci(\fF^{\{i+1,\ldots,N\}})_*(\m^{\{i-1,i,i+1\}})]
=-[{\rm E}^{\{0,\ldots,N\}}].
\end{align*}

That is, after this modification of $\m^{\{0,\ldots,N\}\sm\{i\}}$, we now have $[{\rm E}^{\{0,\ldots,N\}}]=0$. Proposition \ref{ai10prop2} with $\ti M=N$ and $\ti I=\{0,1,\ldots,N\}$ now implies that we can change the definition of $\m^{\{0,\ldots,N\}}$ (which was an arbitrary choice above) to make \eq{ai7eq19} and \eq{ai8eq5} hold with $\ti I=\{0,1,\ldots,N\}$, for all $n\ge 0$ and $x_0,\ldots,x_n\in\sO^0$. This proves $(*)_i^N$ for all $N\ge 4$ and~$0<i<N$.

When $N=2$, we do not use the proof above, as Definition \ref{ai10def3} assumes $\ti M\ge 3$, and does not apply when $\ti I=\{0,1,2\}$. Instead we define $\m^{\{0,2\}}=\m^{\{1,2\}}\ci\m^{\{0,1\}}$ by composition of filtered $A_\iy$-functors in Definition \ref{ai5def9}, and we set $\m_0^{\{0,1,2\}}=\boo_{\m^{\{0,2\}}}:\m^{\{0,2\}}\Ra\m^{\{0,2\}}$ to be the identity filtered $A_\iy$-natural isomorphism in \eq{ai4eq19}. It is easy to check these satisfy all the conditions.

When $N=3$, the proof above for $N\ge 4$ works and produces choices of $\m^{\{0,\ldots,3\}}$ and $\m^{\{0,\ldots,3\}\sm\{i\}}$ such that \eq{ai7eq19} and \eq{ai8eq5} hold with $I=\{0,\ldots,3\},\ab\{0,\ldots,3\}\sm\{i\}$. It remains to show that $[\m^{\{0,1,2,3\}\sm\{i\}}]$ is an isomorphism. Observe that \eq{ai8eq5} with $I=\{0,\ldots,3\}$ may be rewritten in the form
\begin{align*}
&-\n_1\bigl(\m^{\{0,1,2,3\}}\bigr)-\n_2\bigl((\fF^{01})^*(\m^{\{1,2,3\}}),\m^{\{0,1,3\}}\bigr)\\
&+\n_2\bigl((\fF^{23})_*(\m^{\{0,1,2\}}),\m^{\{0,2,3\}}\bigr)=0.
\end{align*}
Passing to cohomology, this implies that
\e
[\m^{\{0,1,3\}}]\ci (\fF^{01})^*([\m^{\{1,2,3\}}])=
[\m^{\{0,2,3\}}]\ci 
(\fF^{23})_*([\m^{\{0,1,2\}}]).
\label{ai10eq40}
\e

When $i=1$, by assumption $[\m^{\{0,1,2\}}],[\m^{\{0,1,3\}}],[\m^{\{1,2,3\}}]$ are isomorphisms, so $(\fF^{01})^*([\m^{\{1,2,3\}}]),(\fF^{23})_*([\m^{\{0,1,2\}}])$ are isomorphisms. We see from \eq{ai10eq40} that $[\m^{\{0,2,3\}}]$ is an isomorphism, with inverse
\begin{equation*}
[\m^{\{0,2,3\}}]^{-1}=(\fF^{23})_*([\m^{\{0,1,2\}}])\ci(\fF^{01})^*([\m^{\{1,2,3\}}])^{-1}\ci[\m^{\{0,1,3\}}]^{-1}.
\end{equation*}
This proves $(*)_1^3$, and $(*)_2^3$ is similar. This completes the proof of Theorem \ref{ai8thm2} for $\AfCat^{\un,\iy}_R$. The remainder of Theorems \ref{ai6thm3}, \ref{ai7thm3}, \ref{ai8thm2} follow as at the beginning of~\S\ref{ai104}.

\subsection{Proofs of Theorems \ref{ai6thm4}, \ref{ai7thm4}, and \ref{ai8thm3}}
\label{ai105}

The proofs are extensions of those of Theorems \ref{ai6thm3}, \ref{ai7thm3}, and \ref{ai8thm2} in \S\ref{ai104}. We will prove Theorem \ref{ai8thm3} for $(\hAfCat^{\un,\iy}_R,\cT^\rfi)$. Working with $\hAfCat^{\un,\iy}_R$ involves using the $A_\iy$-category $F^0\Fun_{A^\rfi_\iy}^\un(\sA^{\ti\imath_0},\ab\sA^{\ti\imath_{\ti M}})$ over $F^0R$ from Definition \ref{ai5def10}(f), with morphisms $F^0\Hom_{A^\rfi_\iy}(\fF,\fG)$. The proof for $(\hAwCat^{\un,\iy}_R,\cT^\wf)$ is essentially the same, except that we omit $F^0(\cdots)$ throughout.

Theorem \ref{ai7thm4} is a special case of Theorem \ref{ai8thm3} for $(\hAfCat^{\un,\iy}_R,\cT^\rfi)$ with $\Ga=\{0,\iy\}$, as in Example \ref{ai5ex2}. Theorem \ref{ai6thm4} follows from Theorem \ref{ai7thm4} by restricting to $A_\iy$-categories $\sA$ with~$\sO=\{*\}$. Let $(\Ga,\pre)$ be a po-monoid and $R$ a complete $\Ga$-filtered ring. Define $\hAfCat^{\un,\iy}_R$ as in Definition \ref{ai8def1}, and $\AfCat^{\un,\iy}_R\subset\hAfCat^{\un,\iy}_R$ as in Definition~\ref{ai8def3}.

We will prove that $(\hAfCat^{\un,\iy}_R,\cT^\rfi)$ satisfies each of Definition \ref{ai2def21}(a)--(d) in \S\ref{ai1051}--\S\ref{ai1054}. These use the very long proof of Theorem \ref{ai8thm2} in \S\ref{ai104}. Thus $(\hAfCat^{\un,\iy}_R,\cT^\rfi)$ is an $\iy$-bicategory. The final claim that $\AfCat^{\un,\iy}_R=\Core(\hAfCat^{\un,\iy}_R,\cT^\rfi)$ is immediate from Definitions \ref{ai2def22} and~\ref{ai8def3}.

\subsubsection{$(\hAfCat^{\un,\iy}_R,\cT^\rfi)$ satisfies Definition \ref{ai2def21}(a)}
\label{ai1051}

Suppose $t:\De^2\ra \AfCat^{\un,\iy}_R\subset\hAfCat^{\un,\iy}_R$ is a $2$-simplex in $\AfCat^{\un,\iy}_R$.
We must show that for $N\ge 3$ and $0<i<N$, all horn-filling problems \eq{ai2eq26} in $\hAfCat^{\un,\iy}_R$ can be solved. Now \S\ref{ai104} proved that all horn-filling problems \eq{ai2eq5} in $\AfCat^{\un,\iy}_R$ can be solved. Consider running the proof of \S\ref{ai104} in $\hAfCat^{\un,\iy}_R$ rather than in $\AfCat^{\un,\iy}_R$, in particular, the part at the end of \S\ref{ai104} when we prove $(*)_i^N$ for~$N\ge 4$. 

The crucial point is around \eq{ai10eq39}. To fill a horn $f:\La^N_i\ra\AfCat^{\un,\iy}_R$, we use the fact that $\bigl[(\fF^{\{0,\ldots,i-1\}})^*\ci(\fF^{\{i+1,\ldots,N\}})_*(\m^{\{i-1,i,i+1\}})\bigr]$ is an isomorphism, since the 2-simplex $f\ci\De^{\{i-1,i,i+1\}}:\De^2\ra\AfCat^{\un,\iy}_R$ satisfies the invertibility condition in Definition \ref{ai8def3}. This is the only point in the proof of filling $f$ where the invertibility condition is used. Thus, the same proof shows we can fill a horn $f:\La^N_i\ra\hAfCat^{\un,\iy}_R$ in $\hAfCat^{\un,\iy}_R$ for $N\ge 4$, provided $t=f\ci\De^{\{i-1,i,i+1\}}$ maps $\De^2$ to $\AfCat^{\un,\iy}_R\subset\hAfCat^{\un,\iy}_R$, as in \eq{ai2eq26}. When $N=3$, this proof also works, as the extra argument in \S\ref{ai104} showing that $\ti f$ maps $\De^3$ to $\AfCat^{\un,\iy}_R\subset\hAfCat^{\un,\iy}_R$ is not relevant. 

\subsubsection{$(\hAfCat^{\un,\iy}_R,\cT^\rfi)$ satisfies Definition \ref{ai2def21}(b)}
\label{ai1052}

A morphism $f:\La^2_i\ra\hAfCat^{\un,\iy}_R$ is equivalent to unital filtered $A_\iy$-categories $\sA^0,\sA^1,\sA^2$ and unital filtered $A_\iy$-functors $\fF^{01}:\sA^0\ra\sA^1$, $\fF^{12}:\sA^1\ra\sA^2$. As in Example \ref{ai8ex1}(c), extending this to $\ti f:\De^2\ra\hAfCat^{\un,\iy}_R$ as in \eq{ai2eq27} is equivalent to also choosing a unital filtered $A_\iy$-functor $\fF^{02}:\sA^0\ra\sA^2$ and a unital filtered $A_\iy$-natural transformation $\fH^{012}:\fF^{12}\ci\fF^{01}\Ra\fF^{02}$, and $\ti f$ is thin if $\fH^{012}$ is a filtered $A_\iy$-natural isomorphism.

We can choose $\fF^{02}=\fF^{12}\ci\fF^{01}$ and $\fH^{012}=\id_{\fF^{12}\ci\fF^{01}}$, which is a filtered $A_\iy$-natural isomorphism. Thus Definition \ref{ai2def21}(b) holds.

\subsubsection{$(\hAfCat^{\un,\iy}_R,\cT^\rfi)$ satisfies Definition \ref{ai2def21}(c)}
\label{ai1053}

Suppose $f:\De^4\ra \hAfCat^{\un,\iy}_R$ is a morphism and $f\ci\De^{\{i,j,k\}}:\De^2\ra \hAfCat^{\un,\iy}_R$ is thin for $\{i,j,k\}=\{0,2,4\},\{1,2,3\},\{0,1,3\},\{1,3,4\},\{0,1,2\}$. As usual, we consider the 0-simplices of $f(\De^4)$ to correspond to unital filtered $A_\iy$-categories $\sA^0,\ldots,\sA^4$, and the nondegenerate 1-simplices of $f(\De^4)$ to correspond to unital filtered $A_\iy$-functors $\fF^{ij}:\sA^i\ra\sA^j$ for $0\le i<j\le 4$, and the nondegenerate 2-simplices of $f(\De^4)$ to correspond to unital filtered $A_\iy$-natural transformations $\fH^{ijk}:\fF^{jk}\ci\fF^{ij}\Ra\fF^{ik}$ for $0\le i<j<k\le 4$.

We regard $\fF^{ij}$ as an object in the $A_\iy$-category $F^0\Fun_{A^\rfi_\iy}^\un(\sA^i,\sA^j)$, and $\fH^{ijk}$ as a morphism in $F^0\Fun_{A^\rfi_\iy}^\un(\sA^i,\sA^k)$, which induces a morphism $[\fH^{ijk}]$ in the cohomology category $H^0(F^0\Fun_{A^\rfi_\iy}^\un(\sA^i,\sA^k))$. Suppose $0\le i<j<k<l\le 4$. Then by the analogue of Example \ref{ai6ex3}(d) for the 3-simplex $f\ci\De^{\{i,j,k,l\}}(\De^3)$, in morphisms in the cohomology category $H^0(F^0\Fun_{A^\rfi_\iy}^\un(\sA^i,\sA^k))$, we have
\e
[\fH^{ijl}]\ci[\fH^{jkl}*\id_{\fF^{ij}}]=[\fH^{ikl}]\ci[\id_{\fF^{kl}}*\fH^{ijk}]:\fF^{kl}\ci\fF^{jk}\ci\fF^{ij}\longra\fF^{il}.
\label{ai10eq41}
\e

Consider the commutative diagram in $H^0(F^0\Fun_{A^\rfi_\iy}^\un(\sA^0,\sA^4))$:
\e
\begin{gathered}
\!\!\!\xymatrix@C=45pt@R=17pt{
*+[r]{\fF^{34}\!\ci\!((\fF^{23}\!\ci\!\fF^{12})\!\ci\!\fF^{01})} \ar[rr]_(0.65){\raisebox{-9pt}{$\st[\id_{\fF^{34}}*\fH^{123}*\id_{\fF^{01}}]$}} \ar@{=}[d] && \fF^{34}\!\ci\!(\fF^{13}\!\ci\!\fF^{01}) \ar[r]_(0.45){\raisebox{-9pt}{$\st[\id_{\fF^{34}}*\fH^{013}]$}} &
*+[l]{\fF^{34}\!\ci\!\fF^{03}} \ar[d]_{[\fH^{034}]}
\\
*+[r]{(\fF^{34}\!\ci\!\fF^{23})\!\ci\!(\fF^{12}\!\ci\!\fF^{01})} \ar[rr]^(0.7){[\fH^{234}*\fH^{012}]} && \fF^{24}\!\ci\!\fF^{02} \ar[r]^(0.45){[\fH^{024}]} & *+[l]{\fF^{04}.\!}
}\!\!\!
\end{gathered}
\label{ai10eq42}
\e
Using \eq{ai10eq41} we can show that this commutes. Since $\fH^{123},\fH^{013},\fH^{234},\fH^{012},\fH^{024}$ are filtered $A_\iy$-natural isomorphisms by assumption, the horizontal morphisms in \eq{ai10eq42} are isomorphisms by Definition \ref{ai5def12}(c). Hence the right hand morphism $[\fH^{034}]$ is an isomorphism, so $\fH^{034}$ is a filtered $A_\iy$-natural isomorphism, and $f\ci\De^{\{0,3,4\}}:\De^2\ra \hAfCat^{\un,\iy}_R$ is thin, as we have to prove.

To show $f\ci\De^{\{0,1,4\}}:\De^2\ra \hAfCat^{\un,\iy}_R$ is thin, we instead use the diagram
\begin{equation*}
\xymatrix@C=55pt@R=17pt{
*+[r]{(\fF^{34}\!\ci\!(\fF^{23}\!\ci\!\fF^{12}))\!\ci\!\fF^{01}} \ar[rr]_(0.62){\raisebox{-9pt}{$\st[\id_{\fF^{34}}*\fH^{123}*\id_{\fF^{01}}]$}} \ar@{=}[d] && (\fF^{34}\!\ci\!\fF^{13})\!\ci\!\fF^{01} \ar[r]_(0.45){\raisebox{-9pt}{$\st[\fH^{134}*\id_{\fF^{01}}]$}} &
*+[l]{\fF^{14}\!\ci\!\fF^{01}} \ar[d]_{[\fH^{014}]}
\\
*+[r]{(\fF^{34}\!\ci\!\fF^{23})\!\ci\!(\fF^{12}\!\ci\!\fF^{01})} \ar[rr]^(0.7){[\fH^{234}*\fH^{012}]} && \fF^{24}\!\ci\!\fF^{02} \ar[r]^(0.45){[\fH^{024}]} & *+[l]{\fF^{04}.\!} }
\end{equation*}

\subsubsection{$(\hAfCat^{\un,\iy}_R,\cT^\rfi)$ satisfies Definition \ref{ai2def21}(d)}
\label{ai1054}

Suppose $N\ge 3$, and as in \eq{ai2eq28}, with $s_{1,0}:\De^2\ra\De^1$ acting on vertices by $0\mapsto 0$, $1\mapsto 0$, $2\mapsto 1$, we are given a commutative diagram of morphisms `$\ra$'
\e
\begin{gathered}
\xymatrix@C=70pt@R=15pt{ *+[r]{\La^N_{0_{\phantom{(}}}} \ar@{^{(}->}[d]^\inc \ar[drr]_(0.3)f 
& \De^2 \ar[l]_{\De^{\{0,1,N\}}} \ar[r]^{s_{1,0}} & \De^1 \ar[d]_g
\\
*+[r]{\De^N} \ar@{..>}[rr]^(0.45){\ti f} && \hAfCat^{\un,\iy}_R. }
\end{gathered}
\label{ai10eq43}
\e
We must show it can completed with a morphism $\ti f$ as shown. As in \S\ref{ai1051}, we do this by adapting the argument in \S\ref{ai104}, but this time we must fill an {\it outer\/} horn $\La^N_0$, which were not considered in~\S\ref{ai104}.

The next part of the proof is similar to the last part of \S\ref{ai104}, starting from $(*^N_i)$. The map $f:\La^N_0\ra\hAfCat^{\un,\iy}_R$ gives us data $(\sO^i)_{i=0}^N,\ab(f^{ij})_{0\le i\le j\le N},\ab(\sH^i)_{i=0}^N$ and $\m^I$ as in $(*^N_0)$. (Section \ref{ai104} only considered $(*^N_i)$ for $0<i<N$, but $(*^N_0)$ still makes sense.) As in \S\ref{ai104}, for consistency with Definition \ref{ai10def1} (which assumed $\m^I$ is defined for all $I$), we first make arbitrary choices of $\m^I$ for $I=\{0,\ldots,N\},\{1,\ldots,N\}$ (e.g.\ $\m^I=0$), which satisfy \eq{ai7eq19}, but which initially need not satisfy \eq{ai8eq5}. We will then explain how to modify these choices to make \eq{ai8eq5} hold for~$I=\{0,\ldots,N\},\{1,\ldots,N\}$.

Observe that the hypotheses of Definition \ref{ai10def3} hold with $\ti M=N$ and $\ti I=\{0,1,
\ldots,N\}$. Thus for $a=1,\ldots,N-1$ we have cohomology classes
\begin{equation*}
\bigl[{\rm E}^{\{0,\ldots,N\}\sm\{a\}}\bigr]\quad\text{in}\quad H^{4-N}\bigl(F^0\Hom_{A^\rfi_\iy}^\un(\fF^{\{0,\ldots,N\}\sm\{a\}},\fF^{0N}),\n_1\bigr),
\end{equation*}
and we also have a cohomology class
\begin{equation*}
\bigl[{\rm E}^{\{1,\ldots,N\}}\bigr]\quad\text{in}\quad H^{4-N}\bigl(F^0\Hom_{A^\rfi_\iy}^\un(\fF^{\{1,\ldots,N\}},\fF^{1N}),\n_1\bigr),
\end{equation*}
and Proposition \ref{ai10prop3}(c) shows that \eq{ai10eq33} holds.

Now if $a=1,\ldots,N-1$ then the hypotheses of Definition \ref{ai10def2} hold with $\ti M=N-1$ and $\ti I=\{0,1,\ldots,N\}\sm\{a\}$, and as \eq{ai8eq5} holds for $I=\{0,\ldots,N\}\sm\{a\}$, Proposition \ref{ai10prop2} implies that $\bigl[{\rm E}^{\{0,\ldots,N\}\sm\{a\}}\bigr]=0$. Therefore \eq{ai10eq33} reduces to
\e
(\fF^{01})^*([{\rm E}^{\{1,2,\ldots,N\}}])\ci[\m^{\{0,1,N\}}]=0.
\label{ai10eq44}
\e
Since $f\ci \De^{\{0,1,N\}}(\De^2)$ is degenerate as in \eq{ai10eq47}, we have $\m^{\{0,1,N\}}=\id_{\fF^{0N}}$, so $[\m^{\{0,1,N\}}]$ is the identity. Thus \eq{ai10eq44} implies that $(\fF^{01})^*([{\rm E}^{\{1,2,\ldots,N\}}])=0$. But by degeneracy in \eq{ai10eq47}, we have $\sA^0=\sA^1$ and $\fF^{01}=\id_{\sA^0}$, so $(\fF^{01})^*$ is the identity. Hence $[{\rm E}^{\{1,2,\ldots,N\}}]=0$.

Proposition \ref{ai10prop2} with $\ti M=N-1$ and $\ti I=\{1,\ldots,N\}$ now implies that we can change the definition of $\m^{\{1,\ldots,N\}}$ (which was an arbitrary choice above, so can be changed) to make \eq{ai7eq19} and \eq{ai8eq5} hold with $I=\{1,2,\ldots,N\}$, for all $n\ge 0$ and~$x_0,\ldots,x_n\in\sO^1$.

We have now arranged that \eq{ai7eq17}--\eq{ai7eq19} and \eq{ai8eq4}--\eq{ai8eq5} hold for all $\es\ne I\subseteq\{0,\ldots,N\}$ with $I\ne\{0,\ldots,N\}$ and $n\ge 0$. Hence the hypotheses of Definition \ref{ai10def2} hold with $\ti M=N$ and $\ti I=\{0,1,\ldots,N\}$, so Proposition \ref{ai10prop2} shows that $\n_1\bigl({\rm E}^{\{0,\ldots,N\}}\bigr)=0$, and we have a cohomology class
\begin{equation*}
[{\rm E}^{\{0,\ldots,N\}}]\quad\text{in}\quad H^{3-N}\bigl(F^0\Hom_{A^\rfi_\iy}^\un(\fF^{\{0,\ldots,N\}},\fF^{0N}),\n_1\bigr).
\end{equation*}

We now make a further modification of $\m^{\{1,2,\ldots,N\}}$. As in Proposition \ref{ai10prop2}(a), to make \eq{ai8eq5} hold with $I=\{1,\ldots,N\}$, we had to choose $\m^{\{1,\ldots,N\}}$ such that $\n_1(\m^{\{1,\ldots,N\}})+{\rm E}^{\{1,\ldots,N\}}=0$. Thus we may change $\m^{\{1,\ldots,N\}}$ by adding on an element of $\Ker(\n_1)$, that is, a representative of a cohomology class in $H^{3-N}\bigl(F^0\Hom_{A^\rfi_\iy}^\un(\fF^{\{1,\ldots,N\}},\fF^{1N}),\n_1\bigr)$, and \eq{ai8eq5} with $I=\{1,\ldots,N\}$ will still hold.

Because of the degeneracy in \eq{ai10eq43} we have $\fF^{\{0,\ldots,N\}}=\fF^{\{1,\ldots,N\}}$ and $\fF^{0N}=\fF^{1N}$. Let us change $\m^{\{1,\ldots,N\}}$ by adding on a representative of the class $-[{\rm E}^{\{0,\ldots,N\}}]$. Since $[\m^{\{0,1,N\}}]$ and $(\fF^{01})^*$ are both the identity from above, by Proposition \ref{ai10prop3}(b), after this modification of $\m^{\{1,\ldots,N\}}$, we now have $[{\rm E}^{\{0,\ldots,N\}}]=0$. Proposition \ref{ai10prop2} with $\ti M=N$ and $\ti I=\{0,1,\ldots,N\}$ now implies that we can change the definition of $\m^{\{0,\ldots,N\}}$ (which was an arbitrary choice above) to make \eq{ai7eq19} and \eq{ai8eq5} hold with $I=\{0,1,\ldots,N\}$, for all $n\ge 0$ and $x_0,\ldots,x_n\in\sO^0$. We have now constructed all the data needed for the map $\ti f$ in \eq{ai10eq43}. This completes the proof of Theorem \ref{ai8thm3} for~$(\hAfCat^{\un,\iy}_R,\cT^\rfi)$.

\subsection{Proofs of Propositions \ref{ai6prop2} and \ref{ai7prop3}}
\label{ai106}

We prove Proposition \ref{ai7prop3}, which implies Proposition \ref{ai6prop2}. It is easy to see that $\bigl((\bar \sO^i)_{i=0}^{\bar N},(\bar  f^{ij})_{0\le i\le j\le\ti N},(\bar\sH{}^i)_{i=0}^{\bar N},\ab(\bar\m{}^{\bar I})_{\es\ne\bar I\subseteq\{0,\ldots,\bar N\}}\bigr)$ in Definition \ref{ai7def8}(b) satisfies all of Definition \ref{ai7def8}(a)(i)--(v), except perhaps equation \eq{ai7eq25}. So we must show that the $\bar\m{}^{\bar I;x_0,\ldots,x_k}_k$ satisfy \eq{ai7eq25} for all $\es\ne \bar I\subseteq\{0,\ldots,\bar N\}$, which we write as $\bar I=\{\bar\imath_0,\ab \bar\imath_1,\ab\ldots,\ab \bar\imath_{\bar M}\}$ for $\bar\imath_0<\bar\imath_1<\cdots<\bar\imath_{\bar M}$ with $\bar M\ge 0$, and $n\ge 0$, $x_0,\ldots,x_n\in\bar\sO^{\bar\imath_0}$. Write $\bar\d_{n,k;y_0,\ldots,y_k}^{\bar I;x_0,\ldots,x_n}$ for the operators $\d_{n,k;y_0,\ldots,y_k}^{I;x_0,\ldots,x_n}$ in Definition \ref{ai7def2} defined using $\bigl((\bar \sO^i),(\bar f^{ij}),(\bar\sH{}^i),\ab(\bar\m{}^{\bar I})\bigr)$. First we claim that
\e
\begin{gathered}
\bar\d_{n,k;y_0,\ldots,y_k}^{\bar I;x_0,\ldots,x_n}=0 \quad\text{if $\al\vert_{\bar I}$ is not injective, unless $k=n$ and $\bar I=\{\bar\imath_0,\bar\imath_1\}$} \\
\text{with $\bar\imath_0\!\ne\! \bar\imath_1$, $\al(\bar\imath_0)\!=\!\al(\bar\imath_1)\!=\!i$, and $y_j=x_j$ for $i=0,\ldots,n$,} \\
\text{and then $\bar\d_{n,n;x_0,\ldots,x_n}^{\{\bar\imath_0,\bar\imath_1\};x_0,\ldots,x_n}=\id_{\sH^i(x_0,x_1)}\ot\cdots \ot\id_{\sH^i(x_{n-1},x_n)}$.}
\end{gathered}
\label{ai10eq45}
\e

To see this, note that if $\bar M\ge 2$ and $\al(\bar\imath_j)=\al(\bar\imath_{j+1})$ then for each term $\bs\Ga_1,\ldots,\bs\Ga_k$ in \eq{ai7eq9}, by considering the least $i=1,\dots,k$ such that one of $\bar\imath_j$ or $\bar\imath_{j+1}$ does not lie in $R(\bs\Ga_i)$ we can show that $\bs\Ga_i$ has a vertex $v$ with $\bar\imath_j,\bar\imath_{j+1}\in I_v$ and $\md{I_v}\ge 2$, which forces $\bar\m_k^{\bar I;w_0,\ldots,w_k}=0$ by \eq{ai7eq26}, so $C(\bs\Ga_i)=0$. Thus $\bar\d_{n,k;y_0,\ldots,y_k}^{\bar I;x_0,\ldots,x_n}=0$, as this holds for all terms in \eq{ai7eq9}. If $\bar M=1$ then $\bar\d_{n,k;y_0,\ldots,y_k}^{\{\bar\imath_0,\bar\imath_1\};x_0,\ldots,x_n}$ is as in \eq{ai7eq11}. Since $\bar\m_{n_a}^{\{\bar\imath_0,\bar\imath_1\};w_0,\ldots,w_{n_a}}=0$ if $n_a>1$ and $\bar\m_{n_a}^{\{\bar\imath_0,\bar\imath_1\};w_0,w_1}=\id_{\sH^i(w_0,w_1)}$ if $n_a=1$ by \eq{ai7eq26}, equation \eq{ai10eq45} follows.

Now to prove \eq{ai7eq25} for $\es\ne \bar I\subseteq\{0,\ldots,\bar N\}$, divide into cases:
\begin{itemize}
\setlength{\itemsep}{0pt}
\setlength{\parsep}{0pt}
\item[(a)] $\al\vert_{\bar I}$ is injective.
\item[(b)] $\bar M=1$ and $\al(\bar\imath_0)=\al(\bar\imath_1)$.
\item[(c)] $\bar M>1$ and $\al(\bar\imath_0)=\al(\bar\imath_1)$.
\item[(d)] $\bar M>1$ and $\al(\bar\imath_{\bar M-1})=\al(\bar\imath_{\bar M})$.
\item[(e)] There exists $0<j<\bar M-1$ such that $\al(\bar\imath_j)=\al(\bar\imath_{j+1})$.
\end{itemize}
These exhaust all cases, as $\al:\{0,\ldots,\bar N\}\ra\{0,\ldots,N\}$ is nondecreasing, so can only fail to be injective on $\bar I$ if $\al(\bar\imath_j)=\al(\bar\imath_{j+1})$ for~$0\le j<\bar M$.

In case (a), by \eq{ai7eq26} we have $\bar\m_k^{\bar J;w_0,\ldots,w_l}=\m_k^{\al(\bar J);w_0,\ldots,w_l}$ for all $\es\ne \bar J\subset \bar I$, so equation \eq{ai7eq25} for $\bigl((\bar \sO^i),(\bar  f^{ij}),(\bar\sH{}^i),\ab(\bar\m{}^{\bar I})\bigr)$ and $\bar I$ is equivalent to \eq{ai7eq25} for $\bigl((\sO^i),(f^{ij}),(\sH^i),(\m^I)\bigr)$ and~$I=\al(\bar I)$.

In case (b), the right hand side of \eq{ai7eq25} for $\bigl((\bar \sO^i),(\bar  f^{ij}),(\bar\sH{}^i),\ab(\bar\m{}^{\bar I})\bigr)$ is zero as $\bar M=1$, and the left hand side becomes
\begin{align*}
&\bar\m_1^{\{\bar\imath_0,\bar\imath_1\};x_0,x_n}\ci\bar\d_{n,1;x_0,x_n}^{\{\bar\imath_0\};x_0,\ldots,x_n}-\bar\m_n^{\{\bar\imath_0\};x_0,\ldots,x_n}\ci\bar\d_{n,n;x_0,\ldots,x_n}^{\{\bar\imath_0,\bar\imath_1\};x_0,\ldots,x_n}\\
&=\id_{\sH^{\al(\bar\imath_0)}(x_0,x_n)}\ci\m_n^{\{\al(\bar\imath_0)\};x_0,\ldots,x_n}\\
&\qquad -\m_n^{\{\al(\bar\imath_0)\};x_0,\ldots,x_n}\ci\bigl(\id_{\sH^{\al(\bar\imath_0)}(x_0,x_1)}\ot\cdots \ot\id_{\sH^{\al(\bar\imath_0)}(x_{n-1},x_n)}\bigr)=0,
\end{align*}
using equations \eq{ai7eq10}, \eq{ai7eq26} and \eq{ai10eq45}.

In case (c), equation \eq{ai7eq25} for $\bigl((\bar \sO^i),(\bar  f^{ij}),(\bar\sH{}^i),\ab(\bar\m{}^{\bar I})\bigr)$ reduces to
\begin{equation*}
-\bar\m_n^{\{\bar\imath_1,\ldots,\bar\imath_{\bar M}\};x_0,\ldots,x_n}\ci\bar\d_{n,n;x_0,\ldots,x_n}^{\{\bar\imath_0,\bar\imath_1\};x_0,\ldots,x_n}=-\bar\m_n^{\{\bar\imath_0,\bar\imath_2,\ldots,\bar\imath_{\bar M}\};x_0,\ldots,x_n},\end{equation*}
as on the left hand side the only potentially nonzero term is $j=1$ and $k=n$ since $\bar\m_k^{\{\bar\imath_j,\ldots,\bar\imath_{\bar M}\};y_0,\ldots,y_k}=0$ for $j=0$ and $\bar\d_{k,n;y_0,\ldots,y_k}^{\{\bar\imath_0,\ldots,\bar\imath_j\};x_0,\ldots,x_n}=0$ if either $j>1$ or $j=1$ and $k<n$, and on the right hand side the only potentially nonzero term is $j=1$ as otherwise $\al$ is not injective on $I\sm\{\bar\imath_j\}$. But $\bar\d_{n,n;x_0,\ldots,x_n}^{\{\bar\imath_0,\bar\imath_1\};x_0,\ldots,x_n}=\id$ by \eq{ai10eq45} and $\bar\m_n^{\{\bar\imath_1,\ldots,\bar\imath_{\bar M}\};x_0,\ldots,x_n}=\bar\m_n^{\{\bar\imath_0,\bar\imath_2,\ldots,\bar\imath_{\bar M}\};x_0,\ldots,x_n}$ as $\al(\bar\imath_0)=\al(\bar\imath_1)$.

In case (d), similarly \eq{ai7eq25} for $\bigl((\bar \sO^i),(\bar  f^{ij}),(\bar\sH{}^i),\ab(\bar\m{}^{\bar I})\bigr)$ reduces to
\begin{align*}
&(-1)^{\bar M-1}\bar\m_1^{\{\bar\imath_{\bar M-1},\bar\imath_{\bar M}\};f^{\al(\{\bar\imath_0,\ldots,\bar\imath_{\bar M-1}\})}(x_0),
f^{\al(\bar\imath_0)\al(\bar\imath_{\bar M-1})}(x_n)}\ci{}\\
&\bar\d_{n,1;f^{\al(\{\bar\imath_0,\ldots,\bar\imath_{\bar M-1}\})}(x_0),f^{\al(\bar\imath_0)\al(\bar\imath_{\bar M-1})}(x_n)}^{\{\bar\imath_0,\ldots,\bar\imath_{\bar M-1}\};x_0,\ldots,x_n}=(-1)^{\bar M-1}\bar\m_n^{\{\bar\imath_0,\ldots,\bar\imath_{\bar M-2},\bar\imath_{\bar M}\};x_0,\ldots,x_n},
\end{align*}
as the only potentially nonzero terms on each side have $j=\bar M-1$ and $k=1$, and then $\bar\m_1^{\{i_{\bar M-1},\bar\imath_{\bar M}\};\ldots}=\id$ by \eq{ai7eq26}, and $\bar\d_{n,1;\ldots}^{\{\bar\imath_0,\ldots,\bar\imath_{\bar M-1}\};\ldots}=\bar\m_n^{\{\bar\imath_0,\ldots,\bar\imath_{\bar M-1}\};\ldots}$ by \eq{ai7eq10}, and $\bar\m_n^{\{\bar\imath_0,\ldots,\bar\imath_{\bar M-1}\};\ldots}=\bar\m_n^{\{\bar\imath_0,\ldots,\bar\imath_{\bar M-2},\bar\imath_{\bar M}\};\ldots}$ as~$\al(\bar\imath_{\bar M-1})=\al(\bar\imath_{\bar M})$.

In case (e), equation \eq{ai7eq25} for $\bigl((\bar \sO^i),(\bar  f^{ij}),(\bar\sH{}^i),\ab(\bar\m{}^{\bar I})\bigr)$ reduces to
\begin{equation*}
0=(-1)^j\bar\m_n^{\bar I\sm\{\bar\imath_j\};x_0,\ldots,x_n}+(-1)^{j+1}\bar\m_n^{\bar I\sm\{\bar\imath_{j+1}\};x_0,\ldots,x_n},
\end{equation*}
which holds as $\al(\bar I\sm\{\bar\imath_j\})=\al(\bar I\sm\{\bar\imath_{j+1}\})$.

\subsection{Proofs of Theorems \ref{ai6thm5}, \ref{ai7thm5}, and \ref{ai8thm4}}
\label{ai107}

We will prove Theorem \ref{ai8thm4}. Theorem \ref{ai7thm5} is a special case of this with $\Ga=\{0,\iy\}$, as in Example \ref{ai5ex2}. Theorem \ref{ai6thm5} follows from Theorem \ref{ai7thm5} by restricting to $A_\iy$-categories $\sA$ with $\sO=\{*\}$. The only differences in the proofs of Theorems \ref{ai7thm5} and \ref{ai8thm4} are in the ranges of variables in the sums in \eq{ai10eq46}, and that in the filtered case some sums are infinite and convergent, but in the non-filtered case there are only finitely many nonzero terms. The proof is similar to that of Theorem \ref{ai8thm1} in~\S\ref{ai103}.

Observe that when $k=0$ equation \eq{ai8eq12} follows from Definition \ref{ai7def2}(a), where both sides are zero except when $n=0$ and $M=2$, when both sides are $-\id_R$. Also when $k=1$ equation \eq{ai8eq12} follows from \eq{ai7eq10} and \eq{ai8eq11}.

Thus it remains to prove \eq{ai8eq12} when $M\ge 0$ and $k\ge 2$. We do this by a double induction on $M$ and $k$. For given $\ti M\ge 0$ and $\ti k\ge 2$, our inductive hypothesis is
\begin{itemize}
\setlength{\itemsep}{0pt}
\setlength{\parsep}{0pt}
\item[$(*)^{\ti M}_{\ti k}$] Equation \eq{ai8eq12} holds for $M\ge 0$ and $k,n\ge 0$ provided either $M<\ti M$, or $M=\ti M$ and $k\le\ti k$.
\end{itemize}
Let $\ti M\ge 0$ and $\ti k\ge 2$ be given. If $\ti k>2$ we suppose that $(*)^{\ti M}_{\ti k-1}$ holds. If $\ti k=2$ we suppose that $(*)^{\ti M-1}_k$ holds for all $k\ge 2$, unless $\ti M=0$, when we make no assumption. For any $\ti n\ge 0$, we will now prove \eq{ai8eq12} with $\ti M,\ti k,\ti n$ in place of $M,k,n$, in the long equation
\ea
&\sum_{\ti\jmath=0}^{\ti M}\sum_{\ti l\ge 0, \; w_0,\ldots,w_{\ti l}\in\sO^{i_{\ti\jmath}}}(-1)^{\ti\jmath}\d_{\ti l,\ti k;y_0,\ldots,y_{\ti k}}^{\{i_{\ti\jmath},\ldots,i_{\ti M}\};w_0,\ldots,w_{\ti l}}\ci\d_{\ti n,\ti l;w_0,\ldots,w_{\ti l}}^{\{i_0,\ldots,i_j\};x_0,\ldots,x_{\ti n}}
\label{ai10eq46}\\
\allowdisplaybreaks
&=\!\!\sum_{\ti\jmath=1}^{\ti M-1}(-1)^{\ti\jmath}\!\!\!\!\!\!\!\!\!\!\!\!\!\!
\begin{aligned}[t]&\sum_{\begin{subarray}{l} \ti l\ge 0, \; w_0,\ldots,w_{\ti l}\in\sO^{i_{\ti\jmath}} \\ 
0<m'\le \ti M-\ti\jmath, \\
l',l''\ge 0,\; \ti l=l'+l'', \\
\ti\jmath=j'_0<j'_1<\cdots<j'_{m'}=\ti M, \\
k'_0,\ldots,k'_{m'}\ge 0, \; k'_0=l', \\ k'_{m'}=\ti k-1, \\
w_0^{\prime a},\ldots,w_{k'_a}^{\prime a}\in\sO^{i_{j'_a}},\\ 0<a<m'
\end{subarray}}\Bigl(\begin{aligned}[t]&
(-1)^{\ti k-1+m'+\ti M-\ti\jmath+l'(1+m'+1)+\sum_{a=1}^{m'-1}(j'_a-\ti\jmath-a)}\\
&\bigl(\d_{k'_{m'-1},k'_{m'};y_0,\ldots,y_{\ti k-1}}^{\{i_{j'_{m'-1}},\ldots,i_{j'_{m'}}\};w^{\prime m-1}_0,\ldots,w^{\prime m-1}_{k'_{m-1}}}\ci{}\\
&\d_{k'_{m'-2},k'_{m'-1};w^{\prime m'-1}_0,\ldots,w^{\prime m'-1}_{k'_{m'-1}}}^{\{i_{j'_{m'-2}},\ldots,i_{j'_{m'-1}}\};w^{\prime m'-2}_0,\ldots,w^{\prime m'-2}_{k'_{m'-2}}}\ci\cdots\ci{}\\
&\d_{k'_0,k'_1;w^{\prime 1}_0,\ldots,w^{\prime 1}_{k'_1}}^{\{i_{j'_0},\ldots,i_{j'_1}\};w_0,\ldots,w_{l'}}\bigr)\ot\d_{l'',1;y_{\ti k-1},y_{\ti k}}^{\{i_{j'_0},i_{j'_1},\ldots,i_{j'_{m'}}\};w_{l'},\ldots,w_{\ti l}}\Bigr)\end{aligned}\\
&\!\!\!\!\!\!\!\!\!\!\!\!{}\ci\biggl(\sum_{\begin{subarray}{l} 
0<m''\le\jmath, \\ n',n''\ge 0, \; \ti n=n'+n'', \\
0=j''_0<j''_1<\cdots<j''_{m''}=\ti\jmath, \\
k''_0,\ldots,k''_{m''}\ge 0, \; k''_0=n', \\ k''_{m''}=l', \\
w_0^{\prime\prime a},\ldots,w_{k''_a}^{\prime\prime a}\in\sO^{i_{j''_a}},\\ 0<a<m''
\end{subarray}}\begin{aligned}[t]&
(-1)^{l''(l'+m''+\ti\jmath)+n'(l''+m''+1)+\sum_{a=1}^{m''-1}(j''_a-a)}\\
&\bigl(\d_{k''_{m''-1},k''_{m''};w_0,\ldots,w_{l'}}^{\{i_{j''_{m''-1}},\ldots,i_{j''_{m''}}\};w^{\prime\prime m''-1}_0,\ldots,w^{\prime\prime m''-1}_{k''_{m''-1}}}\ci{}\\
&\d_{k''_{m''-2},k''_{m''-1};w^{\prime\prime m''-1}_0,\ldots,w^{\prime\prime m''-1}_{k''_{m''-1}}}^{\{i_{j''_{m''-2}},\ldots,i_{j''_{m''-1}}\};w^{\prime\prime m''-2}_0,\ldots,w^{\prime\prime m''-2}_{k''_{m''-2}}}\ci\cdots\ci{}\\
&\d_{k''_0,k''_1;w^{\prime\prime 1}_0,\ldots,w^{\prime\prime 1}_{k''_1}}^{\{i_{j''_0},\ldots,i_{j''_1}\};x_0,\ldots,x_{n'}}\bigr)\ot\d_{n'',l'';w_{l'},\ldots,w_{\ti l}}^{\{i_{j''_0},i_{j''_1},i_{j''_2},\ldots,i_{j''_{m''}}\};x_{n'},\ldots,x_{n''}}\biggr)\end{aligned}
\end{aligned}
\allowdisplaybreaks
\nonumber\\
&+\begin{aligned}[t]&\sum_{\begin{subarray}{l} 
\ti l\ge 0, \; w_0,\ldots,w_{\ti l}\in\sO^{i_0} \\
0<m\le \ti M, \;
l',l''\ge 0, \; \ti l=l'+l'', \\
0=j_0<j_1<\cdots<j_m=\ti M, \\
k_0,\ldots,k_m\ge 0, \; k_0=l', \; k_m=\ti k-1, \\
w_0^a,\ldots,w_{k_a}^a\in\sO^{i_{j_a}},\; 0<a<m
\end{subarray}}\Bigl(\begin{aligned}[t]&
(-1)^{\ti k-1+m+\ti M+l'm+\sum_{a=1}^{m-1}(j_a-a)}\\
&\bigl(\d_{k_{m-1},k_m;y_1,\ldots,y_{\ti k-1}}^{\{i_{j_{m-1}},\ldots,i_{j_m}\};w^{m-1}_0,\ldots,w^{m-1}_{k_{m-1}}}\ci{}\\
&\d_{k_{m-2},k_{m-1};w^{m-1}_0,\ldots,w^{m-1}_{k_{m-1}}}^{\{i_{j_{m-2}},\ldots,i_{j_{m-1}}\};w^{m-2}_0,\ldots,w^{m-2}_{k_{m-2}}}\ci\cdots\ci{}\\
&\d_{k_0,k_1;w^1_0,\ldots,w^1_{k_1}}^{\{i_{j_0},\ldots,i_{j_1}\};w_0,\ldots,w_{l'}}\bigr)\ot\d_{l'',1;y_{\ti k-1},y_{\ti k}}^{\{i_{j_0},i_{j_1},i_{j_2},\ldots,i_{j_m}\};w_{l'},\ldots,w_{\ti l}}\Bigr)\end{aligned}\\
&\ci\Bigl((-1)^{l''(l'+\ti n)}\de_{\begin{subarray}{l} x_{i+\ti n-\ti l}=w_i,\\ l'\le i\le \ti l\end{subarray}}\d_{\ti n-l'',l';w_0,\ldots,w_{l'}}^{\{i_0\};x_0,\ldots,x_{\ti n-l''}}\ot\id^{\ot^{l''}}\\
&\qquad +(-1)^{l'}\de_{\begin{subarray}{l} x_i=w_i,\\ 0\le i\le l'\end{subarray}}\id^{\ot^{l'}}\ot\d_{\ti n-l',l'';w_{l'},\ldots,w_{\ti l}}^{\{i_0\};x_{l'},\ldots,x_{\ti n}}
\Bigr)\end{aligned}
\allowdisplaybreaks
\nonumber\\
&+\sum_{\begin{subarray}{l} \ti l\ge 0, \\ w_0,\ldots,w_{\ti l}\in\sO^{i_{\ti M}}\end{subarray}}\begin{aligned}[t]&
(-1)^{\ti M}\Bigl((-1)^{\ti k-1+\ti l}\de_{w_{\ti l-i}=y_{\ti k-i},\; i=0,1}\d_{\ti l-1,\ti k-1;y_0,\ldots,y_{\ti k-1}}^{\{i_{\ti M}\};w_0,\ldots,w_{\ti l-1}}\ot\id\Bigr)\\
&{}\ci\biggl(\sum_{\begin{subarray}{l} 
0<m\le \ti M, \\ n',n''\ge 0, \; \ti n=n'+n'', \\
0=j_0<j_1<\cdots<j_m=\ti M, \\
k_0,\ldots,k_m\ge 0, \; k_0=n', \\ k_m=\ti l-1, \\
w_0^a,\ldots,w_{k_a}^a\in\sO^{i_{j_a}},\\ 0<a<m
\end{subarray}}\begin{aligned}[t]&
(-1)^{\ti l-1+m+\ti M+n'm+\sum_{a=1}^{m-1}(j_a-a)}\\
&\bigl(\d_{k_{m-1},k_{m};w_0,\ldots,w_{\ti l-1}}^{\{i_{j_{m-1}},\ldots,i_{j_m}\};w^{m-1}_0,\ldots,w^{m-1}_{k_{m-1}}}\ci{}\\
&\d_{k_{m-2},k_{m-1};w^{m-1}_0,\ldots,w^{m-1}_{k_{m-1}}}^{\{i_{j_{m-2}},\ldots,i_{j_{m-1}}\};w^{m-2}_0,\ldots,w^{m-2}_{k_{m-2}}}\ci\cdots\ci{}\\
&\d_{k_0,k_1;w^1_0,\ldots,w^1_{k_1}}^{\{i_{j_0},\ldots,i_{j_1}\};x_0,\ldots,x_{n'}}\bigr)\\
&\ot\d_{n'',1;w_{\ti l-1},w_{\ti l}}^{\{i_{j_0},i_{j_1},i_{j_2},\ldots,i_{j_m}\};x_{n'},\ldots,x_{\ti n}}\biggr)\end{aligned}
\end{aligned}
\allowdisplaybreaks
\nonumber\\
&+\sum_{\begin{subarray}{l} \ti l\ge 0, \\ w_0,\ldots,w_{\ti l}\in\sO^{i_{\ti M}}\end{subarray}}\begin{aligned}[t]&
(-1)^{\ti M}\Bigl((-1)^{\ti k-1}\de_{\begin{subarray}{l} w_i=y_i, \\ 0\le i\le\ti k-1\end{subarray}}
\id^{\ot^{\ti k-1}}\ot\d_{\ti l-\ti k+1,1;y_{\ti k-1},y_{\ti k}}^{\{i_{\ti M}\};w_{\ti k-1},\ldots,w_{\ti l}}\Bigr)\\
&{}\ci\biggl(\sum_{\begin{subarray}{l} 
0<m\le \ti M, \\ n',n''\ge 0, \; \ti n=n'+n'', \\
0=j_0<j_1<\cdots<j_m=\ti M, \\
k_0,\ldots,k_m\ge 0, \; k_0=n', \\ k_m=\ti k-1, \\
w_0^a,\ldots,w_{k_a}^a\in\sO^{i_{j_a}},\\ 0<a<m
\end{subarray}}\begin{aligned}[t]&
(-1)^{\begin{subarray}{l}(\ti l-\ti k+1)(\ti k-1+m+\ti M) \\ +n'(\ti l-\ti k+m)+\sum_{a=1}^{m-1}(j_a-a)\end{subarray}}\\
&\bigl(\d_{k_{m-1},k_{m};w_0,\ldots,w_{\ti k-1}}^{\{i_{j_{m-1}},\ldots,i_{j_m}\};w^{m-1}_0,\ldots,w^{m-1}_{k_{m-1}}}\ci{}\\
&\d_{k_{m-2},k_{m-1};w^{m-1}_0,\ldots,w^{m-1}_{k_{m-1}}}^{\{i_{j_{m-2}},\ldots,i_{j_{m-1}}\};w^{m-2}_0,\ldots,w^{m-2}_{k_{m-2}}}\ci\cdots\ci{}\\
&\d_{k_0,k_1;w^1_0,\ldots,w^1_{k_1}}^{\{i_{j_0},\ldots,i_{j_1}\};x_0,\ldots,x_{n'}}\bigr)\\
&\ot\d_{n'',\ti l-\ti k+1;w_{\ti k-1},\ldots,w_{\ti l}}^{\{i_{j_0},i_{j_1},i_{j_2},\ldots,i_{j_m}\};x_{n'},\ldots,x_{\ti n}}\biggr)\end{aligned}
\end{aligned}
\allowdisplaybreaks
\nonumber\\
&=\!\!\!\sum_{\begin{subarray}{l} 
0<l\le \ti M, \\ n',n''\ge 0, \; n=n'+n'', \\
0=j_0<j_1<\cdots<j_l=\ti M, \\
k_0,\ldots,k_l\ge 0, \; k_0=n', \\ k_l=\ti k-1, \\
0\le m\le l,\; l''\ge 1 \\
w_0^a,\ldots,w_{k_a}^a\in\sO^{i_{j_a}}, \\ 0<a<l, \\
z_0,\ldots,z_{l''}\in\sO^{i_{j_m}}
\end{subarray}}
\begin{aligned}[t]&
(-1)^{j_m}(-1)^{\begin{subarray}{l} \ti k-1+l-m+\ti M-j_m+k_m(l-m) \\ +\sum_{a=m+1}^{l-1}(j_a-j_m-a+m)\end{subarray}}\\
&(-1)^{l''(k_m+m+j_m)+n'(l''+m+1)+\sum_{a=1}^{m-1}(j_a-a)}\\
& (-1)^{(l+l''+m)(j_m+k_m+m+n')}
\\
&\bigl(\d_{k_{l-1},k_l;y_0,\ldots,y_{\ti k-1}}^{\{i_{j_{l-1}},\ldots,i_{j_l}\};w^{l-1}_0,\ldots,w^{l-1}_{k_{l-1}}}\ci{}\\
&\d_{k_{l-2},k_{l-1};w^{l-1}_0,\ldots,w^{l-1}_{k_{l-1}}}^{\{i_{j_{l-2}},\ldots,i_{j_{l-1}}\};w^{l-2}_0,\ldots,w^{l-2}_{k_{l-2}}}\ci\cdots\ci{}\\
&\d_{k_0,k_1;w^1_0,\ldots,w^1_{k_1}}^{\{i_{j_0},\ldots,i_{j_1}\};x_0,\ldots,x_{n'}}\bigr)\\
&{}\ot\bigl(\d_{l'',1;y_{\ti k-1},y_{\ti k}}^{\{i_{j_m},i_{j_{m+1}},\ldots,i_{j_l}\};z_0,\ldots,z_{l''}}\ci\d_{n'',l'';z_0,\ldots,z_{l''}}^{\{i_{j_0},i_{j_1},i_{j_2},\ldots,i_{j_m}\};x_{n'},\ldots,x_{\ti n}}\bigr)
\end{aligned}
\allowdisplaybreaks
\nonumber\\
&+\sum_{\begin{subarray}{l} \ti l\ge 0, \;
0<m\le \ti M, \\ l',l''\ge 0, \; \ti l=l'+l'', \\
0=j_0<j_1<\cdots<j_m=\ti M, \\
k_0,\ldots,k_m\ge 0, \; k_0=l', \\ k_m=\ti k-1, \\
w_0^a,\ldots,w_{k_a}^a\in\sO^{i_{j_a}},\\ 0\le a<m \end{subarray}}
\begin{aligned}[t]&
(-1)^{\ti k-1+m+\ti M+l'm+\sum_{a=1}^{m-1}(j_a-a)}\\
&{}\cdot(-1)^{l''(l'+\ti n)}\cdot (-1)^{(\ti l+\ti n+1)(l''+m')}
\\
&{}\cdot\bigl(\d_{k_{m-1},k_m;y_1,\ldots,y_{\ti k-1}}^{\{i_{j_{m-1}},\ldots,i_{j_m}\};w^{m-1}_0,\ldots,w^{m-1}_{k_{m-1}}}\ci{}\\
&\d_{k_{m-2},k_{m-1};w^{m-1}_0,\ldots,w^{m-1}_{k_{m-1}}}^{\{i_{j_{m-2}},\ldots,i_{j_{m-1}}\};w^{m-2}_0,\ldots,w^{m-2}_{k_{m-2}}}\ci\cdots\ci{}\\
&\d_{k_0,k_1;w^1_0,\ldots,w^1_{k_1}}^{\{i_{j_0},\ldots,i_{j_1}\};w^0_0,\ldots,w^0_{k_0}}\ci\d_{\ti n-l'',l';w^0_0,\ldots,w^0_{k_0}}^{\{i_0\};x_0,\ldots,x_{\ti n-l''}}\bigr)\\
&\ot\d_{l'',1;y_{\ti k-1},y_{\ti k}}^{\{i_{j_0},i_{j_1},\ldots,i_{j_m}\};x_{\ti n-l''},\ldots,x_{\ti n}}\end{aligned}
\allowdisplaybreaks
\nonumber\\
&+\sum_{\begin{subarray}{l} 
\ti l\ge 0, \; 0<m\le \ti M, \; n',n''\ge 0, \\ \ti n=n'+n'', \\
0=j_0<j_1<\cdots<j_m=\ti M, \\
k_0,\ldots,k_m\ge 0, \; k_0=n', \; k_m=\ti l-1, \\
w_0^a,\ldots,w_{k_a}^a\in\sO^{i_{j_a}},\\ 0<a\le m
\end{subarray}}\begin{aligned}[t]&
(-1)^{\ti k+m+n'm+\sum_{a=1}^{m-1}(j_a-a)}\\
&{}\cdot\bigl(\d_{\ti l-1,\ti k-1;y_0,\ldots,y_{\ti k-1}}^{\{i_{\ti M}\};w^m_0,\ldots,w^m_{k_m}}\ci\d_{k_{m-1},k_{m};w^m_0,\ldots,w^m_{k_m}}^{\{i_{j_{m-1}},\ldots,i_{j_m}\};w^{m-1}_0,\ldots,w^{m-1}_{k_{m-1}}}\\
&{}\ci\d_{k_{m-2},k_{m-1};w^{m-1}_0,\ldots,w^{m-1}_{k_{m-1}}}^{\{i_{j_{m-2}},\ldots,i_{j_{m-1}}\};w^{m-2}_0,\ldots,w^{m-2}_{k_{m-2}}}\\
&{}\ci\cdots\ci\d_{k_0,k_1;w^1_0,\ldots,w^1_{k_1}}^{\{i_{j_0},\ldots,i_{j_1}\}}\bigr)\ot\d_{n'',1;y_{\ti k-1},y_{\ti k}}^{\{i_{j_0},i_{j_1},\ldots,i_{j_m}\};x_{n'},\ldots,x_{\ti n}}\end{aligned}
\allowdisplaybreaks
\nonumber\\
&=\!\!\!\sum_{\begin{subarray}{l} 
0<l\le \ti M, \; n',n''\ge 0, \\ n=n'+n'', \\
0=j_0<j_1<\cdots<j_l=\ti M, \\
k_0,\ldots,k_l\ge 0, \; k_0=n', \\ k_l=\ti k-1, \\
0<m<l \\
w_0^a,\ldots,w_{k_a}^a\in\sO^{i_{j_a}},\\ 0<a<l
\end{subarray}}
\begin{aligned}[t]&
(-1)^{\ti k+(l+1)(n'+1)+\ti M+\sum_{a=1}^{l-1}(j_a-a)}\cdot(-1)^m\\
&\bigl(\d_{k_{l-1},k_l;y_0,\ldots,y_{\ti k-1}}^{\{i_{j_{l-1}},\ldots,i_{j_l}\};w^{l-1}_0,\ldots,w^{l-1}_{k_{l-1}}}\ci\d_{k_{l-2},k_{l-1};w^{l-1}_0,\ldots,w^{l-1}_{k_{l-1}}}^{\{i_{j_{l-2}},\ldots,i_{j_{l-1}}\};w^{l-2}_0,\ldots,w^{l-2}_{k_{l-2}}}\\
&{}\ci\cdots\ci\d_{k_0,k_1;w^1_0,\ldots,w^1_{k_1}}^{\{i_{j_0},\ldots,i_{j_1}\};x_0,\ldots,x_{n'}}\bigr)\\
&{}\ot\d_{n'',1;y_{\ti k-1},y_{\ti k}}^{\{i_{j_0},i_{j_1},\ldots,i_{j_{m-1}},i_{j_{m+1}},\ldots,i_{j_l}\};x_{n'},\ldots,x_{\ti n}}
\end{aligned}
\nonumber\\
&+(\text{second and third terms from previous step})
\allowdisplaybreaks
\nonumber\\
&=\!\!\!\sum_{\begin{subarray}{l} 
0<\hat l< \ti M, \\ n',n''\ge 0, \; \ti n=n'+n'', \\
0=\hat\jmath_0<\hat\jmath_1<\cdots<\hat\jmath_{\hat l}=\ti M, \\
\hat k_0,\ldots,\hat k_{\hat l}\ge 0, \; \hat k_0=n', \\ \hat k_{\hat l}=\ti k-1, \\
0<m\le\hat l,\; \hat\jmath_{m-1}\le\ti\jmath\le\hat\jmath_m \\
w_0^a,\ldots,w_{k_a}^a\in\sO^{i_{\hat\jmath_a}},\\ 0<a<\hat l, \\
\ti l\ge 0,\; z_0,\ldots,z_{\ti l}\in\sO^{i_{\ti\jmath}} \end{subarray}}
\begin{aligned}[t]&
(-1)^{\ti k+\hat l(n'+1)+\ti M+\sum_{a=1}^{\hat l-1}(\hat\jmath_a-a)}\cdot(-1)^{\hat\jmath_{m-1}+m}\\
&\Bigl(\d_{\hat k_{\hat l-1},\hat k_{\hat l};y_0,\ldots,y_{\ti k-1}}^{\{i_{\hat\jmath_{\hat l-1}},\ldots,i_{\hat\jmath_{\hat l}}\};w^{\hat l-1}_0,\ldots,w^{\hat l-1}_{\hat k_{\hat l-1}}}\ci\cdots\ci\d_{\hat k_m,\hat k_{m+1};w^{m+1}_0,\ldots,w^{m+1}_{\hat k_{m+1}}}^{\{i_{\hat\jmath_m},\ldots,i_{\hat\jmath_{m+1}}\};w^m_0,\ldots,w^m_{\hat k_m}}\ci{}\\
&\Bigl[(-1)^{\ti\jmath-\hat\jmath_{m-1}}\d_{\ti l,\hat k_m;w^m_0,\ldots,w^m_{\hat k_m}}^{\{i_{\ti\jmath},\ldots,i_{\hat\jmath_m}\};z_0,\ldots,z_{\ti l}}\ci\d_{\hat k_{m-1},\ti l;z_0,\ldots,z_{\ti l}}^{\{i_{\hat\jmath_{m-1}},\ldots,i_{\ti\jmath}\};w^{m-1}_0,\ldots,w^{m-1}_{\hat k_{m-1}}}\Bigr]\\
&{}\ci\d_{\hat k_{m-2},\hat k_{m-1};w^{m-1}_0,\ldots,w^{m-1}_{\hat k_{m-1}}}^{\{i_{\hat\jmath_{m-2}},\ldots,i_{\hat\jmath_{m-1}}\};w^{m-2}_0,\ldots,w^{m-2}_{\hat k_{m-2}}}\ci\cdots\ci\d_{\hat k_0,\hat k_1;w^1_0,\ldots,w^1_{\hat k_1}}^{\{i_{\hat\jmath_0},\ldots,i_{\hat\jmath_1}\};x_0,\ldots,x_{n'}}\Bigr)\\
&\ot\d_{n'',1;y_{\ti k-1},y_{\ti k}}^{\{i_{\hat\jmath_0},i_{\hat\jmath_1},\ldots,i_{\hat\jmath_{\hat l}}\};x_{n'},\ldots,x_{\ti n}}
\end{aligned}
\allowdisplaybreaks
\nonumber\\
&=\!\!\!\sum_{\begin{subarray}{l} 
0<\hat l< \ti M, \\ n',n''\ge 0, \; \ti n=n'+n'', \\
0=\hat\jmath_0<\hat\jmath_1<\cdots<\hat\jmath_{\hat l}=\ti M, \\
\hat k_0,\ldots,\hat k_{\hat l}\ge 0, \; \hat k_0=n', \\ \hat k_{\hat l}=\ti k-1, \\
0<m\le\hat l,\; \hat\jmath_{m-1}<\ti\jmath<\hat\jmath_m, \\
w_0^a,\ldots,w_{k_a}^a\in\sO^{i_{\hat\jmath_a}},\\ 0<a<\hat l \end{subarray}}
\begin{aligned}[t]&
(-1)^{\ti k+\hat l(n'+1)+\ti M+\sum_{a=1}^{\hat l-1}(\hat\jmath_a-a)}\cdot(-1)^{\hat\jmath_{m-1}+m}\\
&\Bigl(\d_{\hat k_{\hat l-1},\hat k_{\hat l};y_0,\ldots,y_{\ti k-1}}^{\{i_{\hat\jmath_{\hat l-1}},\ldots,i_{\hat\jmath_{\hat l}}\};w^{\hat l-1}_0,\ldots,w^{\hat l-1}_{\hat k_{\hat l-1}}}\ci\cdots\ci\d_{\hat k_m,\hat k_{m+1};w^{m+1}_0,\ldots,w^{m+1}_{\hat k_{m+1}}}^{\{i_{\hat\jmath_m},\ldots,i_{\hat\jmath_{m+1}}\};w^m_0,\ldots,w^m_{\hat k_m}}\\
&\ci\Bigl[(-1)^{\ti\jmath-\hat\jmath_{m-1}}\d_{\hat k_{m-1},\hat k_m;w^m_0,\ldots,w^m_{\hat k_m}}^{\{i_{\hat\jmath_{m-1}},\ldots,i_{\ti\jmath-1},i_{\ti\jmath+1},\ldots,i_{\hat\jmath_m}\};w^{m-1}_0,\ldots,w^{m-1}_{\hat k_{m-1}}}\Bigr]\ci{}\\
&\d_{\hat k_{m-2},\hat k_{m-1};w^{m-1}_0,\ldots,w^{m-1}_{\hat k_{m-1}}}^{\{i_{\hat\jmath_{m-2}},\ldots,i_{\hat\jmath_{m-1}}\};w^{m-2}_0,\ldots,w^{m-2}_{\hat k_{m-2}}}\ci\cdots\ci\d_{\hat k_0,\hat k_1;w^1_0,\ldots,w^1_{\hat k_1}}^{\{i_{\hat\jmath_0},\ldots,i_{\hat\jmath_1}\};x_0,\ldots,x_{n'}}\Bigr)\\
&\ot\d_{n'',1;y_{\ti k-1},y_{\ti k}}^{\{i_{\hat\jmath_0},i_{\hat\jmath_1},\ldots,i_{\hat\jmath_{\hat l}}\};x_{n'},\ldots,x_{\ti n}}
\end{aligned}
\nonumber\\
&=\sum_{\ti\jmath=1}^{\ti M-1}(-1)^{\ti\jmath}\,\d_{\ti n,\ti k;y_0,\ldots,y_{\ti k}}^{\ti I\sm\{i_{\ti\jmath}\};x_0,\ldots,x_{\ti n}}.
\nonumber
\ea
Here in the first step we apply \eq{ai7eq12}--\eq{ai7eq13} twice, first to $\d_{\ti l,\ti k;y_0,\ldots,y_{\ti k}}^{\{i_{\ti\jmath},\ldots,i_{\ti M}\};w_0,\ldots,w_{\ti l}}$ with $k'=\ti k-1$ and $k''=1$, which yields a sum over splittings $\ti l=l'+l''$, and then to $\d_{\ti n,\ti l;w_0,\ldots,w_{\ti l}}^{\{i_0,\ldots,i_j\};x_0,\ldots,x_{\ti n}}$ splitting $\ti l=l'+l''$. The first term, with $\sum_{\ti\jmath=1}^{\ti M-1}\cdots$, uses \eq{ai7eq13} twice. The second term is for $\ti\jmath=0$, and uses \eq{ai7eq13} first and then \eq{ai7eq12}. The third and fourth terms are for $\ti\jmath=\ti M$, and use \eq{ai7eq12} first and then~\eq{ai7eq13}. 

In the second step we combine the sums in the first term into one with the notation change
\begin{align*}
l=m'+m'', \quad m&=m'', \quad l-m=m',\quad \quad j_m=\ti\jmath, \quad k_m=l'\\
(j_0,\ldots,j_l)&=(j_0'',j_1'',\ldots,j_{m''}''=\ti\jmath=j_0',j_1',\ldots,j'_{m'}),\\
(k_0,\ldots,k_l)&=(k_0'',k_1'',\ldots,k_{m''}''=l'=k_0',k_1',\ldots,k'_{m'}).
\end{align*}
The extra sign $(-1)^{(l+l''+m)(j_m+k_m+m+n')}$ comes from exchanging the order of $\ot$ and $\ci$ using \eq{ai3eq2}. The first term in the first step corresponds to the terms with $0<m<l$ in the first term of the second step. The second half of the second term in the first step corresponds to the terms with $m=0$ in the first term of the second step, with the notation change $n'=l'$ and $n''=\ti n-l'$. The fourth term in the first step corresponds to the terms with $m=l$ in the first term of the second step, with the notation change $l''=\ti l-\ti k+1$.

The second and third terms in the second step correspond to the first half of the second term and the third term in the first step, with extra signs coming from exchanging the order of $\ot$ and $\ci$ using \eq{ai3eq2}.

In the third step of \eq{ai10eq46} we first note that the product of signs in the first term of the second step may be simplified to
\begin{equation*}
(-1)^{\ti k+(l+1)(n'+1)+\ti M+\sum_{a=1}^{l-1}(j_a-a)}\cdot(-1)^m.
\end{equation*}
Then we rewrite the sum $\sum_{m,l''}(-1)^m\d_{l'',1;\ldots}^{\{i_{j_m},i_{j_{m+1}},\ldots,i_{j_l}\};\ldots}\ci\d_{n'',l'';\ldots}^{\{i_{j_0},i_{j_1},i_{j_2},\ldots,i_{j_m}\};\ldots}$ using \eq{ai8eq12}, which holds as above since~$k=1$.

In the fourth step we rewrite the first sum in the third step using the notation
\begin{align*}
\hat l&=l-1,\quad \ti\jmath=j_m, \quad \ti l=k_m,
\\
(\hat\jmath_0,\ldots,\hat\jmath_{\hat l})&=(j_0,j_1,\ldots,j_{m-1},j_{m+1},\ldots,j_l),\\
(\hat k_0,\ldots,\hat k_{\hat l})&=(k_0,k_1,\ldots,k_{m-1},k_{m+1},\ldots,k_l).
\end{align*}
Observe that the condition $j_{m-1}<j_m<j_{m+1}$ in the third step should transform to $\hat\jmath_{m-1}<\ti\jmath<\hat\jmath_m$ in the sum in the fourth step, but we have instead written $\hat\jmath_{m-1}\le\ti\jmath\le\hat\jmath_m$. This is valid because for $m=1,\ldots,\hat l-1$, the term in the sum with $\ti\jmath=\hat\jmath_m$, which includes a factor $\d_{\ti l,\hat k_m;\ldots}^{\{i_{\hat\jmath_m}\};\ldots}$, cancels with the term in the sum with $m+1$ in place of $m$ and $\ti\jmath=\hat\jmath_{(m+1)-1}$, which includes a factor $\d_{\hat k_{(m+1)-1},\ti l;\ldots}^{\{i_{\hat\jmath_m}\};\ldots}$, and we exchange the values of $\hat k_m,\ti l$. This leaves two terms uncancelled, the case $m=1$ and $0=\hat\jmath_0=\ti\jmath$, and the case $m=\hat l$ and $\ti\jmath=\hat\jmath_{\hat l}=\ti M$. These cancel with the second and third terms in the third step, respectively, which is why these have disappeared in the fourth step.

The fifth step rewrites $\sum_{\ti\jmath=\hat\jmath_{m-1}}^{\hat\jmath_m}\bigl[(-1)^{\ti\jmath-\hat\jmath_{m-1}}\d_{\ti l,\hat k_m;\ldots}^{\{i_{\ti\jmath},\ldots,i_{\hat\jmath_m}\};\ldots}\ci\d_{\hat k_{m-1},\ti l;\ldots}^{\{i_{\hat\jmath_{m-1}},\ldots,i_{\ti\jmath}\};\ldots}\bigr]$ in the fourth step using \eq{ai8eq12}. This is valid by induction as either $\hat\jmath_m-\hat\jmath_{m-1}<\ti M$, or $m=\hat l=1$ and $\hat k_m=\ti k-1<\ti k$. The final step uses \eq{ai7eq13} with $k'=\ti k-1$, $k''=1$. To match up the signs, note that $\hat\jmath_0,\ldots,\hat\jmath_{m-1},\hat\jmath_{m+1},\ldots,\hat\jmath_{\hat l}$ should be renumbered as $j_0,j_1,\ldots,j_{m-1},j_m,\ldots,j_{\hat l-1}$ in \eq{ai7eq13}, so the signs $(-1)^{\sum_{a=1}^{\hat l-1}(\hat\jmath_a-a)}$ in the fifth step and $(-1)^{\sum_{a=1}^{\hat l-2}(j_a-a)}$ in \eq{ai7eq13} differ by $(-1)^{\hat\jmath_m-m}\cdot(-1)^{\hat l-m}$. This proves \eq{ai10eq46}, and Theorem \ref{ai8thm4} follows by induction.

\subsection{Proofs of Theorems \ref{ai6thm6}, \ref{ai7thm6}, and \ref{ai8thm5}}
\label{ai108}

We will prove Theorem \ref{ai8thm5} for $\AfCat^{\un,\iy}_{R,K,L}$, where $K,L\in\{0,1,2,\ldots,\iy\}$ with $K\le L$ and $L\ge 3$. The proof for $\AwCat^{\un,\iy}_{R,K,L}$ is the same, as the different filtration conditions play no part in the proof. Theorem \ref{ai7thm6} is a special case of Theorem \ref{ai8thm5} with $K=L$ and $\Ga=\{0,\iy\}$, as in Example \ref{ai5ex2}. Theorem \ref{ai6thm6} follows from Theorem \ref{ai7thm6} by restricting to $A_\iy$-categories $\sA$ with $\sO=\{*\}$. 

Compare the definitions of an $N$-simplex $\bigl((\sO^i)_{i=0}^N,(f^{ij})_{0\le i\le j\le N},(\sH^i)_{i=0}^N,\ab(\m^I)_{\es\ne I\subseteq\{0,\ldots,N\}}\bigr)$ in $\AfCat^{\un,\iy}_{R,K,L}$ from Definitions \ref{ai8def1} and \ref{ai8def4}, and an $N$-simplex in $\dAfCat^{\un,\iy}_{R,K,L}$ from Definition \ref{ai8def6}. There are only two differences:
\begin{itemize}
\setlength{\itemsep}{0pt}
\setlength{\parsep}{0pt}
\item[(a)] Suppose $I=\{i_0<i_1<i_2\}\subseteq\{0,\ldots,N\}$ and $x_0\in\sO^{i_0}$. Then in $\AfCat^{\un,\iy}_{R,K,L}$ we have $\m_0^{I;x_0}(1)=\boo^{i_2}_{f^{i_0i_2}(x_0)}$ if $K\ge 3$ and $\m_0^{I;x_0}(1)-\boo^{i_2}_{f^{i_0i_2}(x_0)}\in F^{\ssu 0}\Hom_{\sA^{i_2}}(f^{i_0i_2}(x_0),f^{i_0i_2}(x_0))$ if~$K\le 2$.

However, in $\dAfCat^{\un,\iy}_{R,K,L}$ we have $\m_0^{I;x_0}(1)=0$ if $K\ge 3$ and $\m_0^{I;x_0}(1)\in F^{\ssu 0}\Hom_{\sA^{i_2}}(f^{i_0i_2}(x_0),f^{i_0i_2}(x_0))$ if~$K\le 2$.
\item[(b)] In $\AfCat^{\un,\iy}_{R,K,L}$ the $\m_k^{J;w_0,\ldots,w_k}$ must satisfy equations \eq{ai8eq4}--\eq{ai8eq5} for all $\es\ne I\subseteq\{0,\ldots,N\}$ and $n\ge 0$, but in $\dAfCat^{\un,\iy}_{R,K,L}$ the $\m_k^{J;w_0,\ldots,w_k}$ must satisfy equation \eq{ai8eq11} for all $\es\ne I\subseteq\{0,\ldots,N\}$ and $n\ge 0$.
\end{itemize}

Let $N\ge 0$, and suppose we are given data $\bigl((\sO^i),(f^{ij}),(\sH^i),\ab(\m^I)\bigr)$ satisfying all the conditions to be an $N$-simplex in $\AfCat^{\un,\iy}_{R,K,L}$, except that we do not require \eq{ai8eq4}--\eq{ai8eq5} to hold for all $I,n,x_0,\ldots,x_n$. So in particular $\m_0^{I;x_0}(1)=\boo^{i_2}_{f^{i_0i_2}(x_0)}$ if $K\ge 3$ and $\m_0^{I;x_0}(1)-\boo^{i_2}_{f^{i_0i_2}(x_0)}\in F^{\ssu 0}\Hom_{\sA^{i_2}}(f^{i_0i_2}(x_0),f^{i_0i_2}(x_0))$ if~$K\le 2$.

Define data $\dot\m^{I;x_0,\ldots,x_k}_k$ as in \eq{ai8eq16}. Then $\bigl((\sO^i),(f^{ij}),(\sH^i),\ab(\dot\m^I)\bigr)$ satisfies all the conditions to be an $N$-simplex in $\dAfCat^{\un,\iy}_{R,K,L}$, except that we do not require \eq{ai8eq11} to hold for all $I,n,x_0,\ldots,x_n$. In particular, $\dot\m_0^{I;x_0}(1)=0$ if $K\ge 3$ and $\dot\m_0^{I;x_0}(1)\in F^{\ssu 0}\Hom_{\sA^{i_2}}(f^{i_0i_2}(x_0),f^{i_0i_2}(x_0))$ if~$K\le 2$.

We will show below that \eq{ai8eq4}--\eq{ai8eq5} hold for the $\m_k^{J;w_0,\ldots,w_k}$ for given $I,\ab n,\ab x_0,\ab\ldots,x_n$ if and only if \eq{ai8eq11} holds for the $\dot\m_k^{J;w_0,\ldots,w_k}$ for the same $I,n,x_0,\ldots,x_n$. This implies that $I^\un$ in \eq{ai8eq14}--\eq{ai8eq16} is a 1-1 correspondence between $N$-simplices in $\AfCat^{\un,\iy}_{R,K,L}$ and $N$-simplices in $\dAfCat^{\un,\iy}_{R,K,L}$. It is easy to see that $I^\un$ commutes with the maps $\AfCat^{\un,\iy}_{R,K,L}(\al)$ in Definition \ref{ai8def1}(b) and $\dAfCat^{\un,\iy}_{R,K,L}(\al)$ in Definition \ref{ai8def6}(b) for all morphisms $\al:[\bar N]\ra[N]$ in $\De$. Thus $I^\un$ is a strict isomorphism of simplicial sets, and so of $\iy$-categories, and the proof will be complete.

Let $I=\{i_0<i_1<\cdots<i_M\}\subseteq\{0,\ldots,N\}$, $n\ge 0$, and $x_0,\ldots,x_n\in\sO^{i_0}$. We must show that \eq{ai8eq4}--\eq{ai8eq5} for the $\m_k^J$ for $I,n,x_0,\ldots,x_n$ are equivalent to \eq{ai8eq11} for the $\dot\m_k^J$ for $I,n,x_0,\ldots,x_n$. This is immediate for $M=0,1$, as then $\m_k^J=\dot\m_k^J$ for all terms in the equations, and \eq{ai8eq4}--\eq{ai8eq5} and \eq{ai8eq11} agree. So suppose~$M\ge 2$. 

We will show that if we rewrite \eq{ai8eq5} in terms of the $\dot\m_k^J$, and cancel terms using \eq{ai7eq17}--\eq{ai7eq19}, it becomes equivalent to \eq{ai8eq11}. To do this we must understand all terms involving $\m_0^{\{i_0,i_1,i_2\};w_0}$ which can make a nonzero contribution to \eq{ai8eq5}. The first term $\sum_{k,y_0,\ldots,y_k}\m_k^{I;y_0,\ldots,y_k}\ci\d_{n,k;y_0,\ldots,y_k}^{\{i_0\};x_0,\ldots,x_n}$ in \eq{ai8eq5} does not involve $\m_0^{\{i_0,i_1,i_2\};w_0}$. In the second term $(-1)^M\sum_{k, y_0,\ldots,y_k}\m_k^{\{i_M\};y_0,\ldots,y_k}\ci\d_{n,k;y_0,\ldots,y_k}^{I;x_0,\ldots,x_n}$, the expansion \eq{ai7eq9} of $\d_{n,k;y_0,\ldots,y_k}^{I;x_0,\ldots,x_n}$ may involve~$\m_0^{\{i_0,i_1,i_2\};w_0}$.

An argument we gave in the proof of Proposition \ref{ai7prop1} in \S\ref{ai102} shows that for a term $\pm C(\bs\Ga_1)\ot\cdots\ot C(\bs\Ga_k)$ in \eq{ai7eq9} to have a nonzero term involving an $\m_0^{\{i_0,i_1,i_2\};w_0}$, it is necessary that one tree $\bs\Ga_i$ should be of the form \eq{ai10eq8}, and then $C(\bs\Ga_i)=\boo^{i_M}_{y_j}$. But by \eq{ai7eq17}, substituting this $\boo^{i_M}_{y_j}$ into $\m_k^{\{i_M\};y_0,\ldots,y_k}$ gives zero unless $k=2$. Thus, the only nonzero terms involving an $\m_0^{\{i_0,i_1,i_2\};w_0}$ in \eq{ai8eq5} come from $(-1)^M\m_2^{\{i_M\};y_0,y_1,y_2}\ci\d_{n,2;y_0,y_1,y_2}^{I;x_0,\ldots,x_n}$, in a term $\pm C(\bs\Ga_1)\ot C(\bs\Ga_2)$ in \eq{ai7eq9} for $\d_{n,2;y_0,y_1,y_2}^{I;x_0,\ldots,x_n}$ in which either $\bs\Ga_1$ or $\bs\Ga_2$ is of the form~\eq{ai10eq8}.

Terms in \eq{ai7eq9} with $\bs\Ga_1$ of the form \eq{ai10eq8} look like this for~$0<j<M$:
\e
\begin{gathered}
\xymatrix@C=19pt@R=-4pt{
\overset{\!\!\!\!\!\!\{i_{j-1},i_j,i_{j+1}\}\!\!\!\!\!\!}{\sq} 
\ar[r] & 
\overset{\!\!\!\!\!\!\{i_{j+1},i_{j+2}\}\!\!\!\!\!\!}{\bu}  \ar[r] & \overset{\!\!\!\!\!\!\{i_{j+2},i_{j+3}\}\!\!\!\!\!\!}{\bu}  \ar[r] &  \cdots \ar[r] & 
\overset{\!\!\!\!\!\!\{i_{M-1},i_M\}\!\!\!\!\!\!}{\bu} \ar[r] & \tr & \!\!\!\!\!\!\!\!\!\!\!\!\bs\Ga_1, \\
&&& \!\!\!\!{\begin{subarray}{l} \ts 1\; \ci \\[-3pt] \ts \kern .91em\vdots \\[2pt] \ts n\; \ci\end{subarray}}
\ar[r] \ar@<1.5ex>[r]\ar@<-1.5ex>[r]
& \overset{I\sm\{i_j\}}{\bu} \ar[r] & \tr & \!\!\!\!\!\!\!\!\!\!\!\!\bs\Ga_2. }
\end{gathered}
\label{ai10eq47}
\e
These give an extra contribution to \eq{ai8eq5}:
\e
\sum_{j=1}^{M-1}(-1)^M(-1)^{M+j+1}\dot\m_n^{I\sm\{i_j\};x_0,\ldots,x_n}.
\label{ai10eq48}
\e

Terms in \eq{ai7eq9} with $\bs\Ga_2$ of the form \eq{ai10eq8} look like this for~$0<j<M$:
\begin{equation*}
\xymatrix@C=21pt@R=-4pt{
&& {\begin{subarray}{l} \ts 1\vphantom{\{}\; \ci \\[-3pt] \ts \kern .91em\vdots \\[2pt] \ts n\vphantom{\{}\; \ci\end{subarray}} \ar@<2.25ex>[r]  \ar@<.75ex>[r] \ar@<-.75ex>[r] \ar@<-2.25ex>[r] & \d^{\{i_0,\ldots,i_j\};x_0,\ldots,x_n}_{n,k;y_0,\ldots,y_k} 
\ar@<1.5ex>[r]  \ar[r] \ar@<-1.5ex>[r] 
& \overset{\{i_j,i_{j+1},\ldots,i_M\}}{\bu} \ar[r] & \tr & \!\!\!\!\!\!\!\!\!\!\!\!\bs\Ga_1, \\
&&&& \overset{\!\!\!\!\{i_0,i_j,i_M\}\!\!\!\!}{\sq} \ar[r] & \tr & \!\!\!\!\!\!\!\!\!\!\!\!\bs\Ga_2. }
\end{equation*}
Here we mean that the shown vertex $\bu$ in $\bs\Ga_1$ has $k$ incoming edges, and in $\d^{\{i_0,\ldots,i_j\};x_0,\ldots,x_n}_{n,k;y_0,\ldots,y_k}$ we substitute all the possible trees $\bs\Ga_1,\ldots,\bs\Ga_k$ in the sum for $\d^{\{i_0,\ldots,i_j\};x_0,\ldots,x_n}_{n,k;y_0,\ldots,y_k}$ in \eq{ai7eq9}. These give an extra contribution to~\eq{ai8eq5}: 
\e
\sum_{j=1}^{M-1}\sum_{k\ge 0, \; y_0,\ldots,y_k\in\sO^{i_j}\!\!\!\!\!\!\!\!\!\!\!\!\!\!\!\!\!\!}(-1)^M(-1)^{M+j}\,\dot\m_k^{\{i_j,\ldots,i_M\};y_0,\ldots,y_k}\ci\dot\d_{n,k}^{\{i_0,\ldots,i_j\};y_0,\ldots,y_k}.
\label{ai10eq49}
\e

In conclusion, if we rewrite \eq{ai8eq5} in terms of the $\dot\m_k^{J;w_0,\ldots,w_k}$ rather than the $\m_k^{J;w_0,\ldots,w_k}$, the extra terms coming from replacing all $\m_0^{\{i_0,i_1,i_2\};w_0}(1)$ by $\dot\m_0^{\{i_0,i_1,i_2\};w_0}+\boo^{i_2}_{f^{i_0i_2}(w_0)}$ are \eq{ai10eq48} and \eq{ai10eq49}. Thus the new equation satisfied by the $\dot\m_k^J$ is the sum of \eq{ai8eq5}, \eq{ai10eq48} and \eq{ai10eq49}. But this sum is \eq{ai8eq11}. Hence \eq{ai8eq4}--\eq{ai8eq5} for the $\m_k^{J;w_0,\ldots,w_k}$ are equivalent to \eq{ai8eq11} for the $\dot\m_k^{J;w_0,\ldots,w_k}$, and the proof of Theorem \ref{ai8thm5} is complete.

\subsection{Proofs of Theorems \ref{ai6thm7} and \ref{ai7thm7}}
\label{ai109}

We now prove Theorem \ref{ai7thm7}, covering the four cases of equation \eq{ai7eq30} in \S\ref{ai1091}--\S\ref{ai1094}. Theorem \ref{ai6thm7} follows from Theorem \ref{ai7thm7} by restricting to $A_\iy$-categories $\sA$ with~$\sO=\{*\}$. 

\subsubsection{Proof that $\ACat^\iy_{R,K}$ is an $\iy$-category}
\label{ai1091}

We first show that $\ACat^\iy_{R,K}$ is an $\iy$-category for $K\in\{3,4,5,\ldots,\iy\}$. Definition \ref{ai7def10} shows that $\ACat^\iy_{R,K}$ is a simplicial set. By Definition \ref{ai2def4}, to prove that $\ACat^\iy_{R,K}$ is an $\iy$-category we must show that every map of simplicial sets $f:\La^N_{i'}\ra \ACat^\iy_{R,K}$ can be extended to an $N$-simplex $\ti f:\De^N\ra\ACat^\iy_{R,K}$ for all~$0<i'<N$. 

Unwinding the definitions, this means:
\begin{itemize}
\setlength{\itemsep}{0pt}
\setlength{\parsep}{0pt}
\item[$(*)_{i'}^N$] Suppose we are given data $(\sO^i)_{i=0}^N,(f^{ij})_{0\le i\le j\le N},\ab(\sH^i)_{i=0}^N$ and $\m^I$ for all $\es\ne I\subseteq\{0,\ldots,N\}$ with $I\ne\{0,\ldots,N\},\{0,\ldots,N\}\sm\{i'\}$ (note this exclusion) as in Definitions \ref{ai7def8}(a)(i)--(iv) and \ref{ai7def10}, such that equation \eq{ai7eq25} holds for all $\es\ne I\subseteq\{0,\ldots,N\}$ with $I\ne\{0,\ldots,N\},\{0,\ldots,N\}\sm\{i'\}$ and $x_0,\ldots,x_n\in\sO^{\min I}$ for $n\ge 0$. We must show that there exist (not necessarily unique) choices of $\m^I$ for $I=\{0,\ldots,N\},\{0,\ldots,N\}\sm\{i'\}$, such that \eq{ai7eq25} also holds for $I=\{0,\ldots,N\},\{0,\ldots,N\}\sm\{i'\}$ and all $x_0,\ldots,x_n\in\sO^0$ for~$n\ge 0$.
\end{itemize}

Work in the situation of $(*)_{i'}^N$ for some $0<i'<N$. Choose arbitrary $R$-linear maps $\m_k^{\{0,\ldots,N\};x_0,\ldots,x_k}$ as in \eq{ai7eq24} graded of degree $2-N-k$ for all $k\ge 0$ and $x_0,\ldots,x_k\in\sO^0$, for example we could just set $\m_k^{\{0,\ldots,N\};x_0,\ldots,x_k}=0$. When $k=0$ we also require $\m_0^{\{0,\ldots,N\};x_0}=0$ if $K\ge N+1$, as in Definition \ref{ai7def8}(a)(iv). Consider \eq{ai7eq25} for $I=\{0,\ldots,N\}$ and $x_0,\ldots,x_n\in\sO^0$ for $n\ge 0$. This contains a term $(-1)^{i'}\m_n^{\{0,\ldots,N\}\sm\{i'\};x_0,\ldots,x_n}$ from $j=i'$ on the right hand side. Thus we can rearrange \eq{ai7eq25} for $I=\{0,\ldots,N\}$ and $x_0,\ldots,x_n$ to give
\ea
\m_n^{\{0,\ldots,N\}\sm\{i'\};x_0,\ldots,x_n}&=\sum_{\begin{subarray}{l} 0\le j\le N, \; 
0\le k\le\max(n,n+j-1) \\ y_0,\ldots,y_k\in\sO^{i_j} \end{subarray}
\!\!\!\!\!\!\!\!\!\!\!\!\!\!\!\!\!\!\!\!\!\!\!\!\!\!\!\!\!\!\!\!\!\!\!\!\!\!\!\!\!\!\!\!\!\!\!\!\!\!\!\!\!\!\!\!\!\!}
(-1)^{j-i'}\,\m_k^{\{j,j+1,\ldots,N\};y_0,\ldots,y_k}\ci\d_{n,k;y_0,\ldots,y_k}^{\{0,1,\ldots,j\};x_0,\ldots,x_n}
\nonumber\\
&-\sum_{0<j<M: \; j\ne i'}(-1)^{j-i'}\m_n^{\{0,\ldots,N\}\sm\{j\};x_0,\ldots,x_n}.
\label{ai10eq50}
\ea

We would like to take \eq{ai10eq50} as the {\it definition\/} of $\m_n^{\{0,\ldots,N\}\sm\{i'\};x_0,\ldots,x_n}$. (Note that $j=0$ in the first sum includes $\m_k^{\{0,\ldots,N\};y_0,\ldots,y_k}$, and when $j=N$ in the first sum, $\d_{n,k;y_0,\ldots,y_k}^{\{0,\ldots,N\};x_0,\ldots,x_n}$ can also involve $\m_l^{\{0,\ldots,N\};w_0,\ldots,w_l}$, which we chose arbitrarily above, so $\m_n^{\{0,\ldots,N\}\sm\{i'\};x_0,\ldots,x_n}$ will depend on the choices of $\m_l^{\{0,\ldots,N\};w_0,\ldots,w_l}$.) However, there is a problem: in the case $j=N$ in the first sum in \eq{ai10eq50}, the operator $\d_{n,k;y_0,\ldots,y_k}^{\{0,\ldots,N\};x_0,\ldots,x_n}$ from \eq{ai7eq9} can involve other operators $\m_l^{\{0,\ldots,N\}\sm\{i'\};w_0,\ldots,w_l}$.  To get round this problem we prove:

\begin{lem}
\label{ai10lem1}
If\/ $0<i'<N,$ the definition \eq{ai7eq9} of\/ $\d_{n,k;y_0,\ldots,y_k}^{\{0,\ldots,N\};x_0,\ldots,x_n}$ only involves\/ $\m_l^{\{0,\ldots,N\}\sm\{i'\};w_0,\ldots,w_l}$ for\/~$l<n$.
\end{lem}

\begin{proof}
In the definition of $\d_{n,k;y_0,\ldots,y_k}^{\{0,\ldots,N\};x_0,\ldots,x_n}$ in \eq{ai7eq9}, for an allowed choice of labelled trees $\bs\Ga_1,\ldots,\bs\Ga_k$ giving a nonzero term to contain an internal vertex or bud $v$ with $I_v=\{0,\ldots,N\}\sm\{i'\}$, for some $1\le i<j\le k$ we must have:
\begin{itemize}
\setlength{\itemsep}{0pt}
\setlength{\parsep}{0pt}
\item[(i)] $\bs\Ga_1,\ldots,\bs\Ga_{i-1}$ have all internal vertices $v$ with $I_v\!=\!\{m-1,m\}$ for $1\!\le\! m\!\le\!N$.
\item[(ii)] $\bs\Ga_i$ has exactly one internal vertex $w$ with $I_w=\{i'-1,i',i'+1\}$, and all other vertices $v$ have $I_v=\{i'-1,i'+1\}$ or $I_v=\{m-1,m\}$ for~$1\le m\le N$.

Note that $\bs\Ga_i$ cannot contain a bud $b$, as then $\m_0^{I_b;w_0}=0$ by Definition \ref{ai7def8}(a)(iv) as $\md{I_b}=2$ or 3. Hence $\bs\Ga_i$ has at least one leaf, and~$n_i\ge 1$.
\item[(iii)] $\bs\Ga_{i+1},\ldots,\bs\Ga_{j-1}$ have all internal vertices $v$ with $I_v=\{i'-1,i'+1\}$ or $I_v=\{m-1,m\}$ for $1\le m\le N$ with $m\ne i',i'+1$.
\item[(iv)] $\bs\Ga_j$ has exactly one internal vertex or bud $v$, with~$I_v=\{0,\ldots,N\}\sm\{i'\}$.
\item[(v)] $\bs\Ga_{j+1},\ldots,\bs\Ga_k$ have exactly one internal vertex $v$, which has~$I_v=\{0,N\}$.
\end{itemize}
The case $N=2$, $i'=1$ is degenerate as (iii),(iv),(v) coincide, but the following argument still works. We have $n=n_1+\cdots+n_k$, with $n_l\ge 0$ the number of leaves of $\bs\Ga_l$. As $n_i\ge 1$ we have $n_j<n$. But the vertex $v$ in (iv) has $n_j$ incoming edges, and corresponds to the operator $\m_{n_j}^{\{0,\ldots,N\}\sm\{i'\};w_0,\ldots,w_{n_j}}$. The lemma follows.
\end{proof}

We now use equation \eq{ai10eq50} to define $\m_n^{\{0,\ldots,N\}\sm\{i'\};x_0,\ldots,x_n}$ by induction on $n=0,1,\ldots.$ The right hand side of \eq{ai10eq50} depends on the $\m^{I;w_0,\ldots,w_k}_k$ for $I\!\ne\!\{0,\ldots,N\},\ab\{0,\ldots,N\}\sm\{i'\}$ given in $(*)_{i'}^N$, and on the arbitrary values of $\m_k^{\{0,\ldots,N\};w_0,\ldots,w_k}$ chosen above, and on $\m_l^{\{0,\ldots,N\}\sm\{i'\};w_0,\ldots,w_l}$ for $l<n$, which is defined by induction. Thus, we define unique $\m_n^{\{0,\ldots,N\}\sm\{i'\};x_0,\ldots,x_n}$ for all $n\ge 0$ and $x_0,\ldots,x_n\in\sO^0$ such that \eq{ai10eq50} holds for all $n\ge 0$ and $x_0,\ldots,x_n\in\sO^0$. Equivalently, equation \eq{ai7eq25} holds for $I=\{0,\ldots,N\}$ and all $n\ge 0$ and~$x_0,\ldots,x_n\in\sO^0$.

To satisfy Definition \ref{ai7def8}(a)(iv), when $n=0$ we need $\m_0^{\{0,\ldots,N\}\sm\{i'\};x_0}=0$ if $K\ge N$. From \eq{ai10eq50} and Definition \ref{ai7def2}(a) we find that
\ea
&\m_0^{\{0,\ldots,N\}\sm\{i'\};x_0}
\label{ai10eq51}\\
&=(-1)^{1-i'}
\m_0^{\{1,2,\ldots,N\};x_0}
-\sum_{0<j<M: \; j\ne i'\!\!\!\!\!\!\!\!\!\!\!\!\!}(-1)^{j-i'}\m_0^{\{0,\ldots,N\}\sm\{j\};x_0}.
\nonumber
\ea
But if $K\ge N$ each term on the right hand side is zero, so~$\m_0^{\{0,\ldots,N\}\sm\{i'\};x_0}=0$.

It remains to show that \eq{ai7eq25} holds for $I=\{0,\ldots,N\}\sm\{i'\}$ and $x_0,\ldots,x_n\in\sO^0$. We do this by induction on $n=0,1,\ldots.$ Let $\ti n\ge 0$, and suppose by induction that \eq{ai7eq25} holds for $I=\{0,\ldots,N\}\sm\{i'\}$ and all $0\le n<\ti n$ and $x_0,\ldots,x_n\in\sO^0$. The initial step $\ti n=0$ is trivial. To prove the inductive step, for $x_0,\ldots,x_{\ti n}\in\sO^0$, consider the equation
\ea
&\sum_{ 0\le j<i', \; k\ge 0,\; y_0,\ldots,y_k\in\sO^j
\!\!\!\!\!\!\!\!\!\!\!\!\!\!\!\!\!\!\!\!\!\!\!\!\!\!\!\!\!\!\!\!\!\!\!\!\!\!\!\!\!\!\!\!\!\!\!\!\!}
(-1)^j\,\m_k^{\{j,\ldots,N\}\sm\{i'\};y_0,\ldots,y_k}\ci\d_{\ti n,k;y_0,\ldots,y_k}^{\{0,1,\ldots,j\};x_0,\ldots,x_{\ti n}}
\label{ai10eq52}\\
&-\sum_{i'<j\le N, \; k\ge 0,\; y_0,\ldots,y_k\in\sO^j
\!\!\!\!\!\!\!\!\!\!\!\!\!\!\!\!\!\!\!\!\!\!\!\!\!\!\!\!\!\!\!\!\!\!\!\!\!\!\!\!\!\!\!\!\!\!\!\!\!\!}
(-1)^j\,\m_k^{\{j,\ldots,N\};y_0,\ldots,y_k}\ci\d_{\ti n,k;y_0,\ldots,y_k}^{\{0,1,\ldots,j\}\sm\{i'\};x_0,\ldots,x_{\ti n}}
\allowdisplaybreaks
\nonumber\\
&=\sum_{0\le j<i', \; k\ge 0,\; y_0,\ldots,y_k\in\sO^j
\!\!\!\!\!\!\!\!\!\!\!\!\!\!\!\!\!\!\!\!\!\!\!\!\!\!\!\!\!\!\!\!\!\!\!\!\!\!\!\!\!\!\!\!\!\!\!\!\!\!\!}
(-1)^j\m_k^{\{j,\ldots,N\}\sm\{i'\};y_0,\ldots,y_k}\ci \d_{\ti n,k;y_0,\ldots,y_k}^{\{0,\ldots,j\};x_0,\ldots,x_{\ti n}}
\allowdisplaybreaks
\nonumber\\
&-\sum_{i'<j\le N,\; k\ge 0,\; y_0,\ldots,y_k\in\sO^j
\!\!\!\!\!\!\!\!\!\!\!\!\!\!\!\!\!\!\!\!\!\!\!\!\!\!\!\!\!\!\!\!\!\!\!\!\!\!\!\!\!\!\!\!\!\!\!\!\!\!\!\!}
(-1)^j\m_k^{\{j,\ldots,N\};y_0,\ldots,y_k}\ci\d_{\ti n,k;y_0,\ldots,y_k}^{\{0,\ldots,j\}\sm\{i'\};x_0,\ldots,x_{\ti n}}
\allowdisplaybreaks
\nonumber\\
&-\sum_{\begin{subarray}{l} 0\le j\le N \end{subarray}}\begin{aligned}[t]&(-1)^{i'}\biggl[\sum_{\begin{subarray}{l} l\ge 0, \\ z_0,\ldots,z_l\in\sO^j \end{subarray}}\begin{aligned}[t]& 
\biggl(\sum_{j\le j'\le N, \; k\ge 0, \; y_0,\ldots,y_k\in\sO^{j'} \!\!\!\!\!\!\!\!\!\!\!\!\!\!\!\!\!\!\!\!\!\!\!\!\!\!\!\!\!\!\!\!\!\!\!\!\!\!\!\!\!\!\!\!\!\!\!\!\!\!\!\!\!\!\!}
(-1)^{j'-j}\m_k^{\{j',\ldots,N\};y_0,\ldots,y_k}\ci\d_{l,k;y_0,\ldots,y_k}^{\{j,\ldots,j'\};z_0,\ldots,z_l}\biggr)\\
&\ci \d_{\ti n,l;z_0,\ldots,z_l}^{\{0,\ldots,j\};x_0,\ldots,x_{\ti n}} 
\end{aligned}\\ 
&-(-1)^j\sum_{\begin{subarray}{l} k\ge 0, \\ y_0,\ldots,y_k\in\sO^j \end{subarray}}\begin{aligned}[t]& \m_k^{\{j,\ldots,N\};y_0,\ldots,y_k}\ci{}\\
&\biggl(\sum_{0\le j'\le j,\; l\ge 0, \; z_0,\ldots,z_l\in\sO^{j'}  
\!\!\!\!\!\!\!\!\!\!\!\!\!\!\!\!\!\!\!\!\!\!\!\!\!\!\!\!\!\!\!\!\!\!\!\!\!\!\!\!\!\!\!\!\!\!\!\!\!}
(-1)^{j'}\d_{l,k;y_0,\ldots,y_k}^{\{j',\ldots,j\};z_0,\ldots,z_l}\ci\d_{\ti n,l;z_0,\ldots,z_l}^{\{0,\ldots,j'\};x_0,\ldots,x_{\ti n}}\biggr)\biggr]
\end{aligned}\end{aligned}
\allowdisplaybreaks
\nonumber\\
&=\sum_{0\le j<i', \; k\ge 0,\; y_0,\ldots,y_k\in\sO^j
\!\!\!\!\!\!\!\!\!\!\!\!\!\!\!\!\!\!\!\!\!\!\!\!\!\!\!\!\!\!\!\!\!\!\!\!\!\!\!\!\!\!\!\!\!\!\!\!\!\!\!\!}(-1)^j\m_k^{\{j,\ldots,N\}\sm\{i'\};y_0,\ldots,y_k}\ci \d_{\ti n,k;y_0,\ldots,y_k}^{\{0,\ldots,j\};x_0,\ldots,x_{\ti n}}
\allowdisplaybreaks
\nonumber\\
&-\sum_{i'<j\le N,\; k\ge 0,\; y_0,\ldots,y_k\in\sO^j
\!\!\!\!\!\!\!\!\!\!\!\!\!\!\!\!\!\!\!\!\!\!\!\!\!\!\!\!\!\!\!\!\!\!\!\!\!\!\!\!\!\!\!\!\!\!\!\!\!\!\!\!}(-1)^j\m_k^{\{j,\ldots,N\};y_0,\ldots,y_k}\ci\d_{\ti n,k;y_0,\ldots,y_k}^{\{0,\ldots,j\}\sm\{i'\};x_0,\ldots,x_{\ti n}}
\allowdisplaybreaks
\nonumber\\
&-\sum_{\begin{subarray}{l} 0\le j\le N \end{subarray}}\begin{aligned}[t]&(-1)^{i'}\biggl[\sum_{\begin{subarray}{l} l\ge 0, \\ z_0,\ldots,z_l\in\sO^j \end{subarray}\!\!\!\!\!\!\!\!}\biggl(\sum_{j<j'<N\!\!\!\!\!\!}(-1)^{j'-j}\m_l^{\{j,\ldots,N\}\sm\{j'\};z_0,\ldots,z_l}\biggr)\ci \d_{\ti n,l;z_0,\ldots,z_l}^{\{0,\ldots,j\};x_0,\ldots,x_{\ti n}} \\ 
&-(-1)^j\sum_{\begin{subarray}{l} k\ge 0, \\ y_0,\ldots,y_k\in\sO^j \end{subarray}\!\!\!\!\!\!\!\!}\m_k^{\{j,\ldots,N\};y_0,\ldots,y_k}\ci\biggl(\sum_{0<j'<j\!\!\!\!\!\!}(-1)^{j'}\d_{\ti n,k;y_0,\ldots,y_k}^{\{0,\ldots,j\}\sm\{j'\};x_0,\ldots,x_{\ti n}}\biggr)\biggr]\end{aligned}
\allowdisplaybreaks
\nonumber\\
&=\sum_{\begin{subarray}{l} 0<j'<N \\ j'\ne i'\end{subarray}}(-1)^{i'+j'+1}\begin{aligned}[t]&\biggl(\sum_{0\le j<j', \; l\ge 0,\; z_0,\ldots,z_l\in\sO^j \!\!\!\!\!\!\!\!\!\!\!\!\!\!\!\!\!\!\!\!\!\!\!\!\!\!\!\!\!\!\!\!\!\!\!\!\!\!\!\!\!}
(-1)^j\m_l^{\{j,\ldots,N\}\sm\{j'\};z_0,\ldots,z_l}\ci \d_{\ti n,l;z_0,\ldots,z_l}^{\{0,\ldots,j\};x_0,\ldots,x_{\ti n}} \\ 
&-\sum_{j'<j\le N,\; k\ge 0,\; y_0,\ldots,y_k\in\sO^j \!\!\!\!\!\!\!\!\!\!\!\!\!\!\!\!\!\!\!\!\!\!\!\!\!\!\!\!\!\!\!\!\!\!\!\!\!\!\!\!\!}(-1)^j\m_k^{\{j,\ldots,N\};y_0,\ldots,y_k}\ci\d_{\ti n,k;y_0,\ldots,y_k}^{\{0,\ldots,j\}\sm\{j'\};x_0,\ldots,x_{\ti n}}\biggr)\end{aligned}
\allowdisplaybreaks
\nonumber\\
&=\sum_{\begin{subarray}{l} 0<j'<N \\ j'\ne i'\end{subarray}\!\!\!\!\!}(-1)^{i'+j'+1}\biggl(\begin{aligned}[t]& \sum_{0<j<j'\!\!\!}(-1)^j\m_{\ti n}^{\{0,\ldots,N\}\sm\{j,j'\};x_0,\ldots,x_{\ti n}}\\
&-\sum_{j'<j<N\!\!\!}(-1)^j\m_{\ti n}^{\{0,\ldots,N\}\sm\{j',j\};x_0,\ldots,x_{\ti n}}\biggr)\end{aligned}
\allowdisplaybreaks
\nonumber\\
&=\sum_{0<j<i'}(-1)^j\m_{\ti n}^{\{0,\ldots,N\}\sm\{j,i'\};x_0,\ldots,x_{\ti n}}-\sum_{i'<j<N}(-1)^j\m_{\ti n}^{\{0,\ldots,N\}\sm\{i',j\};x_0,\ldots,x_{\ti n}}.
\nonumber
\ea

Here in the first step the final sum $\sum_{0\le j\le N}(-1)^{i'}[\cdots]$ is zero, as the two inner sums cancel on swapping $j,j'$. In the second step we use \eq{ai7eq25} for $I=\{j,\ldots,N\}$ and $n=l$ in the second line, and \eq{ai7eq27} for $I=\{0,\ldots,j\}$ and $n=\ti n$ in the third. Note that when $j=0$, equation \eq{ai7eq25} for $I=\{0,\ldots,N\}$ holds as above by construction. 

\begin{rem}
\label{ai10rem1}
When $j=N$, we just used \eq{ai7eq27} for $I=\{0,\ldots,N\}$. This was proved in Theorem \ref{ai7thm5} (which also holds for $\ACat^\iy_{R,K}$) assuming $\bigl((\sO^i),\ab(f^{ij}),\ab(\sH^i),\ab(\m^I)\bigr)$ is an $N$-simplex in $\ACat^\iy_{R,K}$, so in particular assuming that \eq{ai7eq25} holds for both $I=\{0,\ldots,N\}$, which is true as above, and {\it also\/} for $I=\{0,\ldots,N\}\sm\{i'\}$, which is what we are trying to prove. In fact we can show using Lemma \ref{ai10lem1} that the proof of \eq{ai7eq27} in \S\ref{ai107} for $I=\{0,\ldots,N\}$ and $n=\ti n$ in Theorem \ref{ai7thm5} involved only \eq{ai7eq25} for $I=\{0,\ldots,N\}\sm\{i'\}$ and $n<\ti n$, which holds by induction, so the argument is not circular. 
\end{rem}

In the third step of \eq{ai10eq52} we rearrange the sums. The fourth step holds by \eq{ai7eq25} for $I=\{0,\ldots,N\}\sm\{j'\}$ for $j'\ne i'$. The fifth step holds as the terms in $\m_n^{\{0,\ldots,N\}\sm\{j,j'\};x_0,\ldots,x_{\ti n}}$ for $j\ne i'\ne j'$ cancel in pairs. Thus equation \eq{ai10eq52} proves \eq{ai7eq25} for $I=\{0,\ldots,N\}\sm\{i'\}$ and $n=\ti n$. This completes the inductive step, so \eq{ai7eq25} holds for $I=\{0,\ldots,N\}\sm\{i'\}$ and all $n\ge 0$ and $x_0,\ldots,x_n\in\sO^0$. Therefore $(*)_{i'}^N$ holds, and $\ACat^\iy_{R,K}$ is an $\iy$-category.

\subsubsection{Proof that $\ACat^{\cu,\iy}_{R,K}$ is an $\iy$-category}
\label{ai1092}

For $\ACat^{\cu,\iy}_{R,K}$, we modify \S\ref{ai1091} as follows: to prove $(*)_{i'}^N$ for $N>2$ we suppose the initial data $(\sO^i),(f^{ij}),\ab(\sH^i)$ and $\m^I$ for $I\ne\{0,\ldots,N\},\{0,\ldots,N\}\sm\{i'\}$ is chosen as in Definition \ref{ai7def9}(v$)'$,(vi$)'$, so that $\sA^i$ and $\fF^{ij}$ are cohomologically unital. Then the proof in \S\ref{ai1091} works without change. When $N=2$, to prove $(*)_1^2$ we also need to show we can choose $\m^{\{0,2\}}$ such that $\fF^{02}$ is cohomologically unital. We do this by taking $\m^{\{0,1,2\}}=0$ and~$\fF^{02}=\fF^{12}\ci\fF^{01}$.

\subsubsection{Proof that $\dACat^{\un,\iy}_{R,K}$ is an $\iy$-category}
\label{ai1093}

For $\dACat^{\un,\iy}_{R,K}$, we modify \S\ref{ai1091} as follows: when $N>2$ we require the initial data $(\sO^i),(f^{ij}),\ab(\sH^i)$ and $\m^I$ for $I\ne\{0,\ldots,N\},\ab\{0,\ldots,N\}\sm\{i'\}$ to satisfy Definition \ref{ai7def3}(a)(v), as in Definition \ref{ai7def8}(a)(v). We also choose the data $\m_k^{\{0,\ldots,N\};x_0,\ldots,x_k}$ above to satisfy \eq{ai7eq19}. Then \S\ref{ai1091} constructs unique data $\m_n^{\{0,\ldots,N\}\sm\{i'\};x_0,\ldots,x_n}$ such that \eq{ai7eq25} holds for $I=\{0,\ldots,N\},\{0,\ldots,N\}\sm\{i'\}$ and all $n\ge 0$ and $x_0,\ldots,x_n\in\sO^0$.

We need show that this data $\m_n^{\{0,\ldots,N\}\sm\{i'\};x_0,\ldots,x_n}$ also satisfies \eq{ai7eq19}. We do this by induction on $n$, so for $n\ge 0$ given, suppose by induction that the $\m_k^{\{0,\ldots,N\}\sm\{i'\};x_0,\ldots,x_k}$ satisfy \eq{ai7eq19} for $k<n$. For $1\le m\le n$, suppose $x_0,\ldots,x_n\in\sO^0$ with $x_{m-1}=x_m$. We must prove that
\begin{equation*}
\m_n^{\{0,\ldots,N\}\sm\{i'\};x_0,\ldots,x_n}\ci\bigl(\id^{\ot^{m-1}}\ot\boo^0_{x_m}\ot\id^{\ot^{n-m}}\bigr)=0.
\end{equation*}
We do this by showing that the composition of the right hand side of \eq{ai10eq50} with $\id^{\ot^{m-1}}\ot\boo^0_{x_m}\ot\id^{\ot^{n-m}}$ is zero. In more detail:
\smallskip

\noindent(a) In the term $\m_n^{\{0,\ldots,N\}\sm\{j\};x_0,\ldots,x_n}$ on the right hand side of \eq{ai10eq50}, we must have $N\ge 3$, as if $N=2$ then $i'=1$ and the condition $0<j<2$, $j\ne i'$ has no solutions for $j$. Hence the composition of $\m_n^{\{0,\ldots,N\}\sm\{j\};x_0,\ldots,x_n}$ on the right hand side of \eq{ai10eq50} with $\id^{m-1}\ot\boo^0_{x_m}\ot\id^{n-m}$ is zero by Definition~\ref{ai7def3}(a)(v).
\smallskip

\noindent(b) Consider when the composition of $\m_k^{\{j,\ldots,N\};y_0,\ldots,y_k}\!\ci\d_{n,k;y_0,\ldots,y_k}^{\{0,\ldots,j\};x_0,\ldots,x_n}$ in \eq{ai10eq50} with $\id^{m-1}\ot\boo^0_{x_m}\ot\id^{n-m}$ can be nonzero. Divide into cases:
\begin{itemize}
\setlength{\itemsep}{0pt}
\setlength{\parsep}{0pt}
\item[(i)] $j=0$ and $1<m<n$. Then $\d_{n,k;y_0,\ldots,y_k}^{\{0\};x_0,\ldots,x_n}$ is defined by \eq{ai7eq8}, with conditions on $y_0,\ldots,y_k$ to get a nonzero answer. If $\boo^0_{x_m}$ is inserted into an $\id_{\sH^0(x_{i-1},x_i)}$ in \eq{ai7eq8} then $\boo^0_{x_m}$ is passed on as an input to $\m_k^{\{0,\ldots,N\};y_0,\ldots,y_k}$, and so gives zero as $\m_k^{\{0,\ldots,N\};y_0,\ldots,y_k}$ satisfies \eq{ai7eq19}. 

If $\boo^0_{x_m}$ is inserted into $\m_{n-k+1}^{\{0\};x_j,\ldots,x_{j+n-k+1}}$ in \eq{ai7eq8} then by \eq{ai7eq17}--\eq{ai7eq19} we get a nonzero term only if $n-k+1=2$, so $k=n-1$. Then the sign in \eq{ai7eq8} simplifies to $(-1)^j$. To insert $\boo^0_{x_m}$ into $\m_2^{\{0\};y_0,y_1,y_2}$ we must take $j=m-2$ or $j=m-1$ in \eq{ai7eq8}, and by \eq{ai7eq17} the compositions of these terms in $\sum_{y_0,\ldots,y_{n-1}}\d_{n,n-1;y_0,\ldots,y_{n-1}}^{\{0\};x_0,\ldots,x_n}$ with $\id^{m-1}\ot\boo^0_{x_m}\ot\id^{n-m}$ are $(-1)^{m-2}\id^{n-1}$ and $(-1)^{m-1}\id^{n-1}$, which cancel. Thus this case gives~zero.
\item[(ii)] $j=0$ and $m=1$. Then the case $j=m-2$ in (i) does not occur, so the composition of $\sum_{y_0,\ldots,y_{n-1}}\d_{n,n-1;y_0,\ldots,y_{n-1}}^{\{0\};x_0,\ldots,x_n}$ with $\id^{m-1}\ot\boo^0_{x_m}\ot\id^{n-m}$ is $\id^{n-1}$, and we get a contribution of $(-1)^{i'}\m_{n-1}^{\{0,\ldots,N\};x_1,\ldots,x_n}$ to the composition of \eq{ai10eq50} with $\id^{m-1}\ot\boo^0_{x_m}\ot\id^{n-m}$.
\item[(iii)] $j=0$ and $m=n$. Then the case $j=m-1$ in (i) does not occur, and we get a contribution of $(-1)^{i'+n}\m_{n-1}^{\{0,\ldots,N\};x_0,\ldots,x_{n-1}}$ to the composition.
\item[(iv)] $0<j<N-1$. Then as $j>0$, $\d_{n,k;y_0,\ldots,y_k}^{\{0,\ldots,j\};x_0,\ldots,x_n}$ is defined by \eq{ai7eq9}, and it involves only operators $\m_l^{I;w_0,\ldots,w_l}$ for $I\subseteq\{0,1,\ldots,j\}$ with $\md{I}\ge 2$. Thus by \eq{ai7eq18}--\eq{ai7eq19}, the only way to get a nonzero composition of $\d_{n,k;y_0,\ldots,y_k}^{\{0,\ldots,j\};x_0,\ldots,x_n}$ with $\id^{m-1}\ot\boo^0_{x_m}\ot\id^{n-m}$ is by inserting the $\boo^0_{x_m}$ into a tensor factor $\m_1^{\{i_{l-1},i_l\};w_{l-1},w_{l-1}}\ci\m_1^{\{i_{l-2},i_{l-1}\};w_{l-2},w_{l-2}}\ci\cdots\ci\m_1^{\{i_0,i_1\};x_m,x_m}$ in a term in $\d_{n,k;y_0,\ldots,y_k}^{\{0,\ldots,j\};x_0,\ldots,x_n}$ for $0=i_0<i_1<\cdots<i_l=j$. By \eq{ai7eq18}, the composition of this with $\boo^0_{x_m}$ is $\boo^j_{w_l}$, which then appears as a input in $\m_k^{\{j,\ldots,N\};y_0,\ldots,y_k}$. As $j<N-1$ this gives zero by \eq{ai7eq19}.
\item[(v)] $j=N-1$, and either $n>1$ or $N>2$. As in (iv), to get a nonzero term we must insert $\boo^0_{x_m}$ into a tensor factor $\m_1^{\{i_{l-1},i_l\};w_{l-1},w_{l-1}}\ci\cdots\ci\m_1^{\{i_0,i_1\};x_m,x_m}$, giving an input $\boo^j_{w_l}$ to $\m_k^{\{N-1,N\};y_0,\ldots,y_k}$. By \eq{ai7eq18} this can only be nonzero if $k=1$. But then for a term in $\d_{n,k;y_0,\ldots,y_k}^{\{0,\ldots,N-1\};x_0,\ldots,x_n}$ to have a tensor factor $\m_1^{\{i_{l-1},i_l\};w_{l-1},w_{l-1}}\ci\cdots\ci\m_1^{\{i_0,i_1\};x_m,x_m}$ we must have $n=1$, so we get zero if $n>1$. If $n=1$ and $N>2$ then \eq{ai7eq10} says that $\sum_{y_0,y_1}\d_{1,1;y_0,y_1}^{\{0,\ldots,N-1\};x_0,x_1}=\m_1^{\{0,\ldots,N-1\};x_0,x_1}$, contradicting the tensor factor $\m_1^{\{i_{l-1},i_l\};w_l,w_l}\ci\cdots\ci\m_1^{\{i_0,i_1\};x_m,x_m}$ as $N>2$. So this case gives zero.
\item[(vi)] $N=2$, $i'=1$, $m=n=1$ and $j=1$. As in (v) we get $\sum_{y_1}\d_{1,1;y_1,y_1}^{\{0,1\};x_1,x_1}=\m_1^{\{0,1\};x_1,x_1}$, and $\m_1^{\{1,2\};y_1,y_1}\ci\m_1^{\{0,1\};x_1,x_1}(\boo_{x_1}^0)=\boo^2_{z_1}$ by \eq{ai7eq18} for $z_1=f^{12}\ci f^{01}(x_1)$, so this case contributes $\boo^2_{z_1}$ to \eq{ai10eq50} composed with~$\boo^0_{x_1}$.
\item[(vii)] $j=N$, and $1<m<n$. As in (iv), to get a nonzero term we must insert $\boo^0_{x_m}$ into a tensor factor $\m_1^{\{i_{l-1},i_l\};w_{l-1},w_{l-1}}\ci\cdots\ci\m_1^{\{i_0,i_1\};x_m,x_m}$ in a term in $\d_{n,k;y_0,\ldots,y_k}^{\{0,1,\ldots,N\};x_0,\ldots,x_n}$, giving an input $\boo^N_z$ to $\m_k^{\{N\};y_0,\ldots,y_k}$. By \eq{ai7eq17}--\eq{ai7eq19} this can only be nonzero if $k=2$. However, in the definition of $\d_{n,k;y_0,\ldots,y_k}^{\{0,\ldots,N\};x_0,\ldots,x_n}$ in \eq{ai7eq9}, we need at least one graph $\bs\Ga_i$ to cover leaves $1,\ldots,m-1$, and one graph $\bs\Ga_i$ for the factor $\m_1^{\{i_{l-1},i_l\};w_{l-1},w_{l-1}}\ci\cdots\ci\m_1^{\{i_0,i_1\};x_m,x_m}$ on leaf $m$, and at least one graph $\bs\Ga_i$ to cover leaves $m+1,\ldots,n$, giving $k\ge 3$. So this term contributes zero.
\item[(viii)] $j=N$, and $m=1$. As in (vii) we must insert $\boo^0_{x_m}$ into a tensor factor $\m_1^{\{i_{l-1},i_l\};w_{l-1},w_{l-1}}\ci\cdots\ci\m_1^{\{i_0,i_1\};x_m,x_m}$ in leaf 1 in $\d_{n,k}^{\{0,1,\ldots,N\}}$, and we must have $k=2$. The only possible term in $\d_{n,2;y_0,y_1,y_2}^{\{0,\ldots,N\};x_0,\ldots,x_n}$ giving a nonzero contribution is $(-1)^{N+1}(\m_1^{\{N-1,N\};w_{N-1},w_{N-1}}\ci\m_1^{\{N-2,N-1\};w_{N-2},w_{N-2}}\ci\cdots\ci\m_1^{\{0,1\};x_1,x_1})\ot\m_{n-1}^{\{0,\ldots,N\};x_1,\ldots,x_n}$. Composing this with $\boo^0_{x_1}\ot\id^{n-1}$ and $\sum_{y_i}(-1)^{N+i'}\m_2^{\{N\};y_0,y_1,y_2}$ gives a contribution of $(-1)^{1+i'}\m_{n-1}^{\{0,\ldots,N\};x_1,\ldots,x_n}$ to the composition of \eq{ai10eq50} with $\id^{m-1}\ot\boo^0_{x_m}\ot\id^{n-m}$.
\item[(ix)] $j=N$, and $m=n$. As in (vii)--(viii) we must insert $\boo^0_{x_m}$ into a tensor factor $\m_1^{\{i_{l-1},i_l\};w_{l-1},w_{l-1}}\ci\cdots\ci\m_1^{\{i_0,i_1\};x_m,x_m}$ in leaf $n$ in $\d_{n,k;y_0,\ldots,y_k}^{\{0,\ldots,N\};x_0,\ldots,x_n}$, and we must have $k=2$. The only possible term in $\d_{n,2;y_0,y_1,y_2}^{\{0,\ldots,N\};x_0,\ldots,x_n}$ giving a nonzero contribution is $(-1)^{n+N+1}\m_{n-1}^{\{0,\ldots,N\};x_0,\ldots,x_{n-1}}\ot \m_1^{\{0,N\};x_n,x_n}$. Composing this with $\id^{n-1}\ot\boo^0_{x_n}$ and $(-1)^{i'+N}\m_2^{\{N\}}$ gives a contribution of $(-1)^{1+i'+n}\m_{n-1}^{\{0,\ldots,N\};x_0,\ldots,x_{n-1}}$ to the composition of \eq{ai10eq50} with $\id^{m-1}\ot\boo^0_{x_m}\ot\id^{n-m}$.
\end{itemize}

Observe that when $m=1$ cases (ii) and (viii) cancel, and when $m=n$ cases (iii) and (ix) cancel. The only other nonzero case is (vi) when $N=2$, $i'=1$, and $n=1$, which gives \eq{ai7eq18} for $\m_n^{\{0,\ldots,N\}\sm\{i'\};x_0,\ldots,x_n}$. In all other cases the composition of $\m_n^{\{0,\ldots,N\}\sm\{i'\};x_0,\ldots,x_n}$ with $\id^{m-1}\ot\boo^0_{x_m}\ot\id^{n-m}$ is zero, as in \eq{ai7eq18}--\eq{ai7eq19}. This completes the inductive step. Hence $\m_n^{\{0,\ldots,N\}\sm\{i'\};x_0,\ldots,x_n}$ satisfies \eq{ai7eq19} when $N>2$.

Thus $\dACat^{\un,\iy}_{R,K}$ satisfies the extension condition for maps $f:\La^N_{i'}\ra \dACat^{\un,\iy}_{R,K}$ provided $N>2$. When $N=2$, to prove $(*)_1^2$, as in \S\ref{ai1092} we take $\m^{\{0,1,2\}}=0$ and $\fF^{02}=\fF^{12}\ci\fF^{01}$, so that $\fF^{02}$ is unital.

\subsubsection{Proof that $\ACat^{\un,\iy}_{R,K}$ is an $\iy$-category}
\label{ai1094}

Theorem \ref{ai7thm6} shows $\ACat^{\un,\iy}_{R,K}$ and $\dACat^{\un,\iy}_{R,K}$ are isomorphic as simplicial sets, without assuming they are $\iy$-categories. Thus $\dACat^{\un,\iy}_{R,K}$ an $\iy$-category in \S\ref{ai1093} implies that $\ACat^{\un,\iy}_{R,K}$ is an $\iy$-category.

This completes the proof of Theorem~\ref{ai7thm7}.

\subsection{Proof of Theorem \ref{ai8thm6}}
\label{ai1010}

Theorem \ref{ai8thm6} is the generalization of Theorem \ref{ai7thm7}, proved in \S\ref{ai109}, to the filtered case. We start by showing that $\AfCat^\iy_{R,K,L}$ is an $\iy$-category for $K,L\in\{0,1,2,\ldots,\iy\}$ with $K\le L$ and $L\ge 3$, generalizing the proof for $\ACat^\iy_{R,K}$ with $K\ge 3$ in \S\ref{ai1091}. An important difference with \S\ref{ai1091} is that Lemma \ref{ai10lem1} no longer holds in the filtered case, since part (ii) of the proof used the fact that in $\ACat^\iy_{R,K}$ for $K\ge 3$ we have $\m_0^{I_b;w_0}=0$ if $\md{I_b}=2$ or 3. However, in $\AfCat^\iy_{R,K,L}$ when $K<3$ we can have $\m_0^{I_b;w_0}\ne 0$ if $\md{I_b}=2$ or~3.

This is used in the proof of Lemma \ref{ai10lem1} to show that $n_i\ge 1$, so that $\m_l^{\{0,\ldots,N\}\sm\{i'\};w_0,\ldots,w_l}$ has $l=n_j\le n-1$. If instead we allow $n_i=0$ then $l=n_j\le n$. Because of this, the appropriate generalization of Lemma \ref{ai10lem1} is:

\begin{lem}
\label{ai10lem2}
Working in $\AfCat^\iy_{R,K,L},$ if\/ $0<i'<N,$ the definition \eq{ai7eq9} of\/ $\d_{n,k;y_0,\ldots,y_k}^{\{0,\ldots,N\};x_0,\ldots,x_n}$ only involves\/ $\m_l^{\{0,\ldots,N\}\sm\{i'\};w_0,\ldots,w_l}$ for\/ $l<n,$ plus terms in $\m_n^{\{0,\ldots,N\}\sm\{i'\};x_0,\ldots,x_n}$ (with the same $x_0,\ldots,x_n$) and at least one factor $\m_0^{I_b;z_0}$ with\/ $\md{I_b}=2$ or $3$ and\/ $\md{I_b}>K$. Here $z_0$ must be of the form $\fF^J(x_0)$ for $\es\ne J\subseteq\{0,\ldots,N\}$ with\/ $\min J=0$ and\/ $\max J=\min I_b$.
\end{lem}

There are only finitely many possible choices for $I_b,z_0$ in Lemma \ref{ai10lem2}, and all have $\m_0^{I_b;z_0}(1)\in F^{\ssu 0}\Hom_{\sA^*}(*,*)$ by Definition \ref{ai8def4}(i). Hence we may choose $\ep\in\Ga$ with $\ep\ssu 0$, such that $\m_0^{I_b;z_0}(1)\in F^\ep\Hom_{\sA^*}(*,*)$ for all such $I_b,z_0$.

Observe that we can rewrite \eq{ai10eq50} schematically as
\ea
&\m_n^{\{0,\ldots,N\}\sm\{i'\};x_0,\ldots,x_n}=
\label{ai10eq53}\\
&\bigl(\text{terms in r.h.s.\ of \eq{ai10eq50} with $\m_n^{\{0,\ldots,N\}\sm\{i'\};x_0,\ldots,x_n}$ and at least one $\m_0^{I_b;z_0}$}\bigr)
\nonumber\\
&+\bigl(\text{terms in r.h.s.\ of \eq{ai10eq50} not involving $\m_n^{\{0,\ldots,N\}\sm\{i'\};x_0,\ldots,x_n}$}\bigr).
\nonumber
\ea

In the proof of $(*)_{i'}^N$ in \S\ref{ai1091}, we first chose data $\m_k^{\{0,\ldots,N\};x_0,\ldots,x_k}$ for all $k,x_0,\ldots,x_k$ subject to the obvious necessary conditions, but otherwise arbitrary ($\m_k^{\{0,\ldots,N\};x_0,\ldots,x_k}=0$ will do), and then we defined data $\m_n^{\{0,\ldots,N\}\sm\{i'\};x_0,\ldots,x_n}$ by induction on $n=0,1,\ldots$ by \eq{ai10eq50}, noting by Lemma \ref{ai10lem1} that the right hand side of \eq{ai10eq50} only depended on $\m_l^{\{0,\ldots,N\}\sm\{i'\};w_0,\ldots,w_l}$ for $l<n$.

For the filtered case, we modify this method as follows. In the inductive step in which we define $\m_n^{\{0,\ldots,N\}\sm\{i'\};x_0,\ldots,x_n}$, using \eq{ai10eq53} we first define a sequence
\ea
&\bigl(\m_{n\,(k)}^{\{0,\ldots,N\}\sm\{i'\};x_0,\ldots,x_n}\bigr)_{k=0}^\iy,\quad \text{where $\m_{n\,(k)}^{\{0,\ldots,N\}\sm\{i'\};x_0,\ldots,x_n}$ is a map as in \eq{ai8eq3},}
\nonumber\\
&\text{by $\m_{n\,(0)}^{\{0,\ldots,N\}\sm\{i'\};x_0,\ldots,x_n}=0$, and by induction for $k=0,1,\ldots,$}
\nonumber\\
&\m_{n\,(k+1)}^{\{0,\ldots,N\}\sm\{i'\};x_0,\ldots,x_n}=
\label{ai10eq54}\\
&\bigl(\text{terms in r.h.s.\ of \eq{ai10eq50} with $\m_{n\, (k)}^{\{0,\ldots,N\}\sm\{i'\};x_0,\ldots,x_n}$ and at least one $\m_0^{I_b;z_0}$}\bigr)
\nonumber\\
&+\bigl(\text{terms in r.h.s.\ of \eq{ai10eq50} not involving $\m_n^{\{0,\ldots,N\}\sm\{i'\};x_0,\ldots,x_n}$}\bigr).
\nonumber
\ea

Subtracting \eq{ai10eq54} for $k+1$ and $k$ shows that
\begin{align*}
&\m_{n\,(k+2)}^{\{0,\ldots,N\}\sm\{i'\};x_0,\ldots,x_n}-\m_{n\,(k+1)}^{\{0,\ldots,N\}\sm\{i'\};x_0,\ldots,x_n}=\\
&\bigl(\text{terms in r.h.s.\ of \eq{ai10eq50} with $\m_{n\,(k+1)}^{\{0,\ldots,N\}\sm\{i'\};x_0,\ldots,x_n}-\m_{n\, (k)}^{\{0,\ldots,N\}\sm\{i'\};x_0,\ldots,x_n}$}\\
&\text{in place of $\m_{n}^{\{0,\ldots,N\}\sm\{i'\};x_0,\ldots,x_n}$, and at least one $\m_0^{I_b;z_0}$}\bigr).
\end{align*}
Since $\m_{n\,(k)}^{\{0,\ldots,N\}\sm\{i'\};x_0,\ldots,x_n}\in F^0(\cdots)$ for $k=0,1$, and $\m_0^{I_b;z_0}(1)\in F^\ep(\cdots)$ from above, we see by induction on $k=0,1,\ldots$ that
\e
\m_{n\,(k+1)}^{\{0,\ldots,N\}\sm\{i'\};x_0,\ldots,x_n}-\m_{n\,(k)}^{\{0,\ldots,N\}\sm\{i'\};x_0,\ldots,x_n}\in F^{k\ep}(\cdots).
\label{ai10eq55}
\e
As $\ep\ssu 0$ we have $k\ep\ra\iy$ in $\Ga$ as $k\ra\iy$. Thus $\bigl(\m_{n\,(k)}^{\{0,\ldots,N\}\sm\{i'\};x_0,\ldots,x_n}\bigr)_{k=0}^\iy$ is a Cauchy sequence in the filtered sense, and there is a unique $\m_n^{\{0,\ldots,N\}\sm\{i'\};x_0,\ldots,x_n}$ with $\m_{n\,(k)}^{\{0,\ldots,N\}\sm\{i'\};x_0,\ldots,x_n}\ra\m_n^{\{0,\ldots,N\}\sm\{i'\};x_0,\ldots,x_n}$ as $k\ra\iy$. Taking the limit $k\ra\iy$ in \eq{ai10eq54} shows that $\m_n^{\{0,\ldots,N\}\sm\{i'\};x_0,\ldots,x_n}$ satisfies \eq{ai10eq50}. It is easy to see that this is the {\it unique\/} $\m_n^{\{0,\ldots,N\}\sm\{i'\};x_0,\ldots,x_n}$ satisfying \eq{ai10eq50}, noting that both sides of \eq{ai10eq50} now involve $\m_n^{\{0,\ldots,N\}\sm\{i'\};x_0,\ldots,x_n}$ as above.

Hence, by induction on $n$ we define $\m_n^{\{0,\ldots,N\}\sm\{i'\};x_0,\ldots,x_n}$ satisfying \eq{ai10eq50}, and therefore \eq{ai10eq29}, for all $n\ge 0$ and~$x_0,\ldots,x_n\in\sO^0$.

In \S\ref{ai1091} we used \eq{ai10eq51} to show that $\m_0^{\{0,\ldots,N\}\sm\{i'\};x_0}=0$ if $K\ge N$. For $\AfCat^\iy_{R,K,L}$, the same argument shows that $\m_0^{\{0,\ldots,N\}\sm\{i'\};x_0}=0$ if $K\ge N$ and $\m_0^{\{0,\ldots,N\}\sm\{i'\};x_0}(1)\in F^{\ssu 0}(\cdots)$ if $K<N\le L$, as in Definition~\ref{ai8def6}(a)(iv).

Next, in \S\ref{ai1091}, we proved that \eq{ai7eq25} holds for $I=\{0,\ldots,N\}\sm\{i'\}$ by induction on $n=0,1,\ldots$ in \eq{ai10eq52}. However, as explained in Remark \ref{ai10rem1}, equation \eq{ai10eq52} used \eq{ai7eq27} for $I=\{0,\ldots,N\}$ and $n=\ti n$, and we applied Lemma \ref{ai10lem1} to show this only depends on \eq{ai7eq25} for $n<\ti n$. Since Lemma \ref{ai10lem1} is false in the filtered case, this argument needs changing.

We resolve this issue by introducing an extra induction over $k=0,1,\ldots$ in the proof of \eq{ai10eq52}, with inductive hypothesis that \eq{ai7eq25} holds for $I=\{0,\ldots,N\}\sm\{i'\}$ and $x_0,\ldots,x_{\ti n}$ modulo terms in $F^{k\ep}(\cdots)$. At the point in the proof in \S\ref{ai107} of the filtered case of \eq{ai7eq27} for $I=\{0,\ldots,N\}$ and $x_0,\ldots,x_{\ti n}$ that involves \eq{ai7eq25} for $I=\{0,\ldots,N\}\sm\{i'\}$ and $x_0,\ldots,x_{\ti n}$, as in Lemma \ref{ai10lem2} it also involves at least one factor $\m_0^{I_b;z_0}$ with $\md{I_b}=2$ or 3, so that $\m_0^{I_b;z_0}(1)\in F^\ep\Hom_{\sA^*}(*,*)$ as above. Therefore the inductive hypothesis implies that \eq{ai7eq27} for $I=\{0,\ldots,N\}$ and $x_0,\ldots,x_{\ti n}$ holds modulo terms in $F^{(k+1)\ep}(\cdots)$. So following \eq{ai10eq52} we show that \eq{ai7eq25} holds for $I=\{0,\ldots,N\}\sm\{i'\}$ and $x_0,\ldots,x_{\ti n}$ modulo terms in $F^{(k+1)\ep}(\cdots)$, proving the inductive step. Then taking the limit $k\ra\iy$ shows that \eq{ai7eq25} holds for $I=\{0,\ldots,N\}\sm\{i'\}$ and~$x_0,\ldots,x_{\ti n}$.

This proves that $\AfCat^\iy_{R,K,L}$ is an $\iy$-category.

Next consider how to modify the proof above to show $\AwCat^\iy_{R,K,L}$ is an $\iy$-category. The important difference is that, as in Definitions \ref{ai8def7} and \ref{ai8def8}, the maps $\m_k^{I;w_0,\ldots,w_k}$ need not preserve filtrations for $\md{I}\ge 4$ and $k\ge 1$, or if $\md{I}>L$ and $k=0$, but instead map $F^\ga(\cdots)\ra F^{\ga+\de}(\cdots)$ for all $\ga\in\Ga$. Note that in Definition \ref{ai8def4}(iv) we assumed, by hand, that $\m_k^{I;w_0,\ldots,w_k}$ preserves filtrations if $\md{I}=3$ and $k\ge 1$. This is necessary for the next part of the proof.

Because of this, in the argument with $(\m_{n\,(k)}^{\{0,\ldots,N\}\sm\{i'\};x_0,\ldots,x_n})_{k=0}^\iy$ above, we do not know that $\m_{n\,(k)}^{\{0,\ldots,N\}\sm\{i'\};x_0,\ldots,x_n}\in F^0(\cdots)$ for $k=0,1$, but only that $\m_{n\,(k)}^{\{0,\ldots,N\}\sm\{i'\};x_0,\ldots,x_n}\in F^\de(\cdots)$ for $k=0,1$ and some $\de\in\Ga$. This only affects the third line of \eq{ai10eq53}: the second line of \eq{ai10eq53} involves $\m_n^{\{0,\ldots,N\}\sm\{i'\};x_0,\ldots,x_n}$ together with $\m_k^{I;w_0,\ldots,w_k}$ for $\md{I}\le 3$, which are filtration-preserving. 

Hence \eq{ai10eq55} should be replaced by
\begin{equation*}
\m_{n\,(k+1)}^{\{0,\ldots,N\}\sm\{i'\};x_0,\ldots,x_n}-\m_{n\,(k)}^{\{0,\ldots,N\}\sm\{i'\};x_0,\ldots,x_n}\in F^{\de+k\ep}(\cdots).
\end{equation*}
This still implies that $(\m_{n\,(k)}^{\{0,\ldots,N\}\sm\{i'\};x_0,\ldots,x_n})_{k=0}^\iy$ is Cauchy, as we want.

Similarly, in the part of the proof involving \eq{ai10eq52}, the inductive hypothesis should be that \eq{ai7eq25} holds for $I=\{0,\ldots,N\}\sm\{i'\}$ and $x_0,\ldots,x_{\ti n}$ modulo terms in $F^{\de+k\ep}(\cdots)$ for $k=0,1,\ldots,$ for some sufficiently negative $\de\in\Ga$. The rest of the proof goes through without change.

We prove that $\dAfCat^{\un,\iy}_{R,K,L},\dAwCat^{\un,\iy}_{R,K,L}$ are $\iy$-categories by modifying the proofs above that $\AfCat^\iy_{R,K,L},\AwCat^\iy_{R,K,L}$ are $\iy$-categories to include unitality, by the method of~\S\ref{ai1093}.

Theorem \ref{ai8thm5}, proved in \S\ref{ai108}, shows that $\AfCat^{\un,\iy}_{R,K,L},\AwCat^{\un,\iy}_{R,K,L}$ are isomorphic to $\dAfCat^{\un,\iy}_{R,K,L},\ab\dAwCat^{\un,\iy}_{R,K,L}$ as simplicial sets. So, as in \S\ref{ai1094}, $\AfCat^{\un,\iy}_{R,K,L},\AwCat^{\un,\iy}_{R,K,L}$ are also $\iy$-categories. Theorem \ref{ai8thm6} follows.

\section{Proofs of results in Chapter \ref{ai9}}
\label{ai11}

\subsection{Proof of Theorem \ref{ai9thm1}}
\label{ai111}

For Theorem \ref{ai9thm1}(a), we begin with two long examples. The first computes left mapping spaces in $\hACat^{\un,\iy}_R$ and~$\ACat^{\un,\iy}_R$.

\begin{ex}
\label{ai11ex1}
Let $R$ be a commutative ring, and let $\sA,\sB$ be unital $A_\iy$-categories over $R$, thought of as objects in the $\iy$-category $\ACat^{\un,\iy}_R$ from Definition \ref{ai7def4}. We will explicitly compute the simplicial set of left morphisms $\Hom_{\ACat^{\un,\iy}_R}^L(\sA,\sB)$ from Definition \ref{ai2def8}, which is a Kan complex, and one model for the $\iy$-category mapping space~$\Hom_{\ACat^{\un,\iy}_R}(\sA,\sB)$.

In fact we will first compute the simplicial set $\Hom_{\hACat^{\un,\iy}_R}^L(\sA,\sB)$ for $\hACat^{\un,\iy}_R$ in Definition \ref{ai7def3}, and then describe $\Hom_{\ACat^{\un,\iy}_R}^L(\sA,\sB)$ as a simplicial subset of~$\Hom_{\hACat^{\un,\iy}_R}^L(\sA,\sB)$.

By Definition \ref{ai2def8}, an $N$-simplex in $\Hom_{\hACat^{\un,\iy}_R}^L(\sA,\sB)$ is a morphism $z:\De^{N+1}\ra\hACat^{\un,\iy}_R$ such that $z\vert_{\De^{\{0\}}}$ is the 0-simplex $\sA$, and $z\vert_{\De^{\{1,\ldots,N+1\}}}$ is a constant $N$-simplex at $\sB$. From Definition \ref{ai7def3}, $z$ is equivalent to data $\bigl((\sO^i)_{i=0}^{N+1},(f^{ij})_{0\le i\le j\le N+1},(\sH^i)_{i=0}^{N+1},(\m^I)_{\es\ne I\subseteq\{0,\ldots,N+1\}}\bigr)$ satisfying Definition \ref{ai7def3}(a),
such that
\begin{itemize}
\setlength{\itemsep}{0pt}
\setlength{\parsep}{0pt}
\item[(i)] $\sO^0=\Obj\sA$ and $\sO^i=\Obj\sB$ for~$i=1,\ldots,N+1$.
\item[(ii)] $f^{00}=\id_{\Obj_\sA}$, and $f^{ij}=\id_{\Obj_\sB}$ for~$1\le i\le j\le N+1$.
\item[(iii)] $\sH^0(x_0,x_1)=\Hom_\sA(x_0,x_1)$ for all $x_0,x_1\in\Obj_\sA$, and $\sH^i(y_0,y_1)\ab=\ab\Hom_\sB(y_0,y_1)$ for all $1\le i\le N+1$ and~$y_0,y_1\in\Obj_\sB$.
\item[(iv)] $\m_k^{\{0\};x_0,\ldots,x_k}=\m_k^{\sA;x_0,\ldots,x_k}$ for all $k\ge 0$ and $x_0,\ldots,x_k$ in $\Obj\sA$, and $\m_k^{\{i\};y_0,\ldots,y_k}=\m_k^{\sB;y_0,\ldots,y_k}$ for all $k\ge 0$, $1\le i\le N+1$ and $y_0,\ldots,y_k$ in~$\Obj\sB$.
\item[(v)] $\m_1^{\{i,j\};y_0,y_1}=\id_{\Hom_\sB(y_0,y_1)}$ for all $1\le i<j\le N+1$ and $y_0,y_1\in\Obj\sB$.
\item[(vi)] $\m_0^{\{i,j,k\};y_0}:\la\mapsto\boo^\sB_{y_0}$ for all $1\le i<j<k\le N+1$ and $y_0\in\Obj\sB$.
\item[(vii)] If $\es\ne I\subseteq\{1,\ldots,N+1\}$, $k\ge 0$ and $y_0,\ldots,y_k\in\Obj\sB$ then in all cases not covered in (iv)--(vi), we have~$\m_k^{I;y_0,\ldots,y_k}=0$.
\item[(viii)] Observe that (iv)--(vii) prescribe $\m_k^{I;x_0,\ldots,x_k}$ in all cases {\it except\/} when $I=\{0,j_0,\ldots,j_k\}$ for $k\ge 0$ and $1\le j_0<j_1<\cdots<j_k\le N+1$. We impose no conditions on these, except that they form part of an $(N+1)$-simplex in $\hACat^{\un,\iy}_R$. We generally write $j_a=i_a+1$, so that $0\le i_0<i_1<\cdots<i_k\le N$.
\end{itemize}

Let $\fF^i:\sA\ra\sB$ for $i=0,\ldots,N$ be the unital $A_\iy$-functor on edge $[0,i+1]$, so that $\fF_\Obj^i=f^{0(i+1)}$ and
\e
\fF_k^{i;x_0,\ldots,x_k}=\m_k^{\{0,i+1\};x_0,\ldots,x_k}.
\label{ai11eq1}
\e
For all $0\le i_0<\cdots<i_k\le N$ for $k\ge 1$ and $l\ge 0$, $x_0,\ldots,x_l\in\Obj\sA$, define
\e
\fT_{k;l}^{i_0,\ldots,i_k;x_0,\ldots,x_l}=\m_l^{\{0,i_0+1,i_1+1,\ldots,i_k+1\};x_0,\ldots,x_l}.
\label{ai11eq2}
\e
This is an $R$-linear map graded of degree $1-k-l$
\ea
\fT_{k;l}^{i_0,\ldots,i_k;x_0,\ldots,x_l}&:\Hom_\sA(x_0,x_1)\ot_R\Hom_\sA(x_1,x_2)\ot_R\cdots \ot_R\Hom_\sA(x_{l-1},x_l)
\nonumber\\
&\longra\Hom_\sB(\fF^{i_0}(x_0),\fF^{i_k}(x_l)).
\label{ai11eq3}
\ea
Here the target of \eq{ai11eq3} is correct as by (ii) and $\fF_\Obj^i=f^{0(i+1)}$\begin{align*}
f^{(i_{k-1}+1)(i_k+1)}\!\ci\! f^{(i_{k-2}+1)(i_{k-1}+1)}\!\ci\!\cdots \!\ci\! f^{(i_0+1)(i_2+1)}\!\ci\! f^{0(i_0+1)}(x_0)&=\fF^{i_0}(x_0),\\
f^{0(i_k+1)}(x_l)&=\fF^{i_k}(x_l).
\end{align*}
Write $\fT_k^{i_0,\ldots,i_k}=(\fT_{k;l}^{i_0,\ldots,i_k;x_0,\ldots,x_l})_{l\ge 0}^{x_0,\ldots,x_l\in\Obj\sA}$.

Then the data $\bigl((\fF^i)_{i=0}^N,(\fT_k^{i_0,\ldots,i_k})_{k\ge 1}^{0\le i_0<\cdots<i_k\le N}\bigr)$ is all the data in $\bigl((\sO^i)_{i=0}^{N+1},\ab(f^{ij})_{0\le i\le j\le N+1},(\sH^i)_{i=0}^{N+1},(\m^I)_{\es\ne I\subseteq\{0,\ldots,N+1\}}\bigr)$ not already determined by (i)--(vii), so it is equivalent to $z$. Taking the $\fF^i$ to be unital $A_\iy$-functors, the remaining conditions in Definition \ref{ai7def3}(a)(i)--(v) which are not automatically satisfied are equations \eq{ai7eq16} and \eq{ai7eq19} for $\es\ne I\subseteq \{0,\ldots,N+1\}$ with $0\in I$ and $\md{I}\ge 3$. Writing $I=\{0,i_0+1,i_1+1,\ldots,i_k+1\}$ for $0\le i_0<\cdots<i_k\le N$ and $k\ge 1$, for objects $x_0,\ldots,x_n$ in $\sA$, equation \eq{ai7eq16} becomes
\ea
&0=\sum_{\begin{subarray}{l} l=1,\ldots,n \\
y_0,\ldots,y_l\in\Obj\sA \end{subarray}}\fT_{k;l}^{i_0,\ldots,i_k;y_0,\ldots,y_l}\ci\d_{n,l;y_0,\ldots,y_l}^{\sA;x_0,\ldots,x_n}
\label{ai11eq4}\\
&+(-1)^{k+1}\sum_{\begin{subarray}{l} q=1,\ldots,k+n,\\
y_0,\ldots,y_q\in\Obj\sB \end{subarray}}\m_q^{\sB;y_0,\ldots,y_q}\ci\d_{n,q;y_0,\ldots,y_q}^{\{0,i_0+1,i_1+1,\ldots,i_k+1\};x_0,\ldots,x_n},
\nonumber
\ea
where $\d_{n,l;y_0,\ldots,y_l}^{\{0,i_0+1,i_1+1,\ldots,i_k+1\};x_0,\ldots,x_n}$ should be written in terms of the $\fF_k^{i;0,\ldots,w_k}$ and $\fT_{k;l}^{J;z_0,\ldots,z_l}$ following \eq{ai7eq7}--\eq{ai7eq9}. Also \eq{ai7eq19} becomes 
\ea
\begin{aligned}
&\fT_{k;l+1}^{i_0,\ldots,i_k;x_0,\ldots,x_j,x_j,\ldots,x_l}(a_1\ot\cdots\ot a_j\ot \boo^{i_0}_{x_j}\ot a_{j+1}\ot \cdots \ot a_l)=0\\
&\text{for all $k\ge 1$, $0\le i_0<\cdots<i_k\le N$, $0\le j\le l$,}\\
&\text{$x_0,\ldots,x_l\in\Obj\sA$, and $a_m\in\Hom_\sA(x_{m-1},x_m)$, $m=1,\ldots,l$.}
\end{aligned}
\label{ai11eq5}
\ea

We will show that we can write \eq{ai11eq4} more explicitly as, for all $1\le k\le N$ and $0\ab\le\ab i_0\ab<\ab i_1\ab<\ab\cdots\ab<\ab i_k\ab\le N,$ $n\ge 0$ and objects $x_0,\ldots,x_n$ in~$\sA$
\ea
&0=\sum_{\begin{subarray}{l} p,r\ge 0,\; q\ge 1:\\ p+q+r=n \end{subarray}}\begin{aligned}[t]&(-1)^{p+qr}\fT^{i_0,\ldots,i_k;x_0,\ldots,x_p,x_{p+q},\ldots,x_n}_{k;p+r+1}\ci{}\\
&\bigl(\id_{\Hom_\sA(x_0,x_1)}\ot\cdots\ot\id_{\Hom_\sA(x_{p-1},x_p)}\\
&{}\ot\m_q^{\sA;x_p,\ldots,x_{p+q}}\ot \id_{\Hom_\sA(x_{p+q},x_{p+q+1})}\ot\cdots\ot\id_{\Hom_\sA(x_{n-1},x_n)}\bigr) \end{aligned}
\allowdisplaybreaks
\nonumber\\
&+\sum_{\begin{subarray}{l} 
1\le p\le q,\\
0=l_0<l_1<\\
\cdots <l_p=k, \\
0=n_0\le n_1\le\\
\cdots\le n_q=n,\\
1\le a_1<a_2< \\
\cdots<a_p\le q: \\
n_j>n_{j-1}, \\ 
j\ne a_i
\end{subarray}}
\begin{aligned}[t]&(-1)^{\sum_{i=1}^p(l_i-l_{i-1}+1)(p-i)}\cdot{}\\
&(-1)^{p+q}\cdot
(-1)^{\sum_{i=1}^p\sum_{j=1}^{a_i-1}(n_j-n_{j-1}+1)}
\cdot(-1)^{\sum_{i=1}^q(n_i-n_{i-1}+1)(q-i)}\cdot{}\\
&\m_q^{\begin{subarray}{l} \sB;\fF^{l_0}(x_{n_0}),\fF^{l_0}(x_{n_1}),\fF^{l_0}(x_{n_2}),\ldots,\fF^{l_0}(x_{n_{a_1-1}}),\fF^{l_1}(x_{n_{a_1}}),\fF^{l_1}(x_{n_{a_1+1}}),\\
\ldots,\fF^{l_1}(x_{n_{a_2-1}}),\ldots,\fF^{l_p}(x_{n_{a_p}}),\fF^{l_p}(x_{n_{a_p+1}}),\ldots,\fF^{l_p}(x_{n_q})
\end{subarray}}\ci\\
&\bigl(\fF_{n_1-n_0}^{l_0;x_{n_0},\ldots,x_{n_1}}\ot\fF_{n_2-n_1}^{l_0;x_{n_1},\ldots,x_{n_2}}\ot\cdots \\
&\ot\fF_{n_{a_1-1}-n_{a_1-2}}^{l_0;x_{n_{a_1-2}},\ldots,x_{n_{a_1-1}}}\ot
\fT_{l_1-l_0;n_{a_1}-n_{a_1-1}}^{i_{l_0},\ldots,i_{l_1};x_{n_{a_1-1}},\ldots,x_{n_{a_1}}}\\
&\ot\fF_{n_{a_1+1}-n_{a_1}}^{l_1;x_{n_{a_1}},\ldots,x_{n_{a_1+1}}}\!\ot\!\cdots\!\ot\!\fF_{n_{a_2-1}-n_{a_2-2}}^{l_1;x_{n_{a_2-2}},\ldots,x_{n_{a_2-1}}} \\
&{}\ot
\fT_{l_2-l_1;n_{a_2}-n_{a_2-1}}^{i_{l_1},\ldots,i_{l_2};x_{n_{a_2-1}},\ldots,x_{n_{a_2}}}\ot \cdots \ot\fT_{l_p-l_{p-1};n_{a_p}-n_{a_p-1}}^{i_{l_{p-1}},\ldots,i_{l_p};x_{n_{a_p-1}},\ldots,x_{n_{a_p}}} {}\\
&\ot\fF_{n_{a_p+1}-n_{a_p}}^{l_p;x_{n_{a_p}},\ldots,x_{n_{a_p+1}}}\ot\cdots\ot\fF_{n_q-n_{q-1}}^{l_p;x_{n_{q-1}},\ldots,x_{n_q}}\bigr)
\end{aligned}
\nonumber\\
&+\sum_{0<j<k}(-1)^j\fT_{k-1;n}^{i_0,i_1,\ldots,i_{j-1},i_{j+1},\ldots,i_k;x_0,\ldots,x_n}.
\label{ai11eq6}
\ea

Here the first sum in \eq{ai11eq6} agrees with the first sum in \eq{ai11eq4} by \eq{ai7eq8}. We will show by expanding out $\d_{n,q;y_0,\ldots,y_q}^{\{0,i_0+1,i_1+1,\ldots,i_k+1\};x_0,\ldots,x_n}$ in \eq{ai11eq4} using \eq{ai7eq9} that the second sum in \eq{ai11eq4} is the total of the second and third sums in \eq{ai11eq6}, where the third sum in \eq{ai11eq6} is a special case which arises when~$q=2$.

Now $\d_{n,q;y_0,\ldots,y_q}^{\{0,i_0+1,i_1+1,\ldots,i_k+1\};x_0,\ldots,x_n}$ is defined in terms of the $\m_k^{I;x_0,\ldots,x_k}$ by \eq{ai7eq9}, since $q\ge 1$, as a sum of terms $\pm C(\bs\Ga_1)^{x_{n_0},\ldots,x_{n_1}}\ot C(\bs\Ga_2)^{x_{n_1},\ldots,x_{n_2}}\ot\ab\cdots\ab\ot\ab C(\bs\Ga_q)^{x_{n_{q-1}},\ldots,x_{n_q}}$. For most $I,k$ we have $\m_k^{I;x_0,\ldots,x_k}=0$ by (vii) above. The others are given by (v),(vi), \eq{ai11eq1} and \eq{ai11eq2}. We first describe which labelled trees $\bs\Ga_a$ can give nonzero terms $C(\bs\Ga_a)^{x_{n_{a-1}},\ldots,x_{n_a}}$. There are three cases:
\begin{itemize}
\setlength{\itemsep}{0pt}
\setlength{\parsep}{0pt}
\item[(a)] For $0\le l_0<l_1<\cdots<l_r\le N$, $0<s\le r$, and $0\le n_{a-1}\le n_a\le n$, $\bs\Ga_a$ is
\end{itemize}
\begin{align*}
\xymatrix@C=13pt{
{\begin{subarray}{l} \!\!\!\!\!\!\!\!\ts n_{a-1}\vphantom{\{}\; \ci \\[-3pt] \ts \kern 1.5em\vdots \\[2pt] \ts n_a\vphantom{\{}\; \ci\end{subarray}}
 \ar@<3ex>[r]  \ar@<1ex>[r] \ar@<-1ex>[r] \ar@<-3ex>[r] &
\overset{\!\!\{0,i_{l_0}+1,i_{l_1}+1,\ldots,i_{l_s}+1\}\!\!\!\!}{\bu}
\ar[r] & \overset{\!\!\!\!\{i_{l_s}+1,i_{l_{s+1}}+1\}\!\!\!\!}{\bu}
\ar[r] & \cdots \ar[r] & \overset{\!\!\!\!\{i_{l_{r-1}}+1,i_{l_r}+1\}\!\!\!\!}{\bu} \ar[r] & \tr & \!\!\!\!\!\!\!\!\bs\Ga_a. }
\end{align*}
\begin{itemize}
\setlength{\itemsep}{0pt}
\setlength{\parsep}{0pt}
\item[] Then $C(\bs\Ga_a)^{x_{n_{a-1}},\ldots,x_{n_a}}=\fT_{r-s;n_a-n_{a-1}}^{i_{l_0},\ldots,i_{l_s};x_{n_{a-1}},\ldots,x_{n_a}}$. This is because the vertices $\{i_{l_{b-1}}\!+\!1,i_{l_b}\!+\!1\} \; \bu$ contribute maps $\m_1^{\{i_{l_{b-1}}+1,i_{l_b}+1\};y,y'}\!=\!\id_{\Hom_\sB(y,y')}$ by (v) above, so $C(\bs\Ga_a)^{x_{n_{a-1}},\ldots,x_{n_a}}=\m_{n_a-n_{a-1}}^{\{0,i_{l_0}+1,i_{l_1}+1,\ldots,i_{l_s}+1\};x_{n_{a-1}},\ldots,x_{n_a}}$.
\item[(b)] For $0\le l_0<l_1<\cdots<l_r\le N$, $r\ge 0$, and $0\le n_{a-1}<n_a\le n$, $\bs\Ga_a$ is
\end{itemize}
\begin{equation*}
\xymatrix@C=17pt{
{\begin{subarray}{l} \!\!\!\!\!\!\!\!\ts n_{a-1}\vphantom{\{}\; \ci \\[-3pt] \ts \kern 1.5em\vdots \\[2pt] \ts n_a\vphantom{\{}\; \ci\end{subarray}}
 \ar@<3ex>[r]  \ar@<1ex>[r] \ar@<-1ex>[r] \ar@<-3ex>[r] &
\overset{\!\!\{0,i_{l_0}+1\}\!\!\!\!}{\bu}
\ar[r] & \overset{\!\!\!\!\{i_{l_0}+1,i_{l_1}+1\}\!\!\!\!}{\bu}
\ar[r] & \cdots \ar[r] & 
\overset{\!\!\!\!\{i_{l_{r-1}}+1,i_{l_r}+1\}\!\!\!\!}{\bu}  \ar[r] & \tr & \!\!\!\!\!\!\!\!\bs\Ga_a. }
\end{equation*}
\begin{itemize}
\setlength{\itemsep}{0pt}
\setlength{\parsep}{0pt}
\item[] Then $C(\bs\Ga_a)^{x_{n_{a-1}},\ldots,x_{n_a}}=\fF_{n_a-n_{a-1}}^{i_{l_0};x_{n_{a-1}},\ldots,x_{n_a}}$. Note that (a) required that $s>0$, but this is in effect the case $s=0$. We require $n_{a-1}<n_a$ as $\fF_0^{i_{l_0};x_{n_a}}=0$ when $n_{a-1}=n_a$.
\item[(c)] For $0\le l_0<l_1<\cdots<l_r\le N$, $r\ge 2$, $\Ga_a$ is
\end{itemize}
\begin{equation*}
\xymatrix@C=14pt{
n_a \ar@{..>}[r] & \overset{\!\!\!\!\!\!\{i_{l_0}+1,i_{l_1}+1,i_{l_2}+1\}\!\!\!\!\!\!}{\sq} 
\ar[r] & 
\overset{\!\!\!\!\!\!\{i_{l_2}+1,i_{l_3}+1\}\!\!\!\!\!\!}{\bu}  \ar[r] & \overset{\!\!\!\!\!\!\{i_{l_3}+1,i_{l_4}+1\}\!\!\!\!\!\!}{\bu}  \ar[r] &  \cdots \ar[r] & 
\overset{\!\!\!\!\!\!\{i_{l_{r-1}}+1,i_{l_r}+1\}\!\!\!\!\!\!}{\bu} \ar[r] & \tr & \!\!\!\!\!\!\!\!\!\!\!\!\bs\Ga_a. }
\end{equation*}
\begin{itemize}
\setlength{\itemsep}{0pt}
\setlength{\parsep}{0pt}
\item[] Here the initial `$n_a\dashra$' is not part of $\bs\Ga_a$, but indicates that the bud $\{i_{l_0}+1,i_{l_1}+1,i_{l_2}+1\}\; \sq$ is labelled by the object $x_{n_a}$ in Definition~\ref{ai7def1}. 

Then $C(\bs\Ga_a)^{x_{n_a}}=\boo^\sB_{\fF^{i_j}(x_{n_a})}$ for some $0\le j\le l_0$ depending on the prefix of $\{i_{l_0}+1,i_{l_1}+1,i_{l_2}+1\}\; \sq$. This is because $\m_0^{\{i_{l_0},i_{l_1},i_{l_2}\};\fF^{i_j}(x_{n_a})}=\boo^\sB_{\fF^{i_j}(x_{n_a})}$ by (vi) above, and then the $\m_1^{\{i_{l_{b-1}},i_{l_b}\};y,y}$ act on this as the identity by~(v).
\end{itemize}

These are the only possibilities to get $C(\bs\Ga_a)^{x_{n_{a-1}},\ldots,x_{n_a}}\ne 0$. We exclude other cases involving buds $\{i_{l_0}+1,i_{l_1}+1,i_{l_2}+1\}\; \sq$ as in (c), because when $\m_0^{\{i_{l_0}+1,i_{l_1}+1,i_{l_2}+1\};y}=\boo^\sB_y$ is input into any $\m_k^{I;w_0,\ldots,w_k}$ with $\md{I}\ge 2$ other than $\m_1^{\{i_{l_{b-1}},i_{l_b}\};y,y}$ as in (c), it gives zero by~\eq{ai7eq19}.

We can now show the second sum in \eq{ai11eq6} is $\sum_{q\ge 1,\; y_0,\ldots,y_q}\m_q^{\sB;y_0,\ldots,y_q}\ci -$ applied to the sum of all terms $\pm C(\bs\Ga_1)^{x_{n_0},\ldots,x_{n_1}}\ot\ab\cdots\ab\ot\ab C(\bs\Ga_q)^{x_{n_{q-1}},\ldots,x_{n_q}}$ in \eq{ai7eq9} for $\d_{n,q;y_0,\ldots,y_q}^{\{0,i_0+1,i_1+1,\ldots,i_k+1\};x_0,\ldots,x_n}$ with all $\bs\Ga_a$ of type (a) or (b), where $\bs\Ga_a$ is of type (a) for $a=a_j$ for some $j=1,\ldots,p$, and $\bs\Ga_a$ is of type (b) otherwise. 

This leaves the terms $\pm C(\bs\Ga_1)^{x_{n_0},\ldots,x_{n_1}}\ot\ab\cdots\ab\ot\ab C(\bs\Ga_q)^{x_{n_{q-1}},\ldots,x_{n_q}}$ with some $\bs\Ga_a$ of type (c). Then $C(\bs\Ga_a)^{x_{n_a}}=\boo^\sB_{\fF^{i_j}(x_{n_a})}$ is an input in $\m_q^{\sB;y_0,\ldots,y_q}$. By \eq{ai7eq17}--\eq{ai7eq19} this gives zero unless $q=2$. We can then show that the only possibility is for $\bs\Ga_1$ to be of type (c) and $\bs\Ga_2$ of type (a), and for some $0<j<k$ to be given by the diagram
\begin{equation*}
\xymatrix@C=14pt@R=-4pt{
0 \ar@{..>}[r] & \overset{\!\!\!\!\!\!\{i_{j-1}+1,i_{j}+1,i_{j+1}+1\}\!\!\!\!\!\!}{\sq} 
\ar[r] & 
\overset{\!\!\!\!\!\!\{i_{j+1}+1,i_{j+2}+1\}\!\!\!\!\!\!}{\bu}  \ar[r] &  \cdots \ar[r] & 
\overset{\!\!\!\!\!\!\{i_{k-1}+1,i_k+1\}\!\!\!\!\!\!}{\bu} \ar[r] & \tr & \!\!\!\!\!\!\!\!\!\!\!\!\bs\Ga_1 \\
&&& \!\!\!\!{\begin{subarray}{l} \ts 1\; \ci \\[-3pt] \ts \kern .91em\vdots \\[2pt] \ts n\; \ci\end{subarray}}
\ar[r] \ar@<1ex>[r]\ar@<-1ex>[r]
& \overset{\!\!\{0,i_0+1,\ldots,i_k+1\}\sm\{i_j+1\}\!\!}{\bu} \ar[r] & \tr & \!\!\!\!\!\!\!\!\!\!\!\!\bs\Ga_2. 
}
\end{equation*}
The sum over $0<j<k$ of this term composed with $\m_2^{\sB;y_0,y_1,y_2}$ gives the third sum in \eq{ai11eq6}, proving~\eq{ai11eq6}.

To summarize our work so far: an $N$-simplex in $\Hom_{\hACat^{\un,\iy}_R}^L(\sA,\sB)$ is equivalent to the choice of unital $A_\iy$-functors $\fF^i:\sA\ra\sB$ for $i=0,\ldots,N$, and for all $0\le i_0<\cdots<i_k\le N$ for $k\ge 1$ and $l\ge 0$, $x_0,\ldots,x_l\in\Obj\sA$, $R$-linear maps \eq{ai11eq3} graded of degree $1-k-l$,  satisfying equations \eq{ai11eq5} and \eq{ai11eq6}. The morphisms between simplices induced by morphisms $\al:[\bar N]\ra[N]$ in $\De$, as in \eq{ai7eq20}, are straightforward to write down, but we will not do this.

The simplicial subset $\Hom_{\ACat^{\un,\iy}_R}^L(\sA,\sB)$ in $\Hom_{\hACat^{\un,\iy}_R}^L(\sA,\sB)$ consists of those $N$-simplices as above coming from a morphism $z:\De^{N+1}\ra\ACat^{\un,\iy}_R\subset\hACat^{\un,\iy}_R$. By Definition \ref{ai7def4}, this imposes a condition on $\m^{\{i_0,i_1,i_2\}}$ for all $0\le i_0<i_1<i_2\le N+1$. If $i_0\ge 1$ this condition is automatic as $z\vert_{\De^{\{1,2,\ldots,N+1\}}}$ is a constant $N$-simplex at $\sB$. When $i_0=0$, replacing $i_0,i_1,i_2$ by $0,i_0+1,i_1+1$, it becomes the condition that for all $0\le i_0<i_1\le N$ and objects $x_0$ in $\sA$, we have a closed element $\fT_{1,0}^{i_0,i_1;x_0}(1)\in\Hom_\sB^0(\fF^{i_0}(x_0),\fF^{i_1}(x_0))$, and we require that
\begin{equation*}
[\fT_{1,0}^{i_0,i_1;x_0}(1)]:\fF^{i_0}(x_0)\longra\fF^{i_1}(x_0)
\end{equation*}
should be an isomorphism in the cohomology category~$H^0(\sB)$.

As $\ACat^{\un,\iy}_R$ is an $\iy$-category by Theorem \ref{ai7thm3}, $\Hom_{\ACat^{\un,\iy}_R}^L(\sA,\sB)$ is a Kan complex. It is the {\it maximal Kan complex\/} in $\Hom_{\hACat^{\un,\iy}_R}^L(\sA,\sB)$, since one can show that the condition on $\fT_{1,0}^{i_0,i_1;x_0}(1)$ above is necessary for any Kan subcomplex in $\Hom_{\hACat^{\un,\iy}_R}^L(\sA,\sB)$, as it is needed for the horn-filling conditions for $\La^2_0$ and $\La^2_2$ to hold, as in Definition~\ref{ai2def4}.
\end{ex}

Next we describe $\sN_{A_\iy}(\Fun_{A_\iy}^\un(\sA,\sB))$ as a simplicial set in a similar way.

\begin{ex}
\label{ai11ex2}
Let $R$ be a commutative ring, and let $\sA,\sB$ be unital $A_\iy$-categories over $R$. Definitions \ref{ai4def7}--\ref{ai4def8} define the $A_\iy$-category $\Fun_{A_\iy}^\un(\sA,\sB)$ of $A_\iy$-functors $\fF:\sA\ra\sB$, and Definition \ref{ai4def20} defines the simplicial nerve $\sN_{A_\iy}(\Fun_{A_\iy}^\un(\sA,\sB))$ as an $\iy$-category. Corollary \ref{ai4cor4} describes the $N$-sim\-pl\-ices in $\sN_{A_\iy}(\Fun_{A_\iy}^\un(\sA,\sB))$ and their morphisms explicitly.

By Corollary \ref{ai4cor4}(a), an $N$-simplex in $\sN_{A_\iy}(\Fun_{A_\iy}^\un(\sA,\sB))$ for $N\ge 0$ is equivalent to data $\bigl((\fF_i)_{i=0}^N,(\fT_k^{i_0,\ldots,i_k})^{\bs i}_{k\ge 1}\bigr),$ where $\fF_i:\sA\ra\sB$ is a unital $A_\iy$-functor for $i=0,\ldots,N$ and for all $k\ge 1$ and $0\ab\le\ab i_0\ab<\ab i_1\ab<\ab\cdots\ab<\ab i_k\ab\le N,$ we are given $\fT_k^{i_0,\ldots,i_k}\in\Hom_{\Fun_{A_\iy}^\un(\sA,\sB)}^{1-k}(\fF_{i_0},\fF_{i_k}),$ such that for all $1\le k\le N$ and $0\ab\le\ab i_0\ab<\ab i_1\ab<\ab\cdots\ab<\ab i_k\ab\le N,$ as in \eq{ai4eq38} we have
\ea
&0=\sum_{\begin{subarray}{l} p\ge 1,\; 0=l_0<l_1<\cdots <l_p=k\end{subarray}}
\begin{aligned}[t]&(-1)^{\sum_{i=1}^p(l_i-l_{i-1}+1)(p-i)}\m^{\Fun^\un_{A_\iy}(\sA,\sB);\fF_{l_0},\fF_{l_1},\ldots,\fF_{l_p}}_p\ci{}\\
&\bigl(\fT_{l_1-l_0}^{i_{l_0},\ldots,i_{l_1}}\ot \fT_{l_2-l_1}^{i_{l_1},\ldots,i_{l_2}}\ot\cdots\ot \fT_{l_p-l_{p-1}}^{i_{l_{p-1}},\ldots,i_{l_p}})\end{aligned} 
\nonumber\\
&+\sum_{0<j<k}(-1)^j\fT_{k-1}^{i_0,i_1,\ldots,i_{j-1},i_{j+1},\ldots,i_k}.
\label{ai11eq7}
\ea

By Definitions \ref{ai4def7}--\ref{ai4def8}, $\fT_k^{i_0,\ldots,i_k}=(\fT_{k;l}^{i_0,\ldots,i_k;x_0,\ldots,x_l})_{l\ge 0}^{x_0,\ldots,x_l\in\Obj\sA}$ is the data of an $R$-linear map $\fT_{k;l}^{i_0,\ldots,i_k;x_0,\ldots,x_l}$ as in \eq{ai11eq3}, graded of degree $1-k-l$, for all $l\ge 0$ and objects $x_0,\ldots,x_l$ in $\sA$. The condition that $\fT_k^{i_0,\ldots,i_k}$ lies in $\Fun^\un_{A_\iy}(\sA,\sB)\subseteq\Fun_{A_\iy}(\sA,\sB)$ is~\eq{ai11eq5}.

We can show using \eq{ai4eq15} and \eq{ai4eq17} that equation \eq{ai11eq7} is equivalent to \eq{ai11eq6} for all $n\ge 0$ and objects $x_0,\ldots,x_n$ in $\sA$. Here the second sum in \eq{ai11eq7} obviously becomes the third sum in \eq{ai11eq6}. Now \eq{ai4eq15} and \eq{ai4eq17} define $\m^{\Fun^\un_{A_\iy}(\sA,\sB)}_p$ for $p=1$ and $p>1$ respectively. The first sum in \eq{ai4eq15} is \eq{ai4eq17} with $p=1$, and the second sum is an extra term for $p=1$ only. The first sum in \eq{ai11eq6} is the second sum in \eq{ai4eq15} for the term $p=1$ in the first sum in \eq{ai11eq7} with $\fF^{i_0},\fF^{i_k}$ in place of $\fF,\fG$. The second sum in \eq{ai11eq6} is the first sum in \eq{ai4eq15} for $p=1$, and is \eq{ai4eq17} for $p>1$, with $\fF^{i_{l_0}},\ldots,\fF^{i_{l_p}}$ in place of~$\fF^0,\ldots,\fF^p$.

We have proved that an $N$-simplex of $\sN_{A_\iy}(\Fun_{A_\iy}^\un(\sA,\sB))$ corresponds exactly to an $N$-simplex of $\Hom_{\hACat^{\un,\iy}_R}^L(\sA,\sB)$, as described in Example \ref{ai11ex1}. It is easy to show that the morphisms between simplices induced by morphisms $\al:[\bar N]\ra[N]$ in $\De$, as in \eq{ai7eq20}, also correspond. This proves equation \eq{ai9eq2}. As $\Hom_{\ACat^{\un,\iy}_R}^L(\sA,\sB)\subset\Hom_{\hACat^{\un,\iy}_R}^L(\sA,\sB)$ is the maximal Kan complex by Example \ref{ai11ex1}, we deduce equation~\eq{ai9eq1}.
\end{ex}

For Theorem \ref{ai9thm1}(b), we can extend Examples \ref{ai11ex1}--\ref{ai11ex2} to the $\iy$-categories $\AfCat^{\un,\iy}_R,\ab\AwCat^{\un,\iy}_R$ of unital filtered $A_\iy$-categories over $R$ in \S\ref{ai81}, for $(\Ga,\pre)$ a po-monoid and $R$ a $\Ga$-filtered ring. The main differences are:
\begin{itemize}
\setlength{\itemsep}{0pt}
\setlength{\parsep}{0pt}
\item[(i)] In equation \eq{ai11eq6} in Example \ref{ai11ex1}, in the first sum we should allow $q\ge 0$, and in the second sum we should allow $n_j\ge n_{j-1}$ for all $j$. Then the second sum becomes an infinite sum, which converges in the filtered sense.
\item[(ii)] The objects (unital filtered $A_\iy$-categories over $R$) and 1-morphisms (unital filtered $A_\iy$-functors) in $\AfCat^{\un,\iy}_R,\ab\AwCat^{\un,\iy}_R$ are the same. But when $\md{I}\ge 3$, the maps $\m_k^{I;x_0,\ldots,x_k}$ are filtration-preserving in $\AfCat^{\un,\iy}_R$, and can decrease filtrations by a bounded amount in $\AwCat^{\un,\iy}_R$.

Because of this, $\fT_k^{i_0,\ldots,i_k}$ in Example \ref{ai11ex1} should be identified in Example \ref{ai11ex2} with an element of $\smash{\Hom_{\Fun_{A^\rfi_\iy}^\un(\sA,\sB)}^{1-k}(\fF_{i_0},\fF_{i_k})}$ for $\AwCat^{\un,\iy}_R$, but an element of
\begin{equation*}
F^0\Hom_{\Fun_{A^\rfi_\iy}^\un(\sA,\sB)}^{1-k}(\fF_{i_0},\fF_{i_k})=\Hom_{F^0\Fun_{A^\rfi_\iy}^\un(\sA,\sB)}^{1-k}(\fF_{i_0},\fF_{i_k})
\end{equation*}
for $\AfCat^{\un,\iy}_R$. This explains the right hand sides of \eq{ai9eq3}--\eq{ai9eq4}.
\end{itemize}

\subsection{Proof of Theorem \ref{ai9thm2}}
\label{ai112}

First we prove Theorem \ref{ai9thm2}(c). For unital $A_\iy$-categories $\sA,\sB$ over $R$ and $N\ge 0$ we have canonical bijections
\ea
&\Hom_{\ACat_R^\un}\bigl(\sA,\Fun_{A_\iy}^\un(\sC^\De_{[N]},\sB)\bigr)\cong 
\Hom_{\ACat_R^\un}^{(2)}(\sA,\sC^\De_{[N]},\sB)
\nonumber\\
&\qquad\cong\Hom_{\ACat_R^\un}^{(2)}(\sC^\De_{[N]},\sA,\sB)
\cong\Hom_{\ACat_R^\un}\bigl(\sC^\De_{[N]},\Fun_{A_\iy}^\un(\sA,\sB)\bigr)
\nonumber\\
&\qquad\cong \sN_{A_\iy}(\Fun_{A_\iy}^\un(\sA,\sB))([N]),
\label{ai11eq8}
\ea
where $\Hom_{\ACat_R^\un}^{(2)}(\sA_1,\sA_2,\sB)$ is the set of {\it unital\/ $A_\iy$-bifunctors\/} $\sA_1\t\sA_2\ra\sB$ defined and studied in Bespalov--Lyubashenko--Manzyuk \cite[\S 8]{BLM}. They give a notion of $A_\iy$-multifunctor $\sA_1\t\cdots\t\sA_k\ra\sB$, and give a multifunctor version of the $A_\iy$-category $\Fun_{A_\iy}^\un(\sA,\sB)$ in \S\ref{ai43}. We will use the canonical bijections $\Hom_{\ACat_R^\un}\bigl(\sA_1,\Fun_{A_\iy}^\un(\sA_2,\sB)\bigr)\cong\Hom_{\ACat_R^\un}^{(2)}(\sA_1,\sA_2,\sB)$ and $\Hom_{\ACat_R^\un}^{(2)}(\sA_1,\sA_2,\sB)\cong\Hom_{\ACat_R^\un}^{(2)}(\sA_2,\sA_1,\sB)$, together with functoriality of these in $\sA_1,\sA_2,\sB$. The last step of \eq{ai11eq8} uses Definition~\ref{ai4def20}.

We now claim there is a natural bijection `$\dashra$' in the commuting diagram:
\e
\begin{gathered}
\!\!\!\xymatrix@C=145pt@R=15pt{
*+[r]{\;\begin{subarray}{l}\ts F^\un_{S_\sB}(\sA)([N])= \\ \ts \Hom_{\ACat_R^\un}\bigl(\sA,S_\sB([N])\bigr)_{\vphantom{(}} \end{subarray}} \ar@{..>}[r]_(0.47)\cong \ar@<2ex>@{^{(}->}[d] & *+[l]{\sN_{A_{\iy_{\vphantom{(}}}}(\Fun_{A_\iy}^\un(\sA,\sB))^\Kan([N]) } \ar@<-3ex>@{^{(}->}[d]
\\
*+[r]{\Hom_{\ACat_R^\un}\bigl(\sA,\Fun_{A_\iy}^\un(\sC^\De_{[N]},\sB)\bigr)} \ar[r]^(0.55){\eq{ai11eq8}}_(0.55)\cong &
*+[l]{\sN_{A_\iy}(\Fun_{A_\iy}^\un(\sA,\sB))([N]).\!} }\!\!\!\!\!\!\!\!\!\!\!\!\!\!
\end{gathered}
\label{ai11eq9}
\e
Here $S_\sB([N])$ is the full $A_\iy$-subcategory of $\Fun_{A_\iy}^\un(\sC^\De_{[N]},\sB)$ with objects those $\fF:\sC^\De_{[N]}\ra\sB$ with $\fF(\om_{ij})$ a homotopy equivalence in $\sB$ for $0\le i<j\le N$, as in Definition \ref{ai9def1}, and $\sN_{A_\iy}(\Fun_{A_\iy}^\un(\sA,\sB))^\Kan$ is the maximal Kan complex in $\sN_{A_\iy}(\Fun_{A_\iy}^\un(\sA,\sB))$, which explains the vertical inclusion morphisms.

As in \S\ref{ai48}, $\sN_{A_\iy}(\Fun_{A_\iy}^\un(\sA,\sB))$ is an $\iy$-category. By Land \cite[Def.~1.2.24]{Land}, an $\iy$-category $\sC$ is called an $\iy$-{\it groupoid\/} if every morphism in $\sC$ is an equivalence. If $\sC$ is an $\iy$-category, as in \cite[Cor.~1.12.28]{Land} the maximal $\iy$-subgroupoid $\sC'\subseteq\sC$ is the simplicial subset of $\sC$ consisting of $N$-simplices in $\sC$ all of whose edges are equivalences. It is a nontrivial theorem \cite[Lem.~1.2.29 \& Cor.~2.1.12]{Land} that $\iy$-groupoids and Kan complexes are the same thing. Hence $\sN_{A_\iy}(\Fun_{A_\iy}^\un(\sA,\sB))^\Kan$ is the simplicial subset of $\sN_{A_\iy}(\Fun_{A_\iy}^\un(\sA,\sB))$ consisting of $N$-simplices all of whose edges are equivalences.

Now \eq{ai11eq8}--\eq{ai11eq9} are functorial under morphisms $f:[N]\ra[M]$ in $\De^\rop$ used to define $S_\sB:\De^\rop\ra\ACat^\un_R$ and to define $\sN_{A_\iy}(\Fun_{A_\iy}^\un(\sA,\sB))$ as a simplicial set. So to prove \eq{ai11eq9}, we must show that taken in the context of the whole simplicial sets, the restriction that $\fF(\om_{ij})$ a homotopy equivalence for $0\le i<j\le N$ in $S_\sB([N])$ corresponds to the condition that an $N$-simplex in $\sN_{A_\iy}(\Fun_{A_\iy}^\un(\sA,\sB))$ has edges equivalences, under the isomorphism~\eq{ai11eq8}.

By \S\ref{ai43}--\S\ref{ai44} and \S\ref{ai48}, an edge $e$ in the simplicial set $\sN_{A_\iy}(\Fun_{A_\iy}^\un(\sA,\sB))$ is equivalent to unital $A_\iy$-functors $\fF_0,\fF_1:\sA\ra\sB$, and a unital $A_\iy$-natural transformation $\fT:\fF_0\Ra\fF_1$. Also $e$ is an equivalence if and only if $\fT$ is an $A_\iy$-natural isomorphism. Note that as in  Proposition \ref{ai4prop3}, the condition of being an $A_\iy$-natural isomorphism only depends on a small part of the full data $(\fT_k^{x_0,\ldots,x_k})_{k,x_0,\ldots,x_k}$ in $\fT$: it requires that $[\fT_0^x(1)]:\fF_0(x)\ra \fF_1(x)$ should be an isomorphism in the cohomology category $H^0(\sB)$ for each object $x$ in~$\sA$. 

The effect of \eq{ai11eq8} is to translate $(\fF_0,\fF_1,\fT)$ into a unital $A_\iy$-functor $\fG:\sA\ra\Fun_{A_\iy}^\un(\sC^\De_{[1]},\sB)$. Under this translation, an object $x$ in $\sA$ is mapped to an object $\fH$ in $\Fun_{A_\iy}^\un(\sC^\De_{[1]},\sB)$ which is equivalent to the objects $\fF_0(x),\fF_1(x)$ in $\sB$, which are $\fH(0),\fH(1)$, and a morphism $\fT_0^x(1):\fF_0(x)\ra\fF_1(x)$ in $\sB$, which is $\fH(\om_{01})$. The condition that $\fG$ maps to $S_\sB([1])\subset \Fun_{A_\iy}^\un(\sC^\De_{[1]},\sB)$ is that for every object $x$ in $\sA$, this morphism $\fH(\om_{01})=\fT_0^x(1)$ should be an equivalence in $\sB$, that is, $[\fT_0^x(1)]$ should be an isomorphism in $H^0(\sB)$. This proves~\eq{ai11eq9}.

Equation \eq{ai11eq9} gives the isomorphism of simplicial sets \eq{ai9eq6} at $[N]$ in $\De^\rop$. As \eq{ai11eq9} is functorial in $[N]$ in $\De^\rop$, taken over all $N\ge 0$ it defines the isomorphism \eq{ai9eq6}. The identification \eq{ai9eq7} follows from the fact that the bijections \eq{ai9eq8} are functorial in $\sA$. This completes part~(c).

For (a), we must show that $(S_\sB,i_\sB)$ satisfies Definition \ref{ai2def18}(i),(ii). For (i), first let $N\ge 0$ and define morphisms $\pi_N:[N]\ra[0]$ and $\io_N:[0]\ra[N]$ in $\De$ by $\pi_N:i\mapsto 0$ and $\io_N:0\mapsto 0$. Definition \ref{ai4def20} now gives unital $A_\iy$-functors $\sC^\De_{\pi_N}:\sC^\De_{[N]}\ra \sC^\De_{[0]}$ and $\sC^\De_{\io_N}:\sC^\De_{[0]}\ra\sC^\De_{[N]}$, so composing gives $A_\iy$-functors
\e
\begin{split}
&\sC^\De_{\io_N\ci\pi_N}=\sC^\De_{\io_N}\ci\sC^\De_{\pi_N}:\sC^\De_{[N]}\longra\sC^\De_{[N]},\\
&\Id_{\sC^\De_{[0]}}=\sC^\De_{\id_{[0]}}=\sC^\De_{\pi_N\ci\io_N}=\sC^\De_{\pi_N}\ci\sC^\De_{\io_N}:\sC^\De_{[0]}\longra\sC^\De_{[0]}.
\end{split}
\label{ai11eq10}
\e
Define a unital $A_\iy$-natural transformation $\fK:\sC^\De_{\io_N\ci\pi_N}\Ra\Id$ of unital $A_\iy$-functors $\sC^\De_{[N]}\ra\sC^\De_{[N]}$ by
\e
\fK_0^i(1)=\om_{0i}\in\Hom_{\sC^\De_{\io_N}}(0,i),\; i=0,\ldots,N, \quad \fK_k^{x_0,\ldots,x_k}=0,\quad k>0.
\label{ai11eq11}
\e
Here $\fK_k^{x_0,\ldots,x_k}=0$ for $k>0$ is necessary for degree reasons as $\deg\fK_k^{x_0,\ldots,x_k}=-k$ and $\sC^\De_{[N]}$ is concentrated in degree 0. It is easy to check that \eq{ai11eq11} does define a unital $A_\iy$-natural transformation.

Applying Definitions \ref{ai4def9} and \ref{ai4def16}(a) to $\Fun_{A_\iy}^\un(-,\sB)$ and \eq{ai11eq10}--\eq{ai11eq11} gives unital $A_\iy$-functors
\begin{align*}
(\sC^\De_{\pi_N})^*\ci(\sC^\De_{\io_N})^*&:\Fun_{A_\iy}^\un(\sC^\De_{[N]},\sB)
\longra\Fun_{A_\iy}^\un(\sC^\De_{[N]},\sB),\\
\Id_{\Fun_{A_\iy}^\un(\sC^\De_{[0]},\sB)}=(\sC^\De_{\io_N})^*\ci(\sC^\De_{\pi_N})^*&:\Fun_{A_\iy}^\un(\sC^\De_{[0]},\sB)\longra\Fun_{A_\iy}^\un(\sC^\De_{[0]},\sB),
\end{align*}
and a unital $A_\iy$-natural transformation
\begin{equation*}
\fK^*:(\sC^\De_{\pi_N})^*\ci(\sC^\De_{\io_N})^*\Longra\Id_{\Fun_{A_\iy}^\un(\sC^\De_{[N]},\sB)}.
\end{equation*}
By Definition \ref{ai9def1} we can restrict these to the $A_\iy$-subcategories $S_\sB([i])\subseteq\Fun_{A_\iy}^\un(\sC^\De_{[i]},\sB)$ to get unital $A_\iy$-functors
\e
\begin{split}
S_\sB(\pi_N)\ci S_\sB(\io_N)&:S_\sB([N])\longra S_\sB([N]),\\
\Id_{S_\sB([0])}=S_\sB(\io_N)\ci S_\sB(\pi_N)&:S_\sB([0])\longra S_\sB([0]),
\end{split}
\label{ai11eq12}
\e
and a unital $A_\iy$-natural transformation
\e
\fK^*\vert_{S_\sB([N])}:S_\sB(\pi_N)\ci S_\sB(\io_N)\Longra\Id_{S_\sB([N])}.
\label{ai11eq13}
\e

We now claim that $\fK^*\vert_{S_\sB([N])}$ is an $A_\iy$-natural isomorphism. By Proposition \ref{ai4prop3}, it is enough to show that for each object $\fF:\sC^\De_{[N]}\ra\sB$ in $S_\sB([N])$,
\e
\bigl[\bigl(\fK^*\vert_{S_\sB([N])}\bigr){}_0^\fF\bigr]:S_\sB(\pi_N)\ci S_\sB(\io_N)(\fF)\ra \fF
\label{ai11eq14}
\e
is an isomorphism in the cohomology category $H^0(S_\sB([N]))$. But $\bigl(\fK^*\vert_{S_\sB([N])}\bigr){}_0^\fF$ acts by composition with $\fF(\om_{0i})$ in $\sB$ for $i=0,1,\ldots,N$, and the $\fF(\om_{0i})$ are homotopy equivalences in $\sB$ by definition of $S_\sB$, so that $[\fF(\om_{0i})]$ are isomorphisms in $H^0(\sB)$. From this we can show that \eq{ai11eq14} is an isomorphism. Thus \eq{ai11eq12}--\eq{ai11eq13} imply that $S_\sB(\pi_N),S_\sB(\io_N)$ are unital homotopy equivalences.

Now let $f:[N]\ra[M]$ be a morphism in $\De$. Then $\pi_N\ci f=\pi_M$, so applying $S_\sB$ gives $S_\sB(f)\ci S_\sB(\pi_N)=S_\sB(\pi_M)$. As $S_\sB(\pi_N),S_\sB(\pi_M)$ are unital homotopy equivalences, $S_\sB(f)$ is one too. This gives Definition~\ref{ai2def18}(i).

To prove $(S_\sB,i_\sB)$ satisfies Definition \ref{ai2def18}(ii), we must show that if $\fF:\sA_1\ra\sA_2$ is a unital homotopy equivalence in $\ACat_R^\un$, then $F^\un_{S_\sB}(\fF):F^\un_{S_\sB}(\sA_2)\ra F^\un_{S_\sB}(\sA_1)$ is a weak homotopy equivalence of simplicial sets. By part (b), this is equivalent to showing that \eq{ai9eq7} is a weak homotopy equivalence. By Theorem \ref{ai9thm1}(a), we may identify \eq{ai9eq7} with the morphism of mapping spaces 
\e
\fF^*:\Hom_{\ACat^{\un,\iy}_R}^L(\sA_2,\sB)\longra \Hom_{\ACat^{\un,\iy}_R}^L(\sA_1,\sB)
\label{ai11eq15}
\e
in $\ACat^{\un,\iy}_R$. Since $\ACat^{\un,\iy}_R$ is an $\iy$-category by Theorem \ref{ai7thm3}, and unital homotopy equivalences are equivalences in $\ACat^{\un,\iy}_R$, it follows from properties of mapping spaces in $\iy$-categories that \eq{ai11eq15} is a weak homotopy equivalence. Thus Definition \ref{ai2def18}(ii) holds, proving~(a).

For (b), Holstein \cite[\S 3]{Hols} gives an explicit construction of simplicial resolutions of a dg-category $\sB$ in the model category $\DGCat_R$. He first takes a fibrant replacement $\sB'$ of $\sB$. Then he writes down \cite[Def.~3.2]{Hols} an explicit functor $\De^\rop\ra\DGCat_R$. We can check that this is strictly isomorphic to our $S_{\sB'}$ in Definition \ref{ai9def1}. Since we assume $\sB$ is fibrant we can take $\sB'=\sB$. Thus (b) follows from~\cite[\S 3]{Hols}.

Part (d) is immediate. 

For (e), first note that if $\sB$ is a dg-category then $\Fun_{A_\iy}^\un(\sC^\De_{[N]},\sB)$ is a dg-category by Remark \ref{ai4rem12}, and so $S_\sB([N])\subset\Fun_{A_\iy}^\un(\sC^\De_{[N]},\sB)$ is a dg-category. Hence $S_\sB$ maps to $\DGCat_R^{A_\iy}\subset \ACat^\un_R$, proving the claim for~$\DGCat_R^{A_\iy}$.

By Definition \ref{ai2def18} and Theorem \ref{ai9thm2}(a), for any $N\ge 0$, taking $\pi:[N]\ra[0]$ to be the projection in $\De$, then $S_\sB(\pi):\sB\cong S_\sB([0])\ra S_\sB([N])$ is a unital homotopy equivalence. Thus, for any property $\bs P$ of unital $A_\iy$-categories which is invariant under unital homotopy equivalence, if $\sB$ satisfies $\bs P$, then so does $S_\sB([N])$ for all $N\ge 0$.

By Propositions \ref{ai4prop5} and \ref{ai4prop8} and Definition \ref{ai4def22}, having split units, and being h-projective, and being pretriangulated, are all invariant under unital homotopy equivalence. Thus we see that if $\sB$ lies in any of the full subcategories $\cC$ in \eq{ai9eq9}, then $S_\sB$ maps to $\cC$. This proves the first part of (e). Hence mapping spaces in $\cC[(\HE_\un\cap\cC)^{-1}]^\iy$ are given by \eq{ai9eq8} as in (c),(d). 

The inclusion $\cC[(\HE_\un\cap\cC)^{-1}]^\iy\ra \ACat^\un_R[\HE_\un^{-1}]^\iy$ also induces homotopy equivalences on mapping spaces, as the computation of mapping spaces using simplicial resolution is functorial under this inclusion. So the inclusion is fully faithful by Definition \ref{ai2def9}. This completes part (e), and Theorem~\ref{ai9thm2}.

\subsection{Proof of Theorem \ref{ai9thm3}}
\label{ai113}

Here is equation \eq{ai9eq10} with morphisms labelled, to help with the proof:
\e
\begin{gathered}
\xymatrix@C=150pt@R=15pt{
*+[r]{\DGCat_R^\hs[\HE_{\dg,\hs}^{-1}]^\iy}
\ar[r]^\simeq_{\Al} \ar@{.>}[d]^\Be
& *+[l]{\DGCat_R[\QE_\dg^{-1}]^\iy} \ar[d]^\simeq_{\Ga}
\\
*+[r]{\DGCat^{A_\iy,\hs}_R[\HE_{\dg,A_\iy,\hs}^{-1}]^\iy} \ar@{.>}[r]^\De \ar@{.>}[d]^\Ep
& *+[l]{\DGCat^{A_\iy}_R[\QE_{\dg,A_\iy}^{-1}]^\iy} \ar[d]^\simeq_{\Ze}
\\
*+[r]{\ACat^{\un,\hs}_R[\HE_{\un,\hs}^{-1}]^\iy} \ar@{.>}[r]^\Eta &
*+[l]{\ACat_R^\un[\QE_\un^{-1}]^\iy.}
}
\end{gathered}
\label{ai11eq16}
\e
We already know that $\Al,\Ga,\Ze$ are Joyal equivalences by \eq{ai4eq44} and \eq{ai4eq46}. We must show that $\Be,\De,\Ep,\Eta$ are Joyal equivalences. Since \eq{ai11eq16} commutes, it is sufficient to prove that $\Be,\Ep$ are Joyal equivalences.

By Theorem \ref{ai2thm1}, we must show $\Be,\Ep$ are {\it essentially surjective}, as in Definition \ref{ai2def5}, and {\it fully faithful}, as in Definition \ref{ai2def9}. Here $\Be$ is trivially essentially surjective, as $\DGCat_R^\hs$ and $\DGCat^{A_\iy,\hs}_R$ have the same objects, and $\Ep$ is essentially surjective by Theorem \ref{ai4thm4} and Proposition \ref{ai4prop8}. Also $\Ep$ is fully faithful by Theorem \ref{ai9thm2}(e), since the computation of mapping spaces in Theorem \ref{ai9thm2}(d) applies in both $\DGCat^{A_\iy,\hs}_R[\HE_{\dg,A_\iy,\hs}^{-1}]^\iy$ and $\ACat^{\un,\hs}_R[\HE_{\un,\hs}^{-1}]^\iy$, and $\Ep$ identifies these mapping spaces. It remains to prove that $\Be$ is fully faithful.

By Canonaco--Ornaghi--Stellari \cite[Prop.~3.5]{COS2} (see also Pascaleff \cite{Pasc} and Faonte \cite{Faon2}), there is an adjunction
\e
\xymatrix@C=70pt{
\ACat^\un_R \ar@<.5ex>[r]^U & \DGCat_R, \ar@<.5ex>[l]^\inc }
\label{ai11eq17}
\e
such that writing $\rho:\id_{\ACat^\un_R}\Ra\inc\ci U$ and $\si:U\ci\inc\Ra\id_{\DGCat_R}$ for the unit and counit of the adjunction, then $\rho_\sA:\sA\ra U(\sA)$ and $\si_\sB:U(\sB)\ra\sB$
are quasi-isomorphisms, and $\rho_\sA:\sA\ra U(\sA)$ is a homotopy equivalence if $\sA$ has split units. 

Here $U$ is called the {\it cobar-bar resolution functor}, as it is the composition of the `bar construction' mapping $A_\iy$-categories to dg-cocategories, and the `cobar construction' mapping dg-cocategories to dg-categories. Also $U\ci\inc:\DGCat_R\ra\DGCat_R$ is called the {\it enveloping dg-category functor\/} \cite[Def.~2.10 \& Rem.~2.11]{Faon2}, and is a cofibrant replacement functor for the Tabuada model structure on $\DGCat_R$ from Theorem~\ref{ai4thm1}. 

Suppose $\sA,\sB$ lie in $\DGCat^\hs_R\subset\ACat^{\un,\hs}_R,$ and $\sB$ is fibrant in the Tabuada model structure on $\DGCat_R$ from Theorem \ref{ai4thm1}. As $\rho_\sA:\sA\ra U(\sA)$ is a homotopy equivalence, $U(\sA)$ also lies in $\DGCat^\hs_R$ by Propositions \ref{ai4prop5} and \ref{ai4prop8}, and it is cofibrant. We consider the action of $\Ep$ on mapping spaces:
\ea
\Ep_{U(\sA),\sB}&:\Hom_{\DGCat_R^\hs[\HE_{\dg,\hs}^{-1}]^\iy}(U(\sA),\sB)\cong
\Hom_{\DGCat_R[\QE_{\dg}^{-1}]^\iy}(U(\sA),\sB)
\nonumber\\
&\longra\Hom_{\DGCat^{A_\iy,\hs}_R[\HE_{\dg,A_\iy,\hs}^{-1}]^\iy}(U(\sA),\sB).
\label{ai11eq18}
\ea
Here the top equivalence holds by Theorem \ref{ai4thm11}.

Now Theorem \ref{ai9thm2} allows us to compute the mapping spaces in the top and bottom lines of \eq{ai11eq18} explicitly, up to homotopy equivalence, using (two different kinds of) simplicial resolutions. For the top line, as $\sB$ is fibrant, Theorem \ref{ai9thm2}(b) says that $(S_\sB,i_\sB)$ is a simplicial resolution of $\sB$ in the model category $\DGCat_R$ in the sense of Remark \ref{ai2rem4}(c) and Hirschhorn \cite[\S 16]{Hirs}. Since $U(\sA)$ is cofibrant, by the theory of simplicial resolutions, the mapping space $\Hom_{\DGCat_R[\QE_{\dg}^{-1}]^\iy}(U(\sA),\sB)$ is represented by the simplicial set
\e
[N]\longmapsto \Hom_{\DGCat_R}(U(\sA),S_\sB([N])).
\label{ai11eq19}
\e
But by the adjunction \eq{ai11eq17}, this is strictly isomorphic to the simplicial set
\e
[N]\longmapsto \Hom_{\ACat^\un_R}(\sA,S_\sB([N]))=F^\un_{S_\sB}(\sA)([N]),
\label{ai11eq20}
\e
where $F^\un_{S_\sB}(\sA)$ is as in Definition \ref{ai2def18}.

By Theorem \ref{ai9thm2}(a),(e), $(S_\sB,i_\sB)$ is a simplicial resolution of $\sB$ in the category with weak equivalences $(\DGCat^{A_\iy,\hs}_R,\HE_{\dg,A_\iy,\hs})$, in the sense of  Definition \ref{ai2def18}. Thus by Theorem \ref{ai2thm5}, $\Hom_{\DGCat^{A_\iy,\hs}_R[\HE_{\dg,A_\iy,\hs}^{-1}]^\iy}(U(\sA),\sB)$ is represented by the simplicial set
\e
[N]\longmapsto \Hom_{\ACat^\un_R}(U(\sA),S_\sB([N]))=F^\un_{S_\sB}(U(\sA))([N]).
\label{ai11eq21}
\e

It is now clear, from the theory of simplicial resolutions, that up to homotopy, the morphism $\Ep_{U(\sA),\sB}$, as a map between the representatives \eq{ai11eq19} and \eq{ai11eq21}, is induced by the map between the right hand sides of \eq{ai11eq19} and \eq{ai11eq21} induced by the inclusion $\inc:\DGCat_R\hookra\ACat^\un_R$. So by the relation between \eq{ai11eq19} and \eq{ai11eq20} using the adjunction \eq{ai11eq17}, we see that $\Ep_{U(\sA),\sB}$, as a map between the representatives \eq{ai11eq20} and \eq{ai11eq21} is induced by composition with $\si_\sA:U(\sA)\ra\sA$.

We now have a diagram of simplicial sets commuting up to homotopy:
\begin{equation*}
\xymatrix@C=140pt@R=1pt{
*+[r]{\Hom_{\DGCat_R^\hs[\HE_{\dg,\hs}^{-1}]^\iy}(U(\sA),\sB)} \ar[dd]^{} 
\ar[dr]_\simeq \\
& *+[l]{\sN_{A_\iy}(\Fun_{A_\iy}^\un(\sA,\sB))^\Kan} \ar[dd]_{\sN_{A_\iy}(\si_\sA^*)^\Kan} \\
*+[r]{\Hom_{\DGCat^{A_\iy,\hs}_R[\HE_{\dg,A_\iy,\hs}^{-1}]^\iy}(U(\sA),\sB)} \ar[dr]_\simeq
\\
& *+[l]{\sN_{A_\iy}(\Fun_{A_\iy}^\un(U(\sA),\sB))^\Kan.\! }
}
\end{equation*}
Here the diagonal morphisms combine the representatives \eq{ai11eq20}--\eq{ai11eq21} with equation \eq{ai9eq6} in Theorem \ref{ai9thm2}(b), and the right hand morphism is as in~\eq{ai9eq7}.

As above, $\si_\sA$ is a quasi-equivalence, and thus a homotopy equivalence of unital $A_\iy$-categories by Corollary \ref{ai4cor3}(b), since $\sA,U(\sA)$ are h-projective with split units. Thus by the argument after \eq{ai11eq15}, $\sN_{A_\iy}(\si_\sA^*)^\Kan$ is a homotopy equivalence of Kan complexes, so $\Ep_{U(\sA),\sB}$ is also a homotopy equivalence.

To prove that $\Ep$ is fully faithful, as in Definition \ref{ai2def9}, we must show that $\Ep_{\sA,\sB}$ is a homotopy equivalence of Kan complexes for all objects $\sA,\sB$ in $\DGCat_R^\hs$. We have proved that $\Ep_{U(\sA),\sB}$ is a homotopy equivalence of Kan complexes provided $\sB$ is fibrant. Since $\sA$ is homotopy equivalent to $U(\sA)$ in $\DGCat_R^\hs$, and any $\sB$ in $\DGCat_R^\hs$ is homotopy equivalent to a fibrant $\sB'$, we see that $\Ep$ is fully faithful. This completes the proof of Theorem~\ref{ai9thm3}.

\subsection{Proof of Theorem \ref{ai9thm4}}
\label{ai114}

As in \eq{ai9eq5} we have a homotopy commutative diagram of $\iy$-functors
\e
\begin{gathered}
\xymatrix@C=80pt@R=15pt{
*+[r]{\sN(\ACat^\un_R)} \ar[rr]_(0.55){\Pi_\loc^\iy} \ar@/_.5pc/[drr]_{\iinc^\un_R} &  \drtwocell_{}\omit^{}\omit{^{}} & *+[l]{\ACat^\un_R[\HE_\un^{-1}]^\iy} \ar[d]^{\Phi} \\
&& *+[l]{\ACat^{\un,\iy}_R.} }
\end{gathered}
\label{ai11eq22}
\e
We will prove that $\Phi$ is a Joyal equivalence. By Theorem \ref{ai2thm1}, it is enough to show $\Phi$ is {\it essentially surjective}, as in Definition \ref{ai2def5}, and {\it fully faithful}, as in Definition \ref{ai2def9}. In \eq{ai11eq22}, $\sN(\ACat^\un_R),\ACat^\un_R[\HE_\un^{-1}]^\iy$ and $\ACat^{\un,\iy}_R$ have the same objects, unital $A_\iy$-categories over $R$, and $\iinc^\un_R$ and $\Pi_\loc^\iy$ are the identity on objects. Hence $\Phi$ is also the identity on objects, so it is trivially essentially surjective.

To prove $\Phi$ is fully faithful, let $\sA,\sB$ be unital $A_\iy$-categories over $R$, and consider the homotopy commutative diagram in the $\iy$-category $\sSpc$ from~\S\ref{ai25}:
\e
\begin{gathered}
\xymatrix@C=140pt@R=15pt{
*+[r]{\Hom^L_{\ACat^\un_R[\HE_\un^{-1}]^\iy}(\sA,\sB)} \ar@<2ex>[d]_\simeq^{\eq{ai9eq8}} \ar[r]_(0.54){\Phi^L_{\sA,\sB}} & *+[l]{\Hom^L_{\ACat^{\un,\iy}_R}(\sA,\sB)} \ar@<-2ex>[d]^\simeq_{\eq{ai9eq1}}  \\
*+[r]{\sN_{A_\iy}(\Fun_{A_\iy}^\un(\sA,\sB))^\Kan}  \ar@{..>}[r]^{\Phi_{\sA,\sB}^\sh}
& *+[l]{\sN_{A_\iy}(\Fun_{A_\iy}^\un(\sA,\sB))^\Kan.\!}  }
\end{gathered}
\label{ai11eq23}
\e
Here $\Phi^L_{\sA,\sB}$ is the action of $\Phi$ on mapping spaces, as in Definition \ref{ai2def9}. The columns are the equivalences \eq{ai9eq1} and \eq{ai9eq8}, and $\Phi_{\sA,\sB}^\sh$ is the morphism making \eq{ai11eq23} commute up to homotopy. To prove $\Phi$ is essentially surjective, we must show $\Phi^L_{\sA,\sB}$ is a homotopy equivalence for all $\sA,\sB$. We will do this by showing that $\Phi_{\sA,\sB}^\sh$ is homotopic to the identity.

Definition \ref{ai9def1} defines a functor $S_\sB:\De^\rop\ra\ACat^{\un,\iy}_R$ and an isomorphism $i_\sB:\sB\ra S_\sB([0])$, and Theorem \ref{ai9thm2}(a) shows that $(S_\sB,i_\sB)$ is a simplicial resolution of $\sB$ in $(\ACat^\un_R,\HE_\un)$, in the sense of Definition \ref{ai2def18}. Applying nerves gives an $\iy$-functor $\sN(S_\sB):\sN(\De^\rop)\ra\sN(\ACat^{\un,\iy}_R)$.

Recall that for an $\iy$-category $\sC$ and an object $X\in\sC$, Example \ref{ai2ex7} defined a left fibration $p_\sC:\sC_{X/}\ra\sC$. Applying the straightening/unstraightening correspondence in Theorems \ref{ai2thm2} and \ref{ai2thm3} to $p_\sC$ yields an $\iy$-functor $F_X:\sC\ra\sSpc$, the straightening of $p_\sC$, such that $F_X(Y)$ is homotopy equivalent to $\Hom^L_\sC(X,Y)$. We apply this to the three $\iy$-categories in \eq{ai11eq22} with object $X=\sA$, and then pull the left fibrations back to $\sN(\De^\rop)$ using the $\iy$-functors $\sN(S_\sB)$, $\Pi_\loc^\iy\ci\sN(S_\sB)$, $\iinc^\un_R\ci\sN(S_\sB)$. In this way we obtain a homotopy commutative diagram of $\iy$-categories
\e
\begin{gathered}
\xymatrix@C=135pt@R=20pt{
*+[r]{\sN(S_\sB)^*(\sN(\ACat^\un_R)_{\sA/})} \ar[r]_(0.55)a \ar[d]^p \ar[dr]_(0.3)c  & *+[l]{\begin{subarray}{l} \ts (\Pi_\loc^\iy\ci\sN(S_\sB))^* \\ \ts (\ACat^\un_R[\HE_\un^{-1}]^\iy)_{\sA/})\end{subarray}\;\>} \ar[d]_(0.6)b \ar[dl]^(0.6)q \\
*+[r]{\sN(\De^\rop)} & *+[l]{(\iinc^\un_R\ci\sN(S_\sB))^*((\ACat^{\un,\iy}_R)_{\sA/}).} \ar[l]_(0.7)r }\!\!\!\!\!\!\!\!\!\!\!
\end{gathered}
\label{ai11eq24}
\e

Here $p,q,r$ are left fibrations, as they are pullbacks of $p_\sC$ in Example \ref{ai2ex7} for $\sC=\sN(\ACat^\un_R)$, $\ACat^\un_R[\HE_\un^{-1}]^\iy$, $\ACat^{\un,\iy}_R$, and $a,b,c$ are induced by the morphisms in \eq{ai11eq22}. The top left and bottom left triangles in \eq{ai11eq24} commute strictly, and the rest commutes up to homotopy since \eq{ai11eq22} does. The fibre of \eq{ai11eq24} over $[N]\in\sN(\De^\rop)$ is
\e
\begin{gathered}
\xymatrix@C=140pt@R=20pt{
*+[r]{\Hom^L_{\sN(\ACat^\un_R)}(\sA,S_\sB([N]))} \ar[r]^(0.47){\raisebox{8pt}{$\st(\Pi_\loc^\iy)^L_{\sA,S_\sB([N])}$}} \ar[d] \ar[dr]_(0.25){(\iinc^\un_R)^L_{\sA,S_\sB([N])}}  & *+[l]{\Hom^L_{\ACat^\un_R[\HE_\un^{-1}]^\iy}(\sA,S_\sB([N]))} \ar[d]_{\Phi^L_{\sA,S_\sB([N])}} \ar[dl] \\
*+[r]{*} & *+[l]{\Hom^L_{\ACat^{\un,\iy}_R}(\sA,S_\sB([N])),} \ar[l] }\!\!\!\!\!
\end{gathered}
\label{ai11eq25}
\e
with top right triangle the action of \eq{ai11eq22} on mapping spaces.

As $p,q,r$ are left fibrations, we can apply `straightening' in Theorems \ref{ai2thm2}(b) and \ref{ai2thm3} to \eq{ai11eq24}, to get a diagram of $\iy$-functors $\sN(\De^\rop)\ra\sSpc$, which we regard as a homotopy commutative triangle in the presheaf $\iy$-category $\PSh(\sN(\De))=\Fun(\sN(\De^\rop),\sSpc)$ from Definition \ref{ai2def15}:
\e
\begin{gathered}
\xymatrix@C=80pt@R=15pt{
*+[r]{\Ga_{\sA,\sB}} \ar[rr]_(0.55){\la_{\sA,\sB}} \ar@/_.5pc/[drr]_{\nu_{\sA,\sB}} && *+[l]{\De_{\sA,\sB}} \ar[d]_{\mu_{\sA,\sB}} \\
&& *+[l]{\Ep_{\sA,\sB}.} }
\end{gathered}
\label{ai11eq26}
\e

Conversely, equation \eq{ai11eq24} is up to homotopy equivalence the `unstraightening' of \eq{ai11eq26}. Here the actions of the $\iy$-functors $\Ga_{\sA,\sB},\De_{\sA,\sB},\Ep_{\sA,\sB}$ and $\iy$-natural transformations $\la_{\sA,\sB},\mu_{\sA,\sB},\nu_{\sA,\sB}$ at an object $[N]\in\sN(\De^\rop)$ are homotopy equivalent to the top right triangle in \eq{ai11eq25}, that is
\e
\begin{gathered}
\xymatrix@C=130pt@R=20pt{
*+[r]{\begin{subarray}{l} \ts \Ga_{\sA,\sB}([N])\simeq \\
\ts \Hom^L_{\sN(\ACat^\un_R)}(\sA,S_\sB([N])) \end{subarray}} \ar[r]^(0.47){\raisebox{8pt}{$\st \la_{\sA,\sB}([N])\simeq(\Pi_\loc^\iy)^L_{\sA,S_\sB([N])}$}} \ar@/_1pc/[dr]_(0.25){\nu_{\sA,\sB}([N])\simeq(\iinc^\un_R)^L_{\sA,S_\sB([N])}}  & *+[l]{\begin{subarray}{l} \ts \qquad\qquad\qquad \De_{\sA,\sB}([N])\simeq \\
\ts \Hom^L_{\ACat^\un_R[\HE_\un^{-1}]^\iy}(\sA,S_\sB([N])) \end{subarray}} \ar[d]_{\mu_{\sA,\sB}([N])\simeq\Phi^L_{\sA,S_\sB([N])}} \\
& *+[l]{\begin{subarray}{l} \ts \Ep_{\sA,\sB}([N])\simeq \\
\ts  \Hom^L_{\ACat^{\un,\iy}_R}(\sA,S_\sB([N])). \end{subarray}} }\!\!\!\!\!\!\!\!
\end{gathered}
\label{ai11eq27}
\e
In particular, when $N=0$ we can use the isomorphism $i_\sB:\sB\ra S_\sB([0])$ to identity \eq{ai11eq27} with
\e
\begin{gathered}
\xymatrix@C=130pt@R=20pt{
*+[r]{\begin{subarray}{l} \ts \Ga_{\sA,\sB}([0])\simeq \\
\ts \Hom^L_{\sN(\ACat^\un_R)}(\sA,\sB) \end{subarray}} \ar[r]^(0.47){\raisebox{8pt}{$\st \la_{\sA,\sB}([0])\simeq(\Pi_\loc^\iy)^L_{\sA,\sB}$}} \ar@/_1pc/[dr]_(0.25){\nu_{\sA,\sB}([0])\simeq(\iinc^\un_R)^L_{\sA,\sB}}  & *+[l]{\begin{subarray}{l} \ts \qquad\qquad\qquad \De_{\sA,\sB}([0])\simeq \\
\ts \Hom^L_{\ACat^\un_R[\HE_\un^{-1}]^\iy}(\sA,\sB) \end{subarray}} \ar[d]_{\mu_{\sA,\sB}([0])\simeq\Phi^L_{\sA,\sB}} \\
& *+[l]{\begin{subarray}{l} \ts \Ep_{\sA,\sB}([0])\simeq \\
\ts  \Hom^L_{\ACat^{\un,\iy}_R}(\sA,\sB). \end{subarray}} }\!\!\!\!\!\!\!\!
\end{gathered}
\label{ai11eq28}
\e

The effect of a 1-morphism $\al:[N]\ra[N']$ in $\sN(\De^{\rm op})$ on the spaces in \eq{ai11eq27} is up to homotopy the action of $S_\sB(\al):S_\sB([N])\ra S_\sB([N'])$ on mapping spaces $\Hom^L_{\cdots}(\sA,-)$. Since $S_\sB(\al)$ is a homotopy equivalence in $\ACat^\un_R$ by Definition \ref{ai2def18}(i) and Theorem \ref{ai9thm1}(a), it is an equivalence in both $\ACat^\un_R[\HE_\un^{-1}]^\iy$ and $\ACat^{\un,\iy}_R$. Thus $\al$ induces homotopy equivalences $\De_{\sA,\sB}([N])\ra \De_{\sA,\sB}([N'])$ and $\Ep_{\sA,\sB}([N])\ra \Ep_{\sA,\sB}([N'])$. (In general it does {\it not\/} induce homotopy equivalences $\Ga_{\sA,\sB}([N])\ra \Ga_{\sA,\sB}([N'])$.) As $\sN(\De^\rop)$ is contractible, this means that $\De_{\sA,\sB},\Ep_{\sA,\sB}$ are homotopic to constant $\iy$-functors, which by \eq{ai11eq28} we may take to have constant values $\Hom^L_{\ACat^\un_R[\HE_\un^{-1}]^\iy}(\sA,\sB)$ and~$\Hom^L_{\ACat^{\un,\iy}_R}(\sA,\sB)$.

Let us now take homotopy colimits in $\sSpc$ over the diagram $\sN(\De^\rop)$ of equation \eq{ai11eq26}. (For colimits in $\iy$-categories, see Lurie \cite[\S 4]{Luri2}, \cite[\S 7]{Luri5}; all colimits exist in $\sSpc$ by \cite[Cor.~7.4.3.3]{Luri5}.) As $\sN(\De^\rop)$ is contractible, and $\De_{\sA,\sB},\Ep_{\sA,\sB}$ are homotopic to constant with values $\Hom^L_{\ACat^\un_R[\HE_\un^{-1}]^\iy}(\sA,\sB)$ and $\Hom^L_{\ACat^{\un,\iy}_R}(\sA,\sB)$, this induces a homotopy commutative diagram in~$\sSpc$:
\e
\begin{gathered}
\xymatrix@C=90pt@R=15pt{
*+[r]{\colim\limits_{\sN(\De^\rop)}\Ga_{\sA,\sB}} \ar[rr]_(0.45){\colim\la_{\sA,\sB}}^(0.45)\simeq \ar@/_1pc/[drr]_{\colim\nu_{\sA,\sB}} && *+[l]{\begin{subarray}{l} \ts \colim\limits_{\sN(\De^\rop)}\De_{\sA,\sB} \simeq \\ \ts \Hom^L_{\ACat^\un_R[\HE_\un^{-1}]^\iy}(\sA,\sB) \end{subarray}} \ar[d]_{\colim\mu_{\sA,\sB}\simeq\Phi^L_{\sA,\sB}} \\
&& *+[l]{\begin{subarray}{l} \ts  \colim\limits_{\sN(\De^\rop)}\Ep_{\sA,\sB}\simeq \\ \ts  \Hom^L_{\ACat^{\un,\iy}_R}(\sA,\sB). \end{subarray}} }\!\!\!\!\!\!\!
\end{gathered}
\label{ai11eq29}
\e
Since $(S_\sB,i_\sB)$ is a simplicial resolution of $\sB$ in $(\ACat^\un_R,\HE_\un)$, as in Definition \ref{ai2def18} and Theorem \ref{ai2thm5}, it follows from Arakawa--Carranza--Kapulkin \cite[\S 2]{ACK} that the top morphism $\colim\la_{\sA,\sB}$ in \eq{ai11eq29} is an equivalence.

The next proposition will be proved in \S\ref{ai1141}--\S\ref{ai1149}.

\begin{prop}
\label{ai11prop1}
There is a diagram of explicit simplicial sets and morphisms
\e
\begin{gathered}
\xymatrix@C=135pt@R=17pt{
*+[r]{\sN(S_\sB)^*(\sN(\ACat^\un_R)_{\sA/})} \ar[r]_(0.45){a'} \ar[d]^p \ar[dr]_(0.3)c  & *+[l]{\sN(\De^\rop)\!\t\!\Hom^L_{\ACat^{\un,\iy}_R}(\sA,\sB)} \ar[d]_{b'} \ar[dl]^(0.62){\Pi_{\sN(\De^\rop)}} \\
*+[r]{\sN(\De^\rop)} & *+[l]{(\iinc^\un_R\ci\sN(S_\sB))^*((\ACat^{\un,\iy}_R)_{\sA/}).} \ar[l]_(0.7)r }\!\!\!\!\!\!\!\!\!\!\!
\end{gathered}
\label{ai11eq30}
\e
Here the bottom left triangle is as in {\rm\eq{ai11eq24},} so that\/ $p,r$ are left fibrations. The projection to the first factor $\Pi_{\sN(\De^\rop)}:\sN(\De^\rop)\t\Hom^L_{\ACat^{\un,\iy}_R}(\sA,\sB)\ra\sN(\De^\rop)$ is also a left fibration. Furthermore:
\begin{itemize}
\setlength{\itemsep}{0pt}
\setlength{\parsep}{0pt}
\item[{\bf(a)}] The morphism $a'$ is strictly isomorphic to the `last vertex map' construction of Lurie {\rm\cite[\S 3.3.4]{Luri5},} for the simplicial set\/ $\Hom^L_{\ACat^{\un,\iy}_R}(\sA,\sB)$.
\item[{\bf(b)}] The bottom left, top left, and bottom right triangles in \eq{ai11eq30} strictly commute. The top right triangle in \eq{ai11eq30} does not strictly commute. However, we can construct an explicit morphism of simplicial sets 
\e
\begin{split}
H&:\sN(S_\sB)^*(\sN(\ACat^\un_R)_{\sA/})\t\De^1\\
&\qquad\qquad \longra (\iinc^\un_R\ci\sN(S_\sB))^*((\ACat^{\un,\iy}_R)_{\sA/}),
\end{split}
\label{ai11eq31}
\e
which satisfies
\e
\begin{split}
&H\vert_{\sN(S_\sB)^*(\sN(\ACat^\un_R)_{\sA/})\t\De^{\{0\}}}=c, \\
&H\vert_{\sN(S_\sB)^*(\sN(\ACat^\un_R)_{\sA/})\t\De^{\{1\}}}=b'\ci a',\\
&r\ci H=p\ci\Pi_{\sN(S_\sB)^*(\sN(\ACat^\un_R)_{\sA/})}.
\end{split}
\label{ai11eq32}
\e
Thus, the top right triangle commutes up to homotopy over~$\sN(\De^\rop)$.
\item[{\bf(c)}] As for {\rm\eq{ai11eq25},} the fibre of\/ \eq{ai11eq30} over $[0]\in\sN(\De^\rop)$ is
\end{itemize}
\e
\begin{gathered}
\xymatrix@C=140pt@R=20pt{
*+[r]{\Hom^L_{\sN(\ACat^\un_R)}(\sA,S_\sB([0]))} \ar[r]_(0.55){a'\vert_{[0]}} \ar[d] \ar[dr]_(0.25){(\iinc^\un_R)^L_{\sA,S_\sB([0])}}  & *+[l]{\Hom^L_{\ACat^{\un,\iy}_R}(\sA,\sB)} \ar[d]_{b'\vert_{[0]}} \ar[dl] \\
*+[r]{*} & *+[l]{\Hom^L_{\ACat^{\un,\iy}_R}(\sA,S_\sB([0])).} \ar[l] }\!\!\!\!\!
\end{gathered}
\label{ai11eq33}
\e
\begin{itemize}
\setlength{\itemsep}{0pt}
\setlength{\parsep}{0pt}
\item[] Then $b'\vert_{[0]}$ is a strict isomorphism of simplicial sets, induced by composition with the isomorphism $i_\sB:\sB\ra S_\sB([0])$.
\end{itemize}
\end{prop}

As for the construction of \eq{ai11eq26} from \eq{ai11eq25} we can apply `straightening' in Theorems \ref{ai2thm2}(b) and \ref{ai2thm3} to \eq{ai11eq30}, to get a homotopy commutative triangle in $\PSh(\sN(\De))$, where we take $\Ga_{\sA,\sB},\Ep_{\sA,\sB},\nu_{\sA,\sB}$ to be as in~\eq{ai11eq26}:
\e
\begin{gathered}
\xymatrix@C=80pt@R=15pt{
*+[r]{\Ga_{\sA,\sB}} \ar[rr]_(0.55){\la'_{\sA,\sB}} \ar@/_.5pc/[drr]_{\nu_{\sA,\sB}} && *+[l]{\De'_{\sA,\sB}} \ar[d]_{\mu'_{\sA,\sB}} \\
&& *+[l]{\Ep_{\sA,\sB}.} }
\end{gathered}
\label{ai11eq34}
\e
The conditions on $H$ in Proposition \ref{ai11prop1}(b) ensure this commutes up to homotopy. We may also take $\De'_{\sA,\sB}$ to be the constant $\iy$-functor $\sN(\De^\rop)\ra\sSpc$ with value~$\Hom^L_{\ACat^{\un,\iy}_R}(\sA,\sB)$.

As for \eq{ai11eq29}, we take homotopy colimits in $\sSpc$ over the diagram $\sN(\De^\rop)$ of equation \eq{ai11eq34} to get a homotopy commutative diagram in~$\sSpc$:
\e
\begin{gathered}
\xymatrix@C=90pt@R=15pt{
*+[r]{\colim\limits_{\sN(\De^\rop)}\Ga_{\sA,\sB}} \ar[rr]_(0.45){\colim\la'_{\sA,\sB}}^(0.45)\simeq \ar@/_1pc/[drr]_{\colim\nu_{\sA,\sB}} && *+[l]{\begin{subarray}{l} \ts \colim\limits_{\sN(\De^\rop)}\De'_{\sA,\sB} \simeq \\ \ts \Hom^L_{\ACat^{\un,\iy}_R}(\sA,\sB) \end{subarray}} \ar[d]_{\colim\mu'_{\sA,\sB}\simeq \id_{\Hom^L_{\ACat^{\un,\iy}_R}(\sA,\sB)}} \\
&& *+[l]{\begin{subarray}{l} \ts  \colim\limits_{\sN(\De^\rop)}\Ep_{\sA,\sB}\simeq \\ \ts  \Hom^L_{\ACat^{\un,\iy}_R}(\sA,\sB). \end{subarray}} }\!\!\!\!\!\!\!
\end{gathered}
\label{ai11eq35}
\e
Since $a'$ is strictly isomorphic to a `last vertex map' by Proposition \ref{ai11prop1}(a), with straightening $\la'_{\sA,\sB},$ by Lurie \cite[Prop.~3.3.4.8]{Luri5} $\colim_{\sN(\De^\rop)}\la'_{\sA,\sB}$ is a homotopy equivalence. 

As $\De'_{\sA,\sB},\Ep_{\sA,\sB}$ are constant up to homotopy, $b'\vert_{[0]}$ in \eq{ai11eq33} represents $\colim\mu'_{\sA,\sB}$ up to homotopy, so Proposition \ref{ai11prop1}(c) implies that $\colim\mu'_{\sA,\sB}$ is a homotopy equivalence (in fact it is homotopic to $\id_{\Hom^L_{\ACat^{\un,\iy}_R}(\sA,\sB)}$). Now comparing \eq{ai11eq29} and \eq{ai11eq35} shows that $\Phi^L_{\sA,\sB}$ is a homotopy equivalence (and in fact $\Phi_{\sA,\sB}^\sh$ in \eq{ai11eq23} is homotopic to $\id_{\sN_{A_\iy}(\Fun_{A_\iy}^\un(\sA,\sB))^\Kan}$). As this holds for all $\sA,\sB$, we see that $\Phi$ in \eq{ai9eq11} is fully faithful, and thus a Joyal equivalence. This proves the first part of Theorem~\ref{ai9thm4}.

For the second part of Theorem \ref{ai9thm4}, the $\iy$-functors analogous to $\Phi$ for each of the full subcategories $\cC$ of $\ACat^\un_R$ in \eq{ai9eq9} are trivially essentially surjective, as they act by the identity on objects. And the proof above shows that they are fully faithful, since by Theorem \ref{ai9thm2}(e) if $\sB\in\cC\subset\ACat^\un_R$ then $(S_\sB,i_\sB)$ is a simplicial resolution of $\sB$ in $(\cC,\HE_\un\cap\cC)$. So they are Joyal equivalences. Given Proposition \ref{ai11prop1}, this completes the proof of Theorem \ref{ai9thm4}. The rest of the section proves Proposition~\ref{ai11prop1}.

\subsubsection{A map $\Xi_\sB:\cN(\De^\rop)\t\De^1\t\De^1\ra\ACat^{\un,\iy}_R$}
\label{ai1141}

We begin the proof of Proposition \ref{ai11prop1} by defining a morphism $\Xi_\sB$, which depends on the unital $A_\iy$-category $\sB$ over $R$ fixed above.

\begin{dfn}
\label{ai11def1}
We will define a map of simplicial sets
\e
\Xi_\sB:\cN(\De^\rop)\t\De^1\t\De^1\longra\ACat^{\un,\iy}_R.
\label{ai11eq36}
\e
Writing $0,1$ for the vertices of $\De^1$, so that $(0,0),(1,0),(0,1),(1,1)$ are the vertices of $\De^1\t\De^1$, this will have the following properties:
\begin{itemize}
\setlength{\itemsep}{0pt}
\setlength{\parsep}{0pt}
\item[(a)] $\Xi_\sB\vert_{\cN(\De^\rop)\t(i,j)}=\iinc^\un_R\ci\cN(S_\sB):\cN(\De^\rop)\ra\ACat^{\un,\iy}_R$ for $(i,j)=(0,0),(0,1),(1,1)$, where $S_\sB$ is as in Definition~\ref{ai9def1}.
\item[(b)] $\Xi_\sB\vert_{\cN(\De^\rop)\t(1,0)}$ is the constant map $\cN(\De^\rop)\!\ra\!\ACat^{\un,\iy}_R$ with value~$\sB$.
\item[(c)] $\Xi_\sB\vert_{\cN(\De^\rop)\t\{0\}\t\De^1}=\iinc^\un_R\ci\cN(S_\sB)\ci\Pi_{\cN(\De^\rop)}$ as a map $\cN(\De^\rop)\t\De^1\ra\ACat^{\un,\iy}_R$, where $\Pi_{\cN(\De^\rop)}$ is projection to the first factor.
\item[(d)] $\Xi_\sB\vert_{\cN(\De^\rop)\t\De^1\t\{1\}}=\iinc^\un_R\ci\cN(S_\sB)\ci\Pi_{\cN(\De^\rop)}$.\end{itemize}

Let $N\ge 0$, and for $a=0,\ldots,N+1$ define $\si_a:\De^N\ra\De^1$ to take vertex $i$ in $\De^N$ to vertex 0 in $\De^1$ if $i<a$, and to take vertex $i$ to vertex 1 if $i\ge a$. Then an $N$-simplex in $\cN(\De^\rop)\t\De^1\t\De^1$ is of the form $(S,\si_a,\si_b)$ for $S$ an $N$-simplex in $\cN(\De^\rop)$ and $a,b=0,\ldots,N+1$. The $N$-simplex $S$ in $\cN(\De^\rop)$ is equivalent to a sequence of morphisms in~$\De^\rop$:
\begin{equation*}
\xymatrix@C=30pt{ [M_0] \ar[r]^{\al_1} & [M_1] \ar[r]^{\al_2} & \,\cdots\,  \ar[r]^{\al_N} & [M_N], }
\end{equation*}
where $M_0,\ldots,M_N\ge 0$ and $\al_k:[M_{k-1}]\ra[M_k]$ is a morphism in $\De^\rop$ (equivalently, $\al_k:[M_k]\ra[M_{k-1}]$ is a morphism in $\De$). If $0\le i<j\le N$ we write $\al_{ij}=\al_j\ci\al_{j-1}\ci\cdots\ci\al_{i+1}:[M_i]\ra[M_j]$, and we write~$\al_{ii}=\id_{[M_i]}$.

We will define an $N$-simplex $T=\Xi_\sB(S,\si_a,\si_b)$ in $\ACat^{\un,\iy}_R$, using Definitions \ref{ai7def3} and \ref{ai7def4}. First we define $T$ on all 0- and 1-subsimplices of~$\De^N$:
\begin{itemize}
\setlength{\itemsep}{0pt}
\setlength{\parsep}{0pt}
\item[(i)] For $i=0,\ldots,N$, the unital $A_\iy$-category $\sA^i$ at vertex $i$ of $\De^N$ is $\sB$ if $a\le i<b$, and is $S_\sB([M_i])$ otherwise. This is necessary by (a),(b) above.
\item[(ii)] If $0\le i\le j\le N$ with $a\le i\le j<b$ then the unital $A_\iy$-functor $\fF^{ij}:\sA^i\ra\sA^j$ at edge $i\ra j$ of $\De^N$ is $\id_\sB:\sB\ra\sB$. This is necessary by (b) above.
\item[(iii)] If $0\le i\le j\le N$ with either $j<a$, or $b\le i$, then $\fF^{ij}:\sA^i\ra\sA^j$ is $S_\sB(\al_{ij}):S_\sB([M_i])\ra S_\sB([M_j])$. This is necessary by (c),(d) above.
\item[(iv)] If $0\le i\le j\le N$ with $i\le a<j\le b$, so that edge $i\ra j$ of $\De^N$ maps to edge $(0,0)\mapsto(1,0)$ of $\De^1\t\De^1$, then $\fF^{ij}:\sA^i\ra\sA^j$ is given by the commutative diagram
\e
\begin{gathered}
\xymatrix@C=130pt@R=15pt{
*+[r]{\A^i=S_\sB([M_{i_{\vphantom{(}}}])} \ar[r]_{\fF^{ij}} \ar@{^{(}->}[d]^{\text{ Definition \ref{ai9def1}}} & *+[l]{\sA^j=\sB} \\
*+[r]{\Fun_{A_\iy}^\un(\sC^\De_{[M_i]},\sB)} \ar[r]^{i_{M_i}^*}  & *+[l]{\Fun_{A_\iy}^\un(\sP^\un,\sB),\!} \ar[u]^{\text{Proposition \ref{ai4prop2} }}_\cong }
\end{gathered}
\label{ai11eq37}
\e
where $\sP^\un$ and $\sC^\De_{[M_i]}$ are as in Definitions \ref{ai4def11} and \ref{ai4def20}, and $i_{M_i}:\sP^\un\ra\sC^\De_{[M_i]}$ is the unique unital $A_\iy$-functor acting on objects by $*\mapsto M_i$. 
\item[(v)] If $0\le i\le j\le N$ with $a<i\le b<j$, so that edge $i\ra j$ of $\De^N$ maps to edge $(1,0)\mapsto(1,1)$ of $\De^1\t\De^1$, then $\fF^{ij}:\sA^i\ra\sA^j$ is
given by the commutative diagram
\e
\begin{gathered}
\xymatrix@C=130pt@R=15pt{
*+[r]{\A^i=\sB} \ar[r]_{\fF^{ij}} \ar[d]^{\text{ Proposition \ref{ai4prop2}}}_\cong & *+[l]{\sA^j=S_\sB([M_{j_{\vphantom{(}}}])} \ar@{^{(}->}[d]_{\text{Definition \ref{ai9def1} }} \\
*+[r]{\Fun_{A_\iy}^\un(\sP^\un,\sB)} \ar[r]^{\pi_{M_j}^*}  & *+[l]{\Fun_{A_\iy}^\un(\sC^\De_{[M_j]},\sB),\!} }
\end{gathered}
\label{ai11eq38}
\e
where $\pi_{M_j}:\sC^\De_{[M_j]}\ra\sP^\un$ is the unital $A_\iy$-functor acting on objects by $k\mapsto *$, and on morphisms by $\om_{kl}\mapsto\id_*$. Then $\pi_{M_j}^*$ maps to $S_\sB([M_j])\subset \Fun_{A_\iy}^\un(\sC^\De_{[M_j]},\sB)$, as for any object $F:\sC^\De_{[M_j]}\ra\sB$ in the image of $\pi_{M_j}^*$, $F(\om_{kl})$ is an identity in $\sB$, and so a homotopy equivalence.
\item[(vi)] If $0\le i\le j\le N$ with $i\le a,b<j$, so that edge $i\ra j$ of $\De^N$ maps to edge $(0,0)\mapsto(1,1)$ of $\De^1\t\De^1$, then $\fF^{ij}:\sA^i\ra\sA^j$ is given by the commutative diagram
\e
\begin{gathered}
\xymatrix@C=130pt@R=15pt{
*+[r]{\A^i=S_\sB([M_{i_{\vphantom{(}}}])} \ar[r]_{\fF^{ij}} \ar@{^{(}->}[d]^{\text{ Definition \ref{ai9def1}}} & *+[l]{\sA^j=S_\sB([M_{j_{\vphantom{(}}}])} \ar@{^{(}->}[d]_{\text{Definition \ref{ai9def1} }} \\
*+[r]{\Fun_{A_\iy}^\un(\sC^\De_{[M_i]},\sB)} \ar[r]^{(i_{M_i}\ci\pi_{M_j})^*}  & *+[l]{\Fun_{A_\iy}^\un(\sC^\De_{[M_j]},\sB).\!} }
\end{gathered}
\label{ai11eq39}
\e
\end{itemize}
Here (i)--(vi) determine the data $(\sO^i),(f^{ij}),(\sH^i)$ in $T$ in Definition \ref{ai7def3}, and $\m^I$ for $\es\ne I\subseteq\{0,\ldots,N\}$ with $\md{I}=1,2$. It remains to define the $\m^I$ for~$\md{I}\ge 3$.

Define $\m^I=0$ if $\md{I}\ge 4$. When $\md{I}=3$, let $I=\{i_0,i_1,i_2\}$ with $0\le i_0<i_1<i_2\le N$. We divide into cases according to where $a,b$ sit in this sequence, or equivalently, depending on~$F_{a,b}=(\si_a,\si_b)\ci\De^{\{i_0,i_1,i_2\}}:\De^2\ra\De^1\t\De^1$:
\begin{itemize}
\setlength{\itemsep}{0pt}
\setlength{\parsep}{0pt}
\item[(vii)] Suppose that $F_{a,b}$ maps $\De^2$ to any of the four vertices in $\De^1\t\De^1$, or surjectively onto any of the five edges in $\De^1\t\De^1$ {\it except\/} the edge $(0,0)\ra(1,0)$, or surjectively onto the 2-simplex with vertices $(0,0),(1,0),(1,1)$. In this case we can show from (ii)--(vi) that $\fF^{i_0i_2}=\fF^{i_1i_2}\ci\fF^{i_0i_1}$, and we define $\m^I$ to correspond to the $A_\iy$-natural isomorphism $\id_{\fF^{i_0i_2}}:\fF^{i_1i_2}\ci\fF^{i_0i_1}\Ra\fF^{i_0i_2}$, as in Example \ref{ai7ex2}(c). (Here if $F_{a,b}$ maps to the 2-simplex with vertices $(0,0),(1,0),(1,1)$ then in \eq{ai11eq37}--\eq{ai11eq39} we identify $\fF^{i_0i_1}\cong i_{M_{i_0}}^*$, $\fF^{i_1i_2}\cong \pi_{M_{i_2}}^*$, and~$\fF^{i_0i_2}\cong (i_{M_{i_0}}\ci\pi_{M_{i_2}})^*=\pi_{M_{i_2}}^*\ci i_{M_{i_0}}^*$.)
\item[(viii)] Suppose that $i_0<a\le i_1<i_2<b$, so that $F_{a,b}$ maps $\De^2$ to the edge $(0,0)\ra(1,0)$. Then $\fF^{i_0i_1}=\fF^{i_0i_2}\cong i_{M_{i_0}}^*$ and $\fF^{i_1i_2}=\id_\sB$. Thus $\fF^{i_0i_2}=\fF^{i_1i_2}\ci\fF^{i_0i_1}$, and we define $\m^I$ to correspond to~$\id_{\fF^{i_0i_2}}$.
\item[(ix)] Suppose that $i_0<i_1<a\le i_2<b$, so that $F_{a,b}$ again maps $\De^2$ to the edge $(0,0)\ra(1,0)$. Let us identify $\A^{i_2}=\sB$ with $\Fun_{A_\iy}^\un(\sP^\un,\sB)$ using Proposition \ref{ai4prop2}. Then 
\begin{align*}
\fF^{i_0i_1}&=S_\sB(\al_{i_0i_1})=(\sC^\De_{\al_{i_0i_1}})^*\vert_{S_\sB([M_{i_0}])}, &
\fF^{i_1i_2}&=i_{M_{i_1}}^*\vert_{S_\sB([M_{i_1}])},\\
\fF^{i_1i_2}\ci\fF^{i_0i_1}&=(\sC^\De_{\al_{i_0i_1}}\ci i_{M_{i_1}})^*\vert_{S_\sB([M_{i_0}])}, &
\fF^{i_0i_2}&=i_{M_{i_0}}^*\vert_{S_\sB([M_{i_0}])}.
\end{align*}

Here $i_{M_{i_0}}$ and $\sC^\De_{\al_{i_0i_1}}\ci i_{M_{i_1}}$ are unital $A_\iy$-functors $\sP^\un\ra\sC^\De_{[M_{i_0}]}$, where $i_{M_{i_0}}$ maps $*\mapsto M_{i_0}$ and $\sC^\De_{\al_{i_0i_1}}\ci i_{M_{i_1}}$ maps $*\mapsto \al_{i_0i_1}(M_{i_1})$ on objects. There is a unital $A_\iy$-natural transformation $\eta^{i_0i_1i_2}:\sC^\De_{\al_{i_0i_1}}\ci i_{M_{i_1}}\Ra i_{M_{i_0}}$ with $(\eta^{i_0i_1i_2})_0^*=\om_{\al_{i_0i_1}(M_{i_1}),M_{i_0}}:\al_{i_0i_1}(M_{i_1})\ra M_{i_0}$, and $(\eta^{i_0i_1i_2})_k^{*,\ldots,*}=0$ for $k>0$, using the notation of Definition \ref{ai4def20}. We regard $\eta^{i_0i_1i_2}$ as a morphism in~$\Fun_{A_\iy}^\un(\sP^\un,\sC^\De_{[M_{i_0}]})$. 

For $k>0$ define $\m_k^{\{i_0,i_1,i_2\};\fG_0,\ldots,\fG_k}=0$ for all objects $\fG_0,\ldots,\fG_k$ in $S_\sB([M_{i_0}])$. When $k=0$, define
\begin{align*}
&\m_0^{\{i_0,i_1,i_2\};\fG_0}:R\longra \Hom_{\Fun_{A_\iy}^\un(\sP^\un,\sB)}^0\bigl(\fG_0\ci \sC^\De_{\al_{i_0i_1}}\ci i_{M_{i_1}},\fG_0\ci i_{M_{i_0}}\bigr) \\
&\text{by}\quad \m_0^{\{i_0,i_1,i_2\};\fG_0}(1)=(\fG_0)_*(\eta^{i_0i_1i_2}).
\end{align*}
We can now check that this defines a 2-simplex in $\hACat^{\un,\iy}_R$. Moreover, as $\fG_0$ lies in $S_\sB([M_{i_0}])\subset\Fun_{A_\iy}^\un(\sC^\De_{[M_{i_0}]},\sB)$ it maps $\om_{\al_{i_0i_1}(M_{i_1}),M_{i_0}}$ to a homotopy equivalence in $\sB$. This implies the condition in Definition \ref{ai7def4} for this to define a 2-simplex in $\ACat^{\un,\iy}_R\subset\hACat^{\un,\iy}_R$.
\item[(x)] Suppose that $i_0<b\le i_1<a\le i_2$, so that $F_{a,b}$ maps $\De^2$ surjectively onto the 2-simplex with vertices $(0,0),(0,1),(1,1)$. Then
\begin{align*}
\fF^{i_0i_1}&=S_\sB(\al_{i_0i_1})=(\sC^\De_{\al_{i_0i_1}})^*\vert_{S_\sB([M_{i_0}])},\\
\fF^{i_1i_2}&=S_\sB(\al_{i_1i_2})=(\sC^\De_{\al_{i_1i_2}})^*\vert_{S_\sB([M_{i_1}])},\\
\fF^{i_1i_2}\ci\fF^{i_0i_1}&=S_\sB(\al_{i_0i_2})=(\sC^\De_{\al_{i_0i_2}})^*\vert_{S_\sB([M_{i_0}])},\\
\fF^{i_0i_2}&=(i_{M_{i_0}}\ci\pi_{M_{i_2}})^*\vert_{S_\sB([M_{i_0}])}.
\end{align*}

Here $\sC^\De_{\al_{i_0i_2}}$ and $i_{M_{i_0}}\ci\pi_{M_{i_2}}$ are unital $A_\iy$-functors $\sC^\De_{[M_{i_2}]}\ra\sC^\De_{[M_{i_0}]}$, where in fact $i_{M_{i_0}}\ci\pi_{M_{i_2}}=\sC^\De_\be$ for $\be:[M_{i_2}]\ra[M_{i_0}]$ the morphism in $\De$ mapping $\be:j\mapsto M_{i_0}$ for all $j$. There is a unital $A_\iy$-natural transformation $\ze^{i_0i_1i_2}:\sC^\De_{\al_{i_0i_2}}\Ra i_{M_{i_0}}\ci\pi_{M_{i_2}}$ with $(\eta^{i_0i_1i_2})_0^{j_0}=\om_{\al_{i_0i_2}(j_0),M_{i_0}}:\al_{i_0i_2}(j_0)\ra M_{i_0}$ for $j_0=0,1,\ldots,M_{i_2}$, and $(\ze^{i_0i_1i_2})_k^{j_0,\ldots,j_k}=0$ for $k>0$. We regard $\ze^{i_0i_1i_2}$ as a morphism in~$\Fun_{A_\iy}^\un(\sC^\De_{[M_{i_2}]},\sC^\De_{[M_{i_0}]})$. 

For $k>0$ define $\m_k^{\{i_0,i_1,i_2\};\fG_0,\ldots,\fG_k}=0$ for all objects $\fG_0,\ldots,\fG_k$ in $S_\sB([M_{i_0}])$. When $k=0$, define
\begin{align*}
&\m_0^{\{i_0,i_1,i_2\};\fG_0}:R\longra \Hom_{\Fun_{A_\iy}^\un(\sC^\De_{[M_{i_2}]},\sB)}^0\bigl(\fG_0\ci \sC^\De_{\al_{i_0i_2}},\fG_0\ci i_{M_{i_0}}\ci\pi_{M_{i_2}}\bigr) \\
&\text{by}\quad \m_0^{\{i_0,i_1,i_2\};\fG_0}(1)=(\fG_0)_*(\ze^{i_0i_1i_2}).
\end{align*}
The argument in (ix) shows this gives a 2-simplex in~$\ACat^{\un,\iy}_R$.
\end{itemize}

This defines all the data in $T=\Xi_\sB(S,\si_a,\si_b)$. To show $T$ is an $N$-simplex in $\hACat^{\un,\iy}_R$, we must show \eq{ai7eq15}--\eq{ai7eq19} hold. Here \eq{ai7eq15} holds as $\sA^i$ is an $A_\iy$-category, and \eq{ai7eq17}--\eq{ai7eq19} hold as $\sA^i,\fF^{ij},\eta^{i_0i_1i_2},\ze^{i_0i_1i_2}$ are unital and $\m^I=0$ if $\md{I}\ge 4$. Also \eq{ai7eq16} holds for $M=1$ as the $\fF^{ij}$ are $A_\iy$-functors, and for $M=2$ as $(\fG_0)_*(\eta^{i_0i_1i_2}),(\fG_0)_*(\ze^{i_0i_1i_2})$ are $A_\iy$-natural transformations.

Thus we need to prove \eq{ai7eq16} holds for $M\ge 3$. Then the first line of \eq{ai7eq16} is zero as $\m_k^{I;y_0,\ldots,y_k}=0$. For the second line, $\d_{n,k;y_0,\ldots,y_k}^{I;x_0,\ldots,x_n}$ is defined in \eq{ai7eq9} as a sum over trees $\bs\Ga_1,\ldots,\bs\Ga_k$. Since $\m^{I;w_0,\ldots,w_k}_k=0$ if either $\md{I}\ge 4$, or if $\md{I}=3$ and $k>0$, any such $\bs\Ga_1,\ldots,\bs\Ga_k$ can only yield a nonzero term if $\md{I_v}=2$ or 3 for all internal vertices $v$ in each $\bs\Ga_i$, and if $\md{I_v}=3$ then $v$ has no incoming edges. Also, in each of (ii)--(vi) above we can show that if $\md{I}=2$ then $\m^{I;w_0,\ldots,w_k}_k=0$ unless $k=1$. Thus, if $\md{I_v}=2$ then $v$ must have exactly one incoming edge. Therefore any $\bs\Ga_i$ appearing in a nonzero term must be either of the form
\e
\xymatrix@C=25pt{ \overset{\{i_{a_0},i_{a_1},i_{a_2}\}}{\bu}  \ar[r]  &  \overset{\{i_{a_2},i_{a_3}\}}{\bu}  \ar[r] & \cdots \ar[r] & \overset{\{i_{a_{k-1}},i_{a_k}\}}{\bu}  \ar[r] & \tr }
\label{ai11eq40}
\e
for $0\le a_0<a_1<\cdots<a_k=M$, or of the form
\e
\xymatrix@C=25pt{ \ci  \ar[r]  &  \overset{\{i_{a_0},i_{a_1}\}}{\bu}  \ar[r] & \cdots \ar[r] & 
\overset{\{i_{a_{k-1}},i_{a_k}\}}{\bu}  \ar[r] & \tr }
\label{ai11eq41}
\e
for $0=a_0<a_1<\cdots<a_k=M$. Furthermore, the conditions on $L(\bs\Ga_i),R(\bs\Ga_i)$ in \eq{ai7eq9} imply there are exactly $M-1$ $\bs\Ga_i$ of the form \eq{ai11eq40}, so that~$k\ge M-1$.

By considering the last two $\bs\Ga_i$ of the form \eq{ai11eq40} in the list $\bs\Ga_1,\ldots,\bs\Ga_k$, we can see that for some $0<c<d<M$, either $\{i_0,i_c,i_d\}$ and $\{i_0,i_d,i_M\}$ appear as $I_v,I_w$ for vertices $v,w$ in $\bs\Ga_1,\ldots,\bs\Ga_k$, or $\{i_c,i_d,i_M\}$ and $\{i_0,i_c,i_M\}$ appear as $I_v,I_w$. By considering how $c,d$ may be ordered in relation to $a,b$, we can show that in both cases, at least one of the two must be of the form (vii) or (viii) above, and so correspond to an identity 2-morphism $\id_{\fF^{i_0i_2}}$. It follows that some $\bs\Ga_i$ must have $C(\bs\Ga_i)^{x_r}=\boo^{i_M}_{y_i}$ for some object $y_i$ in $\sA^{i_M}$. If $k\ne 2$ then composing with $\m_k^{\{i_M\};y_0,\ldots,y_k}$ in \eq{ai7eq16} gives zero by~\eq{ai7eq19}.

Therefore $\m_k^{\{i_M\};y_0,\ldots,y_k}\ci\d_{n,k;y_0,\ldots,y_k}^{I;x_0,\ldots,x_n}=0$ in \eq{ai7eq16} unless $k=2$. Since $k\ge M-1$ for $\d_{n,k;y_0,\ldots,y_k}^{I;x_0,\ldots,x_n}$ to have nonzero terms, this proves \eq{ai7eq16} if $M\ge 4$. This leaves the case $M=3$ and $k=2$. Then there are only two possibilities for $\bs\Ga_1,\bs\Ga_2$, firstly:
\e
\begin{gathered}
\xymatrix@!0@C=25pt@R=25pt{ 
&& \overset{\{i_0,i_1,i_2\}}{\bu}  \ar[rr]  &&  \overset{\{i_2,i_3\}}{\bu}  \ar[rr] && \tr & *+[r]{\bs\Ga_1,} \\
&&&& \overset{\{i_0,i_2,i_3\}}{\bu} \ar[rr] && \tr & *+[r]{\bs\Ga_2,} 
 }
\end{gathered}
\label{ai11eq42}
\e
and secondly:
\e
\begin{gathered}
\xymatrix@!0@C=25pt@R=25pt{ 
&&&& \overset{\{i_1,i_2,i_3\}}{\bu} \ar[rr] && \tr & *+[r]{\bs\Ga_1,} \\
&&&& \overset{\{i_0,i_1,i_3\}}{\bu} \ar[rr] && \tr & *+[r]{\bs\Ga_2.}  }
\end{gathered}
\label{ai11eq43}
\e

We claim that the contributions of \eq{ai11eq42} and \eq{ai11eq43} to \eq{ai7eq16} cancel, so that \eq{ai7eq16} holds with $M=3$. To see this, first note that the signs in \eq{ai7eq16} from \eq{ai11eq42} and \eq{ai11eq43} are opposite. This is because the signs are the same except for the contributions $(-1)^{N(v,w)}$ in \eq{ai7eq9}, for $N(v,w)$ as in \eq{ai6eq12}, where in \eq{ai11eq42} there is a contribution $-1$ to $(-1)^{N(v,w)}$ from $v$ the vertex in $\bs\Ga_1$ with $I_v=\{i_0,i_1,i_2\}$ and $w$ the vertex in $\bs\Ga_2$ with $I_w=\{i_0,i_2,i_3\}$, but there is no corresponding sign in~\eq{ai11eq43}.

Then we show that \eq{ai11eq42} and \eq{ai11eq43} make opposite contributions to \eq{ai7eq16} by a case by case analysis depending on the orders of $a,b$ relative to $i_0<i_2<i_2<i_3$, which determine the $\m_v^{I_v}$ for vertices $v$ in $\bs\Ga_1,\bs\Ga_2$ in \eq{ai11eq42} and \eq{ai11eq43} by (ii)--(x) above. We will explain the case $i_1\le a<i_2\le b<i_3$ as an example. Then \eq{ai11eq42} becomes
\e
\begin{gathered}
\xymatrix@!0@C=25pt@R=25pt{ 
&& \overset{\eta^{i_0i_1i_2}}{\bu}  \ar[rr]  &&  \overset{\pi_{M_3}^*}{\bu}  \ar[rr] && \tr & *+[r]{\bs\Ga_1,} \\
&&&& \overset{\id_{S_\sB(\al_{i_0i_3})}}{\bu} \ar[rr] && \tr & *+[r]{\bs\Ga_2,} }
\end{gathered}
\label{ai11eq44}
\e
and \eq{ai11eq43} becomes
\e
\begin{gathered}
\xymatrix@!0@C=25pt@R=25pt{ 
&&&& \overset{\id_{S_\sB(\al_{i_1i_3})}}{\bu} \ar[rr] && \tr & *+[r]{\bs\Ga_1,
}
\\
&&&& \overset{\id_{S_\sB(\al_{i_0i_3})}}{\bu} \ar[rr] && \tr & *+[r]{\bs\Ga_2.} }
\end{gathered}
\label{ai11eq45}
\e
Here in \eq{ai11eq44}, since $\pi_{M_3}$ maps $\om_{\al_{i_0i_1}(M_{i_1}),M_{i_0}}$ in the definition of $\eta^{i_0i_1i_2}$ in (ix) to $\id_*$, the composition $\pi_{M_3}^*\ci \eta^{i_0i_1i_2}$ is $\id_{S_\sB(\al_{i_0i_3})}$, which in fact agrees with $\id_{S_\sB(\al_{i_1i_3})}$ in \eq{ai11eq45} because of the conventions on `prefixes' in \S\ref{ai61} and \S\ref{ai71}. Hence \eq{ai11eq44}--\eq{ai11eq45} give cancelling terms in this case.

This proves that $T=\Xi_\sB(S,\si_a,\si_b)$ is an $N$-simplex in $\hACat^{\un,\iy}_R$. Since all the 2-simplices in (vii)--(x) lie in $\ACat^{\un,\iy}_R\subset\hACat^{\un,\iy}_R$, $T$ is an $N$-simplex in $\ACat^{\un,\iy}_R$. It is easy to see from the definitions that $T$ commutes with the action of morphisms $\al:[N]\ra[N']$ in $\De$. Thus we have constructed a map of simplicial sets $\Xi_\sB$ as in \eq{ai11eq36}. Parts (a)--(d) above are evident from the construction.
\end{dfn}

\subsubsection{A unique horn-filling type result in $\ACat^{\un,\iy}_R$}
\label{ai1142}

\begin{dfn}
\label{ai11def2}
Define $C_N$ to be the simplicial set illustrated in Figure \ref{ai11fig1}.
\begin{figure}[htb]
\centering
\includegraphics[width=4in]{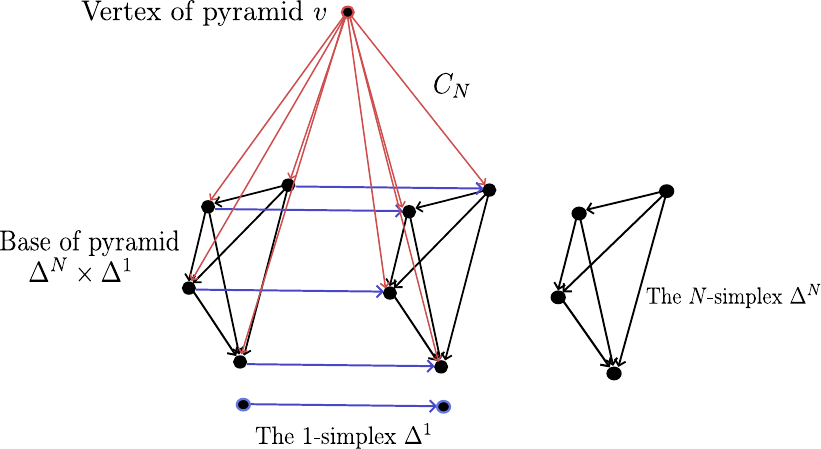}    
\caption{The simplicial set $C_N$}
\label{ai11fig1}
\end{figure}
It is a pyramid with base $\De^N\t\De^1$, and vertex $v$. Each $k$-simplex $f:\De^k\ra\De^N\t\De^1$ lifts to a $(k+1)$-simplex $\ti f:\De^{k+1}\ra C_N$ with $\ti f\ci\De^{\{0\}}=v$ and $\ti f\ci\De^{\{1,\ldots,k+1\}}=f$. Write $a_N:\De^{N+1}\ra C_N$ for the $(N+1)$-simplex lifting $\id_{\De^N}\t\De^{\{0\}}:\De^N\ra\De^N\t\De^1$ as above. Write $b_N:\De^{N+1}\ra C_N$ for the $(N+1)$-simplex lifting $\id_{\De^N}\t\De^{\{1\}}:\De^N\ra\De^N\t\De^1$. Write $c_N:\De^N\t\De^1\ra C_N$ for the inclusion of the base into $C_N$, so that $c_N=a_N\ci\De^{\{1,\ldots,N+1\}}$. Write $\pi_N:C_N\ra\De^{N+1}$ for the morphism mapping $v\mapsto 0$ and mapping the base $\De^N\t\De^1$ by the composition $\De^N\t\De^1\,{\buildrel\Pi_{\De^N}\over\longra}\,\De^N\,{\buildrel\De^{\{1,\ldots,N+1\}}\over\longra}\,\De^{N+1}$. 

Let $f:\De^k\ra C_N$ be a $k$-simplex in $C_N$ containing the vertex $v$. We call $f$ {\it flat\/} if $\pi_N\ci f:\De^k\ra\De^{N+1}$ is a degenerate $k$-simplex in~$\De^{N+1}$.
\end{dfn}

\begin{lem}
\label{ai11lem1}
Suppose we are given a commutative diagram of morphisms `$\,\ra$' of simplicial sets
\e
\begin{gathered}
\xymatrix@C=80pt@R=15pt{
*+[r]{\De^N} \ar[r]^{\De^{\{1,2,\ldots,N+1\}}} \ar[d]^{\id_{\De^N}\t\De^{\{0\}}} & \De^{N+1} \ar[d]_{a_N} \ar@/^.5pc/[ddr]^U
\\
*+[r]{\De^N\t\De^1} \ar@/_.5pc/[drr]_V \ar[r]^(0.6){c_N} & C_N \ar@{..>}[dr]^W
\\
&& \ACat^{\un,\iy}_R. }
\end{gathered}
\label{ai11eq46}
\e
Then there is a unique morphism $W:C_N\ra\ACat^{\un,\iy}_R$ making \eq{ai11eq46} commute, such that if\/ $f:\De^k\ra C_N$ is a \begin{bfseries}flat\end{bfseries} $k$-simplex in $C_N$ containing the vertex $v$ for $k\ge 2,$ then $W\ci f$ is a \begin{bfseries}special\end{bfseries}\/ $k$-simplex in $\ACat^{\un,\iy}_R,$ in the sense of Definition\/~{\rm\ref{ai7def5}}.
\end{lem}

\begin{proof}
Since $W$ is defined on the base $\De^N\t\De^1$ by $V$, it is sufficient to define $W$ on {\it nondegenerate\/} $(k+1)$-simplices in $C_N$ containing $v$, for $k=0,\ldots,N$. These correspond to nondegenerate $k$-simplices in the base $\De^N\t\De^1$. A nondegenerate $k$-simplex $g:\De^k\ra\De^N\t\De^1$ corresponds to choices of $0\le i\le k+1$ and $0\le a_0<a_1<\cdots<a_{i-1}\le a_i<a_{i+1}<\cdots<a_k\le N$, such that $g$ maps vertex $j$ of $\De^k$ for $k=0,\ldots,k$ to vertex $(a_j,0)$ of $\De^N\t\De^1$ if $0\le j<i$ and to vertex $(a_j,1)$ of $\De^N\t\De^1$ if $i\le j\le k$. Write $g_{\bs a}^i$ for this $k$-simplex, with $\bs a=(a_0,a_1,\ldots,a_k)$. Write $h_{\bs a}^i:\De^{k+1}\ra C_N$ for the nondegenerate $(k+1)$-simplex with $h_{\bs a}^i\ci\De^{\{0\}}=v$ and $h_{\bs a}^i\ci\De^{\{1,\ldots,k+1\}}=g_{\bs a}^i$. Then $h_{\bs a}^i$ is {\it flat\/} if and only if $a_{i-1}=a_i$, since then the vertices $(a_{i-1},0)$, $(a_i,1)$ in $g_{\bs a}^i(\De^k)\subset\De^N\t\De^1\subset C_N$ are mapped by $\pi_N$ to the same vertex $a_i+1$ in~$\De^{N+1}$.

We will construct $W$ on nondegenerate $(k+1)$-simplices including $v$ by induction on $k=0,1,\ldots,N$. Our inductive hypothesis, for $k=0,1,\ldots,N$, is:
\begin{itemize}
\setlength{\itemsep}{0pt}
\setlength{\parsep}{0pt}
\item[$(*)_k$] Suppose we have chosen $W\ci f$ uniquely on all nondegenerate $l$-simplices $f:\De^l\ra C_N$ with $l\le k$ including vertex $v$ if $f$ is not flat, and with $l\le k+1$ if $f$ is flat, and that these choices are consistent on restricting to subsimplices $\De^{l'}\hookra\De^l$, and with the values given by $U$ or $V$ in \eq{ai11eq46} on restricting to subsimplices $\De^{l'}\hookra\De^l$ which factor via $a_N$ or~$c_N$.
\end{itemize}
For the first step $(*)_0$ we have $W(v)=U\ci\De^{\{0\}}$, where $v$ is the only 0-simplex containing $v$, and there are no nondegenerate flat 1-simplices.

For the inductive step, suppose $(*)_k$ holds for some $k=0,\ldots,N$. To prove $(*_k)$, suppose $0\le a_0<a_1<\cdots<a_k\le N$, and write $\bs a=(a_0,\ldots,a_k)$. Also, for $i=0,\ldots,k$, define $\bs a_i=(a_0,\ldots,a_{i-1},a_i,a_i,a_{i+1},\ldots,a_k)$. Then we have non-flat $(k+1)$-simplices $h_{\bs a}^i:\De^{k+1}\ra C_N$ for $i=0,\ldots,k+1$, and flat $(k+2)$-simplices $h_{{\bs a}_i}^{i+1}:\De^{k+2}\ra C_N$ for $i=0,\ldots,k$. Observe that:
\begin{itemize}
\setlength{\itemsep}{0pt}
\setlength{\parsep}{0pt}
\item[(i)] When $i=k+1$, the map $h_{\bs a}^{k+1}:\De^{k+1}\ra C_N$ factors via $a_N:\De^{N+1}\ra C_N$. Thus $W\ci h_{\bs a}^{k+1}$ is determined uniquely by $W\ci a_N=U$ in~\eq{ai11eq46}.
\item[(ii)] Consider the action of $h_{{\bs a}_i}^{i+1}:\De^{k+2}\ra C_N$ on the $k+3$ boundary faces $\De^{k+1}$ in $\De^{k+2}$. The boundary face $\De^{\{0,\ldots,i,i+2,\ldots,k+2\}}(\De^{k+1})$, omitting vertex $i+1$ in $\De^{k+2}$, has $h_{{\bs a}_i}^{i+1}\ci \De^{\{0,\ldots,i,i+2,\ldots,k+2\}}=h_{\bs a}^{i+1}$. The boundary face $\De^{\{0,\ldots,i+1,i+3,\ldots,k+2\}}(\De^{k+1})$, omitting vertex $i+2$ in $\De^{k+2}$, has $h_{{\bs a}_i}^{i+1}\ci \De^{\{0,\ldots,i,i+2,\ldots,k+2\}}=h_{\bs a}^i$. The bottom boundary face $\De^{\{1,\ldots,k+2\}}(\De^{k+1})$, omitting vertex 0, has $W\ci h_{{\bs a}_i}^{i+1}\ci\De^{\{1,\ldots,k+2\}}$ determined by $V$.

All other boundary faces $h_{{\bs a}_i}^{i+1}\ci \De^{\{0,\ldots,j,j+2,\ldots,k+2\}}$ for $j\ne i+1,i+2$ are flat, as they contain vertices $i+1,i+2$ of $\De^{k+2}$ which are mapped to the same point by $\pi_N$. Thus, $W\ci h_{{\bs a}_i}^{i+1}\ci \De^{\{0,\ldots,j,j+2,\ldots,k+2\}}$ for $j\ne i+1,i+2$ is determined by the inductive hypothesis $(*)_k$.
\end{itemize}

We now define $W\ci h_{\bs a}^i$ and $W\ci h_{{\bs a}_i}^{i+1}$ by reverse induction on $i=k,k-1,\ldots,0$. The inductive hypothesis is that for some $k\ge j>0$ we have defined $W\ci h_{\bs a}^i$ and $W\ci h_{{\bs a}_i}^{i+1}$ for all $i>j$. When $j=k$ this is true by (i). For the inductive step, consider the diagram
\e
\begin{gathered}
\xymatrix@C=90pt@R=15pt{ \La^{k+2}_{j+1_{\vphantom{(}}} \ar@{^{(}->}[d]^\inc \ar@/^.5pc/[dr]^{W\ci h_{{\bs a}_j}^{j+1}\ci\inc} \\
\De^{k+2} \ar@{..>}[r]^{W\ci h_{{\bs a}_j}^{j+1}}  & \ACat^{\un,\iy}_R. }
\end{gathered}
\label{ai11eq47}
\e
Here $h_{{\bs a}_j}^{j+1}(\La^{k+2}_{j+1})\subset C_N$ is the union of $\De^{\{0,\ldots,l,l+2,\ldots,k+2\}}(\De^{k+1})$ for $l\ne j+1$. By (ii) this consists of $h_{\bs a}^{j+1}(\De^{k+1})$, on which $W$ is defined by induction, the base $(k+1)$-simplex on which $W$ is determined by $V$, and $k+1$ other $(k+1)$-simplices which are flat, so that $W$ is defined on them by $(*)_k$. Thus, we have already defined $W$ on $h_{{\bs a}_j}^{j+1}(\La^{k+2}_{j+1})$, so $W\ci h_{{\bs a}_j}^{j+1}\ci\inc$ in \eq{ai11eq47} is defined. Proposition \ref{ai7prop2} now says that we can fill \eq{ai11eq47} with a unique $(k+2)$-simplex $W\ci h_{{\bs a}_j}^{j+1}$ which is special. This defines $W$ on the flat $(k+2)$-simplex $h_{{\bs a}_j}^{j+1}(\De^{k+2})$ in~$C_N$.

Since $h_{{\bs a}_j}^{j+1}\ci \De^{\{0,\ldots,j,j+2,\ldots,k+2\}}=h_{\bs a}^j$, defining $W\ci h_{{\bs a}_j}^{j+1}$ restricts to a definition of $W\ci h_{\bs a}^j=W\ci h_{{\bs a}_j}^{j+1}\ci \De^{\{0,\ldots,j,j+2,\ldots,k+2\}}$. This completes the inductive step. Hence we can uniquely define $W\ci h_{\bs a}^i$ and $W\ci h_{{\bs a}_i}^{i+1}$ for all $i$. All nondegenerate, non-flat $(k+1)$-simplices in $C_N$ including $v$ are of the form $h_{\bs a}^i$, and all nondegenerate, flat $(k+2)$-simplices in $C_N$ including $v$ are of the form $h_{{\bs a}_i}^{i+1}$. Thus we have proved $(*)_{k+1}$. So by induction, $(*)_k$ holds for all $k=0,1,\ldots,N$.

We have now defined $W$ on all nondegenerate simplices in $C_N$, compatibly with restriction to subsimplices. Thus we have defined a morphism $W$, which makes \eq{ai11eq46} commute by $(*)_k$. Each choice we made in the construction was forced upon us uniquely by the condition that \eq{ai11eq46} commutes and $W$ takes flat simplices to special simplices. This proves Lemma~\ref{ai11lem1}.
\end{proof}

\subsubsection{Outline of the proof of Proposition \ref{ai11prop1}}
\label{ai1143}

We can now sketch the main idea in the proof of Proposition \ref{ai11prop1}, which we will carry out in \S\ref{ai1144}--\S\ref{ai1149}. Section \ref{ai1141} defined $\Xi_\sB:\cN(\De^\rop)\t\De^1\t\De^1\longra\ACat^{\un,\iy}_R$. Here is a diagram of $\De^1\t\De^1$ in the domain of $\Xi_\sB$:
\e
\begin{gathered}
\xymatrix@C=80pt@R=15pt{
*+[r]{\overset{(0,1)}{\bu}} \ar[r]_H & *+[r]{\overset{(1,1)}{\bu}} \\
*+[r]{\overset{(0,0)}{\bu}} \ar[r]^{a'} \ar[u]_c & *+[r]{\overset{(1,0)}{\bu}} \ar[u]^{b'}
}
\end{gathered}
\label{ai11eq48}
\e
Given an $N$-simplex $(S,T):\De^N\ra \sN(S_\sB)^*(\sN(\ACat^\un_R)_{\sA/})$, so that $p\ci(S,T)=S:\De^N\ra\cN(\De^\rop)$, we will define $a'(S,T)$ by applying Lemma \ref{ai11lem1} with $U$ determined by $T$ and $V$ determined by the restriction of $\Xi_\sB$ to $S(\De^N)\t\De^1\t\De^{\{0\}}$, giving a morphism $W$ in \eq{ai11eq46}, and then $a'(S,T)$ is determined by~$W\ci b_N$. 

Thus $a'$ is constructed by a kind of `parallel transport' of $N$-simplices along the edge $(0,0)\ra(1,0)$ of $\De^1\t\De^1$, using $\Xi_\sB$ and Lemma \ref{ai11lem1}. Similarly, $b'$ is constructed by `parallel transport' along the edge $(1,0)\ra(1,1)$, and $c$ by `parallel transport' along the edge~$(0,0)\ra(0,1)$.

We then define $H$ in \eq{ai11eq31} in a different way. Given an $N$-simplex $(S,T):\De^N\ra \sN(S_\sB)^*(\sN(\ACat^\un_R)_{\sA/})$, we will construct a map $K:\De^N\t\De^1\t\De^1\longra(\ACat^{\un,\iy}_R)_{\sA/}$. Here the values of $K$ over the three edges $(0,0)\ra(1,0)$, $(1,0)\ra(1,1)$, $(0,0)\ra(0,1)$ of $\De^1\t\De^1$ are determined by the constructions of $a'(S,T),b'\ci a'(S,T)$, and $c(S,T)$. We fill in $K$ over the rest of $\De^1\t\De^1$ uniquely using special $(N+1)$-simplices in $\ACat^{\un,\iy}_R$. The restriction of $K$ to the edge $(0,1)\ra(1,1)$ of $\De^1\t\De^1$ then determines~$H$.

\subsubsection{Definition of $a'$ in \eq{ai11eq30}}
\label{ai1144}

Let $(S,T)$ be an $N$-simplex in $\sN(S_\sB)^*(\sN(\ACat^\un_R)_{\sA/})$ for $N\ge 0$, so that $S$ is an $N$-simplex in $\sN(\De^\rop)$ and $T$ is an $N$-simplex in $\sN(\ACat^\un_R)_{\sA/}$ with $p_{\sN(\ACat^\un_R)}(T)=\sN(S_\sB)(S)$. By Definition \ref{ai2def8}, $T$ corresponds to $\ti T:\De^{N+1}\ra\ACat^{\un,\iy}_R$ such that $\ti T\vert_{\De^{\{0\}}}$ is the the 0-simplex $\sA$, and $\ti T\vert_{\De^{\{1,2,\ldots,N+1\}}}=\sN(S_\sB)\ci S$. 

Consider the following diagram of simplicial sets:
\e
\begin{gathered}
\xymatrix@C=80pt@R=15pt{
*+[r]{\De^N} \ar[r]^{\De^{\{1,2,\ldots,N+1\}}} \ar[d]^{\id_{\De^N}\t\De^{\{0\}}} & \De^{N+1} \ar[d]_{a_N} \ar@/^.5pc/[ddr]^{\ti T}
\\
*+[r]{\De^N\t\De^1} \ar@/_.5pc/[drr]_{\Xi_\sB\ci(S\t \id_{\De^1}\t\De^{\{0\}})}
 \ar[r]^(0.6){c_N} & C_N \ar@{..>}[dr]^W
\\
&& \ACat^{\un,\iy}_R. }
\end{gathered}
\label{ai11eq49}
\e
This is \eq{ai11eq46} with $U=\ti T$ and $V=\Xi_\sB\ci(S\t \id_{\De^1}\t\De^{\{0\}})$. Here for $V$ we first map $\De^N\t\De^1\ra \sN(\De^\top)\t\De^1\t\De^1$ by the product of $S:\De^N\ra \sN(\De^\top)$ and $\id_{\De^1}\t\De^{\{0\}}:\De^1\ra\De^1\t\De^1$ which includes $\De^1$ as the edge $(0,0)\ra(1,0)$ in $\De^1\t\De^1$, as in \eq{ai11eq48}, and then we compose with $\Xi_\sB$ in~\eq{ai11eq36}.

By definition of $\ti T$ and $\Xi_\sB$, the two compositions $\De^N\ra\ACat^{\un,\iy}_R$ round the outside of \eq{ai11eq49} are $S_\sB\ci S$. Hence the diagram of morphisms `$\ra$' in \eq{ai11eq49} commutes, so Lemma \ref{ai11lem1} gives a unique morphism $W$ as shown making \eq{ai11eq49} commute and taking flat $k$-simplices to special $k$-simplices. 

Define $\ti F=W\ci b_N$, for $b_N$ as in Definition \ref{ai11def2}. Then $\ti F:\De^{N+1}\ra\ACat^{\un,\iy}_R$ is such that $\ti F\vert_{\De^{\{0\}}}$ is the the 0-simplex $\sA$, and $\ti T\vert_{\De^{\{1,2,\ldots,N+1\}}}$ is the constant $N$-simplex at $\sB$, by Definition \ref{ai11def1}(b). Hence $\ti F$ corresponds to a map $F:\De^N\ra \Hom^L_{\ACat^{\un,\iy}_R}(\sA,\sB)$. Set $a'(S,T)=(S,F)$. This defines $a'$ on all $N$-simplices $(S,T)$ in $\sN(S_\sB)^*(\sN(\ACat^\un_R)_{\sA/})$. One can show that the definition of $a'(S,T)$ commutes with the action of morphisms in $\De^\rop$, so that $a'$ is a morphism of simplicial sets. 

\subsubsection{Alternative description of $a'$, proof of Proposition \ref{ai11prop1}(a)}
\label{ai1145}

We will define $a'$ in \S\ref{ai1144} in a different way. Use the notation of Definition \ref{ai9def1}, Theorem \ref{ai9thm2}, and Definition \ref{ai2def18} for $(S_\sB,i_\sB)$. For $[N]\in\sN(\De^\rop)$ we have
\ea
p^{-1}([N])&=\Hom^L_{\sN(\ACat^\un_R)}(\sA,S_\sB([N]))\cong\Hom_{\ACat^\un_R}(\sA,S_\sB([N]))
\nonumber\\
&=F^\un_{S_\sB}(\sA)([N])\cong \sN_{A_\iy}(\Fun_{A_\iy}^\un(\sA,\sB))^\Kan([N])
\nonumber\\
&\cong \Hom^L_{\ACat^{\un,\iy}_R}(\sA,\sB)([N]).
\label{ai11eq50}
\ea
Here we use Theorem \ref{ai9thm2}(c) in the fourth step and Theorem \ref{ai9thm1}(a) in the fifth.

Write $\ul\De:\De^\rop\ra\SSet^\rop$ for the (opposite of the) usual functor which maps $[N]\mapsto\De^N$, and $\sN(\ul\De):\sN(\De^\rop)\ra\sN(\SSet^\rop)$ for its nerve. As for $\sC_{X/}$ in Example \ref{ai2ex7}, the {\it undercategory\/} of the object $\Hom^L_{\ACat^{\un,\iy}_R}(\sA,\sB)$ in $\sN(\SSet^\rop)$ is a left fibration
\begin{equation*}
p_{\sN(\SSet^\rop)}:\sN(\SSet^\rop)_{\Hom^L_{\ACat^{\un,\iy}_R}(\sA,\sB)/}\longra\sN(\SSet^\rop).
\end{equation*}
The pullback of this to $\sN(\De^\rop)$ under $\sN(\ul\De)$ is a left fibration
\begin{equation*}
p'\!=\!\sN(\ul\De)^*(p_{\sN(\SSet^\rop)})\!:\!\sN(\ul\De)^*(\sN(\SSet^\rop)_{\Hom^L_{\ACat^{\un,\iy}_R}(\sA,\sB)/})\!\ra\!\sN(\De^\rop).
\end{equation*}
Then for $[N]\in\sN(\De^\rop)$ we have 
\ea
(p')^{-1}([N])&\cong\Hom_{\SSet^\rop}(\Hom^L_{\ACat^{\un,\iy}_R}(\sA,\sB),\ul\De([N]))
\nonumber\\
&\cong\Hom_{\SSet}(\De^N,\Hom^L_{\ACat^{\un,\iy}_R}(\sA,\sB))
\nonumber\\
&\cong\Hom^L_{\ACat^{\un,\iy}_R}(\sA,\sB)([N])\cong p^{-1}([N]),
\label{ai11eq51}
\ea
using \eq{ai11eq50} in the last step.

By showing that all the isomorphisms in \eq{ai11eq50}--\eq{ai11eq51} are functorial in $[N]\in\sN(\SSet^\rop)$, we can construct an isomorphism of simplicial sets in a strictly commutative diagram:
\begin{equation*}
\xymatrix@C=120pt@R=17pt{
*+[r]{\sN(S_\sB)^*(\sN(\ACat^\un_R)_{\sA/})} \ar@/^1pc/[dr]^(0.6)d_(0.6)\cong \ar[d]^p   \\
*+[r]{\sN(\De^\rop)} & *+[l]{\sN(\ul\De)^*(\sN(\SSet^\rop)_{\Hom^L_{\ACat^{\un,\iy}_R}(\sA,\sB)/}).} \ar[l]_(0.75){p'}  }
\end{equation*}
Next we define a morphism of simplicial sets
\e
e\!:\!\sN(\ul\De)^*(\sN(\SSet^\rop)_{\Hom^L_{\ACat^{\un,\iy}_R}(\sA,\sB)/})\!\ra\! \Hom^L_{\ACat^{\un,\iy}_R}(\sA,\sB),
\label{ai11eq52}
\e
using the `last vertex map' construction of Lurie \cite[\S 3.3.4]{Luri5}. An $N$-simplex in the left hand side of \eq{ai11eq52} is a commutative diagram in $\SSet^\rop$:
\ea
\label{ai11eq53}\\[-21pt]
\xymatrix@!0@C=75pt@R=25pt{ && \Hom^L_{\ACat^{\un,\iy}_R}(\sA,\sB) \ar@/_.4pc/[dll]_(0.7){\si_0} \ar[dl]_(0.7){\si_1} \ar@<-1.5ex>@{}[d]^(0.6){\cdots} \ar@/^.4pc/[drr]^(0.7){\si_N} \ar[dr]^(0.7){\si_{N-1}} \\
 \De^{M_0} \ar[r]^(0.6){\al_1} & \De^{M_1} \ar[r]^(0.55){\al_2} & \,\,\cdots\,\, \ar[r]^(0.45){\al_{N-1}} & \De^{M_{N-1}} \ar[r]^(0.4){\al_N} & \De^{M_N},  }
\nonumber
\ea
where $M_0,\ldots,M_N\ge 0$.  Equivalently, it is a commutative diagram in $\SSet$:
\ea
\xymatrix@!0@C=75pt@R=25pt{  \De^{M_N} \ar@/_.4pc/[drr]_(0.3){\si_N} \ar[r]_(0.6){\al_N} & \De^{M_{N-1}} \ar[dr]_(0.3){\si_{N-1}} \ar[r]_(0.55){\al_{N-1}} & \,\,\cdots\,\, \ar@<-1.5ex>@{}[d]^(0.4){\cdots} \ar[r]_(0.45){\al_2} & \De^{M_1} \ar[r]_(0.4){\al_1} \ar[dl]^(0.3){\si_1} & \De^{M_0} \ar@/^.4pc/[dll]^(0.3){\si_0} \\
&& \Hom^L_{\ACat^{\un,\iy}_R}(\sA,\sB),  }
\nonumber\\[-21pt]
\label{ai11eq54}
\ea
where $\al_i:\De^{M_i}\ra\De^{M_{i-1}}$ is equivalent to a morphism $\al_i:[M_i]\ra[M_{i-1}]$ in $\De$. As in \cite[Ex.~3.3.4.7]{Luri5}, we now define $e$ to map the $N$-simplex \eq{ai11eq54} to the composition $\De^N\,{\buildrel f\over\longra}\,\De^{M_0}\,{\buildrel \si_0\over\longra}\,\Hom^L_{\ACat^{\un,\iy}_R}(\sA,\sB)$, where $f$ maps vertex $i$ in $\De^N$ to the image of the last vertex $M^{N-i}$ in $\De^{M_{N-i}}$ under $\al_1\ci\al_2\ci\cdots\ci\al_{N-i}$. 

We now claim that $a'=(p',e)\ci d$, where $(p',e)$ is the direct product of $p'$ and $e$. To see this, note that an $N$-simplex $T$ in $\sN(S_\sB)^*(\sN(\ACat^\un_R)_{\sA/})$ is equivalent to a commutative diagram in~$\ACat^\un_R$:
\ea
\label{ai11eq55}\\[-21pt]
\xymatrix@!0@C=73pt@R=25pt{ && \sA \ar@/_.4pc/[dll]_(0.7){\tau_0} \ar[dl]_(0.7){\tau_1} \ar@<-1.5ex>@{}[d]^(0.6){\cdots} \ar@/^.4pc/[drr]^(0.7){\tau_N} \ar[dr]^(0.7){\tau_{N-1}} \\
S_\sB(\De^{M_0}) \ar[r]_(0.5){S_\sB(\al_1)} & S_\sB(\De^{M_1}) \ar[r]_(0.55){S_\sB(\al_2)} & \,\,\cdots\,\, \ar[r]_(0.35){S_\sB(\al_{N-1})} & S_\sB(\De^{M_{N-1}}) \ar[r]_(0.5){S_\sB(\al_N)} & S_\sB(\De^{M_N}).  }
\nonumber
\ea
This is mapped by $d$ to an $N$-simplex as in \eq{ai11eq53}--\eq{ai11eq54}, where morphisms $\si_i:\De^{M_i}\ra\Hom^L_{\ACat^{\un,\iy}_R}(\sA,\sB)$ in \eq{ai11eq54} are identified with $\tau_i:\sA\ra S_\sB(\De^{M_i})$ in \eq{ai11eq55} by~\eq{ai11eq50}.

We can now check that for $a'(T)=(S,F)$ in \S\ref{ai1144}, $F=e\ci d(T)$, where $d(T)$ is the $N$-simplex in \eq{ai11eq54}. Since $a'$ is defined using $V=\Xi_\sB\ci(S\t \id_{\De^1}\t\De^{\{0\}})$, it involves `parallel transport' along 1-simplices lying over the edge $(0,0)\mapsto(1,0)$ of $\De^1\t\De^1$, as in Definition \ref{ai11def1}(iv), which are defined using $i_{M_i}:\sP^\un\ra\sC^\De_{[M_i]}$ acting on objects by taking $*$ to the `last vertex' $M_i$ in $\sC^\De_{[M_i]}$, and this leads to the `last vertex map'.

This proves Proposition \ref{ai11prop1}(a), since $a'$ is the composition of a strict isomorphism $d$ with the `last vertex map' of Lurie~\cite[\S 3.3.4]{Luri5}.

\subsubsection{Definition of $b'$ in \eq{ai11eq30}}
\label{ai1146}

Next we define $b'$ in \eq{ai11eq30}. Let $N\ge 0$, and $(S,F)$ be an $N$-simplex in $\sN(\De^\rop)\t\Hom^L_{\ACat^{\un,\iy}_R}(\sA,\sB)$. We will construct an $N$-simplex $(S,G)$ in $(\iinc^\un_R\!\ci\sN(S_\sB))^*((\ACat^{\un,\iy}_R)_{\sA/})$, i.e.\ $G$ is an $N$-simplex in $(\ACat^{\un,\iy}_R)_{\sA/}$ with $p_{\ACat^{\un,\iy}_R}(G)=\iinc^\un_R\ci\sN(S_\sB)(S)$, and set~$b'(S,F)=(S,G)$.

By Definition \ref{ai2def8}, $F$ corresponds to $\ti F:\De^{N+1}\ra\ACat^{\un,\iy}_R$ such that $\ti F\vert_{\De^{\{0\}}}$ is the the 0-simplex $\sA$, and $\ti F\vert_{\De^{\{1,2,\ldots,N+1\}}}$ is the constant $N$-simplex at $\sB$. Similarly, $G$ corresponds to $\ti G:\De^{N+1}\ra\ACat^{\un,\iy}_R$ such that $\ti G\vert_{\De^{\{0\}}}$ is the the 0-simplex $\sA$, and $\ti G\vert_{\De^{\{1,2,\ldots,N+1\}}}$ is the $N$-simplex~$\iinc^\un_R\ci\sN(S_\sB)(S)$.

Consider the following diagram of simplicial sets:
\e
\begin{gathered}
\xymatrix@C=80pt@R=15pt{
*+[r]{\De^N} \ar[r]^{\De^{\{1,2,\ldots,N+1\}}} \ar[d]^{\id_{\De^N}\t\De^{\{0\}}} & \De^{N+1} \ar[d]_{a_N} \ar@/^.5pc/[ddr]^{\ti F}
\\
*+[r]{\De^N\t\De^1} \ar@/_.5pc/[drr]_{\Xi_\sB\ci(S\t \De^{\{1\}}\t\id_{\De^1})}
 \ar[r]^(0.6){c_N} & C_N \ar@{..>}[dr]^{\hat W}
\\
&& \ACat^{\un,\iy}_R. }
\end{gathered}
\label{ai11eq56}
\e
This is \eq{ai11eq46} with $U=\ti F$ and $V=\Xi_\sB\ci(S\t \De^{\{1\}}\t\id_{\De^1})$. Here for $V$ we first map $\De^N\t\De^1\ra \sN(\De^\top)\t\De^1\t\De^1$ by the product of $S:\De^N\ra \sN(\De^\top)$ and $\De^{\{1\}}\t\id_{\De^1}:\De^1\ra\De^1\t\De^1$ which includes $\De^1$ as the edge $(1,0)\ra(1,1)$ in $\De^1\t\De^1$, as in \eq{ai11eq48}, and then we compose with $\Xi_\sB$ in~\eq{ai11eq36}.

By definition of $\ti F$ and $\Xi_\sB$, the two compositions $\De^N\ra\ACat^{\un,\iy}_R$ round the outside of \eq{ai11eq56} are the constant $N$-simplex at $\sB$. Hence the diagram of morphisms `$\ra$' in \eq{ai11eq56} commutes, so Lemma \ref{ai11lem1} gives a unique morphism $\hat W$ as shown making \eq{ai11eq56} commute and taking flat $k$-simplices to special $k$-simplices. Define $\ti G=\hat W\ci b_N$, for $b_N$ as in Definition \ref{ai11def2}. This satisfies the conditions above, and so is equivalent to an $N$-simplex $G$ in $(\ACat^{\un,\iy}_R)_{\sA/}$ with $(S,G)$ an $N$-simplex in $(\iinc^\un_R\ci\sN(S_\sB))^*((\ACat^{\un,\iy}_R)_{\sA/})$. Set $b'(S,F)=(S,G)$. This defines $b'$ in \eq{ai11eq30} on all $N$-simplices $(S,F)$ in $\sN(\De^\rop)\t\Hom^L_{\ACat^{\un,\iy}_R}(\sA,\sB)$. One can show that the definition of $b'(S,F)$ commutes with the action of morphisms in $\De^\rop$, so that $b'$ is a morphism of simplicial sets. 

\subsubsection{Commuting triangles in \eq{ai11eq30}, for Proposition \ref{ai11prop1}(b)}
\label{ai1147}

The bottom left triangle of \eq{ai11eq30} commutes as in \eq{ai11eq24}. The top left triangle commutes since for an $N$-simplex $T$ in $\sN(S_\sB)^*(\sN(\ACat^\un_R)_{\sA/})$ we have $a'(S,T)=(S,F)$ with $S=p(S,T)$. The bottom right commutes as for an $N$-simplex $(S,F)$ in $\sN(\De^\rop)\t\Hom^L_{\ACat^{\un,\iy}_R}(\sA,\ab\sB)$ we have $b'(S,F)=(S,G)$ with~$\Pi_{\sN(\De^\rop)}(S,F)=S=r(S,G)$.

\subsubsection{Definition of $H$ in \eq{ai11eq31}, proof of Proposition \ref{ai11prop1}(b)}
\label{ai1148}

Here are analogues of Definition \ref{ai11def2} and Lemma \ref{ai11lem1}, proved in the same way.

\begin{dfn}
\label{ai11def3}
Define $D_N$ to be the simplicial set defined as for $C_N$ in Definition \ref{ai11def2}, but with base $\De^N\t\De^1\t\De^1$ rather than $\De^N\t\De^1$. Write $d_N:C_N\hookra D_N$ for the inclusion which takes vertex $v$ to vertex $v$, and takes the base $\De^N\t\De^1$ of $C_N$ to $\De^N\t\De^1\t\De^{\{0\}}$ in the base $\De^N\t\De^1\t\De^1$ of $D_N$. Similarly, define $e_N:C_N\hookra D_N$ to take $v$ to $v$ and take $\De^N\t\De^1$ to $\De^N\t\De^1\t\De^{\{1\}}$ in $\De^N\t\De^1\t\De^1$. Define $f_N:C_N\hookra D_N$ to take $v$ to $v$ and take $\De^N\t\De^{\{0\}}\t\De^1$ to $\De^N\t\De^1$ in $\De^N\t\De^1\t\De^1$. Define $g_N:C_N\hookra D_N$ to take $v$ to $v$ and take $\De^N\t\De^{\{1\}}\t\De^1$ to $\De^N\t\De^1$ in $\De^N\t\De^1\t\De^1$. Write $h_N:\De^N\t\De^1\t\De^1\ra D_N$ for the inclusion of the base. Write $\Pi_N:D_N\ra C_N$ for the morphism mapping $v\mapsto v$ and mapping the base $\De^N\t\De^1\t\De^1$ to $\De^N\t\De^1$ by projection to the first and second factors. 

Let $f:\De^k\ra D_N$ be a $k$-simplex in $D_N$ containing $v$. We call $f$ {\it flat\/} if $\Pi_N\ci f:\De^k\ra C_N$ is a degenerate $k$-simplex in~$C_N$.
\end{dfn}

\begin{lem}
\label{ai11lem2}
Suppose we are given a commutative diagram of morphisms `$\,\ra$' of simplicial sets
\e
\begin{gathered}
\xymatrix@C=80pt@R=15pt{
*+[r]{\De^N\t\De^1} \ar[r]^{c_N} \ar[d]^{\id_{\De^N}\t\id_{\De^1}\t\De^{\{0\}}} & C_N \ar[d]_{d_N} \ar@/^.5pc/[ddr]^{U'}
\\
*+[r]{\De^N\t\De^1\t\De^1} \ar@/_.5pc/[drr]_{V'} \ar[r]^(0.6){h_N} & D_N \ar@{..>}[dr]^{W'}
\\
&& \ACat^{\un,\iy}_R. }
\end{gathered}
\label{ai11eq57}
\e
Then there is a unique morphism $W':D_N\ra\ACat^{\un,\iy}_R$ making \eq{ai11eq57} commute, such that if\/ $f:\De^k\ra D_N$ is a \begin{bfseries}flat\end{bfseries} $k$-simplex in $D_N$ containing the vertex $v$ for $k\ge 2,$ then $W\ci f$ is a \begin{bfseries}special\end{bfseries}\/ $k$-simplex in $\ACat^{\un,\iy}_R,$ in the sense of Definition\/~{\rm\ref{ai7def5}}.
\end{lem}

We can now construct $H$ in \eq{ai11eq31}. Let $(S,T)$ be an $N$-simplex in $\sN(S_\sB)^*\ab (\sN(\ACat^\un_R)_{\sA/})$ for $N\ge 0$, so that $S$ is an $N$-simplex in $\sN(\De^\rop)$ and $T$ is an $N$-simplex in $\sN(\ACat^\un_R)_{\sA/}$ with $p_{\sN(\ACat^\un_R)}(T)=\sN(S_\sB)(S)$. Then $T$ corresponds to $\ti T:\De^{N+1}\ra\ACat^{\un,\iy}_R$ such that $\ti T\vert_{\De^{\{0\}}}=\sA$ and $\ti T\vert_{\De^{\{1,2,\ldots,N+1\}}}=\sN(S_\sB)\ci S$. In the definition of $a'(S,T)$ in \S\ref{ai1144}, we constructed a map $W:C_N\ra \ACat^{\un,\iy}_R$ in \eq{ai11eq49}. Consider the following diagram:
\e
\begin{gathered}
\xymatrix@C=80pt@R=15pt{
*+[r]{\De^N\t\De^1} \ar[r]^{c_N} \ar[d]^{\id_{\De^N}\t\id_{\De^1}\t\De^{\{0\}}} & C_N \ar[d]_{d_N} \ar@/^.5pc/[ddr]^W
\\
*+[r]{\De^N\t\De^1\t\De^1} \ar@/_.5pc/[drr]_{\Xi_\sB\ci(S\t \id_{\De^1}\t\id_{\De^1})} \ar[r]^(0.6){h_N} & D_N \ar@{..>}[dr]^{W'}
\\
&& \ACat^{\un,\iy}_R. }
\end{gathered}
\label{ai11eq58}
\e
This is \eq{ai11eq57} with $U'=W$ and $V'=\Xi_\sB\ci(S\t \id_{\De^1}\t\id_{\De^1})$. By \eq{ai11eq43}
\begin{equation*}
W\ci c_N\!=\!\Xi_\sB\ci(S\!\t\! \id_{\De^1}\!\t\!\De^{\{0\}})\!=\!\bigl(\Xi_\sB\ci(S\!\t\! \id_{\De^1}\!\t\!\id_{\De^1})\bigr)\!\ci\!\bigl(\id_{\De^N}\!\t\!\id_{\De^1}\!\t\!\De^{\{0\}}\bigr).
\end{equation*}
Hence the diagram of morphisms `$\ra$' in \eq{ai11eq58} commutes, so Lemma \ref{ai11lem2} gives a unique morphism $W'$ as shown making \eq{ai11eq58} commute and taking flat $k$-simplices to special $k$-simplices. 

Now $W'\ci e_N:C_N\ra\ACat^{\un,\iy}_R$ satisfies $W'\ci e_N(v)=\sA$, and
\begin{align*}
W'\ci e_N\ci c_N&=\Xi_\sB\ci(S\t \id_{\De^1}\t\De^{\{1\}})\\
&=\iinc^\un_R\ci\cN(S_\sB)\ci S\ci\Pi_{\De^N}:\De^N\t\De^1\longra \ACat^{\un,\iy}_R,
\end{align*}
using Definition \ref{ai11def1}(d) in the second step. Hence $W'\ci e_N$ is equivalent to a morphism $L:\De^N\t\De^1\ra (\ACat^{\un,\iy}_R)_{\sA/}$ with
\e
p_{\ACat^{\un,\iy}_R}\ci L=\iinc^\un_R\ci\cN(S_\sB)\ci S\ci\Pi_{\De^N}. 
\label{ai11eq59}
\e
Therefore we have a morphism
\e
(S\ci\pi_{\De^N},L):\De^N\t\De^1\longra(\iinc^\un_R\ci\sN(S_\sB))^*((\ACat^{\un,\iy}_R)_{\sA/}).
\label{ai11eq60}
\e

We will show that
\ea
(S\ci\pi_{\De^N},L)\vert_{\De^N\t\De^{\{0\}}}&=(S,\iinc^\un_R(T))=c(S,T),
\label{ai11eq61}\\
(S\ci\pi_{\De^N},L)\vert_{\De^N\t\De^{\{1\}}}&=(S,G)=b'(S,F)=b'\ci a'(S,T),
\label{ai11eq62}
\ea
where $F$ is as in \S\ref{ai1144} with $a'(S,T)=(S,F)$, and $G$ as in \S\ref{ai1146} with $b'(S,F)=(S,G)$. To prove \eq{ai11eq61}, consider the diagram
\e
\begin{gathered}
\xymatrix@C=80pt@R=15pt{
*+[r]{\De^N} \ar[r]^{\De^{\{1,2,\ldots,N+1\}}} \ar[d]^{\id_{\De^N}\t\De^{\{0\}}} & \De^{N+1} \ar[d]_{a_N} \ar@/^.5pc/[ddr]^{W\ci a_N=\ti T}
\\
*+[r]{\De^N\t\De^1} \ar@/_.5pc/[drr]_{\Xi_\sB\ci(S\t \De^{\{0\}}\t\id_{\De^1})}
 \ar[r]^(0.6){c_N} & C_N \ar@{..>}[dr]^{W'\ci f_N}
\\
&& \ACat^{\un,\iy}_R. }
\end{gathered}
\label{ai11eq63}
\e
This is essentially the restriction of \eq{ai11eq58} to $\De^N\t\De^{\{0\}}\t\De^1\subset \De^N\t\De^1\t\De^1$. Because a $k$-simplex $f:\De^k\ra C_N$ is flat in $C_N$ if and only if $f_N\ci f:\De^k\ra D_N$ is flat in $D_N$, we see that $W'\ci f_N$ takes flat $k$-simplices in $C_N$ to special $k$-simplices in $\ACat^{\un,\iy}_R$, and therefore $W'\ci f_N$ is the unique morphism given by Lemma \ref{ai11lem1}, comparing \eq{ai11eq46} and \eq{ai11eq63} (note that $W'$ came from Lemma~\ref{ai11lem2}).

As $\Xi_\sB\ci(S\t \De^{\{0\}}\t\id_{\De^1})=\iinc^\un_R\ci\cN(S_\sB)\ci\Pi_{\cN(\De^\rop)}$ by Definition \ref{ai11def1}(c), we see that $W'\ci f_N=\ti T\ci\pi_N$, since $\ti T\ci\pi_N$ also makes \eq{ai11eq63} commute and takes flat $k$-simplices to degenerate, hence special, $k$-simplices. But then
\begin{equation*}
W'\ci e_N\ci a_N=W'\ci f_N\ci b_N=\ti T\ci\pi_N\ci b_N=\ti T.
\end{equation*}
As $W'\ci e_N$ corresponds to $L:\De^N\t\De^1\ra (\ACat^{\un,\iy}_R)_{\sA/}$, $W'\ci e_N\ci a_N$ corresponds to $L\vert_{\De^N\t\De^{\{0\}}}$, and so $L\vert_{\De^N\t\De^{\{0\}}}=\iinc^\un_R(T)$, proving~\eq{ai11eq61}.
 
Similarly, to prove \eq{ai11eq62}, consider the diagram
\begin{equation*}
\xymatrix@C=80pt@R=15pt{
*+[r]{\De^N} \ar[r]^{\De^{\{1,2,\ldots,N+1\}}} \ar[d]^{\id_{\De^N}\t\De^{\{0\}}} & \De^{N+1} \ar[d]_{a_N} \ar@/^.5pc/[ddr]^{W\ci b_N=\ti F}
\\
*+[r]{\De^N\t\De^1} \ar@/_.5pc/[drr]_{\Xi_\sB\ci(S\t \De^{\{1\}}\t\id_{\De^1})}
 \ar[r]^(0.6){c_N} & C_N \ar@{..>}[dr]^{W'\ci g_N}
\\
&& \ACat^{\un,\iy}_R. }
\end{equation*}
Comparing this with \eq{ai11eq56} in the definition of $b'(S,F)$, we see that $W'\ci g_N=\hat W$. Hence
\begin{equation*}
W'\ci e_N\ci b_N=W'\ci g_N\ci b_N=\hat W\ci b_N=\ti G.
\end{equation*}
As $W'\ci e_N$ corresponds to $L:\De^N\t\De^1\ra (\ACat^{\un,\iy}_R)_{\sA/}$, $W'\ci e_N\ci b_N$ corresponds to $L\vert_{\De^N\t\De^{\{1\}}}$, and so $L\vert_{\De^N\t\De^{\{1\}}}=G$, proving~\eq{ai11eq62}.

Using \eq{ai11eq60}, define $H$ in \eq{ai11eq31} to satisfy
\e
H\ci((S,T)\t\id_{\De^1})=(S\ci\pi_{\De^N},L): \De^N\t\De^1\longra (\ACat^{\un,\iy}_R)_{\sA/}.
\label{ai11eq64}
\e
This defines $H$ on $(S,T)(\De^N)\t\De^1$ for every $N$-simplex $(S,T)$ in $\sN(S_\sB)^*\ab(\sN(\ACat^\un_R)_{\sA/})$. One can show that the definition of $H\ci((S,T)\t\id_{\De^1})$ commutes with the action of morphisms in $\De^\rop$ on $(S,T)$, so that $H$ is a morphism of simplicial sets. The first and second equations of \eq{ai11eq32} follow from \eq{ai11eq61}--\eq{ai11eq62}. The third equation of \eq{ai11eq32} follows from \eq{ai11eq59} and \eq{ai11eq64}. This proves Proposition~\ref{ai11prop1}(b).

\subsubsection{Proof of Proposition \ref{ai11prop1}(c)}
\label{ai1149}

Consider the restriction of the construction of $b'$ in \S\ref{ai1146} to $\Pi_{\sN(\De^\rop)}^{-1}([0])$. Then we must take the $N$-simplex $S$ in $\sN(\De^\rop)$ to be constant at $[0]$. So by Definition \ref{ai11def2} the morphism $\Xi_\sB\ci(S\t \De^{\{1\}}\t\id_{\De^1})$ in \eq{ai11eq56} is the composition $\De^N\t\De^1\,{\buildrel\Pi_{\De^1}\over\longra}\,\De^1\,{\buildrel i_\sB\over\longra}\,\ACat^{\un,\iy}_R$, regarding the morphism $i_\sB:\sB\ra S_\sB([0])$ in $\ACat^\un_R$ as a 1-simplex in~$\ACat^{\un,\iy}_R$.

Thus the construction of $(S,G)$ from $(S,F)$ in \S\ref{ai1146} is a kind of `parallel transport' of $F$ along the morphism $i_\sB$, using special simplices. Since $i_\sB$ is an isomorphism, this process is reversible: we can reconstruct $(S,F)$ from $(S,G)$ by `parallel transport' along $i^{-1}_\sB:S_\sB([0])\ra\sB$. As these two constructions are inverse, $b'\vert_{[0]}$ in \eq{ai11eq33} induces bijections on $N$-simplices for all $N\ge 0$, and so is a strict isomorphism of simplicial sets.

This completes the proofs of Proposition \ref{ai11prop1} and Theorem~\ref{ai9thm4}.

\subsection{Proof of Theorem \ref{ai9thm5}}
\label{ai115}

For (a), first we show that $p:\NCat^{\un,\iy}_R\ra \ACat^{\un,\iy}_R$ is an {\it inner fibration\/} in the sense of Definition \ref{ai2def12}. This means that as in \eq{ai2eq8}, for all $0<i<N$, a commutative square of maps of simplicial sets `$\ra$' as follows can be completed with a map `$\dashra$' in a commutative diagram
\begin{equation*}
\xymatrix@C=80pt@R=15pt{
*+[r]{\La^N_{i_{\phantom{(}}}} \ar[r] \ar@{^{(}->}[d] & *+[l]{\NCat^{\un,\iy}_R} \ar[d]_p \\
*+[r]{\De^N} \ar@{..>}[ur] \ar[r] & *+[l]{\ACat^{\un,\iy}_R.\!} }
\end{equation*}
Using Definition \ref{ai9def2} we translate this into a filling problem in $\hACat^{\un,\iy}_R$:
\e
\begin{gathered}
\xymatrix@C=70pt@R=11pt{
*+[r]{\De^0\,\,\,} \ar@{_{(}->}[r]_{\De^{\{0\}}} \ar@<1ex>@/^1pc/[rr]^(0.3){\sP^\un} & {C(\La^N_{i_{\phantom{(}}})} \ar[r] \ar@{^{(}->}[d] & *+[l]{\hACat^{\un,\iy}_R}  \\
& {\De^{N+1}} \ar@{..>}[ur] & \\
*+[r]{\La^N_i\,\,\,} \ar@{^{(}->}[uur] \ar@{^{(}->}[r] & {\De^{N^{\phantom{(}}}} \ar@{^{(}->}[u]_{\De^{\{1,2,\ldots,N+1\}}} 
 \ar[r] & *+[l]{\ACat^{\un,\iy^{\phantom{(}}}_R.\!} \ar@{^{(}->}[uu] }
\end{gathered}
\label{ai11eq65}
\e
Here $C(\La^N_i)$ is the `cone on $\La^N_i$', and we can think of it as the simplicial subcomplex of $\De^{N+1}$ which has $\La^N_i$ in the $N$-simplex face $\De^{\{1,\ldots,N+1\}}$ in $\De^{N+1}$, and all simplices linking $\La^N_i$ to the remaining vertex 0.

Now $\La^N_i$ is the union of the $(N-1)$-simplex faces $\De^{\{0,\ldots,N\}\sm\{j\}}$ in $\De^N$ for $0\le j\le N$ with $j\ne i$. Therefore $C(\La^N_i)$ is the union of the $N$-simplex faces $\De^{\{0,\ldots,N+1\}\sm\{j\}}$ in $\De^{N+1}$ for $1\le j\le N+1$ with $j\ne i+1$. Taking the union of this with the $N$-simplex face $\De^{\{1,\ldots,N+1\}}$, glued to $C(\La^N_i)$ along $\La^N_i$, gives the inner horn $\La^{N+1}_{i+1}$ in $\De^{N+1}$. Thus \eq{ai11eq65} is equivalent to the lifting problem
\e
\begin{gathered}
\xymatrix@C=55pt@R=9pt{
*+[r]{\De^0\,\,\,} \ar@{_{(}->}[r]_{\De^{\{0\}}} \ar@<1ex>@/^1pc/[rrr]^(0.3){\sP^\un} & \La^{N+1}_{i+1} \ar[rr] \ar@{_{(}->}[dr] && *+[l]{\hACat^{\un,\iy}_R}  \\
&& {\De^{N+1}} \ar@{..>}[ur] & \\
 & {\De^N} \ar@{^{(}->}[uu]_(0.3){\De^{\{1,2,\ldots,N+1\}}} 
 \ar[rr] && *+[l]{\ACat^{\un,\iy^{\phantom{(}}}_R.\!} \ar@{^{(}->}[uu] }
\end{gathered}
\label{ai11eq66}
\e
This is simply a filling condition for an inner horn $h:\La^{N+1}_{i+1}\ra\hACat^{\un,\iy}_R$, with the additional conditions on $h$ that $h(\De^{\{0\}})=\sP^\un$ and $h(\De^{\{1,\ldots,N+1\}})\ab\subset\ab\ACat^{\un,\iy}_R\subset\hACat^{\un,\iy}_R$.

Let us enlarge \eq{ai11eq66} to the commutative diagram:
\e
\begin{gathered}
\xymatrix@C=55pt@R=9pt{
*+[r]{\De^0\,\,\,} \ar@{_{(}->}[r]_{\De^{\{0\}}} \ar@<1ex>@/^1pc/[rrr]^(0.3){\sP^\un} & \La^{N+1}_{i+1} \ar[rr] \ar@{_{(}->}[dr] && *+[l]{\hACat^{\un,\iy}_R}  \\
&& {\De^{N+1}} \ar@{..>}[ur] & \\
*+[r]{\De^2} \ar[r]^(0.55){\De^{\{i-1,i,i+1\}}} \ar[uur]^(0.45){\De^{\{i,i+1,i+2\}}} & {\De^N} \ar@{^{(}->}[uu]_(0.3){\De^{\{1,2,\ldots,N+1\}}} 
 \ar[rr] && *+[l]{\ACat^{\un,\iy^{\phantom{(}}}_R.\!} \ar@{^{(}->}[uu] }
\end{gathered}
\label{ai11eq67}
\e
Consider the composition $t:\De^2\ra \hACat^{\un,\iy}_R$, as a 2-simplex in $\hACat^{\un,\iy}_R$. Since it factors via $\ACat^{\un,\iy}_R$, it is a 2-simplex in $\ACat^{\un,\iy}_R\!\subset\!\hACat^{\un,\iy}_R$. Thus Theorem \ref{ai7thm4} says it is a {\it thin\/} 2-simplex in $(\hACat^{\un,\iy}_R,\cT)$, and it satisfies Definition \ref{ai2def21}(a). So comparing equations \eq{ai2eq26} and \eq{ai11eq67}, we see that we can fill in the morphisms `$\dashra$' to make \eq{ai11eq66} and \eq{ai11eq67} commute. Therefore $p:\NCat^{\un,\iy}_R\ra \ACat^{\un,\iy}_R$ is an inner fibration.

A morphism $f:\De^1\ra \NCat^{\un,\iy}_R$ is $p$-{\it cocartesian\/} in the sense of Definition \ref{ai2def12} if as in \eq{ai2eq9}, we can fill the following commutative diagrams for all~$N\ge 2$
\begin{equation*}
\xymatrix@C=80pt@R=15pt{
*+[r]{\De^{\{0,1\}}} \ar@<1ex>@/^.8pc/[rr]^f  \ar@{^{(}->}[r] & \La^N_{0_{\phantom{(}}} \ar[r] \ar@{^{(}->}[d] & *+[l]{\NCat^{\un,\iy}_R} \ar[d]_p \\
& *+[r]{\De^N} \ar@{..>}[ur] \ar[r] & *+[l]{\ACat^{\un,\iy}_R.\!} }
\end{equation*}
As for \eq{ai11eq65}--\eq{ai11eq66}, this translates into the filling problem
\e
\begin{gathered}
\xymatrix@C=55pt@R=9pt{
*+[r]{\De^2\,\,\,} \ar@{_{(}->}[r]_{\De^{\{0,1,2\}}} \ar@<.5ex>@/^.8pc/[rrr]^(0.2){f'} & \La^{N+1}_1 \ar[rr] \ar@{_{(}->}[dr] && *+[l]{\hACat^{\un,\iy}_R}  \\
&& {\De^{N+1}} \ar@{..>}[ur] & \\
 & {\De^N} \ar@{^{(}->}[uu]^{\De^{\{1,2,\ldots,N+1\}}} 
 \ar[rr] && *+[l]{\ACat^{\un,\iy^{\phantom{(}}}_R,\!} \ar@{^{(}->}[uu] }
\end{gathered}
\label{ai11eq68}
\e
where $f':\De^2\ra\hACat^{\un,\iy}_R$ corresponds to $f:\De^1\ra \NCat^{\un,\iy}_R$. Theorem \ref{ai7thm4} says that if the 2-simplex $f'$ is filled by an $A_\iy$-natural isomorphism then it is thin in the $\iy$-bicategory $(\hACat^{\un,\iy}_R,\cT)$, so $f'$ satisfies Definition \ref{ai2def21}(a). Thus, comparing \eq{ai2eq26} and \eq{ai11eq68}, we see that we can fill in the morphism `$\dashra$' to make \eq{ai11eq68} commute. That is, if $f'$ is filled by an $A_\iy$-natural isomorphism, then $f$ is $p$-cocartesian.

To show $p$ is a cocartesian fibration, as in \eq{ai2eq10} we must show we can fill the following commutative diagram
\e
\begin{gathered}
\xymatrix@C=80pt@R=15pt{
*+[r]{\De^0_{\phantom{(}}} \ar[r]^g \ar@{^{(}->}[d]_{\De^{\{0\}}} & *+[l]{\NCat^{\un,\iy}_R} \ar[d]_p \\
*+[r]{\De^1} \ar@{..>}[ur]^(0.4){\ti f} \ar[r]^f & *+[l]{\ACat^{\un,\iy}_R,\!} }
\end{gathered}
\label{ai11eq69}
\e
such that $\ti f$ is $p$-cocartesian. As for \eq{ai11eq65} we translate this into a filling problem
\e
\begin{gathered}
\xymatrix@C=70pt@R=15pt{
*+[r]{\De^0\,\,\,} \ar@{_{(}->}[r]_{\De^{\{0\}}} \ar@<1ex>@/^1pc/[rr]^(0.3){\sP^\un} & {\De^1_{\phantom{(}}} \ar[r]^(0.3){g'} \ar@{^{(}->}[d]^{\De^{\{0,1\}}} & *+[l]{\hACat^{\un,\iy}_R}  \\
& {\De^2} \ar@{..>}[ur]_(0.45){\ti f'} & \\
*+[r]{\De^0\,\,\,} \ar@{^{(}->}[uur]^{\De^{\{1\}}} \ar@{^{(}->}[r]^{\De^{\{0\}}} & {\De^{1^{\phantom{(}}}} \ar@{^{(}->}[u]_{\De^{\{1,2\}}} \ar@/_1pc/[uur]_(0.6)f
 \ar[r] & *+[l]{\ACat^{\un,\iy^{\phantom{(}}}_R.\!} \ar@{^{(}->}[uu] }
\end{gathered}
\label{ai11eq70}
\e
Here $g'\!:\!\sP^\un\!\ra\!\sA$ and $f\!:\!\sA\!\ra\!\sB$ are morphisms (1-simplices) in $\ACat^{\un,\iy}_R$, as all 1-simplices in $\hACat^{\un,\iy}_R$ lie in $\ACat^{\un,\iy}_R$, so they are unital $A_\iy$-functors. The source $\sA$ of $f$ and target $\sA$ of $g'$ agree because of the commutative diamond starting at the bottom left $\De^0$ in~\eq{ai11eq70}.

Define $\ti f':\De^2\ra\hACat^{\un,\iy}_R$ to have unital $A_\iy$-categories $\sP^\un,\sA,\sB$ at vertices $0,1,2$, and to have unital $A_\iy$-functors $g',f,f\ci g'$ on the edges $\{0,1\},\{1,2\},\{0,2\}$. Define the data $\m_k^{\{0,1,2\}}$ in $\ti f'(\De^2)$ by $\m_0^{\{0,1,2\};*,*}:\la\mapsto \la\boo^\sB_{f\ci g'(*)}$ and $\m_k^{\{0,1,2\};*,\ldots,*}=0$ for $k\ge 1$. It is easy to show this defines a 2-simplex $\ti f'$ of $\hACat^{\un,\iy}_R$, which corresponds to an edge $\ti f:\De^1\ra \NCat^{\un,\iy}_R$ making \eq{ai11eq69} commute. As $[\boo^\sB_{f\ci g'(*)}]$ is an isomorphism in $H^0(\sB)$, $\ti f'$ is filled by an $A_\iy$-natural isomorphism, and so $\ti f$ is $p$-cocartesian from above. This proves $p$ is a cocartesian fibration. The rest of Theorem \ref{ai9thm5}(a) is immediate.

For (b), note from Definition \ref{ai9def2} that $p^{-1}(\sA)=\Hom_{\hACat^{\un,\iy}_R}^L(\sP^\un,\sA)$. Thus by \eq{ai9eq2} we have a strict isomorphism $p^{-1}(\sA)\cong \sN_{A_\iy}(\Fun_{A_\iy}^\un(\sP^\un,\sA))$. But Proposition \ref{ai4prop2} gives a strict isomorphism $\sA\cong\Fun^\un_{A_\iy}(\sP^\un,\sA)$. Hence $p^{-1}(\sA)\cong\sN_{A_\iy}(\sA)$, as we have to prove. The second part of (b) is immediate.

For (c), we first construct a strict isomorphism of simplicial sets $\Phi$ in the commutative diagram
\e
\begin{gathered}
\xymatrix@C=150pt@R=15pt{
*+[r]{\NCat^{\un,\iy}_R\vert_{\sN(\ACat^\un_R)}} \ar@{..>}[r]_(0.55)\Phi^(0.55)\cong \ar[d]^{p\vert_{\cdots}} & *+[l]{N_{\sN_{A_\iy}}(\ACat^\un_R)} \ar[d]_\pi \\
*+[r]{\sN(\ACat^\un_R)} \ar@{=}[r] &  *+[l]{\sN(\ACat^\un_R).\!} }
\end{gathered}
\label{ai11eq71}
\e
Here we use $\iinc^\un_R$ to identify $\sN(\ACat^\un_R)$ with its image as an $\iy$-subcategory of $\ACat^{\un,\iy}_R$. The left hand column of \eq{ai11eq71} is the restriction of $p:\NCat^{\un,\iy}_R\ra\ACat^{\un,\iy}_R$ in Definition \ref{ai9def2} to $\sN(\ACat^\un_R)$. The right hand column is the {\it relative simplicial nerve\/} from Example \ref{ai2ex6} of $\sN_{A_\iy}:\ACat^\un_R\ra\Cat^1_\iy$ from Definition~\ref{ai4def20}.

Fix an $N$-simplex $S$ in $\sN(\ACat^\un_R)\subset\ACat^{\un,\iy}_R$ for the next two examples. For convenience of notation, we number the vertices of $S$ by $1,2,\ldots,N+1$ rather than $0,1,\ldots,N$. Thus $S$ is equivalent to a choice of unital $A_\iy$-categories $\sA^i$ over $R$ for $i=1,\ldots,N+1$, and unital $A_\iy$-functors $\fF^{ij}:\sA^i\ra\sA^j$ for $1\le i\le j\le N+1$, where $\fF^{ii}=\id_{\sA^i}$ and $\fF^{ik}=\fF^{jk}\ci\fF^{ij}$ for $1\le i\le j\le k\le N+1$. The next two long examples compute the sets of $N$-simplices in $\NCat^{\un,\iy}_R\vert_{\sN(\ACat^\un_R)},N_{\sN_{A_\iy}}(\ACat^\un_R)$ mapped to $S$ by~$p,\pi$.

\begin{ex}
\label{ai11ex3}
An $N$-simplex $T$ in $\NCat^{\un,\iy}_R\vert_{\sN(\ACat^\un_R)}$ with $p(T)=S$ is, by Definition \ref{ai9def2}, an $(N+1)$-simplex $f:\De_{N+1}\ra\hACat^{\un,\iy}_R$ such that $f(\De^{\{0\}})=\sP^\un$ and $f(\De^{\{1,2,\ldots,N+1\}})=S$. By Definition \ref{ai7def3}, the data in $f$ which is not determined by $f(\De^{\{0\}})\!=\!\sP^\un$ and $f(\De^{\{1,2,\ldots,N+1\}})\!=\!S$ consists~of:
\begin{itemize}
\setlength{\itemsep}{0pt}
\setlength{\parsep}{0pt}
\item[(i)] Unital $A_\iy$-functors $\fG^{0i}:\sP^\un\ra\sA^i$ for $i=1,\ldots,N+1$, which live on the edge $[0,i]$ of $\De^{N+1}$. By Definition \ref{ai4def11}, $\fG^{0i}$ is equivalent to the choice of an object $X^i$ in $\sA^i$.
\item[(ii)] Let $\bar I=\{0,i_0,\ldots,i_M\}$ for $0\le M\le N$ and $1\le i_0<i_1<\cdots<i_M\le N+1$, and write $I=\{i_0,\ldots,i_M\}$, so that $\es\ne I\subseteq\{1,\ldots,N+1\}$ and $\bar I=\{0\}\amalg I$. Then $f$ contains a map $\m_k^{\bar I;*,\ldots,*}$ as in \eq{ai7eq14} for all such $\bar I$ and $k\ge 0$, where $x_0,\ldots,x_k$ must be the unique object $*$ in $\sP^\un$. As $\sH^0(*,*)=R\cdot\boo_*$, equations \eq{ai7eq18}--\eq{ai7eq19} imply that if $M=1$ and $k=1$ then $\m_k^{\{0,i_0\};*,*}:R\cdot\boo_*\ra\Hom_{\sA^{i_0}}(X^{i_0},X^{i_0})$ maps $\boo_*\mapsto\boo_{X^{i_0}}$, and if either $k>1$, or $M\ne 1$ and $k=1$, then $\m_k^{\bar I;*,\ldots,*}=0$.

Thus the only free data in the maps $\m_k^{\bar I;*,\ldots,*}$ is when $k=0$. As $f^{\bar I}(*)=\fF^{i_{M-1}i_M}\ci\cdots\ci\fF^{i_0i_1}(X^{i_0})=\fF^{i_0i_M}(X^{i_0})$, and $f^{0,i_M}(*)=X^{i_M}$, this is a degree $1-M$ map
\begin{equation*}
\m_0^{\bar I;*}:R\cdot\boo_*\longra \Hom_{\sA^{i_M}}(\fF^{i_0i_M}(X^{i_0}),X^{i_M}),
\end{equation*}
which is determined by $\la^I=\m_0^{\bar I;*}(\boo_*)\in\Hom_{\sA^{i_M}}^{1-M}(\fF^{i_0i_M}(X^{i_0}),X^{i_M})$. Since $\m^{\bar I;*}_0=0$ when $\md{\bar I}=1,2$, we have $\la^I=0$ when $M=0$.
\end{itemize}

Then $S$ and $(X^i)_{i=1}^{N+1},(\la^I)_{\es\ne I\subseteq\{1,\ldots,N+1\}}$ determine the data $\bigl((\sO^i),\ab(f^{ij}),\ab(\sH^i),\ab(\m^I)\bigr)$ in the $(N+1)$-simplex $f:\De_{N+1}\ra\hACat^{\un,\iy}_R$ uniquely. For a given choice of $(X^i),(\la^I)$ to come from an $(N+1)$-simplex $f$, we need the corresponding $(\m^I)$ to satisfy \eq{ai7eq15}--\eq{ai7eq19}. This holds by construction on $\De^{\{0\}}$ and $\De^{\{1,\ldots,N+1\}}$, and the extra $\m_k^{\bar I;*,\ldots,*}$ satisfy \eq{ai7eq18}--\eq{ai7eq19} by~(ii). 

It remains to ensure \eq{ai7eq16} holds with $\bar I$ in place of $I$, for $\bar I=\{0,i_0,\ldots,i_M\}$ as in (ii). This involves a choice of $n\ge 0$ and $x_0,\ldots,x_n$ which must be the unique object $*$ in $\sP^\un$, and then \eq{ai7eq16} is an equation in maps
\begin{equation*}
(R\cdot\boo_*)^{\ot^n}\longra \Hom_{\sA^{i_M}}^{2-M}(\fF^{i_0i_M}(X^{i_0}),X^{i_M}).
\end{equation*}
If $n>0$ this holds automatically, as substituting $\boo_*$ into any $\m_l^{\bar J;*,\ldots,*}$ is prescribed by \eq{ai7eq17}--\eq{ai7eq19}. So the only nontrivial equation is \eq{ai7eq16} with $\bar I$ in place of $I$ and $n=0$. The first term on the right hand side of \eq{ai7eq16} vanishes as $\d_{0,k;y_0,\ldots,y_k}^{\{i_0\};*}=0$. So omitting the $(-1)^M$, equation \eq{ai7eq16} becomes
\e
\sum_{\begin{subarray}{l} k=1,\ldots,M,\\
y_0,\ldots,y_k\in\Obj\sA^{i_M} \end{subarray}}\m_k^{\sA^{i_M};y_0,\ldots,y_k}\ci\d_{0,k;y_0,\ldots,y_k}^{\bar I;*}=0.
\label{ai11eq72}
\e

Now $\d_{0,k;y_0,\ldots,y_k}^{\bar I;*}$ is defined in \eq{ai7eq9} as a sum over labelled trees $\bs\Ga_1,\ldots,\bs\Ga_k$. Thus we may rewrite \eq{ai11eq72} as
\e
0=\sum_{\begin{subarray}{l} 
k=1,\ldots,M, \\
\text{isomorphism classes of $\bs\Ga_1,\ldots,\bs\Ga_k$:}\\ 
\text{$\bs\Ga_i$ is a balanced $\bar I$-labelled tree} \\ 
\text{with no leaves, $1\le i\le k$,} \\
I=I_0\supseteq I_1 \supseteq \cdots\supseteq I_k=\{i_M\}, \\
L(\bs\Ga_i)=\{0\}\amalg I_{i-1}, \;
R(\bs\Ga_i)=\{0\}\amalg I_i, \; 1\le i\le k, \\
y_i=\fF^{(\min I_i)i_M}(X^{\min I_i}), \; 1\le i\le k \end{subarray}
\!\!\!\!\!\!\!\!\!\!\!\!\!\!}
\begin{aligned}[t]
&(-1)^\dag\m_k^{\sA^{i_M};y_0,\ldots,y_k}\ci{}\\
& \bigl(C(\bs\Ga_1)\!\ot\! C(\bs\Ga_2)\!\ot\!\cdots\!\ot\! C(\bs\Ga_k)\bigr),
\end{aligned}
\label{ai11eq73}
\e
where $(-1)^\dag$ is the sign in \eq{ai7eq9}. Here the trees $\bs\Ga_i$ in \eq{ai11eq73} may have internal vertices $v$ of two kinds:
\begin{itemize}
\setlength{\itemsep}{0pt}
\setlength{\parsep}{0pt}
\item[(a)] If $0\in I_v$ then $I_v=\{0\}\amalg J$ for $\es\ne J\subseteq I$, and $v$ must be a bud, and $\m^{I_v;*}_0$ used in Definition \ref{ai7def1} to define $C(\bs\Ga_i)$ maps $r\mapsto r\la^J$ for $r\in R$.
\item[(b)] If $0\notin I_v$ then $\m_l^{I_v;w_0,\ldots,w_l}$ in $C(\bs\Ga_i)$ comes from the $N$-simplex $S$. Thus
\e
\m_l^{I_v;w_0,\ldots,w_l}=\begin{cases}
(\fF^{ij})_l^{w_0,\ldots,w_l}, & I_v=\{i,j\},\; 1\le i<j\le N+1, \\[4pt]
\boo_{\fF^{hj}(X^h)}, & {\begin{subarray}{l} \ts l=0,\; I_v=\{h,i,j\},\\ \ts 1\le h<i<j\le N+1,\end{subarray}} \\[4pt]
0, & \text{otherwise.}
\end{cases}
\label{ai11eq74}
\e
\end{itemize}

Now by the argument used in the proof of Theorem \ref{ai7thm6} in \S\ref{ai108}, terms involving $\boo_{\fF^{hj}(X^h)}$ in \eq{ai11eq74} can only make a nonzero contribution to \eq{ai11eq73} if $k=2$ and $\bs\Ga_1,\bs\Ga_2$ are as in \eq{ai10eq47} with $n=0$. Taking these terms out of the sum, we may rewrite \eq{ai11eq73} as
\ea
&0=\sum_{0<j<M}(-1)^j\la^{I\sm\{i_j\}}
\label{ai11eq75}\\
&+\sum_{\begin{subarray}{l} 
k=1,\ldots,M, \\
\text{isomorphism classes of $\bs\Ga_1,\ldots,\bs\Ga_k$:}\\ 
\text{$\bs\Ga_i$ is a balanced $\bar I$-labelled tree} \\ 
\text{with no leaves, $1\le i\le k$,} \\
I=I_0\supseteq I_1 \supseteq \cdots\supseteq I_k=\{i_M\}, \\
L(\bs\Ga_i)=\{0\}\amalg I_{i-1}, \;
R(\bs\Ga_i)=\{0\}\amalg I_i, \; 1\le i\le k, \\
y_i=\fF^{(\min I_i)i_M}(X^{\min I_i}), \; 1\le i\le k, \\
\text{if $v$ is an internal vertex in $\bs\Ga_i$ then} \\
\text{either $0\in I_v$ or $\md{I_v}=2$} 
\end{subarray}
\!\!\!\!\!\!\!\!\!\!\!\!\!\!}
\begin{aligned}[t]
&(-1)^\dag\m_k^{\sA^{i_M};y_0,\ldots,y_k}\ci{}\\
& \bigl(C(\bs\Ga_1)\!\ot\! C(\bs\Ga_2)\!\ot\!\cdots\!\ot\! C(\bs\Ga_k)\bigr),
\end{aligned}
\nonumber
\ea
where the first line is the analogue of the second line of~\eq{ai7eq25}.

To summarize: we have shown $N$-simplices $T$ in $\NCat^{\un,\iy}_R\vert_{\sN(\ACat^\un_R)}$ with $p(T)=S$ are equivalent to data $(X^i)_{i=1}^{N+1},(\la^I)_{\es\ne I\subseteq\{1,\ldots,N+1\}}$ as above, satisfying equation \eq{ai11eq75} for all~$\es\ne I\subseteq\{1,\ldots,N+1\}$.
\end{ex}

\begin{ex}
\label{ai11ex4}
By Example \ref{ai2ex6}, an $N$-simplex $U$ in $N_{\sN_{A_\iy}}(\ACat^\un_R)$ with $\pi(U)=S$ is $(S,(\tau(I))_{\es\ne I\subseteq\{1,\ldots,N+1\}})$ where $\tau(I):\De^I\ra\sN_{A_\iy}(\sA^{\max I})$ is a map of simplicial sets for all $\es\ne I\subseteq\{1,\ldots,N+1\}$, such that for all $\es\ne J\subseteq I\subseteq\{1,\ldots,N+1\}$, as in \eq{ai2eq12} the following diagram commutes:
\e
\begin{gathered}
\xymatrix@C=150pt@R=15pt{
*+[r]{\De^J} \ar[r]_(0.45){\tau(J)} \ar[d]^{\De^{\inc(J\hookra I)}} & *+[l]{\sN_{A_\iy}(\sA^{\max J})} \ar[d]_{\sN_{A_\iy}(\fF^{(\max J)(\max I)})} \\
*+[r]{\De^I} \ar[r]^(0.45){\tau(I)} & *+[l]{\sN_{A_\iy}(\sA^{\max I}).\!}
}
\end{gathered}
\label{ai11eq76}
\e

Generalizing $\sC^\De_{[N]}$ in Definition \ref{ai4def20}, for each $\es\ne I\subseteq\{1,\ldots,N+1\}$ define a unital $A_\iy$-category $\sC^\De_I$ over $R$ to have objects $\Obj\sC^\De_I=I$ and Hom spaces
\begin{equation*}
\Hom(i,j)=\begin{cases} R\cdot\om_{ij}, & i\le j, \\
0, & i>j,
\end{cases}
\end{equation*}
with $\deg\om_{ij}=0$, and units $\boo_i=\om_{ii}$, and maps $\m_k^{i_0,\ldots,i_k}$ as in \eq{ai4eq37}, so $\sC^\De_I\cong\sC^\De_{[\md{I}-1]}$. Then $\tau(I)$ above is equivalent to a unital $A_\iy$-functor $\up(I):\sC^\De_I\ra\sA^{\max I}$, and \eq{ai11eq76} is equivalent to the following diagram commuting
\e
\begin{gathered}
\xymatrix@C=150pt@R=15pt{
*+[r]{\sC^\De_J} \ar[r]_(0.45){\up(J)} \ar[d]^{\cC^\De_{\inc(J\hookra I)}} & *+[l]{\sA^{\max J}} \ar[d]_{\fF^{(\max J)(\max I)}} \\
*+[r]{\sC^\De_I} \ar[r]^(0.45){\up(I)} & *+[l]{\sA^{\max I},\!}
}
\end{gathered}
\label{ai11eq77}
\e
where $\cC^\De_{\inc(J\hookra I)}:\sC^\De_J\hookra\sC^\De_I$ is the obvious inclusion of $A_\iy$-categories.

We will rewrite $(\up(I))_{\es\ne I\subseteq[N]}$ in terms of simpler data. First consider how the $A_\iy$-functors $\up(I)$ act on objects. For $i=1,\ldots,N+1$, $\up(\{i\}):\De^{\{i\}}\cong\De^0\ra\sA^i$ is equivalent to the choice of an object $X^i$ in $\sA^i$. These will correspond to $X^1,\ldots,X^{N+1}$ in Example \ref{ai11ex3}.
For all $\es\ne I\subseteq\{0,\ldots,N\}$, if $i\in I$, taking $J=\{i\}$ in \eq{ai11eq77} implies that $\tau(I)$ acts on objects by ~$\tau(I):i\mapsto \fF^{i\max I}(X^i)$. 

Let $\es\ne J\subseteq I\subseteq\{1,\ldots,N+1\}$, and write $J=\{j_0<j_1<\cdots<j_L\}$ and $I=\{i_0<i_1<\cdots<i_M\}$. Then $\up(I)$ includes a degree $1-L$ map
\begin{equation*}
\up(I)_L^{j_0,j_1,\ldots,j_L}:(R\cdot\om_{j_0j_1})\!\ot\!\cdots\!\ot\!(R\cdot\om_{j_{L-1}j_L})\!\longra\!\Hom_{\sA^{j_L}}(\fF^{j_0j_L}(X^{j_0}),X^{j_L}).
\end{equation*}
Define $\up^{J,I}\in \Hom_{\sA^{j_L}}^{1-L}(\fF^{j_0j_L}(X^{j_0}),X^{j_L})$ by
\begin{equation*}
\up^{J,I}=\up(I)_L^{j_0,j_1,\ldots,j_L}(\om_{j_0j_1}\ot\cdots\ot\om_{j_{L-1}j_L}).
\end{equation*}
When $J=I$, set $\la^I=\up^{I,I}$. These will correspond to $\la^I$ in Example \ref{ai11ex3}. When $M=0$, as $\up(\{i_0\})_0^{i_0}=0$ we have~$\up^{\{i_0\},\{i_0\}}=\la^{\{i_0\}}=0$.

We claim that $X^1,\ldots,X^{N+1}$ and $\up^{J,I}$ for $\es\ne J\subseteq I$ determine the unital $A_\iy$-functor $\up(I):\sC^\De_I\ra\sA^{\max I}$. This holds as in the maps $\up(I)_k^{x_0,\ldots,x_k}$ for $x_0,\ldots,x_k$ objects in $\sC^\De_I$, so that $x_0,\ldots,x_k\in I$, if we do not have $x_0\ab\le\ab x_1\ab\le\ab\cdots\ab\le x_k$ then $\up(I)_k^{x_0,\ldots,x_k}=0$ as its domain is zero, and if $x_{i-1}=x_i$ then the domain of $\up(I)_k^{x_0,\ldots,x_k}$ includes a factor $R\cdot\om_{x_ix_i}=R\cdot\boo_{x_i}$, and so $\up(I)_k^{x_0,\ldots,x_k}$ is determined by \eq{ai4eq11}--\eq{ai4eq12}. Hence the only free data is when $x_0<\cdots<x_k$, and then $\up(I)_k^{x_0,\ldots,x_k}$ maps $r(\om_{x_0x_1}\ot\cdots\ot \om_{x_{k-1}x_k})\mapsto r\up^{J,I}$ for~$J=\{x_0,\ldots,x_k\}$.

Let $\es\ne A\subseteq J\subseteq I\subseteq \{1,\ldots,N+1\}$, and write $A=\{a_0<a_1<\cdots<a_n\}$ and $\max J=j_L$, $\max I=i_M$. Then evaluating the condition for \eq{ai11eq77} to commute on objects $a_0,\ldots,a_n$ in $\sC^\De_J$ and on $\om_{a_0a_1}\ot\cdots\ot\om_{a_{n-1}a_n}$ in its domain, and using \eq{ai4eq10} to expand the composition $\fF^{j_Li_M}\ci\up(J)$, yields the identity
\ea
&\up^{A,I}=
\label{ai11eq78}\\
&\sum_{\begin{subarray}{l} 0=l_0<l_1<\cdots<l_k=n
\end{subarray}}
\begin{aligned}[t]& (-1)^{\sum_{i=1}^k(l_i-l_{i-1}+1)(k-i)}\cdot{}\\
&(\fF^{j_Li_M})_k^{\fF^{a_{l_0}j_L}(X^{a_{l_0}}),\fF^{a_{l_1}j_L}(X^{a_{l_1}}),\ldots,\fF^{a_{l_k}j_L}(X^{a_{l_k}})}\ci{} \\ 
&\bigl(\up^{\{a_{l_0},\ldots,a_{l_1}\},J}\ot\up^{\{a_{l_1},\ldots,a_{l_2}\},J}\ot\cdots\ot\up^{\{a_{l_{k-1}},\ldots,a_{l_k}\},J}).
\end{aligned}
\nonumber
\ea
If $j_L=i_M$ then $\fF^{j_Li_M}=\id_{\sA^{i_M}}$, and $(\fF^{j_Li_M})_k^{\ldots}=0$ unless $k=1$, so that 
\eq{ai11eq78} reduces to
\e
\up^{A,I}=\up^{A,J}\qquad\text{if $\max J=\max I$.}
\label{ai11eq79}
\e

Note that if $\md{A}<\md{I}$ then taking $J=A$, equation \eq{ai11eq78} determines $\up^{A,I}$ in terms of $\up^{B,J}$ for $\md{J}<\md{I}$. So by induction on $\md{I}$, using \eq{ai11eq78} we can iteratively determine $\up^{A,I}$ for all $\es\ne A\subseteq I$ in terms of $\la^J=\up^{J,J}$ for all $\es\ne J\subseteq I$. In fact we can give an explicit formula. We claim that
\e
\up^{J,I}=\sum_{\begin{subarray}{l} 
\text{isomorphism classes of $\bs\Ga$: $\bs\Ga$ is a}\\ 
\text{balanced $\bar I$-labelled tree with no leaves, } \\
L(\bs\Ga)=\{0\}\amalg J\cup\{i_M\}, \; R(\bs\Ga)=\{0,j_L,i_M\} \\
\text{if $v$ is an internal vertex in $\bs\Ga$ then} \\
\text{either $0\in I_v$ or $\md{I_v}=2$} 
\end{subarray}
\!\!\!\!\!\!\!\!\!\!\!\!\!\!}
(-1)^\dag C(\bs\Ga).
\label{ai11eq80}
\e

Here $(-1)^\dag$ is the sign in \eq{ai7eq9} for $n=0$, $k=1$, and $\bs\Ga_1=\bs\Ga$. To define $C(\bs\Ga)$ the vertices $v$ of $\bs\Ga$ are labelled with maps $\m_l^{I_v;w_0,\ldots,w_l}$ as in \S\ref{ai71} according to the following rules, which should be compared to Example \ref{ai11ex3}(a),(b):
\begin{itemize}
\setlength{\itemsep}{0pt}
\setlength{\parsep}{0pt}
\item[(a)] If $0\in I_v$ then $I_v=\{0\}\amalg J$ for $\es\ne J\subseteq I$, and $v$ must be a bud, and $\m^{I_v;*}_0$ used in Definition \ref{ai7def1} to define $C(\bs\Ga_i)$ maps $r\mapsto r\la^J$ for $r\in R$.
\item[(b)] If $0\notin I_v$ then $\md{I_v}=2$ by \eq{ai11eq80}, so $I_v=\{i,j\}$ for $1\le i<j\le N+1$, and~$\m_l^{I_v;w_0,\ldots,w_l}=(\fF^{ij})_l^{w_0,\ldots,w_l}$.
\end{itemize}

Before explaining the proof of \eq{ai11eq80}, observe that if $j_L<b_1<\cdots<b_d<i_M$ for $d\ge 0$ then \eq{ai11eq80} is equivalent to 
\e
\up^{J,I}=\sum_{\begin{subarray}{l} 
\text{isomorphism classes of $\bs\Ga$: $\bs\Ga$ is a}\\ 
\text{balanced $\bar I$-labelled tree with no leaves, } \\
L(\bs\Ga)=\{0\}\amalg J\cup\{b_1,\ldots,b_d,i_M\}, \; R(\bs\Ga)=\{0,j_L,b_1,\ldots,b_d,i_M\} \\
\text{if $v$ is an internal vertex in $\bs\Ga$ then} \\
\text{either $0\in I_v$ or $\md{I_v}=2$} 
\end{subarray}
\!\!\!\!\!\!\!\!\!\!\!\!\!\!}
(-1)^\dag C(\bs\Ga).
\label{ai11eq81}
\e
By including $b_1,\ldots,b_d$ in $L(\bs\Ga),R(\bs\Ga)$, vertices $v$ with $I_v=\{j_L,i_M\}$ in $\bs\Ga$ in \eq{ai11eq80} get replaced by subgraphs $\bs\Ga'$ of graphs $\bs\Ga$ in \eq{ai11eq81} with vertices $v'$ with $I_{v'}=\{j_L,b_1\},\{b_1,b_2\},\ldots,\{b_{d-1},b_d\},\{b_d,i_M\}$. To get from \eq{ai11eq81} to \eq{ai11eq80}, we collapse each such subgraph $\bs\Ga'$ to the single vertex $v$ with $I_v=\{j_L,i_M\}$. We can then show that the right hand sides of \eq{ai11eq80} and \eq{ai11eq81} are equal using~$\fF^{j_Li_M}=\fF^{b_di_M}\ci\fF^{b_{d-1}b_d}\ci\cdots\ci\fF^{j_Lb_1}$.

We prove \eq{ai11eq80} by induction on $i_M=\max I=1,2,\ldots,N+1$. The initial step $i_M=1$ is trivial. To prove the inductive step, note that if $j_L=i_M$ then there is only one graph $\bs\Ga$ in the sum in \eq{ai11eq80}, which has one vertex $v$ with $I_v=\{0\}\amalg J$ plus the root vertex, and \eq{ai11eq80} reduces to $\up^{J,I}=\la^J$, which is true as $\up^{J,I}=\up^{J,J}=\la^J$ by~\eq{ai11eq79}.

If $j_L<i_M$ then every tree $\bs\Ga$ in the sum in \eq{ai11eq80} must have a unique vertex $v$ connected to the root vertex with $I_v=\{j_L,j_M\}$, which has $k\ge 1$ incoming edges. We `saw' $\bs\Ga$ at level $j_L$, as in Definition \ref{ai6def7}. This divides $\bs\Ga$ into a tree $\bs\Ga_0$ consisting of $v$ with $k$ leaves, corresponding to $(\fF^{j_Li_M})_k^{\ldots}$ in \eq{ai11eq78}, and further trees $\bs\Ga_1,\ldots,\bs\Ga_k$ substituted into the leaves of $\bs\Ga_0$. These $\bs\Ga_1,\ldots,\bs\Ga_k$ can be interpreted as terms in the sums in \eq{ai11eq81} (which holds by induction as $j_L<i_M$) for $\up^{\{a_{l_0},\ldots,a_{l_1}\},J},\ab\up^{\{a_{l_1},\ldots,a_{l_2}\},J},\ab\ldots,\ab\up^{\{a_{k-1},\ldots,a_k\},J}$ in \eq{ai11eq78}, for some added elements $b_1,\ldots,b_d$ in each case. The difference in the sign $(-1)^\dag$ in \eq{ai11eq80} for $\bs\Ga$ and the product of $(-1)^\dag$ for $\bs\Ga_1,\ldots,\bs\Ga_k$ is the sign in \eq{ai11eq78}. Hence \eq{ai11eq78} implies the inductive step.

We proved \eq{ai11eq78}--\eq{ai11eq81} using the condition that \eq{ai11eq77} commutes for $\es\ne J\subseteq I$. However, we have not yet discussed the conditions on the $\up^{J,I}$ for $\up(I):\sC^\De_I\ra\sA^{\max I}$ to be a unital $A_\iy$-functor. Let $\es\ne J\subseteq I\subseteq\{1,\ldots,N+1\}$, and write $J=\{j_0<j_1<\cdots<j_L\}$ and $I=\{i_0<i_1<\cdots<i_M\}$. Consider equation \eq{ai4eq8} for $\up(I)$ with $n=L$ and $(x_0,\ldots,x_n)=(j_0,\ldots,j_L)$, evaluated on $\om_{j_0j_1}\ot\cdots\ot\om_{j_{L-1}j_L}$. This becomes
\ea
&\sum_{0<a<L}(-1)^{a-1}\up^{J\sm\{j_a\},I}
\label{ai11eq82}\\
&=\sum_{\begin{subarray}{l} 0=l_0<l_1<\cdots<l_k=L
\end{subarray}}
\begin{aligned}[t]& (-1)^{\sum_{i=1}^k(l_i-l_{i-1}+1)(k-i)}\cdot{}\\
&\m_k^{\sA^{i_M};\fF^{j_{l_0}i_M}(X^{j_{l_0}}),\fF^{j_{l_1}i_M}(X^{j_{l_1}}),\ldots,\fF^{j_{l_k}i_M}(X^{j_{l_k}})}\ci{} \\ 
&\bigl(\up^{\{j_{l_0},\ldots,j_{l_1}\},I}\ot\up^{\{j_{l_1},\ldots,j_{l_2}\},I}\ot\cdots\ot\up^{\{j_{l_{k-1}},\ldots,j_{l_k}\},I}).
\end{aligned}
\nonumber
\ea
Here the term $(-1)^{a-1}\up^{J\sm\{j_a\},I}$ on the left hand side of \eq{ai11eq82} comes from the case $i=a-1$, $j=2$ and $k=L-a-1$ on the left hand side of \eq{ai4eq8}, and uses $\m_2^{\sC^\De_I;j_{a-1},j_a,j_{a+1}}(\om_{j_{a-1}j_a}\ot\om_{j_aj_{a+1}})=\om_{j_{a-1}j_{a+1}}$. Terms on the left hand side of \eq{ai4eq8} with $j\ne 2$ are zero, as $\m_j^{\sC^\De_I;\ldots}=0$ for~$j\ne 2$.

Equation \eq{ai11eq82} for all $\es\ne J\subseteq I$ is a sufficient condition for the data $\up^{J,I}$ to come from a unital $A_\iy$-functor $\up(I):\sC^\De_I\ra\sA^{\max I}$, since for all other choices of objects $(x_0,\ldots,x_n)$ in \eq{ai4eq8}, either we do not have $x_0\ab\le\ab x_1\ab\le\ab\cdots\ab\le x_n$, and then the domain of \eq{ai4eq8} is zero, or $x_{i-1}=x_i$ for some $1\le i\le n$, and then the domain of \eq{ai4eq8} includes a factor $R\cdot\om_{x_ix_i}=R\cdot\boo_{x_i}$, and \eq{ai4eq8} follows by unitality.

We can now substitute \eq{ai11eq80}, or \eq{ai11eq81}, into \eq{ai11eq82} to write \eq{ai11eq82} as a condition on the data $\la^A$ rather $\up^{A,I}$. First consider the case in $J=I$, so that $L=M$ and $j_a=i_a$. This yields
\ea
&\sum_{0<a<M}(-1)^{a-1}\la^{I\sm\{i_a\}}
\label{ai11eq83}\\
&=\sum_{\begin{subarray}{l} 
k=1,\ldots,M, \\
\text{isomorphism classes of $\bs\Ga_1,\ldots,\bs\Ga_k$:}\\ 
\text{$\bs\Ga_i$ is a balanced $\bar I$-labelled tree} \\ 
\text{with no leaves, $1\le i\le k$,} \\
0=l_0<l_1<\cdots<l_k=M, \\
L(\bs\Ga_a)=\{0,i_{l_{a-1}},i_{l_{a-1}}+1,\ldots,i_M\},\\
R(\bs\Ga_a)=\{0,i_{l_a},i_{l_a}+1,\ldots,i_M\},\; 1\le a\le k,\\
y_a=\fF^{i_{l_a}i_M}(X^{i_{l_a}}), \; 0\le a\le k, \\
\text{if $v$ is an internal vertex in $\bs\Ga_i$ then} \\
\text{either $0\in I_v$ or $\md{I_v}=2$} 
\end{subarray}}
\begin{aligned}[t]
&(-1)^\dag\m_k^{\sA^{i_M};y_0,\ldots,y_k}\ci{}\\
& \bigl(C(\bs\Ga_1)\!\ot\! C(\bs\Ga_2)\!\ot\!\cdots\!\ot\! C(\bs\Ga_k)\bigr).
\end{aligned}
\nonumber
\ea
The top line is correct since $\up^{I\sm\{i_a\},I}=\up^{I\sm\{i_a\},I\sm\{i_a\}}=\la^{I\sm\{i_a\}}$ by \eq{ai11eq79}. Here for the term $\up^{\{i_{l_{a-1}},\ldots,i_{l_a}\},I}$ in \eq{ai11eq82}, we have substituted \eq{ai11eq81} with $(b_1,\ldots,b_d)=(i_{l_a+1},i_{l_a+2},\ldots,i_{M-1})$. This explains the conditions $L(\bs\Ga_a)=\{0,i_{l_{a-1}},i_{l_{a-1}}+1,\ldots,i_M\}$, $R(\bs\Ga_a)=\{0,i_{l_a},i_{l_a}+1,\ldots,i_M\}$.

Now \eq{ai11eq83} is equivalent to \eq{ai11eq75}, where the two correspond by taking $I_a$ in \eq{ai11eq75} to be $\{i_{l_a},i_{l_a}+1,\ldots,i_M\}$ in \eq{ai11eq83}. This works as the condition in \eq{ai11eq75} that either $0\in I_v$ or $\md{I_v}=2$ excludes $\bs\Ga_1,\ldots,\bs\Ga_k$ which do not have $I_0,\ldots,I_k$ of this form.

We also claim that \eq{ai11eq82} for $J=I$, taken over all $\es\ne I\subseteq\{1,\ldots,N+1\}$, together with \eq{ai11eq80}--\eq{ai11eq81}, imply \eq{ai11eq82} for all $\es\ne J\subseteq I\subseteq\{1,\ldots,N+1\}$. If $\max J=\max I$ this is immediate by \eq{ai11eq79}. If $\max J<\max I$, we can show that \eq{ai11eq82} for $J\subseteq I$ can be rewritten as a sum of compositions of $\fF^{(\max J)(\max I)}$ with \eq{ai11eq82} with $\{j_{l_{a-1}},\ldots,j_{l_a}\},J$ in place of $J,I$, in a similar way to \eq{ai11eq78}. Then we can prove by induction on $\max I$ that \eq{ai11eq82} for $J=I$, all $I$, implies \eq{ai11eq82} for all~$\es\ne J\subseteq I$.

To summarize: we have shown $N$-simplices $U$ in $N_{\sN_{A_\iy}}(\ACat^\un_R)$ with $\pi(U)=S$ are equivalent to data $(X^i)_{i=1}^{N+1},(\la^I)_{\es\ne I\subseteq\{1,\ldots,N+1\}}$ as above, satisfying \eq{ai11eq75} for all $\es\ne I\subseteq\{1,\ldots,N+1\}$. This concludes Example~\ref{ai11ex4}.
\end{ex}

Examples \ref{ai11ex3}--\ref{ai11ex4} show the $N$-simplices in $\NCat^{\un,\iy}_R\vert_{\sN(\ACat^\un_R)}$ and $N_{\sN_{A_\iy}}(\ACat^\un_R)$ mapped to $S$ by $p,\pi$ are in natural bijection, since both can be identified with data $(X^i)_{i=1}^{N+1},(\la^I)_{\es\ne I\subseteq\{1,\ldots,N+1\}}$ satisfying \eq{ai11eq75}. Taken over all $N$-simplices $S$ in $\sN(\ACat^\un_R)$, this defines a bijection
\begin{equation*}
\Phi([N]):\NCat^{\un,\iy}_R\vert_{\sN(\ACat^\un_R)}([N])\,{\buildrel\cong\over\longra}\, N_{\sN_{A_\iy}}(\ACat^\un_R)([N]),
\end{equation*}
such that \eq{ai11eq71} commutes at level $[N]$. It is not difficult to show that these maps $\Phi([N])$ for all $N\ge 0$ are compatible with the actions of morphisms $\al:[N]\ra[\bar N]$ in $\De^\rop$, and so define a strict isomorphism of simplicial sets $\Phi$ in the commuting diagram \eq{ai11eq71}. This proves the final part of Theorem~\ref{ai9thm5}(c).

By Example \ref{ai2ex6}, we have shown that $(\iinc^\un_R)^*(p)$ is isomorphic to a model for the unstraightening of $I^{\sCat_\iy}\ci\sN(\sN_{A_\iy})$ in Theorem \ref{ai2thm2}(c). Thus by Theorem \ref{ai2thm2}, $I^{\sCat_\iy}\ci\sN(\sN_{A_\iy})$ is a model for the straightening of $(\iinc^\un_R)^*(p)$. But as $\sN^\iy_{A_\iy}$ is the straightening of $p$, $\sN^\iy_{A_\iy}\ci\iinc^\un_R$ is the straightening of $(\iinc^\un_R)^*(p)$. Therefore \eq{ai9eq15} commutes up to homotopy, since straightenings are unique up to homotopy. This completes the proof of Theorem~\ref{ai9thm5}.

\subsection{Proof of Theorem \ref{ai9thm7}}
\label{ai116}

We will prove that $\hFun_{A_\iy}^{\un,\iy}$ in \eq{ai9eq19} is an $\iy$-functor. First we will show that $\hFun_{A_\iy}^{\un,\iy}(S,T)$ in \eq{ai9eq20} is an $N$-simplex in $\hACat^{\un,\iy}_R$ in Definition \ref{ai7def3}. To do this we must show that the maps $\o_p^{I;\fH^0,\ldots,\fH^p}$ in \eq{ai9eq21} satisfy \eq{ai7eq15}--\eq{ai7eq16}. We write $\bar\d{}_{p,q;\fI^0,\ldots,\fI^q}^{I;\fH^0,\ldots,\fH^p}$ for the maps constructed in Definition \ref{ai7def2} from the data in $\hFun_{A_\iy}^{\un,\iy}(S,T),$ so that $\bar\d{}_{p,q;\fI^0,\ldots,\fI^q}^{I;\fH^0,\ldots,\fH^p}$ is a sum of multiple compositions of maps $\o_r^{J;\fK^0,\ldots,\fK^r}$ in \eq{ai9eq21}. The next proposition can be deduced from the (complicated) definitions, and $q-1$ nested applications of Theorem~\ref{ai7thm1}.

\begin{prop}
\label{ai11prop2}
Suppose $I=\{i_0<\cdots<i_M\},$ $p,q\ge 0,$ $\fH^0,\ldots,\fH^p\in\sQ^{i_0},$ $\fJ^0,\ldots,\fJ^q\in\sQ^{i_M},$ $\fT^j\in\sJ^{i_0,d_j}(\fH^{j-1},\fH^j)$ for $j=1,\ldots,p,$ and\/ $w_j^i\in\sO^{N-i_M}$ for\/ $i=1,\ldots,q$ and\/ $j=0,\ldots,k_i$ with\/ $k_1,\ldots,k_q\ge 0$. Then if\/ $M=0$ then $\bar\d{}_{p,q;\fI^0,\ldots,\fI^q}^{\{i_0\};\fH^0,\ldots,\fH^p}=0$ if\/ $p<q,$ and if\/ $p\ge q$ then
\ea
&\bigl(\bar\d{}_{p,q;\fI^0,\ldots,\fI^q}^{\{i_0\};\fH^0,\ldots,\fH^p}(\fT^1\ot\cdots\ot\fT^p)\bigr)^{w_0^1,\ldots,w_{k_1}^1;\ldots;w_0^q,\ldots,w_{k_q}^q}_{k_1,\ldots,k_q}
\label{ai11eq84}\\
&=\!\!\!\!\sum_{\begin{subarray}{l} 
1\le r\le q,\\
p-q\le m\le k_r+p-q,\\
y_0,\ldots,y_m\in\sP^{i_0}, \\
\fJ^i=\fH^i,\; 0\le i<r,\\
\fJ^i=\fH^{i+p-q},\\ r\le i\le q
\end{subarray}}
\begin{aligned}[t]&
(-1)^\dag(\fT^1)^{w_0^1,\ldots,w_{k_1}^1}_{k_1}\!\ot\!\cdots\!\ot\!(\fT^{r-1})^{w_0^{r-1},\ldots,w_{k_{r-1}}^{r-1}}_{k_{r-1}}\!\ot\!\bigl(\n_m^{\sB,\{i_0\};y_0,\ldots,y_m} \\
&{}\ci\bD_{k_r,m;y_0,\ldots,y_m}^{p-q+1;w_0^r,\ldots,w_{k_r}^r}(\fH^{r-1},\ldots,\fH^{r+p-q},\fT^r\!\ot\!\cdots\!\ot\!\fT^{r+p-q})\bigr)\\
&{}\ot(\fT^{r+p-q+1})^{w_0^{r+1},\ldots,w_{k_{r+1}}^{r+1}}_{k_{r+1}}\ot\cdots\ot(\fT^p)^{w_0^q,\ldots,w_{k_q}^q}_{k_q}
\end{aligned}
\nonumber\\
&+\!\!\!\sum_{\begin{subarray}{l} 
1\le r\le q,\\
0\le l\le k_r,\\
x_0,\ldots,x_l\in\sO^{N-i_0}, \\
\fJ^i=\fH^i,\; 0\le i<r,\\
\fJ^i=\fH^{i+p-q},\\ r\le i\le q \end{subarray}}
\begin{aligned}[t]&
(-1)^{d_1+\cdots+d_p}(-1)^\ddag
(\fT^1)^{w_0^1,\ldots,w_{k_1}^1}_{k_1}\ot\cdots\ot(\fT^{r-1})^{w_0^{r-1},\ldots,w_{k_{r-1}}^{r-1}}_{k_{r-1}}\ot{}\\
&\bigl(\bD_{l,1;\fH^{r-1}(x_0),\fH^{r+p-q}(x_l)}^{p-q;x_0,\ldots,x_l}
(\fH^{r-1},\ldots,\fH^{r+p-q},\fT^r\!\ot\!\cdots\!\ot\!\fT^{r+p-q})\\
&\ci\d^{\sA,\{N-i_0\};w_0^r,\ldots,w_{k_r}^r}_{k_r,l;x_0,\ldots,x_l}
\bigr)\\
&{}\ot(\fT^{r+p-q+1})^{w_0^{r+1},\ldots,w_{k_{r+1}}^{r+1}}_{k_{r+1}}\ot\cdots\ot(\fT^p)^{w_0^q,\ldots,w_{k_q}^q}_{k_q},
\end{aligned}
\nonumber
\ea
and if\/ $M>0$
\ea
&\bigl(\bar\d{}_{p,q;\fI^0,\ldots,\fI^q}^{I;\fH^0,\ldots,\fH^p}(\fT^1\ot\cdots\ot\fT^p)\bigr)^{w_0^1,\ldots,w_{k_1}^1;\ldots;w_0^q,\ldots,w_{k_q}^q}_{k_1,\ldots,k_q}
\label{ai11eq85}\\
&=\sum_{\begin{subarray}{l} 
p=p_1+\cdots+p_q,\; p_i\ge 0, \\
r_i=\sum_{j=0}^ip_j,\; 0\le i\le q, \\
0\le l_i\le k_i,\\ p_i\le m_i\le l_i+p_i,\; 1\le i\le q,\\
x^i_j\in\sO^{N-i_0}, \\ 1\le i\le q,\; 0\le j\le l_i,\\
y^i_j\in\sP^{i_0}, \\ 1\le i\le q,\; 0\le j\le m_i, \\
I=A_0\supseteq A_1\supseteq\cdots \\
\supseteq A_q=\{i_0,i_M\},\\
\md{A_i}=1+n_i, \\
A_i=\{a_0^i<a_1^i<\cdots<a_{n_i}^i\},\\
a^i_j=a^{i-1}_{b^i_j},\; 1\le i\le q,\\
0\le j\le n_i,\\
\fJ^i= \fG^{A_i}\ci\fH^{r_i}\\
\ci\fF^{(N-i_M)(N-i_0)}, \\
0\le i\le q, \\
\text{{\rm`$*$'} mean matching} \\
\text{choices of objects}
\end{subarray}}
\begin{aligned}[t]&
(-1)^{\star} \Bigl[\bigot_{i=1}^q
\bigl(\n_{*}^{\{a^{i-1}_{b^i_{n_i-1}},a^{i-1}_{b^i_{n_i-1}+1},\ldots,a^{i-1}_{b^i_{n_i}}\};*}\\
&{}\ci\d_{*}^{\sB,\{a^{i-1}_{b^i_{n_i-2}},a^{i-1}_{b^i_{n_i-2}+1},\ldots,a^{i-1}_{b^i_{n_i-1}}\};*}\\
&{}\ci\cdots\ci\d_{m_i,*}^{\sB,\{a^{i-1}_{b^i_0},a^{i-1}_{b^i_0+1},\ldots,a^{i-1}_{b^i_1}\};y^i_0,\ldots,y^i_{m_i}}\bigr)
\Bigr]\ci{} \\
&\bigl(\bD_{l_1,m_1;y^1_0,\ldots,y^1_{m_1}}^{p_1;x^1_0,\ldots,x^1_{l_1}}(\fH^{r_0},\ldots,\fH^{r_1},\fT^{r_0+1}\ot\cdots\ot\fT^{r_1})\ot\cdots \\
&{}\ot\bD_{l_q,m_q;y^q_0,\ldots,y^q_{m_q}}^{p_q;x^q_0,\ldots,x^q_{l_q}}(\fH^{r_{q-1}},\ldots,\fH^{r_q},\fT^{r_{q-1}+1}\ot\cdots\ot\fT^{r_q})\bigr)
\\
&{}\ci\bigl(\d^{\sA,\{N-i_M,N-i_0\};w^1_0,\ldots,w^1_{k_1}}_{k_1,l_1;x^1_0,\ldots,x^1_{l_1}}\ot\cdots \\
&{}\ot\d^{\sA,\{N-i_M,N-i_0\};w^q_0,\ldots,w^q_{k_q}}_{k_q,l_q;x^q_0,\ldots,x^q_{l_q}}\bigr).
\end{aligned}
\nonumber
\ea

Here to reduce notation we have omitted the explicit choices of objects in the operators $\d_{\ldots}^{\sB,\ldots},$ which are like the choices of $k_i\ge 0$ and\/ $w_0^i,\ldots,w_{k_i}^i\in\sO^i$ in {\rm\eq{ai7eq13},} and have just written {\rm`$*$'} instead. By this we mean that we should sum over all possible choices of strings of objects such that the operators $\d_{\ldots}^{\sB,\ldots}$ are well defined and compose consistently. The unspecified signs $(-1)^\dag,(-1)^\ddag,(-1)^{\star}$ in\/ \eq{ai11eq85} are the products of signs in the $q-1$ applications of\/~{\rm\eq{ai7eq12}--\eq{ai7eq13}}.
\end{prop}

\begin{rem}
\label{ai11rem2}
The details of \eq{ai11eq84}--\eq{ai11eq85} -- for example, the precise signs $(-1)^\dag,(-1)^\ddag,(-1)^{\star}$ -- will not really matter to us, but only their general structure. This is because when we apply them in \eq{ai11eq86} and \eq{ai11eq88}, we will immediately invert the $q-1$ nested applications of \eq{ai7eq12}--\eq{ai7eq13}, which will simplify the equations and cancel the signs~$(-1)^\dag,(-1)^\ddag,(-1)^{\star}$.
\end{rem}

We will now show that the $\o_p^{I;\fH^0,\ldots,\fH^p}$ in \eq{ai9eq21} satisfy \eq{ai7eq15}--\eq{ai7eq16}. Suppose $I=\{i_0,\ldots,i_M\}\subseteq\{0,1,\ldots,N\}$ with $i_0<i_1<
\cdots<i_M$, $p\ge 0$, and $\fH^0,\ldots,\fH^p\in\sQ^{i_0}$, and $\fT^j\in\sJ^{i_0,d_j}(\fH^{j-1},\fH^j)$ for $j=1,\ldots,p$, and $k\ge 0$ and $w_0,\ldots,w_k\in\sO^{N-i_M}$. Then if $M=0$ we have
\ea
&\biggl[\sum_{\begin{subarray}{l} q=1,\ldots,p, \\
\fI^0,\ldots,\fI^q\in\sQ^{i_0} \end{subarray}}
\o_q^{\{i_0\};\fI^0,\ldots,\fI^q}\ci
\bar\d{}_{p,q;\fI^0,\ldots,\fI^q}^{\{i_0\};\fH^0,\ldots,\fH^p}(\fT^1\ot\cdots\ot\fT^p)\biggr]_k^{w_0,\ldots,w_k}
\label{ai11eq86}\\
&=\sum_{\begin{subarray}{l} 
p\le l\le k+p,\\
1\le m\le l,\\
x_0,\ldots,x_l\in\sP^{i_0}, \\
y_0,\ldots,y_m\in\sP^{i_0}
\end{subarray}}
\begin{aligned}[t]&
\n_m^{\{i_0\};y_0,\ldots,y_m}\ci\bigl\llbracket\d^{\sB,\{i_0\};x_0,\ldots,x_l}_{l,m;y_0,\ldots,y_m}\ci{} \\
&\bD_{k,l;x_0,\ldots,x_l}^{p;w_0,\ldots,w_k}(\fH^0,\ldots,\fH^p,\fT^1\ot\cdots\ot\fT^p)\bigr\rrbracket
\end{aligned}
\allowdisplaybreaks
\nonumber\\
&+(-1)^{d_1+\cdots+d_p}\sum_{\begin{subarray}{l} 
1\le l\le k, \\
p\le m\le l+p,\\
x_0,\ldots,x_l\in\sO^{N-i_0}, \\
y_0,\ldots,y_m\in\sP^{i_0}
\end{subarray}}
\begin{aligned}[t]&
\n_m^{\{i_0\};y_0,\ldots,y_m}\ci\bigl\llbracket\bD_{l,m;y_0,\ldots,y_m}^{p;x_0,\ldots,w_l}(\fH^0,\ldots,\fH^p,\\
&\fT^1\ot\cdots\ot\fT^p)\ci\d^{\sA,\{N-i_0\};w_0,\ldots,w_k}_{k,l;x_0,\ldots,x_l}\bigr\rrbracket
\end{aligned}
\allowdisplaybreaks
\nonumber\\
&+(-1)^{d_1+\cdots+d_p+1}\sum_{\begin{subarray}{l} 
1\le l\le k, \\
p\le m\le l+p,\\
x_0,\ldots,x_l\in\sO^{N-i_0}, \\
y_0,\ldots,y_m\in\sP^{i_0}
\end{subarray}}
\begin{aligned}[t]&
\bigl\llbracket\n_m^{\{i_0\};y_0,\ldots,y_m}\ci\bD_{l,m;y_0,\ldots,y_m}^{p;x_0,\ldots,w_l}(\fH^0,\ldots,\fH^p,\\
&\fT^1\ot\cdots\ot\fT^p)\bigr\rrbracket\ci\d^{\sA,\{N-i_0\};w_0,\ldots,w_k}_{k,l;x_0,\ldots,x_l}
\end{aligned}
\allowdisplaybreaks
\nonumber\\
&-\sum_{\begin{subarray}{l} 
0\le m\le l\le k, \\
x_0,\ldots,x_l\in\sO^{N-i_0}, \\
y_0,\ldots,y_m\in\sO^{N-i_0}
\end{subarray}}
\begin{aligned}[t]&
\bigl\llbracket\bD_{m,1;\fH^0(y_0),\fH^p(y_m)}^{p;y_0,\ldots,y_m}(\fH^0,\ldots,\fH^p,\fT^1\ot\cdots\ot\fT^p)\\
&{}\ci\d^{\sA,\{N-i_0\};x_0,\ldots,x_l}_{l,m;y_0,\ldots,y_m}\bigr\rrbracket\ci\d^{\sA,\{N-i_0\};w_0,\ldots,w_k}_{k,l;x_0,\ldots,x_l}
\end{aligned}
\nonumber\\
&=0.
\nonumber
\ea

Here in the first step we combine the expressions \eq{ai9eq23} for $\o_q^{\{i_0\};\fI^0,\ldots,\fI^q}$ and \eq{ai11eq84} for $\bar\d{}_{p,q;\fI^0,\ldots,\fI^q}^{\{i_0\};\fH^0,\ldots,\fH^p}(\fT^1\ot\cdots\ot\fT^p)$, and rewrite. The brackets $\llbracket\cdots\rrbracket$ indicate terms originating from $\bar\d{}_{p,q;\fI^0,\ldots,\fI^q}^{\{i_0\};\fH^0,\ldots,\fH^p}(\fT^1\ot\cdots\ot\fT^p)$, so for example, the third line of \eq{ai11eq86} combines the second term on the right hand side of \eq{ai9eq23} and the first term on the right hand side of \eq{ai11eq84}, but the fourth line of \eq{ai11eq86} combines the first term on the right hand side of \eq{ai9eq23} and the second term on the right hand side of~\eq{ai11eq84}.

The rough idea in the simplification in the first step of \eq{ai11eq86} is that
\e
\begin{split}
&\ts\sum_{r=1}^q\bD_{\ldots}^{q,\ldots}\bigl(\fH^0,\ldots,\fH^{r-1},\fH^{r+p-q},\ldots,\fH^p;\fT^1\ot\cdots\ot\fT^{r-1}\ot
\\
&\bD^{p-q+1}(\fH^{r-1},\ldots,\fH^{r+p-q},\fT^r\!\ot\!\cdots\!\ot\!\fT^{r+p-q})\ot \fT^{r+p-q+1}\ot\cdots\ot \fT^p\bigr)\\
&\quad \approx\bD_{\ldots}^{p;\ldots}(\fH^0,\ldots,\fH^p,\fT^1\ot\cdots\ot\fT^p).
\end{split}
\label{ai11eq87}
\e
The two sides of \eq{ai11eq87} are not equal on the nose: on the left hand side there is extra information in the choice of $r$ and in where we choose to split the sequences of factors $\fH_{\ldots}^{r-1;\ldots}$ and $\fH_{\ldots}^{r+p-q;\ldots}$ in \eq{ai9eq25}, which are not present on the right hand side. But this extra information is transferred to choices of $j$ in equation \eq{ai7eq8} for the other terms $\d^{\sA,\{N-i_0\};\ldots}_{\ldots}$ and $\d^{\sB,\{i_0\};\ldots}_{\ldots}$ in~\eq{ai11eq86}.

For the signs in the second--fifth lines in \eq{ai11eq86}, the signs $(-1)^\dag$ and $(-1)^\ddag$ in \eq{ai11eq84} are absorbed into signs in the definition \eq{ai7eq8} of the terms $\d^{\sA,\{N-i_0\};\ldots}_{\ldots}$ and $\d^{\sB,\{i_0\};\ldots}_{\ldots}$. In the third line the sign $(-1)^{d_1+\cdots+d_p}$ comes from the corresponding sign in \eq{ai11eq84}, where $\deg\fT^i=d_i$. In the fourth line the sign $(-1)^{d_1+\cdots+d_p+1}$ corresponds to $(-1)^{d_1+\cdots+d_p}$ in \eq{ai9eq23}, but noting that we are applying $\o_q^{I;\fI^0,\ldots,\fI^q}$ not to $\fT^1\ot\cdots\ot\fT^p$ with $\deg(\fT^1\ot\cdots\ot\fT^p)=d_1+\cdots+d_p$, but to $\bar\d{}_{p,q;\fI^0,\ldots,\fI^q}^{\{i_0\};\fH^0,\ldots,\fH^p}(\fT^1\ot\cdots\ot\fT^p)$ with $\deg\bigl(\bar\d{}_{p,q;\fI^0,\ldots,\fI^q}^{\{i_0\};\fH^0,\ldots,\fH^p}(\fT^1\ot\cdots\ot\fT^p)\bigr)=d_1+\cdots+d_p+1$. In the fifth line the minus sign is the product of these signs $(-1)^{d_1+\cdots+d_p}$ and~$(-1)^{d_1+\cdots+d_p+1}$.

For the second step of \eq{ai11eq86}, note that the second line is zero by \eq{ai7eq15} for $T$, the third and fourth lines cancel, and the fifth line is zero by \eq{ai7eq21} for $S$. This proves \eq{ai11eq86}, and shows that when $M=0$ the $\o_p^{I;\fH^0,\ldots,\fH^p}$ satisfy~\eq{ai7eq15}.

Similarly, if $M>0$ we can show that
\ea
&\biggl[\sum_{\begin{subarray}{l} q=1,\ldots,p \\
\fI^0,\ldots,\fI^q\in\sQ^{i_0}  \end{subarray}}
\o_q^{I;\fI^0,\ldots,\fI^q}\ci
\bar\d{}_{p,q;\fI^0,\ldots,\fI^q}^{\{i_0\};\fH^0,\ldots,\fH^p}(\fT^1\ot\cdots\ot\fT^p)
\label{ai11eq88}\\
&+(-1)^M\sum_{\begin{subarray}{l} q=1,\ldots,p+M-1,\\
\fI^0,\ldots,\fI^q\in\sQ^{i_M} \end{subarray}}
\o_q^{\{i_M\};\fI^0,\ldots,\fI^q}\ci
\bar\d{}_{p,q;\fI^0,\ldots,\fI^q}^{I;\fH^0,\ldots,\fH^p}(\fT^1\ot\cdots\ot\fT^p)\biggr]_k^{w_0,\ldots,w_k}
\allowdisplaybreaks
\nonumber\\
&=\sum_{\begin{subarray}{l} 
1\le l\le k, \;
p\le m\le l+p,\;
1\le n\le m,\\
x_0,\ldots,x_l\in\sO^{N-i_0}, \\
y_0,\ldots,y_m\in\sP^{i_0}, \\
z_0,\ldots,z_n\in\sP^{i_0}
\end{subarray}}
\begin{aligned}[t]&
\n_m^{I;z_0,\ldots,z_n}\ci\bigl\llbracket\d^{\sB,\{i_0\};y_0,\ldots,y_m}_{m,n;z_0,\ldots,z_n}\ci{} \\
&\bD_{l,m;y_0,\ldots,y_m}^{p;x_0,\ldots,x_l}(\fH^0,\ldots,\fH^p,\fT^1\ot\cdots\ot\fT^p)\bigr\rrbracket\\
&\ci\d^{\sA,\{N-i_M,N-i_0\};w_0,\ldots,w_k}_{k,l;x_0,\ldots,x_l}
\end{aligned}
\allowdisplaybreaks
\nonumber\\
&+(-1)^{d_1+\cdots+d_p}\!\!\!\!\!\sum_{\begin{subarray}{l} 
1\le m\le l\le k, \;
p\le n\le m+p,\\
x_0,\ldots,x_l\in\sO^{N-i_0}, \\
y_0,\ldots,y_m\in\sO^{N-i_0}, \\
z_0,\ldots,z_n\in\sP^{i_0}
\end{subarray}}
\begin{aligned}[t]&
\n_n^{I;z_0,\ldots,z_n}\ci\bigl\llbracket\bD_{m,n;z_0,\ldots,z_n}^{p;y_0,\ldots,y_m}(\fH^0,\ldots,\fH^p,\\
&\fT^1\ot\cdots\ot\fT^p)\ci\d^{\sA,\{N-i_0\};x_0,\ldots,x_l}_{l,m;y_0,\ldots,y_m}\bigr\rrbracket\\
&\ci\d^{\sA,\{N-i_M,N-i_0\};w_0,\ldots,w_k}_{k,l;x_0,\ldots,x_l}
\end{aligned}
\allowdisplaybreaks
\nonumber\\
&+(-1)^M\sum_{\begin{subarray}{l} 
1\le l\le k, \;
p\le m\le l+p,\; 
1\le n\le m, \\
x_0,\ldots,x_l\in\sO^{N-i_0}, \\
y_0,\ldots,y_m\in\sP^{i_0}, \\
z_0,\ldots,z_n\in\sP^{i_M}
\end{subarray}}
\begin{aligned}[t]&
\n_n^{\{i_M\};z_0,\ldots,z_n}\ci\bigl\llbracket\d^{\sB,I;y_0,\ldots,y_m}_{m,n;z_0,\ldots,z_n}\ci\\
&\bD_{l,m;y_0,\ldots,y_m}^{p;x_0,\ldots,x_l}(\fH^0,\ldots,\fH^p,\fT^1\ot\cdots\ot\fT^p)\\
&\ci\d^{\sA,\{N-i_M,N-i_0\};w_0,\ldots,w_k}_{k,l;x_0,\ldots,x_l}
\bigr\rrbracket
\end{aligned}
\nonumber\\
&+(-1)^M(-1)^{d_1+\cdots+d_p+1-M}\cdot{}
\nonumber\\
&\quad\sum_{\begin{subarray}{l} 
1\le m\le l\le k, \; p\le n\le m+p,\\
x_0,\ldots,x_l\in\sO^{N-i_M}, \\
y_0,\ldots,y_m\in\sO^{N-i_0}, \\
z_0,\ldots,z_n\in\sP^{i_0}
\end{subarray}}
\begin{aligned}[t]&
\bigl\llbracket\n_n^{I;z_0,\ldots,z_n}\ci\bD_{m,n;z_0,\ldots,z_n}^{p;y_0,\ldots,y_m}(\fH^0,\ldots,\fH^p,\\
&\fT^1\ot\cdots\ot\fT^p)
\ci\d^{\sA,\{N-i_M,N-i_0\};x_0,\ldots,x_l}_{l,m;y_0,\ldots,y_m}
\bigr\rrbracket\\
&\ci\d^{\sA,\{N-i_M\};w_0,\ldots,w_k}_{k,l;x_0,\ldots,x_l}
\end{aligned}
\nonumber\\
&=0.
\nonumber
\ea

The first step is proved as for \eq{ai11eq86}, but using \eq{ai9eq23}--\eq{ai9eq24} and \eq{ai11eq84}--\eq{ai11eq85}. In the sixth line, the sign $(-1)^{d_1+\cdots+d_p+1-M}$ comes from $(-1)^{d_1+\cdots+d_p}$ in \eq{ai9eq23}, but noting that $\deg\bigl(\bar\d{}_{p,q;\fI^0,\ldots,\fI^q}^{I;\fH^0,\ldots,\fH^p}(\fT^1\ot\cdots\ot\fT^p)\bigr)=d_1+\cdots+d_p+1-M$. The second step follows as the third and fifth lines sum to zero by \eq{ai7eq16} for $T$, and the fourth and sixth lines sum to zero by \eq{ai7eq21} for $S$ with $\{N-i_M,N-i_0\}$ in place of $I$. This proves \eq{ai11eq88}, and shows that when $M>0$ the $\o_p^{I;\fH^0,\ldots,\fH^p}$ satisfy \eq{ai7eq16}. Hence $\hFun_{A_\iy}^{\un,\iy}(S,T)$ in \eq{ai9eq20} is an $N$-simplex in $\hACat^{\un,\iy}_R$ in Definition~\ref{ai7def3}. 

Next we show $\hFun_{A_\iy}^{\un,\iy}(S,T)$ is an $N$-simplex in $\ACat^{\un,\iy}_R$. If $0\le i_0<i_1<i_2\le N$ the corresponding 2-simplex in $\hFun_{A_\iy}^{\un,\iy}(S,T)$ may be written
\begin{equation*}
\xymatrix@C=50pt@R=15pt{
\drrtwocell_{}\omit^{}\omit{_{\,\,\,\,\,\,\,\,\,\,\,\,\,\,\,\,\,\,\,\,\,\,\,\,\,\,\,\,\,\,\,\,\,\,\,\,\,\, (\fH^{i_0i_1i_2})_*}} & \Fun_{A_\iy}^\un(\sA^{N-i_1},\sB^{i_1})  \ar@/^1.4pc/[dr]^(0.8){\!\!\!\!\!\!\!\!\!\!\!\!\!\!\!\!\!\!\!\!\!\!\!\!\!\!\!\!\!\!\!\!(\fF^{(N-i_2)(N-i_1)})^*\ci(\fG^{i_1i_2})_*} \\
*+[r]{\Fun_{A_\iy}^\un(\sA^{N-i_0},\sB^{i_0})} \ar@/^1.4pc/[ur]^(0.2){(\fF^{(N-i_1)(N-i_0)})^*\ci(\fG^{i_0i_1})_*\!\!\!\!\!\!\!\!\!\!\!\!\!\!\!\!\!\!\!\!\!\!\!\!\!\!\!\!\!\!\!\!} \ar[rr]_{(\fF^{(N-i_2)(N-i_0)})^*\ci(\fG^{i_0i_2})_*} && *+[l]{\Fun_{A_\iy}^\un(\sA^{N-i_2},\sB^{i_2}).}
}
\end{equation*}
Here $\fH^{i_0i_1i_2}:\fG^{i_2i_2}\ci\fG^{i_0i_1}\Ra\fG^{i_0i_2}$ is the unital $A_\iy$-natural transformation in the $i_0i_1i_2$ 2-simplex in $T$, as in Example \ref{ai7ex2}(c) and Definition \ref{ai7def4}(iii), and $(\fH^{i_0i_1i_2})_*$ is the induced unital $A_\iy$-natural transformation in $\hFun_{A_\iy}^{\un,\iy}(S,T)$. The important point is that as $\fH^{i_0i_1i_2}$ is an $A_\iy$-natural isomorphism, so is $(\fH^{i_0i_1i_2})_*$, and thus $\hFun_{A_\iy}^{\un,\iy}(S,T)$ lies in $\ACat^{\un,\iy}_R$ by Definition~\ref{ai7def4}.

Our goal is to prove $\hFun_{A_\iy}^{\un,\iy}$ in \eq{ai9eq19} is an $\iy$-functor, a morphism of simplicial sets. That is, we must show it is a natural transformation of functors $\De^\rop\ra\Set$. So far we know that for each $[N]\in\De^\rop$, $\hFun_{A_\iy}^{\un,\iy}([N])$ does map $(\sN(\ACat^\un_R)^\rop\t\ACat^{\un,\iy}_R)([N])$ to $\ACat^{\un,\iy}_R([N])$. It remains to prove these $\hFun_{A_\iy}^{\un,\iy}([N])$ commute with the action of morphisms $\al:[N]\ra[\bar N]$ in $\De^\rop$ (equivalently, $\al:[\bar N]\ra[N]$ in $\De$), where $\al$ acts as in Definition~\ref{ai7def3}(b).

It is clear that the $\al$-actions on the data $(\sO^i),(f^{ij}),(\sH^i)$, $(\sP^i),(g^{ij}),(\sI^i)$ and $(\sP^i),(h^{ij}),(\sJ^i)$ commute. We will show the $\al$-actions on $(\m^I),(\n^I),(\o^I)$ commute. The data $\o^I$ is defined in \eq{ai9eq22}--\eq{ai9eq25}, and the action of $\al:[\bar N]\ra[N]$ in $\De$ on this is defined in four cases in \eq{ai7eq20}. The first case, that $\al\vert_{\bar I}$ is injective, holds because $\al$ is the identity on the data attached to $\bar I$ in $S,T$, so it is also the identity on the data attached to $\bar I$ in~$\hFun_{A_\iy}^{\un,\iy}(S,T)$.

For the second case, with $p=1$ and $\bar I=\{\bar\imath_0,\bar\imath_1\}$ with $\bar\imath_0\ne \bar\imath_1$, $\al(\bar\imath_0)=\al(\bar\imath_1)=i$, in the definition of $\o_1^{\bar I}$ in \eq{ai9eq24}, from the action of $\al$ on $S$ we see that
\begin{equation*}
\d^{\sA,\{\bar N-\bar\imath_1,\bar N-\bar\imath_0\};w_0,\ldots,w_k}_{k,l;x_0,\ldots,x_l}=\begin{cases} \id, & \text{$k=l$ and $(w_0,\ldots,w_k)=(x_0,\ldots,x_l)$,} \\
0, & \text{otherwise.}\end{cases}
\end{equation*}
Also from the action of $\al$ on $T$ we have
\begin{equation*}
\n_m^{\{\bar\imath_0,\bar\imath_1\};y_0,\ldots,y_m}=\begin{cases}\id, & m=1, \\ 0, & \text{otherwise.} \end{cases}
\end{equation*}
Hence \eq{ai9eq24} reduces to
\begin{equation*}
\bigl(\o_1^{\{\bar\imath_0,\bar\imath_1\}}(\fT^1)\bigr)_k^{w_0,\ldots,w_k}=
\sum_{\begin{subarray}{l} 
y_0,y_1\in\sP^{i_0}
\end{subarray}}
\bD_{k,1;y_0,y_1}^{1;w_0,\ldots,w_k}(\fH^0,\fH^1,\fT^1)
=(\fT^1)_k^{w_0,\ldots,w_k},
\end{equation*}
using \eq{ai9eq25}. So $\o_1^{\{\bar\imath_0,\bar\imath_1\}}=\id$, as we want in \eq{ai7eq20}.

For the third case, with $p=0$, $\bar I=\{\bar\imath_0<\bar\imath_1<\bar\imath_2\}$ and $\al\vert_{\bar I}$ not injective, in \eq{ai9eq24} for $\o_0^{\bar I}$, from the action of $\al$ on $T$ we have
\begin{equation*}
\n_m^{\bar I;y_0,\ldots,y_m}=\begin{cases} \boo^{\al(\bar\imath_2)}_{g^{\al(\bar\imath_0)\al(\bar\imath_2)}(y_0)}, & m=0, \\ 0, & m>0.
\end{cases}
\end{equation*}
When $m=0$ we see that
\begin{equation*}
\bD_{l,0;y_0}^{0;x_0,\ldots,x_l}(\fH^0)=\begin{cases}
\id_R, & \text{$l=0$ and $y_0=\fH^0(x_0)$,} \\ 0, & \text{otherwise.}
\end{cases}
\end{equation*}
And when $l=0$ we have 
\begin{equation*}
\d^{\sA,\{\bar N-\bar\imath_2,\bar N-\bar\imath_0\};w_0,\ldots,w_k}_{k,0;x_0}=
\begin{cases}
\id_R, & \text{$k=0$ and $x_0=f^{(N-\al(\bar\imath_2))(N-\al(\bar\imath_0))}(w_0)$,} \\ 0, & \text{otherwise.}
\end{cases}
\end{equation*}
Combining these shows that
\begin{equation*}
(\o_0^{\bar I})_k^{w_0,\ldots,w_k}=\begin{cases}
\boo^{\al(\bar\imath_2)}_{g^{\al(\bar\imath_0)\al(\bar\imath_2)}\ci\fH^0\ci f^{(N-\al(\bar\imath_2))(N-\al(\bar\imath_0))}(w_0)}, & k=0, \\
0, & \text{otherwise,}
\end{cases}
\end{equation*}
which is what we want.

For the fourth case, note that we are considering $\o_p^{\bar I}$ with $\bar I,p$ satisfying the case conditions in \eq{ai7eq20}, but \eq{ai9eq24} involves $\n_m^{\bar I;y_0,\ldots,y_m}$ where we need not have $m=p$. So there are two cases in which we could have $\n_m^{\bar I;y_0,\ldots,y_m}\ne 0$ in~\eq{ai9eq24}:
\begin{itemize}
\setlength{\itemsep}{0pt}
\setlength{\parsep}{0pt}
\item[(i)] $\bar I=\{\bar\imath_0,\bar\imath_1\}$ with $\bar\imath_0\ne \bar\imath_1$, $\al(\bar\imath_0)=\al(\bar\imath_1)=i$, and $p\ne 1$ but $m=1$.
\item[(ii)] $\bar I=\{\bar\imath_0<\bar\imath_1<\bar\imath_2\}$ and $\al\vert_{\bar I}$ not injective, and $m=0$ but $p>0$.
\end{itemize}
Neither can give a nonzero term in \eq{ai9eq24}, as $\bD_{l,m;y_0,\ldots,y_m}^{p;x_0,\ldots,x_l}(\fH^0,\ldots,\fH^p,\fT^1\ot\cdots\ot\fT^p)=0$ if either $m=1$ and $p\ne 1$, or if $m=0$ and $p>0$. So \eq{ai9eq24} gives $\o_p^{\bar I}=0$. This completes the proof that $\hFun_{A_\iy}^{\un,\iy}$ is an $\iy$-functor.

Next we show that the top rectangle of \eq{ai9eq28} strictly commutes. Recall that the $\iy$-functor $\iinc^\un_R:\sN(\ACat^\un_R)\hookra\ACat^{\un,\iy}_R$ in Definition \ref{ai7def11} identifies $\sN(\ACat^\un_R)$ with the simplicial subset of $\ACat^{\un,\iy}_R$ whose $N$-simplices have $\m_0^{\{i_0,i_1,i_2\};x_0}=\boo^{i_2}_{f^{i_0i_2}(x_0)}$ if $\fF^{i_1i_2}\ci\fF^{i_0i_1}=\fF^{i_0i_2}$, and $\m_k^{I;x_0,\ldots,x_k}=0$ in all other cases with $\md{I}\ge 3$.

We can show from \eq{ai9eq24} that if $T$ (as well as $S$) lies in $\sN(\ACat^\un_R)\subset\ACat^{\un,\iy}_R$ then $\hFun_{A_\iy}^{\un,\iy}(S,T)$ lies in $\sN(\ACat^\un_R)$. Hence $\hFun_{A_\iy}^{\un,\iy}$ restricts to an $\iy$-functor $\sN(\ACat^\un_R)^\rop\t\sN(\ACat^\un_R)\ra\sN(\ACat^\un_R)$, which must be the nerve of a functor $(\ACat^\un_R)^\rop\t\ACat^\un_R\ra\ACat^\un_R$. Considering the action of $\hFun_{A_\iy}^{\un,\iy}$ on 0- and 1-simplices, we see that this functor is $\Fun_{A_\iy}^\un$. Thus the top rectangle of \eq{ai9eq28} commutes. The bottom rectangle then commutes up to homotopy by comparing \eq{ai9eq18} and~\eq{ai9eq28}.

This completes the parts of Theorem \ref{ai9thm7} concerning $\hFun_{A_\iy}^{\un,\iy}$ in \eq{ai9eq19}. It remains to prove the analogues for $\cFun_{A_\iy}^{\un,\iy}$ in \eq{ai9eq26}. The proofs are very similar, so we will only outline the differences.

In the analogue of Proposition \ref{ai11prop2} for $\cFun_{A_\iy}^{\un,\iy}$, equation \eq{ai11eq84} is unchanged, but \eq{ai11eq85} should be replaced by
\ea
&\bigl(\bar\d{}_{p,q;\fI^0,\ldots,\fI^q}^{I;\fH^0,\ldots,\fH^p}(\fT^1\ot\cdots\ot\fT^p)\bigr)^{w_0^1,\ldots,w_{k_1}^1;\ldots;w_0^q,\ldots,w_{k_q}^q}_{k_1,\ldots,k_q}
=(-1)^{(M+1)(d_1+\cdots+d_p)}\cdot{}
\nonumber\\
&\sum_{\begin{subarray}{l} 
p=p_1+\cdots+p_q,\; p_i\ge 0, \\
r_i=\sum_{j=0}^ip_j,\; 0\le i\le q, \\
0\le l_i\le k_i,\\ p_i\le m_i\le l_i+p_i,\; 1\le i\le q,\\
x^i_j\in\sO^{N-i_0}, \\ 1\le i\le q,\; 0\le j\le l_i,\\
y^i_j\in\sP^{i_0}, \\ 1\le i\le q,\; 0\le j\le m_i, \\
I=A_0\supseteq A_1\supseteq\cdots \\
\supseteq A_q=\{i_0,i_M\},\\
\md{A_i}=1+n_i, \\
A_i=\{a_0^i<a_1^i<\cdots<a_{n_i}^i\},\\
a^i_j=a^{i-1}_{b^i_j},\; 1\le i\le q,\\
0\le j\le n_i,\\
\fJ^i= \fG^{i_0i_M}\ci\fH^{r_i}\\
\ci\fF^{N-A_i}, \;
0\le i\le q, \\
\text{{\rm`$*$'} mean matching} \\
\text{choices of objects}
\end{subarray}}
\begin{aligned}[t]&
(-1)^{\star}\bigl(\n^{\{i_0,i_M\};y^1_0,\ldots,y^1_{m_1}}_{m_1}\ot\cdots\ot\n^{\{i_0,i_M\};y^q_0,\ldots,y^q_{m_q}}_{m_q}\bigr)\\
&{}\ci\bigl(\bD_{l_1,m_1;y^1_0,\ldots,y^1_{m_1}}^{p_1;x^1_0,\ldots,x^1_{l_1}}(\fH^{r_0},\ldots,\fH^{r_1},\fT^{r_0+1}\ot\cdots\ot\fT^{r_1})\ot\cdots \\
&{}\ot\bD_{l_q,m_q;y^q_0,\ldots,y^q_{m_q}}^{p_q;x^q_0,\ldots,x^q_{l_q}}(\fH^{r_{q-1}},\ldots,\fH^{r_q},\fT^{r_{q-1}+1}\ot\cdots\ot\fT^{r_q})\bigr)
\\
&{}\ci\Bigl[\bigot_{i=1}^q
\bigl(\d_{*,l_i;x^i_0,\ldots,x^i_{l_i}}^{\sA,N-\{a^{i-1}_{b^i_0},a^{i-1}_{b^i_0+1},\ldots,a^{i-1}_{b^i_1}\};*}\\
&{}\ci\d_{*}^{\sA,N-\{a^{i-1}_{b^i_1},a^{i-1}_{b^i_1+1},\ldots,a^{i-1}_{b^i_2}\};*}\\
&{}\ci\cdots\ci\d_{k_i,*}^{\sA,N-\{a^{i-1}_{b^i_{n_i-1}},a^{i-1}_{b^i_{n_i-1}+1},\ldots,a^{i-1}_{b^i_{n_i}}\};w^i_0,\ldots,w^i_{k_i}}\bigr)
\Bigr],
\end{aligned}
\nonumber\\[-19pt]
\label{ai11eq89}
\ea
where if $J\subseteq\{0,\ldots,N\}$ then $N-J=\{N-j:j\in J\}$, and the sign $(-1)^{(M+1)(d_1+\cdots+d_p)}$ comes from~\eq{ai9eq27}.

To prove that the $\o_p^{I;\fH^0,\ldots,\fH^p}$ for $\cFun_{A_\iy}^{\un,\iy}$ satisfy \eq{ai7eq15}--\eq{ai7eq16}, equation \eq{ai11eq86} needs no change, but \eq{ai11eq88} should be replaced by
\begin{align*}
&\biggl[\sum_{\begin{subarray}{l} q=1,\ldots,p \\
\fI^0,\ldots,\fI^q\in\sQ^{i_0}  \end{subarray}}
\o_q^{I;\fI^0,\ldots,\fI^q}\ci
\bar\d{}_{p,q;\fI^0,\ldots,\fI^q}^{\{i_0\};\fH^0,\ldots,\fH^p}(\fT^1\ot\cdots\ot\fT^p)
\\
&+(-1)^M\sum_{\begin{subarray}{l} q=1,\ldots,p+M-1,\\
\fI^0,\ldots,\fI^q\in\sQ^{i_M} \end{subarray}}
\o_q^{\{i_M\};\fI^0,\ldots,\fI^q}\ci
\bar\d{}_{p,q;\fI^0,\ldots,\fI^q}^{I;\fH^0,\ldots,\fH^p}(\fT^1\ot\cdots\ot\fT^p)\biggr]_k^{w_0,\ldots,w_k}
\allowdisplaybreaks\\
&=(-1)^{(M+1)(d_1+\cdots+d_p+1)}\cdot{}\\
&\quad\sum_{\begin{subarray}{l} 
1\le l\le k, \;
p\le m\le l+p,\;
1\le n\le m,\\
x_0,\ldots,x_l\in\sO^{N-i_0}, \\
y_0,\ldots,y_m\in\sP^{i_0}, \\
z_0,\ldots,z_n\in\sP^{i_0}
\end{subarray}}
\begin{aligned}[t]&
\n_m^{\{i_0,i_M\};z_0,\ldots,z_n}\ci\bigl\llbracket\d^{\sB,\{i_0\};y_0,\ldots,y_m}_{m,n;z_0,\ldots,z_n}\ci{} \\
&\bD_{l,m;y_0,\ldots,y_m}^{p;x_0,\ldots,x_l}(\fH^0,\ldots,\fH^p,\fT^1\ot\cdots\ot\fT^p)\bigr\rrbracket\\
&\ci\d^{\sA,N-I;w_0,\ldots,w_k}_{k,l;x_0,\ldots,x_l}
\end{aligned}
\allowdisplaybreaks\\
&+(-1)^{(M+1)(d_1+\cdots+d_p+1)}\cdot(-1)^{d_1+\cdots+d_p}\cdot{}\\
&\quad\sum_{\begin{subarray}{l} 
1\le m\le l\le k, \;
p\le n\le m+p,\\
x_0,\ldots,x_l\in\sO^{N-i_0}, \\
y_0,\ldots,y_m\in\sO^{N-i_0}, \\
z_0,\ldots,z_n\in\sP^{i_0}
\end{subarray}}
\begin{aligned}[t]&
\n_n^{\{i_0,i_M\};z_0,\ldots,z_n}\ci\bigl\llbracket\bD_{m,n;z_0,\ldots,z_n}^{p;y_0,\ldots,y_m}(\fH^0,\ldots,\fH^p,\\
&\fT^1\ot\cdots\ot\fT^p)\ci\d^{\sA,\{N-i_0\};x_0,\ldots,x_l}_{l,m;y_0,\ldots,y_m}\bigr\rrbracket\\
&\ci\d^{\sA,N-I;w_0,\ldots,w_k}_{k,l;x_0,\ldots,x_l}
\end{aligned}
\allowdisplaybreaks\\
&+(-1)^M(-1)^{(M+1)(d_1+\cdots+d_p)}\cdot{}\\
&\quad\sum_{\begin{subarray}{l} 
1\le l\le k, \;
p\le m\le l+p,\; 
1\le n\le m, \\
x_0,\ldots,x_l\in\sO^{N-i_0}, \\
y_0,\ldots,y_m\in\sP^{i_0}, \\
z_0,\ldots,z_n\in\sP^{i_M}
\end{subarray}}
\begin{aligned}[t]&
\n_n^{\{i_M\};z_0,\ldots,z_n}\ci\bigl\llbracket\d^{\sB,\{i_0,i_M\};y_0,\ldots,y_m}_{m,n;z_0,\ldots,z_n}\ci\\
&\bD_{l,m;y_0,\ldots,y_m}^{p;x_0,\ldots,x_l}(\fH^0,\ldots,\fH^p,\fT^1\ot\cdots\ot\fT^p)\\
&\ci\d^{\sA,N-I;w_0,\ldots,w_k}_{k,l;x_0,\ldots,x_l}
\bigr\rrbracket
\end{aligned}
\allowdisplaybreaks\\
&+(-1)^M(-1)^{(M+1)(d_1+\cdots+d_p)}(-1)^{d_1+\cdots+d_p+1-M}\cdot{}\\
&\quad\sum_{\begin{subarray}{l} 
1\le m\le l\le k, \; p\le n\le m+p,\\
x_0,\ldots,x_l\in\sO^{N-i_M}, \\
y_0,\ldots,y_m\in\sO^{N-i_0}, \\
z_0,\ldots,z_n\in\sP^{i_0}
\end{subarray}}
\begin{aligned}[t]&
\bigl\llbracket\n_n^{\{i_0,i_M\};z_0,\ldots,z_n}\ci\bD_{m,n;z_0,\ldots,z_n}^{p;y_0,\ldots,y_m}(\fH^0,\ldots,\fH^p,\\
&\fT^1\ot\cdots\ot\fT^p)
\ci\d^{\sA,N-I;x_0,\ldots,x_l}_{l,m;y_0,\ldots,y_m}
\bigr\rrbracket\\
&\ci\d^{\sA,\{N-i_M\};w_0,\ldots,w_k}_{k,l;x_0,\ldots,x_l}
\end{aligned}\\
&=0.
\end{align*}
Here the additional signs $(-1)^{(M+1)(d_1+\cdots+d_p[+1])}$ compared to \eq{ai11eq88} come from those in \eq{ai9eq27} and \eq{ai11eq89}. The rest of the proof is essentially the same.

\medskip

\noindent{\small\textsc{The Mathematical Institute, Radcliffe
Observatory Quarter, Woodstock Road, Oxford, OX2 6GG, U.K.}

\noindent\textsc{E-mail:} {\tt dominic.joyce@maths.ox.ac.uk.}}

\end{document}